\documentclass[a4paper,12pt,onecolumn,openany]{memoir}

\pdfoutput=1

\usepackage[utf8]{inputenc} 
\usepackage[T1]{fontenc}    %
\usepackage[english]{babel} 
\usepackage[final]{microtype} 
\usepackage{imakeidx} 
\makeindex[intoc]
\usepackage{tikz}
\usetikzlibrary{external} 
\usepackage{pgf}
\pgfmathsetseed{\number\pdfrandomseed}  

\usepackage{amsmath,amssymb,amsthm,mathtools} 

\usepackage{color}
\usepackage{esint}
\usepackage{graphicx}
\usepackage{MnSymbol}
\usepackage[colorlinks=true, pdfstartview=FitV, linkcolor=blue, citecolor=blue, urlcolor=blue,pagebackref=false]{hyperref}
\usepackage{microtype}

\usepackage{bm}
\usepackage{scalerel} 
\usepackage{dsfont}

\renewcommand{\hat}{\widehat}

\setlrmarginsandblock{0.15\paperwidth}{*}{1} 
\setulmarginsandblock{0.2\paperwidth}{*}{1}  
\checkandfixthelayout

\maxsecnumdepth{subsection} 
\setsecnumdepth{subsection}

\makeatletter %
\makechapterstyle{standard}{
  \setlength{\beforechapskip}{0\baselineskip}
  \setlength{\midchapskip}{1\baselineskip}
  \setlength{\afterchapskip}{8\baselineskip}
  
  \renewcommand{\chapnamefont}{\centering\normalfont\Large}
  \renewcommand{\printchaptername}{\chapnamefont \@chapapp}

}
\makeatother

\chapterstyle{standard}

\setsecheadstyle{\normalfont\large\bfseries}
\setsubsecheadstyle{\normalfont\normalsize\bfseries}
\setparaheadstyle{\normalfont\normalsize\bfseries}
\setparaindent{0pt}\setafterparaskip{0pt}

\makeatletter                  
\renewcommand\fps@figure{htbp} %
\renewcommand\fps@table{htbp}  %
\makeatother                   %

\captiondelim{\space } 
\captionnamefont{\small\bfseries} 
\captiontitlefont{\small\normalfont} 

\changecaptionwidth          
\captionwidth{1\textwidth} %

\makepagestyle{standard} 

\makeatletter                 
\makeevenfoot{standard}{}{}{} %
\makeoddfoot{standard}{}{}{}  %
\makeevenhead{standard}{\bfseries\thepage\normalfont\qquad\small\leftmark}{}{}
\makeoddhead{standard}{}{}{\small\rightmark\qquad\bfseries\thepage}
\makeatother                  %

\makeatletter
\makepsmarks{standard}{
\createmark{chapter}{both}{shownumber}{\@chapapp\ }{ \quad }
\createmark{section}{right}{shownumber}{}{ \quad }
\createplainmark{toc}{both}{\contentsname}
\createplainmark{lof}{both}{\listfigurename}
\createplainmark{lot}{both}{\listtablename}
\createplainmark{bib}{both}{\bibname}
\createplainmark{index}{both}{\indexname}
\createplainmark{glossary}{both}{\glossaryname}
}
\makeatother                               %

\makepagestyle{chap} 

\makeatletter
\makeevenfoot{chap}{}{\small\bfseries\thepage}{} 
\makeoddfoot{chap}{}{\small\bfseries\thepage}{}  %
\makeevenhead{chap}{}{}{}   %
\makeoddhead{chap}{}{}{}    %
\makeatother

\nouppercaseheads
\aliaspagestyle{chapter}{chap} %

\newtheorem{theorem}{Theorem}
\newtheorem{proposition}[theorem]{Proposition}
\newtheorem{lemma}[theorem]{Lemma}
\newtheorem{corollary}[theorem]{Corollary}

\theoremstyle{definition}
\newtheorem{exercise}{Exercise}
\newtheorem{remark}[theorem]{Remark}
\newtheorem{example}[theorem]{Example}
\newtheorem{definition}[theorem]{Definition}

\numberwithin{equation}{chapter}
\numberwithin{theorem}{chapter}
\numberwithin{exercise}{chapter}

\newcommand{\Z}{\mathbb{Z}}
\newcommand{\N}{\mathbb{N}}
\newcommand{\R}{\mathbb{R}}

\newcommand{\A}{\mathcal{A}}
\newcommand{\Po}{\mathcal{P}}
\newcommand{\Ahom}{\overline{\A}} 
\newcommand{\Phom}{\overline{\Po}}

\newcommand{\E}{\mathbb{E}}
\renewcommand{\P}{\mathbb{P}}
\newcommand{\F}{\mathcal{F}}
\newcommand{\Zd}{\mathbb{Z}^d}
\newcommand{\Rd}{{\mathbb{R}^d}}

\newcommand{\ep}{\varepsilon}
\newcommand{\eps}{\varepsilon}

\renewcommand{\a}{\mathbf{a}}
\renewcommand{\b}{\mathbf{b}}
\newcommand{\Abf}{{\mathbf{A}}}

\newcommand*\xxbar[1]{%
   \hbox{%
     \vbox{%
       \hrule height 0.7pt 
       \kern0.25ex
       \hbox{%
         \kern+0em
         \ensuremath{#1}%
         \kern-0.1em
       }%
     }%
   }%
}
\newcommand{\ahom}{\xxbar{\mathbf{a}}\,}

\newcommand{\Abfh}{{\overbracket[.8pt][-1pt]{{\mathbf{A}}}}}  

\renewcommand{\subset}{\subseteq}

\newcommand{\cu}{{\scaleobj{1.2}{\square}}}
\newcommand{\cut}{{\scaleobj{1.2}{\boxbox}}}

\renewcommand{\fint}{\strokedint}

\newcommand{\Ll}{\left}
\newcommand{\Rr}{\right}

\DeclareMathOperator{\dist}{dist}

\DeclareMathOperator*{\esssup}{ess\,sup}

\DeclareMathOperator*{\osc}{osc}
\DeclareMathOperator{\var}{var}
\DeclareMathOperator{\cov}{cov}
\DeclareMathOperator{\diam}{diam}

\DeclareMathOperator{\supp}{supp}
\DeclareMathOperator{\spn}{span}

\newcommand{\Tr}{\mathsf{Tr}}
\newcommand{\Ext}{\mathsf{Ext}}

\newcommand{\X}{\mathcal{X}}  
\newcommand{\Y}{\mathcal{Y}} 

\let\shortbar=\bar
\renewcommand{\bar}{\overline}
\renewcommand{\tilde}{\widetilde}
\newcommand{\td}{\widetilde}
\newcommand{\indc}{\mathds{1}}
\newcommand{\1}{\mathds{1}}
\newcommand{\Id}{\mathsf{Id}}

\newcommand{\mcl}{\mathcal}
\newcommand{\msf}{\mathsf}
\newcommand{\al}{\alpha}
\newcommand{\be}{\beta}
\newcommand{\CC}{\mathbf{C}}
\newcommand{\cc}{\mathbf{c}}
\newcommand{\nH}{{\underline H}}
\newcommand{\nL}{{\underline L}}

\renewcommand{\O}{\mathcal{O}}
\renewcommand{\S}{\mathcal{S}}

\newcommand{\de}{\delta}
\newcommand{\Jd}{{J_1^{(\delta)}}}
\newcommand{\tJd}{{\widetilde{J}_1^{(\delta)}}}
\newcommand{\per}{\mathsf{per}}

\newcommand{\Add}{\mathsf{Add}}
\newcommand{\Fluc}{\mathsf{Fluc}}

\newcommand{\Loc}{\mathsf{Loc}}

\newcommand{\W}{\mathsf{W}}
\newcommand{\Proj}{\mcl P_{\shortbar \a}}

\newcommand{\Lso}{L^2_{\mathrm{sol},0}}
\newcommand{\Ls}{L^2_{\mathrm{sol}}}
\newcommand{\Lpot}{L^2_{\mathrm{pot}}}
\newcommand{\Lpoot}{L^2_{\mathrm{pot},0}}

\newcommand\irregularcircle[2]{
  +(0:#1+rand*#2)
  \foreach \a in {10,20,...,350}{
    -- +(\a:#1+rand*#2)
  } -- cycle
}

\usepgflibrary{arrows} 
\usetikzlibrary{arrows}
\usetikzlibrary{arrows.meta}

\makeatletter 
\newcommand{\negphantom}{\v@true\h@true\negph@nt} 
\newcommand{\neghphantom}{\v@false\h@true\negph@nt} 
\newcommand{\negph@nt}{\ifmmode\expandafter\mathpalette 
  \expandafter\mathnegph@nt\else\expandafter\makenegph@nt\fi} 
\newcommand{\makenegph@nt}[1]{%
  \setbox\z@\hbox{\color@begingroup#1\color@endgroup}\finnegph@nt} 
\newcommand{\finnegph@nt}{%
  \setbox\tw@\null 
  \ifv@ \ht\tw@\ht\z@\dp\tw@\dp\z@\fi \ifh@\wd\tw@-\wd\z@\fi\box\tw@} 
\newcommand{\mathnegph@nt}[2]{%
  \setbox\z@\hbox{$\m@th #1{#2}$}\finnegph@nt} 
\makeatother 

\newcommand{\Hminus}{\hat{\phantom{H}}\negphantom{H}H^{-1}}
\newcommand{\Hminusul}{\hat{\phantom{H}}\negphantom{H}\underline{H}^{-1}}

\newcommand{\g}{\mathbf{g}}
\newcommand{\f}{\mathbf{f}}
\newcommand{\h}{\mathbf{h}}

\newcommand{\s}{\mathbf{s}}
\newcommand{\m}{\mathbf{m}}

\renewcommand{\leq}{\leqslant}
\renewcommand{\le}{\leqslant}
\renewcommand{\geq}{\geqslant}
\renewcommand{\ge}{\geqslant}

\newcommand{\pa}{\mathrm{par}}

\maxtocdepth{subsection} 
\settocdepth{subsection}

\AtEndDocument{\addtocontents{toc}{\par}} 

\usepackage{hyperref}   
\hypersetup{
pdfborder={0 0 0},      
pdfauthor={I am the Author} 
}
\usepackage{memhfixc}   %

\author{Scott Armstrong \and Tuomo Kuusi \and Jean-Christophe Mourrat}
\date{Preliminary version of \today}

\title{Quantitative stochastic homogenization and large-scale regularity}

\makeindex

\begin{document}

\frontmatter

\maketitle

\bigskip

\begin{figure*}
\centering
\includegraphics[width=\linewidth]{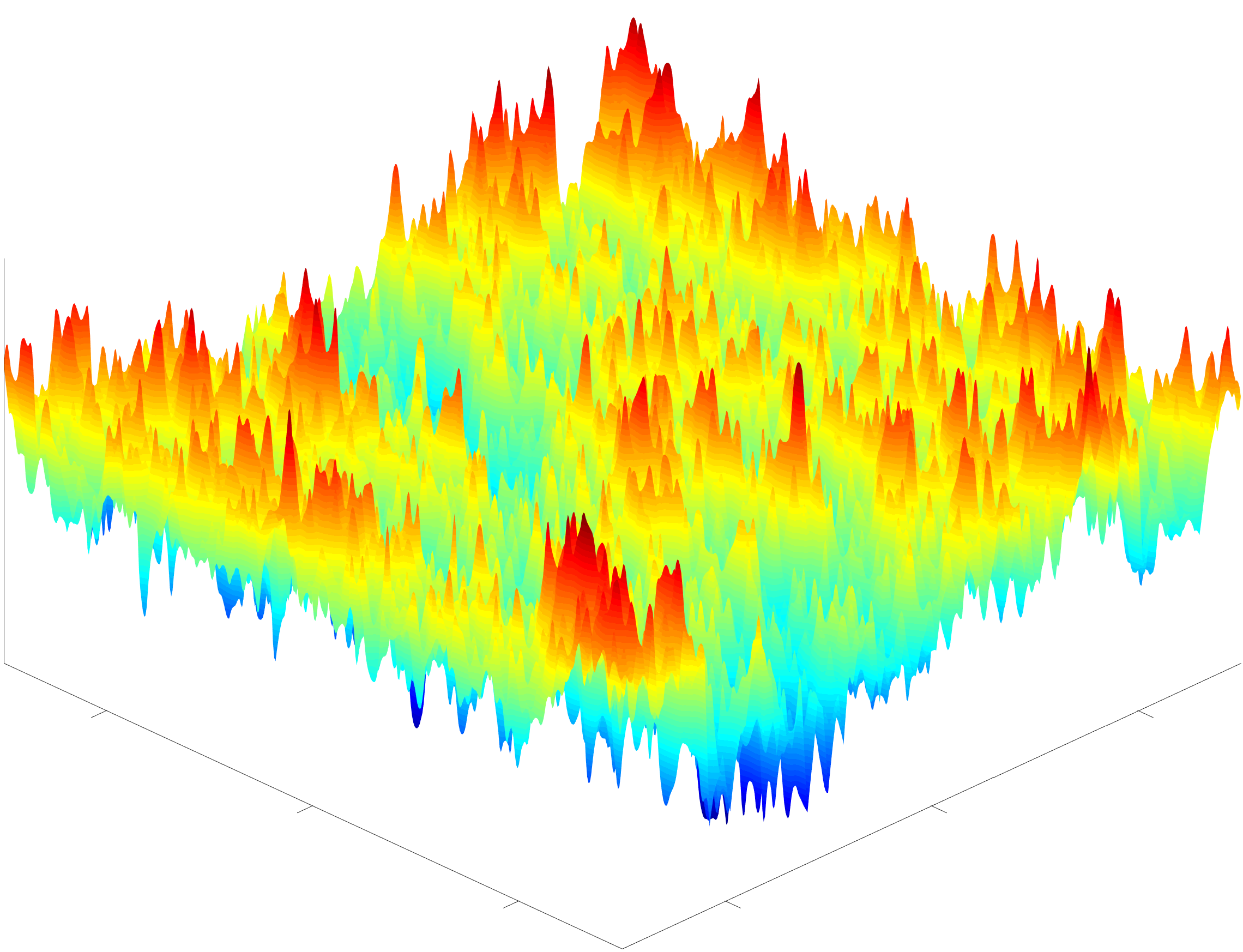}
\end{figure*}

\thispagestyle{empty}

\clearpage

\vspace{12cm}

\tableofcontents*

\clearpage



\chapter{Preface}

Many microscopic models lead to partial differential equations with rapidly oscillating coefficients. A particular example, which is the main focus of this book, is the scalar, uniformly elliptic equation
\begin{equation*} \label{}
-\nabla \cdot \left( \a(x) \nabla u \right) = f,
\end{equation*}
where the interest is in the behavior of the solutions on length scales much larger than the unit scale (the \emph{microscopic} scale on which the coefficients are varying). The coefficients are assumed to be valued in the positive definite matrices, and may be periodic, almost periodic, or stationary random fields. Such equations arise in a variety of contexts such as heat conduction and electromagnetism in heterogeneous materials, or through their connection with stochastic processes. 

\smallskip

To emphasize the highly heterogeneous nature of the problem, it is customary to introduce a parameter $0 < \ep\ll1$ to represent the ratio of the microscopic and macroscopic scales. The equation is then rescaled as
\begin{equation*} \label{}
-\nabla \cdot \left( \a\left(\tfrac x\ep\right) \nabla u^\ep \right) = f,
\end{equation*}
with the problem reformulated as that of determining the asymptotic behavior of~$u^\ep$, subject to appropriate boundary conditions, as~$\ep \to 0$. 

\smallskip

It has been known since the early 1980s that, under very general assumptions, the solution~$u^\ep$ of the heterogeneous equation converges in $L^2$ to the solution $u$ of a constant-coefficient equation
\begin{equation*} \label{}
-\nabla \cdot \left( \ahom \nabla u \right) = f.
\end{equation*}
We call this the \emph{homogenized equation} and the coefficients the \emph{homogenized} or \emph{effective coefficients}. The matrix $\ahom$ will depend on the coefficients~$\a\left( \cdot \right)$ in a very complicated fashion: there is no simple formula for~$\ahom$ except in dimension~$d=1$ and some special situations in~$d=2$. However, if one is willing to perform the computational work of approximating the homogenized coefficients and to tolerate the error in replacing~$u^\ep$ by~$u$, then there is a potentially huge payoff to be gained in terms of a reduction of the complexity of the problem. Indeed, up to a change of variables, the homogenized equation is simply the Poisson equation, which can be numerically computed in linear time and memory and is obviously independent of~$\ep>0$. In contrast, the cost of computing the solution to the heterogeneous equation explodes as $\ep$ becomes small, and can be considered out of reach. 

\smallskip

There is a vast and rich mathematical literature on homogenization developed in the last forty years and already many good expositions on the topic (see for instance the books~\cite{Allaire,BLP,BraidesDe,CD,DalMaso,JKO,KLO,PavStuart,Tartar}). Most of these works are focused on qualitative results, such as proving the existence of a homogenized equation which characterizes the limit as~$\ep \to0$ of solutions. The need to develop efficient methods for determining $\ahom$ and for estimating the error in the homogenization approximation (e.g., $\| u^\ep - u \|_{L^2}$) motivates the development of a \emph{quantitative} theory of homogenization. However, until recently, nearly all of the quantitative results were confined to the rather restrictive case of periodic coefficients. The main reason for this is that quantitative homogenization estimates in the periodic case are vastly simpler to prove than under essentially any other hypothesis (even the almost periodic case). Indeed, the problem can be essentially reduced to one on the torus and compactness arguments then yield optimal estimates. In other words, in the periodic setting, the typical arguments of qualitative homogenization theory can be made quantitative in a relatively straightforward way.

\smallskip

This book is concerned with the quantitative theory of homogenization for nonperiodic coefficient fields, focusing on the case in which $\a(x)$ is a stationary random field satisfying quantitative ergodicity assumptions. 
This is a topic which has undergone a rapid development since its birth at the beginning of this decade, with new results and more precise estimates coming at an ever accelerating pace. Very recently, there has been a convergence toward a common philosophy and set of core ideas, which have resulted in a complete and optimal theory. The purpose of this book is to give this theory a complete and self-contained presentation.

\smallskip

We have written it with several purposes and audiences in mind. Experts on the topic will find new results as well as arguments which have been greatly simplified compared to the previous state of the literature. Researchers interested in stochastic homogenization will hopefully find a useful reference to the main results in the field and a roadmap to the literature. Our approach to certain topics, such as the construction of the Gaussian free field or the relation between Sobolev norms and the heat kernel, could be of independent interest to certain segments of the probability and analysis communities. We have written the book with newcomers to homogenization in mind and, most of all, graduate students and young researchers. In particular, we expect that readers with a basic knowledge of probability and analysis, but perhaps without expertise in elliptic regularity, the Gaussian free field, negative and fractional Sobolev spaces, etc, should not have difficulty following the flow of the book. These topics are introduced as they arise and are developed in a mostly self-contained way.

\smallskip

Before we give a summary of the topics we cover and the approach we take, let us briefly recall the historical and mathematical context. In the case of stationary random coefficients, there were very beautiful, soft arguments given independently in the early 1980s by Kozlov~\cite{K1}, Papanicolaou and Varadhan~\cite{PV1} and Yurinski{\u\i}~\cite{Y0} which give proofs of qualitative homogenization under very general hypotheses. A few years later, Dal Maso and Modica~\cite{DM1,DM2} extended these results to nonlinear equations using variational arguments inspired by~$\Gamma$-convergence. Each of the proofs in these papers relies in some way on an application of the ergodic theorem applied to the gradient (or energy density) of certain solutions of the heterogeneous equation. In order to obtain a convergence rate for the limit given by the ergodic theorem, it is necessary to verify quantitative ergodic conditions on the underlying random sequence or field. It is therefore necessary and natural to impose such a condition on the coefficient field~$\a(x)$. However, even under the strongest of mixing assumptions (such as the finite range of dependence assumption we work with for most of this book), one faces the difficulty of transferring the quantitative ergodic information contained in these strong mixing properties from the coefficients to the solutions, since the ergodic theorem is applied to the latter. This is difficult because, of course, the solutions depend on the coefficient field in a very complicated, nonlinear and nonlocal way.

\smallskip

Gloria and Otto~\cite{GO1,GO2} were the first to address this difficulty in a satisfactory way in the case of coefficient fields that can be represented as functions of countably many independent random variables. They used an idea from statistical mechanics, previously introduced in the context of homogenization by Naddaf and Spencer~\cite{NS},
of viewing the solutions as functions of these independent random variables and applying certain general concentration inequalities such as the Efron-Stein or logarithmic Sobolev inequalities. If one can quantify the dependence of the solutions on a resampling of each independent random variable, then these inequalities immediately give bounds on the fluctuations of solutions. Gloria and Otto used this method to derive estimates on the first-order correctors which are optimal in terms of the ratio of length scales (although not optimal in terms of stochastic integrability). 

\smallskip

The point of view developed in this book is different and originates in works of Armstrong and Smart~\cite{AS}, Armstrong and Mourrat~\cite{AM}, and the authors~\cite{AKM1,AKM2}. Rather than study solutions of the equation directly, the main idea is to focus on certain energy quantities, which allow us to implement a progressive coarsening of the coefficient field and capture the behavior of solutions on large---but finite---length scales. The approach can thus be compared with renormalization group arguments in theoretical physics. The core of the argument is to establish that on large scales, these energy quantities are in fact essentially local, additive functions of the coefficient field. It is then straightforward to optimally transfer the mixing properties of the coefficients to the energy quantities and then to the solutions.

\smallskip

The quantitative analysis of the energy quantities is the focus of the first part of the book. After a first introductory chapter, the strategy naturally breaks into several distinct steps: 

\begin{itemize}

\item Obtaining an algebraic rate of convergence for the homogenization limits, using the subadditive and convex analytic structure endowed by the variational formulation of the equation (Chapter~\ref{c.two}). Here the emphasis is on obtaining estimates with optimal stochastic integrability, while the exponent representing the scaling of the error is suboptimal. 

\item Establishing a large-scale regularity theory: it turns out that solutions of an equation with stationary random coefficients are much more regular than one can show from the usual elliptic regularity for equations with measurable coefficients (Chapter~\ref{c.regularity}). We prove this by showing that the extra regularity is inherited from the homogenized equation by approximation, using a Campanato-type iteration and the quantitative homogenization results obtained in the previous chapter. 

\item Implementing a modification of the renormalization scheme of Chapter~\ref{c.two}, with the major additional ingredient of the large-scale regularity theory, to improve the convergence of the energy quantities to the optimal rate predicted by the scaling of the central limit theorem. Consequently, deriving optimal quantitative estimates for the first-order correctors (Chapter~\ref{c.A1}). 

\item Characterizing the fluctuations of the energy quantities by proving convergence to white noise and consequently obtaining the scaling limit of the first-order correctors to a modified Gaussian free field (Chapter~\ref{c.gff}). 

\item Combining the optimal estimates on the first-order correctors with classical arguments from homogenization theory to obtain optimal estimates on the homogenization error, and the two-scale expansion, for Dirichlet and Neumann boundary value problems (Chapter~\ref{c.twoscale}). 

\end{itemize}

These six chapters represent, in our view, the essential part of the theory. The first four chapters should be read consecutively (Sections~\ref{s.boundaryreg} and~\ref{s.regularity.optimal} can be skipped), while the Chapters~\ref{c.gff} and~\ref{c.twoscale} are independent of each other. 

\smallskip

Chapter~\ref{c.CZ} complements  the regularity theory of Chapter~\ref{c.regularity} by developing local and global gradient $L^p$ estimates ($2<p<\infty$) of Calder\'on-Zygmund-type for equations with right-hand side. Using these estimates, in Section~\ref{s.twoscale.p} we extend the results of Chapter~\ref{c.twoscale} by proving optimal quantitative bounds on the error of the two-scale expansion in $W^{1,p}$-type norms. Except for the last section, which requires the optimal bounds on the first-order correctors proved in Chapter~\ref{c.A1}, this chapter can be read after Chapter~\ref{c.regularity}. 

\smallskip

Chapter~\ref{c.parabolic} extends the analysis to the  time-dependent parabolic equation
\begin{equation*}
\partial_t u - \nabla\cdot \a\nabla u = 0. 
\end{equation*}
The main focus is on obtaining a suboptimal error estimate for the Cauchy-Dirichlet problem and a parabolic version of the large-scale regularity theory. 
Here the coefficients $\a(x)$ depend only on space, and the arguments in the chapter rely on the estimates on first-order correctors obtained in Chapters~\ref{c.two} and~\ref{c.regularity} in addition to some relatively routine deterministic arguments. In Chapter~\ref{c.parabolic} we also prove decay estimates on the elliptic and parabolic Green functions as well as on their derivatives, homogenization error and two-scale expansions. 

\smallskip

In Chapter~\ref{c.semigroup}, we study the decay,  as $t\to \infty$, of the solution $u(t,x)$ of the parabolic initial-value problem
\begin{equation*}
\left\{
\begin{aligned}
& \partial_t u -\nabla \cdot \left( \a \nabla u \right) = 0 &  \mbox{in} & \ (0,\infty) \times \Rd, \\
& u(0,\cdot) = \nabla\cdot \g & \mbox{on} & \ \Rd,
\end{aligned}
\right.
\end{equation*}
where~$\g$ is a bounded, stationary random field with a unit range of dependence. We show that the solution $u$ decays to zero at the same rate as one has in the case $\a = \Id$. As an application, we upgrade the quantitative homogenization estimates for the parabolic and elliptic Green functions to the optimal scaling (see Theorem~\ref{t.PtoPbar} and Corollary~\ref{c.GtoGbar}).

\smallskip

In Chapter~\ref{c.subadd-fitz}, we show how the variational methods in this book can be adapted to non-self adjoint operators, in other words, linear equations with nonsymmetric coefficients. In Chapter~\ref{c.nonlinear} we give a generalization to the case of nonlinear equations. In particular, in both of these chapters we give a full generalization of the results of Chapters~\ref{c.one} and~\ref{c.two} to these settings as well as the large-scale $C^{0,1}$ estimate of Chapter~\ref{c.regularity}. 

\smallskip

This version of the manuscript is essentially complete and, except for small changes and corrections and a modest expansion of Chapter~\ref{c.semigroup}, we expect to publish it in close to its present form. 

\smallskip

We would like to thank several of our colleagues and students for their helpful comments, suggestions, and corrections: Alexandre Bordas, Sanchit Chaturvedi, Paul Dario, Sam Ferguson, Chenlin Gu, Jan Kristensen, Jules Pertinand, Christophe Prange, Armin Schikorra, Charlie Smart, Tom Spencer, Stephan Wojtowytsch, Wei Wu and Ofer Zeitouni. We particularly thank Antti Hannukainen for his help with the numerical computations that generated Figure~\ref{f.correctors}. SA was partially supported by NSF Grant DMS-1700329. TK was partially supported by the Academy of Finland and he thanks Giuseppe Mingione for the invitation to give a graduate course at the University of Parma. JCM was partially supported by the ANR grant LSD (ANR-15-CE40-0020-03).

\smallskip

There is no doubt that small mistakes and typos remain in the manuscript, and so we encourage readers to send any they may find, as well as any comments, suggestions and criticisms, by email. Until the manuscript is complete, we will keep the latest version on our webpages. After it is published as a book, we will also maintain a list of typos and misprints found after publication. 

\bigskip

\begin{flushright}
Scott Armstrong, New York \par
Tuomo Kuusi, Helsinki \par
Jean-Christophe Mourrat, Paris \par
\end{flushright}



\chapter{Assumptions and examples}

We state here the assumptions which are in force throughout most of the book, and present some concrete examples of coefficient fields satisfying them.  

\subsection*{Assumptions}

Except where specifically indicated otherwise, the following standing assumptions are in force throughout the book. 

\smallskip

We fix a constant $\Lambda>1$ called the \emph{ellipticity constant}, and a dimension $d\geq 2$.

\smallskip

We let $\Omega$ denote the set of all measurable maps $\a(\cdot)$ from~$\Rd$ into the set of symmetric~$d\times d$ matrices, denoted by~$\R^{d\times d}_{\mathrm{sym}}$, which satisfy the uniform ellipticity and boundedness condition
\begin{equation} 
\label{e.ue}
 \left| \xi \right|^2 \leq \xi\cdot \a(x) \xi \leq \Lambda \left| \xi \right|^2, \quad \forall \xi \in\Rd.
\end{equation}
That is, 
\begin{equation} 
\label{e.omega}
\Omega:= 
\left\{ \a \,:\, \a \ \mbox{is a Lebesgue measurable map from $\Rd$ to $\R^{d\times d}_{\mathrm{sym}}$ satisfying~\eqref{e.ue}} \right\}. 
\end{equation}
The entries of an element $\a\in \Omega$ are written as $\a_{ij}$, $i,j\in\{1,\ldots,d\}$. 

\smallskip

We endow $\Omega$ with a family of $\sigma$-algebras $\left\{ \F_U \right\}$ indexed by the family of Borel subsets $U\subseteq \Rd$, defined by
\begin{equation}
 \label{e.FU.def}
\begin{aligned}
\F_U & :=  \mbox{the $\sigma$-algebra generated by the following family:}  \\
& \quad \quad \left\{ \a \mapsto \int_{\Rd} \a_{ij} (x) \varphi(x)\,dx \, : \, \varphi \in C^\infty_c(U), \ i,j\in \{1,\ldots,d\} \right\}.
\end{aligned}
\end{equation}
The largest of these $\sigma$-algebras is denoted $\F:= \F_{\Rd}$. For each $y\in\Rd$, we let $T_y:\Omega \to\Omega$ be the action of translation by $y$,
\begin{equation} 
\label{e.Ty}
\left( T_y\a\right)(x):= \a(x+y),
\end{equation}
and extend this to elements of~$\F$ by defining $T_y E:= \left\{ T_y\a\,:\, \a\in E \right\}$. 

\smallskip

Except where indicated otherwise, we assume throughout the book that~$\P$ is a probability measure on the measurable space $(\Omega,\F)$ satisfying the following two important assumptions:

\begin{itemize}

\item \emph{Stationarity with respect to $\Zd$-translations:}
\begin{equation} 
\label{e.assumption.stationarity}
\P \circ T_z = \P \quad \mbox{for every} \ z\in\Zd. 
\end{equation}

\item \emph{Unit range of dependence:}
\begin{equation}
\label{e.assumption.independence}
\begin{aligned}
& \F_U \ \mbox{and} \ \F_V \ \  \mbox{are $\P$-independent for every pair $U,V\subseteq \Rd$}  
\\
& \mbox{of Borel subsets satisfying $\dist(U,V)\geq 1$.}
\end{aligned}
\end{equation}
\end{itemize}
We denote the expectation with respect to $\P$ by $\E$. That is, if $X:\Omega \to \R$ is an $\F$-measurable random variable, we write
\begin{equation} 
\label{e.E.def}
\E \left[ X \right] := \int_\Omega X(\a) \, d\P(\a). 
\end{equation}
While all random objects we study in this text are functions of $\a\in\Omega$, we do not typically display this dependence explicitly in our notation. We rather use the symbol~$\a$ or~$\a(x)$ to denote the \emph{canonical coefficient field} with law~$\P$.

\subsection*{Examples satisfying the assumptions}

The simplest way to construct explicit examples satisfying the assumptions of uniform ellipticity \eqref{e.ue}, stationarity \eqref{e.assumption.stationarity} and \eqref{e.assumption.independence} is by means of a ``random checkerboard'' structure: we pave the space by unit-sized cubes and color each cube either white or black independently at random. Each color is then associated with a particular value of the diffusivity matrix. More precisely, let $(b(z))_{z \in \Zd}$ be independent random variables such that for every $z \in \Zd$,
\begin{equation*}  
\P[b(z) = 0] = \P[b(z) = 1] = \frac 1 2,
\end{equation*}
and fix two matrices $\a_0, \a_1$ belonging to the set
\begin{equation}  
\label{e.ue.set}
\Ll\{\td \a \in \R^{d \times d}_{\mathrm{sym}} \ : \ \forall \xi \in\Rd, \quad \left| \xi \right|^2 \leq \xi\cdot \td \a \xi \leq \Lambda \left| \xi \right|^2 \Rr\}.
\end{equation}
We can then define a random field $x \mapsto \a(x)$ satisfying \eqref{e.ue} and with a law satisfying \eqref{e.assumption.stationarity} and \eqref{e.assumption.independence} by setting, for every $z \in \Zd$ and  $x \in z + \Ll[ -\frac 1 2, \frac 1 2  \Rr)^d$,
\begin{equation*}  
\a(x) = \a_{b(z)}.
\end{equation*}
\begin{figure}[tb]
\centering
\begin{tikzpicture}[scale=0.267]
\draw[lightgray] (0,0) rectangle (30,30);
\foreach \n in {1,2,...,621}
\pgfmathsetmacro{\x}{random(0,29)}
\pgfmathsetmacro{\y}{random(0,29)}
\draw[fill = black] (\x,\y) rectangle (\x+1,\y+1);
\end{tikzpicture}
\caption{\small{A piece of a sample of a random checkerboard. The conductivity matrix is equal to $\a_0$ in the black region, and $\a_1$ in the white region.}}
\label{f.checkerboard}
\end{figure}
This example is illustrated on Figure~\ref{f.checkerboard}. It can be generalized as follows: we give ourselves a family $(\a(z))_{z \in \Zd}$ of independent and identically distributed (i.i.d.)\ random variables taking values in the set \eqref{e.ue.set}, and then extend the field $z \mapsto \a(z)$ by setting, for every $z \in \Zd$ and $x \in z + \Ll[ -\frac 1 2, \frac 1 2  \Rr)^d$,
\begin{equation*}  
\a(x) := \a(z). 
\end{equation*}

Another class of examples can be constructed using homogeneous Poisson point processes. We recall that a Poisson point process on a measurable space $(E,\mcl E)$ with (non-atomic, $\sigma$-finite) intensity measure $\mu$ is a random subset $\Pi$ of $E$ such that the following properties hold (see also \cite{kingman}):
\begin{itemize}
\item 
For every measurable set $A \in \mcl E$, the number of points in $\Pi \cap A$, which we denote by $N(A)$, follows a Poisson law of mean $\mu(A)$; 
\item
For every pairwise disjoint measurable sets $A_1,\ldots, A_k \in \mcl E$, the random variables $N(A_1)$, \ldots, $N(A_k)$ are independent.
\end{itemize}
Let $\Pi$ be a Poisson point process on $\Rd$ with intensity measure given by a multiple of the Lebesgue measure. Fixing two matrices $\a_0, \a_1$ belonging to the set \eqref{e.ue.set}, we may define a random field $x \mapsto \a(x)$ by setting, for every $x \in \Rd$,
\begin{equation}  
\label{e.easy.bubbles}
\a(x) := 
\Ll\{
\begin{aligned}
 \a_0 & \quad \text{if } \dist(x,\Pi) \le \frac 1 2, \\
 \a_1 & \quad \text{otherwise}.
\end{aligned}
\Rr.
\end{equation}

\begin{figure}[tb]
\centering
\begin{tikzpicture}[scale=0.8]
\clip (0,0) rectangle (10,10);
\draw[lightgray] (0,0) rectangle (10,10);
\foreach \n in {1,2,...,275}
\pgfmathsetmacro{\x}{(1+rand)*5.5-0.5}
\pgfmathsetmacro{\y}{(1+rand)*5.5-0.5}
\draw[fill = black] (\x,\y) circle (0.2);
\end{tikzpicture}
\caption{
\small{
A sample of the coefficient field defined in~\eqref{e.easy.bubbles} by the Poisson point cloud. The matrix~$\a$ is equal to~$\a_0$ in the black region and to~$\a_1$ in the white region.
}
}
\label{f.easy.bubbles}
\end{figure}

\begin{figure}[tb]
\centering
\begin{tikzpicture}[scale=0.8]
\clip (0,0) rectangle (10,10);
\draw[lightgray] (0,0) rectangle (10,10);
\foreach \n in {1,2,...,275}
\pgfmathsetmacro{\X}{(1+rand)*5.5-0.5}
\pgfmathsetmacro{\Y}{(1+rand)*5.5-0.5}
\draw[fill=black,smooth] (\X,\Y) \irregularcircle{0.2}{0.05};	
\end{tikzpicture}
\caption{
\small{
This coefficient field is sampled from the same distribution as in~Figure~\ref{f.easy.bubbles}, except that the balls have been replaced by random shapes.
}
}
\label{f.hard.bubbles}
\end{figure}

\begin{figure}
\centering
\includegraphics[width=.8\linewidth]{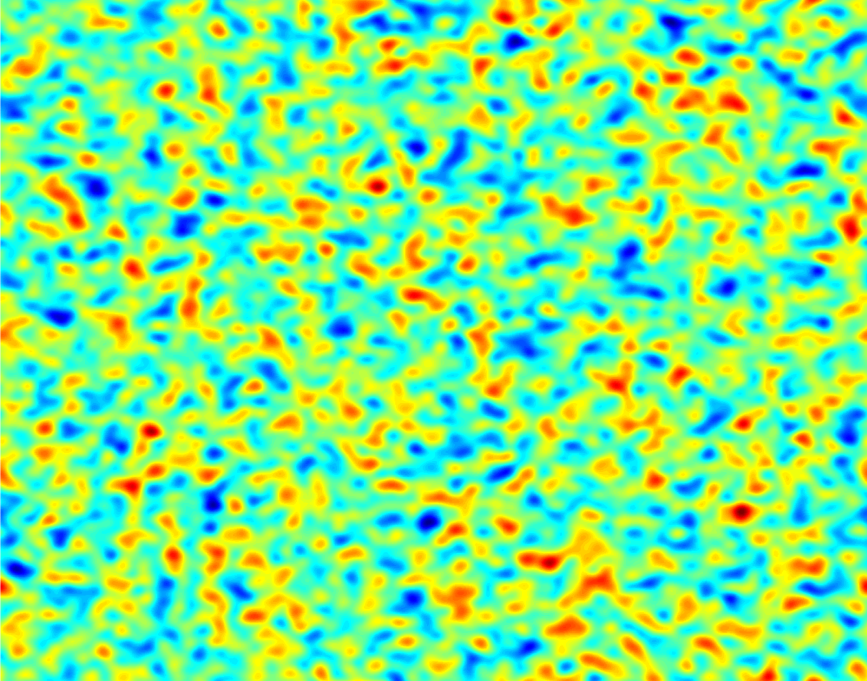} 
\caption{\small The figure represents the convolution of white noise with a smooth function of compact support, using a color scale. This scalar field can be used to construct a matrix field $x \mapsto \a(x)$ satisfying our assumptions, see \eqref{e.def.local.white}.}
\label{f.whitenoise}
\end{figure}

This example is illustrated on Figure~\ref{f.easy.bubbles}. More complicated examples can be constructed using richer point processes. For instance, in the construction above, each point of $\Pi$ imposes the value of $\a(x)$ in a centered ball of radius $1/2$; we may wish to construct examples where this radius itself is random. In order to do so, let $\lambda > 0$, let $\mu$ denote a probability measure on $\Ll[ 0,\frac 1 2 \Rr]$ (the law of the random radius), and let $\Pi$ be a Poisson point process on $\Rd \times \R$ with intensity measure $\lambda dx \otimes \mu$ (where $dx$ denotes the Lebesgue measure on $\Rd$). We then set, for every $x \in \Rd$,
\begin{equation*}  
\a(x) := 
\Ll\{
\begin{aligned}
 \a_0 & \quad \text{if there exists } (z,r) \in \Pi  \text{ such that }  |x-z| \le r , \\
 \a_1 & \quad \text{otherwise}.
\end{aligned}
\Rr.
\end{equation*}
Minor variants of this example allow for instance to replace balls by random shapes,  to allow the conductivity matrix to take more than two values, etc. See Figure~\ref{f.hard.bubbles} for an example. 

\smallskip

Yet another class of examples can be obtained by defining the coefficient field $x \mapsto \a(x)$ as a local function of a white noise field. We refer to Definition~\ref{d.scalar.noise} and Proposition~\ref{p.white.noise.exists} for the definition and construction of white noise. For instance, given a scalar white noise $\W$, we may fix a smooth function $\phi \in C^\infty_c(\Rd)$ with support in $B_{1/2}$, a smooth function $F$ from $\Rd$ into the set \eqref{e.ue.set}, and define
\begin{equation}  
\label{e.def.local.white}
\a(x) = F((\W \ast \phi)(x)).
\end{equation}
See Figure~\ref{f.whitenoise} for a representation of the scalar field $x \mapsto (\W \ast \phi)(x)$.



\chapter{Frequently asked questions}

\subsection*{Where is the independence assumption used?} 
The unit range of dependence assumption~\eqref{e.assumption.independence} is obviously very important, and to avoid diluting its power we use it sparingly. We list here all the places in the book where it is invoked:
\begin{itemize}
\item The proof of Proposition~\ref{p.nuconv} (which is made redundant by the following one).
\item The proof of Lemma~\ref{l.spatavg} (and the generalizations of this lemma appearing in Chapters~\ref{c.subadd-fitz} and~\ref{c.nonlinear}). This lemma lies at the heart of the iteration argument in Chapter~\ref{c.two}, as it is here that we obtain our first estimate on the correspondence between spatial averages of gradients and fluxes of solutions. Notice that the proof does not use the full strength of the independence assumption, it actually requires only a very weak assumption of correlation decay. 

\item The last step of the proof of Theorem~\ref{t.subadd} (and the generalizations of this theorem appearing in Chapters~\ref{c.subadd-fitz} and~\ref{c.nonlinear}). Here independence is used very strongly to obtain homogenization estimates with optimal stochastic integrability. 
\item The proof of Proposition~\ref{p.improvefluc} in Section~\ref{s.fluctuations}, where we control the fluctuations of the quantity~$J_1$ inside the bootstrap argument. 
\item The proof of Proposition~\ref{p.d2vengeance} in Section~\ref{s.twopoint}, where we prove sharper bounds on the first-order correctors in dimension~$d=2$.  
\item In Section~\ref{s.CLT}, where we prove the central limit theorem for the quantity~$J_1$. This can be considered a refinement of Proposition~\ref{p.improvefluc}. 

\item In Section~\ref{s.sgp.opt} in the proofs of Lemmas~\ref{l.SG.opt1} and~\ref{l.SG.opt4}.
\end{itemize}
In particular,  all of the results of Chapters~\ref{c.two} and~\ref{c.regularity} are obtained with only two very straightforward applications of independence.

\subsection*{Can the independence assumption be relaxed?}

Yes. One of the advantages of the approach presented in this book is that \emph{the independence assumption is applied only to sums of local random variables.} Any reasonable decorrelation condition or mixing-type assumption will give estimates regarding the stochastic cancellations of sums of local random variables (indeed, this is essentially a tautology). Therefore, while the statements of the theorems may need to be modified for weaker assumptions (for instance, the strong stochastic integrability results we obtain under a finite range of dependence assumption may have to be weakened), the proofs will only require straightforward adaptations. In fact, since we have only used independence in a handful of places in the text, enumerated above, it is not a daunting task to perform these adaptations. This is in contrast to alternative approaches in quantitative stochastic homogenization which use nonlinear concentration inequalities and therefore are much less robust to changes in the hypotheses. 

\smallskip

The reason for formalizing the results under the strongest possible mixing assumption (finite range of dependence) rather than attempting write a very general result is, therefore, not due to a limitation of the arguments. It is simply because we favor clarity of exposition over generality.

\subsection*{Can the uniform ellipticity assumption be relaxed?}

One of the principles of this book is that one should avoid using small-scale or pointwise properties of the solutions or of the equation and focus rather on large-scale, averaged information. In particular, especially in the first part of the book, we concentrate on the energy quantities~$\nu$,~$\nu^*$ and~$J_1$ which can be thought of as ``coarsened coefficients'' in analogy to a renormalization scheme (see Remark~\ref{r.coarsened}). The arguments we use adhere to this philosophy rather strictly. In particular, they are adaptable to situations in which the matrix $\a(x)$ is not necessarily uniformly positive definite, provided we have some quantitative information regarding the law of its condition number. This is because such assumptions can be translated into quantitative information about~$J_1$ which suffices to run the renormalization arguments of Chapter~\ref{c.two}. A demonstration of the robustness of these methods can also be found in~\cite{AD}, which adapted Chapters~\ref{c.two} and~\ref{c.regularity} of this manuscript to obtain the first quantitative homogenization results on supercritical percolation clusters (a particularly extreme example of a degenerate environment).

\subsection*{Do the results in this book apply to elliptic systems?}

Since the notation for elliptic systems is a bit distracting, we have decided to use scalar notation. However, throughout most of the book, we use exclusively arguments which also work for systems of equations (satisfying the uniform ellipticity assumption). 
The only exceptions are the last two sections of Chapter~\ref{c.parabolic} and Chapter~\ref{c.semigroup}, where we do use some scalar estimates (the De Giorgi-Nash $L^\infty$ bound and variations) which make it easier to work with Green functions. 
In particular, we claim that all of the statements and proofs appearing in this book,
with the exception of those appearing in those two chapters, 
can be adapted to the case of elliptic systems with easy and straightforward modifications to the notation.

\subsection*{This book is written for equations in the continuum. Do the arguments apply to finite difference equations on $\Zd$?}

The techniques developed in the book are robust to the underlying structure of the environment on the microscopic scale. What is the important is that the ``geometry'' of the macroscopic medium is like that of~$\Rd$  in the sense that certain functional inequalities (such as the Sobolev inequality) hold on large scales. In the case that~$\Rd$ is replaced by~$\Zd$, the modifications are relatively straightforward: besides changes to the notation, there is just the slight detail that the boundary of a large cube has a nonzero volume, which creates an additional error term in Chapter~\ref{c.two} causing no harm. If one has a more complicated microstructure like a random graph, such as a supercritical Bernoulli percolation cluster, it is necessary to first establish the ``geometric regularity'' of the graph in the sense that Sobolev-type inequalities hold on large scales. The techniques described in this book can then be readily applied: see~\cite{AD}. 

\subsection*{Can I find a simple proof of \emph{qualitative} homogenization somewhere here?}

The arguments in Chapter~\ref{c.one} only need to be slightly modified in order to obtain a more general qualitative homogenization result valid in the case that the unit range of dependence assumption is relaxed to mere \emph{ergodicity}. In other words, in place of~\eqref{e.assumption.independence} we assume instead that 
\begin{equation}
\label{e.onlyergodic}
\begin{aligned}
\mbox{if $A\in \F$ satisfies $T_zA=A$ for all $z\in\Zd$, then $\P\left[A \right] \in \{ 0,1\}$.} 
\end{aligned}
\end{equation}
In fact, the only argument that needs to be modified is the proof of Proposition~\ref{p.nuconv}, since it is the only place in the chapter where independence is used. Moreover, that argument is essentially a proof of the subadditive ergodic theorem in the special case of the unit range of dependence assumption~\eqref{e.assumption.independence}. In the general ergodic case~\eqref{e.onlyergodic}, one can simply directly apply the subadditive ergodic theorem (see for instance~\cite{AK}) to obtain, in place of~\eqref{e.preL1conv}, the estimate
\begin{equation*} 
\P \left[ \limsup_{n\to \infty}  \left| \a(\cu_n) - \ahom  \right| =0  \right] 
=1.
\end{equation*}
The rest of the arguments in that chapter are deterministic and imply that the random variable~$\mathcal{E}'(\ep)$ in Theorem~\ref{t.DP.blackbox} satisfies~$\P \left[ \limsup_{\ep \to 0} \mathcal{E}'(\ep) = 0 \right] =1$. 

\subsection*{What do we learn about reversible diffusions in random environments?}

Just as we learn about Brownian motion from properties of harmonic functions (and conversely), the study of divergence-form operators gives us information about the associated diffusion processes. To start with, De Giorgi-Nash-Aronson estimates recalled in \eqref{e.upperP.NA} and Proposition~\ref{p.par.reg.NA} can be used together with the classical Kolmogorov extension and continuity theorems (see \cite[Theorem~36.2]{billingsley} and \cite[Theorem~I.2.1]{revuz-yor}) to define these stochastic processes. Denoting by $\mathbf{P}_x^\a$ the probability law of the diffusion process starting from $x\in\Rd$, and by $(X(t))_{t \ge 0}$ the canonical process, we have by construction that, for every $\a\in\Omega$, Borel measurable set $A \subset \Rd$ and $(t,x)\in (0,\infty)\times\Rd$, 
\begin{equation}  
\label{e.local.clt}
\mathbf{P}_x^\a \Ll[ X_t \in A \Rr] = \int_A P(t,x,y) \, dy,
\end{equation}
where $P(t,x,y)$ is the parabolic Green function defined in Proposition~\ref{p.GFexistence.NA}. 
The statement 
\begin{equation*}  
\text{for every } x \in \Rd, \quad t^{\frac d 2} P(t,0, t^\frac 1 2  x) \xrightarrow[t\to \infty]{\text{a.s.}}\bar P(1,0,x),
\end{equation*}
where $\bar P$ is the parabolic Green function for the homogenized operator, can thus be interpreted as a (quenched) local central limit theorem for the diffusion process. Seen in this light, Theorem~\ref{t.PGF.basecase} gives us a first quantitative version of this local central limit theorem. The much more precise Theorem~\ref{t.PtoPbar} gives an optimal rate of convergence for this statement, and can thus be interpreted as analogous to the classical Berry-Esseen result \index{Berry-Esseen theorem} on the rate of convergence in the central limit theorem for sums of independent random variables (see \cite[Theorem 5.5]{Petrov}).



\chapter{Notation}

\subsection*{Sets and Euclidean space}

The set of nonnegative integers is denoted by~$\N:=\{0,1,2,\ldots\}$. The set of real numbers is written~$\R$. When we write $\R^m$ we implicitly assume that $m\in\N\setminus\{0\}$. For each  $x,y\in\R^m$, the scalar product of $x$ and $y$ is denoted by $x\cdot y$, their tensor product by $x\otimes y$ and the Euclidean norm on $\R^m$ is~$\left| \,\cdot\, \right|$. The canonical basis of $\R^m$ is written as $\{e_1,\ldots,e_m\}$. We let~$\mathcal{B}$ denote the Borel $\sigma$--algebra on~$\R^m$. A \emph{domain} is an open connected subset of~$\R^m$. The notions of~\emph{$C^{k,\alpha}$ domain} and \emph{Lipschitz domain} are defined in Definition~\ref{def.Lipdomain}. The boundary of~$U\subseteq \R^m$ is denoted by $\partial U$ and its closure by~$\overline{U}$. The open ball of radius $r>0$ centered at $x\in\R^m$ is $B_r(x):= \left\{ y\in\R^m\,:\, |x-y|<r \right\}$. The distance from a point to a set $V\subseteq \R^m$ is written $\dist(x,V):= \inf\left\{ |x-y|\,:\, y\in V\right\}$. For $r>0$ and $U \subseteq\R^m$, we define
\begin{equation}
\label{e.Ur.def}
U_r:= \left\{ x\in U\,:\, \dist(x,\partial U)>r \right\} \quad \mbox{and} \quad U^r := \left\{ x\in \R^m\,:\, \dist(x,U) < r \right\}.
\end{equation}
For $\lambda >0$, we set $\lambda U:= \left\{ \lambda x \,:\, x\in U\right\}$. 
If $m,n\in\N\setminus\{0\}$, the set of $m\times n$ matrices with real entries is denoted by $\R^{m\times n}$. We typically denote an element of $\R^{m\times n}$ by a boldfaced latin letter, such as $\mathbf{m}$, and its entries by $(\mathbf{m}_{ij})$. The subset of $\R^{n\times n}$ of symmetric matrices is written $\R^{n\times n}_{\mathrm{sym}}$ and the set of $n$-by-$n$ skew-symmetric matrices is~$\R^{n\times n}_{\mathrm{skew}}$. The identity matrix is denoted~$\Id$. If $r,s\in \R$ then we write $r\vee s := \max\{ r,s \}$ and $r\wedge s := \min\{r,s\}$. We also denote $r_+:=r \vee 0$ and $r_-:= -(r\wedge 0)$. 

\smallskip

We use triadic cubes throughout the book. For each $m\in\N$, we denote 
\begin{equation} 
\label{e.triad.def}
\cu_m := \left( -\frac12 3^m, \frac12 3^m \right)^d\subseteq \Rd. 
\end{equation}
Observe that, for each $n\in\N$ with $n\leq m$, the cube $\cu_m$ can be partitioned (up to a set of zero Lebesgue measure) into exactly $3^{d(m-n)}$ subcubes which are $\Zd$-translations of $\cu_n$, namely $\left\{ z+ \cu_n\,:\, z\in 3^n\Zd\cap \cu_m \right\}$.

\subsection*{Calculus}

If $U\subseteq\Rd$ and $f: U \to \R$, we denote the partial derivatives of $f$ by $\partial_{x_i} f$ or simply $\partial_i f$, which unless otherwise indicated, is understood in the sense of distributions. The gradient of $f$ is denoted by $\nabla f:= \left( \partial_1 f ,\ldots, \partial_d f \right)$. The Hessian of $f$ is denoted by $\nabla^2f:= \left( \partial_{i} \partial_j f \right)_{i,j\in\{1,\ldots,d\}}$, and higher derivatives are denoted similarly:
\begin{equation*} \label{}
\nabla^kf := \left( \partial_{i_1}\cdots\partial_{i_k} f \right)_{i_1,\ldots,i_k\in\{1,\ldots,d\}}
\end{equation*}
A vector field on~$U\subseteq \Rd$, typically  denoted by a boldfaced latin letter, is a function~$\mathbf{f}: U \to \Rd$. The divergence of~$\f$ is~$\nabla \cdot \f = \sum_{k=1}^d \partial_i f_i$, where the $(f_i)$ are the entries of $\f$, i.e., $\f = \left( f_1,\ldots f_d\right)$. 

\smallskip

\subsection*{H\"older and Lebesgue spaces}

For $k\in\N\cup\{\infty\}$, the set of functions $f:U \to \R$ which are $k$ times continuously differentiable in the classical sense is denoted by~$C^k(U)$. We denote by $C^k_c(U)$ the collections of $C^k(U)$ functions with compact support in $U$. For $k\in\N$ and $\alpha \in (0,1]$, we denote the classical H\"older spaces by $C^{k,\alpha}(U)$, which are the functions $u\in C^k(U)$ for which the norm 
\begin{equation*} \label{}
\left\| u \right\|_{C^{k,\alpha}(U)} :=
\sum_{n=0}^k \sup_{x\in U} \left| \nabla^n u(x) \right|
+ \left[\nabla^ku\right]_{C^{0,\alpha}(U)}
\end{equation*}
is finite, where~$\left[\,\cdot\,\right]_{C^{0,\alpha}(U)}$ is the seminorm defined by
\begin{equation*} \label{}
\left[ u \right]_{C^{0,\alpha}(U)} 
:=  \sup_{x,y\in U,\,x\neq y} \frac{ \left|  u(x) - u(y) \right|}{|x-y|^\alpha}.
\end{equation*}

\smallskip

For every Borel set $U \in \mathcal{B}$, we denote by $|U|$ the Lebesgue measure of $U$. For an integrable function $f:U \to \R$, we may denote the integral of $f$ in a compact notation by
\begin{equation*} \label{}
\int_U f:=  \int_U f(x)\,dx. 
\end{equation*}
For $U\subseteq \Rd$ and $p\in[1,\infty]$, we denote by~$L^p(U)$ the Lebesgue space on~$U$ with exponent~$p$, that is, the set of measurable functions $f:U \to \R$ satisfying
\begin{equation*} \label{}
\left\| f \right\|_{L^p(U)} := \left( \int_U \left| f \right|^p \right)^{\frac1p} < \infty. 
\end{equation*}
The vector space of functions on $U$ which belong to $L^p(V)$ whenever $V$ is bounded and $\overline{V} \subseteq U$ is denoted by $L^p_{\mathrm{loc}}(U)$. If $|U| < \infty$ and $f \in L^1(U)$, then we write 
\begin{equation*}  
\fint_U f := \frac 1 {|U|} \int_U f.
\end{equation*}
The average of a function $f\in L^1(U)$ on~$U$ is also sometimes denoted by
\begin{equation*} \label{}
\left( f \right)_U:= \fint_Uf.
\end{equation*}
To make it easier to keep track of scalings, we very often work with rescaled versions of $L^p$ norms: for every $p \in [1,\infty)$ and $f \in L^p(U)$, we set
\begin{equation*}  
\|f\|_{\nL^p(U)} := \Ll( \fint_U |f|^p \Rr) ^{\frac 1 p} = | U |^{-\frac 1p} \left\| f \right\|_{L^p(U)}.
\end{equation*}
For convenience, we may also use the notation $\|f\|_{\nL^\infty(U)} := \|f\|_{L^\infty(U)}$. If $X$ is a Banach space, then $L^p(U;X)$ denotes the set of measurable functions $f:U \to X$ such that $x\mapsto \left\| f(x) \right\|_{X} \in L^p(U)$. We denote the corresponding norm by~$\left\| f \right\|_{L^p(U;X)}$. 
By abuse of notation, we will sometimes write $\f \in L^p(U)$ if $\f:U \to \R^m$ is a vector field such that $|\f|\in L^p(U)$ and denote $\| \f \|_{\underline{L}^p(U)}:= \| \f  \|_{\underline{L}^p(U;\R^m)} = \| |\f|  \|_{\underline{L}^p(U)}$. For $f\in L^p(\Rd)$ and $g\in L^{p'}(\Rd)$ with $\frac1p+\frac1{p'}=1$, we denote the convolution of~$f$ and~$g$ by
\begin{equation*} \label{}
( f\ast g )(x) := \int_{\Rd} f(x-y)g(y)\,dy. 
\end{equation*}

\subsection*{Special functions}

For $p\in\Rd$, we denote the affine function with slope $p$ passing through the origin by
\begin{equation*} \label{}
\ell_p(x):= p\cdot x. 
\end{equation*}
Unless otherwise indicated, $\zeta \in C^\infty_c(\Rd)$ denotes the standard mollifier
\index{standard mollifier~$\zeta$}
\begin{equation}
 \label{e.standardmollifier}
\zeta(x) : = \left\{ 
\begin{aligned}
& c_d \, \exp\left( - (1-|x|^2)^{-1} \right) & \mbox{if} & \ |x|<1,\\
& 0 & \mbox{if} & \ |x| \geq 1,
\end{aligned}
\right.
\end{equation}
with the multiplicative constant~$c_d$ chosen so that $\int_{\R^d} \zeta = 1$. We denote, for~$\delta>0$,
\begin{equation} 
\label{e.standardmollifer.delta}
\zeta_\delta(x):= \delta^{-d} \zeta\left( \frac x\delta \right). 
\end{equation}
The standard heat kernel\index{heat kernel~$\Phi$} is denoted by
\begin{equation*} \label{}
\Phi(t,x):= \left( 4\pi t\right)^{-\frac d2} \exp\left( -\frac{|x|^2}{4t} \right)
\end{equation*}
and define, for each $z\in\Rd$ and $r>0$, 
\begin{equation*} \label{}
\Phi_{z,r} (x):= \Phi(r^2, x-z) \quad \mbox{and} \quad \Phi_r := \Phi_{0,r}. 
\end{equation*}
We also denote by~$\mathcal{P}_k$ the set of real polynomials on~$\Rd$ of order at most~$k$. 

\subsection*{Sobolev and fractional Sobolev spaces}

For $k \in \N$ and $p \in [1,\infty]$, we denote by $W^{k,p}(U)$ the classical Sobolev space, see Definition~\ref{def.Sobolev} in Appendix~\ref{a.sobolev}. The corresponding norm is~$\left\| \,\cdot \, \right\|_{W^{k,p}(U)}$.  For $\alpha\in (0,\infty)\setminus \N$, the space $W^{\alpha,p}(U)$ is the fractional Sobolev space introduced in~Definition~\ref{def.sobolev.fractional}, with norm $\left\| \,\cdot\,\right\|_{W^{\alpha,p}(U)}$. For $\alpha\in (0,\infty)$, we denote by $W^{\alpha,p}_0(U)$ the closure of $C^\infty_c(U)$ in $W^{\alpha,p}(U)$. We also use the shorthand notation
\begin{equation}
\label{e.sobolev.hilbert}
H^\alpha(U) := W^{\alpha,2}(U) \quad \text{and} \quad H^\alpha_0(U) := W^{\alpha,2}_0(U)
\end{equation}
with corresponding norms~$\left\| \,\cdot \, \right\|_{H^{\alpha}(U)}=\left\| \,\cdot \, \right\|_{W^{\alpha,p}(U)}$.

\smallskip

As for $L^p$ spaces, it is useful to work with normalized and scale-invariant versions of the Sobolev norms. We define the rescaled~$W^{k,p}(U)$ norm of a function~$u\in W^{k,p}(U)$ by 
\begin{equation*} \label{}
\left\| u \right\|_{\underline{W}^{k,p}(U)} 
:= \sum_{j=0}^k \left| U \right|^{\frac{j-k} d} \left\| \nabla^j u \right\|_{\underline{L}^p(U)} .
\end{equation*}
Observe that 
\begin{equation} 
\label{e.notation.Sobolev.scaling}
\left\| u\left( \tfrac{\cdot}{\lambda} \right) \right\|_{\underline{W}^{k,p}(\lambda U)}  
= \lambda^{-k} \left\| u \right\|_{\underline{W}^{k,p}(U)} .
\end{equation}
We also set $\left\| u \right\|_{\underline{H}^1(U)} := \left\| u \right\|_{\underline{W}^{1,2}(U)}$. We use the notation~$W^{\alpha,p}_{\mathrm{loc}}(U)$ and $H^\alpha_{\mathrm{loc}}(U)$ for the spaces defined analogously to~$L^p_{\mathrm{loc}}(U)$.

\subsection*{Negative Sobolev spaces}

For $\alpha\in (0,\infty)$, $p\in [1,\infty]$ and $p':= \frac{p}{p-1}$ the H\"older conjugate exponent of $p$, the space $W^{-\al,p}(U)$ is the space of distributions~$u$ such that the following norm is finite: 
\begin{equation}  
\label{e.first.def.w-alpha.norm}
\|u\|_{W^{-\al,p}(U)} := \sup  \Ll\{ \int_U u v \ : \ v \in C^\infty_c(U), \ \|v\|_{W^{\al,p'}(U)} \le 1 \Rr\}.
\end{equation}
We also set $H^{-\al}(U) := W^{-\alpha,2}(U)$. 
For $p > 1$, the space $W^{-\alpha,p}(U)$ is the space dual to $W^{\al,p'}_0(U)$, and we have
\begin{equation}
\label{e.second.def.w-alpha.norm}
\|u\|_{W^{-\al,p}(U)} = \sup  \Ll\{ \int_U u v \ : \ v \in W^{\al,p'}_0(U), \ \|v\|_{W^{\al,p'}(U)} \le 1 \Rr\}.
\end{equation}
We refer to Definition~\ref{def.negsob} and Remark~\ref{r.non.separable.sobolev} for details. The rescaled $W^{-\alpha,p}(U)$ norm is defined by
\begin{equation} 
\label{e.first.def.rescaled.w-alpha.norm}
\|u\|_{\underline W^{-\al,p}(U)} := \sup  \Ll\{ \fint_U u v \ : \ v \in C^\infty_c(U), \ \|v\|_{\underline W^{\al,p'}(U)} \le 1 \Rr\}.
\end{equation}
We also set $\left\| \,\cdot\,\right\|_{\underline{H}^{-1}(U)} = \left\| \,\cdot\,\right\|_{\underline{W}^{-1,2}(U)}$.
These rescaled norms behave under dilations in the following way:
\begin{equation}  \label{e.notation.Sobolev.scaling-}
\left\| u \left( \tfrac{\cdot}{\lambda} \right) \right\|_{\underline{W}^{-\alpha,p}(\lambda U)} = \lambda^{\alpha} \left\| u  \right\|_{\underline{W}^{-\alpha,p}(U)}.
\end{equation}

\smallskip

Note that we have abused notation in \eqref{e.first.def.w-alpha.norm}-\eqref{e.second.def.w-alpha.norm}, denoting by $\int_U uv$ the duality pairing between $u$ and $v$. In other words, ``$\int_U u v$'' denotes the duality pairing that is normalized in such a way that if $u, v \in C^\infty_c(U)$, then the notation agrees with the usual integral. In \eqref{e.first.def.rescaled.w-alpha.norm}, we understand that $\fint_U uv = |U|^{-1} \int_U uv$.
\smallskip

It is sometimes useful to consider the slightly different space $\Hminus(U)$ which is the dual space to $H^1(U)$. We denote its (rescaled) norm by
\begin{equation}
\label{e.def.H-1}
\|u\|_{\Hminusul(U)} := \sup\left\{  \fint_U uv  \,:\, v \in H^1(U) 
\ \mbox{and} \ 
\left\| v \right\|_{\underline{H}^1(U)} \leq 1
\right\}.
\end{equation}
It is evident that $\Hminus(U) \subseteq H^{-1}(U)$ and that we have
\begin{equation} 
\label{e.H.vs.hatH}
\|u\|_{\nH^{-1}(U)} 
\leq 
\|u\|_{\Hminusul(U)}.
\end{equation}

\subsection*{Solenoidal and potential fields}

We let $\Lpot(U)$ and $\Ls(U)$ denote the closed subspaces of $L^2(U;\Rd)$ consisting respectively of \emph{potential} (gradient) and \emph{solenoidal} (divergence-free) vector fields. These are defined by
\begin{equation}
\label{e.Lpot.not}
\Lpot(U):= \left\{ \nabla u \,:\,  u \in H^1(U) \right\}
\end{equation}
and
\begin{equation}
\label{e.Ls.not}
\Ls(U):= \left\{ \g \in L^2(U;\Rd) \,:\, \nabla \cdot \g = 0 \right\}.
\end{equation}
Here we interpret the condition~$\nabla \cdot \g = 0$ in the sense of distributions. In this case, it is equivalent to the condition
\begin{equation*} \label{}
\forall \phi\in H^1_0(U), \ \int_U \g\cdot \nabla \phi = 0.
\end{equation*}
We can also make sense of the condition $\nabla\cdot \g=0$ if the entries of the vector field~$\g$ are distributions by restricting the condition above to $\phi \in C^\infty_c(U)$. 
We also set 
\begin{equation} 
\label{e.def.Lp0}
\Lpoot(U):= \left\{ \nabla u \,:\, u\in H^1_0(U) \right\}
\end{equation}
and
\begin{equation}
\label{e.def.Ls0}
\Lso(U):= \left\{ \g \in L^2(U;\Rd) \,:\, \forall \phi\in H^1(U), \ \int_U \g\cdot \nabla \phi = 0 \right\}.
\end{equation}
Notice that from these definitions we immediately have the Helmholtz-Hodge orthogonal decompositions
\begin{equation} 
\label{e.helmoltz-hodge}
L^2(U;\Rd) = \Lpoot(U) \oplus \Ls(U) = \Lpot(U) \oplus \Lso(U),
\end{equation}
which is understood with respect to the usual inner product on~$L^2(U;\Rd)$. Finally, we denote 
\begin{equation*} \label{}
L^2_{\mathrm{pot},\,\mathrm{loc}}(U):= \left\{ \nabla u \,:\, u \in H^1_{\mathrm{loc}}(U) \right\}
\end{equation*}
and
\begin{equation*} \label{}
L^2_{\mathrm{sol},\,\mathrm{loc}}(U):= \left\{ \g \in L^2_{\mathrm{loc}}(U;\Rd) \,:\, \nabla \cdot \g = 0 \right\}.
\end{equation*}

\subsection*{The probability space}

As explained previously in the assumptions section, except where explicitely indicated to the contrary,~$(\Omega,\F)$ refers to the pair defined in~\eqref{e.omega} and~\eqref{e.FU.def} which is endowed by the probability measure~$\P$ satisfying the assumptions~\eqref{e.assumption.stationarity} and~\eqref{e.assumption.independence}. The constants~$\Lambda\geq 1$ and $d\in\N$ in the definition of~$\Omega$ (i.e., the ellipticity constant and the dimension) are fixed throughout the book. 

A \emph{random element} on a measurable space~$(G,\mathcal{G})$ is an $\F$--measurable map~$X:\Omega\to G$. If $G$ is also a topological space and we do not specify $\mathcal{G}$, it is understood to be the Borel~$\sigma$-algebra. In the case that~$G=\R$ or~$G=\R^m$, we say respectively that~$X$ is a \emph{random variable} or \emph{vector}. If $G$ is a function space on $\Rd$, then we typically say that $X$ is a \emph{random field}. If $X$ is a random element taking values in a topological vector space, then we  define $\E\left[ X \right]$ to be the expectation with respect to~$\P$ via the formula~\eqref{e.E.def}. Any usage of the expression ``almost surely,'' also abbreviated ``a.s.''~or ``$\P$--a.s.,''~is understood to be interpreted with respect to~$\P$. For example, if $\{ X, X_1, X_2, \ldots \}_{n\in\N}$ is a sequence of random variables, then we may write either ``$X_n \to X$ as $n\to \infty$, $\P$--a.s.,'' or 
\begin{equation*} \label{}
X_n \xrightarrow[n\to \infty]{\text{a.s.}}  X
\end{equation*}
in place of the statement~$\P \left[ \limsup_{n\to \infty} \left|X_n-X \right| = 0 \right] = 1$.

\smallskip

We denote by~$\a$ the \emph{canonical element} of~$\Omega$ which we consider to be a random field by identifying it with the identity map~$\a \mapsto \a$. Therefore, all random elements can be considered as functions of~$\a$ and conversely, any object which we can define as a function of~$\a$ in an~$\F$--measurable way. Since this will be the case of nearly every object in the entire book, we typically do not display the dependence on~$\a$ explicitly. 

\smallskip

We note that this convention is in contrast to the usual one made in the literature on stochastic homogenization, where it is more customary to let $(\Omega,\F)$ be an abstract probability space, denote a typical element of $\Omega$ by~$\omega$, and consider the coefficient field $\a=\a(x,\omega)$ to be a matrix-valued mapping on $\Rd \times \Omega$. The~$\omega$ then appears everywhere in the notation, since all objects are functions of~$\omega$, which in our opinion makes the notation unnecessarily heavy without bringing any additional clarity. We hope that, for readers who find this notation more familiar, it is nevertheless obvious that it is equivalent to ours via the identification of the abstract $\Omega$ with our~$\Omega$ by the map $\omega\mapsto \a(\cdot,\omega)$.

\subsection*{Solutions of the PDE}

Given a coefficient field $\a\in\Omega$ and an open subset $U\subseteq \Rd$, the set of weak solutions of the equation
\begin{equation*} \label{}
-\nabla \cdot \left( \a\nabla u \right)  = 0 \quad \mbox{in} \ U
\end{equation*}
is denoted by 
\begin{equation} \label{e.def.A(U)}
\A(U) := \left\{ u\in H^1_{\mathrm{loc}}(U)\,:\,\a\nabla u \in L^2_{\mathrm{sol},\,\mathrm{loc}}(U) \right\}.
\end{equation}
We sometimes display the spatial dependence by writing
\begin{equation*} \label{}
-\nabla \cdot \left( \a(x)\nabla u \right)  = 0.
\end{equation*}
This is to distinguish it from the equation we obtain by introducing a small parameter~$\ep>0$ and rescaling, which we also consider, namely
\begin{equation*} \label{}
-\nabla \cdot \left( \a\left(\tfrac x\ep \right)\nabla u \right)  = 0.
\end{equation*}
Notice that these are \emph{random} equations and~$\mathcal{A}(U)$ is a \emph{random} linear vector space since of course it depends on the coefficient field~$\a$. Throughout the book,~$\ahom$ is the homogenized matrix defined in Definition~\ref{d.ahom}. We denote the set of solutions of the homogenized equation in~$U$ by 
\begin{equation*} \label{}
\Ahom(U) := \left\{ u\in H^1_{\mathrm{loc}}(U)\,:\, \ahom\nabla u \in L^2_{\mathrm{sol},\,\mathrm{loc}}(U)  \right\}.
\end{equation*}
For each $k\in\N$, we denote by $\Ahom_k$ and $\A_k$ the subspaces of $\Ahom(\Rd)$ and $\A(\Rd)$, respectively, which grow like $o\left(|x|^{k+1}\right)$ as measured in the~$L^2$ norm:
\begin{equation} 
\label{e.Ahom000.predef}
\Ahom_k := \left\{ u \in \Ahom(\Rd) \, : \,  \limsup_{r \to \infty} r^{-(k+1)} \left\|  u \right\|_{\underline{L}^2(B_r)} = 0  \right\}
\end{equation}
and
\begin{equation} 
\label{e.A000.predef}
\A_k := \left\{ u \in \A(\R^d) \, : \, \limsup_{r \to \infty} r^{-(k+1)} \left\|  u \right\|_{\underline{L}^2(B_r)} = 0  \right\}.
\end{equation}

\subsection*{The $\O_s(\cdot)$ notation}

Throughout the book, we use the following notation to control the size of random variables: for every random variable $X$ and $s, \theta \in (0,\infty)$, we write
\begin{equation}
\label{e.Os.not}
X \leq \O_s(\theta)
\end{equation}
to denote the statement 
\begin{equation}
\label{e.Os2.not}
\E \Ll[ \exp\left(\left(\theta^{-1} X_+\right)^s\right) \Rr] \le 2.
\end{equation}
We likewise write 
\begin{equation*}
X\leq Y +  \O_s(\theta)  \quad \iff \quad X-Y \leq \O_s(\theta)
\end{equation*}
and
\begin{equation*}
X = Y + \O_s(\theta) \quad\iff \quad X-Y \leq \O_s(\theta) \ \ \mbox{and} \ \ Y-X \leq \O_s(\theta). 
\end{equation*}
This notation allows us to write bounds for random variables conveniently, in the spirit of the standard ``big-$O$'' notation it evokes, enabling us to compress many computations that would otherwise take many lines to write. The basic properties of this notation are collected in Appendix~\ref{a.bigO}. 

\subsection*{Convention for the constants $C$}

Throughout, the symbols~$c$ and~$C$ denote positive constants which may vary from line to line, or even between multiple occurrences in the same line, provided the dependence of these constants is clear from the context. Usually, we use~$C$ for large constants (those we expect to belong to~$[1,\infty)$) and~$c$ for small constants (those we expect to belong to~$(0,1]$). When we wish to explicitly declare the dependence of $C$ and $c$ on the various parameters, we do so using parentheses~$(\cdots)$. For example, the phrase \emph{``there exists~$C(p,d,\Lambda)<\infty$ such that...''} is short for \emph{``there exists a constant~$C\in [1,\infty)$, depending on~$p$,~$d$ and~$\Lambda$, such that...''} We use the same convention for other symbols when it does not cause confusion, for instance the small positive exponents~$\alpha$ appearing in Chapter~\ref{c.two}.

\mainmatter



\chapter{Introduction and qualitative theory}
\label{c.one}

We start this chapter by giving an informal introduction to the problem studied in this book. A key role is played by a \emph{subadditive quantity}, denoted by $\nu(U,p)$, which is the energy per unit volume of the solution of the Dirichlet problem in a bounded Lipschitz domain $U\subseteq\Rd$ with affine boundary data with slope~$p\in\Rd$. We give a simple argument to show that $\nu(U,p)$ converges in $L^1(\Omega,\P)$ to a deterministic limit if $U$ is a cube with side length growing to infinity. Then, using entirely deterministic arguments, we show that this limit contains enough information to give a qualitative homogenization result for quite general Dirichlet boundary value problems. Along the way, we build some intuition for the problem and have a first encounter with some ideas playing a central role in the rest of the book.

\section{A brief and informal introduction to homogenization}
\label{s.intro}

The main focus of this book is the study of the elliptic equation
\begin{equation}  
\label{e.main.eq}
-\nabla \cdot \Ll( \a(x) \nabla u \Rr)  = f \qquad \text{in } U,
\end{equation}
where $U \subset \Rd$ is a bounded Lipschitz domain.
We assume throughout that $x \mapsto \a(x)$ is a random vector field taking values in the space of symmetric matrices and satisfying the assumptions \eqref{e.ue} to \eqref{e.assumption.independence}, namely uniform ellipticity, stationarity with respect to $\Z^d$-translations, and unit range of dependence. Our goal is to describe the large-scale behavior of solutions of \eqref{e.main.eq}. 

\smallskip

In dimension $d = 1$, the solution $u_\eps$ of 
\begin{equation*}  
\Ll\{
\begin{aligned}
& -\nabla \cdot \Ll( \a \Ll( \tfrac{x}{\eps} \Rr)  \nabla u_\eps \Rr) = 0 & \ \text{in } (0,1), \\
& u_\eps(0) = 0, & \\
& u_\eps(1) = 1 &
\end{aligned}
\Rr.
\end{equation*}
can be written explicitly as
\begin{equation*}  
u_\eps(x) = m_\eps^{-1} \int_0^x \a^{-1} \Ll( \tfrac y \eps \Rr) \, dy
\qquad \mbox{where} \qquad
m_\eps := \int_0^1 \a^{-1} \Ll( \tfrac y \eps \Rr) \, dy.
\end{equation*}
One can then verify using the ergodic theorem (which in this case can be reduced to the law of large numbers for sums of independent random variables) that, almost surely (henceforth abbreviated ``a.s.''), the solution~$u_\eps$ converges to the linear function $x \mapsto x$ in $L^2(0,1)$ as $\eps\to 0$. Moreover, observing once again by the ergodic theorem, or law of large numbers, that 
\begin{equation*}  
m_\eps \rightarrow m_0 := \E \Ll[ \int_0^1 \a^{-1}(y) \, dy \Rr] 
\quad 
\mbox{a.s. as} \ \ep \to 0, 
\end{equation*}
we also obtain the \emph{weak} convergences in $L^2(0,1)$ of
\begin{equation*}  
\nabla u_\eps \rightharpoonup 1 = \nabla u 
\qquad \text{and} \qquad 
\a \Ll( \tfrac \cdot \eps \Rr) \nabla u_\eps \rightharpoonup m_0^{-1} = m_0^{-1} \nabla u
\quad
\mbox{a.s. as} \ \ep \to 0.
\end{equation*}
Note that we do \emph{not} expect these limits to hold in the (strong) topology of $L^2(0,1)$: see Figure~\ref{f.twoscale}. It is natural to define the homogenized coefficients~$\ahom$ so that the flux $\a \Ll( \tfrac \cdot \eps \Rr) \nabla u_\eps$ of the heterogeneous solution weakly converges to the flux of the homogenized solution. This leads us to define
\begin{equation}  
\label{e.formula.oned}
\ahom := m_0^{-1} = \E \Ll[ \int_0^1 \a^{-1}(y) \, dy \Rr]^{-1}.
\end{equation}
In other words, $\ahom$ is the harmonic mean of $\a$. This definition is also consistent with the convergence of the energy
\begin{equation*}  
\frac 1 2  \int_0^1 \nabla u_\eps(x) \cdot \a \Ll( \tfrac{x}{\eps} \Rr) \nabla u_\eps(x) \, dx \rightarrow \frac 1 2 \int_0^1  \nabla u(x) \cdot \ahom \nabla u(x) \, dx
\quad
\mbox{a.s. as} \ \ep \to 0.
\end{equation*}

\smallskip

In dimensions $d \ge 2$, one still expects that the randomness of the coefficient field averages out, in the sense that there exists a \emph{homogenized matrix} $\ahom$ such that the solution $u_\eps$ of
\begin{equation}  
\label{e.heterogeneous.nobdy}
-\nabla \cdot \Ll( \a \Ll( \tfrac{x}{\eps} \Rr)  \nabla u_\eps \Rr) = 0 \qquad \text{in } U
\end{equation}
with suitable boundary condition converges to the solution $u$ of 
\begin{equation}  
\label{e.homogeneous.nobdy}
-\nabla \cdot \Ll( \ahom  \nabla u \Rr) = 0 \qquad \text{in } U
\end{equation}
with the same boundary condition. However, there will be no longer any simple formula such as \eqref{e.formula.oned} for the homogenized matrix. Indeed, in dimension $d \ge 2$, the flux can circumvent regions of small conductivity which are surrounded by regions of high conductivity, and thus $\ahom$ must incorporate subtle geometric information about the law of the coefficient field. To make this point clear, consider the example displayed on Figure~\ref{f.barriers}.

\begin{figure}[tb]
\centering
\begin{tikzpicture}[scale=1]
\clip (0,0) rectangle (10,10);
\draw[lightgray] (0,0) rectangle (10,10);
\foreach \n in {1,2,...,2000}
\pgfmathsetmacro{\x}{(1+rand)*5.1-0.05}
\pgfmathsetmacro{\y}{(1+rand)*7-2}
\pgfmathsetmacro{\t}{rand*5}
\pgfmathsetmacro{\H}{(rand+5)*0.1}
\pgfmathsetmacro{\P}{rand*0.03}
\pgfmathsetmacro{\Q}{rand*0.015}
\draw[thick,cm={cos(\t) ,-sin(\t) ,sin(\t) ,cos(\t) ,(0 cm,0 cm)}] plot [smooth] coordinates{  (\x,\y)  (\x+\P*\H,\y+\H/3)  (\x+\P*\H/2 +\Q*\H,\y+2*\H/3)  (\x,\y+\H) };
\end{tikzpicture}
\caption{
\small{
If this image is a sample of a composite material which is a good conductor in embedded thin (black) fibers and a good insulator in the (white) ambient material, then the effective conductivity will be larger in the~$e_2$ direction than the~$e_1$ direction. In particular, in this situation we have $\a(x) = a(x)\Id$ at every point~$x$, yet $\overline{\a}$ is not a scalar matrix. This shows us that in~$d\geq 2$ we should not expect to find a simple formula for~$\overline{\a}$ which extends~\eqref{e.formula.oned} for~$d=1$.
}
}
\label{f.barriers}
\end{figure}

\smallskip

\begin{figure}
\centering
\begin{tikzpicture}[xscale=0.08,yscale=0.10,domain=(-87.673):(68.488),samples=800]
\draw[blue,smooth,variable=\x] 
plot ({\x},{8.4-sin(deg(\x/5))+.004*\x*\x+0.4*cos(deg(\x/5))+ 9.9*sin(deg(\x/10)) });
\draw[red,smooth,variable=\x] 
plot ({\x},{ 
9-sin(deg(\x/5))+.004*\x*\x+0.4*cos(deg(\x/5))+ 9.9*sin(deg(\x/10)) 
+1.2*cos(deg((2*\x-20)*cos(deg(3*\x-4))/7+101) )
-sin(deg(\x+47+2*cos(deg(2.2*\x+1)+87))) - 1.75*exp(-(\x*\x/(80*80)))*cos(deg(\x/20))
+0.5*sin(deg(3*\x+2))+ 1.2*sin(deg(2.5*\x))-0.6*cos(deg(2.6*\x))-0.6*(sin(deg(2*\x)))*(sin(deg(2*\x))) +0.005*\x
});
\draw[-] (5,10)--(15,10)--(15,21)--(5,21)--(5,10);
\draw[dashed] (5,21) -- (-49.593,-8);
\draw[dashed] (5,10) -- (-49.593,-50); 
\draw[dashed] (15,10) -- (30.408,-50); 
\draw[dashed] (15,21) -- (30.408,-8); 
\draw[-] (-49.593,-50)-- (30.408,-50)--(30.408,-8)--(-49.593,-8)--(-49.593,-50);
\draw[blue]  (-49.593,-42)--(30.408,-16);
\draw[red,smooth,domain=(-49.593):30.408,variable=\x] 
plot({\x},{
0.325*\x-26 - 1.6*cos(deg(1+\x+sin(deg(\x/11)))) + sin(deg(2+ cos(deg(\x/4)))) -0.7 +\x/45 + 2.9*cos(deg((\x/4+42)*cos(deg(\x-4))/25) )
});
\end{tikzpicture}
\caption{
\small{
The red and blue curves represent the solutions of equations \eqref{e.heterogeneous.nobdy} and \eqref{e.homogeneous.nobdy} respectively; the rectangle at the bottom is a (schematic) close-up of the rectangle on top. On mesoscopic scales, the blue curve is essentially affine, and the red curve is close to the solution to a Dirichlet problem with affine boundary condition (with slope given by the local gradient of the homogeneous solution).
}
}
\label{f.twoscale}
\end{figure}

Since there is no explicit formula for $\ahom$ in dimension $d \ge 2$, we need to identify quantities which will allow us to track the progressive homogenization of the equation \eqref{e.main.eq} as we move to larger and larger scales. Before doing so, we first argue that understanding the homogenization phenomenon for simple domains such as balls or cubes, and with affine boundary condition, should be sufficient; it should be possible to deduce homogenization results for more complicated domains and boundary conditions (and possibly non-zero right-hand side) a posteriori. The idea is that, since the solution of the homogeneous equation is smooth, it will be well-approximated by an affine function on scales smaller than the macroscopic scale. On scales intermediate between the microscopic and macroscopic scales, the behavior of the solution of the equation with rapidly oscillating coefficients should thus already be typical of homogenization, while tracking an essentially affine function. In other words, for $u_\ep$ and $u$ the solutions to \eqref{e.heterogeneous.nobdy} and \eqref{e.homogeneous.nobdy} respectively, with the same boundary condition, we expect that for $z \in U$ and for scales $r$ such that $\eps \ll r \ll 1$, 
\begin{equation}  
\label{e.first.twoscale}
\|\nabla u_\ep - \nabla \td u_{\eps,z,r}\|_{\underline L^2(B_r(z))} \ll 1,
\end{equation}
where $\td u_{\eps,z,r}$ solves 
\begin{equation}  
\label{e.first.corrector}
\Ll\{
\begin{aligned}
& -\nabla \cdot \Ll( \a \Ll( \tfrac{x}{\eps} \Rr)  \nabla \td u_{\eps,z,r} \Rr) = 0 & \ \text{in } B_r(z), \\
& \td u_{\ep,z,r}(x) = \nabla u(z) \cdot x  & \ \text{on } \partial B_r(z). 
\end{aligned}
\Rr.
\end{equation}
See Figure~\ref{f.twoscale} for a cartoon visualization of this idea.

\smallskip

These considerations motivate us to focus on understanding the homogenization of problems such as \eqref{e.first.corrector}, that is, Dirichlet problems on simple domains with affine boundary data. The approach taken up here is inspired by earlier work of Dal Maso and Modica \cite{DM1,DM2}, who introduced, for every $p \in \Rd$ (and for more general nonlinear equations), the quantity
\begin{equation*}
\nu(U,p) := \inf_{v \in \ell_p + H^1_0(U)} \fint_U \frac 1 2 \nabla v \cdot \a \nabla v,
\end{equation*}
where $\ell_p$ denotes the affine function $x \mapsto p\cdot x$. Note that the minimizer $v(\cdot,U,p)$ in the definition of $\nu$ is the solution of the Dirichlet problem 
\begin{equation*}  
\Ll\{
\begin{aligned}
& -\nabla \cdot \Ll( \a   \nabla v(\cdot,U,p) \Rr) = 0 & \ \text{in } U, \\
& v(\cdot,U,p) = \ell_p  & \ \text{on } \partial U,
\end{aligned}
\Rr.
\end{equation*}
which should be compared with~\eqref{e.first.corrector}. The first key observation of Dal Maso and Modica is that $\nu(\cdot,p)$ is \emph{subadditive}\footnote{Note that our use of the term \emph{subadditive} is not standard: it is usually the unnormalized quantity $U\mapsto |U|\nu(U,p)$ which is called subadditive.}: if the domain $U$ is partitioned into subdomains $U_1, \ldots, U_k$, up to a Lebesgue null set, then
\begin{equation*}  
\nu(U,p) \le \sum_{i = 1}^k \frac{|U_i|}{|U|} \nu(U_i,p). 
\end{equation*}
Indeed, the minimizers for each $\nu(U_i,p)$ can be glued together to create a minimizer candidate for the minimization problem in the definition of~$\nu(U,p)$. The true minimizer cannot have more energy, which yields the claimed inequality. By an appropriate version of the ergodic theorem (found for instance in~\cite{AK}), we deduce the convergence 
\begin{equation*}  
\nu((-r,r)^d,p) \xrightarrow[r\to \infty]{\text{a.s.}} \frac 1 2 p \cdot \ahom p.
\end{equation*}
That the limit can be written in the form above is a consequence of the fact that $p \mapsto \nu(U,p)$ is a quadratic form; we take this limit as the definition of the effective matrix $\ahom$. Dal Maso and Modica then observed (even in a more general, nonlinear setting) that this convergence suffices to imply qualitative  homogenization. 

\smallskip

In this chapter, we will carry out the program suggested in the previous paragraph and set the stage for the rest of the book. In particular, in Theorem~\ref{t.DP.blackbox} we will obtain a fairly general (although at this stage, still only qualitative) homogenization result for Dirichlet problems. We begin in the next two sections with a proof of the convergence of the quantity~$\nu(U,p)$.

\section{The subadditive quantity \texorpdfstring{$\nu$}{nu} and its basic properties}

In this section, we review the basic properties of the quantity~$\nu(U,p)$ introduced in the previous section, which is defined for each bounded Lipschitz domain~$U\subseteq \Rd$ and~$p \in \Rd$ by
\index{subadditive quantity!$\nu$}
\begin{equation}
\label{e.def.nu}
\nu(U,p) 
:= \inf_{v \in \ell_p + H^1_0(U)}  \fint_U \frac 1 2 \nabla v \cdot \a \nabla  v
= \inf_{w \in H^1_0(U)}  \fint_U \frac 1 2\left( p + \nabla w \right)\cdot \a\left( p + \nabla w \right).
\end{equation}
Recall that $\ell_p(x):= p\cdot x$ is the affine function of slope~$p$. We denote the (unique) minimizer of the optimization problem in the definition of $\nu(U,p)$ by
\begin{equation} 
\label{e.vdeff}
v(\cdot,U,p) := \mbox{unique $v\in \ell_p+H^1_0(U)$ minimizing} \  \fint_U \frac 1 2\nabla v \cdot \a \nabla  v.
\end{equation}
The uniqueness of the minimizer is immediate from the uniform convexity of the integral functional, which is recalled in Step~1 of the proof of Lemma~\ref{l.basicnu} below. The existence of a minimizer follows from the weak lower semicontinuity of the integral functional (cf.~\cite[Chapter 8]{Evans}), a standard fact from the calculus of variations, the proof of which we do not give here.

\smallskip

The quantity $\nu(U,p)$ is the energy (per unit volume) of its minimizer $v(\cdot,U,p)$ which, as we will see below, is the solution of the Dirichlet problem
\begin{equation*} \label{}
\left\{
\begin{aligned}
& -\nabla \cdot \left( \a(x) \nabla u \right) = 0 & \mbox{in} & \ U, \\
& u = \ell_p & \mbox{on} & \ \partial U. 
\end{aligned}
\right.
\end{equation*}
We next explore some basic properties of $\nu(U,p)$.

\begin{lemma}[Basic properties of $\nu$]
\label{l.basicnu}
Fix a bounded Lipschitz domain $U\subseteq \Rd$. 
The quantity $\nu(U,p)$ and its minimizer $v(\cdot,U,p)$ satisfy the following properties:

\begin{itemize}

\item \emph{Representation as quadratic form.} The mapping $p\mapsto \nu(U,p)$ is a positive quadratic form, that is, there exists a symmetric matrix $\a(U)$ such that 
\begin{equation} \label{e.numatrixbounds}
\Id \leq \a(U) \leq \Lambda \Id
\end{equation}
and 
\begin{equation} \label{e.nurepresentation}
\nu(U,p) = \frac12 p\cdot \a(U) p. 
\end{equation}

\item \emph{Subadditivity.} Let $U_1, \ldots, U_N \subset U$ be bounded Lipschitz domains that form a partition of~$U$, in the sense that $U_i \cap U_j = \emptyset$ if $i\neq j$ and 
\begin{equation*} \label{}
\left| U \setminus \bigcup_{i=1}^N U_i \right| = 0. 
\end{equation*}
Then, for every $p\in\Rd$,
\begin{equation}
\label{e.nusubadd}
\nu(U,p) \leq \sum_{i=1}^N \frac{\left|U_i\right|}{|U|} \nu(U_i,p). 
\end{equation}

\item \emph{First variation.} For each $p\in\Rd$, the function $v(\cdot,U,p)$ is characterized as the unique solution of the Dirichlet problem
\begin{equation} 
\label{e.DP}
\left\{
\begin{aligned}
& -\nabla \cdot \left( \a \nabla v \right) = 0 & \mbox{in} & \ U, \\
& v = \ell_p & \mbox{on} & \ \partial U. 
\end{aligned}
\right. 
\end{equation}
The precise interpretation of~\eqref{e.DP} is
\begin{equation*} \label{}
\mbox{$v$ solves~\eqref{e.DP}} 
\iff
v\in \ell_p+H^1_0(U) \ \mbox{and}  \ \forall w\in H^1_0(U), \ 
\fint_U \nabla w \cdot \a\nabla v = 0. 
\end{equation*}

\item \emph{Quadratic response.} For every $w\in \ell_p + H^1_0(U)$, 
\begin{equation} 
\label{e.nuquadresponse}
\frac12 \fint_U \left| \nabla w - \nabla v(\cdot,U,p) \right|^2
  \leq
\fint_U \frac12\nabla w \cdot \a\nabla w - \nu(U,p)
\leq
\frac{\Lambda}2 \fint_U \left| \nabla w - \nabla v(\cdot,U,p) \right|^2.
\end{equation}

\end{itemize}
\end{lemma}
\begin{proof}

\emph{Step 1.} We first derive the first variation of the minimization problem in the definition of $\nu$. We write $v:= v(\cdot,U,p)$ for short. Fix $w\in H^1_0(U)$, $t\in\R$ and compare the energy of $\tilde{v}_t:= v + tw$ to the energy of $v$:
\begin{align*} \label{}
\fint_U \frac12 \nabla v \cdot \a\nabla v 
= \nu(U,p) 
 & 
\leq \fint_U \frac12 \nabla \tilde{v}_t \cdot \a\nabla  \tilde{v}_t 
\\ & 
=  \fint_U \frac12 \nabla v \cdot \a\nabla v  + 2t \fint_U \frac12 \nabla w \cdot \a\nabla v + t^2 \fint_U \frac12 \nabla w \cdot \a \nabla w. 
\end{align*}
Rearranging this and dividing by $t$, we get
\begin{equation*} \label{}
 \fint_U  \nabla w \cdot \a\nabla v 
 \geq 
 - \frac12 t \fint_U  \nabla w \cdot \a \nabla w.
\end{equation*}
Sending $t\to 0$ gives 
\begin{equation*} \label{}
 \fint_U  \nabla w \cdot \a\nabla v \geq 0.  
\end{equation*}
Applying this inequality with both $w$ and $-w$, we deduce that, for every $w\in H^1_0(U)$,  
\begin{equation*} \label{}
 \fint_U  \nabla w \cdot \a\nabla v = 0. 
\end{equation*}
This confirms that $v$ is a solution of~\eqref{e.DP}. To see it is the unique solution, we assume that $\hat v$ is another solution and test the equation for $v$ and for $\hat v$ with $v-\hat v$ and subtract the results to obtain
\begin{equation*} \label{}
\fint_U \frac12\left| \nabla v - \nabla \hat v\right|^2 \leq \fint_U \frac12 \left( \nabla v - \nabla \hat v\right) \cdot \a\left( \nabla v - \nabla \hat v\right) = 0. 
\end{equation*}

\smallskip

\emph{Step 2.} We show that 
\begin{equation} 
\label{e.nubounds}
\frac12|p|^2 \leq \nu(U,p) \leq \frac{\Lambda}2|p|^2. 
\end{equation}
The upper bound is immediate from testing the definition of $\nu(U,p)$ with $\ell_p$:
\begin{equation*} \label{}
\nu(U,p) 
\leq 
\fint_{U} \frac12 \nabla \ell_p \cdot \a \nabla \ell_p
 = \fint_U \frac12 p\cdot \a p
 \leq \frac{\Lambda}2|p|^2. 
\end{equation*}
The lower bound comes from Jensen's inequality: for every $w\in H^1_0(U)$, 
\begin{equation*} \label{}
\fint_U \frac12(p+\nabla w) \cdot \a (p+\nabla w) \geq \fint_U \frac12 \left|p+\nabla w \right|^2 \geq   \frac12 \left|p+\fint_U \nabla w \right|^2 = \frac12|p|^2. 
\end{equation*}
Taking the infimum over $w\in H^1_0(U)$ yields the lower bound of~\eqref{e.nubounds}. 

\smallskip

\emph{Step 3.} We show that $\nu(U,\cdot)$ is a positive quadratic form as in~\eqref{e.nurepresentation} satisfying bounds in~\eqref{e.numatrixbounds}. Observe first that a consequence of the characterization~\eqref{e.DP} of the minimizer $v(\cdot,U,p)$ is that 
\begin{equation} 
\label{e.minimizerlinearity}
p\mapsto v(\cdot,U,p) \quad \mbox{is a linear map from $\Rd$ to $H^1(U)$.}
\end{equation}
Indeed, the formulation~\eqref{e.DP} makes linearity immediate. 
Moreover, since 
\begin{equation} 
\label{e.obviousformula}
\nu(U,p) = \fint_U \frac 1 2 \nabla v(\cdot,U,p) \cdot \a\nabla v(\cdot,U,p)
\end{equation}
we deduce that 
\begin{equation} 
\label{e.nuquadratic}
p\mapsto \nu(U,p) \quad \mbox{is a quadratic form.} 
\end{equation}
That is, there exists a symmetric matrix $\a(U) \in \R^{d\times d}$ as in~\eqref{e.nurepresentation}. 
The inequalities in~\eqref{e.nubounds} can thus be rewritten as
\begin{equation*} \label{}
\frac12 |p|^2 \leq \frac12 p\cdot \a(U) p \leq \frac{\Lambda}2 |p|^2,
\end{equation*}
which gives~\eqref{e.numatrixbounds}.  

\smallskip

\emph{Step 4.} We next prove~\eqref{e.nuquadresponse}, the quadratic response of the energy around the  minimizer. This is an easy consequence of the first variation: in fact, we essentially proved it already in Step~1. 

\smallskip

We fix $w\in \ell_p+H^1_0(U)$ and compute
\begin{align*} \label{}
\fint_U \frac12\nabla w \cdot \a\nabla w - \nu(U,p)
& 
= \fint_U \frac12\nabla w \cdot \a\nabla w - \fint_U \frac12\nabla v(\cdot,U,p) \cdot \a\nabla v(\cdot,U,p) 
\\ & 
=  \fint_U \frac12\left( \nabla w -\nabla v(\cdot,U,p) \right) \cdot \a\left( \nabla w -\nabla v(\cdot,U,p) \right) 
\\ & \qquad 
+ \fint_U \left( \nabla w -\nabla v(\cdot,U,p) \right)\cdot \a \nabla v(\cdot,U,p).
\end{align*}
Noting that $w - v \in H^1_0(U)$, we see that the last term on the right side of the previous display is zero, by the first variation. Thus 
\begin{equation*} \label{}
\fint_U \frac12\nabla w \cdot \a\nabla w - \nu(U,p) 
= \fint_U \frac12\left( \nabla w -\nabla v(\cdot,U,p) \right) \cdot \a\left( \nabla w -\nabla v(\cdot,U,p) \right),
\end{equation*}
which is a more precise version of~\eqref{e.nuquadresponse}. 

\smallskip

\emph{Step 5.} The proof of subadditivity. We glue together the minimizers $v(\cdot,U_i,p)$ of~$\nu$ in the subdomains~$U_i$ and compare the energy of the result to~$v(\cdot,U,p)$. We first need to argue that the function defined by 
\begin{equation} 
\label{e.gluingfun}
\tilde{v}(x):= v(x,U_i,p), \quad x \in U_i,
\end{equation}
belongs to $\ell_p+H^1_0(U)$. To see this, we observe that each $v(\cdot,U_i,p)$ can be approximated in $H^1$ by the sum of $\ell_p$ and a $C^\infty$ function with compact support in $U_i$. We can glue these functions together to get a smooth function in $\ell_p+H^1_0(U)$ which clearly approximates $\tilde{v}$ in $H^1(U)$. Therefore $\tilde{v}\in \ell_p+H^1_0(U)$. This allows us to test the definition of $\nu(U,p)$ with $\tilde{v}$ and yields
\begin{align*} \label{}
\nu(U,p) 
\leq  \fint_U \frac12 \nabla \tilde{v} \cdot \a\nabla \tilde{v}
 & 
= \frac{1}{|U|} \sum_{i=1}^N \int_{U_i} \frac12 \nabla \tilde{v} \cdot \a\nabla \tilde{v}
\\ & 
= \frac{1}{|U|} \sum_{i=1}^N \left| U_i \right| \fint_{U_i} \frac12 \nabla v(\cdot,U_i,p) \cdot \a\nabla v(\cdot,U_i,p) 
\\ & 
= \sum_{i=1}^N \frac{\left| U_i \right|}{|U|} \nu(U_i,p). 
\end{align*}
This completes the proof of~\eqref{e.nusubadd} and therefore of the lemma. 

\smallskip

For future reference, however, let us continue by recording the slightly more precise estimate that the argument for~\eqref{e.nusubadd} gives us. The above computation can be rewritten as
\begin{equation*} 
\sum_{i=1}^N \frac{\left| U_i \right|}{|U|} \nu(U_i,p) - \nu(U,p)
=  \fint_U \frac12 \nabla \tilde{v} \cdot \a\nabla \tilde{v} - \nu(U,p).
\end{equation*}
Quadratic response thus implies 
\begin{equation}
\label{e.gluingerror0}
\frac12 \fint_U \left| \nabla v(\cdot,U,p) - \nabla \tilde{v} \right|^2
\leq
\sum_{i=1}^N \frac{\left| U_i \right|}{|U|} \nu(U_i,p) - \nu(U,p)
\leq \frac{\Lambda}2 \fint_U \left| \nabla v(\cdot,U,p) - \nabla \tilde{v} \right|^2.
\end{equation}
This can be written as 
\begin{multline}
\label{e.gluingerror}
\frac 1 2 \sum_{i=1}^N \frac{\left| U_i \right|}{|U|} \fint_{U_i} \left| \nabla v(\cdot,U,p) - \nabla v(\cdot,U_i,p) \right|^2
\\
\leq
\sum_{i=1}^N \frac{\left| U_i \right|}{|U|} \nu(U_i,p) - \nu(U,p)
\leq \frac \Lambda 2\sum_{i=1}^N \frac{\left| U_i \right|}{|U|} \fint_{U_i} \left| \nabla v(\cdot,U,p) - \nabla v(\cdot,U_i,p) \right|^2.
\end{multline}
In other words, the strictness of the subadditivity inequality is proportional to the weighted average of the $L^2$ differences between $\nabla v(\cdot,U,p)$ and $\nabla v(\cdot,U_i,p)$ in the subdomains $U_i$.
\end{proof}

\section{Convergence of the subadditive quantity}

In order to study the convergence of $\nu(U,p)$ as the domain $U$ becomes large, it is convenient to work with the family of triadic cubes $\left\{ x+ \cu_n\,:\, n\in\N, \,x\in \Zd \right\}$ defined in~\eqref{e.triad.def}. Recall that for each $n\in\N$ with $n\leq m$, up to a set of zero Lebesgue measure, the cube $\cu_m$ can be partitioned into exactly $3^{d(m-n)}$ subcubes which are $\Zd$-translations of $\cu_n$, namely $\left\{ z+ \cu_n\,:\, z\in 3^n\Zd\cap \cu_m \right\}$. 

\smallskip

An immediate consequence of subadditivity and stationarity is the monotonicity of $\E \left[ \nu(\cu_m,p)\right]$: for every $m\in\N$ and $p\in\Rd$, 
\begin{equation}
\label{e.nuEmonotone}
\E \left[ \nu(\cu_{m+1},p)\right] \leq \E \left[ \nu(\cu_{m},p)\right]. 
\end{equation}
To see this, we first apply the subadditivity property with respect to the partition $\{ z+ \cu_m \,:\, z\in \{ -3^m,0,3^m\}^d \}$ of $\cu_{m+1}$ into its $3^d$ largest triadic subcubes, to get
\begin{align*}
\nu(\cu_{m+1},p) 
\leq 
\sum_{z\in \{ -3^m,0,3^m \}^d} \frac{\left|z+\cu_m\right|}{\left|\cu_{m+1}\right|} \nu(z+\cu_m,p)
= 3^{-d} \sum_{z\in \{ -3^m,0,3^m \}^d} \nu(z+\cu_m,p).
\end{align*}
Stationarity tells us that, for every $z\in\Zd$, the law of $\nu(z+\cu_m,p)$ is the same as the law of $\nu(\cu_m,p)$. Thus they have the same expectation, and so taking the expectation of the previous display gives 
\begin{equation*} \label{}
\E \left[ \nu(\cu_{m+1},p)  \right]
\leq
3^{-d} \sum_{z\in \{ -3^m,0,3^m \}^d}\E\left[ \nu(z+\cu_m,p) \right]
= \E \left[ \nu(\cu_m,p) \right],
\end{equation*}
which is~\eqref{e.nuEmonotone}.

\smallskip

Therefore, for each $p\in\Rd$, the sequence $\left\{ \E \left[ \nu(\cu_m,p) \right] \right\}_{m\in\N}$ is bounded by~\eqref{e.nubounds} and nonincreasing by~\eqref{e.nuEmonotone}. It therefore has a limit, which we denote by
\begin{equation}
\label{e.defnubar}
\overline{\nu}(p):= \lim_{m\to \infty} \E \left[ \nu(\cu_m,p) \right] = \inf_{m\in\N} \E \left[ \nu(\cu_m,p) \right].
\end{equation}
In Lemma~\ref{l.basicnu} we found (cf.~\eqref{e.nurepresentation}) that
\begin{equation} 
\label{e.nuquadratic2}
p\mapsto \nu(U,p) \quad \mbox{is quadratic.} 
\end{equation}
It follows that $p\mapsto \E\left[ \nu(U,p) \right]$ is also quadratic, and hence
\begin{equation*} \label{}
p\mapsto \overline\nu(p) \quad \mbox{is quadratic.} 
\end{equation*}
It is clear in view of~\eqref{e.numatrixbounds},~\eqref{e.nurepresentation} and~\eqref{e.defnubar} that we have
\begin{equation} 
\label{e.nubarbounds}
\frac12|p|^2 \leq \overline{\nu}(p) \leq \frac\Lambda2|p|^2. 
\end{equation}
The deterministic object $\overline\nu$ allows us to identify the homogenized coefficients and motivates the following definition. 

\begin{definition}[{Homogenized coefficients~$\ahom$}]
\index{homogenized coefficients}
\label{d.ahom}
We denote by $\ahom\in \R^{d\times d}$ the unique symmetric matrix satisfying
\begin{equation*} \label{}
\forall p\in\Rd,  \quad \overline{\nu}(p) = \frac12 p\cdot \ahom p. 
\end{equation*}
We call $\ahom$ the \emph{homogenized coefficients}. By~\eqref{e.nubarbounds}, we see that $\ahom$ is a positive definite matrix and satisfies the bounds
\begin{equation}
\label{e.ahombounds}
\Id \leq \ahom \leq \Lambda \Id. 
\end{equation}
\end{definition}

\begin{exercise}
\label{ex.isotropyforahom}
Show that if the coefficient field~$\a$ is \emph{isotropic in law} in the sense that~$\P$ is invariant under any linear isometry which maps the union of the coordinate axes to itself, then $\ahom$ is a multiple of the identity matrix.
\end{exercise}

\begin{exercise}
\label{ex.voigt-reiss1}
The \emph{Voigt-Reiss bounds}\index{Voigt-Reiss bounds} for the effective coefficients assert that 
\begin{equation} 
\label{e.voigt-reiss}
\E \left[  \int_{\cu_0} \a^{-1} (x) \,dx\right]^{-1}  
\leq \ahom 
\leq \E \left[ \int_{\cu_0} \a(x)\,dx \right].  
\end{equation}
Show that the second inequality of~\eqref{e.voigt-reiss} follows from our definitions of~$\nu$ and~$\ahom$, stationarity and the subadditivity of~$\nu$. (See Exercise~\ref{ex.voigt-reiss2} for the other inequality.)
\end{exercise}

\begin{exercise}
\label{ex.CS}
Assume that, for some $\rho \in (0,1]$,
\begin{equation*} \label{}
\left\| \a - \Id \right\|_{L^\infty(\Rd)} \leq \rho, \quad \mbox{$\P$--a.s.}
\end{equation*}
Using the inequalities of~\eqref{e.voigt-reiss}, show that, for a constant $C(d,\Lambda)<\infty$, 
\begin{equation*} \label{}
\left| \ahom - \E \left[ \int_{\cu_0} \a(x)\,dx \right] \right| \leq C \rho^2. 
\end{equation*}
In other words, the homogenized matrix coincides with the average of the coefficients at first order in the regime of small ellipticity contrast.
\end{exercise}

\smallskip

We show in the next proposition that, using the independence assumption, we can upgrade the limit~\eqref{e.defnubar} from convergence of the expectations to convergence in $L^1(\Omega,\P)$. That is, we prove that, for each $p\in\Rd$,
\begin{equation*} \label{e.nuconv}
\E \left[ \left| \nu(\cu_n,p) - \overline{\nu}(p) \right| \right] \to 0 \quad \mbox{as} \ n \to \infty. 
\end{equation*}
In the process, we will try to extract as much quantitative information about the rate of this limit as we are able to at this stage. For this purpose, we introduce the modulus $\omega$ which governs the rate of the limit~\eqref{e.defnubar}, uniformly in $p\in B_1$:
\begin{equation}
 \label{e.defomega}
\omega(n):= \sup_{p\in B_1} \left(  \E \left[ \nu(\cu_n,p) \right] - \overline{\nu}(p) \right), \quad n\in\N.
\end{equation} 
Since we have taken a supremum over $p\in B_1$, to ensure that $\omega(n) \to 0$ we need an argument. For this it suffices to note that, since $p\mapsto \E \left[ \nu(\cu_n,p) \right] - \overline{\nu}(p)$ is a quadratic form with corresponding matrix $ \E \left[ \a(\cu_n) \right] - \ahom$ which is nonnegative by~\eqref{e.defnubar}, there are constants $C(d)<\infty$ such that 
\begin{equation*} \label{}
\omega(n) 
\leq
C \left| \E \left[ \a(\cu_n) \right] - \ahom \right| 
\leq
C \sum_{i=1}^d \left|  \E \left[ \nu(\cu_n,e_i) \right] - \overline{\nu}(e_i) \right| 
\to 0 \quad \mbox{as} \ n\to \infty. 
\end{equation*}
Notice also that $\omega(n)$ is nonincreasing in $n$. 


\smallskip

In order to apply the independence assumption, we  require the observation that 
\begin{equation}
\label{e.locality}
\nu(U,p) \quad \mbox{is $\F(U)$--measurable.}
\end{equation}
This is immediate from the definition of $\nu(U,p)$ since the latter depends only on $p$ and the restriction of the coefficient field~$\a$ to~$U$. That is, $\nu(U,p)$ is a \emph{local quantity.} 
 
 \smallskip

\begin{proposition}
\label{p.nuconv}
There exists $C(d,\Lambda)<\infty$ such that, for every $m\in\N$,
\begin{equation} 
\label{e.preL1conv}
\E \left[ \left| \a(\cu_m) - \ahom  \right| \right] \leq C 3^{-\frac d4m} + C \omega\left( \left\lceil  \tfrac m2 \right\rceil \right).
\end{equation}
\end{proposition}
\begin{proof}
This argument is a simpler variant of the one in the proof of the subadditive ergodic theorem (we are thinking of the version proved in~\cite{AK}). Compared to the latter, the assumptions here are stronger (finite range of dependence instead of a more abstract ergodicity assumption) and we prove less (convergence in $L^1$ versus almost sure convergence). The idea is relatively simple: if we wait until the expectations have almost converged, then the subadditivity inequality will be almost sharp (at least in expectation). That is, we will have almost additivity in expectation. On the other hand, for sums of independent random variables, it is of course very easy to show improvement in the scaling of the variance. 

\smallskip

To begin the argument, we fix $p\in B_1$, $m\in\N$ and choose a mesoscopic scale given by $n\in\N$ with $n<m$. 

\smallskip

\emph{Step 1.} 
We show that there exists $C(d,\Lambda)<\infty$ such that 
\begin{equation} 
\label{e.indycontr}
\E\left[\left( \nu(\cu_m,p) - \E \left[ {\nu}(\cu_n,p)\right] \right)_+^2\right]
\leq C 3^{-d(m-n)},
\end{equation}
where we write $x_+ = \max(x,0)$. 
Using the subadditivity of $\nu$ with respect to the partition $\left\{ z+\cu_n\,:\, z\in 3^n\Zd\cap \cu_m \right\}$ of $\cu_m$ into triadic subcubes of size $3^n$, we get
\begin{equation*} \label{}
\nu(\cu_m,p) \leq 3^{-d(m-n)} \sum_{z\in 3^n\Zd\cap \cu_m } \nu(z+\cu_n,p). 
\end{equation*}
Thus
\begin{align*} \label{}
\left( \nu(\cu_m,p) - \E \left[ {\nu}(\cu_n,p)\right] \right)_+^2
&
\leq  \left( 3^{-d(m-n)} \sum_{z\in 3^n\Zd\cap \cu_m } \nu(z+\cu_n,p) -\E \left[ {\nu}(\cu_n,p)\right] \right)_+^2.
\end{align*}
By stationarity, we have that 
\begin{align*} \label{}
\lefteqn{
\E\left[ \left( 3^{-d(m-n)} \sum_{z\in 3^n\Zd\cap \cu_m } \nu(z+\cu_n,p) - \E \left[ \nu(\cu_n,p) \right] \right)_+^2 \right]
}
\quad & 
\\ & 
=
\E\left[ \left( 3^{-d(m-n)} \sum_{z\in 3^n\Zd\cap \cu_m } \nu(z+\cu_n,p) - \E \left[ 3^{-d(m-n)} \sum_{z\in 3^n\Zd\cap \cu_m } \nu(z+\cu_n,p) \right] \right)_+^2\right]
\\ &
\leq
\var\left[ 3^{-d(m-n)} \sum_{z\in 3^n\Zd\cap \cu_m } \nu(z+\cu_n,p) \right]
\\ & 
= 3^{-2d(m-n)} \sum_{z,z'\in 3^n\Zd\cap \cu_m } \cov\left[\nu(z+\cu_n,p),\,\nu(z'+\cu_n,p)  \right]. 
\end{align*}
By the unit range of dependence assumption and~\eqref{e.locality}, we have that 
\begin{equation*} \label{}
\dist(z+\cu_n,z'+\cu_n)\geq 1\implies  \cov\left[\nu(z+\cu_n,p),\nu(z'+\cu_n,p)  \right]=0. 
\end{equation*}
Each subcube $z+\cu_n$ has at most $3^d-1$ neighboring subcubes, those which satisfy $\dist(z+\cu_n,z'+\cu_n)< 1$. There are $3^{d(m-n)}$ subcubes in $\cu_m$, which means that there are at most $C3^{d(m-n)}$ pairs of neighboring subcubes. For neighboring subcubes, we give up and use H\"older's inequality to estimate the covariance, which gives
\begin{align*} \label{}
\cov\left[\nu(z+\cu_n,p),\nu(z'+\cu_n,p)  \right] 
& 
\leq \left( \var\left[\nu(z+\cu_n,p)\right] \cdot \var\left[\nu(z'+\cu_n,p)\right] \right)^{\frac12} 
\\ & 
= \var\left[\nu(\cu_n,p)\right] 
\leq C,
\end{align*}
where in the last line we used~\eqref{e.nubounds}. Putting all this together, we obtain
\begin{equation*} \label{}
\sum_{z,z'\in 3^n\Zd\cap \cu_m } \cov\left[\nu(z+\cu_n,p),\nu(z'+\cu_n,p)  \right] \leq C 3^{d(m-n)}. 
\end{equation*}
Combining this result with the previous displays above, we get~\eqref{e.indycontr}.

\smallskip

\emph{Step 2.} We now use~\eqref{e.indycontr} to obtain convergence in $L^1(\Omega,\P)$. Observe that, by the triangle inequality, the fact that $|r| = 2 r_+ - r$ for any $r \in \R$, the monotonicity~\eqref{e.nuEmonotone} of~$\E\left[ \nu(\cu_m,p)\right]$ and H\"older's inequality, we get, for every $m,n\in\N$ with $n < m$, 
\begin{align*}
\E \left[ \left| \nu(\cu_m,p) - \overline{\nu}(p) \right| \right] 
& 
\leq
\E \left[ \left| \nu(\cu_m,p) - \E\left[ \nu(\cu_m,p)\right] \right| \right] 
+ \omega(m) 
\\ & 
= 2\E \left[ \left( \nu(\cu_m,p) - \E\left[ \nu(\cu_m,p)\right] \right)_+ \right] + \omega(m)
\\ & 
\leq 
2\E \left[ \left( \nu(\cu_m,p) - \E\left[ \nu(\cu_n,p)\right] \right)_+ \right] + \omega(m)+ 2\omega(n)
\\ &
\leq 2\E \left[ \left( \nu(\cu_m,p) -  \E\left[ \nu(\cu_n,p)\right]  \right)_+^2 \right]^{\frac12} +  \omega(m) + 2\omega(n)
\\ & 
\leq C3^{-\frac{d}{2}(m-n)} + 3\omega(n).
\end{align*}
The crude choice of the mesoscale $n:= \left\lceil \frac m2 \right\rceil$ gives us 
\begin{equation}
\label{e.cruderate}
\E \left[ \left| \nu(\cu_m,p) - \overline{\nu}(p) \right| \right] \leq C 3^{-\frac d4m} + C \omega\left( \left\lceil  \tfrac m2 \right\rceil \right).
\end{equation}
In view of~\eqref{e.nusubadd}, taking the supremum over~$p\in B_1$ yields the proposition. 
\end{proof}

The previous argument gives more information than the limit $\nu(\cu_n,p) \to \overline \nu(p)$ in $L^1(\Omega,\P)$. Namely, it provides an \emph{explicit, quantitative convergence rate} for the limit, up to the knowledge of the speed of convergence of the expectations in~\eqref{e.defnubar}. This motivates us to estimate the modulus~$\omega(n)$. Unfortunately, the qualitative argument for the limit~\eqref{e.defnubar}, which was a one-line soft argument based on the monotonicity of $\E\left[ \nu(\cu_n,p) \right]$ in $n$, does not tell us how to obtain a quantitative rate of convergence. The task of estimating $\omega$ turns out to be rather more subtle and will be undertaken in Chapter~\ref{c.two}.

\smallskip

We next demonstrate the convergence of the minimizers $v(\cdot,\cu_m,p)$ to the affine function $\ell_p$ in $L^2(\cu_m)$. For qualitative convergence, what we should expect is that the $\underline{L}^2(\cu_m)$ norm of the difference~$v(\cdot,\cu_m,p)-\ell_p$ is much smaller than the $\underline{L}^2(\cu_m)$ norm of $\ell_p$ itself, which is $\asymp 3^m |p|$. In other words, we should show that, in some appropriate sense, 
\begin{equation*} \label{}
3^{-m}  \left\| v(\cdot,\cu_m,p) - \ell_p \right\|_{\underline{L}^2(\cu_m)} \to 0 \quad \mbox{as} \ m\to \infty. 
\end{equation*}
For now, we will prove this convergence in $L^1(\Omega,\P)$ with an explicit rate depending only on the modulus~$\omega$. 

\begin{proposition}
\label{p.vellpconv}
There exists $C(d,\Lambda)>0$ such that, for every $m\in\N$ and $p\in B_1$, 
\begin{equation} 
\label{e.vellpconv}
\E \left[ 3^{-2m}  \left\| v(\cdot,\cu_m,p) - \ell_p \right\|_{\underline{L}^2(\cu_m)}^2
\right]
\leq 
C 3^{-\frac m4  }+ C \omega\left( \left\lceil \tfrac m2 \right\rceil \right). 
\end{equation}
\end{proposition}
\begin{proof}
Let us first summarize the rough idea underlying the argument. By quadratic response, the expected squared $L^2$ difference of the gradients of the minimizer $v(\cdot,\cu_m,p)$ and the function obtained by gluing the minimizers $v(\cdot,z+\cu_n,p)$ for $z\in 3^n\Zd\cap\cu_m$ is controlled by the difference in their expected energies. We have encountered this fact already in~\eqref{e.gluingerror0}. With the help of the Poincar\'e inequality, this tells us that the $L^2$ difference between these functions is appropriately small. But because the glued function is equal to the affine function $\ell_p$ on the boundary of each subcube, it cannot deviate much from $\ell_p$ when viewed from the larger scale. In other words, because we can use the Poincar\'e inequality in each smaller subcube, we gain from the scaling of the constant in the Poincar\'e inequality.

\smallskip

Fix $m\in\N$ and $p\in B_1$. We denote, for every $n\in\N$ with $n<m$,
\begin{equation} 
\label{e.defZn}
\mathcal{Z}_n:= 3^n\Zd \cap \cu_m
\end{equation}
so that $\left\{ z+\cu_n\,:\, z\in\mathcal{Z}_n \right\}$ is a partition of $\cu_m$. 

\emph{Step 1.} We show that, for every $n\in\N$ with $n<m$, we have
\begin{multline} 
\label{e.downascale}
3^{-2m}  \left\| v(\cdot,\cu_m,p) - \ell_p \right\|_{\underline{L}^2(\cu_m)}^2
\\
\leq 
C3^{2(n-m)} + \frac{C}{|\mathcal{Z}_n|} \sum_{z\in \mathcal{Z}_n} \left( \nu(z+\cu_n,p) - \nu(\cu_m,p)\right).
\end{multline}
Let $\tilde{v}$ be the function defined in~\eqref{e.gluingfun} for the partition $\left\{ z+\cu_n \,:\, z\in 3^n\Zd\cap \cu_m \right\}$ of $\cu_m$. That is, $\tilde{v}\in \ell_p+ H^1_0(\cu_m)$ satisfies 
\begin{equation*} \label{}
\tilde{v} = v(\cdot,z+\cu_n,p) \quad \mbox{in} \ z+\cu_n. 
\end{equation*}
By the triangle inequality and the Poincar\'e inequality, 
\begin{align*} \label{}
\left\| v(\cdot,\cu_m,p) - \ell_p \right\|_{\underline{L}^2(\cu_m)}^2 
&
 \leq
 2\left\| v(\cdot,\cu_m,p) - \tilde{v} \right\|_{\underline{L}^2(\cu_m)}^2 
 +
 2\left\| \tilde{v} - \ell_p \right\|_{\underline{L}^2(\cu_m)}^2 
 \\ &
 \leq 
C3^{2m}  \left\| \nabla v(\cdot,\cu_m,p) - \nabla \tilde{v} \right\|_{\underline{L}^2(\cu_m)}^2 
 +
 2\left\| \tilde{v} - \ell_p \right\|_{\underline{L}^2(\cu_m)}^2. 
\end{align*}
Applying the first inequality of~\eqref{e.gluingerror0} to the partition $\mathcal{Z}_n$ of $\cu_m$ yields 
\begin{align} 
\label{e.vvtildesnap}
\left\| \nabla v(\cdot,\cu_m,p) - \nabla \tilde{v} \right\|_{\underline{L}^2(\cu_m)}^2
\leq 
\frac{2}{|\mathcal{Z}_n|} \sum_{z\in \mathcal{Z}_n} \left( \nu(z+\cu_n,p) - \nu(\cu_m,p) \right).
\end{align}
Meanwhile, it is clear from~\eqref{e.nubounds},~\eqref{e.obviousformula} and $|p|\leq 1$ that 
\begin{align*} \label{}
\left\| \tilde{v} - \ell_p \right\|_{\underline{L}^2(\cu_m)}^2 
& 
= \frac{1}{|\mathcal{Z}_n|} \sum_{z\in\mathcal{Z}_n} \left\| v(\cdot,z+\cu_n,p) - \ell_p \right\|_{\underline{L}^2(z+\cu_n)}^2 
\\ &  
\leq  \frac{C}{|\mathcal{Z}_n|} \sum_{z\in\mathcal{Z}_n} \left( 3^{2n} \left\| \nabla v(\cdot,z+\cu_n,p) - p \right\|_{\underline{L}^2(z+\cu_n)}^2 \right) 
\\ &  
\leq  \frac{C}{|\mathcal{Z}_n|} \sum_{z\in\mathcal{Z}_n} 3^{2n} \left(   \left\| \nabla v(\cdot,z+\cu_n,p) \right\|_{\underline{L}^2(z+\cu_n)}^2 +|p|^2 \right) 
\leq C3^{2n}.
\end{align*}
Combining the above yields~\eqref{e.downascale}.

\smallskip

\emph{Step 2.} The conclusion. Taking the expectation of~\eqref{e.downascale} and using stationarity, we obtain, for any $m,n\in\N$ with $n<m$,
\begin{align*} \label{}
\E \left[ 3^{-2m}  \left\| v(\cdot,\cu_m,p) - \ell_p \right\|_{\underline{L}^2(\cu_m)}^2
\right]
\leq 
C3^{2(n-m)} + C \E \left[ \nu(\cu_n,p) - \nu(\cu_m,p)\right].
\end{align*}
Taking $n:= \left\lceil \frac m2 \right\rceil$ and using~\eqref{e.nuEmonotone}, we obtain
\begin{equation*} \label{}
\E \left[ 3^{-2m}  \left\| v(\cdot,\cu_m,p) - \ell_p \right\|_{\underline{L}^2(\cu_m)}^2
\right]
\leq
C 3^{-m} + C 3^{-\frac d8 m} + C \omega\left( \left\lceil \tfrac m2 \right\rceil \right). 
\end{equation*}
This yields~\eqref{e.vellpconv}. 
\end{proof}

The previous proposition can be seen as an $L^2$ estimate for the error in homogenization for the Dirichlet problem in a cube with affine boundary data. To see this, set $\ep_m:=3^{-m}$ and notice that the function
\begin{equation*} \label{}
w_m (x):= \ep_m v\left( \frac x{\ep_m}, \cu_m,p\right)
\end{equation*}
is the solution of the Dirichlet problem
\begin{equation*} \label{}
\left\{
\begin{aligned}
& -\nabla \cdot\left( \a\left( \tfrac \cdot{\ep_m} \right) \nabla w_m \right) = 0 & \mbox{in} & \ \cu_0,\\
& w_m = \ell_p & \mbox{on} & \ \partial \cu_0. 
\end{aligned}
\right.
\end{equation*}
Obviously, the solution of
\begin{equation*} \label{}
\left\{
\begin{aligned}
& -\nabla \cdot\left( \ahom \nabla w_{\mathrm{hom}} \right) = 0 & \mbox{in} & \ \cu_0, \\
& w_{\mathrm{hom}} = \ell_p & \mbox{on} & \ \partial \cu_0, 
\end{aligned}
\right.
\end{equation*}
is the function $w_{\mathrm{hom}} = \ell_p$. Moreover, notice by changing variables that
\begin{equation*} \label{}
\left\| w_m - \ell_p \right\|_{L^2(\cu_0)} = 3^{-m} \left\| v( \cdot, \cu_m,p) - \ell_p \right\|_{\underline{L}^2(\cu_m)}.
\end{equation*}
Therefore the limit~\eqref{e.vellpconv} shows homogenization in~$L^2(\cu_0)$ along a subsequence of~$\ep$'s for this specific Dirichlet problem. One may consider this demonstration to be ``cheating'' because it provides no evidence that we have chosen the correct $\ahom$ (any choice of $\ahom$ will give the same solution of the Dirichlet problem with affine boundary data). This is a valid objection, but evidence that $\ahom$ has been chosen correctly and that the quantity $\nu(U,p)$ is capturing information about the homogenization process will be given in the next section. Recall that we encountered a similar phenomenon in the one dimensional case when deciding how the homogenized coefficients should be defined: see~\eqref{e.formula.oned} and the discussion there.

\section{Weak convergence of gradients and fluxes}

In the previous section, we proved the convergence in $L^1(\Omega, \P)$ of the limits 
\begin{equation}
 \label{e.naivelimits}
\left\{
\begin{aligned}
& \fint_{\cu_m} \frac12 \nabla v(\cdot,\cu_m,p) \cdot \a\nabla v(\cdot,\cu_m,p) \to \frac12 p\cdot \ahom p, \\
& 3^{-2m} \fint_{\cu_m} \left| v(\cdot,\cu_m,p) - \ell_p \right|^2 \to 0
\end{aligned}
\right. \qquad \mbox{as} \ m\to \infty. 
\end{equation}
We also showed that an explicit quantitative bound for the modulus $\omega(m)$ for the limit~\eqref{e.defnubar} would give us a rate of convergence for these limits as well. In this section, we push this analysis a bit further by proving some more precise results in the same spirit. We are particularly interested in obtaining results which quantify (up to bounds on~$\omega$) of the following \emph{weak} limits:
\begin{equation} 
\label{e.lessnaivelimits}
\left\{
\begin{aligned}
& \frac12 \nabla v(3^m \cdot,\cu_m,p) \cdot \a(3^m\cdot) \nabla v(3^m\cdot,\cu_m,p) \rightharpoonup \frac12 p\cdot \ahom p, \\
& \a\left( 3^m\cdot\right) \nabla v(3^m\cdot,\cu_m,p) \rightharpoonup \ahom p, \\
& \nabla v(3^m\cdot,\cu_m,p) \rightharpoonup p, \\
\end{aligned}
\right. \qquad \mbox{as} \ m\to \infty.
\end{equation}
That is, we want to address the weak convergence of the energy density, flux, and gradient of the rescaled minimizers $x\mapsto 3^{-m} v(3^mx,\cu_m,p)$. Notice that the limits~\eqref{e.lessnaivelimits} are more precise than~\eqref{e.naivelimits}, in the sense that the former implies the latter. Moreover, while one may first be tempted to focus on the~$L^2$ convergence of solutions, the structure of the problem in fact gives much more direct access to information on the gradients, fluxes and energies of solutions, and it is therefore more efficient to focus our attention on the weak convergence of these quantities, and to derive the $L^2$ convergence of solutions a posteriori. In other words, \emph{homogenization is about the weak convergence of gradients, fluxes and energy densities of solutions. Other convergence results are just consequences of these, but not the main point.}

\smallskip

As in the previous section, our desire is to train ourselves for future chapters by obtaining some crude quantitative bound for the limits~\eqref{e.lessnaivelimits} in terms of~$\omega(m)$. But what is the right way to quantify weak convergence? While there are many ways to do it and the ``right'' way may depend on the context, one very natural choice is to use a negative Sobolev norm like $H^{-1}$. We continue with an informal discussion around this point before getting into the analysis of proving~\eqref{e.lessnaivelimits}. 

\smallskip

Let us say that we have a bounded sequence of functions $\{ f_m \}_{m\in\N} \subseteq L^2(U)$ which weakly converges to some $f\in L^2(U)$. This means that
\begin{equation*} \label{}
\forall g\in L^2(U), \quad \lim_{m\to \infty}\left| \fint_{U} (f_m-f)g \right| = 0. 
\end{equation*}
If we want to quantify weak convergence, then we obviously have to quantify this limit. 
Now, if this limit is uniform over $\| g \|_{\underline{L}^2(U)} \leq 1$ then we have another name for this, which is \emph{strong convergence in $L^2(U)$}, and then perhaps we should quantify this instead! Therefore it makes sense to assume we are in a situation in which convergence does not happen uniformly in~$g$. 
Yet, we hope that the convergence rate depends in some natural way on $g$: perhaps we can get a uniform convergence rate for all smooth $g$ with some derivatives under control. To be more concrete, perhaps we can hope that the limit is uniform over the set of $g$'s with unit $H^1$ norm? This leads us to consider trying to prove a convergence rate for
\begin{equation*} \label{}
\sup\left\{ \left| \fint_U (f_m-f)g \right| \,:\, g\in H^1(U), \,\left\| g \right\|_{\underline{H}^1(U)}  \leq 1 \right\} =  \left\| f_m - f \right\|_{\Hminusul(U)}. 
\end{equation*}
We thus see that the $\Hminusul(U)$ norm (defined in~\eqref{e.def.H-1}) appears in a very natural way (and that other choices are possible as well). If we succeed, we will get an explicit convergence rate for any given $g\in L^2(U)$ in terms of how well $g$ can be approximated by an $H^1(U)$ function. Indeed, for every $M>0$ and $h \in H^1(U)$ with $\left\| h \right\|_{\underline{H}^1(U)}\leq M$,
\begin{align*} \label{}
\left| \fint_U (f_m-f)g \right|
&
\leq 
\left| \fint_U (f_m-f)h \right|
+
\left| \fint_U (f_m-f)(g-h) \right|
\\ & 
\leq 
M \left\| f_m - f \right\|_{\Hminusul(U)} 
+ \left\| f_m -f \right\|_{\underline{L}^2(U)} \left\| g - h \right\|_{\underline{L}^2(U)}.
\end{align*}
Thus, for any $g\in L^2(U)$ and $M>0$,
\begin{multline*} \label{}
\left| \fint_U (f_m-f)g \right| 
\\
\leq
M \left\| f_m - f \right\|_{\Hminusul(U)} 
+  \sup_{m\in\N} \left( \left\| f_m \right\|_{\underline{L}^{2}(U)} + \left\| f \right\|_{\underline{L}^{2}(U)}  \right)
 \inf_{ \left\| h \right\|_{\underline{H}^1(U)} \leq M } \left\| g- h \right\|_{\underline{L}^2(U)}.
\end{multline*}
We can then optimize over the parameter $M$ to make the right side as small as possible (the choice of $M$ will naturally depend on how well $g$ can be approximated by $H^1(U)$ functions). 

\begin{exercise}
Assuming that the sequence $\left\{ f_m \right\}_{m\in\N}\subseteq L^2(U)$ is uniformly bounded, show that
\begin{equation*} \label{}
f_m \rightharpoonup 0 \quad \mbox{weakly in $L^2(U)$} 
\iff 
\left\| f_m  \right\|_{\Hminusul(U)} \rightarrow 0.
\end{equation*}
\end{exercise}

The above discussion motivates the following result, which we split into a pair of propositions. We define, for each $m\in\N$, the random variable 
\begin{equation} 
\label{e.Em.def}
\mathcal{E}(m)
:=  \left( \sum_{n=0}^{m} 3^{n-m}    \left(\frac{1}{|\mathcal{Z}_n|} \sum_{z\in \mathcal{Z}_n}  \left| \a(z+\cu_n) - \ahom \right| \right)^{\frac12} \right)^2.
\end{equation}
This random variable monitors the convergence of the subadditive quantities $\nu$ on the mesoscale grid $\{ z+\cu_n\,:\, z\in \mathcal{Z}_n\}$ of subcubes of $\cu_m$, for every $p\in\Rd$. 

\begin{proposition}
\label{p.qualweakconv}
There exists $C(d,\Lambda)<\infty$ such that, for every $m\in\N$ and $p\in B_1$, 
\begin{multline} 
\label{e.weakminconv}
3^{-2m} \left\| \nabla v (\cdot,\cu_m,p)- p \right\|_{\Hminusul (\cu_m)}^2
+ 3^{-2m}  \left\| \a \nabla v (\cdot,\cu_m,p) - \ahom p \right\|_{\Hminusul (\cu_m)}^2
\\
\leq 
C 3^{-2m} 
+ 
 C \mathcal{E}(m).
\end{multline}
\end{proposition}

\begin{proposition}
\label{p.Em.bound}
There exists $C(d,\Lambda)<\infty$ such that, for every $m\in\N$, 
\begin{equation} 
\label{e.EEm}
\E \left[ \mathcal{E}(m) \right] 
\leq 
C  3^{-\left( \frac d4 \wedge 1\right) m}
+
C \sum_{n=0}^{m} 3^{n-m} \omega\left( \left\lceil  \tfrac n2 \right\rceil \right).
\end{equation}
In particular, $\E\left[ \mathcal{E}(m) \right] \to 0$ as $m\to \infty$. 
\end{proposition}

We could have chosen to state these two propositions as one and hidden the random variable $\mathcal{E}(m)$ inside the proof. The reasons for splitting it like this are threefold: (i) it better reveals the structure of the proof, which comes in two very independent pieces; (ii) the random variable $\mathcal{E}(m)$ will appear on the right side of several other interesting estimates below; (iii) the bound for~$\mathcal{E}(m)$ in~\eqref{e.EEm} has very weak stochastic integrability (only $L^1(\Omega,\P)$) and will be replaced by a much stronger, more explicit bound in the next chapter. 

\smallskip

The random variable~$\mathcal{E}(m)$ may be compared to the right side of~\eqref{e.downascale}, and the proof of Proposition~\ref{p.Em.bound} is actually a simple exercise. We postpone it to later in the section and focus first on the proof of Proposition~\ref{p.qualweakconv}, which is a purely deterministic argument that does not use the stochastic structure of the problem. 

\smallskip

We begin by presenting a functional inequality which gives an estimate of the~$\Hminusul(\cu_m)$ norm of an~$L^2(\cu_m)$ function in terms of its spatial averages on all triadic subcubes of~$\cu_m$. In some sense, this result is just a better formalization of the ``crude'' separation of scales arguments from the previous section, where we used only one mesoscale~$n$. Here we will use \emph{all} the scales, rather than just one. Since we need to use this argument many times in the following chapters, we take care to state a more refined version of it here. 

\begin{proposition}[Multiscale Poincar\'e inequality]
\label{p.mspoin}
\index{Poincar\'e inequality!multiscale}
\index{multiscale Poincar\'e inequality|see Poincar\'e inequality!multiscale}
Fix $m\in\N$ and, for each $n\in\N$ with $n\leq m$, define $\mcl Z_n:= 3^n\Zd \cap \cu_m$. There exists a constant $C(d)<\infty$ such that, for every $f \in L^2(\cu_m)$,
\begin{equation}
\label{e.multiscalepoincare}
 \|f \|_{\Hminusul(\cu_m)}
\leq C \left\| f  \right\|_{\underline{L}^2(\cu_m)}  
+  C \sum_{n=0}^{m-1} 3^n    \left( \left| \mcl Z_n \right|^{-1} \sum_{y\in \mcl Z_n}\left| \left( f \right)_{y+\cu_{n} } \right|^2 \right)^{\frac12}.
\end{equation}
\end{proposition}
\begin{proof}
In view of the definition of the normalized $\Hminusul$ norm given in \eqref{e.def.H-1}, we consider a function $g \in H^1(\cu_m)$ such that 
\begin{equation*} \label{}
3^{-2m} \left| (g)_{\cu_m} \right|^2 + \|\nabla g\|_{\nL^2(\cu_m)}^2 \leq 1
\end{equation*}
and aim to control $\fint_{\cu_m} f g$. We may assume without loss of generality that $ (g)_{\cu_m}  = 0$. Indeed, this follows from the inequalities
\begin{equation*} \label{}
\left| \fint_{\cu_m} fg \right| 
 \leq \left| \fint_{\cu_m} f \cdot \left( g -(g)_{\cu_m} \right) \right| + \left| (g)_{\cu_m} \right| \left| (f)_{\cu_m} \right|
\end{equation*}
and
\begin{equation*} \label{}
\left| (g)_{\cu_m} \right| \left| (f)_{\cu_m} \right| 
\leq
3^{m} \left| (f)_{\cu_m} \right|
\leq 
C 3^{m-1} \left( \left| \mathcal{Z}_{m-1} \right|^{-1} \sum_{y\in \mathcal{Z}_{m-1}} \left| \left( f \right)_{y+\cu_{m-1}} \right|^2 \right)^{\frac12}. 
\end{equation*}
Now, observe that, for every $n\in \{ 0,\ldots,m-1\}$ and $z\in \mcl Z_{n+1}$, we have
\begin{multline*}
\int_{z+\cu_{n+1}} f \cdot \left( g - \left( g\right)_{z+\cu_{n+1}} \right) 
= \sum_{y\in \mcl Z_{n} \cap (z+\cu_{n+1})} \int_{y+\cu_{n}} f \cdot \left( g - \left( g\right)_{y+\cu_{n}} \right) \\
+ |\cu_{n}|\sum_{y\in \mcl Z_{n} \cap (z+\cu_{n+1})} \left(  \left( g\right)_{y+\cu_{n} }  -  \left( g\right)_{z+\cu_{n+1}} \right) \cdot \left( f \right)_{y+\cu_{n}}.
\end{multline*}
Summing over $z \in \mcl Z_{n+1}$ and using H\"older's inequality, we get
\begin{multline*}  
\sum_{z\in \mcl Z_{n+1}} \int_{z+\cu_{n+1}} f \cdot \left( g - \left( g\right)_{z+\cu_{n+1}} \right)  \leq \sum_{y\in \mcl Z_{n}} \int_{y+\cu_{n}} f \cdot \left( g - \left( g\right)_{y+\cu_{n}} \right) \\
+|\cu_{n}| \Bigg(\sum_{\substack{z \in \mcl Z_{n+1} \\ y \in \mcl Z_{n} \cap (z+\cu_{n+1})}} \Ll|   \left( g\right)_{y+\cu_{n} } -  \left( g\right)_{z+\cu_{n+1}}\Rr|^2\Bigg)^{\frac 1 2} \Ll( \sum_{y \in \mcl Z_{n}}\left| \left( f\right)_{y+\cu_{n} } \right|^2  \Rr)^\frac 1 2.
\end{multline*}
By the Jensen and Poincar\'e inequalities, for each $z \in \mcl Z_{n+1}$,
\begin{align*}
\sum_{y\in \mcl Z_{n} \cap (z+\cu_{n+1})} \left|  \left( g\right)_{y+\cu_{n} } -  \left( g\right)_{z+\cu_{n+1}}  \right|^{2}
& 
=  \sum_{y\in \mcl Z_{n} \cap (z+\cu_{n+1})} \left|  \fint_{y+\cu_n} \left( g - \left( g \right)_{z+\cu_{n+1}} \right) \right|^2
\\ & 
\leq  \sum_{y\in \mcl Z_{n} \cap (z+\cu_{n+1})}   \fint_{y+\cu_n} \left| g - \left( g \right)_{z+\cu_{n+1}} \right|^2 
\\ & 
= 3^d \fint_{z+\cu_{n+1}} \left| g - \left( g \right)_{z+\cu_{n+1}} \right|^2
\\ &
\leq C3^{2n}  \fint_{z+\cu_{n+1}} \left| \nabla g \right|^{2} .
\end{align*}
Since $\|\nabla g\|_{\nL^2(\cu_m)} \leq 1$ and $|\mcl Z_n| = |\cu_m|/|\cu_n|$, combining the last two displays yields
\begin{multline*}
\sum_{z\in \mcl Z_{n+1}} \int_{z+\cu_{n+1}} f \cdot \left( g - \left( g\right)_{z+\cu_{n+1}} \right) \\ \leq \sum_{y\in \mcl Z_{n}} \int_{y+\cu_{n}} f \cdot \left( g - \left( g\right)_{y+\cu_{n}} \right)
+C\, |\cu_m| \, 3^{n} \left( |\mcl Z_{n}|^{-1}\sum_{y\in \mcl Z_{n}}\left| \left( f\right)_{y+\cu_{n} } \right|^2 \right)^{\frac12}.
\end{multline*}
Iterating this inequality and using that $(g)_{\cu_m} = 0$, we get
\begin{equation*}
\int_{\cu_m} f g \leq \sum_{z\in \mcl Z_0} \int_{z+\cu_0} f \cdot \left( g - \left( g\right)_{z+\cu_{0}} \right) + C \, \left| \cu_m\right| \, \sum_{n=0}^{m-1} 3^n    \left( \left| \mcl Z_n \right|^{-1} \sum_{y\in \mcl Z_n}\left| \left( f \right)_{y+\cu_{n} } \right|^2 \right)^{\frac12}.
\end{equation*}
By the Poincar\'e inequality,
$$
\sum_{z \in \mcl Z_0} \int_{z + \cu_0} |g - (g)_{z + \cu_0}|^{2}  \le C \int_{\cu_m} |\nabla g|^{2}  = C\, |\cu_m|.
$$
Hence, by H\"older's inequality, 
\begin{equation*}
\int_{\cu_m} f  g  \le C|\cu_m|^{\frac1{2}} \left(  \int_{\cu_m} \left| f\right|^2 \right)^{\frac12}  
+ C \, \left| \cu_m\right| \, \sum_{n=0}^{m-1} 3^n    \left( \left| \mcl Z_n \right|^{-1} \sum_{y\in \mcl Z_n}\left| \left( f \right)_{y+\cu_{n} } \right|^2 \right)^{\frac12}.
\end{equation*}
Dividing by $|\cu_m|$ and taking the supremum over all such $g$ yields the result.
\end{proof}

We coined the previous proposition the ``multiscale Poincar\'e inequality'' without explanation. Let us give one here. The usual Poincar\'e inequality for the cube $\cu_m$ states that there exists a constant $C(d)<\infty$ such that, for every $u\in H^1(\cu_m)$,
\begin{equation}
\label{e.usualPoincare}
\left\| u - \left( u \right)_{\cu_m} \right\|_{\underline{L}^2(\cu_m)} \leq C 3^{m} \left\| \nabla u \right\|_{\underline{L}^2(\cu_m)}. 
\end{equation}
We can see that the scaling of the factor $3^m$ is sharp by considering any nonconstant affine function. However, we may wonder if we can do better for a function whose gradient exhibits large-scale cancellations. 
Applying the previous inequality to $f = \nabla u$ gives us the bound 
\begin{equation*} \label{}
\left\| \nabla u \right\|_{\Hminusul(\cu_m)} 
\leq
C \left\|  \nabla u \right\|_{\underline{L}^2(\cu_m)}  
+  C \sum_{n=0}^{m-1} 3^n    \left( \left| \mcl Z_n \right|^{-1} \sum_{y\in \mcl Z_n}\left| \left(  \nabla u \right)_{y+\cu_{n} } \right|^2 \right)^{\frac12}.
\end{equation*}
If~$\nabla u$ is ``canceling itself out'' in the sense that its spatial averages on large scale triadic subcubes are much smaller than the size of the gradient itself, then the right side of the previous inequality will be much smaller than that of~\eqref{e.usualPoincare}, because we pay only a factor~$3^n$, for each triadic scale $n<m$, against the spatial average of~$\nabla u$, not its absolute size in~$L^2$. Moreover, while this gives us an estimate of the~$\Hminusul(\cu_m)$ norm of~$\nabla u$, the latter actually bounds the~$L^2(\cu_m)$ oscillation of~$u$ itself. This is the content of the next lemma. 
\begin{lemma}
\label{l.integrate.H-1}
There exists a $C(d) < \infty$ such that, for every $m \in \N$ and $u \in H^1(\cu_m)$, 
\begin{equation}
\label{e.suff.mP.gradient}
\left\| u - \left( u \right)_{\cu_m} \right\|_{\underline{L}^2(\cu_m)} 
\leq 
C  \|\nabla u\|_{\Hminusul(\cu_m)}
\end{equation}
and, similarly, for every $v\in H^1_0(\cu_m)$,
\begin{equation}
\label{e.suff.mP.gradient.0}
\left\| v \right\|_{\underline{L}^2(\cu_m)} 
\leq 
C  \| \nabla v \|_{\Hminusul(\cu_m)}.
\end{equation}
\end{lemma}

We will very often use Proposition~\ref{p.mspoin} and Lemma~\ref{l.integrate.H-1} in combination, and therefore state the combined estimate as an immediate corollary.
\begin{corollary}  
\label{c.mP.gradient}
There exists a constant $C(d) < \infty$ such that for every $m \in \N$ and $u \in H^1(\cu_m)$,
\begin{equation} 
\label{e.mpcor}
 \|u - (u)_{\cu_m}\|_{\underline{L}^2(\cu_m)}
\leq C \left\| \nabla u  \right\|_{\underline{L}^2(\cu_m)}  
+  C \sum_{n=0}^{m} 3^n    \left( \left| \mcl Z_n \right|^{-1} \sum_{y\in \mcl Z_n}\left| \left( \nabla u \right)_{y+\cu_{n} } \right|^2 \right)^{\frac12}.
\end{equation}
\end{corollary}
\begin{proof}[Proof of Lemma~\ref{l.integrate.H-1}]
We start with the proof of the first estimate~\eqref{e.suff.mP.gradient}. 
We may assume that $(u)_{\cu_m} = 0$. 
Let $w \in H^1(\cu_m)$ be the unique (up to an additive constant) solution of the Neumann boundary-value problem 
\begin{equation}
\label{e.mP.Neumann}
\left\{ 
\begin{aligned}
& -\Delta w = u &  \mbox{in} & \ \cu_m, \\
&  \mathbf{n} \cdot \nabla w = 0 & \mbox{on} & \ \partial \cu_m.
\end{aligned} 
\right.
\end{equation}
By Lemma~\ref{l.H2.est}, we have $w \in H^2(\cu_m)$ and, for some $C(d) < \infty$, 
\begin{equation}
\label{e.mP.H2}
\left\| \nabla^2 w \right\|_{\underline{L}^2(\cu_m)} 
\leq 
C \left\| u \right\|_{\underline{L}^2(\cu_m)}.
\end{equation}
In order to estimate $\left| \left( \nabla w \right)_{\cu_m} \right|$, we fix 
$$
p = \frac{(\nabla w)_{\cu_m}}{|(\nabla w)_{\cu_m}|} \in B_1
$$ 
and test the equation \eqref{e.mP.Neumann} against the affine function $x \mapsto p \cdot x$ to obtain
\begin{equation*}  
\left|(\nabla w)_{\cu_m}\right| = \fint_{\cu_m} p \cdot \nabla w = \fint_{\cu_m} u(x) \, p\cdot x \, dx \le C 3^{m} \|u\|_{\underline L^2(\cu_m)}.
\end{equation*}
Testing the equation \eqref{e.mP.Neumann} with $u \in H^1(\cu_m)$ and using~\eqref{e.mP.H2} yields
\begin{align*}
\left\| u \right\|_{\underline{L}^2(\cu_m)}^2 
 = 
\fint_{\cu_m} u^2  
& = \fint_{\cu_m} \nabla u\cdot \nabla w \\
& \leq \left\| \nabla u \right\|_{\Hminusul(\cu_m)} \left( 3^{-2m} \left|(\nabla w)_{\cu_m}\right|^2 + \left\| \nabla \nabla w \right\|_{\underline{L}^2(\cu_m)}^2 \right)^{\frac12}  \\
& \leq C \left\| \nabla u \right\|_{\Hminusul(\cu_m)} \left\| u \right\|_{\underline{L}^2(\cu_m)} ,
\end{align*}
which completes the proof of~\eqref{e.suff.mP.gradient}.

\smallskip

The second estimate~\eqref{e.suff.mP.gradient.0} is a consequence of~\eqref{e.suff.mP.gradient} and the following estimate: for every $v\in H^1_0(\cu_m)$, 
\begin{equation*}
\left| \left( v \right)_{(\cu_m)} \right|
\leq
C\left\| \nabla v \right\|_{\Hminusul(\cu_m)}.
\end{equation*}
To see this, fix a unit direction $e$ and let $\mathbf{f}$ denote the vector field $\mathbf{f}(x) = (e\cdot x) e$ so that $\nabla \cdot \mathbf{f} = 1$ and the components $\mathbf{f}_i$ of $\mathbf{f}$ belong to $H^1(\cu_m)$ and satisfy the bound $\left\| \mathbf{f}_i \right\|_{\underline{H}^1(\cu_m)} \leq C$. Then we find that 
\begin{equation*}
\left| \left( u \right)_U \right| 
= 
\left| \fint_U \nabla u \cdot \mathbf{f} \right| 
\leq 
C \left\| \nabla u \right\|_{\Hminusul(\cu_m)}.
\end{equation*}
This completes the proof. 
\end{proof}

We are now ready to give the proof of Proposition~\ref{p.qualweakconv}. 

\begin{proof}[{Proof of Proposition~\ref{p.qualweakconv}}]

\emph{Step 1.} We fix $p\in B_1$, $m\in\N$ and denote $v:= v(\cdot,\cu_m,p)$.  We estimate the difference between~$\nabla v$ and $p$ in $\hat H^{-1}(\cu_m)$ by mimicking the proof of Proposition~\ref{p.vellpconv}. Let $\mathcal{Z}_n$ be denoted by~\eqref{e.defZn}. By the multiscale Poincar\'e inequality (Proposition~\ref{p.mspoin}),
\begin{equation*}
\left\| \nabla v - p \right\|_{\Hminusul (\cu_m)}^2
\leq 
C\left\| \nabla v - p \right\|_{\underline{L}^{2} (\cu_m)}^2
+ 
 C \left( \sum_{n=0}^{m-1} 3^n    \left( \frac1{\left| \mcl Z_n \right|} \sum_{z\in \mcl Z_n}\left| \left( \nabla v - p \right)_{z+\cu_{n} } \right|^2 \right)^{\frac12} \right)^2.
\end{equation*}
The first term on the right side is estimated easily by~\eqref{e.nubounds} and~\eqref{e.obviousformula}:
\begin{equation} 
\label{e.easyboundgradweak}
\left\| \nabla v - p \right\|_{\underline{L}^{2} (\cu_m)}^2 
\leq 
2\left\| \nabla v  \right\|_{\underline{L}^{2} (\cu_m)}^2 + 2|p|^2
\leq C|p|^2 \leq C. 
\end{equation}
To estimate the expectation of the second term, we
fix $n\in\N$ and let $\tilde{v}_n$ be the function $\tilde{v}$ described in Step~1 of the proof of Proposition~\ref{p.vellpconv} which is obtained by gluing the minimizers $v(\cdot,z+\cu_n,p)$ for $z\in 3^n\Zd\cap \cu_m$. 
Using that $\left( \nabla \tilde v_n \right)_{z+\cu_n} = p$ for every $z\in 3^n\Zd\cap \cu_m$, we obtain from Jensen's inequality that 
\begin{equation*} \label{}
 \frac{1}{|\mathcal{Z}_n|} \sum_{z\in \mathcal{Z}_n} 
 \left| \left( \nabla v \right)_{z+\cu_n}  - p \right|^2 
  \leq 
   \left\| \nabla v(\cdot,\cu_m,p) - \nabla \tilde{v}_n \right\|_{\underline{L}^2(\cu_m)}^2.
\end{equation*}
According to~\eqref{e.vvtildesnap}, we have 
\begin{align*} 
\left\| \nabla v(\cdot,\cu_m,p) - \nabla \tilde{v} \right\|_{\underline{L}^2(\cu_m)}^2
\leq 
 \frac{2}{|\mathcal{Z}_n|} \sum_{z\in \mathcal{Z}_n} \left( \nu(z+\cu_n,p) - \nu(\cu_m,p) \right).
\end{align*}
Thus
\begin{equation*} \label{}
 \frac{1}{|\mathcal{Z}_n|} \sum_{z\in \mathcal{Z}_n} 
 \left| \left( \nabla v \right)_{z+\cu_n}  - p \right|^2 
  \leq 
\frac{2}{|\mathcal{Z}_n|} \sum_{z\in \mathcal{Z}_n} \left( \nu(z+\cu_n,p) - \nu(\cu_m,p) \right).  
\end{equation*}
Combining the above displays yields 
\begin{equation}
3^{-2m} \left\| \nabla v - p \right\|_{\Hminusul (\cu_m)}^2
\leq 
C 3^{-2m} 
+ 
 C \mathcal{E}(m).
\end{equation}
This completes the proof of the estimate for the first term on the left side of~\eqref{e.weakminconv}. 

\smallskip

The proof of the $H^{-1}$ estimate for the fluxes is the focus of the next three steps. We begin with an  identity which gives a relationship between the spatial average of the flux of the minimizer $v(\cdot,U,p)$ and the gradient of the quantity $\nu(U,p)$. This is a first encounter with an important idea that will reappear in later chapters.

\smallskip

\emph{Step 2.} We show that, for every bounded Lipschitz domain $U\subseteq \Rd$ and $p\in\Rd$,  
\begin{equation}
\label{e.nufluxidentity}
 \fint_{U} \a\nabla v (\cdot,U,p) = \a(U) p.
\end{equation}
By the definition of $\a(U)$ in \eqref{e.nurepresentation}, we have, for every $p \in \Rd$,
\begin{equation*}  
\fint_U \nabla v(\cdot,U,p) \cdot \a \nabla v(\cdot,U,p) = p \cdot \a(U) p.
\end{equation*}
By polarization, we deduce that for every $p,q \in \Rd$,
\begin{align*}  
q \cdot \a(U) p  & = \frac 1 4 \Ll( (p+q) \cdot \a(U)(p+q) - (p-q) \cdot \a(U) (p-q) \Rr) \\
& = \fint_U \nabla v(\cdot,U,q) \cdot \a \nabla v(\cdot,U,p).
\end{align*}
By the first variation (i.e., using that $v(\cdot,U,p)$ is a weak solution of the equation \eqref{e.DP} and testing it with $v(\cdot,U,q) - \ell_q \in H^1_0(U)$), we obtain that
\begin{equation*} \label{}
q \cdot \a(U) p = \fint_U  q\cdot \a\nabla v(\cdot,U,p),
\end{equation*}
and thus~\eqref{e.nufluxidentity} holds. 

\smallskip

\emph{Step 3.} We argue that, for every $n\in\N$,
\begin{equation}
\label{e.nufluxcaptcha}
\sup_{p\in B_1}
\left| \fint_{\cu_n} \a\nabla v(\cdot,\cu_n,p) - \ahom p \right|^2
 \leq 
 C \sup_{p \in B_1}  \left| \nu(\cu_n,p) -  \overline{\nu}(p)\right|.
\end{equation}
By~\eqref{e.nufluxidentity}, we have 
\begin{equation*} \label{}
\left| \fint_{\cu_n} \a\nabla v(\cdot,\cu_n,p) - \ahom p \right|  
= \left| \a(U)p - \ahom p \right|,
\end{equation*}
and moreover, by \eqref{e.numatrixbounds} and \eqref{e.ahombounds},
\begin{equation*} \label{}
\sup_{p\in B_1} \left|  (\a(U) - \ahom)p \right|^2 
\leq C \sup_{p\in B_1}  \left|\frac12p\cdot (\a(U) - \ahom)p \right| 
= C\sup_{p\in B_1} \left|\nu(\cu_n,p) -  \overline{\nu}(p)\right|.
\end{equation*}
Combining these two displays yields~\eqref{e.nufluxcaptcha}.

\smallskip

\emph{Step 4.} We complete the proof of~\eqref{e.weakminconv}. We once again fix $p\in B_1$ and denote $v:= v(\cdot,\cu_m,p)$. By the multiscale Poincar\'e inequality, 
\begin{multline}
 \label{e.qualweakconvgradmspoin}
\left\| \a \nabla v - \ahom p \right\|_{\Hminusul (\cu_m)}^2
\\
\leq 
C\left\| \a \nabla v - \ahom p \right\|_{\underline{L}^{2} (\cu_m)}^2
+ 
 C \left( \sum_{n=0}^{m-1} 3^n    \left( \frac1{\left| \mcl Z_n \right|} \sum_{z\in \mcl Z_n}\left| \left( \a \nabla v - \ahom p \right)_{z+\cu_{n} } \right|^2 \right)^{\frac12} \right)^2.
\end{multline}
By the triangle inequality, 
\begin{multline*} \label{}
\frac1{\left| \mcl Z_n \right|} \sum_{z\in \mcl Z_n}\left| \left( \a \nabla v - \ahom p \right)_{z+\cu_{n} } \right|^2
\\
\leq 
\frac2{\left| \mcl Z_n \right|} \sum_{z\in \mcl Z_n} \left( \left\| \a \nabla v - \a\nabla v(\cdot,z+\cu_n,p) \right\|_{\underline{L}^2(z+\cu_n)}^2 +  \left| \fint_{z+\cu_n} \a\nabla v(\cdot,z+\cu_n,p) - \ahom p \right|^2 \right).
\end{multline*}
Arguing as in Step~1, we obtain the bounds
\begin{equation*} \label{}
\left\| \a \nabla v - \ahom p \right\|_{\underline{L}^{2} (\cu_m)}^2
\leq C,
\end{equation*}
and
\begin{equation*} \label{}
\frac1{\left| \mcl Z_n \right|} \sum_{z\in \mcl Z_n} \left\| \a \nabla v - \a\nabla v(\cdot,z+\cu_n,p) \right\|_{\underline{L}^2(z+\cu_n)}^2
\leq 
\frac{C}{|\mathcal{Z}_n|} \sum_{z\in \mathcal{Z}_n} \left| \nu(z+\cu_n,p) - \nu(\cu_m,p) \right|.
\end{equation*}
The previous display and~\eqref{e.nufluxcaptcha} yield
\begin{multline*} \label{}
 \frac1{\left| \mcl Z_n \right|} \sum_{z\in \mcl Z_n} \left( \left\| \a \nabla v - \a\nabla v(\cdot,z+\cu_n,p) \right\|_{\underline{L}^2(z+\cu_n)}^2 +  \left| \fint_{z+\cu_n} \a\nabla v(\cdot,z+\cu_n,p) - \ahom p \right|^2 \right) 
\\
\leq
C  \frac{1}{|\mathcal{Z}_n|} \sum_{z\in \mathcal{Z}_n} \left( \left|   \nu(z+\cu_n,p) - \nu(\cu_m,p) \right| + \sup_{e\in B_1} \left|  \nu(z+\cu_n,e) -  \overline{\nu}(e) \right| \right). 
\end{multline*}
Combining the above, we obtain
\begin{align*} \label{}
\lefteqn{ 
\left\| \a \nabla v - \ahom p \right\|_{\Hminusul (\cu_m)}^2
} \qquad & \\ & 
\leq C +  C3^{2m} \mathcal{E}(m) + C  \left( \sum_{n=0}^m 3^n \left( \frac{1}{|\mathcal{Z}_n|} \sum_{z\in \mathcal{Z}_n} \sup_{e\in B_1} \left| \nu(z+\cu_n,e) - \overline{\nu}(e) \right|   \right)^{\frac12}  \right)^2
\\ & 
\leq C + C3^{2m}\mathcal{E}(m) . 
\end{align*}
This completes the proof of~\eqref{e.weakminconv}.
\end{proof}

We now give the proof of Proposition~\ref{p.Em.bound}.

\begin{proof}[{Proof of Proposition~\ref{p.Em.bound}}]
By H\"older's inequality, 
\begin{equation} 
\mathcal{E}(m)
\leq C\sum_{n=0}^{m} 3^{n-m} \frac{1}{|\mathcal{Z}_n|} \sum_{z\in \mathcal{Z}_n}  \left| \a(z+\cu_n) - \ahom \right|.
\end{equation}
Taking the expectation of both sides and using stationarity and~\eqref{e.preL1conv}, we find that
\begin{align*} \label{}
\E \left[ \mathcal{E}(m) \right]
&
\leq 
C \sum_{n=0}^{m} 3^{n-m} \left( 3^{-\frac d4n} + C \omega\left( \left\lceil  \tfrac n2 \right\rceil \right) \right)
\\ & 
\leq 
C  3^{-\left( \frac d4 \wedge 1\right) m}
+
C \sum_{n=0}^{m} 3^{n-m} \omega\left( \left\lceil  \tfrac n2 \right\rceil \right).
\end{align*}
To see that $\E \left[ \mathcal{E}(m) \right] \to 0$ as $m\to \infty$, observe that we can estimate the right side of~\eqref{e.EEm} crudely by  
\begin{align*} \label{}
 \sum_{n=0}^{m} 3^{n-m} \omega(n)
&
\leq 
C  \sum_{n=0}^{\left\lceil \frac m2\right\rceil } 3^{n-m}
+  \sum_{n=\left\lceil \frac m2\right\rceil }^{m-1} \left( 3^{-\frac d4n} + C \omega\left( \left\lceil  \tfrac n2 \right\rceil \right) \right) 
\\ & 
\leq C3^{-\frac m2} + C 3^{-\frac d8m} + 
 C \omega\left( \left\lceil \tfrac m4 \right\rceil \right) \to 0 \quad \mbox{as} \ m\to \infty.\qedhere
\end{align*}
\end{proof}

We next present a $W^{-2,q}(\cu_m)$ estimate which quantifies the weak convergence of the energy densities of $v(\cdot,\cu_m,p)$. Note that the energy density is a priori only an $L^1(\cu_m)$ function and $L^1$ embeds into $W^{-2,q}$ for every $q\in \left[1,\frac{d}{d-1}\right)$ by Sobolev embedding and duality. Therefore it is natural to quantify weak convergence in $L^1$ by strong convergence in $W^{-2,q}$ for this range of~$q$. The fact that we can obtain $q>1$ is not very important here: in the next section we will just use $q=1$ to make the statement simpler. Also observe that we can get an estimate in $W^{-1,1}$, if desired, by interpolating between $W^{-2,1}$ and $L^1$.

\begin{proposition}
\label{p.qualweakconv.energy}
Suppose that $q \in \left[1,\frac{d}{d-1}\right)$.
There exists $C(q,d,\Lambda)<\infty$ such that, for every $m\in\N$ and $p\in B_1$, 
\begin{equation} 
\label{e.weakminconv.energy}
3^{-2m} \left\| \tfrac12 \nabla v (\cdot,\cu_m,p) \cdot \a\nabla v (\cdot,\cu_m,p) - \tfrac12  p\cdot \ahom p \right\|_{\underline{W}^{-2,q}(\cu_m)}
\leq  
C 3^{-2m}  +  C \mathcal{E}(m).
\end{equation}
\end{proposition}

Proposition~\ref{p.qualweakconv.energy} is an immediate consequence of Proposition~\ref{p.qualweakconv}, Lemma~\ref{l.integrate.H-1}, and the following quantitative version of the so-called\index{div-curl lemma|(}
 ``div-curl lemma'' of Murat-Tartar (cf.~\cite[Chapter 7]{Tartar}). Recall that~$\Ls(U)$ is the space of~$L^2$ solenoidal vector fields on~$U$ defined in~\eqref{e.Ls.not};  the negative Sobolev norms are defined in~\eqref{e.first.def.w-alpha.norm}. 

\begin{lemma}[Div-curl lemma]
\label{l.divcurl}
Suppose that $q \in \left[1,\frac{d}{d-1}\right)$. Let $u_1,u_2\in H^1(U)$ and $\g_1,\g_2\in \Ls(U)$. Then there exists $C(U,q,d)<\infty$ such that 
\begin{multline} 
\label{e.divcurl}
\left\| \nabla u_1 \cdot \g_1 - \nabla u_2 \cdot \g_2 \right\|_{W^{-2,q}(U)}
\\
\leq C \left( \left\| \nabla u_1 \right\|_{L^2(U)}\left\| \g_1-\g_2 \right\|_{H^{-1}(U)} 
+\left\| \g_2 \right\|_{L^2(U)} \left\| u_1-u_2\right\|_{L^2(U)}\right).
\end{multline}
\end{lemma}
\begin{proof}
We may assume that $(u_1)_U =(u_2)_U = 0$. 
Fix $q\in \left[1,\frac{d}{d-1}\right)$. Let $p := q'$ and observe that $p\in (d,\infty]$. 
Select $\varphi\in W^{2,p}_0(U)$ and compute
\begin{align*} \label{}
\int_U \varphi \left( \nabla u_1 \cdot \g_1 - \nabla u_2 \cdot \g_2 \right) 
&
= - \int_U \nabla \varphi \cdot \left( u_1 \g_1-u_2\g_2\right) 
\\ &
= - \int_U \nabla \varphi \cdot\left( u_1 (\g_1-\g_2) + (u_1-u_2)\cdot \g_2 \right).
\end{align*}
Thus 
\begin{multline*} \label{}
\left| \int_U \varphi \left( \nabla u_1 \cdot \g_1 - \nabla u_2 \cdot \g_2 \right)  \right| 
\\
\leq
\left\| u_1 \nabla \varphi  \right\|_{H^1(U)} \left\| \g_1-\g_2 \right\|_{H^{-1}(U)} 
+ \left\| \g_2 \nabla \varphi \right\|_{L^2(U)} \left\| u_1-u_2\right\|_{L^2(U)}.
\end{multline*}
Using the H\"older and Sobolev-Poincar\'e inequalities, we estimate
\begin{align*} \label{}
\left\| u_1 \nabla \varphi  \right\|_{H^1(U)} 
& \leq C \left( \left\| \nabla u_1 \nabla \varphi  \right\|_{L^2(U)} +\left\|  u_1 \nabla \nabla \varphi  \right\|_{L^2(U)} \right)
\\ & 
\leq 
C \left( \left\| \nabla \varphi \right\|_{L^\infty(U)} \left\| \nabla u_1 \right\|_{L^2(U)}  
+ \left\| \nabla \nabla \varphi \right\|_{L^p(U)} \left\| u_1 \right\|_{L^{\frac{2p}{p-2}}(U)}  \right) 
\\ & 
\leq C \left\| \nabla \nabla \varphi \right\|_{L^p(U)} \left\| \nabla u_1 \right\|_{L^2(U)},
\end{align*}
where we used $p > d$ for the first term, and $\frac{2p}{p-2} < \frac{2d}{d-2}$ for the second term. Similarly,
\begin{equation*} \label{}
\left\| \g_2 \nabla \varphi \right\|_{L^2(U)}
 \leq \left\| \g_2 \right\|_{L^2(U)} \left\| \nabla \varphi \right\|_{L^\infty(U)} 
 \leq \left\| \g_2 \right\|_{L^2(U)} \left\| \nabla\nabla  \varphi \right\|_{L^p(U)}. 
\end{equation*}
Combining the above, we deduce that 
\begin{multline*} \label{}
\left| \int_U \varphi \left( \nabla u_1 \cdot \g_1 - \nabla u_2 \cdot \g_2 \right)  \right| 
\\
\leq
C  \left\| \nabla \nabla \varphi \right\|_{L^p(U)} 
\left( 
\left\| \nabla u_1 \right\|_{L^2(U)}\left\| \g_1-\g_2 \right\|_{H^{-1}(U)} 
+\left\| \g_2 \right\|_{L^2(U)} \left\| u_1-u_2\right\|_{L^2(U)}
 \right).
\end{multline*}
Taking the supremum over $\varphi\in W^{2,p}_0(U)$ completes the proof. 
\end{proof}\index{div-curl lemma|)}

\section{Homogenization of the Dirichlet problem}
\label{s.homogDPintro}

The purpose of this section is to prove an estimate of the homogenization error for a quite general class of Dirichlet problems. The statement is given in Theorem~\ref{t.DP.blackbox}, below, which is the culmination of the theory developed in this chapter. Similarly to the previous section, the estimate is stated in terms of the convergence of the subadditive quantity~$\nu$ in triadic cubes and gives a rate of homogenization in $L^1(\Omega,\P)$ which essentially depends on the modulus~$\omega(m)$ defined in~\eqref{e.defomega}. While we do not estimate $\omega(m)$ explicitly until the next chapter, the estimate in Theorem~\ref{t.DP.blackbox} suffices to yield a qualitative homogenization result.

\smallskip

It is with an eye toward an important future application of the result here (see Chapter~\ref{c.regularity}) that we allow for relatively rough boundary conditions. 
We therefore consider the Dirichlet problem in Lipschitz domains with boundary data slightly better than $H^1(U)$, namely in $W^{1,2+\delta}(U)$ for some $\delta >0$. 
For the same reason, we also provide an error estimate whose random part is uniform in the boundary condition. The argument presented here will not lead to sharp error estimates, even for smooth boundary conditions, but it foreshadows similar arguments in Chapter~\ref{c.twoscale} which do yield sharp estimates. 

\smallskip

We follow very natural and by-now standard arguments introduced in periodic homogenization in the 1970s and 80s, which use the first-order correctors. The argument here can be compared for instance to those of~\cite[Chapter 1.5.5]{BLP} or of~\cite[Theorem 2.6]{Allaire0}. A central idea in homogenization theory is that one can reduce more complicated problems---such as obtaining the convergence of the Dirichlet problem for general boundary data---to simpler ones, like the construction of ``affine-like'' solutions. (We have already informally discussed this idea, see the paragraph containing~\eqref{e.first.corrector}.) 

\smallskip

This reduction is sometimes called the \emph{oscillating test function method} and the affine-like solutions which are usually used in the argument are called the \emph{(first-order) correctors}. However, the particular choice of ``affine-like solution'' is not very important and since we have not yet introduced the first-order correctors (this will wait until Section~\ref{s.correctorsdefined}), their role is played by the ``finite-volume correctors'' defined for each $n\in\N$ and $e\in\Rd$ by
\begin{equation} 
\label{e.FVC}
\phi_{n,e} (x):= v(x,\cu_n,e) - e\cdot x.
\end{equation}
In other words,~$\phi_{n,e}$ is the difference of the solution of the Dirichlet problem in~$\cu_n$ with boundary data~$\ell_e(x):=e\cdot x$ and the affine function~$\ell_e$ itself. 

\smallskip

The forthcoming argument consists of deriving estimates for the homogenization error for a general Dirichlet problem in terms of the following quantity: for each $\ep\in(0,1)$, we choose $m:= \left\lfloor \left| \log \ep \right| / \log 3 \right\rfloor \in\N$ so that $\ep \in \left[ 3^{-m}, 3^{-m+1} \right)$ and then set
\begin{equation} 
\label{e.E'm}
\mathcal{E}'(\ep):= \sum_{k=1}^d \left( \ep \left\|  \phi_{m,e_k}\left( \tfrac \cdot \ep \right) \right\|_{L^2(\ep\cu_m)}
+
\left\| \a\left( \tfrac \cdot\ep\right)\left( e_k + \nabla \phi_{m,e_k}\left(\tfrac\cdot\ep\right) \right) - \ahom e_k \right\|_{H^{-1}(\ep\cu_m)}  \right)^2.
\end{equation}
Notice that this quantity measures: (i) how flat (affine-like) the functions $x\mapsto e\cdot x + \phi_e(x)$ are, and (ii) the rate at which their fluxes are weakly converging to the homogenized flux. Also observe that it is just a rescaled version of the quantity estimated in Proposition~\ref{p.qualweakconv}. Indeed, the latter result and Lemma~\ref{l.integrate.H-1} imply that 
\begin{equation} 
\label{e.linktopast}
\mathcal{E}'(\ep)
\leq C\left(  \ep^2  +  \mathcal{E}(m) \right), \quad \mbox{for} \ m = \left\lfloor \left| \log \ep \right| / \log 3 \right\rfloor.
\end{equation}
In order to avoid confusion coming from the rescaling of the $H^{-1}$ norm in the application of Proposition~\ref{p.qualweakconv}, we recall from~\eqref{e.notation.Sobolev.scaling-} that, if we take a distribution~$f \in H^{-1} ( \cu_m)$ and we rescale it by setting~$f_\ep (x):= f\left( \tfrac x\ep \right)$, then 
\begin{equation*} \label{}
\left\| f_\ep \right\|_{\underline{H}^{-1}(\ep \cu_m)}  
= 
\ep \left\| f \right\|_{\underline{H}^{-1}(\cu_m)},
\end{equation*} 
and similarly with $\Hminus$ in place of~$\underline{H}^{-1}$. 

\smallskip

The main result of this section is the following theorem, which in view of the above discussion provides a direct link between the convergence of the subadditive quantity~$\nu$ and the error in homogenization of the Dirichlet problem. 

\begin{theorem}
\label{t.DP.blackbox}
Fix $\delta>0$, a bounded Lipschitz domain $U\subseteq \cu_0$, $\ep\in (0,1]$ and $f\in W^{1,2+\delta}(U)$. Let $m\in\N$ be such that $3^{m-1} < \ep^{-1} \leq 3^m$, and let $u^\ep, u \in f+H^1_0(U)$ respectively denote the solutions of the Dirichlet problems 
\begin{equation}
\label{e.uep.and.ubar.bb}
\left\{
\begin{aligned}
& -\nabla \cdot \left( \a\left(\tfrac x\ep\right) \nabla u^\ep \right) = 0 &  \mbox{in} & \ U, \\
& u^\ep = f & \mbox{on} & \ \partial U,
\end{aligned}
\right.
\quad \mbox{and} \quad 
\left\{
\begin{aligned}
& -\nabla \cdot \left( \ahom \nabla  u  \right) = 0 &  \mbox{in} & \ U, \\
&  u  = f & \mbox{on} & \ \partial U.
\end{aligned}
\right.
\end{equation}
There exist $\beta(\delta, d,\Lambda)>0$ and $C(\delta,U,d,\Lambda)<\infty$ such that for every $r\in (0,1)$,
\begin{multline}
\label{e.DPestimates.bb}
\left\| u^\ep -   u  \right\|_{L^{2}(U)} + 
\left\| \nabla u^\ep - \nabla  u  \right\|_{\Hminus(U)}  
+ \left\| \a\left(\tfrac\cdot\ep\right) \nabla u^\ep - \ahom \nabla  u  \right\|_{\Hminus(U)} 
\\
\leq C\left\| \nabla f \right\|_{L^{2+\delta}(U)} 
\left(  r^{\beta}+ \frac{1}{r^{3+d/2}} \mathcal{E}'(\ep)^{\frac12} \right),
\end{multline}
and
\begin{equation*}
\left\| \tfrac12 \nabla u^\ep \cdot  \a\left(\tfrac\cdot\ep\right) \nabla u^\ep - \tfrac12 \nabla  u  \cdot \ahom \nabla  u  \right\|_{W^{-2,1}(U)}
\leq
C\left\| \nabla f \right\|_{L^{2+\delta}(U)}^2
\left(  r^{\beta}+ \frac{1}{r^{3+d/2}}\mathcal{E}'(\ep)^{\frac12} \right).
\end{equation*}
\end{theorem}
\begin{proof}
The idea of the proof is to compare $u^\ep$ to the function 
\begin{equation*} \label{}
w^\ep(x):=  u (x) + \ep \zeta_r (x)  \sum_{k=1}^d \partial_k  u (x) \phi_{m,e_k} \left( \tfrac x\ep \right).
\end{equation*}
where $\zeta_r \in C^\infty_c(U)$ is a smooth cutoff function satisfying, for every $k\in\N$,
\begin{equation} 
\label{e.zetarcutoffch1}
0\leq \zeta_r \leq 1,  \quad 
\zeta_r = 1  \ \mbox{in} \ U_{2 r }, \quad 
\zeta_r \equiv 0  \  \mbox{in} \ U \setminus U_{ r }, \quad  
\left| \nabla^k \zeta_r \right| \leq C(k,d,U) r^{-k }. 
\end{equation}
Recall from~\eqref{e.Ur.def} that $U_r:= \left\{ x\in U\,:\, \dist(x,\partial U)>r \right\}$. 
The purpose of $\zeta_r$ is to cut off a boundary layer with the thickness $ r >0$, which is a given parameter. The function $w^\ep$ is a primitive version of the \emph{two-scale expansion} of $ u $.\index{two-scale expansion} The plan is to plug $w^\ep$ into the equation for $u^\ep$ and estimate the error. This will allow us to get an estimate on $\left\| u^\ep - w^\ep \right\|_{H^1(U)}$. Separately, we make direct estimates of~$\left\| \nabla w^\ep -\nabla u \right\|_{\Hminus(U)}$ and $\left\| \a\left(\tfrac\cdot\ep \right)\nabla w^\ep - \ahom\nabla u \right\|_{\Hminus(U)}$ and then conclude by the triangle inequality. 

\smallskip

Notice that if the cutoff function $\zeta_r$ was not present and the gradient of $u$ were constant (i.e., $u$ is affine), then $w^\ep$ would be an exact solution of the equation for $u^\ep$. Therefore all of the error terms will involve~$\zeta_r$ (which is only active in the boundary layer $U \setminus U_{2r}$, a set of small measure that will be ``H\"oldered away'') and higher derivatives of~$u$. 

\smallskip

Since $m$ is fixed throughout, we drop the dependence of~$\phi_{m,e_j}$ on $m$ and henceforth just write~$\phi_{e_j}$. 

\smallskip

\emph{Step 0.} Before we begin the proof, we record some standard estimates for constant-coefficient equations which are needed. First, we need the pointwise interior estimates for $\ahom$-harmonic functions (cf.~\eqref{e.pointwiseharm} in Chapter~\ref{c.regularity}) which give us, for every $k\in\N$, the existence of a constant $C(d,k,\Lambda)<\infty$ such that 
\begin{equation} 
\label{e.pointwiseharmch1}
\left\| \nabla^k u \right\|_{L^\infty(U_{r})} \leq \frac{C}{r^{k+d/2}} \left\| u - (u)_{U} \right\|_{L^2(U)} \leq \frac{C}{ r^{k+d/2}} \left\| f \right\|_{H^1(U)}. 
\end{equation}
The second estimate we need is the Meyers estimate: there exists $\delta_0(d,\Lambda)>0$ and $C(U,d,\Lambda)<\infty$ such that, if $\delta \in \left[0,\delta_0\right]$, then 
\begin{equation} 
\label{e.meyersch1}
\left\| \nabla u \right\|_{L^{2+\delta}(U)} 
\leq 
C\left\| \nabla f\right\|_{L^{2+\delta}(U)}. 
\end{equation}
This estimate is proved in Appendix~\ref{a.meyers}, see Theorem~\ref{t.Meyers appendix global}.
We may assume without loss of generality that the exponent~$\delta$ in the assumption of the theorem belongs to the interval~$(0,\delta_0]$, so that~\eqref{e.meyersch1} is in force. The only way~\eqref{e.meyersch1} enters the proof is that it allows us to estimate the $L^2$ norm of $\nabla u$ in the boundary layer. By H\"older's inequality, we have:
\begin{equation} 
\label{e.nablaubndrylayer}
\left\| \nabla u \right\|_{L^2(U\setminus U_{2 r })} \leq C \left| U \setminus U_{2 r }\right|^{\frac{\delta}{4+2\delta}} \left\| \nabla u \right\|_{L^{2+\delta} (U)} \leq C  r ^{\frac{\delta}{4+2\delta}} \left\| \nabla f \right\|_{L^{2+\delta} (U)}.
\end{equation}
Finally, we make one reduction, which is to notice that without loss of generality we may suppose that $\left( f \right)_U = 0$. Otherwise we may subtract~$\left( f \right)_U $ from each of the functions $u^\ep$,~$u$ and~$f$. By the Poincar\'e inequality, this gives us the bound
\begin{equation*} \label{}
\| f \|_{H^1(U)}
\leq C \| \nabla f \|_{L^{2}(U)} 
\leq C \| \nabla f \|_{L^{2+\delta}(U)}.
\end{equation*}

\smallskip

\emph{Step 1.} We show that
\begin{equation} 
\label{e.meltdownH}
\left\|  \nabla \cdot \left( \a\left( \tfrac \cdot\ep \right)\nabla {w}^\ep  - \ahom \nabla u\right)  \right\|_{\Hminus(U)}
\leq 
C \left\|\nabla  f \right\|_{L^{2+\delta} (U)} \left( r ^{\frac{\delta}{4+2\delta}}+ \frac{1}{ r^{3+d/2}} \mathcal{E}'(\ep)^{\frac12} \right).
\end{equation}
We compute first the gradient of $w^\ep$:
\begin{align*} 
\nabla w^\ep & = \nabla u + \zeta_r \sum_{j=1}^d \nabla \phi_{e_j}\left( \tfrac \cdot\ep \right)  \partial_j u  +  \ep \sum_{j=1}^d \phi_{e_j}\left( \tfrac \cdot \ep \right) \nabla \left( \zeta_r\partial_j u\right)  \\ 
\notag & =  \sum_{j=1}^d \left( \zeta_r \partial_j u \left( e_j   +  \nabla \phi_{e_j}\left( \tfrac \cdot\ep \right) \right) + (1-\zeta_r) \partial_j u  e_j  + 
  \ep  \phi_{e_j}\left( \tfrac \cdot \ep \right) \nabla \left( \zeta_r\partial_j u\right) \right).
\end{align*}
By the definition of $\phi_{e_j}$, we have that 
\begin{equation*} 
e_j   +  \nabla \phi_{e_j}\left( \tfrac \cdot\ep \right) = \nabla v \left( \tfrac{\cdot}{\ep} , \cu_m , e_j \right), 
\end{equation*}
and thus, since $v \left( \tfrac{\cdot}{\ep} , \cu_m , e_j \right)$ is a weak solution, we have in the weak sense the identity 
\begin{multline*} 
\nabla \cdot \left(\a\left( \tfrac \cdot\ep \right) \nabla w^\ep \right)   = 
\sum_{j=1}^d \nabla \left( \zeta_r \partial_j u \right) \cdot \a\left( \tfrac \cdot\ep \right)\left( e_j   +  \nabla \phi_{e_j}\left( \tfrac \cdot\ep \right) \right)  
\\ + \nabla \cdot \left(\a\left( \tfrac \cdot\ep \right) \left(  (1-\zeta_r) \nabla u +  \ep  \sum_{j=1}^d \phi_{e_j}\left( \tfrac \cdot \ep \right) \nabla \left( \zeta_r\partial_j u\right)  \right)  \right).
\end{multline*}
On the other hand, since $u$ is $\ahom$-harmonic,
\begin{align*} 
\nabla (\zeta_r \partial_j u ) \cdot \ahom e_j & = \nabla \cdot \left( \zeta_r  \ahom  \nabla u\right) = 
- \nabla \cdot \left( (1-\zeta_r)  \ahom  \nabla u\right)   ,
\end{align*}
and hence, in the weak sense, 
\begin{multline*} 
\nabla \cdot \left(\a\left( \tfrac \cdot\ep \right) \nabla w^\ep 
\right)   = 
\sum_{j=1}^d \nabla \left( \zeta_r \partial_j u \right) \cdot \left( \a\left( \tfrac \cdot\ep \right)\left( e_j   +  \nabla \phi_{e_j}\left( \tfrac \cdot\ep \right) \right) - \ahom e_j  \right)
\\ + \nabla \cdot \left( (1-\zeta_r)  \left( \a\left( \tfrac \cdot\ep \right) -\ahom \right)   \nabla u  \right)   
+ \nabla \cdot \left( \ep  \sum_{j=1}^d \phi_{e_j} \left( \tfrac \cdot \ep \right) \a\left( \tfrac \cdot\ep \right) \nabla \left( \zeta_r \partial_j u\right)  \right) .
\end{multline*}
It follows that 
\begin{align*} 
\lefteqn{\left\| \nabla \cdot \left( \a\left( \tfrac \cdot\ep \right)\nabla {w}^\ep - \ahom \nabla u \right) \right\|_{\Hminus(U)} } \quad &
\\ & \leq \sum_{j=1}^d \left\| \nabla \left( \zeta_r \partial_j u \right) \right\|_{W^{1,\infty}(U)} \left\| \a\left( \tfrac \cdot\ep \right)\left( e_j + \nabla \phi_{e_j}\left( \tfrac \cdot\ep \right) \right) - \ahom e_j \right\|_{\Hminus(\ep \cu_m)} 
\\ & \quad + \left\|  (1-\zeta_r)  \left( \a\left( \tfrac \cdot\ep \right) -\ahom \right)   \nabla u  \right\|_{L^2(U)} 
\\ & \quad +  \sum_{j=1}^d  \left\|  \ep  \phi_{e_j} \left( \tfrac \cdot \ep \right) \a\left( \tfrac \cdot\ep \right) \nabla \left( \zeta_r \partial_j u\right) \right\|_{L^2(U)} 
\\ & =: T_1 + T_2 + T_3. 
\end{align*}
To bound the terms on the right, observe first that,  by~\eqref{e.zetarcutoffch1} and~\eqref{e.pointwiseharmch1},
\begin{align} \label{e.nabla(zeta_r u_j)}
\left\| \nabla \left( \zeta_r \partial_j u \right) \right\|_{W^{1,\infty}(U)} 
\leq \frac{C}{ r ^{3+d/2}} \left\| f \right\|_{H^1(U)},
\end{align}
and hence, by H\"older's inequality and the definition of $\mathcal{E}'(\ep)$,
\begin{equation*} 
T_1 + T_3 \leq \frac{C}{r ^{3+d/2}} \left\| f \right\|_{H^1(U)} \mathcal{E}'(\ep)^{\frac 12}.
\end{equation*}
For the second term $T_2$, we use~\eqref{e.nablaubndrylayer} and the fact that $1-\zeta_r$ is supported in $U\setminus U_{2 r }$ to get, again by H\"older's inequality,
\begin{align*} \label{}
T_2 
\leq C \left\| \nabla u \right\|_{L^{2}(U\setminus U_{2 r } )}
\leq 
C \left| U \setminus U_{2r} \right|^{\frac{\delta}{4+2\delta}}
\left\| \nabla u \right\|_{L^{2+\delta}(U)}
\leq
C r ^{\frac{\delta}{4+2\delta}} \left\| \nabla f \right\|_{L^{2+\delta}(U)}.
\end{align*}
Combining the estimates for $T_1$, $T_2$, and $T_3$ finishes the proof of~\eqref{e.meltdownH}.
\smallskip

\emph{Step 2.} We next deduce that 
\begin{align} 
\label{e.meltdownH1}
\left\| u^\ep - w^\ep \right\|_{H^1(U)}
\leq
C \left\|\nabla  f \right\|_{L^{2+\delta} (U)} \left( r ^{\frac{\delta}{4+2\delta}}+ \frac{1}{r^{3+d/2}} \mathcal{E}'(\ep)^{\frac12} \right).
\end{align}
Indeed, testing~\eqref{e.meltdownH} with $u^\ep - w^\ep$, which belongs to $H^1_0(U)$, yields 
\begin{equation*} \label{}
\left| \int_U \nabla (u^\ep - w^\ep)(x) \cdot \a\left( \tfrac x\ep \right) \nabla w^\ep(x)\,dx \right| 
\leq \left\| u^\ep - w^\ep \right\|_{H^1(U)} \left\|  \nabla \cdot \left( \a\left( \tfrac \cdot\ep \right)\nabla {w}^\ep \right)  \right\|_{H^{-1}(U)},
\end{equation*}
and meanwhile testing the equation for $u^\ep$ with $u^\ep -w^\ep$ gives
\begin{equation*} \label{}
\int_U \nabla (u^\ep - w^\ep)(x) \cdot \a\left( \tfrac x\ep \right) \nabla u^\ep(x)\,dx = 0. 
\end{equation*}
Combining these and using the Poincar\'e inequality, we obtain
\begin{align*} \label{}
\left\| \nabla u^\ep - \nabla w^\ep \right\|_{L^2(U)}^2
&
\leq C \int_U \nabla (u^\ep - w^\ep)(x) \cdot \a\left( \tfrac x\ep \right) \nabla (u^\ep - w^\ep)(x)\,dx
\\ &
\leq C \left\| u^\ep - w^\ep \right\|_{H^1(U)} \left\|  \nabla \cdot \left( \a\left( \tfrac \cdot\ep \right)\nabla {w}^\ep \right)  \right\|_{H^{-1}(U)} 
\\ & 
\leq C \left\| \nabla u^\ep - \nabla w^\ep \right\|_{L^2(U)} \left\|  \nabla \cdot \left( \a\left( \tfrac \cdot\ep \right)\nabla {w}^\ep \right)  \right\|_{H^{-1}(U)}.
\end{align*}
Thus 
\begin{equation*} \label{}
\left\| \nabla u^\ep - \nabla w^\ep \right\|_{L^2(U)}
\leq 
C  \left\|  \nabla \cdot \left( \a\left( \tfrac \cdot\ep \right)\nabla {w}^\ep \right)  \right\|_{H^{-1}(U)},
\end{equation*}
and so we obtain~\eqref{e.meltdownH1} from~\eqref{e.meltdownH} (in view of~\eqref{e.H.vs.hatH}) and another application of the Poincar\'e inequality.

\smallskip

\emph{Step 3.} We prove that 
\begin{multline} 
\label{e.weptou.bb}
\left\| \nabla w^\ep - \nabla  u  \right\|_{\Hminus(U)}
+ \left\| \a\left(\tfrac\cdot\ep\right) \nabla w^\ep - \ahom \nabla  u  \right\|_{\Hminus(U)}
\\
\leq 
C \left\|\nabla  f \right\|_{L^{2+\delta} (U)} \left( r ^{\frac{\delta}{4+2\delta}}+ \frac{1}{ r^{2 + d/2}}\mathcal{E}'(\ep)^{\frac12} \right).
\end{multline}
Here we just have to estimate the size of the second term in the definition of~$w^\ep$, which should be small in the appropriate norms if the~$\phi_{e_k}$ are, the latter being controlled by the random variable~$\mathcal{E}(m)$. This follows from some fairly straightforward computations. 

\smallskip

We have that 
\begin{equation*} \label{}
\nabla w^\ep - \nabla  u
= \nabla \left(  \ep \zeta_r   \sum_{k=1}^d \partial_k  u  \phi_{e_k} \left( \tfrac \cdot \ep \right) \right),
\end{equation*}
and therefore 
\begin{align*}
\left\| \nabla w^\ep - \nabla  u \right\|_{\Hminus(U)}
&
\leq 
C \left\|  \ep \zeta_r   \sum_{k=1}^d \partial_k  u  \phi_{e_k} \left( \tfrac \cdot \ep \right) \right\|_{L^2(U)}
\\ & 
\leq
C \left\| \nabla u \right\|_{L^\infty(U_{ r })} \sum_{k=1}^d \ep \left\|\phi_{e_k} \left( \tfrac \cdot \ep \right) \right\|_{L^2(\ep \cu_m)}.
\end{align*}
This yields by~\eqref{e.pointwiseharmch1} that
\begin{equation*} \label{}
\left\| \nabla w^\ep - \nabla  u \right\|_{\Hminus(U)} 
\leq 
C \frac1{r^{1+d/2}} \left\| \nabla f \right\|_{L^{2+\delta}(U)} \mathcal{E}'(\ep)^{\frac12}.
\end{equation*}
We turn to the estimate for $\left\| \a\left(\tfrac\cdot\ep\right) \nabla w^\ep - \ahom \nabla  u  \right\|_{\Hminus(U)}$. 
We have
\begin{multline*} \label{}
\a\left(\tfrac \cdot\ep\right) \nabla w^\ep - \ahom \nabla u
= 
\zeta_r \sum_{j=1}^d \partial_j u \left( \a\left(\tfrac \cdot\ep\right)  \left( e_j + \nabla \phi_{e_j}\left( \tfrac \cdot\ep \right)\right)  - \ahom e_j  \right)
\\ + \left( 1-\zeta_r  \right) \left( \a\left(\tfrac \cdot\ep\right) - \ahom \right) \nabla u
+ \a\left( \tfrac \cdot\ep \right)\sum_{j=1}^d \nabla \left( \zeta_r \partial_j u \right)  \ep \phi_{e_j}\left( \tfrac \cdot\ep \right).
\end{multline*}
The $\Hminus(U)$ norm of the last two terms on the right side can be estimated analogously to~$T_2$ and~$T_3$ in Step~1. For the first term, we compute
\begin{multline*} 
\left\| \zeta_r \sum_{j=1}^d \partial_j u \left( \a\left(\tfrac \cdot\ep\right)  \left( e_j + \nabla \phi_{e_j}\left( \tfrac \cdot\ep \right)\right)  - \ahom e_j  \right) \right\|_{\Hminus(U)} \\ 
\leq C \left\|\zeta_r \nabla u\right\|_{W^{1,\infty}(U)} \left\|  \a\left(\tfrac \cdot\ep\right) \left( e_j + \nabla \phi_{e_j}\left( \tfrac \cdot\ep \right)\right)  - \ahom e_j   \right\|_{\Hminus(\ep \cu_m)} ,
\end{multline*}
and then use~\eqref{e.linktopast} and~\eqref{e.nabla(zeta_r u_j)} to obtain the desired estimate. 

\smallskip

\emph{Step 4.} The conclusion. We deduce from the triangle inequality and Steps~2 and~3 of the proof that 
\begin{align*} \label{}
\left\| \nabla u^\ep - \nabla u \right\|_{\Hminus(U)}
& 
\leq
\left\| \nabla w^\ep - \nabla u \right\|_{\Hminus(U)}
+\left\| \nabla u^\ep - \nabla w^\ep \right\|_{L^{2}(U)}
\\ &
\leq C \left\|\nabla  f \right\|_{L^{2+\delta} (U)} 
\left( r ^{\frac{\delta}{4+2\delta}} 
+ \frac{1}{r^{3+d/2}}\mathcal{E}'(\ep)^{\frac12} \right).
\end{align*}
The bound for the fluxes follows in a similar manner. The bound for the energy densities then follows from these, using the easy bound 
\begin{equation*} \label{}
\left\| \nabla u^\ep \right\|_{L^2(U)} + \left\| \nabla u \right\|_{L^2(U)}
\leq C \left\| \nabla f\right\|_{L^2(U)} 
\end{equation*}
and Lemma~\ref{l.divcurl}.
Finally, the bound for $\left\| u^\ep - u \right\|_{L^2(U)}$ is obtained from the~$\Hminus(U)$  bound for the gradients and~Lemma~\ref{l.integrate.H-1}. While the latter lemma is stated for cubes only, we can apply it here since $u_\eps - u \in H^1_0(U)$ can be extended to $\cu_0$ by setting it to be zero on $\cu_0 \setminus U$. 
\end{proof}

\begin{exercise}
\label{ex.makeEgoaway}
Show that 
\begin{equation*} \label{}
\E\left[ \inf_{r\in(0,1)} \left(  r^{\beta}+ \frac{1}{r^{3+d/2}} \mathcal{E}'(\ep)^{\frac12} \right) \right] 
\to 0 \quad \mbox{as} \ m \to \infty. 
\end{equation*}
In fact, show using Proposition~\ref{p.Em.bound} that the rate of this limit can be estimated in terms of the modulus~$\omega(m)$ defined in~\eqref{e.defomega}. 
\end{exercise}

\section*{Notes and references}

The first (qualitative) stochastic homogenization results were obtained independently, all around the same time, by Kozlov~\cite{K1}, Papanicolaou and Varadhan~\cite{PV1} and Yurinski{\u\i}~\cite{Y0}. Many of the ideas presented in this chapter originate in the variational perspective introduced by Dal Maso and Modica~\cite{DM1,DM2}, which has its roots in the work of De Giorgi and Spagnolo~\cite{DeGS}, although the arguments here are simpler since the setting is less general. In particular, the use of the subadditivity of the quantity $\nu(U,p)$ to prove qualitative homogenization in the stochastic case first appeared in~\cite{DM1,DM2}.



\chapter{Convergence of the subadditive quantities}
\label{c.two}

In the previous chapter, we introduced a subadditive quantity~$\nu(U,p)$ and showed that its convergence on large cubes implies a very general homogenization result for the Dirichlet problem. In fact, we saw that an estimate of the modulus $\omega(n)$ defined in~\eqref{e.defomega}, which represents the rate of convergence of the means $\E\left[\nu(\cu_n,p)\right]$ to their limit~$\overline{\nu}(p)$, would imply an explicit convergence rate for the homogenization limits. The goal of this chapter is to obtain an estimate showing~$\omega(n)$ is at most a negative power of the length scale:~$\omega(n) \leq C 3^{-n\alpha}$ for an exponent~$\alpha(d,\Lambda)>0$.

\smallskip

Unfortunately, the argument from the previous chapter for the convergence of the limit $\lim_{n\to \infty} \E\left[\nu(\cu_n,p)\right]=\overline{\nu}(p)$ provides no clues for how we should obtain a convergence rate. Indeed, the limit was obtained from the monotonicity of the sequence $n\mapsto \E\left[\nu(\cu_n,p)\right]$, and obviously a monotone sequence may converge at any speed. This is a central issue one must face when trying to quantify a limit obtained from the subadditive ergodic theorem. 

\smallskip

To estimate the speed of convergence of a subadditive quantity to its limit, it is natural to search for a \emph{superadditive} quantity which is close to the original subadditive one. If a good estimate of the difference between the two can be obtained, then an estimate for the speed of convergence of both quantities would follow since each gives us monotonicity on one side, and the limit would thus be squeezed in between. In this chapter we accomplish this by introducing a new subadditive ``dual'' quantity which we denote by~$\nu^*(U,q)$. We then compare the family of subadditive quantities $\nu(U,p)$ indexed by $p\in\Rd$ to the family of superadditive quantities~$p\cdot q - \nu^*(U,q)$ indexed by~$q\in\Rd$. We will show that, for each $p$, there is a~$q$ such that these two quantities are close in expectation (in fact, the correspondence between $p$ and $q$ is linear and gives us another way to define the homogenized coefficients~$\ahom$). An equivalent way of saying this in the terminology of convex analysis is that we show that \emph{$\nu(U,\cdot)$ and $\nu^*(U,\cdot)$ are close to being convex dual functions.} We will prove that this is so by constructing an iteration scheme which establishes the decay of the defect in their convex dual relationship.

\smallskip


We begin in the first section by giving the definition of~$\nu^*$ and exploring its relationship to~$\nu$. The iteration procedure giving the convergence rate for~$\omega(n)$ is then carried out in Section~\ref{s.iteration}. At the end of that section, we upgrade the stochastic integrability of our estimates, making use of subadditivity and using independence in a stronger way than we did in the previous chapter. This allows us to obtain very strong bounds on the tail of the distribution of the random variable $\mathcal{E}(m)$ encountered in Proposition~\ref{p.qualweakconv} and Theorem~\ref{t.DP.blackbox}, thereby giving us useful quantitative versions of these qualitative homogenization results: see Theorems~\ref{t.Jmaximizers} and~\ref{t.DP} in Section~\ref{s.Emwhipped}. 

\section{The dual subadditive quantity \texorpdfstring{$\nu^*$}{nu*}}
\label{s.dualsubadd}

\index{subadditive quantity!$\nu^*$}

The dual subadditive quantity is defined for each $q\in\Rd$ and bounded Lipschitz domain $U\subseteq \Rd$ by
\begin{equation}
\label{e.def.nu*}
\nu^*(U,q):= \sup_{u\in H^1(U)} \fint_{U} \left( -\frac12 \nabla u \cdot \a\nabla u + q\cdot \nabla u \right).
\end{equation}
Note that we impose no restriction on the trace of~$u$ on the boundary of~$U$ in~\eqref{e.def.nu*}. 
It is routine to check that the supremum on the right of~\eqref{e.def.nu*} is finite and that a maximizer exists which is unique up to an additive constant. Although we make no explicit use of this fact, we mention that it is characterized as the weak solution of the Neumann problem
\begin{equation} \label{e.Neumann}
\left\{
\begin{aligned}
& -\nabla \cdot \left( \a\nabla u \right) = 0 & \mbox{in} & \ U, \\
& \mathbf{n} \cdot \a\nabla u = \mathbf{n} \cdot q  & \mbox{on} & \ \partial U,
\end{aligned}
\right.
\end{equation}
where $\mathbf{n}$ denotes the outward-pointing unit normal vector to $\partial U$. Recall that we denote by $\A(U)$ the linear space of $\a$-harmonic functions in $U$, see \eqref{e.def.A(U)}. Thus, since the maximizer belongs to $\A(U)$, we may restrict the supremum in \eqref{e.def.nu*} to~$u \in \A(U)$ without changing the value of $\nu^*(U,q)$: 
\begin{equation} 
\label{e.def.nu*2}
\nu^*(U,q) = \sup_{u\in \A(U)} \fint_{U} \left( -\frac12 \nabla u \cdot \a\nabla u + q\cdot \nabla u \right).
\end{equation}
We see immediately that this quantity is subadditive by taking a partition $\{U_j\}$ of $U$ and testing the definition of each $\nu^*(U_j,q)$ with the restrictions to $U_j$ of the maximizer of $\nu(U,q)$. See the proof of Lemma~\ref{l.basicJ} below for details.

\smallskip

\index{Legendre transform|(}
The definition of~\eqref{e.def.nu*} itself should remind the reader of the Legendre transform. Recall that if $L(p)$ is a uniformly convex function of $p$, then the \emph{Legendre transform} 
of $L$ is the function $\displaystyle{L^{*}(q)}$ defined by
\begin{equation*} \label{}
L^*(q) := \sup_{p\in\Rd} \left( -L(p) + p\cdot q\right). 
\end{equation*}
Notice that~\eqref{e.def.nu*} resembles the formula for $L^*$, but there is an integral and the supremum is over~$H^1(U)$ instead of~$\Rd$. We can compare $\nu^*(U,q)$ to the other subadditive quantity $\nu(U,p)$ defined in the previous chapter (see~\eqref{e.def.nu}) by testing the minimizer of $\nu(U,p)$ in the definition of $\nu^*(U,q)$. This gives
\begin{equation} 
\label{e.obvioustesting}
\nu(U,p) + \nu^*(U,q) \geq p\cdot q. 
\end{equation}
Indeed, letting $v(\cdot,U,p)$ denote the minimizer in the definition of~$\nu(U,p)$, as in~\eqref{e.vdeff}, and using that $p = \fint_U \nabla v(\cdot,U,p)$ since $v(\cdot,U,p)\in \ell_p+ H^1_0(U)$, we find that 
\begin{align*} \label{}
\nu^*(U,q) 
\geq 
\fint_U \left( -\frac12 \nabla v(\cdot,U,p) \cdot \a \nabla v(\cdot,U,p) + q\cdot \nabla v(\cdot,U,p) \right)
= -\nu(U,p) + q\cdot p.
\end{align*}

\smallskip

The inequality~\eqref{e.obvioustesting} reminds us of convexity duality, for it is equivalent to the statement that $\nu^*(U,q)$ and $\nu(U,p)$ are each bounded below by the Legendre-Fenchel transform of the other. We emphasize however that, despite our notation, $\nu^*(U,\cdot)$ is \emph{not} in general equal to the Legendre-Fenchel transform of $\nu(U,\cdot)$. Indeed, the set of solutions of the Neuman problem~\eqref{e.Neumann} with affine data is not in general equal to the set of solutions of the Dirichlet problem~\eqref{e.DP} with affine data (unless the coefficients are constant). Therefore, a precise convex dual relationship should be expected to hold only in the large-scale limit. We do expect such duality to hold in the large-scale limit: to see why, imagine that for a fixed $q$, the maximizer of $\nu^*(U,q)$ is relatively ``flat'' for a large domain~$U$ (meaning that it is well-approximated by an affine function in a quantitative sense). In that case, the inequality~\eqref{e.obvioustesting} should be nearly an equality for the right choice of $p$ (which should be the slope of the affine function approximating $v(\cdot,U,0,q)$). Furthermore, since $-\nu^*(U,q) + q\cdot p$ is \emph{superadditive} and lower bounds $\nu(U,p)$, its monotonicity complements that of $\nu(U,p)$. This gives us a chance to quantify the convergence rate for both quantities by studying the sharpness of the inequality~\eqref{e.obvioustesting}. In view of \eqref{e.defnubar}, it also suggests that 
\begin{equation}
\label{e.nu*.conv}
\lim_{n \to \infty} \E[\nu^*(\cu_n,q)] = \frac 1 2 q \cdot \ahom^{-1} q.
\end{equation}

\index{Legendre transform|)}

\smallskip

These considerations motivate us to study the quantity \index{subadditive quantity!$J$}
\begin{equation} 
\label{e.Jsplitting}
J(U,p,q):= \nu(U,p)+\nu^*(U,q) - p \cdot q,
\end{equation}
which monitors the ``defect'' in the convex duality relationship between $\nu$ and $\nu^*$. Notice that $J(U,p,q)$ contains both $\nu$ and $\nu^*$ separately since $J(U,p,0) = \nu(U,p)$ and $J(U,0,q) =\nu^*(U,q)$, so we can think of it as the ``master'' subadditive quantity. We expect the following limit to hold:
\begin{equation} 
\label{e.Jlimit}
\lim_{n\to \infty} J(\cu_n,p,\ahom p) = 0. 
\end{equation}
The goal is to prove~\eqref{e.Jlimit} with an explicit convergence rate for the limit and to show that this gives a convergence rate for $\nu(\cu_n,p)$ and $\nu^*(\cu_n,q)$ separately. The main result is stated in Theorem~\ref{t.subadd}, below.  

\smallskip

We begin our analysis by giving a variational representation for $J(U,p,q)$. As a first step, we look for a formula for~$\nu(U,p)$ which resembles that of~$\nu^*(U,q)$ in~\eqref{e.def.nu*2}. In fact, as we will prove below in Lemma~\ref{l.Jsplitting}, we have
\begin{equation}
\label{e.formulafornuUp}
\nu(U,p) =  \sup_{u\in \A(U)} \fint_{U} \left( -\frac12 \nabla u \cdot \a\nabla u + p\cdot \a \nabla u \right).
\end{equation}
Notice that the comparison between this expression and~\eqref{e.def.nu*2} provides compelling evidence that~$\nu^*(U,q)$ is a natural quantity to consider. 

\smallskip

The formula~\eqref{e.formulafornuUp} can be derived from a geometric (or least-squares) interpretation of~$\nu(U,p)$, which we briefly summarize as follows. Recall from \eqref{e.def.Lp0} that the space of gradients of functions in $H^1_0(U)$ is denoted by $\Lpoot(U)$. The gradient of the minimizer for $\nu(U,p)$ is the element of $p + \Lpoot(U)$ with smallest ``$\a(\cdot)$-norm'', where the $\a(\cdot)$-norm is the square root of the quadratic form
\begin{equation*}  
\nabla v \mapsto  \fint_U \frac 1 2 \nabla v \cdot \a \nabla v.
\end{equation*}
It is therefore the orthogonal projection of the constant vector field $p$ onto the orthogonal complement of $\Lpoot(U)$ relative to~$\Lpot(U)$, where the notion of orthogonality is defined according to the inner product  implicit in the quadratic form above. This linear space is precisely the set~$\nabla (\A(U))$ of gradients of $\a$-harmonic functions in $U$, see \eqref{e.def.A(U)}. Hence, the gradient of the minimizer for $\nu(U,p)$ is also the projection of $p$ onto $\nabla(\A(U))$ with respect to the $\a(\cdot)$-inner product, that is, the element of~$\nabla(\A(U))$ which  minimizes the $\a(\cdot)$-distance to $p$. It can therefore be written as the maximizer over $v\in \A(U)$ of the functional 
\begin{equation*}  
v \mapsto - \fint_U \frac 1 2 (\nabla v - p) \cdot \a (\nabla v - p) = \fint_U \Ll( -\frac 1 2 \nabla v \cdot \a \nabla v + p \cdot \a \nabla v - \frac 1 2 p \cdot \a p \Rr).
\end{equation*}
Since the term~$\frac 1 2 p\cdot \a p$ on the right side does not depend on~$v$, it can be dropped without changing the identity of the maximizer. It turns out that the value of the minimum of the resulting functional coincides with $\nu(U,p)$.

\smallskip

We next give the formal details of the argument sketched above and actually prove a general formula, valid for the quantity~$J$ defined in~\eqref{e.Jsplitting}, which combines both~\eqref{e.def.nu*2} and~\eqref{e.formulafornuUp}.

\begin{lemma}
\label{l.Jsplitting}
For every bounded Lipschitz domain $U\subseteq \Rd$ and $p,q\in\Rd$,
\begin{equation} 
\label{e.formulafornuJpq}
J(U,p,q)= \sup_{w\in \A(U)} \fint_U \left( -\frac12 \nabla w \cdot \a\nabla w - p\cdot \a\nabla w + q\cdot \nabla w \right).
\end{equation}
Moreover, the maximizer $v(\cdot,U,p,q)$ is the difference between the maximizer of $\nu^*(U,q)$ in \eqref{e.def.nu*} and the minimizer of $\nu(U,p)$ in \eqref{e.def.nu}.
\end{lemma}
\begin{proof}
Let $v \in \ell_p + H^1_0(U)$ denote the minimizer in the definition of $\nu(U,p)$. For every $u \in \A(U)$, we have
\begin{multline}
\label{e.mu.nu.nu}
\nu(U,p) + \fint_U  \left( -\frac 1 2 \nabla u \cdot \a\nabla u + q\cdot \nabla u \right) - p \cdot q \\
 = \fint_U \Ll( \frac 1 2 \nabla v \cdot \a \nabla v - \frac 1 2 \nabla u \cdot \a\nabla u + q\cdot \nabla u \right) - p \cdot q.
\end{multline} 
Since $v \in \ell_p + H^1_0(U)$, we have 
\begin{equation*}  
p  = \fint_U \nabla v.
\end{equation*}
Since $u  \in \A(U)$, we also have
\begin{equation*}  
\fint_U \nabla u \cdot \a \nabla v = \fint_U \nabla u \cdot \a p,
\end{equation*}
and this last identity holds true in particular for~$u = v$. 
We deduce that the left side of \eqref{e.mu.nu.nu} equals
\begin{equation*}  
\fint_U \Ll( -\frac 1 2 \Ll( \nabla u - \nabla v \Rr) \cdot \a \Ll( \nabla u - \nabla v \Rr) 
- p \cdot \a \Ll( \nabla u - \nabla v \Rr) + q \cdot \Ll( \nabla u - \nabla v \Rr) \Rr) .
\end{equation*}
Comparing this result with~\eqref{e.def.nu*2},~\eqref{e.Jsplitting} and the right side of~\eqref{e.formulafornuJpq}, we obtain the announced result.
\end{proof}

We next collect some further properties of~$J$ in the following lemma, which extends Lemma~\ref{l.basicnu}. 

\begin{lemma}[Basic properties of $J$]
\label{l.basicJ}
Fix a bounded Lipschitz domain $U\subseteq \Rd$. 
The quantity $J(U,p,q)$ and its maximizer $v(\cdot,U,p,q)$ satisfy the following properties:

\begin{itemize}

\item \emph{Representation as quadratic form.} The mapping $(p,q) \mapsto J(U,p,q)$ is a quadratic form and there exist matrices $\a(U)$ and $\a_*(U)$ such that 
\begin{equation}
\label{e.aastar.bounds}
\Id \leq  \a_*(U)  \leq \a(U) \leq \Lambda \Id 
\end{equation}
and
\begin{equation}
\label{e.Jrepresentation}
J(U,p,q) = \frac12p\cdot \a(U) p + \frac12q\cdot \a_*^{-1}(U) q - p\cdot q.
\end{equation}
The matrices $\a(U)$ and $\a_*(U)$ are characterized by the following relations, valid for every $p,q\in\Rd$:
\begin{equation} 
\label{e.Jderp}
\a(U) p  = -\fint_U \a \nabla v(\cdot,U,p,0),
\end{equation}
\begin{equation} 
\label{e.Jderq}
\a_*^{-1}(U)q  =  \fint_U  \nabla v(\cdot,U,0,q).
\end{equation}


\item \emph{Subadditivity.} Let $U_1, \ldots, U_N \subset U$ be bounded Lipschitz domains that form a partition of~$U$, in the sense that $U_i \cap U_j = \emptyset$ if $i\neq j$ and 
\begin{equation*} \label{}
\left| U \setminus \bigcup_{i=1}^N U_i \right| = 0. 
\end{equation*}
Then, for every $p,q\in\Rd$,
\begin{equation}
\label{e.Jsubadd}
J(U,p,q) \leq \sum_{i=1}^N \frac{\left|U_i\right|}{|U|} J(U_i,p,q). 
\end{equation}

\item \emph{First variation for $J$.} For $p,q\in\Rd$, the function $v(\cdot,U,p,q)$ is characterized as the unique element of $\A(U)$ which satisfies 
\begin{equation}
\label{e.Jfirstvar} 
\fint_U \nabla w \cdot \a\nabla v(\cdot,U,p,q) = \fint_U \left( - p\cdot \a\nabla w + q\cdot \nabla w \right),
\quad 
\forall\, w \in \A(U). 
\end{equation}

\item \emph{Quadratic response.} For every $p,q\in\Rd$ and $w\in \A(U)$, 
\begin{multline} 
\label{e.Jquadresponse}
\fint_U \frac12 \left(  \nabla w - \nabla v(\cdot,U,p,q)\right) \cdot\a\left(  \nabla w - \nabla v(\cdot,U,p,q)\right)
\\
=
J(U,p,q) - \fint_U \left( -\frac12\nabla w \cdot \a\nabla w  -p\cdot\a \nabla w + q\cdot \nabla w  \right).
\end{multline}
\end{itemize}
\end{lemma}
\begin{remark}
\label{r.coarsened}
From \eqref{e.Jsplitting} and \eqref{e.Jrepresentation}, we see that for every $q \in \Rd$,
\begin{equation*}  
\nu^*(U,q) = \frac 12 q \cdot \a_*^{-1}(U) q.
\end{equation*}
The choice of notation for $\a_*(U)$ is justified by the fact that we expect the quantity $\nu^*(U,q)$ to converge to $\frac 1 2 q \cdot \ahom^{-1} q$, see \eqref{e.nu*.conv}. We consider~$\a(U)$ and $\a_*(U)$ to be ``best guess approximations for~$\ahom$ by looking only at~$\a(\cdot)$ restricted to~$U$,'' from the point of view of the quantities $\nu$ and~$\nu^*$ respectively. The main point of the arguments in this chapter is to show that these two notions are almost the same for $U=\cu_n$ and large~$n\in\N$ and that the matrices~$\a(\cu_n)$ and~$\a_*(\cu_n)$ are both close to~$\ahom$: see Theorem~\ref{t.subadd} below. The matrices $\a(U)$ and $\a_*(U)$ should be thought of as ``coarsened coefficients'' with respect to~$U$.
\index{coarsened coefficient field}
\end{remark}
\begin{proof}[Proof of Lemma~\ref{l.basicJ}]
\emph{Step 1.} We derive the first variation and prove the quadratic response~\eqref{e.Jquadresponse} inequalities. Fix $p,q\in\Rd$ and set $v:= v(\cdot,U,p,q)$. Also fix~$w\in \A(U)$ and~$t\in\R$ test the expression~\eqref{e.Jsplitting} for $J(U,p,q)$ with $\tilde{v}_t: = v + tw$:
\begin{align*}
\lefteqn{
\fint_U \left( -\frac12 \nabla v \cdot\a\nabla v -p\cdot\a\nabla v+q\cdot \nabla v \right) 
} \qquad & \\ &
\geq 
\fint_U \left( -\frac12 \nabla \tilde{v}_t \cdot\a\nabla \tilde{v}_t -p\cdot\a\nabla \tilde{v}_t+q\cdot \nabla \tilde{v}_t \right) 
\\ &
=
 \fint_U \left( -\frac12 \nabla v \cdot\a\nabla v -p\cdot\a\nabla v+q\cdot \nabla v \right) 
+ t \fint_U \left( -\nabla v \cdot \a\nabla w -p\cdot\a\nabla w+q\cdot \nabla w \right) 
\\ & \qquad 
+ t^2 \fint_U -\frac12 \nabla w \cdot\a\nabla w.
\end{align*}
Rearranging this inequality, dividing by $t$ and sending $t\to 0$ yields 
 \begin{equation*} \label{}
\fint_U \left( -\nabla v \cdot \a\nabla w -p\cdot\a\nabla w+q\cdot \nabla w \right) =0. 
\end{equation*}
This is~\eqref{e.Jfirstvar}. This calculation also shows that any function satisfying~\eqref{e.Jfirstvar} will be a maximizer for $J(U,p,q)$. In fact, taking $t=1$, writing $\tilde{v}:=\tilde{v}_1$ and using the previous two displays, we see that 
\begin{equation*} \label{}
J (U,p,q)  - \fint_U \left( -\frac12 \nabla \tilde{v} \cdot\a\nabla \tilde{v} -p\cdot\a\nabla \tilde{v}+q\cdot \nabla \tilde{v} \right) 
=  \fint_U \frac12 \nabla w \cdot\a\nabla w.
\end{equation*}
This is~\eqref{e.Jquadresponse}, which yields, in particular, that the maximizer is unique. By~\eqref{e.Jfirstvar}, we have that 
\begin{equation}
\label{e.Jmaximizerlinear}
(p,q) \mapsto v(\cdot,U,p,q) \quad \mbox{is linear.}
\end{equation}

\smallskip

\emph{Step 2.} We prove that $J(U,p,q)$ is quadratic and establish the representation formula~\eqref{e.Jrepresentation} for matrices $\a(U)$ and $\a_*(U)$ given by~\eqref{e.Jderp} and~\eqref{e.Jderq}, respectively. First, we observe that we have the following identity:
\begin{equation}
 \label{e.firstvarJ2}
J(U,p,q) = \fint_U \frac12 \nabla v(\cdot,U,p,q)\cdot \a  \nabla v(\cdot,U,p,q).
\end{equation}
Indeed, this is immediate from~\eqref{e.formulafornuJpq} and~\eqref{e.Jfirstvar}. This together with~\eqref{e.Jmaximizerlinear} implies
\begin{equation} 
\label{e.Jisquadratic}
(p,q)\mapsto J(U,p,q)\quad \mbox{is quadratic.}
\end{equation}
In view of the formula~\eqref{e.Jsplitting}, this implies in particular that $q \mapsto \nu^*(U,q)$ is quadratic. We define $\a_*(U)$ to be the symmetric matrix such that for every $q \in \Rd$,
\begin{equation*}  
\nu^*(U,q) = \frac 1 2 q \cdot \a_*^{-1}(U)q.
\end{equation*}
The representation \eqref{e.Jrepresentation} is then immediate from \eqref{e.Jsplitting}. Moreover, for every $q \in \Rd$,
\begin{equation*}  
q \cdot \a_*^{-1}(U) q = \fint_U \nabla v(\cdot,U,0,q) \cdot \a \nabla v (\cdot,U,0,q). 
\end{equation*}
This and \eqref{e.Jmaximizerlinear} imply that, for every $p,q \in \Rd$,
\begin{equation*}  
p \cdot \a_*^{-1}(U) q = \fint_U \nabla v(\cdot,U,0,p) \cdot \a \nabla v (\cdot,U,0,q).
\end{equation*}
By \eqref{e.Jfirstvar}, we deduce that
\begin{equation*}  
p \cdot \a_*^{-1}(U) q = \fint_U p \cdot \nabla v (\cdot,U,0,q),
\end{equation*}
and this is \eqref{e.Jderq}. The relation \eqref{e.Jderp} follows from \eqref{e.nufluxidentity} and Lemma~\ref{l.Jsplitting}.

\smallskip

%

\smallskip

\emph{Step 3.}  We prove the bounds~\eqref{e.aastar.bounds}. The upper bound for $\a(U)$ was proved already in Lemma~\ref{l.basicnu}. The bound $\a_*(U) \leq \a(U)$ is immediate from~\eqref{e.obvioustesting} and~\eqref{e.Jrepresentation} by taking $q = \a_*(U) p$ in~\eqref{e.Jrepresentation}.  
To obtain the lower bound for $\a_*(U)$, we first use \eqref{e.Jderq} to write
\begin{equation*}
\nu^*(U,q) = \frac 12 q \cdot \a_*^{-1}(U) q = \frac12 q\cdot \fint_U  \nabla v(\cdot,U,0,q).
\end{equation*}
By Young's inequality, 
\begin{equation*}
\fint_U q\cdot  \nabla v(\cdot,U,0,q)
\leq 
\fint_U \left( \frac12 q\cdot \a^{-1} q + \frac12 \nabla v(\cdot,U,0,q)\cdot \a \nabla v(\cdot,U,0,q) \right).
\end{equation*}
and thus, by \eqref{e.firstvarJ2},
\begin{equation*}  
\frac 12 q \cdot \a^{-1}(U) q \le \frac12 q\cdot \left( \fint_U \a^{-1} \right) q
\leq \frac12 \left| q \right|^2,
\end{equation*}
which gives the desired lower bound for $\a_*(U)$. 

\smallskip

\emph{Step 4.} We verify the subadditivity of~$J$. Testing~\eqref{e.Jsplitting} for~$J(U_i,p,q)$ with $v=v(\cdot,U,p,q)$, we get
\begin{equation*} \label{}
J(U_i,p,q) \geq \fint_{U_i} \left(  -\frac12 \nabla v \cdot\a\nabla v -p\cdot\a\nabla v+q\cdot \nabla v \right). 
\end{equation*}
Multiplying this by $|U_i|/|U|$ and summing over $i$, we get
\begin{align*} \label{}
J(U,p,q) & = \fint_U \left(  -\frac12 \nabla v \cdot\a\nabla v -p\cdot\a\nabla v+q\cdot \nabla v \right) \\
& = \sum_{i=1}^N \frac{|U_i|}{|U|} \fint_{U_i} \left(  -\frac12 \nabla v \cdot\a\nabla v -p\cdot\a\nabla v+q\cdot \nabla v \right) \\
& \leq \sum_{i=1}^N \frac{|U_i|}{|U|} J(U_i,p,q).
\end{align*}
This completes the proof of the lemma. 
\end{proof}

The main result of this chapter is the following theorem, which gives quantifies convergence of~$\a(\cu_n)$ and~$\a_*(\cu_n)$ to~$\ahom$. We remind the reader that the $\O_s(\cdot)$ notation is defined in~\eqref{e.Os.not} and its basic properties are given in Appendix~\ref{a.bigO}. 
\index{subadditive quantity!convergence of}

\begin{theorem}
\label{t.subadd}
Fix $s\in (0,d)$. 
There exist $\alpha(d,\Lambda)\in \left(0,\frac12\right]$ and $C(s,d,\Lambda)<\infty$ such that, for every $n\in\N$,  
\begin{equation}
\label{e.subadderror}
\left| \a(\cu_n) - \ahom \right| + \left| \a_*(\cu_n) - \ahom \right|
\leq C 3^{-n \alpha(d-s)} + \O_1\left( C3^{-ns} \right). 
\end{equation}
\end{theorem}

The right side of~\eqref{e.subadderror} conveniently breaks the error into a deterministic part with a relatively small exponent of $\alpha(d-s)$, and a random part with a relatively large exponent $s<d$ (we will typically take $s$ very close to $d$). Thus, while we are controlling the total size of the error in a very mild way, we are controlling the tail of the distribution of the error very strongly (in fact, as strongly as is possible, see Remark~\ref{r.Od.optimal} below). In this sense,~\eqref{e.subadderror} is close to being a deterministic estimate.

\begin{remark}[{Optimality of stochastic integrability}]
\label{r.Od.optimal}
By Chebyshev's inequality, we may deduce from~\eqref{e.subadderror} that, for every $t>0$,
\begin{equation*} 
\P \left[ 
\left| \a(\cu_n) - \ahom \right| + \left| \a_*(\cu_n) - \ahom \right|
\geq  C 3^{-n \alpha(d-s)} + t 
\right]
\leq
C\exp\left( -3^{ns} t \right). 
\end{equation*}
In particular, we find that the probability of finding a very large error of size~$O(1)$ is smaller than~$O(\exp\left( - c3^{ns} \right))$ for every $s<d$. This estimate is optimal, in the sense that the probability of this event cannot be made smaller than $O\left( \exp\left( -A \, 3^{nd}  \right)\right)$ for a large constant~$A\gg1$. To see why this is so, we think of the random checkerboard example, in which ``white'' squares have $\a(x)\equiv \Id$ and ``black'' squares have $\a(x)\equiv 2\Id$, and a fair coin is tossed independently to determine whether each square $z+ \cu_0$ with $z\in\Zd$ is white or black. The probability of having all the squares in $\cu_n$ be of a same given color is exactly 
\begin{equation*} \label{}
\left(\frac{1}{2}\right)^{3^{nd}} = \exp\left( -c\, 3^{nd} \right).
\end{equation*}
These two outcomes (all white versus all black) will necessarily have $J(\cu,p,q)$'s which are $O(1)$ apart for $p,q\in B_1$. So the probability of an error of size~$O(1)$ must be at least $O\left( \exp\left( -c \, 3^{nd}  \right)\right)$ in general.
\end{remark}

\begin{remark} 
\label{r.arbitrary}
The choice of the stochastic integrability exponent $1$ in the term $\O_1(C 3^{-ns})$ appearing in \eqref{e.subadderror} is somewhat arbitrary; the only important point is that the range of allowed values for the product of this exponent and $s$ is $(0,d)$. Indeed, since for every $p,q \in B_1$ and $n \in \N$, we have $J(\cu_n,p,q) \le C(d,\Lambda)$, an application of Lemma~\ref{l.change-s}
yields that for every $r \ge 1$, the right side of \eqref{e.subadderror} can be replaced by
\begin{equation}  
\label{e.alternative.rhs}
C 3^{-n\al(d-s)} + \O_r \Ll( C 3^{-ns/r} \Rr) .
\end{equation}
Conversely, if for some $r\ge 1$ and $\eps, \theta > 0$, a random variable $X$ satisfies
\begin{equation*}  
X \le \ep + \O_r \Ll( \theta \Rr) ,
\end{equation*}
then by Young's inequality, for every $\rho > 0$, 
\begin{equation}  
\label{e.convese.of.change-s}
X \le \rho + \rho^{1-r} X^r \le \rho + \rho^{1-r} \ep + \O_1 \Ll( \rho^{1-r} \theta^r \Rr) . 
\end{equation}
In particular, applying this result with $\rho = 3^{-\beta (d-s) n}$, we get
\begin{multline}  
\label{e.impl.converse.of}
X \le  C 3^{-n\al(d-s)} + \O_r \Ll(  C 3^{-ns/r} \Rr) \\
\implies  X \le  C 3^{-\beta n(d-s)} + C 3^{-n(\al-(r-1)\be)(d-s)} + \O_1 \Ll(  C 3^{-ns + (r-1)\be n(d-s) } \Rr).
\end{multline}
Hence, if Theorem~\ref{t.subadd} holds with the right side of \eqref{e.subadderror} replaced by \eqref{e.alternative.rhs}, we recover the original formulation of this theorem by applying \eqref{e.impl.converse.of} with $\beta(\alpha,r) > 0$ sufficiently small.
\end{remark}

In view of~\eqref{e.Jsplitting}, the estimate~\eqref{e.subadderror} can be thought of as two separate estimates, one for each of the subadditive quantities $\nu$ and $\nu^*$. However, it is better to think of it as just one estimate for the ``convex duality defect'' $J(U,p,\ahom p)$, as we discussed above, because this is how we prove it. Since this simple observation is central to our strategy, we formalize it next. 

\begin{lemma}
\label{l.minimalset}
There exists a constant~$C(d,\Lambda)<\infty$ such that, for every symmetric matrix $\tilde{\a}\in \R^{d\times d}$ satisfying
\begin{equation*} \label{}
 \Id \leq \tilde{\a} \leq \Lambda \Id
\end{equation*}
and every bounded Lipschitz domain $U\subseteq\Rd$, we have 
\begin{equation}
\label{e.minimalset}
\left| \a(U)   -  \tilde{\a}\right| + \left| \a_*(U)   -  \tilde{\a}\right|
 \leq C \sup_{p \in B_1} \Ll(J (U,p,\tilde{\a}p)\Rr)^\frac 1 2 . 
\end{equation}
\end{lemma}
\begin{proof}
Denote
\begin{equation*}
\delta := \sup_{p\in B_1} J(U,p,\tilde{\a} p).
\end{equation*}
%
%
In view of~\eqref{e.Jrepresentation}, for each fixed $p\in\Rd$, the uniformly convex function
\begin{equation}
\label{e.qmap}
q\mapsto J(U,p,q) = \frac12 p \cdot \a(U) p + \frac12 q \cdot \a_*^{-1}(U) q - p\cdot q
\end{equation}
achieves its minimum at the point $\tilde{q} :=  \a_*(U) p$. By~\eqref{e.aastar.bounds} and~\eqref{e.Jrepresentation}, for every $q \in \Rd$, we have that
\begin{multline*} 
 \frac12 (q-\tilde{q}) \cdot \a_*^{-1}(U) (q-\tilde{q}) 
 \\ = \frac12 p \cdot \a(U) p + \frac12 q \cdot \a_*^{-1}(U) q - p\cdot q  + \frac12 p\cdot (\a_*(U) - \a(U)) p \leq J(U,p,q).
\end{multline*}
Choosing $q= \tilde{\a}  p$, we deduce, after taking the supremum over~$p \in B_1$ and recalling the definition of~$\delta$,  that 
\begin{equation*}
\left| \a_*(U)   -  \tilde{\a}\right|^2  \leq C\delta  .
\end{equation*}
A similar argument, fixing $q \in \Rd$ and minimizing $p \mapsto J(U,p,q)$, gives that 
\begin{equation*}
\left| \a(U)   -  \tilde{\a}\right|^2  \leq C\delta ,
\end{equation*}
and hence~\eqref{e.minimalset} follows, completing the proof. 
%
%
\end{proof}
%

\begin{exercise}
\label{ex.voigt-reiss2}
Recall the~\emph{Voigt-Reiss} estimates\index{Voigt-Reiss bounds} for the homogenized coefficients:
\begin{equation} 
\label{e.voigt-reiss.again}
\E \left[  \int_{\cu_0} \a^{-1} (x) \,dx\right]^{-1}  
\leq \ahom 
\leq \E \left[ \int_{\cu_0} \a(x)\,dx \right].  
\end{equation}
In Exercise~\ref{ex.voigt-reiss1}, we proved the second inequality of~\eqref{e.voigt-reiss.again}. Assuming that~\eqref{e.nu*.conv} holds, use properties of $\nu^*$ and the inequality
\begin{equation*} \label{}
 \frac12p\cdot \a(x) p + \frac12 q\cdot \a^{-1}(x)q \geq p\cdot q
\end{equation*}
to prove the first inequality of~\eqref{e.voigt-reiss.again}. 
\end{exercise}

\begin{exercise}
\label{ex.2dgeometricmean}
Suppose~$d = 2$ and let~$R \in \R^{2 \times 2}$ denote the rotation by~$\frac \pi 2$ about the origin: 
\begin{equation*}  
R = 
\begin{pmatrix}  
0 & 1 \\ -1 & 0
\end{pmatrix}
.
\end{equation*}
Assume that there exists a constant $\sigma > 0$ such that the random field~$x \mapsto \a(x)$ has the same law as the field~$x \mapsto \sigma \a^{-1}(Rx)$. The purpose of this exercise is to show, in this case, that~$\ahom$ is given explicitly by the following formula (called \emph{Dykhne formula}~\cite{dykhne}): 
\index{Dykhne formula}
\begin{equation}
\label{e.Dykhne}
\ahom = \sqrt{\sigma} \, I_2.
\end{equation}
Here is an outline of the proof:
\begin{enumerate}  
\item 
Recalling the definifion of the space $\Ls(U)$ in \eqref{e.Ls.not} and \eqref{e.formulafornuUp}, show that
\begin{equation*}  
\nu(U,p) = \sup_{\g \in \Ls(U)} \int_U \Ll( -\frac 1 2 \g \cdot \a^{-1} \g + p \cdot \g \Rr) .
\end{equation*}

\item Show that $\g \in \Ls(U)$ if and only if there exists $u \in H^1(U)$ such that 
\begin{equation*}  
\g = 
\begin{pmatrix}  
\partial_{x_2} u \\ - \partial_{x_1} u
\end{pmatrix}
=: \nabla^\perp u.
\end{equation*}

\item
Using that $x \mapsto \sigma \a^{-1}(Rx)$ has the same law as $x \mapsto \a(x)$ and \eqref{e.def.nu*}, show that for every $n \in \N$,
\begin{equation*}  
\E \Ll[\nu(\cu_n,p)\Rr] = \sigma^{-1} \E \Ll[ \nu^*(\cu_n,\sigma p) \Rr].
\end{equation*}
Assuming \eqref{e.nu*.conv}, show that \eqref{e.Dykhne} holds. Deduce that for the $2$-dimensional random checkerboard displayed on Figure~\ref{f.checkerboard} (page \pageref{f.checkerboard}) with $\a_0 = \al I_2$ and $\a_1 = \beta I_2$, we have $\ahom = \sqrt{\alpha \beta} \, I_2$. 
\end{enumerate}
\end{exercise}

\section{Quantitative convergence of the subadditive quantities}
\label{s.iteration}

In this section, we present the proof of Theorem~\ref{t.subadd}. The main step in the argument and the focus of most of the section is to obtain the following proposition. 

\begin{proposition}
\label{p.subaddE}
There exist $\alpha(d,\Lambda)\in \left(0,\frac 12 \right]$ and $C(d,\Lambda)<\infty$ such that, for every $n\in\N$, 
\begin{equation}
\label{e.subaddE}
\sup_{p\in B_1} \E \left[ J(\cu_n,p,\ahom p)\right] \leq C3^{-n\alpha}. 
\end{equation}
\end{proposition}

The quantity $U \mapsto J(U,p,q)$ is nonnegative and subadditive. We expect it to converge to zero for $q=\ahom p$ and as $U$ becomes large, and indeed, what~\eqref{e.subaddE} gives is a rate of convergence for this limit in $L^1(\Omega,\P)$ in large cubes, asserting that the expectation is at most the $C(3^n)^{-\alpha}|q|^2$, where $3^n$ is the side length of the cube. The basic idea for proving this estimate is to show that $J(\cu_n,p,\ahom p)$ must contract by a factor less than~$1$ as we pass from triadic scale~$n$ to~$n+1$. That is, we must show that, for some $\theta<1$,
\begin{equation}
\label{e.wts1}
\E \left[ J(\cu_{n+1},p,\ahom p)\right]
\leq \theta \, \E \left[ J(\cu_n,p,\ahom p)\right].
\end{equation}
If we could prove~\eqref{e.wts1}, then a simple iteration would give~\eqref{e.subaddE}. The inequality~\eqref{e.wts1} can be rewritten as 
\begin{equation*}
\E \left[ J(\cu_{n+1},p,\ahom p)\right] 
\leq C \left( \E \left[ J(\cu_{n},p,\ahom p)\right] - \E \left[ J(\cu_{n+1},p,\ahom p)\right]  \right).
\end{equation*}
We do not prove precisely this, but rather a slightly weaker inequality in this spirit which is still strong enough to yield~\eqref{e.subaddE} after iteration. We handle all values of the parameter~$p$ at the same time, and we replace the term in parentheses on the right side of the previous inequality with 
\begin{align*}  
\tau_n  := & \sup_{p,q \in B_1} \Ll(  \E \Ll[ J(\cu_{n},p,q) \Rr] - \E \Ll[ J(\cu_{n+1},p,q) \Rr] \Rr) \\
 = & \sup_{p\in B_1} \Ll(  \E \Ll[ \nu(\cu_{n},p) \Rr] - \E \Ll[ \nu(\cu_{n+1},p) \Rr] \Rr) + \sup_{q\in B_1} \Ll(  \E \Ll[ \nu^*(\cu_{n},q) \Rr] - \E \Ll[ \nu^*(\cu_{n+1},q) \Rr] \Rr).
\end{align*}
This is morally close to the same, but the resulting inequality is a bit weaker. Notice that~$\tau_n$ measures the strictness in the subadditivity relation between scales~$n$ and~$n=1$. We sometimes call it the \emph{additivity defect}. This leaves us with the task of proving
\begin{equation}
\label{e.wts2}
\E \left[ J(\cu_{n+1},p,\ahom p)\right]
\leq C \tau_n. 
\end{equation}
In view of Lemma~\ref{l.minimalset}, we can roughly summarize~\eqref{e.wts2} in words as the statement that ``if the expectations of both of the subadditive quantities do not change much as we change the scale from $n$ to $n+1$ (that is, if the additivity defect is small) then the subadditive quantities have already converged.'' Actually, what we prove is still slightly weaker than~\eqref{e.wts2}, because we need to monitor the additivity defect across many scales rather than just the largest scale. The precise version of~\eqref{e.wts2} is stated in Lemma~\ref{l.iterstep}, below, and most of this section is focused on its proof.

\smallskip

We begin the proof of Proposition~\ref{p.subaddE} by putting the formula~\eqref{e.Jquadresponse} in a more useful form for the analysis in this section. It states that we can control the difference of the optimizers for both $\nu$ and $\nu^*$ on two scales by the strictness of the subadditivity. We have essentially encountered this inequality already in Chapter~\ref{c.one} in the case $q=0$: see~\eqref{e.gluingerror0}.

\begin{lemma}
\label{l.quadresponse}
Fix a bounded Lipschitz domain $U\subseteq \Rd$ and let $\{U_1,\ldots,U_k\}$ be a partition of $U$ into smaller Lipschitz domains, up to a set of measure zero. Then, for every $p,q\in \Rd$, 
\begin{multline} 
\label{e.quadresponse}
\sum_{j=1}^k  \frac{|U_j|}{|U|}  \frac12 \left\| \a^{\frac12} \left( \nabla v (\cdot,U,p,q) - \nabla v(\cdot,U_j,p,q) \right) \right\|_{\underline{L}^2(U_j)}^2 
\\
= \sum_{j=1}^k \frac{|U_j|}{|U|}\left( J(U_j,p,q) -  J(U,p,q) \right).
\end{multline}
\end{lemma}
\begin{proof}
Write $v:= v(\cdot,U,p,q)$ and $v_j:= v(\cdot,U_j,p,q)$. According to~\eqref{e.Jquadresponse}, for every $j \in\{1,\ldots,k\}$,
\begin{multline*}
 \fint_{U_j} \frac12 \left( \nabla v - \nabla v_j \right)\cdot \a \left( \nabla v - \nabla v_j \right)
\\
= J(U_j,p,q) - \fint_{U_j} \left( -\frac12\nabla v \cdot \a\nabla v  -p\cdot\a \nabla v + q\cdot \nabla v  \right).
\end{multline*}
Summing this over $j \in \{1,\ldots,k\}$ gives the lemma. 
\end{proof}

We next turn to one of the key steps in the proof of Proposition~\ref{p.subaddE}, which is to obtain some control of the variance of the spatial average of the gradient of the maximizer of $\nu^*$. By~\eqref{e.Jderq}, this is the same as controlling the variance of the gradient of $J$ in the~$q$ variable, that is, the matrix $\a_*$. 

\smallskip

Note that the following lemma is obvious when $q = 0$, since $v(\cdot,U,p,0) \in \ell_p+H^1_0(U)$, which implies that
\begin{equation}
\label{e.nablavobvious}
\fint_{U} \nabla v(\cdot,U,p,0) = \fint_{U} \nabla \ell_p(\cdot,U,p,0) = p. 
\end{equation}
Thus the spatial average of $\nabla v(\cdot,U,p,0)$ over $U$ is actually deterministic. All of the interest in the lemma therefore concerns the maximizers of~$\nu^*$. We note that its proof is the only place where we use the unit range of dependence assumption in the proof of Proposition~\ref{p.subaddE}. We also note that we do not need the full power of independence here: essentially any mild decorrelation assumption would do. 

\smallskip

It is convenient to extend the notation $\var[X]$ to random \emph{vectors} $X$ (taking values in~$\Rd$) by setting 
$\var[X] := \E \left[ \left| X- \E\left[ X \right] \right|^2 \right]$.

\begin{lemma}
\label{l.spatavg}
There exist $\kappa(d)>0$ and $C(d,\Lambda)<\infty$ such that, for every $p,q\in B_1$ and $m\in\N$,
\begin{equation} 
\label{e.spatavg}
\var \left[ \fint_{\cu_m} \nabla v(\cdot,\cu_m,p,q)  \right] 
\leq C3^{-m\kappa} + C \sum_{n=0}^m 3^{-\kappa(m-n)} \tau_n. 
\end{equation}
\end{lemma}
\begin{proof}
\emph{Step 1.} 
To ease the notation, denote
\begin{equation}
\label{e.nots}
v := v(\cdot,\cu_{m},p,q), \quad
v_z:= v(\cdot,z+\cu_{n},p,q). 
\end{equation}
We claim that, for every $p,q\in B_1$,
\begin{equation} 
\label{e.vardownscales}
\var^{\frac12} \left[  \fint_{\cu_m} \nabla v\right]
\leq \var^{\frac12} \left[ 3^{-d(m-n)}\sum_{z\in 3^n\Zd\cap \cu_m}  \fint_{z+\cu_n} \nabla v_z\right] 
+ C \left( \sum_{k=n}^{m-1} \tau_k \right)^{\frac12}.
\end{equation}
Using the identity 
\begin{equation*} \label{}
\fint_{\cu_m} \nabla v = 3^{-d(m-n)} \sum_{z\in 3^n\Zd \cap \cu_m} \fint_{z+\cu_n} \left( \nabla v - \nabla v_z \right)
+ 3^{-d(m-n)} \sum_{z\in 3^n\Zd\cap \cu_m} \fint_{z+\cu_n} \nabla v_z,
\end{equation*}
we find that 
\begin{multline} 
\label{e.varsplitting}
\var^{\frac12} \left[ \fint_{\cu_m} \nabla v \right]
\leq 
 \E \left[ 3^{-d(m-n)} \sum_{z\in 3^n\Zd \cap \cu_m} \left| \fint_{z+\cu_n} \left( \nabla v - \nabla v_z \right) \right|^2 \right]^{\frac12}   
\\
+  \var^{\frac12} \left[ 3^{-d(m-n)}\sum_{z\in 3^n\Zd\cap \cu_m} \fint_{z+\cu_n} \nabla v_z \right]. 
\end{multline}
By Lemma~\ref{l.quadresponse},
\begin{align}
\notag
\sum_{z\in 3^n\Zd\cap \cu_m} \left| \fint_{z+\cu_n} \left(  \nabla v - \nabla v_z \right) \right|^2
& \leq 
\sum_{z\in 3^n\Zd\cap \cu_m}
\fint_{z+\cu_n} \left|  \nabla v - \nabla v_z \right|^2
\\ & 
\leq C \sum_{z\in 3^n\Zd\cap \cu_m} \left(J(z+\cu_n,p,q) - J(\cu_m,p,q) \right).
 \label{e.downscales}
\end{align}
Taking the expectation of~\eqref{e.downscales} yields, by stationarity, 
\begin{align*} \label{}
\E \left[ 3^{-d(m-n)} \sum_{z\in 3^n\Zd \cap \cu_m} \left| \fint_{z+\cu_n} \left( \nabla v - \nabla v_z \right) \right|^2 \right]
&
\leq C \E \left[ J(\cu_n,p,q) - J(\cu_m,p,q) \right] 
\\ & 
= C \sum_{k=n}^{m-1} \tau_k. 
\end{align*}
Combining the previous display and~\eqref{e.varsplitting} yields~\eqref{e.vardownscales}.

\smallskip

\emph{Step 2.} We use independence to obtain the existence of $\theta(d) \in (0,1)$ such that, for every $p,q\in B_1$ and $n\in\N$,
\begin{equation} 
\label{e.varsumcontract}
\var^{\frac12}  \left[ 3^{-d}\sum_{z\in 3^n\Zd\cap \cu_{n+1}} \fint_{z+\cu_n} \nabla v(\cdot,z+\cu_n,p,q)\right] 
\leq 
\theta \var^{\frac12}  \left[\fint_{\cu_n} \nabla v(\cdot,\cu_n,p,q) \right]. 
\end{equation}
To ease the notation, we set, for each $z\in 3^n\Zd$,
\begin{equation*} \label{}
X_z := \fint_{z+\cu_n} \nabla v(\cdot,z+\cu_n,p,q) - \E \left[ \fint_{z+\cu_n} \nabla v(\cdot,z+\cu_n,p,q)\right] .
\end{equation*}
We expand the variance by writing
\begin{equation*} \label{}
\var\left[ 3^{-d}\sum_{z\in 3^n\Zd\cap\cu_{n+1}} X_z \right] 
\\
= 3^{-2d}\sum_{z,z'\in 3^n\Zd\cap \cu_{n+1}}  \cov\left[ X_z,X_{z'} \right].
\end{equation*}
Since $X_z$ is $\F(z+\cu_n)$-measurable, 
by the independence assumption we have that $\cov\left[ X_z,X_{z'} \right]= 0$ in the case that the subcubes corresponding to $z,z'\in 3^n\Zd\cap \cu_{n+1}$ are not neighbors: that is, $\dist(z+\cu_n,z'+\cu_n) \neq 0$ (which, since $n\geq 0$, implies that $\dist(z+\cu_n,z'+\cu_n) \geq1$). For $z,z'\in 3^n\Zd\cap \cu_{n+1}$ corresponding to neighboring subcubes, we use the Cauchy-Schwarz inequality to get
\begin{equation*} \label{}
\left| \cov\left[ X_z,X_{z'} \right] \right| 
\leq \var^{\frac12}\left[X_z \right]\var^{\frac12}\left[X_{z'} \right] 
= \var\left[X_0 \right].
\end{equation*}
The number of such pairs of cubes is clearly at most $3^{2d}-1$, since it suffices to find a single pair of subcubes which are not neighbors (for instance, opposite corners). 
We therefore obtain 
\begin{equation*} \label{}
\var\left[ 3^{-d}\sum_{z\in 3^n\Zd\cap\cu_{n+1}} X_z \right] 
\leq \left( \frac{3^{2d}-1}{3^{2d}}\right) \var\left[ X_0 \right]. 
\end{equation*}
This implies the claim for $\theta := \left( 1-3^{-2d} \right)^{\frac12}$.

\smallskip

\emph{Step 3.} Iteration and conclusion. Fix $p,q\in\Rd$ and denote 
\begin{equation*} \label{}
\sigma_n^2 := \var  \left[\fint_{\cu_n} \nabla v(\cdot,\cu_n,p,q) \right].
\end{equation*}
Combining~\eqref{e.vardownscales} and~\eqref{e.varsumcontract}, we deduce the existence of~$\theta(d) \in (0,1)$ such that, for every $n\in\N$,
\begin{equation} 
\sigma_{n+1}
\leq 
\theta \sigma_n 
+ C \tau_n^{\frac12}. 
\end{equation}
An iteration yields
\begin{equation*} \label{}
\sigma_m \leq \theta^m \sigma_0 + C \sum_{n=0}^m \theta^{m-n} \tau_n^{\frac12}.
\end{equation*}
Squaring this and using $\sigma_0 \leq C$, we get 
\begin{align*} \label{}
\sigma_m^2 
\leq 
2 \theta^{2m}\sigma_0^2 + C \left( \sum_{n=0}^m \theta^{m-n} \tau_n^{\frac12} \right)^2
\leq
C \theta^{2m} + C \sum_{n=0}^m \theta^{m-n} \tau_n.
\end{align*}
Setting $\kappa:= \log \theta / \log 3$, we get 
\begin{equation*} \label{}
\sigma_m^2 \leq C 3^{-m\kappa} + C \sum_{n=0}^m 3^{-\kappa(m-n)} \tau_n.
\end{equation*}
The proof is complete.
\end{proof}

\begin{definition}
\label{def.ahomU}
We define the deterministic matrix $\ahom_{U}$ by
\begin{equation} \label{e.ahomUdef}
\ahom_U
:= 
\E \left[ \a_*^{-1} (U) \right]^{-1} .
\end{equation}
We also denote $\ahom_n:= \ahom_{\cu_n}$ for short.
\end{definition}

By~\eqref{e.aastar.bounds}, we have 
\begin{equation} 
\label{e.ahomUbounds}
\Id \leq \ahom_U \leq \Lambda \Id.
\end{equation}
Moreover, by~\eqref{e.Jderq}, we note that for every $q\in\Rd$,
\begin{equation*} \label{}
\ahom_U^{-1} q =  \E \left[ \fint_{U} \nabla v(\cdot,U,0,q) \right],
\end{equation*}
which is closely related to the quantity we encountered in the statement of Lemma~\ref{l.spatavg}. Indeed, since $(p,q) \mapsto v(\cdot,U,p,q)$ is a linear mapping, and since by Lemma~\ref{l.Jsplitting}, the function $v(\cdot,U,p,0)$ is minus the minimizer of $\nu(U,p)$, which has spatial average $p$, we have
\begin{equation}  
\label{e.esp.spat.av.v}
\E \Ll[ \fint_U \nabla v(\cdot,U,p,q) \Rr] = \ahom_U^{-1} q - p.
\end{equation}

\smallskip

We note for future reference that, for every $q\in B_1$ and $m,n\in\N$ with $m<n$,
\begin{equation}
\label{e.ahommn}
\left| \ahom_n^{-1} q - \ahom_{m}^{-1} q \right|^2
\leq C \sum_{k=m}^{n-1} \tau_k. 
\end{equation}
Indeed, by stationarity and Lemma~\ref{l.quadresponse}, we find
\begin{align*}
\lefteqn{
\left| \ahom_n^{-1} q - \ahom_{m}^{-1} q \right|^2
} \quad & \\
& = \left| \E \left[  \fint_{\cu_{n}} \nabla v(x,\cu_{n},0,q)\,dx- 3^{-d(n-m)} \sum_{z\in 3^m\Zd\cap \cu_{n}} \fint_{z+\cu_m} \nabla v(x,\cu_m,0,q)\,dx  \right]\right|^2 \\
& \leq \E \left[ 3^{-d(n-m)} \sum_{z\in 3^m\Zd\cap \cu_{n}}  \fint_{z+\cu_m} \left| \nabla v(x,z+\cu_{m},0,q) - \nabla v(x,\cu_{n},0,q) \right|^2 \,dx \right] \\
& \leq C\left( \E \left[ J(\cu_m,0,q)\right] - \E\left[ J(\cu_{n},0,q) \right] \right) \\
& \leq C \sum_{k=m}^{n-1} \tau_k. 
\end{align*}
We also note that the uniformly convex function
\begin{equation}
\label{e.mumap}
q\mapsto \, \E \left[ J\left(U, p,q \right) \right] = \frac12 p \cdot \E \left[ \a(U) \right] p + \frac12 q \cdot \E \left[ \a_*^{-1}(U) \right] q-  p\cdot q
\end{equation}
achieves its minimum at the point $\tilde{q} = \E \left[ \a_*^{-1}(U) \right]^{-1} p$, which by definition is $\tilde{q} =\ahom_U p$. It therefore follows from~\eqref{e.Jrepresentation} that, for $C(\Lambda)<\infty$, 
\begin{equation}
\label{e.quadresponseahom}
\E \left[ J\left(U, p,\ahom_Up \right) \right]
 \leq \E \left[ J\left(U, p,q \right) \right] 
 \leq \E \left[ J\left(U, p,\ahom_Up \right) \right] + C \left| q - \ahom_U p \right|^2. 
\end{equation}

\smallskip

We next use the previous lemma and the multiscale Poincar\'e inequality to control the expected flatness of $v(\cdot,\cu_{n+1},p,q)$ in terms of the sequence $\tau_1,\ldots,\tau_n$. Note that the statement of the lemma can be compared to Proposition~\ref{p.vellpconv} and in particular essentially generalizes that proposition with an independent proof.

\begin{lemma}
\label{l.flatness}
There exist $\kappa(d)>0$ and $C(d,\Lambda)<\infty$ such that, for every $n\in\N$ and $p,q \in B_1$,
\begin{multline}
\label{e.flatness}
\E \left[\fint_{\cu_{n+1}} \left| v(x,\cu_{n+1},p,q) - \Ll(\ahom_n^{-1} q - p\Rr)  \cdot x  \right|^2\,dx \right]  \\
 \leq 
 C 3^{2n} \left( 3^{-\kappa n}+ \sum_{m=0}^n 3^{-\kappa (n-m)} \tau_m  \right).
\end{multline}
\end{lemma}
\begin{proof}
Fix $p,q\in B_1$ and denote $Z_m:= 3^m\Zd \cap \cu_{n+1}$. 

\smallskip

\emph{Step 1.} Application of the multiscale Poincar\'e inequality. Corollary~\ref{c.mP.gradient} gives
\begin{align}
\lefteqn{
 \fint_{\cu_{n+1}} \left| v(x,\cu_{n+1},p,q) - (\ahom_n^{-1} q - p) \cdot x  \right|^2\,dx
 } \quad & \notag \\
& \leq C  \fint_{\cu_{n+1}} \left| \nabla v(x,\cu_{n+1},p,q) -  \ahom_n^{-1} q + p \right|^2\,dx \notag \\
& \quad 
+C\left(  \sum_{m=0}^n 3^m \left( |Z_{m}|^{-1}\sum_{y\in Z_{m}} \left| \fint_{y+\cu_m} \nabla v(x,\cu_{n+1},p,q) \,dx -  \ahom_n^{-1}q + p \right|^2 \right)^{\frac12}  \right)^2. \label{e.multiscaleapp}
\end{align}
The first term on the right side is almost surely bounded:
\begin{multline}
\label{e.dumbbound}
 \fint_{\cu_{n+1}} \left| \nabla v(\cdot,\cu_{n+1},p,q) -  \ahom_n^{-1} q  + p \right|^2 \\
 \leq  2 \left| \ahom_n^{-1} q - p \right|^2 +  2\fint_{\cu_{n+1}} \left| \nabla v(\cdot,\cu_{n+1},p,q) \right|^2  \leq C.
\end{multline}
This leaves us with the task of bounding the expectation of the difference between the spatial average of $\nabla v(\cdot,\cu_{n+1},p,q)$ and the deterministic slope $\ahom_n^{-1} q-p$ in all triadic subcubes of $\cu_{n+1}$ down to the unit scale.

\smallskip

\emph{Step 2.} For $\kappa(d)>0$ as in Lemma~\ref{l.spatavg}, we show that, for every $m\in \{ 0,\ldots,n\}$,
\begin{multline}
\label{e.SAdownscales}
|Z_m|^{-1} \sum_{y\in Z_m} \E \left[ \left| \fint_{y+\cu_m} \nabla v(\cdot,\cu_{n+1},p,q)\,dx  -   \ahom_n^{-1} q + p \right|^2 \right] \\
\leq 
C  \left( 3^{-\kappa m} + \sum_{k=0}^m 3^{\kappa(k-m)}\tau_k + \sum_{k=m}^n \tau_k\right).
\end{multline}
We derive~\eqref{e.SAdownscales} by applying the previous lemma and using the triangle inequality to go down to smaller scales. 

\smallskip

By Lemma~\ref{l.quadresponse}, we have, for every $q\in\Rd$, 
\begin{multline}
\label{e.descentthescales}
|Z_m|^{-1} \sum_{y\in Z_m}\fint_{y+\cu_m}  \left| \nabla v(\cdot,\cu_{n+1},p,q) - \nabla v(\cdot,y+\cu_{m},p,q)   \right|^2 
\\  
\leq C|Z_m|^{-1} \sum_{y\in Z_m} \left(  J(y+\cu_m,p,q) - J(\cu_{n+1},p,q) \right) .
\end{multline}
Taking expectations and using stationarity, we obtain
\begin{multline*} \label{}
|Z_m|^{-1} \sum_{y\in Z_m}\E \left[ \fint_{y+\cu_m}  \left| \nabla v(\cdot,\cu_{n+1},p,q) - \nabla v(\cdot,y+\cu_{m},p,q)   \right|^2 \right]
\\
\leq C \left( \E \left[ J(\cu_{m},p,q)\right] - \E\left[ J(\cu_{n+1},p,q) \right] \right) 
= C \sum_{k=m}^n \tau_k.
\end{multline*}
The triangle inequality, the previous display, Lemma~\ref{l.spatavg}, \eqref{e.esp.spat.av.v} and~\eqref{e.ahommn} yield
\begin{align*}
\lefteqn{
|Z_m|^{-1} \sum_{y\in Z_m} \E \left[\left|   \fint_{y+\cu_m}  \nabla v(\cdot,\cu_{n+1},p,q) -  \ahom^{-1}_{n}q + p  \right|^2  \right]
} \qquad & \\
& \leq 
3|Z_m|^{-1} \sum_{y\in Z_m}\E \left[  \fint_{y+\cu_m} \left|  \nabla v(\cdot,\cu_{n+1},p,q) - \nabla v(\cdot,y+\cu_{m},p,q) \right|^2    \right] \\
& \qquad + 3 |Z_m|^{-1} \sum_{y\in Z_m} \E \left[\left|   \fint_{y+\cu_m}  \nabla v(\cdot,y+\cu_{m},p,q) - \ahom^{-1}_{m}q + p  \right|^2  \right] \\
& \qquad + 3  \left|  \ahom^{-1}_{m}q  - \ahom^{-1}_{n}q  \right|^2  \\
& \leq 
C   \sum_{k=m}^n \tau_k
+ C \left( 3^{-\kappa m} + \sum_{k=0}^m 3^{\kappa(k-m)}\tau_k \right).
\end{align*}
This is~\eqref{e.SAdownscales}. 

\smallskip

\emph{Step 3.} We complete the proof of~\eqref{e.flatness}. To summarize what we have shown so far, we combine~\eqref{e.multiscaleapp},~\eqref{e.dumbbound} and~\eqref{e.SAdownscales} into the estimate
\begin{equation}
\label{e.summary}
 \fint_{\cu_{n+1}} \left| v(x,\cu_{n+1},p,q) - \Ll(\ahom_n^{-1} q -p \Rr) \cdot x  \right|^2\,dx
 \leq C\left( 1 +\left( \sum_{m=0}^n 3^m X_m^{\frac12} \right)^2 \right),
\end{equation}
where the random variable 
\begin{equation*}
X_m:= |Z_m|^{-1} \sum_{y\in Z_m}  \left| \fint_{y+\cu_m} \nabla v(\cdot,\cu_{n+1},p,q)  - \left( \ahom_n^{-1} q - p\right) \right|^2
\end{equation*}
satisfies
\begin{equation}
\label{e.EXm}
\E \left[ X_m \right] \leq C \left( 3^{-\kappa m} + \sum_{k=0}^m 3^{\kappa(k-m)}\tau_k + \sum_{k=m}^n \tau_k\right).
\end{equation}
Next we use H\"older's inequality to see that 
\begin{equation*}
\left( \sum_{m=0}^n 3^m X_m^{\frac12} \right)^2 
\leq \left( \sum_{m=0}^n 3^{m} \right) \left( \sum_{m=0}^n 3^m X_m\right) 
\leq C 3^n \sum_{m=0}^n 3^m X_m
\end{equation*}
and then take the expectation of this and apply~\eqref{e.EXm} to get
\begin{align*}
\E\left[ \left( \sum_{m=0}^n 3^m X_m^{\frac12} \right)^2 \right]
& \leq C  3^n \sum_{m=0}^n \left(  3^m \left( 3^{-\kappa m} + \sum_{k=0}^m 3^{\kappa(k-m)}\tau_k + \sum_{k=m}^n \tau_k\right)\right) \\
& \leq C 3^{2n} \left( 3^{-\kappa n}+ \sum_{k=0}^n 3^{-\kappa (n-k)} \tau_k + \sum_{k=0}^n 3^{-(n-k)} \tau_k  \right)
\\ 
& \leq 
C 3^{2n} \left( 3^{-\kappa n}+ \sum_{k=0}^n 3^{-\kappa (n-k)} \tau_k  \right).
\end{align*}
Combining the above yields~\eqref{e.flatness} and completes the proof of the lemma. 
\end{proof}

Now that we have control on the flatness of the maximizers of~$J(\cu_n,p,q)$, we can estimate~$J(\cu_n,p,\ahom_np)$ with the help of the Caccioppoli inequality.

\begin{lemma}
\label{l.iterstep}
There exist $\kappa(d) > 0$ and $C(d,\Lambda)<\infty$ such that, for every $n\in\N$ and $p \in B_1$,
\begin{equation}
\label{e.itersteprealz}
\E \left[  J(\cu_n, p,\ahom_n p)\right]  
\leq 
C \left( 3^{-\kappa n}+ \sum_{m=0}^n 3^{-\kappa (n-m)} \tau_m  \right).
\end{equation}
\end{lemma}
\begin{proof}
Fix $p\in B_1$. Lemma~\ref{l.flatness} asserts that 
\begin{equation*}
\E \left[ \fint_{\cu_{n+1}} \left| v(\cdot,\cu_{n+1},p,\ahom_np) \right|^2 \right] 
\leq 
C 3^{2n} \left( 3^{-\kappa n}+ \sum_{m=0}^n 3^{-\kappa (n-m)} \tau_m  \right).
\end{equation*}
Applying the Caccioppoli inequality (cf. Lemma~\ref{l.Caccioppoli appendix}), we find that
\begin{equation*}
\E \left[ \fint_{\cu_{n}} \left| \nabla v(\cdot,\cu_{n+1},p,\ahom_np)  \right|^2  \right] 
\leq 
C \left( 3^{-\kappa n}+ \sum_{m=0}^n 3^{-\kappa (n-m)} \tau_m \right).
\end{equation*}
To conclude, it suffices to show that, for every $p,q\in B_1$,
\begin{equation}
\label{e.toconclude}
\E \left[  J(\cu_n, p,q)\right] 
\leq 
C\E \left[ \fint_{\cu_{n}} \left| \nabla v(\cdot,\cu_{n+1},p,q)  \right|^2  \right] 
+ C\tau_n. 
\end{equation}
To prove~\eqref{e.toconclude}, fix $p,q\in B_1$ and denote 
\begin{equation*} \label{}
v:=  v(\cdot,\cu_n,p,q)
\quad \mbox{and} \quad
w:=  v(\cdot,\cu_{n+1},p,q).
\end{equation*}
Using~\eqref{e.firstvarJ2}, we find that 
\begin{align*}
J(\cu_n,p,q)
=
\fint_{\cu_n} \frac12 \nabla v  \cdot \a \nabla v  
= 
\fint_{\cu_n} \frac12 \nabla w  \cdot \a \nabla w
-
\fint_{\cu_n} \frac12 \left( \nabla v - \nabla w \right) \cdot \a\left( \nabla v + \nabla w \right).
\end{align*}
By H\"older's inequality, 
\begin{align*} \label{}
\fint_{\cu_n} \left( \nabla v - \nabla w \right) \cdot \a\left( \nabla v + \nabla w \right)
&
\leq 
C\left\| \nabla v - \nabla w \right\|_{\underline{L}^2(\cu_n)} 
\left(  \left\| \nabla v \right\|_{\underline{L}^2(\cu_n)} 
+\left\| \nabla w \right\|_{\underline{L}^2(\cu_{n})} \right)
\end{align*}
and thus, taking expectations, we find that
\begin{align*}
\lefteqn{ 
\E \left[ \fint_{\cu_n} \frac12 \left( \nabla v - \nabla w \right) \cdot \a\left( \nabla v + \nabla w \right) \right]
} \quad &
\\ &
\leq C \E \left[  \left\| \nabla v - \nabla w \right\|_{\underline{L}^2(\cu_n)}^2 \right]^{\frac12} \left( \E \left[\left\| \nabla v \right\|_{\underline{L}^2(\cu_n)}^2 \right] + \E \left[ \left\| \nabla w \right\|_{\underline{L}^2(\cu_{n+1})}^2\right]  \right)^{\frac12}
\\ & 
\leq C\tau_n^{\frac12} \left( \E \left[ J(\cu_n,p,q)\right] + \E \left[ J(\cu_{n+1},p,q)\right] \right)^{\frac12}
\\ & 
\leq C\tau_n^{\frac12} \E \left[ J(\cu_n,p,q)\right]^{\frac12}. 
\end{align*}
Combining these and using Young's inequality, we obtain
\begin{align*} \label{}
\E \left[ J(\cu_n,p,q) \right] 
& \leq 
\E \left[ \fint_{\cu_n} \frac12 \nabla w\cdot \a \nabla w \right] 
+ C\tau_n^{\frac12} \E \left[ J(\cu_n,p,q)\right]^{\frac12}
\\ & 
\leq C \E \left[ \fint_{\cu_n} \left| \nabla w\right|^2 \right] + \frac12 \E \left[ J(\cu_n,p,q)\right] + C\tau_n. 
\end{align*}
Rearranging this gives~\eqref{e.toconclude}. 
This completes the proof of the lemma. 
\end{proof}

We are now ready to complete the proof of Proposition~\ref{p.subaddE}  which, as previously indicated,  is accomplished by an iteration of the result of the previous lemma. 

\begin{proof}[{Proof of Proposition~\ref{p.subaddE}}]
Define the quantity
\begin{equation*} \label{}
E_n:= \sup_{p\in B_1}  \E \left[ J(\cu_n, p,\ahom_n p)\right].
\end{equation*}
Our first goal, which is close to the conclusion of the proposition, is to prove the estimate
\begin{equation} \label{e.Engoal}
E_n \leq C3^{-n\alpha}. 
\end{equation}
To show this, it is natural to attempt to prove the bound $E_{n+1} \leq \theta E_n$ for some constant $\theta(d,\Lambda) \in (0,1)$. This could be iterated to immediately yield~\eqref{e.Engoal}. However, we cannot show exactly this directly because the supremum in the definition of $E_n$ is inconvenient to work with. Instead we consider the quantity defined for $n\in\N$ by
\begin{equation} \label{e.defDn}
F_n : =  \sum_{i=1}^d  \E \left[ J(\cu_n, e_i,\ahom_n e_i)\right]  .
\end{equation}
It is clear that $F_n$ is equivalent to~$E_n$ 
in the sense that
\begin{equation}
\label{e.DnEn}
cE_n \leq F_n \leq C E_n. 
\end{equation}
The reason $F_n$ is more favorable to consider than $E_n$ is because of the bound
\begin{equation}
\label{e.increment}
F_{n} - F_{n+1} \geq c\tau_n. 
\end{equation}
Indeed, to prove~\eqref{e.increment}, we use~\eqref{e.quadresponseahom}, which gives us that
\begin{equation} 
\label{e.Dneasymin}
F_{n+1} \leq   \sum_{i=1}^d 
\E \left[ J(\cu_{n+1}, e_i,\ahom_n e_i)\right],
\end{equation}
and therefore that 
\begin{align*}
\lefteqn{
 F_{n}  - F_{n+1} 
} \  & \\ &
\geq \sum_{i=1}^d \left( \E \left[ J(\cu_n,e_i,\ahom_n e_i) \right]  - \E \left[ J(\cu_{n+1},e_i,\ahom_n e_i) \right] \right)
\\ & 
= \sum_{i=1}^d  \left( \E \left[ \nu(\cu_n,e_i) \right]  - \E \left[ \nu(\cu_{n+1},e_i) \right] +
\E \left[ \nu^*(\cu_n,\ahom_n e_i) \right]  - \E \left[ \nu^*(\cu_{n+1},\ahom_n e_i) \right] \right) 
\\ &
\geq 
c\left( \sup_{p\in B_1} \Ll(  \E \Ll[ \nu(\cu_{n},p) \Rr] - \E \Ll[ \nu(\cu_{n+1},p) \Rr] \Rr) + \sup_{q\in B_1} \Ll(  \E \Ll[ \nu^*(\cu_{n},q) \Rr] - \E \Ll[ \nu^*(\cu_{n+1},q) \Rr] \Rr) \right)
\\ &
= c\tau_n. 
\end{align*}
(Note that we also used the inequality~\eqref{e.ahomUbounds} to get the fourth line in the previous display.)
We can thus prove $F_{n+1} \leq \theta  F_{n} $ by showing that $F_{n} \leq C\tau_n$, since the latter together with~\eqref{e.increment} give us that 
\begin{equation*} \label{}
F_{n} \leq C\tau_n \implies  F_{n}  \leq C( F_{n}  - F_{n+1}) \iff F_{n+1}\leq \left( 1 - \frac1C\right)  F_{n} ,
\end{equation*}
The inequality $F_{n} \leq C\tau_n$ is close to what Lemma~\ref{l.iterstep} gives us, but unfortunately it is actually slightly stronger, since the right side of~\eqref{e.itersteprealz} is more complicated and involves a sum of $\tau_m$ over all scales smaller than $n$ with smaller scales discounted by the exponent~$\kappa=\kappa(d)$ from that lemma. Therefore, rather than $ F_{n} $, it is natural to work instead with the sequence $\tilde{F}_n$ defined by 
\begin{equation*} \label{}
\tilde{F}_n:= 3^{-\frac \kappa 2n} \sum_{m=0}^n 3^{\frac \kappa2m}  F_m .
\end{equation*}
Notice that~$\tilde{F}_n$ is essentially a weighted average of~$F_m $ over all scales smaller than~$n$, and is in particular a stronger quantity than~$ F_{n} $ in the sense that~$ F_{n}  \leq \tilde{F}_n$.
Therefore it suffices to control~$\tilde{F}_n$. 

\smallskip

We next show that there exist $\theta(d,\Lambda) \in \left[\frac12 ,1 \right)$ and $C(d,\Lambda)<\infty$ such that, for every $n\in\N$, 
\begin{equation}
\label{e.tildeiterstep}
\tilde{F}_{n+1} \leq \theta \tilde{F}_n + C3^{-\frac \kappa 2n}. 
\end{equation}
Using~\eqref{e.increment} and $F_0  \leq C$, we have
\begin{equation*} \label{}
\tilde{F}_n - \tilde{F}_{n+1} 
\geq 
3^{-\frac \kappa 2n} \sum_{m=0}^n 3^{\frac \kappa 2m} \left( F_m  - F_{m+1} \right) - C3^{-\frac \kappa 2n} 
\geq  
c \left( 3^{-\frac \kappa 2n} \sum_{m=0}^n 3^{\frac \kappa 2m} \tau_m - C3^{-\frac \kappa 2n} \right). 
\end{equation*}
Lemma~\ref{l.iterstep} gives
\begin{align} \notag
\tilde{F}_{n+1} 
\leq \tilde{F}_n 
= 3^{-\frac \kappa 2n} \sum_{m=0}^n 3^{-\frac \kappa 2m} F_m  
& \leq C 3^{-\frac \kappa 2n} \sum_{m=0}^n 3^{\frac \kappa 2m} \left( 3^{-\kappa m}+ \sum_{k=0}^m 3^{-\kappa (m-k)} \tau_k  \right)  \\ \notag
& \leq C3^{-\frac \kappa 2n} + C 3^{-\frac \kappa 2n}  \sum_{m=0}^n \sum_{k=0}^m 3^{\frac\kappa2( 2k-m)} \tau_k
\\ \notag
& =  C3^{-\frac \kappa 2n} + C3^{-\frac \kappa 2n}  \sum_{k=0}^n \sum_{m=k}^n 3^{\frac\kappa2( 2k-m)} \tau_k
\\
& \leq C 3^{-\frac \kappa 2n} + C3^{-\frac \kappa 2n}  \sum_{k=0}^n 3^{\frac \kappa 2 k} \tau_k.
\label{e.octopusalign}
\end{align}
Comparing the previous two displays yields
\begin{equation*} \label{}
\tilde{F}_{n+1} \leq C \left( \tilde{F}_n - \tilde{F}_{n+1}  \right) + C3^{-\frac \kappa2n}.
\end{equation*}
A rearrangement of the previous inequality is~\eqref{e.tildeiterstep}. 

\smallskip

Now an iteration of~\eqref{e.tildeiterstep} gives
\begin{equation*} \label{}
\tilde{F}_{n} \leq \theta^n \tilde{F}_0 + C \sum_{k=0}^n \theta^k3^{-\frac \kappa2(n-k)}.
\end{equation*}
This implies (after making $\theta$ closer to $1$ if necessary) that
\begin{equation*} \label{}
 \sum_{k=0}^n \theta^k3^{-\frac\kappa2(n-k)} \leq Cn \left( \theta \vee 3^{-\frac\kappa2} \right)^n \leq C \theta^n.
\end{equation*}
Using this and the fact that $\tilde{F}_0 = F_0  \leq C$, we obtain
\begin{equation*} \label{}
\tilde{F}_{n}  \leq C\theta^n. 
\end{equation*}
Taking $\alpha(d,\Lambda):= \log 3 / \left| \log \theta \right|>0$ so that $\theta = 3^{-\alpha}$ yields the bound $\tilde{F}_n \leq C3^{-n\alpha}$ and thereby completes the proof of~\eqref{e.Engoal}.

\smallskip

To complete the proof of the proposition,
 we need to obtain the same estimate as~\eqref{e.Engoal} after replacing $\ahom^{-1}_{n}$ in the definition of $E_n$ by $\ahom^{-1}$. Notice that~\eqref{e.Engoal},~\eqref{e.DnEn} and~\eqref{e.increment} imply that 
\begin{equation*}
\tau_n \leq C 3^{-n\alpha}. 
\end{equation*}
Thus, by~\eqref{e.ahommn}, for every $m,n\in\N$ with $m\geq n$, we have
\begin{equation*} \label{}
\left| \ahom_{n} - \ahom_{m} \right|^2 \leq  C3^{-n\alpha}. 
\end{equation*}
It follows that there exists a matrix $\tilde\a\in \R^{d\times d}$ such that 
\begin{equation*} \label{}
\left| \ahom_n - \tilde\a \right|^2 \leq  C3^{-n\alpha}. 
\end{equation*}
The previous inequality,~\eqref{e.quadresponseahom} and~\eqref{e.Engoal} imply that 
\begin{equation} 
\label{e.subaddE.tilde}
\sup_{p\in B_1} \E \left[ J(\cu_n,p,\tilde\a p)\right] \leq C3^{-n\alpha}. 
\end{equation}
In view of Lemma~\ref{l.minimalset}, we obtain that 
\begin{equation*} \label{}
\left|\E \left[ \a(\cu_n)\right]    -  \tilde{\a}\right| + \left| \E \left[ \a_*(\cu_n)\right]   -  \tilde{\a}\right|
\leq C3^{-n\alpha/2}.
\end{equation*}
Comparing this to Definition~\ref{d.ahom}, we see that~$\tilde \a = \ahom$. Therefore~\eqref{e.subaddE.tilde} is the same as~\eqref{e.subaddE}. This completes the proof. 
\end{proof}

We did not use the full strength of the unit range of dependence assumption in the proof of Proposition~\ref{p.subaddE}, and in fact the proof applies essentially verbatim if this assumption is replaced by a milder decorrelation condition (we just need an algebraic decay of correlations to obtain Lemma~\ref{l.spatavg}). We next use independence in a much stronger way to complete the proof of Theorem~\ref{t.subadd} by upgrading the stochastic integrability of the result of Proposition~\ref{p.subaddE}. Here we leverage strongly on the fact that the subadditive quantities are bounded and localized. 

\begin{proof}[{Proof of Theorem~\ref{t.subadd}}]
We first apply the unit range of dependence assumption. For every $m,n\in\N$ with $m\leq n$, denote $Z_m:=3^m\Zd\cap \cu_n$. Note that $J(U,p,\ahom p)$ is $\F(U)$--measurable and $J(U,p,\ahom p) \leq C|p|^2$ by~\eqref{e.aastar.bounds}. Therefore an application of Lemma~\ref{l.barO.boxes} yields, for every $m,n\in\N$ with $m\leq n$,
\begin{equation*}
J(\cu_{n},p,\ahom p)
\leq 
\frac{1}{|Z_m|} \sum_{z\in Z_m} 
J(z+\cu_m,p,\ahom p) 
\leq 
\E \left[ J(\cu_m,p,\ahom p) \right] 
+
\O_2\left( C3^{-\frac d2(n-m)} \right).  
\end{equation*}
As $J(\cu_m,p,\ahom p) \leq C|p|^2$, we obtain from Lemma~\ref{l.change-s} that 
\begin{equation*}
J(\cu_{n},p,\ahom p) - \E \left[ J(\cu_m,p,\ahom p) \right]
\leq 
\O_1\left( C3^{-  d(n-m)} \right).
\end{equation*}
By Proposition~\ref{p.subaddE}, we thus obtain, for some~$\alpha(d,\Lambda)\in \left(0,\tfrac 12\right]$, 
\begin{equation}
\label{e.Jestwithm}
J(\cu_{n},p,\ahom p) 
\leq 
C3^{-m\alpha}
+
\O_1\left( C3^{-d(n-m)} \right).
\end{equation}
To conclude, we use the previous display, Lemma~\ref{l.minimalset} and the elementary inequality
\begin{equation} 
\label{e.sqrtelem}
(a+b)^{\frac12} \leq a^{\frac12} + \frac12 a^{-\frac12} b \quad \forall a,b > 0
\end{equation}
to obtain, for every $m,n\in\N$ with $m\leq n$, 
\begin{equation*} \label{}
\left| 
\a(\cu_n) - \ahom
\right|  + \left| 
\a_*(\cu_n) - \ahom
\right|
 \leq C 3^{-m\alpha/2} + \O_1\left( C3^{m\alpha/2-d(n-m)} \right).
\end{equation*}
Given~$s\in (0,d)$, the choice of $m$ to be the closest integer to $(d-s)n/(d+\alpha/2)$ and then shrinking~$\alpha$ yields the desired inequality~\eqref{e.subadderror}. This completes the proof of the theorem. 
\end{proof}

\section{Quantitative homogenization for the Dirichlet problem}
\label{s.Emwhipped}

In this section, we record some estimates on the convergence of the maximizers of $J(\cu_m,p,q)$ and on the homogenization of the Dirichlet problem which will be useful to us later in the book. We deduce the following two theorems by combining Propositions~\ref{p.qualweakconv},~\ref{p.qualweakconv.energy} and Theorem~\ref{t.DP.blackbox}, proved in the previous chapter, with Theorem~\ref{t.subadd}. 

\begin{theorem}
\label{t.Jmaximizers}
Fix $s\in (0,d)$. 
There exist $\alpha(d,\Lambda)\in \left(0,\frac1d\right]$ and $C(s,d,\Lambda)<\infty$ such that, for every $p,q\in B_1$ and $m\in\N$, 
\begin{multline} 
\label{e.Jmaximizers.conv}
\left\| \nabla v(\cdot,\cu_m,p,q) - \left( \ahom^{-1} q - p \right) \right\|_{\Hminusul (\cu_m)}
+ \left\| \a \nabla v(\cdot,\cu_m,p,q) - \left( q -  \ahom p \right) \right\|_{\Hminusul (\cu_m)}
\\
\leq 
C 3^{m -m \alpha(d-s)} + \O_1\left( C3^{m-ms} \right).
\end{multline}
\end{theorem}
\begin{proof}
Using Proposition~\ref{p.qualweakconv}, and Lemma~\ref{l.Jsplitting} stating that $-v(\cdot,U,p,0)$ is the minimizer of $\nu(U,p)$, we have the bound
\begin{equation*} 
\left\| \nabla v(\cdot,\cu_m,p,0) + p \right\|_{\Hminusul (\cu_m)}^2 + 
\left\| \a \nabla v(\cdot,\cu_m,p,0) + \ahom p \right\|_{\Hminusul (\cu_m)}^2
\leq C +  C3^{2m} \mathcal{E}(m).
\end{equation*}
Here $\mathcal{E}(m)$ is the random variable defined in~\eqref{e.Em.def}. 
On the other hand, by the triangle inequality, we get
\begin{multline*} 
\left\| \nabla v(\cdot,\cu_m,0,q) - \ahom^{-1}q \right\|_{\Hminusul (\cu_m)} \\
\leq \left\| \nabla v(\cdot,\cu_m,\ahom^{-1}q,0) + \ahom^{-1}q \right\|_{\Hminusul (\cu_m)} + \left\| \nabla v(\cdot,\cu_m,\ahom^{-1}q ,q)  \right\|_{\Hminusul (\cu_m)},
\end{multline*}
and similarly 
\begin{multline*} 
\left\| \a \nabla v(\cdot,\cu_m,0,q) - q \right\|_{\Hminusul (\cu_m)} \\
\leq \left\| \a \nabla v(\cdot,\cu_m,\ahom^{-1}q,0) + q \right\|_{\Hminusul (\cu_m)} + \left\| \a \nabla v(\cdot,\cu_m,\ahom^{-1}q ,q)  \right\|_{\Hminusul (\cu_m)}.
\end{multline*}
We have that
\begin{multline*} 
3^{-m} \left\| \nabla v(\cdot,\cu_m,\ahom^{-1}q ,q) \right\|_{\Hminusul (\cu_m)} \\ \leq C 
\left\| \nabla v(\cdot,\cu_m,\ahom^{-1}q ,q)  \right\|_{\underline{L}^{2} (\cu_m)}
 \leq C  \left( J(\cu_m,\ahom^{-1}q ,q) \right)^{\frac 12} 
 ,
\end{multline*}
and applying Theorem~\ref{t.subadd} with $\frac{s+d}{2}$ instead of $s$, together with the representation of $J$ in~\eqref{e.Jrepresentation}, gives
\begin{multline*} 
J(\cu_m,\ahom^{-1}q ,q) = \frac12 \ahom^{-1} q  \cdot \left( \a(\cu_m) - \ahom \right) \ahom^{-1} q  +  \frac12 \ahom^{-1} q \cdot \left( \ahom  - \a_*(\cu_m)\right)
\a_*^{-1}(\cu_m)  q 
\\ \leq C 3^{-m \alpha \frac{d-s}{2}} + \O_1\left( C3^{-m\frac{s+d}2} \right),
\end{multline*}
and then applying~\eqref{e.sqrtelem} and reducing $\al > 0$, if necessary, yields
\begin{equation*}  
3^{-m} \left\| \nabla v(\cdot,\cu_m,\ahom^{-1}q ,q) \right\|_{\Hminusul (\cu_m)} \le C 3^{-m \alpha(d-s)} + \O_1\left(C 3^{-ms} \right).
\end{equation*}
A similar bound also holds for the flux $\a \nabla v(\cdot,\cu_m,\ahom^{-1}q ,q)$. We obtain, for $p,q\in B_1$, 
\begin{multline}
\label{e.lockedloaded}
\left\| \nabla v(\cdot,\cu_m,p,q) - \left( \ahom^{-1} q - p \right) \right\|_{\Hminusul (\cu_m)}^2
+ \left\| \a \nabla v(\cdot,\cu_m,p,q) - \left( q -  \ahom p \right) \ \right\|_{\Hminusul (\cu_m)}^2
\\
\leq 
C \left( 1 +  3^{2m} \mathcal{E}(m)\right) +  C 3^{m-m \alpha(d-s)} + \O_1\left( C3^{m-ms} \right).
\end{multline}
We are therefore left with the task of estimating~$\mathcal{E}(m)$. The claim is that 
\begin{equation}
\label{e.Em.bigbound2}
\mathcal{E}(m)^{\frac12} \leq 
3^{ -m \alpha(d-s)/2} + \O_{1}\left( C3^{-ms} \right).
\end{equation}
Recall that the definition of~$\mathcal{E}(m)$ in~\eqref{e.Em.def} is
\begin{equation*} 
\mathcal{E}(m)
\leq C \left( \sum_{n=0}^{m} 3^{n-m}    \left(\frac{1}{|\mathcal{Z}_n|} \sum_{z\in \mathcal{Z}_n} \left| \a(z+ \cu_n) -  \ahom  \right| \right)^{\frac12} \right)^2.
\end{equation*}
By Theorem~\ref{t.subadd}, for every $z \in \mathcal{Z}_n$,
\begin{equation*} 
\left( \left| \a(z+ \cu_n) -  \ahom  \right|  - C 3^{-n \alpha (d-s)} \right)_+ \leq  \O_1(C 3^{-n(s+d)/2}).
\end{equation*}
We then use Lemma~\ref{l.sum-O} to obtain
\begin{equation*} 
\left(   \left(\frac{1}{|\mathcal{Z}_n|} \sum_{z\in \mathcal{Z}_n} \left| \a(z+ \cu_n) -  \ahom  \right| \right)^{\frac12}  - C 3^{-n \alpha(d-s)/2} \right)_+ 
\leq  \O_2\left( C3^{-n(s+d)/4} \right). 
\end{equation*}
Since the left side is also bounded by $C$, we have by Lemma~\ref{l.change-s} that it is bounded by~$\O_{d}\left( C3^{-n(s+d)/2d}\right)$. As $(s+d)/2d<1$, we can  use Lemma~\ref{l.sum-O} again to obtain
\begin{equation*} 
 \sum_{n=0}^{m} 3^{n-m}  
 \left(   \left(\frac{1}{|\mathcal{Z}_n|} \sum_{z\in \mathcal{Z}_n} \left| \a(z+ \cu_n) -  \ahom  \right| \right)^{\frac12}  - C 3^{-n \alpha(d-s)/2} \right)_+ 
 \leq
\O_{d}\left( C3^{-m(s+d)/2d}\right).
\end{equation*}
Squaring this and rearranging, we get 
\begin{equation} \label{e.Em.bigbound} 
\mathcal{E}(m) \leq C 3^{ -m \alpha(d-s)} + \O_{d/2}\left( C3^{-m(s+d)/d} \right).
\end{equation}
Applying Remark~\ref{r.arbitrary}---see~\eqref{e.impl.converse.of}---and then using~\eqref{e.sqrtelem}, we obtain~\eqref{e.Em.bigbound2}. 
We then obtain~\eqref{e.Jmaximizers.conv} 
by combining~\eqref{e.lockedloaded} and~\eqref{e.Em.bigbound2} and reducing $\alpha$ again, if necessary. The proof is complete.
\end{proof}

We conclude this chapter with the following quantitative homogenization result for the Dirichlet problem.

\begin{theorem}
\label{t.DP}
Fix $s\in (0,d)$, a bounded Lipschitz domain $U\subseteq B_1$ and an exponent $\delta>0$. There exist an exponent $\beta(\delta,d,\Lambda)>0$, a constant $C(s,U,\delta,d,\Lambda)<\infty$ and a random variable $\X_s$ satisfying
\begin{equation}
\label{e.sizeX}
\X_s = \O_1\left(C \right)
\end{equation}
such that the following statement holds. For each $\ep\in (0,1]$ and $f\in W^{1,2+\delta}(U)$, let $u^\ep, u \in f+H^1_0(U)$ respectively denote the solutions of the Dirichlet problems 
\begin{equation}
\label{e.uep.and.ubar}
\left\{
\begin{aligned}
& -\nabla \cdot \left( \a\left(\tfrac x\ep\right) \nabla u^\ep \right) = 0 &  \mbox{in} & \ U, \\
& u^\ep = f & \mbox{on} & \ \partial U,
\end{aligned}
\right.
\quad \mbox{and} \quad 
\left\{
\begin{aligned}
& -\nabla \cdot \left( \ahom \nabla  u  \right) = 0 &  \mbox{in} & \ U, \\
&  u  = f & \mbox{on} & \ \partial U.
\end{aligned}
\right.
\end{equation}
Then we have the estimate
\begin{multline}
\label{e.DPestimates}
\left\| u^\ep -   u  \right\|_{L^{2}(U)}^2
+ \left\| \nabla u^\ep - \nabla  u  \right\|_{\Hminus(U)}^2 
+ \left\| \a\left(\tfrac\cdot\ep\right) \nabla u^\ep - \ahom \nabla  u  \right\|_{\Hminus(U)}^2 
\\
+ \left\| \tfrac12 \nabla u^\ep \cdot  \a\left(\tfrac\cdot\ep\right) \nabla u^\ep - \tfrac12 \nabla  u  \cdot \ahom \nabla  u  \right\|_{W^{-2,1}(U)}
\\
\leq
C \left\| \nabla f \right\|_{L^{2+\delta}(U)}^2 \left( 
\ep^{\beta(d-s)}
+ \X_s \ep^s
\right).  
\end{multline}
\end{theorem}
\begin{proof}
By Theorem~\ref{t.DP.blackbox}, inequality \eqref{e.linktopast} and the fact that $\mcl E(m) \le C(d,\Lambda)$, there exists $\gamma(\delta,d,\Lambda)>0$ such that, for each $\ep\in (0,1)$, if we choose $m_\ep\in\N$ so that $3^{m_\ep-1}<\ep^{-1} \leq 3^{m_\ep}$, then we have
\begin{multline}
\label{e.DPestimates.bb.app}
\left\| u^\ep -   u  \right\|_{L^{2}(U)}^2+ 
\left\| \nabla u^\ep - \nabla  u  \right\|_{\Hminus(U)}^2 
+ \left\| \a\left(\tfrac\cdot\ep\right) \nabla u^\ep - \ahom \nabla  u  \right\|_{\Hminus(U)}^2 
\\
+ \left\| \tfrac12 \nabla u^\ep \cdot  \a\left(\tfrac\cdot\ep\right) \nabla u^\ep - \tfrac12 \nabla  u  \cdot \ahom \nabla  u  \right\|_{W^{-2,q}(U)}
\\
\leq C\left\| \nabla f \right\|_{L^{2+\delta}(U)}^2 
\inf_{r\in (0,1)} \left(  r^{2\gamma}+ \frac{1}{r^{6+d}} \left(  3^{-2m_\ep}  +  \mathcal{E}(m_\ep) \right)^{\frac12} \right),
\end{multline}
Fix $s\in (0,d)$. It remains to show that, for some $\beta>0$ and $\X_s$ as in the statement of the theorem,  
\begin{equation*} \label{}
\inf_{r\in (0,1)} \left(  r^{2\gamma}+ \frac{1}{r^{6+d}} \left(  3^{-2m_\ep}  +  \mathcal{E}(m_\ep) \right)^{\frac12}  \right) 
\leq C\ep^{\beta(d-s)} + \X \ep^s.
\end{equation*}
In other words, we need to show that, for some $\beta(\delta,d,\Lambda)>0$
\begin{equation} 
\label{e.sumthisupO}
\sup_{m\in\N } 3^{ms} \inf_{r\in (0,1)} \left(   r^{2\gamma}+ \frac{1}{r^{6+d}} \left(  3^{-2m}  +  \mathcal{E}(m) \right)^{\frac12}  - C3^{-m\beta(d-s)} \right)_+ \leq \O_1(C). 
\end{equation}
We have already seen in~\eqref{e.Em.bigbound2} that, for some $\alpha(d,\Lambda)>0$ and with $s':= \frac12(d+s)$, we have the bound
\begin{equation*}
\left( 3^{-2m}+ \mathcal{E}(m) \right)^{\frac12}
 \leq C 3^{ -m \alpha(d-s)} + \O_1\left( C3^{-ms'} \right).
\end{equation*}
We next select $r \in (0,1)$ to be as small as possible, but still large enough that
\begin{equation*} \label{}
\frac1{r^{6+d}} 3^{-m\alpha(d-s)} \leq 3^{-m\alpha(d-s)/2}
\quad \mbox{and} \quad
\frac1{r^{6+d}} 3^{-ms'/2} \leq 3^{-ms/2}. 
\end{equation*}
We can take for example 
$r:= 3^{-m(s'-s)/ (12+2d)} \vee 3^{-m\alpha(d-s)/(12+2d)}.$
We obtain, for some $\beta(\gamma,d,\Lambda)>0$, 
\begin{equation*} \label{}
r^{2\gamma} + \frac1{r^{6+d}} \left( 3^{-2m}+ \mathcal{E}(m) \right) 
 \leq C3^{ -m \beta(d-s)} + \O_1\left( C3^{-ms -m \beta (d-s)} \right).
\end{equation*}
Thus, by Lemma~\ref{l.sum-O},
\begin{align*} \label{}
\sum_{m\in \N} 3^{ms} \left( \inf_{r>0} \left( r^{2\gamma} + \frac1{r^{6+d}} \left( 3^{-2m}+ \mathcal{E}(m) \right) \right) - C3^{- m\beta (d-s)} \right)_+
& \leq \O_1\left(C \sum_{m\in\N} 3^{-m\beta (d-s)} \right) 
\\ & 
\leq \O_1(C). 
\end{align*}
This completes the proof of~\eqref{e.sumthisupO}.
\end{proof}

\section*{Notes and references}

The first quantitative estimate for the homogenization error was obtained by Yurinski{\u\i}~\cite{Y101}, who proved an algebraic rate of convergence in $L^2(\Omega,\P)$ using different methods. 
The dual subadditive quantity~$\nu^*(U,q)$ was introduced in~\cite{AS} and this chapter is a simplified version of the arguments which appeared in that paper. In particular,~\cite{AS} proved Theorem~\ref{t.subadd} and essentially contained a proof of Theorem~\ref{t.DP}. Most of the improvements to the argument we make here are due to the use of the multiscale Poincar\'e inequality (Proposition~\ref{p.mspoin}) which was introduced later in~\cite{AKM1}. This allowed for an easier proof of Lemma~\ref{l.flatness}, which was accomplished in~\cite{AS} by a long calculation involving Helmholtz-Hodge projections. 

\smallskip

It was also shown in~\cite{AKM1}, by a renormalization-type (bootstrap) argument, that the scaling of the error in Theorem~\ref{t.subadd} can be improved to the statement that, for every $\alpha<1$, 
\begin{equation*} \label{}
\left| J(\cu_n,p,q) - \left( \frac12 p\cdot \ahom p + \frac12 q\cdot \ahom^{-1} q - p\cdot q \right) \right|
\leq \O_2\left( C3^{-n\alpha} \right). 
\end{equation*}
This estimate cannot be improved beyond the exponent~$\alpha=1$ due to boundary layers effects, as explained in~\cite{AKM1}. 

\smallskip

The ideas and techniques in this chapter are robust in the sense that they can be extended beyond the case of linear, uniformly elliptic equations. See~\cite{AD} for an extension to the very degenerate case of supercritical percolation clusters and~\cite{AS,AM} for general nonlinear equations and integral functionals. While the techniques are based on variational methods, they also work in the case of linear equations with non-symmetric coefficients (see Chapter~\ref{c.subadd-fitz}).
In fact, they extend to the setting of general uniformly monotone operators~\cite{AM} and even to parabolic equations with space-time coefficients~\cite{ABM}. 

\smallskip

We remark that the philosophy underlying the analysis in this chapter, as well as that of Chapter~\ref{c.A1}, is similar to the ``variational approach to coarse-graining'' described in~\cite{DLPS}.

\addtocontents{toc}{\protect\newpage}



\chapter{Regularity on large scales}
\label{c.regularity}

In this chapter, we show that quantitative homogenization implies that solutions are typically much more regular, on large scales, than the deterministic regularity theory for uniformly elliptic equations with rough coefficients would predict.

\smallskip

Such improved regularity results are of fundamental importance to the quantitative theory of homogenization. Indeed, to get finer estimates for the rate of homogenization than what we have obtained thus far, we need better a priori estimates on the gradients of solutions. We know from the very definition of weak solution that a solution~$u$ of our equation will satisfy~$\nabla u \in L^2$. By the Meyers estimate (see \eqref{e.display.Meyers} below), this can be slightly improved to~$\nabla u\in L^{2+\delta}$ for some small $\delta >0$. Our first goal in this chapter is to improve this to something close to~$\nabla u \in L^\infty$, showing roughly that the energy density of a solution cannot concentrate on small scales. Since the Meyers estimate is the best gradient estimate for an arbitrary elliptic equation, such a result will necessarily require some probabilistic ingredient. It turns out that the ingredient we need is the conclusion of Theorem~\ref{t.DP}, which gives us \emph{$\ahom$-harmonic approximation}.
 
\smallskip

The basic idea, borrowed from the classical proof of the Schauder estimates, is to implement a $C^{1,\gamma}$-type iteration to obtain \index{excess decay iteration} \emph{decay of the excess}\footnote{The \emph{excess} is the name sometimes given to the quantity that measures how close a function~$u$ is to an affine function at a given scale. The excess of~$u$ in $B_r$ for example could be defined by $\frac1r \inf_{\ell \in \mathcal{P}_1} \left\| u - \ell \right\|_{\underline{L}^p(B_r)}$ for any convenient~$p\in[1,\infty]$.} by comparing a given solution, in each large ball, to the~$\ahom$-harmonic function with the same boundary values. Homogenization ensures that these two functions are very close to each other and, since the error is at most algebraic, it sums as we sum over a sequence of dyadic balls. On the other hand, the regularity of $\ahom$-harmonic functions gives an improvement in polynomial approximation as we pass from~$B_R$ to $B_{R/2}$. The main difference from the proof of the Schauder estimates is that here the error in comparing to $\ahom$-harmonic functions is controlled by homogenization rather than by freezing the coefficients. This idea was previously used in the context of periodic homogenization by Avellaneda and Lin~\cite{AL1}, and the analysis here can be thought of as a quantitative version of their compactness, blow-down argument.

\smallskip

The main results of the chapter are the $C^{0,1}$-type estimate in Theorem~\ref{t.Lipschitz} and its generalization, Theorem~\ref{t.regularity}, which gives higher $C^{k,1}$-type regularity estimates for every $k\in\N$ and related Liouville results. We set up the context for these results in the following section. 

\section{Brief review of classical elliptic regularity}
\label{s.regularity_brief}

Let $u\in H^1(B_R)$ be an $\a(\cdot)$-harmonic function in $B_R$, that is, a  weak solution of the linear, uniformly elliptic equation
\begin{equation}
\label{e.divformeq}
-\nabla \cdot \left( \a(x) \nabla u \right) = 0 \quad \mbox{in} \ B_R. 
\end{equation}
Here we make no further assumptions on the regularity of the coefficient field beyond measurability or any structural assumption on $\a(\cdot)$ beyond uniform ellipticity. What is the best regularity we can expect for the solution $u$?

\begin{itemize}

\item In terms of H\"older spaces, the best regularity is given by the \emph{De Giorgi-Nash estimate} (cf.~\cite[Chapter 3]{HL}): there exists $\alpha(d,\Lambda)>0$ and a constant $C(d,\Lambda)<\infty$ such that $u \in C^{0,\alpha}_{\mathrm{loc}}\left(B_{R}\right)$ and
\begin{equation}
\label{e.DeGiorgi}
\left\| u \right\|_{L^\infty(B_{R/2})} + R^\alpha \left[ u \right]_{C^{0,\alpha}(B_{R/2})} \leq C \left\| u \right\|_{\underline{L}^2(B_R)}. 
\end{equation}
This estimate only holds for scalar equations. It is false, in general, for linear elliptic systems (there is a counterexample due to De Giorgi~\cite{DG2}). \index{De Giorgi-Nash estimate}

\item In terms of Sobolev spaces, the best regularity is given by the \emph{Meyers estimate}: there exists $\delta(d,\Lambda)>0$ and a constant $C(d,\Lambda) < \infty$ such that $\nabla u \in L^{2+\delta}_{\mathrm{loc}}(B_R)$ and 
\begin{equation} \label{e.display.Meyers}
\left\| \nabla u \right\|_{\underline{L}^{2+\delta}(B_{R/2})} \leq \frac{C}R \left\| u \right\|_{\underline{L}^2(B_R)}.
\end{equation}
This improvement of integrability of the gradient also extends to the case of linear elliptic systems. A complete proof of this estimate can be found in Appendix~\ref{a.meyers}, see Theorem~\ref{t.Meyers appendix}. \index{Meyers estimate}

\end{itemize}

Lest the reader imagine that it is possible to improve either of these estimates, or even to show that $\alpha$ or $\delta$ can be  bounded away from zero independently of the ellipticity constant $\Lambda$, we present the following example due to Meyers~\cite{Me}. Another example is given later in Section~\ref{s.regularity.optimal}. 

\begin{example}[Optimality of De Giorgi-Nash \& Meyers estimates]
\label{ex.meyers}
\index{De Giorgi-Nash estimate!optimality of|(}
\index{Meyers estimate!optimality of|(}
Fix an exponent $\alpha \in (0,1]$ and consider the function $u : \R^2 \to \R$ defined by
\begin{equation*} \label{}
u(x) := |x|^{\alpha-1} x_1. 
\end{equation*}
For every $r\in (0,1]$, we have 
\begin{equation*}
u \not \in C^{0,\beta} (B_r) \quad \mbox{for any}\  \beta>\alpha
\quad \mbox{and} \quad 
\nabla u \not \in L^p(B_r) \quad \mbox{for any} \ p \geq d + \frac{d\alpha}{1-\alpha}.
\end{equation*}
Consider the coefficient matrix~$\a(x)$ given by 
\begin{equation*}
\a(x):= \left( I - \frac{x\otimes x}{|x|^2} \right) + \Lambda \left( \frac{x\otimes x}{|x|^2}\right). 
\end{equation*}
Notice that~$\a(x)$ has eigenvalue~$1$ with multiplicity~$d-1$ (the eigenvectors are orthogonal to~$x$) and eigenvalue~$\Lambda$ with multiplicity one (the eigenvector is~$x$). See Figure~\ref{f.meyerscounterexample} for a visualization in $d=2$. 

\smallskip

We now compute 
\begin{equation*}
\nabla u(x) = |x|^{\alpha-1} \left( e_1 + (\alpha-1)|x|^{-2} x_1 x\right).
\end{equation*}
Therefore
\begin{equation*}
\a(x) \nabla u(x) = |x|^{\alpha-1} e_1 + \left( \Lambda \alpha-1 \right) |x|^{\alpha-3} x_1 x. 
\end{equation*}
Taking the divergence of this, we find 
\begin{align*}
(\nabla \cdot \a\nabla u)(x)
&
= |x|^{\alpha-3} x_1 \left( \alpha-1 + (\alpha\Lambda-1)(\alpha+d-2) \right)
\\ & 
= |x|^{\alpha-3} x_1\Lambda \left( \alpha^2 + (d-2)\alpha - \frac{d-1}{\Lambda} \right). 
\end{align*}
This vanishes therefore if 
\begin{equation*}
\alpha^2 + (d-2)\alpha - \frac{d-1}{\Lambda} = 0. 
\end{equation*}
Therefore $u$ satisfies
\begin{equation*}
-\nabla \cdot \a\nabla u = 0 \quad \mbox{in} \ \Rd
\end{equation*}
if we choose $\alpha$ by the formula
\begin{equation}
\label{e.Meyersexamchoosealpha}
\alpha := \frac12 \left( -(d-2) + \sqrt{(d-2)^2 + \frac{4(d-1)}{\Lambda}} \right).
\end{equation}
Notice that in $d=2$ we have $\alpha = \Lambda^{-1/2}$ and 
\begin{equation*}
\alpha \sim \left( \frac{d-1}{d-2}\right)\Lambda^{-1} \quad \mbox{as} \ \Lambda \to \infty \quad \mbox{in} \ d>2. 
\end{equation*}
In particular, in every dimension, the best exponent  in the De Giorgi-Nash estimate is no larger than the  $\alpha$ in~\eqref{e.Meyersexamchoosealpha}. For the Meyers estimate in $d=2$, the exponent $\delta$ can be no larger than $2\alpha(1-\alpha)$ with $\alpha$ in~\eqref{e.Meyersexamchoosealpha}. In higher dimensions, however, the exponent in the Meyers estimate cannot be improved because we can just add dummy variables to the two-dimensional example.

\smallskip

Also notice that, for $u$ as above, 
\begin{equation}
\label{e.nocoarselipschitz}
\fint_{B_r} \left| \nabla u(x) \right|^2\,dx =  r^{2(\alpha-1)} \fint_{B_1} \left| \nabla u(x) \right|^2\,dx.
\end{equation}
In other words, the energy density of the solution~$u$ is very large in~$B_1$ compared to average of the energy density of~$u$ in $B_r$ for large~$r$. This should be compared with the statement of Theorem~\ref{t.Lipschitz}, below. 
\begin{figure}[tb]
\centering
\begin{center}
\begin{tikzpicture}[scale=8]
\draw[<->,black,very thin] (-0.5,0) -- (0.5,0);
\draw[<->,black,very thin] (0,-0.5) -- (0,0.5);
\def\N{12};
\def\M{24};
\def\Q{25};
\def\S{ 0.004 };
\foreach \i in {1,...,\M}
{
\foreach \j in {1,...,\M}
{
\def\x{ (\i-\N-0.5) / \Q};
\def\y{ (\j-\N-0.5) / \Q };
\def\Le{ sqrt(   (\x)*(\x) + (\y)*(\y)   ) } ;
\def\dx{ (\S*\x / \Le) };
\def\dy{ (\S*\y / \Le ) };
\draw[-] ( {\x + 4*\dx}, {\y + 4*\dy} ) -- ({\x - 4*\dx},{\y - 4*\dy});
\draw[-] ( {\x + \dy}, {\y -\dx} ) -- ({\x - \dy},{\y + \dx});
}
}
\end{tikzpicture}
\end{center}
\caption{The coefficient field $\a$ in Example~\ref{ex.meyers} with~$d=2$ and~$\Lambda=4$. The matrix $\a(x)$ at each grid point $x$ is represented by two lines which indicate the direction of the eigenvectors of~$\a(x)$ with lengths proportional to the corresponding eigenvalue. Notice that the eigenvector with the largest eigenvalue points toward the origin, where the singularity occurs in the solution.}
\label{f.meyerscounterexample}
\end{figure}

\smallskip

Incidentally, it is known that in $d=2$ solutions of~\eqref{e.divformeq} belong to $C^{0,\Lambda^{-1/2}}(B_{1/2})$ and therefore~$\Lambda^{-1/2}$ is the sharp exponent for H\"older regularity in dimension~$d=2$: see~\cite{PS} for a short and beautiful proof. In~$d>2$, the best exponent is not known. In fact, a problem posed by De Giorgi (which is still open) is to show that the exponent~$\alpha(d,\Lambda)$ in the De Giorgi-Nash estimate satisfies a lower bound of the form $\alpha(d,\Lambda) \geq c \Lambda^{-p}$ for some~$c(d)>0$ and~$p(d)<\infty$. \qed
\index{De Giorgi-Nash estimate!optimality of|)}
\index{Meyers estimate!optimality of|)}
\end{example}

The relatively weak regularity of solutions of~\eqref{e.divformeq} with oscillating coefficients stands in sharp contrast to the very strong regularity possessed by solutions of a constant-coefficient equation
\begin{equation}
\label{e.divformeqhom}
-\nabla \cdot \left( \ahom \nabla u \right) = 0 \quad \mbox{in} \ B_R
\end{equation}
which in the case of equations, of course, up to the affine change of variables $x\mapsto \ahom^{-\frac 12}x$, is just Laplace's equation
\begin{equation*} \label{}
-\Delta u(\ahom^{\frac12}\, \cdot) = 0 \quad \mbox{in} \ \ahom^{-\frac12}B_R. 
\end{equation*}
We call solutions of~\eqref{e.divformeqhom} $\ahom$-\emph{harmonic.}
These $\ahom$-harmonic functions possess the same regularity as harmonic functions (cf.~\cite{Evans} or~\cite{HL}): for every $k\in\N$, there exists $C=C(k,d,\Lambda)<\infty$ such that, for every~$R>0$ and~$\ahom$-harmonic function $u$ on $B_R$, 
\begin{equation}
\label{e.pointwiseharm}
\left\| \nabla^k u \right\|_{L^\infty(B_{R/2})}  \leq \frac{C}{R^{k}} \left\| u \right\|_{\underline{L}^2(B_{R})}. 
\end{equation}
Moreover, the constant $C_k$ grows in $k$ at a slow enough rate to guarantee the \emph{real analyticity} of $\ahom$-harmonic functions (see~\cite{Evans}). This is to say that an $\ahom$-harmonic function in $B_r$ has the property that at any point~$z \in B_r$, the Taylor series of $u$ at $z$ is convergent with a certain convergence radius. The proof of~\eqref{e.pointwiseharm} is based on the fact that any partial derivative of $u$ is still $\ahom$-harmonic, and in this way one gets strong $H^k$-estimates. These $H^k$-estimates can be turned into $C^{k,1}$-estimates using  Sobolev (Morrey-type) inequalities. 

\smallskip

\index{Schauder estimates|(}

Each of the estimates discussed so far is \emph{scale-invariant}. In other words, while we have stated them in~$B_R$ for general~$R>0$, it is actually enough to consider~$R=1$ since the general statement can be recovered by scaling. Indeed, rescaling the coefficients under either assumption (measurable in~$x$ or constant) does not change the assumption. There are certain elliptic estimates, however, which are not scale-invariant. In these cases, the assumptions on the coefficients~$\a(\cdot)$ always refer to an implicit length scale. A good example of this is the family of \emph{Schauder estimates}. 

\smallskip

For divergence-form equations, one version of the Schauder estimates states that the De Giorgi-Nash and Meyers estimates may be improved if we assume that the coefficients are H\"older continuous. Precisely, for every $\beta\in (0,1)$, if we assume in addition that $\a\in C^{0,\beta}(B_R)$, then any weak solution $u\in H^1(B_1)$ of~\eqref{e.divformeq} satisfies $u \in C^{1,\beta}_{\mathrm{loc}}(B_R)$ and, if we define the length scale
\begin{equation*} \label{}
r_0(\a):= \left[ \a \right]_{C^{0,\beta}(B_{R/2})}^{-\frac1\beta},
\end{equation*}
then there exists $C(d,\Lambda)<\infty$ such that, for every $r \in (0,r_0(\a)]$ and $x\in B_{R/4}$,
\begin{equation} 
\label{e.schauder}
\left\| \nabla u \right\|_{L^\infty(B_r(x))} + r^\beta \left[ \nabla u \right]_{C^{0,\beta}(B_r(x))} \leq C r^{-1}\| u \|_{\underline{L}^2(B_{2r}(x))}. 
\end{equation}
In other words, it is only on length scales small enough that the coefficients have their properly scaled H\"older seminorm at most of order $O(1)$ that the quantitative estimate holds. In the typical way this estimate is stated, the dependence on the seminorm $\left[ \a \right]_{C^{0,\beta}(B_R)}$ is hidden in the constant $C$ and the restriction $r\leq r_0(\a)$ is removed. Even though the two statements are obviously equivalent, we prefer our formulation, in which the estimate is only valid on certain length scales. 

\smallskip

Let us say a few words about the proof of the Schauder estimate~\eqref{e.schauder}, which can be found in~\cite[Chapter 3]{HL}. Like most estimates which are not scale-invariant, its proof is \emph{perturbative}. The idea is a very classical one: since the coefficients are continuous (an assumption made quantitative by the introduction of $r_0(\a)$), on small scales they are close to being constant. It is therefore natural to wonder to what extent, in this situation, one can ``borrow'' the very strong regularity of solutions of constant-coefficient equations (harmonic functions) by approximation. 

\smallskip

We present an outline of how the perturbation argument works in this classical framework to obtain $C^{0,1-}$ regularity of solutions under the assumption that the coefficients are uniformly continuous:  for a bounded domain $U \subset  \R^d$,\index{oscillation decay iteration|(}
\begin{equation}  \label{e.unifcont}
\lim_{r \to 0} \sup_{x \in U} \| \a(\cdot) - \a(x) \|_{L^\infty(U \cap B_r(x))} = 0. 
\end{equation}
Proving H\"older regularity involves understanding the pointwise behavior of a function. However, a priori a weak solution is only an $H^1$ function and therefore is not necessarily even defined at every point; likewise, the definition of weak solution is a statement involving integrals, not pointwise information. It is therefore helpful to think of H\"older continuity in terms of the Campanato characterization (see Exercise~\ref{ex.campanato} below), as a quantitative statement concerning the $L^2$-oscillation of a function on small scales. We will proceed by trying to show that the $L^2$ oscillation of a weak solution must improve as we zoom in to smaller and smaller scales. 

\smallskip

Let $u \in H^1(U)$ be a solution in $U$ of $-\nabla \cdot(\a(x) \nabla u(x)) = 0$ and $B_r(x_0) \subset U$. By translation, we may assume $x_0=0$ to ease the notation. One now takes $w \in u + H_0^1(B_{r/2})$ solving $-\nabla \cdot(\a(x_0) \nabla w(x)) = 0$ and subtracts the equations of $u$ and $w$. Testing the resulting equation with $u-w \in H_0^1(B_{r/2})$, one arrives, after an application of Young's inequality and ellipticity of $\a(\cdot)$, at
\begin{equation*} 
\fint_{B_{r/2}} \left|\nabla u - \nabla w\right|^2  \leq \left\| \a(\cdot)-\a(x_0)\right\|_{L^\infty(B_{r/2})}^2 \fint_{B_{r/2}} \left|\nabla u \right|^2  \,.
\end{equation*}
By the Poincar\'e and Caccioppoli inequalities (cf. Lemma~\ref{l.Caccioppoli appendix}), we deduce 
\begin{equation*} 
\fint_{B_{r/2}} \left| u - w\right|^2  \leq C \left\| \a(\cdot)-\a(x_0)\right\|_{L^\infty(B_r)}^2 \fint_{B_{r}} \left| u - (u)_{B_{r}}\right|^2 .
\end{equation*}
Using~\eqref{e.pointwiseharm} for $k=1$ and the Poincar\'e inequality, it follows that, for every  $t \in \left(0,\frac14r\right]$, 
\begin{equation*} 
\left\| w - (w)_{B_t} \right\|_{\underline{L}^2(B_{t})} 
\leq Ct \left\| \nabla w \right\|_{\underline{L}^2(B_t)}
\leq C \left(\frac{t}{r}\right)\left\| w -(w)_{B_{r/2}}  \right\|_{\underline{L}^2(B_{r/2})}.
\end{equation*}
Combining the above displays and using the triangle inequality yields 
\begin{equation*} 
\left\| u - (u)_{B_t} \right\|_{\underline{L}^{2}(B_{t})} \leq C\left( \frac{t}{r}  + \left(\frac{r}{t}\right)^{\frac{d}{2}} \left\| \a(\cdot)-\a(x_0)\right\|_{L^\infty(B_r)} \right) \left\| u - (u)_{B_{r}}\right\|_{\underline{L}^2(B_{r})}.
\end{equation*}
Now, for any $\delta \in (0,1)$ we may take $t = \theta r$ with $\theta(\delta,d,\Lambda) \in (0,1)$ small enough, and subsequently find, using~\eqref{e.unifcont}, $r_0(\a)$  small enough so that $r \leq r_0(\a)$ implies
\begin{equation*} 
\left\| u - (u)_{B_{\theta r}} \right\|_{\underline{L}^{2}(B_{\theta r})}  \leq \theta^{1-\delta} \left\| u - (u)_{B_{r}}\right\|_{\underline{L}^2(B_{r})}.
\end{equation*}
This can be iterated. Indeed, we set $r_{j+1} := \theta r_j$ and choose $r_0 \in  \left( 0 , \frac12 r_0(\a) \right]$ to get
\begin{equation*} 
\left\| u - (u)_{B_{r_j}} \right\|_{\underline{L}^{2}(B_{r_j})}  \leq  \left(\frac{r_j}{r_0}\right)^{1-\delta} \left\| u - (u)_{B_{r_0}}\right\|_{\underline{L}^2(B_{r_0})}.
\end{equation*}
This gives, for any $r \in \left( 0 , r_0(\a) \wedge \dist(x_0,\partial U)  \right]$, 
\begin{equation} \label{e.regcamp000}
\sup_{t \in (0,r)} t^{\delta-1} \left\| u - (u)_{B_{t}(x_0)} \right\|_{\underline{L}^{2}(B_{t}(x_0))}  \leq C(\delta,d,\Lambda) r^{\delta-1} \left\| u - (u)_{B_{r}(x_0)} \right\|_{\underline{L}^{2}(B_{r}(x_0))},
\end{equation}
The definition of Campanato spaces (see e.g.~\cite[Section 2.3]{Giu} or~\cite[Chapter 3]{HL}) stems from the inequality as above, and standard iteration techniques used in such spaces lead to full $C^{0,1-}$-regularity, that is, $C^{0,1-\delta}$-regularity for every $\delta \in (0,1)$. We leave the details as an exercise. 
\begin{exercise} 
\label{ex.campanato}
Let $\gamma \in (0,1]$, $f \in L_{\textrm{loc}}^2(\R^d)$. Suppose that, for some $R>0$ and a domain $U\subseteq\R^d$, 
\begin{equation*} 
\sup_{x \in U} \sup_{r \in (0,R)} r^{-\gamma} \left\| f - (f)_{B_{r}(x)} \right\|_{\underline{L}^{2}(B_{r}(x))} < \infty.
\end{equation*}
Show that $f \in C^{0,\gamma}(U)$. Hint: see~\cite[Theorem 3.1]{HL}.
 \end{exercise}
\index{Schauder estimates|)}
 \index{oscillation decay iteration|)}

\smallskip

In subsequent sections, we see that this relatively simple idea can be pushed much further. We emphasize that the perturbation argument above does not use any properties of the function~$u$ apart from \emph{harmonic approximation}: that at each given scale, there exists a harmonic function approximating the original 
function in $L^2$. In particular, apart from this property, we do not even use that~$u$ is a solution of an equation. In the particular setting of homogenization studied here, in which $u$ is a solution of the equation~\eqref{e.divformeq} with coefficients~$\a(\cdot)$ sampled by~$\P$, this property is essentially implied by Theorem~\ref{t.DP}. To make this more transparent, we present in Proposition~\ref{p.harmonicapproximation} below a reformulation of the latter result.

\smallskip

A conceptually important point is that homogenization allows us to prove the harmonic approximation only on length scales larger than a certain random length scale~$\X_s$ which is roughly of size~$O(1)$, the order of the correlation length scale. Since we cannot expect homogenization to give us information on scales smaller than~$O(1)$, we cannot expect harmonic approximation to hold on these scales. Thus, in contrast to the Schauder theory---which controls small scales but not larger scales---here we are able to control large scales but not small scales.

\begin{proposition}[{Harmonic approximation}]
\label{p.harmonicapproximation}
Fix $s\in (0,d)$. There exist constants $\alpha(d,\Lambda)>0$, $C(s,d,\Lambda)<\infty$,  and a random variable $\X_s:\Omega \to [1,\infty]$ satisfying 
the estimate
\begin{equation}
\label{e.Xharmonicapproximation}
\X_s = \O_s\left( C \right),
\end{equation}
such that the following holds: for every $R\geq \X_s$ and weak solution $u\in H^1(B_{R})$ of
\begin{equation} 
\label{e.wksrt}
-\nabla \cdot \left( \a\nabla u \right) = 0 \quad \mbox{in} \ B_{R},
\end{equation}
there exists a solution $\overline{u} \in H^1(B_{R/2})$ satisfying 
\begin{equation*} \label{}
-\nabla \cdot \left( \ahom \nabla \bar u \right) = 0 \quad \mbox{in} \ B_{R/2}
\end{equation*}
and
\begin{equation}
\label{e.harmonicapproximation}
\left\| u - \overline{u}  \right\|_{\underline{L}^2(B_{R/2})} \leq CR^{-\alpha(d-s)} \left\|  u - \left( u \right)_{B_{R}} \right\|_{\underline{L}^2(B_{R})}.
\end{equation}
\end{proposition}
\begin{proof} 
We first apply the interior Meyers estimate (Theorem~\ref{t.Meyers appendix}), which states that there exist $\delta(d,\Lambda)>0$ and $C(d,\Lambda)>0$ such that, for every $u\in H^1(B_{R})$ satisfying~\eqref{e.wksrt}, we have that $\nabla u \in L^{2+\delta}(B_{R/2})$ and 
\begin{equation}
\label{e.meyersrt}
\left\| \nabla u \right\|_{\underline{L}^{2+\delta}(B_{R/2})} \leq \frac{C}{R} \left\|  u - \left( u \right)_{B_{R}} \right\|_{\underline{L}^2(B_{R})}. 
\end{equation}
For each $s\in (0,d)$, we let $\tilde{\X}_s$ denote the random variable in the statement of Theorem~\ref{t.DP} with $\delta$ the exponent in the Meyers estimate, as above. \index{Meyers estimate}Recall from~\eqref{e.sizeX} that $\tilde{\X}_s \leq \O_1(C)$. 
We now define 
\begin{equation*} \label{}
\X_s := \tilde{\X}_{(s+d)/2}^{1/s} \vee 1
\end{equation*}
It is clear that~$\X_s = \O_s(C)$. Applying Theorem~\ref{t.DP} with~$\ep = R^{-1}$ and~$\tilde s = \frac{s+d}{2}$ gives us the desired conclusion for every~$R\geq \X_s$ if we take~$\overline{u}$ to be the solution of the Dirichlet problem for the homogenized equation with Dirichlet boundary condition~$u$ on~$\partial B_{R/2}$. 
\end{proof}

\section{\texorpdfstring{$C^{0,1}$}{C0,1}-type estimates}
\label{s.lipschitz}

In this section, we prove a~$C^{0,1}$-type estimate for~$\a(\cdot)$-harmonic functions which holds on \emph{large} length scales, that is, on length scales larger than some \emph{random scale}, which we denote by~$\X_s$. We also prove sharp estimates on the distribution of~$\X_s$ with respect to~$\P$. The arguments are here completely deterministic: the only probabilistic ingredient is Proposition~\ref{p.harmonicapproximation}. What we present here is a real variables argument for any~$L^2$-function allowing a good harmonic approximation, similar to what we saw in the previous section for Schauder estimates, but now valid only above a certain length scale.

\smallskip

The main result of this section is the following theorem.

\index{$C^{0,1}$ estimate}
\begin{theorem}
[{Quenched $C^{0,1}$-type estimate}]
\label{t.Lipschitz}
Fix $s\in (0,d)$. There exist a constant  $C(s,d,\Lambda)<\infty$ and a random variable $\X_s :\Omega \to [1,\infty]$ satisfying
\begin{equation}
\label{e.XLipschitz}
\X_s = \O_s\left( C \right),
\end{equation}
such that the following holds: for every $R\geq \X_s$ and weak solution $u\in H^1(B_{R})$ of
\begin{equation} 
\label{e.wksrt2}
-\nabla \cdot \left( \a\nabla u \right) = 0 \quad \mbox{in} \ B_{R},
\end{equation}
we have, for every  $r\in \left[ \X_s , R\right]$, the estimates 
\begin{equation}
\label{e.Lipschitz}
 \frac1{r} \left\| u - \left( u \right)_{B_{r}} \right\|_{\underline{L}^2(B_{r})}
\leq
\frac{C}{R} \left\| u - \left( u \right)_{B_{R}} \right\|_{\underline{L}^2(B_{R})}
\end{equation}
and
\begin{equation} 
\label{e.Linfty}
\left\| u  \right\|_{\underline{L}^2(B_{r})} \leq C  \left\| u \right\|_{\underline{L}^2(B_{R})}.
\end{equation} 
\end{theorem}

Observe that, by the Caccioppoli and Poincar\'e inequalities, we can write~\eqref{e.Lipschitz} equivalently (up to changing the constant $C$) as: 
\begin{equation}
\label{e.Lipschitz2}
\sup_{r\in \left[ \X_s , R\right]} \left\| \nabla u \right\|_{\underline{L}^2(B_{r})}
\leq
C \left\| \nabla u\right\|_{\underline{L}^2(B_{R})}.
\end{equation}
In other words, the energy density of the solution in a small ball centered at the origin is controlled by the average energy density of the solution on the largest scale if the radius of the ball is larger than~$\X_s$. This can be contrasted with Example~\ref{ex.meyers}, where we found the reverse inequality in~\eqref{e.nocoarselipschitz}. We may also give up a volume factor in the previous inequality to get an estimate in~$B_1$:
\begin{equation*} \label{}
\left\| \nabla u \right\|_{\underline{L}^2(B_{1})} \leq C \X_s^{\frac d2} \left\| \nabla u\right\|_{\underline{L}^2(B_{R})}.
\end{equation*}
Therefore Theorem~\ref{t.Lipschitz} gives very strong control of the gradients of solutions.

\smallskip

It may seem strange to call Theorem~\ref{t.Lipschitz} a ``$C^{0,1}$-type'' estimate since the conclusion involves~$L^2$--type bounds. The reason is that, from the point of view of homogenization, the quantities~$ \left\| \nabla u \right\|_{\underline{L}^2(B_{\X_s})}$ or~$\left\| \nabla u \right\|_{\underline{L}^2(B_{1})} $ should be considered as equivalent to $| \nabla u(0) |$. Similarly,~\eqref{e.Linfty} should be thought of as a bound for~$|u(0)|$. Indeed, on the one hand, homogenization theory is concerned only with the behavior of solutions on scales larger than the correlation length scale and cannot be expected to provide useful estimate on smaller scales: without further assumptions, our model could be a checkerboard with each cell containing a translate of the matrix in Example~\ref{ex.meyers}. On the other hand, if we impose the (quite reasonable) additional assumption that the coefficients~$\a(\cdot)$ are uniformly H\"older continuous, $\P$-a.s., then we may combine the Schauder estimates with~\eqref{e.Lipschitz2} to obtain the true pointwise estimate
\begin{equation} 
\label{e.schaudertopointwise}
| \nabla u(0) | 
\leq 
C \X_s^{\frac d2} \left\| \nabla u\right\|_{\underline{L}^2(B_{\X_s})}
\leq 
C\X_s^{\frac d2}  \left\| \nabla u\right\|_{\underline{L}^2(B_{R})}.
\end{equation}
However, we prefer to consider the problem of controlling the small scales to be a separate matter which does not concern us here, because a theory of homogenization should not ask questions about such small scales. In any case, we will discover that a microscopic-scale estimate like~\eqref{e.Lipschitz2} will serve us just as well as a pointwise bound like~\eqref{e.schaudertopointwise}.

\smallskip

Theorem~\ref{t.Lipschitz} is a consequence of Proposition~\ref{p.harmonicapproximation} and the following lemma.

\begin{lemma}
\label{l.harmapproxlipschitz}
Fix $\alpha \in (0,1]$,  $K\geq 1$ and $X\geq 1$. 
Let $R\geq 2X$ and $u\in L^2(B_{R})$ have the property that, for every $r\in \left[ X, \frac12R \right]$, there exists $w\in H^1(B_{r})$ which is a weak solution of 
\begin{equation*} \label{}
-\nabla \cdot \left( \ahom \nabla w \right) = 0 \quad \mbox{in} \ B_{r},
\end{equation*}
satisfying 
\begin{equation} 
\label{e.wharmapproxit}
\left\| u - w \right\|_{\underline{L}^2(B_{r/2})} \leq K r^{-\alpha} \left\| u - \left( u \right)_{B_{r}} \right\|_{\underline{L}^2(B_{r})}. 
\end{equation}
Then there exists $C(\alpha,K,d,\Lambda)<\infty$ such that, for every $r\in \left[ X, R\right]$, 
\begin{equation}
\label{e.elipsc}
 \frac1{r} \left\| u - \left( u \right)_{B_{r}} \right\|_{\underline{L}^2(B_{r})}
\leq
\frac{C}{R} \left\| u - \left( u \right)_{B_{R}} \right\|_{\underline{L}^2(B_{R})}.
\end{equation}
\end{lemma}


Before presenting the proof of Lemma~\ref{l.harmapproxlipschitz}, we give some further comments on the regularity of $\ahom$-harmonic functions. As already mentioned before, $\ahom$-harmonic functions are smooth; indeed, they are even real analytic. For our purposes, it is convenient to quantify the smoothness via $L^2$-integrals. We introduce the linear vector space
\begin{equation} \label{e.Ahom000}
\Ahom_k := \left\{ u \in H_{\textrm{loc}}^1(\R^d) \, : \, -\nabla \cdot\left( \ahom \nabla  u \right) = 0, \ \  \lim_{r \to \infty} r^{-(k+1)} \left\|  u \right\|_{\underline{L}^2(B_r)} = 0  \right\}.
\end{equation}
The $C^{k,1}$ estimate~\eqref{e.pointwiseharm} for $\ahom$-harmonic functions implies that, for every $u \in \Ahom_k$,
\begin{align*} 
\left\| \nabla^{k+1} u \right\|_{L^\infty(\R^d)} \leq \limsup_{R \to \infty} \left\| \nabla^{k+1} u \right\|_{L^\infty(B_{R})} \leq C_{k+1} \limsup_{R \to \infty} R^{-k-1} \left\| u \right\|_{\underline{L}^2(B_{R})} = 0. 
\end{align*}
Thus $\Ahom_k$ coincides with the set of $\ahom$-harmonic polynomials of degree at most~$k$. 
\index{harmonic polynomial}
Conversely,~\eqref{e.pointwiseharm}  also asserts that the local behavior of an arbitrary $\ahom$-harmonic function can be described in terms of~$\Ahom_k$. To see this, let~$w$ be $\ahom$-harmonic in $B_r$. Then, for any $k\in \N$, there exists a constant $C_k(d,\Lambda)\in [1,\infty)$ and $p \in \Ahom_k$ such that, for every $t \in (0,\frac12 r]$, 
\begin{equation} 
\label{e.Ahomharm}
\left\|w -  p \right\|_{L^\infty(B_t)} \leq C_k \left(\frac{t}{r}\right)^{k+1}  \left\|w -  p \right\|_{\underline{L}^2(B_r)}.
\end{equation}
We leave the proof of this fact from~\eqref{e.pointwiseharm} as an exercise. 
\begin{exercise}
Let $w \in H^1(B_R(y))$ be $\ahom$-harmonic. Prove using~\eqref{e.pointwiseharm} that, for every $k \in \N$, there exists $p \in \Ahom_k$ and a constant $C_k$ such that~\eqref{e.Ahomharm} holds. Hint: Prove that the Taylor expansion of $w$ at $y$ is $\ahom$-harmonic by an induction on the order of the expansion, using blow-ups, homogeneity and~\eqref{e.pointwiseharm}, and then conclude with the aid of~\eqref{e.pointwiseharm}.
\end{exercise}

We next extract from the pointwise regularity estimate~\eqref{e.pointwiseharm} some properties of~$\ahom$-harmonic polynomials which will be useful below. Fix $k \in \N$ and $p \in \Ahom_k$. Then~\eqref{e.pointwiseharm} gives, for every $r\geq 1$,
\begin{equation} \label{e.nicecontrolonp}
\sum_{m=0}^k r^m \left\| \nabla^m p \right\|_{L^\infty(B_r)} \leq C \left\| p \right\|_{\underline{L}^2(B_{2r})} \leq 2^k C  \left\| p \right\|_{\underline{L}^2(B_{r})}.
\end{equation}
In what follows, we need some control on the homogeneous parts of $p$. For this, let $\pi_{m,y} p$ stand for the Taylor expansion of $p$ at $y$ of degree $m$, that is 
\begin{equation}  \label{e.projection pi}
(\pi_{m,y} p)(x) := \sum_{j=0}^m \frac1{j!} \nabla^jp(y) (x-y)^{\otimes j} =  \sum_{j=0}^m \sum_{|\alpha| = j} \frac1{\alpha!} \partial^\alpha p(y) (x-y)^\alpha\,.
\end{equation}
We have by~\eqref{e.pointwiseharm} that
\begin{equation*} 
\left\| (\pi_{m,y} - \pi_{m-1,y})p \right\|_{\underline{L}^2(B_{r}(y))} \leq r^{m} | \nabla^m p(y) | \leq  C_m  \left\| p \right\|_{\underline{L}^2(B_{r}(y))}.
\end{equation*}
This implies the equivalence between the following two norms of $\Ahom_k$:
\begin{equation} \label{e.equivalenceofnorms}
\left\| p \right\|_{\underline{L}^2(B_{r}(y))} \leq \sum_{m=0}^k \left\| (\pi_{m,y} - \pi_{m-1,y})p \right\|_{\underline{L}^2(B_{r}(y))} \leq C_k \left\| p \right\|_{\underline{L}^2(B_{r}(y))}.
\end{equation}
We remark that, in the case $\ahom=\Id$ and for any~$r\geq 1$, the constant $C_k$ in~\eqref{e.equivalenceofnorms} can be taken to be~$1$ due to the orthogonality, with respect to the usual inner product of $L^2(B_r)$, of homogeneous harmonic polynomials of different degrees.

\smallskip

We next show that the strong regularity described in~\eqref{e.Ahomharm} can be transferred to a function~$u$ 
which possesses good harmonic approximation. While in this section we are interested only in the case $k=1$, we state an estimate for general $k\in\N$ since it will be useful in the characterization of~$\A_k$ in the next section. 

\begin{lemma}
\label{l.u:decayestimate}
Fix $\alpha \in [0,1]$,  $K\geq 1$ and $X\geq 1$. 
Let $R\geq 2X$ and $u\in L^2(B_{R})$ have the property that, for every $r\in \left[ X, R \right]$, there exists $w_r \in H^1(B_{r/2})$ which is a weak solution of 
\begin{equation*} \label{}
-\nabla \cdot \left( \ahom \nabla w_r \right) = 0 \quad \mbox{in} \ B_{r/2},
\end{equation*}
satisfying 
\begin{equation} 
\label{e.wharmapproxit02}
\left\| u - w_r \right\|_{\underline{L}^2(B_{r/2})} \leq K r^{-\alpha} \left\| u - \left( u \right)_{B_{r}} \right\|_{\underline{L}^2(B_{r})}. 
\end{equation}
Then, for every $k \in \N$, there exists $\theta(\alpha,k,d,\Lambda) \in (0,\frac12)$  and $C(\alpha,k,d,\Lambda)<\infty$ such that, for every $r\in \left[ X, R\right]$, 
\begin{equation}
\label{e.u:decayestimate}
\inf_{p \in \Ahom_k } \left\| u - p \right\|_{\underline{L}^2(B_{\theta r})}
\leq
\frac14 \theta^{k+1 - \alpha/2} \inf_{p \in \Ahom_k} \left\| u - p  \right\|_{\underline{L}^2(B_{r})}   + C K r^{-\alpha} \left\| u - (u)_{B_r} \right\|_{\underline{L}^2(B_{r})} .
\end{equation}
\end{lemma}
\begin{proof} Let $t \in \left(0,\frac14 r \right]$.  Then, by the triangle inequality and~\eqref{e.Ahomharm}, we get
\begin{align*}  
\lefteqn{\inf_{p \in \Ahom_k } \left\| u - p \right\|_{\underline{L}^2(B_{t})}} \qquad  & 
\\ & \leq  
 \inf_{p \in \Ahom_k}  \left\|w_r -  p \right\|_{\underline{L}^2(B_t)}  +  \left\|w_r -  u \right\|_{\underline{L}^2(B_t)} \\
& \leq C_k \left(\frac{t}{r}\right)^{k+1}  \inf_{p \in \Ahom_k}  \left\|w_r -  p \right\|_{\underline{L}^2(B_{r/2})} +  \left(\frac{r}{t} \right)^{\frac d2} \left\|w_r -  u \right\|_{\underline{L}^2(B_{r/2})}  \\
 & \leq C_k \left(\frac{t}{r}\right)^{k+1} \inf_{p \in \Ahom_k}  \left\|u -  p \right\|_{\underline{L}^2(B_{r/2})}   +   
\left(C_k \left(\frac{t}{r}\right)^{k+1}  + \left(\frac{r}{t} \right)^{\frac d2} \right)\left\|w_r -  u \right\|_{\underline{L}^2(B_{r/2})}.
\end{align*}
We choose $t = \theta r$ with $\theta  = (C_k 2^{2+\frac d2})^{-2/\alpha}$, and use~\eqref{e.wharmapproxit02} and an elementary fact
\begin{equation*} 
\inf_{p \in \Ahom_k}  \left\|u -  p \right\|_{\underline{L}^2(B_{r/2})} \leq 2^{\frac d2} \inf_{p \in \Ahom_k}   \left\|u -  p \right\|_{\underline{L}^2(B_{r})} 
\end{equation*}
to obtain the statement.
\end{proof}

We are now ready to give the proof of~Lemma~\ref{l.harmapproxlipschitz}. 

\begin{proof}[Proof of Lemma~\ref{l.harmapproxlipschitz}]
It is typically difficult to obtain a $C^{0,1}$ estimate by performing a decay of oscillation estimate. Instead we try to perform a $C^{1,\beta}$-type ``improvement of flatness'' iteration. In other words, we measure the distance between our solution and the closest \emph{affine} function rather than the closest constant function (as in the previous section). Because we are actually interested in the oscillation of the function and not its flatness (since we are \emph{not} trying to get a $C^{1,\beta}$ estimate), we must keep track of the size of the slopes of the approximating affine functions. This makes the argument a little more technical.

\smallskip

We assume that $R\geq R_0 := X \vee H$, where $H$ is a large constant to be selected below and which will depend on~$(\alpha,K,d,\Lambda)$. Fix $u\in L^2(B_{R})$ satisfying the hypothesis of the lemma. It is convenient to rephrase the claim of the lemma in terms of the~\emph{maximal function}:
\begin{equation*} 
M := \sup_{r \in [R_0,R]} \frac1r \left\| u - (u)_{B_r} \right\|_{\underline{L}^2(B_r)}.
\end{equation*}
 Our goal is to prove the bound
\begin{equation} \label{e.Mgoal}
M \leq \frac{C}{R} \left\| u - \left( u \right)_{B_{R}} \right\|_{\underline{L}^2(B_{R})}.
\end{equation}
Indeed, if $X\geq H$, then~\eqref{e.Mgoal} is the same as~\eqref{e.elipsc}. On the other hand, if $X \leq H$, then we have $R_0 = H$ and 
\begin{equation*} 
\sup_{r \in [X,H]} \frac1r \left\| u - (u)_{B_r} \right\|_{\underline{L}^2(B_r)} \leq 
H^{1+\frac{d}{2}} M,
\end{equation*}
which implies~\eqref{e.elipsc} since $H^{1+\frac{d}{2}} \leq C$. It therefore suffices to prove~\eqref{e.Mgoal}. Notice that we have here appealed to the following elementary fact: for every $0<r<t<\infty$, 
\begin{equation*} 
\left\| u - (u)_{B_r} \right\|_{\underline{L}^2(B_r)} \leq \left\| u - (u)_{B_t} \right\|_{\underline{L}^2(B_r)} \leq \left(\frac tr\right)^{\frac d2} \left\| u - (u)_{B_t} \right\|_{\underline{L}^2(B_t)},
\end{equation*}
where the first inequality follows by the fact that $(u)_{B_r} $ is realizing the infimum in $\inf_{a \in \R} \left\| u - a \right\|_{\underline{L}^2(B_r)}$ and the second one follows by giving up a volume factor. 

\smallskip

\emph{Step 1.}
We first set up the argument.  Using the definition of $M$,~\eqref{e.wharmapproxit} implies that, for every $r \in [R_0,R]$, there exists an $\ahom$-harmonic function $w_r$ such that 
\begin{equation}  \label{e.w_r}
\left\| u - w_r \right\|_{\underline{L}^2(B_{r/2})}  \leq  K r^{1-\alpha} M\,.
\end{equation} 
Denote the $L^2$-flatness of $u$ in $B_r$ by 
\begin{equation} 
\label{e.Fflatness}
E_1(r):= \frac1{r} \inf_{p \in \Ahom_1} \left\| u - p \right\|_{\underline{L}^2(B_r)} .
\end{equation}
Then Lemma~\ref{l.u:decayestimate} yields, for $k=1$, 
\begin{equation}
\label{e.E_1(r)}
E_1(\theta r) \leq \frac12 E_1(r)   + C K r^{-\alpha} M,
\end{equation}
where both $\theta$ and $C$ depend on $\alpha,d,\Lambda$. Define, for $j \in \N$, 
\begin{equation*} \label{}
r_j  := \theta^j R
\end{equation*}
so that the previous inequality reads, for every $r_j\geq 2X$, as
\begin{equation*} 
E_1(r_{j+1}) \leq \frac12 E_1(r_j) +  C r_j^{-\alpha} M.
\end{equation*}
Denote by $J$ the largest integer such that $r_J \geq R_0$. Summing over $j$ then leads to 
\begin{equation*} 
\sum_{j=1}^{J+1} E_1(r_{j}) \leq \frac12 \sum_{j=0}^{J}  E_1(r_j) +   C H^{-\alpha} M \sum_{j=0}^{J} \theta^{j \alpha}.
\end{equation*}
We may reabsorb the sum on the right and get
\begin{equation} \label{e.E_1 sum}
\sum_{j=0}^{J+1} E_1(r_{j}) \leq  2 E_1(R) +  \frac{2 C}{1-\theta^\alpha} H^{-\alpha} M \leq 2 E_1(R) +  C H^{-\alpha} M.
\end{equation}
The small number $H^{-\alpha}$ will allow us to reabsorb $M$ later in the proof.

\smallskip

\emph{Step 2.} Let $p_j \in \Ahom_1$ denote (the unique) affine function realizing the infimum in $E_1(r_j)$.  It is easy to see that $p_j$ must be of the form $x\mapsto (u)_{B_{r_j}}  + \nabla p_j \cdot x$, and hence 
\begin{equation*}  \label{e.affinebounds}
c |\nabla p_j| \leq  \frac1{r_j} \left\|  p_j -  (u)_{B_{r_j}}  \right\|_{\underline{L}^2(B_{r_j})} \leq C |\nabla p_j|.
\end{equation*}
It follows that
\begin{equation*} 
\frac1{r_j} \left\| u -  (u)_{B_{r_j}}  \right\|_{\underline{L}^2(B_{r_j})} \leq E_1(r_j) + C  |\nabla p_j|.
\end{equation*}
This gives us a way to estimate the maximal function $M$:
\begin{equation*} 
M = \sup_{r \in [R_0,R]} r^{-1} \left\| u - (u)_{B_r} \right\|_{\underline{L}^2(B_r)}
\leq \theta^{-(\frac{d}{2}+1)} \max_{j \in \{0,\ldots,J\}} \left( E_1(r_j) + |\nabla p_j|   \right).
\end{equation*}
By~\eqref{e.E_1 sum}, we have
\begin{equation*} 
\max_{j \in \{0,\ldots,J\}} E_1(r_j) \leq \sum_{j=0}^J E_1(r_j) \leq 2 E_1(R) +  C H^{-\alpha} M.
\end{equation*}
Thus it remains to estimate $|\nabla p_j|$. For this, observe that 
\begin{equation*} 
\left| \nabla p_j - \nabla p_{j+1}\right| \leq \frac{C}{r_{j+1}} \left\| p_j - p_{j+1} \right\|_{\underline{L}^2(B_{r_{j+1}})} \leq \theta^{-(\frac{d}{2}+1)} \left(E_1(r_{j+1}) + E_1(r_{j}) \right).
\end{equation*}
Summing the previous inequality over $j\in\{ 0,\ldots,J\}$ yields
\begin{equation*} 
\max_{j \in \{0,\ldots,J\}} \left| \nabla p_j \right| 
\leq \left| \nabla p_0 \right|  +  2 \theta^{-(\frac{d}{2}+1)} \sum_{j=0}^{J+1}  E_1(r_{j})
\leq \left| \nabla p_0 \right|  +  C\sum_{j=0}^{J+1}  E_1(r_{j})
\end{equation*}
where the last term on the right can again be estimated with the aid of~\eqref{e.E_1 sum}. Our initial effort to estimate the sum on the right stems precisely from the inequality above. The term $\left| \nabla p_0 \right|$ on the right can be easily treated with the aid of the triangle inequality as follows:
\begin{align*} 
|\nabla p_0|  & \leq \frac CR \left\|  p_0 -  (u)_{B_{R}}  \right\|_{\underline{L}^2(B_{R})}
\\ & \leq \frac CR \left\|  u - p_0\right\|_{\underline{L}^2(B_{R})}   
+ \frac CR \left\|  u - (u)_{B_{R}}\right\|_{\underline{L}^2(B_{R})}
\\ & \leq CE_1(R) + \frac{C}{R} \left\| u - \left( u \right)_{B_{R}} \right\|_{\underline{L}^2(B_{R})} 
\\ & \leq \frac{C}{R} \left\| u - \left( u \right)_{B_{R}} \right\|_{\underline{L}^2(B_{R})}.
\end{align*}
Combining the above estimates leads to
\begin{equation*} 
M\leq \frac{C}{R} \left\| u - \left( u \right)_{B_{R}} \right\|_{\underline{L}^2(B_{R})}  + C H^{-\alpha} M.
\end{equation*}
Choosing now $H: = (2C)^{1/\alpha}$ allows to reabsorb the last term and yields~\eqref{e.Mgoal}.

\smallskip

\emph{Step 3.} We finally demonstrate~\eqref{e.Linfty}. From~\eqref{e.Lipschitz} it follows that, for $N \in \N$ such that  $2^{-N}R < r \leq 2^{-N+1}R$,
\begin{equation*} 
\sum_{j=0}^N \left\| u - (u)_{B_{2^{-j} R}}\right\|_{\underline{L}^2(B_{2^{-j} R})} \leq C \left\| u - (u)_{B_{R}}\right\|_{\underline{L}^2(B_{ R})}.
\end{equation*}
Using Jensen's inequality we get 
\begin{equation*} 
\left| (u)_{B_{2^{-j-1} R}} - (u)_{B_{2^{-j} R}} \right|  \leq 2^{\frac d2} \left\| u - (u)_{B_{2^{-j} R}}\right\|_{\underline{L}^2(B_{2^{-j} R})},
\end{equation*}
and thus by the two previous displays we obtain
\begin{equation*} 
\sum_{j=0}^{N-1} \left| (u)_{B_{2^{-j-1} R}} - (u)_{B_{2^{-j} R}} \right| \leq C \left\| u \right\|_{\underline{L}^2(B_{ R})}.
\end{equation*}
Consequently,
\begin{equation*} 
\sup_{j \in \{0,\ldots,N\}} \left| (u)_{B_{2^{-j}R}}\right|  \leq C \left\| u \right\|_{\underline{L}^2(B_{ R})},
\end{equation*}
from which~\eqref{e.Linfty} follows easily using~\eqref{e.Lipschitz}. The proof is complete.
\end{proof}

\begin{exercise}
\label{ex.C01.RHS}
Generalize Theorem~\ref{t.Lipschitz} to allow for right-hand sides by proving the following statement:
for every $s\in (0,d)$, there exist~$C(s,d,\Lambda)<\infty$  and a random variable $\X_s :\Omega \to [1,\infty]$ satisfying
\begin{equation}
\label{e.XLipschitz000}
\X_s = \O_s\left( C_1 \right),
\end{equation}
such that the following holds: for every $R\geq \X$, $f \in H^{-1}(B_{2R})$, and weak solution $u\in H^1(B_{2R})$ of
\begin{equation} 
\label{e.wksrt2000}
-\nabla \cdot \left( \a\nabla u \right) = f \quad \mbox{in} \ B_{2R},
\end{equation}
we have, for $r \geq \X_s$, 
\begin{equation}
\label{e.Lipschitz000}
 \sup_{t \in (r ,R )} \frac1{t} \left\| u - \left( u \right)_{B_{t}} \right\|_{\underline{L}^2(B_{t})}
 \leq \frac{C}{R}  \left\| u - \left( u \right)_{B_{R}} \right\|_{\underline{L}^2(B_{R})} + C \int_r^R \frac1s \left\| f \right\|_{\underline{H}^{-1}(B_{s})}\,ds.
\end{equation}

\noindent Hint: Obtain a harmonic approximation statement for $u$ by just throwing away the $f$. In other words, show first that, in every ball $B_r$ with $r \leq R$,
\begin{equation*} \label{}
\left\| \nabla u - \nabla \tilde{u}_r \right\|_{\underline{L}^2(B_r)} 
\leq
C \left\| f \right\|_{\underline{H}^{-1}(B_r)}
\end{equation*}
where $\tilde{u}_r$ satisfies
\begin{equation*} \label{}
\left\{
\begin{aligned}
& -\nabla \cdot \left( \ahom \nabla \tilde{u}_r \right) = 0 & \mbox{in} & \ B_r, \\
& u_r = u & \mbox{on} & \ \partial B_r. 
\end{aligned}
\right. 
\end{equation*}
We then pick up an extra error term in our harmonic approximation, namely $C \left\| f \right\|_{\underline{H}^{-1}(B_r)}$. We then rerun the iteration in the proof of Lemma~\ref{l.harmapproxlipschitz} and observe that, as we sum this over the scales, these extra terms are bounded by the second term on the right side of~\eqref{e.Lipschitz000}. 

\smallskip

Notice also that, for every $p\in(d,\infty]$, there exists~$C(p,d)<\infty$ such that 
\begin{equation*} \label{}
 \int_r^R \frac1s \left\| f \right\|_{\underline{H}^{-1} (B_{s})}\,ds \leq C R \left\| f \right\|_{\underline{L}^p(B_R)}.
\end{equation*}

\end{exercise}

\section{Higher-order regularity theory and Liouville theorems}
\label{s.higherreg}

A classical fact in partial differential equations is that the regularity theory can be used to classify entire solutions (solutions in the whole $\R^d$) which have at most polynomial growth. We encountered this idea for $\ahom$-harmonic functions in the previous section in the characterization of~$\Ahom_k$. These classification theorems are usually called \emph{Liouville theorems}.\index{Liouville theorem}

\smallskip

The purpose of this section is to show that appropriate versions of such Liouville theorems and of~\eqref{e.Ahomharm} continue to hold (with high $\P$ probability) in the case of heterogeneous equations. For this purpose, it is necessary to generalize the notion of $\Ahom_k$. It is natural therefore to define the vector space of solutions of~\eqref{e.divformeq} which grow at most like $o(|x|^{k+1})$ as $|x| \to \infty$:
\begin{equation} \label{e.A000}
\A_k := \left\{ u \in \A(\R^d) \, : \, \limsup_{r \to \infty} r^{-(k+1)} \left\|  u \right\|_{\underline{L}^2(B_r)} = 0  \right\}.
\end{equation}
The next theorem is the main result of this chapter. In particular, it supersedes Theorem~\ref{t.Lipschitz}. It identifies the dimension of the space $\A_k$, relates it to $\Ahom_k$ and gives an appropriate version of~\eqref{e.Ahomharm}. 
\index{$C^{k,1}$ estimate}

\begin{theorem}[Higher regularity theory]
\label{t.regularity}
\index{Liouville theorem}
Fix $s \in (0,d)$. There exist an exponent $\delta(s,d,\Lambda)\in \left( 0, \frac12 \right]$ and a random variable $\X_s$ satisfying the estimate
\begin{equation}
\label{e.X}
\X_s \leq \O_s\left(C(s,d,\Lambda)\right)
\end{equation}
such that the following statements hold, for every $k\in\N$:
\begin{enumerate}
\item[{$\mathrm{(i)}_k$}] There exists $C(k,d,\Lambda)<\infty$ such that, for every $u \in \A_k$, there exists $p\in \overline{\A}_k$ such that, for every $R\geq \X_s$,
\begin{equation} \label{e.liouvillec}
\left\| u - p \right\|_{\underline{L}^2(B_R)} \leq C R^{-\delta} \left\| p \right\|_{\underline{L}^2(B_R)}.
\end{equation}

\item[{$\mathrm{(ii)}_k$}]For every $p\in \overline{\A}_k$, there exists $u\in \A_k$ satisfying~\eqref{e.liouvillec} for every $R\geq \X_s$. 

\item[{$\mathrm{(iii)}_k$}]
There exists $C(k,d,\Lambda)<\infty$ such that, for every $R\geq \X_s$ and $u\in \A(B_R)$, there exists $\phi \in \A_k$ such that, for every $r \in \left[ \X_s,  R \right]$, we have the estimate
\begin{equation}
\label{e.intrinsicreg}
\left\| u - \phi \right\|_{\underline{L}^2(B_r)} \leq C \left( \frac r R \right)^{k+1} 
\left\| u \right\|_{\underline{L}^2(B_R)}.
\end{equation}
\end{enumerate}
In particular, $\P$-almost surely, we have, for every $k\in\N$,
\begin{equation} 
\label{e.dimensionofAk}
\dim(\A_k) = \dim(\Ahom_k) =  \binom{d+k-1}{k} + \binom{d+k-2}{k-1}.
\end{equation}
\end{theorem}

Theorem~\ref{t.regularity} allows us to think of the finite-dimensional space $\A_k$ in the same way that we think of harmonic polynomials. The third statement says that we can approximate an arbitrary solution of our equation by these ``heterogeneous polynomials'' to the same precision that an analytic function can be approximated by true polynomials. This fact has important consequences and, from a computational perspective, greatly reduces the complexity of the space of solutions. We think of~{$\mathrm{(iii)}_k$} as a ``large-scale $C^{k,1}$ estimate.'' 

\smallskip

The three basic tools we use to prove Theorem~\ref{t.regularity} are Theorem~\ref{t.DP}, Proposition~\ref{p.harmonicapproximation} and Lemma~\ref{l.u:decayestimate}. 
Throughout we fix $s \in (0,d)$ and 
\begin{equation} 
\label{e.r0.def}
r_0  := H \vee \X_s,
\end{equation}
where the random variable $\X_s$ is the maximum of the respective random variables from Theorem~\ref{t.DP} and Proposition~\ref{p.harmonicapproximation}, with the given $s \in (0,d)$. Observe, in particular, that when $r\geq r_0$, we may apply both of them. The large constant $H$ will be fixed later on. Notice that it is in fact enough to prove the statement for $r \geq r_0$, since if $\X_s \leq H$, we may simply enlarge the constant in each statement by the prefactor $H^{\frac d2}$ to obtain the result. We will give several conditions on $H$ during the proof, but in the end it can be chosen to depend only on $s,k,d,\Lambda$. The parameter $\delta$ in the statement is fixed to be
\begin{equation}
\label{e.deltaexplicit}
\delta := \alpha(d-s) \wedge \frac12,
\end{equation}
where~$\alpha(d,\Lambda)>0$ is as in the statement of Proposition~\ref{p.harmonicapproximation}.

\smallskip

The proof of Theorem~\ref{t.regularity} is an induction on~$k\in\N$. As will be explained below, Theorem~\ref{t.Lipschitz} is actually the base case, that is, it is equivalent to the statement of Theorem~\ref{t.regularity} for~$k=0$.

\smallskip

To set up the induction argument, it is necessary to formulate a slightly weaker version of $\mathrm{(iii)}_k$. The reason has to do with the trouble of getting on top of the exponent $\beta=1$ in a $C^{k,\beta}$ estimate, as discussed in the previous section in the case $k=0$. We therefore denote by $\mathrm{(iii')}_k$ the statement that 
\begin{enumerate}
\item[{$\mathrm{(iii')}_k$}] There exists $C(s, k,d,\Lambda)<\infty$ such that, for every $R\geq r_0$ and $u\in \A(B_R)$, there exists $\phi \in \A_k(\Rd)$ such that, for every $r \in \left[r_0, R \right]$, we have the estimate
\begin{equation}
\label{e.intrinsicregprime}
\left\| u - \phi \right\|_{\underline{L}^2(B_r)} \leq C \left( \frac r R \right)^{k+1 - \delta/2} 
\left\| u \right\|_{\underline{L}^2(B_R)}.
\end{equation}
\end{enumerate}
Theorem~\ref{t.regularity} is a consequence of Theorem~\ref{t.Lipschitz} and the following four implications: for each $k\in\N$,
\begin{gather}
\mbox{$\mathrm{(i)}_{k-1}$,\, $\mathrm{(ii)}_{k-1}$ \ \mbox{and} \ $\mathrm{(iii')}_{k-1} \implies \mathrm{(ii)}_{k}$ } \label{e.findmeap} \\
\mbox{$\mathrm{(i)}_{k-1}$ \ \mbox{and} \ $\mathrm{(ii)}_{k}  \implies \mathrm{(i)}_{k}$ } 
\label{e.findmeau} \\
\mbox{$\mathrm{(i)}_{k}$ \ \mbox{and} \ $\mathrm{(ii)}_{k}  \implies \mathrm{(iii')}_{k}$ }
\label{e.findmeaphi}
\end{gather}
and
\begin{equation} 
\label{e.climbontop}
\mbox{$\mathrm{(i)}_{k+1}$,\, $\mathrm{(ii)}_{k+1}$ \ \mbox{and} \ $\mathrm{(iii')}_{k+1} \implies \mathrm{(iii)}_{k}$. }
\end{equation}
We present the proof of each of these four implications separately, followed by the proof of~\eqref{e.dimensionofAk}. We first argue that Theorem~\ref{t.Lipschitz} does indeed give us the base case.

\begin{proof}[{Proof of~$\mathrm{(i)}_{0}$,~$\mathrm{(ii)}_{0}$ and~$\mathrm{(iii)}_{0}$.}]
Observe that~$\mathrm{(ii)}_{0}$ is trivial, since $\Ahom_0$ contains only constants, and hence $\Ahom_0 \subset \A_0$. Property~$\mathrm{(iii)}_{0}$, on the other hand, is a restatement of Theorem~\ref{t.Lipschitz}. Finally, the property~$\mathrm{(i)}_{0}$ follows also from Theorem~\ref{t.Lipschitz} since, for every~$u\in \A_0$ and~$r \geq r_0$, 
\begin{equation*} 
\frac1r \left\| u - (u)_{B_r}\right\|_{\underline{L}^2(B_r)} \leq C \limsup_{R\to \infty} \frac1R \left\| u - (u)_{B_R}\right\|_{\underline{L}^2(B_R)} = 0.
\end{equation*}
We deduce that $\A_0 = \Ahom_0$ is the set of constant functions.
\end{proof}

\smallskip

We now proceed with the proofs of the assertions~\eqref{e.findmeap},~\eqref{e.findmeau},~\eqref{e.findmeaphi} and~\eqref{e.climbontop}.
We first establish a basic decay estimates used in the proof. Recall that~$\delta(s,d,\Lambda)>0$ is defined in~\eqref{e.deltaexplicit}. For each $k\in\N$ and $r>0$, we define
\begin{equation*} \label{}
\label{e.Dk.def}
E_k(r):= \frac1{r^k} \inf_{p \in \Ahom_k} \left\| u - p \right\|_{\underline{L}^2(B_r)} .
\end{equation*}
Proposition~\ref{p.harmonicapproximation} and Lemma~\ref{l.u:decayestimate} yield the existence of~$\theta(s,k,d,\Lambda) \in (0,1)$ such that, for every $r \geq r_0$, $u\in \A(B_r)$, $k \in \N$,
\begin{equation} \label{e.basicreg}
E_{k}(\theta r)  
 \leq \frac14 \theta^{1-\delta/2} E_{k}(r)  + C_{s,k} r^{-(k+\delta)} \left\| u \right\|_{\underline{L}^2(B_{r})}
\end{equation}
and
\begin{equation} \label{e.basicreg10}
E_{k+1}(\theta r)  
 \leq \frac14  \theta^{1-\delta/2} E_{k+1}(r)  + C_{s,k} r^{-(k+1+\delta)} \left\| u \right\|_{\underline{L}^2(B_{r})}.
\end{equation}
For each  $j\in\N$, we define
\begin{equation*}
r_{j} := \theta^{-j} r_0.
\end{equation*} 
Recall that~$r_0$ is the random variable defined in~\eqref{e.r0.def}.
Note that~$r_j$ is defined differently here compared to the previous section.

\smallskip

\begin{proof}[{Proof of~\eqref{e.findmeap}}]
Due to $\mathrm{(ii)}_{k-1}$ and~\eqref{e.equivalenceofnorms}, we may assume that  $q \in \Ahom_k$ is a homogeneous polynomial of degree $k$. By this we mean that $(\pi_{k}-\pi_{k-1}) q = q$, where $\pi_k = \pi_{k,0}$ is defined in~\eqref{e.projection pi}. Fix thus such $q \in \Ahom_k$. 

\smallskip

\emph{Step 1.}
First, by Theorem~\ref{t.DP}, after appropriate rescaling ($x \mapsto x/\ep$), we find, for each $m \in \N$, a solution $u_{m} \in \A(B_{r_m})$ such that 
\begin{equation} 
\label{e.u_m local}
\left\| u_m - q \right\|_{\underline{L}^2(B_{r_m})} 
\leq C r_m^{-\delta} \left\| q \right\|_{\underline{L}^2(B_{r_m})}.
\end{equation}
Indeed, we may take $u_m$ to be the solution of the Dirichlet problem 
\begin{equation*} \label{}
\left\{
\begin{aligned}
& - \nabla \cdot \left( \a(x) \nabla u_m \right) = 0 & \mbox{in} & \ B_{r_m}, \\
& u_m = q & \mbox{on} & \ \partial B_{r_m}. 
\end{aligned}
\right.
\end{equation*}
Theorem~\ref{t.DP} applied to $\tilde{u}_m:= u_m(r_m\cdot)$ for $\ep := r_m^{-1}$ yields the bound~\eqref{e.u_m local}. 

\smallskip

Now let $w_m := u_{m+1} - u_{m}$ and $\phi_m \in \A_{k-1}(\R^d)$ be given by $\mathrm{(iii')}_{k-1}$ for~$w_m$. Then~$\mathrm{(iii')}_{k-1}$ together with the triangle inequality and the previous display imply, for every $j \leq m$,
\begin{align*} 
\left\| w_m -\phi_m  \right\|_{\underline{L}^2(B_{r_j})}
 \leq C \left(\theta^{k-\delta/2}\right)^{m-j} \left\| w_m  \right\|_{\underline{L}^2(B_{r_m})} 
& \leq C r_{m}^{-\delta} \left(\theta^{k-\delta/2}\right)^{m-j} \left\| q  \right\|_{\underline{L}^2(B_{r_m})} 
\\ &
= C r_j^{-\delta} \left(\theta^{\delta/2} \right)^{m-j}  \left\| q  \right\|_{\underline{L}^2(B_{r_j})}.
\end{align*}
By $\mathrm{(i)}_{k-1}$, there exists $p_{m} \in \Ahom_{k-1}$ such that, for every $n\in\N$,
\begin{equation*} \label{}
\left\| p_{m} - \phi_m \right\|_{\underline{L}^2(B_{r_n})} \leq C r_n^{-\delta}  \left\| p_{m}  \right\|_{\underline{L}^2(B_{r_n})}. 
\end{equation*}
If $H$ in~\eqref{e.r0.def} is taken sufficiently large, depending only on $(\delta,k,d,\Lambda)$, then $Cr_0^{-\delta} \leq \frac12$. Then the triangle inequality and the two previous displays gives
\begin{align*} 
\left\| p_{m}   \right\|_{\underline{L}^2(B_{r_m})} 
& 
\leq 2 \left\| {\phi_m}   \right\|_{\underline{L}^2(B_{r_m})} 
\\ & 
\leq  2 \left\| u_{m+1} - u_{m}  \right\|_{\underline{L}^2(B_{r_m})} 
+ 2\left\| w_m -\phi_m  \right\|_{\underline{L}^2(B_{r_m})} 
\leq  C r_{m}^{-\delta} \left\| q  \right\|_{\underline{L}^2(B_{r_m})} \,. 
\end{align*}
Therefore we have, for every $n,m\in\N$ with $n>m$,
\begin{align*} 
\left\| \phi_m  \right\|_{\underline{L}^2(B_{r_n})} & \leq C  \left(\frac{r_n}{r_m} \right)^{k-1} \left\| p_{m}  \right\|_{\underline{L}^2(B_{r_m})} 
\\ & \leq C \theta^{(m-n)(k-1)} r_{m}^{-\delta} \left\| q  \right\|_{\underline{L}^2(B_{r_m})}  = C \theta^{(n-m)(1-\delta)}  r_n^{-\delta} \left\| q  \right\|_{\underline{L}^2(B_{r_n})}\,.
\end{align*}
Set $v_n := u_n - \sum_{m=1}^{n-1} \phi_m$, so that $v_n - u_j = \sum_{m=j}^{n-1} (w_m-\phi_m) - \sum_{m=1}^{j-1} \phi_m$. Then, for every $j\in\N$ with $j < n$, 
\begin{align} \label{e.v_n vs q}
\left\| v_n - q  \right\|_{\underline{L}^2(B_{r_j})}  & \leq \left\| u_j - q  \right\|_{\underline{L}^2(B_{r_j})} + \sum_{m=j}^{n-1} \left\| w_m -\phi_m \right\|_{\underline{L}^2(B_{r_j})} + \sum_{m=1}^{j-1} \left\| \phi_m \right\|_{\underline{L}^2(B_{r_j})}
\\ \nonumber & \leq C r_j^{-\delta} \left\| q  \right\|_{\underline{L}^2(B_{r_j})} \left(1 + \sum_{m=0}^{n-1-j} \left(\theta^{\delta/2} \right)^{m}   +  \sum_{m=1}^{j-1} \left( \theta^{1-\delta}\right)^{m} \right)
\\ \nonumber & \leq C r_j^{-\delta} \left\| q  \right\|_{\underline{L}^2(B_{r_j})}.
\end{align}
In the next step we will show that the previous estimate is enough to conclude. 

\smallskip

\emph{Step 2.}
We will show that, up to a subsequence, $u = \lim_{n\to \infty} v_n$ exists, belongs to $\A_k(\R^d)$, and satisfies~\eqref{e.liouvillec} for $r\geq r_0$. First, by the Caccioppoli estimate (cf. Lemma~\ref{l.Caccioppoli appendix}) we have that
\begin{equation*} 
\left\| \nabla v_n \right\|_{\underline{L}^2(B_{r_j/2})} \leq Cr_j^{-1} \left\| v_n \right\|_{\underline{L}^2(B_{r_j})}  \leq Cr_j^{-1}\left\| q  \right\|_{\underline{L}^2(B_{r_j})}. 
\end{equation*}
Weak compactness in $H^1$ and Rellich's theorem imply that there exists $\tilde v_j \in H^1(B_{r_j/2})$ and a subsequence $\{n_m\} \subset \N$ such that, as $m\to \infty$, 
\begin{equation*} 
\left\{ 
\begin{aligned}
& \nabla v_{n_m} \rightharpoonup  \nabla \tilde v_j  \quad \mbox{weakly in } L^2(B_{r_j/2}), \ \mbox{and} 
\\
& v_{n_m} \to \tilde v_j  \quad \mbox{strongly in } L^2(B_{r_j/2}).
\end{aligned}
\right.
\end{equation*}
The above convergence guarantees that $\tilde v_j \in \A(B_{r_{j}/2})$, and by~\eqref{e.v_n vs q} it moreover satisfies by the strong convergence in $L^2$, for every $h < j$, 
 \begin{equation*} 
\left\| \tilde v_j - q \right\|_{\underline{L}^2(B_{r_h})}  \leq Cr_h^{-\delta}\left\| q  \right\|_{\underline{L}^2(B_{r_h})}. 
\end{equation*}
Taking the above subsequence $\{n_m\}$, and repeating the argument in $B_{r_{j+1}}$ with this subsequence, gives $\tilde v_{j+1}$, and by the uniqueness of weak limits, $\tilde v_j = \tilde v_{j+1}$ in $B_{r_{j}/2}$. We can repeat the argument for every $j$, and appeal to a diagonal argument to obtain a subsequence $\{ n_j \}_{j\in\N}$ such that $u = \lim_{j\to \infty} v_{n_j}$ belongs to $\A_k$ and satisfies, for every $r\geq R_0$, 
\begin{equation*} 
\left\| u - q \right\|_{\underline{L}^2(B_{r})} \leq  C r^{-\delta} \left\| q \right\|_{\underline{L}^2(B_{r})}.
\end{equation*}
This completes the proof of~\eqref{e.findmeap} for general $k\in\N$. 
\end{proof}

\begin{proof}[{Proof of~\eqref{e.findmeau}}]  The proof of this implication is the most challenging one, and we divide it into several steps. Let us give the outline of the argument. The starting point is~\eqref{e.basicreg10} with $k+1$ degree polynomial approximations. We will use it together with the fact that $u \in \A_k$, implying that $\lim_{r \to \infty} E_{k+1}(r) = 0$.  Using this it is possible to show that $k+1$ degree derivatives of approximating polynomials are quantitatively small, see~\eqref{e.nabla^(k+1)p_m in control}. In Step 3 we will prove that the lower degree parts of approximating polynomials are small as well, and identify the $k$th degree part of the approximating polynomial via a suitable Cauchy sequence. Finally, in Step 4, we conclude using $\mathrm{(i)}_{k-1}$ and  $\mathrm{(ii)}_{k}$. 

\smallskip

\emph{Step 1.} We first initialize the argument.  Assume that $u \in \A_{k}$, that is, $u \in \A(\R^d)$ and
\begin{equation} \label{e.uinA_k}
\lim_{j \to \infty} r_j^{-(k+1)}\|u\|_{\underline{L}^2(B_{r_j})} = 0.
\end{equation}
Let $p_j\in \Ahom_{k+1}$ be the unique minimizer appearing in $E_{k+1}(r_j)$. Then, by the triangle inequality,
\begin{equation*} 
 r_j^{-(k+1)}\left\| u \right\|_{\underline{L}^2(B_{r_j})} \leq  E_{k+1}(r_j) +  C r_j^{-(k+1)} \left\| p_j\right\|_{\underline{L}^2(B_{r_j})}. 
\end{equation*}
Hence, if we choose the parameter $H$ to satisfy a condition 
\begin{equation} \label{e.H cond 1}
H \geq \left( 4 C_{s,k}   \right)^{\frac{1}{\delta}}  \implies C_{s,k} r_{j}^{-\delta} \leq \frac14 \,,
\end{equation}
we obtain from~\eqref{e.basicreg10} that
\begin{equation} \label{e.basicreg2}
E_{k+1}(r_j)  
 \leq \frac12 E_{k+1}(r_{j+1})  + C r_{j+1}^{-(k+1+\delta)} \left\| p_{j+1} \right\|_{\underline{L}^2(B_{r_{j+1}})}.
\end{equation}
From~\eqref{e.uinA_k} it follows that $\lim_{j \to \infty }E_{k+1}(r_j) = 0$, and therefore we may sum over $j$ in~\eqref{e.basicreg2}, and obtain after reabsorption that, for every $m\in\N$,
\begin{equation}  \label{e.infinite sum Dk 00}
\sum_{j=m}^\infty E_{k+1}(r_j) \leq  C \sum_{j=m+1}^\infty r_{j}^{-(k+1+\delta)} \left\| p_{j} \right\|_{\underline{L}^2(B_{r_{j}})}.
\end{equation}
Observe also that~\eqref{e.uinA_k} and~\eqref{e.nicecontrolonp} imply
\begin{equation} \label{e.nabla^(k+1) p_j to zero}
\limsup_{j \to \infty} 
\left| \nabla^{k+1} p_j \right|
\leq C \limsup_{j \to \infty}
r_j^{-(k+1)} \left\| p_{j} \right\|_{\underline{L}^2(B_{r_{j}})} \leq C
\lim_{j \to \infty} r_j^{-(k+1)}\|u\|_{\underline{L}^2(B_{r_j})}  = 0  .
\end{equation}

\smallskip

\emph{Step 2.} In this step we will show that $\nabla^{k+1} p_j$ is small for large $j$ with a quantified rate. 
First, the triangle inequality and~\eqref{e.nicecontrolonp} yield
\begin{equation*} 
|\nabla^{k+1} (p_j - p_{j+1})| \leq C r_j^{-(k+1)} \left\| p_{j} - p_{j+1} \right\|_{\underline{L}^2(B_{r_{j+1}})} \leq C\left( E_{k+1}(r_j) + E_{k+1}(r_{j+1})  \right),
\end{equation*}
and thus by~\eqref{e.nabla^(k+1) p_j to zero} and~\eqref{e.infinite sum Dk 00},
\begin{equation*} 
|\nabla^{k+1} p_m | \leq C \sum_{j=m}^\infty E_{k+1}(r_j) \leq C \sum_{j=m+1}^\infty r_{j}^{-(k+1+\delta)} \left\| p_{j} \right\|_{\underline{L}^2(B_{r_{j}})}.
\end{equation*}
Furthermore, denoting $\pi_k = \pi_{k,0}$,  we have that 
\begin{equation*} 
r_{j}^{-(k+1)}  \left\| p_{j} \right\|_{\underline{L}^2(B_{r_{j}})} 
\leq 
r_{j}^{-(k+1)}  \left\| \pi_k p_{j} \right\|_{\underline{L}^2(B_{r_{j}})} + |\nabla^{k+1} p_j |,
\end{equation*}
and therefore
\begin{equation*} 
|\nabla^{k+1} p_m | \leq C \sum_{j=m+1}^\infty r_{j}^{-(k+1+\delta)} \left\| \pi_k p_{j} \right\|_{\underline{L}^2(B_{r_{j}})} + C \sum_{j=m+1}^\infty r_j^{-\delta} |\nabla^{k+1} p_j |.
\end{equation*}
It follows that 
\begin{equation*} 
\sup_{j \geq m} |\nabla^{k+1} p_j| \leq C \sum_{j=m+1}^\infty r_{j}^{-(k+1+\delta)} \left\| \pi_k p_{j} \right\|_{\underline{L}^2(B_{r_{j}})}  + C  \sup_{j\geq m} |\nabla^{k+1} p_j |  \sum_{j=m+1}^\infty r_j^{-\delta}.
\end{equation*}
The last sum can be estimated, with large enough $H$ in the definition of $r_0$, as
\begin{equation} \label{e.H cond 2}
C \sum_{j=m}^\infty r_j^{-\delta} \leq \frac{C}{1-\theta^\delta} H^{-\delta} \leq \frac12.
\end{equation}
Therefore we obtain, after reabsorption,
\begin{equation} \label{e.nabla^(k+1)p_m in control}
 |\nabla^{k+1} p_m | \leq C \sum_{j=m}^\infty r_{j}^{-(k+1+\delta)} \left\| \pi_k p_{j} \right\|_{\underline{L}^2(B_{r_{j}})}.
\end{equation}

\emph{Step 3.}
We then proceed to identify the limit of $(\pi_k - \pi_{k-1}) p_j$ as $j \to \infty$. Setting 
\begin{equation*} 
\tilde E_k(r_m) := r_m^{-k} \left\| u - \pi_k p_{m} \right\|_{\underline{L}^2(B_{r_{m}})},
\end{equation*}
we deduce from~\eqref{e.infinite sum Dk 00} and~\eqref{e.nabla^(k+1)p_m in control} that 
\begin{equation}  \label{e.tilde E_k}
\tilde E_k(r_m) \leq r_m E_{k+1}(r_m) + Cr_m |\nabla^{k+1} p_m | \leq C r_m  \sum_{j=m}^\infty r_{j}^{-(k+1+\delta)} \left\| \pi_k p_{j} \right\|_{\underline{L}^2(B_{r_{j}})}.
\end{equation}
To obtain suitable estimates for $(\pi_k - \pi_{k-1}) p_j$, we need to analyze the sum of $\tilde E_k(r_m)$'s. Rearrange the last sum in~\eqref{e.tilde E_k} as
\begin{equation*} 
r_m  \sum_{j=m}^\infty r_{j}^{-(k+1+\delta)} \left\| \pi_k p_{j} \right\|_{\underline{L}^2(B_{r_{j}})} \leq r_{m}^{-(k+\delta)} \omega(m),
\end{equation*}
where we denote 
\begin{equation*} 
\omega(m) = \sum_{h=0}^k \omega_h(k) \,, \qquad  \omega_h(m) := \sum_{j=m}^\infty \left( \theta^{k+1+\delta} \right)^{j-m}  \left\| (\pi_h - \pi_{h-1}) p_{j} \right\|_{\underline{L}^2(B_{r_{j}})} ,
\end{equation*}
with a convention that $\pi_{-1} = 0$. Since~$(\pi_h-\phi_{h-1})q$ is a homogeneous polynomial of degree $h$, we have in particular that, for every~$q \in \Ahom_k$ and $r,t>0$, 
\begin{equation} \label{e.q grwoth 11}
\left\| (\pi_h - \pi_{h-1}) q \right\|_{\underline{L}^2(B_r)} = \left(\frac{r}{t}\right)^h  \left\| (\pi_h - \pi_{h-1}) q \right\|_{\underline{L}^2(B_t)} .
\end{equation}
Moreover, we have by the equivalence of norms described in~\eqref{e.equivalenceofnorms} that, for $q \in \Ahom_k$,
\begin{equation} \label{e.q grwoth 22}
\left\| q \right\|_{\underline{L}^2(B_r)}  
\leq \sum_{h=0}^k \left\| (\pi_h - \pi_{h-1}) q \right\|_{\underline{L}^2(B_r)}  
\leq C_k \left\| q \right\|_{\underline{L}^2(B_r)} .
\end{equation}
With this notation we get from~\eqref{e.tilde E_k} that 
\begin{equation} \label{e.tilde E_k 2}
\tilde E_k(r_m) \leq C r_m^{-(k+\delta)} \omega(m). 
\end{equation}
We next estimate the growth of $\omega(m)$. A first, crude bound is
\begin{equation}
\label{e.crude omega}
\omega(m+1) \leq \theta^{-(k+1+\delta)} \omega(m).
\end{equation}
By the triangle inequality,~\eqref{e.q grwoth 11} and~\eqref{e.q grwoth 22}, we get
\begin{align*} 
\omega_h (m+1) & =  \sum_{j=m+1}^\infty \left( \theta^{k+1+\delta} \right)^{j-(m+1)}  \left\| (\pi_h - \pi_{h-1}) p_{j} \right\|_{\underline{L}^2(B_{r_{j}})} 
\\ & = \sum_{j=m}^\infty \left( \theta^{k+1+\delta} \right)^{j-m}  \left\| (\pi_h - \pi_{h-1}) p_{j+1} \right\|_{\underline{L}^2(B_{r_{j+1}})} 
\\ & =  \theta^{-h} \sum_{j=m}^\infty \left( \theta^{k+1+\delta} \right)^{j-m} \left\| (\pi_h - \pi_{h-1}) p_{j+1} \right\|_{\underline{L}^2(B_{r_{j}})} 
\\ & \leq \theta^{-h} \omega_h(m) + C \sum_{j=m}^\infty \left( \theta^{k+1+\delta} \right)^{j-m} \left\| \pi_h(p_j - p_{j+1} )\right\|_{\underline{L}^2(B_{r_{j}})} .
\end{align*}
Let us estimate the last term. By~\eqref{e.tilde E_k 2} and~\eqref{e.crude omega} we have
\begin{multline*} 
\left\| \pi_h(p_j - p_{j+1} )\right\|_{\underline{L}^2(B_{r_{j}})}  \leq \left\| \pi_k(p_j - p_{j+1} )\right\|_{\underline{L}^2(B_{r_{j}})} \\ \leq
Cr_j^{k} \left(\tilde E_k(r_j) + \tilde E_k(r_{j+1})\right)
\leq C r_j^{-\delta} \omega(j) 
\leq  C r_j^{-\delta} \left( \theta^{k+1+\delta} \right)^{m-j} \omega(m)
\end{multline*}
and hence
\begin{equation*} 
\sum_{j=m}^\infty \left( \theta^{k+1+\delta} \right)^{j-m} \left\| \pi_k(p_j - p_{j+1} )\right\|_{\underline{L}^2(B_{r_{j}})}  
\leq C \omega(m)  \sum_{j=m}^\infty r_j^{-\delta} \leq C r_m^{-\delta} \omega(m).
\end{equation*}
Therefore we obtain 
\begin{equation} \label{e.omega(m) iter 0}
\omega_h(m+1) \leq \theta^{-h} \omega_h(m) + C r_m^{-\delta}  \omega(m).
\end{equation}
By choosing $H$ in the definition of $r_0$ larger, if necessary, we have that $Cr_0^{-\delta} \leq 1$, and therefore iterating the previous inequality, we get, for every $n \geq m$,
\begin{equation} \label{e.omega(m) bounded 1}
\omega(n) \leq \omega(m) \prod_{j=m}^n (\theta^{-k} + \theta^{\delta j}) \leq C \left(\frac{r_n}{r_m}\right)^k \omega(m).
\end{equation}
In particular,
\begin{equation} \label{e.omega(m) bounded 2}
\sup_{m\in \N} r_m^{-k} \omega(m) \leq C r_0^{-k} \omega(0) < \infty.
\end{equation}
Furthermore, defining
\begin{equation*} 
\tilde \omega_h(m) := r_{m}^{\delta  -k} \omega_h(m), 
\end{equation*}
we have by~\eqref{e.omega(m) iter 0} that, for $h \leq k-1$,
\begin{equation*} 
\tilde \omega_h(m+1) \leq \theta^{1-\delta} \tilde \omega_h(m)  + C r_m^{-k} \omega(m).
\end{equation*}
By taking supremum over $m$ and reabsorbing,
\begin{equation*} 
\sup_{m \in \N} r_m^{\delta  -k} \omega_h(m)  < \infty.
\end{equation*}
As a consequence, we obtain from the definition of $\omega_h(m)$ that
\begin{equation}  \label{e.pi_(k-1) p_j bounded}
\sup_{m \in \N} r_m^{\delta  -k}\left\|\pi_{k-1} p_j \right\|_{\underline{L}^2(B_{r_j})} < \infty.
\end{equation}
Returning to~\eqref{e.tilde E_k 2} and using~\eqref{e.omega(m) bounded 1}, we obtain, after summation,
\begin{equation} \label{e.tilde E_k 3}
\sum_{n=m}^\infty \tilde E_k(r_n) \leq C \sum_{n=m}^\infty  r_n^{-(k+\delta)} \omega(n) \leq C r_m^{-(k+\delta)}\omega(m).
\end{equation}
This implies that $\{\nabla^k \pi_k p_j\}_j$ is a Cauchy-sequence: for every $n \geq m$, 
\begin{align*} 
\left|\nabla^k \pi_k(p_n - p_m)\right| \leq C \sum_{j=m}^n  \tilde E_k(r_n) \leq C r_m^{-(k+\delta)}\omega(m) \to 0 \quad \mbox{as} \ m\to \infty
\end{align*}
Therefore there exists a polynomial~$q \in \Ahom_k$, which is homogeneous of degree~$k$ ($\pi_{k-1} q =0$)  and, for every~$m \in \N$, 
\begin{equation}
\label{e.tilde E_k 4}
\left|\nabla^k (\pi_k p_m - q) \right| \leq C r_m^{-(k+\delta)}\omega(m).
\end{equation}

\emph{Step 4.} We conclude the argument.  By $\mathrm{(ii)}_k$ we find $\phi_q \in \A_k$ such that, for every~$r \geq r_0$, 
\begin{equation*} 
\left\| \phi_q - q \right\|_{\underline{L}^2(B_{r})} \leq C r^{-\delta} \left\| q \right\|_{\underline{L}^2(B_{r})} .
\end{equation*}
We will then show that $u - \phi_q \in \A_{k-1}$. By the triangle inequality, we have 
\begin{equation*} 
r_j^{-k} \left\| u - \phi_q \right\|_{\underline{L}^2(B_{r_j})} \leq r_j^{-k}  \left\| u - q \right\|_{\underline{L}^2(B_{r_j})} + C r_j^{-(k+\delta)} \left\| q \right\|_{\underline{L}^2(B_{r_j})} 
\end{equation*}
and, using also~\eqref{e.tilde E_k 3} and~\eqref{e.tilde E_k 4}, 
\begin{align*} 
r_j^{-k}  \left\| u - q \right\|_{\underline{L}^2(B_{r_j})} & \leq \tilde E_{k}(r_j) + C \left| \nabla^k(p_j - q) \right|+ r_j^{-k}  \left\| \pi_{k-1} p_j \right\|_{\underline{L}^2(B_{r_j})} 
\\ & \leq   C r_{j}^{-\delta } \sup_{m \in \N} \left( r_m^{-k} \omega(m) + r_m^{\delta -k}  \left\| \pi_{k-1} p_j \right\|_{\underline{L}^2(B_{r_j})} \right).
\end{align*}
Then, by~\eqref{e.omega(m) bounded 2} and~\eqref{e.pi_(k-1) p_j bounded}, we obtain
\begin{equation*} 
\lim_{j \to \infty} r_j^{-k} \left\| u - \phi_q \right\|_{\underline{L}^2(B_{r_j})} = 0.
\end{equation*}
Therefore, by the definition, $u - \phi_q \in \A_{k-1}$. By $\mathrm{(i)}_{k-1}$ we consequently find $\tilde q \in \Ahom_{k-1}$ such that, for every $r\geq r_0$, 
\begin{equation*} 
\left\| u - \phi_q - \tilde q \right\|_{\underline{L}^2(B_{r})} \leq C r^{-\delta} \left\| \tilde q \right\|_{\underline{L}^2(B_{r})}. 
\end{equation*}
 Finally, setting $p = q+\tilde q \in \Ahom_k$, we have by the triangle inequality and~\eqref{e.equivalenceofnorms} that, for every~$r \geq r_0$, 
 \begin{equation*} 
\left\| u - p \right\|_{\underline{L}^2(B_{r})} \leq C r^{-\delta} \left(\left\| q \right\|_{\underline{L}^2(B_{r})} +   \left\| \tilde q \right\|_{\underline{L}^2(B_{r})} \right) \leq C r^{-\delta} \left\| p \right\|_{\underline{L}^2(B_{r})} ,
\end{equation*}
which was to be proven. 
\end{proof}

\begin{proof}[{Proof of~\eqref{e.findmeaphi}}]
Assuming that $\mathrm{(i)}_k$ and $\mathrm{(ii)}_k$  hold, we will show that also $\mathrm{(iii')}_k$ is true. Let $n \in \N$ be such that $R \in [r_n, r_{n+1})$. In view of~\eqref{e.basicreg}, by imposing the condition on~$H$ that 
\begin{equation} 
\label{e.H cond 3}
H \geq \left( 4 C_{s,k} \theta^{-1}  \right)^{\frac{1}{\delta}},
\end{equation}
which in view of the definition of~$r_0$ also implies that 
\begin{equation*}
C_{s,k} r_{j}^{-\delta} \leq \frac14 \theta, 
\end{equation*}
we deduce that, for every~$u \in \A(B_{r_{j+1}})$, there exists~$p \in \Ahom_{k}$ such that 
\begin{equation*}
\left\| u - p  \right\|_{\underline{L}^2(B_{r_j})} \leq \frac12 \theta^{k+1- \delta/2}  \left\| u  \right\|_{\underline{L}^2(B_{r_{j+1}})}.
\end{equation*}
Define sequences $\{ u_j\}$, $\{ p_j\}$ and $\{ \phi_j\}$ recursively by setting $u_n :=u$ and, for every $j \in \{1,\ldots,n\}$, selecting $p_{j-1}\in \Ahom_k$ by way of the previous display to satisfy
\begin{equation*} \label{e.pickpolypj}
 \left\| u_j - p_{j-1} \right\|_{\underline{L}^2(B_{r_{j-1}})} 
\leq \frac12\theta^{k+1 - \delta/2} \left\| u_j  \right\|_{\underline{L}^2(B_{r_j}) }  \quad \mbox{and} \quad \left\| p_{j-1} \right\|_{L^\infty(B_{r_{j}})} \leq C\left\| u_j  \right\|_{\underline{L}^2(B_{r_j}) }.
\end{equation*}
Then pick $\phi_{j-1} \in \A_k(\Rd)$ using the assumption $\mathrm{(ii)}_{k}$ satisfying
\begin{equation*} 
\left\| p_{j-1} - \phi_{j-1} \right\|_{\underline{L}^2(B_{r_{j-1}})} \leq C r_j^{-\delta} \left\| p_{j-1} \right\|_{\underline{L}^2(B_{r_{j-1}})}\,,
\end{equation*}
and set $u_{j-1}:= u_j - \phi_{j-1}$. The triangle inequality and the above estimates imply
\begin{equation*}
\left\| u_{j-1}  \right\|_{\underline{L}^2(B_{r_{j-1}})} \leq \left( \frac{1}{2} \theta^{k+1 - \delta/2} + C r_j^{-\delta} \right) \left\| u_j \right\|_{\underline{L}^2(B_{r_j})}.
\end{equation*}
As before, the condition on $H$ guarantees that $C r_j^{-\delta} \leq \frac12 \theta^{k+1 - \delta/2}$, and thus the previous inequality gives after iteration that, for every $j\in \{ 0,\ldots,n\}$,
\begin{equation*}
\left\| u_{j}  \right\|_{\underline{L}^2(B_{r_j})} \leq  \left( \frac{r_j}{r_n}\right)^{k+1-\delta/2} \left\| u \right\|_{\underline{L}^2(B_{r_n})} \leq C \left( \frac{r_j}{R}\right)^{k+1-\delta/2} \left\| u \right\|_{\underline{L}^2(B_{R})}.
\end{equation*}
Since
\begin{equation*} 
u_j = u - \sum_{m=j}^{n-1} \phi_m ,
\end{equation*}
we have hence shown that, for every $r\in\left[ r_0 ,R \right]$, 
\begin{equation} \label{e.iiiww}
\inf_{\phi \in \A_k(\Rd) } \left\| u - \phi \right\|_{\underline{L}^2(B_r)} \leq C\left( \frac rR \right)^{k+1 - \delta/2} \left\| u \right\|_{\underline{L}^2(B_{R})}.
\end{equation}
To complete the proof of $\mathrm{(iii')}_k$, we need to check that we can select $\phi\in \A_k(\Rd)$ independent of the radius $r$. Let $\phi_r\in \A_k(\Rd)$ achieve the infimum on the left side of~\eqref{e.iiiww} for $r \in [r_0,R]$. Then by the triangle inequality, for $r\in\left[ r_0 , \frac12 R \right]$,
\begin{equation*}
\left\| \phi_r - \phi_{2r} \right\|_{\underline{L}^2(B_r)} \leq C\left( \frac rR \right)^{k+1 - \delta/2} \left\| u \right\|_{\underline{L}^2(B_{R})}.
\end{equation*}
The statements $\mathrm{(i)}_k$ and $\mathrm{(ii)}_k$ imply that every $\phi \in \A_k(\R^d)$ satisfies, for every $s\geq r$, 
\begin{equation*} \label{}
\left\| \phi \right\|_{\underline{L}^2(B_s)} \leq C \left( \frac sr \right)^k \left\| \phi \right\|_{\underline{L}^2(B_r)}.
\end{equation*}
Combining the previous two displays, we get that, for every $s\geq r\geq r_0$,
\begin{equation*}
\left\| \phi_r - \phi_{2r} \right\|_{\underline{L}^2(B_s)} \leq C\left( \frac sr \right)^k \left( \frac rR \right)^{k+1-\delta/2} \left\| u \right\|_{\underline{L}^2(B_{R})} = C\left( \frac rs \right)^{1-\delta/2} \left( \frac sR \right)^{k+1-\delta/2} \left\| u \right\|_{\underline{L}^2(B_{R})}. 
\end{equation*}
Summing the previous inequality over dyadic radii yields, for every $s\geq r\geq r_0$,
\begin{equation*}
\left\| \phi_r - \phi_{s} \right\|_{\underline{L}^2(B_s)}  \leq C \left( \frac sR \right)^{k+\alpha} \left\| u \right\|_{\underline{L}^2(B_{R})}. 
\end{equation*}
In particular, if we take $\phi:= \phi_{r_0}$, then we obtain, for every $r \geq r_0$,
\begin{equation*}
\left\| u - \phi \right\|_{\underline{L}^2(B_r)} \leq \left\| u - \phi_r \right\|_{\underline{L}^2(B_r)} + \left\| \phi - \phi_r \right\|_{\underline{L}^2(B_r)}\leq C\left( \frac rR \right)^{k+1-\delta/2} \left\| u \right\|_{\underline{L}^2(B_{R})}.
\end{equation*}
This completes the proof of $\mathrm{(iii')}_k$.
\end{proof}

\begin{proof}[{Proof of~\eqref{e.climbontop}}]
We finish the proof of the theorem by showing that $\mathrm{(i)}_{k+1}$, $\mathrm{(ii)}_{k+1}$ and $\mathrm{(iii')}_{k+1}$ imply $\mathrm{(iii)}_k$. Fix $R\geq 2\X_s$ and $u\in \A(B_R)$. Select first 
$\psi \in \A_{k+1}(\Rd)$ such that, for every $r \in \left[ \X_s, R \right]$, 
\begin{equation}
\label{e.kk11}
\left\| u - \psi \right\|_{\underline{L}^2(B_r)} \leq C \left( \frac r R \right)^{k+2-\delta/2} 
\left\| u \right\|_{\underline{L}^2(B_R)}.
\end{equation}
By $\mathrm{(i)}_{k+1}$, we can choose $p_\psi \in \Ahom_{k+1}$ to satisfy
\begin{equation*} \label{}
\left\|\psi - p_\psi \right\|_{\underline{L}^2(B_r)} \leq Cr^{-\delta} \left\| p_\psi \right\|_{\underline{L}^2(B_r)} .
\end{equation*}
By~$\mathrm{(ii)}_{k+1}$, we may also take $\tilde \psi\in \A_{k+1}$ such that 
\begin{equation*} \label{}
\left\|\tilde \psi - (\pi_{k+1}-\pi_k) p_\psi \right\|_{\underline{L}^2(B_r)} \leq Cr^{-\delta} \left\| (\pi_{k+1}-\pi_k) p_\psi \right\|_{\underline{L}^2(B_r)}.
\end{equation*}
It is clear that $\tilde \psi$ has growth of degree $k+1$: by~\eqref{e.equivalenceofnorms}, we have
\begin{align} 
\notag 
\left\| \tilde \psi \right\|_{\underline{L}^2(B_r)} 
 \leq 2 \left\| (\pi_{k+1}-\pi_k) p_\psi \right\|_{\underline{L}^2(B_r)}  
& = 2 \left(\frac rR \right)^{k+1} \left\| (\pi_{k+1}-\pi_k) p_\psi \right\|_{\underline{L}^2(B_R)} 
\\  & \notag
\leq C \left(\frac rR \right)^{k+1} \left\| p_\psi \right\|_{\underline{L}^2(B_R)}
\\ &
\leq C \left(\frac rR \right)^{k+1} \left\| u \right\|_{\underline{L}^2(B_R)}.
\label{e.kk1111}
\end{align}
The function $\phi := \psi - \tilde \psi$ is now our candidate in~$\mathrm{(iii)}_k$. To see that $\phi \in \A_k$, we get by~\eqref{e.equivalenceofnorms}  that
\begin{align*} 
\notag 
 \left\| \phi \right\|_{\underline{L}^2(B_r)} 
 & \leq  \left\| \psi -  p_\psi \right\|_{\underline{L}^2(B_r)}   
 + \left\| \tilde \psi - (\pi_{k+1}-\pi_k) p_\psi \right\|_{\underline{L}^2(B_r)}  
 + \left\| \pi_k p_\psi \right\|_{\underline{L}^2(B_r)}   
 \\ & \leq 
 C r^{-\delta} \left\| p_\psi \right\|_{\underline{L}^2(B_r)} + \left\| \pi_k p_\psi \right\|_{\underline{L}^2(B_r)} 
 \\ & \leq C r^{k+1-\delta}  \left\| p_\psi \right\|_{\underline{L}^2(B_1)} ,
\end{align*}
and therefore $\lim_{r\to \infty} r^{-(k+1)}  \left\| \phi \right\|_{\underline{L}^2(B_r)}  = 0$ and $\phi \in \A_k$. 
Furthermore, by~\eqref{e.kk11},~\eqref{e.kk1111} and the triangle inequality, we find  that 
\begin{align*} 
\left\| u - \phi \right\|_{\underline{L}^2(B_r)} &  \leq C\left\| \tilde \psi \right\|_{\underline{L}^2(B_r)} +  C\left\| u - \psi \right\|_{\underline{L}^2(B_r)}  
 \\ &  \leq C \left(\left(\frac rs \right)^{k+1} +   \left( \frac r R \right)^{k+2 -  \delta/2} \right) \left\| u \right\|_{\underline{L}^2(B_R)} \leq C \left(\frac rs \right)^{k+1} \left\| u \right\|_{\underline{L}^2(B_R)}  .
\end{align*}
The proof is now complete. 
\end{proof}

\begin{proof}[{Proof of~\eqref{e.dimensionofAk}}]
An explicit expression for~$\dim \left( \Ahom_k \right)$ in terms of $(k,d)$ is 
\begin{equation*} \label{}
\dim \left( \Ahom_k \right)
=
\binom{d+k-1}{k} + \binom{d+k-2}{k-1}.
\end{equation*}
This formula holds also for $k=0$ if we interpret $\binom{d-2}{-1}=0$. 
A proof of this fact can be found in~\cite[Corollary 2.1.4]{Armitage}. 

\smallskip

We now argue by induction on~$k$ that 
\begin{equation} 
\label{e.samedimensions}
\dim(\A_k) = \dim(\Ahom_k).
\end{equation}
A proof of~\eqref{e.samedimensions} in the case $k=0$ has already been given above in the argument for~$\mathrm{(i)}_{0}$,~$\mathrm{(ii)}_{0}$ and~$\mathrm{(iii)}_{0}$, where we saw that $\A_0 = \Ahom_0$ is the set of constant functions. We next observe that ~$\mathrm{(i)}_{k}$ and~$\mathrm{(ii)}_{k}$ give us a canonical isomorphism between the spaces $\A_k / \A_{k-1}$ and $\Ahom_k / \Ahom_{k-1}$. Thus, for every $k\in\N$ with $k\geq 1$, 
\begin{equation*} \label{}
\dim(\A_k/\A_{k-1}) = \dim(\Ahom_k/\Ahom_{k-1}).
\end{equation*}
We now easily obtain~\eqref{e.samedimensions} by induction on~$k$. 
\end{proof}

Theorem~\ref{t.regularity} may be generalized to allow for polynomial right-hand sides, essentially replacing~$\Ahom_k$ by~$\Phom_k$ in each of the assertions, where $\Phom_k$ is the set of polynomials of degree at most~$k$. To give the precise statement, we consider, for each $k\in\N$ and $\Phom_k$, the equation
\begin{equation*} \label{}
-\nabla \cdot \left( \a(x) \nabla u \right) = -\nabla \cdot \left( \ahom \nabla p \right) \quad \mbox{in} \ U
\end{equation*}
and denote the set of solutions by 
\begin{equation*} \label{}
\A\left[ p \right] (U):= \left\{ u \in H^1_{\mathrm{loc}}(U)\, :\, \forall w \in C^\infty_c(U), \ \int_{U} \nabla w \cdot \left( \a \nabla u -\ahom \nabla p\right) = 0\right\}.
\end{equation*}
We then define, for each $k\in\N$ and $p\in \Phom_k$, the affine space 
\begin{equation*} \label{}
\A_k\left[ p \right] := \left\{ u \in \A\left[ p \right](\Rd)  \,:\, \limsup_{r\to \infty} r^{-(k+1)} \left\| u \right\|_{\underline{L}^2(B_r)} = 0 \right\},
\end{equation*}
and set
\begin{equation*} \label{}
\Po_k:= \bigcup_{p\in\Phom_k} \A_k\left[ p \right]. 
\end{equation*}
In order to extend Theorem~\ref{t.regularity} to allow for polynomial right-hand sides, it is necessary only to extend the statement $\mathrm{(ii)}_k$. This is done in the following lemma. 

\begin{lemma}
\label{l.morepees}
For each $s\in (0,d)$, there exists $\delta(s,d,\Lambda)>0$ and a random variable $\X_s$ satisfying~\eqref{e.X} such that, for every $k\in \N$ and $p\in \mathcal{P}_k$, there exists $u \in \A_k\left[ p \right]$ such that, for every $R\geq \X_s$, 
\begin{equation*} \label{}
\left\| u - p \right\|_{\underline{L}^2(B_R)} \leq CR^{-\delta} \left\| p \right\|_{\underline{L}^2(B_R)}. 
\end{equation*}
\end{lemma}
\begin{proof}
The proof is almost identical to that of~\eqref{e.findmeap}, except that it requires a version of Theorem~\ref{t.DP} with (smooth) right-hand sides. The required modifications are straightforward and left as an exercise to the reader.  
\end{proof}

\begin{corollary}[Higher regularity theory, generalized]
\label{cor.regularity}
\index{Liouville theorem}
Fix $s \in (0,d)$. There exist an exponent $\delta(s,d,\Lambda)\in \left( 0, \frac12 \right]$ and a random variable $\X_s$ satisfying the estimate
\begin{equation}
\label{e.X.p}
\X_s \leq \O_s\left(C(s,d,\Lambda)\right)
\end{equation}
such that the following statements hold, for every $k\in\N$:
\begin{enumerate}
\item[{$\mathrm{(i)}_k$}] There exists $C(k,d,\Lambda)<\infty$ such that, for every $u \in \Po_k$, there exists $p\in \Phom_k$ such that $u \in \A_k\left[p\right]$ and, for every $R\geq \X_s$,
\begin{equation} \label{e.liouvillec.p}
\left\| u - p \right\|_{\underline{L}^2(B_R)} \leq C R^{-\delta} \left\| p \right\|_{\underline{L}^2(B_R)}.
\end{equation}

\item[{$\mathrm{(ii)}_k$}]For every $p\in \Phom_k$, there exists $u\in \A_k\left[ p \right]$ satisfying~\eqref{e.liouvillec} for every $R\geq \X_s$. 

\item[{$\mathrm{(iii)}_k$}]
There exists $C(k,d,\Lambda)<\infty$ such that, for every $R\geq \X_s$, $p\in \Phom_k$ and solution $u\in H^1(B_R)$ of the equation
\begin{equation*} \label{}
-\nabla \cdot \left( \a\nabla u \right) = -\nabla \cdot \left( \ahom \nabla p \right) =: f \quad \mbox{in} \ B_R,
\end{equation*}
there exists $\phi \in \A_k\left[ p \right]$ such that, for every $r \in \left[ \X_s,  R \right]$,
\begin{multline}
\label{e.intrinsicreg.p}
\left\| u - \phi \right\|_{\underline{L}^2(B_r)}
\leq C \left( \frac r R \right)^{k+1} 
\left\| u \right\|_{\underline{L}^2(B_R)}
\\
+ C\left( \frac r R \right)^{k+1}\left( \sum_{j=0}^{k-1} R^j \left\| \nabla^j f \right\|_{L^\infty(B_R)} +R^{k-1+\beta} \left[ \nabla^{k-1} f \right]_{C^{0,\beta}(B_R) } \right).
\end{multline}
\end{enumerate}
In particular, $\P$-almost surely, we have, for every $k\in\N$ and $p\in \Phom_k$,
\begin{equation} 
\label{e.dimensionofAk.p}
\dim\left(\A_k\left[ p \right]\right) 
= 
\dim(\Ahom_k) 
= \binom{d+k-1}{k} + \binom{d+k-2}{k-1}.
\end{equation}
\end{corollary}
\begin{proof}
The corollary is immediate from Theorem~\ref{t.regularity} and Lemma~\ref{l.morepees}. Indeed, we may use the lemma and linearity to reduce the statement of the corollary to Theorem~\ref{t.regularity}. 
\end{proof}

\begin{exercise}
\label{ex.fullslate}
Show that for every $k\in\N$ and $p\in \mathcal{P}_k$, there exists  $q \in \mathcal{P}_{k+2}$ such that $p = -\nabla \cdot \left( \ahom \nabla q\right)$. Hint: consider the solution $v\in C^\infty(B_1)$ of the Dirichlet problem
\begin{equation*} \label{}
\left\{ 
\begin{aligned}
& -\nabla \cdot \left( \ahom \nabla v \right) = p & \mbox{in} & \ B_1, \\
& v = 0 & \mbox{on} & \ \partial B_1,
\end{aligned}
\right.
\end{equation*}
and define 
\begin{equation*} \label{}
q(x):= \sum_{j=0}^{k-2} \nabla^j v(0) x^{\otimes j}. 
\end{equation*}
Deduce~$p = -\nabla \cdot \left( \ahom \nabla q\right)$ by a scaling argument and the~$C^{k-1,\alpha}$ regularity of~$v$.
\end{exercise}

\begin{exercise}
\label{ex.regularity}
Further generalize Corollary~\ref{cor.regularity}$\mathrm{(iii)}_k$ to allow for non-polynomial right-hand sides. The precise statement is:
for every $s\in (0,d)$ and $\beta\in (0,1)$, there exist~$C(s,\beta,d,\Lambda)<\infty$ and a random variable $\X_s :\Omega \to [1,\infty]$ satisfying
\begin{equation}
\label{e.XLipschitz040}
\X_s = \O_s\left( C \right),
\end{equation}
such that the following holds: for every $R\geq \X_s$, $k\in\N$ with $k\geq 1$, $f \in C^{k-1,\alpha}(B_{2R})$, and weak solution $u\in H^1(B_{2R})$ of
\begin{equation} 
\label{e.wksrt2040}
-\nabla \cdot \left( \a\nabla u \right) = f \quad \mbox{in} \ B_{2R},
\end{equation}
if we use Exercise~\ref{ex.fullslate} to find $p\in \Phom_k$ such that  $-\nabla \cdot \left( \ahom \nabla p \right) (x)= \sum_{j=0}^{k-2} \nabla^j f(x) x^{\otimes j}$, then there exists $\phi\in\A_k[p]$ satisfying, for every $r\in \left[\X_s,R \right]$, 
\begin{multline}
\label{e.intrinsicreg2040}
\left\| u - \phi \right\|_{\underline{L}^2(B_r)} 
\leq C \left( \frac r R \right)^{k+1} 
\left\| u \right\|_{\underline{L}^2(B_R)} 
\\
+ C\left( \frac r R \right)^{k+1}\left( \sum_{j=0}^{k-1} R^j \left\| \nabla^j f \right\|_{L^\infty(B_R)} +R^{k-1+\beta} \left[ \nabla^{k-1} f \right]_{C^{0,\beta}(B_R) } \right).
\end{multline}
Note that last factor of the last term in parentheses on the right side is just the properly scaled $C^{k-1,\beta}(B_R)$ norm of $f$. 
\end{exercise}

\begin{remark}[{Stationary extension of $\X_s$}]
\label{r.Lipminscale}
The random variable~$\X_s$ in the statements of Theorems~\ref{t.Lipschitz},~\ref{t.regularity} and Corollary~\ref{cor.regularity} controls the regularity in balls centered at the origin. We obtain similar statements for balls centered at $x\in\Rd$ by replacing with $\X_s$ by $T_{x}\X_s$, where we recall $T_x\X_s$ is the translation of $\X_s$ given by
\begin{equation*} \label{}
\left( T_x \X_s\right)(\a):= \X_s \left( T_x\a\right),
\end{equation*}
and $T_x$ acts on the coefficient field by $(T_x\a)(y):= \a(y+x)$. It is convenient therefore to think of $\X_s$ as a \emph{stationary random field} \index{stationary random field}
(see Definition~\ref{def.stat.field} below) and to denote its spatial dependence by writing
\begin{equation*} \label{}
\X_s(x):= T_x\X_s. 
\end{equation*}
It will be helpful to have some regularity in $x$ of $\X_s(x)$. We will now argue that we may modify $\X_s$ slightly, without changing any of the conclusions in the above theorems, so that, for every $x,y\in\Rd$,
\begin{equation} 
\label{e.Xs.lip}
\left| \X_s(x) - \X_s(y) \right| \leq 2 |x-y|. 
\end{equation}
To see this, we  replace~$\X_s(x)$ by $\X'_s(x)$, defined by
\begin{equation*} \label{}
\X'_s(x):= 2 + \inf_{y\in\Rd} \left( 2|y-x| + \X_s(y) \right). 
\end{equation*}
Observe that, for every $x\in\Rd$,  $\X'_s(x) \leq 2 + \sqrt{d}+ \X_s\left(\left\lceil x\right\rceil\right)$, where $\left\lceil x\right\rceil\in\Zd$ is the lattice point with coordinates $\left(\left\lceil x_i \right\rceil\right)$. Thus $\X'_s(x)$ satisfies the estimate~\eqref{e.X.p} after enlarging the constant~$C$. It is clear that $x\mapsto \X'_s(x)$ is $2$-Lipschitz, that is, 
\begin{equation*} \label{}
\left| \X'_s(x) - \X'_s(y) \right| \leq 2 |x-y|. 
\end{equation*}
Finally, we obtain all the same statements $\mathrm{(i)}_k$--$\mathrm{(iii)}_k$ in Theorem~\ref{t.regularity} and Corollary~\ref{cor.regularity} (and likewise in Theorem~\ref{t.Lipschitz}) after inflating the constant $C$. To see this, fix $x$ and find $y$ such that, with $R:=|x-y| \vee 1$, we have 
\begin{equation*} \label{}
\X_s(y) \leq \X'_s(x) - 2 R. 
\end{equation*}
Then $r \geq \X'_s(x)$ implies $r\geq 2R$ and hence we have both $B_r(x) \subseteq B_{2r}(y)$ and $B_r(y) \subseteq B_{2r}(x)$  for every $r\geq \X_s'(y)$. Thus for every $f\in L^2(B_R(x))$ and $r\in \left[ \X'_s(x), \tfrac 14R \right]$, 
\begin{equation*} \label{}
\left\| f \right\|_{\underline{L}^2(B_r(x))}
\leq 
2^d \left\| f \right\|_{\underline{L}^2(B_{2r}(y))}
\quad \mbox{and} \quad
\left\| f \right\|_{\underline{L}^2(B_r(y))}
\leq 
2^d \left\| f \right\|_{\underline{L}^2(B_{2r}(x))}
\end{equation*}
This observation and the triangle inequality gives us the statements $\mathrm{(i)}_k$--$\mathrm{(iii)}_k$ in Theorem~\ref{t.regularity} and Corollary~\ref{cor.regularity} with $\X'_s(x)$ in place of $\X_s$ and with the balls centered at $x$ rather than the origin. 
\end{remark}

\section{The first-order correctors}
\label{s.correctorsdefined}

In this section, we identify the vector space $\A_1/\R$ with a family of $\Zd$--stationary potential field which we call the \emph{first-order correctors}. We then obtain some preliminary estimates on these functions. 

\smallskip

The name \emph{corrector} comes from the fact that $\phi_\xi$ represents the difference between the element of $\A_1$ tracking the affine function $x\to \xi\cdot x$ and the affine function itself. We should therefore think of~$\A_1$ as an analogue in the whole space of the solutions to the Dirichlet problem in a large domain with affine boundary data. In fact, these two sets of solutions will be converging to each other in the sense that,~for some $\delta>0$ and every $r$ larger than a random minimal scale, 
\begin{equation*} \label{}
\left\| \phi_\xi - \left( \phi_\xi \right)_{B_r} \right\|_{\underline{L}^2(B_r)} 
\leq 
C r^{1-\delta}.
\end{equation*}
We can therefore use the estimates we have already proved in the previous chapters to get some quantitative bounds on the correctors in addition to the bound on the~$L^2$-oscillation in~$B_r$. There are two advantages to studying $\phi_\xi$ rather than the finite-volume analogues in previous chapters: one is that its gradient $\nabla \phi_\xi$ is \emph{statistically stationary} (a concept explained below), and the second is that there are no boundaries to perturb its behavior. There are also disadvantages, such as the fact that the correctors are not localized: they depend on the coefficients in the whole space.

\smallskip

We begin by introducing the concept of a stationary random field. 
\index{stationary random field|(}

\smallskip

\begin{definition}[Stationary random field]
\label{def.stat.field}
Let $(S,\S)$ be a measurable space. An \emph{$S$--valued $\Zd$--stationary random field} is a $\mathcal{B}\times \F$--measurable mapping
\begin{equation*} \label{}
F:\Rd \times \Omega \to S
\end{equation*}
with the property that, for every $z\in\Zd$ and $\a\in\Omega$,
\begin{equation*} \label{}
F(x+z,\a) = F(x,T_z\a) \quad \mbox{for a.e.} \ x\in\Rd \ \ \mbox{and} \ \ \mbox{$\P$--a.s. in} \ \a\in\Omega.
\end{equation*}
If $S = \R$ then we typically call $F$ a \emph{stationary random field} or just \emph{stationary field}. If $S=\Rd$ then we may say that $F$ is a \emph{stationary vector field}.
\end{definition}

If $F$ is an $S$--valued stationary random field, we typically keep the dependence of $F$ on $\a$ implicit in our notation and just write $F(x)$. Of course, the canonical process~$\a \mapsto \a$ is itself the most obvious example of a stationary random field, as is any function of~$\a$ which commutes with translations.  

\begin{definition}[Stationary random potential field]
\label{d.stat.l2pot}
\emph{A stationary random potential field} is an $\Rd$--valued stationary random field $F$ such that,~$\P$--a.s.~in~$\a\in\Omega$, the vector field $F(\cdot,\a)$ is the gradient of a function in $H^1_{\mathrm{loc}}(\Rd)$.
In this case we can write $F(x,\a) = \nabla u(x,\a)$. We denote the collection of stationary random potential fields by $\mathbb{L}^2_{\mathrm{pot}}$. 
\end{definition}

We can also think of an element of $\mathbb{L}^2_{\mathrm{pot}}$ as a mapping from $\Omega$ to $L^2_{\mathrm{pot},\,\mathrm{loc}}(\Rd)$. Note that an element $\nabla u$ of $\mathbb{L}^2_{\mathrm{pot}}$ need not be the spatial gradient of a stationary random field, since the potential $u$ is only well-defined up to a constant. 
\index{stationary random field|)}

\begin{lemma}
\label{l.identifyA1}
There exists a family $\left\{ \nabla \phi_e \,:\, e\in\Rd\right\} \subseteq \mathbb{L}^2_{\mathrm{pot}}$ such that 
\begin{equation*} \label{}
\nabla \A_1 : = \left\{ \nabla u \,:\, u \in \A_1 \right\} = \left\{ e + \nabla \phi_e \,:\, e\in\Rd \right\}
\end{equation*}
and
\begin{equation*} \label{}
\limsup_{r\to\infty} \frac1r \left\| \phi_e - \left( \phi_e \right)_{B_r} \right\|_{\underline{L}^2(B_r)} = 0 \quad \mbox{$\P$--a.s.}
\end{equation*}
\end{lemma}
\begin{proof}
According to Theorem~\ref{t.regularity}, to every $u \in \A_1$, there corresponds a unique---up to additive constants---element of $\Ahom_1$. Conversely, every $p\in \Ahom_1$ is associated to an element of $\A_1$ which is again unique up to additive constants. Therefore there exists a canonical bijection between $\nabla \A_1$ and $\nabla \Ahom_1$. Identifying the latter with~$\Rd$, we obtain an isomorphism $\xi \mapsto (\xi+\nabla \phi_\xi)$ from $\Rd$ to $\nabla \A_1$, where the potential field $\xi +\nabla \phi_\xi$ is uniquely determined, $\P$--a.s., by membership in $\nabla\A_1$ and the condition
\begin{equation*} \label{}
\limsup_{R \to \infty} \frac1R \left\| \phi_\xi \right\|_{\underline{L}^2(B_{R})} = 0.
\end{equation*}
We next argue that the potential field $\nabla \phi_\xi$ is $\Z^d$-stationary. Indeed, for each fixed $z \in \Z^d$, consider the potential field $\xi + \nabla \phi'_\xi$ obtained through the canonical bijection associated with the shifted coefficient field $T_z\a:= \a(\cdot+z)$. We have that the function~$x \mapsto \xi \cdot x+ \phi'_\xi(x-z)$ belongs to~$\A_1$ and, since ($\P$--a.s.)
\begin{equation*}  
\limsup_{R \to \infty} \frac1R \|\phi'_\xi\|_{\underline L^2(B_R)} = 0,
\end{equation*}
the triangle inequality yields ($\P$--a.s.)
\begin{equation*}  
\limsup_{R \to \infty} \frac1R \|\phi'_\xi(\cdot \, - z)\|_{\underline L^2(B_R)} = 0.
\end{equation*}
We deduce that $\nabla \phi_\xi = \nabla \phi'_\xi(\cdot \, - z)$, $\P$--a.s. In other words,
$\nabla \phi_\xi(x+z,\a) = \nabla \phi_\xi(\cdot,T_z\a)$, $\P$--a.s. in $\a\in\Omega$, therefore $\nabla \phi_\xi \in \mathbb{L}^2_{\mathrm{pot}}$. 
\end{proof}

The potential functions~$\phi_e$ for the gradient fields in the family $\{ \nabla \phi_e \}_{e\in\Rd}$ are called the \emph{first-order correctors}. \index{corrector!first-order~$\phi_e$|(}
We next record an important consequence of Theorem~\ref{t.regularity} for the first-order correctors, giving us bounds on~$\| \nabla \phi_e \|_{\underline L^2(B_r)}$ for~$r$ larger than a minimal scale. 

\begin{lemma}
\label{l.correctorminbounds}
Fix $s\in (0,d)$ and let $\X_s$ be as in the statement of Theorem~\ref{t.regularity}. There exist~$\beta(d,\Lambda)>0$ and $C(d,\Lambda)<\infty$ such that, for every $r\geq \X_s$ and $e\in \partial B_1$, 
\begin{multline}
\label{e.correctorminbounds0}
\left\| \phi_{e} - \left( \phi_e \right)_{B_r} \right\|_{\underline{L}^2(B_r)} 
+
\left\| \nabla \phi_e  \right\|_{\hat{\underline{H}}^{-1}(B_{r})}
+ \left\| \a \left( e + \nabla \phi_e \right)  - \ahom e \right\|_{\hat{\underline{H}}^{-1}(B_{r})} 
\\
+ \left\| \tfrac12 \left( e + \nabla \phi_e \right)  \cdot  \a\left( e + \nabla \phi_e \right)  - \tfrac12 e \cdot \ahom e \right\|_{\underline{W}^{-1,1}(B_{r})}
\leq Cr^{1-\beta(d-s)}
\end{multline}
and
\begin{equation} 
\label{e.correctorminbounds}
\left\| \nabla \phi_{e} \right\|_{\underline{L}^2(B_r)} 
\leq C. 
\end{equation}
In particular, there exists $\ep(d,\Lambda)>0$ and $C(d,\Lambda)<\infty$ such that
\begin{equation} 
\label{e.correctorgradbound}
\sup_{e\in\partial B_1} \left\| \nabla \phi_{e} \right\|_{\underline{L}^2(B_1)} 
\leq 
\O_{2+\ep}\left( C \right).
\end{equation}
\end{lemma}
\begin{proof}
According to statement $\mathrm{(i)}_1$ of Theorem~\ref{t.regularity}, for every $r\geq \X_s$, we have
\begin{equation} 
\label{e.glafhbrah}
 \left\| \phi_e - \left( \phi_e \right)_{B_r} \right\|_{\underline{L}^2(B_r)} \leq Cr^{1-\delta}.
\end{equation}
We next prove the estimate for the other three terms on the right side of~\eqref{e.correctorminbounds0}. 
Let~$\ell_e$ denote the affine function $\ell_e(x):= e\cdot x$. For each $r>0$, let $w_r \in H^1(B_r)$ denote the solution of the Dirichlet problem
\begin{equation*} \label{}
\left\{ 
\begin{aligned}
& -\nabla \cdot \left( \a\nabla w_r \right) = 0 & \mbox{in} & \ B_r, \\
& w_r = \ell_e & \mbox{on} & \ \partial B_r. 
\end{aligned}
\right.
\end{equation*}
According to Theorem~\ref{t.DP} and the proof of Proposition~\ref{p.harmonicapproximation}, there exists an exponent~$\delta(s,d,\Lambda)>0$ and a random variable $\X_s$ such that $\X_s \leq \O_s(C)$ and, for every $r\geq \X_s$, 
\begin{multline} 
\label{e.BrDPbounds}
 \left\| w_r - \ell_e \right\|_{\underline{L}^2(B_r)}^2 
+ \left\| \nabla w_r - e \right\|_{\hat{\underline{H}}^{-1}(B_r)}^2 
+ \left\| \a \nabla w_r - \ahom e \right\|_{\hat{\underline{H}}^{-1}(B_r)}^2 
\\
+ \left\| \tfrac12 \nabla w_r \cdot  \a\nabla w_r - \tfrac12 e \cdot \ahom e \right\|_{\underline{W}^{-2,1}(B_r)}
\leq Cr^{2(1-\delta)}. 
\end{multline}
Combining~\eqref{e.glafhbrah}, the previous line, the triangle inequality, the fact that $w_r - ( \ell + \phi_e)\in \A(B_r)$ and the Caccioppoli inequality, we find that, for every $r\geq \X_s$,
\begin{equation} 
\label{e.snapcorrector}
\left\| \nabla w_r - \left( e + \nabla \phi_e \right) \right\|_{\underline{L}^2(B_{r/2})} \leq Cr^{-\delta}.
\end{equation}
The triangle inequality,~\eqref{e.BrDPbounds} and~\eqref{e.snapcorrector} imply that, for every $r\geq \X_s$,
\begin{multline} 
\label{e.BrDPboundscorr}
 \left\| \nabla \phi_e  \right\|_{\hat{\underline{H}}^{-1}(B_{r/2})}^2 
+ \left\| \a \left( e + \nabla \phi_e \right)  - \ahom e \right\|_{\hat{\underline{H}}^{-1}(B_{r/2})}^2 
\\
+ \left\| \tfrac12 \left( e + \nabla \phi_e \right)  \cdot  \a\left( e + \nabla \phi_e \right)  - \tfrac12 e \cdot \ahom e \right\|_{\underline{W}^{-2,1}(B_{r/2})}
\leq Cr^{2(1-\delta)}. 
\end{multline}
Using~\eqref{e.correctorminbounds}, which we will prove next, and interpolation we can improve the norm for the last term from $W^{-2,1}$ to ${W}^{-1,1}$, after shrinking the exponent~$\delta$. 
This completes the proof of~\eqref{e.correctorminbounds0}. 

\smallskip

We turn to the proof of~\eqref{e.correctorminbounds}.  
Set $\overline{\phi}_e:= \ell_e + \phi_e$ and observe that, for every $r\geq \X_s$,
\begin{equation*} \label{}
 \left\| \overline{\phi}_e - \left( \phi_e \right)_{B_r} \right\|_{\underline{L}^2(B_r)}
\leq \left\| \ell_e \right\|_{\underline{L}^2(B_r)} 
+  \left\|  \phi_e - \left( \phi_e \right)_{B_r} \right\|_{\underline{L}^2(B_r)}
\leq Cr.
\end{equation*}
As~$\overline{\phi}_e\in\A_1 \subseteq \A$, we apply the Caccioppoli inequality to obtain, for every~$r\geq \X_s$, 
\begin{equation*} \label{}
\left\| \nabla \overline{\phi}_e\right\|_{\underline{L}^2(B_r)} 
\leq 
\frac{C}{r}  \left\| \overline{\phi}_e -  \left( \phi_e \right)_{B_r} \right\|_{\underline{L}^2(B_{2r})}
\leq C. 
\end{equation*}
As $\left\| \nabla {\phi}_e\right\|_{\underline{L}^2(B_r)}  \leq 1 + \left\| \nabla \overline{\phi}_e\right\|_{\underline{L}^2(B_r)}$, we obtain~\eqref{e.correctorminbounds}. 

\smallskip

To prove~\eqref{e.correctorgradbound}, we first notice that we can just give up a volume factor to get 
\begin{equation*}
\sup_{e\in\partial B_1} \left\| \nabla \phi_e \right\|_{\underline{L}^2(B_1)} \leq \X_s^{\frac d2} \sup_{e\in\partial B_1} \left\| \nabla \phi_e \right\|_{\underline{L}^2(B_{\X_s} )} \leq C\X_s^{\frac d2} \leq \O_{2s/d}(C). 
\end{equation*}
The right side is slightly worse than $\O_2(C)$, but we can do slightly better than give up a volume factor by using H\"older's inequality and the Meyers estimate:
\begin{align*}
 \left\| \nabla \overline{\phi}_e \right\|_{\underline{L}^2(B_1)} 
&
\leq \left( \X_s^{\frac d2} \right)^{\frac{2}{2+\delta}} \left\| \nabla \overline{\phi}_e \right\|_{\underline{L}^{2+\delta}(B_{\X_s/2})} \leq C \left( \X_s^{\frac d2} \right)^{\frac{2}{2+\delta}} \left\| \nabla \overline{\phi}_e \right\|_{\underline{L}^{2}(B_{\X_s})} \leq C \X_s^{\frac d{2+\delta}}.
\end{align*}
Now the right side is~$\leq \O_{(2+\delta)s/d}(C)$ and so we can take~$s$ close enough to~$d$ to obtain the result. 
\end{proof}

\begin{remark}[Pointwise bounds for scalar equations]
\label{r.scalar.Linfty}
Using the De Giorgi-Nash H\"older estimate, we can upgrade~\eqref{e.correctorminbounds0}  from~$L^2$ to~$L^\infty$ on the small scales. The claim is that, in the notation of Lemma~\ref{l.correctorminbounds} and after shrinking the exponent $\beta(s,d,\Lambda)>0$, we have, for every~$e\in\partial B_1$ and $r\geq \X_s$, the estimate
\begin{equation} 
\label{e.correctorminbounds.Linfty}
\sup_{e\in\partial B_1} \left\| \phi_{e} - \left( \phi_e \right)_{B_r} \right\|_{L^\infty(B_r)} 
\leq Cr^{1-\beta(d-s)}.
\end{equation}
This estimate is only valid for scalar equations---since the De Giorgi-Nash estimate is unavailable for elliptic systems. We will not use it to prove any of the results in this chapter (or, except for a few exceptions, in the rest of the book).

\smallskip

By the interpolation inequality~\eqref{e.interp.LinftyL2Calpha} stated below and the sublinear estimate~\eqref{e.correctorminbounds0}, the bound~\eqref{e.correctorminbounds.Linfty} reduces to the following claim: there exist~$\alpha(d,\Lambda)>0$ and~$C(s,d,\Lambda)<\infty$ such that, for every $r\geq 2 \X_s$, 
\begin{equation} 
\label{e.correctorminCalpha}
r^{\alpha} \left[ \phi_e \right]_{C^{0,\alpha}(B_{r/2})} \leq Cr. 
\end{equation}
To obtain the latter, we apply the De Giorgi-Nash estimate~\eqref{e.DeGiorgi} in the ball $B_r$ to the function $w(x):= x\mapsto e\cdot x+\phi_{e}(x)$, which satisfies
\begin{equation*} \label{}
-\nabla \cdot \left( \a\nabla w \right) = 0 \quad \mbox{in} \ \Rd, 
\end{equation*}
to obtain, by the triangle inequality and~\eqref{e.correctorminbounds0},
\begin{align*} \label{}
r^{\alpha}  \left[ \phi_e \right]_{C^{0,\alpha}(B_{r/2})} 
&
\leq r^{\alpha} + r^{\alpha} \left[ w \right]_{C^{0,\alpha}(B_{r/2})}
\\ & 
\leq r^{\alpha} + C\left\|  w - (w)_{B_r} \right\|_{\underline{L}^2(B_{r})}
\\ & 
\leq r^{\alpha} + C \left( r + \left\|  \phi_e - (\phi_e)_{B_r} \right\|_{\underline{L}^2(B_{r})} \right) 
\\ & 
\leq Cr. 
\end{align*}
This gives~\eqref{e.correctorminCalpha} and completes the proof of~\eqref{e.correctorminbounds.Linfty}. 
\end{remark}

\begin{exercise}
\label{ex.interp}
Prove the following interpolation inequality for~$L^\infty$ between~$L^2$ and~$C^{0,\alpha}$: 
for every $r\in (0,\infty)$ and~$\alpha \in (0,1]$, there exists $C(\alpha,d)<\infty$ such that, for every function~$u \in C^{0,\alpha}(B_r)$, 
\begin{equation} 
\label{e.interp.LinftyL2Calpha}
\left\| u \right\|_{L^\infty(B_r)} 
\leq 
C \left\| u \right\|_{\underline{L}^2(B_r)}^{\frac{2\alpha}{d+2\alpha}} 
\left( r^\alpha \left[ u \right]_{C^{0,\alpha}(B_r)} \right) ^{\frac{d}{d+2\alpha}}.
\end{equation}
Sketch of proof: By scaling we may assume $r=1$ and $\| u \|_{L^\infty(B_1)} = 1$. Fix $x_0\in B_1$ such that $u(x_0) > \frac 12$. Deduce that 
\begin{equation*} \label{}
u(x) > \frac12 - \left[ u \right]_{C^{0,\alpha}(B_1)} |x-x_0|^\alpha =:w(x)
\end{equation*}
Notice that $w(x) \geq \frac 14$ in a ball $B_r(x_0)$ with $r = c \left[ u \right]_{C^{0,\alpha}(B_1)}^{-\frac1\alpha}$ and conclude that
\begin{equation*} \label{}
\left\| u \right\|_{L^2(B_1)} \geq \left\| w \right\|_{L^2(B_r(x_0))} \geq c |B_r|^{\frac12} = c \left[ u \right]_{C^{0,\alpha}(B_1)}^{-\frac d{2\alpha}}.
\end{equation*}
\end{exercise}

The main focus of the next chapter is on obtaining sharp estimates of the gradient, flux and energy density of the first-order correctors in negative Sobolev spaces. We prepare for the analysis there by giving a suboptimal version of these bounds, with a small exponent $\beta>0$, but which is optimal in stochastic integrability. 

\smallskip

For reasons that will become clear below, it is very convenient to formulate these estimates in terms of spatial averages of these functions against the standard heat kernel. Here and throughout the rest of this chapter and the next one, we use the notation, for a given $\Psi\in L^1(\Rd)$,
\begin{equation} 
\int_\Psi f  = \int_{\Psi} f(x) \, dx := \int_{\Rd} \Psi(x) f(x)\,dx. 
\end{equation}
Recall that we denote the standard heat kernel by 
\begin{equation} 
\label{e.HK.1}
\Phi(t,x):= \left( 4\pi t\right)^{-\frac d2} \exp\left( -\frac{|x|^2}{4t} \right)
\end{equation}
and define, for each $z\in\Rd$ and $r>0$, 
\begin{equation} 
\label{e.HK.2}
\Phi_{z,r} (x):= \Phi(r^2,x-z) \quad \mbox{and} \quad \Phi_r := \Phi_{0,r}. 
\end{equation}

\begin{proposition}
\label{p.correctorsbasecase}
For each $s\in(0,d)$, there exist $C(s,d,\Lambda)<\infty$, $\beta(s,d,\lambda)>0$ and a random variable $\X_s$ satisfying 
\begin{equation} \label{e.Ysbasecase}
\X_s(x) \leq \O_s(C)
\end{equation}
such that, for every $r\geq \X_s$ and $e\in\partial B_1$, we have 
\begin{equation}
\label{e.gradientbasecase}
\left| \int_{\Phi_{r}} \nabla \phi_e \right| \leq Cr^{-\beta},
\end{equation}
\begin{equation}
\label{e.fluxbasecase}
\left| \int_{\Phi_{r}} \a (e +  \nabla \phi_e) - \ahom e \right| \leq Cr^{-\beta},
\end{equation}
and
\begin{equation}
\label{e.energydensitybasecase}
\left| \int_{\Phi_{r}} \frac12 \left( e+ \nabla \phi_e\right) \cdot \a\left( e+ \nabla \phi_e\right)  - \frac12 e\cdot \ahom e \right| \\
 \leq Cr^{-\beta}.
\end{equation}
\end{proposition}
\begin{proof}
We apply Lemma~\ref{l.layercake} below with $k \in \{1,2\}$, $p=2$, $\beta = 2d$, and $\psi(x) = \exp\left(-\tfrac{|x|^2}{4}\right)$. Clearly~\eqref{e.layercakeweight} is valid with $A(d) < \infty$, and thus by Lemmas~\ref{l.correctorminbounds} and~\ref{l.layercake} we get, for every $r\geq \X_s$, with $\X_s$ as in Lemma~\ref{l.correctorminbounds},
\begin{equation*} 
\left| \int_{\Phi_{r}} \nabla \phi_e \right| \leq  C  
r^{-1} \int_{1}^\infty  t^{-d -1} \left\| \nabla \phi_e \right\|_{\underline{H}^{-1}(B_{tr})} \, \frac{dt}{t} \leq C r^{-1 }\int_1^\infty t^{-d-1}  (t r)^{1-\delta} \, \frac{dt}{t} \leq C r^{-\delta}.
\end{equation*}
Proofs of~\eqref{e.fluxbasecase} and~\eqref{e.energydensitybasecase} follow similarly. 
\end{proof}

We next prove the lemma used in the previous argument.
 
\begin{lemma}
\label{l.layercake}
Fix $r>0$, $k \in \N \cup \{0\}$, $p \in [1,\infty]$, $A>0$ and $\beta>0$. Suppose that $\psi \in W_{\mathrm{loc}}^{k,\frac{p}{p-1}}(\R^d)$ satisfies
\begin{equation} \label{e.layercakeweight}
\left\| \psi \right\|_{\underline{W}^{k,\frac{p}{p-1}}(B_1)} +   \sup_{t>1} \left( t^{\beta + k} \left\| \psi \right\|_{\underline{W}^{k,\frac{p}{p-1}}(B_{2t} \setminus B_t)} \right)
\leq 
A . 
\end{equation}
Then there exists $C(d)<\infty$ such that 
\begin{equation*} 
\left| \int_{\R^d} f(x) r^{-d} \psi\left( \tfrac{x}{r} \right) \, dx \right| 
\leq  
C A  
r^{-k} \int_{1}^\infty  t^{d - \beta -k} \left\| f \right\|_{\underline{W}^{-k,p}(B_{tr})} \, \frac{dt}{t}.
\end{equation*}
\end{lemma}
\begin{proof}
Let $\eta \in C_0^\infty(B_1)$ be such that $\eta = 1$ in $B_{1/2}$. Set $\eta_0  = 0$ and, for $m \in \N$, $\eta_m = \eta_0(2^{-m}\cdot)$.  Changing variables and decomposing $\psi$ using $\eta_m$'s, we get 
\begin{align} \notag 
\int_{\R^d} f(x) r^{-d} \psi\left(\tfrac{x}{r}\right) \, dx & = \int_{\R^d} f(r x) \psi\left(x \right) \, dx 
\\ & = \sum_{m=0}^\infty \int_{\R^d} f(r x) \psi\left(x \right) \left(\eta_{m+1}(x) - \eta_m(x) \right) \, dx .
\end{align}
Now, 
\begin{multline*} 
\left|  \int_{\R^d} f(r x) \psi\left(x \right) \left(\eta_{m+1}(x) - \eta_m(x) \right) \, dx \right| 
\\ 
\leq C2^{md} r^{-k} \left\| f \right\|_{\underline{W}^{-k,p}\left(B_{2^{m+1}r}\right)} \left\|\psi \left(\eta_{m+1} - \eta_m \right)  \right\|_{\underline{W}^{k,\frac{p}{p-1}}\left(B_{2^{m+1}}\right)} .
\end{multline*}
By assumption~\eqref{e.layercakeweight}, 
\begin{equation*} 
\left\|\psi \left(\eta_{k+1} - \eta_k \right)  \right\|_{\underline{W}^{k,\frac{p}{p-1}}\left(B_{2^{m+1}}\right)} \leq C A 2^{-m(\beta+k)} \left\| \eta_0 \right\|_{W_0^{k,\infty}(B_1)},
\end{equation*}
and hence we arrive at
\begin{equation*} 
\left|  \int_{\R^d} f(r x) \psi\left(x \right)  \, dx \right| \leq C r^{-k} \sum_{m=0}^\infty  2^{m(d-\beta-k)} \left\| f \right\|_{\underline{W}^{-k,p}\left(B_{2^{m+1}r}\right)} \left\| \eta_0 \right\|_{W_0^{k,\infty}(B_1)}.
\end{equation*}
The result follows from this. 
\end{proof}

The previous proposition and the stationarity of $\nabla \phi_e$ imply strong bounds in expectation for the same quantities, a fact which will be useful to us later.

\begin{lemma}
\label{l.correctorbias}
There exist $c(d,\Lambda)>0$ and $C(d,\Lambda)<\infty$ such that, for every $e\in\partial B_1$, $x\in\Rd$ and $r\geq 1$,
\begin{equation} 
\label{e.biascorrectorsgrad}
\left| \E \left[ \int_{\Phi_{x,r}} \nabla \phi_e \right] \right| \leq C\exp\left( -cr^2 \right),
\end{equation}
\begin{equation} 
\label{e.biascorrectorsflux}
\left| \E \left[ \int_{\Phi_{x,r}} \a(e+\nabla \phi_e) \right] -\ahom e\right| \leq C\exp\left( -cr^2 \right)
\end{equation}
and 
\begin{equation} 
\label{e.biascorrectorsenergy}
\left| \E \left[ \int_{\Phi_{x,r}} \frac12(e+\nabla \phi_e)\cdot \a(e+\nabla \phi_e) \right] -\frac12 e\cdot \ahom e\right| \leq C\exp\left( -cr^2 \right).
\end{equation}
\end{lemma}
\begin{proof}
We prove only the first estimate, since the argument for the others are similar. By the stationarity of $\nabla\phi_e$ and the fact that 
\begin{equation*} \label{}
\fint_{[0,1]^d} \left| \E \left[ \nabla \phi_e(x) \right] \right|^2\,dx \leq \E \left[ \fint_{[0,1]^d} \left| \nabla \phi_e(x) \right|^2\,dx\right] \leq C,
\end{equation*}
we deduce that the map
\begin{equation*} \label{}
x\mapsto \E \left[ \nabla \phi(x) \right]
\end{equation*}
is a $\Zd$--periodic function belonging to $L^2_{\mathrm{loc}}(\Rd)$. In view of~\eqref{e.gradientbasecase}, it therefore suffices to show that, for any $\Zd$--periodic function $f\in L^1_{\mathrm{loc}}(\Rd)$, we have
\begin{equation} 
\label{e.periodicheatexpdecay}
\left\| f \ast \Phi(t,\cdot) - \left( f \right)_{[0,1]^d} \right\|_{L^\infty(\Rd)} 
\leq C  \exp\left( - ct \right)\left\| f - \left( f \right)_{[0,1]^d} \right\|_{L^1([0,1]^d)} . 
\end{equation}
We may assume that $\left( f \right)_{[0,1]^d}= 0$. We see then that the claim is equivalent to the exponential decay in $L^\infty$ of a solution of the heat equation on the torus with mean zero initial data. This is a classical fact that can be seen in several ways and we outline one such proof in Exercise~\ref{ex.heatdecay} below.
\end{proof}

\begin{exercise}
\label{ex.heatdecay}
Prove~\eqref{e.periodicheatexpdecay}: that is, there exist constants $C(d) <\infty$ and $c(d) >0$ such that, for every~$\Zd$-periodic function $f\in L^1_{\mathrm{loc}}(\Rd)$ and $t\in [1,\infty)$, 
\begin{equation} 
\label{e.periodicheatflowdecay}
\sup_{x\in\Rd} \left| \left( f \ast \Phi(t,\cdot) \right)(x) - \left( f \right)_{[0,1]^d} \right| \leq C \left\| f - \left( f \right)_{[0,1]^d} \right\|_{L^1([0,1]^d)} \exp (-ct). 
\end{equation}
Hint: Assume $\left( f \right)_{[0,1]^d}=0$, let $f(t,x) = \left( f \ast \Phi(t,\cdot) \right)(x)$, which is $\Zd$--periodic in~$x$ and solves the heat equation, and argue that  
\begin{equation*} \label{}
\partial_t \int_{[0,1]^d} \frac12 \left| f(t,x) \right|^2\,dx = - \int_{[0,1]^d}  \left| \nabla f(t,x) \right|^2\,dx. 
\end{equation*}
Deduce by the Poincar\'e inequality (note $\left( f(t,\cdot) \right)_{[0,1]^d}=0$ since heat flow preserves mass) that
\begin{equation} 
\label{e.snapoincca}
\partial_t \int_{[0,1]^d} \frac12 \left| f(t,x) \right|^2\,dx 
\leq -c \int_{[0,1]^d} \frac12 \left| f(t,x) \right|^2\,dx. 
\end{equation}
Deduce that the quantity $\int_{[0,1]^d} \frac12 \left| f(t,x) \right|^2\,dx$ decays exponentially. Now use pointwise bounds for the heat equation (or the semigroup property again) to conclude.
\end{exercise}

\begin{exercise}
Show that the optimal Poincar\'e constant for $\Zd$-periodic, mean-zero  functions is $1/2\pi$. That is, for every $\Zd$-periodic $f\in L^2_{\mathrm{loc}}(\Rd)$ with $(f)_{[0,1]^d} =0$,
\begin{equation*} \label{}
\int_{[0,1]^d}  \left| f(x) \right|^2\,dx 
\leq 
\frac{1}{4\pi^2} \int_{[0,1]^d}  \left| \nabla f(x) \right|^2\,dx.
\end{equation*}
Deduce that the best constant $c$ in~\eqref{e.snapoincca} is~$8\pi^2$, and therefore the best constant~$c$ in~\eqref{e.periodicheatflowdecay} is~$4\pi^2$. (In particular, it can be taken independently of~$d$.) 
\end{exercise}

\section{Boundary regularity}
\label{s.boundaryreg}

In this section, we extend some of the interior regularity estimates presented in previous sections by giving global versions for Dirichlet boundary conditions. In particular, we prove a~$C^{0,1}$-type estimate in~$C^{1,\gamma}$ domains with boundary data in~$C^{1,\gamma}$. We also give Liouville-type theorems by characterizing the set of solutions in a half space which have polynomial traces on the boundary plane and grow at most polynomially at infinity. This allows us to give a~$C^{k,1}$-type regularity result for arbitrary solutions in half-balls. 

\smallskip

The first result is an extension of Theorem~\ref{t.Lipschitz} to a neighborhood of a boundary point of a $C^{1,\gamma}$ domain $U$ for solutions which satisfy a given $C^{1,\gamma}$ Dirichlet boundary condition on~$\partial U$. Since the domain~$U$ imposes a macroscopic length scale, it is natural to denote the microscopic scale by a small parameter~$\ep>0$.
\index{$C^{0,1}$ estimate}

\begin{theorem}
[{Quenched $C^{0,1}$-type estimate, boundary version}]
\label{t.Lipschitz.boundary}
Fix $s\in (0,d)$, $\gamma \in (0,1]$ and a $C^{1,\gamma}$ domain~$U$ with $0\in \partial U$. 
There exist~$C(s,\gamma,U,d,\Lambda)<\infty$, and a random variable $\X_s :\Omega \to [1,\infty]$ satisfying
\begin{equation}
\label{e.XLipschitz.global}
\X_s = \O_s\left( C \right),
\end{equation}
such that the following holds: for every $\ep\in \left( 0,\frac12 \right]$, $R\geq 2\ep \X_s$ and functions~$g\in C^{1,\gamma}(\partial U)$ and $u\in H^1(B_R \cap U)$ satisfying
\begin{equation}
\label{e.wksrt2.global}
\left\{
\begin{aligned}
& -\nabla \cdot \left( \a\left( \tfrac x\ep \right)\nabla u \right) = 0 & \mbox{in} & \ B_R \cap U,\\
& u = g & \mbox{on} & \ B_R \cap \partial U,
\end{aligned}
\right.
\end{equation}
we have, for every $r\in \left[ \ep \X_s ,\frac12 R\right]$, the estimate 
\begin{multline}
\label{e.Lipschitz.global}
 \frac1{r} \left\| u -g(0) \right\|_{\underline{L}^2(B_{r}\cap U)} \leq
\frac{C}{R} \left\| u - g(0) \right\|_{\underline{L}^2(B_{R}\cap U)}
\\
+ CR^\gamma \left[ \nabla g \right]_{C^{0,\gamma}(B_R \cap U)} + C \left\| \nabla g \right\|_{L^\infty(B_R\cap U)} .
\end{multline}
\end{theorem}

\begin{remark}
\label{r.gradbound.global}
By the global Caccioppoli inequality (cf. Lemma~\ref{l.Caccioppoli appendix glob}), we can write the estimate~\eqref{e.Lipschitz.global}, for every  $r\in \left[ \ep \X_s , R\right]$, as
\begin{equation}
\label{e.Lipschitz.global.grad}
\left\| \nabla u \right\|_{L^2(B_r\cap U)} 
\leq 
C \left( \left\| \nabla u \right\|_{L^2(B_R\cap U)} +R^\gamma \left[ \nabla g \right]_{C^{0,\gamma}(B_R\cap \partial U)} +  \left\| \nabla g \right\|_{L^\infty(B_R\cap  \partial U)} \right). 
\end{equation}
Indeed, one can test the weak formulation of $u$ with $(u-g)\phi^2$, where  $g$ is extended to $B_R \cap \partial U$ and $\phi$ is a cut-off function in $B_{2r}$, and then get 
\begin{equation*} 
\left\| \nabla u \right\|_{L^2(B_r\cap U)} \leq C \left( \frac1r \left\| u - g(0) \right\|_{L^2(B_{2r} \cap U)}  +  \left\| \nabla g \right\|_{L^\infty(B_{2r} \cap U)} \right).
\end{equation*}
On the other hand, by the triangle inequality and Poincar\'e's inequality, 
\begin{align*} 
\frac1R \left\| u - g(0) \right\|_{L^2(B_{R} \cap U)} & \leq \frac1R \left\| u - g \right\|_{L^2(B_{R} \cap U)} +  \frac1R \left\| g - g(0) \right\|_{L^2(B_{R} \cap  \partial  U)} 
\\ & \leq C \left\| \nabla u - \nabla g \right\|_{L^2(B_{R} \cap U)}  +  \frac1R \left\| g - g(0) \right\|_{L^2(B_{R} \cap  \partial U)}  \\
\notag & \leq C \left( \left\| \nabla u\right\|_{L^2(B_{R} \cap U)}  +  \left\| \nabla g \right\|_{L^\infty(B_{R} \cap  \partial  U)} \right).
\end{align*}
Combining the above two displays with~\eqref{e.Lipschitz.global} leads to~\eqref{e.Lipschitz.global.grad}. 
\end{remark}

\begin{remark}
\label{r.local+global}
If we relax the assumption that $0 \in \partial U$ to the assumption that $B_R \cap \partial U \neq \emptyset$, one may combine Theorems~\ref{t.Lipschitz} and~\ref{t.Lipschitz.boundary} to obtain
\begin{multline}
\label{e.Lipschitz.local+global}
 \frac1{r} \left\| u - m_r \right\|_{\underline{L}^2(B_{r}\cap U)} \leq
\frac{C}{R} \left\| u - g(0) \right\|_{\underline{L}^2(B_{R}\cap U)}
\\
+ CR^\gamma \left[ \nabla g \right]_{C^{0,\gamma}(B_R \cap U)} + C \left\| \nabla g \right\|_{L^\infty(B_R\cap U)} ,
\end{multline}
where 
\begin{equation*} 
m_r := \left\{
\begin{aligned}
& g(0) & \mbox{if } & \dist(0,\partial U) \leq r,\\
& (u)_{B_r} & \mbox{if } &  \dist(0,\partial U) > r.
\end{aligned}
\right.\end{equation*} 
\end{remark}

\smallskip

The proof of Theorem~\ref{t.Lipschitz.boundary} closely follows the one of Theorem~\ref{t.Lipschitz}, and we only need to modify the argument slightly. The classical boundary regularity estimate for~$\ahom$-harmonic functions which we need is stated as follows. Let~$\gamma \in (0,1]$, $U$ be a~$C^{1,\gamma}$ domain and~$u$ an $\ahom$-harmonic function in $U$ such that $g:=u\vert_{\partial U} \in C^{1,\gamma}(\partial U)$. Then there exists $C(d,\beta,U)<\infty$ such that, for every $x_0 \in \partial U$ and $r \in \left(0,\tfrac 12 \diam(U) \right)$, 
\begin{multline} 
\label{e.classicalbndryreg}
\left\| \nabla u \right\|_{L^\infty\left(B_{r/2}(x_0) \cap U \right)} + r^\gamma \left[ \nabla u \right]_{C^{0,\gamma}\left(B_{r/2}(x_0) \cap U \right)}   
\\
\leq 
C \left( \left\| \nabla u \right\|_{\underline{L}^2\left(B_{r}(x_0) \cap U\right)} + \left\| \nabla g \right\|_{L^\infty\left(B_r\cap \partial U\right)} + r^\gamma \left[ \nabla g \right]_{C^{0,\gamma}\left(B_r\cap \partial U\right)}   
\right).
\end{multline}
The estimate can be easily obtained by iterating the result of the following lemma.

\begin{lemma} \label{l.classicalbndryreg1}
Let $\gamma\in(0,1]$ and $R>0$. Suppose that $U$ is a $C^{1,\gamma}$-domain with $0\in \partial U$ and $g \in C^{1,\gamma}(\partial B_{R} \cap U)$. Let $u$ solve
\begin{equation} 
\label{e.wksrt2.global-pre2}
\left\{
\begin{aligned}
& -\nabla \cdot \left( \ahom \nabla u \right) = 0 & \mbox{in} & \ B_R \cap U,\\
& u = g & \mbox{on} & \ B_R \cap \partial U,
\end{aligned}
\right.
\end{equation}
Then there exists a constant $C(\gamma,U,d,\Lambda)<\infty$ such that, for every $\theta \in \left( 0, \frac18 \right]$, 
\begin{align} 
\label{e.classicalbndryreg1}
\inf_{p\in \mathcal{L}_g }  \left\| u  - p  \right\|_{\underline{L}^2(B_{\theta R} \cap U)}  & \leq C\theta^2 \inf_{p\in \mathcal{L}_g } \left\| u  - p  \right\|_{\underline{L}^2(B_{R} \cap U)} 
\\ \notag & \quad +  C \theta^{-\frac d2} R^\gamma  \left( \left\| u - \ell_g \right\|_{\underline{L}^2(B_{R} \cap U)} + R \left[ \nabla g \right]_{C^{0,\gamma}(\partial B_{R} \cap U)}   
\right),
\end{align}
where 
\begin{equation} \label{e.L_g}
\ell_g(x) := g(0) + \nabla_T g(0) \cdot x  \quad \mbox{and} \quad \mathcal{L}_g := \left\{ \ell_g(x) + q \cdot x \, : \, q \in \R^d \right\}.
\end{equation}
Above $\nabla_T g(0)$ stands for the tangential gradient of $g$. 
\end{lemma}

Before presenting the proof, let us prove yet another preliminary version in the special case when the solution is zero on a hyperplane. 
We denote, for each $e\in\partial B_1$ and $r >0$,  the half space
\begin{equation*} \label{}
H_e := \left\{ x\in\Rd \,:\, x\cdot e > 0\right\}, \quad   E_r :=  B_r \cap H_e  \quad \mbox{and} \quad E_r^0 := B_r \cap \partial H_e. 
\end{equation*}

\begin{lemma} \label{l.classicalbndryreg2}
Let $R>0$ and let $u$ solve
\begin{equation} 
\label{e.wksrt2.global-pre3}
\left\{
\begin{aligned}
& -\nabla \cdot \left( \ahom \nabla u \right) = 0 & \mbox{in} & \  E_R ,\\
& u = 0 & \mbox{on} & \ E_R^0 .
\end{aligned}
\right.
\end{equation}
Then, for $k \in \N$, there exists $C(k,d,\Lambda)<\infty$ such that 
\begin{equation} \label{e.classicalbndryreg2}
\left\|\nabla^k u \right\|_{L^\infty(E_{R/2} )} \leq \frac{C}{R^k}  \left\| u \right\|_{\underline{L}^2(E_{R})} .
\end{equation}
\end{lemma}

\begin{proof}
After rotation we may assume that $e=e_d$. Let us denote by $\nabla' $ the gradient with respect to the $d-1$ first variables. Note that, for every $m \in \N$, each component of $\mathbf{u}_m := (\nabla')^m u$ solves the same equation as~$u$ and 
$\mathbf{u}_m = 0$ on $E_R^0$. The Caccioppoli estimate (cf. Lemma~\ref{l.Caccioppoli appendix glob}) then implies that, for $r \in (0,R)$ and $\sigma \in \left[\tfrac12,1\right)$, 
\begin{equation*} 
\left\| \nabla \mathbf{u}_m \right\|_{\underline{L}^2(E_{\sigma r} )} \leq \frac{C_m}{(1-\sigma)} \frac1r \left\|  \mathbf{u}_m \right\|_{\underline{L}^2(E_{r} ) }.
\end{equation*}
This gives iteratively that 
\begin{equation} \label{e.Cacc.hsbnd}
\left\| \nabla \mathbf{u}_{m-1} \right\|_{\underline{L}^2(E_{r/2} )} \leq \frac{C_m}{r^{m}} \ \left\|  u \right\|_{\underline{L}^2(E_{r} ) }.
\end{equation}
By the equation, we deduce the existence of a bounded matrix $\mathbf{A}_m$ such that 
\begin{equation*} 
\partial_{d}^2 \mathbf{u}_{m-1}(x) = \mathbf{A}_m \nabla \mathbf{u}_m(x).
\end{equation*}
Indeed, this follows from the positivity of $\ahom_{dd}$ or, in the case of elliptic systems, from the positive definitiveness of the matrix $\left(\ahom_{dd}^{\alpha \beta}\right)_{\alpha\beta}$.  Using this formula repeatedly it is easy to see that there exists a bounded matrix $\mathbf{B}_{m}$ such that 
\begin{equation*} 
\nabla^m u = \mathbf{B}_m \nabla \mathbf{u}_{m-1}. 
\end{equation*}
Together with~\eqref{e.Cacc.hsbnd} we get
\begin{equation*} 
\left\| \nabla^m u \right\|_{\underline{L}^2(E_{R/2} )} \leq \frac{C_m}{R^{m}} \left\|  u \right\|_{\underline{L}^2(E_{R} ) }.
\end{equation*}
Applying  then Morrey's inequality for sufficiently large $m$ finishes the proof.
\end{proof}

\begin{proof}[Proof of Lemma~\ref{l.classicalbndryreg1}]
Without loss of generality we may assume that the normal vector to~$\partial U$ at the origin is~$e_d$. Moreover, since $u - g(0) - \nabla g(0) \cdot x$ still solves the same equation, we may assume that both $g(0)=0$ and $\nabla g(0) = 0$. 

\smallskip

\emph{Step 1.} By the assumption that $U$ is a $C^{1,\gamma}$-domain, we can find $r_0$ small enough so that the boundary of $U$ is  the graph $x_d  = \psi(x' ) $, where we denote $x' = (x_1,\ldots,x_{d-1})$. 
The new coordinate system is given by $y = \Theta(x) = (x',x_d - \psi(x'))$ so that the boundary of $U$ in $B_{r_0}$ in the new coordinates is $E_{r_0}^0 = B_{r_0} \cap \partial H_{e_d}$. Clearly $\Theta$ is invertible, since $\nabla_T \psi(0) = 0$, $\psi \in C^{1,\gamma}(E_{r_0}^0)$, and $r_0$ is small. Moreover,
\begin{equation} 
\label{e.Theta vs Id}
\sup_{x \in B_{r} \cap U} \left| \Theta(x) - x \right| \leq \sup_{x \in B_{r} \cap U} \left| \psi(x ) \right| \leq r^{1+\gamma} \left[\nabla \psi\right]_{C^{0,\gamma}} ,
\end{equation}
and
\begin{equation} 
\label{e.Theta vs Id 2}
\sup_{x \in B_{r} \cap U} \left| \nabla \Theta(x) - \Id \right| \leq \sup_{x \in B_{r} \cap U} \left| \nabla \psi(x) \right| \leq r^{\gamma} \left[\nabla \psi\right]_{C^{0,\gamma}} ,
\end{equation}
again by the fact that $\psi(0) = 0$ and $\nabla \psi(0) = 0$. Consider
\begin{equation*} 
\tilde u = u \circ \Theta^{-1} \quad \mbox{and} \quad \tilde g = g \circ \Theta^{-1},
\end{equation*}
so that $\tilde u = \tilde g$ on $E_{r_0}^0$. 

\smallskip

\emph{Step 2.}  Observe that $\tilde u$ satisfies the equation 
\begin{equation*} 
\nabla_y \cdot \left(\mathbf{b}(y) \nabla_y \tilde u(y)\right) = 0 \quad \mbox{in } \; E_r ,
\end{equation*}
where we define
\begin{multline*} 
\mathbf{b}_{ij} (y) := 
\ahom_{ij} 
- \delta_{id} \sum_{k=1}^{d-1} \partial_{k} \psi(y') \ahom_{kj} 
- \delta_{jd} \sum_{k=1}^{d-1} \partial_{k} \psi(y') \ahom_{ik} 
\\ + \delta_{id} \delta_{jd} \sum_{k,m=1}^{d-1} \partial_{k} \psi(y') \partial_{m} \psi(y') \ahom_{km}.
\end{multline*}
In view of the fact that $\nabla \psi(0) = 0$, we have, for every $r \in (0,r_0]$,
\begin{equation} 
\label{e.b close to ahom}
\left\| \mathbf{b}(\cdot) - \ahom\right\|_{L^\infty(E_{r})} \leq C r^\gamma \left[\nabla \psi\right]_{C^{0,\gamma}(E_{r})}.
\end{equation}
Therefore, for $r_0$ small enough, the ellipticity of~$\mathbf{b}$ is dictated by~$\ahom$. 
Since~$\nabla_T g(0) = 0$, we can extend $\tilde g$ to $B_{r_0}$ so that, for $r \in (0,r_0]$,   
\begin{equation} \label{e.tilde g est}
\left\| \nabla \tilde g \right\|_{L^\infty(B_r)} \leq C r^{\gamma} \left[\nabla g\right]_{C^{0,\gamma}(E_{r_0}^0)} .
\end{equation}

\smallskip

\emph{Step 3.} 
We next prove a suitable comparison estimate. Solve 
\begin{equation} 
\label{e.bnd comp sol 1}
\left\{
\begin{aligned}
& -\nabla \cdot \left(\mathbf{b}(\cdot) \nabla \tilde v \right) = 0  & \mbox{in} & \ E_{r/2},\\
& \tilde v = \tilde u - \tilde g   & \mbox{on} & \ \partial B_{r} \cap E_r, \\
& \tilde v = 0 & \mbox{on} & \ E_{r}^0,
\end{aligned}
\right.
\end{equation}
so that, by~\eqref{e.tilde g est},
\begin{equation} \label{e.bnd u vs v}
\fint_{E_r} \left|\nabla \tilde u(x) - \nabla \tilde v(x) \right|^2 \, dx \leq C \left\| \nabla \tilde g \right\|_{L^\infty(E_r)}^2 \leq  
C r^{2\gamma} \left[\nabla g\right]_{C^{0,\gamma}(E_{r_0}^0)} .
\end{equation}
On the other hand, solving 
\begin{equation} 
\label{e.bnd comp sol 2}
\left\{
\begin{aligned}
& -\nabla \cdot \left(\ahom \nabla \tilde w \right) = 0   & \mbox{in} & \ E_{r/2},\\
& \tilde w = \tilde v    & \mbox{on} & \ \partial  E_{r},
\end{aligned}
\right.
\end{equation}
leads to
\begin{equation} \label{e.bnd w vs v}
\fint_{E_r} \left|\nabla \tilde w(x) - \nabla \tilde v(x) \right|^2 \, dx 
\leq C \left[\psi \right]_{C^{0,\gamma}(E_{r_0}^0 )} \fint_{E_r} \left|\nabla \tilde v(x) \right|^2 \, dx.
\end{equation}
Therefore, combining~\eqref{e.bnd u vs v} and~\eqref{e.bnd w vs v} and allowing~$C$ to depend also on $\psi$ (which depends on~$\partial U$), we obtain
\begin{equation}  
\label{e.bnd w vs u}
\fint_{E_r} \left|\nabla \tilde w(x) - \nabla \tilde u(x) \right|^2 \, dx 
\\
\leq C r^{2\gamma} \left( \left[\nabla g\right]_{C^{0,\gamma}(E_{r_0}^0)} + \fint_{E_r} \left|\nabla \tilde u(x) \right|^2 \, dx   \right).
\end{equation}
Testing the equation of $\tilde u$ with $(u-\tilde g)\phi^2$, where $\phi$ is a cut-off function, we also obtain
\begin{equation*} 
\fint_{E_r} \left|\nabla \tilde u(x) \right|^2 \, dx \leq  \frac{C}{r^2} \fint_{E_{2r}} \left|\tilde u(x) \right|^2 \, dx  + C \left\| \nabla \tilde g \right\|_{L^\infty(E_r)}^2.
\end{equation*}
Consequently, since $\tilde u = \tilde w + \tilde g$ on $\partial  E_r $, we get by Poincar\'e's inequality that 
\begin{equation}  \label{e.bnd w vs u 2}
\left\|  \tilde w - \tilde u \right\|_{\underline{L}^2(E_r )} 
\\
\leq C r^{\gamma} \left( r \left[ \nabla g\right]_{C^{0,\gamma}(E_{2r}^0)} + \left\| \tilde u  \right\|_{\underline{L}^2(E_{2r} )}    \right).
\end{equation}

\smallskip

\emph{Step 4.} We now prove~\eqref{e.classicalbndryreg1}.  Since $\tilde w = 0$ on $\partial E_{r}^0$, we may use Lemma~\ref{l.classicalbndryreg2} to obtain
\begin{equation*} 
\inf_{p\in \mathcal{L}_g }  \left\|  \tilde w - p  \right\|_{\underline{L}^2(E_{\theta r})} 
\leq C\theta^2 \inf_{p\in \mathcal{L}_g } \left\|  \tilde w - p  \right\|_{\underline{L}^2(E_{r})} ,
\end{equation*}
where $ \mathcal{L}_g $ is defined in~\eqref{e.L_g}. Now~\eqref{e.classicalbndryreg1} follows using the triangle inequality,~\eqref{e.bnd w vs u 2} and changing variables to obtain the estimate for $u$ instead of $\tilde u$. 
\end{proof}

We next formulate the harmonic approximation lemma near the boundary of a $C^{1,\gamma}$ domain~$U$.

\begin{lemma}[Harmonic approximation, boundary version]
\label{l.harmonicapproximation.boundary}
Fix $s\in (0,d)$. Suppose $\gamma\in (0,1]$ and $U$ is a $C^{1,\gamma}$ domain with $0\in \partial U$. There exist constants $\alpha(d,\Lambda)>0$, $C(s,\gamma,U,d,\Lambda)<\infty$, and a random variable $\X_s:\Omega \to [1,\infty]$ satisfying
\begin{equation}
\label{e.Xharmonicapproximation.boundary}
\X_s = \O_s\left( C \right)
\end{equation}
such that the following holds: for every~$\ep\in(0,1)$, $R\geq \X_s$ and $u\in H^1(B_{\ep R}\cap U)$ satisfying
\begin{equation} 
\label{e.wksrt.boundary}
-\nabla \cdot \left( \a\left(\tfrac x\ep \right) \nabla u \right) = 0 \quad \mbox{in} \ B_{\ep R} \cap U,
\end{equation}
there exists a solution $\overline{u} \in H^1(B_{\ep R/2}\cap U)$ satisfying 
\begin{equation*} \label{}
\left\{ 
\begin{aligned}
&-\nabla \cdot \left( \ahom \nabla \bar u \right) = 0 & \mbox{in} & \ B_{\ep R/2}\cap U, \\
&  \quad \bar u = u & \mbox{on} & \ B_{\ep R/2} \cap \partial U, 
\end{aligned} 
\right.
\end{equation*}
and
\begin{equation}
\label{e.harmonicapproximation.boundary}
\left\| u - \overline{u}  \right\|_{\underline{L}^2(B_{\ep R/2} \cap U)} \leq 
CR^{-\alpha(d-s)} \left( \left\|  u - g(0) \right\|_{\underline{L}^2(B_{\ep R} \cap U)} + \ep R \left\|  \nabla g \right\|_{L^\infty(B_{\ep R} \cap \partial U)}  \right).
\end{equation}
\end{lemma}
\begin{proof}
The proof is almost the same as that of Proposition~\ref{p.harmonicapproximation}, the only difference being the use of the global (rather than interior) Meyers estimate, which is proved in Theorem~\ref{t.Meyers appendix global}. Indeed, we have that 
\begin{equation*} 
 \left\|  \nabla u \right\|_{\underline{L}^{2+\delta}(B_{\ep R/2} \cap U)} \leq C \left( \left\|  \nabla u \right\|_{\underline{L}^{2}(B_{ \frac 34 \ep R} \cap U)} + \left\|  \nabla g \right\|_{L^\infty(B_{\ep R} \cap \partial U)}   \right),
\end{equation*}
and the Caccioppoli estimate (test with $(u-g)\phi^2$) yields 
\begin{equation*} 
\left\|  \nabla u \right\|_{\underline{L}^{2}(B_{ \frac 34 \ep R} \cap U)} \leq C \left( \frac{1}{\ep R} \left\| u - g(0) \right\|_{\underline{L}^{2}(B_{ \ep R} \cap U)} + \left\|  \nabla g \right\|_{L^\infty(B_{\ep R} \cap \partial U)}   \right).
\end{equation*}
The rest of the proof is left to the reader. 
\end{proof}

We next formulate a version of Lemma~\ref{l.harmapproxlipschitz}, which holds up to the boundary. 

\begin{lemma}
\label{l.harmapproxlipschitz.boundary}
Fix $\alpha \in (0,1]$,  $K\geq 1$ and $X\geq 1$. Suppose that $\gamma>0$ and $U$ is a $C^{1,\gamma}$ domain with $0\in \partial U$. Also fix~$g\in C^{1,\gamma}(\overline{U})$. 
Let $\ep \in (0,1)$, $R\geq 2\ep X$ and $u\in L^2(B_{R})$ have the property that, for every $r\in \left[ X \ep, \frac12R \right]$, there exists $w\in H^1(B_{r})$ which is a weak solution of 
\begin{equation} 
\label{e.wksrt2.global-pre}
\left\{
\begin{aligned}
& -\nabla \cdot \left( \ahom \nabla w_r \right) = 0 & \mbox{in} & \ B_r \cap U,\\
& w_r = g & \mbox{on} & \ B_r \cap \partial U,
\end{aligned}
\right.
\end{equation}
satisfying 
\begin{equation} 
\label{e.wharmapproxit.boundary2}
\left\| u - w_r \right\|_{\underline{L}^2(B_{r/2} \cap U)} \leq K \left( \frac{\ep}{r} \right)^{\alpha} \left( \left\| u - g(0) \right\|_{\underline{L}^2(B_{r} \cap U)} + r \left\|  \nabla g \right\|_{L^\infty(B_{r} \cap \partial U)}  \right). 
\end{equation}
Then there exists $C(\alpha,K,U,d,\Lambda)<\infty$ such that, for every $r\in \left[ \ep X, R\right]$, 
\begin{multline}
\label{e.elipsc.boundary}
 \frac1{r} \left\| u - g(0) \right\|_{\underline{L}^2(B_{r} \cap U)}
\leq \frac{C}{R} \left\| u - g(0) \right\|_{\underline{L}^2(B_{R} \cap U)}
\\ + C R^{\gamma} \left[ \nabla g \right]_{C^{0,\gamma}(\partial B_{R} \cap U)}  + C \left\| \nabla g \right\|_{L^{\infty}(\partial B_{R} \cap U)} .
\end{multline}
\end{lemma}

\begin{proof}
We first deduce an appropriate decay estimate for $u$. Use~\eqref{e.wharmapproxit.boundary2}  to get
\begin{equation}
\label{e.harmonicapproximation.boundary-appl}
\left\| u - w_r \right\|_{\underline{L}^2(B_{r/2} \cap U)} 
\leq C\left( \frac{\ep}{r} \right)^{\alpha} \left( \left\| u - g(0) \right\|_{\underline{L}^2(B_{r} \cap U)} +  r \left\|  \nabla g \right\|_{L^\infty( \partial B_{r} \cap U)} \right).
\end{equation}
Apply then~\eqref{e.classicalbndryreg1} to $w_r$: for any $\theta \in \left(0,\frac 18\right]$ we have
\begin{align} 
\label{e.classicalbndryreg3}
\inf_{p\in \mathcal{L}_g }  \left\| w_r  - p  \right\|_{\underline{L}^2(U \cap  B_{\theta r})}  & \leq C\theta^2 \inf_{p\in \mathcal{L}_g } \left\| w_r  - p  \right\|_{\underline{L}^2(B_{r/2} \cap U)} 
\\ \notag & \quad +  C \theta^{-\frac d2} r^\gamma  \left( \left\| w_r - \ell_g \right\|_{\underline{L}^2(B_{r/2} \cap U)} + r \left[ \nabla g \right]_{C^{0,\gamma}(\partial B_{r} \cap U)}   
\right),
\end{align}
where $\ell_g = g(0) + \nabla_T g(0) \cdot x$ and $\mathcal{L}_g := \left\{ \ell_g(x) + q \cdot x \, : \, q \in \R^d \right\}$.  We now choose $\theta$ so that $C \theta = \frac12$ and set
\begin{equation*} 
E_1(r) := \frac1r \inf_{p \in \mathcal{L}_g} \left\| u - p \right\|_{\underline{L}^2(B_{r} \cap U)} \quad \mbox{and} \quad E_0(r) := \frac 1r \left\| u - \ell_g \right\|_{\underline{L}^2( B_{r} \cap U )}.
\end{equation*}
Rearranging terms after using the triangle inequality yields, for $r \in (\ep X, \ep R]$, 
\begin{align*} \label{e.u bndr iter 1}
E_1(\theta r) & \leq \frac12  E_1(r) + C\left( r^{\gamma}  + \left( \frac{\ep}{r} \right)^{\alpha}  \right) E_0(r) 
\\ & \quad + Cr^{\gamma} \left[ \nabla g \right]_{C^{0,\gamma}(\partial B_{R} \cap U)} 
+ C \left( \frac{\ep}{r} \right)^{\alpha}  \left\| \nabla g \right\|_{L^{\infty}(\partial B_{R} \cap U)} .
\end{align*}
The rest of the proof is very similar to Lemma~\ref{l.harmapproxlipschitz}. Indeed, for $r_0 \in  [\ep X,R]$ and $r_1 \in [\ep X,r_0]$, set
\begin{equation*} 
M := \max_{r \in [r_1,r_0]} E_0(r). 
\end{equation*}
Let $n \in \N$ be such that $\ep X \in (\theta^{n+1} r_0,\theta^n r_0]$ and take $k\in \N$ and $r_1 = \theta^{n-k} r_0 $. Summation then gives 
\begin{align*} \label{e.u bndr iter 2}
\sum_{j=0}^{n-k} E_1(\theta^j r) & \leq 2 E_1(r_0) + C \left( r_0^\gamma  +\theta^{k \alpha(d-s)} \right) M 
\\ & \quad + C r_0^{\gamma} \left[ \nabla g \right]_{C^{0,\gamma}(\partial B_{R} \cap U)}  + C X^{-\alpha}  \left\| \nabla g \right\|_{L^{\infty}(\partial B_{R} \cap U)} .
\end{align*}
Letting $p_j \in \mathcal{L}_g$ be the minimizing affine function in $E_1(\theta^j r)$ and $p_0 = \ell_g$, we then obtain 
\begin{align*} 
\max_{j \in \{1,\ldots,n-k\} } \left | \nabla p_j \right| & \leq 2 E_0(r_0) + C \left( r_0^\gamma  +\theta^{k \alpha(d-s)} \right) M 
\\ & \quad + C r_0^{\gamma} \left[ \nabla g \right]_{C^{0,\gamma}(\partial B_{R} \cap U)}  + C X^{-\alpha}  \left\| \nabla g \right\|_{L^{\infty}(\partial B_{R} \cap U)} .
\end{align*}
Thus we also get 
\begin{align*} 
M & \leq  C \theta^{-\frac d2}\left( r_0^\gamma  +\theta^{k \alpha(d-s)} \right) M 
\\ & \quad + C \theta^{-\frac d2} \left(E_0(r_0) + r_0^{\gamma} \left[ \nabla g \right]_{C^{0,\gamma}(\partial B_{R} \cap U)}  +  X^{-\alpha}  \left\| \nabla g \right\|_{L^{\infty}(\partial B_{R} \cap U)} \right) . 
\end{align*}
We now choose $k$ large enough and $r_0$ small enough so that
\begin{equation*} 
C \theta^{-\frac d2}\left( r_0^\gamma  +\theta^{k \alpha(d-s)} \right)  \leq \frac12.
\end{equation*}
Therefore we obtain after reabsorption that
\begin{equation*} 
M  \leq   C \theta^{-\frac d2} \left(E_0(r_0) + r_0^{\gamma} \left[ \nabla g \right]_{C^{0,\gamma}(\partial B_{R} \cap U)}  +  X^{-\alpha}  \left\| \nabla g \right\|_{L^{\infty}(\partial B_{R} \cap U)} \right) . 
\end{equation*}
From this~\eqref{e.Lipschitz.global} follows easily.
\end{proof}

Theorem~\ref{t.Lipschitz.boundary} is now a consequence of Lemmas~\ref{l.harmonicapproximation.boundary} and~\ref{l.harmapproxlipschitz.boundary}.

\begin{exercise}
\label{ex.C01.RHS.global}
State and prove a global version of Exercise~\ref{ex.C01.RHS}. 
\end{exercise}

We next present an analogue of Theorem~\ref{t.regularity} in half spaces with Dirichlet boundary data. In particular, we classify all solutions in the upper half space 
\begin{equation*} \label{}
\Rd_+:=\left\{ x_d > 0 \right\}
\end{equation*}
which have at most polynomial growth in~$\Rd_+$ and a polynomial trace on~$\partial\Rd_+$. 

\smallskip

\index{harmonic polynomial|(}
For the precise statement, see Theorem~\ref{t.Liouville.halfspace} below. We first discuss such a characterization for~$\ahom$-harmonic functions. Let us consider first the scalar case and assume that~$\ahom = \Id$, so that we are discussing harmonic functions. The first question we should address is (identifying~$\R^{d-1}$ with~$\partial \Rd_+$) whether every polynomial in~$d-1$ variables is the trace of a \emph{harmonic} polynomial in~$\Rd$. The answer is affirmative: for every polynomial $p$ on $\R^{d-1}$, we first extend it to $\Rd$ and then define a harmonic polynomial $q$ by
\begin{equation} 
\label{e.magicq}
q(x):= \sum_{n=0}^{\left\lceil\deg(p)/2\right\rceil} \frac{(-1)^n}{(2n)!} x_d^{2n}\Delta^n p(x).
\end{equation}
It is straightforward to check that $q$ is harmonic and that $q = p$ on~$\partial\Rd_+$. Observe also that $q$ has the same degree as $p$, and that if $p$ is homogeneous then so is~$q$. The mapping $p\mapsto q$ given by~\eqref{e.magicq} defines a linear injection from $\mathcal{P}(\R^{d-1})$ to $\Ahom(\Rd)$, but it is not a surjection: there are many nonzero harmonic polynomials which have zero trace on~$\partial \Rd_+$, for instance $x\mapsto x_d p(x)$, where $p$ is any harmonic polynomial of the first~$d-1$ variables. 

\smallskip 

Motivated by the previous discussion, we extend the result to a general hyperplane. Recall the definition of the half space, for each $e\in\partial B_1$, 
\begin{equation*} \label{}
H_e := \left\{ x\in\Rd \,:\, x\cdot e > 0\right\}  \quad \mbox{and} \quad E_r := H_e \cap B_r.
\end{equation*}
The purpose of Exercises~\ref{ex.liouville.bndr1} and~\ref{ex.liouville.bndr2} is to show that 
\begin{equation*} \label{}
\Ahom_k \left( H_e \right):= \left\{ u \in \Ahom(\Rd)\,:\, \limsup_{r\to \infty} r^{-(k+1)} \left\| u \right\|_{\underline{L}^2(E_r)} = 0, \ u \,\vert_{\partial H_e} = p \in \mathcal{P}_k(\Rd) \right\}
\end{equation*}
is simply the set of restrictions of elements of~$\Ahom_k=\Ahom_k(\Rd)$ to~$H_e$. In particular, traces of $\ahom$-harmonic polynomials exhaust the set of polynomials on $\partial H_e$. 

\index{harmonic polynomial|)}

\begin{exercise}
\label{ex.liouville.bndr1}
Verify the claim that
\begin{equation} 
\label{e.identifyAhomk+}
\Ahom_k \left( H_e \right)
=
\left\{ u\,\vert_{H_e} \,:\, u\in \Ahom_k \right\},
\end{equation}
using the following argument (which does not work in the case of elliptic systems): given~$u\in \Ahom_k\left(\Rd_+\right)$ with~$u\,\vert_{\Rd_+}=p\in \mathcal{P}_k$, check that $w\in \Ahom_k$, where~$w$ is the odd reflection of~$u-q$ and~$q\in \Ahom_k$ is defined in~\eqref{e.magicq}. Notice then that~$u=w+q$ in~$\Rd_+$. 
\end{exercise}

In the following lemma, we give a quantitative version of the Liouville principle in Exercise~\ref{ex.liouville.bndr1}, using an argument which extends to the case of elliptic systems.

\smallskip

\begin{lemma}
\label{l.oddreflection.gen}
Fix $k\in\N$, $e \in \partial B_1$, $\gamma \in (0,1]$, $r>0$ and $g \in C^{k,\gamma}(B_r \cap \partial H_e)$. Let $u$ be an $\ahom$-harmonic function $u \in H^1(E_r)$ satisfying $u=g$ on $B_r \cap \partial H_e$. Then there exist $C(k,d,\Lambda)<\infty$ and $p\in \Ahom_k$ such that, for every $s\in \left(0, r \right]$,
\begin{equation} 
\label{e.Ahomharm+}
\left\| u - p \right\|_{\underline{L}^2(E_s)} 
\leq 
C \left( \frac sr \right)^{k+1} \left\| u \right\|_{\underline{L}^2(E_r)} + C \left( \frac rs \right)^{\frac d2} r^{k+\gamma} \left[ \nabla^k g \right]_{C^{0,\gamma}(B_r \cap \partial H_e)}. 
\end{equation}
\end{lemma}

\begin{proof}
Letting $q \in \Ahom_k$ be an $\ahom$-harmonic polynomial agreeing with the $k$-th degree Taylor series of $g$ at origin, we see that, for $s \in (0,r]$,
\begin{equation*} 
\left\|q - g \right\|_{L^\infty(B_s \cap \partial H_e)} \leq C s^{k+\gamma} \left[ \nabla^k g \right]_{C^{0,\gamma}(B_s \cap \partial H_e)}.
\end{equation*}
 Thus, since $u - q$ is still $\ahom$-harmonic, we may without loss of generality assume that~$q = 0$. Extending $g$ to $E_r$ satisfying
 \begin{equation*} 
\sum_{j=0}^k r^j \left\| \nabla^j \Ext(g) \right\|_{L^\infty(E_r)}  \leq Cr^{k+\gamma} \left[ \nabla^k g \right]_{C^{0,\gamma}(B_r \cap \partial H_e)},
\end{equation*}
and taking~$v_r \in H^1(E_r)$ to be the~$\ahom$-harmonic function belonging to~$u - g + H_0^1(E_r)$,
we obtain by testing that 
\begin{equation*} 
\left\|\nabla u - \nabla v_r \right\|_{\underline{L}^2(E_r)} \leq C \left\| \nabla \Ext(g) \right\|_{\underline{L}^2(E_r)} \leq C r^{k-1+\gamma} \left[ \nabla^k g \right]_{C^{0,\gamma}(B_r \cap \partial H_e)}.
\end{equation*}
Poincar\'e's inequality hence gives
\begin{equation*} 
\left\|u - v_r \right\|_{\underline{L}^2(E_r)} \leq C r^{k+\gamma} \left[ \nabla^k g \right]_{C^{0,\gamma}(B_r \cap \partial H_e)}.
\end{equation*}
Since $v_r$ has zero boundary values on $E_r \cap \partial H_e$, Lemma~\ref{l.classicalbndryreg2} and the previous display imply the result by the triangle inequality. 
\end{proof}

\begin{exercise}
\label{ex.liouville.bndr2}
Verify~\eqref{e.identifyAhomk+} using Lemma~\ref{l.oddreflection.gen}.
\end{exercise}

\begin{lemma}
\label{l.u:decayestimate.bndr}
Fix $\alpha \in [0,1]$,  $K\geq 1$ and $X\geq 1$. 
Let $R\geq 2X$ and $u\in L^2(B_{R})$ have the property that, for every $r\in \left[ X, R \right]$, there exists $w_r \in H^1(B_{r/2})$ which is a weak solution of 
\begin{equation*} \label{}
\left\{
\begin{aligned}
& -\nabla \cdot \left( \ahom \nabla w_r \right) = 0 & \mbox{in} & \ E_r ,\\
& w_r = g & \mbox{on} & \  B_r \cap \partial H_e,
\end{aligned}
\right.
\end{equation*}
satisfying 
\begin{equation} 
\label{e.wharmapproxit.boundary3}
\left\| u - w_r \right\|_{\underline{L}^2(E_{r/2} )} \leq K r^{-\alpha} \left( \left\| u - g(0) \right\|_{\underline{L}^2(E_{r})} + r \left\|  \nabla_{e^\perp} g \right\|_{L^\infty(E_r\cap \partial H_e)}  \right). 
\end{equation}
Then, for every $k \in \N$, there exists $\theta(\alpha,k,d,\Lambda) \in (0,\frac12)$  and $C(\alpha,k,d,\Lambda)<\infty$ such that, for every $r\in \left[ X, R\right]$, 
\begin{multline}
\label{e.u:decayestimate.bndr}
\inf_{p \in \Ahom_k } \left\| u - p \right\|_{\underline{L}^2(E_{\theta r})}
 \leq
\frac14 \theta^{k+1 - \alpha/2} \inf_{p \in \Ahom_k} \left\| u - p \right\|_{\underline{L}^2(E_{r})} + C r^{k+\gamma} 
\left[ \nabla^k g \right]_{C^{0,\gamma}(B_r \cap  \partial H_e)}
\\ + C K r^{-\alpha} \left( \left\| u - g(0) \right\|_{\underline{L}^2(E_{r})} + r \left\| \nabla_{e^\perp} g \right\|_{L^\infty(E_r\cap \partial H_e)}\right)  .
\end{multline}
\end{lemma}

\begin{proof} Fix $r\in \left[ X, R \right]$ and let $w_r$ be as in the statement satisfying~\eqref{e.wharmapproxit.boundary3}. An application of Lemma~\ref{l.oddreflection.gen} yields
\begin{equation*} 
\inf_{p \in \Ahom_k } \left\| w_r - p \right\|_{\underline{L}^2(B_{t})} \leq C \left( \frac{t}{r}\right)^{k+1} \inf_{p \in \Ahom_k }  \left\| w_r - p \right\|_{\underline{L}^2(B_{t})} + C r^{k+\gamma} \left[ \nabla^k g \right]_{C^{0,\gamma}(B_r \cap  \partial H_e)}.
\end{equation*}
Choosing $t = \theta r$ with $\theta  = (C 2^{2+\frac d2})^{-2/\alpha}$ and using~\eqref{e.wharmapproxit.boundary3} together with the triangle inequality then yields the result.
\end{proof}

We define
\begin{equation*} \label{}
\A_k\left( H_e \right):= \left\{ u\in \A(H_e) \,:\, \limsup_{r\to \infty} r^{-(k+1)} \left\| u \right\|_{\underline{L}^2(B_r\cap H_e)} = 0, \ u \vert_{\partial H_e} = p \in \mathcal{P}_k(\Rd) \right\}.
\end{equation*}
The main result of this section is the following characterization of $\A_k\left( H_e \right)$, which is a version of Theorem~\ref{t.regularity} for half spaces. 

\index{$C^{k,1}$ estimate}
\begin{theorem}[Higher regularity theory in half spaces]
\label{t.Liouville.halfspace}
\index{Liouville theorem}
Fix $s \in (0,d)$ and $e\in\partial B_1$. There exist an exponent $\delta(s,d,\Lambda)\in \left( 0, \frac12 \right]$ and a random variable $\X_s$ satisfying the estimate
\begin{equation}
\label{e.X.liouville}
\X_s \leq \O_s\left(C(s,d,\Lambda)\right)
\end{equation}
such that the following statements hold, for every $k\in\N$:
\begin{enumerate}
\item[{$\mathrm{(i)}_k$}] There exists $C(k,d,\Lambda)<\infty$ such that, for every $u \in \A_k(H_e)$, there exists $p\in \overline{\A}_k(H_e)$ such that, for every $R\geq \X_s$,
\begin{equation} \label{e.liouvillec.halfspace}
\left\| u - p \right\|_{\underline{L}^2(B_R\cap H_e)} \leq C R^{-\delta} \left\| p \right\|_{\underline{L}^2(B_R\cap H_e)}.
\end{equation}

\item[{$\mathrm{(ii)}_k$}]For every $p\in \overline{\A}_k(H_e)$, there exists $u\in \A_k(H_e)$ satisfying~\eqref{e.liouvillec} for every $R\geq \X_s$. 

\item[{$\mathrm{(iii)}_k$}]
Fix $\gamma\in (0,1]$. There exists $C(k,\gamma,d,\Lambda)<\infty$ such that, for every $R\geq \X_s$ and $u\in \A(B_R\cap H_e)$ satisfying $u = g \in C^{k,\gamma}(B_R \cap \partial H_e)$ on $B_R\cap \partial H_e$, there exists $\phi \in \A_k(H_e)$ such that, for every $r \in \left[ \X_s,  R \right]$, we have the estimate
\begin{equation}
\label{e.intrinsicreg.halfspace}
\left\| u - \phi \right\|_{\underline{L}^2(B_r \cap H_e)} 
\leq 
C \left( \frac r R \right)^{k+1} 
\left\| u \right\|_{\underline{L}^2(B_R \cap H_e)}
+
CR^{k+\gamma} \left[ \nabla^{k} g \right]_{C^{0,\gamma}(B_R\cap \partial H_e)}
\end{equation}
\end{enumerate}
In particular, $\P$-almost surely, we have, for every $k\in\N$,
\begin{equation} 
\label{e.dimensionofAk.halfspace}
\dim(\A_k(H_e)) = \dim(\Ahom_k(H_e)) =  \binom{d+k-1}{k} + \binom{d+k-2}{k-1}.
\end{equation}
\end{theorem}

Having Lemmas~\ref{l.harmonicapproximation.boundary} and~\ref{l.u:decayestimate.bndr} in hand, the proof of Theorem~\ref{t.Liouville.halfspace} follows closely that of Theorem~\ref{t.regularity}, and we leave the details as an exercise. 

\begin{exercise}
\label{ex.liouville.bndr3}
Prove Theorem~\ref{t.Liouville.halfspace} using Lemmas~\ref{l.harmonicapproximation.boundary} and~\ref{l.u:decayestimate.bndr}.
\end{exercise}

\section{Optimality of the regularity theory}
\label{s.regularity.optimal}

In this section we will demonstrate that the estimate~\eqref{e.XLipschitz} for the random minimal scale~$\X_s$ in Theorem~\ref{t.Lipschitz} is optimal by giving an explicit example in $d=2$ for which 
\begin{equation} 
\label{e.no.cannot.have.d}
\P \left[ \X_s > t \right] \gtrsim \exp\left( -ct^2 \right). 
\end{equation}
It is straightforward to generalize the example below to dimensions $d>2$, to obtain~\eqref{e.no.cannot.have.d} with $t^d$ in place of $t^2$ on the right side, but for simplicity we consider only the case $d=2$ here. 

\smallskip

We work with a Bernoulli random checkerboard example. 
Define two matrices
\begin{equation*} \label{}
\a_1:= \frac12
\begin{pmatrix}
5 & 3 \\
3 & 5
\end{pmatrix}
\quad \mbox{and} \quad 
\a_2:= \frac12
\begin{pmatrix}
5 & -3 \\
-3 & 5
\end{pmatrix}.
\end{equation*}
These matrices are simultaneously diagonalizable, with unit eigenvectors
\begin{equation*} \label{}
\mathbf{v}_1 := \frac1{\sqrt{2}}
\begin{pmatrix}
1\\1
\end{pmatrix}
\quad \mbox{and} \quad 
\mathbf{v}_2 := \frac1{\sqrt{2}}
\begin{pmatrix}
-1 \\ 1
\end{pmatrix}.
\end{equation*}
Note that $\a_1v_1=4v_1$, $\a_1v_2=v_2$, $\a_2v_1=v_1$, $\a_2v_2=4v_2$. Moreover, since a reflection of the $x_2$-axis flips~$\mathbf{v}_1$ and~$\mathbf{v}_2$, it also maps~$\a_1$ to~$\a_2$. We have chosen~$1$ and~$4$ for the eigenvalues in order to make the computations as simple as possible (the square roots of~$\a_1$ and~$\a_2$ and their inverses have a pleasant form, see below). In particular, we could take the eigenvalues to be~$1$ and~$1+\delta$ for any~$\delta>0$ and the details below would be only slightly different. 

\smallskip

We let~$\P_0$ be the probability measure on~$\{ \a_1,\a_2\}$ which assigns probability~$\frac12$ to each matrix. We take~$\tilde{\P}$ to be~$\P_0^{\Z^2}$, which is the law of the i.i.d.~ensemble~$\{ \a(z) \}_{z\in\Z^2}$. We assign to each $\a(z)$ a coefficient field~$\a(x)$ defined for all~$x\in\R^2$ by extending~$\a(\cdot)$ to be constant on the cubes of the form~$z+[0,1)^2$, and the pushforward of this extension is denoted by~$\P$, which is the law of~$\a(x)$.

\smallskip

Notice that $\P$ satisfies our assumptions with $\Lambda=4$ but with range of dependence equal to $\sqrt{2}$ instead of~$1$. (We could modify the example to have unit range of dependence by performing a dilation, but this makes the notation messy so we will not do so.)

\smallskip

We define the ``bad matrix'' by
\begin{equation} 
\label{e.badmatrix}
\tilde{\a} = \left\{
\begin{aligned}
& \a_1 & \mbox{in} & \ \left\{ (x_1,x_2) \in \R^2\,:\, x_1x_2 > 0 \right\}, \\
& \a_2 & \mbox{in} & \ \left\{ (x_1,x_2) \in \R^2\,:\, x_1x_2 < 0 \right\}.
\end{aligned}
\right.
\end{equation}
In other words, $\tilde{\a}$, restricted to $Q_L$, is equal to~$\a_1$ in the first and third quadrant and equal to $\a_2$ in the second and fourth quadrants (see Figure~\ref{f.atilde}). Observe that, for each $L\in\N$, the probability that the coefficient field sampled by $\P$ is equal to $\tilde{\a}$ in the cube $Q_L:=[-L,L)^2$ is
\begin{equation} 
\label{e.hahagotcha}
\P \left[ \a = \tilde{\a} \ \mbox{in} \ Q_L \right]
=
\left(\frac12\right)^{(2L)^2} = \exp\left( -( \log 16) L^2 \right). 
\end{equation}
Let $\mathcal{B}_L$ denote the ``bad event'' $\mathcal{B}_L:= \left\{ \a\,:\, \a = \tilde{\a} \ \mbox{a.e. in} \ Q_L \right\}$. We will check in the following lemma that any minimal radius~$\X$ for which the statement of Theorem~\ref{t.Lipschitz} holds must be larger than $L$ on the event $\mathcal{B}_L$ holds, at least for large enough~$L$. 

\begin{figure}[tb]
\centering
\begin{center}
\begin{tikzpicture}[scale=10]
\draw[<->,black,very thin] (0,1/2) -- (1,1/2);
\draw[<->,black,very thin] (1/2,0) -- (1/2,1);
\def\N{12};
\def\S{1/ (24*\N) };
\foreach \i in {1,...,\N}
{
\foreach \j in {1,...,\N}
{
\def\x{ \i/(2*\N+1) };
\def\y{ \j/(2*\N+1) };
\draw[-] ( {\x + 4*\S}, {\y + 4*\S}) -- ({\x - 4*\S},{\y -4*\S});
\draw[-] ( {\x - \S}, {\y + \S}) -- ({\x + \S},{\y -\S});
}
}
\foreach \i in {1,...,\N}
{
\foreach \j in {1,...,\N}
{
\def\x{ (\i+\N )/(2*\N+1) };
\def\y{ \j/(2*\N+1) };
\draw[-] ( {\x + \S}, {\y + \S}) -- ({\x - \S},{\y -\S});
\draw[-] ( {\x - 4*\S}, {\y + 4*\S}) -- ({\x + 4*\S},{\y -4*\S});
}
}
\foreach \i in {1,...,\N}
{
\foreach \j in {1,...,\N}
{
\def\x{ (\i+\N )/(2*\N+1) };
\def\y{ (\j+\N) /(2*\N+1) };
\draw[-] ( {\x + 4*\S}, {\y + 4*\S}) -- ({\x - 4*\S},{\y -4*\S});
\draw[-] ( {\x - \S}, {\y + \S}) -- ({\x + \S},{\y -\S});
}
}
\foreach \i in {1,...,\N}
{
\foreach \j in {1,...,\N}
{
\def\x{ \i /(2*\N+1) };
\def\y{ (\j+\N) /(2*\N+1) };
\draw[-] ( {\x + \S}, {\y + \S}) -- ({\x - \S},{\y -\S});
\draw[-] ( {\x - 4*\S}, {\y + 4*\S}) -- ({\x + 4*\S},{\y -4*\S});
}
}
\end{tikzpicture}
\end{center}
\caption{The coefficient field $\tilde{\a}$ defined in~\eqref{e.badmatrix}. The matrix at each grid point is represented by two lines which indicate the direction of the eigenvectors of the matrix with lengths proportional to the corresponding eigenvalue. Like the coefficient field in Figure~\ref{f.meyerscounterexample}, the eigenvector with the largest eigenvalue generally points in the direction of the origin.
}
\label{f.atilde}
\end{figure}

\begin{lemma}
There exists a solution $u\in H^1_{\mathrm{loc}}(\R^2) \cap L^\infty_{\mathrm{loc}}(\R^2)$ of the equation
\begin{equation*} \label{}
-\nabla \cdot \left( \tilde{\a}\nabla u\right) = 0 \quad \mbox{in} \ \R^2
\end{equation*}
such that, with $\alpha:=  \pi /  2\arccos \left( -\frac35 \right) \approx 0.7094\ldots < 1$, we have 
\begin{equation*} \label{}
u(\lambda x) = \lambda^\alpha u(x)
\end{equation*}
\end{lemma}
\begin{proof}
Let $\mathcal{Q}\subseteq \R^2$ be the first quadrant, that is, 
\begin{equation*} \label{}
\mathcal{Q}:= \left\{ (x_1,x_2)\in\R^2\,:\, x,y>0 \right\}. 
\end{equation*}
By symmetry, it suffices to find a solution $u\in H^1_{\mathrm{loc}}(\mathcal{Q}) \cap L^\infty_{\mathrm{loc}}(\mathcal{Q})$ satisfying 
\begin{equation} 
\label{e.quadrant.pde}
\left\{
\begin{aligned}
& -\nabla \cdot \left( \a_1 \nabla u\right) = 0 
& \mbox{in} & \ \mathcal{Q}, \\
& u = 0 & \mbox{on} & \ \partial \mathcal{Q} \cap \{ x_2=0\}, \\
& e_1\cdot \a_1 \nabla u = 0 & \mbox{on} & \ \partial \mathcal{Q} \cap \{ x_1=0\},
\end{aligned}
\right. 
\end{equation}
and which is~$\alpha$-homogeneous. We then obtain the desired solution by extending~$u$ by reflecting~$u$ across the Neumann boundary~$\{x_1=0\}$ to get a solution in~$\{ x_2>0\}$, and then performing a negative reflection across~$\{x_2=0\}$ to obtain a solution in~$\R^2$.

\smallskip

In order to solve~\eqref{e.quadrant.pde}, we can perform an affine change of variables so that the equation becomes the Laplace equation. The change of variables will map $\mathcal{Q}$ to a sector~$\mathcal{S}$ and we will find a radial function in~$\mathcal{S}$ which solves the equation. Since the angle of $\mathcal{S}$ will turn out to be larger than $\frac{\pi}{2}$, this will allow us to select $\alpha < 1$. 

\smallskip

The change of variables is $x\mapsto \a_1^{1/2}x$.
Given a function
$u\in H^1_{\mathrm{loc}}(\mathcal{Q}) \cap L^\infty_{\mathrm{loc}}(\mathcal{Q})$, 
we define
\begin{equation*} \label{}
v(x) := u\left( \a_1^{1/2}x \right), \quad x\in \mathcal{S}:= \a_1^{-1/2} \mathcal{Q}. 
\end{equation*}
For the sake of computations, we note that 
\begin{equation*} \label{}
\a_1^{1/2} := \frac12
\begin{pmatrix}
3 & 1 \\
1 & 3
\end{pmatrix}
\quad \mbox{and} \quad
\a_1^{-1/2} := \frac14
\begin{pmatrix}
3 & -1 \\
-1 & 3
\end{pmatrix}.
\end{equation*}
It is easy to check that $u$ solves~\eqref{e.quadrant.pde} if and only if $v$ solves the problem 
\begin{equation} 
\label{e.sector.pde}
\left\{
\begin{aligned}
& -\Delta v = 0 
& \mbox{in} & \ \mathcal{S}, \\
& v = 0 & \mbox{on} & \ \partial \mathcal{S} \cap \a_1^{-1/2}\{ x_2=0\}, \\
& n \cdot \nabla u = 0 & \mbox{on} & \ \partial \mathcal{S} \cap \a_1^{-1/2} \{ x_1=0\}. 
\end{aligned}
\right. 
\end{equation}
In other words, $v$ should be a harmonic function in the sector~$\mathcal{S}$ satisfying a zero Dirichlet condition on one boundary line and a zero Neumann condition on the other (here the vector $n=\a_1^{1/2} e_1$ is the outer-pointing unit normal to $\partial \mathcal{S}$ on the second boundary line $\partial \mathcal{S} \cap \a_1^{-1/2} \{ x_1=0\}$). We denote the the angle of opening of $\mathcal{S}$ (the angle between the two boundary lines) by $\theta_*$, and notice by an easy computation that 
\begin{equation*} \label{}
\theta_* = \arccos\left( -\frac35 \right). 
\end{equation*}
Notice that $\theta_* \approx 2.2143{\hfil.\hfil.\hfil.}$, in particular, $\frac\pi2 < \theta_* < \pi$. 

\smallskip

To solve~\eqref{e.sector.pde}, by symmetry we can perform a rotation and reflection so that the boundary line with the Dirichlet condition lies on the positive $x_1$-axis, and the second boundary line for the Neumann condition lies in the second quadrant. We may then search for a solution taking the following the form, in polar coordinates:
\begin{equation*} \label{}
v(r,\theta) = r^\alpha \sin\left( \frac{\pi}{2\theta_*}  \theta \right) 
\end{equation*}
Since this function clearly satisfies the appropriate boundary conditions, we can just ask for which values of~$\alpha$ the function~$v$ is harmonic. By an easy computation we see that we should take $\alpha = \pi/2\theta_* < 1$. This completes the proof of the lemma. 
\end{proof}

Let us now suppose that we are given a large constant $\mathsf{C}>1$ and minimal scale~$\X$ which satisfies the result of Theorem~\ref{t.Lipschitz} with $C=\mathsf{C}$. This means that, for each $L\geq \frac12\X$ and every solution $w\in H^1(B_L)$ of 
\begin{equation*} \label{}
-\nabla \cdot \left( \a\nabla w \right) = 0 \quad \mbox{in} \ B_L, 
\end{equation*}
we have, for every $r\in \left[ \X, \frac 12L \right]$, 
\begin{equation*} \label{}
\left\| \nabla w  \right\|_{\underline{L}^2(B_r)} 
\leq \mathsf{C} \left\| \nabla w  \right\|_{\underline{L}^2(B_L)}. 
\end{equation*}
However, by homogeneity, the function~$u$ in the statement of the previous lemma satisfies
\begin{equation*} \label{}
\left\| \nabla u  \right\|_{\underline{L}^2(B_r)} 
= L^{1-\alpha} r^{\alpha-1} \left\| \nabla u  \right\|_{\underline{L}^2(B_L)}. 
\end{equation*}
Therefore if $\a\in \mathcal{B}_L$ (i.e. $\a = \tilde\a$ in $Q_L$) and $L > \frac12\X$, then taking $r=\X$ yields 
\begin{equation*} \label{}
\X \geq \mathsf{C}^{-\frac{1}{1-\alpha}} L. 
\end{equation*}
In view of~\eqref{e.hahagotcha}, we deduce that 
\begin{equation*} \label{}
\P \left[ 
\X \geq \left( 2 \vee \mathsf{C}^{-\frac{1}{1-\alpha}} \right) L
\right] 
\geq  \exp\left( -( \log 16) L^2 \right).
\end{equation*}
Therefore we conclude that, no matter how large the constant~$\mathsf{C}$ is chosen, any minimal scale~$\X$ satisfying the statement of Theorem~\ref{t.Lipschitz} will not be bounded by $\O_s(C)$ for any $s>d$ and, moreover, if we take $c$ small enough then 
\begin{equation*} \label{}
\E \left[ \exp\left( c\X^d \right) \right] = \infty.
\end{equation*}

\section*{Notes and references}

The regularity theory for homogenization was introduced in the celebrated work of Avellaneda and Lin~\cite{AL1,AL2,AL3} in the~1980s for periodic, smooth coefficients. They proved $C^{0,1}$ estimates, higher $C^{k,1}$-type estimates and Liouville theorems as well as Calder\'on-Zygmund estimates. Their original arguments are similar in spirit to the ones here, based on $\ahom$-harmonic approximation, although their softer approach using a compactness blow-down argument and qualitative homogenization does not extend beyond the non-compact (i.e.\ non-periodic) setting. 

\smallskip

The arguments in Section~\ref{s.lipschitz} are taken from~\cite{AS}, where the $C^{0,1}$-type estimates was first proved in the stochastic setting and the idea of using a quantitative excess decay or $C^{1,\gamma}$-type iteration was introduced. It was that paper which introduced the concept of a \emph{random scale} above which stronger regularity holds and gave a recipe for obtaining sharp stochastic integrability of this random scale. Following~\cite{AS}, the regularity theory was expanded and completed rather rapidly in a series of papers by several authors~\cite{ASh,GNO3,AM,FO,AKM1,AKM2}, where the higher regularity theory and related Liouville theorems were established. Prior to~\cite{AS}, some weaker $C^{0,1-\ep}$-type regularity results and qualitative Liouville theorems were obtained in~\cite{BDKY,MO} by rather different methods. 

\smallskip

The arguments in Section~\ref{s.higherreg} are taken mostly from~\cite[Section 3]{AKM2}, where the full statement of Theorem~\ref{t.regularity} first appeared. Example~\ref{ex.meyers} is due to Meyers~\cite[Section 5]{Me}. Boundary regularity results appeared previously in~\cite{ASh,FR}.



\chapter{Quantitative description of first-order correctors}
\label{c.A1}

The regularity theory, namely Theorem~\ref{t.regularity}$\mathrm{(iii)}_1$, tells us that an arbitrary solution of our equation can be approximated on scales smaller than the macroscopic scale, up to additive constants and at leading order, by an element of the~$d$-dimensional space~$\A_1/\R$. This result gives us a very fine quantitative understanding of the microscopic- and mesoscopic-scale behavior of solutions and can greatly reduce the complexity of many natural homogenization problems---provided that we have a good quantitative understanding of~$\A_1$ itself. 

\smallskip

The elements of $\A_1/\R$ are, modulo constants, the set of solutions in the whole space~$\Rd$ which exhibit linear growth at infinity. We know from Theorem~\ref{t.regularity} that each of these solutions is close to some affine function in the sense of~\eqref{e.liouvillec}: roughly, it diverges from an affine function like $O(|x|^{1-\delta})$. Indeed, every $u\in \A_1$ can be written for some $p\in \Ahom_1$ as
\begin{equation*} \label{}
u(x) = p(x) + \phi_{\nabla p} (x),
\end{equation*}
where $\{ \phi_\xi \,:\, \xi\in\Rd\}$ is the set of \emph{first-order correctors} (defined in Section~\ref{s.correctorsdefined} above). 

\smallskip

The goal in this chapter is to prove estimates on quantities such as the $L^2$-oscillation of the first-order correctors and estimates on the gradients, fluxes and energy densities in negative Sobolev spaces. Unlike the theory developed in previous chapters (see for instance Lemma~\ref{l.correctorminbounds} and Proposition~\ref{p.correctorsbasecase}), our intention here is to prove \emph{optimal} quantitative bounds on the first-order correctors, that is, we wish to discover the actual scaling of the errors.

\smallskip

We next present one of the main quantitative results for the first-order correctors, stated in terms of convolutions against the heat kernel. Most of this chapter is focused on its proof, which is finally completed in Section~\ref{s.correctors} and where we also present the version of the estimates in terms of~$W^{-\alpha,p}$ norms.
Here and throughout the rest of this chapter, we use the notation, for a given $\Psi\in L^1(\Rd)$,
\begin{equation} 
\label{e.def.int.against.psi}
\int_\Psi f  = \int_{\Psi} f(x) \, dx := \int_{\Rd} \Psi(x) f(x)\,dx. 
\end{equation}
Recall that the standard heat kernel is denoted by 
\begin{equation*} 
\Phi(t,x):= \left( 4\pi t\right)^{-\frac d2} \exp\left( -\frac{|x|^2}{4t} \right).
\end{equation*}
We also define, for each $z\in\Rd$ and $r>0$, 
\begin{equation*} 
\Phi_{z,r} (x):= \Phi(r^2,x-z) \quad \mbox{and} \quad \Phi_r := \Phi_{0,r}. 
\end{equation*}

\begin{theorem}
\label{t.correctors}
Fix $s<2$. There exists $C(s,d,\Lambda)<\infty$ such that, for every $r \geq 1$, $x\in\Rd$ and $e\in B_1$,
\begin{equation}
\label{e.gradient}
\left| \int_{\Phi_{x,r}} \nabla \phi_e \right| \leq \O_s\left( Cr^{-\frac d2} \right),
\end{equation}
\begin{equation}
\label{e.flux}
\left| \int_{\Phi_{x,r}} \a (e +  \nabla \phi_e) - \ahom e \right| \leq \O_s\left( Cr^{-\frac d2} \right),
\end{equation}
and
\begin{equation}
\label{e.energydensity}
\left| \int_{\Phi_{x,r}} \frac12 \left( e+ \nabla \phi_e\right) \cdot \a\left( e+ \nabla \phi_e\right)  - \frac12 e\cdot \ahom e \right| \\
 \leq \O_s\left( Cr^{-\frac d2} \right).
\end{equation}
Moreover, in dimensions $d>2$, the corrector~$\phi_e$ exists as a $\Zd$-stationary function which is identified uniquely by the choice $\E \left[ \fint_{\cu_0} \phi_e \right] = 0$ and there exist $\delta(d,\Lambda)>0$ and~$C(d,\Lambda)<\infty$ such that 
\begin{equation} 
\label{e.phioscd>2}
\left\| \phi_e \right\|_{\underline{L}^2(\cu_0)} 
\leq 
\O_{2+\delta}(C). 
\end{equation}
In~$d=2$, there exists $C(s,d,\Lambda)<\infty$ such that, for every $r\in [2,\infty)$, $R\in [r,\infty)$ and $x,y\in\Rd$, we have
\begin{equation}
\label{e.phioscd=2}
\left\|  \phi_e - \left( \phi_e\ast \Phi_r\right)(0) \right\|_{\underline{L}^2(B_r)}  \leq \O_{s} \left( C \log^{\frac12} r \right)
\end{equation}
and
\begin{equation}
\label{e.phid=2}
\left| \left( \phi_e\ast \Phi_r \right)(x) - \left( \phi_e\ast \Phi_R \right)(y) \right|
\leq
\O_s\left( C\log^{\frac12} \left(2+\frac{R+|x-y|}r\right) \right).
\end{equation}
\end{theorem}

If $f(x)$ is a stationary random field with a unit range of dependence and $|f|\leq 1$ almost surely, then $\int_{\Phi_{r}} f(x)\,dx$ should have roughly the same scaling as the average of $Cr^d$ many independent random variables---the fluctuations of which are of order $(Cr^{d})^{-\frac12} = Cr^{-\frac d2}$ with Gaussian tails. That is, the average of $Cr^{d}$ many bounded and independent random variables has fluctuations bounded by $\O_2\left( Cr^{-\frac d2}\right)$. Therefore, the conclusion of Theorem~\ref{t.correctors} can be interpreted as the statement that the spatial averages of the gradient, flux and energy density of the correctors behave in the same way as a stationary random field with finite range of dependence, or a sum of independent random variables. In particular, these estimates are optimal, both in the scaling $\simeq r^{-\frac d2}$ of the error as well as in stochastic integrability $\lesssim \O_{2-}$ (strictly speaking, they are only \emph{almost} optimal in stochastic integrability, since it is probably possible to prove the estimates with $\O_2$ replacing $\O_s$ for $s<2$). The estimates~\eqref{e.phioscd>2},~\eqref{e.phioscd=2} and~\eqref{e.phid=2} on the sublinear growth of the correctors are also optimal, in every dimensions. In particular, the logarithmic divergence in $d=2$ is intrinsic to the problem and not an artifact of our proof (see Theorem~\ref{t.gff} and Exercise~\ref{ex.gff.div.2d}). Note that the choice of the heat kernel in~\eqref{e.phioscd=2} and~\eqref{e.phid=2} is not important and can be replaced by any other $H^2$ function: see Remark~\ref{r.Phidoesntmatter}. 

\index{corrector!first-order~$\phi_e$|)}

\smallskip

The proof of Theorem~\ref{t.correctors} is accomplish by a \emph{bootstrap} argument, 
using the suboptimal bounds proved in previous chapters as a starting point. We will show, as we pass to larger and larger length scales, that we can improve the exponent from the tiny~$\alpha>0$ or~$\delta>0$ found in previous chapters to the optimal exponent. This bootstrap argument will use more sophisticated versions of ideas already present in Chapter~\ref{c.two}. We will identify a variant of the energy quantity~$J$, which we denote by~$J_1$, and show that it is becoming more \emph{additive} on larger and larger scales. The regularity theory will allow us to \emph{localize} this quantity. The combination of additivity and locality will then allow us to apply the independence assumption in a very simple and straightforward way to control the fluctuations of~$J_1$ by essentially reducing it to a sum of independent random variables (or more precisely a convolution of a random field with finite range of dependence). 

\smallskip

This bootstrap argument, which is the focus of most of the chapter, is outlined in Section~\ref{s.boots}. We begin in the next section by introducing the energy quantity~$J_1$. 

\section{The energy quantity \texorpdfstring{$J_1$}{J1} and its basic properties}

In order to study the first-order correctors and prove Theorem~\ref{t.correctors}, we introduce a variant of the quantity~$J$ studied in Chapter~\ref{c.two} which is adapted to~$\A_1$ and which we denote by~$J_1$. There is another important difference between $J_1$ and the quantity~$J$ considered in previous chapters: since our goal is to prove very precise estimates, the analysis requires us to relax the concept of \emph{domain} on which our quantity is defined to allow for smoother objects. The integral $\fint_U f(x) \,dx$ can of course be written as $\int_{\Rd} |U|^{-1} \indc_{U}(x) f(x) \,dx$; we will replace $ |U|^{-1} \indc_{U}$ in our definition of $J$ by allowing for an arbitrary nonnegative $L^1$ function $\Psi$ with unit mass. Usually we will take $\Psi=\Phi_r$. Recall that we use the notation in \eqref{e.def.int.against.psi} to denote integrals against $\Psi$. 
We also write $\| w \|_{L^p(\Psi)} := \left( \int_\Psi |w|^p \right)^{\frac1p}$ for $p\in [1,\infty)$, and
\begin{equation*} \label{}
\A_k(\Psi):=\left\{ u \in \A_k \,:\, \left\| \nabla u \right\|_{L^2(\Psi)}\leq 1 \right\}.
\end{equation*}
We define, for every nonnegative $\Psi\in L^1(\Rd)$ with $\int_{\Rd} \Psi=1$ and $p,q\in\Rd$,
\index{subadditive quantity!$J_1$}
\begin{equation} 
\label{e.def.J1}
J_1 (\Psi,p,q) := 
\sup_{w\in \A_1} 
\int_\Psi 
\left( -\frac12\nabla w \cdot \a\nabla w - p\cdot \a\nabla w + q\cdot \nabla w \right).
\end{equation}
As we will see, the quantity $J_1(\Psi,p,q)$ shares many properties possessed by its predecessor, the quantity $J(U,p,q)$ analyzed in Chapter~\ref{c.two}. One notable exception is that there is no version of the splitting formula~\eqref{e.Jsplitting} for $J_1$. Indeed, the connection between~\eqref{e.Jsplitting} and~\eqref{e.formulafornuJpq} relies on an integration by parts to determine the spatial average of the gradient of $v(\cdot,U,p,0)$, and this computation is destroyed by the presence of the weight $\Psi$. Thus there is no version of Lemma~\ref{l.minimalset} and no way to control the fluctuations of $J_1$ from its behavior for $(p,q)$ with $q=\ahom p$. This makes our analysis more complicated, compared to Chapter~\ref{c.two}, by not allowing us to make much use of subadditivity, forcing us to considerably alter our strategy. The advantage we have here, compared with Chapter~\ref{c.two}, is that the regularity theory is at our disposal as well as the suboptimal bounds summarized in Proposition~\ref{p.basecase} below, which are used as the base case of our bootstrap argument. Both of these crucial ingredients are actually consequences of the analysis in Chapter~\ref{c.two}.

\smallskip

Related to the loss of the splitting formula is the (relatively minor, technical) problem that $J_1(\Psi,p,q)$ is \emph{not necessarily} uniformly convex in $p$ and/or $q$. This leads us to introduce the notion of a nondegenerate mask for $J_1$. 
\begin{definition}[{Nondegeneracy at $\Psi$}]
We say that \emph{$J_1$ is nondegenerate at $\Psi$} if, for every $p,q\in\Rd$, 
\begin{equation}
\label{e.nondegeneratePsi}
J(\Psi,p,0) \geq \frac{1} 4 p\cdot\ahom p
\quad \mbox{and} \quad 
J(\Psi,0,q) \geq \frac{1} 4 q\cdot\ahom^{-1}q.
\end{equation}
\end{definition}
Since $J_1$ is quadratic, nondegeneracy at $\Psi$ implies that it is uniformly convex in $p$ and $q$ separately, as we will see in Lemma~\ref{l.basicJ1} below. 

\smallskip

The maximizer of $J_1(\Psi,p,q)$ is denoted by
\begin{equation*} \label{}
v(\cdot,\Psi,p,q) := \mbox{the element of $\A_1$ achieving the supremum in the definition of $J_1$.}
\end{equation*}
Note that such a maximizer exists since $\nabla \A_1$ is a $d$--dimensional vector subspace of $L^2(\Rd,\Psi(x)\,dx)$ by Theorem~\ref{t.regularity}. Furthermore, it is unique---up to additive constants---by the strict concavity of the integral functional in the definition of~$J_1$.

\smallskip

We next use Proposition~\ref{p.correctorsbasecase} to get quantitative convergence of $J_1$ on scales larger than a random minimal scale. 

\begin{proposition}
\label{p.basecase}
Fix $s \in (0,d)$. Then there exist an exponent $\beta(s,d,\Lambda)>0$, a random variable $\X_s$ and a constant $C(s,d,\Lambda)<\infty$ such that 
\begin{equation*} \label{}
\X_s = \O_s(C)
\end{equation*}
and, for every $r\geq \X_s$ and $p,q\in \R^d$, 
\begin{equation} 
\label{e.J1basecase}
\left| J_1(\Phi_{r},p,q) - \left( \frac 12 p\cdot \ahom p + \frac12 q\cdot \ahom^{-1}q - p\cdot q\right) \right| \leq Cr^{-\beta}\left( |p|^2+|q|^2 \right).
\end{equation}
In particular, $J_1$ is nondegenerate at $\Phi_r$ for every $r\geq \X_s$. 
\end{proposition}
\begin{proof}
The idea is to evaluate $J_1$ in terms of the integral of the energy densities, fluxes and gradients of the correctors $\{\phi_\xi\}_{\xi\in\Rd}$ which reduces, up to a suitable error, to a deterministic computation by Proposition~\ref{p.correctorsbasecase}. 

\smallskip

We take $\X_s$ to be as in Proposition~\ref{p.correctorsbasecase} and fix $r\geq \X_s$ and $p,q\in\Rd$. Then Proposition~\ref{p.correctorsbasecase} yields, for every $\xi \in \Rd$, 
\begin{multline} 
\label{e.determineJ1}
\bigg| \int_{\Phi_r}  
\left( -\frac12( \xi + \nabla \phi_\xi) \cdot \a ( \xi + \nabla \phi_\xi) - p\cdot \a  ( \xi + \nabla \phi_\xi) + q\cdot ( \xi + \nabla \phi_\xi) \right) 
\\
- \left( -\frac12 \xi\cdot \ahom \xi - p \cdot \ahom \xi + q\cdot \xi \right)  \bigg|
\leq C\left( |\xi|^2 + |\xi|(|p|+|q|) \right) r^{-\beta}. 
\end{multline}
By Lemma~\ref{l.identifyA1}, we may write $J_1(\Phi_r,p,q)$ as  
\begin{multline} \label{e.J1xidef}
J_1(\Phi_r,p,q) \\
= \sup_{\xi \in\Rd} \int_{\Phi_r}  
\left( -\frac12( \xi + \nabla \phi_\xi) \cdot \a ( \xi + \nabla \phi_\xi) - p\cdot \a ( \xi + \nabla \phi_\xi) + q\cdot  ( \xi + \nabla \phi_\xi) \right).
\end{multline}
We deduce that
\begin{equation*} \label{}
 \sup_{\xi \in\Rd} \left( -\frac12 \xi\cdot \ahom \xi - p \cdot \ahom \xi + q\cdot \xi -Cr^{-\beta}\left( |\xi|^2 + |\xi|(|p|+|q|) \right) \right)
 \leq 
J_1(\Phi_r,p,q).
\end{equation*}
Using the identity
\begin{equation*} \label{}
\sup_{\xi \in\Rd} \left( -\frac12 \xi\cdot \ahom \xi - p \cdot \ahom \xi + q\cdot \xi \right) =  \frac 12 p\cdot \ahom p + \frac12 q\cdot \ahom^{-1}q - p\cdot q 
\end{equation*}
and the fact that the maximizing $\xi$ is $\xi = -p+\ahom^{-1}q$ which satisfies $|\xi| \leq C(|p|+|q|)$, we obtain that
\begin{equation*} \label{}
\frac 12 p\cdot \ahom p + \frac12 q\cdot \ahom^{-1}q - p\cdot q  \leq J_1(\Phi_r,p,q) + Cr^{-\beta} \left( |p|^2+|q|^2 \right).
\end{equation*}
This is one half of~\eqref{e.J1basecase}. To ensure that $J_1$ is nondegenerate at $\Phi_r$ for $r\geq \X_s$, we note that 
\begin{equation*} \label{}
 \frac 12 p\cdot \ahom p + \frac12 q\cdot \ahom^{-1}q - p\cdot q = \frac12 \left( p-\ahom^{-1}q\right)\cdot \ahom \left( p-\ahom^{-1}q\right)
\end{equation*}
and replace $\X_s$ with $\X_s\vee C$ for a suitable constant~$C$ so that $Cr^{-\beta}\leq C\X_s^{-\beta} \leq \frac{1}4$. 

\smallskip

To get the other half of~\eqref{e.J1basecase}, we argue similarly, first using the nondegeneracy of~$J_1$ at~$\Phi_r$ to get that the parameter $\xi$ achieving the supremum in~\eqref{e.J1xidef} is bounded by $C(|p|+|q|)$. Then we apply the estimate~\eqref{e.determineJ1} to obtain
\begin{equation*} \label{}
J_1(\Phi_r,p,q) \leq \frac 12 p\cdot \ahom p + \frac12 q\cdot \ahom^{-1}q - p\cdot q + Cr^{-\beta} \left( |p|^2+|q|^2 \right).
\end{equation*}
This completes the proof. 
\end{proof}

We next summarize some basic properties of $J_1$ in the following lemma, which can be compared with Lemma~\ref{l.basicJ}. In what follows, we use~$D_p$ and~$D_q$, respectively, to denote the derivatives of $J_1$ in the variables $p$ and~$q$ (the symbol~$\nabla$ is reserved for spatial variables).

\begin{lemma}[Basic properties of $J_1$]
\label{l.basicJ1}
The quantity $J_1(\Psi,p,q)$ and its maximizer $v(\cdot,\Psi,p,q)$ satisfy the following properties:

\begin{itemize}

\item \emph{First variation for $J_1$.} For $p,q\in\Rd$, the function $v(\cdot,\Psi,p,q)$ is characterized as the unique element of $\A_1$ which satisfies 
\begin{equation}
\label{e.J1firstvar} 
\int_\Psi \nabla w \cdot \a\nabla v(\cdot,\Psi,p,q) = \int_\Psi \left( - p\cdot \a\nabla w + q\cdot \nabla w \right),
\quad 
\forall\, w \in \A_1. 
\end{equation}
In particular, 
\begin{equation} 
\label{e.Jexpressbyv}
J_1(\Psi,p,q) = \int_{\Psi} \frac12 \nabla v(\cdot,\Psi,p,q) \cdot  \a  \nabla v(\cdot,\Psi,p,q).
\end{equation}

\item \emph{Quadraticity, linearity and boundedness.} The mapping $(p,q) \mapsto \nabla v(\cdot,\Psi,p,q)$ is linear, $(p,q) \mapsto J_1(\Psi,p,q)$ is quadratic, and there exists $C(\Lambda)<\infty$ such that 
\begin{equation}
\label{e.Jvbounded}
\left\| \nabla v(\cdot,\Psi,p,q) \right\|_{L^2(\Psi)}^2 \leq J_1(\Psi,p,q) \leq C \left( |p|^2 + |q|^2 \right). 
\end{equation}

\item \emph{Polarization.} For each $p,q,p',q' \in \R^d$, 
\begin{multline} 
\label{e.J polarization}
J_1(\Psi,p+p',q+q') \\ = J_1(\Psi,p,q) + J_1(\Psi,p',q') + \int_\Psi \left( - \a p' + q' \right) \cdot \nabla v(\cdot,\Psi,p,q) .
\end{multline}

\item \emph{Uniformly convexity for nondegenerate $\Psi$.} Suppose that $J_1$ is nondegenerate at $\Psi$. Then there exists $C(\Lambda)<\infty$ such that, for every  $p_1,p_2,q\in\Rd$,
\begin{multline} 
\label{e.J1unifconvp}
\frac1C \left| p_1-p_2\right|^2
\leq
\frac12 J_1\left( \Psi,p_1,q \right) + \frac12 J_1\left( \Psi,p_2,q \right)  - J_1\left( \Psi,\frac12p_1+\frac12 p_2 ,q\right)
\end{multline}
and, for every $q_1,q_2,p\in\Rd$, 
\begin{multline} 
\label{e.J1unifconvq}
\frac1C \left| q_1 - q_2\right|^2
\leq
\frac12 J_1\left( \Psi,p,q_1 \right) + \frac12 J_1\left( \Psi,p,q_2 \right)  - J_1\left( \Psi,p,\frac12q_1+\frac12 q_2 \right).
\end{multline}

\item \emph{Second variation and quadratic response.} For every $p,q\in\Rd$ and $w\in \A_1$, 
\begin{multline} 
\label{e.J1quadresponse}
 \int_\Psi \frac12 \left( \nabla w - \nabla v(\cdot,\Psi,p,q) \right)\cdot \a \left( \nabla w - \nabla v(\cdot,\Psi,p,q) \right)
\\
  =
J_1(\Psi,p,q) - \int_\Psi \left( -\frac12\nabla w \cdot \a\nabla w  -p\cdot\a \nabla w + q\cdot \nabla w  \right).
\end{multline}

\item \emph{Formulas for derivatives of $J_1$.} For every $p,q\in\Rd$,
\begin{equation} 
\label{e.J1derp}
 D_p  J_1(\Psi,p,q)  = -\int_\Psi \a \nabla v(\cdot,\Psi,p,q)
\end{equation}
and
\begin{equation} 
\label{e.J1derq}
 D_q  J_1(\Psi,p,q)  =  \int_\Psi  \nabla v(\cdot,\Psi,p,q).
\end{equation}
\end{itemize}
\end{lemma}
\begin{proof}

Let us begin with the first and second variations for $J_1$. Let $p,q\in\Rd$ and $u,v \in \A_1$. Test $J_1(\Psi,p,q)$ with $w = v(\cdot,\Psi,p,q) + t u$, $t\in \R \setminus \{0\}$, to obtain
\begin{align*} 
0 & \leq  J_1(\Psi,p,q) - \int_\Psi \left(-\frac12\nabla w \cdot \a\nabla w - p\cdot \a\nabla w + q\cdot \nabla w \right)
\\ \notag & =  - t \int_\Psi \left( -  \a  \nabla v(\cdot,\Psi,p,q)   - \a p  + q  \right) \cdot  \nabla u
+ \frac{t^2}{2} \int_\Psi \nabla u \cdot \a\nabla u .
\end{align*}
Dividing by $t$ and then sending $t \to 0$ gives~\eqref{e.J1firstvar}. The second variation formula~\eqref{e.J1quadresponse} is a consequence of~\eqref{e.J1firstvar} and the previous display with the choice $t=1$. The first variation formula immediately implies
\begin{equation*} 
(p,q) \mapsto \nabla v(\cdot,\Psi,p,q) \quad \mbox{is a linear mapping from \quad $\Rd \times \Rd \to \nabla \A_1$,}
\end{equation*}
as well as the formula~\eqref{e.Jexpressbyv}. 
By Young's inequality, we deduce  
\begin{equation} \label{e.J rough bound}
J_1(\Psi,p,q) \leq C \left\| \Psi \right\|_{L^1(\R^d)}\left(|p| + |q| \right)^2 = C\left(|p| + |q| \right)^2
\end{equation}
and hence~\eqref{e.Jvbounded}.
By~\eqref{e.J1firstvar} and~\eqref{e.Jexpressbyv}, we also obtain
\begin{align} 
\label{e.J polar1}
\lefteqn{J_1(\Psi,p+p',q+q')} \quad &
\\  \notag & = \int_\Psi \left(\frac12\nabla v(\cdot,\Psi,p+p',q+q') \cdot \a\nabla v(\cdot,\Psi,p+p',q+q') \right)
\\ \notag &  =  \int_\Psi \left(\frac12\nabla v(\cdot,\Psi,p,q) \cdot \a\nabla v(\cdot,\Psi,p,q) + \frac12\nabla v(\cdot,\Psi,p',q') \cdot \a\nabla v(\cdot,\Psi,p',q') \right) 
\\ \notag &  \qquad +  \int_\Psi \nabla v(\cdot,\Psi,p,q) \cdot \a\nabla v(\cdot,\Psi,p',q')
\\ \notag &  = J_1(\Psi,p,q) + J_1(\Psi,p',q') + \int_\Psi \left( - \a p' + q' \right) \cdot \nabla v(\cdot,\Psi,p,q).
\end{align}
This is~\eqref{e.J polarization}.
Together with~\eqref{e.J rough bound}, this yields both~\eqref{e.J1derp} and~\eqref{e.J1derq}. Indeed,  for $h \in \R^d$  and $t\in \R \setminus \{0\}$, the formulas for gradients follow by taking $p' = t h$ and $q'=0$ or $q=t h$ and $p'=0$ in the above formula, and then dividing by $t$ and finally sending $t\to 0$, noticing that $|J(\Psi,p',q')| \leq C|h|^2 t^2$ in both cases by~\eqref{e.J rough bound}. 

\smallskip

We are left to prove both~\eqref{e.J1unifconvp} and~\eqref{e.J1unifconvq}. But they are now simple consequences of~\eqref{e.J polar1} and the bounds~\eqref{e.nondegeneratePsi} and~\eqref{e.J rough bound}. Indeed, taking $p = p_1$ and $p' = - \frac12 (p_1 -p_2)$, and $p = p_2$ and $p' = \frac12 (p_1 -p_2)$, together with choices $q'=0$ and $q$, we easily deduce~\eqref{e.J1unifconvp} using linearity,~\eqref{e.J polar1},~\eqref{e.nondegeneratePsi} and~\eqref{e.J rough bound}. The estimate~\eqref{e.J1unifconvq} is showed analogously. The proof is complete. 
\end{proof}

\begin{definition}[{Minimal scale~$\Y_s$}]
\label{d.Ys}
For the rest of this chapter, we let $\Y_s$ denote, for each $s\in (0,d)$, the maximum of the random variables in Theorem~\ref{t.regularity} and Propositions~\ref{p.correctorsbasecase} and~\ref{p.basecase}. In particular, we have that 
\begin{equation*} \label{}
\Y_s = \O_s\left(C(s,d,\Lambda)\right)
\end{equation*}
and that $r\geq \Y_s$ implies each of the conclusions in these results. We also let $\Y_s(x)$ denote, for each $x\in\Rd$, the random variable $\a\mapsto \Y_s(T_x\a)$ where $T_x$ is the translation operator on $\Omega$ defined in~\eqref{e.Ty}. It follows that $\Y_s(x)$ satisfies the same estimate as~$\Y_s$ and that similar statements hold for $r\geq \Y_s(x)$, appropriately translated by~$x$. 
\end{definition}

Since it will be used many times below, we record here the observation that, for every $x\in\Rd$, $s\in (0,d)$ and $t>0$, we have the estimate
\begin{equation} 
\label{e.indcYsbound}
\indc_{\{ \Y_s(x) \geq r \}} \leq \left( \frac{\Y_s}{r} \right)^{t} \leq \O_{s/t} \left( Cr^{-t} \right).
\end{equation}
We also require a useful consequence of~\eqref{e.correctorminbounds} and~\eqref{e.gradientbasecase}, which is that 
\begin{equation} 
\label{e.correctorsballz}
r \geq \Y_t(x) \implies
\left\{ \overline{\phi}_\xi \,:\, \left|\xi\right|\leq c \right\} \subseteq \nabla \A_1(\Phi_{x,r})  \subseteq \left\{ \overline{\phi}_\xi \,:\, \left|\xi\right|\leq C \right\},
\end{equation}
where here and throughout the rest of the section we denote
\begin{equation}
\label{e.correctorbarzz}
\overline{\phi}_e:= \ell_e + \phi_e.
\end{equation}
Also, by Proposition~\ref{p.basecase} and~\eqref{e.Jexpressbyv},
\begin{equation} 
\label{e.vqnondegen}
r\geq \Y_t(x) \implies
c|q| \leq \left\| \nabla v(\cdot,\Phi_{x,r},0,q) \right\|_{L^2(\Phi_{x,r})} \leq C|q|. 
\end{equation}
It will also be useful to have a deterministic scale above which the quadratic function $(p,q)\mapsto \E\left[ J(\Psi_r,p,q) \right]$ behaves nicely. By Proposition~\ref{p.basecase}, there exists $r_0(d,\Lambda) \in [1,\infty)$ such that, for every $x\in\Rd$, 
\begin{equation} 
\label{e.def.r0}
r\geq r_0
 \implies 
 \E \left[ J_1(\Phi_{x,r},p,0) \right] \geq \frac14p\cdot \ahom p 
\quad  \mbox{and} \quad
 \E \left[ J_1(\Phi_{x,r},0,q) \right] \geq \frac14q\cdot\ahom^{-1} q. 
\end{equation}

\section{The additive structure of \texorpdfstring{$J_1$}{J1}}
\label{s.boots}

In this section, we state the main quantitative estimates on~$J_1$ and give an overview of the structure of the proofs.
\index{subadditive quantity!convergence of}

\begin{theorem}[{The additive structure of $J_1$}]
\label{t.additivity}
For every $s < 2$, there exists a constant $C(s,d,\Lambda) < \infty$ such that the following statements hold.
\begin{enumerate}

\item[{(i)}] \emph{Additivity}. For every $R > r \ge 1$, $z\in\Rd$ and $p, q \in B_1$,
\begin{equation*}
\Ll| J_1(\Phi_{z,R},p,q) 
- \int_{\Phi_{z,\sqrt{R^2 - r^2}}} J_1(\Phi_{x,r},p,q) \, dx \Rr| \le \O_{s/2}\Ll(C  r^{-d}\Rr).
\end{equation*}

\smallskip

\item[(ii)] \emph{Control of the expectation}. For every $r \ge 1$, $z\in\Rd$ and $p,q \in B_1$,
\begin{equation*}
\left|
\E \left[J_1(\Phi_{z,r},p,q)\right] 
- \left( \frac12 p\cdot \ahom p + \frac12 q\cdot \ahom^{-1} q - p\cdot q \right) 
\right|
\leq C r^{-d}.
\end{equation*}

\smallskip

\item[(iii)] \emph{CLT scaling of the fluctuations}. For every $r \ge 1$ and $p,q \in B_1$,
\begin{equation*}
\big|J_1 (\Phi_{z,r},p,q)  - \E\left[ J_1(\Phi_{z,r},p,q)\right]\big| \le \O_{s}\left(C r^{-\frac d 2}\right).
\end{equation*}

\smallskip

\item[(iv)] \emph{Localization}. For every $\delta, \eps > 0$, there exist $C(\delta,\eps,s,d,\Lambda)<\infty$ and, for every $r\geq 1$, $z\in\Rd$ and $p,q\in B_1$, an $\F(B_{r^{1+\delta}}(z))$-measurable random variable $J_1^{(\delta)}(z,r,p,q)$ such that, for $\gamma := \frac d{s}\wedge \left( \frac d2(1+\delta)+\delta\right) -\ep$,
\begin{equation*} 
\left| J_1(\Phi_{z,r},p,q) - J_1^{(\delta)}(z,r,p,q)\right| \leq \O_{s}\left( Cr^{-\gamma}\right).
\end{equation*}
\end{enumerate}
\end{theorem}

We will see in Section~\ref{s.correctors} that Theorem~\ref{t.correctors} is a corollary of Theorem~\ref{t.additivity}. Therefore most of this chapter is focused on the proof of the latter, which is via a bootstrap on the exponent $\alpha$ for which the following statements hold. 

\begin{definition}[$\Add(s,\al)$]  
For each $s, \al \in (0,\infty)$, we let $\Add(s,\al)$ denote the statement that there exists a constant $C(s,\al,d,\Lambda) < \infty$ such that, for every $z \in \Rd$, $1 \le r \le R$ and $p,q \in B_1$,
\begin{equation*}  
J_1(\Phi_{z,R},p,q) = \int_{\Phi_{z,\sqrt{R^2 - r^2}}} J_1(\Phi_{x,r},p,q)\,dx + \O_s \Ll( C r^{-\al} \Rr) .
\end{equation*}
\end{definition}
\begin{definition}[$\Fluc(s,\al)$]  
For each $s, \al \in (0,\infty)$, we let $\Fluc(s,\al)$ denote the statement that there exists a constant $C(s,\al,d,\Lambda) < \infty$ such that, for every $z \in \Rd$, $r \ge 1$ and $p,q \in B_1$,
\begin{equation*}  
J_1(\Phi_{z,r},p,q) = \E \Ll[ J_1(\Phi_{z,r},p,q) \Rr] + \O_s \Ll( C r^{-\al} \Rr) .
\end{equation*}
\end{definition}
\begin{definition}[$\Loc(s,\delta,\alpha)$]  
For each $s, \delta, \alpha \in (0,\infty)$, we let $\Loc(s,\delta,\alpha)$ denote the statement that there exists a constant $C(s,\al,\de,d,\Lambda) < \infty$ and, for every $z \in \Rd$, $r \ge 1$ and $p,q \in B_1$, an $\F(B_{r^{1+\de}}(z))$-measurable random variable $\Jd(z,r,p,q)$ such that
\begin{equation*}  
J_1(\Phi_{z,r},p,q) = \Jd(z,r,p,q)  + \O_s \Ll( C r^{-\al} \Rr) .
\end{equation*}
\end{definition}

The bootstrap argument is split into three main subclaims: (i) an  \emph{improvement of additivity}, (ii) \emph{localization estimates} and (iii) \emph{fluctuation estimates}. 
These are proved in the following three sections, and their arguments are mostly independent of each other. In the rest of this section, we give their precise statements and show that Theorem~\ref{t.additivity} is a consequence of them.

\begin{proposition}[{Improvement of additivity}]
\label{p.improveadditivity}
For every $\ep>0$, there exists an exponent~$\eta(\ep,s,d,\Lambda)>0$ such that, for every 
$s\in (0,2]$ and $\alpha \in \left( 0, \left( \frac ds -\ep\right) \wedge \frac d2 \right]$,
\begin{equation} 
\label{e.add.improvement}
\Fluc(s,\alpha)   \implies   \Add\left( s,\alpha + \eta\right). 
\end{equation}
Moreover, for every $s\in (0,2]$ and $\alpha \in \left( 0,\frac ds \right) \cap \left(0,\frac d2\right]$,
\begin{equation}
\label{e.add.doubling}
\Fluc(s,\alpha)   \implies    \Add\left( \tfrac s2,2\alpha \right). 
\end{equation}
\end{proposition}

\begin{proposition}[{Localization estimates}]
\label{p.improvelocalization}
For every $s\in (0,\infty)$, $\alpha \in \left(0,\tfrac ds \right)$ and $\delta ,\eta>0$, 
\begin{equation*} \label{}
\Fluc(s,\alpha) 
\implies
\Loc\left(s,\delta,\left( \alpha(1+\delta)+\delta\right)\wedge \tfrac ds-\eta\right).
\end{equation*}
\end{proposition}

\begin{proposition}[{Fluctuation estimates}]
\label{p.improvefluc} For every $s \in (1,2]$, $\be \in \Ll(0,\frac d 2\Rr]$ and $\al,\de > 0$ satisfying $\al > \be(1+\de)$, we have
\begin{equation*}  
\Add(s,\al) \; \; \text{and} \; \;   \Loc(s,\delta,\alpha)  \implies  \Fluc(s,\be).
\end{equation*}
\end{proposition}

We conclude this section by demonstrating that Theorem~\ref{t.additivity} follows from the previous three propositions.

\begin{proof}[{Proof of Theorem~\ref{t.additivity}}]

\emph{Step 1.} We claim that, for each $\ep>0$ and $\alpha \in \left(0,\frac d2-\ep\right]$, there exists $\tau(\ep,d,\Lambda)>0$ such that 
\begin{equation*} \label{}
\Fluc(2,\alpha) \implies \Fluc(2,\alpha+\tau).
\end{equation*}
Let $\eta = \eta(\ep,d,\Lambda)>0$ be as in Proposition~\ref{p.improveadditivity} with $s=2$. By shrinking $\eta$, we may also assume that $\eta \leq \frac12 \ep$. Next, select $\delta(\ep,d,\Lambda)>0$ by $\delta:= \eta / (1+\frac d2)$ and observe that 
\begin{equation*} \label{}
\alpha (1+\delta) + \delta \leq \alpha + \eta \leq \frac d2. 
\end{equation*}
According to Propositions~\ref{p.improveadditivity} and~\ref{p.improvelocalization}, we have that 
\begin{equation*} \label{}
\Fluc(2,\alpha) 
\implies 
\Add\left(2,\alpha (1+\delta) + \tfrac12 \delta\right) 
\ \mbox{and} \
\Loc\left(2,\alpha(1+\delta)+\tfrac12 \delta \right).
\end{equation*}
Taking $\tau(\ep,d,\Lambda)>0$ to be defined by $\tau:= \delta / (2+2\delta)$ and applying Proposition~\ref{p.improvefluc}, we see that 
\begin{equation*} \label{}
\Add\left(2,\alpha (1+\delta) + \tfrac12 \delta\right) 
\ \mbox{and} \
\Loc\left(2,\alpha(1+\delta)+\tfrac12 \delta \right)
\implies
\Fluc(2,\alpha+\tau). 
\end{equation*}
This completes the proof of the claim. 

\smallskip

\emph{Step 2.} We show that 
\begin{equation} 
\label{e.flucbelowd2}
\forall \alpha < \tfrac d2, \quad \Fluc(2,\alpha) \  \mbox{holds.}
\end{equation}
According to Proposition~\ref{p.basecase} and Lemma~\ref{l.change-s}, there exists $\beta(d,\Lambda)>0$ such that $\Fluc(t,\beta/t)$ holds for every $t\geq 1$. Thus the result of Step~1 yields~\eqref{e.flucbelowd2}. 

\smallskip

\emph{Step 3.} We show that
\begin{equation} \label{}
\forall s\in(1,2), \quad \Fluc\left(s,\tfrac d2\right) \  \mbox{holds.}
\end{equation}
Fix $s\in (1,2)$ and let $\eta(s,d,\Lambda)>0$ be as in Proposition~\ref{p.improveadditivity} for $\ep:= \frac ds - \frac d2$. Next, select $\delta(s,d,\Lambda)>0$ by $\delta:= \eta / (d+1)$ so that  $\frac d2 (1+\delta) + \frac 12\delta = \frac d2 + \frac12 \eta$. Using~\eqref{e.flucbelowd2} with $\alpha$ sufficiently close to $\frac d2$, we may apply Propositions~\ref{p.improveadditivity} and~\ref{p.improvelocalization} to obtain 
\begin{equation*} \label{}
\Add\left(s,\tfrac d2 (1+\delta) + \tfrac12 \delta\right) 
\ \mbox{and} \
\Loc\left(s,\delta, \tfrac d2 (1+\delta)+\tfrac12 \delta \right) 
\quad \mbox{hold.}
\end{equation*}
Now we apply Proposition~\ref{p.improvefluc} once more to obtain that~$\Fluc(s,\tfrac d2)$ holds. 

\smallskip

\emph{Step 4.} The conclusion. Now that we have proved~$\Fluc(s,\tfrac d2)$ for every $s<2$, we may apply once more Propositions~\ref{p.improveadditivity} and~\ref{p.improvelocalization} to obtain that $\Add\left(\tfrac s2, d\right)$ and $\Loc\left(s,\delta,  \tfrac d2 (1+\delta)+\delta -\eta \right)$ hold for every $s<2$ and $\delta,\eta >0$. This completes the proofs of statements (i), (iii) and (iv) of the theorem. 

\smallskip

We have left to prove statement (ii). We postpone this to the next section, since we obtain the result as a byproduct of the analysis in the proof of Proposition~\ref{p.improveadditivity}, see Remark~\ref{r.expectationsconverge}. 
\end{proof}

\section{Improvement of additivity}
\label{s.additivity}

The purpose of this section is to prove Proposition~\ref{p.improveadditivity}. The argument here can be compared with that of Chapter~\ref{c.two}, where we compare two solutions (maximizers of $J$) by first estimating the spatial average of their gradients. However, the mechanism by which we gain control over these spatial averages is somewhat different here. 

\smallskip

We begin by using the $\Fluc(s,\alpha)$ assumption together with the identities~\eqref{e.J1derp} and~\eqref{e.J1derq} to obtain control of the spatial averages of the gradients and fluxes of $v(\cdot,\Phi_r,p,q)$. Since the span of $\nabla v(\cdot,\Phi_r,0,q)$ over $q\in \{ e_1,\ldots,e_d\}$ is $\nabla \A_1$, at least for $r\geq \Y_s$, it is reasonable to expect that we can gain control of the spatial averages of gradients and fluxes of arbitrary elements of~$\A_1$. To put it another way, we are able to show that, for an arbitrary $u\in\A_1$, the convolution $u \ast \Phi_r$ is very close to an $\ahom$-harmonic function. Since $\ahom$-harmonic functions are very smooth, we can use elliptic regularity to rule out the presence of ``wiggles'' in $u \ast \Phi_r$. This allows us to match $\nabla v(\cdot,\Phi_r,p,q)$ to the appropriate corrector field $\xi+\nabla\phi_\xi$, \emph{with $\xi$ deterministic and independent of the scale~$r$}, up to an error of $\O_s(Cr^{-\alpha})$. As a consequence, we deduce that~$\nabla v(\cdot,\Phi_{x,r},p,q)$ is close to $\nabla v(\cdot,\Phi_R,p,q)$ in $L^2(\Phi_{x,r})$ for $R>r$, with an error of $\O_s(Cr^{-\alpha})$. By quadratic response, we then deduce that $J_1(\Phi_R,p,q)$ and $\int_{\Phi_{\sqrt{R^2-r^2}}} J_1(\Phi_{x,r},p,q)\,dx$ are within $\O_{s}(Cr^{-\alpha} )^2 = \O_{s/2}(Cr^{-2\alpha})$ of each other, which is precisely $\Add(s/2,2\alpha)$. 

\smallskip

Throughout this section, we fix parameters 
\begin{equation} 
\label{e.salphaparams}
s\in (0,2] 
\quad \mbox{and} \quad 
\alpha \in \left( 0, \tfrac{d}{s} \right) \cap \left(0,\tfrac d2\right] 
\end{equation}
and suppose that 
\begin{equation} 
\label{e.fluc.ass2}
\Fluc(s,\alpha) \quad \mbox{holds.}
\end{equation}

We first use the $\Fluc(s,\alpha)$ assumption to show a correspondence between the spatial averages of gradients and fluxes of elements of $\A_1$.

\begin{lemma}
\label{l.commutator} 
There exists a constant $C(s,\alpha,d,\Lambda)<\infty$ such that, for every $z\in \Rd$,
\begin{equation} 
\label{e.commutator}
\sup_{u \in \A_1(\Phi_{z,r})} \left|  \int_{\Phi_{z,r}} \left( \a(x) - \ahom\right)\nabla u(x)\,dx  \right| = \O_s\left(Cr^{-\alpha}\right). 
\end{equation}
\end{lemma}
\begin{proof}
\emph{Step 1.} We first show that, for every $x\in\Rd$ and $r\geq r_0$, we can implicitly define a matrix $\ahom_{x,r}$ to satisfy, for every $q\in\Rd$, 
\begin{equation}
 \label{e.abarr}
\E \left[ \int_{\Phi_{x,r}} \left( \a - \ahom_{x,r} \right) \nabla  v(\cdot,\Phi_{x,r},0,q)\right] 
=
0.
\end{equation}
By~\eqref{e.def.r0} and the fact that $(p,q) \mapsto \E\left[ J_1(\Phi_{x,r},p,q) \right]$ is quadratic, we have that, for every $r\geq r_0$, 
\begin{equation*} \label{}
\mbox{the linear map}\quad  q \mapsto \E \left[ \int_{\Phi_{x,r}} \nabla  v(\cdot,\Phi_{x,r},0,q)\right] 
\quad \mbox{is invertible.}
\end{equation*}
Denote the matrix representing the inverse of this map by $\tilde{\a}_{x,r}$. Also let~$\mathbf{b}_{x,r}$ denote the matrix representing the linear map $q\mapsto  \E \left[ -  D_p   J(\Phi_{x,r},0,q) \right]$. We may now define the matrix~$\ahom_{x,r} := \mathbf{b}_{x,r}\tilde{\a}_{x,r}$ and observe, by~\eqref{e.J1derp}, that~\eqref{e.abarr} holds. Moreover, we have 
\begin{equation} 
\label{e.abarrahom}
\left| \ahom - \ahom_{x,r} \right| \leq C r^{-\ep}. 
\end{equation}

\smallskip

\emph{Step 2.} We use the $\Fluc(s,\alpha)$ assumption and polarization to show that, for every $x\in\Rd$, $q\in B_1$ and $r\geq r_0$,
\begin{equation} \label{e.fluxmapflucJ}
\left|  \int_{\Phi_{x,r}} \left( \a - \ahom_{x,r} \right) \nabla  v(\cdot,\Phi_{x,r},0,q) \right| = \O_s\left( C r^{-\alpha} \right). 
\end{equation}
To verify~\eqref{e.fluxmapflucJ}, it suffices by~\eqref{e.abarr} to show that 
\begin{equation} \label{e.fluxmappolar1}
\int_{\Phi_{x,r}} \nabla  v(\cdot,\Phi_{x,r},0,q)
=
\E \left[ \int_{\Phi_{x,r}} \nabla  v(\cdot,\Phi_{x,r},0,q)\right] + \O_s\left( Cr^{-\alpha} \right)
\end{equation}
and
\begin{equation} \label{e.fluxmappolar2}
\int_{\Phi_{x,r}} \a \nabla  v(\cdot,\Phi_{x,r},0,q)
=
\E \left[ \int_{\Phi_{x,r}} \a \nabla  v(\cdot,\Phi_{x,r},0,q)\right] + \O_s\left( Cr^{-\alpha} \right).
\end{equation}
We check only~\eqref{e.fluxmappolar2} since the proof of~\eqref{e.fluxmappolar1} is similar. 
We use~\eqref{e.J polarization} and the $\Fluc(s,\alpha)$ assumption to find that, for every $p' \in B_1$, 
\begin{align*} \label{}
\lefteqn{
p' \cdot \int_{\Phi_{x,r}} \a \nabla  v(\cdot,\Phi_{x,r},0,q)
} \qquad & 
\\ & 
= -J(\Phi_{x,r},p',q) + J(\Phi_{x,r},0,q) + J(\Phi_{x,r},p',0) 
\\ &
= \E\left[ -J(\Phi_{x,r},p',q) + J(\Phi_{x,r},0,q) + J(\Phi_{x,r},p',0) \right] + \O_s\left( Cr^{-\alpha} \right)
\\ &
= \E\left[ p' \cdot \int_{\Phi_{x,r}} \a \nabla  v(\cdot,\Phi_{x,r},0,q) \right] +  \O_s\left( Cr^{-\alpha} \right).
\end{align*}
Thus 
\begin{align*} \label{}
\lefteqn{
\left| \int_{\Phi_{x,r}} \a \nabla  v(\cdot,\Phi_{x,r},0,q) - \E \left[ \int_{\Phi_{x,r}} \a \nabla  v(\cdot,\Phi_{x,r},0,q)\right] \right|
} \qquad & \\
& = 
\sup_{p' \in B_1} \left( p' \cdot \int_{\Phi_{x,r}} \a \nabla  v(\cdot,\Phi_{x,r},0,q) - p'\cdot \E \left[ \int_{\Phi_{x,r}} \a \nabla  v(\cdot,\Phi_{x,r},0,q)\right] \right) 
\\ &
\leq \sum_{i=1}^d \left| e_i \cdot \left(  \int_{\Phi_{x,r}} \a \nabla  v(\cdot,\Phi_{x,r},0,q) - \E \left[ \int_{\Phi_{x,r}} \a \nabla  v(\cdot,\Phi_{x,r},0,q)\right] \right)  \right|
\\ & 
\leq \O_s\left( Cr^{-\alpha} \right). 
\end{align*}
This completes the proof of~\eqref{e.fluxmappolar2} and thus of~\eqref{e.fluxmapflucJ}.

\smallskip

\emph{Step 3.} We show that, for every $r\geq r_0$,
\begin{equation} 
\label{e.fluxmapA1}
\sup_{w\in \A_1(\Phi_{x,r})} 
\left|  \int_{\Phi_{x,r}} \left( \a - \ahom_{x,r} \right) \nabla w \right| \leq \O_s\left( C r^{-\alpha} \right). 
\end{equation}
By~\eqref{e.fluxmapflucJ} and the linearity of~$q\mapsto \nabla  v(\cdot,\Phi_{x,r},0,q)$, we have that 
\begin{equation} 
\label{e.fluxmapflucJsupq}
\sup_{q\in B_1} \left|  \int_{\Phi_{x,r}} \left( \a - \ahom_{x,r} \right) \nabla  v(\cdot,\Phi_{x,r},0,q) \right| = \O_s\left( C r^{-\alpha} \right). 
\end{equation}
Fix $t\in (0,d)$. By~\eqref{e.vqnondegen}, we have that 
\begin{multline*} \label{}
\sup_{u \in \A_1(\Phi_{x,r})} \left|  \int_{\Phi_{x,r}} \left( \a - \ahom_{x,r} \right) \nabla  u \right| \indc_{\{ r \geq \Y_t (x)\}} 
\\
\leq 
C\sup_{q\in B_1} \left|  \int_{\Phi_{x,r}} \left( \a - \ahom_{x,r} \right) \nabla  v(\cdot,\Phi_{x,r},0,q) \right| = \O_s\left( C r^{-\alpha} \right). 
\end{multline*}
On the other hand, taking $t:=\alpha s < d$ and using~\eqref{e.indcYsbound}, we find
\begin{equation*} \label{}
\sup_{u \in \A_1(\Phi_{x,r})} \left|  \int_{\Phi_{x,r}} \left( \a - \ahom_{x,r} \right) \nabla  u \right| \indc_{\{ r < \Y_t(x) \}}
\leq C \indc_{\{ r < \Y_t(x) \}} 
\leq \O_{s} \left( C r^{-\frac ts} \right) 
= \O_{s} \left( C r^{-\alpha} \right).
\end{equation*}
Combining the previous two displays yields~\eqref{e.fluxmapA1}. 

\smallskip

\emph{Step 4.} In view of~\eqref{e.fluxmapA1}, to complete the proof of the lemma it remains to check that, for every $x\in\Rd$ and $r\geq r_0$,
\begin{equation} 
\label{e.abarrahomalpha}
\left| \ahom_{x,r} - \ahom \right|  \leq Cr^{-\alpha}. 
\end{equation}
Indeed, note that we trivially have the estimate~\eqref{e.commutator} for every $r\in [ 1,r_0)$ by~\eqref{e.Jvbounded} and~$r_0\leq C$. 
Fix $x\in\Rd$. By  Lemma~\ref{l.correctorbias}, namely~\eqref{e.biascorrectorsgrad} and~\eqref{e.biascorrectorsflux}, we have
\begin{equation} 
\label{e.aasnap}
\sup_{e\in\partial B_1} \left| \E \left[ \int_{\Phi_{x,r}} \left( \a(x) - \ahom \right) \left( e + \nabla \phi_e \right) \right] \right| 
\leq C \exp\left( -cr^2 \right). 
\end{equation}
On the other hand, by Lemma~\ref{l.correctorminbounds}, namely~\eqref{e.correctorminbounds}, 
\begin{equation*} \label{}
r \geq \Y_t(x)
\implies 
\sup_{e\in\partial B_1} \left\| e + \nabla \phi_e \right\|_{L^2(\Phi_{x,r})} \leq C,
\end{equation*}
and therefore, by~\eqref{e.fluxmapA1},
\begin{equation*} \label{}
\sup_{e\in\partial B_1} \left| \E \left[ \int_{\Phi_{x,r}} \left( \a(x) - \ahom_{x,r} \right) \left( e + \nabla \phi_e \right) \indc_{\{ r \geq \Y_t(x)\}} \right] \right| 
\leq C r^{-\alpha}. 
\end{equation*}
On the other hand, by~\eqref{e.correctorgradbound} and~\eqref{e.indcYsbound},
\begin{equation*} \label{}
\sup_{e\in\partial B_1} \left| \int_{\Phi_{x,r}} \left( \a(x) - \ahom_{x,r} \right) \left( e + \nabla \phi_e \right)\right|  \indc_{\{ r < \Y_t(x)\}}
\leq
C \O_{2}(C) \cdot \O_{2t/d}(Cr^{-\frac{d}2}).
\end{equation*}
We deduce from the previous two displays (and the assumption $\alpha \leq \frac d2$) that 
\begin{equation} 
\label{e.aaxrsnap}
\sup_{e\in\partial B_1} \left| \E \left[ \int_{\Phi_{x,r}} \left( \a(x) - \ahom_{x,r} \right) \left( e + \nabla \phi_e \right)  \right] \right| 
\leq C r^{-\alpha}. 
\end{equation}
The triangle inequality,~\eqref{e.aasnap},~\eqref{e.aaxrsnap} and~\eqref{e.biascorrectorsgrad} give us that \begin{equation*} \label{}
\left| \ahom_{x,r} - \ahom \right|  \leq 
\sup_{e\in\partial B_1} \left| \left(\ahom_{x,r} - \ahom\right) e \right| \leq Cr^{-\alpha},
\end{equation*}
as desired. This completes the proof. 
\end{proof}

We next rephrase the previous lemma as the statement that $u\ast \Phi_r$ is ``almost $\ahom$-harmonic.'' To see why this should be so, notice that if the right side of~\eqref{e.almostahomharm} below identically vanished, then~\eqref{e.almostahomharm} would simply be the weak formulation of the equation $-\nabla \cdot \left(\ahom \nabla (u\ast \Phi_r) \right)=0$. We can think of~\eqref{e.almostahomharm} as roughly asserting that
\begin{equation*} \label{}
\left\| -\nabla \cdot \left(\ahom \nabla (u\ast \Phi_r) \right) \right\|_{H^{-1}(\Rd)} \lesssim Cr^{-\alpha},
\end{equation*}
although this estimate would be a bit too strong since the random variables $\mathcal{H}_r(x)$ introduced in the lemma are not necessarily uniformly bounded in~$x$. 
\begin{lemma}
\label{l.coarsenedequation}
There exists a constant $C(s,\alpha,d,\Lambda)<\infty$ and, for every $x\in \Rd$ and $r\geq 1$, a nonnegative random variable $\mathcal{H}_r(x)$ satisfying 
\begin{equation} \label{e.Hrbound}
\mathcal{H}_r(x) \leq C \wedge \O_s\left( Cr^{-\alpha} \right)
\end{equation}
such that, for every $u \in \A_1$ and $\psi \in H^1_c(\Rd)$, 
\begin{equation} 
\label{e.almostahomharm}
\left| \int_{\Rd} \nabla \psi(x)\cdot \ahom \nabla \left(u\ast\Phi_r\right)(x) \,dx \right| \leq \int_{\Rd} \left| \nabla \psi(x) \right| \left\| \nabla u \right\|_{L^2(\Phi_{x,r})} \mathcal{H}_r(x)\,dx.
\end{equation}
\end{lemma}
\begin{proof}
We define 
\begin{equation*} \label{}
\mathcal{H}_r(x)
:= C \sup_{u \in \A_1\left(\Phi_{x,r}\right)} \left| \int_{\Phi_{x,r}} \left( \a-\ahom\right)\nabla u \right|.
\end{equation*}
By Lemma~\ref{l.commutator}, we see that~\eqref{e.Hrbound} is satisfied. Fix $u\in \A_1$, $\psi\in H^1_c(\Rd)$ and compute
\begin{align*} \label{}
\lefteqn{
\left| \int_{\Rd} \nabla \psi(x) \cdot \ahom \nabla \left( u\ast \Phi_r \right)(x) \,dx \right|
} \qquad & \\
& = \left| \int_{\Rd} \nabla \psi(x) \cdot \left( \int_{\Phi_{x,r}} \ahom \nabla u \right)\,dx \right|
\\ &
\leq \left| \int_{\Rd} \nabla \psi(x) \cdot \left( \int_{\Phi_{x,r}} \a \nabla u \right)\,dx \right| +  \int_{\Rd} \left| \nabla \psi(x) \right| \left\| \nabla u \right\|_{L^2(\Phi_{x,r})} \mathcal{H}_r(x)\,dx
\\ & 
=  \left| \int_{\Rd} \nabla \left( \psi \ast \Phi_{r} \right) (x) \cdot \a(x) \nabla u(x) \,dx \right| +  \int_{\Rd} \left| \nabla \psi(x) \right| \left\| \nabla u \right\|_{L^2(\Phi_{x,r})} \mathcal{H}_r(x)\,dx.
\end{align*}
The first term on the last line is zero, by the equation. To make the integration by parts rigorous, we observe that $\psi \ast \Phi_r$ decays faster than any polynomial at infinity and $u$ has at most linear growth. This completes the proof of the lemma. 
\end{proof}

The next lemma is the backbone of this section. It asserts that the spatial average of elements of $\A_1$ are close to constants on every scale. This will eventually allow us to compare maximizers on different scales.

\begin{lemma}
\label{l.harmonicapproximation}
There exists $C(s,\alpha,d,\Lambda)<\infty$ such that, for every $r\in \left[1,R/\sqrt{2} \right]$, 
\begin{equation} 
\label{e.harmapproxflat}
\sup_{w\in \A_1(\Phi_R)}  \inf_{\xi \in \Rd} 
\left\|  \nabla w \ast \Phi_{r} - \xi  \right\|_{L^2\left( \Phi_{\sqrt{R^2-r^2}} \right)}^2 
\leq
\O_{s/2}\left( Cr^{-2\alpha} \right). 
\end{equation}
\end{lemma}

\begin{proof}
Throughout, we fix $T:= \sqrt{R^2-r^2}$.  

\smallskip

\emph{Step 1.}  
We first reduce ourselves to analyzing maximizers of $J_1$. For convenience, let $X$ denote the random variable
\begin{equation*}
X :=  \sup_{w \in \A_1(\Phi_R)} \inf_{\xi \in \R^d} \| \nabla w  \ast \Phi_{r} - \xi \|_{L^2(\Phi_{T })}^2.
\end{equation*}
Note that $0 \leq X \leq 1$ by definition. Using $s \alpha < d$ and~\eqref{e.indcYsbound}, we find that 
\begin{align*}
X  \indc_{\{\Y_s \geq T \}} 
\leq \indc_{\{\Y_s \geq T \}} 
\leq \O_{s/2}\left(C T^{-2\alpha}\right) 
\leq \O_{s/2}\left(C r^{-2\alpha}\right)\,.
\end{align*}
On the other hand, we see from~\eqref{e.vqnondegen} and $R \geq T$ that, for some $C(\Lambda)<\infty$,
\begin{equation*} \label{}
T > \Y_s \implies
\nabla \A_1(\Phi_T) 
\subseteq
\left\{ \nabla v(\cdot,\Phi_{R},0,q) \,:\, q\in B_C \right\}.
\end{equation*}
Therefore 
\begin{equation} 
\label{e.tildescriptXbound}
X \indc_{\{ \Y_s < T \}} \leq C \sup_{q\in B_1} \inf_{\xi\in\Rd} \left\|  \nabla v(\cdot,\Phi_{R},0,q) \ast \Phi_{r} - \xi \right\|_{L^2(\Phi_{T })}^2\indc_{\{ \Y_s < T \}}.
\end{equation}
By linearity, the desired inequality~\eqref{e.harmapproxflat} follows from~\eqref{e.tildescriptXbound} and the following statement:
\begin{multline}
\label{e.harmonicapprox'}
\sum_{j=1}^d \inf_{\xi \in \R^d} \left\| \nabla v(\cdot,\Phi_R,0,e_j) \ast \Phi_{r} - \xi \right\|_{L^2(\Phi_{T })}^2 \indc_{\{ \Y_s < T \}}
\\ 
\leq \max\left\{\O_{s/2}\left(C r^{-2\alpha}\right)\,, \O_{\frac{s \alpha}{d+2\alpha} }\left(Cr^{-(d+2\alpha) } \right) \right\}\,.
\end{multline}
Indeed, by Lemma~\ref{l.change-s}(i),
\begin{equation*} \label{}
\O_{\frac{s \alpha}{d+2\alpha} }\left(Cr^{-(d+2\alpha) } \right) \wedge 1 \le \O_{s/2}\Ll( C r^{-2\al} \Rr).
\end{equation*}
Thus, since $X  \in [0,1]$, it suffices to prove~\eqref{e.harmonicapprox'}. For the rest of the argument, we fix $j\in \{1,\ldots,d\}$ and set
\begin{equation*} \label{}
v := v(\cdot,\Phi_R,0,e_{j}).
\end{equation*}
To simplify the notation, we also assume that $T > \Y_s$, so that we may drop $\indc_{\{ T > \Y_s\}}$ from all the expressions.

\smallskip

\emph{Step 2.} Harmonic approximation of $v \ast \Phi_{r}$ and iteration. For convenience, we denote 
\begin{equation*}  
w := v \ast \Phi_r.
\end{equation*}
For each $S\geq T $, we introduce an~$\ahom$-harmonic approximation of~$w$ in $B_S$, which we denote by~$h_S$. We take $h_S$ to be the unique element of $(w + H^1_0(B_S))\cap \Ahom(B_S)$. It follows from Lemma~\ref{l.coarsenedequation} (choosing $\psi = h_S-w$ and applying the Cauchy-Schwarz inequality) that
\begin{align*}
\fint_{B_S} \left| \nabla w(x) - \nabla h_S(x) \right|^2\,dx 
& \leq 
 \fint_{B_S} \left\| \nabla v \right\|_{L^2(\Phi_{x,r})}^2 \left( \mathcal{H}_r(x) \right)^2  \,dx,
\end{align*}
where $\mathcal{H}_r(x) = C \wedge \O_s(C r^{-\alpha})$ is as in Lemma~\ref{l.coarsenedequation}.
Therefore, for every $\theta  \in (0,1]$,
\begin{equation}
\label{e.harmapprox}
\fint_{B_{\theta S}} \left| \nabla w(x) - \nabla h_S(x) \right|^2\,dx 
 \leq \theta^{-d}  \fint_{B_S}  \left\| \nabla v \right\|_{L^2(\Phi_{x,r})}^2  \left( \mathcal{H}_r(x) \right)^2  \,dx.
\end{equation}
By the regularity of $\ahom$-harmonic functions, we find, for every $\eta \in \R^d$, 
\begin{align}
\label{e.tildep}
\inf_{\xi \in \R^d}  \left\|\nabla h_S- \xi \right\|_{\underline L^2(B_{\theta S})}   
& \leq C \theta   \left\| \nabla h_S- \eta \right\|_{\underline L^2(B_S)}. 
\end{align}
By the triangle inequality and the previous two displays, denoting
\begin{equation*}
\omega(\varrho):= \varrho^{-\frac12} \inf_{\xi\in\Rd} \left\| \nabla w- \xi \right\|_{\underline{L}^2(B_\varrho)},
\end{equation*}
we obtain, for $\theta = (2C)^{-\frac12}$,
\begin{equation*}
\omega(\theta S) \leq \frac12 \omega(S) + C S^{-\frac12} 
\left( \fint_{B_S}  \left\| \nabla v \right\|_{L^2(\Phi_{x,r})}^2  \left( \mathcal{H}_r(x) \right)^2 \,dx \right)^{\frac12}.
\end{equation*}
Setting $S_j := \theta^{-j} T $ and summing over all the scales, using also the fact that $\omega(\varrho) \to 0$ as $\varrho \to \infty$ on the event $\{ T > \Y_s \}$, yields
\begin{equation}
\label{e.omegaRiter}
 \sum_{j=0}^\infty \omega(S_j) 
 \leq  C T ^{-\frac12}\overline{\mathcal{H}}  \,,
\end{equation}
where we have defined
\begin{equation} \label{e.Y=sum}
\overline{\mathcal{H}} := \sum_{n=1}^\infty \theta^{\frac n2} \left( \fint_{B_{\theta^{-n} T }}  \left\| \nabla v \right\|_{L^2(\Phi_{x,r})}^2  \left( \mathcal{H}_r(x) \right)^2 \,dx \right)^{\frac12}.
\end{equation}
Letting $\xi_j \in \R^d$ be the minimizer appearing in the definition of $\omega(S_j)$, we obtain by the triangle inequality that, for $m>j$,
\begin{align*} 
\left|\xi_j - \xi_{j+1} \right| & \leq CS_j^{\frac12}   \left(\omega(S_j) + \omega(S_{j+1})\right) \,.
\end{align*}
Therefore we get, again by the triangle inequality, 
\begin{align*} 
\left|\xi_m - \xi_{0} \right|  & \leq  \sum_{j=0}^{m-1} \left|\xi_j - \xi_{j+1} \right| \leq \sum_{j=0}^m S_j^{\frac12}  \omega(S_j) \leq C\left( \frac{S_m}{T } \right)^{\frac12} \overline{\mathcal{H}}  \,.
\end{align*}
Using now the decay properties of $\Phi_{T }$ we conclude that
\begin{align*} 
\left\| \nabla w - \xi_0 \right\|_{L^2(\Phi_{T })}^2
 &  = \sum_{m=-\infty}^\infty \int_{\{ S_m \leq |y| \leq S_{m+1} \} } \Phi_{T }(y) |\nabla w(y) - \xi_0|^2 \, dy
\\ &    \leq C  \sum_{m=0}^\infty \exp\left( - c \theta^{-2m}\right) \theta^{-dm} \left( \left | \xi_m {-} \xi_0 \right|^2 +  \left\| \nabla w {-} \xi_m \right\|_{\underline{L}^2(B_{S_m})}^2 \right) 
\\ &   \leq C  \overline{\mathcal{H}} \sum_{m=0}^\infty \exp\left( - c \theta^{-2m}\right) \theta^{-(d + 1)m}   
\\ &   \leq C  \overline{\mathcal{H}} \,.
\end{align*}
The rest of the proof is thus devoted to estimating the random variable $\overline{\mathcal{H}} $. 

\smallskip

\emph{Step 3.}  We use the Lipschitz estimate to pull the $\left\| \nabla v \right\|_{L^2(\Phi_{x,r})}^2 $ term outside the integrals appearing in $\overline{\mathcal{H}} $ defined in~\eqref{e.Y=sum}. The claim is that 
\begin{multline}
\label{e.lippullout}
 \fint_{B_{\theta^{-n} T }}  \left\| \nabla v \right\|_{L^2(\Phi_{x,r})}^2  \left( \mathcal{H}_r(x) \right)^2 \,dx  \\
  \leq C  \fint_{B_{\theta^{-n}T } } \left( 1+\indc_{\{r\leq \Y_s(x)\}} \left(  \frac{\X(x)}{r} \right)^{d} \right) \left(\mathcal{H}_r(x) \right)^2 \,dx \,.
\end{multline}
First, since we are working on the event $\{ T  >  \Y_s\}$, we have that there is $p_v \in \R^d$ such that, for all $S \geq T $, 
$$
\left\|v \right\|_{\underline{L}^2(B_S)} \leq C S \left | p_v \right| \leq C \left(\frac{S}{T }\right) \left\| v \right\|_{\underline{L}^2(B_{T })}.
$$
Without loss of generality we may assume that $(v)_{B_{T }} = 0$. The  Poincar\'e inequality and the fact that we have $ \left\| \nabla v \right\|_{L^2(\Phi_R)} \leq C$ yield that  
\begin{equation} \label{e.vgrowthharmaprox}
\left\|v \right\|_{\underline{L}^2(B_S)} \leq C S  .
\end{equation}
Suppose next that $r \leq S \leq  |x| + T $ and $\Y_s(x) \leq |x|+T $. The Caccioppoli inequality and the $C^{0,1}$-type estimate (Theorem~\ref{t.Lipschitz}) then yield
\begin{align*} 
\left\| \nabla v \right\|_{\underline{L}^2(B_S(x))} 
& \leq \left(\frac{\Y_s(x) \vee S}{S}\right)^{\frac d2} \left\| \nabla v \right\|_{\underline{L}^2(B_{\Y_s(x) \vee S}(x))} 
\\
& \leq C \left(\frac{\Y_s(x) \vee S}{S}\right)^{\frac d2} 
\left( \Y_s(x) \vee S \right)^{-1}
\left\| v - (v)_{B_{2(\Y_s(x) \vee S )}}\right\|_{\underline{L}^2(B_{2(\Y_s(x) \vee S )}(x))}
\\& \leq C \left(\frac{\Y_s(x) \vee r}{r}\right)^{\frac d2} \,. 
\end{align*}
On the other hand, if $ r\leq S \leq |x| + T $ and $\Y_s(x) \geq |x|+T $, then we get 
\begin{equation*} 
\left\| \nabla v \right\|_{\underline{L}^2(B_S(x))} \leq \left(\frac{|x|+T }{S}\right)^{\frac d2} \left\| \nabla v \right\|_{\underline{L}^2(B_{|x| + T  }(x))} \leq 
C \left(\frac{\Y_s(x) }{r}\right)^{\frac d2} \,,
\end{equation*}
and finally if $ S \geq |x| + T $, then directly $\left\| \nabla v \right\|_{\underline{L}^2(B_S(x))}  \leq C $ by the Caccioppoli inequality and~\eqref{e.vgrowthharmaprox}. Using these gives
\begin{align*} 
\left\| \nabla v \right\|_{L^{2}(\Phi_{x,r})}^2 & \leq C \int_r^\infty \left(\frac{S}{r} \right)^{d} \exp \left(- c \frac{S^2}{r^2} \right) \fint_{B_{S}(x)} \left|\nabla v(y)\right|^2 \, dy\, \frac{dS}{S}
\\ & \leq C \left(\frac{\Y_s(x) \vee r}{r}\right)^{d}  \int_r^{|x|+T } \left(\frac{S^2}{r^2} \right)^{\frac d2} \exp \left(- c \frac{S^2}{r^2} \right) \, \frac{dS}{S}
\\ & \qquad + C \int_{|x|+T }^\infty \left(\frac{S}{r} \right)^{d } \exp \left(- c \frac{S^2}{r^2} \right) \, \frac{dS}{S}
\\ & \leq C \left(\frac{\Y_s(x) \vee r}{r}\right)^{d}.
\end{align*}
We deduce that 
\begin{align*}
 \fint_{B_{\theta^{-n} T }}  \left\| \nabla v \right\|_{L^2(\Phi_{x,r})}^2  \left( \mathcal{H}_r(x) \right)^2  \,dx 
 & \leq \fint_{B_{\theta^{-n} T }}  \left( 1+\indc_{\{r\leq \Y_s(x)\}} \left(  \frac{\Y_s(x)}{r} \right)^{d} \right)  \left( \mathcal{H}_r(x) \right)^2 \,dx , 
\end{align*}
which completes the proof of~\eqref{e.lippullout}. 

\smallskip

\emph{Step 4.} We next estimate the term in~\eqref{e.lippullout}.  The claim is that there exists $C(\alpha,s,d,\Lambda)<\infty$ such that
\begin{equation}
\label{e.Xbigweirdo}
\left( 1 + \indc_{\{r\leq \Y_s(x)\}} \left(  \frac{\Y_s(x)}{r} \right)^{d} \right) \left( \mathcal{H}_r(x) \right)^2 
\leq \O_{\frac{s \alpha}{d + 2\alpha}}\left(C r^{-(d + 2\alpha)} \right) \vee \O_{s/2}(Cr^{-2\alpha}).
\end{equation}
First, note that~\eqref{e.Hrbound} gives 
\begin{equation*} \label{}
\left( \mathcal{H}_r(x) \right)^2 = \O_{s/2}(C r^{-2\alpha}).
\end{equation*}
Moreover, since $ s \alpha < d$,  we have that
\begin{equation}
\label{e.Xweirdo}
\indc_{\{r\leq\Y_s(x)\}} \left(  \frac{\Y_s(x)}{r} \right)^{d}  \le \left(  \frac{\Y_s(x)}{r} \right)^{d}  = \O_{s \alpha/d} \left(Cr^{-d }\right).
\end{equation}
By Lemma~\ref{l.change-s}(ii), applied with $s_1 = \frac s2$, $s_2 = \frac{s \alpha}{d}$, $\theta_1 = Cr^{-2\alpha}$ and $\theta_2 = Cr^{-d}$, we therefore obtain
$$
\indc_{\{r\leq \Y_s(x)\}} \left(  \frac{\Y_s(x)}{r} \right)^{d} \left( \mathcal{H}_r(x) \right)^2  = \O_{\frac{s \alpha}{d + 2\alpha}}\left(C r^{-(d+ 2\alpha)} \right).
$$
This proves~\eqref{e.Xbigweirdo}.  

\smallskip

\emph{Step 5.} We complete the proof. Combining~\eqref{e.omegaRiter},~\eqref{e.lippullout} and~\eqref{e.Xbigweirdo}, we get by Lemma~\ref{l.sum-O} that
\begin{align*}
 \overline{\mathcal{H}}
& \leq C\sum_{n=0}^\infty \theta^{\frac n2}  \left( \fint_{B_{\theta^{-n} R}}  \left\| \nabla v \right\|_{L^2(\Phi_{x,r})}^2  \left( \mathcal{H}_r(x) \right)^2 \,dx \right)^{\frac12} \\
& \leq C\sum_{n=0}^\infty \theta^{\frac n2}  \left( 
 \fint_{B_{\theta^{-n} R}}  \left( 1 + \indc_{\{r\leq \Y_s(x)\}} \left(  \frac{\Y_s(x)}{r} \right)^{d} \right) \left( \mathcal{H}_r(x) \right)^2   \,dx \right)^{\frac12}  \\
& \leq \left( \O_{\frac{s \alpha}{d + 2\alpha}}\left(C \sum_{n=0}^\infty \theta^{\frac n2}  r^{-(d + 2\alpha)} \right) \vee \O_{s/2}\left(C \sum_{n=0}^\infty \theta^{\frac n2}  r^{-2\alpha} \right) \right)^{\frac12}    
 \\
 & = \left( \O_{\frac{s \alpha}{d + 2\alpha}}\left(C r^{-(d + 2\alpha)} \right) \vee \O_{s/2}(Cr^{-2\alpha}) \right)^{\frac12}.
 \end{align*}
This combined with the result of Step~1 yields~\eqref{e.harmonicapprox'} and completes the proof by the discussion in the first step.
 \end{proof}

By summing the result of the previous lemma over the scales, we can control the spatial averages of the gradients of the correctors, up to an error of $\O_{s}(Cr^{-\alpha})$.

\begin{lemma}
\label{l.okcupid}
There exists $C(s,d,\Lambda)<\infty$ such that, for every $\xi \in B_1$, $x\in\Rd$, $R\geq 2$ and $r\in \left[ 1 , R/\sqrt{2} \right]$,
\begin{equation*} \label{}
\left\| \nabla \phi_\xi \ast \Phi_r \right\|_{L^2\left( \Phi_{x,\sqrt{R^2-r^2}} \right)} = \O_{s}\left(Cr^{-\alpha}\right). 
\end{equation*}
\end{lemma}
\begin{proof}
Fix $\xi\in B_1$. By Lemma~\ref{l.harmonicapproximation}, for each $r\geq 1$, there exists a random vector $\xi(r)\in\Rd$ satisfying 
\begin{equation} 
\label{e.findslope}
\left\| \nabla \phi_\xi \ast \Phi_r - \xi(r) \right\|_{L^2\left( \Phi_{2r } \right)}
\leq 
\O_s\left( Cr^{-\alpha} \right).
\end{equation}
It follows from the triangle inequality, the fact that $\Phi_{r}\leq C \Phi_{4r}$, the semigroup property and the previous display that
\begin{align*} \label{}
\left| \xi(r) - \xi(2r) \right| 
& \leq 
\left\| \left(\nabla \phi_\xi - \xi(r) \right)  \ast \Phi_{2r}  \right\|_{L^2\left( \Phi_{r} \right)}
+
\left\| \left(\nabla \phi_\xi - \xi(2r) \right)  \ast \Phi_{2r} \right\|_{L^2\left( \Phi_{r} \right)}
\\ 
& =
\left\| \left(\nabla \phi_\xi \ast \Phi_{r} - \xi(r) \right) \ast \Phi_{\sqrt{3r^2}}   \right\|_{L^2\left( \Phi_{r} \right)}
+
\left\| \nabla \phi_\xi \ast \Phi_{2r}  - \xi(2r)   \right\|_{L^2\left( \Phi_{r} \right)}
\\ & 
\leq 
\left\| \nabla \phi_\xi \ast \Phi_{r} - \xi(r)   \right\|_{L^2\left( \Phi_{2r} \right)}
+
C\left\| \nabla \phi_\xi \ast \Phi_{2r} - \xi(2r) \right\|_{L^2\left( \Phi_{4r} \right)} 
\\ &
\leq \O_s\left( Cr^{-\alpha} \right).
\end{align*}
 By summing over the scales, we obtain that $\{ \xi(2^nr) \}_{n\in\N}$ is Cauchy and its limit $\xi(\infty)$ satisfies
\begin{equation*} \label{}
\left| \xi(r) - \xi(\infty) \right| \leq \O_s\left( Cr^{-\alpha} \right).
\end{equation*}
Since $\nabla \phi_\xi$ is a stationary function with zero mean, the ergodic theorem implies that, in fact, $\xi(\infty) = 0$, $\P$-almost surely. Combining this with~\eqref{e.findslope} yields 
\begin{equation*} \label{}
\left\| \nabla \phi_\xi \ast \Phi_r \right\|_{L^2\left( \Phi_{2r} \right)}
\leq 
\O_s\left( Cr^{-\alpha} \right).
\end{equation*}
Since $\Phi_{x,r} \leq C \Phi_{2r}$ for $x\in B_r$, we deduce that, for every $x\in B_r$, 
\begin{equation*} \label{}
\left\| \nabla \phi_\xi \ast \Phi_r \right\|_{L^2\left( \Phi_{x,r} \right)}
\leq 
\O_s\left( Cr^{-\alpha} \right).
\end{equation*}
By stationarity, we have the same estimate for every $x\in\Rd$. The desired estimate now follows from the semigroup property of the heat kernel and Lemma~\ref{l.sum-O}. 
\end{proof}

Using the previous lemma, we can ``match the parameters''  by identifying the maximizer $v(\cdot,\Phi_{x,r},p,q)$ as the corrector $\bar \phi_\xi$ with $\xi = -p + \ahom^{-1}q$, up to an error of~$\O_{s} \left( Cr^{-\alpha} \right)$. Recall that~$\bar\phi_\xi$ is defined in~\eqref{e.correctorbarzz}. 

\begin{lemma}
\label{l.shortcuts}
There exists $C(s,\alpha,d,\Lambda)<\infty$ such that, for every $r\geq 1$, $x\in\Rd$ and $p,q\in B_1$, 
\begin{equation}
\label{e.matchJtophi}
\left\| \nabla v(\cdot,\Phi_{x,r},p,q) - \nabla\overline{\phi}_{(-p + \shortbar{\a}^{\,-1}q)} \right\|_{L^2\left(\Phi_{x,r}\right)} 
\leq 
\O_{s} \left( Cr^{-\alpha} \right).
\end{equation}
\end{lemma}

\begin{proof}
\emph{Step 1.}
We first show that, for some $C(d,\Lambda)<\infty$,
\begin{equation}
\label{e.L2byspatavg}
\sup_{w\in \A_1(\Phi_{x,r})} \left( \left\| \nabla w \right\|_{L^2(\Phi_{x,r})} - C\left| \left(  \nabla w \right)_{\Phi_{x,r}}  \right| \right) \leq \O_s\left( Cr^{-\alpha} \right). 
\end{equation}
By~Lemma~\ref{l.okcupid}, for each $\xi \in B_C$,
\begin{align}
\label{e.followthrough}
\left| \int_{\Phi_{x,r}} 
\nabla \overline{\phi}_{\xi}  -\xi \right| 
= \left| \int_{\Phi_{x,r}} 
\nabla {\phi}_{\xi} \right| 
\le \left\| \nabla \phi_\xi \ast \Phi_{r/\sqrt{2}} \right\|_{L^2\left( \Phi_{x,r/\sqrt{2}} \right)} 
\leq \O_{s} \left( Cr^{-\alpha} \right).
\end{align}
It follows from this and~\eqref{e.correctorsballz} that 
\begin{equation*} \label{}
\sup_{w\in \A_1(\Phi_{x,r})} \left( \left\| \nabla w \right\|_{L^2(\Phi_{x,r})} - C\left| \left(  \nabla w \right)_{\Phi_{x,r}} \right| \right)_+ \indc_{\{ r \geq \Y_s(x)\}} \leq \O_s\left( Cr^{-\alpha} \right). 
\end{equation*}
On the other hand, by~\eqref{e.indcYsbound} and $s\alpha < d$, 
\begin{equation*} \label{}
\sup_{w\in \A_1(\Phi_{x,r})} \left\| \nabla w \right\|_{L^2(\Phi_{x,r})}  \indc_{\{ r \leq \Y_s(x)\}} 
\leq  \indc_{\{ r \le \Y_s(x)\}}  \leq \O_s\left( Cr^{-\alpha} \right).
\end{equation*}
This completes the proof of~\eqref{e.L2byspatavg}.

\smallskip

\emph{Step 2.} We match each~$\nabla v(\cdot,\Phi_{x,r},p,q)$ to~$\nabla \overline{\phi}_\xi$ for a deterministic~$\xi$. Let $L_{x,r}:\Rd\times \Rd \to \Rd$ denote the deterministic linear map 
\begin{equation*} \label{}
(p,q) \mapsto \E \left[ \int_{\Phi_{x,r}} 
\nabla v(\cdot,\Phi_{x,r},p,q) \right].
\end{equation*}
Fix $p,q\in B_1$ and set $\xi:= L_{x,r}(p,q)$. Then $|\xi| \leq C$ by~\eqref{e.Jvbounded} and, by~\eqref{e.fluxmappolar1},
\begin{equation} 
\label{e.giveittomefluc}
\left| \int_{\Phi_{x,r}} \nabla v(\cdot,\Phi_{x,r},p,q) - \xi \right| \leq \O_s \left( Cr^{-\alpha} \right).
\end{equation}
By the previous line and~\eqref{e.followthrough},
\begin{equation*} \label{}
\left| \int_{\Phi_{x,r}} \left( \nabla v(\cdot,\Phi_{x,r},p,q) - \nabla \overline{\phi}_{\xi}\right) \right| 
 \leq 
\O_{s} \left( Cr^{-\alpha} \right).
\end{equation*}
By~\eqref{e.L2byspatavg}, we deduce that 
\begin{equation} 
\label{e.yescupid}
\left\|  \nabla v(\cdot,\Phi_{x,r},p,q) - \nabla \overline{\phi}_{\xi} \right\|_{L^2(\Phi_{x,r})} 
\leq \O_{s} \left( Cr^{-\alpha} \right).
\end{equation}
%
%
%
%
\emph{Step 3.} We identify $L_{x,r}(p,q)$ up to a suitable error. The claim is that 
\begin{equation} 
\label{e.identifyLxr}
\left| L_{x,r}(p,q) - \left( -p + \ahom^{-1}q \right) \right| \leq Cr^{-\alpha}.
\end{equation}
As before, we denote $\xi:= L_{x,r}(p,q)$. By~\eqref{e.yescupid} and quadratic response,
\begin{equation*} \label{}
 \left| J_1(\Phi_{x,r},p,q) 
- 
 \int_{\Phi_{x,r}} \left( -\frac12 \nabla \overline{\phi}_{\xi} \cdot \a \nabla \overline{\phi}_{\xi}  - p\cdot \a\nabla \overline{\phi}_{\xi} +q\cdot\nabla \overline{\phi}_{\xi}  \right)\right| \\
\leq \O_{s/2} \left( Cr^{-2\alpha} \right). 
 \end{equation*}
Taking expectations and applying Lemma~\ref{l.correctorbias} yields
\begin{align*} \label{}
\E \left[ J_1(\Phi_{x,r},p,q) \right] 
&
\leq \E \left[  \int_{\Phi_{x,r}} \left( -\frac12 \nabla \overline{\phi}_{\xi} \cdot \a \nabla \overline{\phi}_{\xi}  - p\cdot \a\nabla \overline{\phi}_{\xi} +q\cdot\nabla \overline{\phi}_{\xi}  \right)\right] 
+C r^{-2\alpha} 
\\ &
\leq \left( -\frac12 \xi \cdot \ahom \xi - p\cdot \ahom \xi +q\cdot\xi  \right)
+C r^{-2\alpha}.
\end{align*}
On the other hand, it is clear by testing the definition of $J_1$ with $\overline{\phi}_{\xi'}$ that 
\begin{align*} \label{}
\E \left[ J_1(\Phi_{x,r},p,q) \right] 
&
\geq \sup_{\xi'\in\Rd}  \E \left[  \int_{\Phi_{x,r}} \left( -\frac12 \nabla \overline{\phi}_{\xi'} \cdot \a \nabla \overline{\phi}_{\xi'}  - p\cdot \a\nabla \overline{\phi}_{\xi'} +q\cdot\nabla \overline{\phi}_{\xi'}  \right)\right] 
\\ & 
\geq  \sup_{\xi'\in\Rd} \left( -\frac12 \xi' \cdot \ahom \xi' - p\cdot \ahom \xi' +q\cdot\xi'  \right) - C\exp(-cr^2).
\end{align*}
Thus
\begin{equation*} \label{}
\left( -\frac12  \xi \cdot \ahom \xi - p\cdot \ahom \xi +q\cdot\xi  \right) 
\geq \sup_{\xi'\in\Rd} \left( -\frac12  \xi' \cdot \ahom \xi' - p\cdot \ahom \xi' +q\cdot\xi'  \right)
- Cr^{-2\alpha}. 
\end{equation*}
The supremum on the right side is attained at $\xi_* := -p + \ahom^{-1}q$ and hence by quadratic response, we get $\left|\xi - \xi_*\right|\leq Cr^{-\alpha}$. This is~\eqref{e.identifyLxr}. 

\smallskip

\emph{Step 4.} The conclusion. With $\xi_* := -p + \ahom^{-1}q$ and $\xi:=L_{x,r}(p,q)$ as above, we apply the triangle inequality,~\eqref{e.yescupid} and~\eqref{e.identifyLxr} to find that
\begin{align*} \label{}
\lefteqn{
\left\| \nabla v(\cdot,\Phi_{x,r},p,q) - \nabla\overline{\phi}_{\xi_*} \right\|_{L^2\left(\Phi_{x,r}\right)}
} \qquad & \\
& 
\leq 
\left\|  \nabla v(\cdot,\Phi_{x,r},p,q) - \nabla \overline{\phi}_{\xi} \right\|_{L^2(\Phi_{x,r})} 
+ \left\| \nabla \overline{\phi}_{\xi} - \nabla \overline{\phi}_{\xi_*} \right\|_{L^2(\Phi_{x,r})}
\\ & 
\leq \left\| \nabla \overline{\phi}_{\xi-\xi_*}  \right\|_{L^2(\Phi_{x,r})} + \O_s\left( Cr^{-\alpha} \right). 
\\ &
\leq Cr^{-\alpha}  \sup_{e\in \partial B_1}\left\| \nabla \overline{\phi}_{e}  \right\|_{L^2(\Phi_{x,r})}+ \O_s\left( Cr^{-\alpha} \right).
\end{align*}
On the other hand, we see from~\eqref{e.correctorgradbound} and Lemma~\ref{l.sum-O} that 
\begin{equation*} \label{}
\sup_{e\in \partial B_1}\left\| \nabla \overline{\phi}_{e}  \right\|_{L^2(\Phi_{x,r})}
\leq \O_{2+\ep} \left( C \right) \leq \O_s(C). 
\end{equation*}
This completes the proof of~\eqref{e.matchJtophi}.
\end{proof}

The ``parameter matching'' of the previous lemma allows us, in particular, to compare the maximizers of $J_1$ on different scales by matching them to the same corrector. By quadratic response, this implies the additivity statement---with double the exponent and half the stochastic integrability. We then interpolate this result with the base case to obtain the additivity statement with no loss of stochastic integrability but maintaining still a slight improvement of the exponent. We conclude this section by presenting the details of this argument.

\begin{proof}[{Proof of Proposition~\ref{p.improveadditivity}}]
The triangle inequality and~\eqref{e.matchJtophi} imply
\begin{equation} 
\label{e.matchtwoscales}
\int_{\Phi_{\sqrt{R^2-r^2}}} \left\| \nabla v\left(\cdot,\Phi_R,p,q\right)  -  \nabla v(\cdot,\Phi_{x,r},p,q) \right\|_{L^2(\Phi_{x,r})}^2 \,dx \\
\leq \O_{s/2}\left(Cr^{-2\alpha}\right).
\end{equation}
By~\eqref{e.matchtwoscales} and quadratic response, we have, for every $p,q\in B_1$ and $1\leq r < R$,
\begin{equation} 
\label{e.add2alpha}
J(\Phi_R,p,q) 
= \int_{\Phi_{\sqrt{R^2-r^2}}} J(\Phi_{x,r},p,q)\,dx 
+ \O_{s/2}\left(Cr^{-2\alpha}\right).
\end{equation}
This is $\Add(s/2,2\alpha)$. 

\smallskip

Suppose now that $\alpha < \frac{d}s-\ep$ for a given $\ep>0$. We next argue that $\Add(s/2,2\alpha)$ and the base case estimate Proposition~\ref{p.basecase} imply that $\Add(s,\alpha+\eta)$ holds for some $\eta(\ep,s,d,\Lambda)>0$, by an interpolation between these statements. Fix $t\in (\alpha s,d)$ to be selected. We first observe that, by Proposition~\ref{p.basecase} and the triangle inequality, we have the~$\P$--a.s.~bound
\begin{equation} 
\label{e.adddumbbeta}
\left| 
 \int_{\Phi_{\sqrt{R^2-r^2}}} \left(J(\Phi_R,p,q) - J(\Phi_{x,r},p,q)\right) \indc_{\{ \Y_t(x) \leq r \}} \,dx \right| \indc_{\{ \Y_t \leq R \}}
\leq Cr^{-\beta}.
\end{equation}
By~\eqref{e.Jvbounded} and Lemma~\ref{l.sum-O}, for every $s'>0$,
\begin{multline} 
\label{e.killindcr}
\int_{\Phi_{\sqrt{R^2-r^2}}} \left( \left| J(\Phi_R,p,q)\right|+\left|J(\Phi_{x,r},p,q)\right| \right) \indc_{\{ \Y_t(x) \geq r \}} \,dx
\\
\leq C \int_{\Phi_{\sqrt{R^2-r^2}}} \indc_{\{ \Y_t(x) \geq r \}} \,dx \leq \O_{s'/2}\left( Cr^{-\frac {2t}{s'}} \right)
\end{multline}
in view of the fact that, for every $x\in \Rd$, 
\begin{equation*} \label{}
\indc_{\{ \Y_t(x) \geq r \}} \leq \left( \frac{\Y_t(x)}{r}\right)^{\frac {2t}{s'}} \leq \O_{s'/2}\left( Cr^{-\frac {2t}{s'}} \right).
\end{equation*}
Taking $s'=s$ and using $2t/s>2\alpha$, the previous inequality and~\eqref{e.add2alpha} imply that 
\begin{equation} 
\label{e.add2alphaindc}
\left| 
 \int_{\Phi_{\sqrt{R^2-r^2}}} \left(J(\Phi_R,p,q) - J(\Phi_{x,r},p,q)\right) \indc_{\{ \Y_t(x) \leq r \}} \,dx \right| \indc_{\{ \Y_t \leq R \}}
\leq \O_{s/2}\left( Cr^{-2\alpha} \right). 
\end{equation}
The estimates~\eqref{e.add2alphaindc},~\eqref{e.adddumbbeta} and Lemma~\ref{l.change-s}(i) applied to the random variable
\begin{equation*} \label{}
X:=cr^{\beta} \left| \int_{\Phi_{\sqrt{R^2-r^2}}} \left(J(\Phi_R,p,q) - J(\Phi_{x,r},p,q)\right) \indc_{\{ \Y_t(x) \leq r \}}  \,dx \right| \indc_{\{ \Y_t \leq R \}}
\end{equation*}
with parameters $s':=s$ and $s/2$ in place of $s$, yield the bound 
\begin{equation*} \label{}
\left| 
 \int_{\Phi_{\sqrt{R^2-r^2}}} \left(J(\Phi_R,p,q) - J(\Phi_{x,r},p,q)\right) \indc_{\{ \Y_t(x) \leq r \}}  \,dx \right| \indc_{\{ \Y_t \leq R \}}
\leq \O_{s} \left( Cr^{-\alpha-\frac12 \beta}  \right).
\end{equation*}
Using the previous line together with the triangle inequality and~\eqref{e.killindcr} again, this time with $s'=2s$ and observing that $t/s > \alpha$, we obtain that  
\begin{equation*} \label{}
\left| 
 \int_{\Phi_{\sqrt{R^2-r^2}}} \left(J(\Phi_R,p,q) - J(\Phi_{x,r},p,q)\right)  \,dx \right| \indc_{\{ \Y_t \leq R \}}
\leq \O_{s} \left( Cr^{-\alpha-\eta}  \right).
\end{equation*}
where 
\begin{equation*} \label{}
\eta := \frac12 \beta  \wedge \left( \frac{t}{s} - \alpha \right).
\end{equation*}
We take $t:= \frac12\left( \alpha s+d \right)$ to get that $\eta>0$ depends only on~$(\ep,s,d,\Lambda)$. 
Finally, using~\eqref{e.indcYsbound} and~\eqref{e.Jvbounded} again we see that 
\begin{multline*} \label{}
\left| 
 \int_{\Phi_{\sqrt{R^2-r^2}}} \left(J(\Phi_R,p,q) - J(\Phi_{x,r},p,q)\right)  \,dx \right| \indc_{\{ \Y_t > R \}}
 \\
 \leq C \indc_{\{ \Y_t > R \}} 
 \leq \O_s \left( CR^{-\frac ts} \right) 
 \leq \O_s \left( Cr^{-\alpha-\eta} \right). 
\end{multline*}
Combining the above yields $\Add(s,\alpha+\eta)$ and completes the argument. 
\end{proof}

\begin{remark}
\label{r.expectationsconverge}
Using the result of Lemma~\ref{l.shortcuts} and quadratic response to estimate $J_1(\Phi_{x,r},p,q)$, and then taking the expectation of the result and combining it with Lemma~\ref{l.correctorbias} yields that 
\begin{equation} 
\label{e.Jbiasestimate}
\left|
 \E\left[ J(\Phi_{x,r},p,q) \right] 
- \left( \frac12 p\cdot \ahom p + \frac12 q\cdot \ahom^{-1} q - p\cdot q \right) 
\right| 
\leq 
Cr^{-2\alpha}. 
\end{equation}
Thus the assumption of~$\Fluc(s,\alpha)$ implies~\eqref{e.Jbiasestimate}. 
\end{remark}

\section{Localization estimates}
\label{s.localization}

In this section, we prove the localization estimate stated in Proposition~\ref{p.improvelocalization}. The proposition roughly states that, under the assumption that $\Fluc(s,\alpha)$ holds, we can compute $J_1(\Phi_r,p,q)$ up to an error of $\O_s(Cr^{-\gamma})$, for some $\gamma > \alpha(1+\delta)$, with only the knowledge of the coefficients restricted to~$B_{r^{1+\delta}}$. To prove this result, we are therefore confronted with what can be thought of as a computational question: 
\begin{quote} 
\emph{What is the best way to compute $J_1(\Phi_r,p,q)$?}
\end{quote}
Since $J_1(\Phi_r,p,q)$ can be computed straightforwardly once~$\A_1$ is known, we can rephrase the task into a more specific problem:
\begin{quote}
\emph{Find an algorithm which accurately computes the correctors in $B_r$ by only looking at the coefficients in~$B_{r^{1+\delta}}$.}
\end{quote}
We have already shown in the proof of Proposition~\ref{p.correctorsbasecase} that the first-order correctors are close to the solutions of the Dirichlet problem with affine data. Therefore we could compute the latter. A similar idea is to consider all solutions in a large ball and choose the ones which are closest in $L^2$ to affine functions. These proposals can be shown to work as long as $\alpha < 1$. However, when $\alpha > 1$ the precision they give is not enough, because we cannot expect the normalized $L^2$ difference between an element of~$\A_1$ and the corresponding affine function to be better than $O(1)$. Therefore, at length scale~$r$, we would not be able to show that the gradient of our approximation and that of the true corrector are closer than $\simeq O(r^{-1})$.

\smallskip

To do better, rather than look for the solution minimizing the distance to an affine function in the $L^2$ norm, we measure the distance in a weaker norm. As in the rest of the chapter, we find it convenient to formulate the ``weak norm'' in terms on spatial averages of the heat kernel. We find the solution $u$ such that the difference in the $L^2$ norm between a given affine function and $u \ast \Phi_{cr}$ is minimized, where $c>0$ is some appropriately small constant. Using some lemmas from the previous section---which assert that the $\Fluc(s,\alpha)$ assumption implies that elements of~$\A_1$ are within  order~$r^{1-\alpha}$ of affine functions after convolution by $\Phi_r$---we are able to show roughly that the resulting solution is within order~$r^{1-\alpha}$ of a true corrector in the ball of radius $r^{1+\delta}$ in this weak norm. Separately, we can use the higher regularity theory, Theorem~\ref{t.regularity}, and a version of the multiscale Poincar\'e inequality (adapted for heat kernels) to show that the difference in the weak norm actually controls the difference in the $L^2$ norm. We can then use the higher regularity theory a second time to improve the error in the approximation as we zoom in to~$B_r$. 

\smallskip

Before we proceed with the proof of Proposition~\ref{p.improvelocalization}, we give a version of the multiscale Poincar\'e inequality which controls the~$L^2$ norm of a function~$w$ with respect to an exponential weight~$\Psi_R$ at scale~$R$ by the $L^2$ norms of $\nabla w \ast \Phi_{\sqrt{t}}$ over all the scales $0 < \sqrt{t} < \sigma R$. To compare the second term on the right side of~\eqref{e.mspoincare.masks} with that of~\eqref{e.multiscalepoincare} for $f=\nabla w(0,\cdot)$, use the change of variables $t = r^2$ to write the former as
\begin{equation*} \label{}
 \int_{0}^{(\sigma R)^2} \int_{\Psi_R} \left| \nabla w(t,y) \right|^2\,dy\,dt
 = 2 \int_0^{\sigma R} r \int_{\Psi_R} \left| \int_{\Phi_{x,r}} \nabla w(y,0)\,dy \right|^2\, dx\,dr.
\end{equation*}
Then observe that a Riemann sum approximation of the right side (using intervals of the form $[3^n,3^{n+1}]$) would look very much like the square of the right side of~\eqref{e.multiscalepoincare}, with averages over triadic cubes replaced by averages with respect to heat kernels with roughly the same length scale.

\begin{lemma}[{Multiscale Poincar\'e inequality, heat kernel version}]
\label{l.mspoincare.masks}
\index{Poincar\'e inequality!multiscale}
Suppose that $w=w(t,x)$ is a solution of the heat equation
\begin{equation*}
\partial_t w - \Delta w = 0 \quad \mbox{in} \ \Rd \times (0,\infty)
\end{equation*}
satisfying $w(0,\cdot) \in L^2(\Psi_R)$ and 
\begin{equation*}
\int_{0}^{R^2} \int_{\Psi_R} \left| w(t,y) \right|^2\,dy\,dt < \infty,
\end{equation*}
where $\Psi_R$ is the function
\begin{equation}
\label{e.PsiR}
\Psi_R(x):= R^{-d} \exp\left( - \frac{|x|}{R} \right).
\end{equation}
Then there exists $C(d,\Lambda)< \infty$ such that, for every $\sigma \in (0,1]$,
\begin{multline}
\label{e.mspoincare.masks}
 \int_{\Psi_R}  \left| w(0,y) \right|^2\,dy \\
 \leq C\int_{\Psi_R} \left| w((\sigma R)^2,y) \right|^2\,dy 
 +C \int_{0}^{(\sigma R)^2} \int_{\Psi_R} \left| \nabla w(t,y) \right|^2\,dy\,dt.
 \end{multline}
\end{lemma}
\begin{proof}
We compute, for any $\ep>0$
\begin{align*}
& \left| \partial_t \int_{\Psi_R} \frac12 \left| w(t,y) \right|^2\,dy \right| 
 = \left| \int_{\Rd} \nabla \left( \Psi_Rw(t,\cdot)\right)(y) \cdot \nabla w(t,y)\,dy \right| \\
& \qquad  \qquad \leq C \int_{\Rd} \Psi_R(y) \left( \frac1\ep\left| \nabla w(t,y) \right|^2 + \ep \left(\frac{\left| \nabla \Psi_R(y) \right| }{\Psi_R(y)}\right)^2 \left|w(t,y) \right|^2\right)\,dy\\
& \qquad \qquad \leq \frac{C}{\ep} \int_{\Psi_R} \left| \nabla w(t,y) \right|^2\,dy + \frac{C\ep}{R^{2}} \int_{\Psi_R}  \left| w(t,y) \right|^2\,dy.
\end{align*}
Integrating with respect to $t$ yields
\begin{align*}
\lefteqn{
\sup_{t\in [0,(\sigma R)^2]}\int_{\Psi_R}  \left| w(\cdot,y) \right|^2\,dy 
} \qquad & \\
& \leq \int_{\Psi_R} \left| w((\sigma R)^2,y) \right|^2\,dy + \int_{0}^{(\sigma R)^2} \left| \partial_t \int_{\Psi_R} \left| w(t,y) \right|^2\,dy \right| \,dt \\
 & \leq \int_{\Psi_R}\left| w((\sigma R)^2,y) \right|^2\,dy +  \frac{C}{\ep}\int_{0}^{(\sigma R)^2} \int_{\Psi_R} \left| \nabla w(t,y) \right|^2\,dy\,dt
  \\ & \qquad 
+ \int_0^{(\sigma R)^2} \frac{C\ep}{R^{2}} \int_{\Psi_R}  \left| w(t,y) \right|^2\,dy \,dt.
\end{align*}
Now taking $\ep =c$ sufficiently small, we can absorb the last term on the right side to obtain
\begin{align*}
 \int_{\Psi_R}  \left| w(0,y) \right|^2\,dy 
 & \leq \sup_{t\in [0,(\sigma R)^2] } \int_{\Psi_R}  \left| w(t,y) \right|^2\,dy \\
 & \leq C \int_{\Psi_R} \left| w\left( (\sigma R)^2),y\right) \right|^2\,dy +  C \int_{0}^{(\sigma R)^2} \int_{\Psi_R} \left| \nabla w(t,y) \right|^2\,dy\,dt. 
\end{align*}
This completes the proof of the lemma. 
\end{proof}

We next specialize the previous lemma to elements of~$\A_k(\Rd)$. We obtain that elements of $\A_k(\Rd)$ behave like normal polynomials in the sense that their spatial averages bound their oscillation. We need this result for $\A_k(\Rd)$ for very large $k\in\N$ in order to transfer this property, by approximation using Theorem~\ref{t.regularity}, to arbitrary solutions.

\begin{lemma} 
\label{l.mspoincareAk}
Fix $s\in (0,d)$, $k\in\N$, and let $\X_s$ denote the random variable in Theorem~\ref{t.regularity}. There exist $\sigma_0(k,d,\Lambda) \in \left(0,\frac12 \right]$ and, for every $\sigma \in (0,\sigma_0]$, 
constants $C(\sigma,k,d,\Lambda)<\infty$ and $\theta(\sigma,k,d,\Lambda) \in \left(0,\frac12 \right]$ such that, for  every $v \in \A_k$ and $r \ge \X_s$,
\begin{equation}  
\label{e.mspoincareAk}
\left\| v \right\|_{L^2(\Psi_r)}^2  \leq C \fint_{B_{r/\theta} } \left| \int_{\Phi_{y,\sigma r}} v(z) \, dz  \right|^2 \, dy \,.
\end{equation}
\end{lemma}
\begin{proof}
Fix $k \in \N$ and $v \in \A_k$. Define
\begin{equation*}
w (t,y):= \int_{\Phi_{y,\sqrt{t}}} v(z)\,dz  \, ,
\end{equation*}
which is the solution of
\begin{equation*}
\left\{
\begin{aligned}
& \partial_t w   - \Delta w = 0 & \mbox{in} & \ \Rd \times (0,\infty), \\
& w  = v & \mbox{on} & \ \Rd \times \{0 \}.
\end{aligned}
\right.
\end{equation*}
By Theorem~\ref{t.regularity}, for every  $r \geq\X_s$, there exists a unique polynomial $q \in \Ahom_k$ such that, for every $S \geq R \geq r$, we have 
\begin{equation} \label{e.upolygrowth1}
\left\| v \right\|_{\underline{L}^2(B_S)} \leq C \left\| q \right\|_{\underline{L}^2(B_S)} \leq C\left(\frac{S}{R}\right)^k \left\| q \right\|_{\underline{L}^2(B_R)} 
\leq C\left(\frac{S}{R}\right)^k \left\| v \right\|_{\underline{L}^2(B_R)} \,.
\end{equation}
In particular,~$v^2(x) \,dx$ is a doubling measure and $v$ has at most polynomial growth. 

\smallskip

Our starting point is Lemma~\ref{l.mspoincare.masks}, which gives, for every $\sigma \in (0,1]$,  
\begin{equation}
\label{e.intintime00}
 \int_{\Psi_{r}}  \left| v (y) \right|^2\,dy  
 \leq C \int_{\Psi_{r}} \left| w((\sigma r)^2,y) \right|^2 \,dy 
 +C \int_{0}^{(\sigma r)^2} \int_{\Psi_{r}} \left| \nabla w(t,y) \right|^2\,dy\,dt.
\end{equation}

\smallskip

\emph{Step 1.} We show that there exists a small $\sigma_0 = \sigma_0(k,d,\Lambda) \in (0,1]$ such that
\begin{equation} \label{e.additivityb100}
C \int_{0}^{(\sigma_0 r)^2} \int_{\Psi_{r}} \left| \nabla w(y,t) \right|^2\,dy\,dt \leq \frac14  \int_{\Psi_r}  \left| v (y) \right|^2\,dy\,,
\end{equation}
so that by absorbing it back onto the left side,~\eqref{e.intintime00} can be improved to 
\begin{equation} \label{e.intintime10}
\int_{\Psi_{r}}  \left| v (y) \right|^2\,dy  \\
 \leq C \int_{\Psi_r } \left| w\left((\sigma r)^2,y\right) \right|^2\,dy 
\end{equation}
for all $\sigma \in (0,\sigma_0]$. 
To prove~\eqref{e.additivityb100}, we first get by H\"older's inequality that 
\begin{equation*} 
\int_{\Psi_{r}} \left| \nabla w(t,y) \right|^2\,dy
  \leq \int_{\Psi_{r}} \int_{\Phi_{y,\sqrt{t}}} \left|\nabla v(z) \right|^2 \, dz \, dy\,.
\end{equation*}
We then notice that the right side can be rewritten with the aid of Fubini's theorem, taking into account the definitions of $\Psi_r$ and $\Phi_{\sqrt{t}}$,  for all $\ep>0$ as
\begin{multline*} 
 \int_{0}^{\ep r^2} \int_{\Psi_{r}} \int_{\Phi_{y,\sqrt{t}}} \left|\nabla v(z) \right|^2  \, dz \, dy \,dt
\\  = \int_{\Psi_{r}} \left|\nabla v (z) \right|^2 \int_{0}^{\ep r^2} \int_{\Phi_{z,\sqrt{t}} } \exp\left(\frac{|z|}{r} - \frac{|y|}{r}\right)   \, dy \, dt \, dz\,.
\end{multline*}
We analyze the integral in the middle. By the triangle inequality,
we obtain
\begin{align*} 
\int_{0}^{\ep r^2} \int_{\Phi_{z,\sqrt{t}} } \exp\left(\frac{|z|}{r} - \frac{|y|}{r}\right)   \, dy \, dt 
&  \leq \int_{0}^{\ep r^2} \int_{ \R^d } \exp\left(\frac{|z-y|}{r}\right) \Phi_{\sqrt{t}}(z-y)  \, dy \, dt 
\\ &   \leq C \int_{0}^{\ep r^2} t^{-\frac d2} \int_{ \R^d }  \exp\left(\frac{|y|}{r} - \frac{|y|^2}{t} \right) \, dy \, dt 
\\ &   \leq C \int_{0}^{\ep r^2} t^{-\frac d2} \left| B_{\frac{2t}{r}}\right| \, dt +   C \int_{0}^{\ep r^2}  \int_{ \Phi_{\sqrt{2t}}}\, dy \, dt 
\\ &  \leq C \ep r^2 \,.
\end{align*}
Therefore, combining the above three displays, we arrive at
\begin{equation} \label{e.mspoincare201}
\int_{0}^{\ep r^2} \int_{\Psi_{r}} \left| \nabla w(t,y) \right|^2\,dy \, dt \leq C\ep r^2 \int_{\Psi_{r}} \left|\nabla v(z) \right|^2 \, dz\,.
\end{equation}
Furthermore, we have the layer-cake formula, for any $g \in L^1(\Psi_r)$ and $V\subset \R^d$, 
\begin{equation} \label{e.psilayer}
\int_{V} \Psi_{r}(z) g(z) \, dz = \frac1{r^{d+1}} \int_{0}^\infty \exp\left(-\frac{\lambda}{r} \right) \int_{V \cap  B_\lambda } g(z) \, dz \, d\lambda\,.
\end{equation}
Using the Caccioppoli estimate and the doubling property~\eqref{e.upolygrowth1}, we deduce that
\begin{equation*} 
 \int_{ B_{R}} \left|\nabla v(y) \right|^2 \, dy  
 \leq  \frac{C}{R^2} \int_{  B_{2R}} \left| v(y) \right|^2 \, dy  
  \leq   \frac{C}{r^2}  \int_{ B_{R}} \left| v(y) \right|^2 \, dy 
\end{equation*}
for any $R \geq r$. Thus the layer-cake formula~\eqref{e.psilayer} yields
\begin{equation*} 
r^2 \int_{\Psi_{r}} \left|\nabla v(z) \right|^2 \, dz  \leq C \int_{\Psi_{r}} \left| v(z) \right|^2 \, dz  \,.
\end{equation*}
Now our claim~\eqref{e.additivityb100} follows from~\eqref{e.intintime00},~\eqref{e.mspoincare201} and the above display provided we take $\ep =\sigma_0^2$ sufficiently small. 

\smallskip

\emph{Step 2.}
We next show that for any $\ep \in (0,1)$ and $\sigma \in (0,\sigma_0]$, there are constants $C(\ep,\sigma,k,d,\Lambda)< \infty$ and $\theta(\ep,\sigma,k,d,\Lambda) \in (0,1)$ such that
\begin{equation} \label{e.psitail00}
 \int_{\Psi_{r}} \left| w((\sigma r)^2,y) \right|^2\,dy 
 \leq  C \fint_{B_{r/\theta}} \left| w((\sigma r)^2,y) \right|^2\,dy +  \ep \int_{\Psi_{r}} \left| v(y) \right|^2 \, dy\,.
\end{equation}
This together with~\eqref{e.intintime10} proves our claim by taking small enough $\ep$, which then also fixes the parameter $\theta$. 
We first decompose the integral on the left as  
\begin{multline}  \label{e.psitail10a} 
 \int_{\Psi_{r}} \left| w((\sigma r)^2,y) \right|^2\,dy  
 \\  \leq C_\theta  \fint_{B_{r/\theta}} \left| w((\sigma r)^2,y) \right|^2\,dy + \int_{\R^d \setminus B_{r/\theta}}  \Psi_{r}(y) \left| \int_{\Phi_{y,\sigma r}} v(z) \, dz \right|^2 \, dy \,.
\end{multline}
As in Step 1, with the aid of the triangle and H\"older's inequalities we obtain
\begin{align}  \label{e.psitail10} 
& \int_{\R^d \setminus B_{r/\theta}}   \Psi_{r}(y) \left| \int_{\Phi_{y,c r}} v(z) \, dz \right|^2 \, dy 
\\  \nonumber & \qquad  
\leq \int_{\R^d \setminus B_{r/\theta}} \Psi_{r}(z) \left| v(z) \right|^2  \int_{\Phi_{z,\sigma r}} \exp\left(\frac{|y-z|}{r}\right)\,dy \, dz 
\\  \nonumber & \qquad  
\leq  C \int_{\R^d \setminus B_{r/\theta}} \Psi_{r}(z) \left| v(z) \right|^2  \, dz \,.
\end{align}
Now the layer-cake formula~\eqref{e.psilayer} and the polynomial growth in~\eqref{e.upolygrowth1} imply
\begin{align} 
\label{e.psitail20}   
\int_{\R^d \setminus B_{r/\theta}} \Psi_{r}(z) \left| v(z) \right|^2  \, dz  & \leq \frac1{r^{d+1}} \int_{r/\theta}^\infty \exp\left(-\frac{\lambda}{r} \right) \int_{ B_\lambda} |v(z)|^2 \, dz \, d\lambda 
\\\nonumber & \leq \frac{C}{r} \fint_{B_r} |v(z)|^2 \, dz  \int_{r/\theta}^\infty \exp\left(-\frac{\lambda}{r} \right)  \left(\frac{\lambda}{r}\right)^k \, d\lambda 
\,.
\end{align}
For any given $\ep \in (0,1)$, we may choose $\theta(\ep,k,d,\Lambda) > 0$ so small that
\begin{equation*} \label{}
\frac1r \int_{r/\theta}^\infty \exp\left(-\frac{\lambda}{r} \right)  \left(\frac{\lambda}{r}\right)^k \, d\lambda = \int_{1/\theta}^\infty \exp\left(-\lambda \right)  \lambda^k \, d\lambda  = \ep.
\end{equation*}
Combining~\eqref{e.psitail10}  and~\eqref{e.psitail20}, we get 
\begin{equation*} \label{}
\int_{\R^d \setminus B_{r/\theta}}   \Psi_{r}(y) \left| \int_{\Phi_{y,\delta r}} v(z) \, dz \right|^2 \, dy  \leq C  \ep \fint_{B_{r} } |v(z)|^2 \, dz
\end{equation*}
with $C$ independent of $\ep$. Inserting this into~\eqref{e.psitail10a}, we obtain that~\eqref{e.psitail00} holds.
\end{proof}

We now proceed with the details of the proof of Proposition~\ref{p.improvelocalization}. As in the previous section, here we fix parameters 
\begin{equation} 
\label{e.salphaparams2}
\delta > 0, \quad 
s\in (0,2], \quad  
\alpha \in \left( 0, \tfrac{d}{s} \right) \cap \left(0,\tfrac d2\right]
\quad \mbox{and} \quad
\gamma \in \left( 0, \left( \alpha(1+\delta)+\delta \right) \wedge \tfrac ds  \right). 
\end{equation}
and suppose that 
\begin{equation} 
\label{e.fluc.ass}
\Fluc(s,\alpha) \quad \mbox{holds.}
\end{equation}
The random variable~$\Y_s$ is defined, for each $s\in (0,d)$, in Definition~\ref{d.Ys}. 

\smallskip

We begin with the construction, for each fixed $r\geq 2$, of a (random) vector subspace of~$\A(B_{r^{1+\delta}})$ denoted by~$\mathcal{V}^\delta_r$ which should well approximate~$\A_1$ in~$B_{r}$ in a sense to be made precise. 

\smallskip

Fix $r\geq 2$ and set $T:=r^{1+\delta}$. Also fix a large integer $m\in\N$ to be selected below---it is chosen in~\eqref{e.choosem} in Step~2 of the proof of Lemma~\ref{l.Vdeltar} and depends only on~$(\delta,\alpha,d)$---and let $\sigma$ and $\theta$ be the constants in the statement of Lemma~\ref{l.mspoincareAk} for $k=m$. We also define an intermediate scale
\begin{equation*} \label{}
S:=r^{1+\delta - \ep},
\end{equation*}
where $\ep>0$ is defined by the relation
\begin{equation*} \label{}
\gamma = \alpha(1+\delta)+\delta - (1+\alpha)\ep.
\end{equation*}
Note that $r < S < T$. The reason for introducing the scale~$S$ is that we have to apply the regularity theory twice: the first time (for $k=m$ large) to transfer the result of Lemma~\ref{l.mspoincareAk} to more general solutions and the second time (with $k=1$) to improve the scaling in the localization estimate. 

\smallskip

Motivated by the right side of~\eqref{e.mspoincareAk}, we define $I:H^1(B_T) \to \R_+$ by
\begin{equation*} \label{}
I\left[  w  \right] := \fint_{B_{S/\theta}} \left| \int_{\Phi_{y,\sigma S}} w(x) \indc_{B_T}(x)\,dx \right|^2 \, dy.
\end{equation*}
For each $i \in \{ 1,\ldots,d\}$, let $w_i \in \A(B_T)$ be the minimizer of the optimization problem 
\begin{equation*} \label{}
\mbox{minimize} \  I\left[ w - \ell_{e_i} \right] \quad \mbox{among} \ w \in \mathcal{K},
\end{equation*}
where $\ell_e(x):=e\cdot x$ and we take the admissible class $\mathcal{K}$ to be the subset of $\A(B_T)$ consisting of functions with gradients under some control:
\begin{equation*} \label{}
\mathcal{K}:= \left\{ w \in \A(B_T) \,:\, \left\| \nabla w \right\|_{\underline{L}^2(B_T)} \leq \overline{C} \right\},
\end{equation*}
where $\overline{C}(s,d,\Lambda)$ is chosen large enough compared to the constants $C$ on the right side of~\eqref{e.correctorminbounds} and of~\eqref{e.vqnondegen} to ensure that 
\begin{equation} 
\label{e.admissibles}
r \geq \Y_s \implies \left\{ v(\cdot,\Phi_{r},0,\ahom e) \,:\, e\in \partial B_1 \right\} \cup \left\{ \overline{\phi}_e\,:\, e\in\partial B_1 \right\} \subseteq \mathcal{K}.
\end{equation}
The functional $I[\cdot]$ is clearly weakly continuous on the convex set $\mathcal{K}$ with respect to the~$H^1(B_T)$-norm, and the existence of a minimizer follows. The uniqueness of the minimizer is not important, in the case of multiple minimizers we can for instance select the one with the smallest $L^2(B_T)$-norm. It is not hard to show that, by~\eqref{e.gradientbasecase} and~\eqref{e.admissibles}, that 
\begin{equation} 
\label{e.basecaseLI}
r\geq \Y_s \implies \left\{ w_1,\ldots,w_d \right\} \ \  \mbox{is linearly independent.}
\end{equation}
Finally, we define
\begin{equation*} \label{}
\mathcal{V}^\delta_r:= \mathrm{span}\left\{ w_1,\ldots,w_d \right\}.
\end{equation*}
It is clear that each of the functions $w_i$ is $\F(B_T)$--measurable, and thus 
\begin{equation*} \label{}
\mathcal{V}^\delta_r 
\quad \mbox{is $\F(B_T)$--measurable.}
\end{equation*}
The main part of the argument for Proposition~\ref{p.improvelocalization} is formalized in the following lemma, which states roughly that every element of $\A_1$ can be approximated by an element of $\mathcal{V}^\delta_r$, up to an error of order~$\simeq r^{-\gamma}$. We let~$\Phi_r^{\delta}$ denote the truncated heat kernel given by
\begin{equation} 
\label{e.truncatedheatkernel}
\Phi_r^{\delta}(x):= \Phi_r(x) \indc_{x\in B_T}. 
\end{equation}

\begin{lemma}
There exists a constant $C(\gamma,\delta, s, \alpha, d,\Lambda)<\infty$ such that 
\label{l.Vdeltar}
\begin{equation} 
\label{e.Vdeltar}
\sup_{u \in \A_1(\Phi_r)} 
\inf_{w \in \mathcal{V}^{\delta}_r}
\left\| \nabla u - \nabla w \right\|_{L^2(\Phi_{r}^{\delta})}
+
\sup_{w \in \mathcal{V}^{\delta}_r(\Phi_r) } 
\inf_{u \in \A_1}
\left\| \nabla u - \nabla w \right\|_{L^2(\Phi_{r}^{\delta})}
\leq 
\O_s\left( Cr^{-\gamma} \right).
\end{equation}
\end{lemma}
\begin{proof}
Since $\gamma < \frac ds$, we have by~\eqref{e.indcYsbound} that 
\begin{multline*} \label{}
\sup_{u \in \A_1(\Phi_r)} 
\inf_{w \in \mathcal{V}^{\delta}_r}
\left\| \nabla u - \nabla w \right\|_{L^2(\Phi_{r}^{\delta})} \indc_{ \{ \Y_s \ge r\}}
\\
\leq \sup_{u \in \A_1(\Phi_r)} \left\| \nabla u \right\|_{L^2(\Phi_{r}^{\delta})} \indc_{ \{ \Y_s \ge r\}}
\leq \indc_{ \{ \Y_s \ge r\}} 
\leq  \O_s\left( Cr^{-\gamma} \right),
\end{multline*}
and a similar inequality holds if we switch the roles of $\A_1$ and $\mathcal{V}^\delta_r$. Therefore for the rest of the argument we work on the event~$\{ \Y_s < r\}$ while dropping the~$\indc_{\{ \Y_s < r \}}$ from the expressions to lighten the notation. For each~$j\in\{ 1,\ldots,d\}$, we denote~$v_j:= v(\cdot,\Phi_r,0,\ahom e_j)$.

\smallskip

\emph{Step 1.} We show that 
\begin{equation} 
\label{e.mintocorrect}
\sup_{j\in\{1,\ldots,d\}} I \left[ w_j -  v_j \right] \leq \O_{s/2} \left( CS^{2-2\alpha} \right).
\end{equation}
Fix $j\in\{1,\ldots,d\}$. Observe that~\eqref{e.admissibles} and the definition of $w_j$ imply that
\begin{equation} 
\label{e.comparemintocorrect}
I\left[ w_j - \ell_{e_j} \right] \leq I\left[ v_j - \ell_{e_j} \right]. 
\end{equation}
We may assume that the additive constant for $v_j$ is chosen so that
\begin{equation} 
\label{e.normforpoinc}
\fint_{B_{S/\theta}} v_j \ast \Phi_{\sigma S} = 0. 
\end{equation}
We now use the assumption that $\Fluc(s,\alpha)$ holds, which allows us to apply Lemma~\ref{l.shortcuts}---actually it is best to use directly~\eqref{e.giveittomefluc} and~\eqref{e.identifyLxr} from its proof---which, together with the Poincar\'e inequality and~\eqref{e.normforpoinc}, yields
\begin{equation} 
\label{e.correctorsareclose}
I \left[ v_j  - \ell_{e_j} \right] 
\leq \O_{s/2} \left( CS^{2-2\alpha} \right).
\end{equation}
Indeed, we compute
\begin{align*}
I \left[ v_j - \ell_{e_j} \right]
 & = \fint_{B_{S/\theta}} \left| \int_{\Phi_{y,\sigma S}} \left( v_j(x) - \ell_{e_j}(x) \right) \,dx \right|^2 \,dy \\
&  \leq CS^2 \fint_{B_{S/\theta}} \left| \int_{\Phi_{y,\sigma S}} \left( \nabla v_j(x) - e_j \right) \,dx \right|^2 \,dy \\
& \leq CS^2 \int_{\Phi_{\sqrt{S^2-(\sigma S)^2}}}  \left|  \int_{\Phi_{y,\sigma S}} \left( \nabla v_j(x) - e_j \right) \,dx \right|^2 \,dy \\
& \leq \O_{s/2}\left( C  S^{2-2\alpha}\right)\,.
\end{align*}
We deduce from~\eqref{e.comparemintocorrect} and~\eqref{e.correctorsareclose} that 
\begin{equation} 
\label{e.minimizersarecloser}
I \left[ w_j - \ell_{e_j} \right]  \leq \O_{s/2} \left( CS^{2-2\alpha} \right).
\end{equation}
The triangle inequality then gives us
\begin{equation*} \label{}
I \left[ w_j -  v_j \right] \leq \O_{s/2} \left( CS^{2-2\alpha} \right).
\end{equation*}
Summing over $j\in \{1,\ldots,d\}$ yields~\eqref{e.mintocorrect}.

\smallskip

\emph{Step 2.} We upgrade~\eqref{e.mintocorrect} to convergence in the strong $L^2$-norm by using the regularity theory and~Lemma~\ref{l.mspoincareAk}. The precise claim is that 
\begin{equation} 
\label{e.mintocorrectL2}
\sup_{j\in\{1,\ldots,d\}} \left\| w_j -  v_j \right\|_{\underline{L}^2(B_{S/\theta})} \leq \O_{s} \left( CS^{1-\alpha} \right).
\end{equation}
According to $w_j \in \mathcal{K}$ and $r\geq \Y_s$, we have that, for every $r' \in [r,T]$,
\begin{equation*} \label{}
\inf_{u \in \A_m} \left\| w_j - u \right\|_{\underline{L}^2(B_{r'})} 
\leq  C\left( \frac{r'}{T} \right)^{m+1}  \left\| w_j  \right\|_{\underline{L}^2(B_{r'})} 
\leq   C\left( \frac{r'}{T} \right)^{m+1}T. 
\end{equation*}
This implies that
\begin{equation*} \label{}
\inf_{u \in \A_m} \left\| \nabla w_j - \nabla u \right\|_{\underline{L}^2(B_{S})} 
\leq   C\left( \frac{S}{T} \right)^{m} = C r^{-\ep m} \leq CS^{-\ep m/ (1+\delta)}. 
\end{equation*}
We now select 
\begin{equation} 
\label{e.choosem}
m:=\left\lceil d  (1+\delta) \ep^{-1} \right\rceil
\end{equation}
so that $\ep m /  (1+\delta) \geq d > \alpha$. Observe that $m$ depends only on~$(\delta,\alpha,d)$, as claimed above. Taking $\tilde{u}_j\in \A_m$ to achieve the infimum in the previous display, we have shown that 
\begin{equation*} \label{}
 \left\| \nabla w_j - \nabla \tilde{u}_j \right\|_{\underline{L}^2(B_{S})} 
\leq CS^{-\alpha}. 
\end{equation*}
It follows from the triangle inequality and~\eqref{e.mintocorrect} that 
\begin{equation*} \label{}
I\left[ \tilde{u}_j - v_j \right] \leq \O_{s/2} \left( CS^{2-2\alpha} \right). 
\end{equation*}
Since $\tilde{u}_j - v_j  \in \A_m$, we may apply Lemma~\ref{l.mspoincareAk} to obtain
\begin{equation*} \label{}
\left\| \tilde{u}_j - v_j \right\|_{L^2(B_S)}^2 \leq C \left\| \tilde{u}_j - v_j \right\|_{L^2(\Psi_S)}^2  \leq CI\left[ \tilde{u}_j - v_j \right] \leq \O_{s/2} \left( CS^{2-2\alpha} \right). 
\end{equation*}
By the triangle inequality again, we obtain~\eqref{e.mintocorrectL2}.

%

\smallskip

\emph{Step 3.} We apply the regularity theory to obtain better estimates on balls $B_{\rho}$ with $\rho \in \left[ r , S \right]$. The claim is that, for every $j \in\{1,\ldots,d\}$,
\begin{equation} 
\label{e.approxdownscales}
\inf_{u \in \A_1} \sup_{\rho \in  \left[ r , S/8\theta \right]} 
\frac{1}{\rho}  \left\| \nabla w_j - \nabla u \right\|_{\underline{L}^2(B_{\rho})}
\leq 
\O_s\left( CS^{-1-\alpha} \right). 
\end{equation}
Fix $j$. Applying statement~$\mathrm{(iii)}_1$ of Theorem~\ref{t.regularity} and using that $r > \Y_s$, we deduce the existence of $u_j \in \A_1$ satisfying, for every $\rho \in \left[ r , S/8\theta \right]$ and $a\in\R$,
\begin{equation*} \label{}
\left\| w_j - u_j \right\|_{\underline{L}^2(B_{2\rho})} 
\leq C \left( \frac{\rho}{S} \right)^2 \left\| w_j - v_j - a\right\|_{\underline{L}^2(B_{S/\theta})}.
\end{equation*}
By the Caccioppoli and Poincar\'e inequalities, this yields
\begin{equation*} \label{}
\left\| \nabla w_j - \nabla u_j \right\|_{\underline{L}^2(B_{\rho})} 
\leq C \left( \frac{\rho}{S} \right) \left\| \nabla w_j - \nabla  v_j \right\|_{\underline{L}^2(B_{S/\theta})}.
\end{equation*}
Applying~\eqref{e.mintocorrectL2}, we deduce that, for every $\rho\in [r,S]$, 
\begin{equation*} \label{}
\frac{1}{\rho}  \left\| \nabla w_j - \nabla u_j \right\|_{\underline{L}^2(B_{\rho})} 
\leq \frac{C}S  \left\| \nabla w_j - \nabla  v_j \right\|_{\underline{L}^2(B_{S/\theta})}
\leq \O_s\left( CS^{-1-\alpha} \right). 
\end{equation*}
Taking the supremum over $\rho\in [r,S/8\theta]$, we obtain~\eqref{e.approxdownscales}.

\smallskip

\emph{Step 4.} We show that, for every $j\in \{1,\ldots,d\}$, 
\begin{equation} 
\label{e.approxdownscalesPhi}
 \left\| \nabla w_j - \nabla u_j \right\|_{{L}^2(\Phi^{\delta}_r)}
\leq 
\O_s\left( Cr^{-\gamma} \right). 
\end{equation}
By~\eqref{e.approxdownscales} and a layer-cake approximation,
\begin{align}
\left\| \nabla w_j - \nabla u_j \right\|_{\Phi^{\delta}_r}^2
& 
\leq C \int_r^T \left( \frac\rho r \right)^d \exp\left( -c \left( \frac{\rho}{r} \right)^2 \right) \left\| \nabla w_j - \nabla u_j \right\|_{\underline{L}^2(B_\rho)}^2 \frac{d\rho}\rho 
\\ & 
\leq \O_{s/2}\left( C\left( \frac{r}{S} \right)^2 S^{-2\alpha} \right). 
\end{align}
Using the definition of $S$, we see that 
\begin{equation*} \label{}
\left( \frac{r}{S} \right)^2 S^{-2\alpha} = r^{-2\gamma}.
\end{equation*}
This completes the proof of~\eqref{e.approxdownscalesPhi}. 

\smallskip

\emph{Step 5.} The conclusion. Let $\mathcal{T}: \mathcal{V}^\delta_r/\R \to \A_1/\R$ be the linear mapping for which $\mathcal{T}(w_j) = u_j$ for every $j\in\{1,\ldots,d\}$. By linearity, we deduce from~\eqref{e.approxdownscalesPhi} that, for any $w\in \mathcal{V}^\delta_r$, 
\begin{equation*} \label{}
\inf_{u\in\A_1}
 \left\| \nabla w - \nabla u \right\|_{{L}^2(\Phi^{\delta}_r)}
 \leq
 \left\| \nabla w - \nabla \mathcal{T}(w) \right\|_{{L}^2(\Phi^{\delta}_r)}
\leq 
\O_s\left( C \left\| \nabla w \right\|_{ {L}^2(\Phi^{\delta}_r)} r^{-\gamma} \right). 
\end{equation*}
This gives us an estimate for the second term on the left of~\eqref{e.Vdeltar}. 
Since by~\eqref{e.basecaseLI}, we have that $r\geq \Y_s$ implies  $\dim( \mathcal{V}^\delta_r/\R) = \dim( \A_1/\R) = d$, we deduce that $\mathcal{T}$ is within $\O_s(Cr^{-\gamma})$ of an isometry. If a linear operator on finite dimensional spaces is within~$h$ of an isometry, its inverse is within $O(h)$ of an isometry as well. This observation yields the estimate for the first term on the left of~\eqref{e.Vdeltar}. 
\end{proof}

Now that we have approximated~$\A_1$ by a local vector space $\mathcal{V}^\delta_r$, we define the localized version of $J(\Phi_r,p,q)$ by
\begin{equation*} \label{}
J^{(\delta)}_1( 0,r,p,q)
:= 
\sup_{w\in \mathcal{V}^\delta_r} 
\int_{\Phi_r^{\delta} }
\left( -\frac12\nabla w \cdot \a\nabla w - p\cdot \a\nabla w + q\cdot \nabla w \right),
\end{equation*}
where $\Phi_r^{\delta}$ is the truncated heat kernel defined in~\eqref{e.truncatedheatkernel}. 
For future reference, we note that, for every $r\geq 1$ and $s' < d$, 
\begin{equation} 
\label{e.tailchopping}
\sup_{u \in \A_1(\Phi_r)} \int_{\Phi_r - \Phi_r^{\delta}}\left| \nabla u \right|^2
\leq C \wedge \O_{s'} \left( \exp\left( -cr^{-2\delta} \right)  \right).
\end{equation}
This is easy to show from~\eqref{e.correctorminbounds},~\eqref{e.gradientbasecase} and~\eqref{e.correctorsballz}. It is clear that 
\begin{equation*} \label{}
J^{(\delta)}_1( 0,r,p,q) 
\quad
\mbox{is $\F(B_T)$-measurable.}
\end{equation*}
To complete the proof of Proposition~\ref{p.improvelocalization}, it remains to prove the following lemma. 
\begin{lemma}
\label{l.J1approx}
Recall that $\gamma \in \left( 0, \left( \alpha(1+\delta)+\delta \right) \wedge \frac ds \right)$. There exists a constant $C(\gamma,\delta, s, \alpha, d,\Lambda)<\infty$ such that, for every $p,q\in B_1$,
\begin{equation*} \label{}
\left| J_1(\Phi_{r},p,q) - J^{(\delta)}_1(0,r,p,q) \right|
\leq
\O_s\left(Cr^{-\gamma} \right).
\end{equation*}
\end{lemma}
\begin{proof}
The main principle underlying the proof is that if we compute the minimum of a uniformly convex function on a finite dimensional subspace and then perturb this subspace by~$h$ and recompute the minimum, we find that the minimum is perturbed by~$O(h)$.

\smallskip

Fix $p,q\in B_1$. Denote by $v^{(\delta)}(\cdot,0,r,p,q)\in \mathcal{K}$ the maximizer in the definition of $J_1^{(\delta)}(0,r,p,q)$. Write $v= v(\cdot,\Phi_r,p,q)$ and $v^{(\delta)} = v^{(\delta)}(\cdot,0,r,p,q)$ for short. By Lemma~\ref{l.Vdeltar}, we can find $u \in \A_1$ and $w \in V^{\delta}_r$ such that 
\begin{equation*} \label{}
\left\| \nabla v^{(\delta)} - \nabla u \right\|_{L^2(\Phi_{r}^{\delta})}
+\left\| \nabla v - \nabla w \right\|_{L^2(\Phi_{r}^{\delta})}
\leq \O_{s}\left( Cr^{-\gamma} \right).
\end{equation*}
Here we used~\eqref{e.Jvbounded} and the fact that a similar bound holds for $v^{(\delta)}$ by the same argument. Using also~\eqref{e.tailchopping}, we deduce that 
\begin{align*}
\lefteqn{
J_1(\Phi_r,p,q) 
} \qquad & 
\\ & 
\geq \int_{\Phi_r} 
\left( -\frac12\nabla u \cdot \a\nabla u - p\cdot \a\nabla u + q\cdot \nabla u \right)
\\ & 
\geq  \int_{\Phi_r^{\delta} } 
\left( -\frac12\nabla u \cdot \a\nabla u - p\cdot \a\nabla u + q\cdot \nabla u \right) 
-  \O_{s} \left( \exp\left( -cr^{-2\delta} \right)\right)
\\ &
\geq  \int_{\Phi_r^{\delta} } 
\left( -\frac12\nabla v^{(\delta)} \cdot \a\nabla v^{(\delta)} - p\cdot \a\nabla v^{(\delta)} + q\cdot \nabla v^{(\delta)} \right) 
- \O_{s} \left( Cr^{-\gamma} \right)
\\ & 
= J_1^{(\delta)} (0,r,p,q) 
- \O_{s} \left( Cr^{-\gamma} \right)
\end{align*}
and, similarly,
\begin{align*}
J_1^{(\delta)} (0,r,p,q) 
& 
\geq \int_{\Phi_r^{\delta} } 
\left( -\frac12\nabla w \cdot \a\nabla w - p\cdot \a\nabla w + q\cdot \nabla w \right)
\\ & 
\geq  \int_{\Phi_r^{\delta} } 
\left( -\frac12\nabla v \cdot \a\nabla v - p\cdot \a\nabla v + q\cdot \nabla v \right) 
- \O_{s} \left( Cr^{-\gamma} \right)
\\ &
\geq  \int_{\Phi_r} 
\left( -\frac12\nabla v \cdot \a\nabla v - p\cdot \a\nabla v + q\cdot \nabla v \right) 
- \O_{s} \left( Cr^{-\gamma} \right)
\\ & 
= J(\Phi_r,p,q) 
- \O_{s} \left( Cr^{-\gamma} \right).
\end{align*}
This completes the proof. 
\end{proof}

\section{Fluctuation estimates}
\label{s.fluctuations}

In this section, we prove Proposition~\ref{p.improvefluc}. The argument is a variation of a classical proof of the central limit theorem for sums of bounded and 
i.i.d.~(independent and identically distributed) random variables, which we now summarize. For a random variable $X$, we call the mapping 
\index{Laplace transform}
\begin{equation*}  
\lambda \mapsto \log \E \Ll[ \exp \Ll( \lambda X \Rr)  \Rr] 
\end{equation*}
its \emph{Laplace transform}. We first observe that the Laplace transform of a single centered random variable is close to a parabola near $\lambda = 0$, and then use this fact to monitor the behavior of the Laplace transform of the rescaled sums. Passing to larger and larger scales, these Laplace transforms approach the parabola on larger and larger intervals. This implies that the rescaled sums converge in law to the unique distribution whose Laplace transform is this parabola, that is, a Gaussian. (See also the paragraph below \eqref{e.def.mclN} for a more precise discussion.)

\smallskip

It is convenient to introduce a new notation for \emph{centered} random variables.
For every random variable $X$, $s \in (1,2]$ and $\theta \in (0,\infty)$, we write
\begin{equation}
\label{e.def.Ost}
X = \bar \O_s(\theta)
\end{equation}
to mean that
\begin{equation*}
\forall \lambda \in \R, \quad \log \E \left[ \exp\left( \lambda \theta^{-1}  X \right) \right] \leq  \lambda^2 \vee   | \lambda|^{\frac s {s-1}}  .
\end{equation*}
It is classical to verify that the assumption of \eqref{e.def.Ost} implies that $X$ is centered. Moreover, there exists~$C(s) < \infty$ such that, for  every random variable $X$, 
\begin{equation*}  
X = \O_s(1)  \ \  \text{and}  \ \  \E[X] = 0 
\ \  \implies  \ \
X = \bar \O_s(C).
\end{equation*}
See Lemma~\ref{l.bigO.barO} for a proof of this fact.
The key ingredient of the proof of Proposition~\ref{p.improvefluc} is the simple observation that a sum of $k$ independent $\bar \O_s(\theta)$ random variables is $\bar \O_s(\sqrt{k} \, \theta)$, in agreement with the scaling of the central limit theorem: see Lemmas~\ref{l.barO} and~\ref{l.barO.boxes} in Appendix~\ref{a.bigO}. 

\smallskip

We begin with a weaker form of Proposition~\ref{p.improvefluc}, which is sufficient to show a weaker form of Theorem~\ref{t.additivity} with \emph{almost optimal} (rather than \emph{optimal}) exponents. Although we will not use this lemma, its proof is a simpler version of the argument used below to prove Proposition~\ref{p.improvefluc} and thus easier to understand on a first reading.

\begin{lemma}
\label{l.weak.fluc}
For every $\al \in (0,\infty)$, $\beta \in \Ll( 0, \al \wedge \frac d 2 \Rr)$ and $\de \in \Ll( 0, \frac {(\al - \be)(d-2\be)}{\be d} \Rr)$, we have
\begin{equation*}  
\Add(s,\al) \ \ \text{and} \ \ \Loc(s,\al,\de) \quad \implies \quad \Fluc(s,\be).
\end{equation*}
\end{lemma}
\begin{proof}
Throughout the proof, the value of the constant $C(\al,s,\de,\be,d,\Lambda) < \infty$ may change from place to place. We fix $p,q \in B_1$, and aim to show that for every $z \in \Rd$ and $r \ge 1$,
\begin{equation}  
\label{e.weak.fluc1}
J_1(\Phi_{z,r},p,q) = \E[J_1(\Phi_{z,r},p,q)] + \O_s \Ll( C r^{-\be} \Rr) .
\end{equation}
Denote by $\Jd(z,r,p,q)$ the random variables provided by the assumption of $\Loc(s,\al,\de)$, and set
\begin{equation}  
\label{e.def.tJd}
\tJd(z,r) := \Jd(z,r,p,q) - \E \Ll[ \Jd(z,r,p,q) \Rr] ,
\end{equation}
where we dropped the dependency on the fixed vectors $p,q, \in \Rd$ in the notation for concision. By the construction of $\Jd$, in order to prove \eqref{e.weak.fluc1}, it suffices to show that for every $z \in \Rd$ and $r \ge 1$,
\begin{equation}
\label{e.weak.fluc2}
\tJd(z,r) = \O_s \Ll( C r^{-\be} \Rr) .
\end{equation}
By the assumption of $\Add(s,\al)$ and $\Loc(s,\al,\de)$, we have, for every $z \in \Rd$ and $R \ge r \ge 1$,
\begin{equation}  
\label{e.additivity.unchopped}
\tJd(z,R) = \int_{\Phi_{z,\sqrt{R^2 - r^2}}} \tJd(\cdot,r) + \O_s \Ll( C r^{-\al} \Rr) .
\end{equation}
For $R_1, \CC \ge 1$, we denote by $\msf A(R_1,\CC)$ the statement that for every $z \in \Rd$ and $r \in [1,R_1]$,  
\begin{equation*}  
\tJd(z,r) = \O_s \Ll( \CC r^{-\be} \Rr) .
\end{equation*}
Since $0 \le J_1(\Phi_{z,r},p,q) \le C$, we can assume $|\tJd(z,r)| \le C$, and thus for every $R_1 \in [1,\infty)$, we can find a constant $\CC \in [1,\infty)$ such that $\msf A(R_1,\CC)$ holds. In order to show \eqref{e.weak.fluc2} and complete the proof, it thus suffices to show that for every $\CC$ sufficiently large and $R_1$ sufficiently large,
\begin{equation}
\label{e.fluct-induction}
\msf A(R_1, \CC) \quad \implies \quad  \msf A(2R_1, \CC).
\end{equation}
Indeed, this yields the existence of a constant $\CC$ such that $\msf A(R_1,\CC)$ holds for every $R_1$ sufficiently large, which means that \eqref{e.weak.fluc2} holds.

We thus assume $\msf A(R_1,\CC)$, fix $R \in (R_1,2R_1]$ and $z \in \Rd$ and aim to show that
\begin{equation*}  
\tJd(z,R) = \bar \O_s \Ll( C \CC R^{-\be - \eps} \Rr)  + \bar \O_s \Ll( C R^{-\be} \Rr) ,
\end{equation*}
for some exponent $\eps(\al,\be, d) > 0$ (recall that the constant $C$ is \emph{not} allowed to depend on $\CC$). 
Without loss of generality, we assume that $z = 0$, fix 
\begin{equation}
\label{e.choice.of.r}
r := R^{\frac \be \al},
\end{equation}
and rewrite the additivity property \eqref{e.additivity.unchopped} as
\begin{equation}  
\label{e.additivity-again}
\tJd(0,R) = \int_{\Phi_{\sqrt{R^2 - r^2}}} \tJd(\cdot,r) + \bar \O_s \Ll( C R^{-\be} \Rr) ,
\end{equation} 
We decompose the integral on the right side of \eqref{e.additivity-again} as
\begin{equation*}
\sum_{y \in r^{1+\de} \Z^d} \int_{\Phi_{\sqrt{R^2 - r^2}}} \tJd(\cdot,r) \1_{\cu_{r^{1+\de}}(y)}.
\end{equation*}
For $R_1$ sufficiently large, our choice of $r$ in \eqref{e.choice.of.r} ensures that $r \in [1,R_1]$, and therefore, by the induction hypothesis,
\begin{equation*}  
\forall x \in \Rd, \quad \tJd(x,r) = \bar \O_s \Ll( \CC r^{-\be} \Rr) .
\end{equation*}
By Lemmas~\ref{l.sum.barO} and \ref{l.barO.boxes}, we deduce that
\begin{equation*}  
\sum_{y \in r^{1+\de} \Z^d} \int_{\Phi_{\sqrt{R^2 - r^2}}} \tJd(\cdot,r) \1_{\cu_{r^{1+\de}}(y)} = \bar \O_s \Ll(C \CC r^{-\be} \Ll[\sum_{y \in \mcl Z} \Ll(\int_{\Phi_{\sqrt{R^2 - r^2}}} \1_{\cu_{r^{1+\de}}(y)} \Rr)^2 \Rr]^\frac 1 2\Rr) .
\end{equation*}
Moreover,
\begin{equation*}  
\sum_{y \in r^{1+\de} \Z^d} \Ll(\int_{\Phi_{\sqrt{R^2 - r^2}}} \1_{\cu_{r^{1+\de}}(y)} \Rr)^2  \le C {r^{d(1+\de)}}{R^{-d}},
\end{equation*}
so we have shown that
\begin{equation*}  
\int_{\Phi_{\sqrt{R^2 - r^2}}} \tJd(\cdot,r) = \bar \O_s \Ll( C \CC r^{-\be + \frac d 2(1+\de)}R^{-\frac d 2}  \Rr) .
\end{equation*}
Recalling our choice of $r$ in \eqref{e.choice.of.r}, we see that
\begin{equation*}  
r^{-\be + \frac d 2(1+\de) }R^{-\frac d 2}  = R^{-\Ll[\frac d 2 \Ll( 1-\frac \be \al \Rr) + \be \, \frac \be \al  \Rr] + \frac d 2 \de},
\end{equation*}
and since $\be < \frac d 2$, the exponent between square brackets is larger than $\be$.
It thus suffices to check that $\de$ is sufficiently small to obtain the desired result. That is, we need to verify that
\begin{equation*}  
\frac d 2 \Ll( 1 - \frac \be \al \Rr) + \be \, \frac \be \al  - \frac  d 2 \de > \be.
\end{equation*}
This condition is equivalent to our assumption $\de < \frac {(\al - \be)(d-2\be)}{\be d}$, so the proof is complete.
\end{proof}
We now turn to the proof of Proposition~\ref{p.improvefluc}, which is a refined version of Lemma~\ref{l.weak.fluc} above. To guide our intuition, we can compare the behavior of $J_1(\Phi_{z,r},p,q)$ with that of a bounded random field with finite range of dependence, integrated against $\Phi_{z,r}$. The additivity and localization properties of the former can only be less good than those of the latter. For the latter quantity, the significant overlap between the heat kernel masks makes it difficult to reach the CLT scaling, if we insist on manipulating only averages of the random field against $\Phi_{z,r}$. If we want to reach the CLT scaling for this simpler quantity via a renormalization-type argument, it would be more appropriate to consider averages of the original random fields against masks that are constant on a large cube, and then have a much smaller transition layer to $0$ occupying an asymptotically negligible volume. This would indeed allow us to decompose larger-scale quantities into a sum of independent quantities on a smaller scale, plus negligible boundary regions.

\smallskip

In our context, we implement this idea of using ``masks with a mesoscopic transition layer'' by considering quantities of the form
\begin{equation}  
\label{e.meso.mask}
\fint_{\cu_R} J_1(\Phi_{z,r},p,q) \, dz
\end{equation}
for $r \le R$ relatively close to $R$ (on a logarithmic scale). For the toy example mentioned above, this would amount to considering convolutions of the original random field against the function
\begin{equation*}  
\Phi(r,\cdot) \ast \frac{\1_{\cu_R}}{|\cu_R|},
\end{equation*}
which indeed is essentially equal to $1$ in $\cu_R$ away from the boundary, and essentially vanishes to $0$ over a width of the order of $r$ around the boundary of $\cu_R$. In our context, it is of course essential that we also make sure that $r$ is sufficiently large in \eqref{e.meso.mask}, since additivity is not exact.

\begin{proof}[Proof of Proposition~\ref{p.improvefluc}]
Since we assume $\al > \be(1+\de)$, there exist exponents $\eta_1 < \eta_2 < 1$ such that 
\begin{equation}
\label{e.cond.eta}
\eta_1 \al > \be \quad \text{and} \quad \eta_2(1+\de) < 1.
\end{equation}
These exponents will serve to describe the interval of allowed values for $r$ in quantities of the form \eqref{e.meso.mask}. Throughout the proof, the value of the constant $C(\al,s,\be,\de,\eta_1,\eta_2,d,\Lambda) < \infty$ and of the exponent $\eps(\al,\be,\de, \eta_1,\eta_2,d) > 0$ may change from place to place without further mention. As in the proof of Lemma~\ref{l.weak.fluc}, we fix $p,q \in B_1$ and define $\tJd(z,r)$ by \eqref{e.def.tJd}. Since $\al > \beta$, in order to prove $\Fluc(s,\be)$, it suffices to show that for every $R \ge 1$ and $z \in \Rd$, we have
\begin{equation}  
\label{e.suffice.fluc}
\tJd(z,R) = \O_s \Ll( C R^{-\be} \Rr) .
\end{equation}
Recall from \eqref{e.additivity.unchopped} and Lemma~\ref{l.bigO.barO} that for every $z \in \Rd$ and $R \ge r \ge 1$, we have the additivity property
\begin{equation}  
\label{e.additivity.tJd}
\tJd(z,R) = \int_{\Phi_{z,\sqrt{R^2 - r^2}}} \tJd(\cdot,r) + \bar \O_s \Ll( C r^{-\al} \Rr) .
\end{equation}
For $R_1, \CC \ge 1$, we let $\msf A(R_1,\CC)$ denote the statement that for every $R \in [1,R_1]$, $r \in [R^{\eta_1}, R^{\eta_2}]$, $z \in \Rd$, and every (deterministic) function $f \in L^\infty(\Rd)$ satisfying $\|f\|_{L^\infty(\Rd)}\le 1$, we have
\begin{equation*}  
\fint_{\cu_R(z)} f \,\tJd(\cdot,r) = \bar \O_s \Ll( \CC R^{-\beta} \Rr) .
\end{equation*}
We focus on showing that
\begin{equation}
\label{e.induc.scale}
\mbox{there exists $\CC < \infty$ such that, for every $R \ge 1$, $\msf A(R,\CC)$ holds.}
\end{equation}
Since $0 \le J(z,r,p,q) \le C$, we can assume that $|\tJd(z,r)| \le C$. It is therefore clear that for every given $R_1$, we can find $\CC$ sufficiently large that $\msf A(R_1,\CC)$ holds. In order to prove \eqref{e.induc.scale}, it therefore suffices to show that for every $\CC$ and $R_1$ sufficiently large,
\begin{equation}
\label{e.induction-fluct}
\msf A(R_1, \CC) \quad \implies \quad \msf A(2R_1, (1+R_1^{-\eps})\CC).
\end{equation}
Hence, we assume that $\msf A(R_1,\CC)$ holds for possibly large values of $R_1$ and $\CC$, give ourselves $R \in (R_1,2R_1]$, $r \in [R^{\eta_1}, R^{\eta_2}]$, $z \in \Rd$ and a function $f \in L^\infty(\Rd)$ satisfying $\|f\|_{L^\infty(\Rd)} \le 1$, and aim to show that 
\begin{equation}
\label{e.rewrite.induc}
\fint_{\cu_R(z)} f \, \tJd(\cdot,r) = \bar \O_s \Ll( \CC (1+R^{-\eps})R^{-\beta} \Rr) + \bar \O_s \Ll( C \CC R^{-\beta - \eps} \Rr)  .
\end{equation}
For notational convenience, we fix $z = 0$ from now on. Roughly speaking, in order to reduce the proof of \eqref{e.rewrite.induc} to an application of the induction hypothesis, we will first lower the value of $r$ to
\begin{equation}
\label{e.def.r1}
r_1 := \Ll(\frac R 2 \Rr)^{\eta_1}
\end{equation}
in the expression on the left side of \eqref{e.rewrite.induc}, using additivity; and then lower the value of $R$ to $\frac R 2$, by chopping the integral. We do each of these operations in the next two steps.

\smallskip

\emph{Step 1.}
For $r_1$ as defined in \eqref{e.def.r1}, define
\begin{equation*}  
g := \Ll(f \1_{\cu_R} \Rr)\ast \Phi\Ll(\cdot, {r^2 - r_1^2}\Rr).
\end{equation*}
Note that $\|g\|_{L^\infty(\Rd)} \le 1$. In this step, we show that
\begin{equation}  
\label{e.replace.r}
\fint_{\cu_R} f \, \tJd(\cdot,r) = R^{-d} \int_{\Rd} g \, \tJd(\cdot,r_1) + \bar \O_s \Ll( C R^{-\be - \eps} \Rr) .
\end{equation}
By \eqref{e.additivity.tJd}, we have, for every $x \in \Rd$,
\begin{equation*}  
\tJd(x,r) = \int_{\Phi_{x,\sqrt{r^2 - r_1^2}}} \tJd(\cdot,r_1) + \bar \O_s \Ll( C r_1^{-\al} \Rr) .
\end{equation*}
Recalling the choice of $\eta_1$ in \eqref{e.cond.eta}, we see that the $\bar \O$ term above can be replaced by $\bar \O_s \Ll( C R^{-\be - \eps} \Rr) $.
Multiplying the last display by $f(x)$, integrating over $\cu_R$ and using Lemma~\ref{l.sum.barO}, we get
\begin{multline*}  
\fint_{\cu_R} f \, \tJd(\cdot,r) \\
= R^{-d} \int_{\cu_R} \int_{\Rd} f(x) \Phi\Ll(y-x, {r^2 - r_1^2}\Rr) \tJd(y,r_1) \, dy \, dx + \bar \O_s \Ll( C R^{-\be - \eps} \Rr).
\end{multline*}
This is \eqref{e.replace.r}.

\smallskip

\emph{Step 2.}
We prove \eqref{e.rewrite.induc}. We first cover the box of side length $R + R^{1-2\eps}$ with $2^d$ boxes of side length $\frac R 2$ at distance at least $R^{1-2\eps}$ from one another, and a set of boxes of side length $R^{1-2\eps}$ for the remainder. More precisely, we define
\begin{equation*}  
\mcl Z := \Ll( \frac {R} 4  + \frac{R^{1-2\eps}}{2} \Rr) \{-1,1\}^d,
\end{equation*}
a subset of $\Rd$ with $2^d$ elements. The set
\begin{equation*}  
\cu_{R + R^{1-2\eps}} \setminus \bigcup_{z \in \mcl Z} \cu_{\frac R 2}(z)
\end{equation*}
can be covered (up to a set of null measure) with the non-overlapping boxes
\begin{equation*}  
\bigcup_{z \in \mcl Z'} \cu_{R^{1-2\eps}}(z),
\end{equation*}
for some $\mcl Z' \subset R^{1-2\eps} \Zd$ such that 
\begin{equation}
\label{e.card.Z'}
|\mcl Z'| \le C R^{2\eps(d-1)}.
\end{equation}
Up to a set of null measure, the complementary set
\begin{equation*}  
\Rd \setminus \Ll( \bigcup_{z \in \mcl Z} \cu_{\frac R 2}(z) \cup \bigcup_{z \in \mcl Z'} \cu_{R^{1-2\eps}}(z) \Rr) 
\end{equation*}
is contained in $ \Rd \setminus \cu_{R + R^{1-2\eps}}$. Since 
\begin{equation}  
\label{e.upper.bound.g}
x \in \Rd \setminus \cu_R \quad \implies \quad R^{-d} g(x) \le C \exp \Ll( -\frac{|x|^2}{C r_1^2} \Rr) ,
\end{equation}
and  $|\tJd(x,r_1)| \le C$, we have, for $\eps > 0$ sufficiently small that $1-2\eps > \eta_1$, 
\begin{equation*}  
\Ll|R^{-d} \int_{\Rd \setminus \Ll( \bigcup_{z \in \mcl Z} \cu_{\frac R 2}(z) \cup \bigcup_{z \in \mcl Z'} \cu_{R^{1-2\eps}}(z) \Rr)}  g \, \tJd(\cdot,r_1)\Rr| \le C R^{- 100d} 
\end{equation*}
almost surely. We treat the remaining parts of the integral separately, starting with $\bigcup_{z \in \mcl Z} \cu_{\frac R 2}(z)$. By the definition of $\tJd(z,r_1)$, the random variables
\begin{equation*}  
\Ll( \int_{\cu_{\frac R 2}(z)} g \, \tJd(\cdot,r_1) \Rr)_{z \in \mcl Z}
\end{equation*}
are independent. By the induction hypothesis, each satisfies
\begin{equation*}  
\Ll(\frac{R}2\Rr)^{-d} \int_{\cu_{\frac R 2}(z)} g \, \tJd(\cdot,r_1) = \bar \O_s \Ll( \CC \Ll(\frac R 2\Rr)^{-\beta}  \Rr). 
\end{equation*}
Hence, by Lemma~\ref{l.barO}, we deduce that
\begin{equation*}  
 R^{-d} \sum_{z \in \mcl Z} \int_{\cu_{\frac R 2}(z)} g \, \tJd(\cdot,r_1) = \bar \O_s \Ll( 2^{\beta-\frac d 2} \CC R^{-\beta} \Rr) .
\end{equation*}
Since $\beta \le \frac d 2$, this implies
\begin{equation*}
\label{e.the.tight.place}
R^{-d} \sum_{z \in \mcl Z} \int_{\cu_{\frac R 2}(z)} g \, \tJd(\cdot,r_1) = \bar \O_s \Ll( \CC R_1^{-\beta} \Rr).
\end{equation*}
It remains to analyze the contribution of the integral over $\bigcup_{z \in \mcl Z'} \cu_{R^{1-2\eps}}(z)$. For $\eps > 0$ sufficiently small (and $R_1$ sufficiently large), we have
\begin{equation*}  
r_1 = \Ll(\frac R 2 \Rr)^{\eta_1} \in [R^{(1-2\eps)\eta_1}, R^{(1-2\eps)\eta_2}].
\end{equation*}
We can thus apply the induction hypothesis and obtain
\begin{equation*}  
R^{-d(1-2\eps)} \int_{\cu_{R^{1-2\eps}}(z)} g \, \tJd(\cdot,r_1) = \bar \O_s \Ll( \CC R^{-\beta(1-2\eps)} \Rr) .
\end{equation*}
Using Lemma~\ref{l.barO.boxes} and \eqref{e.card.Z'}, we deduce that
\begin{equation*}  
R^{-d} \sum_{z \in \mcl Z'} \int_{\cu_{R^{1-2\eps}}(z)} g \, \tJd(\cdot,r_1) = \bar \O_s \Ll( C \CC R^{-2d\eps + (d-1)\eps - \beta(1-2\eps)} \Rr) .
\end{equation*}
The exponent on the right side above can be rewritten as
\begin{equation*}  
-\Ll[\beta + \eps  + \eps\Ll( d-2\beta \Rr) \Rr].
\end{equation*}
Since we assume $\beta \le \frac d 2$, this completes the proof of \eqref{e.rewrite.induc}, and therefore of \eqref{e.induc.scale}.

\smallskip

\emph{Step 3.} We have thus shown that for every $R \ge 1$, $r \in [R^{\eta_1}, R^{\eta_2}]$, $z \in \Rd$, and every function $f \in L^\infty(\Rd)$ satisfying $\|f\|_{L^\infty(\Rd)} \le 1$, we have
\begin{equation}  
\label{e.you.may.recall}
\fint_{\cu_R(z)} f \,\tJd(\cdot,r) = \bar \O_s \Ll( C R^{-\beta} \Rr) .
\end{equation}
We now use additivity once more to show that this implies \eqref{e.suffice.fluc}. We fix $r = R^{\eta_1}$ and apply \eqref{e.additivity.tJd} to get
\begin{equation*}  
\tJd(0,R) = \int_{\Phi_{\sqrt{R^2 - r^2}}} \tJd(\cdot,r) + \O_s \Ll( C r^{-\al} \Rr) .
\end{equation*}
By the choice of $\eta_1$ in \eqref{e.cond.eta}, the $\O$ term on the right side can be replaced by $\O_s \Ll( C R^{-\beta} \Rr)$. We then decompose the integral on the right side above as
\begin{equation*}  
\sum_{z \in R \Z^d} \int_{\cu_R(z)} \Phi_{\sqrt{R^2 - r^2}} \, \tJd(\cdot,r).
\end{equation*}
By \eqref{e.you.may.recall}, we have, for every $z \in R \Z^d$,
\begin{equation*}  
R^{-d} \int_{\cu_R(z)} \Phi_{\sqrt{R^2 - r^2}} \, \tJd(\cdot,r) =  \O_s \Ll( C \|\Phi_{\sqrt{R^2 - r^2}}\|_{L^\infty(\cu_R(z))} R^{-\beta} \Rr) .
\end{equation*}
Moreover, by a Riemann sum argument, we verify that
\begin{equation*}  
\sum_{z \in R \Z^d} \|\Phi_{\sqrt{R^2 - r^2}}\|_{L^\infty(\cu_R(z))} \le C R^{-d}.
\end{equation*}
By Lemma~\ref{l.sum-O}, we therefore obtain
\begin{equation*}  
\int_{\Phi_{\sqrt{R^2 - r^2}}} \tJd(\cdot,r) = \O_s \Ll( C R^{-\beta} \Rr) .
\end{equation*}
This is \eqref{e.suffice.fluc}, so the proof is complete.
\end{proof}

\section{Corrector estimates in weak norms}
\label{s.correctors}

\index{corrector!first-order~$\phi_e$|(}

In this section, we demonstrate that Theorem~\ref{t.additivity} implies Theorem~\ref{t.correctors}, with the exception of the final estimate~\eqref{e.phid=2}, which is proved in Section~\ref{s.twopoint}. Much of the work is already implicit in Section~\ref{s.additivity}, namely in the estimate~\eqref{e.matchJtophi}. This estimate and Theorem~\ref{t.additivity} imply the first three estimates of Theorem~\ref{t.correctors} on the spatial averages of the gradient, flux and energy density of the first-order correctors. The estimates~\eqref{e.phioscd>2} and~\eqref{e.phioscd=2} on the oscillation of the correctors themselves are obtained from these, as expected, by an appropriate version of the multiscale Poincar\'e inequality, Lemma~\ref{l.mspoincare.masks} in this case. The final ``almost pointwise'' estimate~\eqref{e.phid=2} on the oscillation of~$\phi_e$ in dimension $d=2$ requires a more subtle analysis which is postponed to Section~\ref{s.twopoint}. 

\smallskip

We will actually prove something more than the statement of Theorem~\ref{t.correctors} by allowing for spatial averages against more general masks~$\Psi$ than the heat kernel $\Phi_r$. In particular, we would like to be able to take compactly supported functions but, to have the best estimates, it is necessary to take them sufficiently smooth. It turns out to be convenient to take Sobolev functions (with possibly a fractional order). It is therefore natural to write the estimates we seek in terms of the $W^{-\alpha,p}(B_1)$ norms of rescalings of the gradients, fluxes and energy densities of the correctors. The main result is stated in the following theorem. Recall that~$\zeta_\ep$ is the standard mollifier defined in~\eqref{e.standardmollifier}. 


\begin{theorem}[{Corrector estimates in $W^{-\alpha,p}$}]
\label{t.HKtoSob}
Fix~$s\in (0,2)$, $\alpha\in (0,\infty)$ and $p\in [2,\infty)$.
There exist $C(s,\alpha,p,d,\Lambda)<\infty$ and $\delta(\alpha,d,\Lambda)>0$ such that, for every $r \in[1,\infty)$, $e\in\partial B_1$ and $\ep \in \left(0,\tfrac12\right]$, 
\begin{multline}
\label{e.SobHd2}
\left\| \nabla \phi_e  \left( \tfrac \cdot \ep\right)  \ast\zeta_\ep \right\|_{W^{-\alpha,p}(B_1)}  
+
\left\| \left(  \a \left( \tfrac \cdot \ep\right) \left(e +  \nabla \phi_e \left( \tfrac \cdot \ep\right)\right) \right)  \ast\zeta_\ep - \ahom e \right\|_{W^{-\alpha,p}(B_1)}
\\
+ \left\| \left( \tfrac12 \left(e +  \nabla \phi_e \left( \tfrac \cdot \ep\right)\right) \cdot  \a \left( \tfrac \cdot \ep\right) \left(e +  \nabla \phi_e \left( \tfrac \cdot \ep\right)\right) \right) \ast\zeta_\ep - \tfrac12 e\cdot \ahom e \right\|_{W^{-\alpha,p}(B_1)}
\\
\leq
\left\{ 
\begin{aligned}
&  \O_{2+\delta}\left( C\ep^{\alpha}\right) & \mbox{if} & \ \alpha < \tfrac{d}{2}, \\
& \O_s\left( C \ep^{\frac d2} \left| \log \ep \right|^{\frac12} \right) & \mbox{if} & \ \alpha = \tfrac{d}{2}, \\
& \O_s\left( C\ep^{\frac d2} \right) & \mbox{if} & \ \alpha> \tfrac d2.
\end{aligned}
\right. 
\end{multline}
\end{theorem}

The reason for convolving the small scales with the standard mollifier $\zeta_\ep$ is that, without further assumptions on the coefficients, we cannot expect to control the regularity of the correctors on scales smaller than the correlation length scale (which in this case is~$\ep$) beyond what deterministic elliptic regularity estimates provide. In particular, we cannot improve the integrability of their gradients beyond what the Meyers estimate gives. As we have explained in other contexts, these small scales do not concern homogenization and it is therefore natural to ignore them. The mollifier allows us to do just that. In the case $p=2$ (or $p$ slightly larger than~$2$, by the Meyers estimate), we can control the first two terms on the right side of~\eqref{e.SobHd2} without the mollifier. An estimate we have without the mollifier is, for instance:
\begin{multline}
\label{e.SobHd2nozeta}
\left\| \nabla \phi_e \left( \tfrac \cdot \ep\right) \right\|_{H^{-\alpha}(B_1)}  
+
\left\|  \a \left( \tfrac \cdot \ep\right) \left(e +  \nabla \phi_e \left( \tfrac \cdot \ep\right)\right) - \ahom e \right\|_{H^{-\alpha}(B_1)}
\\
\leq
\left\{ 
\begin{aligned}
&  \O_{2+\delta}\left( C\ep^{\alpha}\right) & \mbox{if} & \ \alpha < \tfrac{d}{2}, \\
& \O_s\left( C \ep^{\frac d2} \left| \log \ep \right|^{\frac12} \right) & \mbox{if} & \ \alpha = \tfrac{d}{2}, \\
& \O_s\left( C\ep^{\frac d2} \right) & \mbox{if} & \ \alpha> \tfrac d2.
\end{aligned}
\right. 
\end{multline}
Similarly,~\eqref{e.SobHd2} holds for every $p\in[2,\infty)$ without $\zeta_\ep$ provided we work under the additional assumption that the coefficient field $\a(x)$ is uniformly continuous (with the constant~$C$ depending additionally on the modulus of continuity). In the course of the proof of Theorem~\ref{t.HKtoSob}, we will indicate where the argument should be modified to obtain these results.

\smallskip

Theorem~\ref{t.HKtoSob} is obtained from Theorem~\ref{t.correctors} and some functional inequalities which transfer estimates on the spatial averages of a function $f$ against translates of the heat kernel into estimates on~$\| f \|_{W^{-\alpha,p}(B_1)}$ for general $\alpha>0$ and $p\in (1,\infty)$.  

\smallskip

We proceed now with the proof of most of Theorem~\ref{t.correctors}. 

\begin{proof}[{Proof of Theorem~\ref{t.correctors}, with the exception of~\eqref{e.phid=2}}]

\emph{Step 1.} We prove~\eqref{e.gradient},~\eqref{e.flux} and~\eqref{e.energydensity}.
These estimates are a straightforward consequence of Theorem~\ref{t.additivity} and the basic properties of~$J_1$. 

\smallskip

Fix $s\in (0,2)$. 
According to Theorem~\ref{t.additivity},  
\begin{equation} 
\label{e.Jcaptchagg}
\left| J_1(\Phi_r,p,q) - \left( \frac12 p\cdot \ahom p + \frac12 q \cdot \ahom^{-1}q - p\cdot q \right) \right| 
\leq \O_s\left( Cr^{-\frac d2} \right). 
\end{equation}
Theorem~\ref{t.additivity} also implies that~$\Fluc\left(s,\tfrac d2\right)$ holds and thus we may use Lemma~\ref{l.shortcuts} to obtain that, for every $r\geq 1$ and $p,q\in B_1$,
\begin{equation} 
\label{e.L2paramatching}
\left\| \nabla v(\cdot,\Phi_r,p,q) - \nabla\overline{\phi}_{(-p + \shortbar{\a}^{\,-1}q)} \right\|_{L^2\left(\Phi_{r }\right)} 
\leq 
\O_{s} \left( Cr^{-\frac d2} \right).
\end{equation}
Recall that, by~\eqref{e.J polarization}, for every $e,p,q\in\Rd$, 
\begin{equation} \label{e.L2paramatching000}
p\cdot \int_{\Psi} \a\nabla v(\cdot,\Psi,-e,0)
=
J_1(\Psi,-e+p,0) - J_1(\Psi,-e,0) - J_1(\Psi,p,0) 
\end{equation}
and
\begin{equation} \label{e.L2paramatching001}
q\cdot \int_{\Psi} \nabla v(\cdot,\Psi,-e,0)
=
J_1(\Psi,-e,q) - J_1(\Psi,-e,0) + J_1(\Psi,0,q).
\end{equation}
We also have that, by~\eqref{e.Jexpressbyv},
\begin{equation*} \label{}
J_1(\Psi,p,q) = \int_{\Phi_r} \frac12  \nabla v(\cdot,\Psi,p,q)\cdot \a \nabla v(\cdot,\Psi,p,q).
\end{equation*}
The above identities,~\eqref{e.Jcaptchagg},~\eqref{e.L2paramatching} and the triangle inequality now give the desired bounds~\eqref{e.gradient},~\eqref{e.flux} and~\eqref{e.energydensity}. For the latter we also use the upper bound~\eqref{e.Jvbounded}. 

\smallskip

\emph{Step 2.} We prove~\eqref{e.phioscd>2} and~\eqref{e.phioscd=2} by combining the bound~\eqref{e.gradient} for the spatial averages of the gradient of the corrector with the multiscale Poincar\'e inequality, Lemma~\ref{l.mspoincare.masks}. Applying the latter to $\phi_e- (\phi_e\ast \Phi_R)(0)$ for fixed $e\in\partial B_1$ yields
\begin{multline} 
\label{e.bigMSPapp}
 \int_{\Psi_R}  \left| \phi_e(y)- (\phi_e\ast \Phi_R)(0)  \right|^2\,dy 
 \leq C\int_{\Psi_R} \left| \left( \phi_e \ast \Phi_R \right)(y) - (\phi_e\ast \Phi_R)(0) \right|^2\,dy 
\\
 +C \int_{0}^{R^2} \int_{\Psi_R} \left| \left( \nabla \phi_e \ast \Phi_{\sqrt{t}} \right)(y) \right|^2\,dy\,dt.
\end{multline}
In dimensions $d>2$, we use Lemma~\ref{l.change-s} to interpolate between~\eqref{e.correctorgradbound} and~\eqref{e.gradient}. We take $\ep \in (0,\tfrac12)$, set 
\begin{equation*} \label{}
s:= \left\{ 
\begin{aligned}
& 2 & \mbox{if} & \ d=2, \\
& 2+\ep & \mbox{if} & \ d>2,
\end{aligned} 
\right.
\end{equation*}
and obtain, for every $s' < {s}$ and $x\in\Rd$,
\begin{equation} 
\label{e.gradboundforMSP}
\left| \int_{\Phi_{x,r}} \nabla \phi_e \right| \leq \O_{s'} \left( C r^{-\frac ds} \right). 
\end{equation}
This and Lemma~\ref{l.sum-O} imply that, for every $x, y \in \Rd$, 
\begin{align} 
\label{e.averagesdervs}
\left| \phi_e \ast \Phi_{R} (y) - \phi_e \ast \Phi_{R}(x) \right| 
&
= \left| \int_{0}^1 (y-x) \cdot \nabla \phi_e\ast \Phi_R(ty+(1-t)x)\,dt \right| \notag
\\ & \notag
\leq |y-x| \int_{0}^1 \left|  \nabla \phi_e\ast \Phi_R(ty+(1-t)x) \right| \,dt
\\ &  
\le \O_{s'}\left( C |y-x| R^{-\frac ds} \right). 
\end{align}
We thus obtain the following estimate for the first term on the right side of~\eqref{e.bigMSPapp}: 
\begin{align*} \label{}
\int_{\Psi_R} \left| \left( \phi_e \ast \Phi_R \right)(y) - (\phi_e\ast \Phi_R)(0) \right|^2\,dy 
&
\leq \O_{s'/2} \left(CR^{-\frac {2d}s} \int_{\Rd} |y|^2\Psi_R(y)\,dy \right)
\\ &
\leq \O_{s'/2} \left( CR^{2-\frac{2d}{s}} \right)
\\ &
\leq \O_{s'/2}(C).
\end{align*}
We split the second term on the right side of~\eqref{e.bigMSPapp} into small scales and large scales. The small scales we crudely estimate using~\eqref{e.correctorgradbound}:
\begin{equation} 
\label{e.smallscalestuff}
 \int_{0}^{1} \int_{\Psi_R} \left| \left( \nabla \phi_e \ast \Phi_{\sqrt{t}} \right)(y) \right|^2\,dy\,dt
 \leq C \int_{\Psi_R} \left| \nabla \phi_e(y) \right|^2\,dy
  \leq \O_{s'/2}(C).
 \end{equation}
For the large scales, we use~\eqref{e.gradboundforMSP} and Lemma~\ref{l.sum-O}:
\begin{align*} \label{}
 \int_{1}^{R^2} \int_{\Psi_R} \left| \left( \nabla \phi_e \ast \Phi_{\sqrt{t}} \right)(y) \right|^2\,dy\,dt
 &
 \leq \O_{s'/2} \left( C \int_{1}^{R^2} t^{-\frac d2} \,dt\right)
\\ &
 = \left\{ 
 \begin{aligned}
& \O_{s'/2} \left( C\log R \right) & \mbox{if} & \ d=2, \\
& \O_{s'/2} (C) & \mbox{if} & \ d>2. 
 \end{aligned}
 \right.
 \end{align*}
Combining the above inequalities yields, for every $s'<s$, 
\begin{equation}
\label{e.step2mostproof}
 \int_{\Psi_R}  \left| \phi_e(y)- (\phi_e\ast \Phi_R)(0)  \right|^2\,dy
\leq 
\left\{ 
 \begin{aligned}
& \O_{s'/2} \left( C\log R \right) & \mbox{if} & \ d=2, \\
& \O_{s'/2} (C) & \mbox{if} & \ d>2. 
 \end{aligned}
 \right.
\end{equation}
Note that in dimension~$d=2$ we have now completed the proof of~\eqref{e.phioscd=2}. 

\smallskip

To complete the proof in dimensions $d>2$ of~\eqref{e.phioscd>2} and the statement preceding it, we next perform a variant of the above estimates to show that, with $s < \frac 52$ as above, for every $s'<2$, we have 
\begin{equation} 
\label{e.averagesarecauchy}
\left| \phi_e \ast \Phi_{2R} (x) - \phi_e \ast \Phi_{R}(x) \right| \leq \O_{s'}\left( C R^{1-\frac ds} \right). 
\end{equation}
Note that this inequality makes sense even though $\phi_e$ is only well-defined, a priori, up to an additive constant. 
By \eqref{e.averagesdervs} and Lemma~\ref{l.sum-O}, we have
\begin{align*}
\left| \phi_e \ast \Phi_{2R} (x) - \phi_e \ast \Phi_{R}(x) \right|
&
= \left| \int_{\Phi_{\sqrt{3}R}} \left( \phi_e \ast \Phi_R(y+x) - \phi_e \ast \Phi_R(x) \right)\,dy \right|
\\ & 
\leq \int_{\Phi_{\sqrt{3}R}} \left| \phi_e \ast \Phi_R(y+x) - \phi_e \ast \Phi_R(x)  \right|\,dy
\\ & 
\leq \int_{\Phi_{\sqrt{3}R}} \O_s\left( C|y| R^{-\frac ds} \right) \,dy
\\ & 
\leq \O_{s'}\left( C R^{1-\frac ds} \right). 
\end{align*}
Now, in dimensions $d>2$, since $\frac ds > 1$, we find that the sequence $$\{\phi_e(x)  -  \phi_e \ast \Phi_{2^n} (x) \}_{n\in\N}$$
is Cauchy,~$\P$--a.s. We denote its limit by $\phi_e(x) - (\phi_e)_\infty(x)$. By~\eqref{e.averagesdervs}, we see that the quantity $(\phi_e)_\infty(x)$ exists for every $x\in\Rd$, $\P$--a.s., and is independent of~$x$. We denote this deterministic constant by~$(\phi_e)_\infty$. We may therefore select the additive constant for $\phi_e$ in a stationary way by requiring that $(\phi_e)_\infty = 0$ and obtain a stationary random field. Note that, by the estimates above,  
\begin{equation*}
\left| \phi_e\ast \Phi_R(x) \right| 
\leq \sum_{k\in\N} \left| \phi_e \ast \Phi_{2^kR} (x) - \phi_e \ast \Phi_{2^{k+1}R}(x) \right|
\leq \O_{s'}(C).
\end{equation*}
This completes the proof of~\eqref{e.phioscd>2}.
\end{proof}

In the next lemma, we squeeze some further  information from the previous argument by showing that, up to a small error, the spatial averages of the gradient, flux and energy density of the first-order correctors are equivalent to the quantity~$J_1$. 

\begin{lemma}
\label{l.redwedding}
For every~$\alpha\in (0,d]$ and~$s\in\left( 0, 2\wedge \frac{d}{\alpha}\right)$, there exists~$C(\alpha,s,d,\Lambda)<\infty$ such that, for every~$p,q\in B_1$,~$x\in\Rd$ and~$r\geq 1$, 
\begin{align}
\label{e.J1isgradflux.hk}
\left| 
J_1(\Phi_{x,r},p,q) 
-
\frac12 \int_{\Phi_{x,r}} \left( - p\cdot \a\nabla\overline{\phi}_{(-p + \shortbar{\a}^{\,-1}q)} + q\cdot \nabla\overline{\phi}_{(-p + \shortbar{\a}^{\,-1}q)} \right)
\right|
\leq \O_{s}\left( Cr^{-\alpha} \right)
\end{align}
and
\begin{align}
\label{e.J1isenergy.hk}
\left| 
J_1(\Phi_{x,r},p,q) 
-
\frac12 \int_{\Phi_{x,r}}  \nabla\overline{\phi}_{(-p + \shortbar{\a}^{\,-1}q)} \cdot \a\nabla\overline{\phi}_{(-p + \shortbar{\a}^{\,-1}q)}
\right|
\leq \O_{s}\left( Cr^{-\alpha} \right).
\end{align}
\end{lemma}
\begin{proof}
By Theorem~\ref{t.additivity}, $\Fluc\left(s,\tfrac d2\right)$ holds for every $s\in (0,2)$. Therefore  Lemma~\ref{l.shortcuts} holds for $\alpha=\frac d2$ and $s\in (0,2)$ which, in view of~\eqref{e.J1firstvar} and~\eqref{e.Jexpressbyv}, implies the lemma for such~$\alpha$ and~$s$. It therefore suffices by Lemma~\ref{l.change-s} to prove the lemma in the case~$\alpha=d$ and~$s\in (0,1)$.

\smallskip

The second estimate~\eqref{e.J1isenergy.hk} is actually an immediate consequence of~\eqref{e.J1quadresponse} and~\eqref{e.L2paramatching}.
To prove the first estimate~\eqref{e.J1isgradflux.hk}, we recall that, by~\eqref{e.J1firstvar} and~\eqref{e.Jexpressbyv},
\begin{equation}
\label{e.J1lbiden1} 
J_1(\Phi_{x,r},p,q) 
=
\frac12 \int_\Psi \left( - p\cdot \a\nabla v(\cdot,\Phi_{x,r},p,q) + q\cdot \nabla v(\cdot,\Phi_{x,r},p,q) \right).
\end{equation}
To obtain~\eqref{e.J1isgradflux.hk}, it suffices therefore to establish that, for each~$s\in (0,2)$, there exists a constant~$C(s,d,\Lambda)<\infty$ such that, for every~$p,q\in B_1$,~$x\in\Rd$ and~$r\geq 1$, 
\begin{equation}
\label{e.lbgradest1}
\left| \int_{\Phi_{x,r}} \left( \nabla v(\cdot,\Phi_r,p,q) - \nabla\overline{\phi}_{(-p + \shortbar{\a}^{\,-1}q)} \right) \right| 
\leq 
\O_{s/2} \left( Cr^{-d} \right)
\end{equation}
and
\begin{equation}
\label{e.lbgradest2}
\left| \int_{\Phi_{x,r}} \a \left( \nabla v(\cdot,\Phi_r,p,q) - \nabla\overline{\phi}_{(-p + \shortbar{\a}^{\,-1}q)} \right) \right| 
\leq 
\O_{s/2} \left( Cr^{-d} \right).
\end{equation}
As $\nabla v(\cdot,\Phi_r,p,q) - \nabla\overline{\phi}_{(-p + \shortbar{\a}^{\,-1}q)} \in \nabla \mathcal{A}_1$, we have, by~\eqref{e.gradient},~\eqref{e.correctorsballz} and~\eqref{e.L2paramatching},
\begin{align*}
\lefteqn{
\left| \int_{\Phi_{x,r}} \left( \nabla v(\cdot,\Phi_r,p,q) - \nabla\overline{\phi}_{(-p + \shortbar{\a}^{\,-1}q)} \right) \right| \indc_{\{ \Y_{sd/2}(x) \leq r \}}
} \quad & 
\\ &
\leq
\O_{s} \left( Cr^{-\frac d2} \left\| \nabla v(\cdot,\Phi_r,p,q) - \nabla\overline{\phi}_{(-p + \shortbar{\a}^{\,-1}q)} \right\|_{L^2\left(\Phi_{r }\right)} \right)  
\leq 
\O_{s/2} \left( Cr^{-d} \right). 
\end{align*}
On the other hand,~\eqref{e.indcYsbound} and~\eqref{e.correctorsballz} imply that 
\begin{align*}
\lefteqn{
\left| \int_{\Phi_{x,r}} \left( \nabla v(\cdot,\Phi_r,p,q) - \nabla\overline{\phi}_{(-p + \shortbar{\a}^{\,-1}q)} \right) \right| \indc_{\{ \Y_{sd/2}(x) \geq r \}}
} \quad & 
\\ &
\leq
C \left( \left\| \nabla v(\cdot,\Phi_r,p,q) \right\|_{L^2\left(\Phi_{x,r}\right)} 
+
\left\| \nabla\overline{\phi}_{(-p + \shortbar{\a}^{\,-1}q)} \right\|_{L^2\left(\Phi_{x,r}\right)}\right) \indc_{\{ \Y_s(x) \geq r \}}
\\ & 
\leq C \left( 1 + \left( \frac{\Y_{sd/2}(x)}{r} \right)^{\frac{d}2}\left\| \nabla\overline{\phi}_{(-p + \shortbar{\a}^{\,-1}q)} \right\|_{L^2\left(\Phi_{x,\Y_{sd/2}(x)}\right)} \right)\indc_{\{ \Y_{sd/2}(x) \geq r \}}
\\ & 
\leq 
C \left( \frac{\Y_{sd/2}(x)}{r} \right)^{\frac{d}2} \indc_{\{ \Y_s(x) \geq r \}}
\\ & 
\leq \O_{s/2} \left( Cr^{-d} \right).
\end{align*}
Combining these yields~\eqref{e.lbgradest1}. The argument for~\eqref{e.lbgradest2} is almost identical. The proof of~\eqref{e.J1isgradflux.hk} is complete.
\end{proof}

\begin{remark}
\index{Liouville theorem|(}
\label{r.canonical}
Theorem~\ref{t.regularity} asserts that~$\A_{k+1} / \A_k$ and~$\Ahom_{k+1} / \Ahom_k$ are \emph{canonically} isomorphic vector spaces, for every~$k\in\N$, in the sense that an isomorphism can be chosen which is invariant under translation of the coefficients. The final paragraph of the proof of Theorem~\ref{t.correctors}, which proves that~$\phi_e$ exists as a stationary random field, also proves the essentially equivalent statement that this isomorphism in the case $k=0$ can be extended to a canonical isomorphism between~$\A_1$ and~$\Ahom_1$ in dimensions~$d>2$. It is an open problem to identify a canonical isomorphism for the space $\A_{k + \Ll\lceil\frac d 2 \Rr\rceil}/\A_k$.
\end{remark}

\index{Liouville theorem|)}

We now give the proof of Theorem~\ref{t.HKtoSob} by combining Theorem~\ref{t.correctors} with the functional inequalities in Appendix~\ref{a.MSP}.

\begin{proof}[{Proof of Theorem~\ref{t.HKtoSob}}]
We give only the proof of the estimates for the gradient of the correctors. The arguments for the bounds for the flux and energy density are almost identical. It is convenient to denote $\Psi(x):= \exp \left( -|x| \right)$.  

\smallskip

Fix~$\alpha\in (0,\infty)$, $p\in [2,\infty)$ and $\ep \in \left(0,\frac12\right]$. 
Denote 
\begin{equation*} \label{}
\mathbf{f}_\ep:= 
\nabla \phi_e\left(\tfrac\cdot\ep\right)
\ast \zeta_\ep. 
\end{equation*}
We intend to apply the inequalities of Remark~\ref{r.MSPlocal} to $\mathbf{f}_\ep$. In the case $\alpha\in\N$ we use~\eqref{e.triebel_loc}, while if $\alpha\not\in\N$ we apply~\eqref{e.MSP.Walpha_loc}. In order to prepare for the application of these, we first use Theorem~\ref{t.correctors} to estimate $\left| \left( \mathbf{f}_\ep \ast \Phi(t, \cdot) \right)(x)\right|$.
What the theorem gives us is that, 
for every $x\in \Rd$ and $t\geq \ep^2$, 
\begin{equation} 
\label{e.tcorrectorsapp}
\left| \left(  \mathbf{f}_\ep \ast \Phi(t, \cdot) \right)(x) \right| 
\leq 
\O_s\left( C \left( \ep t^{-\frac12} \right)^{\frac d2} \right). 
\end{equation}
Indeed, by Lemma~\ref{l.sum-O},
\begin{align*}
\left| \left(  \mathbf{f}_\ep \ast \Phi(t, \cdot) \right)(x) \right| 
& =
\left| 
\int_{\Rd}\int_{\Rd} \Phi(t,x-y) \nabla \phi_e \left( \frac {y-z}\ep \right) \zeta_\ep(z) \,dz\,dy
\right|
\\ & 
= \left| 
\int_{\Rd}\int_{\Rd} \Phi\left(\frac t{\ep^2}, x-\frac y\ep   \right) \nabla \phi_e \left( \frac {y}\ep - z \right) \zeta(z) \,dz\,dy
\right|
\\ & 
\leq 
\int_{\Rd} \zeta(z) \left| \int_{\Rd} \Phi\left(\frac t{\ep^2}, x-\frac y\ep  \right) \nabla \phi_e \left( \frac {y}\ep - z \right) \,dy \right|\,dz
\\ & 
\leq 
\O_s\left( C \left( \ep t^{-\frac12} \right)^{\frac d2} \right). 
\end{align*}
We will use this inequality for the larger scales ($t\geq \ep^2$) but we still need something to control the small scales. These are handled rather crudely using~\eqref{e.correctorgradbound} and Young's inequality for convolutions. The claim is that there exists~$\delta(d,\Lambda)>0$ such that, for every $t>0$, 
\begin{align}
\label{e.smallscalecaptcha}
 \left\| \mathbf{f}_\ep \ast \Phi(t, \cdot) \right\|_{L^p(\Psi)}
\leq \O_{2+\delta} \left( C \right).
\end{align}
To obtain this, we simply use Young's inequality for convolutions to get 
\begin{equation*} \label{}
\left\| \mathbf{f}_\ep \ast \Phi(t, \cdot) \right\|_{L^p(\Psi)}
\leq \left\| \mathbf{f}_\ep \right\|_{L^{p}(\Psi)}  \left\| \Phi \left(t, \cdot\right) \right\|_{L^1(\Psi)} 
\leq \left\| \mathbf{f}_\ep \right\|_{L^{p}(\Psi)}.
\end{equation*}
The right side is controlled by~\eqref{e.correctorgradbound} and Lemma~\ref{l.sum-O}:
\begin{align*}
\left\| \mathbf{f}_\ep \right\|_{L^{p}(\Psi)}^p
 = 
\int_{\Rd} \Psi(x) \left| \int_{\Rd} \nabla \phi \left( \frac {y-x}\ep \right) \zeta_\ep(y) \,dy \right|^p\,dx
& \leq
\int_{\Rd} \Psi(x) \left\| \nabla \phi  \right\|_{L^2\left(B_1(x/\ep)\right)}^p  \,dx
\\
& \leq 
\O_{(2+\delta)/p} (C). 
\end{align*}
It is here, in the proof of~\eqref{e.smallscalecaptcha}, that we used the mollifier~$\zeta_\ep$. If we were working in the case~$p=2$ or under the assumption that the coefficients are uniformly continuous, we could obtain the same estimate without~$\zeta_\ep$ in the definition of~$\mathbf{f}_\ep$. 

\smallskip

Now we are ready to use the multiscale Poincar\'e-type functional inequalities. We must split the argument into two cases, depending on whether $\alpha \in\N$ or not. We first assume that $\alpha\in\N$ and apply~\eqref{e.triebel_loc}, which gives us 
\begin{align*}
\left\| \mathbf{f}_\ep \right\|_{W^{-\alpha,p}(B_1)} 
& 
\leq C \left( \int_{\R^d} \Psi(x) \left( \int_0^1 \left( t^{\frac \alpha2} \left| \mathbf{f}_\ep  \ast \left(\Phi(t, \cdot)  \right)(x)\right| \right)^2 \, \frac{dt}{t}  \right)^{\frac p2}  \,dx \right)^{\frac1p}
\\ & 
\leq 
C \left( \int_{\R^d} \Psi(x) \left( \int_{\ep^2}^1 \left( t^{\frac \alpha2} \left| \mathbf{f}_\ep  \ast \left(\Phi(t, \cdot)  \right)(x)\right| \right)^2 \, \frac{dt}{t}  \right)^{\frac p2}  \,dx \right)^{\frac1p}
\\ & \qquad 
+ C \left( \int_{\R^d} \Psi(x) \left( \int_{0}^{\ep^2} \left( t^{\frac \alpha2} \left| \mathbf{f}_\ep  \ast \left(\Phi(t, \cdot)  \right)(x)\right| \right)^2 \, \frac{dt}{t}  \right)^{\frac p2}  \,dx \right)^{\frac1p}.
\end{align*}
To estimate the first term, we use~\eqref{e.tcorrectorsapp} and Lemma~\ref{l.sum-O} to get 
\begin{multline*}
\int_{\ep^2}^1 \left( t^{\frac \alpha2} \left| \mathbf{f}_\ep  \ast \left(\Phi(t, \cdot)  \right)(x)\right| \right)^2 \, \frac{dt}{t}
\\
 \leq \int_{\ep^2}^1 t^\alpha \O_{s/2} \left( C \left( \ep t^{-\frac12} \right)^{d} \right)\,\frac{dt}{t}
\leq 
\left\{ 
\begin{aligned}
&  \O_{(2+\delta)/2}\left( C\ep^{2\alpha}\right) & \mbox{if} & \ \alpha < \tfrac{d}{2}, \\
& \O_{s/2}\left( C \ep^{d} \left| \log \ep \right| \right) & \mbox{if} & \ \alpha = \tfrac{d}{2}, \\
& \O_{s/2}\left( C\ep^{d} \right) & \mbox{if} & \ \alpha> \tfrac d2.
\end{aligned}
\right. 
\end{multline*}
We remark that, to get the extra $\delta>0$ in the first line above, we used $\alpha<\tfrac d2$,~\eqref{e.smallscalecaptcha} and $\O_s$--interpolation, i.e., Lemma~\ref{l.change-s}(ii). Using Lemma~\ref{l.sum-O} again, we obtain
\begin{multline*}
\left( \int_{\R^d} \Psi(x) \left( \int_{\ep^2}^1 \left( t^{\frac \alpha2} \left| \mathbf{f}_\ep  \ast \left(\Phi(t, \cdot)  \right)(x)\right| \right)^2 \, \frac{dt}{t}  \right)^{\frac p2}  \,dx \right)^{\frac1p}
\\
\leq 
\left\{ 
\begin{aligned}
&  \O_{(2+\delta)}\left( C\ep^{\alpha}\right) & \mbox{if} & \ \alpha < \tfrac{d}{2}, \\
& \O_{s}\left( C \ep^{\frac d2} \left| \log \ep \right|^{\frac12}  \right) & \mbox{if} & \ \alpha = \tfrac{d}{2}, \\
& \O_{s}\left( C\ep^{\frac d2} \right) & \mbox{if} & \ \alpha> \tfrac d2.
\end{aligned}
\right. 
\end{multline*}
The second term is controlled by~\eqref{e.smallscalecaptcha}, Jensen's inequality and Lemma~\ref{l.sum-O}:
\begin{align*}
\lefteqn{
\int_{\R^d} \Psi(x) \left( \int_{0}^{\ep^2} \left( t^{\frac \alpha2} \left| \mathbf{f}_\ep  \ast \left(\Phi(t, \cdot)  \right)(x)\right| \right)^2 \, \frac{dt}{t}  \right)^{\frac p2}  \,dx
} \quad & \\ & 
\leq
C\ep^{p\alpha }\int_{\Rd} \int_0^{\ep^2} \Psi(x)   \left| \mathbf{f}_\ep  \ast \left(\Phi(t, \cdot)  \right)(x)\right|^p \,\frac {dt}{t^{1-\alpha} \ep^{2\alpha}} \,dx 
\\ & 
\leq \O_{(2+\delta)/p}\left( C \ep^{p\alpha} \right)
\end{align*}
Combining the previous two displays yields the theorem in the case $\alpha \in\N$. 

\smallskip

We next consider the case $\alpha \in (0,\infty) \setminus \N$. An application of~\eqref{e.MSP.Walpha_loc} gives us
\begin{align*}
\left\| \mathbf{f}_\ep \right\|_{W^{-\alpha,p}(B_1)} 
\leq 
C\left( \int_0^1 \left( t^{\frac\alpha2} \left\| \mathbf{f}_\ep \ast \Phi(t, \cdot) \right\|_{L^p(\Psi)} \right)^p \,\frac{dt}t \right)^{\frac1p}.
\end{align*}
As before, we split up the integral over time (which represents the square of the length scale) so that we may handle large scales and small scales differently. Using~\eqref{e.smallscalecaptcha} and Lemma~\ref{l.sum-O}, we have
\begin{align*}
\int_0^{\ep^2}  \left( t^{\frac\alpha2} \left\| \mathbf{f}_\ep \ast \Phi(t, \cdot) \right\|_{L^p(\Psi)} \right)^p \,\frac{dt}t 
& 
\leq 
\int_0^{\ep^2}  t^{p\alpha/2-1} \O_{(2+\delta)/p}(C)\,dt
\\ & 
\leq 
\O_{(2+\delta)/p} \left( \int_0^{\ep^2}  t^{p\alpha/2-1}\,dt \right) 
\\ & 
= \O_{(2+\delta)/p} \left( C \ep^{p\alpha} \right). 
\end{align*}
For the large scales, we use~\eqref{e.tcorrectorsapp} and Lemma~\ref{l.sum-O} (twice) to find that 
\begin{align*}
\int_{\ep^2}^1 \left( t^{\frac\alpha2} \left\| \mathbf{f}_\ep \ast \Phi(t, \cdot) \right\|_{L^p(\Psi)} \right)^p \,\frac{dt}t 
& = 
\int_{\ep^2}^1 t^{p\alpha/2-1} \int_{\Rd} \Psi(x)\left| \left(  \mathbf{f}_\ep \ast \Phi(t, \cdot) \right)(x) \right|^p  \,dx  \,dt 
\\ & 
\leq \int_{\ep^2}^1  \O_{s/p} \left( C \left( \ep t^{-\frac12} \right)^{\frac {pd}2} \right)  \,dt
\\ & 
\leq \O_{s/p} \left( C \ep^{\frac{pd}{2}}  \int_{\ep^2}^1 t^{\frac p2 \left( \alpha-\frac d2 \right) -1} \,dt  \right).
\end{align*}
Taking the $p$th root yields
\begin{align*}
\left( \int_{\ep^2}^1 \left( t^{\frac\alpha2} \left\| \mathbf{f}_\ep \ast \Phi(t, \cdot) \right\|_{L^p(\Psi)} \right)^p \,\frac{dt}t \right)^{\frac1p} 
\leq
\left\{ 
\begin{aligned}
&  \O_{(2+\delta)}\left( C\ep^{\alpha}\right) & \mbox{if} & \ \alpha < \tfrac{d}{2}, \\
& \O_{s}\left( C \ep^{\frac d2} \left| \log \ep \right|^{\frac12}  \right) & \mbox{if} & \ \alpha = \tfrac{d}{2}, \\
& \O_{s}\left( C\ep^{\frac d2} \right) & \mbox{if} & \ \alpha> \tfrac d2.
\end{aligned}
\right. 
\end{align*}
As above, we note that the extra $\delta>0$ in the first line was obtained by using $\alpha<\tfrac d2$,~\eqref{e.smallscalecaptcha} and $\O_s$--interpolation, Lemma~\ref{l.change-s}(ii).
Combining this with the estimate for the small scales gives us the desired estimate in the case $\alpha \not\in\N$. This completes the proof of the theorem.
\end{proof}

\section{Corrector oscillation estimates in two dimensions}
\label{s.twopoint}

In this section we complete the proof of Theorem~\ref{t.correctors} by obtaining the estimate~\eqref{e.phid=2} in dimension~$d=2$. Recall that this estimates asserts, for each $s\in(0,2)$, the existence of~$C(s,\Lambda)<\infty$ such that, for every $x,y \in\R^2$ and $r,R\in[2,\infty)$ with $r\leq R$, we have
\begin{equation}
\label{e.phid=2.again}
\left| \left( \phi_e\ast \Phi_r \right)(x) - \left( \phi_e\ast \Phi_R \right)(y) \right|
\leq
\O_s\left( C\log^{\frac12} \left(2+\frac{R+|x-y|}r\right) \right).
\end{equation}

\smallskip

Throughout we fix $s\in (0,2)$. 
By a variant of the argument leading to~\eqref{e.step2mostproof} in the previous section, we have, in~$d=2$, a constant $C(s,\Lambda)<\infty$ such that, for every $s\in (0,2)$, $r,R\in [2,\infty)$ with $r\leq R$, 
\begin{equation}
\label{e.d2Rroscbound}
 \int_{\Psi_R}  \left| (\phi_e\ast \Phi_r)(y)- (\phi_e\ast \Phi_R)(0)  \right|^2\,dy
\leq 
\O_{s/2} \left( C\log \left(2+\frac{R}{r}\right) \right).
\end{equation}
Moreover, we have already shown in~\eqref{e.averagesdervs} that for some $C(s,\Lambda)<\infty$, we have, for every $x, y \in \R^2$ and $R\geq 2$, 
\begin{align} 
\label{e.d2averagesdervs}
\left| \phi_e \ast \Phi_{R} (y) - \phi_e \ast \Phi_{R}(x) \right| 
\leq \O_{s}\left( C |y-x| R^{-1} \right). 
\end{align} 
To obtain~\eqref{e.phid=2.again}, it therefore suffices by the triangle inequality and stationarity to prove the inequality in the case $x=y=0$ and $2r\leq R$. In other words, we will prove the statement of the following proposition.

\begin{proposition}
\label{p.d2vengeance}
Suppose that $d=2$. For each $s\in (0,2)$, there exists~$C(s,\Lambda)<\infty$ such that, for every $r,R\in[2,\infty)$ with $2r\leq R$, 
\begin{equation}
\label{e.phid=2.reduction}
\left| \left( \phi_e\ast \Phi_r \right)(0) - \left( \phi_e\ast \Phi_R \right)(0) \right|
\leq
\O_s\left( C\log^{\frac12} \left(\frac Rr\right) \right).
\end{equation}
\end{proposition}

The proof of Proposition~\ref{p.d2vengeance} is delicate because the scaling of the estimate is in some sense critical. Before embarking on its proof, let us see that a slightly suboptimal version, in which we give up a square root of the logarithm, has a very easy proof. Thes suboptimal estimate states that, for some~$C(s,\Lambda)<\infty$ and every $r,R\in[2,\infty)$ with $2r\leq R$, 
\begin{equation}
\label{e.phid=2.reduction.easy}
\left| \left( \phi_e\ast \Phi_r \right)(0) - \left( \phi_e\ast \Phi_R \right)(0) \right|
\leq
\O_s\left( C\log  \left(\frac Rr\right) \right).
\end{equation}
This follows from the identity 
\begin{align*}
\left( \phi_e\ast \Phi_r \right)(0) - \left( \phi_e\ast \Phi_R \right)(0)
&
= \int_{r^2}^{R^2}\!\!\int_{\R^2} \nabla \Phi(t,x) \cdot \nabla \phi_e(x)\,dx\,dt 
\notag
\\ &  
= \int_{r^2}^{R^2}\!\!\int_{\R^2} \nabla \Phi\left(\tfrac t2,y\right) \cdot \int_{\R^2} \Phi \left(\tfrac t2,x-y\right)\nabla \phi_e(x)\,dx \,dy\,dt,
\end{align*}
which we can crudely estimate using~\eqref{e.gradient}: we obtain
\begin{align*}
& 
\left| \left( \phi_e\ast \Phi_r \right)(0) - \left( \phi_e\ast \Phi_R \right)(0) \right| 
\\ & \qquad
\leq
\int_{r^2}^{R^2}\int_{\R^2} \left| \nabla \Phi\left(\tfrac t2,y\right) \right|  \left|  \int_{\R^2} \Phi \left(\tfrac t2,x-y\right)\nabla \phi_e(x)\,dx \right| \,dy\,dt
\\ & \qquad
\leq 
\int_{r^2}^{R^2} \O_s\left( Ct^{-1} \right) \,dt
= \O_s \left( C \log\left( \frac Rr \right)\right). 
\end{align*}
This is~\eqref{e.phid=2.reduction.easy}.

\smallskip

If we believe that~\eqref{e.phid=2.reduction} is correct, then we should ask ourselves where our proof of~\eqref{e.phid=2.reduction.easy} was inefficient. The first identity split our quantity into a sum (or integral) of spatial averages of~$\nabla \phi_e$, and we showed that each of the dyadic scales between $r$ and $R$ contributed~$O_s(C)$ (where a ``dyadic scale'' corresponds, roughly, to~$t$ in an interval $[2^kr,2^{k+1}r)$). Since there are $\log \left(\frac Rr\right)$ many dyadic scales,  this gives the estimate~\eqref{e.phid=2.reduction.easy} by the triangle inequality. If we hypothesize, however, that the contributions of these dyadic scales are \emph{decorrelated}, we may expect that the \emph{variance} of the quantity is proportional to~$\log \left(\frac Rr\right)$, giving us an improvement by a square root of the logarithm. This is exactly what we will show. 

\begin{proof}[{Proof of Proposition~\ref{p.d2vengeance}}]
Fix~$s\in (1,2)$ and~$r,R\in [2,\infty)$ with~$2r\leq R$. We will use the following rescaled version of~\eqref{e.SobHd2nozeta}: for every $L\geq 1$ and $f\in H^2_0\left(B_{L}\right)$, 
\begin{equation}
\label{e.rescalebtch}
\left| \int_{\R^2} f(x) \nabla \phi_e(x) \,dx\right| 
\leq 
\O_s\left(CL^{3} \left\| f \right\|_{\underline{H}^2(B_L)} \right). 
\end{equation}
As in the proof of the suboptimal bound~\eqref{e.phid=2.reduction.easy}, we start from the formula
\begin{equation}
\label{e.twopt.PGFidentity}
(\phi_e\ast\Phi_r)(0) - (\phi_e\ast \Phi_R)(0)
=
\int_{\R^2} \nabla H(x) \cdot \nabla \phi_e(x)\,dx
\end{equation}
where
\begin{equation*}
H (x) : = \int_{r^2}^{R^2} \Phi(t,x)\,dt. 
\end{equation*}

\smallskip

\emph{Step 1.} We split the integral into dyadic annuli. We write
\begin{align}
\label{e.Psidecomp}
\int_{\R^2} \nabla H(x) \cdot \nabla \phi_e(x)\,dx
=
\sum_{k=0}^\infty 
\int_{\R^2} \eta_k(x) \nabla {H}(x) \cdot \nabla \phi_e(x)\,dx,
\end{align}
where $\left\{ \eta_k \right\}_{k\geq0} \subseteq C^\infty_c(\R^2)$ is a partition of unity subordinate to $\{A_k\}_{k\geq 0}$, with
\begin{equation*}
A_0:= B_r 
\quad \mbox{and} \quad 
A_k : = B_{r2^k} \setminus B_{r 2^{k-2}}, 
\quad \forall  k\geq 1, 
\end{equation*}
with the following properties:
\begin{equation*}
0 \leq \eta_k \leq 1 \quad \mbox{in} \ \R^2, \quad
\sum_{k=0}^\infty \eta_k \equiv 1 \quad \mbox{in} \ \R^2, \quad 
\eta_k \equiv 0 \ \mbox{in} \ \R^2 \setminus A_k,
\end{equation*}
as well as the regularity bounds
\begin{equation}
\label{e.annulusetak.ests}
\left\| \nabla^m \eta_k \right\|_{L^\infty(A_k)} \leq C_m\left( r2^{k} \right)^{-m} \ \forall m\in\N,
\end{equation}
and the dilation symmetry property
\begin{equation}
\label{e.etakdilationsymmetry}
\eta_k(x) = \eta_1(2^{1-k} x),\quad \forall k\geq 2. 
\end{equation}
Such a partition of unity can be constructed explicitly. With $\zeta_\ep$ the standard mollifier defined in~\eqref{e.standardmollifier}, we take 
\begin{equation*}
\eta_0 := 
\indc_{B_{\frac34r}}\ast \zeta_{\frac14 r}
\qquad \mbox{and} \qquad 
\eta_1(x)
:=
\left\{
\begin{aligned}
& \left( \indc_{\R^2 \setminus B_{\frac 34r} } \ast \zeta_{\frac14 r}\right)(x), & \mbox{if} &\ x\in B_r,\\
& \left( \indc_{B_{\frac32 r}}\ast \zeta_{\frac12r} \right)(x), & \mbox{if} &\ x\in \R^2\setminus B_r.
\end{aligned}
\right. 
\end{equation*}
Then for each $k\geq 2$ we let $\eta_k$ be defined by~\eqref{e.etakdilationsymmetry}.

\smallskip

\emph{Step 2.}
We show that, for every $k\in\N$, 
\begin{align}
\label{e.2pt.remover}
\left| \int_{\R^2} \eta_k(x) \nabla {H}(x) \cdot \nabla \phi_e(x)\,dx \right|
\leq \O_s\left(C \exp\left( -\frac{\left(r2^k\right)^2}{CR^2}\right) \right). 
\end{align}
By~\eqref{e.rescalebtch} it suffices to check that 
\begin{equation*}
(2^kr)^3 \left\| \eta_k \nabla {H} \right\|_{\underline{H}^2(B_r)}
\leq C \exp\left( -\frac{\left(r2^k\right)^2}{CR^2}\right).
\end{equation*}
This follows from straightforward computations. Indeed, for every~$x\in A_k$ and $m\in\N$, we have
\begin{align*}
\left| \nabla^m {H}(x) \right| 
\leq \int_{r^2}^{R^2} \left| \nabla^m \Phi(t,x) \right| \,dt 
&
\leq 
C \int_{r^2}^{R^2} t^{-\frac m2-\frac d2} \exp\left( -\frac{r^22^{2k}}{Ct} \right)\,dt
\\ &
\leq C \left(r2^k \right)^{-m} \exp\left( -\frac{\left(r2^k\right)^2}{CR^2}\right).
\end{align*}
Combining this with~\eqref{e.annulusetak.ests}, we obtain, for each $m\in\N$, 
\begin{equation}
\label{e.etakPsibouns}
\left\| \nabla^m \left( \eta_k \nabla {H} \right) \right\|_{L^\infty(A_k)} 
\leq C\left(r2^k\right)^{-(m+1)}\exp\left( -\frac{\left(r2^k\right)^2}{CR^2}\right),
\end{equation}
Thus in particular $\left\| \eta_k \nabla {H} \right\|_{\underline{W}^{2,\infty}\left(B_{r2^k}\right)} \leq C\left( r2^k \right)^{-3}\exp\left( -\frac{\left(r2^k\right)^2}{CR^2}\right)$. 

\smallskip

\emph{Step 3.} Let $k_*\in\N$ be chosen so that $R < r 2^{k_*} \leq 2R$ and let $G(x)$ denote the elliptic Green function for $-\Delta$, that is, 
\begin{equation}
\label{e.GtoPhid2}
G(x) : = -\frac1{2\pi} \log |x| = \int_0^\infty \Phi(t,x)\,dt. 
\end{equation}
The goal of this step is to show that there exists~$C(s,\Lambda)<\infty$ such that
\begin{equation}
\label{e.PsitoG}
\sum_{k=1}^{k_*} \left| \int_{\R^2} \eta_k(x) \left( \nabla {H}(x) - \nabla G(x) \right)\cdot \nabla \phi_e(x) \,dx \right| 
\leq \O_s(C). 
\end{equation}
As an intermediate step, we first show that,
for every $x\in B_{2R} \setminus B_r$, 
\begin{align}
\label{e.Psi.toGF}
\left| \nabla^m G(x) - \nabla^m {H}(x) \right|
\leq
Cr^{-m} \exp\left( -\frac{|x|^2}{Cr^2} \right)
+
CR^{-m}.
\end{align}
In view of~\eqref{e.GtoPhid2}, we have that 
\begin{equation*}
\left| \nabla^m G(x) - \nabla^m {H}(x) \right|
\leq
\int_0^{r^2} \left| \nabla^m\Phi(t,x) \right|\,dt 
+ \int_{R^2}^\infty \left| \nabla^m\Phi(t,x) \right|\,dt.
\end{equation*}
For each $|x|>r$, we have 
\begin{align*}
\int_0^{r^2} \left| \nabla^m\Phi(t,x) \right|\,dt 
&
\leq 
C\int_0^{r^2} t^{-1 - \frac m2} \exp\left( -\frac{|x|^2}{Ct}\right)\,dt
\leq
Cr^{-m} \exp\left( -\frac{|x|^2}{Cr^2} \right). 
\end{align*}
For each $|x|<2R$, we have 
\begin{align*}
\int_{R^2}^\infty \left| \nabla^m\Phi(t,x) \right|\,dt
\leq
C\int_{R^2}^\infty t^{-1-\frac m2} \,dt = CR^{-m}. 
\end{align*}
Combining these yields~\eqref{e.Psi.toGF}. Using~\eqref{e.annulusetak.ests} and~\eqref{e.Psi.toGF}, we obtain, for each $k \in \{ 1,\ldots,k_*\}$,
\begin{equation*}
\left\| \nabla^m \left( \eta_k\left( \nabla G-  \nabla {H} \right) \right) \right\|_{L^\infty(A_k)} 
\leq C\left(r2^k\right)^{-(m+1)}\exp\left( -c2^{2k}\right) + C R^{-1} \left(r2^k\right)^{-m}.
\end{equation*}
In particular, for each $k \in \{ 1,\ldots,k_*\}$,
\begin{align*}
\left(r2^k\right)^3 \left\|  \eta_k\left( \nabla G-  \nabla {H} \right) \right\|_{\underline{H}^{2}\left(B_{r2^k} \right)}
&
\leq
C\left(r2^k\right)^3 \left\|  \eta_k\left( \nabla G-  \nabla {H} \right) \right\|_{\underline{W}^{2,\infty}(A_k)}
\\ &
\leq 
C\left( R^{-1} r2^k +  \exp\left( -c2^{2k}\right)\right).  
\end{align*}
Applying~\eqref{e.rescalebtch}, we obtain, for each $k \in \{ 1,\ldots,k_*\}$,
\begin{equation*}
\left| \int_{\R^2} \eta_k(x)\left( \nabla G(x) -  \nabla {H} (x) \right) \cdot \nabla\phi_e(x) \,dx\right| 
\leq 
\O_s \left( C\left( R^{-1}r2^k +  \exp\left( -c2^{2k}\right)\right) \right). 
\end{equation*}
Using Lemma~\ref{l.sum-O} to sum the previous display over $k\in\{1,\ldots,k_*\}$, we obtain
\begin{align*}
\sum_{k=1}^{k_*} \left| \int_{\R^2} \eta_k(x) \left( \nabla {H}(x) - \nabla G(x) \right)\cdot \nabla \phi_e(x) \,dx \right| 
\leq \O_s (C). 
\end{align*}
This completes the proof of~\eqref{e.PsitoG}.

\smallskip

\emph{Step 4.} We summarize the results of previous steps. The claim is that 
\begin{align}
\label{e.gotossdyadic}
\left| 
(\phi_e\ast\Phi_r)(0) - (\phi_e\ast \Phi_R)(0)
- 
\sum_{k=1}^{k_*}
\int_{\R^2} \eta_k(x) \nabla G(x) \cdot \nabla \phi_e(x)\,dx
\right|
= 
\O_s\left(C \right). 
\end{align}
This follows from~\eqref{e.twopt.PGFidentity},~\eqref{e.Psidecomp},~\eqref{e.2pt.remover},~\eqref{e.PsitoG}
and the triangle inequality.

\smallskip

\emph{Step 5.}
We show that, for each $k\in\N$ and $T\in [1,\infty)$, 
\begin{align}
\label{e.CKish}
\lefteqn{
\int_{\R^2} \eta_k(x) \nabla G(x) \cdot \nabla \phi_e(x)\,dx
} \qquad & 
\notag
\\ &
=
\int_{\R^2}
\eta_k(x) \nabla G(x) \cdot
\left( \Phi\left(T,\cdot\right) \ast \nabla \phi_e\right)(x) \,dx
+ \O_s\left(C T^{\frac12} (r2^k)^{-1}\right). 
\end{align}
Set $h_{k,j}:= \eta_k \partial_{x_j} G$. We have
\begin{equation}
\label{e.hkjidentity}
h_{k,j} = \int_0^\infty \left( \Delta h_{k,j} \ast \Phi(t,\cdot)\right)(x)\,dt. 
\end{equation}
Fix $T\in [1,\infty)$ and compute 
\begin{align*}
\lefteqn{
\int_{\R^2} \eta_k(x) \nabla G(x) \cdot \nabla \phi_e(x)\,dx
} \quad & 
\\ & 
= \int_0^{T}
\int_{\R^2}
\Delta h_{k,j}(x)
\left(  \Phi(t,\cdot)\ast
\partial_{x_j} \phi_e\right)(x)
\,dt \,dx
\\ & \quad 
+ 
\int_{T}^\infty
\int_{\R^2}
\left( \Delta h_{k,j} \ast \Phi\left(t-T,\cdot\right) \right)(x)
\left( \Phi\left(T,\cdot\right) \ast \partial_{x_j} \phi_e\right)(x)
\,dt \,dx
\end{align*}
By~\eqref{e.hkjidentity}, the second integral on the right side of the previous display is equal to 
\begin{equation*}
\int_{\R^2}
\eta_k(x) \nabla G(x) \cdot
\left( \Phi\left(T,\cdot\right) \ast \nabla \phi_e\right)(x) \,dx,
\end{equation*}
and thus to establish~\eqref{e.CKish} it suffices to show that 
\begin{equation}
\label{e.CKish.wts}
\int_0^{T}
\int_{\R^2}
\Delta h_{k,j}(x)
\left(  \Phi(t,\cdot)\ast
\partial_{x_j} \phi_e\right)(x)
\,dt \,dx
=
\O_s\left(C T^{\frac12} (r2^k)^{-1}\right).
\end{equation}
We compute, using~\eqref{e.gradient},~\eqref{e.etakPsibouns} and Lemma~\ref{l.sum-O}, 
\begin{align*}
\lefteqn{
\left|
\int_0^{T}
\int_{\R^2}
\Delta h_{k,j}(x)
\left(  \Phi(t,\cdot)\ast
\partial_{x_j} \phi_e\right)(x)
\,dt \,dx
\right| 
} \qquad & 
\\ &
\leq 
\int_0^{T}
\int_{\R^2}
\left| \Delta h_{k,j}(x) \right| 
\left| \left(  \Phi(t,\cdot)\ast
\partial_{x_j} \phi_e\right)(x)\right| 
\,dt \,dx
\\ & 
\leq 
\O_s\left( 
\left| A_k \right|
\left\| \Delta h_{k,j}(x) \right\|_{L^\infty(A_k)}  \int_0^{T} t^{-\frac12} \,dt \right)
=
\O_s\left(C T^{\frac12} (r2^k)^{-1}\right).
\end{align*}
This completes the proof of~\eqref{e.CKish}. 

\smallskip

\emph{Step 6.} We localize to prepare for the use of the independence assumption. For each $k\geq1$, denote 
\begin{equation*}
\tilde{A}_k:= B_{r2^{k+1}} \setminus B_{r 2^{k-3}}.
\end{equation*}
We will construct, for each $k\in\N$, an $\F(\tilde{A}_k)$--measurable random variable $\mathsf{S}_k$ satisfying, for some $\delta(s)>0$ and $C(s,\Lambda)<\infty$, 
\begin{equation}
\label{e.yeslocal2dbtch}
\left| 
\mathsf{S}_k
-
\int_{\R^2} \eta_k(x) \nabla G(x) \cdot \nabla \phi_e(x)\,dx
\right| 
\leq \O_s\left( C (r2^k)^{-\delta} \right).
\end{equation}
We aim at applying~\eqref{e.CKish} with~$T = T_k := \left(r2^{k-3}\right)^{\frac{2}{1+\ep}}$ and~$\ep>0$ to be selected below. Define
\begin{align}
\label{e.tildeJ1}
\tilde{J}_{1}\left(\Psi,e,e_{j} \right) 
&
:= J_1\left(\Psi,-e,e_{j}\right) 
- J_1\left(\Psi,-e,0\right) 
+ J_1\left(\Psi,0,e_{j}\right) - e_{j} \cdot e.
\end{align}
By~\eqref{e.L2paramatching001} and Lemma~\ref{l.redwedding}, there exists~$\delta_1(s)\in \left(0,\tfrac12\right]$ such that
\begin{equation} 
\label{e.phiconvvstildeJ1}
\left| \left( \Phi(T_k,\cdot) \ast \partial_{x_{j}} \phi_{e} \right)(x)
-
\tilde{J}_{1}\left(\Phi(T_k,\cdot-x),e,e_{j}\right)  \right| \leq \O_{s} \left(C T_k^{- \frac12 (1+\delta_1)}  \right). 
\end{equation}
By Theorem~\ref{t.additivity}, taking $\delta_2(s) \in \left(0,\frac12 \delta_{1}\right]$ sufficiently small so that 
\begin{equation*} 
1+ \frac32 \delta_2 \leq \frac 2{s},
\end{equation*}
we obtain, for each $t\in [1,\infty)$, $x\in\Rd$, $e\in \partial B_1$ and $j\in\{1,2\}$, a random variable $\tilde{J}^{(\delta_2)}_1\left(t,x , e,e_{j} \right)$ satisfying
\begin{equation*}
\tilde J^{(\delta_2)}_1\left(t,x ,e,e_{j}\right) 
\quad 
\mbox{is $\F\left(B_{t^{\frac12 (1+\delta_2)}}(x)\right)$--measurable}
\end{equation*}
and
\begin{equation} \label{e.tildeJ1loc}
\left| \tilde{J}_{1}\left(\Phi(t,\cdot-x),e,e_{j}\right) 
- 
\tilde{J}^{(\delta_2)}_1 \left( t,x,e,e_{j}\right) 
\right| 
\leq \O_s\left( Ct^{-\frac12 \left( 1+ \frac32 \delta_2\right)}\right).  
\end{equation}
Choose now $\ep := \delta_2$ in the definition of $T_{k}$. Since $\delta_{2} \leq \frac12 \delta_{1}$, we obtain
\begin{equation*} 
T_k^{-\frac12 \left( 1+\delta_{1} \right)}  + T_k^{- \frac12 \left(1+\frac32 \delta_2 \right)} \leq C \left(r2^{k}\right)^{- \left( 1 + \frac{\delta_{2}}{4} \right) }
\end{equation*}
and
\begin{equation*} 
T_k^{\frac12 (1+\delta_2)} = r2^{k-3}. 
\end{equation*}
Define
\begin{equation*}
\mathsf{S}_k
:=
\sum_{j=1}^2
\int_{\R^2}
\eta_k(x) \partial_{x_j} G(x) \cdot
\tilde{J}^{(\delta_{2})}_1 \left( T_{k},x,e,e_{j}\right)  \, dx.
\end{equation*}
As the function $\eta_k$ is supported in $A_k$ and $A_k+B_{r2^{k-3}} \subseteq \tilde{A}_k$, we have that indeed
\begin{equation*}
\mathsf{S}_k
\quad \mbox{is $\F(\tilde{A}_k)$--measurable.}
\end{equation*}
Using~\eqref{e.phiconvvstildeJ1} and~\eqref{e.tildeJ1loc}, we obtain by Lemma~\ref{l.sum-O} that, for every $k\in\N$,
\begin{align*}
\lefteqn{
\left| 
\mathsf{S}_k
-
\int_{\R^2}
\eta_k(x) \nabla G(x) \cdot
\left( \Phi\left( T_{k} ,\cdot \right) \ast \nabla \phi_e\right)(x) \,dx
\right|
} \qquad & 
\\ & 
\leq
\sum_{j=1}^2
\int_{A_k}
\left| \nabla G(x) \right| 
\left| \left(  \Phi \left( T_{k} ,\cdot \right) \ast \partial_{x_j} \phi_e \right)(x) 
-
\tilde{J}^{(\delta_{2})}_1 \left( T_{k},x,e,e_{j}\right)  
\right|
\,dx
\\ & 
\leq
\O_s\left( C(r2^k)^{-\frac{\delta_{2}}{4}}\right).
\end{align*}
Combining this with~\eqref{e.CKish} yields~\eqref{e.yeslocal2dbtch} with $\delta = \frac14 \delta_{2}$. 

\smallskip

\emph{Step 7.} The conclusion. We have by~\eqref{e.gotossdyadic} and~\eqref{e.yeslocal2dbtch} that 
\begin{align}
\label{e.twopt.Sksum}
\left| 
(\phi_e\ast\Phi_r)(0) - (\phi_e\ast \Phi_R)(0)
- 
\sum_{k=1}^{k_*}
\mathsf{S}_k
\right|
= 
\O_s\left(C \right). 
\end{align}
Notice that~\eqref{e.2pt.remover} and~\eqref{e.yeslocal2dbtch} also imply that 
\begin{equation*}
\mathsf{S}_k = \O_s(C). 
\end{equation*}
By Lemma~\ref{l.bigO.barO}, we have that 
\begin{equation*}
\overline{\mathsf{S}}_k := \mathsf{S}_k - \E \left[ \mathsf{S}_k \right]
= \overline{\O}_s\left( C \right). 
\end{equation*}
By the unit range of dependence assumption, the random variables~$\overline{\mathsf{S}}_k$ and~$\overline{\mathsf{S}}_l$ are independent for every $k,l\in\N$ satisfying $|k-l|\geq 5$. Therefore, by Lemma~\ref{l.barO.boxes}, 
\begin{equation}
\label{e.Sbarksumsup}
\left| \sum_{k=1}^{k_*}
\overline{\mathsf{S}}_k  \right| 
\leq
\sum_{l=0}^4
\left| \sum_{k=1}^{\left\lfloor (k_*-l)/5 \right\rfloor}
\overline{\mathsf{S}}_{5k+l} \right| 
= \O_s\left( C k_*^{\frac12} \right)
= \O_s\left( C + C\log^{\frac12} \left( \frac Rr\right) \right).
\end{equation}
We also have from~\eqref{e.twopt.Sksum},~\eqref{e.twopt.PGFidentity} and~\eqref{e.biascorrectorsgrad} that 
\begin{align*}
\left| \sum_{k=1}^{k_*}
\E \left[ \mathsf{S}_k \right] \right| 
&
\leq 
C + 
\left| \E \left[ (\phi_e\ast\Phi_r)(0) - (\phi_e\ast \Phi_R)(0) \right] \right| 
\\ & 
= C + 
\left| \E \left[
\int_{r^2}^{R^2}\int_{\R^2}  \nabla \Phi\left( \tfrac t2,x\right)  \cdot \left( \Phi\left( \tfrac t2,\cdot\right) \ast \nabla \phi_e \right) (x) \,dx \,dt 
\right]
\right|
\\ & 
\leq
C + 
\int_{r^2}^{R^2}\int_{\R^2} \left| \nabla \Phi\left( \tfrac t2,x\right)  \right| \left| \E \left[ \left( \Phi\left( \tfrac t2,\cdot\right) \ast \nabla \phi_e \right) (x) \right| \,dx \,dt 
\right]
\\ & 
\leq C + C \exp\left( -cr\right) \leq C. 
\end{align*}
Combining this with~\eqref{e.Sbarksumsup} yields
\begin{equation*}
\left| \sum_{k=1}^{k_*}
\mathsf{S}_k  \right|  
\leq 
\O_s\left( C + C\log^{\frac12} \left( \frac Rr\right) \right).
\end{equation*}
In view of~\eqref{e.twopt.Sksum}, the proof of the proposition is now complete.   
\end{proof}

\begin{remark}
\label{r.Phidoesntmatter}
The choice of the heat kernel in the statement of the estimates~\eqref{e.phioscd=2} and~\eqref{e.phid=2} in Theorem~\ref{t.correctors} is not crucial. Let~$f\in H^2(\R^2)$ with $\int_{\R^2} f = 1$ and denote, for every $r\geq1$, $f_r(x):= r^{-2} f(x/r)$. The claim is that, for every $s\in (0,2)$, there exists $C(f,s,\Lambda)<\infty$ such that 
\begin{equation}
\label{e.Phi.doesntmatter}
\left| \left( \phi_e\ast \Phi_r \right)(x) - \left( \phi_e\ast f_r \right)(x) \right| 
\leq \O_s \left(C \right).
\end{equation}
To prove~\eqref{e.Phi.doesntmatter}, we argue as in the proof of Proposition~\ref{p.d2vengeance}. We use the identity
\begin{align*}
\lefteqn{
\left( \phi_e\ast \Phi_r \right)(0) - \left( \phi_e\ast f_r \right)(0)
} \quad & 
\\ & 
=
\int_0^\infty \int_{\R^2}\int_{\R^2}
\left( \Delta f_r(y) - \Delta \Phi_r(y) \right) \Phi(t,x-y) \phi_e(x) 
\,dx\,dy\,dt
\\ & 
=
- \int_0^\infty \int_{\R^2}\int_{\R^2}
\left( \nabla f_r(y) - \nabla \Phi_r(y) \right) \Phi(t,x-y) \nabla \phi_e(x) 
\,dx\,dy\,dt.
\end{align*}
We estimate
\begin{align*}
\lefteqn{
\left| \int_0^{r^2} \int_{\R^2}\int_{\R^2}
\left( \nabla f_r(y) - \nabla \Phi_r(y) \right) \Phi(t,x-y) \nabla \phi_e(x) 
\,dx\,dy\,dt \right|
} \quad & 
\\ & 
\leq  
\int_0^{r^2} \int_{\R^2}
\left| \nabla f_r(y) - \nabla \Phi_r(y) \right|
\left| \int_{\R^2} \Phi(t,x-y) \nabla \phi_e(x) 
\,dx \right| 
\,dy\,dt
\\ & 
\leq \O_s\left( C \left( \left\| \nabla f_r \right\|_{L^1(\R^2)} 
+
\left\| \nabla \Phi_r \right\|_{L^1(\R^2)}
 \right) \int_0^{r^2} t^{-\frac12}\,dt \right) 
\\ &
=
\O_s \left( C \left( 1 + \left\| \nabla f \right\|_{L^1(\R^2)} \right) \right). 
\end{align*}
To estimate the integral over~$({r^2},\infty)$, we let $w(t,x) := \left( (f_r - \Phi_r) \ast\Phi(t,\cdot ) \right)(x)$ and use the fact that $\int_{\R^2} (f_r - \Phi_r)=0$ implies that $w(t,x)$ has faster decay to zero by factor of $rt^{-\frac12}$: there exists~$C$ depending~on $f$ such that  
\begin{equation*}
\left| w(t,x) \right| \leq Crt^{-\frac 32} \exp\left( -\frac{|x|^2}{Ct} \right).
\end{equation*}
(We leave the proof of this estimate as an exercise to the reader, but see the discussion around~\eqref{e.decayfaster} for a hint.) Using this estimate leads to 
\begin{align*}
\lefteqn{ 
\left| \int_{r^2}^\infty \int_{\R^2}\int_{\R^2}
\left( \nabla f_r(y) - \nabla \Phi_r(y) \right) \Phi(t,x-y) \nabla \phi_e(x) 
\,dx\,dy\,dt \right|
} \quad & 
\\ & 
= 
\left| \int_{r^2}^{\infty} \int_{\R^2} \nabla w(t,x) \cdot \nabla \phi_e(x)\,dx\,dt \right| 
\\ &
\leq
\int_{r^2}^{\infty} \int_{\R^2} \left| w\left(\tfrac t2,y \right) \right|
\left| \int_{\R^2} \nabla \Phi\left(\tfrac t2,x-y \right)  \cdot \nabla \phi_e(x)\,dx\,dt \right| 
\\ & 
= \O_s\left( C \right). 
\end{align*}
Combining the above yields~\eqref{e.Phi.doesntmatter}. 

\smallskip

Note also that, by density, we only require that $f\in W^{1,1}$ in the argument above. In fact, it can also be shown by a variant of this computation that~\eqref{e.Phi.doesntmatter} is still valid if the regularity assumption on~$f$ is relaxed to~$f\in W^{\alpha,1}(\R^2)$ for some $\alpha>0$. In particular, we may take $f:= |B_1|^{-1} \indc_{B_1}$. We leave the verification of this claim as an exercise to the reader. 
\end{remark}

\index{corrector!first-order~$\phi_e$|)}

\section*{Notes and references}

Most of the statements of Theorems~\ref{t.correctors} and~\ref{t.additivity} were first proved in~\cite{AKM2}, and this chapter is based on that paper. 
Another, independent proof appeared later in~\cite{GO6} based on a quantity they call the ``homogenization commutator'' \index{homogenization commutator}
which is essentially equivalent to~$J_1$ (as can be readily seen from Lemma~\ref{l.redwedding} above) and originates from the ideas and philosophy of~\cite{AS,AKM1,AKM2}. The first successful use of this kind of bootstrap or renormalization argument to prove quantitative results in stochastic homogenization appeared previously in~\cite{AKM1}. These approaches are related to the strategy of Chapter~\ref{c.two}, which can be regarded as a more primitive form of this idea. In~\cite{AKM1}, the bootstrap was formalized for the quantity~$J$ from Chapter~\ref{c.two} rather than~$J_1$ and this produced a suboptimal estimate ($\alpha=1$ rather than $\alpha=\frac d2$), due to boundary layer effects.

\smallskip

Quantitative estimates for the first-order correctors with the same scaling~$O(r^{-\frac d2} )$ as in Theorem~\ref{t.correctors}, but with much weaker stochastic integrability, were proved several years earlier in~\cite{GO1,GO2,GNO,GO3} for a certain class of random environments close to the random checkerboard model. The methods in these papers were based not on renormalization but rather on nonlinear concentration inequalities---such as the so-called \emph{spectral gap} or \emph{logarithmic Sobolev} inequalities---an idea originating in statistical physics previously introduced into stochastic homogenization by Naddaf and Spencer~\cite{NS2,NS}. Optimal quantitative bounds on first-order correctors were extended to the setting of  supercritical percolation clusters in~\cite{Dario} by combining the methods described in this book with these nonlinear concentration inequalities.

\smallskip

Some months after a preliminary version of this book was posted to the arXiv, Bella, Giunti and Otto~\cite{BGO2} reproved Remark~\ref{r.canonical} above (by obtaining the full statement of our Theorem~\ref{t.regularity}) and extended the canonical isomorphism to one between $\A_{k+2} / \A_k$ and~$\Ahom_{k+2} / \Ahom_k$ in dimensions~$d>2$.



\chapter{Scaling limits of first-order correctors}
\label{c.gff}

In the previous chapter, we proved that the spatial averages of the gradients and fluxes of the correctors display stochastic cancellations similar to those of a stationary random field with a finite range of correlations. 
In this chapter, we take this analysis a step further by identifying the next-order behavior of correctors. 
We will show here that, as~$r\to\infty$, the random field~$r^{\frac d 2} \, \nabla \phi_e \Ll( r \,\cdot \Rr)$
converges in law to a correlated random field, which turns out to be the gradient of a Gaussian free field. Our analysis is inspired by the previous chapters, and we focus first on describing the next-order behavior of the quantity~$J_1$. Contrary to the gradient of the corrector, this quantity has short-range correlations even at the critical scaling, and is therefore more amenable to direct analysis. We show that~$J_1$ converges in law to a convolution of white noise, and consequently deduce, by a deterministic argument, the convergence in law of the rescaled corrector to a Gaussian free field. 

\smallskip

We begin the chapter by defining the notions of \emph{white noise} and of \emph{Gaussian free field} and providing explicit constructions of these random fields in Sections~\ref{s.gff.defs} and~\ref{s.gff.constr}. We then present an informal, heuristic argument for the convergence in law of the corrector in Section~\ref{s.heuristics}. In the final two sections, we make this heuristic derivation rigorous. We prove a central limit theorem for $J_1$ in Section~\ref{s.CLT}, and then deduce the convergence in law of the corrector as a consequence in Section~\ref{s.conv.gff}.

\section{White noise and the Gaussian free field}
\label{s.gff.defs}

In this section, we introduce the concepts of white noise and Gaussian free field. A white noise can be thought of as the limit in law of a centered random field with finite range of dependence, under the scaling of the central limit theorem. Indeed, we present in the next section an explicit construction of white noise based on this idea. A Gaussian free field can be seen as a generalization of Brownian motion where the time variable becomes multi-dimensional. An explicit construction of the Gaussian free field is also presented in the next section.

\smallskip

\index{white noise|(}
Informally, we think of a scalar white noise $(\W(x))_{x \in \Rd}$ of variance $\sigma^2 \in [0,\infty)$ as being a centered Gaussian random field such that, for every $x, y \in \Rd$,
\begin{equation*}  
\E[\W(x) \W(y)] = \sigma^2 \delta(x-y),
\end{equation*}
where $\delta$ is the $d$-dimensional Dirac distribution. Unfortunately, even if we could construct such a random field, it would not have Lebesgue measurable realizations. Scalar white noise is simply too rough to be a function. Let us continue however with formal computations based on this idea: setting, for each~$f \in L^2(\Rd)$, 
\begin{equation}
\label{e.def.Wf}
\W(f) := \int_{\Rd} \W(x) f(x) \, dx
\end{equation}
we obtain that
\begin{equation}  
\label{e.real.def}
\E[(\W(f))^2] = \sigma^2 \|f\|_{L^2(\Rd)}^2.
\end{equation}
This motivates the following definition.

\index{white noise!scalar}
\begin{definition}[Scalar white noise]
\label{d.scalar.noise}
A \emph{scalar white noise} over $\Rd$ of variance $\sigma^2 \in [0,\infty)$ is a family of real random variables $(\W(f), f \in L^2(\Rd))$ such that 
\begin{itemize}
\item for every $\lambda \in \R$ and $f, g \in L^2(\Rd)$, we have
\begin{equation}  
\label{e.W.linear}
\W(\lambda f + g) = \lambda \W(f) + \W(g) \qquad \text{a.s.},
\end{equation}
\item for every $f \in L^2(\Rd)$, the random variable $\W(f)$ is a centered Gaussian with variance $\sigma^2 \|f\|_{L^2(\Rd)}^2$. 
\end{itemize}
\end{definition}
Figure~\ref{f.whitenoise} (page \pageref{f.whitenoise}) illustrates the convolution of white noise by a smooth function with compact support. 

\smallskip

For notational convenience, in this section we denote by $\P$ the underlying probability measure, although it is unrelated with the law of the coefficient field $\a(x)$ appearing in most of the rest of the book.

\smallskip

It is important to observe that Definition~\ref{d.scalar.noise} completely determines the joint law of the vector $(\W(f_1),\ldots, \W(f_n))$, for every $f_1,\ldots,f_n \in L^2(\Rd)$. Indeed, for every $\lambda_1, \ldots, \lambda_n \in \R$, the law of 
\begin{equation*}  
\sum_{i = 1}^n \lambda_i \W(f_i) =  \W \Ll( \sum_{i = 1}^n \lambda_i f_i \Rr)
\end{equation*}
is a centered Gaussian with variance 
\begin{equation*}  
\sigma^2 \sum_{i,j = 1}^n \lambda_i \lambda_j \int_{\Rd} f_i f_j.
\end{equation*}
The claim of uniqueness is therefore a consequence of the second part of the following elementary lemma, which can be proved using Fourier transforms (see e.g.\ \cite[Theorem~29.4]{billingsley}
).
\index{Cram\'er-Wold lemma}
\begin{lemma}[Cram\'er-Wold lemma]
\label{l.cramer.wold}
\emph{(1)} Let $(X_n, n \in \N)$ and $X$ be random variables taking values in $\Rd$. If for each $\lambda \in \Rd$, 
\begin{equation*}  
\lambda \cdot X_n \xrightarrow[n \to \infty]{\mathrm{(law)}} \lambda \cdot X,
\end{equation*}
then 
\begin{equation*}  
X_n \xrightarrow[n \to \infty]{\mathrm{(law)}} X.
\end{equation*}
\emph{(2)} In particular, if $\lambda \cdot X_1$ and $\lambda \cdot X$ have the same law for every $\lambda \in \Rd$, then $X_1$ and $X$ have the same law.
\end{lemma}
\begin{remark}
\label{r.cramer.wold}
If we only know that, for each $\lambda \in \Rd$, the random variable $\lambda \cdot X_n$ converges in law, then it is automatic that the limit can be put in the form $\lambda \cdot X$ for some random variable $X$ not depending on $\lambda$. Indeed, the convergence in law of each of the coordinates of $X_n$ implies that the sequence $(X_n)_{n \in \N}$ is tight. By the second part of Lemma~\ref{l.cramer.wold}, any limit point of the law of $X_n$ is uniquely characterized. Therefore $X_n$ itself converges in law to some random variable $X$, and in particular $\lambda \cdot X_n$ converges in law to $\lambda \cdot X$.
\end{remark}
\begin{remark}  
\label{r.how.you.doing}
The first part of Lemma~\ref{l.cramer.wold} is useful to recover joint convergence in law of vectors of random variables from the individual convergence in law of scalar-valued random variables. In certain situations, we can also obtain joint convergence in law of random variables by first identifying the convergence in law of a subset of random variables, and then checking that other random variables are essentially deterministic functions of the former. More precisely, for arbitrary real or vector-valued random variables $(X_n)$, $(Y_n)$, $X$ and continuous function $g$, we have
\begin{equation}
\label{e.this.should.go}
(X_n,Y_n) \xrightarrow[n \to \infty]{\mathrm{(law)}} (X,g(X))
\end{equation}
if and only if
\begin{equation}  
\label{e.to.the.annals}
X_n \xrightarrow[n \to \infty]{\mathrm{(law)}} X \quad \text{and} \quad |Y_n - g(X_n)| \xrightarrow[n \to \infty]{\mathrm{(prob.)}} 0.
\end{equation}
Indeed, the implication \eqref{e.this.should.go} $\implies$ \eqref{e.to.the.annals} is clear. Conversely, in order to verify \eqref{e.this.should.go}, it suffices to show that, for every real-valued, uniformly continuous and bounded function $(x,y) \mapsto f(x,y)$, we have
\begin{equation*}  
\E[f(X_n,Y_n)] \xrightarrow[n \to \infty]{} \E\Ll[f(X,g(X))\Rr]. 
\end{equation*}
Assuming \eqref{e.to.the.annals}, this is easily verified using the decomposition
\begin{equation*}  
\E[f(X_n,Y_n)] = \E[f(X_n,Y_n) \1_{|Y_n - g(X_n)| \le \de}] + \E[f(X_n,Y_n) \1_{|Y_n - g(X_n)|  > \de}],
\end{equation*}
appealing to the uniform continuity of $f$ and taking $\de > 0$ sufficiently small.
\end{remark}
We deduce from Lemma~\ref{l.cramer.wold} that, for each $f_1,\ldots,f_k \in L^2(\Rd)$, the vector $(\W(f_1),\ldots,\W(f_k))$ is Gaussian with covariance matrix 
\begin{equation}
\label{e.cov.mat}
\Ll(\sigma^2 \int_{\Rd} f_i f_j\Rr)_{i,j \in \{ 1,\ldots k \}}. 
\end{equation}
As was already mentioned, we will provide with an explicit construction of white noise in the next section. One can alternatively appeal to an abstract extension argument to justify the existence of white noise. Indeed, we first observe that the matrix in \eqref{e.cov.mat} is nonnegative definite. For each $f_1,\ldots,f_k \in L^2(\Rd)$, it is therefore possible to construct a centered Gaussian vector $(\W(f_1),\ldots,\W(f_k))$ with covariance matrix given by \eqref{e.cov.mat}. By Kolmogorov's extension theorem (see e.g.\ \cite[Theorem~36.2]{billingsley}), we deduce the existence of a family of random variables $(\W(f), f \in L^2(\Rd))$ such that, for every $k \in \N$ and every $f_1,\ldots, f_k \in L^2(\Rd)$, the random variables $(\W(f_1),\ldots,\W(f_k))$ form a centered Gaussian vector with covariance \eqref{e.cov.mat}. It is then elementary to verify that the difference between the two sides of the identity \eqref{e.W.linear} has null variance, and therefore that this identity holds almost surely.

\smallskip 

We next extend Definition~\ref{d.scalar.noise} to higher dimensions. Given a nonnegative matrix $\msf Q \in \R^{n \times n}$, we would like a vector white noise with covariance matrix $\msf Q$ to be an~$\R^n$-valued, centered Gaussian random field satisfying, for every $i,j \in \{1,\ldots, n\}$ and $x, y \in \Rd$,
\begin{equation}  
\label{e.formal.W}
\E[\W_i(x) \W_j(y)] = \msf Q_{ij} \delta(x-y). 
\end{equation}
As above, for a vector field $\f = (\f_1,\ldots, \f_n) \in L^2(\Rd; \R^n)$, we may set
\begin{equation}  
\label{e.decomp.W}
\W(\f) := \W_1(\f_1) + \cdots + \W_n(\f_n)
\end{equation}
and integrate the formal identity \eqref{e.formal.W}, suggesting that $\W(\f)$ should be a centered Gaussian with variance $\int_{\Rd} \f \cdot \msf Q \f$.  

\index{white noise!vector}
\begin{definition}[Vector white noise]
\label{d.vector.noise}
Let $\msf Q \in \R^{n \times n}$ be a nonnegative symmetric matrix. An $n$-dimensional \emph{vector white noise} on $\Rd$ with covariance matrix $\msf Q$ is a family of real random variables $(\W(\f), \f \in L^2(\Rd;\R^n))$ such that 
\begin{itemize}
\item for every $\lambda \in \R$ and $\f, \g \in L^2(\Rd;\R^n)$, we have
\begin{equation*}  
\W(\lambda \f+ \g) = \lambda \W(\f) + \W(\g) \qquad \text{a.s.},
\end{equation*}
\item for every $\f \in L^2(\Rd;\R^n)$, the random variable $\W(\f)$ is a centered Gaussian with variance 
\begin{equation*}  
\int_{\Rd} \f \cdot \msf Q \f.
\end{equation*}
\end{itemize}
\end{definition}
\begin{exercise}  
Given a nonnegative symmetric matrix $\msf Q \in \R^{n \times n}$, use Kolmogorov's extension theorem to prove that a vector white noise $\W$ with covariance $\msf{Q}$ as in Definition~\ref{d.vector.noise} exists. For such $\W$ and for every $\f^{(1)}, \ldots, \f^{(k)} \in L^2(\Rd;\R^n)$, show that the vector
\begin{equation*}  
\Ll(\W(\f^{(1)}),\ldots, \W(\f^{k})\Rr)
\end{equation*}
is a centered Gaussian vector with covariance matrix
\begin{equation*}  
\Ll( \int_{\Rd} \f^{(i)} \cdot \msf Q \f^{(j)} \Rr)_{i,j \in \{1,\ldots,k\}}.
\end{equation*}
\end{exercise}
\begin{exercise}  
Verify that a vector white noise $\W$ admits the decomposition \eqref{e.decomp.W}, where $\W_1,\ldots,\W_n$ are scalar white noises.
\end{exercise}
\index{white noise|)}

\index{Gaussian free field|(}
We next introduce the concept of a Gaussian free field. Intuitively, we think of the Gaussian free field as a scalar random field on $\R^d$ which generalizes Brownian motion to higher dimensions. One way to characterize Brownian motion (up to an additive constant) is as the process on $\R$ whose derivative is a scalar white noise. With this in mind, a naive attempt to generalize Brownian motion to dimensions~$d>1$ could then be to look for a random field ${\Psi}$ such that $\nabla {\Psi}$ is a vector white noise. This obviously fails, however, since the requirement that $\nabla \Psi$ be a \emph{gradient field} imposes an additional constraint that a vector white noise does not satisfy: a vector white noise does not vanish when tested against smooth and compactly supported divergence-free fields. (In fact, we encounter a similar phenomenon in~$d = 1$ if we attempt for instance to define Brownian motion on the circle rather than on the full line.) A possible revision of our naive attempt could be to ask that~$\nabla {\Psi}$ is ``as close as possible'' to a white noise vector field. As was recalled in \eqref{e.helmoltz-hodge}, the Helmholtz-Hodge decomposition \index{Helmholtz-Hodge projection}
asserts that the space $L^2(\Rd;\Rd)$ is the orthogonal sum of the space of gradient fields and the space of divergence-free fields. Here, we take up a slightly more general point of view, and allow for a symmetric positive definite matrix $\ahom$ to mediate the orthogonality between vector fields. The decomposition then states that every $\f \in L^2(\Rd;\Rd)$ can be written as
\begin{equation}  
\label{e.helm.ahom}
\f = \ahom \nabla u + \g, \quad \mbox{where} \ u \in H^1(\Rd), \, \g \in \Ls(\Rd).
\end{equation}
The two terms~$ \ahom \nabla u$ and~$\g$ in this decomposition are orthogonal with respect to the scalar product 
\begin{equation} 
\label{e.scalarprodforgff}
(\f_1,\f_2) \mapsto \int_{\Rd} \f_1 \cdot \ahom^{-1} \f_2.
\end{equation}
Note that, given $\f\in L^2(\Rd;\Rd)$, the function $u$ in the decomposition can be obtained by solving the equation
\begin{equation}
\label{e.def.Proj2}
-\nabla \cdot \ahom \nabla u = - \nabla \cdot \f 
\quad \mbox{in} \ \Rd. 
\end{equation}
In view of this, and intending $\ahom \nabla {\Psi}$ to be the projection of a vector white noise $\W$ onto the space $\{\ahom \nabla u \ : \ u \in H^1(\Rd)\}$ with respect to the scalar product in~\eqref{e.scalarprodforgff}, we may expect that $\nabla \Psi$ should satisfy the equation
\begin{equation}  
\label{e.eq.gff}
-\nabla \cdot \ahom \nabla {\Psi} = - \nabla \cdot \W.
\end{equation}
Since white noise is a rather singular object (which we have not even made sense of as a random distribution so far), the meaning of this equation remains to be determined. We may formally test~\eqref{e.eq.gff} against $u \in C^\infty_c(\Rd)$ to get
\begin{equation}  
\label{e.tested.gff}
\int_{\Rd} \nabla {\Psi} \cdot \ahom \nabla u =  \int_{\Rd}  \W \cdot \nabla u .
\end{equation}
The quantity on the right side of \eqref{e.tested.gff} is our interpretation for $\W(\nabla u)$, see  \eqref{e.def.Wf} and \eqref{e.decomp.W}. The identity \eqref{e.tested.gff} can therefore be read as a requirement for how $\nabla {\Psi}$ should behave when tested against~$\ahom \nabla u$ with $u\in C^\infty_c(\Rd)$. Moreover, since we intend $\nabla {\Psi}$ to be a gradient, we require $\nabla {\Psi}$ to vanish when tested against sufficiently smooth divergence-free and compactly supported vector fields. By the Helmholtz-Hodge decomposition, these two requirements entirely specify $\nabla {\Psi}$. We can combine them into the single statement that, for $\f$ decomposed as in \eqref{e.helm.ahom},
\begin{equation*}  
\nabla {\Psi}(\f) = \W(\nabla u).
\end{equation*}
It is convenient to denote by $\Proj:L^2(\Rd;\Rd) \to L^2(\Rd;\Rd)$ the gradient of the solution operator for~\eqref{e.def.Proj2}, that is, $\Proj(\f) = \nabla u$ where $u \in H^1(\Rd)$ is the unique solution of~\eqref{e.def.Proj2}. 

\smallskip

The precise definition of gradient Gaussian free field we adopt is as follows.
\index{Gaussian free field!gradient}
\begin{definition}[Gradient GFF]
\label{d.gff}
Let $\ahom \in \R^{d \times d}$ be a symmetric and positive definite matrix and $\msf Q \in \R^{d \times d}$ be symmetric and nonnegative. A $d$-dimensional \emph{gradient Gaussian free field} (or \emph{gradient GFF} for short) on $\Rd$ associated with $\ahom$ and $\msf Q$ is a family of real random variables $(\nabla{\Psi}(\f), \f \in L^2(\Rd;\Rd))$ such that:
\begin{itemize}  
\item for every $\lambda \in \R$ and $\f, \g \in L^2(\Rd;\Rd)$, we have
\begin{equation*}  
\nabla {\Psi}(\lambda \f + \g) = \lambda \nabla {\Psi}(\f) + \nabla {\Psi}(\g) \qquad \text{a.s.}
\end{equation*}
\item for every $\f \in L^2(\Rd;\Rd)$, the random variable $\nabla {\Psi}(\f)$ is a centered Gaussian with variance
\begin{equation*}  
\int_{\Rd} \Ll(\Proj \f\Rr) \cdot \msf Q \Ll(\Proj \f\Rr).
\end{equation*}
\end{itemize}
\end{definition}

\begin{figure}[tb]
\label{fig.GFF}
\centering
\includegraphics[width=.42\linewidth]{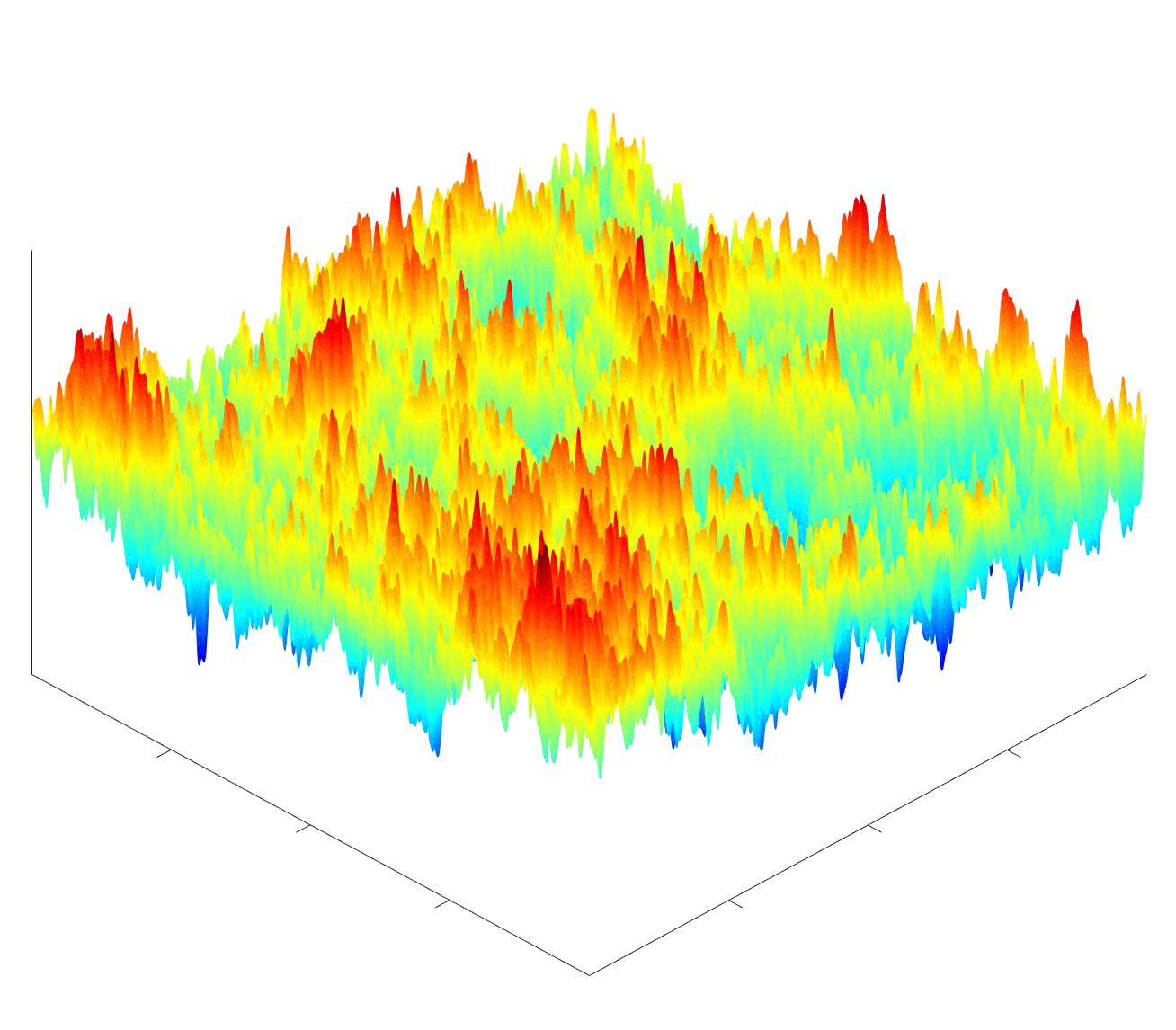}\par \vspace{-1.5em}
\includegraphics[width=.42\linewidth]{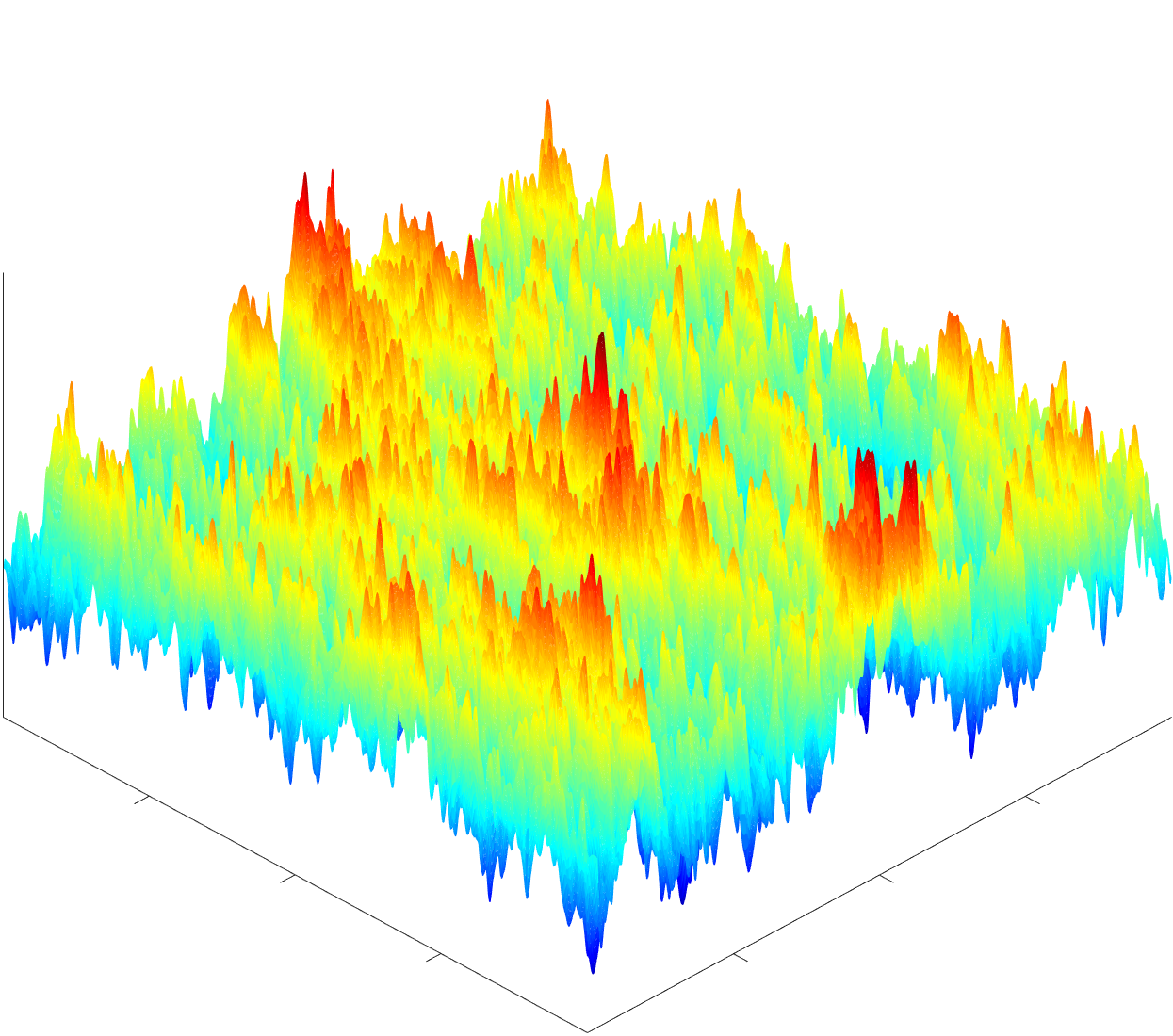} \qquad \
 \includegraphics[width=.42\linewidth]{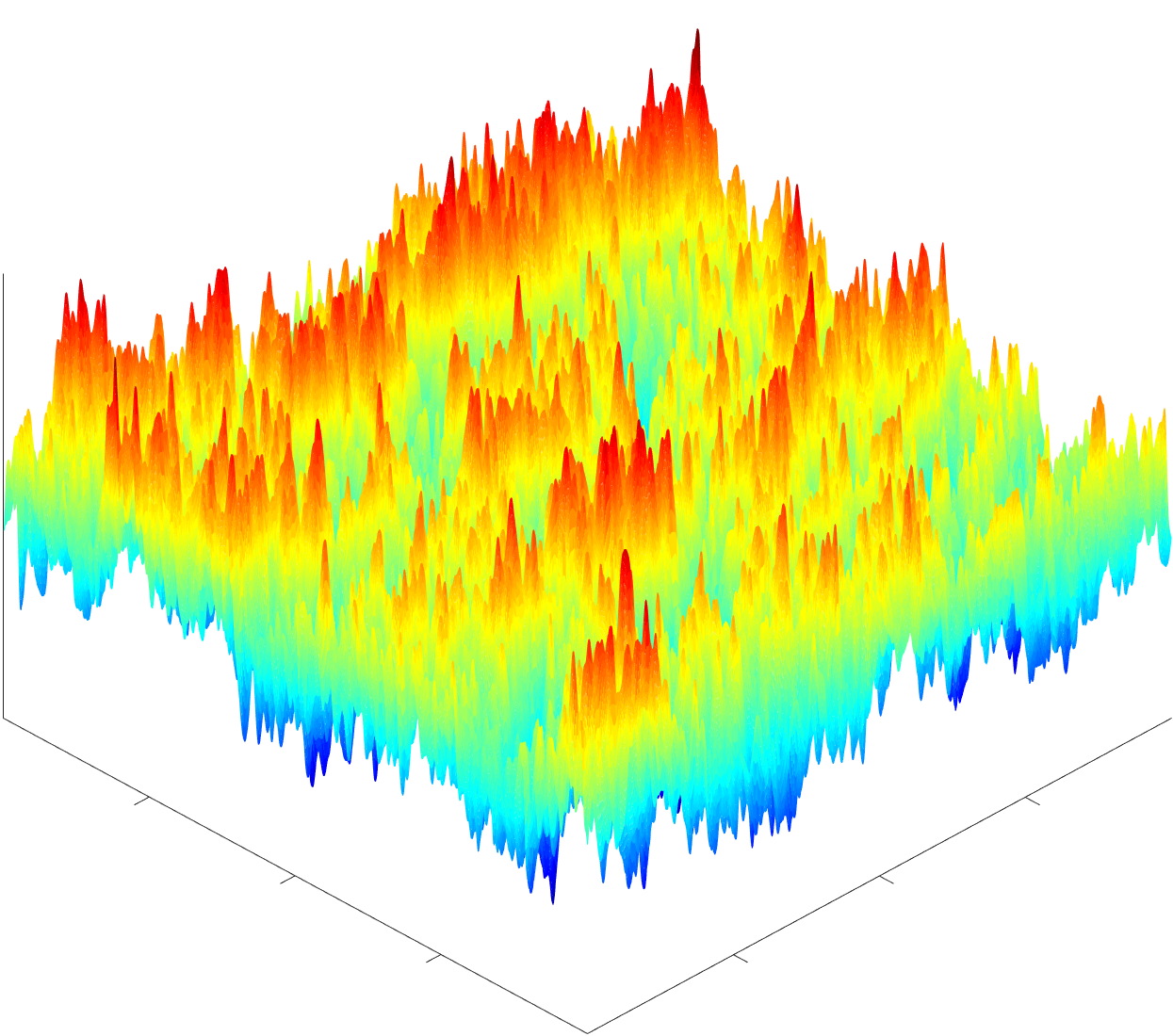}
\caption{\small 
The figures above are samples of the GFF convolved against a fixed smooth and compactly supported function, with $\bar \a = I_2$ and where the matrix~$\msf Q$ has been taken to be 
$I_2$, $e_1\otimes e_1$ and $e_2\otimes e_2$ respectively. Notice that in the first image, the field looks roughly isotropic in the sense that the mountain ranges do not have a preferred orientation. This is not the case in the last two images, where the mountain ranges seem to line up in the $e_1$ and $e_2$ directions respectively.}
\end{figure}

See Figure~\ref{fig.GFF} for samples of Gaussian free fields convolved with a given smooth and compactly supported function, for different values of the matrix $\msf Q$.

\smallskip

As for white noise, one can appeal to Kolmogorov's extension theorem to justify the existence of the gradient Gaussian free field. We provide with a more explicit construction in the next section.
\begin{exercise}  
\label{ex.scaling}
Let $\W$ be a vector white noise and $\lambda > 0$. We define the rescaled field $\W_\lambda(x) := \W(\lambda x)$ by setting, for every $\f \in L^2(\Rd;\R^n)$, 
\begin{equation*}  
\W_\lambda(\f) := \lambda^{-d} \, \W\Ll(x \mapsto \f(\lambda^{-1} x)\Rr).
\end{equation*}
Show that $\W_\lambda$ has the same law as $\lambda^{-\frac d 2} \, \W$. Letting $\nabla \Psi$ be a gradient GFF and defining $(\nabla \Psi)_{\lambda}(x) := (\nabla \Psi)(\lambda x )$ similarly, show that $(\nabla \Psi)_{\lambda}$ has the same law as $\lambda^{-\frac d 2} \, \nabla \Psi$. In particular, we should keep in mind that, as we zoom in to a small length scale $\lambda > 0$, these fields blow up like $\lambda^{-\frac d 2}$.
\end{exercise}
As the name and notation suggest, the gradient GFF $\nabla \Psi$ can be realized as the gradient of a random field $\Psi$. This is a consequence of the following elementary lemma.

\begin{lemma}
\label{l.mean-zero}
For every $f \in C^\infty_c(\Rd)$ of mean zero, there exists $F \in C^\infty_c(\Rd;\Rd)$ such that $\nabla \cdot F = f$.
\end{lemma}
\begin{proof} 
Without loss of generality, we may assume that $\supp f \subset B_1$. 
Let $h \in C^\infty_c(\R)$ be such that $\int_\R h = 1$ and $\supp h \subset B_1$. Writing $x = (x_1, \ldots, x_d) \in \Rd$, we define the functions
\begin{align*}  
g_0(x) & := f(x), \\
g_1(x) & := h(x_1)  \int_{\R}f(y_1,x_2,\ldots,x_d) \, dy_1,  \\
g_2(x) & := h(x_1) h(x_2) \int_{\R^2} f(y_1,y_2,x_3,\ldots,x_d) \, dy_1 \, dy_2, \\
& \vdots \\
g_d(x) & := h(x_1) \, \cdots \, h(x_d) \int_{\R^d} f = 0.
\end{align*}
Observe that $g_0, \ldots,g_d \in C^\infty_c(\Rd)$. For every $i \in \{1, \ldots, d\}$, we set
\begin{align*}  
F_i(x) := \int_{-\infty}^{x_i} (g_{i-1} - g_i)(x_1,\ldots,x_{i-1},y,x_{i+1}, \ldots,x_d)\, dy.
\end{align*}
Since
\begin{equation*}  
\partial_i F_i(x) = g_{i-1}(x) - g_i(x), 
\end{equation*}
we deduce that $\nabla \cdot F = f$. In order to complete the proof, there remains to verify that $F$ is compactly supported. Since $g_0,\ldots,g_d$ are compactly supported, it suffices to check that, for every $i \in \{1,\ldots,d\}$,
\begin{equation}  
\label{e.compact.supp}
x_i > 1 \quad \implies \quad F(x) = 0.
\end{equation}
The fact that 
\begin{equation*}  
x_i > 1 \quad \implies \quad F_j(x) = 0
\end{equation*}
for $j \neq i$ is clear, by the assumption of $\supp f \subset B_1$ (for $j < i$) and of $\supp h \subset B_1$ (for $j > i$). For the remaining case, if $x_i > 1$ then, since $\int_{\R} h = 1$ and $\supp h \subset B_1$, we have
\begin{align*}  
& \int_{-\infty}^{x_i} g_{i-1}(x_1,\ldots,x_{i-1},y,x_{i+1}, \ldots,x_d)\, dy \\
& \qquad = h(x_1) \, \cdots \, h(x_{i-1}) \int_{\R^i} f(y_1,\ldots,y_i,x_{i+1},\ldots,x_d) \, dy_1 \, \cdots \, dy_i \\
& \qquad = \int_{-\infty}^{x_i} g_{i}(x_1,\ldots,x_{i-1},y,x_{i+1}, \ldots,x_d)\, dy .
\end{align*}
This shows \eqref{e.compact.supp} and thereby completes the proof.
\end{proof}

Lemma~\ref{l.mean-zero} allows us to define an action of the Gaussian free field $\Psi$ on the set of mean-zero functions in~$C^\infty_c(\Rd)$ (the smoothness requirement can then be removed by a density argument). Indeed, for every function $f \in C^\infty_c(\Rd)$ of mean zero and $F \in C^\infty_c(\Rd;\Rd)$ such that $\nabla \cdot F = f$, we simply set
\begin{equation}  
\label{e.def.psi.uptoc}
\Psi(f) := - \nabla \Psi(F).
\end{equation}
This definition does not depend on the particular choice of $F \in C^\infty_c(\Rd;\Rd)$ satisfying $\nabla \cdot F = f$, since by definition of the projection operator $\Proj$, 
\begin{equation*}  
\nabla \cdot F = 0 \quad \implies \quad \nabla \Psi(F) = 0 \ \ \text{a.s.}
\end{equation*}
Letting $h \in C^\infty_c(\Rd)$ denote a fixed test function such that $\int_{\Rd} h = 1$, we can then extend this definition by setting, for $f, F$ as in Lemma~\ref{l.mean-zero} and an arbitrary $c \in \R$,
\begin{equation}  
\label{e.def.psi.ch}
\Psi(f + c h) := -\nabla \Psi(F).
\end{equation}
Naturally, every element of $C^\infty_c(\Rd)$ can be represented in the form of $f + ch$ for some $f \in C^\infty_c(\Rd)$ of mean zero and $c \in \R$. The relation \eqref{e.def.psi.ch} thus defines a field whose formal gradient is $\nabla \Psi$, and whose ``additive constant'' is fixed by the requirement that $\Psi$ should evaluate to $0$ when tested against $h$. 

\smallskip

The definition provided by~\eqref{e.def.psi.ch} is however not canonical, due to the arbitrary choice of the function $h$. In particular, the law of $\Psi$ is not invariant under translations: for any fixed $x_0 \in \Rd \setminus \{0\}$, the law of the translated field formally written as $x \mapsto \Psi(x + x_0)$ is different from the law of $\Psi$.

\begin{exercise}
\label{ex.gff.div.2d}
Let $\Psi$ be a Gaussian free field as defined by~\eqref{e.def.psi.uptoc} or~\eqref{e.def.psi.ch}. Fix $f \in C^\infty_c(\Rd)$ and, for every $r > 0$, set $f_r := r^{-d} f(r^{-1} \, \cdot)$. Show that there exists a constant $C(d) < \infty$ such that, for every $r > 0$,
\begin{equation*}  
\Psi(f_{2r}) - \Psi(f_{r}) = \O_2 \Ll( C r^{1-\frac d 2} \Rr) .
\end{equation*}
Deduce that, in the case $d \ge 3$, the sequence $(\Psi(f_{2^n}) - \Psi(f))_{n \in \N}$ is Cauchy with probability one. Conclude that, in~$d\geq 3$, one can construct a field $\Psi$ whose law is invariant under translations. (Hint: compare with the proof that the corrector exists as a stationary random field in Section~\ref{s.correctors}.) When $d = 2$, show that 
\begin{equation*}  
\Ll| \log(r) \Rr|^{-\frac 1 2}  \Ll( \Psi(f_r) - \Psi(f) \Rr) 
\end{equation*}
converges in law to a centered Gaussian, as $r$ tends to $0$ and as $r$ tends to infinity.
\end{exercise}

\begin{exercise}  
\index{Green function|(}
Suppose that $d \ge 3$ and denote by~$\bar G$ the Green function for the operator~$-\nabla \cdot \ahom \nabla$ in~$\Rd$. Equation~\eqref{e.eq.gff} suggests the formal definition
\begin{equation} 
\label{e.heur.psi}
\Psi(x) = \int_{\Rd} \nabla \bar G(x,z) \cdot \W(z) \, dz.
\end{equation}
Using the heuristics of \eqref{e.formal.W}, formally derive from \eqref{e.heur.psi} that, for every $x,y \in \Rd$,
\begin{equation}  
\label{e.correl.Psi}
\E \Ll[ \Psi(x) \Psi(y) \Rr] = \int_{\Rd} \nabla \bar G(x,y) \cdot \msf Q \nabla \bar G(z,y) \, dy.
\end{equation}
Use this property as a starting point for a direct definition of $\Psi$ as a stationary random field, for $d \ge 3$. Assuming further that the matrix $\msf Q$ is a scalar multiple of~$\ahom$, show by integration by parts that \eqref{e.correl.Psi} reduces to
\begin{equation}
\label{e.boring.GFF}
\E \Ll[ \Psi(x) \Psi(y) \Rr] = c \, \bar G(x,y),
\end{equation}
for some constant $c > 0$. 
\end{exercise}
\begin{remark}
\label{r.what.is.gff}
The Gaussian free field is usually defined so that \eqref{e.boring.GFF} holds, with suitable modifications when $d = 2$. Our notion of GFF is more general, and reduces to \eqref{e.boring.GFF} \emph{only} when the matrix $\msf Q$ is a scalar multiple of $\ahom$. One can visually distinguish between samples of the GFFs for different matrices $\msf Q$, see Figure~\ref{fig.GFF}.
\end{remark}

\index{Green function|)}
\index{Gaussian free field|)}

\section{Explicit constructions of random fields}
\label{s.gff.constr}

The definitions of white noise and of the gradient GFF given in the previous section are the weakest possible. Indeed, while Definitions~\ref{d.scalar.noise}, \ref{d.vector.noise} and \ref{d.gff} impose an uncountable number of constraints, each is a constraint for the joint law of the field tested against finitely many test functions. This weak definition allows for a short, abstract proof of existence of these objects, appealing to Kolmogorov's extension theorem as previously mentioned. It however leaves open the question of whether these fields can be realized as random elements valued in the space of distributions. The difficulty lies in the fact that the set of full probability measure over which the identity \eqref{e.W.linear} holds depends on the test functions $f, g \in L^2(\Rd)$, and there are uncountably many such test functions. Moreover, we cannot immediately use a density argument to circumvent this issue, since~$\W(\f)$ will not be bounded uniformly over $\| f \|_{L^2(\Rd)} \leq 1$, as will be explained shortly.

\smallskip

It is possible to complete the construction of white noise (or of the gradient GFF) as a random distribution by appealing to a generalization of Kolmogorov's continuity theorem (see e.g.\ \cite[Theorem~I.2.1]{revuz-yor} for the classical theorem, and \cite[Proposition~2.29]{tightness} for a possible generalization to distributions). However, we prefer to adopt a more direct approach. The goal of this section is to present an explicit construction of white noise, in the spirit of L\'evy's construction of Brownian motion as the limit of continuous piecewise linear approximations. This bypasses appeals to the Kolmogorov extension and continuity theorems, and exposes simple approximations of white noise involving a locally finite number of independent Gaussian random variables. We conclude the section by constructing the gradient GFF as a deterministic function of white noise, in agreement with \eqref{e.eq.gff}.

\smallskip

To which spaces of distributions can we expect a scalar white noise~$\W$ to belong? If~$\W$ has unit variance, then it is immediate from the definition that there exists a universal $C<\infty$ such that, for every $f \in L^2(\Rd)$, 
\begin{equation*}  
\W(f) = \O_2 \Ll( C \|f\|_{L^2(\Rd)} \Rr).
\end{equation*}
A superficial look at this relation may suggest that $\W(f)$ belongs to $L^2(\Rd)$, ``by duality.'' This is false for several reasons, the most obvious one being that since the law of $\W$ is invariant under translation, the quantity $\E[\|\W\|_{L^2(\Rd)}^2]$ can only be infinite. It turns out that $\W$ cannot be realized in $L^2_{\mathrm{loc}}(\Rd)$ either, nor even in $W^{-\al,p}_{\mathrm{loc}}(\Rd)$ for any $\al \le \frac d 2$ and $p\in[1,\infty]$, as will be explained in Lemma~\ref{l.unbounded.W} and Exercise~\ref{ex.unbounded.critical} below. 

\smallskip

Our goal then is to prove the following proposition, which gives the optimal Sobolev space in which scalar white noise can be realized.

\index{white noise!construction}
\begin{proposition}[Construction of white noise]
\label{p.white.noise.exists} 
Let $\al > \frac d 2$. There exists a random element $\W$ of $W^{-\al,\infty}_\mathrm{loc}(\Rd)$ such that, for every $f \in W^{\al,1}(\Rd)$ with compact support, the evaluation of $\W$ against $f$, denoted by $\W(f)$, is a centered Gaussian with variance $\|f\|_{L^2(\Rd)}^2$. Moreover, there exists a constant $C(d) < \infty$ such that 
\begin{equation}
\label{e.white.noise.estimate}
\sup \Ll\{\W(f) \ : \ \|f\|_{W^{\al,1}(\Rd)} \le 1 \ \text{and} \ \supp f \subset B_1 \Rr\} \le \O_2 \Ll( C  \Rr) .
\end{equation}
\end{proposition}
\begin{remark}
\label{r.measurability}
We do not discuss the measurability of the supremum on the left side of \eqref{e.white.noise.estimate}. We simply understand \eqref{e.white.noise.estimate} as the statement that there exists a random variable $\X$ taking values in $[0,+\infty]$ such that the inequality
\begin{equation*}  
\sup \Ll\{\W(f) \ : \ \|f\|_{W^{\al,1}(\Rd)} \le 1 \ \text{and} \ \supp f \subset B_1 \Rr\} \le \X
\end{equation*}
holds everywhere in the probability space, and moreover, 
\begin{equation*}  
\X \le \O_2 \Ll( C  \Rr) .
\end{equation*}
\end{remark}
Since $\W$ from Proposition~\ref{p.white.noise.exists} is a random distribution, it is clear that with probability one,
\begin{multline}  
\label{e.W.is.a.distribution}
\forall f,g \in W^{\al,1}(\Rd) \text{ with compact support}, \ \forall \lambda \in \R, 
\\ 
\W(\lambda f + g) = \lambda \W(f) + \W(g).
\end{multline}
Note the differences between this statement and \eqref{e.W.linear}: the statement above holds with probability one over every allowed test functions simultaneously; but the set of allowed test functions is smaller. 

\smallskip

One can a posteriori extend the definition of $\W(f)$ to every $f \in L^2(\Rd)$, if desired. Indeed, for any given $f \in L^2(\Rd)$, let $(f_n)_{n \in \N}$ be a (deterministic) sequence of compactly supported smooth functions converging to $f$ in $L^2(\Rd)$. Since
\begin{equation*}  
\E \Ll[ (\W(f_n) - \W(f_m))^2 \Rr]  =  \|f_n - f_m\|_{L^2(\Rd)}^2,
\end{equation*}
the sequence $(\W(f_n))_{n \in \N}$ is Cauchy in $L^2(\P)$. It therefore converges to a limit random variable, which we may still denote by $\W(f)$. Once $\W$ is thus extended, the verification of Definition~\ref{d.scalar.noise} is straightforward. Yet, note that the random variable $\W(f)$ is only well-defined up to a set of null probability measure, 
and this set may depend on the function $f$. Because of this, one cannot replace the set of test functions $W^{\al,1}(\Rd)$ in \eqref{e.W.is.a.distribution} by $L^2(\Rd)$. 

\smallskip

The difficulty we are describing here is perhaps best visualized through an orthonormal basis $(\mathbf{e}_k)_{k \in \N}$ of $L^2([0,1])$. Using the properties of white noise and Lemma~\ref{l.cramer.wold}, one can verify that $(\W(\mathbf{e}_k))_{k \in \N}$ is a sequence of independent standard Gaussian variables. Hence, for any fixed sequence $(\lambda_k)_{k \in \N}$ of real numbers satisfying
\begin{equation}  
\label{e.variance.W}
\sum_{k = 0}^\infty \lambda_k^2 < \infty,
\end{equation}
we can define
\begin{equation*}  
\W \Ll( \sum_{k= 0}^{\infty} \lambda_k \mathbf{e}_k \Rr) 
\end{equation*}
as the $L^2(\P)$-convergent series
\begin{equation*}  
\sum_{k = 0}^{\infty}\lambda_k \, \W \Ll( \mathbf{e}_k \Rr) ,
\end{equation*}
which is a Gaussian random variable with variance given by \eqref{e.variance.W}.
On the other hand, we have
\begin{equation*}  
\sup \Ll\{ \sum_{k = 0}^{\infty} \lambda_k \W(\mathbf{e}_k) \ : \ \sum_{k = 0}^{\infty} \lambda_k^2 \le 1 \Rr\} = \Ll(\sum_{k = 0}^{\infty} \W(\mathbf{e}_k)^2\Rr)^\frac 1 2 = +\infty \quad \text{a.s.}
\end{equation*}
In other words, we can always find a \emph{random} choice for the sequence $(\lambda_k)$ such that the left side of \eqref{e.variance.W} is bounded by $1$, and yet
\begin{equation*}  
\sum_{k = 0}^\infty \lambda_k \W(\mathbf e_k) = +\infty \quad \text{a.s.}
\end{equation*}
The next lemma is a more precise version of this argument. 
\begin{lemma}
\label{l.unbounded.W}
Let $\al < \frac d 2$. If $\W$ is a random distribution such that, for every $f \in C^\infty_c(\Rd)$, the evaluation of $\W$ against $f$ is a centered Gaussian with variance $\|f\|_{L^2(\Rd)}^2$, then the event that $\W$ belongs to $W^{-\al,1}_{\mathrm{loc}}(\Rd)$ is of null probability measure.
\end{lemma}

\begin{remark}
Our usage of the word ``event'' for the event that $\W$ belongs to $W^{-\al,1}_{\mathrm{loc}}(\Rd)$ is slightly abusive, since we do not discuss the measurability of this subset of the probability space. We simply understand the conclusion of Lemma~\ref{l.unbounded.W} as stating that this subset of the probability space is contained in a measurable set of null probability measure. Alternatively, one can enlarge the $\sigma$-algebra so that it contains all such subsets. 
\end{remark}

\begin{proof}[Proof of Lemma~\ref{l.unbounded.W}] 
Fix $\al \in \left(0,\frac d 2\right)$ and a nonzero smooth function~$f \in C^\infty_c(\Rd)$ supported in~$B_{1/2}$. For every $z \in \R$ and  $\eps > 0$, we consider 
\begin{equation*}  
f_{\eps,z} := \eps^\frac d 2 f \Ll( \eps^{-1} (\cdot - z) \Rr) . 
\end{equation*}
By scaling, there exists~$C(f,\al,d) < \infty$ such that, for every $z \in \Rd$ and $\eps \in (0,1]$, 
\begin{equation*}  
\|f_{\eps,z}\|_{W^{\al,\infty}(\Rd)} \le C \eps^{\frac d 2 - \al}.
\end{equation*}
Denote
\begin{equation*}  
\mcl P_\eps := \eps \Zd \cap \Ll( -\frac 1 2,\frac 1 2 \Rr)^d.
\end{equation*}
For every $\eps \in (0,1]$, we have
\begin{equation}  
\label{e.s.norm}
\sup_{(s_z) \in \{-1,1\}^{\mcl P_\eps}}\Ll\|\sum_{z \in \mcl P_\eps} s_z f_{\eps,z}\Rr\|_{W^{\al,\infty}(\Rd)} \le C \eps^{\frac d 2 - \al} \xrightarrow[\eps \to 0]{} 0,
\end{equation}
since the supports of the functions $(f_{\eps,z})_{z \in \mcl P_{\eps}}$ are disjoint. For the same reason, 
\begin{equation*}  
\sup_{(s_z) \in \{-1,1\}^{\mcl P_{\eps}}} \W \Ll( \sum_{z \in \mcl P_\eps} s_z f_{\eps,z} \Rr) = \sum_{z \in \mcl P_\eps} \Ll| \W \Ll( f_{\eps,z} \Rr)  \Rr| .
\end{equation*}
On the event that $\W$ belongs to~$W^{-\al,1}_{\mathrm{loc}}(\Rd)$, the quantity above must tend to $0$ almost surely as $\eps$ tends to $0$, by \eqref{e.s.norm}. However, $\W(f_{\eps,z})$ is a centered Gaussian with variance $\eps^{-d} \|f\|_{L^2(\Rd)}^2$. In particular, by the law of large numbers, 
\begin{equation*}  
\lim_{\eps \to 0} \P \Ll[ \sum_{z \in \mcl P_\eps} \Ll| \W \Ll( f_{\eps,z} \Rr)  \Rr| \ge \frac{\|f\|_{L^2(\Rd)}^2}{2} \Rr] = 1.
\end{equation*}
This completes the proof.
\end{proof}
\begin{exercise}  
\label{ex.unbounded.critical}
In this exercise, we extend the result of Lemma~\ref{l.unbounded.W} to the critical case $\al = \frac d 2$. We rely on the characterizations of $W^{-\al,p}(\Rd)$ in terms of spatial averages provided by Appendix~\ref{a.MSP}. However, the case when $\al$ is an integer and $p = 1$ is not covered by the results of the appendix. Thus the goal is to show that, under the assumptions of Lemma~\ref{l.unbounded.W}, the event that $\W$ belongs to $W^{-\frac d 2,p}_{\mathrm{loc}}(\Rd)$ is of null probability measure, where 
\begin{equation*}
\Ll\{
\begin{aligned}  
& p = 1 &  \mbox{if} & \  d\not\in2\N, \\
& p \in (1,\infty) & \mbox{if} & \  d \in 2\N.
\end{aligned}
\Rr.
\end{equation*}

\begin{enumerate}
\item Show that, for every nonnegative random variable $X$ and $\theta \in [0,1]$, we have
\begin{equation}  
\label{e.hint.use.Holder}
\P \Ll[ X > \theta {\E[X]} \Rr] \ge (1-\theta)^2 \frac{\E[X]^2}{\E[X^2]}.
\end{equation}
(This inequality is sometimes called the Paley-Zygmund inequality.)

\item Show that, for every $p \in [1,\infty)$, there exists a constant $c_p < \infty$ such that, for every centered Gaussian random variable $X$,
\begin{equation}  
\label{e.wiener.chaos}
\var \Ll[ |X|^p \Rr] = c_p \, \E \Ll[ |X|^p \Rr]^2.
\end{equation}
\item We fix $\chi \in C^\infty_c(\Rd)$ a smooth function with compact support such that $\chi \equiv 1$ on $B_1$. Show that there exists a constant $c > 0$ such that, for every $x \in B_{1/2}$ and $\eps \in (0,\frac 1 2]$, we have
\begin{equation}  
\label{e.unbounded.odd.d}
\E\Ll[\int_\eps^1 t^\frac d 4 \, \int_{\Rd} |\W(\chi \Phi(t,\cdot - y))| \, d y \, \frac{dt}{t} \Rr] \ge c \log \eps^{-1}
\end{equation}
and
\begin{equation}
\label{e.unbounded.even.d}
\E\Ll[\int_\eps^1 t^\frac d 2 \, |\W(\chi \Phi(t,\cdot - x))|^2 \, \frac{dt}{t} \Rr] \ge c \log \eps^{-1}.
\end{equation}

\item When the dimension $d$ is odd, conclude using Proposition~\ref{p.MSP.Walpha}, Chebyshev's and Minkowski's inequalities, \eqref{e.wiener.chaos} and \eqref{e.unbounded.odd.d}. If $d$ is even, use instead Proposition~\ref{p.MSP_int}, \eqref{e.hint.use.Holder}, Minkowski's inequality, \eqref{e.wiener.chaos} and \eqref{e.unbounded.even.d}.
\end{enumerate}
\end{exercise}

\begin{remark}
One may wonder where white noise lies in finer function spaces, such as Besov spaces. For every $\al \in \R$ and $p,q \in [1,\infty]$, we can define the Besov space $B^\al_{p,q}(\Rd)$ as the space of distributions $f$ such that
\begin{equation*}  
\|f\|_{B^\al_{p,q}(\Rd)} := 
\left( \int_0^1 \left( t^{-\frac\alpha2} \left\| f \ast \Phi(t,\cdot) \right\|_{L^p(\Rd)}\right)^q \,\frac{dt}{t} \right)^{\frac1q}
\end{equation*}
is finite. Note that, in view of Proposition~\ref{p.MSP.Walpha}, for every $\al \in (0,\infty) \setminus \N$ and $p\in (1,\infty)$, we have $W^{-\al,p}(\Rd) = B^{-\al}_{p,p}(\Rd)$. One can verify that white noise belongs to $B^{-\al}_{p,q,\mathrm{loc}}(\Rd)$ if and only if
\begin{equation*}  
\al > \frac d 2 \quad \text{or} \quad \Ll( \al = \frac d 2, \ p < \infty \text{ and } q = \infty \Rr) .
\end{equation*}
In the latter case (and whenever $p$ or $q$ is infinite), the space $B^\al_{p,q}(\Rd)$ we have defined is strictly larger than the completion of the set of smooth functions with respect to the norm $\|\cdot\|_{B^\al_{p,q}(\Rd)}$. White noise belongs only to the local version of the larger of these two sets.
\end{remark}

We now focus on the proof of Proposition~\ref{p.white.noise.exists}, that is, on the construction of white noise. Departing from the rest of the book, we prefer to work here with dyadic instead of tryadic cubes. We define, for every $n \in \N$,
\begin{equation}
\label{e.dyadic}
D_n := \Ll(0, 2^{-n}\Rr)^d \subset \Rd.
\end{equation}
For each~$n\in\N$, the family $(z + D_n, \ z \in 2^{-n} \Z^d)$ is a partition of $\Rd$ up to a Lebesgue null set. Given a family $(X_z)_{z \in 2^{-n} \Z^d}$ of independent standard Gaussian random variables, consider the random function over $\Rd$
\begin{equation*}  
\td \W_n :=  \sum_{z \in 2^{-n} \Z^d} |D_n|^{-\frac 1 2} \,  X_z \, \1_{z + D_n}.
\end{equation*}
This provides us with an approximation of white noise on scale $2^{-n}$. Indeed, for every $f \in L^2(\Rd)$, the random variable
\begin{equation*}
\td \W_n(f) := \int_{\Rd} \td\W_n(x)  f(x) \, d x
\end{equation*}
is a centered Gaussian with variance
\begin{equation}  
\label{e.variance}
\sum_{z \in 2^{-n} \Z^d} |D_n|\Ll(\fint_{z+D_n} f\Rr)^2.
\end{equation}
In particular, this quantity equals $\|f\|_{L^2(\Rd)}^2$ if $f$ is constant over each dyadic cube of the form $z + D_n$, for $z \in 2^{-n} \Z^d$. The point of the proof of Proposition~\ref{p.white.noise.exists} is to construct these approximations consistently for different values of $n$ and then to verify the convergence to a limiting object as $n\to\infty$. 

\begin{proof}[Proof of Proposition~\ref{p.white.noise.exists}]
We break the proof into two steps. In the first step, we construct a consistent sequence of approximations of white noise at scale~$2^{-n}$ and then, in the second step, show that these approximations are Cauchy and estimate the rate of convergence. 

\smallskip

\emph{Step 1.} Define $f^{(0)}, f^{(1)} \in L^2(\R)$ by
\begin{equation*}  
f^{(0)} := \1_{\Ll( 0 , 1 \Rr) } \quad \text{and} \quad f^{(1)} := \1_{\Ll( 0,\frac 1 2 \Rr) } - \1_{\Ll( \frac 1 2,1 \Rr) }.
\end{equation*}
and consider the family of functions on $\Rd$
\begin{equation*}  
x = (x_1,\ldots,x_d) \mapsto \prod_{k = 1}^d f^{(i_k)}(x_k),
\end{equation*}
indexed by all sequences $(i_1,\ldots,i_d) \in \{0,1\}^d$. Denote these functions by $\chi^{(0)}, \ldots$,  $\chi^{(2^d-1)}$, with
\begin{equation*}  
\chi^{(0)} := \1_{D_0}
\end{equation*}
corresponding to the choice of sequence $(i_1,\ldots,i_d) = (0,\ldots,0)$. See Figure~\ref{f.wavelets} for a representation of these functions in two dimensions.

\begin{figure}[tb]
\centering
\begin{tikzpicture}[scale=1]
\foreach \x in {0,3.5,7,10.5}
{
\draw[thick, ->] (\x,0)--(\x+2.5,0);
\draw[thick, ->] (\x,0)--(\x,2.5);
\draw[xshift=\x cm] (0,0) grid (2,2);
}
\draw (0.5,0.5) node[scale=2] {$+$};
\draw (1.5,0.5) node[scale=2] {$+$};
\draw (0.5,1.5) node[scale=2] {$+$};
\draw (1.5,1.5) node[scale=2] {$+$};

\draw (4,0.5) node[scale=2] {$+$};
\draw (5,0.5) node[scale=2] {$+$};
\draw (4,1.5) node[scale=2] {$-$};
\draw (5,1.5) node[scale=2] {$-$};

\draw (7.5,0.5) node[scale=2] {$+$};
\draw (8.5,0.5) node[scale=2] {$-$};
\draw (7.5,1.5) node[scale=2] {$+$};
\draw (8.5,1.5) node[scale=2] {$-$};

\draw (11,0.5) node[scale=2] {$+$};
\draw (12,0.5) node[scale=2] {$-$};
\draw (11,1.5) node[scale=2] {$-$};
\draw (12,1.5) node[scale=2] {$+$};

\end{tikzpicture}
\caption{
\small{
The functions $\chi^{(i)}$, $i \in \{0,1,2,3\}$, in dimension $d = 2$. The plus and minus signs indicate that the function takes the value $1$ and $-1$ respectively.
}
}
\label{f.wavelets}
\end{figure}

The functions $(\chi^{(i)})_{i < 2^d}$ form an orthonormal basis of the subspace of $L^2(\Rd)$ spanned by the indicator functions $\1_{z + D_1}$, for $z$ ranging in $\Ll\{ 0,\frac 1 2 \Rr\}^d$.  For each $n \in \N$ and $i < 2^d$, we define the rescaled functions
\begin{equation}  
\label{e.def.chin}
\chi^{(i)}_n := 2^{\frac {nd}{2}} \chi^{(i)} \Ll( 2^{n} \, \cdot \Rr) .
\end{equation}
The normalization is chosen so that, for each $n$, the family $\Ll( \chi^{(i)}_n  \Rr) _{i < 2^d}$ is orthonormal in $L^2(\Rd)$. As a consequence, the family
\begin{equation}  
\label{e.orth.family}
\Ll\{\chi^{(0)}(\cdot - z), \ z \in \Zd\Rr\} \cup \Ll\{ \chi^{(i)}_k(\cdot - z), \ 0 \le k < n, \ z \in 2^{-k} \Z^d, \ 1 \le i < 2^d\Rr\}
\end{equation}
is an orthonormal basis of the vector space spanned by the functions
\begin{equation*}  
\{\1_{z + D_n}, \ z \in 2^{-n} \Z^d\}.
\end{equation*}
Let $\Ll\{ X^{(i)}_{n,z} \, : \, n \in \N, z \in 2^{-n} \Z^d, 0 \le i < 2^d\Rr\}$ be a family of independent Gaussian random variables. For $n \in \N$, define a random function $\W_n$ on $\Rd$ by
\begin{equation*}  
\W_n :=  \sum_{z \in \Z^d} X^{(0)}_{0,z} \,  \chi^{(0)}(\cdot - z) + \sum_{k = 0}^{n-1} \sum_{\substack{z \in 2^{-k} \Z^d \\ 1 \le i < 2^d}}  X^{(i)}_{k,z} \, \chi^{(i)}_k(\cdot - z).
\end{equation*}
In view of the properties of the family in~\eqref{e.orth.family} reviewed above, we deduce that, for every $f \in L^2(\Rd)$, the random variable 
\begin{equation*}
\W_n(f) := \int_{\Rd} \W_n(x) f(x) \, d x
\end{equation*}
is a centered Gaussian with variance given by \eqref{e.variance}. By the Lebesgue differentiation theorem, we deduce that, for every~$f \in L^2(\Rd)$, the sequence of random variables~$(\W_n(f))_{n \in \N}$ converges in law to a centered Gaussian with variance~$\|f\|_{L^2(\Rd)}^2$. 

\smallskip

\emph{Step 2.} Let $\al \in \Ll( 0,\frac d 2 \Rr) \setminus \N$. In this step, we show that $\W_n$ is a Cauchy sequence in $W^{-\al,\infty}_\mathrm{loc}(\Rd)$, and estimate the rate of convergence. When testing against functions with compact support in $D_0$, we may as well replace $\W_n$ by
\begin{equation*}  
\W_n' := X^{(0)}_{0,0} \, \chi^{(0)} + \sum_{k = 0}^{n-1} \sum_{\substack{z \in Z_k \\ 1 \le i < 2^d}}  X^{(i)}_{k,z}  \, \chi^{(i)}_k(\cdot - z),
\end{equation*}
where we denote
\begin{equation*}  
Z_k := \Ll\{ z \in \Rd \ : \ 2^{k} z \in \{0,\ldots,2^{k}-1\}^d \Rr\}.
\end{equation*}
We therefore focus on proving that $(\W_n')_{n \in \N}$ is, almost surely, a Cauchy sequence in $W^{-\al,\infty}(\Rd)$. More precisely, setting
\begin{equation}
\label{e.def.be.w}
\be := \Ll( \al - \frac d 2 \Rr) \wedge 1,
\end{equation}
we aim to show that there exists $C(\al,d) < \infty$ such that
\begin{equation}  
\label{e.bound.diff.W'}
\|\W'_{n+1} - \W'_n\|_{W^{-\al,\infty}(\Rd)} \le 2^{-\be n } \Ll( C n^{\frac12} + \O_2(C) \Rr).
\end{equation}
By Lemma~\ref{l.sum-O} and the fact that $\be > 0$, this implies that $(\W'_n)_{n \in\N}$ is almost surely a Cauchy sequence in $W^{-\al,\infty}(\Rd)$. Since 
\begin{equation*}  
\|\W'_0\|_{W^{-\al,\infty}(\Rd)} \le \O_2(C),
\end{equation*}
it also implies \eqref{e.white.noise.estimate} for the limit $\W$. Hence, the proof of Proposition~\ref{p.white.noise.exists} is reduced to showing \eqref{e.bound.diff.W'}. 

\smallskip

For each $n \in \N$, let $u_n=u_n(t,x)$ be the solution of the Cauchy problem
\begin{equation*} 
\left\{ 
\begin{aligned}
& \partial_t u_n - \Delta u_n = 0 & \mbox{in}  & \ \R^d \times (0,\infty), \\ 
& u_n(0,\cdot)  = \W'_{n+1} - \W'_n & \mbox{in} & \ \Rd. 
\end{aligned}
\right.
\end{equation*}
By Proposition~\ref{p.MSP.Walpha}, in order to prove \eqref{e.bound.diff.W'}, it suffices to show that
\begin{equation}
\label{e.crit.besov}
\int_0^{1} t^{\frac \al 2 - 1} \|u_n(t,\cdot)\|_{L^{\infty}(\Rd)} \, dt \le 2^{-\be n } \Ll( C n^{\frac12} + \O_2(C) \Rr).
\end{equation}
We decompose the proof of this statement into five substeps.

\smallskip

\emph{Step 2.a.} Denoting 
\begin{equation*}  
\td X_n := \sup \Ll\{ \Ll| X_{n,z}^{(i)} \Rr| \ : \ z \in Z_n, \, 0 \le i < 2^d \Rr\} ,
\end{equation*}
we show that there exists $C(d) < \infty$ such that 
\begin{equation}
\label{e.bound.td.Xn}
\td X_n \le C n^{\frac12} +  \O_2(C).
\end{equation}
Indeed, since $X_{n,z}^{(i)}$ are standard Gaussians, we have, for every $C_0 \in (0,\infty)$ and $z \ge 1$,
\begin{equation*}  
\P \Ll[ X_{n,z}^{(i)} \ge C_0 n^{\frac12} + z \Rr] \le \exp \Ll( -\frac{(C_0 n^{\frac12} + z)^2}{2} \Rr) \le \exp \Ll( -\frac {C_0^2 n + z^2}{2}  \Rr),
\end{equation*}
and therefore, by a union bound,
\begin{equation*}  
\P \Ll[ \td X_n \ge C_0 n^{\frac12} + z \Rr] \le 2^{(n+1)d} \exp \Ll( -\frac {C_0^2 n + z^2}{2}  \Rr) \le  \exp \Ll( -\frac {z^2}{2} \Rr) ,
\end{equation*}
provided that the constant $C_0$ is sufficiently large. By Lemma~\ref{l.bigO.vs.tail}, this shows~\eqref{e.bound.td.Xn}.

\smallskip

\emph{Step 2.b.} 
We show that there exists $C(\al,d) < \infty$ such that 
\begin{equation}
\label{e.crit.small.t}
\int_0^{2^{-2n}} t^{\frac \al 2 - 1} \|u_n(t,\cdot)\|_{L^{\infty}(\Rd)} \, dt \le 2^{-n \Ll( \al - \frac d 2 \Rr) } \Ll( C n^{\frac12} + \O_2(C) \Rr).
\end{equation}
By \eqref{e.bound.td.Xn}, \eqref{e.def.chin} and the maximum principle, we have
\begin{equation}  
\label{e.un.Linfty}
\|u_n(t,\cdot)\|_{L^\infty(\Rd)} \le 2^{\frac{nd}{2}} \Ll( C n^{\frac12} + \O_2(C) \Rr) .
\end{equation}
Using also Lemma~\ref{l.sum-O}, we obtain \eqref{e.crit.small.t}.

\smallskip

\emph{Step 2.c.} We show that there exists $C(d) < \infty$ such that, for every $x \in \Rd$ and $t \ge 2^{-2n}$,
\begin{equation}
\label{e.bound.un.point}
u_n(t,x) = \O_2 \Ll( C 2^{-n}  t^{-\frac {2 + d}{4}} \exp \Ll( -\frac{(|x|_{\infty} - 1)_+^2}{Ct} \Rr)\Rr),
\end{equation}
where we denote, for $x = (x_1, \ldots, x_d)\in\Rd$, the norm $|x|_\infty := \max_{1 \le i \le d} |x_i|$.
By the definition of $u_n$, we have, for every $x \in \Rd$ and $t > 0$,
\begin{equation}  
\label{e.decomp.un}
u_n(t,x) = \sum_{\substack{ z \in Z_n \\ 1 \le i < 2^d}} X_{n,z}^{(i)} \int_{\Rd} \chi_n^{(i)}(y) \Phi(t,x+z-y) \, dy.
\end{equation}
Moreover, recall that, for $i \ge 1$, the function $\chi_n^{(i)}$ is of mean zero. Therefore, for every $t \ge 2^{-2n}$, we have
\begin{align*}  
\Ll|\int_{\Rd} \chi_n^{(i)}(y) \Phi(t,x+z-y) \, dy \Rr| & = \Ll| \int_{\Rd} \chi_n^{(i)}(y) \Ll(\Phi(t,x+z-y) - \Phi(t,x+z)\Rr) \, dy \Rr| \\
& \le C 2^{- \Ll( 1 + \frac d 2 \Rr) n } t^{-\Ll(\frac 1 2 + \frac d 2\Rr)} \exp \Ll( -\frac{|x+z|^2}{Ct} \Rr) .
\end{align*}
By independence of the $(X_{n,z}^{(i)})$, we deduce from this and \eqref{e.decomp.un} that $u_n(t,x)$ is a centered Gaussian with variance bounded by
\begin{equation*}  
C\sum_{\substack{ z \in Z_n \\ 1 \le i < 2^d}} 2^{-(2+d) n} t^{-(1+d)} \exp \Ll( -\frac{|x + z|^2}{Ct} \Rr) \le C 2^{-2n}  t^{-\Ll(1+\frac d 2\Rr)} \exp \Ll( -\frac{(|x|_{\infty} - 1)_+^2}{Ct} \Rr) .
\end{equation*}
Since a Gaussian of unit variance is $\O_2(C)$, this yields \eqref{e.bound.un.point} by homogeneity.

\smallskip

\emph{Step 2.d.} 
We show that there exists~$C(d) < \infty$ such that, for every $t \in[ 2^{-2n},1]$,
\begin{equation}
\label{e.bound.un.Linfty}
\|u_n(t,\cdot)\|_{L^\infty(\Rd)} \le 2^{-n}  t^{-\frac{2+d}{4}} \Ll( C n^{\frac12} +  \O_2 \Ll(C\Rr) \Rr).
\end{equation}
Using elementary heat kernel bounds, we see that there exists~$C(d) < \infty$ such that, for every~$x, y \in \Rd$ and $t > 0$,
\begin{equation*}  
\Ll| u_n(x,t) - u_n(y,t) \Rr| 
\le C t^{-\frac12} |x-y| \Ll\|u_n\Ll(\cdot,0\Rr)\Rr\|_{L^\infty(\Rd)}.
\end{equation*}
Using also \eqref{e.un.Linfty}, we deduce that, for every $\de > 0$ and $t > 0$, 
\begin{equation}  
\label{e.bound.twot}
\|u_n(t,\cdot)\|_{L^\infty(\Rd)} \le \sup_{x \in \de \Z^d} \Ll| u_n(x,t) \Rr|  + C \de t^{-\frac 1 2}2^{\frac{nd}{2}} \Ll( C n^{\frac12} + \O_2(C) \Rr) .
\end{equation}
Taking 
\begin{equation}  
\label{e.choice.of.delta}
\de := 2^{-\Ll( 2 + \frac d 2 \Rr) n},
\end{equation}
we deduce that, for every $t \in (0,1]$,
\begin{equation}  
\label{e.so.you.know.how.to.choose.delta}
\de t^{-\frac 1 2}2^{\frac{nd}{2}} \Ll( C n^{\frac12} + \O_2(C) \Rr) \le \O_2 \Ll(C  2^{- n} t^{-\frac {2+d}{4}} \Rr) .
\end{equation}
We now estimate the first term on the right side of \eqref{e.bound.twot} by a union bound, \eqref{e.bound.un.point} and Lemma~\ref{l.bigO.vs.tail}: for every $z \in (0,\infty)$, we have
\begin{align*}  
\P \Ll[ \sup_{x \in \de \Z^d} \Ll| u_n(x,t) \Rr| \ge z \Rr] & \le \sum_{x \in \de \Z^d} \P \Ll[  \Ll| u_n(x,t) \Rr| \ge z \Rr]  \\
& \le 2 \sum_{x \in \de \Z^d} \exp \Ll( - C^{-1} 2^{2n} t^{1 + \frac d 2} \exp \Ll( \frac{(|x|_{\infty} - 1)_+^2}{Ct} \Rr)  z^2  \Rr) .
\end{align*}
Replacing $z$ by $2^{-n} t^{1 + \frac d 2} \Ll( C_1 n^{\frac12} + z \Rr)$ in the expression above, with $C_1(d)$ sufficiently large, and using \eqref{e.choice.of.delta}, we obtain the existence of~$C(d) < \infty$ such that, for every $t \in (0,1]$,
\begin{equation*}  
\sup_{x \in \de \Z^d} \Ll| u_n(x,t) \Rr| \le 2^{-n} t^{-\frac {2 + d}{4}} \Ll( C n^{\frac12} + \O_2(C) \Rr) .
\end{equation*}
Combining this with \eqref{e.bound.twot} and \eqref{e.so.you.know.how.to.choose.delta} yields \eqref{e.bound.un.Linfty}.

\smallskip

\emph{Step 2.e.} We conclude the proof of \eqref{e.crit.besov}, and therefore of the theorem. By the maximum principle, \eqref{e.bound.un.Linfty} and the definition of $\be$ in \eqref{e.def.be.w}, we have
\begin{align*}  
\int_{2^{-2n}}^{1} t^{\frac \al 2 - 1} \|u_n(t,\cdot)\|_{L^{\infty}(\Rd)} \, dt & \le C \sum_{k = 1}^{2n} 2^{-\frac{\al k}{2}} \Ll\| u_n(\cdot,2^{-k}) \Rr\|_{L^\infty(\Rd)} \\
& \le C \sum_{k = 1}^{2n} 2^{\Ll(1 + \frac d 2 - \al\Rr)\frac k 2 } 2^{-n} \Ll( C n^{\frac12} + \O_2(C) \Rr) \\
& \le 2^{-\be n} \Ll( C n^{\frac12} + \O_2(C) \Rr) .
\end{align*}
This and \eqref{e.crit.small.t} yield \eqref{e.crit.besov} and thus complete the proof.
\end{proof}
\begin{exercise}  
\label{ex.vector.white.noise}
Generalizing Proposition~\ref{p.white.noise.exists}, construct vector white noise as a random element of $W^{-\al,\infty}_{\mathrm{loc}}(\Rd;\R^n)$. That is, show that, for every $n \in \N$, $\al > \frac d 2$ and every symmetric nonnegative matrix $\msf Q \in \R^{n \times n}$, there exists a random element $\W$ of $W^{-\al,\infty}_{\mathrm{loc}}(\Rd;\R^n)$ such that, for every $\f \in W^{\al,1}(\Rd;\R^n)$ with compact support, the evaluation of $\W$  against $f$, denoted by $\W(\f)$, is a centered Gaussian with variance
\begin{equation*}  
\int_{\Rd} \f \cdot \msf Q \f,
\end{equation*}
and moreover, that there exists a constant $C(d,n) < \infty$ such that 
\begin{equation}
\label{e.vector.white.noise.estimate}
\sup \Ll\{\W(\f) \ : \ \|\f\|_{W^{\al,1}(\Rd;\R^n)} \le 1 \ \text{and} \ \supp \f \subset B_1 \Rr\} \le \O_2 \Ll( C |\msf Q|^\frac 1 2 \Rr) .
\end{equation}
\end{exercise}

We now turn to the construction of the gradient Gaussian free field. For every $\alpha \in \R$ and $p\in [1,\infty]$, we denote 
\begin{multline}
\label{e.def.Wsol}
W^{\al,p}_{\mathrm{pot, \, loc}}(\Rd) := \Big\{ G \in W^{\al,p}_{\mathrm{loc}}(\Rd) 
\\
\,:\, \forall \f \in C^\infty_c(\Rd;\Rd), \ 
\nabla \cdot \f = 0 \implies G(\f) = 0 \Big\} .
\end{multline}
Using Lemma~\ref{l.mean-zero}, one may represent any element of~$W^{\al,p}_{\mathrm{pot,\, loc}}(\Rd)$ as the gradient of a distribution, although a canonical choice of the distribution can only be made up to an additive constant, as was explained below this lemma. 

\smallskip

We fix a symmetric positive definite matrix $\ahom \in \R^{d \times d}$, which we assume for convenience to satisfy the uniform ellipticity condition
\begin{equation}
\label{e.gff.ellipticity}
\forall \xi \in \Rd, \quad |\xi|^2 \le \xi \cdot \ahom \xi.
\end{equation}
Under this assumption, we have
\begin{equation}  
\label{e.Proj.contraction}
\|\Proj(\f)\|_{L^2(\Rd)} \le \|\f\|_{L^2(\Rd)}.
\end{equation}
We also remark that, by standard estimates for the Poisson equation,
\begin{equation*} \label{}
\f \in C^\infty_c(\Rd;\Rd) \implies \Proj(\f) \in C^\infty(\Rd;\Rd).
\end{equation*}

\index{Gaussian free field!construction}

\begin{proposition}[Construction of gradient GFF]
\label{p.gff.exists}
Let $\ahom \in \R^{d \times d}$ be a symmetric matrix satisfying \eqref{e.gff.ellipticity}, $\msf Q \in \R^{d \times d}$ be a nonnegative symmetric matrix, and~$\al > \frac d 2$. There exists a random element $\nabla \Psi$ of $W^{-\al,\infty}_\mathrm{pot,\,loc}(\Rd)$ such that the evaluation of $\nabla \Psi$ against $\f \in W^{\al,1}(\Rd;\Rd)$ with compact support, denoted by $\nabla \Psi(\f)$, is a centered Gaussian with variance 
\begin{equation*}  
\int_{\Rd} \Ll(\Proj \f\Rr) \cdot \msf Q \Ll( \Proj \f \Rr) .
\end{equation*}
Moreover, there exists a constant $C(d) < \infty$ such that 
\begin{equation}
\label{e.gff.estimate}
\sup \Ll\{\nabla \Psi(\f) \ : \ \|\f\|_{W^{\al,1}(\Rd;\Rd)} \le 1 \ \text{and} \ \supp \f \subset B_1 \Rr\} \le \O_2 \Ll( C |\msf Q|^\frac 1 2 \Rr) .
\end{equation}
\end{proposition}
\begin{proof}
Let $\W$ be a vector white noise of covariance matrix $\msf Q$, as given by Proposition~\ref{p.white.noise.exists} and Exercise~\ref{ex.vector.white.noise}. Without loss of generality, we assume that $|\msf Q| \le 1$. In view of \eqref{e.eq.gff}, the goal is to justify that the definition of
\begin{equation}  
\label{e.def.grad.gff}
\nabla \Psi(\f) := \W \Ll( \Proj(\f) \Rr) 
\end{equation}
makes sense and satisfies the requirements of Proposition~\ref{p.gff.exists}. The main difficulty is that $\Proj(\f)$ has long tails---in particular, it is not compactly supported. 

\smallskip

We select a cutoff function~$\chi \in C^\infty_c(\Rd)$ such that $\chi \equiv 1$ on $B_{1}$ and $\chi \equiv 0$ outside of $B_2$. For every $r > 0$, we set $\chi_r := \chi(r^{-1} \, \cdot)$, and for every $\f \in C^\infty_c(\Rd;\Rd)$, define
\begin{equation*}  
\nabla \Psi_r(\f):=  \W \Ll( \chi_r \, \Proj(\f) \Rr).
\end{equation*}
We aim to study the convergence, as $r\to\infty$, of the distribution~$\nabla \Psi_r$ in the space~$W^{-\al,\infty}_{\mathrm{pot,\, loc}}(\Rd)$. This distribution is a potential field, that is, we have, almost surely,
\begin{equation*}  
\forall \f \in C^\infty_c(\Rd;\Rd), \quad \nabla \cdot \f = 0 \implies \nabla \Psi_r(\f) = 0.
\end{equation*}
Indeed, this is clear from the definition of the Helmholtz-Hodge projection $\Proj$ and it justifies the use of the notation~$\nabla \Psi_r$. By the translation invariance of the law of~$\nabla \Psi_r$, we may localize the field~$\nabla \Psi_r$ in a neighborhood of the origin by setting
\begin{equation*}  
\msf G_r (\f) := \nabla \Psi_r(\chi \f),
\end{equation*}
and reduce the proof of Proposition~\ref{p.gff.exists} to the demonstration that 
\begin{equation}
\label{e.msfG.r1}
\|\msf G_4\|_{W^{-\al,\infty}(\Rd)} \le \O_2(C),
\end{equation}
and that, for every $r \ge 4$,
\begin{equation}  
\label{e.msfG.2r.r}
\Ll\| \msf G_{2r} - \msf G_r \Rr\|_{W^{-\al,\infty}(\Rd)} \le \O_2 \Ll( C r^{- \frac d 2 } \log^{\frac12} r  \Rr) .
\end{equation}
We decompose the proof of these two statements into three steps.

\smallskip

\emph{Step 1.} We show \eqref{e.msfG.r1}. By Proposition~\ref{p.MSP.Walpha}, it suffices to show that
\begin{equation}  
\label{e.msfG.r1.msp}
\int_0^1 t^{\frac \al 2 - 1} \sup_{x \in \Rd} \Ll| \W \Ll( \chi_4 \Proj \Ll( \chi \Phi(t, \cdot - x) \Rr)  \Rr)  \Rr| \, dt \le \O_2(C).
\end{equation}
By \eqref{e.Proj.contraction}, we have, for every $x \in \Rd$ and $t \in (0,1]$,
\begin{align}  
\notag
\Ll|\W\Ll( \chi_4 \Proj \Ll( \chi \Phi(t, \cdot - x) \Rr)  \Rr)\Rr| & \le \O_2 \Ll( C \Ll\|\chi \Phi(t, \cdot - x)\Rr\|_{L^2(\Rd)} \Rr) \\
\label{e.bound.w.11}
& \le\O_2 \Ll(    C t^{-\frac d 4} \exp \Ll( -\frac{\Ll( |x| - 2\Rr)_+^2  }{C t} \Rr)\Rr) .
\end{align}
We also observe that, for every  $x, y \in \Rd$ satisfying $|y-x| \le 1$ and $t \in (0,1]$,
\begin{align}  
\notag
\lefteqn{
\Ll|\W \Ll( \chi_4 \Proj \Ll( \chi \Ll[\Phi(t, \cdot - y) - \Phi(t, \cdot - x)\Rr] \Rr)  \Rr) \Rr|
}  \qquad & \\
\label{e.W.chi4}
&  \le \X \, \Ll\|\chi \Ll[\Phi(t, \cdot - y) - \Phi(t, \cdot - x)\Rr]\Rr\|_{W^{d,1}(\Rd)} \\
\notag
& \le C \X \,  t^{-\frac {d+1} 2} |x-y| 
,
\end{align}
where the random variable $\X$ satisfies, by Proposition~\ref{p.white.noise.exists},
\begin{equation*}  
\X \le \O_2 \Ll( C \Rr) .
\end{equation*}
We thus deduce that
\begin{equation*}
\sup_{x \in \Rd} \Ll| \W \Ll( \chi_4 \Proj \Ll( \chi \Phi(t, \cdot - x) \Rr)  \Rr)  \Rr| \\
\le \sup_{x \in t^{-d}\Z^d}  \Ll| \W \Ll( \chi_4 \Proj \Ll( \chi \Phi(t, \cdot - x) \Rr)  \Rr)  \Rr| + \O_2 \Ll( C \Rr) .
\end{equation*}
Arguing as in Step 2.a of the proof of Proposition~\ref{p.white.noise.exists}, one can check using \eqref{e.bound.w.11} that the supremum on the right side above is bounded by 
\begin{equation*}  
\O_2 \Ll( C t^{-\frac d 4} \log^{\frac12} (1+t^{-1}) \Rr) ,
\end{equation*}
and therefore the proof of \eqref{e.msfG.r1.msp} is complete.

\smallskip

\emph{Step 2.} Preparing the ground for the proof of \eqref{e.msfG.2r.r}, we show that there exists a constant $C(d)$ such that, for every  $r \ge 4$,
\begin{equation}
\label{e.white.noise.est.annulus}
\sup \Ll\{\W \Ll( (\chi_{2r} - \chi_r) \f \Rr) : \|(\chi_{2r} - \chi_r)\f\|_{W^{d,1}(\Rd)} \le 1\Rr\} \le \O_2 \Ll( C \log^{\frac12} r  \Rr) .
\end{equation}
As in \eqref{e.W.chi4}, the exponent $d$ in the norm $W^{d,1}(\Rd)$ appearing in the above display is so chosen for simplicity, but could be reduced to any exponent larger than $\frac d 2$. In order to show \eqref{e.white.noise.est.annulus}, we use a partition of unity $\theta \in C^\infty_c(\Rd;\Rd)$ such that 
\begin{equation*}  
\sum_{z \in \Z^d} \theta(\cdot - z) = 1,
\end{equation*}
and decompose $\f$ accordingly:
\begin{align*}  
\W \Ll( (\chi_{2r} - \chi_r) \f \Rr) & = \sum_{z \in \Z^d}  \W \Ll( (\chi_{2r} - \chi_r) \f \theta(\cdot - z)\Rr)  \\
& \le \sum_{z \in Z} \X(z) \, \|(\chi_{2r} - \chi_r) \f \theta(\cdot - z)\|_{W^{d,1}(\Rd)},
\end{align*}
where 
\begin{equation*}  
Z := \Ll\{ z \in \Zd \ : \ \supp \theta(\cdot - z) \cap \supp (\chi_{2r} - \chi_r) \neq \emptyset \Rr\} 
\end{equation*}
and 
\begin{equation*}  
\X(z) := \sup \Ll\{ \W(\f) \ : \ \|\f\|_{W^{d,1}(\Rd)} \le 1 \ \text{ and } \ \supp \f \subset \supp \theta(\cdot - z)  \Rr\} .
\end{equation*}
By Proposition~\ref{p.white.noise.exists}, there exists a constant $C$ depending only on our choice of the function $\theta$ and the dimension $d$ such that, for every  $z \in \Z^d$,
\begin{equation*}  
\X(z) \le \O_2(C).
\end{equation*}
Since moreover, we have $|Z| \le C r^{d}$, we can argue as in Step 2.a of the proof of Proposition~\ref{p.white.noise.exists} to obtain that
\begin{equation*}  
\sup_{z \in Z} \X(z) \le \O_2\Ll(C \log^{\frac12} r \Rr),
\end{equation*}
and therefore that
\begin{align*}  
\W \Ll( (\chi_{2r} - \chi_r) \f \Rr) & \le \O_2\Ll(C \log^{\frac12} r \Rr)\sum_{z \in Z} \|(\chi_{2r} - \chi_r) \f \theta(\cdot - z)\|_{W^{d,1}(\Rd)} \\
& \le \O_2\Ll(C \log^{\frac12} r \Rr) \|(\chi_{2r} - \chi_r)\f\|_{W^{d,1}(\Rd)}.
\end{align*}
This completes the proof of \eqref{e.white.noise.est.annulus}.

\smallskip

\emph{Step 3.} We now show \eqref{e.msfG.2r.r}, which by Proposition~\ref{p.MSP.Walpha} is implied by
\begin{equation}
\label{e.msfG.msp}
\int_0^1 t^{\frac \al 2 - 1} \sup_{x \in \Rd} \Ll| \W \Ll( (\chi_{2r} - \chi_r) \Proj \Ll( \chi \Phi(t, \cdot - x) \Rr)  \Rr)   \Rr| \, dt \le \O_2 \Ll( C r^{-\frac d 2}\log^{\frac12} r  \Rr) .
\end{equation}
For each $x \in \Rd$ and $t \in (0,1]$, we denote 
\begin{equation*}  
\nabla u_{x,t} := \Proj \Ll( \chi \Phi(t,\cdot- x)\Rr) .
\end{equation*}
By Green's representation formula, for each $k \in \N$, there exists a constant $C$ depending only on $k, d$ and our choice for the function $\chi$  such that, for every  $z \in \Rd \setminus B_3$ and $t \in (0,1]$,
\begin{align*}  
\Ll|\nabla^k (u_{x,t} - u_{y,t})(z)\Rr| 
& 
\leq C |z|^{-(d+k)} \|\chi \Ll[\Phi(t, \cdot - x) - \Phi(t, \cdot - y)\Rr]\|_{L^1(\Rd)} 
\\ & 
\leq C |z|^{-(d+k)} t^{-\frac {1} 2} \,|x-y|.
\end{align*}
Combining this with \eqref{e.white.noise.est.annulus}, we deduce that, for every $t \in (0,1]$ and $r \ge 4$,
\begin{equation}  
\label{e.decomp.sup.aa}
\sup_{x \in \Rd}  \Ll| \W \Ll( (\chi_{2r} - \chi_r) \nabla u_{x,t}\Rr)  \Rr| \le \sup_{x \in t \Z^d} \Ll| \W \Ll( (\chi_{2r} - \chi_r) \nabla u_{x,t} \Rr)  \Rr|  + \O_2 \Ll( C r^{-\frac d 2} \log^{\frac12} r \Rr) .
\end{equation}
Moreover, using Green's representation formula again, we have, for every $x \in \Rd$, $t \in (0,1]$ and $z \in \Rd \setminus B_3$,
\begin{equation*}  
\Ll|\nabla u_{x,t}(z)\Rr| \le C \, |z|^{-d},
\end{equation*}
and therefore,
\begin{equation*}  
\Ll|\W \Ll( (\chi_{2r} - \chi_r) \nabla u_{x,t}\Rr) \Rr| \le \O_2 \Ll( C \|(\chi_{2r} - \chi_r) \nabla u_{x,t}\|_{L^2(\Rd)} \Rr) \le \O_2 \Ll( C r^{-\frac d 2} \Rr) .
\end{equation*}
From this, we deduce as before that the supremum on the right side of \eqref{e.decomp.sup.aa} is bounded by
\begin{equation*}  
\O_2 \Ll( C r^{-\frac d 2} \log^{\frac12} \Ll(1 + t^{-1} \Rr)  \Rr) .
\end{equation*}
This implies \eqref{e.msfG.msp}, and thus \eqref{e.msfG.2r.r}, thereby completing the proof.
\end{proof}

\begin{exercise}  
\label{ex.unbounded.critical.gff}
Adapt Lemma~\ref{l.unbounded.W} and then Exercise~\ref{ex.unbounded.critical} to the case in which white noise is replaced by a gradient GFF. That is, show that if $\nabla \Psi$ is a random distribution such that, for every  $\f \in C^\infty_c(\Rd)$, the evaluation of $\nabla \Psi$ against~$\f$ is a centered Gaussian with variance 
\begin{equation*}  
\int_{\Rd} \Ll(\Proj \f\Rr) \cdot \msf Q \Ll( \Proj \f \Rr) ,
\end{equation*}
then the event that $\nabla \Psi$ belongs to $W^{-\frac d 2,p}_{\mathrm{loc}}(\Rd)$ has null probability, where $p = 1$ if $d$ is odd, and $p > 1$ if $d$ is even.
\end{exercise}

\section{Heuristic derivation of the scaling limit}
\label{s.heuristics}

We henceforth return to the context of the homogenization problem. In this section, we explain heuristically why we should expect the energy quantity $J_1$ to converge to a convolution of white noise, and why this should imply the convergence of the rescaled corrector to a Gaussian free field.

\smallskip

We showed in Theorem~\ref{t.additivity} that, for every   $s \in (0,2)$, there exists a constant $C(s,d,\Lambda) < \infty$ such that, for every  $r \ge 1$, $z \in \Rd$ and $p,q \in B_1$,
\begin{equation}  
\label{e.recall.our.trophy}
J_1(\Phi_{z,r}, p,q) = \frac 1 2 p \cdot \ahom p + \frac 1 2 q \cdot \ahom^{-1} q - p \cdot q  + \O_s \Ll( C r^{-\frac d 2} \Rr) .
\end{equation}
By \eqref{e.J1derp}--\eqref{e.J1derq}, this implies that, for every  $r \ge 1$, $z \in \Rd$ and $p \in B_1$,
\begin{equation}  
\label{e.recall.spat.av}
\Ll\{
\begin{aligned}
\int_{\Phi_{z,r}} \nabla v(\cdot, \Phi_{z,r}, -p,0) &  = p + \O_s \Ll( C r^{-\frac d 2} \Rr) 
, \quad \text{and} \\
\int_{\Phi_{z,r}} \a\nabla v(\cdot, \Phi_{z,r}, -p,0) &  = \ahom p + \O_s \Ll( C r^{-\frac d 2} \Rr) 
.
\end{aligned}
\Rr.
\end{equation}
In other words, the homogenized matrix identifies the correspondence between the spatial averages on scale $r$ of gradients and fluxes of elements of $\A_1$, up to an error of~$\O_s \Ll(C r^{-\frac d 2}\Rr)$. The natural problem which then arises is to identify the next-order term in the description of this correspondence. In this section, we first rephrase (rigorously) this question in terms of the quantity $J_1$; we then describe in heuristic terms the next-order behavior of $J_1$ and the correctors.

\smallskip

The next-order correction in the correspondence between spatial averages of gradients and fluxes of elements of $\A_1$ is captured by the quantity
\begin{equation}  
\int_{\Phi_{z,r}} p \cdot (\ahom - \a) \nabla v(\cdot,\Phi_{z,r},p',0),
\end{equation}
for $p, p'$ ranging in $B_1$. 
%
Naturally, this quantity can be rewritten in terms of~$J_1$. Indeed, using \eqref{e.J1derp}--\eqref{e.J1derq} and then the fact that $(p,q) \mapsto J_1(\Phi_{z,r},p,q)$ is a quadratic form, we obtain
\begin{align}  
\label{e.identity.br}
\int_{\Phi_{z,r}} p \cdot (\ahom - \a) \nabla v(\cdot,\Phi_{z,r},p',0) & = 
p \cdot  D_p   J_1(\Phi_{z,r},p',0) + \ahom p \cdot  D_q   J_1(\Phi_{z,r},p',0) 
\\
\notag
& = \frac 1 2 \Ll( J_1(\Phi_{z,r}, p + p',\ahom p) - J_1(\Phi_{z,r}, p - p',\ahom p) \Rr) .
\end{align}
This motivates the following definition and lemma, which can be compared with Lemma~\ref{l.commutator}.
\index{coarsened coefficient field}
\begin{definition}[Coarsened coefficients]
\label{d.brz}
For every~$r \ge 1$ and $z \in \Rd$, we let $\b_r(z)$ be the random $d$-by-$d$ matrix such that, for every  $p,p' \in \Rd$,
\begin{equation}  
\label{e.def.bhomrz}
p \cdot \b_r(z) p' = \frac 1 2 \Ll( J_1(\Phi_{z,r}, p + p',\ahom p) - J_1(\Phi_{z,r}, p - p',\ahom p) \Rr).
\end{equation}
\end{definition}

\begin{lemma}  
\label{l.brz}
Fix $s \in (0,2)$. There exists $C(s,d,\Lambda) < \infty$ such that, for every $r \ge 1$ and $z \in \Rd$, we have
\begin{equation}
\label{e.estim.brz}
|\b_r(z)| \le \O_s\Ll( C r^{-\frac d 2} \Rr) ,
\end{equation}
and
\begin{equation}  
\label{e.heres.your.next.order}
\sup_{u \in \A_1(\Phi_{z,r})}\Ll|\int_{\Phi_{z,r}} \a(x) \nabla u(x) \, dx -\Ll( \ahom + \b_r(z) \Rr) \int_{\Phi_{z,r}} \nabla u(x) \, dx \Rr| \le \O_{s/2} \Ll( C r^{-d} \Rr) .
\end{equation}
\end{lemma}
\begin{proof}
We obtain \eqref{e.estim.brz} directly from the definition of $\b_r(z)$ and 
\eqref{e.recall.our.trophy}. The estimate~\eqref{e.heres.your.next.order} then follows from~\eqref{e.def.bhomrz} and Lemma~\ref{l.redwedding}. 
\end{proof}

\smallskip

 By the additivity and localization properties of $J_1$ given in the statement of Theorem~\ref{t.additivity}, it is natural to expect~$J_1$ to display CLT-like fluctuations. In fact, we have already used this intuition in the paragraph preceding \eqref{e.meso.mask}: that $J_1(\Phi_{z,r},p,q)$ should resemble the average against $\Phi_{z,r}$ of a random field with finite range of dependence. Recalling also the estimate of the expectation of $J_1$ in Theorem~\ref{t.additivity}(ii), it seems reasonable to expect that 
\begin{equation}
\label{e.heuristic.J1}
J_1(\Phi_{z,r}, p,q) \simeq \frac 1 2 p \cdot \ahom p + \frac 1 2 q \cdot \ahom^{-1} q - p \cdot q + \int_{\Phi_{z,r}} \W(\cdot,p,q) ,
\end{equation}
where $\W(\cdot,p,q)$ is a white noise with quadratic dependence in~$(p,q)$, and where the integral on the right side of \eqref{e.heuristic.J1} is an informal notation for the evaluation of $\W(\cdot,p,q)$ against the test function $\Phi_{z,r}$. By the basic properties of white noise, this integral is a centered Gaussian, and its fluctuations are of the order of $r^{-\frac d 2}$. The informal statement \eqref{e.heuristic.J1} should be understood as an approximate identity between the \emph{laws} of the quantities on each side, and up to an error of lower order compared with $r^{-\frac d 2}$. Note that this approximate identity is meant to be consistent as we vary the parameters $z \in \Rd$, $r \ge 1$ and $p,q \in B_1$.

\smallskip

In addition, Lemma~\ref{l.brz} strongly suggests that locally averaged elements of $\A_1$ solve a \emph{coarsened equation} with diffusion coefficients given by $\ahom + \b_r$. That is, we expect that the equation
\begin{equation}  
\label{e.coarsened.equation}
-\nabla \cdot (\ahom + \b_r) \Ll[(e+\nabla \phi_e) \ast \Phi_r\Rr] = 0
\end{equation}
holds approximately. 
Rearranging, we get
\begin{equation}
\label{e.coarsened.2}
-\nabla \cdot (\ahom + \b_r) \nabla (\phi_e \ast \Phi_r) = \nabla \cdot ( \b_r e).
\end{equation}
By Lemma~\ref{l.brz}, the coefficients $\b_r$ are of the order of $r^{-\frac d 2}$. Hence, we expect $\nabla (\phi_e \ast \Phi_r)$ to be also of this order of magnitude (in fact this is essentially how we proved \eqref{e.gradient}), and thus the term $\b_r \nabla (\phi_e \ast \Phi_r)$ should be of lower order. Discarding this term, we arrive at the following equation, which should approximately hold for large $r$:
\begin{equation}
\label{e.coarsened.3}
-\nabla \cdot \ahom \nabla  (\phi_e \ast \Phi_r) = \nabla \cdot (\b_r e).
\end{equation}
By the definition of $\b_r$ and \eqref{e.heuristic.J1}, we also expect
\begin{equation}  
\label{e.almost.def.msf.V}
p \cdot \b_r(z) e \simeq \frac 1 2\int_{\Phi_{z,r}} \Ll( \W(\cdot,p+e,\ahom p) - \W(\cdot,p-e,\ahom p) \Rr) ,
\end{equation}
where the approximate identity is understood as in \eqref{e.heuristic.J1}; or equivalently,
\begin{equation*}  
p \cdot \b_r e \simeq \frac 1 2 \Ll( \W(\cdot,p+e,\ahom p) - \W(\cdot,p-e,\ahom p) \Rr) \ast \Phi_{r}.
\end{equation*}
In other words, the law of the vector field $\b_r e$ is close to that of a vector white noise convolved with the heat kernel on scale $r$. With this interpretation in mind, we see that equation~\eqref{e.coarsened.3} has the same form as the defining equation for the gradient Gaussian free field, see \eqref{e.eq.gff}, up to a convolution with the heat kernel on scale $r$. This strongly suggests that the spatial averages of $\nabla \phi_e$ should resemble those of a gradient GFF, and this is indeed what will be proved in Theorem~\ref{t.gff} using an argument inspired by the heuristics in this section.

\section{Central limit theorem for \texorpdfstring{$J_1$}{J1}}
\label{s.CLT}

The goal of this section is to give a next-order description of the energy quantity~$J_1$ introduced in \eqref{e.def.J1}. This is achieved in Theorem~\ref{t.conv.white} below. 
%
%
As explained in the paragraph around \eqref{e.meso.mask}, in order to describe the precise behavior of $J_1$, the main point is to identify the behavior of quantities of the form
\begin{equation}  
\label{e.meso.back}
\fint_{\cu_R} J_1(\Phi_{z,r}, p,q) \, dz,
\end{equation}
where $r \le R$ is a relatively large mesoscale. We have seen in Theorem~\ref{t.additivity}(ii) that 
\begin{equation*}  
\Ll|\fint_{\cu_R} \E\Ll[J_1(\Phi_{z,r}, p,q)\Rr] \, dz -\Ll( \frac 1 2 p \cdot \ahom p + \frac 1 2 q \cdot \ahom^{-1} q - p \cdot q \Rr)\Rr| \le C r^{-d}. 
\end{equation*}
Since we aim for a description of the quantity in \eqref{e.meso.back} only up to a precision of~$o(R^{-\frac d 2})$, the previous estimate gives us a sufficient understanding of its expectation. We focus therefore on analyzing its fluctuations, which is the purpose of the next proposition. Before giving the statement, we introduce the notion of \emph{Gaussian quadratic form}, which \index{Gaussian quadratic form} is a random quadratic form such that the entries of its associated symmetric matrix form a centered Gaussian vector (see also~\eqref{e.def.msfQ.Q} below).
\begin{proposition}
\label{p.conv.white}
There exist~$\eta_1(d) \in \left(\frac 12,1\right)$ and a Gaussian quadratic form~$\msf Q$ on~$\R^d \times \Rd$ such that, for every $\eta \in [\eta_1,1)$ and $z, p, q \in \Rd$,
\begin{equation*} \label{}
\msf Q (p, \ahom p) = 0
\end{equation*}
and
\begin{equation}
\label{e.clt}
r^{\frac d 2} \fint_{\cu_r(z)} \Ll(J_1(\Phi_{x,r^{\eta}},p,q) - \E \Ll[ J_1(\Phi_{x,r^{\eta}},p,q) \Rr] \Rr) dx  \xrightarrow[r \to \infty]{\mathrm{(law)}}  \msf Q(p,q) .
\end{equation}
Moreover the convergence in law holds jointly over $p,q \in \Rd$.
\end{proposition}

Here and throughout, we say that a convergence in law of the form 
\begin{equation}  
\label{e.conv.X}
X_r(p) \xrightarrow[r \to \infty]{\mathrm{(law)}} X(p)
\end{equation}
\emph{holds jointly} over the index $p$ if, for every finite family of such indices $p_1,\ldots,p_n$, we have
\begin{equation*}  
\Ll(X_r(p_1), \ldots, X_r(p_n)\Rr) \xrightarrow[r \to \infty]{\mathrm{(law)}} \Ll(X(p_1), \ldots, X(p_n)\Rr).
\end{equation*}
If we have a second family of random variables $Y_r(q)$ such that 
\begin{equation}  
\label{e.conv.Y}
Y_r(q) \xrightarrow[r \to \infty]{\mathrm{(law)}} Y(q),
\end{equation}
then we say that the convergences in law in \eqref{e.conv.X} and \eqref{e.conv.Y} \emph{hold jointly over $p,q$} if, for every finite sequence of parameters $p_1,\ldots,p_n, q_1,\ldots,q_n$, we have 
\begin{multline*}
\Ll(X_r(p_1), \ldots, X_r(p_n), Y_r(q_1), \ldots, Y_r(q_n)\Rr) \\
\xrightarrow[r \to \infty]{\mathrm{(law)}} \Ll(X(p_1), \ldots, X(p_n),  Y(q_1), \ldots, Y(q_n)\Rr).
\end{multline*}

\smallskip

Just as a Gaussian quadratic form can be seen a particular type of Gaussian vector, we can define a quadratic form white noise as a particular type of vector white noise. More precisely, let $Q$ denote the random $2d$-by-$2d$ matrix such that, for every $p,q \in \Rd$,
\begin{equation}  
\label{e.def.msfQ.Q}
\msf Q(p,q) = \frac 1 2 
\begin{pmatrix}  
p \\ q
\end{pmatrix}
\cdot Q
\begin{pmatrix}  
p \\ q
\end{pmatrix}
,
\end{equation}
where $\msf Q$ is the Gaussian quadratic form appearing in the statement of Proposition~\ref{p.conv.white}. By definition, the family of random variables $(Q_{ij})_{1 \le i,j \le 2d}$ is a centered Gaussian vector, with a certain covariance matrix. We can then define $(W_{ij})_{1 \le i,j \le 2d}$ to be the vector white noise having the same covariance matrix as $(Q_{ij})_{1 \le i,j \le 2d}$, and then, for every $p,q \in \Rd$, consider the random distribution
\begin{equation*}  
\W(\cdot,p,q) := \frac 1 2 
\begin{pmatrix}  
p \\ q
\end{pmatrix}
\cdot W
\begin{pmatrix}  
p \\ q
\end{pmatrix}
,
\end{equation*}
or more explicitly, 
\begin{equation*}  
\W(\cdot,p,q) := \frac 1 2 \sum_{i,j = 1}^d \Ll(p_i p_j W_{ij}(\cdot) + q_iq_j W_{d+i,d+j}(\cdot) +2p_i q_j W_{i,d+j}(\cdot)\Rr).
\end{equation*}
We call such $\W$ a \emph{quadratic form white noise}. Abusing notation slightly, for a sufficiently regular test function $f$, we write
\begin{equation}  
\label{e.f.w.pq}
\int_{\Rd} f(x) \W(x,p,q) \, dx
\end{equation}
to denote the evaluation of $\W(\cdot,p,q)$ against $f$. By the construction of $\W$, the random variable in \eqref{e.f.w.pq} is a centered Gaussian, and has the same law as $\|f\|_{L^2(\Rd)} \, \msf Q(p,q)$. 

\smallskip

With the definition of~$\W$ now in place, we can state the following consequence of Proposition~\ref{p.conv.white}, which is the main result of this section. 

\begin{theorem}  
\label{t.conv.white}
Let $\msf Q$ be as in the statement of Proposition~\ref{p.conv.white} and~$\W$ be as defined above. 
There exists $\eta'(d) \in  \Ll( \frac  1 2, 1 \Rr)$ such that, for every $\eta \in (\eta',1)$, $p,q,z \in \Rd$ and $f \in L^1(\Rd) \cap L^\infty(\Rd)$,
\begin{multline}  
\label{e.clt.mesoJ}
r^{\frac d 2} \int_{\Rd} f \Ll( x \Rr) \Ll( J_1(\Phi_{rx,r^\eta},p,q) - \E \Ll[ J_1(\Phi_{rx,r^\eta},p,q) \Rr] \Rr)  \, dx \\
\xrightarrow[r \to \infty]{\mathrm{(law)}}  \int_{\Rd} f(x) \W(x,p,q) \, dx,
\end{multline}
and 
\begin{equation}  
\label{e.clt.J}
r^{\frac d 2} \Ll( J_1(\Phi_{rz,r},p,q) - \E \Ll[ J_1(\Phi_{rz,r},p,q) \Rr] \Rr) \xrightarrow[r \to \infty]{\mathrm{(law)}}  \int_{\Phi_{z,1}}  \W(\cdot,p,q) .
\end{equation}
Moreover, the convergences in law in \eqref{e.clt.mesoJ} and \eqref{e.clt.J} hold jointly over $p,q,z \in \Rd$ and $f \in L^1(\Rd) \cap L^\infty(\Rd)$.
\end{theorem}

The rest of this section is devoted to the proofs of Proposition~\ref{p.conv.white} and Theorem~\ref{t.conv.white}. 

\smallskip

In order to present the proof of Proposition~\ref{p.conv.white}, it is convenient to introduce the following notation for quantifying how close a random variable~$X$ is to being a Gaussian. For every real random variable $X$ and $\sigma, \lambda_1 ,\cc \ge 0$, we write
$$
X = \mcl N(\sigma^2,\lambda_1, \cc)
$$
to mean that, for every $\lambda \in (-\lambda_1,\lambda_1)$,
\begin{equation}  
\label{e.def.mclN}
\Ll|\log \E\Ll[\exp\Ll(\lambda  X\Rr)\Rr] - \frac{\sigma^2 \lambda^2 }{2} \Rr| \le \cc \lambda^2.
\end{equation}
By the injectivity of the Laplace transform (cf.~\cite[Theorem~26.2]{billingsley}), the statement that $X = \mcl N(\sigma^2,\infty,0)$ is equivalent to the statement that $X$ is a centered Gaussian random variable of variance $\sigma^2$. 

\smallskip

As in the paragraph preceding \eqref{e.meso.mask}, it is useful to think first about the simpler situation where $J_1(\Phi_{z,r},p,q)$ is replaced by the spatial average over $\Phi_{z,r}$ of a random field with finite range of dependence, or even by the normalized sum over boxes of an i.i.d.\ random field indexed by $\Z^d$. We thus seek a renormalization-type argument to prove a central limit theorem for such quantities.\index{central limit theorem} As an initialization step, we observe that, for every centered and sufficiently integrable random variable~$X$ and $\cc > 0$, there exists a possibly very small $\lambda_1 > 0$ such that $X = \mcl N(\E[X^2], \lambda_1,\cc)$. For the induction, if we can decompose the large-scale quantity as the normalized sum of $k$ independent random variables which are each $\mcl N(\sigma^2, \lambda_1,\cc)$, then the large-scale quantity is $\mcl N(\sigma^2, \sqrt{k} \lambda_1,\cc)$.

\smallskip

If these two steps can be justified, then it is clear that over larger and larger scales, the quantity of interest becomes arbitrarily close to $\mcl N(\sigma^2,\infty,0)$. More precisely, for any $\lambda_1 < \infty$ and $\cc > 0$, there exists $\sigma^2 \ge 0$ such that the quantity of interest becomes $\mcl N(\sigma^2,\lambda_1,\cc)$ over sufficiently large scales. For the quantity $J_1$, the difficulty is that the induction step is not exact: additivity is only true up to an error, which becomes negligible as we move to larger and larger scales. The proof consists therefore in initializing the renormalization argument at a sufficiently large scale, so that the cumulative error due to the additivity defect remains small.

\smallskip

The initialization step mentioned above is formalized in the next lemma.
\begin{lemma}
\label{l.centered}
There exists $C < \infty$ such that if a random variable $X$ satisfies 
\begin{equation}  
\label{e.integ.2X}
\E[\exp\left(2|X|\right)] \le 2 \ \text{ and } \ \ \E[X] = 0,
\end{equation}
then for every $\lambda_1 \in (0,1]$,
\begin{equation*}  
X = \mcl N\Ll(\E\left[X^2 \right], \lambda_1, C \lambda_1\Rr). 
\end{equation*}
\end{lemma}
\begin{proof}
This is a direct consequence of Lemma~\ref{l.quad.lapl}.
\end{proof}
As already stated above, the basis of the induction step lies in the following elementary observation.
\begin{lemma}
\label{l.induction}
Let $\sigma,\lambda_1,\cc \ge 0$, and let $X_1, \ldots, X_k$ be independent random variables satisfying, for every $i \in \{1,\ldots,k\}$,
\begin{equation}  
\label{e.l.induction}
X_i = \mcl N \Ll( \sigma^2,\lambda_1,\cc \Rr) .
\end{equation}
We have
\begin{equation*}  
k^{-\frac 1 2} \sum_{i = 1}^k X_i = \mcl N \Ll( \sigma^2, \sqrt{k} \lambda_1, \cc \Rr) .
\end{equation*}
\end{lemma}
\begin{proof}
By the independence assumption, if $|\lambda| <  k^{\frac12} \lambda_1$, then
\begin{equation*}  
\log \E \Ll[ \exp \Ll( \lambda k^{-\frac 1 2} \sum_{i = 1}^k X_i \Rr)  \Rr]  = \sum_{i = 1}^k \log \E \Ll[ \exp \Ll( \lambda k^{-\frac 1 2}  X_i \Rr)  \Rr].
\end{equation*}
The conclusion is then immediate from the definition of \eqref{e.l.induction}, see \eqref{e.def.mclN}.
\end{proof}
For the rest of this section, we fix $s \in (1,2)$, $\de > 0$ sufficiently small and $\al \in \R$ such that
\begin{equation}  
\label{e.fix.s.de.al}
\frac d 2(1+\de) < \al < \frac d s \wedge \Ll[\frac d 2 (1+\de) + \de\Rr].
\end{equation}
Since $\al > \frac d 2(1+\de)$, we can also fix $\frac 1 2 < \eta_1 < \eta_2 < 1$ such that
\begin{equation}
\label{e.cond.eta.2}
\eta_1 \al > \frac d 2 \quad \text{and} \quad \eta_2(1+\de) < 1.
\end{equation}
The exponent $\eta_1$ is that which appears in Proposition~\ref{p.conv.white}, while the exponent $\eta'$ in Theorem~\ref{t.conv.white} is $\eta' = \eta_1(1+\delta) < 1$.

\smallskip

With this in place, we now recall some of the results proved in Chapter~\ref{c.A1}. By Theorem~\ref{t.additivity}(iv), there exists a constant $C(s,\de,\al,d,\Lambda) < \infty$ and, for every $r \ge 1$, $z \in \Rd$ and $p,q \in B_1$, an $\F(B_{r^{1+\de}}(z))$-measurable random variable $J_1^{(\delta)}(z,r,p,q)$ such that 
\begin{equation}
\label{e.w.localization}
\left| J_1(\Phi_{z,r},p,q) - J_1^{(\delta)}(z,r,p,q)\right| \leq \O_{s}\left( Cr^{-\alpha}\right).
\end{equation}
Since $(p,q) \mapsto J_1(\Phi_{z,r},p,q)$ is a quadratic form and is bounded by a deterministic constant on $B_1 \times B_1$, we may and will assume that the same properties also hold for $(p,q) \mapsto J_1^{(\delta)}(z,r,p,q)$. As in \eqref{e.def.tJd}, we center $\Jd$ by setting
\begin{equation}
\label{e.def.tJd.2}
\tJd(z,r,p,q) := \Jd(z,r,p,q) - \E \Ll[ \Jd(z,r,p,q) \Rr].
\end{equation}
Since $\al < \frac d s$, we infer from Theorem~\ref{t.additivity}(i) and Lemma~\ref{l.change-s}(i) that there exists $C(s,\de,\al,d,\Lambda) < \infty$ such that, for every $R > r \ge 1$, $z \in \Rd$ and $p,q \in B_1$,
\begin{equation}  
\label{e.eyes.of.the.tiger}
\Ll|J_1(\Phi_{z,R},p,q) - \int_{\Phi_{z,\sqrt{R^2 - r^2}}} J_1(\Phi_{x,r},p,q) \, dx \Rr| \le \O_s \Ll( C r^{-\al} \Rr) .
\end{equation}
By \eqref{e.w.localization} and the triangle inequality, we deduce the following additivity property for $\tJd$:
\begin{equation}
\label{e.add.tJd.2}
\Ll|\tJd(z,R,p,q) - \int_{\Phi_{z,\sqrt{R^2 - r^2}}} \tJd(\cdot,r,p,q) \Rr| \le \O_s \Ll( C r^{-\al} \Rr) .
\end{equation}
We also recall from \eqref{e.induc.scale} that there exists a constant $C(s,\de,\al,\eta_1,\eta_2,d,\Lambda) < \infty$ such that, for every $R \ge 1$, $r \in [R^{\eta_1}, R^{\eta_2}]$, $z \in \Rd$, $p, q \in B_1$ and every function $f \in L^\infty(\Rd)$ satisfying $\|f\|_{L^\infty(\Rd)} \le 1$, we have
\begin{equation}  
\label{e.basic.estimate.f}
\fint_{\cu_R(z)} f \,\tJd(\cdot,r,p,q) =  \O_s \Ll( C R^{-\frac d 2} \Rr) .
\end{equation}
The restriction $r \le R^{\eta_2}$ was convenient for the proof of \eqref{e.basic.estimate.f}, but can easily be lifted using the additivity property. This is the purpose of the first part of the next lemma; the other parts provide  similar stability estimates for the quantity on the left side of \eqref{e.basic.estimate.f}.
\begin{lemma}
\label{l.estim.error.f}
\emph{(i)}
There exists a constant $C(s,\de,\eta_1,d,\Lambda) < \infty$ such that, for every $R \ge 1$, $r \in [R^{\eta_1},R]$, $z \in \Rd$, $p,q \in B_1$ and every function $f \in L^\infty(\Rd)$ satisfying $\|f\|_{L^\infty(\Rd)} \le 1$, we have
\begin{equation}  
\label{e.estim.error.f.a}
\fint_{\cu_R(z)} f \,\tJd(\cdot,r,p,q) =  \O_s \Ll( C  R^{-\frac d 2} \Rr) .
\end{equation}

\emph{(ii)}
For each $\be \in [0,1]$, there exists a constant $C(\be,s,\de,\eta_1,d,\Lambda) < \infty$ such that, for every $R \ge 1$, $r \in \Ll[R^{\be \eta_1},R^\frac{\beta}{1+\de}\Rr]$, $z \in \Rd$, $p,q \in B_1$ and every function $f \in L^\infty(\Rd)$ satisfying 
\begin{equation}
\label{e.estim.f.boundary}
\|f\|_{L^\infty(\Rd)} \le 1 \quad \text{and} \quad \supp f \subset \cu_R(z)  \setminus \cu_{R - R^\beta}(z),
\end{equation}
we have the improved estimate
\begin{equation}  
\label{e.estim.error.f.b}
\fint_{\cu_R(z)} f \,\tJd(\cdot,r,p,q) =  \O_s \Ll( C  R^{-\frac d 2 - \frac {1-\be}{2}} \Rr) .
\end{equation}

\emph{(iii)} For each $\be \in [\eta_1,1)$, there exists an exponent $\eps(\be,s,\de,\eta_1,d) > 0$ and a constant $C(\be,\eps,s,\de,\eta_1, d, \Lambda) < \infty$ such that, for every $R \ge 1$, $r \in [R^{\eta_1}, R^\be]$, $z \in \Rd$ and $p,q \in B_1$, we have
\begin{equation*}  
\fint_{\cu_R(z)} \tJd(\cdot,r,p,q) = \fint_{\cu_R(z)} \tJd(\cdot,R^{\eta_1},p,q) + \O_s \Ll( C R^{- \frac d 2 - \eps }\Rr) .
\end{equation*}
\end{lemma}
\begin{proof}
We fix $z = 0$ throughout for notational convenience, and start with the proof of part (i). By the additivity property in \eqref{e.add.tJd.2}, we have
\begin{multline}  
\label{e.estim.error.f2}
\fint_{\cu_R} f(x) \,\tJd(x,r,p,q) \, dx \\
= \fint_{\cu_R} f(x) \,\int_{\Phi_{x,\sqrt{r^2 - R^{2\eta_1}}}} \tJd(y,R^{\eta_1},p,q) \, dy \, dx + \O_s \Ll( C R^{-\eta_1 \al } \Rr) .
\end{multline}
By \eqref{e.cond.eta.2}, we have $\eta_1 \al > \frac d2$, and thus the error term is of lower order. We set
\begin{equation}  
\label{e.def.g.f}
g := (f \1_{\cu_R}) \ast \Phi(r^2 - R^{2\eta_1}, \cdot),
\end{equation}
and observe that, for every $x \in \Rd$,
\begin{equation}  
\label{e.g.tail.bound}
g(x) \le C \exp \Ll( -\frac{\Ll( |x|_\infty - R \Rr)_+^2 }{C r^2} \Rr) ,
\end{equation}
where we recall that, for $x = (x_1,\ldots,x_d) \in \Rd$, we write $|x|_{\infty} = \max_{1 \le i \le d} |x_i|$. 
By Fubini's theorem, the double integral on the right side of \eqref{e.estim.error.f2} is equal to 
\begin{equation*}  
R^{-d} \int_{\Rd} g \, \tJd(\cdot,R^{\eta_1},p,q),
\end{equation*}
which we decompose into
\begin{equation*}  
 \sum_{z \in R\Z^d} \fint_{\cu_R(z)} g \, \tJd(\cdot,R^{\eta_1},p,q).
\end{equation*}
The estimate \eqref{e.estim.error.f.a} therefore follows from \eqref{e.basic.estimate.f} and \eqref{e.g.tail.bound}.

\smallskip

We now turn to the proof of part (ii). We decompose the integral on the left side of \eqref{e.estim.error.f.b} into subcubes of side length $R^\be$:
\begin{equation*}  
\int_{\cu_R} f \,\tJd(\cdot,r,p,q) = \sum_{z \in R^\be \Z^d} \int_{\cu_{R^\be}(z)} f \,\tJd(\cdot,r,p,q).
\end{equation*}
By the assumption \eqref{e.estim.f.boundary} on the support of $f$, there are no more than $C R^{(1-\be)(d-1)}$ non-zero summands in the sum above. Applying the result of part (i) to each of these summands yields
\begin{equation*}  
R^{-\be d} \int_{\cu_{R^\be}(z)} f \,\tJd(\cdot,r,p,q) \le \O_s \Ll( C R^{-\be \frac d 2} \Rr) .
\end{equation*}
Moreover, since we impose $r \le R^{\frac \be {1+\de}}$, the quantity on the left side above is $\mcl F(\cu_{2R^{\be}}(z))$-measurable. We can therefore appeal to Lemmas~\ref{l.bigO.barO} and~\ref{l.barO.boxes} to obtain that
\begin{equation*}  
\int_{\cu_R} f \,\tJd(\cdot,r,p,q) \le \O_s \Ll(  R^{(1-\be)\frac{d-1}{2}+\be  \frac  d 2}  \Rr) .
\end{equation*}
Dividing by $R^{-d}$, we obtain \eqref{e.estim.error.f.b}.

\smallskip

We now turn to part (iii). Using the additivity of $\tJd$ as in the beginning of the proof but with $f \equiv 1$, we get that for
\begin{equation*}  
g := \1_{\cu_R} \ast \Phi(r^2 - R^{2\eta_1}, \cdot),
\end{equation*}
we have
\begin{equation*}  
\fint_{\cu_R} \tJd(\cdot,r,p,q) = R^{-d} \int_{\Rd} g \, \tJd(\cdot,R^{\eta_1},p,q) + \O_s \Ll( C R^{-\eta_1 \al} \Rr) ,
\end{equation*}
and we recall that $\eta_1 \al > \frac d 2$. We now decompose the integral on the right side into
\begin{equation*}  
\sum_{z \in R \Z^d} \int_{\cu_R(z)} g \,  \tJd(\cdot,R^{\eta_1},p,q),
\end{equation*}
and proceed to show that, for some small $\eps(\be,s,\de,\eta_1,d) > 0$,
\begin{equation}  
\label{e.ohcomeon1}
\int_{\cu_R} g \,  \tJd(\cdot,R^{\eta_1},p,q) = \int_{\cu_R(z)}   \tJd(\cdot,R^{\eta_1},p,q) + \O_s \Ll( C R^{-\frac d 2 - \eps} \Rr) ,
\end{equation}
while
\begin{equation}  
\label{e.ohcomeon2}
\sum_{z \in R \Z^d \setminus \{0\}} \int_{\cu_R(z)} g \,  \tJd(\cdot,R^{\eta_1},p,q) = \O_s \Ll( C R^{-\frac d 2 - \eps} \Rr) .
\end{equation}
The estimate \eqref{e.ohcomeon1} follows from the observations that $0 \le g \le 1$, 
\begin{equation*}  
\sup\Ll\{|1-g(x)| \ : \ x \in \cu_{R - R^{\frac {1+\be}{2}}}\Rr\} \le C R^{-100d},
\end{equation*}
and an application of parts (i) and (ii) of the lemma. The proof of \eqref{e.ohcomeon2} is obtained similarly, using that, for every $x \in \Rd$,
\begin{equation*}  
0 \le g(x) \le C \exp \Ll( - \frac{(2|x|_\infty-R)_+^2}{C R^{2\be}}   \Rr) .
\end{equation*}
This completes the proof of part (iii) of the lemma.
\end{proof}

Finally, we need a technical lemma justifying that the statement that~$X = \mcl N(\sigma^2,\lambda_1,\cc)$ is stable to perturbations by $\O_s \Ll( \theta \Rr)$-bounded random variables, if $\theta$ is sufficiently small.

\begin{lemma}  
\label{l.mclN.perturb}
Fix $\lambda_2 > 0$. There exists $C(s,\lambda_2)  <\infty$ such that, for every $\sigma, \cc > 0$, $\lambda_1 \in [0,\lambda_2]$ and $\theta \in [0,1]$, if $X_1$, $X_2$ are two centered random variables satisfying
\begin{equation}
\label{e.X1}
X_1 = \mcl N \Ll( \sigma^2, \lambda_1, \cc \Rr) \ \ \text{ and } \ \  X_2 = \O_s \Ll( \theta \Rr) ,
\end{equation}
then 
\begin{equation}  
\label{e.X12}
X_1 + X_2 = \mcl N \Ll( \sigma^2, (1-\sqrt{\theta})\lambda_1, \cc + C \sqrt{\theta}(1+\sigma^2+\cc) \Rr) .
\end{equation}
\end{lemma}
\begin{proof}
Write $X := X_1 + X_2$, and let $\zeta, \zeta'\in (1,\infty)$ be such that $\frac 1 \zeta + \frac 1 {\zeta'} = 1$. By H\"older's inequality,
\begin{equation}
\label{e.mydear1}
\log \E \Ll[ \exp \Ll( \lambda  X   \Rr)  \Rr] \le \frac 1 \zeta \log \E \Ll[ \exp \Ll( \zeta \lambda  X_1 \Rr)  \Rr] + \frac 1 {\zeta'} \log \E \Ll[ \exp \Ll(\zeta' \lambda  X_2 \Rr)  \Rr] ,
\end{equation}
and conversely,
\begin{equation}  
\label{e.mydear2}
\log \E \Ll[ \exp \Ll( \lambda  X_1 \Rr)  \Rr] \le \frac 1 \zeta \log \E \Ll[ \exp \Ll( \zeta \lambda  X \Rr)  \Rr] + \frac 1 {\zeta'} \log \E \Ll[ \exp \Ll(-\zeta' \lambda  X_2 \Rr)  \Rr] .
\end{equation}
By the first part of \eqref{e.X1}, for every $|\lambda| <  \zeta^{-1} \, \lambda_1$,
\begin{equation}  
\label{e.whatcanwedo}
\log \E \Ll[ \exp \Ll( \zeta \lambda  X_1 \Rr)  \Rr] \le \Ll( \frac {\sigma^2} 2 + \cc \Rr) (\zeta \lambda)^2.
\end{equation}
By the second part of \eqref{e.X1} and Lemma~\ref{l.bigO.barO}, there exists $C(s) < \infty$ such that, for every $\lambda \in \R$,
\begin{equation*}  
\log \E \Ll[ \exp \Ll(C^{-1}  \zeta' \lambda  X_2 \Rr)  \Rr] \le  (\zeta' \theta \lambda)^2 \vee |\zeta' \theta\lambda|^{\frac s{s-1}}.
\end{equation*}
Specializing to the case in which $\zeta  =  (1-\theta^{\frac12} )^{-1}$ and $\zeta' = \theta^{-\frac 1 2}$, we deduce  the existence of a constant $C(s,\lambda_2) < \infty$ such that, for every $|\lambda| \le \lambda_2$,
\begin{equation*}  
\log \E \Ll[ \exp \Ll(  \zeta' \lambda  X_2 \Rr)  \Rr]  \le C \Ll(\sqrt{\theta} \lambda\Rr)^2.
\end{equation*}
Using \eqref{e.mydear1} with \eqref{e.whatcanwedo} and the previous display yields, for every $|\lambda| < (1-\theta^{\frac12} ) \lambda_1$,
\begin{equation*}  
\log \E \Ll[ \exp \Ll( \lambda  X   \Rr)  \Rr] \le \Ll( \frac {\sigma^2} 2 + \cc \Rr)\zeta \lambda^2 + C \theta^{\frac12}  \lambda^2.
\end{equation*}
Since
\begin{equation*}  
\zeta - 1 = \left(1-\theta^{\frac12}\right)^{-1} -1 \le C \theta^{\frac12},
\end{equation*}
we obtain one side of the two-sided inequality implicit in \eqref{e.X12}. The other inequality is proved in a similar way, using \eqref{e.mydear2} instead of \eqref{e.mydear1}.
\end{proof}

\begin{proof}[Proof of Proposition~\ref{p.conv.white}]
We first observe that the statement of $\msf Q(p,\ahom p) = 0$ is a consequence of the claimed convergence \eqref{e.clt} and the results of Chapter~\ref{c.A1}. Indeed, by Lemma~\ref{l.shortcuts} and \eqref{e.Jexpressbyv}, there exists $C(s,d,\Lambda) < \infty$ such that, for every $r \ge 1$, $z \in \Rd$ and $p \in B_1$,
\begin{equation*}  
J(\Phi_{z,r},p,\ahom p) \le \O_{s/2} \Ll( C r^{-d} \Rr) ,
\end{equation*}
and therefore $\msf Q(p,\ahom p) = 0$ follows from \eqref{e.clt} and the fact that $\eta_1 > \frac 1 2$.

\smallskip

We turn to the proof of \eqref{e.clt}. By \eqref{e.w.localization}, Lemma~\ref{l.estim.error.f}(iii) and~$\eta_1 \al > \frac d 2$, it suffices to show that
\begin{equation*}
r^{\frac d 2} \fint_{\cu_r(z)} \tJd(\cdot,r^{\eta_1},p,q)  \xrightarrow[r \to \infty]{\mathrm{(law)}}  \msf Q(p,q),
\end{equation*}
with the convergence holding jointly over $p,q \in \Rd$. 
Since $(p,q) \mapsto \tJd(x,r,p,q)$ is a quadratic form, this amounts to showing the joint convergence of the entries of the associated random $2d$-by-$2d$ matrix. By Lemma~\ref{l.cramer.wold}, it suffices to show the convergence in law of arbitrary linear combinations of these entries. These entries can be computed from the bilinear form associated with $\tJd$, that is, 
\begin{equation*}  
(p,q,p',q') \mapsto \frac 1 4 \Ll( \tJd(x,r,p+p',q+q') - \tJd(x,r,p-p',q-q')\Rr) ,
\end{equation*}
for particular choices of $p,q,p',q'$ in the set 
\begin{equation*}  
\mcl B_0 := \{e_i \ : \ 1 \le i \le d\} \cup \{0\} \subset \Rd.
\end{equation*}
We define
\begin{equation}  
\label{e.def.mclB1}
\mcl B_1 := \{v \pm v' \ : \ v,v' \in \mcl B_0\}.
\end{equation}
Using also Remark~\ref{r.cramer.wold}, it therefore suffices to show that, for each fixed finite family of scalars
\begin{equation}  
\label{e.family.scalars}
\Ll\{ \lambda(v,v')\Rr\}_{v,v' \in \mcl B_1} \in [-1,1]^{\mcl B_1^2},
\end{equation}
we have the convergence in law, as $r\to \infty$, of the random variable 
\begin{equation*}  
\sum_{v, v' \in \mcl B_1}\lambda(v, v') \,  r^{\frac d 2} \fint_{\cu_r(z)} \tJd(\cdot,r^{\eta_1},v, v')
\end{equation*}
to a Gaussian random variable. 

\smallskip

We henceforth fix~$\Ll\{ \lambda(v,v')\Rr\}_{v,v' \in \mcl B_1}$ as in \eqref{e.family.scalars} and, for every $x \in \Rd$ and $r \ge 1$, use the shorthand notation
\begin{equation}  
\label{e.lin.comb}
\tJd(x,r) := 
\sum_{v, v' \in \mcl B_1}  \lambda(v, v')\tJd(x,r,v, v').
\end{equation}
Since the number of terms in the sum on the right side is bounded by $C(d) < \infty$ and the coefficients themselves belong to $[-1,1]$, it is clear that the estimates we proved for $\tJd(x,r,p,q)$ with $p,q \in B_1$ transfer to the quantity $\tJd(x,r)$, after adjusting the multiplicative constant.



\smallskip

In order to prove that, for every $z \in \Rd$, the random variable
\begin{equation*}  
r^{\frac d 2} \fint_{\cu_r(z)} \tJd(\cdot, r^{\eta_1})
\end{equation*}
converges in law to a Gaussian as $r$ tends to infinity, it suffices to show that, for any given $\lambda_2 < \infty$ and $\cc > 0$, there exists $\sigma \in [0,\infty)$ such that, for every $z \in \Rd$ and $r$ sufficiently large,
\begin{equation}
\label{e.clt.lambda}
r^{\frac d 2} \fint_{\cu_r(z)} \tJd(\cdot,r^{\eta_1})   = \mcl N\Ll(\sigma^2,  \lambda_2, \cc \Rr).
\end{equation}
Indeed, this statement implies that, for any given bounded interval, the Laplace transform of the random variable on the left side of \eqref{e.clt.lambda} converges to a parabola as $r$ tends to infinity; and this suffices to prove the desired convergence in law. 

\smallskip

Let us therefore fix $\lambda_2 \in [1,\infty)$ and proceed to prove \eqref{e.clt.lambda} by an induction on the scale.
For every $r_1 \ge 2$ and $\sigma \ge 0$, $\lambda_1 \in (0,\lambda_2)$ and $\cc \ge 0$, we denote by $\msf A(r_1,\sigma^2,\lambda_1,\cc)$ the statement that, for every $r \in \Ll[\frac {r_1} 2, r_1\Rr]$ and $z \in \Rd$, we have
\begin{equation}  
\label{e.def.mscA.2}
r^\frac d 2 \fint_{\cu_r(z)} \tJd(\cdot,r^{\eta_1})  = \mcl N(\sigma^2,\lambda_1,\cc).
\end{equation}
Note that by Lemmas~\ref{l.estim.error.f}(i) and \ref{l.bigO.barO}, there exists $C_0(\lambda_2,s,\de,\eta_1,d,\Lambda) < \infty$ such that, for every $r \ge 1$, $z \in \Rd$ and $\lambda \in (0,\lambda_2)$,
\begin{equation*}  
\log \E \Ll[ \exp \Ll( \lambda r^\frac d 2 \fint_{\cu_r(z)} \tJd(\cdot,r^{\eta_1})  \Rr)  \Rr] \le \frac{C_0 \lambda^2}{4}. 
\end{equation*}
We can therefore restrict a priori the range of values of interest for the parameters~$\sigma^2$ and $\cc$ to the interval $[0,C_0]$. 

\smallskip

In the course of the argument, the value of the exponent $\eps(s,\de,\al,\eta_1,\eta_2,d) > 0$ and of the constant $C(s,\de,\al,\eta_1,\eta_2,\lambda_2, d,\Lambda) < \infty$ may change from one occurrence to the next without further mention. 

\smallskip

We decompose the rest of the proof into four steps. 

\smallskip

\emph{Step 1.} In this step, we show that, for every $r$ sufficiently large (depending on the parameters $(s,\delta,\eta_1,d,\Lambda)$), there exists $\sigma^2(r) \in [0,C_0]$ such that, for every $z \in \Rd$ and $\lambda_1 \in \Ll(0,\frac 1 2\Rr]$,
\begin{equation}
\label{e.init0}
r^\frac d 2 \fint_{\cu_r(z)} \tJd(\cdot,r^{\eta_1})  = \mcl N(\sigma^2,\lambda_1,C\lambda_1+ C r^{-\eps}).
\end{equation}
By Lemmas~\ref{l.centered} and \ref{l.estim.error.f}, for each $r \ge 1$, there exists $\sigma^2(r) \in [0,C_0]$ such that, for every $\lambda_1 \in (0,1]$,
\begin{equation*}  
r^{\frac d 2} \fint_{\cu_r} \tJd(\cdot,r^{\eta_1}) = \mcl N\Ll(\sigma^2, \lambda_1, C \lambda_1 \Rr).
\end{equation*}
This shows the claim with the additional constraint of $z = 0$. In order to obtain the full result, we first note that by $\Z^d$-stationarity, it suffices to show that \eqref{e.init0} holds uniformly over $z \in [0,1)^d$. In this case, applying Lemma~\ref{l.estim.error.f}(ii) twice, we can replace integration over $\cu_r$ by integration over $\cu_r \cap \cu_r(z)$ and then over $\cu_r(z)$, at the price of an error of $\O_s \Ll( C r^{-\frac d 2 - \eps} \Rr)$. By Lemma~\ref{l.mclN.perturb}, this completes the proof of \eqref{e.init0}.

\smallskip

\emph{Step 2.} We now complete the initialization of the induction argument by showing that, for every $r_1$ sufficiently large, there exists $\sigma^2(r_1) \in [0,C_0]$ such that 
\begin{equation}
\label{e.init1}
\text{for every }\lambda_1 \in \Ll(0,\tfrac 1 2\Rr], \ \ \msf A(r_1,\sigma^2, \lambda_1, C\lambda_1+ C r_1^{-\eps}) \ \text{holds}.
\end{equation}
Let $\gamma > 1$ be selected sufficiently close to $1$ as a function of $\eta_1$ and $\delta$ so that
\begin{equation}
\label{e.fix.gamma}
\gamma \eta_1 (1+\de) < 1.
\end{equation}
We fix $r_0$ sufficiently large, let $\sigma^2(r_0) \in [0,C_0]$ be as given by the previous step, and aim to show that, for every $r \in \Ll[\frac12 r_0^{\gamma}, r_0^{\gamma} \Rr]$, $z \in \Rd$ and $\lambda_1 \in \Ll(0,\tfrac 1 2\Rr]$, we have
\begin{equation}  
\label{e.init2}
r^{-\frac d 2} \int_{\cu_r(z)} \tJd(\cdot, r^{\eta_1}) = \mcl N \Ll( \sigma^2,\lambda_1,C \lambda_1 + C r^{-\eps} \Rr) .
\end{equation}

\smallskip

The argument for \eqref{e.init2} and the rest of the proof will be completely translation-invariant, so we focus on showing the statement with the additional constraint of $z = 0$ for notational convenience. For every $r \in \Ll[\frac{r_0^{\gamma}}{2}, r_0^{\gamma} \Rr]$, we consider the decomposition
\begin{equation*}  
r^{-\frac d 2} \int_{\cu_{r}} \tJd(\cdot,r^{\eta_1}) =r^{-\frac d 2}  \sum_{z \in r_0 \Zd \cap \cu_r} \int_{\cu_{r_0}(z)} \tJd(\cdot,r^{\eta_1}).
\end{equation*}
We denote by $\cut_r(z)$ the trimmed cube
\begin{equation}  
\label{e.def.trimmed}
\cut_r(z) := \cu_{r - r^{\eta_1(1+\de)}}(z),
\end{equation}
By the choice of the exponent $\gamma$ in \eqref{e.fix.gamma}, we can appeal to Lemma~\ref{l.estim.error.f}(ii) and get that, for every $z \in \Rd$, 
\begin{equation}  
\label{e.trimmed.cube.is.back}
\int_{\cu_{r_0}(z)} \tJd(\cdot,r^{\eta_1}) - \int_{\cut_{r_0}(z)} \tJd(\cdot,r^{\eta_1}) = \O_s \Ll( C r_0^{\frac d 2 - \eps} \Rr) .
\end{equation}
Moreover, the random variable on the left side above is $\mcl F(\cu_{2r_0}(z))$-measurable. By Lemmas~\ref{l.bigO.barO} and \ref{l.barO.boxes}, we deduce that
\begin{equation*}  
r^{-\frac d 2} \sum_{z \in r_0 \Zd \cap \cu_r} \Ll( \int_{\cu_{r_0}(z)} \tJd(\cdot,r^{\eta_1}) - \int_{\cut_{r_0}(z)} \tJd(\cdot,r^{\eta_1}) \Rr)  = \O_s \Ll(C r^{-\frac d 2} r_0^{-\eps} \Rr) .
\end{equation*}
In other words, this difference is of lower order. 
We also note that, by the result of the previous step, \eqref{e.trimmed.cube.is.back} and Lemma~\ref{l.mclN.perturb}, we have
\begin{equation}
\label{e.trimmed.mclN}
r_0^{-\frac d 2} \int_{\cut_{r_0}(z)} \tJd(\cdot,r^{\eta_1}) = \mcl N \Ll( \sigma^2, (1-Cr^{-\eps})\lambda_1, C \lambda_1 + C r^{-\eps} \Rr) .
\end{equation}
Lemma~\ref{l.induction} then yields
\begin{equation*}  
r^{-\frac d 2} \sum_{z \in r_0 \Zd \cap \cu_r}  \int_{\cu_{r_0}(z)} \tJd(\cdot,r^{\eta_1}) = \mcl N\Ll(\sigma^2, 2\lambda_1,C \lambda_1 + C r^{-\eps}\Rr).
\end{equation*}
Another application of Lemma~\ref{l.mclN.perturb} completes the proof of \eqref{e.init2}.

\smallskip

\emph{Step 3.} In this step, we show the induction statement, namely that, for every $r_1$ sufficiently large, $\lambda_1 \in (0,\lambda_2)$ and $\sigma^2, \cc \in [0,C_0]$, we have
\begin{equation}
\label{e.clt.induction}
\msf A \Ll( r_1, \sigma^2, \lambda_1, \cc \Rr)  \quad \implies \quad \msf A \Ll( 2r_1, \sigma^2,  2^{\frac d 2} \lambda_1 (1-Cr_1^{-\eps}), \cc + C r_1^{-\eps}  \Rr) .
\end{equation}
Let $r \in \Ll[ \frac12r_1, r_1 \Rr] $. We decompose 
\begin{equation}  
\label{e.decompose.induction}
(2r)^{-\frac d 2} \int_{\cu_{2r}} \tJd(\cdot, (2r)^{\eta_1}) = (2r)^{-\frac d 2} \sum_{z \in \mcl Z} \int_{\cu_r(z)} \tJd(\cdot,(2r)^{\eta_1}),
\end{equation}
where we denote
\begin{equation*}  
\mcl Z := \Ll\{ -\frac r 2 , \frac r 2 \Rr\} ^d.
\end{equation*}
Recall the definition of the trimmed cube $\cut_r(z)$ in \eqref{e.def.trimmed}. By Lemma~\ref{l.estim.error.f}(ii) and (iii), we have
\begin{equation}  
\label{e.error.induction}
r^{-\frac d 2} \Ll|\int_{\cu_r(z)} \tJd(\cdot,(2r)^{\eta_1}) - \int_{\cut_r(z)} \tJd(\cdot,r^{\eta_1}) \Rr|  = \O_s \Ll( C r^{-\frac d 2 -\eps} \Rr) .
\end{equation}
By the assumption of $\msf A(r_1,\sigma^2,\lambda_1,\cc)$ and Lemma~\ref{l.estim.error.f}(ii), we have
\begin{equation*}  
r^{-\frac d 2} \int_{\cut_r(z)} \tJd(\cdot,r^{\eta_1}) = \mcl N \Ll( \sigma^2, (1-Cr^{-\eps})\lambda_1, \cc + C r^{-\eps} \Rr) .
\end{equation*}
Applying Lemma~\ref{l.induction}, we get
\begin{equation*}  
(2r)^{-\frac d 2} \sum_{z \in \mcl Z} \int_{\cu_r(z)} \tJd(\cdot,(2r)^{\eta_1}) = \mcl N \Ll( \sigma^2, 2^{\frac d 2} \lambda_1 (1-Cr^{-\eps}), \cc + C r^{-\eps} \Rr) .
\end{equation*}
Combining this with \eqref{e.decompose.induction}, \eqref{e.error.induction} and Lemma~\ref{l.mclN.perturb} yields the announced implication \eqref{e.clt.induction}.

\smallskip

\emph{Step 4.} We conclude for the proof of \eqref{e.clt.lambda}. Let $\cc \in (0,1]$. By the result of Step 2, for every $r_1$ sufficiently large, there exists $\sigma^2 \in [0,C_0]$ and $\lambda_1 > 0$ sufficiently small such that
\begin{equation*}  
\msf A \Ll( r_1,\sigma^2 ,\lambda_1, \cc  \Rr) \ \ \text{holds}.
\end{equation*}
Applying the induction step \eqref{e.clt.induction} repeatedly, we obtain that, for every $k$ sufficiently large,
\begin{equation*}  
\msf A \Ll( 2^k r_1, \sigma^2, \lambda_2, \cc + C r_1^{-\eps} \Rr) \ \ \text{holds}.
\end{equation*}
Choosing $r_1$ sufficiently large in terms of $\cc$, it follows that, for every $k$ sufficiently large,
\begin{equation*}  
\msf A \Ll( 2^k r_1, \sigma^2, \lambda_2, 2\cc \Rr) \ \ \text{holds}.
\end{equation*}
This is \eqref{e.clt.lambda}, up to a redefinition of $\cc$. The proof is therefore complete. 
\end{proof}
\begin{remark}  
\label{r.Laplace.converges}
It follows from the proof of Proposition~\ref{p.conv.white} that, for each $\lambda_2 < \infty$ and $\cc > 0$, there exists a sufficiently large scale $r_0$ such that, for every $r \ge r_0$ and $z \in \Rd$, 
\begin{equation*}  
r^{\frac d 2} \fint_{\cu_r(z)} \tJd(\cdot,r^\eta) = \mcl N \Ll( \sigma^2, \lambda_2, \cc \Rr) ,
\end{equation*}
using the notation from \eqref{e.lin.comb}. For a fixed $z \in \Rd$, this can also be proved directly from the statement of Proposition~\ref{p.conv.white} and the a priori stochastic integrability bound given by Lemma~\ref{l.estim.error.f}(i). The statement with arbritrary $z \in \Rd$ can then be recovered by stationarity and Lemma~\ref{l.estim.error.f}(ii).
\end{remark}

We now proceed with the proof of Theorem~\ref{t.conv.white}, which is deduced from Proposition~\ref{p.conv.white} using the approximate additivity and locality of $J_1$. 
\begin{proof}[Proof of Theorem~\ref{t.conv.white}]
The exponent $\eta'$ appearing in the statement is defined by
\begin{equation}  
\label{e.def.eta'.gff}
\eta':= \eta_1(1+\de) < 1.
\end{equation}
We first prove \eqref{e.clt.mesoJ}, and then obtain \eqref{e.clt.J} as a consequence. By \eqref{e.w.localization}, in order to prove \eqref{e.clt.mesoJ}, it suffices to show that
\begin{equation}  
\label{e.clt.mesoJ.2}
r^{\frac d 2} \int_{\Rd} f \Ll( x \Rr) \tJd(rx,r^\eta,p,q)    \, dx \\
\xrightarrow[r \to \infty]{\mathrm{(law)}}  \int_{\Rd} f(x) \W(x,p,q) \, dx,
\end{equation}
and that the convergence in law holds jointly over $f \in L^1(\Rd) \cap L^\infty(\Rd)$ and $p,q \in \Rd$. We fix a family of scalars
$\Ll(\lambda(v,v')\Rr)_{v,v' \in \mcl B_1} \in [-1,1]^{\mcl B_1^2}$, 
with $\mcl B_1$ as in \eqref{e.def.mclB1}, and recall the shorthand notation $\tJd(\cdot,r)$ from \eqref{e.lin.comb}. We know from Proposition~\ref{p.conv.white} that as $r$ tends to infinity, the random variable 
\begin{equation*}  
r^{\frac d 2} \fint_{\cu_r} \tJd(\cdot,r^{\eta_1})
\end{equation*}
converges in law to a centered Gaussian; we denote the variance of this Gaussian by~$\sigma^2$. 
As in the beginning of the proof of Proposition~\ref{p.conv.white}, by appealing to Lemma~\ref{l.cramer.wold}, we see that in order to prove the joint convergence in \eqref{e.clt.mesoJ.2}, it suffices to show that as $r$ tends to infinity, the individual random variable
\begin{equation}  
\label{e.that.should.be.clt}
r^{\frac d 2} \int_{\Rd} f \Ll( x \Rr) \tJd(rx,r^\eta)   \, dx 
\end{equation}
converges in law to a centered Gaussian of variance $\sigma^2 \|f\|_{L^2(\Rd)}^2$. 

We now proceed to prove this convergence in law. Without loss of generality, we assume that
\begin{equation}
\label{e.normalize.f.gff}
\|f\|_{L^1(\Rd)} + \|f\|_{L^\infty(\Rd)} \le 1.
\end{equation}
By the additivity property in~\eqref{e.add.tJd.2}, Lemma~\ref{l.sum-O} and \eqref{e.normalize.f.gff}, we have
\begin{multline*}  
\int_{\Rd} f \Ll(  x  \Rr) \tJd(rx,r^\eta)   \, dx \\
=  \int_{\Rd} f \Ll(  x \Rr) \int_{\Phi_{rx,\sqrt{r^{2\eta} - r^{2\eta_1}}}} \tJd(y,r^{\eta_1}) \, dy   \, dx + \O_s \Ll( C r^{- \eta_1 \al}  \Rr) ,
\end{multline*}
and we recall that $\eta_1 \al > \frac d 2$. The error term is therefore of lower order, and the double integral on the right side above can be rewritten as
\begin{multline*}  
r^{-d} \int_{\Rd} \int_{\Rd}  f \Ll( \tfrac x r \Rr) \Phi \Ll( r^{2\eta} - r^{2\eta_1}, x-y \Rr) \, \tJd(y,r^{\eta_1}) \, dx \, dy \\
= r^{-d} \int_{\Rd} g(y,r) \, \tJd(y,r^{\eta_1}) \, dy,
\end{multline*}
where we set
\begin{equation*}  
g(y,r) := \int_{\Rd} f \Ll( \tfrac x r \Rr) \Phi \Ll(r^{2\eta} - r^{2\eta_1}, x-y \Rr) \, dx.
\end{equation*}
We fix an exponent $\kappa$ function of $\eta_1$, $\eta$ and $\de$ in such a way that
\begin{equation}
\label{e.def.kappa}
 \eta'= \eta_1 (1+\de) < \kappa < \eta,
\end{equation}
and use the decomposition
\begin{equation}  
\label{e.decomp.with.g}
\int_{\Rd} g(y,r) \, \tJd(y,r^{\eta_1}) \, dy = \sum_{z \in r^\kappa \Z^d} \int_{\cu_{r^\kappa(z)}} g(y,r) \, \tJd(y,r^{\eta_1}) \, dy.
\end{equation}
By \eqref{e.normalize.f.gff} and the definition of $g$, there exists a constant $C(d) < \infty$ such that, for every $r \ge 1$ and $y,y' \in \Rd$ satisfying $|y-y'| \le r^\eta$, we have
\begin{equation}  
\label{e.diff.g.yy}
\Ll|g(y,r) - g(y',r) \Rr| \le C\frac{|y-y'|}{r^{(d+1)\eta}}  \int_{\Rd} \Ll| f \Ll( \tfrac x r \Rr)  \Rr| \exp \Ll( -\frac{|x-y|^2}{Cr^{2\eta}} \Rr) \, dx.
\end{equation}
Moreover, by \eqref{e.def.kappa}, each summand on the right side of \eqref{e.decomp.with.g} is $\mcl F(\cu_{2r^{\kappa}})$-measurable. By Lemmas~\ref{l.estim.error.f}(i), \ref{l.bigO.barO} and~\ref{l.barO.boxes}, we get
\begin{multline}  
\label{e.first.approx.clt.jj}
\Ll|\sum_{z \in r^\kappa \Z^d} \int_{\cu_{r^\kappa(z)}} g(y,r) \, \tJd(y,r^{\eta_1}) \, dy - \sum_{z \in r^\kappa \Z^d} g(z,r)\int_{\cu_{r^\kappa(z)}}  \tJd(y,r^{\eta_1}) \, dy \Rr| \\
\le
\O_s \Ll( Cr^{\frac{d\kappa}{2}} \Ll[ \sum_{z \in r^{\kappa} \Z^d} \Ll( \frac{r^{\kappa}}{r^{(d+1)\eta}} \int_{\Rd} \Ll| f \Ll( \tfrac x r \Rr)  \Rr| \exp \Ll( -\frac{|x-z|^2}{Cr^{2\eta}} \Rr)  \, dx \Rr)^2 \Rr]^\frac 1 2   \Rr) .
\end{multline}
By Jensen's inequality, the sum on the right side above is bounded by
\begin{equation*}  
C \sum_{z \in r^{\kappa} \Z^d}  \frac{r^{2\kappa}}{r^{(d+2)\eta}} \int_{\Rd} \Ll| f \Ll( \tfrac x r \Rr)  \Rr|^2 \exp \Ll( -\frac{|x-z|^2}{Cr^{2\eta}} \Rr)  \, dx 
\le
C \frac{r^{2\kappa}}{r^{(d+2)\eta}} r^{d(\eta - \kappa)} \int_{\Rd} \Ll| f \Ll( \tfrac x r \Rr)  \Rr|^2 \, dx.
\end{equation*}
By the normalization \eqref{e.normalize.f.gff}, we have $\|f\|_{L^2(\Rd)} \le 1$, and thus the left side of \eqref{e.first.approx.clt.jj} is bounded by
\begin{equation*}  
\O_s \Ll( C r^{-\frac d 2 - \Ll( \eta - \kappa \Rr)} \Rr) ,
\end{equation*}
and we recall that $\eta - \kappa > 0$. 

\smallskip

By a very similar reasoning, using part (ii) of Lemma~\ref{l.estim.error.f} instead of part (i), we also have
\begin{multline}  
\label{e.la.pelea}
\Ll|\sum_{z \in r^\kappa \Z^d} g(z,r)\int_{\cu_{r^\kappa(z)}}  \tJd(y,r^{\eta_1}) \, dy - \sum_{z \in r^\kappa \Z^d} g(z,r)\int_{\cut_{r^\kappa(z)}}  \tJd(y,r^{\eta_1}) \, dy  \Rr| \\
\le
\O_s \Ll( C r^{-\frac d 2 - \eps} \Rr) ,
\end{multline}
where here we define the trimmed cube to be
\begin{equation*}  
\cut_{r^\kappa}(z) := \cu_{r^{\kappa} - r^{\eta_1(1+\de)}}(z).
\end{equation*}
This definition ensures that the summands in the second sum on the right of \eqref{e.la.pelea} are independent random variables. Moreover, by Proposition~\ref{p.conv.white} (more precisely Remark~\ref{r.Laplace.converges}), for each fixed $\lambda_2 < \infty$ and $\cc > 0$, we have for every $r$ sufficiently large and $z \in \Rd$, 
\begin{equation*}  
r^{-\frac {\kappa d}{2}} \int_{\cut_{r^\kappa}(z)} \tJd(\cdot,r^{\eta_1}) = \mcl N \Ll( \sigma^2, \lambda_2, \cc \Rr) .
\end{equation*}
We deduce that
\begin{equation*}  
r^{-\frac d 2 } g(z,r)\int_{\cut_{r^\kappa}(z)} \tJd(\cdot,r^{\eta_1}) = \mcl N \Ll( r^{-(1-\kappa)d} g^2(z,r) \sigma^2, \lambda_2,r^{-(1-\kappa)d} g^2(z,r) \cc \Rr) ,
\end{equation*}
and by an immediate generalization of Lemma~\ref{l.induction}, 
\begin{multline*}  
r^{-\frac d 2} \sum_{z \in r^\kappa \Z^d} g(z,r)\int_{\cut_{r^\kappa(z)}}  \tJd(y,r^{\eta_1}) \, dy \\
 = \mcl N \Ll( r^{-(1-\kappa)d} \sum_{z \in r^\kappa \Z^d} g^2(z,r) \sigma^2, \lambda_2, r^{-(1-\kappa)d} \sum_{z \in r^\kappa \Z^d}  g^2(z,r) \cc \Rr) .
\end{multline*}
Using \eqref{e.diff.g.yy} again, we verify that
\begin{equation*}  
\Ll|r^{-(1-\kappa)d} \sum_{z \in r^\kappa \Z^d} g^2(z,r) - r^{-d} \int_{\Rd} g^2(\cdot,r) \Rr| \le C r^{-2(\eta - \kappa)}, 
\end{equation*}
and that, for $r$ sufficiently large,
\begin{equation*}  
\Ll|r^{-d} \int_{\Rd} g^2(\cdot,r)  - \|f\|_{L^2(\Rd)}^2\Rr| \le \cc.
\end{equation*}
Using also the elementary observation that, for $|x| \le \cc'$, 
$$
\mcl N(\sigma^2 + x,\lambda_2,\cc) = \mcl N(\sigma^2,\lambda_2,\cc+\cc'),
$$
and that $\|f\|_{L^2(\Rd)} \le 1$, we have thus shown that, for $r$ sufficiently large,
\begin{equation*}  
r^{-\frac d 2 } g(z,r)\int_{\cut_{r^\kappa}(z)} \tJd(\cdot,r^{\eta_1}) = \mcl N\Ll(\sigma^2\|f\|_{L^2(\Rd)}^2 , \lambda_2, 3 \cc\Rr).
\end{equation*}
This completes the proof of the fact that as $r$ tends to infinity, the random variable in \eqref{e.that.should.be.clt} converges to a centered Gaussian of variance $\sigma^2 \|f\|_{L^2(\Rd)}^2$.

\smallskip

In order to complete the proof of Theorem~\ref{t.conv.white}, there remains to show the convergence in \eqref{e.clt.J}. By the additivity property of $J_1$ recalled in \eqref{e.eyes.of.the.tiger}, we have
\begin{multline*}  
\Bigg|   \Ll(J_1(\Phi_{rz,r},p,q) - \E \Ll[ J_1(\Phi_{rz,r},p,q) \Rr] \Rr) \\
 - \int_{\Rd} \Phi(r^{2} - r^{2\eta}, x-rz) \Ll( J_1(\Phi_{x,r^\eta},p,q) - \E \Ll[ J_1(\Phi_{x,r^\eta},p,q) \Rr] \Rr)  \, dx \Bigg| \le \O_s \Ll( C r^{-\eta \al} \Rr) ,
\end{multline*}
with $\eta \al > \frac d 2$. Using Theorem~\ref{t.additivity}(iii), we see that in the integral above, we can replace $\Phi(r^{2} - r^{2\eta}, \cdot - rz)$ by $\Phi(r^2, \cdot - rz) = \Phi \Ll( 1, \frac{\cdot}{r}  - z \Rr) $, up to an error of lower order. Up to a change of variables, the resulting integral is then of the form of the left side of \eqref{e.clt.mesoJ}, with $f := \Phi(1, \cdot - z)$. The announced convergence \eqref{e.clt.J} is thus proved. That the convergences in law in \eqref{e.clt.mesoJ} and \eqref{e.clt.J} hold jointly as stated in the proposition is a consequence of Remark~\ref{r.how.you.doing}.
\end{proof}


\section{Convergence of the correctors to a Gaussian free field}
\label{s.conv.gff}

\begin{figure}
\centering
\includegraphics[width=.8\linewidth]{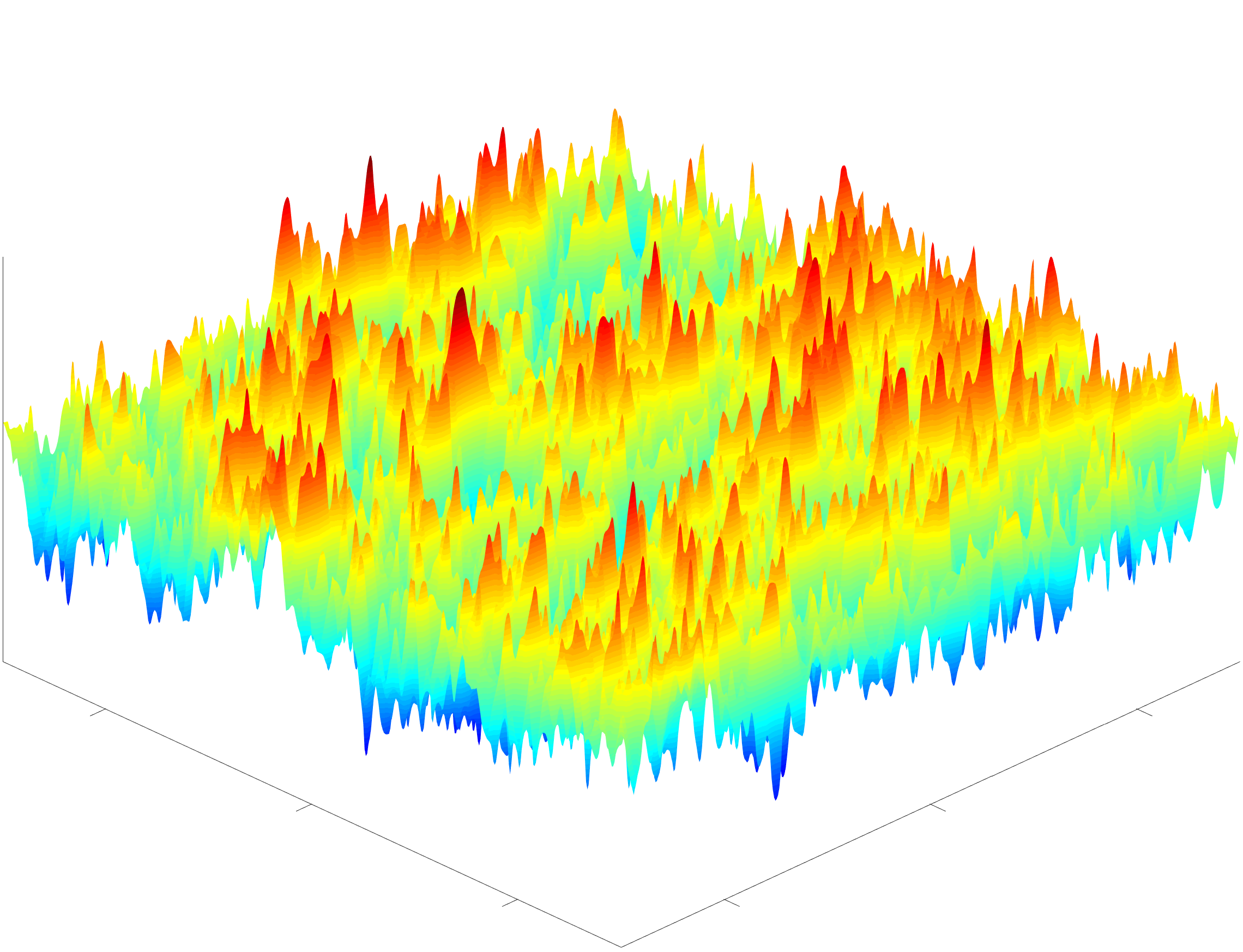} \par\medskip
\includegraphics[width=.8\linewidth]{corr_y_jet.png}
\caption{\small A graph of the first-order correctors, which can be compared to Figure~\ref{fig.GFF}. The coefficient field~$\a$ is a scalar random checkerboard, with checkerboard squares sampled independently from the uniform distribution on~$[1,11]$. In this case, by symmetry we have that~$\overline{\a}$ is also scalar. The top picture is $\phi_{e_1}$ and the bottom picture is of~$\phi_{e_2}$. One can see from the picture that the correlation structures appear to be strongly anisotropic.}
\label{f.correctors}
\end{figure}

In this section, we describe the scaling limit of the corrector. Let $\W(\cdot,p,q)$ be the quadratic form white noise appearing in the statement of Theorem~\ref{t.conv.white}. Inspired by \eqref{e.almost.def.msf.V}, we define, for each $e \in \Rd$, the vector white noise $\msf V(\cdot, e)$ such that, for every $p \in \Rd$, 
\begin{equation*}  
p \cdot \msf V(\cdot, e) = \frac 1 2 \Ll( \W(\cdot,p+e,\ahom p) - \W(\cdot,p-e,\ahom p) \Rr).
\end{equation*}
This definition makes sense since the right side depends linearly on $p \in \Rd$, as follows from the fact that $(p,q) \mapsto \W(\cdot,p,q)$ is a quadratic form. For the same reason, the mapping $e \mapsto \msf V(\cdot,e)$ is linear. Following up on the heuristic argument in \eqref{e.coarsened.3}, we then define the gradient GFF $\nabla \Psi_e$ as solving
\begin{equation}  
\label{e.def.psie}
-\nabla \cdot \ahom \nabla \Psi_e = \nabla \cdot \Ll( \msf V(\cdot, e) \Rr) .
\end{equation}
This is understood as in \eqref{e.eq.gff} and \eqref{e.def.grad.gff}. 
The following theorem, which is the main result of the chapter, shows the correctness of the heuristic argument of Section~\ref{s.heuristics}. See Figure~\ref{f.correctors} for a visualization of the theorem. 

\index{corrector!first-order~$\phi_e$!scaling limit of}

\begin{theorem}[Scaling limit of the corrector]
\label{t.gff}
Let $\gamma > \frac d 2$. For every $e \in \Rd$, we have
\begin{equation}
\label{e.scale.corr}
r^{\frac d 2} \Ll(\nabla \phi_e\Rr)(r \, \cdot) \xrightarrow[r \to \infty]{\mathrm{(law)}} \nabla \Psi_e,
\end{equation}
with respect to the topology of $H^{-\gamma}_{\mathrm{loc}}(\Rd)$. 
Moreover, the convergences in law in \eqref{e.clt.mesoJ}, \eqref{e.clt.J} and \eqref{e.scale.corr} hold jointly over $z,p,q,e \in \Rd$ and $f \in L^1(\Rd) \cap L^\infty(\Rd)$.
\end{theorem}
We recall that the functional convergence in law stated in \eqref{e.scale.corr} means that, for every bounded and continous functional $F : H^{-\gamma}_{\mathrm{loc}}(\Rd) \to \R$, we have
\begin{equation*}  
\E \Ll[ F \Ll( r^{\frac d 2} \Ll(\nabla \phi_e\Rr)(r \, \cdot)  \Rr)  \Rr] \xrightarrow[r \to \infty]{} \E \Ll[ F \Ll( \nabla \Psi_e \Rr)  \Rr] .
\end{equation*}
\begin{exercise}
Deduce from \eqref{e.scale.corr} that, for each compactly supported function $f \in H^{\gamma}(\Rd)$, $\gamma > \frac d 2$, the evaluation of the distribution $r^{\frac d 2} \Ll(\nabla \phi_e\Rr)(r \, \cdot)$ against $f$ converges in law to the evaluation of $\nabla \Psi_e$ against $f$. 
\end{exercise}
\begin{exercise}
Recall that $\zeta$ denotes the standard mollifier. Using Theorems~\ref{t.HKtoSob} and \ref{t.gff}, show that, for each $\gamma > \frac d 2$ and $p \in [2,\infty)$, the convergence in law 
\begin{equation*}
r^{\frac d 2} \Ll(\nabla \phi_e \ast \zeta\Rr)(r \, \cdot) \xrightarrow[r \to \infty]{\mathrm{(law)}} \nabla \Psi_e
\end{equation*}
holds with respect to the topology of $W^{-\gamma,p}_{\mathrm{loc}}(\Rd)$. 
\end{exercise}

We start by restating the central limit theorem of the previous section in terms of the coarsened coefficients $\b_r(z)$. This is immediate from Theorem~\ref{t.conv.white} and the definition of $\b_r(z)$, see Definition~\ref{d.brz}. 
\begin{lemma}
\label{l.b.to.white.noise}
There exists $\eta'(d) \in  \Ll( \frac 1 2, 1 \Rr)$ such that, for every $\eta \in (\eta',1)$, $e \in \Rd$ and $F \in L^1(\Rd;\Rd) \cap L^\infty(\Rd;\Rd)$, we have
\begin{equation}  
\label{e.clt.b}
r^{\frac d 2} \int_{\Rd} F(x) \cdot \b_{r^{\eta}}(rx) e \, dx \xrightarrow[r \to \infty]{\mathrm{(law)}} \int_{\Rd} F(x) \cdot \msf V(x,e) \, dx.
\end{equation}
Moreover, the convergences in law in \eqref{e.clt.mesoJ}, \eqref{e.clt.J} and \eqref{e.clt.b} hold jointly over $e,p,q,z \in \Rd$, $f \in L^1(\Rd) \cap L^\infty(\Rd)$ and $F \in L^1(\Rd;\Rd) \cap L^\infty(\Rd;\Rd)$. 
\end{lemma}

In the next lemma, we give a refinement of Lemma~\ref{l.coarsenedequation} and a rigorous version of~\eqref{e.coarsened.equation} by showing that locally averaged elements of $\A_1$ are approximately~$(\ahom + \b_r)$-harmonic.

\begin{lemma}
\label{l.refined.coarsenedequation}
There exists an exponent $\kappa(s,d) > 0$ and a constant $C(s,d,\Lambda)<\infty$ such that, for every $\psi \in C^1(\Rd)$ satisfying
\begin{equation}  
\label{e.assumption.psi}
\forall x \in \Rd, \quad |x|^{d} |\psi(x)| + |x|^{d+1} \, |\nabla \psi(x)| \le 1,
\end{equation}
and every $e \in B_1$, we have
\begin{equation} 
\label{e.almostahomharm.gff}
\left| \int_{\Rd} \nabla \psi(x)\cdot  (\ahom + \b_r(x))\Ll[e + \nabla (\phi_e \ast \Phi_{r})(x)\Rr] \, dx \right| 
\le \O_s \Ll( C r^{-\frac d 2 - \kappa} \Rr) .
\end{equation}
\end{lemma}
\begin{proof}
The proof follows closely that of Lemma~\ref{l.coarsenedequation}. We define
\begin{equation*}  
\mcl H_r(z) := \sup_{u \in \A_1(\Phi_{z,r})}\Ll|\int_{\Phi_{z,r}} \a(x) \nabla u(x) \, dx -\Ll( \ahom + \b_r(z) \Rr) \int_{\Phi_{z,r}} \nabla u(x) \, dx \Rr|.
\end{equation*}
By the definition of $\A_1(\Phi_{z,r})$, we have $\mcl H_r(x) \le C$. Setting $s' := \frac {s+2}{2}$, we recall from Lemma~\ref{l.brz} that, for every $r \ge 1$ and $z \in \Rd$,
\begin{equation}
\label{e.mclHr.fine.0}
\mcl H_r(z) \le C \wedge \O_{s'/2} \Ll( C r^{-d} \Rr) . 
\end{equation}
By Lemma~\ref{l.change-s}, this implies
\begin{equation}
\label{e.mclHr.fine}
\mcl H_r(z) \le \O_{s} \Ll( C r^{-\frac d 2 - \kappa } \Rr) ,
\end{equation}
for $\kappa(s,d) > 0$ such that $\frac d 2 + \kappa = \frac{ds'}{2s}$. We then compute, for every $u \in \mcl A_1$,
\begin{align*} 
\lefteqn{
\left| \int_{\Rd} \nabla \psi(x) \cdot (\ahom + \b_r(x)) \nabla \left( u\ast \Phi_r \right)(x) \,dx \right|
} \qquad & \\
& = \left| \int_{\Rd} \nabla \psi(x) \cdot (\ahom + \b_r(x)) \Ll(\int_{\Phi_{x,r}} \nabla u \Rr)\,dx \right|
\\ &
\leq \left| \int_{\Rd} \nabla \psi(x) \cdot \left( \int_{\Phi_{x,r}} \a \nabla u \right)\,dx \right| +  \int_{\Rd} \left| \nabla \psi(x) \right| \left\| \nabla u \right\|_{L^2(\Phi_{x,r})} \mathcal{H}_r(x)\,dx
\\ & 
=  \left| \int_{\Rd} \nabla \left( \psi \ast \Phi_{r} \right) (x) \cdot \a(x) \nabla u(x) \,dx \right| +  \int_{\Rd} \left| \nabla \psi(x) \right| \left\| \nabla u \right\|_{L^2(\Phi_{x,r})}  \mathcal{H}_r(x)\,dx.
\end{align*}
The conclusion then follows by specializing to $u : x \mapsto x \cdot e + \phi_e(x)$, provided we can show that
\begin{equation}
\label{e.thats.the.equation}
\int_{\Rd} \nabla \left( \psi \ast \Phi_{r} \right) (x) \cdot \a(x) (e + \nabla\phi_e)(x) \,dx = 0
\end{equation}
and that
\begin{equation}
\label{e.estim.Hr.bis}
\int_{\Rd} \left| \nabla \psi(x) \right| \left\| e + \nabla \phi_e  \right\|_{L^2(\Phi_{x,r})}  \mathcal{H}_r(x)\,dx = \O_s \Ll( C r^{-\frac d 2 - \kappa} \Rr) .
\end{equation}
Let $\chi \in C^\infty_c(\Rd;\R)$ be such that $\chi \equiv 1$ on $B_1$ and $\chi \equiv 0$ outside of $B_2$. The left side of \eqref{e.thats.the.equation} is the limit as $R$ tends to infinity of
\begin{equation}  
\label{e.stare.at.zero}
\int_{\Rd} \nabla \Ll( \chi \Ll( \tfrac \cdot R \Rr) \psi \ast \Phi_r \Rr)(x) \cdot \a(x) (e + \nabla \phi_e)(x) \, dx.
\end{equation}
Indeed, the difference between this term and the left side of \eqref{e.thats.the.equation} can be bounded using \eqref{e.assumption.psi} by a constant times
\begin{equation*}  
R^{-(d+1)} \int_{B_{2R}} |e + \nabla \phi_e|(x) \, dx.
\end{equation*}
and by Theorem~\ref{t.correctors} and Jensen's inequality, this is bounded by
\begin{equation*}  
R^{-(d+1)} \O_s \Ll( C R^d \log^\frac 1 2 R \Rr),
\end{equation*}
which tends to $0$ as $R$ tends to infinity. Moreover, the quantity in \eqref{e.stare.at.zero} is in fact zero for every $R$, since $x \mapsto x \cdot e + \phi_e(x) \in \mcl A_1$. This justifies \eqref{e.thats.the.equation}. 

\smallskip

There remains to show \eqref{e.estim.Hr.bis}. By \eqref{e.assumption.psi} and \eqref{e.mclHr.fine}, we clearly have
\begin{equation*}  
\int_{\Rd} \left| \nabla \psi(x) \right|  \mathcal{H}_r(x)\,dx = \O_s \Ll( C r^{-\frac d 2 - \kappa} \Rr) .
\end{equation*}
We also recall from Theorem~\ref{t.correctors} that
\begin{equation*}  
\|\nabla \phi\|_{L^2(\Phi_{x,r})} \le \O_{s'} \Ll( C r^{-\frac d 2} \Rr) .
\end{equation*}
The conclusion therefore follows from \eqref{e.mclHr.fine.0} and Lemma~\ref{l.change-s}.
\end{proof}
We can now complete the proof of Theorem~\ref{t.gff}.
\begin{proof}[Proof of Theorem~\ref{t.gff}]
We fix $\gamma > \gamma' > \frac d 2$. By Theorem~\ref{t.HKtoSob}, there exists a constant $C(s,\gamma',d,\Lambda) < \infty$ such that, for every $r \ge 1$,
\begin{equation}  
\label{e.bound.Hgamma}
\Ll\| r^{\frac d 2} \Ll(\nabla \phi_e\Rr)(r \, \cdot)  \Rr\|_{H^{-\gamma'}(B_1)}  \le \O_s \Ll( C \Rr) . 
\end{equation}
In particular, by Chebyshev's inequality, we have for every $r, M \ge 1$ that 
\begin{equation*}  
\P \Ll[ \Ll\| r^{\frac d 2} \Ll(\nabla \phi_e\Rr)(r \, \cdot)  \Rr\|_{H^{-\gamma'}(B_1)} \ge  M \Rr] \le \frac C M. 
\end{equation*}
By translation invariance of the law of $\nabla \phi_e$ and the fact that the embedding $H^{-\gamma'}_{\mathrm{loc}}(\Rd) \hookrightarrow H^{-\gamma}_{\mathrm{loc}}(\Rd)$ is compact, this ensures that the family of random distributions $\{r^{\frac d 2} \Ll(\nabla \phi_e\Rr)(r \, \cdot) , r \ge 1\}$ is tight in $H^{-\gamma}_{\mathrm{loc}}(\Rd)$. By Prohorov's theorem (see \cite[Theorem~5.1]{billingsley2}), it thus remains to verify the uniqueness of the limit law. In order to do so, it suffices to show that, for every $F \in C^\infty_c(\Rd;\Rd)$, we have
\begin{equation}
\label{e.conv.F.phi}
r^{\frac d 2} \int_{\Rd} F(x) \cdot (\nabla \phi_e)(r x) \, dx \xrightarrow[r \to \infty]{\mathrm{(law)}} \int_{\Rd} F(x) \cdot (\nabla \Psi_e)(x) \, dx. 
\end{equation}
Since this convergence in law will be deduced deterministically from that of Lemma~\ref{l.b.to.white.noise}, the joint convergence with \eqref{e.clt.mesoJ} and \eqref{e.clt.J} will be clear (see Remark~\ref{r.how.you.doing}).

\smallskip

Let $\kappa(s,d) > 0$ be as given by Lemma~\ref{l.refined.coarsenedequation}, let $\eta'(d) \in (0, 1)$ be as given by Lemma~\ref{l.b.to.white.noise}, and let $\eta(s,d) \in (\eta',1)$ be an exponent to be chosen sufficiently close to $1$ in the course of the proof. We first observe that in order to show \eqref{e.conv.F.phi}, it suffices to show that
\begin{equation}
\label{e.conv.F.phi.ast}
r^{\frac d 2} \int_{\Rd} F(x) \cdot \nabla (\phi_e \ast \Phi_{r^\eta})(r x) \, dx \xrightarrow[r \to \infty]{\mathrm{(law)}} \int_{\Rd} F(x) \cdot (\nabla \Psi_e)(x) \, dx. 
\end{equation}
Indeed, this follows from 
\begin{equation*}  
r^{\frac d 2} \int_{\Rd} F(x) \cdot \nabla (\phi_e \ast \Phi_{r^\eta} - \phi_e)(r x) \, dx 
= r^{\frac d 2} \int_{\Rd} (F \ast \Phi_{r^{\eta - 1}} - F)(x) \cdot (\nabla \phi_e)(r  x) \, dx , 
\end{equation*}
the bound \eqref{e.bound.Hgamma}, and elementary estimates on the convergence of $F \ast \Phi_{r^{\eta - 1}}$ to $F$ in $H^\gamma(\Rd)$ as $r$ tends to infinity.

\smallskip

We first show \eqref{e.conv.F.phi.ast} assuming furthermore that the function $F$ is of zero mean. We introduce the Helmholtz-Hodge decomposition of $F$:
\begin{equation}  
\label{e.helmF}
F = \ahom \nabla h + \g,
\end{equation}
where $h \in L^2(\Rd)$ is the unique solution of
\begin{equation*}  
\nabla \cdot \ahom \nabla h = \nabla \cdot F \quad \text{in } \Rd
\end{equation*}
satisfying $\lim_{|x| \to \infty} |h(x)| = 0$, and $\g := F - \a \nabla h \in \Ls(\Rd)$. The existence of $h$ is a straightforward consequence of the Green representation formula. Using also the fact that $F$ is of zero mean and Lemma~\ref{l.centered}, we infer that
\begin{equation}  
\label{e.h.is.reasonable}
\sup_{x \in \Rd} \Ll( |x|^d |h(x)| + |x|^{d+1} |\nabla h(x)| \Rr) < \infty.
\end{equation}
Recalling \eqref{e.def.grad.gff} or \eqref{e.def.psie}, we see that \eqref{e.conv.F.phi.ast} is equivalent to
\begin{equation}
\label{e.conv.in.law.thisone}
r^{\frac d 2} \int_{\Rd} F(x) \cdot \nabla (\phi_e \ast \Phi_{r^\eta})(r x) \, dx \xrightarrow[r \to \infty]{\mathrm{(law)}} - \int_{\Rd} F(x) \cdot \msf V(x,e) \, dx.
\end{equation}
By Lemma~\ref{l.refined.coarsenedequation} and \eqref{e.h.is.reasonable}, we have
\begin{equation*}  
\int_{\Rd} \nabla h(x) \cdot (\ahom + \b_{r^\eta}(x)))\Ll[e + \nabla (\phi_e \ast \Phi_{r^\eta})(x)\Rr] \, dx 
= \O_s \Ll( C r^{-\eta \Ll(\frac d 2 + \kappa\Rr)} \Rr) .
\end{equation*}
Moreover, an integration by parts gives
\begin{equation*}  
\int_{\Rd} \nabla h(x) \cdot \ahom e \, dx = 0.
\end{equation*}
By Lemma~\ref{l.brz} and Theorem~\ref{t.correctors}, we also have
\begin{equation*}  
\int_{\Rd} \nabla h(x) \cdot\b_{r^\eta}(x)\nabla (\phi_e \ast \Phi_{r^\eta})(x) \, dx \le \O_{s/2} \Ll( C r^{-\eta d} \Rr) ,
\end{equation*} 
and thus
\begin{multline*}  
\Ll|\int_{\Rd} \ahom \nabla h(x) \cdot \nabla (\phi_e \ast \Phi_{r^\eta})(x) \, dx - \int_{\Rd} \nabla h(x) \cdot \b_{r^\eta}(x)  e \, dx \Rr| \\
\le \O_s \Ll(  C r^{-\eta \Ll(\frac d 2 + \kappa\Rr)}\Rr)  + \O_{s/2} \Ll( C r^{-\eta d} \Rr).
\end{multline*}
Provided that we choose $\eta < 1$ sufficiently close to $1$, we therefore deduce from this inequality and Lemma~\ref{l.b.to.white.noise} that
\begin{equation*}  
\int_{\Rd} \ahom \nabla h(x) \cdot \nabla (\phi_e \ast \Phi_{r^\eta})(x) \, dx \xrightarrow[r \to \infty]{\mathrm{(law)}} \int_{\Rd} \nabla h \cdot \msf V(x,e) \, dx.
\end{equation*}
In order to conclude for \eqref{e.conv.in.law.thisone}, it suffices to verify that
\begin{equation*}  
\int_{\Rd} \g \cdot \nabla (\phi_e \ast \Phi_{r^\eta})(x) \, dx = 0.
\end{equation*}
This follows from the fact that $\g \in \Ls(\Rd)$, \eqref{e.h.is.reasonable} and Theorem~\ref{t.correctors}, by arguments very similar to those for \eqref{e.thats.the.equation}. 

\smallskip

We have thus verified that \eqref{e.conv.F.phi} holds for every $F \in C^\infty_c(\Rd;\Rd)$ of zero mean. We now wish to extend the convergence in \eqref{e.conv.F.phi} to a possibly non-zero-mean vector field $F \in C^\infty_c(\Rd;\Rd)$. For this purpose, we introduce the vector field
\begin{equation}  
\label{e.def.Fm}
F_m(x) := F(x) - \Phi_m(x) \int_{\Rd} F.
\end{equation}
While this vector field is of zero mean, it is not compactly supported. However, it has Gaussian tails, and therefore it is elementary to verify that the arguments exposed above also apply to $F_m$ (in particular, the decomposition \eqref{e.helmF} with $F$ replaced by $F_m$, and the estimate \eqref{e.h.is.reasonable}, remain valid). That is, for every $m \ge 1$, we have the convergence 
\begin{equation*}  
r^{\frac d 2} \int_{\Rd} F_m(x) \cdot (\nabla \phi_e)(r x) \, dx \xrightarrow[r \to \infty]{\mathrm{(law)}} \int_{\Rd} F_m(x) \cdot (\nabla \Psi_e)(x) \, dx. 
\end{equation*}
Moreover, since $\|F_m - F\|_{L^2(\Rd)} \to 0$ as $m \to \infty$, we also have
\begin{equation*}  
 \int_{\Rd} F_m(x) \cdot (\nabla \Psi_e)(x) \, dx \xrightarrow[m \to \infty]{\mathrm{(prob)}}  \int_{\Rd} F(x) \cdot (\nabla \Psi_e)(x) \, dx.
\end{equation*}
In order to complete the proof, we thus wish to interchange the limits $m \to \infty$ and $r \to \infty$. By \cite[Theorem~3.2]{billingsley2}, it suffices to verify that, for every $c > 0$,
\begin{equation}  
\label{e.interchange.probas}
\lim_{m \to \infty} \limsup_{r \ge 1} \P \Ll[\Ll|r^{\frac d 2}  \int_{\Rd}(F - F_m) \Ll( x \Rr) \cdot (\nabla \phi_e)(r x) \, dx\Rr| \ge c\Rr] = 0.
\end{equation}
The quantity between absolute values above can be rewritten as
\begin{equation*}  
r^{\frac d 2} \Ll(\int_{\Rd} F \Rr) \cdot \Ll( \int_{\Phi_{mr}} \nabla \phi_e\Rr),
\end{equation*}
and, by Theorem~\ref{t.correctors}, there exists a constant $C(s,d,\Lambda) < \infty$ such that, for every $r \ge 1$ and $m \ge 1$,
\begin{equation*}  
r^{\frac d 2}  \Ll( \int_{\Phi_{mr}} \nabla \phi_e\Rr)  = \O_s \Ll( C m^{-\frac d 2} \Rr) .
\end{equation*}
By Chebyshev's inequality, this implies \eqref{e.interchange.probas}, and thus completes the proof.
\end{proof}
\begin{exercise}
Let $h \in C^\infty_c(\Rd)$ be such that $\supp h \subset B_1$ and $\int_{\Rd} h = 1$. Replace the definition of $F_m$ in \eqref{e.def.Fm} by
\begin{equation*}  
F_m (x) := F(x) - m^{-d} \, h \Ll( m^{-1} x \Rr) \int_{\Rd} F,
\end{equation*}
and show that the rest of the proof can be adapted to yield the same conclusion. As an intermediate step, show using Theorem~\ref{t.HKtoSob} and Lemma~\ref{l.bigO.vs.tail} that, for each $\gamma > \frac d 2$, there exists $C(\gamma,s,d,\Lambda) < \infty$ such that, for every $r \ge 1$ and $m \ge 2$,
\begin{equation*}  
\Ll\| r^{\frac d 2} \Ll(\nabla \phi_e\Rr)(r \, \cdot)  \Rr\|_{H^{-\gamma}(B_m)}  \le \O_s \Ll( C \log^{\frac 1 s} (m)\Rr).
\end{equation*}
\end{exercise}

\section*{Notes and references}

The first non-perturbative proof of convergence of a class of objects to the Gaussian free field was obtained by Naddaf and Spencer in~\cite{NS2}. They studied discrete random fields $(\phi(x))_{x \in \Z^d}$ formally obtained as Gibbs measures with Hamiltonian $\sum_{x \sim y} V(\phi(y) - \phi(x))$, where $V$ is a uniformly convex and $C^{1,1}$ function. Note that the corrector is obtained as the minimizer of a random convex functional, while here the functional is deterministic, but randomness is then added according to the Gibbs principle. One of the main tools in \cite{NS2} is the Helffer-Sj\"ostrand identity (which implies the Efron-Stein or ``spectral gap'' inequality). They used it to map the problem of understanding the large-scale behavior of this model to the homogenization of a differential operator in infinite dimensions. 

\smallskip

The first results showing that the fluctuations of spatial averages of the energy density of the corrector have Gaussian fluctuations were obtained in \cite{N,rossignol,biskup,N2,GN}. The proofs there are based on ``nonlinear'' sensitivity estimates such as the Chatterjee-Stein method of normal approximation \cite{chat1,chat2,nourdin}. The rescaled corrector was shown to converge in law to a Gaussian free field in \cite{MoO,MN} in a discrete setting, using similar tools. It was realized there that the notion of Gaussian free field needed to be extended. Indeed, the definition we adopt in this book is more general than the standard one \cite{Shef}, see Remark~\ref{r.what.is.gff} and \cite{GM2}. The limit law of the fluctuations of the solution $u_\eps$ of
\begin{equation*}  
-\nabla \cdot \a \Ll( \tfrac{\cdot}{\eps} \Rr)  \nabla u_\eps = f \quad \mbox{in } \Rd
\end{equation*}
was identified in \cite{GM,DGO}, still using ``nonlinear'' methods and for discrete equations. The heuristic argument of Section~\ref{s.heuristics} is from \cite{MTalk,GM} and was first made rigorous in~\cite{AKM2}. Sections~\ref{s.CLT} and \ref{s.conv.gff} are based on the approach of \cite{AKM2}. A similar argument was also presented later in~\cite{GO6}.



\chapter{Quantitative two-scale expansions}
\label{c.twoscale}

\index{two-scale expansion|(}

In this chapter, we move beyond estimates on the first-order correctors to consider the homogenization of boundary-value problems. In particular, we consider the convergence, as~$\ep\to0$, of the solution~$u^\ep\in H^1(U)$ of the Dirichlet problem 
\begin{equation}
\label{e.twoscaleDP.pre}
\left\{
\begin{aligned}
& -\nabla \cdot \left(\a^\ep(x) \nabla u^\ep \right) = f & \mbox{in} & \ U, \\
& u^\ep = g & \mbox{on} & \ \partial U,
\end{aligned}
\right.
\end{equation}
where we denote~$\a^\ep:= \a\left(\frac\cdot\ep\right)$, 
to that of the limiting homogenized problem 
\begin{equation}
\label{e.twoscaleDP.prehom}
\left\{
\begin{aligned}
& -\nabla \cdot \left( \ahom \nabla u \right) = f & \mbox{in} & \ U, \\
& u  = g & \mbox{on} & \ \partial U.
\end{aligned}
\right.
\end{equation}
The goal is to not only obtain sharp estimates on the homogenization error $\left\| u^\ep - u \right\|_{L^2(U)}$ and on $\left\| \nabla u^\ep - \nabla u \right\|_{H^{-1}(U)}$, which quantifies the weak convergence of the gradients, but to obtain estimates in the strong~$L^2(U)$ norm of the \emph{corrected gradients}. That is, we study the difference between $u^\ep$ and the \emph{two-scale expansion} of $u$ denoted by
\begin{equation} 
\label{e.tildewepdef}
\tilde{w}^\ep(x):= u(x) + \ep \sum_{k=1}^d \partial_{x_k} u(x)   \phi_{e_k}^\ep (x)
\end{equation}
where $\{ \phi_e \,:\, e\in\partial B_1\}$ is the family of first-order correctors studied in Chapter~\ref{c.A1} and we define  the rescaled corrector~$\phi_e^\ep$, for each $\ep\in \left(0,\tfrac 12\right]$ and $e\in\partial B_1$, by
\begin{equation*} 
\phi^\ep_e(x):= \phi_e \left( \tfrac x\ep \right) - \left(\phi_{e} \ast \Phi_{\ep^{-1}} \right) (0).
\end{equation*}
We have subtracted the constant~$\phi_{e} \ast \Phi_{\ep^{-1}} (0)$ from $\phi_{e}$ since this makes the expression on the right unambiguously defined; recall that~$\phi_e$ is defined a priori only up to a constant. 

\smallskip

We will prove estimates which quantify the limit $\left\| \nabla u^\ep - \nabla \tilde{w}^\ep \right\|_{L^2(U)} \to 0$, justifying the rough expectation that~``$u^\ep$ should look like $u$ plus microscopic-scale wiggles characterized (at first order) by the first-order correctors.'' The main results are stated in Theorems~\ref{t.twoscale} and~\ref{t.L2EE}. We also obtain a result for the Neumann problem in Theorem~\ref{t.twoscale.Neumann}.

\smallskip

While the details of the estimates depend of course on the smoothness of the data (in particular the boundary of the domain~$U$ as well as the regularity of the functions $f$ and $g$), what we roughly show for the Dirichlet problem---neglecting logarithmic corrections appearing in the case~$d=2$ and being vague about stochastic integrability---is that
\begin{equation} 
\label{e.twoscaleresults}
\left\{
\begin{aligned}
& \left\| u^\ep - u \right\|_{L^2(U)} \lesssim O(\ep), \\
& \left\| \nabla u^\ep - \nabla w^\ep \right\|_{L^2(U)} \lesssim O(\ep^{\frac12}). 
\end{aligned}
\right.
\end{equation}
These estimates match the error estimates obtained in periodic homogenization. The larger $O(\ep^{\frac12})$ error in the second estimate is due to \emph{boundary layer} effects. That is, near~$\partial U$, the Dirichlet boundary condition has a stronger effect on the behavior of the solutions than the PDE does, and this unsurprisingly breaks the validity of the two-scale expansion approximation. We also obtain estimates on the size of these boundary layers, showing in particular that in any compactly--contained subdomain $V\subseteq U$, the two-scale expansion error is roughly of order~$\ep$. 

\smallskip

In Section~\ref{s.twoscalenobdry}, we plug the two-scale expansion into the heterogeneous operator and measure the resulting error in $H^{-1}$. It is here that the corrector estimates obtained in Chapter~\ref{c.A1} play a decisive role. Section~\ref{s.twoscale.dirt} is devoted to the proof of the second estimate in~\eqref{e.twoscaleresults}, which is precisely stated in Theorem~\ref{t.twoscale}. Estimates of the boundary layers are the focus of Section~\ref{s.blayers}, where in particular we prove the first estimate of~\eqref{e.twoscaleresults} in Theorem~\ref{t.L2EE}. 
Note that quantitative estimates of the two-scale expansion and of the boundary layer errors are also studied later in Section~\ref{s.twoscale.p}, where even sharper bounds are obtained.

\section{The flux correctors}
\label{s.fluxcorrectors}

We have seen one of the main ideas in this chapter already in Section~\ref{s.homogDPintro}, where we first encountered the two-scale expansion argument. There we showed that quantitative estimates on the sublinearity of the correctors and of their fluxes (in weak norms) imply---by a purely deterministic, analytic argument---estimates on the homogenization error. The goal of this chapter is to refine the arguments in Section~\ref{s.homogDPintro}. In particular, we will be much more careful in how we estimate boundary layer errors (compared to the crude use of H\"older's inequality in~\eqref{e.nablaubndrylayer}) and optimize our assumptions on the regularity of the homogenized solution. 

\smallskip

\index{corrector!flux~$\mathbf{S}_e$|(}

In this section, we prepare for this analysis by introducing the concept of \emph{flux corrector}, which is a vector potential for the flux of the correctors. Roughly speaking, a vector potential for a solenoidal vector field $\g\in \Ls(U)$ is a skew-symmetric matrix-valued function $\mathbf{A}:U\to \R^{d\times d}$ such that, for some constant $\xi\in\Rd$,
\begin{equation*} \label{}
\g = \xi + \nabla \cdot \mathbf{A}. 
\end{equation*}
Here $\mathbf{A}_{ij}$ are the entries of $\mathbf{A}$ and $\nabla \cdot \mathbf{A}$ is the vector field whose $i$th component is 
\begin{equation*} \label{}
\left( \nabla \cdot \mathbf{A} \right)_i := \sum_{j=1}^d \partial_{x_j}\mathbf{A}_{ij}. 
\end{equation*}
There are some situations in which it is useful to refer to the vector potential, just as it is sometimes useful to refer to the scalar potential of a gradient field. As scalar potentials are unique only up to additive constants, vector potentials are not necessarily unique, and compared to scalar potentials there are actually many more degrees of freedom. A canonical choice is to take $\mathbf{A}_{ij}$ to solve the equation
\begin{equation} 
\label{e.VPlaplacian}
-\Delta \mathbf{A}_{ij} = \partial_{x_j}\g_i - \partial_{x_i} \g_j \quad \mbox{in} \ U,
\end{equation}
with suitable boundary conditions on~$\partial U$. 

\smallskip

We want to define a vector potential for the flux of the correctors---in the case that $\g(x)=\a(x)\left( e + \nabla \phi_e(x) \right)-\ahom e$. That is, for each $e\in\Rd$, we seek a stationary random field $\mathbf{S}_e$, taking values in~$\R^{d\times d}$, such that   
\begin{equation*} \label{}
\forall i,j\in\{1,\ldots,d\},\quad \mathbf{S}_{e,ij} = - \mathbf{S}_{e,ji},
\end{equation*}
and
\begin{equation} 
\label{e.yesfluxcorrector}
\a \left( e + \nabla \phi_e \right) - \ahom e 
= 
\nabla \cdot \mathbf{S}_e.
\end{equation}
Once defined, we will call $\mathbf{S}_e$ a \emph{flux corrector} since it plays a similar role for the fluxes of the solutions as the corrector plays for the gradients of solutions. Similar to the usual correctors~$\phi_e$, we will see that $\mathbf{S}_e$ is a priori well-defined only up to additive constants.

\smallskip

As we will see in the next section, writing the flux in divergence form in terms of~$\mathbf{S}_e$ allows us to use the information concerning the weak convergence of the fluxes in a more direct and efficient manner in the two-scale expansion argument. This will lead to two-scale expansion estimates which are optimal in the assumed regularity of the homogenized solution, improving the argument in Section~\ref{s.homogDPintro} by roughly one derivative. 

\smallskip

Proceeding with the definition of $\mathbf{S}_e$, we first let $\g_e$ denote the difference of the heterogeneous and homogenized fluxes of the first-order corrector $\phi_e$:  
\begin{equation*} \label{}
\g_{e}(x) := \a(x)\left( e + \nabla \phi_e(x) \right) - \ahom e. 
\end{equation*}
Motivated by~\eqref{e.VPlaplacian}, we define a family of stationary potential fields 
\begin{equation*} \label{}
\left\{ \nabla \mathbf{S}_{e,ij}  \,:\, e\in\Rd,\, i,j\in\{1,\ldots,d\} \right\}
\end{equation*}
by
\begin{equation} 
\label{e.defnablaSij}
\nabla \mathbf{S}_{e,ij}(x) := \int_0^\infty \int_{\Rd} \left( \g_{e,j}(y) \nabla \partial_{x_i}\Phi(t, y-x) - \g_{e,i}(y) \nabla \partial_{x_j} \Phi(t, y-x) \right)\, dy\,dt,
\end{equation}
where as usual~$\Phi$ is the standard heat kernel. We next check that this definition makes sense and that~\eqref{e.yesfluxcorrector} holds. (Recall  that the space $\mathbb{L}^2_{\mathrm{pot}}$ was introduced in Definition~\ref{d.stat.l2pot}.)


\begin{lemma}
\label{l.fluxcorrector}
For each $e\in\Rd$ and $i,j\in\{1,\ldots,d\}$, we have that~$\nabla \mathbf{S}_{e,ij} \in \mathbb{L}^2_{\mathrm{pot}}$ and there exist $\delta(d,\Lambda)>0$ and $C(d,\Lambda)<\infty$ such that, for every $e\in\partial B_1$, 
\begin{equation} \label{e.fluxcorrectorgradbound}
\left\| \nabla \mathbf{S}_{e,ij} \right\|_{L^2(B_1)} \leq \O_{2+\delta}(C).
\end{equation}
Moreover, for each $s\in(0,2)$, there exists $C(s,d,\Lambda)<\infty$ such that, for each $r\geq 1$ and $x\in\Rd$, 
\begin{equation} 
\label{e.spatavgS}
\left| \left(\nabla \mathbf{S}_{e,ij} \ast \Phi_{r}\right) (x) \right| 
\leq \O_s\left( Cr^{-\frac d2} \right). 
\end{equation}
We also have that, with $\g(x) := \a(x) \left( e+ \nabla \phi_e(x) \right) -\ahom e$, 
\begin{equation} 
\label{e.equationforS}
-\Delta  \mathbf{S}_{e,ij} = \partial_{x_j} \g_i - \partial_{x_i} \g_j \quad \mbox{in} \ \Rd. 
\end{equation}
Finally, the identity~\eqref{e.yesfluxcorrector} holds. 
\end{lemma}
\begin{proof}
Let us first check that the integral in~\eqref{e.defnablaSij} is convergent. We have 
\begin{align*} \label{}
\lefteqn{
\int_{B_1} \left| \nabla \mathbf{S}_{e,ij} (x) \right|^2\,dx
} \qquad & 
\\ &
= \int_{B_1} \left| \int_0^\infty \int_{\Rd} \left( \g_{e,i}(y) \nabla \partial_{x_j}\Phi(t, y-x) - \g_{e,j}(y) \nabla \partial_{x_i} \Phi(t, y-x) \right)\, dy\,dt\right|^2\,dx .
\end{align*}
We aim at proving~\eqref{e.fluxcorrectorgradbound}, and for this, we split the time integral into two pieces and estimate them in the following two steps. 

\smallskip

\emph{Step 1.} Let us first treat the case $t \in (1,\infty)$. According to Theorem~\ref{t.HKtoSob}, for each $x\in\Rd$, $t>0$ and $i ,j \in \{1,\ldots,d\}$,
\begin{equation} 
\label{e.timeslicenablaS}
 \left| \int_{\Rd} \g_{e,i}(y) \nabla \partial_{x_j}\Phi(t, y-x) \, dy \right|
\leq 
 \O_s\left( Ct^{-1-\frac d4} \right) \wedge \O_{2+\delta} \left( Ct^{-1-\frac 14} \right).
\end{equation}
In particular, by the previous inequality and Lemma~\ref{l.sum-O},
\begin{multline*} \label{}
 \int_{B_1} \left| \int_1^\infty \int_{\Rd} \left( \g_{e,i}(y) \nabla \partial_{x_j}\Phi(t, y-x) - \g_{e,j}(y) \nabla \partial_{x_i} \Phi(t, y-x) \right)\, dy\,dt\right|^2\,dx 
 \\
 \leq \O_{1+\delta/2}(C). 
\end{multline*}

\smallskip

\emph{Step 2.} We next analyze the case $t \in (0,1)$, which is more subtle and requires some results from the theory of singular integral operators. Let $i ,j \in \{1,\ldots,d\}$. For each $x \in B_1$, split 
\begin{multline*} 
\int_0^1 \int_{\Rd}  \g_{e,i}(y) \nabla \partial_{x_j}\Phi(t, y-x) \, dy\,dt \\
= \int_0^1 \int_{B_2}  \g_{e,i}(y) \nabla \partial_{x_j}\Phi(t, y-x) \, dy\,dt +  
\int_0^1 \int_{\R^d \setminus B_2}  \g_{e,i}(y) \nabla \partial_{x_j}\Phi(t, y-x) \, dy\,dt .
\end{multline*}
The latter term is easy to control:
\begin{multline*} 
\left| \int_0^1 \int_{\R^d \setminus B_2}  \g_{e,i}(y) \nabla \partial_{x_j}\Phi(t, y-x) \, dy\,dt \right| 
\\ \leq C \sum_{j=1}^\infty  2^{j(d+2)} \exp(-4^j) \fint_{B_{2^{j+1}}\setminus B_{2^j}} | \g_{e,i}(y)| \, dy \leq  \O_{2+\delta}(C). 
\end{multline*}
The first, local term, needs more effort.  Let $G$ be the Green function for the Laplacian operator in~$\Rd$:
 \begin{equation*} 
G(x) = 
\left\{
\begin{aligned}
& -c_2 \log |x| & \mbox{if} & \ d=2,\\
& c_d |x|^{2-d} & \mbox{if} & \ d>2.
\end{aligned}
\right.
\end{equation*}
Then, for $x,y \in B_2$,
\begin{align*} 
\nabla \partial_{x_j}\Phi(t, y-x) & =  - \mathrm{p.v.} \int_{B_3}  \Delta  G(y-z) \nabla \partial_{x_j}\Phi(t, z-x) \, dz \\
 &  = -  \mathrm{\mathrm{p.v.}}  \int_{B_3}  \nabla \partial_{x_j}G(y-z) \Delta \Phi(t, z-x) \, dz  + R(x,z,t)\\ 
 & = - \mathrm{\mathrm{p.v.}} \int_{B_3}  \nabla \partial_{x_j} G(y-z) \partial_t  \Phi(t, z-x) \, dz +   R(x,y,t),
\end{align*}
where the integrals are taken in the sense of principal value, that is $$\mathrm{p.v.} \int_U g(x,y) \, dy = \lim_{\sigma \to 0} \int_{U \setminus B_\sigma(x)} g(x,y) \, dy,$$ and the residual term $R$ comes from boundary terms on $\partial B_3$ when integrating by parts. For small $t$, $R$ is small by means of $t$:
\begin{equation*} 
\left| R(x,y,t) \right| \leq C t^{-2} \exp\left( - \frac 1{4t} \right) . 
\end{equation*}
We then have
\begin{align*} 
\int_0^1 \int_{B_2} \mathbf{g}_{e,i}(x) \nabla \partial_{x_j}\Phi(t, x-z)  \,dx  & =  \int_{B_2} \mathrm{p.v.} \int_{B_3}  \nabla \partial_{x_j} G(y-x) \mathbf{g}_{e,i}(y)  \, dy \, dx    \\ 
\notag &  \quad - \int_{B_2} \int_{B_3}  \nabla \partial_{x_j} G(y-z) \mathbf{g}_{e,i}(x) \Phi(1, z-x) \, dz \, dx \\
\notag & \quad  +  \int_0^1 \int_{B_2}  \g_{e,i}(y)  R(x,y,t) \, dy \, dt .
\end{align*}
For the first part we use the fact that singular integrals as above map $L^2$ to $L^2$. In particular, 
\begin{equation} \label{e.Gfctsingblargh}
 \left\| \mathrm{p.v.}  \int_{B_3}  \nabla \partial_{x_j} G(y - \cdot) \mathbf{g}_{e,i}(y)  \, dy  \right\|_{L^2(B_2)}  \leq  C  \left\| \mathbf{g}_{e,i}  \right\|_{L^2(B_3)} \leq \O_{2+\delta}(C).
\end{equation}
To prove~\eqref{e.Gfctsingblargh}, we apply Plancherel's formula. Indeed, $\left\| \widehat{\partial_{x_i} \partial_{x_j}G} \right\|_{L^\infty(\R^d)} \leq C$, and hence
\begin{multline*} 
 \left\| \mathrm{p.v.}  \int_{B_3}  \nabla \partial_{x_j} G(y - \cdot) \mathbf{g}_{e,i}(y)  \, dy  \right\|_{L^2(B_2)}  
  \\ \leq 
  \left\|  \nabla \partial_{x_j} G \ast \left( \indc_{B_3} \mathbf{g}_{e,i} \right)   \right\|_{L^2(\R^d)} 
  \leq \left\|  \widehat{\nabla \partial_{x_j} G}  \right\|_{L^\infty(\R^d)} \left\| \widehat{ \indc_{B_3} \mathbf{g}_{e,i}}    \right\|_{L^2(\R^d)}  \leq C
  \left\| \mathbf{g}_{e,i}    \right\|_{L^2(B_3)} .
\end{multline*}
The second term can be estimated similarly as
\begin{multline*} 
\left\| \int_{B_2} \int_{B_3} \nabla \partial_{x_j}G(z,y) \mathbf{g}_{e,i}(y) \Phi(1, z-\cdot) \, dz \, dy \right\|_{L^2(B_2)} 
\\ \leq C \left\| \Phi(1, \cdot) \ast \mathbf{g}_{e,i} \right\|_{L^2(B_3)} \leq \O_{2+\delta}(C).
\end{multline*}
The third term can be controlled by the strong bound on $R$: 
\begin{multline*} 
\left\| \int_0^1 \int_{B_2} \g_{e,i}(y) R(\cdot ,y,t) \, dy \, dt \right\|_{L^2(B_1)} 
\\ \leq C \int_{0}^1 t^{-2} \exp\left( - \frac 1{4t} \right) \, dt \left\| \mathbf{g}_{e,i} \right\|_{L^2(B_3)} \leq \O_{2+\delta}(C).
\end{multline*}
Combining, 
\begin{equation*} 
\left\| 
\int_0^1 
\int_{\Rd}  \g_{e,i}(y) \nabla \partial_{x_j}\Phi(t, y-x) 
\, dy\,dt 
 \right\|_{L^2(B_1)} 
\leq \O_{2+\delta}(C),
\end{equation*}
and this together with the first step yields~\eqref{e.fluxcorrectorgradbound}. The fact that $\nabla \mathbf{S}_{e,ij} \in \mathbb{L}^2_{\mathrm{pot}}$ follows from~\eqref{e.fluxcorrectorgradbound} and the stationarity of the vector field~$\g_e $. 

\smallskip

\emph{Step 3.} To prove~\eqref{e.spatavgS}, we use the semigroup property of the heat kernel:
\begin{multline*} \label{}
\left( \nabla \mathbf{S}_{e,ij} \ast \Phi_{r}\right)(x)
\\
= \int_0^\infty \int_{\Rd} \left( \g_{e,j}(y) \nabla \partial_{x_i}\Phi(t+r^2, y-x) - \g_{e,i}(y) \nabla \partial_{x_j} \Phi(t+r^2, y-x) \right)\, dy\,dt.
\end{multline*}
Thus~\eqref{e.spatavgS} follows from~\eqref{e.timeslicenablaS} after integration in~$t$. 

\smallskip

\emph{Step 4.} The identity~\eqref{e.yesfluxcorrector} follows directly from the definition, which gives 
\begin{align*} \label{}
\left( \nabla \cdot \mathbf{S}_e  \right)_i(x)
& 
= \int_0^\infty \int_{\Rd} \left( \g_{e}(y) \cdot \nabla \partial_{x_i}\Phi(t, y-x) - \g_{e,i}(y) \Delta \Phi(t, y-x) \right)\, dy\,dt.
\end{align*}
Indeed, since $\g_e$ is solenoidal and $\partial_{x_j}\Phi(t, \cdot-x)$ has sufficient decay at infinity, the integral of the first term in the integrand vanishes. As for the second term, we use that $\g_e$ has zero mean to deduce that, for almost every $x\in\Rd$, 
\begin{equation*} \label{}
\int_0^\infty \int_{\Rd} \g_{e,i}(y) \Delta \Phi(t, y-x) \, dy\,dt 
= 
\int_0^\infty \partial_t \int_{\Rd} \g_{e,i}(y)  \Phi(t, y-x) \, dy\,dt 
= - \g_{e,i}(x). 
\end{equation*}
This completes the proof. 
\end{proof}

Essentially all of the estimates we have proved for the first-order correctors~$\phi_e$ can also be obtained for the flux correctors~$\mathbf{S}_e$. Instead of integrating the gradients~$\nabla\phi_e$, we instead integrate the (corrected) fluxes $\g_e:=\a \left( e + \nabla \phi_e \right) - \ahom e$. The only estimates we really need for our purposes on~$\mathbf{S}_e $ are the following analogues of~\eqref{e.correctorminbounds0} as well as~\eqref{e.phioscd>2},~\eqref{e.phioscd=2} and~\eqref{e.phid=2}, which are stated in the following proposition.

\begin{proposition}[Flux corrector estimates]
\label{p.fluxcorrectorests}
There exist an exponent~$\beta(d,\Lambda)>0$, a constant~$C(d,\Lambda)<\infty$
and, for each $s\in (0,d)$, a random variable~$\X_s$ satisfying~\eqref{e.X} such that, for every $r\geq \X_s$,
\begin{equation}
\label{e.correctorminbounds0.flux}
\sup_{e\in\partial B_1} \left\| \mathbf{S}_{e} - \left( \mathbf{S}_e \right)_{B_r} \right\|_{\underline{L}^2(B_r)} 
\leq Cr^{1-\beta(d-s)}.
\end{equation}
In dimensions~$d>2$, for every~$e\in\partial B_1$, the flux corrector~$\mathbf{S}_e$ exists as a $\Zd$-stationary $\R^{d\times d}$-valued random field which is identified uniquely by the choice $\E \left[ \fint_{\cu_0} \mathbf{S}_e \right] = 0$; moreover, there exist $\delta(d,\Lambda)>0$ and~$C(d,\Lambda)<\infty$ such that 
\begin{equation} 
\label{e.Soscd>2}
\left\| \mathbf{S}_e\right\|^2_{\underline{L}^2(\cu_0)} 
\leq 
\O_{2+\delta}(C). 
\end{equation}
In dimension~$d=2$, for each $s \in (0,2)$, there exists $C(s,d,\Lambda)<\infty$ such that, for every $e\in\partial B_1$, $1 \leq r < R/\sqrt{2}$ and $x\in\Rd$, 
\begin{equation} 
\label{e.Soscd=2}
\left\| \mathbf{S}_e - \left( \mathbf{S}_e\ast \Phi_r\right)(0) \right\|^2_{\underline{L}^2(B_r)}
\leq \O_{s} \left( C \log^{\frac12} r \right)
\end{equation}
and
\begin{equation}
\label{e.Sd=2}
\left| \left( \mathbf{S}_e\ast \Phi_r \right)(x) - \left( \mathbf{S}_e\ast \Phi_R \right)(y) \right|
\leq
\O_s\left( C\log^{\frac12} \left(2+\frac{R+|x-y|}r\right) \right).
\end{equation}
\end{proposition}
\begin{proof}
The second and third statements follow from~\eqref{e.spatavgS},~\eqref{e.fluxcorrectorgradbound} and Lemma~\ref{l.mspoincare.masks} in exactly the same way as~\eqref{e.phioscd>2} and~\eqref{e.phioscd=2} are proved from~\eqref{e.gradient},~\eqref{e.correctorgradbound} and Lemma~\ref{l.mspoincare.masks}. We omit the details, referring the reader to the argument beginning in the second paragraph of the proof of Theorem~\ref{t.correctors} in Section~\ref{s.correctors}. Similarly, we obtain~\eqref{e.Sd=2} by mirroring the proof of~\eqref{e.phid=2}, in particular analysis in Section~\ref{s.twopoint} in the proof of Proposition~\ref{p.d2vengeance}. 

\smallskip

The bound~\eqref{e.correctorminbounds0.flux} is proved by the same argument as Lemma~\ref{l.fluxcorrector}, except that instead of appealing to Theorem~\ref{t.HKtoSob} we instead use the corrector bounds in Proposition~\ref{p.correctorsbasecase}. 
\end{proof}
\index{corrector!flux~$\mathbf{S}_e$|)}

\section{Quantitative two-scale expansion without boundaries}
\label{s.twoscalenobdry}

Recall that, for each $\ep\in \left(0,\tfrac 12\right]$ and $e\in\partial B_1$, the rescaled corrector~$\phi^\ep_e$ is defined by
\begin{equation} 
\label{e.phiep.def}
\phi^\ep_e(x):= \phi_e \left( \tfrac x\ep \right) - \left(\phi_{e} \ast \Phi_{\ep^{-1}} \right) (0).
\end{equation}
Likewise, we define the rescaled flux correctors by
\begin{equation} 
\label{e.Sep.def}
\mathbf{S}^\ep_{e}:= \mathbf{S}_e \left( \tfrac x\ep \right) - \left( \mathbf{S}_e \ast \Phi_{\ep^{-1}} \right) (0).
\end{equation}
Given~$u\in H^1(\Rd)$ and~$\ep \in\left( 0,\tfrac 12 \right]$, we introduce a function $w^\ep \in H^1(\Rd)$ defined by
\begin{equation} 
\label{e.wepdef}
w^\ep(x) 
:=u(x) + \ep\sum_{k=1}^d \left( \partial_{x_k} u \ast \zeta_\ep \right) (x)  \phi^\ep_{e_k}(x).
\end{equation}
Here as usual we set $\zeta_\ep (x):= \ep^{-d} \zeta(x/\ep)$,  where $\zeta$ is the standard mollifier, see~\eqref{e.standardmollifier}.
We call~$w^\ep$ the \emph{two-scale expansion} of~$u$, and the goal is to show, under appropriate regularity assumptions on~$u$, that it is a good approximation of the solution~$u^\ep$ of the equation
 \begin{equation} 
\label{e.uepRd}
-\nabla \cdot\left( \a^\ep \nabla u^\ep \right) = -\nabla \cdot \left( \ahom \nabla u\right) \quad \mbox{in} \ \Rd. 
\end{equation}
This reflects the expectation that the difference of~$u^\ep$ and~$u$ should be, at leading order, the rescaled corrector with the local slope given by~$\nabla u$. We have already encountered this idea previously in the proof of Theorem~\ref{t.DP.blackbox}. 

\smallskip

Notice that $w^\ep$ is defined slightly differently than the function~$\tilde{w}^\ep$ introduced in~\eqref{e.tildewepdef}, since we have mollified $\nabla u$ on the microscopic scale. This is necessary because otherwise~$\nabla u$ may be too irregular to give a suitable notion of macroscopic (or mesoscopic) slope. Indeed, without regularizing~$\nabla u$, we cannot ensure that~$w^\ep$ is even an~$H^1$ function: if~$u$ is merely an~$H^1$ function, the expression defining~$\tilde{w}^\ep$ may only be slightly better than~$L^1_{\mathrm{loc}}(\Rd)$. On the other hand, if~$u\in W^{2,\infty}(\Rd)$, then this mollification makes essentially no difference, since in this case we have the bound
\begin{equation}
\label{e.weptotildewep}
\left\| w^\ep - \tilde{w}^\ep \right\|_{H^1(B_1)} 
\leq C \left\| u \right\|_{W^{2,\infty} (\Rd)} \cdot
 \left\{ 
\begin{aligned}
& \O_s\left( C\ep\left| \log \ep \right|^{\frac12} \right) & \mbox{if} & \ d=2,\\
& \O_{2+\delta} \left( C \ep \right) & \mbox{if} & \ d>2.
\end{aligned}
\right.
\end{equation}
This follows from~\eqref{e.correctorgradbound},~\eqref{e.phioscd>2}, \eqref{e.phioscd=2} and basic properties of the convolution, which tell us that, for every $p\in [1,\infty]$, 
\begin{equation*} \label{}
\left\| \nabla^2 \left(u\ast \zeta_\ep \right)  \right\|_{L^p(\Rd)} \leq \left\|  \nabla^2 u \right\|_{L^p(\Rd)} 
\end{equation*}
and, for every $\alpha \in (0,1]$, 
\begin{equation} 
\label{e.convolutionpoincare}
\left\| \nabla u - \nabla \left(u\ast \zeta_\ep \right) \right\|_{L^p(\Rd)}
\leq C\ep^\alpha \left\| \nabla u \right\|_{W^{\alpha,p}(\Rd)}.
\end{equation}

\begin{exercise}
Check that~\eqref{e.weptotildewep} is valid for $u\in W^{2,\infty}(\Rd)$.
\end{exercise}

To show that~$w^\ep$ is a good approximation of~$u^\ep$, we plug~$w^\ep$ into the left side of~\eqref{e.uepRd} and estimate how close the result is to the right side, with the error measured in the~$H^{-1}$ norm. This leads to the statement given in the following theorem, which is the main result of this section. It is perhaps the most succinct summary of the principle that good bounds on the correctors can be transferred into good bounds on the homogenization error by a deterministic argument (in fact by a direct computation).

\begin{theorem}[Quantitative two-scale expansion]
\label{t.twoscaleplugveng}
Fix~$\alpha \in (0,1]$, $p\in (2,\infty]$ and $\ep\in \left(0,\tfrac 12\right]$ and define the random variable
\begin{equation} 
\label{e.Xepexpression}
\X_{\ep,p} := \sup_{e\in\partial B_1} \left( \ep^{d} \sum_{z \in  \eps\Z^d \cap B_{1}} 
\left(
\left\|  \mathbf{S}_{e}^\ep  \right\|_{\underline{L}^2(z+2 \ep\cu_0)} + \left\| \phi_{e}^\ep \right\|_{\underline{L}^2(z+2 \ep \cu_0)} \right)^{\frac{2p}{p-2}} 
\right)^{\frac{p-2}{2p}}.
\end{equation}
There exists~$C(\alpha,p,d,\Lambda)<\infty$ such that,
for every~$u\in W^{1+\alpha,p}(\Rd)$, if we define the function~${w}^\ep\in H^1(B_1)$ by~\eqref{e.wepdef}, then we have the estimate
\begin{equation} 
\label{e.twoscaleplugH2}
\left\|  
\nabla \cdot \a^\ep \nabla {w}^\ep - \nabla \cdot \ahom \nabla u 
\right\|_{H^{-1}(B_1)}
\leq
C\ep^\alpha \X_{\ep,p} \|u\|_{W^{1+\alpha,p}(\Rd)}. 
\end{equation}
\end{theorem}

\begin{remark}
\label{r.explicitXep}
For future reference, we record here some estimates on the random variable~$\X_{\ep,p}$ defined in~\eqref{e.Xepexpression} which are consequences of bounds on the correctors and flux correctors we have already proved. 

\smallskip

By Theorem~\ref{t.correctors}, Proposition~\ref{p.fluxcorrectorests} and~Lemma~\ref{l.sum-O}, we have the following estimate: there exist $\delta(d,\Lambda)>0$ and, for every $s\in (0,2)$, a constant $C(s,p,d,\Lambda)<\infty$ such that 
\begin{equation} 
\label{e.Xepp.bound1}
\X_{\ep,p}
\leq 
\left\{ 
\begin{aligned}
& \O_s\left( C \left| \log \ep \right|^{\frac12} \right) & \mbox{if} & \ d=2,\\
& \O_{2+\delta} \left( C \right) & \mbox{if} & \ d>2. 
\end{aligned}
\right.
\end{equation}
This is an optimal bound on the size of~$\X_{\ep,p}$. Combining it with~\eqref{e.twoscaleplugH2}, we obtain, for every $s\in (0,2)$, $\alpha \in (0,1]$ and $p\in (2,\infty]$, a constant $C(s,p,d,\Lambda)<\infty$ such that 
\begin{align} 
\label{e.twoscaleplugH2.explicit}
\left\|  \nabla \cdot \a^\ep \nabla {w}^\ep - \nabla \cdot \ahom \nabla u \right\|_{H^{-1}(B_1)} 
\leq
C 
\left\| u \right\|_{W^{1+\alpha,p} (\Rd)}
\cdot 
\left\{ 
\begin{aligned}
& \O_s\left( \ep^\alpha \left| \log \ep \right|^{\frac12} \right) & \mbox{if} & \ d=2,\\
& \O_{2+\delta} \left( \ep^\alpha \right) & \mbox{if} & \ d>2. 
\end{aligned}
\right. 
\end{align}
%
%
\end{remark}

\begin{remark}
\label{r.extension}
In order to make sense of $w_\ep$ in $B_1$, we need that $u$ be defined over $B_1 + B_\eps$. For simplicity, we assume that $u \in W^{1+\al,p}(\Rd)$ in Theorem~\ref{t.twoscaleplugveng}. Note however that this is not restrictive, since by the Sobolev extension theorem (Proposition~\ref{p.extension}), every $u \in W^{1+\al,p}(B_1)$ admits an extension $\Ext(u) \in W^{1+\al,p}(\Rd)$ which coincides with $u$ in $B_1$ and satisfies
\begin{equation*}  
\|\Ext(u)\|_{W^{1+\al,p}(\Rd)} \le C \|u\|_{W^{1+\al,p}(B_1)}.
\end{equation*}
\end{remark}

We begin the proof of Theorem~\ref{t.twoscaleplugveng} by rewriting the left side of~\eqref{e.twoscaleplugH2} in terms of the correctors and flux correctors.  

\begin{lemma}
\label{l.letswritetheflux}
Fix~$u\in H^1(\Rd)$ and $\ep\in \left(0,\tfrac12\right]$ and let~${w}^\ep\in H^1(B_1)$ be defined by~\eqref{e.wepdef}. Then we have the identity
\begin{equation}
\label{e.letswritetheflux}
\nabla \cdot \left(  \a^\ep \nabla w^\ep 
- \ahom \nabla u \right)
= \nabla \cdot \mathbf{F}^\ep ,
\end{equation}
where the $i$th component of the vector field $\mathbf{F}^\ep$ is given by 
\begin{align} 
\label{e.F_i^ep}
\mathbf{F}_i^\ep(x) 
& :=  
\sum_{j=1}^d  \left(\a_{ij}^\ep(x) - \ahom_{ij} \right)\partial_{x_j} (u - \zeta_\ep \ast u )(x) 
\notag \\ 
& \quad 
+ \ep \sum_{j,k=1}^d \partial_{x_j}\partial_{x_k} \left( \zeta_\ep \ast u\right)(x) \left(  - \mathbf{S}_{e_k,ij}^\ep (x) + \a_{ij}^\ep(x) \phi_{e_k}^\ep (x) \right).
\end{align}
In particular, there is a constant $C(d)<\infty$ such that 
\begin{equation} 
\label{e.twoscalereduce}
\left\| \nabla \cdot  \a^\ep \nabla w^\ep 
- \nabla \cdot \ahom \nabla u  \right\|_{H^{-1}(B_1)} 
\leq 
C\left\| \mathbf{F}^\ep \right\|_{L^2(B_1)} .
\end{equation}
\end{lemma}
\begin{proof}
Denote 
\begin{equation*} 
\hat{w}^\ep(x) := \left( \zeta_\ep \ast u\right)(x) + \ep \sum_{k=1}^d  \partial_{x_k} \left( \zeta_\ep \ast u\right)(x) \phi_{e_k}^\ep \left(x \right).
\end{equation*}
Then
\begin{equation*} 
\nabla \cdot \left( \a^\ep \nabla w^\ep - \ahom \nabla u \right) = \nabla \cdot  (\a^\ep-\ahom) \nabla \left(u - \zeta_\ep \ast u \right) 
+ \nabla \cdot  \left( \a^\ep \nabla \hat{w}^\ep - \ahom \nabla (\zeta_\ep \ast u) \right) ,
\end{equation*}
and in view of the definition of $\mathbf{F}$, we only need to prove the formula for the last term, that is 
\begin{equation}  \label{e.twoscalegoal1}
\nabla \cdot  \left( \a^\ep \nabla \hat{w}^\ep - \ahom \nabla (\zeta_\ep \ast u) \right) 
= 
\sum_{i,j,k=1}^d \partial_{x_i} \left( \ep  \partial_{x_j}\partial_{x_k} \left( \zeta_\ep \ast u\right) \left(  - \mathbf{S}_{e_k,ij}^\ep  + \a_{ij}^\ep \phi_{e_k}^\ep  \right) \right).
\end{equation}
For this, we have that 
\begin{equation*} 
\partial_{x_j}\hat{w}^\ep(x) = \sum_{k=1}^d \left(\delta_{jk}+  \partial_{x_j}\phi_{e_k}^\ep(x) \right)  \partial_{x_k} \left( \zeta_\ep \ast u\right)(x)
+ \ep \sum_{k=1}^d \phi_{e_k}^\ep(x) \partial_{x_j}\partial_{x_k}\left( \zeta_\ep \ast u\right)(x) .
\end{equation*}
Thus 
\begin{align*}
\sum_{j=1}^d  \a_{ij}^\ep\left(x\right) \partial_{x_j}\hat{w}^\ep(x) 
&
= 
\sum_{j,k=1}^d \a_{ij}^\ep\left(x\right)  \left(\delta_{jk}+  \partial_{x_j}\phi_{e_k}^\ep(x)  \right)   \partial_{x_k} \left( \zeta_\ep \ast u\right)(x) 
\\ & \quad 
+ \ep \sum_{j,k=1}^d \a_{ij}^\ep\left(x\right)  \phi_{e_k} \left(\tfrac{x}{\ep}\right)\partial_{x_j}\partial_{x_k}\left( \zeta_\ep \ast u\right)(x)  .
\end{align*}
By the definition of the flux corrector, we have, for every $i,k \in \{1,\ldots,d\}$,
\begin{equation*}
\sum_{j=1}^d \a_{ij}^\ep\left(x\right)  \left(\delta_{jk}+  \partial_{x_j}\phi_{e_k}^\ep (x) \right)  - \sum_{j=1}^d \ahom_{ij} \delta_{jk}= 
\sum_{j=1}^d  \partial_{x_j}\mathbf{S}_{e_k,ij}^\ep\left( x \right) ,
\end{equation*}
and therefore
\begin{align*}
 \nabla \cdot \left( \a^\ep  \nabla \hat{w}^\ep - \ahom \nabla \left( \zeta_\ep \ast u\right) \right)
& =  \sum_{i,j,k=1}^d \partial_{x_i} \left( \partial_{x_j}\mathbf{S}_{e_k,ij}^\ep    \partial_{x_k} \left( \zeta_\ep \ast u\right)  \right) 
\\ & \quad 
+ \ep \sum_{i,j,k=1}^d \partial_{x_i}\left( \a_{ij}^\ep  \phi_{e_k}^\ep \partial_{x_j}\partial_{x_k}\left( \zeta_\ep \ast u\right) \right).
\end{align*}
By the skew-symmetry of $\mathbf{S}_e$, we have
\begin{align*} 
\lefteqn{\sum_{i,j,k=1}^d \partial_{x_i} \left( \partial_{x_j}\mathbf{S}_{e_k,ij}^\ep    \partial_{x_k} \left( \zeta_\ep \ast u\right)  \right) } \quad &   \\ 
\notag & = \sum_{i,j,k=1}^d \left( \partial_{x_i}\partial_{x_j} \left( \ep \mathbf{S}_{e_k,ij}^\ep    \partial_{x_k} \left( \zeta_\ep \ast u\right) \right) - \partial_{x_i} \left(  \ep \mathbf{S}_{e_k,ij}^\ep  \partial_{x_j}\partial_{x_k}\left( \zeta_\ep \ast u\right)  \right) \right)  \\
\notag & = - \sum_{i,j,k=1}^d  \partial_{x_i} \left(  \ep \mathbf{S}_{e_k,ij}^\ep   \partial_{x_j}\partial_{x_k}\left( \zeta_\ep \ast u\right)   \right) ,
\end{align*}
since, for all $i,j,k \in \{1,\ldots,d\}$, 
\begin{equation*} 
\partial_{x_i}\partial_{x_j} \left( \ep \mathbf{S}_{e_k,ij}^\ep  \partial_{x_k} \left( \zeta_\ep \ast u\right)  \right) = -  \partial_{x_i}\partial_{x_j} \left( \ep \mathbf{S}_{e_k,ij}^\ep    \partial_{x_k} \left( \zeta_\ep \ast u\right)  \right) .
\end{equation*}
Combining the previous displays gives us~\eqref{e.twoscalegoal1}, finishing the proof.
\end{proof}

The previous lemma motivates us to estimate~$\left\| \mathbf{F}^\ep \right\|_{L^2(\Rd)}$. Inspecting the definition~\eqref{e.F_i^ep} of $\mathbf{F}^\ep$, we see that this should depend on the assumed regularity of~$u$ and the sizes of~$\mathbf{S}_{e}^\ep$ and~$\phi_e^\ep$, of which we already possess good estimates. The next two lemmas provide the technical labor involved in transferring the regularity assumptions on~$u$ and the bounds on the correctors into an estimate on~$\left\| \mathbf{F}^\ep \right\|_{L^2(\Rd)}$. Since we perform similar computations in the next section and require some flexibility, and because it makes it easier to follow the proof, the statements of these lemmas are somewhat general. 

\begin{lemma} 
\label{l.convolution2}
Fix $1 \leq q\leq p<\infty$ and $0 < \al \le \be \le 1$. There exists $C(\al,\be,p,U,d)<\infty$ such that for every  $g \in W^{\alpha,p}(U + 2\ep \cu_0)$, 
\begin{equation}  \label{e.convolution2}
\left\| g - (g \ast \zeta_\ep) \right\|_{L^q(U)} \leq C|U|^{\frac 1q - \frac1p} \ep^\alpha \| g \|_{W^{\alpha,p}(U+2\ep \cu_0)},
\end{equation}
as well as
\begin{equation}  \label{e.convolution3}
\left\| g \ast \zeta_\ep  \right\|_{W^{\beta,p}(U)} \leq C \ep^{\alpha-\beta} \left\| g \right\|_{W^{\alpha,p}(U+2\ep \cu_0)}.
\end{equation}
\end{lemma}

\begin{proof}
\emph{Step 1.} We first prove~\eqref{e.convolution2}. We only prove it in the case $\alpha \in (0,1)$, the case $\alpha =1$ being similar. 
We first localize as
\begin{equation*} 
\left\| g - (g \ast \zeta_\ep) \right\|_{L^p(U)}^p \leq \sum_{z \in \ep \Z^d \cap U} \left\| g - (g \ast \zeta_\ep) \right\|_{L^p(z + \ep \cu_0)}^p.
\end{equation*}
Since $(a \ast \zeta_\ep) = a$ for all $a \in \R$,  we observe that, for $z \in \ep \Z^d \cap U$, 
\begin{equation*} 
\left\| g - (g \ast \zeta_\ep) \right\|_{L^p(z + \ep \cu_0)} \leq 2 \left\| g - (g)_{z+2\ep \cu_0} \right\|_{L^p(z + 2\ep \cu_0)} .
\end{equation*}
The fractional Poincar\'e inequality (Proposition~\ref{p.sobolevpoincare}) implies
\begin{equation*} 
\left\| g - (g)_{z+2\ep \cu_0} \right\|_{L^p(z + 2\ep \cu_0)}^p  \leq C \ep^{\alpha p} \left[ g \right]_{W^{\alpha,p}(z+2\ep \cu_0)}^p.
\end{equation*}
Moreover, note that
\begin{equation*} 
\left[ g \right]_{W^{\alpha,p}(z+2\ep \cu_0)}^p \leq (1-\alpha) \int_{z+2\ep \cu_0} \int_{U + 2\ep \cu_0} \frac{|g(x) - g(y)|^p}{|x-y|^{d+\alpha p}} \, dx \, dy,
\end{equation*}
and therefore
\begin{equation*} 
\sum_{z \in \ep \Z^d \cap U} \left[ g \right]_{W^{\alpha,p}(z+2\ep \cu_0)}^p \leq C \left[ g \right]_{W^{\alpha,p}(U+2\ep \cu_0)}^p.
\end{equation*}
By H\"older's inequality, we deduce~\eqref{e.convolution2}.

\smallskip

\emph{Step 2.} We now prove \eqref{e.convolution3} assuming $\be < 1$ (the case $\be = 1$ is similar). 
For $x,y \in U$ such that $|x-y| \leq \ep$, we have by the fractional Poincar\'e inequality (Proposition~\ref{p.sobolevpoincare}) that 
\begin{align*} 
\left|(g \ast \zeta_\ep) (y) - (g \ast \zeta_\ep) (x)\right|^p & \leq |y-x|^p \sup_{z \in B_{\ep}(x)} \left| \nabla (g \ast \zeta_\ep)(z)\right|^p
\\ &\leq C \left( \frac{|y-x|}{\ep} \right)^p \fint_{x+ 2\ep \cu_0} \left| g(z) - (g)_{x+2\ep \cu_0} \right|^p \, dz
\\ & \leq C \ep^{p(\alpha-1)-d} |y-x|^p \left[ g \right]_{W^{\alpha,p}(x+2\ep \cu_0)}^p .
\end{align*}
Integration thus gives that 
\begin{equation*} 
\int_{B_\ep(x)} \frac{\left|(g \ast \zeta_\ep) (y) - (g \ast \zeta_\ep) (x)\right|^p}{|x-y|^{d+\beta p}} \, dy  \leq \frac{C}{1-\beta} \ep^{-d + p(\alpha-\beta)}\left[ g \right]_{W^{\alpha,p}(x+2\ep \cu_0)}^p ,
\end{equation*}
and hence
\begin{equation*} 
(1-\beta) \int_U \int_{B_\ep(x)} \frac{\left|(g \ast \zeta_\ep) (y) - (g \ast \zeta_\ep) (x)\right|^p}{|x-y|^{d+\beta p}} \, dy \, dx \leq C \ep^{p(\alpha-\beta)} \left[ g \right]_{W^{\alpha,p}(U+2\ep \cu_0)}^p .
\end{equation*}
The rest of the fractional seminorm can be easily estimated with the aid of the triangle inequality as
\begin{align*} 
\lefteqn{ \int_U \int_{U \setminus B_\ep(x)} \frac{\left|(g \ast \zeta_\ep) (y) - (g \ast \zeta_\ep) (x)\right|^p}{|x-y|^{d+\beta p}} \, dy \, dx }  \qquad & \\ 
& \leq C \int_U \int_{U} \frac{\left|(g \ast \zeta_\ep) (y) - (g \ast \zeta_\ep) (x)\right|^p}{(\ep + |x-y|)^{d+\beta p}} \, dy \, dx \\
& \leq C\ep^{p(\alpha-\beta)}  \left[ g \right]_{W^{\alpha,p}(U)} + C \int_U \int_{U} \frac{\left|(g \ast \zeta_\ep) (x) - g(x)\right|^p}{(\ep + |x-y|)^{d+\beta p}} \, dy \, dx  .
\end{align*}
After first integrating in $y$, Step 1 yields
\begin{multline*} 
 \int_U \int_{U} \frac{\left|(g \ast \zeta_\ep) (x) - g(x)\right|^p}{(\ep + |x-y|)^{d+\beta p}} \, dy \, dx  
 \\ \leq C \ep^{-\beta p} \int_{U} \left|(g \ast \zeta_\ep) (x) - g(x)\right|^p \, dx \leq C \ep^{p(\alpha-\beta)} \left[ g \right]_{W^{\alpha,p}(U+2\ep \cu_0)}^p .
\end{multline*}
Combining the above estimates yields the result. 
\end{proof}
\begin{lemma} 
\label{l.convolution0}
Fix $\al \in (0,1]$ and $p \in (2,\infty)$. There exists $C(\al,p,U,d) < \infty$ such that for every $f \in L^2(U + 2\ep \cu_0)$ and $g \in L^p(U + 2\ep \cu_0)$,
\begin{equation} 
\label{e.convolution0}
\left\| f \left| g \ast \zeta_\ep\right| \right\|_{L^2(U)} 
\leq 
C \left( \ep^{d} \sum_{z \in \eps \Z^d \cap U} 
\left\| f \right\|_{\underline{L}^2(z+2\ep \cu_0)}^{\frac{2p}{p-2}} \right)^{\frac{p-2}{2p}} 
\left\| g \right\|_{L^p(U+2\ep \cu_0)},
\end{equation}
and moreover, for every $f \in L^2(U+2\ep \cu_0)$ and $g \in W^{\alpha,p}(U + 2\ep \cu_0)$, 
\begin{equation} 
\label{e.convolution1}
\left\| f \left| \nabla (g \ast \zeta_\ep) \right| \right\|_{L^2(U)} 
\\ 
\leq 
C \ep^{\alpha-1} \left( \ep^{d} \sum_{z \in \eps \Z^d \cap U} \left\| f \right\|_{\underline{L}^2(z+2\ep \cu_0)}^{\frac{2p}{p-2}} \right)^{\frac{p-2}{2p}} \left\| g \right\|_{W^{\alpha,p}(U+2\ep \cu_0)}.
\end{equation}
The last estimate also holds for $\alpha = 0$ with $\left\| g \right\|_{W^{\alpha,p}(U+2\ep \cu_0)}$ replaced by $\left\| g \right\|_{L^p(U+2\ep \cu_0)}$.
\end{lemma}

\begin{proof}
We will only prove~\eqref{e.convolution1} for $\al \in (0,1)$. The proof of the case $\al  = 1$ and of~\eqref{e.convolution0} are analogous, but simpler. The basic observation is that, for each $z \in \eps \Z^d \cap U$,
\begin{equation*} 
\|g \ast \nabla \zeta_\ep \|_{L^\infty(z+\ep \cu_0)} \leq \|\nabla  \zeta_\ep \|_{L^\infty(\R^d)} \|g  \|_{L^1(z+2\ep \cu_0)} \leq  C \ep^{-1} \ep^{-\frac d 2}  \|g  \|_{L^2(z+2\ep \cu_0)} .
\end{equation*}
Since $\nabla ((g)_{z + 2\ep \cu_0} \ast \zeta_\ep) = 0$, we thus have that
\begin{align*} 
\left\| f \nabla (g \ast \zeta_\ep) \right\|_{L^2(z+\ep \cu_0)}^2 
&  \leq C \ep^{-2} \left\| f \right\|_{\underline{L}^2(z+2\ep \cu_0)}^2 \left\| g - (g)_{z + 2\ep \cu_0} \right\|_{L^2(z+2\ep \cu_0)}^2.
\end{align*}
It follows that
\begin{align*} 
\left\| f \nabla (g \ast \zeta_\ep) \right\|_{L^2(U)}^2 
& \leq  \sum_{z \in \eps \Z^d \cap U } \left\| f \nabla (g \ast \zeta_\ep) \right\|_{L^2(z+\ep \cu_0)}^2 
\\ 
& \leq C  \ep^{-2} \sum_{z \in \eps \Z^d \cap U } \left\| f \right\|_{\underline{L}^2(z+2\ep \cu_0)}^2 \left\| g - (g)_{z + 2\ep \cu_0} \right\|_{L^2(z+2\ep \cu_0)}^2.
\end{align*}
The sum on the right side can be estimated using H\"older's inequality:
\begin{align*} 
\lefteqn{ \sum_{z \in \eps \Z^d \cap U }\left\| f \right\|_{\underline{L}^2(z+2\ep \cu_0)}^2  \left\| g - (g)_{z + 2\ep \cu_0}  \right\|_{L^2(z+2\ep \cu_0)}^2  } \qquad & \\
& \leq \left( \sum_{z \in \eps \Z^d \cap U }\left\| f \right\|_{\underline{L}^2(z+2\ep \cu_0)}^{\frac{2p}{p-2}} \right)^{\frac{p-2}{p}} \left( \sum_{z \in \eps \Z^d \cap U}  \left\| g - (g)_{z + 2\ep \cu_0}  \right\|_{L^2(z+2\ep \cu_0)}^p  \right)^{\frac 2p} \\
& \leq \left( \sum_{z \eps \Z^d \cap U }\left\| f \right\|_{\underline{L}^2(z+2\ep \cu_0)}^{\frac{2p}{p-2}} \right)^{\frac{p-2}{p}} \left( \sum_{z \in \eps \Z^d \cap U }  |2\ep \cu_0|^{\frac{p-2}{2}} \left\| g - (g)_{z + 2\ep \cu_0}  \right\|_{L^p(z+2\ep \cu_0)}^p \right)^{\frac 2p} \\
& \leq 4^{d} \left( \ep^d \sum_{z \in \eps \Z^d \cap U }\left\| f \right\|_{\underline{L}^2(z+2\ep \cu_0)}^{\frac{2p}{p-2}} \right)^{\frac{p-2}{p}} \left( \sum_{z \in \eps \Z^d \cap U } \left\| g - (g)_{z + 2\ep \cu_0}  \right\|_{L^p(z+2\ep \cu_0)}^p \right)^{\frac 2p} .
\end{align*}
Finally, the fractional Poincar\'e inequality (Proposition~\ref{p.sobolevpoincare}) gives us
\begin{equation*} 
\left\| g - (g)_{z + 2\ep \cu_0}  \right\|_{L^p(z+2\ep \cu_0)}^p \leq C \ep^{p \alpha} \left[ g \right]_{W^{\alpha,p}(z+2\ep \cu_0)}^p,
\end{equation*}
and hence we obtain
\begin{equation*} 
\left( \sum_{z \in \eps \Z^d \cap U } \left\| g - (g)_{z + 2\ep \cu_0}  \right\|_{L^p(z+2\ep \cu_0)}^p \right)^{\frac 2p}  \leq \ep^{2 \alpha} \left[ g \right]_{W^{\alpha,p}(U+2\ep \cu_0)}^2.
\end{equation*}
Combining the previous displays yields the lemma. 
\end{proof}

We now complete the proof of Theorem~\ref{t.twoscaleplugveng}.

\begin{proof}[{Proof of Theorem~\ref{t.twoscaleplugveng}}]
By~\eqref{e.twoscalereduce}, it suffices to show that 
\begin{equation} 
\label{e.twoscaleplug.wts}
\left\| \mathbf{F}^\ep \right\|_{L^2(B_1)} 
\leq 
C \ep^\alpha \left\| u \right\|_{W^{1+\alpha,p}(\Rd)} \X_{\ep,p},
\end{equation}
where~$\X_{\ep,p}$ is defined in~\eqref{e.Xepexpression}. By the formula~\eqref{e.F_i^ep} for $ \mathbf{F}^\ep$, we have
\begin{align*}
\left\| \mathbf{F}_i^\ep \right\|_{L^2(B_1)}
&
\leq 
\ep \sum_{k=1}^d
\left\| \, \left| \nabla\nabla \left( \zeta_\ep \ast u\right)\right| \left( - \mathbf{S}_{e_k}^\ep + \a^\ep \phi_{e_k}^\ep \right) \right\|_{L^2(\Rd)}
+ C \left\| \nabla (\zeta_\ep \ast u - u) \right\|_{L^2(\Rd)}.
\end{align*}
According to Lemma~\ref{l.convolution0}, namely~\eqref{e.convolution1}, and the assumption that $u\in W^{1+\alpha,p}(\Rd)$, we have that, for every $e\in \partial B_1$,
\begin{align*}
\lefteqn{
\left\| \, \left| \nabla\nabla \left( \zeta_\ep \ast u\right)\right| \left(  \mathbf{S}_{e}^\ep - \a^\ep \phi_{e}^\ep \right) \right\|_{L^2(B_1)}
} \quad &
\\ & 
\leq 
C \ep^{\alpha-1} 
\left( \ep^{d} \sum_{z \in \eps \Z^d \cap B_1} 
\left\|  \mathbf{S}_{e}^\ep - \a^\ep \phi_{e}^\ep  \right\|_{\underline{L}^2(z+2\ep \cu_0)}^{\frac{2p}{p-2}} 
\right)^{\frac{p-2}{2p}} 
\left\| u \right\|_{W^{1+\alpha,p}(\Rd)}
\\ & 
\leq 
C \ep^{\alpha-1} \X_{\ep,p}
\left\| u \right\|_{W^{1+\alpha,p}(\Rd)}.
\end{align*}
By Lemma~\ref{l.convolution2}, namely~\eqref{e.convolution2}, we have 
\begin{align*}
\left\| \nabla (\zeta_\ep \ast u - u) \right\|_{L^p(B_1)} 
\leq 
C \ep^\alpha \| \nabla u \|_{W^{\alpha,p}(\Rd)} 
\le 
C \ep^\alpha \| u \|_{W^{1+\alpha,p}(\Rd)}.
\end{align*}
Combining the above, we obtain~\eqref{e.twoscaleplug.wts}.
\end{proof}

\begin{exercise}
\label{ex.twoscaleH2}
Prove the following statement, which is a variant of Theorem~\ref{t.twoscaleplugveng} for $p=2$. Fix $\alpha\in (0,1]$, $s \in (0,2)$, $\ep\in \left(0,\tfrac12\right]$, $u\in H^{1+\alpha}(\Rd)$, and let~${w}^\ep \in H^1(B_1)$ be defined by~\eqref{e.wepdef}. 
There exist $C(s,\al,d,\Lambda)<\infty$  and $\delta (d,\Lambda)>0$ such that  
\begin{multline} 
\label{e.twoscaleplugH222}
\left\|  \nabla \cdot \left( \a\left( \tfrac \cdot\ep \right)\nabla {w}^\ep \right)- \nabla \cdot \ahom \nabla u \right\|_{H^{-1}(B_1)}
\\
\leq
\left\{ 
\begin{aligned}
& \O_s\left( C \left\| u \right\|_{H^{1+\alpha} (\Rd)}\ep^\alpha\left| \log \ep \right|^{\frac12} \right) & \mbox{if} & \ d=2,\\
& \O_{2+\delta} \left( C \left\| u \right\|_{H^{1+\alpha} (\Rd)} \ep^\alpha \right) & \mbox{if} & \ d>2.
\end{aligned}
\right.
\end{multline}
\end{exercise}

An important advantage to Theorem~\ref{t.twoscaleplugveng} compared to the result of Exercise~\ref{ex.twoscaleH2}, is that the former gives an estimate which, for every $p>2$, is \emph{uniform} for functions $u\in W^{1+\alpha,p}(\Rd)$ with $\left\| u \right\|_{W^{1+\alpha,p}(\Rd)}\leq 1$. In other words, the random variable defined in~\eqref{e.Xepexpression}, which governs the right side of~\eqref{e.twoscaleplugH2}, clearly does \emph{not} depend on~$u$. In contrast, the random variable on the left side of~\eqref{e.twoscaleplugH222} does depend on~$u$ in the sense that it is \emph{not} bounded uniformly for $\left\| u \right\|_{H^2(\Rd)}\leq 1$. We emphasize this by writing the factor $\left\| u \right\|_{H^2(\Rd)}$ inside the~$\O$ on the right of~\eqref{e.twoscaleplugH222}, rather than outside the~$\O$ as in~\eqref{e.twoscaleplugH2.explicit}.

\section{Two-scale expansions for the Dirichlet problem}
\label{s.twoscale.dirt}

In this section, we present an estimate for the $H^1$ difference between a solution~$u^\ep$ of the Dirichlet problem and the two-scale expansion $w^\ep$ of the solution of the corresponding homogenized Dirichlet problem. This result is an application of Theorem~\ref{t.twoscaleplugveng} and a boundary layer estimate. Under sufficient regularity, the main result~\eqref{e.twoscale} gives a bound of~$\O_2(\ep^{\frac12})$ for this difference (with a logarithmic correction in $d=2$). This estimate is optimal in terms of the exponent of $\ep$, even for smooth data, and agrees with the scaling in the case of periodic homogenization. 

\begin{theorem}[Quantitative two-scale expansion for the Dirichlet problem]
\label{t.twoscale}
\emph{}\\
Fix $s\in (0,2)$, $\alpha \in (0,\infty)$, $p\in (2,\infty]$, a bounded Lipschitz domain $U\subseteq\Rd$ and $\ep\in \left(0,\tfrac 12\right]$. 
There exist $\delta (d,\Lambda)>0$, $C(s,\alpha,p,U,d,\Lambda)<\infty$ and a random variable $\X_\ep$ satisfying
\begin{equation} 
\label{e.XepSI}
\X_\ep \leq 
\left\{ 
\begin{aligned}
& \O_s\left( C \ep^\alpha \left| \log \ep \right|^{\frac12} \right) & \mbox{if} & \ d=2, \ \alpha \in \left( 0,\tfrac12 \right], \\
& \O_s\left( C \ep^{\frac12} \left| \log \ep \right|^{\frac14} \right) & \mbox{if} & \ d=2, \ \alpha \in \left( \tfrac12,\infty \right), \\
& \O_{2+\delta} \left( C  \ep^{\alpha\wedge \frac12} \right) & \mbox{if} & \ d>2, \ \alpha \in (0,\infty),
\end{aligned}
\right.
\end{equation}
such that the following holds: for every $u \in W^{1+\alpha,p}(\Rd)$, if we let $w^\ep\in H^1(U)$ be defined by~\eqref{e.wepdef} and $u^\ep\in H^1(U)$ be the solution of the Dirichlet problem
\begin{equation}
\label{e.twoscaleDP}
\left\{
\begin{aligned}
& -\nabla \cdot \left(\a^\ep(x) \nabla u^\ep \right) = -\nabla \cdot \left( \ahom \nabla u \right) & \mbox{in} & \ U, \\
& u^\ep = u & \mbox{on} & \ \partial U,
\end{aligned}
\right.
\end{equation}
then we have 
\begin{equation} 
\label{e.twoscale}
\left\|  u^\ep - w^\ep  \right\|_{H^{1}(U)}
\leq
\X_\ep \left\| u \right\|_{W^{1+\alpha,p} (\Rd)}.
\end{equation}
\end{theorem}

\begin{remark}
As already pointed out in Remark~\ref{r.extension}, given a function $u \in W^{1+\al,p}(U)$, we can use Proposition~\ref{p.extension} to extend it to a function $u \in W^{1+\al,p}(\Rd)$ with comparable norm, then make sense of $w^{\ep} \in H^1(U)$ and apply Theorem~\ref{t.twoscale}. Therefore there is no restriction in assuming~$u \in W^{1+\alpha,p}(\Rd)$ in Theorem~\ref{t.twoscale}.  
\end{remark}

There are two contributions to the error (the size of the random variables $\X_\ep$) in~\eqref{e.twoscale}: the fact that $w^\ep$ is not an exact solution of the equation (in other words, the error inherited from Theorem~\ref{t.twoscaleplugveng}) and the fact that~$w^\ep$ does not have the same boundary condition as~$u^\ep$ because of the perturbation caused by adding the correctors. We consider these to be separate issues, and remark that it is actually the latter that is responsible for the scaling of the error at leading order. 

\smallskip

To correct the boundary condition, it is natural to introduce~$v^\ep\in H^1(U)$ to be the solution of the Dirichlet problem
\begin{equation} 
\label{e.vepeq}
\left\{
\begin{aligned}
& -\nabla \cdot \left(\a^\ep(x) \nabla v^\ep \right) = 0 & \mbox{in} & \ U, \\
& v^\ep = \ep\sum_{k=1}^d  \phi_{e_k}^\ep  \partial_{x_k} \left( u \ast \zeta_\ep \right) & \mbox{on} & \ \partial U. 
\end{aligned}
\right.
\end{equation}
Note that $w^\ep - u^\ep - v^\ep \in H^1_0(U)$. We are interested therefore in estimating~$\left\| v^\ep \right\|_{H^1(U)}$. This is done by introducing an explicit test function $T_\ep\in v^\ep+H^1_0(U)$ and computing the energy of $T_\ep$. This provides an upper bound for the energy of $v^\ep$ (and hence of~$\left\| v^\ep \right\|_{H^1(U)}$) by the variational principle.

\smallskip

Since the boundary condition for $v^\ep$ is highly oscillating, it is natural to expect that these oscillations will dampen as one moves away from the boundary and that therefore most of the energy of~$v^\ep$ should be confined to a ``boundary layer'', that is, concentrated near the boundary. This expectation is reflected in our definition of~$T_\ep$, which can be thought of as a rough guess for $v^\ep$. It is
\begin{equation} 
\label{e.T_ep}
T_\ep(x):=   \left(\indc_{\R^d \setminus U_{2R(\ep)}} \ast \zeta_{R(\ep)} \right)(x) \sum_{k=1}^d \ep \phi_{e_k}^\ep (x)  \partial_{x_k} \left( u \ast \zeta_\ep \right)(x),
\end{equation}
where we recall that $U_{r} := \left\{ x \in U \, : \, \dist (x,\partial U) > r \right\}$ and the width $R(\ep)$ of the boundary layer is given by
\begin{equation} 
\label{e.Repdef}
R(\ep):=
\left\{
\begin{aligned}
& \ep \left| \log \ep \right|^{\frac12} & \mbox{if} & \ d=2, \\
& \ep & \mbox{if} & \ d>2.
\end{aligned}
\right.
\end{equation}
Observe that the support of $T_\ep$ belongs to $\R^d \setminus U_{3R(\ep)}$.

\smallskip

The next lemma records the estimates of~$T_\ep$ which are needed in the proof of Theorem~\ref{t.twoscale}. 

\begin{lemma} 
\label{l.bndr v}
Let $s \in (0,2)$, $p \in (2,\infty)$, $\alpha \in (0,\infty)$ and 
\begin{equation}  
\label{e.beta bndr v_ep}
\beta \in \left(0,\tfrac 12 \right] \cap \left( 0 ,  \tfrac12  + \al - \tfrac1p   \right).
\end{equation}
There exist $\delta(d,\Lambda) > 0$, $C(s,p,\alpha,\beta,U, d,\Lambda)<\infty$ and, for every $\eps \in \Ll( 0,\frac 1 2 \Rr]$, a random variable $\mathcal{Z}_\ep$ satisfying
\begin{equation}
\label{e.Xeptrippybounds}
\mathcal{Z}_\ep \leq \left\{ 
\begin{aligned}
& \O_s\left( C \ep^{\frac 12} \left| \log \ep \right|^{\frac14} \right) & \mbox{if} & \ d=2  \mbox{ and } \alpha > \tfrac 1p, \\
& \O_s\left( C \ep^{\beta} \right) & \mbox{if} & \ d=2  \mbox{ and } \alpha \leq \tfrac 1p,  \\  
& \O_{2+\delta} \left( C \ep^\beta \right) & \mbox{if} & \ d>2,
\end{aligned}
\right.
\end{equation}
such that for every $u \in W^{1+\alpha,p}(\Rd)$ and $T_\eps$ as defined by \eqref{e.T_ep},
\begin{equation}
\label{e.Tepboundzz}
R(\ep)^{-1}\left\|   T_\ep  \right\|_{L^2(U)}  + \left\| \nabla  T_\ep  \right\|_{L^2(U)} \leq \mathcal{Z}_\ep \left\|  u  \right\|_{W^{1+\alpha,p}(\R^d)}.
\end{equation}
\end{lemma}
\begin{proof}
\emph{Step 1.} Our first goal is to prove the desired estimate for $\left\| \nabla  T_\ep \right\|_{L^2(U)}$. There are two random variables appearing in the proof, namely
\begin{equation}
 \label{e.nastyguys}
\left\{
\begin{aligned} 
\tilde \X_\ep^1 & := 
\sum_{k=1}^d  \left(  \frac{\ep^{d}}{R(\ep)} \sum_{z \in \ep \Z^d \cap (U \setminus U_{3R(\ep)})} \left\| \phi_{e_k}^\ep   \right\|_{\underline{L}^2(z+2\ep \cu_0)}^{\frac{2p}{p-2}} \right)^{\frac{p-2}{2p}} ,
\\ 
 \tilde \X_\ep^2 & := 
 \sum_{k=1}^d  \left( \frac{\ep^{d}}{R(\ep)}  \sum_{z \in \ep \Z^d \cap (U \setminus U_{3R(\ep)})} \left\|\nabla \phi_{e_k} \left( \tfrac\cdot\ep \right)   \right\|_{\underline{L}^2(z+2\ep \cu_0)}^{\frac{2p}{p-2}} \right)^{\frac{p-2}{2p}}.
\end{aligned}
\right. 
\end{equation}
By Theorem~\ref{t.correctors},~\eqref{e.correctorgradbound} and Lemma~\ref{l.sum-O}, we have that 
\begin{equation} \label{e.tilde X1 and tilde X2}
\tilde \X_\ep^1 
\leq 
\left\{ 
\begin{aligned}
& \O_s\left( C  \left| \log \ep \right|^{\frac12} \right) & \mbox{if} & \ d=2,\\
& \O_{2+\delta} \left( C \right) & \mbox{if} & \ d>2,
\end{aligned}
\right. 
\quad \mbox{and} \quad 
\tilde \X_\ep^2 \leq \O_{2+\delta} \left( C \right).
\end{equation}
We next compute
\begin{multline*} \label{}
\left\| \nabla T_\ep \right\|_{L^2(U)}
 \leq \sum_{k=1}^d  \ep \left\| \nabla \left( \left(\indc_{\R^d \setminus U_{2R(\ep)}} \ast \zeta_{R(\ep)} \right) \phi_{e_k}^\ep  \partial_{x_k} \left( u \ast \zeta_\ep \right) \right) \right\|_{L^2(U)}
\\ 
\leq C  \sum_{k=1}^d  \left\| \left( \frac{\ep\left| \nabla  u \ast \zeta_\ep \right| }{R(\ep)} 
 + \ep \left| \nabla \left( \nabla u \ast \zeta_\ep  \right) \right| 
 \right) \left| \phi_{e_k}^\ep  \right| +  \left| \nabla u \ast \zeta_\ep \right|  \left|\nabla \phi_{e_k} \left( \tfrac\cdot\ep \right)  \right|  \right\|_{L^2 \left(U  \setminus U_{3R(\ep)} \right)}.
\end{multline*}
The right side will be estimated with the aid of Lemma~\ref{l.convolution0}. Using~\eqref{e.convolution1}, we get
\begin{equation*} 
\ep  \left\| \, \left| \nabla \left( \nabla u \ast \zeta_\ep  \right) \right|  \left| \phi_{e_k}^\ep  \right| \right\|_{L^2 \left(U  \setminus U_{3R(\ep)} \right)} \\ 
\leq
C R(\ep)^{\frac{p-2}{2p}} \ep^{\alpha} \tilde \X_\ep^1 \left\| \nabla u \right\|_{W^{\alpha,p}(\R^d)}.
\end{equation*}
On the other hand,~\eqref{e.convolution0} implies both 
\begin{equation}  
\label{e.estim.gradu.zeta.phi}
\frac{\ep}{R(\ep)} \left\|  \left| \left( \nabla  u \ast \zeta_\ep  \right) \right| \left| \phi_{e_k}^\ep  \right| \right\|_{L^2 \left(U  \setminus U_{3R(\ep)} \right)}  
     \\
 \leq C
 \frac{\ep}{R(\ep)} R(\ep)^{\frac{p-2}{2p}}  \tilde \X_\ep^1 \left\| \nabla u \right\|_{L^p(V_\ep)}
\end{equation}
and
\begin{equation*} 
\left\| \, \left| ( \nabla u \ast \zeta_\ep ) \right|  \left|\nabla \phi_{e_k} \left( \tfrac\cdot\ep \right)  \right|  \right\|_{L^2 \left(U  \setminus U_{3R(\ep)} \right)}  \leq C R(\ep)^{\frac{p-2}{2p}} \tilde \X_\ep^2 \left\| \nabla u \right\|_{L^p(V_\ep)},
\end{equation*}
where we have set
\begin{equation*}  
V_\eps := U \setminus U_{3R(\ep)} + 2 \ep \cu_0.
\end{equation*}
Lemma~\ref{l.Rellichtype} below asserts that, for each $\sigma >0$, there exists $C(\sigma, \al,p,U,d) < \infty$ such that
\begin{equation}  \label{e.nablau bndr lr} 
\left\| \nabla u \right\|_{L^p(V_\ep)} \leq C R(\ep)^{\frac1p \wedge (\alpha - \sigma)} \left\| \nabla u \right\|_{W^{\alpha,p}(\R^d)}.
\end{equation}
For $\alpha > \tfrac 1p$, we fix $\sigma = \al - \frac 1 p > 0$ and thus obtain
\begin{equation*}  
\left\| \nabla T_\ep \right\|_{L^2(U)}  \le C  R(\ep)^\frac 1 2 \Ll[ \ep^\al R(\ep)^{-\frac 1 p} \tilde \X_\ep^1 + \frac {\ep}{R(\ep)} \tilde \X_\ep^1 + \tilde \X_\ep^2\Rr]\|\nabla u\|_{W^{\al,p}(\Rd)}  .
\end{equation*}
In view of \eqref{e.tilde X1 and tilde X2}, this is the sought-after estimate for $\left\| \nabla T_\ep \right\|_{L^2(U)}$. If instead $\al \le \frac 1 p$, then the estimate becomes
\begin{equation*}  
\left\| \nabla T_\ep \right\|_{L^2(U)}  \le C  R(\ep)^\frac 1 2 \Ll[ \ep^\al R(\ep)^{-\frac 1 p} \tilde \X_\ep^1 + R(\eps)^{\al - \frac 1 p - \sigma} \Ll(\frac {\ep}{R(\ep)}  \tilde \X_\ep^1 + \tilde \X_\ep^2\Rr)\Rr]\|\nabla u\|_{W^{\al,p}(\Rd)} .
\end{equation*}
Noting that
\begin{equation*}  
R(\ep)^\frac 1 2 \Ll[\ep^\al R(\ep)^{- \frac 1 p} + R(\ep)^{\al - \frac 1 p - \sigma} \Ll( \frac \ep {R(\ep)} + 1 \Rr) \Rr] \le C \ep^{\frac 1 2 + \al - \frac 1 p - 2 \sigma},
\end{equation*}
selecting $\sigma > 0$ such that
\begin{equation*}  
\be = \frac 1 2 + \al - \frac 1 p - 3 \sigma,
\end{equation*}
and using \eqref{e.tilde X1 and tilde X2}, we obtain the desired estimate for $\|\nabla T_\ep\|_{L^2(U)}$ in this case as well.

\smallskip

 \emph{Step 2.}
We prove the bound for  $\left\| T_\ep \right\|_{L^2(U)}$.  Notice that
\begin{equation*} 
 \left\| T_\ep \right\|_{L^2(U)} \leq \ep  \sum_{k=1}^d \left\| \phi_{e_k}^\ep (\nabla u \ast \zeta_\ep) \right\|_{L^2(U\setminus U_{3R(\ep)})} .
\end{equation*}
The term on the right side has already been estimated above in~\eqref{e.estim.gradu.zeta.phi}. We therefore obtain an upper bound for $R(\ep)^{-1} \|T_\ep\|_{L^2(U)}$ which matches the upper bound obtained for $\|\nabla T_\ep\|_{L^2(U)}$ in the previous step. This completes the proof.
\end{proof}

Above we made use of the following lemma, which answers a simple question regarding estimates for Sobolev functions in boundary layers: 
given~$\alpha \in (0,1]$ and~$p\in (2,\infty)$ such that $\frac1p<\alpha<\frac dp$ and a function $u \in W^{1+\alpha,p}(U)$, what is the best estimate for 
$\left\| \nabla u \right\|_{L^2(U \setminus U_r)}$ (for small~$r$)? 
 It might be tempting to guess that one should use the fractional Sobolev inequality (Proposition~\ref{p.sobolevpoincare}) to get that $\nabla u\in L^{\frac{dp}{d-p\alpha} }(U)$ and then use this improvement of integrability with H\"older's inequality to take advantage of the small measure of $U\setminus U_r$. This however leads to a suboptimal bound, especially in large dimensions, where we obtain an estimate which scales like a small power of~$r$. But a boundary layer is not any set of small measure, and a better estimate can be found by using the Sobolev trace theorem, Proposition~\ref{p.trace}, which gives the same estimate as for smooth functions, namely $\lesssim O(r^{\frac12})$, provided that $\alpha > \tfrac1p$. Indeed, the question we posed is essentially equivalent to one concerning the image of the trace operator.
\index{trace operator}

\begin{lemma}
\label{l.Rellichtype}
Let $U\subseteq\Rd$ be a bounded Lipschitz domain. Fix $\alpha \in (0,\infty)$, $p\in (1,\infty]$, $q\in [1,p]$, and
$$\beta \in \left( 0, \tfrac1q \right] \cap \left( 0, \tfrac1q - \tfrac 1p + \alpha \right).$$
There exists~$C(\beta,U,d,\alpha,p,q)<\infty$ such that for every $f \in W^{\alpha,p}(\R^d)$ and $r \in (0,1]$,
\begin{equation} 
\label{e.Rellichtype}
\left\| f \right\|_{L^q(\partial U + B_r)} 
\leq
C r^{\beta} \left\| f \right\|_{W^{\alpha,p}(\R^d)}. 
\end{equation}
\end{lemma}
\begin{proof}
\emph{Step 1.} We prove the result in the case $\beta = \tfrac1q$ and $\alpha > \tfrac1p$. First, by H\"older's inequality we have that 
\begin{equation*} 
\left\| f \right\|_{L^q(\partial U + B_r)} \leq C r^{\frac 1q - \frac 1p} \left\| f \right\|_{L^p(\partial U + B_r)}.
\end{equation*}
By the coarea formula we have 
\begin{equation*} 
\left\| f \right\|_{L^p(\partial U + B_r)}^p  = \int_{0}^r \left\| f \right\|_{L^p\left(\{ x \, : \, \dist(x,\partial U) = t\}\right)}^p \, dt ,
\end{equation*}
and  by Proposition~\ref{p.trace} and the assumption of $\al > \frac 1 p$, for each $t \in (0,r)$, we have
\begin{equation*} 
\left\| f \right\|_{L^p\left(\left\{ x \, : \, \dist(x,\partial U) = t\right\}\right)}^p  \leq C \left\| f \right\|_{W^{\alpha,p}\left(\Rd\right)}^p.
\end{equation*}
We thus obtain
\begin{equation*} 
\left\| f \right\|_{L^q(\partial U + B_r)}  \leq Cr^{\frac 1p} \left\| f \right\|_{W^{\alpha,p}(\R^d)},
\end{equation*}
as announced.

\smallskip

\emph{Step 2.} 
We obtain the general statement of the lemma by using the result of Step~1 and an interpolation argument. Let~$\alpha$,~$p$,~$q$ and~$\beta$ be as in the statement of the lemma. We suppose that $\alpha \leq \frac1p$. First, Lemma~\ref{l.convolution2} gives
\begin{equation*} 
 \left\| f  - ( f \ast \zeta_r) \right\|_{L^p(\partial U + B_r)}   \leq C r^\alpha  \left\| f \right\|_{W^{\alpha,p}(\R^d)}. 
\end{equation*}
On the other hand, the result of the first step ensures, for every $\sigma > 0$, the existence of a constant $C(\sigma,p,U,d) < \infty$ such that 
\begin{equation*} 
\left\| f \ast \zeta_r \right\|_{L^p(\partial U + B_r)}   \leq  C r^{\frac 1p} \left\|  f \ast \zeta_r  \right\|_{W^{\sigma + \frac1p,p}(\R^d)}.
\end{equation*}
Lemma~\ref{l.convolution2} implies that
\begin{equation*} 
\left\|  f \ast \zeta_r  \right\|_{W^{\sigma + \frac1p,p}(\R^d)} \leq C r^{\alpha - \frac 1p - \sigma} \left\|  f   \right\|_{W^{\alpha,p}(\R^d)}.
\end{equation*}
By H\"older's inequality and the triangle inequality, we thus have
\begin{equation*} 
\left\|  f \right\|_{L^q(\partial U + B_r)}  \leq C_\sigma r^{\frac 1q - \frac 1p + \alpha -\sigma} \left\|  f   \right\|_{W^{\alpha,p}(\R^d)}.
\end{equation*}
One may now choose $\sigma>0$ so that $\beta = \frac 1q - \frac 1p + \alpha -\sigma$. This completes the argument. 
\end{proof}

We now give the proof of Theorem~\ref{t.twoscale}. 

\begin{proof}[Proof of Theorem~\ref{t.twoscale}] Without loss of generality, we assume that $\al \le 1$. Denoting by $v^\ep$ the boundary layer corrector defined in \eqref{e.vepeq}, we use the triangle inequality to write
\begin{equation}  
\label{e.triangle.boundary.layer}
\left\| \nabla u^\ep - \nabla {w}^\ep \right\|_{L^2(U)} \le  
 \left\| \nabla v^\ep \right\|_{L^2(U)}  + \left\| \nabla (u^\ep +v^\ep - w^\ep )  \right\|_{L^2(U)}.
\end{equation}
We decompose the proof into two steps, each estimating one of the terms above. As will be seen, the dominant error comes from the first term when $\al > \frac 1 2$, and comes from the second term otherwise. 

\smallskip

\emph{Step 1.} We first estimate the contribution of the boundary layer corrector~$v^\ep$. The claim is that, for every 
\begin{equation*}
\beta \in \left(0,\tfrac 12 \right] \cap \left( 0 ,  \tfrac12 - \tfrac1p +  \alpha  \right), 
\end{equation*}
we have
\begin{equation}  
\label{e.vepblest}
\left\| \nabla v^\ep \right\|_{L^2(U)}  
\leq  
C \mathcal{Z}_{\ep}  \left\| u \right\|_{W^{1+\alpha,p}(\Rd)},
\end{equation}
where $\mathcal{Z}_{\ep}$ is the random variable from Lemma~\ref{l.bndr v}.
We use the variational formulation to see that 
\begin{align} \label{e.v_ep var}
\left\| \nabla v^\ep \right\|_{L^2(U)}^2 
& \leq  \int_U  \nabla v^\ep \cdot \a^\ep \nabla v^\ep 
\\ \notag & 
\leq  \inf_{v \in (w^\ep - u)+H^1_0(U)}  \int_U  \nabla v \cdot \a^\ep \nabla v
\leq  C \inf_{v \in (w^\ep - u)+H^1_0(U)}\left\| \nabla v\right\|^2_{L^2(U)} ,
\end{align}
and test it with the function $v = T_\ep$ defined in~\eqref{e.T_ep}. 
Observe that, by the definition of~$w_\ep$, we have $T_\ep \in (w^\ep - u)+H^1_0(U)$. Moreover, Lemma~\ref{l.bndr v} gives us
\begin{equation} \label{e.v bnd est 7.7.7}
R(\ep)^{-1} \left\| T_\ep \right\|_{L^2(U)}  +  \left\| \nabla T_\ep \right\|_{L^2(U)} \leq \mathcal{Z}_\ep \left\| u \right\|_{W^{1+\alpha,p}(\Rd)},
\end{equation}
which combined with~\eqref{e.v_ep var} yields~\eqref{e.vepblest}.

\smallskip

\emph{Step 2.} We complete the proof that 
\begin{equation} 
\label{e.justbreath}
\left\| \nabla u^\ep - \nabla w^\ep \right\|_{L^2(U)} 
\leq
\left\| u \right\|_{W^{1+\alpha,p}(\Rd)}
\cdot
\left\{ 
\begin{aligned}
& \O_s\left( C \ep^\alpha \left| \log \ep \right|^{\frac12} \right) & \mbox{if} & \ d=2, \ \alpha \in \left( 0,\tfrac12 \right], \\
& \O_s\left( C \ep^{\frac12} \left| \log \ep \right|^{\frac14} \right) & \mbox{if} & \ d=2, \ \alpha \in \left( \tfrac12,\infty \right), \\
& \O_{2+\delta} \left( C  \ep^{\alpha\wedge \frac12} \right) & \mbox{if} & \ d>2, \ \alpha \in (0,\infty).
\end{aligned}
\right.
\end{equation} 
Define
\begin{equation*} \label{}
z^\ep:= u^\ep + v^\ep - w^\ep\in H^1_0(U).
\end{equation*}
Testing the equations for $u^\ep$ and $v^\ep$ by $z^\ep$ and summing the results, we get
\begin{equation*} \label{}
\int_{U} \nabla z^\ep(x) \cdot \a^\ep(x) \left( \nabla u^\ep (x) + \nabla v^\ep(x) \right)\,dx 
\\
= \int_{U}  \nabla z^\ep(x)  \cdot \ahom \nabla u(x)\,dx.
\end{equation*}
On the other hand, 
\begin{multline*} \label{}
\left|
\int_{U} \nabla z^\ep(x)
\cdot \a^\ep(x) \nabla w^\ep(x)\,dx 
- \int_U \nabla z^\ep(x) \cdot \ahom \nabla u(x)\,dx
\right|
\\
\leq C  \left\|\nabla z^\ep \right\|_{L^2(U)} \left\| 
\nabla  \cdot \a^\ep \nabla w^\ep -\nabla  \cdot \ahom \nabla u  
 \right\|_{H^{-1}(B_1)}.
\end{multline*}
Combining these yields 
\begin{align*} \label{}
\left\| \nabla z^\ep \right\|_{L^2(U)}^2
&
\leq 
C\int_{U} \nabla z^\ep(x)
\cdot \a^\ep(x) \nabla z^\ep(x)\,dx 
\\ & 
\leq C \left\|\nabla z^\ep \right\|_{L^2(U)} 
\left\| \nabla  \cdot
 \a^\ep \nabla w^\ep - \nabla  \cdot \ahom \nabla u  \right\|_{H^{-1}(B_1)},
\end{align*}
and therefore, by Theorem~\ref{t.twoscaleplugveng},
\begin{align} 
\label{e.uepwepgrad0}
\left\| \nabla ( u^\ep + v^\ep - w^\ep ) \right\|_{L^2(U)}
& 
= \left\| \nabla z^\ep \right\|_{L^2(U)}
\notag \\ & 
\leq
C \left\| 
\nabla \cdot \a^\ep \nabla w^\ep - \nabla \cdot \ahom \nabla u  
 \right\|_{H^{-1}(B_1)}
\notag \\ &
 \leq 
 C \left\| u \right\|_{W^{1+\alpha,p}(\Rd)} \cdot 
 \left\{ 
\begin{aligned}
& \O_s\left( C\ep^\alpha  \left| \log \ep \right|^{\frac 12}\right) & \mbox{if} & \ d=2,\\
& \O_{2+\delta} \left( C\ep^\alpha \right) & \mbox{if} & \ d>2.
\end{aligned}
\right.
\end{align}
With~\eqref{e.triangle.boundary.layer} and the result of the previous step, this yields~\eqref{e.justbreath}.

\smallskip

\emph{Step 3.} We complete the proof of the theorem. The remaining task is to obtain an estimate on~$\left\| u^\ep - w^\ep \right\|_{L^2(U)}$. For this we use the triangle inequality and the Poincar\'e inequality to get
\begin{align*} \label{}
\left\| u^\ep - w^\ep \right\|_{L^2(U)}
&
\leq  \left\| T_\ep  \right\|_{L^2(U)}  + \left\| v^\ep -T_\ep  \right\|_{L^2(U)} + \left\| z^\ep  \right\|_{L^2(U)}
\\ & 
\leq  \left\| T_\ep  \right\|_{L^2(U)}  + \left\| \nabla v^\ep - \nabla T_\ep  \right\|_{L^2(U)} + \left\| \nabla z^\ep  \right\|_{L^2(U)}
\\ & 
\leq  \left\| T_\ep  \right\|_{L^2(U)}  + \left\| \nabla v^\ep \right\|_{L^2(U)} + \left\| \nabla T_\ep  \right\|_{L^2(U)} + \left\| \nabla z^\ep  \right\|_{L^2(U)}.
\end{align*}
Each of the four terms on the right side of the previous inequality has been already estimated above, in~\eqref{e.Tepboundzz},~\eqref{e.vepblest} and~\eqref{e.justbreath}. 
This completes the proof of the theorem. 
\end{proof}

We now present an $H^1$ estimate of the two-scale expansion for the Neumann problem. The proof is an adaptation of the one for Theorem~\ref{t.twoscale} and is left to the reader. 
\index{Neumann problem}

\begin{theorem}
[Quantitative two-scale expansion for the Neumann problem]
\emph{}\\
Fix $s\in (0,2)$, $\alpha \in (0,\tfrac12]$, $p\in (2,\infty]$, a bounded Lipschitz domain $U\subseteq\Rd$ and $\ep\in \left(0,\tfrac 12\right]$. 
There exist $\delta (d,\Lambda)>0$, $C(s,\alpha,p,U,d,\Lambda)<\infty$ and a random variable $\X_\ep$ satisfying~\eqref{e.XepSI}
such that the following holds: for every $u \in W^{1+\alpha,p}(\Rd)$, if we let $w^\ep\in H^1(U)$ be defined by~\eqref{e.wepdef} and $u^\ep\in H^1(U)$ be the solution of the Neumann problem
\label{t.twoscale.Neumann}
\begin{equation*} \label{}
\left\{
\begin{aligned}
& -\nabla \cdot \left(\a^\ep(x) \nabla u^\ep \right) = -\nabla \cdot \left( \ahom \nabla u \right) & \mbox{in} & \ U, \\
& \mathbf{n}\cdot \a \left( \tfrac \cdot\ep\right)\nabla u^\ep = \mathbf{n} \cdot \ahom\nabla u  & \mbox{on} & \ \partial U,
\end{aligned}
\right.
\end{equation*}
then we have the following $H^1$--type estimate for the error of the two-scale expansion:
\begin{equation} 
\label{e.twoscale.Neumann}
\left\|  \nabla u^\ep - \nabla w^\ep  \right\|_{L^2(U)}
\leq
\X_\ep \left\| u \right\|_{W^{1+\alpha,p} (\Rd)}.
\end{equation}
\end{theorem}

\begin{exercise}
Prove Theorem~\ref{t.twoscale.Neumann}. 
\end{exercise}

\section{Boundary layer estimates}
\label{s.blayers}
Reinspecting the proof of Theorem~\ref{t.twoscale}, we see that, for smooth data, there are two sources of error: 
\begin{itemize}

\item The size of the boundary layer solution $v^\ep$. The estimate we have is
\begin{equation*} \label{}
\left\| \nabla v^\ep \right\|_{L^2(U)} \lesssim O\big(\ep^{\frac12}\big),
\end{equation*}
up to the logarithmic correction in $d=2$. This error arises because near the boundary, the homogenization process gives way to the enforcement of the Dirichlet condition---and we have perturbed the boundary condition in our definition of~$w^\ep$. 
The estimate is optimal: indeed, we should expect that, in at least an $O(\ep)$-thick boundary layer, the gradient of $u^\ep$ is essentially tracking the boundary condition, and thus that $|\nabla u^\ep - \nabla w^\ep|$ is of order one there. This already leads to an $L^2$ norm of at least~$O({\ep}^{\frac12})$.

\item The second source of error is the responsibility of homogenization, and it is captured by the quantity~$ \left\|\nabla  \cdot \a^\ep \nabla w^\ep - \nabla  \cdot\ahom \nabla u   \right\|_{H^{-1}(B_1)}$, which we  estimated in Theorem~\ref{t.twoscaleplugveng}. 
This error is much smaller than the boundary layer error: up to logarithmic corrections in $d=2$, it is~$O(\ep)$.
\end{itemize}
While the boundary layer solution~$v^\ep$ is the leading source of the~$H^1$ error, there is good reason to expect it to make a lower-order contribution to the~$L^2$ error. Indeed, a rough guess informed by the maximum principle is that (up to logarithmic corrections in $d=2$), since~$v^\ep(x)$ has boundary values of size~$O(\ep)$, it should be of order~$\ep$ in the whole domain. We would also like to show that the contribution of~$\nabla v^\ep$ to the error in the two-scale expansion is mostly confined to boundary layers and is much smaller than~$O(\ep^{\frac12})$ in the interior of the domain. 

\smallskip

It is the purpose of this section to improve our estimates of~$v^\ep$ and consequently obtain a weighted~$H^1$ estimate for the error in the two-scale expansion as well as an~$L^2$ estimate for the homogenization error. Both of these estimates, stated in the following theorem, are optimal in the scaling of $\ep$ for smooth data. It is convenient to denote 
\begin{equation} 
\label{e.rhoepU}
\rho_U^\ep(x):= \ep \vee \dist(x,\partial U),\quad x\in U. 
\end{equation}

\begin{theorem}[{Optimal $L^2$ homogenization error estimate}]
\label{t.L2EE}
Fix $s\in (0,2)$, $\alpha \in (0,1]$, $p\in (2,\infty]$, a bounded domain $U\subseteq\Rd$ which is either convex or $C^{1,1}$, and $\ep\in \left(0,\tfrac 12\right]$. 
There exist $\delta (d,\Lambda)>0$, $C(s,\alpha,p,U,d,\Lambda)<\infty$ and a random variable $\X_\ep$ satisfying
\begin{equation*} \label{}
\X_\ep \leq 
\left\{ 
\begin{aligned}
& \O_{s/2}\left( C \ep^\alpha \left| \log \ep \right|^{\frac12} \right) & \mbox{if} & \ d=2,\\
& \O_{1+\delta} \left( C  \ep^\alpha \right) & \mbox{if} & \ d>2,
\end{aligned}
\right.
\end{equation*}
such that the following holds: for every $u \in W^{1+\alpha,p}(\Rd)$, if we let $w^\ep$ be given by~\eqref{e.wepdef} and $u^\ep\in H^1(U)$ be the solution of~\eqref{e.twoscaleDP}, then we have
\begin{equation} 
\label{e.weightedtwoscale}
\left\| \left(  \nabla u^\ep - \nabla w^\ep \right) \rho_{U}^\ep \right\|_{L^2(U)} 
\leq 
\X_\ep \left\| u \right\|_{W^{1+\alpha,p}(\Rd)}
\end{equation}
and
\begin{equation} 
\label{e.L2EE}
\left\| u^\ep - u \right\|_{L^2(U)} 
\leq
\X_\ep \left\| u \right\|_{W^{1+\alpha,p}(\Rd)}. 
\end{equation}
\end{theorem}

Notice that~\eqref{e.weightedtwoscale} provides a much stronger estimate than Theorem~\ref{t.twoscale} in the subdomain~$U_r$, for~$\ep^{\frac12}\ll r \leq 1$. In particular, if~$u\in W^{2,p}(\Rd)$ with~$p>2$ and~$\overline{V}\subseteq U$ is a given compact subdomain, then there exists~$C(V,s,p,U,d,\Lambda)<\infty$ such that
\begin{equation} 
\label{e.twoscale.VsubseteqU}
\left\| \nabla u^\ep - \nabla w^\ep \right\|_{L^2(V)} 
\leq 
C \left\| u \right\|_{W^{2,p}(\Rd)} \cdot
\left\{
\begin{aligned}
& \O_{s/2}\left( C \ep \left| \log \ep \right|^{\frac12} \right) & \mbox{if} & \ d=2,\\
& \O_{1+\delta} \left( C  \ep \right) & \mbox{if} & \ d>2.
\end{aligned}
\right.
\end{equation}
However,~\eqref{e.weightedtwoscale} is still not optimal in terms of estimates of the boundary layer errors: 
ignoring the stochastic integrability exponent, the correct estimate should be roughly 
\begin{equation} 
\label{e.boundarylayervengeance}
\left\| \left(  \nabla u^\ep - \nabla w^\ep \right) \left( \rho_{U}^\ep \right)^{\frac12} \right\|_{L^2(U)} 
\leq 
 \left\| u \right\|_{W^{2,p}(\Rd)} \cdot \left\{
\begin{aligned}
& \O\left( C \ep \left| \log \ep \right|^{\frac12} \right) & \mbox{if} & \ d=2,\\
& \O \left( C  \ep \right) & \mbox{if} & \ d>2.
\end{aligned}
\right.
\end{equation}
Compared to~\eqref{t.twoscale} and~\eqref{e.weightedtwoscale}, this estimate gives improved bounds in every subdomain $U_r$ for $\ep \ll r \ll 1$. The proof of~\eqref{e.boundarylayervengeance} is more involved than that of~\eqref{e.weightedtwoscale} and requires estimates for~$\nabla u^\ep -\nabla w^\ep$ in~$L^p(U)$ for~$p>d$.

\smallskip

Regarding the assumptions on $U$: Theorem~\ref{t.L2EE} requires the boundary of~$U$ to have quite a bit more regularity than just that it is Lipschitz, namely that it be convex or~$C^{1,1}$. This can be relaxed: what we really need is that, for some $\delta>0$, the estimate~\eqref{e.psibarw32} holds for solutions of~\eqref{e.psioverlineh}. In general Lipschitz domains, this estimate is borderline, as unfortunately the best estimate---which is known to be sharp---is that~\eqref{e.psibarw32} holds with $\delta=0$ (see~\cite{JeKe}). However, and this will be clear from the proof of Theorem~\ref{t.L2EE}, we can obtain a statement in a general Lipschitz domain by giving up an arbitrarily small power of $\ep$ in the estimates~\eqref{e.weightedtwoscale} and~\eqref{e.L2EE}. We expect this estimate to be true in general $C^{1,\gamma}$ domains for any $\gamma \in (0,1]$, but have not found a proof in the literature. 

\smallskip

We begin the proof of Theorem~\ref{t.L2EE} with an extension of Lemma~\ref{l.Rellichtype} for solutions of the Dirichlet problem~\eqref{e.twoscaleDP}. Like in Chapter~\ref{c.regularity}, the idea is to ``borrow'' the regularity of the homogenized solution given in Lemma~\ref{l.Rellichtype} by using the two-scale expansion result of Theorem~\ref{t.twoscale}.
\index{boundary layer estimate|(}

\begin{proposition}[{Dirichlet boundary layer estimate}]
\label{p.Rellich}
Fix $s\in (0,2)$, $\alpha \in \left(0,\tfrac12\right]$, $p\in (2,\infty]$, a bounded Lipschitz domain $U\subseteq\Rd$, $\ep\in \left(0,\tfrac 12\right]$, and let $R(\ep)$ be given by~\eqref{e.Repdef}.
There exist $\delta (d,\Lambda)>0$, $C(s,\alpha,p,U,d,\Lambda)<\infty$ and, for each $r\in \left[ R(\ep),1\right]$, a random variable $\Y_{r,\ep}$ satisfying
\begin{equation} 
\label{e.Yrep.bound}
\Y_{r,\ep} \leq 
\left\{ 
\begin{aligned}
& \O_s\left( C \right) & \mbox{if} & \ d=2,\\
& \O_{2+\delta} \left( C  \right) & \mbox{if} & \ d>2,
\end{aligned}
\right.
\end{equation}
such that the following holds: for every $u \in W^{1+\alpha,p}(\Rd)$, if we let $u^\ep\in H^1(U)$ be the solution of the Dirichlet problem~\eqref{e.twoscaleDP},
then we have
\begin{equation} 
\label{e.Rellich}
\left\| \nabla u^\ep \right\|_{\underline{L}^2(U\setminus U_r)} 
\leq 
\Y_{r,\ep}
r^{\alpha-\frac12}
\left\| u \right\|_{W^{1+\alpha,p}(\Rd)}.
\end{equation}
\end{proposition}
\begin{proof}
For $w^\eps$ defined in \eqref{e.wepdef}, the triangle inequality gives
\begin{equation*} \label{}
\left\| \nabla u^\ep \right\|_{\underline{L}^2(U\setminus U_r)} 
\leq 
\left\| \nabla w^\ep  \right\|_{\underline{L}^2(U\setminus U_r)} 
+ \left\| \nabla w^\ep - \nabla u^\ep \right\|_{\underline{L}^2(U\setminus U_r)}.
\end{equation*}
By Theorem~\ref{t.twoscale}, with $\X_\ep$ as in the statement of that theorem, we have, for every $r\in [R(\ep),1]$,
\begin{align*} \label{}
\left\| \nabla w^\ep - \nabla u^\ep \right\|_{\underline{L}^2(U\setminus U_r)}
& 
= C \left( \frac1r \int_{U \setminus U_r} \left| \nabla w^\ep - \nabla u^\ep \right|^2   \right)^{\frac12} 
\\ & 
\leq Cr^{-\frac12} \left\| \nabla w^\ep - \nabla u^\ep \right\|_{L^2(U)}
\leq C (R(\ep))^{-\frac12}\X_\ep\left\| u \right\|_{W^{1+\alpha,p}(\Rd)}.
\end{align*}
As $(R(\ep))^{-\frac12}\X_\ep$ is a random variable with the desired integrability, by~\eqref{e.XepSI}, it remains to bound~$\left\| \nabla w^\ep  \right\|_{\underline{L}^2(U\setminus U_r)}$. We break this up by the triangle inequality:
\begin{align}
 \label{e.nablawepbreakup}
\left\| \nabla w^\ep \right\|_{\underline{L}^2(U\setminus U_r)}
& 
\leq  
\left\| \nabla u \right\|_{\underline{L}^2(U\setminus U_r)}
+ 
\sum_{k=1}^d \left\|  \left(\partial_{x_k} u \ast\zeta_\ep \right) \nabla \phi_{e_k} \left( \tfrac\cdot\ep \right)  \right\|_{\underline{L}^2(U\setminus U_r)}
\notag \\ & \qquad 
+ \ep \sum_{k=1}^d \left\| \phi_{e_k}^\ep \nabla   \left(\partial_{x_k} u \ast\zeta_\ep \right) \right\|_{\underline{L}^2(U\setminus U_r)}.
\end{align}
Each term on the right side of~\eqref{e.nablawepbreakup} will be estimated separately. For the remainder of the proof, we fix
\begin{equation*}  
q := \frac 1 2 (p+2) \in (2,p).
\end{equation*}

\smallskip

\emph{Step 1.} We show that, for every $r \in (0,1]$,
\begin{equation} 
\label{e.nablaudumbbound2}
\left\| \nabla u \right\|_{\underline{L}^q(U^r\setminus U_r)}
\leq C r^{\alpha-\frac12} \left\| u \right\|_{W^{1+\alpha,p}(\Rd)}.
\end{equation}
(Recall that~$U^r$ is defined in~\eqref{e.Ur.def}.)
Lemma~\ref{l.Rellichtype} with $\beta =\tfrac1q\wedge \alpha$ yields, for every $r\in (0,1]$,
\begin{equation} 
\label{e.nablaudumbbound1}
\left\| \nabla u \right\|_{\underline{L}^q(U^r\setminus U_r)}
\leq Cr^{\beta - \frac1q} \left\| u \right\|_{W^{1+\alpha,p}(\Rd)}.
\end{equation}
Since~$\beta - \tfrac1q \geq \alpha - \tfrac12$, we deduce that~\eqref{e.nablaudumbbound2} holds. 

\smallskip

\emph{Step 2.} We show that, for every $r \in [\ep,1]$ and $k\in\{1,\ldots,d\}$,
\begin{multline} 
\label{e.nablauXbound1}
\left\|  \left(\partial_{x_k} u \ast\zeta_\ep \right) \nabla \phi_{e_k} \left( \tfrac\cdot\ep \right)  \right\|_{\underline{L}^2(U\setminus U_r)}
\\
\leq 
 C r^{\alpha-\frac12} \left\| u \right\|_{W^{1+\alpha,p}(\Rd)}\left\| \, \left|  \nabla \phi_{e_k} \left( \tfrac\cdot\ep \right)  \right| \ast\zeta_\ep \right\|_{\underline{L}^{q/(q-2)}(U^r\setminus U_{2r})}.
\end{multline}
Fix $r\in [\ep,1]$. Using H\"older's inequality and~\eqref{e.nablaudumbbound2}, we get
\begin{align*}
\left\|  \left(\partial_{x_k} u \ast\zeta_\ep \right) \nabla \phi_{e_k} \left( \tfrac\cdot\ep \right)  \right\|_{\underline{L}^2(U\setminus U_r)}
&
\leq \left\| \, \left|  \nabla  u \right| \left(  \left| \nabla \phi_{e_k} \left( \tfrac\cdot\ep \right)  \right| \ast\zeta_\ep \right) \right\|_{\underline{L}^2(U^r\setminus U_{2r})}
\\ &
\leq \left\| \nabla u \right\|_{\underline{L}^q(U^r\setminus U_{2r})} \left\| \, \left| \nabla \phi_{e_k} \left( \tfrac\cdot\ep \right)  \right| \ast\zeta_\ep \right\|_{\underline{L}^{q/(q-2)}(U^r\setminus U_{2r})}
\\ & 
\leq C r^{\alpha-\frac12} \left\| u \right\|_{W^{1+\alpha,p}(\Rd)}\left\| \, \left|  \nabla \phi_{e_k} \left( \tfrac\cdot\ep \right)  \right| \ast\zeta_\ep \right\|_{\underline{L}^{q/(q-2)}(U^r\setminus U_{2r})}.
\end{align*}
This is~\eqref{e.nablauXbound1}.

\smallskip

\emph{Step 3.} We show that, for every $r \in [\ep,1]$ and $k\in \{1,\ldots,d\}$,
\begin{multline} 
\label{e.nablauXbound2}
 \ep \left\| \phi_{e_k}^\ep \nabla   \left(\partial_{x_k} u \ast\zeta_\ep \right) \right\|_{\underline{L}^2(U\setminus U_r)}
\\
\leq 
 C r^{\alpha-\frac12} \left\| u \right\|_{W^{1+\alpha,p}(\Rd)} \left( \ep \left\| \, \left| \phi_{e_k}^\ep  \right| \ast\zeta_\ep \right\|_{\underline{L}^{q/(q-2)}(U^r\setminus U_{2r})} \right).
\end{multline}
The proof of~\eqref{e.nablauXbound2} is almost identical to that of~\eqref{e.nablauXbound1}, so we omit it.  

\smallskip

\emph{Step 4.} The conclusion. Combining the previous steps yields 
\begin{align*} \label{}
\left\| \nabla u^\ep \right\|_{\underline{L}^2(U\setminus U_r)} 
\leq
\left\| \nabla w^\ep  \right\|_{\underline{L}^2(U\setminus U_r)} 
+ \left\| \nabla w^\ep - \nabla u^\ep \right\|_{\underline{L}^2(U\setminus U_r)}
\leq r^{\alpha-\frac12} \Y_{r,\ep} \left\| u \right\|_{W^{1+\alpha,p}(U)},
\end{align*}
where we define the random variable $\Y_{r,\ep}$ by 
\begin{multline*} \label{}
\Y_{r,\ep}
:=
C(R(\ep))^{-\frac12}\X_\ep
+ C \Big( 1 + \left\| \, \left|  \nabla \phi_{e_k} \left( \tfrac\cdot\ep \right)  \right| \ast\zeta_\ep \right\|_{\underline{L}^{q/(q-2)}(U^r\setminus U_{2r})}
\\ 
+ \ep \left\| \, \left| \phi_{e_k}^\ep  \right| \ast\zeta_\ep \right\|_{\underline{L}^{q/(q-2)}(U^r\setminus U_{2r})}\Big),
\end{multline*}
and where $\X_\ep$ is as in Theorem~\ref{t.twoscale}. 
The estimates~\eqref{e.correctorgradbound},~\eqref{e.phioscd>2},~\eqref{e.phioscd=2} and~\eqref{e.XepSI} imply that~$\Y_{r,\ep}$ satisfies~\eqref{e.Yrep.bound}. This completes the proof. 
\end{proof}

The previous proposition allows us to improve the $L^2$ estimate for the boundary layer solution~$v^\ep$ encountered in the proof of Theorem~\ref{t.twoscale} and obtain a weighted~$L^2$ estimate for its gradient. 

\begin{proposition}[{Estimates for~$v^\ep$}]
\label{p.boundarylayersoln}
Fix exponents~$s \in (0,2)$, $p \in (2,\infty)$, $\alpha \in (0,\infty)$,
\begin{equation}  
\label{e.beta bndr v_ep222}
\beta \in \left(0,1 \right] \cap \left( 0 ,  1 + \al - \tfrac1p   \right),
\end{equation}
a bounded domain $U\subseteq\Rd$ which is either convex or $C^{1,1}$, and $\ep\in \left(0,\tfrac 12\right]$.
There exist $\delta (d,\Lambda)>0$, $C(s,\alpha,p,\beta,U,d,\Lambda)<\infty$ and a random variable $\X_\ep$ satisfying the bound
\begin{equation*} \label{}
\X_\ep \leq 
\left\{ 
\begin{aligned}
& \O_{s/2}\left( C \ep \left| \log \ep \right|^{\frac12} \right) & \mbox{if} & \ d=2 \mbox{ and } \alpha > \tfrac 1p,\\
& \O_{s/2}\left( C \ep^\beta \left| \log \ep \right|^{\frac14} \right) & \mbox{if} & \ d=2  \mbox{ and } \alpha \leq \tfrac 1p,\\
& \O_{1+\delta} \left( C\ep^\beta  \right) & \mbox{if} & \ d>2,
\end{aligned}
\right.
\end{equation*}
such that the following holds: for every $u \in W^{1+\alpha,p}(\Rd)$, if we let 
$v^\ep\in H^1(U)$ be the solution of~\eqref{e.vepeq}, then
\begin{equation} 
\label{e.vepblests}
\left\| v^\ep \right\|_{L^2(U)}
+
\left\| \nabla v^\ep \rho_{U}^\ep \right\|_{L^2(U)} 
\leq
\X_\ep \left\| u \right\|_{W^{1+\alpha,p}(\Rd)}.
\end{equation}
\end{proposition}
\begin{proof}
The main part of the argument concerns the estimate of the first term on the left side of~\eqref{e.vepblests}. The estimate for the second term follows relatively easily from this and Caccioppoli's inequality. 

\smallskip

\emph{Step 1.} The estimate for $\left\| v^\ep \right\|_{L^2(U)}$. We fix $h\in L^2(U)$ and seek to estimate
\begin{equation*} \label{}
\int_U h(x) v^\ep(x) \,dx. 
\end{equation*}
Let $\psi^\ep\in H^1_0(U)$ be the solution of 
\begin{equation}
\label{e.psieph}
\left\{
\begin{aligned}
& -\nabla \cdot \left(\a^\ep(x) \nabla \psi^\ep \right) = h & \mbox{in} & \ U, \\
& \psi^\ep = 0 & \mbox{on} & \ \partial U. 
\end{aligned}
\right.
\end{equation}
We also let $T_\ep$ be as defined in~\eqref{e.T_ep} and $R(\ep)$ be as defined in~\eqref{e.Repdef}. Since $v^\ep - T_\ep\in H^1_0(U)$, we can rewrite the Dirichlet problem for~$v^\ep$ as
\begin{equation} 
\label{e.vepDP2}
\left\{
\begin{aligned}
& -\nabla \cdot \left(\a^\ep(x) \nabla v^\ep \right) = 0 & \mbox{in} & \ U, \\
& v^\ep = T_\ep & \mbox{on} & \ \partial U. 
\end{aligned}
\right.
\end{equation}
Testing~\eqref{e.psieph} with $v^\ep - T_\ep\in H^1_0(U)$ yields
\begin{equation*} \label{}
\int_U h(x) \left( v^\ep(x) - T_\ep(x) \right) \,dx 
=
\int_U \left( \nabla v^\ep(x) - \nabla T_\ep(x) \right) \cdot \a^\ep(x) \nabla \psi^\ep(x)\,dx.
\end{equation*}
Testing~\eqref{e.vepDP2} with $\psi_\ep\in H^1_0(U)$ yields
\begin{equation*} \label{}
\int_U \nabla \psi^\ep(x)\cdot \a^\ep(x) \nabla v^\ep(x)\,dx = 0. 
\end{equation*}
Combining the previous two displays, we get
\begin{equation} 
\label{e.hvepidentity}
\int_U h(x) v^\ep(x)\,dx 
= 
\int_U h(x) T_\ep(x)\,dx 
- \int_U \nabla T_\ep(x) \cdot \a^\ep(x) \nabla \psi^\ep(x)\,dx.
\end{equation}
Since $T_\ep$ is supported in $U \setminus U_{3R(\ep)}$, an application of H\"older's inequality yields
\begin{equation}
\label{e.hvepbound}
\left| \int_U h(x) v^\ep(x)\,dx  \right|
\leq
\left\| h \right\|_{L^2(U)} \left\| T_\ep \right\|_{L^2(U )} 
+ \left\| \nabla \psi^\ep \right\|_{L^2(U\setminus U_{3R(\ep)} )} \left\| \nabla T_\ep \right\|_{L^2(U)} .
\end{equation}
Since we assume~$U$ to be either convex or~$C^{1,1}$, we have by Proposition~\ref{p.H2estimate} that the solution~$\overline{\psi}\in H^1_0(U)$ of the homogenized Dirichlet problem 
\begin{equation}
\label{e.psioverlineh}
\left\{
\begin{aligned}
& -\nabla \cdot \left( \ahom \nabla \overline{\psi} \right) = h & \mbox{in} & \ U, \\
& \overline\psi = 0 & \mbox{on} & \ \partial U,
\end{aligned}
\right.
\end{equation}
satisfies the $H^2$ estimate
\begin{equation} 
\label{e.H2estimateforC11}
\left\| \overline{\psi} \right\|_{H^2(U)} \leq C \left\| h \right\|_{L^2(U)}. 
\end{equation}
In particular, by the fractional Sobolev embedding (cf.~\cite[Theorem 7.58]{Adams}), we obtain that, for an exponent~$\delta(d)>0$, 
\begin{equation} 
\label{e.psibarw32}
\left\| \overline{\psi} \right\|_{W^{\frac32,2+\delta}(U)} \leq C \left\| h \right\|_{L^2(U)}. 
\end{equation}
An application of Proposition~\ref{p.Rellich} therefore yields that, for every $r\in [ R(\ep) ,1]$, 
\begin{equation}
\label{e.nablapsiepbad}
\left\| \nabla \psi_\ep \right\|_{L^2(U\setminus U_r)}
 \leq 
C \Y_{r,\ep} r^{\frac12} \left\| \overline{\psi} \right\|_{W^{\frac32,2+\delta}(U)} 
\leq
C \Y_{r,\ep}  r^{\frac12}\left\| h \right\|_{L^2(U)}. 
\end{equation}
Combining the previous inequality and~\eqref{e.hvepbound} yields 
\begin{align}
\label{e.hvepbound2}
\left| \int_U h(x) v^\ep(x)\,dx  \right|
\leq
\left\| h \right\|_{L^2(U)} 
\left( 
\left\| T_\ep \right\|_{L^2(U )} 
+C\Y_{R(\ep),\ep}  R(\ep)^{\frac12}
\left\| \nabla T_\ep \right\|_{L^2(U )} 
\right).
\end{align}
By Lemma~\ref{l.bndr v} and~\eqref{e.Yrep.bound}, we see that  
\begin{multline*} 
\label{}
\left\| T_\ep \right\|_{L^2(U )} 
+
\Y_{R(\ep),\ep} R(\ep)^{\frac12}
\left\| \nabla T_\ep \right\|_{L^2(U)}
\\
\leq \left\|  u  \right\|_{W^{1+\alpha,p}(\R^d)} \cdot 
\left\{ 
\begin{aligned}
& \O_{s/2}\left( C \ep \left| \log \ep \right|^{\frac12} \right) & \mbox{if} & \ d=2 \mbox{ and } \alpha > \tfrac 1p,\\
& \O_{s/2}\left( C \ep^\beta \left| \log \ep \right|^{\frac12} \right) & \mbox{if} & \ d=2  \mbox{ and } \alpha \leq \tfrac 1p,\\
& \O_{1+\delta} \left( C\ep^\beta  \right) & \mbox{if} & \ d>2.
\end{aligned}
\right.
\end{multline*}
We take $\X_\ep$ to be the random variable implicitly defined in the line above. Now taking the choice $h=v^\ep$ and combining the previous two estimates completes the proof of the estimate for~$\left\| v^\ep \right\|_{L^2(U)}$.

\smallskip

\emph{Step 2.} The estimate for $\left\| \nabla v^\ep \rho^\ep_U \right\|_{L^2(U)}$. The Caccioppoli inequality yields, for every $r\in [\ep,1]$ and $x\in U_{2r}$,
\begin{equation*} \label{}
\int_{B_{r/2}(x)} \left| \nabla v^\ep(y)\right|^2\,dy \leq \frac{C}{r^2} \int_{B_{r}(x)} \left| v^\ep(y) \right|^2\,dy.
\end{equation*}
Using a covering of $U_{2r} \setminus U_{4r}$ with balls of radius $r/2$ and centers in $U_{2r} \setminus U_{4r}$ with the property that any point of~$U_{2r} \setminus U_{4r}$ belongs to at most $C(d) < \infty$ of the balls in the covering, we obtain
\begin{equation*} \label{}
\int_{U_{2r}\setminus U_{4r} } \left| \nabla v^\ep(x)\right|^2\dist^2(x,\partial U)\,dx \leq C\int_{U_{r}\setminus U_{5r} } \left| v^\ep(x) \right|^2\,dx.
\end{equation*}
Summing over $r \in \left\{ r' \,:\, r'=2^{-k} \diam(U) \geq \tfrac14\ep, \, k\in\N \right\}$ and applying the result of Step~1, with~$\X_\ep$ the random variable defined there, we obtain
\begin{equation*} \label{}
\int_{U_{\ep} } \left| \nabla v^\ep (x)\right|^2\rho_U^\ep(x)^2 \,dx
\leq 
C \int_{U} \left| v^\ep(x) \right|^2\,dx
\leq 
C\X_\ep^2 \left\| u \right\|_{W^{1+\alpha,p}(\Rd)}^2.
\end{equation*}
On the other hand, by~\eqref{e.vepblest}, we have 
\begin{equation*} \label{}
\int_{U \setminus U_{\ep} } \left| \nabla v^\ep \right|^2\left(\rho_U^\ep\right)^2 
= \ep^2 \int_{U \setminus U_{\ep} } \left| \nabla v^\ep \right|^2 
\leq \ep^2 \int_{U } \left| \nabla v^\ep \right|^2 
\leq \ep^2 \mathcal{Z}_\ep^2 \left\| u \right\|_{W^{1+\alpha,p}(\Rd)}^2,
\end{equation*}
where $\mathcal{Z}_\ep$ is the random variable given in the statement of Lemma~\ref{l.bndr v} which satisfies~\eqref{e.Xeptrippybounds}. 
Summing the previous two displays and replacing $\X_\ep$ by the random variable~$C\X_\ep + \ep \mathcal{Z}_\ep$, which is less than $C\X_\ep$, completes the proof.
\end{proof}

\index{boundary layer estimate|)}

We now present the proof of Theorem~\ref{t.L2EE}.

\begin{proof}[{Proof of Theorem~\ref{t.L2EE}}]
\emph{Step 1.} The proof of~\eqref{e.weightedtwoscale}. Using the triangle inequality and $\rho_U^\ep \leq C$, we have
\begin{align*} 
\left\| \left(  \nabla u^\ep - \nabla w^\ep \right) \rho_{U}^\ep \right\|_{L^2(U)} 
&
\leq 
\left\| \left(  \nabla u^\ep + \nabla v^\ep - \nabla w^\ep \right) \rho_{U}^\ep \right\|_{L^2(U)} 
 +
 \left\|  \nabla v^\ep  \rho_{U}^\ep \right\|_{L^2(U)} 
 \\ & 
 \leq \left\|  \nabla u^\ep + \nabla v^\ep - \nabla w^\ep  \right\|_{L^2(U)} 
  +
 \left\|  \nabla v^\ep  \rho_{U}^\ep \right\|_{L^2(U)} .
\end{align*}
The estimate now follows from~\eqref{e.uepwepgrad0} and~\eqref{e.vepblests}. 

\smallskip

\emph{Step 2.} The proof of~\eqref{e.L2EE}. Let $w^\ep$ be given by~\eqref{e.wepdef} and~$v^\ep\in H^1(U)$ be the solution of~\eqref{e.vepeq}. Observe that 
\begin{equation*} 
u_\ep - u =  (u^\ep + v^\ep - w^\ep) - v^\ep + (w^\ep - u). 
\end{equation*}
Thus by the triangle inequality, 
\begin{equation*} \label{}
\left\| u^\ep - u \right\|_{L^2(U)} 
\leq
\left\| u^\ep + v^\ep - w^\ep \right\|_{L^2(U)} 
+ \left\|  v^\ep  \right\|_{L^2(U)} 
+ \left\| w^\ep - u \right\|_{L^2(U)} . 
\end{equation*}
By the Poincar\'e inequality and \eqref{e.uepwepgrad0}, we have 
\begin{align*}
\left\| u^\ep + v^\ep - w^\ep \right\|_{L^2(U)} 
& 
\leq 
C \left\| \nabla (u^\ep + v^\ep - w^\ep) \right\|_{L^2(U)}
\\ & 
 \leq 
 C \left\| u \right\|_{W^{1+\alpha,p}(\Rd)} \cdot 
 \left\{ 
\begin{aligned}
& \O_s\left( C\ep^\alpha  \left| \log \ep \right|^{\frac 12}\right) & \mbox{if} & \ d=2,\\
& \O_{2+\delta} \left( C\ep^\alpha \right) & \mbox{if} & \ d>2.
\end{aligned}
\right.
\end{align*}
According to Proposition~\ref{p.boundarylayersoln}, 
\begin{equation*} \label{}
\left\|  v^\ep  \right\|_{L^2(U)} 
\leq 
C \left\| u \right\|_{W^{1+\alpha,p}(\Rd)} \cdot 
\left\{ 
\begin{aligned}
& \O_s\left( C \ep^\al \left| \log \ep \right|^{\frac12} \right) & \mbox{if} & \ d=2,\\
& \O_{2+\delta} \left( C\ep^\al  \right) & \mbox{if} & \ d>2.
\end{aligned}
\right.
\end{equation*}
(Actually, Proposition~\ref{p.boundarylayersoln} gives a stronger estimate for~$\left\|  v^\ep  \right\|_{L^2(U)}$, but it cannot improve the result we prove here.)
Finally, we estimate $w^\ep - u$ using Lemma~\ref{l.convolution0},
\begin{align}
\label{e.dumpthephis}
\Ll\| w^\ep - u \Rr\| _{L^2(U)} 
& 
\leq 
\ep \sum_{k = 1}^d \Ll\| \phi_{e_k}^\ep (\partial_{x_k} u \ast \zeta_\ep) \Rr\|_{L^2(U)} 
\notag \\ & 
\leq 
C 
\ep^\alpha \sum_{k = 1}^d \Ll( \ep^d \sum_{z \in \ep \Z^d \cap U} \Ll\| \phi_{e_k}^\ep  \Rr\|_{L^2(z + 2\ep \cu_0)}^{\frac{2p}{p-2}}  \Rr)^\frac{p-2}{p} \Ll\| \nabla u \Rr\|_{W^{\alpha,p}(\Rd)}.
\end{align}
By Theorem~\ref{t.correctors}, we have
\begin{equation*}
\Ll( \ep^d \sum_{z \in \ep \Z^d \cap U} \Ll\| \phi_{e_k}^\ep \Rr\|_{L^2(z + 2\ep \cu_0)}^{\frac{2p}{p-2}}  \Rr)^\frac{p-2}{p}
\leq  
\left\{ 
\begin{aligned}
& \O_s\left( C \left| \log \ep \right|^{\frac12} \right) & \mbox{if} & \ d=2,\\
& \O_{2+\delta} \left( C  \right) & \mbox{if} & \ d>2.
\end{aligned}
\right.
\end{equation*}
This completes the proof of~\eqref{e.L2EE}.
\end{proof}

\index{two-scale expansion|)}

\section*{Notes and references}

The results in Sections~\ref{s.twoscalenobdry} and~\ref{s.twoscale.dirt} are essentially well-known in the case of periodic coefficient fields (cf.~\cite{BLP}), although the statements are usually somewhat different. Finer estimates on the homogenization error and two-scale expansion have been obtained more recently in~\cite{KS2,KLS1}.

\smallskip

Many of the results in this chapter are new in the stochastic setting, including the optimal~$O(\ep^{\frac12})$ error estimate (up to logarithmic corrections in $d=2$) for the two-scale expansion for the Dirichlet problem, the boundary layers estimates and the~$O(\ep)$ estimate for the homogenization error for the Dirichlet problem. The estimates with respect to fractional Sobolev norms also appear to be new, even in the case of periodic coefficients. Two-scale expansions in stochastic homogenization in situations without boundaries, similar to the results in Section~\ref{s.twoscalenobdry}, have previously appeared in~\cite{GNO2,GO6}. Section~\ref{s.blayers} follows the method introduced by Shen~\cite{Sh} for periodic coefficients; see also~\cite{KLS2}.



\chapter{Calder\'on-Zygmund gradient 
\texorpdfstring{$L^p$}{Lp} estimates} 
\label{c.CZ}

The large-scale regularity results proved in Chapter~\ref{c.regularity} are analogous to Schauder-type pointwise estimates. In this chapter, we give complimentary large-scale regularity results in Sobolev-type norms which can be considered as analogous to the classical Calder\'on-Zygmund estimates. In particular, we are interested in $L^p$-type bounds on the gradients of solutions to equations with a divergence-form right side. As an application, we give an improvement in the (spatial) integrability of the error in the two-scale expansion for the Dirichlet problem proved in Chapter~\ref{c.twoscale}. 

\section{Interior Calder\'on-Zygmund estimates}
\label{s.CZ}

\index{Calder\'on-Zygmund estimate!interior|(}
The classical Calder\'on-Zygmund estimate for the Poisson equation states that for every $p\in (1,\infty)$, there exists a constant $C(p,d)<\infty$ such that, for every $f\in W^{-1,p}(B_1)$ and solution $u\in H^1(B_1)$ of the equation
\begin{equation} 
\label{e.poisson}
-\Delta u = f \quad \mbox{in} \ B_1,
\end{equation}
we have that $\nabla u \in L^p_{\mathrm{loc}}(B_1)$ and the estimate
\begin{equation} 
\label{e.CZest}
\left\| \nabla u \right\|_{L^p(B_{1/2})} \leq C \left(  \left\| f \right\|_{W^{-1,p}(B_{1})} + \left\| \nabla u \right\|_{L^2(B_1)} \right). 
\end{equation}
This estimate is a cornerstone of elliptic regularity theory (we will give a complete proof of it in this section, see Proposition~\ref{p.CZsmallcontrast} below). The purpose of this section is to obtain an analogue of it for the operator $-\nabla \cdot \a\nabla$. The statement is presented below in Theorem~\ref{t.CZ}.

\smallskip

Like the H\"older-type regularity results from Chapter~\ref{c.regularity}, the version of~\eqref{e.CZest} we give is valid ``only on scales larger than the minimal scale,'' which we formalize by mollifying the small scales. As usual, we take~$(\zeta_\delta)_{\delta>0}$ to be the standard mollifier defined in~\eqref{e.standardmollifier} and~\eqref{e.standardmollifer.delta}. 
For $p \in [1,\infty]$ and $f \in L^p(\R^d)$, we define $\zeta_{\delta(\cdot)}\ast f$ pointwise as 
\begin{equation*} 
\left( \zeta_{\delta(\cdot)} \ast f \right)(x) := \left( \zeta_{\delta(x)} \ast f \right)(x)  = \delta(x)^{-d} \int_{\R^d} \zeta\left( \frac {x-y}{\delta(x)} \right) f(y) \, dy.
\end{equation*}
For $f\in L^1(\R^d)$, we define a maximal function 
\index{maximal function}
\begin{equation*} 
\mathsf{M} f(x) :=  \sup_{r > 0} \left\| f \right\|_{\underline{L}^1(B_r(x))}
\end{equation*}
and a truncated version \index{maximal function!coarsened}
\begin{equation*} 
\mathsf{M}_{\delta(\cdot)} f(x) :=  \sup_{ {r} > \delta(x)} \left\| f \right\|_{\underline{L}^1(B_ {r}(x))}.
\end{equation*}
Observe that $\mathsf{M}_{\delta(\cdot)} $ takes into account only scales larger than $\delta(\cdot)$, and thus we refer to it as a \emph{coarsened maximal function}.

\begin{theorem}[Gradient $L^p$-type estimates]
\label{t.CZ}
Fix $s\in (0,2)$, $p \in [2,\infty)$ and $\ep\in\left(0,\tfrac12\right]$. Let $\X(x)$ denote the random variable~$\X_{sd/2}(x)$ defined in Remark~\ref{r.Lipminscale}. There exists~$C(s,p,d,\Lambda)<\infty$ such that, for every $\mathbf{F} \in L^2(B_1;\Rd)$ and solution $u^\ep \in H^1(B_1)$ of the equation
\begin{equation}
 \label{e.theCZeq}
 -\nabla \cdot \left(\a^\ep \nabla u^\ep  \right) 
 = \nabla \cdot \mathbf{F} \quad \mbox{in}  \ B_1,
\end{equation}
we have the estimate
\begin{multline} 
\label{e.CZloc} 
\left\| \mathsf{M}_{\ep \X\left( \frac{\cdot}{\ep}\right) } \left( \indc_{B_1} \left| \nabla u^\ep \right|^2 \right) \right\|_{L^{p/2}(B_{1/2})} 
\\
\leq 
C \left\| \nabla u^\ep \right\|_{L^2(B_{1})}^2
 +  C \left\| \zeta_{\ep \X\left( \frac{\cdot}{\ep}\right)} \ast \left( \indc_{B_1} \left| \mathbf{F} \right|^2  \right)   \right\|_{L^{p/2}(B_{1})}  .
\end{multline}
Moreover,
for every $q \in [2,p)$, there exists $C(q,s,p,d,\Lambda)<\infty$ and a random variable $\Y_{\ep,s}$ satisfying
\begin{equation*} \label{}
\Y_{\ep,s} = \O_s(C) 
\end{equation*}
such that, for every $\mathbf{F} \in L^2(B_1;\Rd)$ and $u^\ep$ solving~\eqref{e.theCZeq}, we have 
\begin{equation} 
\label{e.CZquenched}
\left\| \zeta_{\ep} \ast \left| \nabla u^\ep \right|^2 \right\|_{L^{q/2}(B_{1/2})} 
\leq C \left\| \nabla u^\ep \right\|_{L^2(B_{1})}^2
+  \Y_{\ep,s} \left\| \zeta_{\ep} \ast \left| \mathbf{F} \right|^2  \right\|_{L^{p/2}(B_{1})}.
\end{equation}
\end{theorem}

Recall that the random variable $\X(\cdot)$ in Theorem~\ref{t.CZ} satisfies $\X(x) \leq \O_{\frac{sd}{2}}(C)$, is $2$-Lipschitz continuous (see~\eqref{e.Xs.lip}) and can be used as the random minimal scale in the Lipschitz bound~\eqref{e.Lipschitz}: that is, it satisfies the conclusion of Theorem~\ref{t.Lipschitz}.

\smallskip

Every proof of the Calder\'on-Zygmund estimates uses some kind of measure-theoretic covering argument in order to control the level sets of the function whose $L^p$ norm is being estimated. The original result, stated in terms of the boundedness from~$L^p$ to~$L^p$ of singular integral operators (arising for instance in the analysis of the Poisson equation~\eqref{e.poisson}), was formalized using the Calder\'on-Zygmund decomposition. This measure-theoretic tool was later widely used in many different contexts in harmonic analysis and the analysis of PDEs. Here we prefer to use the Vitali covering lemma, Lemma~\ref{l.vitali}. 

\smallskip

The  Calder\'on-Zygmund-type gradient~$L^p$ estimates for equations with a right-hand side are typically consequences of gradient~$L^\infty$ (Lipschitz) estimates for the same equation with zero right-hand side. There is a more general principle at work here, which roughly asserts that if an $L^2$ function $u$ can be well-approximated in $L^2$, on every scale, by $L^q$ functions (with $q>p$), with local $L^2$ errors controlled by an~$L^p$ function $g$, then this function $u$ must belong to $L^p$ (with an appropriate estimate). We formalize this fact as a pure real analysis lemma. This is where the Vitali covering lemma is used. 

\begin{lemma} 
For each~$p \in ( 2,\infty )$,~$q \in (p,\infty]$ and~$A\geq 1$, there exists $\delta_0(p,q,A,d) > 0$ and $C(p,q,A,d) < \infty$ such that the following holds for every $\delta \in (0,\delta_0]$. Let $f \in L^2(B_1)$ and $g \in L^{p}(B_1)$ be such that for every $x \in B_{\frac12}$ and $r \in \left(0,\tfrac1{4}\right]$, there exists~$f_{x,r} \in L^{q}(B_{r}(x))$ satisfying both
\label{l.CZ-4realz}
\begin{equation} 
\label{e.CZ-v1}
\left\| f_{x,r} \right\|_{\underline{L}^q(B_{r}(x) )} 
\leq 
A  \left\| g \right\|_{\underline{L}^2(B_{2r}(x)  )} + A \left\| f \right\|_{\underline{L}^2(B_{2r}(x)) } 
\end{equation}
and
\begin{equation} \label{e.CZ-v2}
\left\| f - f_{x,r} \right\|_{L^2(B_{r}(x) )} 
 \leq 
A \left\| g \right\|_{L^2(B_{2r}(x)) } + \delta  \left\| f \right\|_{L^2(B_{2r}(x)) }.
\end{equation}
Then $f\in L^{p}(B_{1/2})$ and we have the estimate
\begin{equation}
\label{e.CZ-conc1}
\left\| f \right\|_{L^p(B_{1/2}) } 
\leq C \left( \left\| f \right\|_{L^2(B_1) }  +   \left\| g \right\|_{L^p(B_1) } \right).
\end{equation}
\end{lemma}

\begin{proof}
We fix~$2< p<q\leq \infty$.

\smallskip

For $m\in (0,\infty)$, we denote $f_m:= |f|\wedge m$. By Lemma~\ref{l.simpleiter}, it suffices to show that there exists $C(p,q,A,d)<\infty$ such that, for every $m \in(1,\infty)$ and~$\frac12 \leq s < t \leq 1$, 
\begin{equation}
\label{e.thewts.st}
\left\| f_m \right\|_{L^p(B_{s})}^p \leq \frac12 \left\| f_m \right\|_{L^p(B_{t})}^p + \frac{C}{(t-s)^{\frac{d}{2}}} \left( \left\| f \right\|_{L^2(B_{1})} + \left\| g \right\|_{L^p(B_{1})} \right)^p.
\end{equation}
Indeed, after we prove~\eqref{e.thewts.st}, we can apply the lemma to
obtain 
\begin{equation*} \label{}
\left\| f_m \right\|_{L^p(B_{1/2}) } 
\leq C \left( \left\| f \right\|_{L^2(B_1) }  +   \left\| g \right\|_{L^p(B_1) } \right).
\end{equation*}
We may then send~$m\to\infty$ to obtain~\eqref{e.CZ-conc1}.

\smallskip

We fix $\frac12\leq s < t \leq 1$ for the remainder of the argument. In order to prove~\eqref{e.thewts.st}, we study the measure of the super level sets of~$f$. In fact, it is useful to define the following measures: for each Borel set~$F\subseteq B_1$, we set
\begin{equation*} 
\mu_f(F) := \int_{F \cap B_{t}} |f(x)|^2 \,dx  
\quad \mbox{and} \quad 
\mu_g(F) := \int_{F \cap B_{t}} |g(x)|^2 \,dx.
\end{equation*}
The main step in the proof of~\eqref{e.thewts.st} is to show that there exists~$C(d)<\infty$ such that, for every $\ep\in \left(0,\tfrac{1}{2A}\right)$, $M\in [4A,\infty)$ and $\lambda \in (0,\infty)$ satisfying
\begin{equation*} \label{}
\lambda \geq \lambda_0:=
8^d (t-s)^{-\frac d2} \left(\left\| f \right\|_{\underline{L}^2(B_{t})} + \ep^{-1} \left\| g \right\|_{\underline{L}^2(B_{t})} \right), 
\end{equation*}
we have a good $\lambda$ type inequality
\begin{multline}  
\label{e.CZ-towardsconcl}
\mu_f \left(  \left\{ |f|> M \lambda \right\} \cap B_{s}  \right) 
\leq 
C \left( \left(  \ep A + \delta \right)^2 + \frac{(4A)^{q}}{M^{q-2}} \right) 
\mu_f \left( \left\{ |f|> \tfrac{1}{2}\lambda \right\} \cap B_t \right)\\+ C\ep^{-2}  \left( \left(  \ep A + \delta \right)^2 + \frac{(4A)^{q}}{M^{q-2}} \right) 
\mu_g\left(  \left\{ |g|> \tfrac{\ep}2 \lambda \right\} \cap B_t \right) .
\end{multline}
Let us see that~\eqref{e.CZ-towardsconcl} implies~\eqref{e.thewts.st}: 
using a layer-cake formula and a change of variables (recall that~$p>2$), we see that, with $\ep$ and $M$ as above, for every $m\geq M\lambda_0$, 
\begin{align*} 
\frac{\left\| f_m \right\|_{L^p(B_{s})}^p}{M^{p-2}} 
&  
\leq \int_{B_{s}} \left( \frac{f_m(x)}{M}\right)^{p-2} f(x)^2\,dx
\\ &
= (p-2) \int_{0}^{\infty} \lambda^{p-3} \mu_f \left(  \left\{ f_m > M\lambda \right\} \cap B_{s} \right) \, d\lambda 
\\ & 
\leq 
\mu_f (B_{s}) \lambda_0^{p-2} +  (p-2) \int_{\lambda_0}^{m/M} \lambda^{p-3} \mu_f \left( \left\{ f_m > M\lambda \right\}  \cap B_{s} \right) \, d\lambda .
\end{align*}
We next use the estimate~\eqref{e.CZ-towardsconcl} to bound the second term on the right side:
\begin{align*} 
\lefteqn{ 
(p-2) \int_{\lambda_0}^{m/M} \lambda^{p-3} \mu_f \left( \left\{ f_m > M \lambda \right\} \cap B_{s} \right) \, d\lambda 
} 
\qquad & \\ 
\notag & 
\leq 
C \left( \left(  \ep A + \delta \right)^2 + \frac{(4A)^{q}}{M^{q-2}} \right) (p-2)  \int_{0}^{m/M} \lambda^{p-3}  \mu_f \left(\left\{ |f| > \tfrac{1}{2}\lambda \cap B_t \right\} \right)  \, d\lambda
\\ & \quad 
+ C \ep^{-2} \left( \left(  \ep A + \delta \right)^2 + \frac{(4A)^{q}}{M^{q-2}} \right)  (p-2) \int_{0}^{\infty} \lambda^{p-3}  \mu_g \left( \left\{ |g|> \tfrac{\ep}{2} \lambda \right\} \cap B_t \right)  \, d\lambda 
\\ \notag & 
\leq 
C \left( \left(  \ep A + \delta \right)^2 + \frac{(4A)^{q}}{M^{q-2}} \right) 
\left( \left\| f_m \right\|_{L^p(B_{t})}^p 
+  \ep^{-p} \left\| g \right\|_{L^p(B_{t})}^p \right).
\end{align*}
Putting these together, we obtain 
\begin{align*} 
\left\| f_m \right\|_{L^p\left(B_{s} \right)}^p
&
\leq 
C M^{p-2} \left( \left(  \ep A + \delta \right)^2 + \frac{(4A)^{q}}{M^{q-2}} \right) 
\left( \left\| f_m \right\|_{L^p(B_{t})}^p 
+  \ep^{-p} \left\| g \right\|_{L^p(B_{t})}^p \right).
\\ & \quad 
+  M^{p-2}\mu_f(B_{t}) \lambda_0^{p-2} .
\end{align*}
We next select the parameters so that the prefactor of~$\left\| f_m \right\|_{L^p(B_{t})}^p$ is less than $\frac12$. We first choose $M$ large enough, depending only on~$(p,q,A,d)$, so that~$CM^{p-q}A^q \leq \frac 14$. We then take~$\ep$ and~$\delta_0$ small enough, depending only the choice of~$M$ and~$(p,q,A,d)$, so that $CM^{p-2} (\ep A+\delta_0)^2 \leq \frac 14$. This yields a constant $C(p,q,A,d)<\infty$ such that, for every $\delta\in (0,\delta_0]$,
\begin{align*} 
\left\| f_m \right\|_{L^p\left(B_{s} \right)}^p
&
\leq 
\frac12 \left\| f_m \right\|_{L^p(B_{t})}^p 
+ C \mu_f(B_{t}) \lambda_0^{p-2}  + C \left\| g \right\|_{L^p(B_{t})}^p.
\end{align*}
Next, we note that 
\begin{align*} \label{}
\mu_f(B_{t}) \lambda_0^{p-2} 
&
\leq 
C (t-s)^{-\frac d2} \left\| f \right\|_{L^2(B_t)}^2 
\left(\left\| f \right\|_{\underline{L}^2(B_{t})} + \left\| g \right\|_{\underline{L}^2(B_{t})} \right)^{p-2}
\\ & 
\leq 
C (t-s)^{-\frac d2} \left( \left\| f \right\|_{\underline{L}^2(B_{t})} + \left\| g \right\|_{\underline{L}^2(B_{t})} \right)^{p}.
\end{align*}
Combining these yields~\eqref{e.thewts.st}. 

\smallskip

We therefore focus on the proof of~\eqref{e.CZ-towardsconcl} and fix, for the remainder of the argument, the parameters $\ep \in \left(0,\tfrac{1}{2A}\right)$, $M\in [4A,\infty)$, $\delta >0$ and $\lambda\in (\lambda_0,\infty)$. We define, for every~$x\in B_1$ and~$r\in (0,1)$, 
\begin{equation*} 
E(x,r) := \left\| f \right\|_{\underline{L}^2(B_{r}(x) \cap B_1)} + \ep^{-1} \left\| g \right\|_{\underline{L}^2(B_{r}(x) \cap B_R)}
\end{equation*}
and the set 
\begin{equation*} 
 G  := \left\{ x \in B_{s} \, : \, \lim_{r \to 0} E(x,r) = |f(x)| + \ep^{-1} |g(x)| \right\}. 
\end{equation*}
Notice that, by the Lebesgue differentiation theorem, we have $\left| B_{s} \setminus  G  \right| =0$.

\smallskip

\emph{Step 1.} We show that, for every~$x \in \{|f| > \lambda \} \cap  G$, there exists~$r_x \in \left(0,\frac{t-s}{64}R\right)$ satisfying 
\begin{equation} 
\label{e.CZ-stoppingtime}
E(x,r_x) = \lambda
\quad \mbox{and} \quad 
\sup_{r \in \left(r_x, (t-s) R\right)} E(x,r) \leq \lambda.
\end{equation}
Indeed, we have 
\begin{equation*}
\sup_{r \in \left( \frac{t-s}{64}, (t-s)\right)}
\sup_{x \in B_{s}} E(x,r) 
\leq 
8^d \sup_{x \in B_{s}} E(x,(t-s)) 
= 
\lambda_0
\end{equation*}
and, since $\lambda>\lambda_0$, for every $x\in\{|f| > \lambda\} \cap  G$,   
\begin{equation*}
\lim_{r \to 0} E(x,r) = |f(x)| + \ep^{-1} |g(x)| > \lambda \geq \lambda_0.
\end{equation*}
Therefore~\eqref{e.CZ-stoppingtime} follows from the intermediate value theorem. 

\smallskip

\emph{Step 2.} Fix~$x\in\{|f| > \lambda\} \cap  G$ and let~$r_x \in \left(0,\frac{t-s}{64}\right)$ be as in the previous step.  The stopping time information can be used to transfer information between levels $\lambda$ and radii $r_x$. Indeed, we show that
\begin{equation}  \label{e.|B_r|.CZ}
\lambda^2 \left|B_{r_x} \right| \leq 2 \mu_f\left( \left\{ |f|> \tfrac12 \lambda \right\} \cap B_{r_x}(x) \right) +2 \ep^{-2} \mu_g\left(\left\{ |g|> \tfrac{\ep}{2} \lambda \right\} \cap  B_{r_x}(x) \right) ,
\end{equation}
and letting~$f_x := f_{x,5r_x}$ be as in the hypothesis of the lemma, 
\begin{equation} \label{e.CZ-v1-2}
\left\| f_{x} \right\|_{\underline{L}^q(B_{5r_x}(x))} \leq 2A\lambda
\end{equation}
and 
\begin{equation} \label{e.CZ-v2-2}
\left\| f - f_{x} \right\|_{\underline{L}^2(B_{5r_x}(x) )}  \leq \left( \ep A + \delta \right) \lambda .
\end{equation}
To deduce~\eqref{e.|B_r|.CZ}, rephrasing~\eqref{e.CZ-stoppingtime} by means of $\mu_f$ and $\mu_g$, we get that
\begin{align} \notag 
\lambda^2 \left|B_{r_x} \right| & = \mu_f(B_{r_x}(x))+ \ep^{-2} \mu_g(B_{r_x}(x)) 
\\ \notag &\leq \mu_f\left( \left\{ |f|> \tfrac12 \lambda \right\} \cap B_{r_x}(x) \right) +  \ep^{-2} \mu_g\left(\left\{ |g|> \tfrac{\ep}{2} \lambda \right\} \cap  B_{r_x}(x) \right) + \frac 12 \lambda^2 \left|B_{r_x} \right|,
\end{align}
from which~\eqref{e.|B_r|.CZ} follows. On the other hand,~\eqref{e.CZ-v1-2} is a direct consequence of~\eqref{e.CZ-v1} and~\eqref{e.CZ-stoppingtime}, because
\begin{equation*} 
\left\| f_{x} \right\|_{\underline{L}^q(B_{5r_x}(x))} \leq A \left\| f \right\|_{\underline{L}^2(B_{10r_x}(x))} + A \left\| g \right\|_{\underline{L}^2(B_{10r_x}(x))} \leq A(1+\ep)E(x,10r_x) \leq 2A \lambda,
\end{equation*}
and~\eqref{e.CZ-v2-2} follows by~\eqref{e.CZ-v2} and~\eqref{e.CZ-stoppingtime}:
\begin{align} \notag 
\left\| f - f_{x} \right\|_{\underline{L}^2(B_{5r_x}(x) )} 
& \leq 
A  \left\| g \right\|_{\underline{L}^2(B_{10r_x}(x) ) } + \delta  \left\| f \right\|_{\underline{L}^2(B_{10r_x}(x) ) }
\\ \notag & 
\leq 
\left( \ep A + \delta \right) E\left(x,10r_x\right)
\\ \notag &
\leq \left( \ep A + \delta \right) \lambda .
\end{align}

\smallskip

\emph{Step 3.} Fix~$x\in\{|f| > \lambda\} \cap  G$, $r_x \in \left(0,\frac{t-s}{64}\right)$ and $f_x$ be as in Steps 1 and 2. We show that there exists~$C(d)<\infty$ such that, 
\begin{multline} 
\label{e.CZ-towardsconcl0}
\mu_f \left(  \left\{ |f|> M \lambda \right\} \cap B_{5r_x}(x) \right)
\leq  
C \left( \left(   \ep A  + \delta \right)^2 + \frac{(4A)^{q}}{M^{q-2}} \right)  \mu_f \left(   \left\{ |f|> \tfrac{1}{2}\lambda \right\} \cap B_{r_x}(x) \right) 
\\ 
+  C \ep^{-2} \left( \left(   \ep A  + \delta \right)^2 + \frac{(4A)^{q}}{M^{q-2}} \right) \mu_g\left(  \left\{ |g|> \tfrac {\ep}{2} \lambda \right\} \cap B_{r_x}(x)  \right). 
\end{multline}
We begin by observing that 
\begin{align*} 
\left\{ |f|> M \lambda \right\} \cap B_{5r_x}(x)
\subset 
\left( \left\{ |f_{x}|> \tfrac 12  M\lambda \right\}  \cup  \left\{ |f| \leq  2  |f - f_{x}| \right\} \right) \cap B_{5r_x}(x),
\end{align*}
and thus 
\begin{multline*} 
\mu_f \left( \left\{ |f|> M \lambda \right\} \cap B_{5r_x}(x) \right) 
\\ \leq \mu_f \left( \left\{ |f_{x}|> \tfrac 12 M\lambda \right\} \cap B_{5r_x}(x) \right) + \mu_f \left( \left\{ |f| \leq 2 |f - f_{x}| \right\} \cap B_{5r_x}(x) \right).
\end{multline*}
The second term can be easily estimated using~\eqref{e.CZ-v2-2} as
\begin{align*} 
\mu_f\left( \left\{ |f| \leq 2 \left|f - f_{x}\right| \right\}\cap B_{5r_x}(x) \right)  
& 
\leq 
4 |B_{5 r_x} (x)| \left\| f - f_{x} \right\|_{\underline{L}^2(B_{5r_x}(x) )}^2 
 \\  \notag & \leq  C \lambda^2 |B_{r_x}| \left( \ep A + \delta \right)^2.
\end{align*}
For the first one, on the other hand, we first use the triangle inequality to get
\begin{equation*} 
\mu_f \left(  \left\{ |f_{x}|> \tfrac 12  M\lambda \right\} \cap B_{5r_x}(x) \right) 
\leq  
2  \int_{ \left\{ |f_{x}|> \frac 12  M\lambda \right\} \cap B_{5r_x}(x) } \left(|f_{x}(y)|^2   + |f(y) - f_{x}(y)|^2\right) \, dy, 
\end{equation*}
and to estimate the two terms on the right we apply Chebyshev's inequality and~\eqref{e.CZ-v1-2} to obtain
\begin{align*} 
2 \int_{ \left\{ |f_{x}|> \frac 12  M\lambda \right\}\cap B_{5r_x}(x)  } |f_{x}(y)|^2 \, dy & \leq \left( \tfrac12 M \lambda \right)^{2-q} |B_{5 r_x} (x) |  \left\| f_{x} \right\|_{\underline{L}^q(B_{5r_x}(x) )}^q  \\  \notag & \leq C \lambda^2 \left| B_{r_x}\right| \frac{(4A)^{q}}{M^{q-2}} 
\end{align*}
and, by~\eqref{e.CZ-v2-2}, 
\begin{equation*} 
2  \int_{ \left\{ |f_{x}|> \frac 12  M\lambda \right\} \cap B_{5r_x}(x) } |f(y) - f_{x}(y)|^2 \, dy \leq C  \lambda^2 \left| B_{r_x}\right|  \left( \ep A + \delta \right)^2.
\end{equation*}
Therefore, we obtain that, for $C(d)<\infty$, 
\begin{equation*} 
\mu_f \left( \left\{ |f|> M \lambda \right\} \cap B_{5r_x}(x) \right)  \leq C \lambda^2 \left| B_{r_x}\right|  \left( \left( \ep A + \delta \right)^2 + \frac{(4A)^{q}}{M^{q-2}}  \right).
\end{equation*}
Plugging in~\eqref{e.|B_r|.CZ} then yields~\eqref{e.CZ-towardsconcl0}. 

\smallskip

\emph{Step 4.}
The covering argument and conclusion. Applying the Vitali covering lemma (Lemma~\ref{l.vitali}) to the open cover
\begin{equation*} \label{}
\left\{ B_{r_x}(x) \,:\, x\in \left\{ |f|> M \lambda \right\} \cap G \right\}
\end{equation*}
of~$\left\{ |f|> M \lambda \right\} \cap G$ now yields~\eqref{e.CZ-towardsconcl} by~\eqref{e.CZ-towardsconcl0}, and completes the proof of the lemma. 
\end{proof}

Before proceeding to the proof of Theorem~\ref{t.CZ}, we present, as a warm-up exercise, the derivation of the classical Calder\'on-Zygmund estimate~\eqref{e.CZest} from Lemma~\ref{l.CZ-4realz}. In fact, we give a proof of the following more general statement. 

\begin{proposition}
For each $p \in (2,\infty)$, there exist $\delta_0(p,d,\Lambda) > 0$ and $C(p,d,\Lambda) < \infty$ such that the following holds for every $\delta \in (0,\delta_0]$. Let $\a\in \Omega$ be such that, for some constant matrix $\a_0\in \R^{d\times d}_{\mathrm{sym}}$,
\label{p.CZsmallcontrast}
\begin{equation*} \label{}
\sup_{x\in B_1} \left| \a(x) - \a_0 \right|
\leq \delta. 
\end{equation*}
For every $u\in H^1(B_1)$, we have the estimate
\begin{equation} 
\label{e.CZsmallcontrast}
\left\| \nabla u \right\|_{L^p(B_{1/2})} 
\leq 
C \left(\left\| \nabla u \right\|_{L^2(B_{1})}  +  \left\| \nabla \cdot \left( \a  \nabla u \right)  \right\|_{W^{-1,p}(B_1)} \right). 
\end{equation}
\end{proposition}
\begin{proof}
If $\nabla \cdot \left( \a(\cdot) \nabla u \right)  \in W^{-1,p}(B_1)$, then by the Riesz representation theorem (see Remark~\ref{r.RRT}) there exists a vector field~$\mathbf F \in L^p(B_1;\Rd)$ such that 
\begin{equation*} \label{}
-\nabla \cdot \left( \a(\cdot) \nabla u \right)
= \nabla \cdot \mathbf F 
\quad \mbox{in} \ B_1
\end{equation*} 
and
\begin{equation*} \label{}
\left\| \mathbf F \right\|_{L^p(B_1)}
\leq
\left\| 
\nabla \cdot \left( \a(\cdot) \nabla u \right)
\right\|_{W^{-1,p}(B_1)}.
\end{equation*}
We assume~$\delta_0 \leq \frac12$ so that $\frac12 I_d\leq \a_0 \leq 2\Lambda I_d$. We will show that the assumptions of Lemma~\ref{l.CZ-4realz} are valid with $f:= \nabla u$ and $f_{x,r}:=\nabla u_{x,r}$, where $u_{x,r} \in u+H_0^1(B_{2r}(x))$ is the solution of 
\begin{equation*} \label{}
\left\{
\begin{aligned}
& -\nabla \cdot \left( \a_0 \nabla u_{x,r} \right) = 0 & \mbox{in} & \ B_{2r}(x), \\
& u_{x,r} = u & \mbox{on} & \ \partial B_{2r}(x).
\end{aligned}
\right.
\end{equation*}
Since $\a_0$ is constant, we deduce that $u_{x,r}$ satisfies the estimate, for $C(d,\Lambda)<\infty$,
\begin{equation} 
\label{e.pointwiseharmonic.a0}
\left\| \nabla u_{x,r} \right\|_{L^\infty(B_{r}(x))} \leq C \left\| \nabla u_{x,r} \right\|_{\underline{L}^2(B_{2r}(x))}.
\end{equation}
Indeed, this estimate follows from~\eqref{e.pointwiseharm} after a change of variables. Now, testing the equations for $u$ and $u_{x,r}$ with $u-u_{x,r}$ and subtracting, we obtain
\begin{multline*} \label{}
\int_{B_{2r}(x)} \left( \nabla u - \nabla u_{x,r} \right) \cdot \a   \left( \nabla u - \nabla u_{x,r} \right)
\\
= \int_{B_{2r}(x)} \left( \nabla u - \nabla u_{x,r} \right) \cdot \left( \a_0 - \a \right) \nabla u_{x,r}
 + \int_{B_{2r}(x)} \mathbf F \cdot \left( \nabla u - \nabla u_{x,r} \right).
\end{multline*}
Applying Young's inequality,  we obtain
\begin{equation*} \label{}
 \left\| \nabla u - \nabla u_{x,r} \right\|_{L^2(B_{2r}(x))}^2
\leq C\delta^2 \left\| \nabla u_{x,r} \right\|_{L^2(B_{2r}(x))}^2 + C\left\| \mathbf F \right\|_{L^2(B_{2r}(x))}^2.
\end{equation*}
Using the triangle inequality, we get
\begin{align*} \label{}
\left\| \nabla u_{x,r} \right\|_{L^2(B_{2r}(x))}
& \leq
\left\| \nabla u - \nabla u_{x,r} \right\|_{L^2(B_{2r}(x))}
+
\left\| \nabla u \right\|_{\underline{L}^2(B_{2r}(x))}
\\ &
\leq 
C\delta \left\| \nabla u_{x,r} \right\|_{L^2(B_{2r}(x))}
+ C\left\| \mathbf F \right\|_{L^2(B_{2r}(x))}
+ \left\| \nabla u \right\|_{\underline{L}^2(B_{2r}(x))}.
\end{align*}
Therefore, taking $\delta_0 > 0$ sufficiently small, we may reabsorb the first term on right side, and so we obtain
\begin{equation*} \label{}
\left\| \nabla u_{x,r} \right\|_{L^2(B_{2r}(x))}
\leq 
C \left\| \nabla u \right\|_{\underline{L}^2(B_{2r}(x))} 
+ C \left\| \mathbf F \right\|_{\underline{L}^2(B_{2r}(x))} .
\end{equation*}
Combining this with~\eqref{e.pointwiseharmonic.a0}, we get
\begin{equation*} 
\left\| \nabla u_{x,r} \right\|_{L^\infty(B_{r}(x))} 
\leq 
C \left\| \nabla u \right\|_{\underline{L}^2(B_{2r}(x))} 
+ C \left\| \mathbf F \right\|_{\underline{L}^2(B_{2r}(x))} .
\end{equation*}
We may therefore apply Lemma~\ref{l.CZ-4realz} with $q=\infty$, $f:= \nabla u$, $f_{x,r}:=\nabla u_{x,r}$ and $g = \mathbf F$, to obtain the conclusion of the proposition. 
\end{proof}

\begin{exercise}
\label{ex.CZcontinuous}
Use Proposition~\ref{p.CZsmallcontrast} to obtain the following statement: suppose~$\a\in\Omega$ is uniformly continuous in $B_1$ with modulus $\rho$, that is,
\begin{equation*} \label{}
\left| \a(x) - \a(y) \right| \leq \rho\left( |x-y| \right) \quad \forall x,y\in B_1,
\end{equation*}
for some continuous function $\rho:[0,2) \to \R_+$ with $\rho(0)=0$. Then, for every exponent~$p\in(2,\infty)$, there exists~$C(\rho,p,d,\Lambda)<\infty$ such that~\eqref{e.CZsmallcontrast} holds. 
\end{exercise}

\begin{remark}[{Counterexample for $p=\infty$}]
\label{r.noCZforp=infty}
Proposition~\ref{p.CZsmallcontrast} is false, in general, for $p=\infty$, even in the case of the Laplace operator~$-\Delta$. The standard counterexample (taken from~\cite[Example 10.2]{LL}) in dimension $d = 2$ is, for $x = (x_1,x_2) \in \R^2 \setminus \{0\}$,
\begin{equation*} \label{}
u(x) := x_1x_2 \log \log \frac{1}{|x|}. 
\end{equation*}
One can check that~$-\Delta u \in L^\infty(B_1)$ but~$\partial_{x_1}\partial_{x_2} u \not \in L^\infty(B_r)$ for any~$r>0$. Therefore, if we take $v:= \partial_{x_1} u$, then 
\begin{equation*} \label{}
-\Delta v = \nabla \cdot \mathbf{f} 
\quad \mbox{for} \quad
\mathbf{f} := \begin{pmatrix} 
\Delta u \\
0
\end{pmatrix}.
\end{equation*}
We thus have~$\mathbf{f} \in L^\infty(B_1;\R^2)$, but~$\nabla v\not \in L^\infty(B_r;\R^2)$ for any~$r>0$. 

\smallskip

Interior Calder\'on-Zygmund estimates may also be formulated for $p<2$. For the Laplacian operator, the estimates are true for $p\in (1,2)$, but false for $p=1$: see Exercise~\ref{ex.CZduality.interior} below. 
\end{remark}

We now focus on the proof of Theorem~\ref{t.CZ}. The following lemma  together with Lemma~\ref{l.maximalfct} gives the first statement of Theorem~\ref{t.CZ}. The proof of the lemma is similar to that of Proposition~\ref{p.CZsmallcontrast}: the main difference is that, rather than freezing the coefficients and using pointwise bounds for harmonic functions, we use the $C^{0,1}$-type estimate from Chapter~\ref{c.regularity}. Another difference is technical: since the latter is not a true pointwise estimate, we have to work with the coarsened maximal function.  Here and throughout the rest of this section,~$\X(\cdot)$ is the random variable in the statement of Theorem~\ref{t.CZ}.

\begin{lemma}
\label{l.thetheoremreally} 
For each $p \in [2,\infty)$, there exists a constant $C(d,p) < \infty$ such that for every $\mathbf{F} \in L^2(B_1;\Rd)$ and $u^\ep \in H^1(B_1)$ solution of 
\begin{equation*} \label{}
 -\nabla \cdot \left(\a^\ep \nabla u^\ep  \right) = \nabla \cdot \mathbf{F} \quad \mbox{in}  \ B_1,
\end{equation*}
we have the estimate
\begin{multline} \label{e.CZmaxest} 
\left\| \mathsf{M}_{\ep \X(\frac \cdot\ep)} \left( \indc_{B_1} |\nabla u^\ep|^2 \right) \right\|_{L^{\frac p2}(B_{1/2})} 
\\ \leq C \left\| \nabla u^\ep \right\|_{L^2(B_{1})}^2 + C \left\| \mathsf{M}_{\ep \X(\frac \cdot\ep)} (\indc_{B_1} \left| \mathbf{F} \right|^2) \right\|_{L^{\frac p2}(B_{1})} .
\end{multline}
\end{lemma}

\begin{proof}
\emph{Step 1.} We set up the argument. 
For each $x \in B_{1/2}$, define 
\begin{equation*} 
f(x) := \max_{ r  \in \left[ \ep \X_s \left(\frac x\ep\right) \wedge \frac12 , \frac12 \right]} \left( \frac{1}{ r } \left\| u^\ep - (u^\ep)_{B_ r (x)} \right\|_{\underline{L}^2(B_ r (x) )} \right) + \left\| \nabla u^\ep \right\|_{L^2(B_{1})}.
\end{equation*}
Observe that the Caccioppoli inequality gives
\begin{equation} \label{e.basiccomp000}
 \left\| \nabla u^\ep  \right\|_{\underline{L}^2(B_{ r /2}(x))} \leq  \frac{C}{ r } \left\| u^\ep - (u^\ep)_{B_ r (x)} \right\|_{\underline{L}^2(B_{ r }(x) )}  + C \left\| \mathbf F   \right\|_{\underline L^2(B_{ r }(x))} ,
\end{equation}
and hence we have a pointwise bound
\begin{equation} \label{e.max vs f}
\left( \mathsf{M}_{\ep \X(\frac \cdot\ep)} \left( \indc_{B_1} |\nabla u^\ep|^2 \right)  \right) (x) \leq C f^2(x) + C  \left( \mathsf{M}_{\ep \X(\frac \cdot\ep)} \left(\indc_{B_1} \left| \mathbf{F}\right|^2 \right) \right) (x),
\end{equation}
for every $x \in B_{1/2}$. We set
\begin{equation*} 
g(x) := \left( \mathsf{M}_{\ep \X(\frac \cdot\ep)} \left(\indc_{B_1} \left| \mathbf{F}\right|^2 \right) \right)^\frac12  (x).
\end{equation*}
Now, it is enough to establish $L^p$-integrability for $f$. We aim at applying Lemma~\ref{l.CZ-4realz}, and hence, for any fixed $x$ and $r$, we need to find $f_{x,r}$ as in Lemma~\ref{l.CZ-4realz}. For this, fix $x \in B_{1/2}$ and $r \in \left(0,\tfrac 1{4}\right]$.

\smallskip 

\emph{Step 2.} We claim that, if $\ep \X(\frac x\ep) \geq 4r$, then
\begin{equation} \label{e.locCZ-stoc-f bnd1}
\sup_{z \in B_r(x)} f(z) \leq C \inf_{z \in B_r(x)} f(z) \leq C \left\| f \right\|_{\underline{L}^2(B_r(x))}, 
\end{equation}
and in this case we simply set $f_{x,r} := f$. Since $\X(\cdot)$ is 2-Lipschitz continuous and $\ep \X(\frac x\ep) \geq 4r$, we have that $\inf_{z \in B_r(x)} \ep \X\left( \frac{z}{\ep} \right) \geq 2r$. We obtain, for all $y,z \in B_r(x)$  and $t \in [r,1/20]$,  
\begin{equation*} 
\frac1t \left\| u^\ep - (u^\ep)_{B_t(y)} \right\|_{\underline{L}^2(B_t(y))} \leq  \frac{C}{t+2r} \left\| u^\ep - (u^\ep)_{B_{t+2r}(z)} \right\|_{\underline{L}^2(B_{t+2r}(z))} ,
\end{equation*}
and from this it is easy to verify~\eqref{e.locCZ-stoc-f bnd1}. 

\smallskip

\emph{Step 3.} 
For the rest of the proof, we assume that  $\ep \X(\frac x\ep) \leq 4r$. By the 2-Lipschitz continuity, we also have
$\sup_{z \in B_{2r}(x)} \ep \X\left( \frac{z}{\ep} \right) \leq 8r$. We define
\begin{equation*} 
f_{x,r}(y) := \max_{ r'  \in \left[ \ep \X_s \left(\frac y\ep\right) \wedge  \frac12  ,  \frac12  \right]} \left( \frac{1}{ r' } \left\| v_r^\ep - (v_r^\ep)_{B_{r'} (y)} \right\|_{\underline{L}^2(B_{r'} (y))} \right) +  \left\| \nabla u^\ep \right\|_{L^2(B_{1})}, 
\end{equation*}
where $v_r^\ep \in u^\ep + H_0^1(B_{2r}(x) )$ solves $-\nabla\cdot(\a \nabla v^\ep) = 0$ in $B_{2r}(x) $, and $v_r^\ep$ is extended to be $u^\ep$ outside of $B_{2r}(x)$. 
We will show that 
\begin{equation} \label{e.CZ-v2-verified1pre}
\sup_{y \in B_r(x)}f_{x,r}(y) \leq C \inf_{y \in B_r(x)} f(y) + C\left\| g   \right\|_{\underline{L}^2(B_{2r}(x))}
\end{equation}
and
\begin{equation} \label{e.CZ-v2-verified2pre}
  \left\| f - f_{x,r}   \right\|_{\underline{L}^2(B_{r}(x) )}  \leq C\left\|  g  \right\|_{\underline{L}^2(B_{2r}(x))}. 
\end{equation}

\smallskip

\emph{Step 4.} We show~\eqref{e.CZ-v2-verified1pre}.  By testing the subtracted equations of $u^\ep$ and $v^\ep$ with $u^\ep - v^\ep$, we obtain
\begin{equation} \label{e.basiccomp001}
 \left\| \nabla u^\ep - \nabla v^\ep \right\|_{\underline{L}^2(B_{2r}(x))} \leq C  \left\| \mathbf F   \right\|_{\underline{L}^2(B_{2r}(x))}. 
\end{equation}
Furthermore,  since $\sup_{z \in B_{2r}(x)} \ep \X\left( \frac{z}{\ep} \right) \leq 8r$, we have 
\begin{equation} \label{e.basiccomp002}
\left\| \mathbf F   \right\|_{L^2(B_{2r}(x))} \leq  C \left\| \mathsf{M}_{8r} \left(\indc_{B_1} \left| \mathbf{F}\right|^2 \right) \right\|_{L^1(B_{2r}(x))} \leq C  \left\| g \right\|_{L^2(B_{2r}(x))}. 
\end{equation}
By Theorem~\ref{t.Lipschitz} we have that, for every $y,z \in B_{r}(x)$ and $r' \in \left[ \ep \X_s \left(\tfrac y\ep\right) , r \wedge \tfrac 14 \right] $,
\begin{equation} 
\frac{1}{ r' } \left\| v_r^\ep - (v_r^\ep)_{B_{r'} (y)} \right\|_{\underline{L}^2(B_{r'} (y))} \leq  \frac{C}{ 2r } \left\| v_{r}^\ep - (v_r^\ep)_{B_{2r} (z)} \right\|_{\underline{L}^2(B_{2r} (z))}.
\end{equation}
By the triangle inequality and Poincar\'e's inequality, together with~\eqref{e.basiccomp001} and~\eqref{e.basiccomp002}, recalling that $r \leq \tfrac14$, 
\begin{equation*} 
\frac1{2r}\left\| v_{r}^\ep - (v_r^\ep)_{B_{2r} (z)} \right\|_{\underline{L}^2(B_{2r} (z))} 
 \leq f(z) + C \left\| g \right\|_{\underline{L}^2(B_{2r}(x))},
\end{equation*}
and similarly, for any $z \in B_r(x)$,
\begin{equation*} 
\sup_{r' \in \left[r \wedge \frac14,  \frac 12\right]} \frac1{r'}\left\| v_{r}^\ep - (v_r^\ep)_{B_{r'} (y)} \right\|_{\underline{L}^2(B_{r'} (y))} \leq C f(z) +  \left\| g \right\|_{\underline{L}^2(B_{2r}(x))}.
\end{equation*}
Combining the last two displays yields~\eqref{e.CZ-v2-verified1pre}.

\smallskip

\emph{Step 5.} We show~\eqref{e.CZ-v2-verified2pre}. The triangle inequality yields 
\begin{align*} 
\left|f(y) - f_{x,r}(y) \right| \leq \max_{ r'  \in \left[ (\ep \X_s\left(\frac y\ep\right) \wedge \frac 12), \frac12 \right]} \left( \frac{1}{ r' } \left\| u^\ep - v^\ep - (u^\ep - v^\ep)_{B_{r'} (y)} \right\|_{\underline{L}^2(B_{r'} (y))} \right),
\end{align*}
and hence the Sobolev-Poincar\'e inequality with the exponent 
\begin{equation*} 
2_* := 
\left\{
\begin{aligned}
&\tfrac{2d}{d+2},   &  d>2 , \\
& \tfrac32,   &  d = 2,
\end{aligned}
\right. 
\end{equation*} 
gives, for some constant $C(d)<\infty$, the bound
\begin{align} \label{e.CZ-SP}
\left|f(y) - f_{x,r}(y) \right|^{2_*} & \leq C  \left( \mathsf{M}_{\ep \X\left(\frac \cdot\ep\right)} \left| \nabla u^\ep - \nabla v^\ep \right|^{2_*} \right)(y).
\end{align}
The strong-type $\left(\frac{2}{2_*}, \frac{2}{2_*}\right)$ estimate for maximal functions in Lemma~\ref{l.maximalfct} yields
\begin{multline*} 
\left\|  \mathsf{M}_{\ep \X\left(\frac \cdot\ep\right)} \left| \nabla u^\ep - \nabla v^\ep \right|^{2_*}   \right\|_{L^{\frac{2}{2_*}}(\R^d)} \\ \leq C
\left\| \left| \nabla u^\ep - \nabla v^\ep \right|^{2_*}  \right\|_{L^{\frac{2}{2_*}}(\R^d)} = C \left\| \nabla u^\ep - \nabla v^\ep  \right\|_{L^{2}(B_{2r}(x))}^{2_*}.
\end{multline*}
Combining this with~\eqref{e.basiccomp001},~\eqref{e.basiccomp002} and~\eqref{e.CZ-SP}, we obtain~\eqref{e.CZ-v2-verified2pre}.

\smallskip

\emph{Step 6.} The conclusion. We have verified assumptions~\eqref{e.CZ-v1} and~\eqref{e.CZ-v2} of Lemma~\ref{l.CZ-4realz} with~\eqref{e.CZ-v2-verified1pre} and~\eqref{e.CZ-v2-verified2pre}, respectively, and hence 
we obtain 
\begin{equation*} 
\left\| f \right\|_{L^p(B_{1/2}) }  \leq C \left(  \left\| g   \right\|_{L^p(B_{1}) } +  \left\| f \right\|_{L^2(B_1) }  \right).
\end{equation*}
This together with~\eqref{e.max vs f} finishes the proof. 
\end{proof}

We now complete the proof of Theorem~\ref{t.CZ}.

\begin{proof}[{Proof of Theorem~\ref{t.CZ}}] 
To prove~\eqref{e.CZquenched}, fix $s \in [1,2)$ and take the random variable $\X_{\frac{sd}{2}}$ accordingly. Fix also $p \in (2,\infty)$ and $q \in [2,p)$.  We denote, in short,  $\X = \X_{\frac{sd}{2}}$.  We apply Lemma~\ref{l.thetheoremreally} and start estimating terms appearing in~\eqref{e.CZmaxest}. We have
\begin{align*} 
\mathsf{M}_{\ep \X(\frac \cdot\ep)} \left(\indc_{B_1} \left| \mathbf{F} \right|^2 \right) (x) 
& \leq C \sup_{t>0} \fint_{B_t(x)}  \left( \zeta_{\ep \X(\frac x\ep)} \ast \left(  \indc_{B_1} \left| \mathbf{F} \right|^2  \right) \right)(y) \,dy 
\\ & \leq  C  \mathsf{M} \left( \zeta_{\ep \X(\frac \cdot\ep)} \ast   \left(  \indc_{B_1} \left| \mathbf{F} \right|^2   \right)\right)(x).
\end{align*}
Therefore, the strong $\left(\frac p2, \frac p2\right)$-type estimate for maximal functions in Lemma~\ref{l.maximalfct} yields
\begin{equation*} 
\left\| \mathsf{M}_{\ep \X(\frac \cdot\ep)} \left(\indc_{B_1} \left| \mathbf{F} \right|^2 \right)  \right\|_{L^{\frac p2}(B_1) } \leq C 
\left\|  \zeta_{\ep \X(\frac \cdot\ep)} \ \ast \left(  \indc_{B_1} \left| \mathbf{F} \right|^2   \right) \right\|_{L^{\frac p2}(B_1) }.
\end{equation*}
This implies~\eqref{e.CZloc} by Lemma~\ref{l.thetheoremreally}. Note that a similar argument also gives 
\begin{equation*} 
\left\| \mathsf{M}_{\ep \X(\frac \cdot\ep)} \left(\indc_{B_1} \left| \mathbf{F} \right|^2 \right)  \right\|_{L^{\frac p2}(B_1) } \leq C 
\left\|  \zeta_\ep \ast \left(  \indc_{B_1} \left| \mathbf{F} \right|^2   \right) \right\|_{L^{\frac p2}(B_1) }.
\end{equation*}
Using this we proceed to prove~\eqref{e.CZquenched}, and in view of the previous display and Lemma~\ref{l.thetheoremreally}, we only need to estimate the integral on left in~\eqref{e.CZmaxest}. We first have a pointwise bound 
\begin{equation*} 
\mathsf{M}_{\ep } \left(\indc_{B_1} \left| \nabla u^\ep \right|^2 \right) (x)  
\leq \X^d\left(\tfrac x\ep\right) \mathsf{M}_{\ep \X\left(\frac \cdot\ep\right)} \left(\indc_{B_1} \left| \nabla u^\ep \right|^2 \right) (x).
\end{equation*}
Therefore, H\"older's inequality gives
\begin{equation*} 
\left\| \mathsf{M}_{\ep } \left(\indc_{B_1} \left| \nabla u^\ep \right|^2 \right) \right\|_{L^{\frac q2}(B_{1/2})} 
\leq  \left\| \mathsf{M}_{\ep \X\left(\frac \cdot\ep\right)} \left(\indc_{B_1} \left| \nabla u^\ep \right|^2 \right) \right\|_{L^{\frac p2}(B_{1/2})} 
\left\| \X\left(\tfrac \cdot\ep\right) \right\|_{L^{\frac d2 \frac{qp}{p-q}} (B_{1/2})}^{\frac d2}.
\end{equation*}
By Lemma~\ref{l.sum-O}, we have $\left\| \X_s \left(\tfrac \cdot\ep\right) \right\|_{L^{\frac d2 \frac{qp}{p-q}} (B_{1/2})} \leq \O_{\frac{s d}{2}}(C(p,q,d,\Lambda))$, and hence
\begin{equation*} 
\Y_s := C \left\| \X\left(\tfrac \cdot\ep\right) \right\|_{L^{\frac d2 \frac{qp}{p-q}} (B_{1/2})}^{\frac d2} \leq \O_{s}(C). 
\end{equation*}
 The proof is complete.
 \end{proof}

\index{Calder\'on-Zygmund estimate!interior|)}

%
%
%

%

\section{Global Calder\'on-Zygmund estimates}

\index{Calder\'on-Zygmund estimate!global|(}

In this section, we prove the following global version of the Calder\'on-Zygmund gradient $L^p$-type estimates for Dirichlet boundary conditions. The argument is a variation of the one from the previous section.

\begin{theorem}[Global gradient $L^p$-type estimates]
\label{t.CZ.global}
Fix $s\in (0,2)$, $\gamma \in (0,1]$, $p \in [2,\infty)$, $\ep\in\left(0,\tfrac12\right]$ and a bounded $C^{1,\gamma}$-domain~$U\subseteq\Rd$. 
Let $\X(x)$ denote the random variable~$\X_{sd/2}(x)$ defined in Remark~\ref{r.Lipminscale}.
There exists~$C(s,p,U,d,\Lambda)<\infty$ such that, for every $\mathbf{F} \in L^2(B_1 \cap U;\Rd)$ and solution $u^\ep \in H^1(B_1 \cap U)$ of the equation
\begin{equation}
 \label{e.theCZeq.global}
 \left\{ 
 \begin{aligned}
&  -\nabla \cdot \left(\a^\ep \nabla u^\ep  \right) 
 = \nabla \cdot \mathbf{F}  & \mbox{in} & \  B_1 \cap U,\\
& u^\ep = 0 & \mbox{on} & \ B_1 \cap \partial U,
\end{aligned}
 \right.
\end{equation}
we have the estimate
\begin{multline} 
\label{e.CZloc.global} 
\left\| \mathsf{M}_{\ep \X\left( \frac{\cdot}{\ep}\right) } \left( \indc_{B_1 \cap U} \left| \nabla u^\ep \right|^2 \right) \right\|_{L^{p/2}(B_{1/2} \cap U)} 
\\
\leq 
C \left\| \nabla u^\ep \right\|_{L^2(B_{1} \cap U)}^2
 +  C \left\| \zeta_{\ep \X\left( \frac{\cdot}{\ep}\right)} \ast \left( \indc_{B_1 \cap U} \left| \mathbf{F} \right|^2  \right)   \right\|_{L^{p/2}(B_{1} \cap U)}  .
\end{multline}
\end{theorem}

We begin the proof of Theorem~\ref{t.CZ.global} by giving a global version of Lemma~\ref{l.CZ-4realz}. The main difference is that since the statement is now global, there is no need for a localization argument. The only regularity requirement we impose on~$\partial U$ is a mild geometric condition~\eqref{e.CZ-geomU.global}, which provides a suitable doubling condition.

\begin{lemma} 
\label{l.CZ-4realz.global}
For each~$p \in [ 2,\infty)$,~$q \in (p,\infty]$,~$A\geq 1$, and bounded domain $U\subseteq \Rd$ 
satisfying the geometric condition
\begin{equation} \label{e.CZ-geomU.global}
\sup_{x \in U} \sup_{r>0}\left(\frac{|B_{2r}(x) \cap U|}{|B_{r}(x) \cap U|}\right) \leq  C_U < \infty,
\end{equation}
there exist $\delta_0(p,q,A,C_U,d)>0$ and $C(p,q,A,C_U,d)<\infty$ such that the following holds for every $\delta \in (0,\delta_0]$.  Let $f \in L^2(U)$ and $g \in L^{p}(U)$ be such that
for every $x \in U$ and $r \in \left(0,\tfrac1{8} \diam(U)\right]$, there exists~$f_{x,r} \in L^{q}(B_{r}(x) \cap U)$ satisfying both
\begin{equation} 
\label{e.CZ-v1.global}
\left\| f_{x,r} \right\|_{\underline{L}^q(B_{r}(x)\cap U)} 
\leq 
A  \left\| g \right\|_{\underline{L}^2(B_{2r}(x) \cap U )} + A \left\| f \right\|_{\underline{L}^2(B_{2r}(x)\cap U) } 
\end{equation}
and
\begin{equation} \label{e.CZ-v2.global}
\left\| f - f_{x,r} \right\|_{L^2(B_{r}(x) \cap U )} 
 \leq 
A \left\| g \right\|_{L^2(B_{2r}(x)\cap U) } 
+ \delta  \left\| f \right\|_{L^2(B_{2r}(x)\cap U) }.
\end{equation}
Then $f\in L^p(U)$ and we have the estimate
\begin{equation}
\label{e.CZ-conc1.global}
\left\| f \right\|_{L^p(U) } 
\leq C \left(  \left\| g \right\|_{L^p(U) } +  |U|^{-\frac{p-2}{p}}\left\| f \right\|_{L^2(U) }  \right).
\end{equation}
\end{lemma}
\begin{proof}
Assume without loss of generality that $p>2$ and $q<\infty$. The proof closely follows the lines of the proof of Lemma~\ref{l.CZ-4realz}. Consider $\delta \in (0,1]$ to be a free parameter in the proof. Fix also a parameter $\ep\in \left(0,\tfrac{1}{2A} \right) $ to be selected below and define, for every $x\in U$ and $r>0$, 
\begin{equation*} 
E(x,r) := \left\| f \right\|_{\underline{L}^2(B_{r}(x) \cap U)} + \ep^{-1} \left\| g \right\|_{\underline{L}^2(B_{r}(x)\cap U)}.
\end{equation*}
Let $R:= \diam(U)$ and define 
\begin{equation} \label{e.lambda_0.global}
\lambda_0 := \sup_{x \in U} \left( 16^d C_U^{8} E(x,R)^2 \right)^{\frac12} \leq 8^d C_U^{4} \left(\left\| f \right\|_{\underline{L}^2(U)} + \ep^{-1} \left\| g \right\|_{\underline{L}^2(U)} \right). 
\end{equation}
Proceeding as in the proof of Lemma~\ref{l.CZ-4realz}, by introducing pointwise exit radii as in the proof of Lemma~\ref{l.CZ-4realz} using~\eqref{e.CZ-geomU.global} to control volumes close to the boundary, and finally applying Vitali's covering theorem, one deduces, for $\lambda>\lambda_0$, that
\begin{multline}  \label{e.CZ-towardsconcl.global}
\mu_f \left(  \left\{ |f|> M \lambda \right\}  \right) 
\\ 
\leq 
C  \left( \left(   \ep A  + \delta \right)^2 + \frac{(4A)^{q}}{M^{q-2}} \right) \left( 
\mu_f \left( \left\{ |f|> \tfrac{1}{2}\lambda \right\} \right) + \ep^{-2} \mu_g\left( \left\{ |g|> \tfrac \ep2  \lambda \right\}  \right)\right),
\end{multline}
where, for each Borel set~$F\subseteq U$, we set
\begin{equation*} 
\mu_f(F) := \int_{F \cap U} |f(x)|^2 \,dx  
\quad \mbox{and} \quad 
\mu_g(F) := \int_{F \cap U} |g(x)|^2 \,dx.
\end{equation*}
Using a layer-cake formula, for $p>2$, to get
\begin{align*} 
\frac{\left\| f \right\|_{L^p(U)}^p}{M^{p-2}} &  = (p-2) \int_{0}^{\infty} \lambda^{p-3} \mu_f \left(  \left\{ |f| > M\lambda \right\} \right) \, d\lambda \\
& \leq 
 \mu_f(U) \lambda_0^{p-2} +  (p-2) \int_{\lambda_0}^{\infty} \lambda^{p-3} \mu_f \left( \left\{ |f| > M\lambda \right\} \right) \, d\lambda .
\end{align*}
we then obtain by~\eqref{e.CZ-towardsconcl.global} that
\begin{align*} 
\lefteqn{ (p-2) \int_{\lambda_0}^{\infty} \lambda^{p-3} \mu_f \left( \left\{ |f| > M \lambda \right\} \right) \, d\lambda } \qquad & \\ 
\notag & \leq  C \left( \left(   \ep A  + \delta \right)^2 + \frac{(4A)^{q}}{M^{q-2}} \right) \left( \left\| f \right\|_{L^p(U)}^p + \ep^{-p}\left\| g \right\|_{L^p(U)}^p\right) .
\end{align*}
Thus, we arrive at
\begin{align*} 
 \left\| f \right\|_{L^p(U)}^p & \leq  C M^{p-2} \left( \left(   \ep A  + \delta \right)^2 + \frac{(4A)^{q}}{M^{q-2}} \right)\left( \left\| f \right\|_{L^p(U)}^p + \ep^{-p}  \left\| g \right\|_{L^p(U)}^p \right)
 \\ & \quad +  M^{p-2}\mu(U) \lambda_0^{p-2}  .
\end{align*}
Consequently, we may choose first $M$ large and $\delta$ small, and then $\ep$ accordingly small, all by means of $(p,q,A,C_U,d)$, and then reabsorb the first term on the right. This is possible since $q>p>2$.  Thus, for $C(p,q,A,C_U,d) < \infty$, we have
\begin{equation*} 
\left\| f \right\|_{L^p(U)} \leq C |U|^{-\frac{p-2}{p}} \left\| f \right\|_{L^2(U)} + C \left\| g \right\|_{L^p(U)}.
\end{equation*}
To make the argument rigorous, one repeats the computation of the last step with $f_m := |f| \wedge m$  instead, and after reabsorption pass to the limit $m \to \infty$ using the monotone convergence theorem. The proof is complete. 
\end{proof}

\begin{exercise}
Show that if $|U| = \infty$, then the conclusion of Lemma~\ref{l.CZ-4realz.global} becomes 
\begin{equation*} 
\left\| f \right\|_{L^p(U)} \leq C \left\| g \right\|_{L^p(U)}.
\end{equation*}
\end{exercise}

\begin{exercise}
\label{ex.CZsmallcontrast.global}
By mimicking the proof of Proposition~\ref{p.CZsmallcontrast}, use Lemma~\ref{l.CZ-4realz.global} to prove a global Calder\'on-Zygmund estimate for Dirichlet boundary conditions. The statement is: fix a bounded domain $U\subseteq\Rd$ satisfying the geometric condition~\eqref{e.CZ-geomU.global}, $p\in (2,\infty)$ and $\delta >0$. Suppose that $\a\in \Omega$ satisfies, for some constant matrix $\a_0\in \R^{d\times d}_{\mathrm{sym}}$,
\begin{equation*} \label{}
\sup_{x\in U} \left| \a(x) - \a_0 \right|
\leq \delta. 
\end{equation*}
Then there exists~$\delta_0(U,p,d,\Lambda)>0$ and~$C(U,p,d,\Lambda)<\infty$ such that~$\delta \leq \delta_0$ implies that, for every~$u\in H^1(B_1)$, $f \in W^{-1,p}(U)$ and~$g\in W^{1,p}(U)$ satisfying 
\begin{equation*} \label{}
\left\{
\begin{aligned}
& -\nabla \cdot \left( \a(x) \nabla u \right) = f & \mbox{in} & \ U, \\
& u = g & \mbox{on} & \ \partial U,
\end{aligned}
\right.
\end{equation*}
we have the estimate
\begin{equation} 
\label{e.CZsmallcontrast2}
\left\| \nabla u \right\|_{L^p(U)} 
\leq 
C \left(  \left\| f \right\|_{W^{-1,p}(B_1)} + \left\| \nabla g \right\|_{L^p(U)} \right). 
\end{equation}
\end{exercise}

\begin{exercise}
State and prove a version of the global Calder\'on-Zygmund estimates for Neumann boundary conditions. 
\end{exercise}

\begin{exercise}
State and prove a global version of Exercise~\ref{ex.CZcontinuous} for Dirichlet boundary conditions. 
\end{exercise}

\begin{remark}[Duality for Calder\'on-Zygmund exponents]
\label{r.CZduality}
The statement proved in Exercise~\ref{ex.CZsmallcontrast.global} can be extended to all~$p\in (1,\infty)$. Indeed, the global Calder\'on-Zygmund estimate for $p\in [1,\infty]$ is \emph{equivalent} to the same statement for the H\"older conjugate exponent~$p'$. To see this, we suppose the statement is true for some $p\in [1,\infty]$ and seek to prove the estimate for $p'$. We therefore consider, for a given $\f \in C^\infty(U;\Rd)$, a solution~$u\in H^1(U)$ of the problem 
\begin{equation*} \label{}
\left\{
\begin{aligned}
& -\nabla \cdot \left( \a(x) \nabla u \right) = \nabla \cdot \f & \mbox{in} & \ U, \\
& u = 0 & \mbox{on} & \ \partial U,
\end{aligned}
\right.
\end{equation*}
To obtain a bound on $\left\| \nabla u \right\|_{L^{p'}(U)}$, we argue by duality. 
We fix a vector field $\g \in L^p(U;\Rd)$ and seek to bound the quantity $\left| \int_U \g\cdot \nabla u \right|$. We let $v$ be the solution of 
\begin{equation*} \label{}
\left\{
\begin{aligned}
& -\nabla \cdot \left( \a(x) \nabla v \right) = \nabla \cdot \g & \mbox{in} & \ U, \\
& v = 0 & \mbox{on} & \ \partial U.
\end{aligned}
\right.
\end{equation*}
By assumption, we have 
\begin{equation*} \label{}
\left\| \nabla v \right\|_{L^p(U)} \leq C \left\| \g \right\|_{L^p(U)}. 
\end{equation*}
Testing the equations for $u$ and $v$ with each other and using the H\"older inequality, we get
\begin{align*}
\left| \int_U \g \cdot \nabla u  \right|
=
\left| \int_U \nabla v \cdot \nabla u \right| 
=
\left| \int_U \f \cdot \nabla v \right|
\leq \left\| \f \right\|_{L^{p'}(U)} \left\| \nabla v \right\|_{L^p(U)} 
\leq C \left\| \f \right\|_{L^{p'}(U)} \left\| \g \right\|_{L^p(U)}. 
\end{align*}
By duality, this yields 
\begin{equation*} \label{}
\left\| \nabla u \right\|_{L^{p'}(U)} 
\leq 
C  \left\| \f \right\|_{L^{p'}(U)}
\end{equation*}
and completes the proof of the claim. 

\smallskip

In particular, in view of the example given in Remark~\ref{r.noCZforp=infty}, the global Calder\'on-Zygmund statement proved in Exercise~\ref{ex.CZsmallcontrast.global} is false, even for the Laplacian operator, for the exponents~$p = 1$ and~$p=\infty$. 
\end{remark}

\begin{exercise}
\label{ex.CZduality.interior}
Using the previous remark, extend the interior estimate of Proposition~\ref{p.CZsmallcontrast} to all exponents~$p\in (1,\infty)$. Since~\eqref{e.CZsmallcontrast} is trivial for $p<2$, replace this inequality with 
\begin{equation*} \label{}
\left\| \nabla u \right\|_{L^p(B_{1/2})} 
\leq 
C \left(\left\| u \right\|_{L^p(B_{1})}  +  \left\| \nabla \cdot \left( \a  \nabla u \right)  \right\|_{W^{-1,p}(B_1)} \right). 
\end{equation*}
\end{exercise}

Returning to the context of homogenization and continuing  with the proof of Theorem~\ref{t.CZ.global}, we next give a global analogue of Lemma~\ref{l.thetheoremreally}. 

\begin{lemma}
\label{l.thetheoremreally.global} 
For each $s \in (0,2)$, $\gamma \in (0,1]$, $p \in [2,\infty)$, $\ep\in\left(0,\tfrac12\right]$. and $C^{1,\gamma}$ domain $U \subset \Rd$, there exists a constant $C(s,p,d,U,\Lambda) < \infty$ such that for every $\mathbf{F} \in L^2(B_1 \cap U)$ and solution $u^\ep \in H^1(B_1 \cap U)$ of the equation
\begin{equation*}
 \left\{ 
 \begin{aligned}
&  -\nabla \cdot \left(\a^\ep \nabla u^\ep  \right) 
 = \nabla \cdot \mathbf{F}  & \mbox{in} & \  B_1 \cap U,\\
& u^\ep = 0 & \mbox{on} & \ B_1 \cap \partial U,
\end{aligned}
 \right.
\end{equation*}
we have the estimate
\begin{multline} 
\label{e.CZmaxest.global} 
\left\| \mathsf{M}_{\ep \X(\frac \cdot\ep)} \left( \indc_{B_{1}\cap U} |\nabla u^\ep|^2 \right) \right\|_{L^{\frac p2}(B_{1/2}\cap U)} 
\\ 
\leq 
C \left\| \nabla u^\ep \right\|_{L^2(B_{1}\cap U)}^2 
+ C \left\| \mathsf{M}_{\ep \X(\frac \cdot\ep)} \left(\indc_{B_1\cap U} \left| \mathbf{F} \right|^2\right) \right\|_{L^{\frac p2}(B_{1}\cap U)} .
\end{multline}
\end{lemma}
\begin{proof}
Since the estimate we are looking for is localized in $B_{1/2} \cap U$, we let $\tilde U$ be a $C^{1,\gamma}$ domain such that 
\begin{equation*} 
B_{3/4} \cap U \subset \tilde U \subset B_1 \cap U \quad \mbox{and} \quad B_{3/4} \cap \partial U = B_{3/4} \cap \partial \tilde U .
\end{equation*}

\smallskip

\emph{Step 1.}
It is convenient to decompose $u^\ep = u_1^\ep + u_2^\ep$, where $u_1^\ep \in H_0^1(\tilde U)$ has zero boundary values and solves $-\nabla \cdot \left(\a^\ep \nabla u_1^\ep  \right) = \nabla \cdot \mathbf{F}$ in $\tilde U$, whereas $u_2^\ep \in u^\ep + H_0^1(\tilde U)$ solves the homogeneous problem $-\nabla \cdot \left(\a^\ep \nabla u_2^\ep  \right) = 0 $ with boundary values given by $u^\ep$. Extend $u_1^\ep$ to be zero outside of $\tilde U$. Due to the subadditivity of the maximal function, we have
\begin{multline*} 
\left\| \mathsf{M}_{\ep \X(\frac \cdot\ep)} \left( \indc_{\tilde U} |\nabla u^\ep|^2 \right) \right\|_{L^{\frac p2}(B_{1/2} \cap U)}  
\\ 
\leq  \left\| \mathsf{M}_{\ep \X(\frac \cdot\ep)} \left( \indc_{\tilde U} |\nabla u_1^\ep|^2 \right) \right\|_{L^{\frac p2}(B_{1/2} \cap U)}  
+ \left\| \mathsf{M}_{\ep \X(\frac \cdot\ep)} \left( \indc_{\tilde U} |\nabla u_2^\ep|^2 \right) \right\|_{L^{\frac p2}(B_{1/2}\cap U)}. 
\end{multline*}
We prove that the two maximal functions on the right can be controlled by the right-hand side of~\eqref{e.CZmaxest.global}. The estimate~\eqref{e.CZmaxest.global} then follows easily. We treat the terms separately in the following two steps.

\smallskip

\emph{Step 2.} As in the local case, for a given $x \in  \tilde U$, we define 
\begin{equation*} 
f(x) := \max_{ r'  \in \left[ \ep \X_s \left(\frac x\ep\right) \wedge 1, 1 \right]} \left( \frac{1}{ r' } \left\| u_1^\ep - (u_1^\ep)_{B_{ r' }(x)} \right\|_{\underline{L}^2(B_{r'} (x) )} \right) .
\end{equation*}
Recall that  $u_1^\ep$ is extended to be zero outside of $\tilde U$. 
Then the Caccioppoli estimate gives the pointwise bound
\begin{equation} \label{e.max vs f.global}
\left( \mathsf{M}_{\ep \X(\frac \cdot\ep)} ( \indc_{\tilde U} |\nabla u_1^\ep|^2 )  \right) (x) \leq C f^2(x) + C  \left( \mathsf{M}_{\ep \X(\frac \cdot\ep)} \left(\indc_{B_1} \left| \mathbf{F}\right|^2 \right) \right) (x). 
\end{equation}
Note here that the global Caccioppoli estimate, Lemma~\ref{l.Caccioppoli appendix glob}, is applicable, since if $B_{ r' }(x)$ intersects $\partial \tilde U$, then we simply estimate 
\begin{equation*} 
 \frac{1}{ r' }  \left\| u_1^\ep\right\|_{\underline{L}^2(B_{r'} (x))}  \leq  \frac{C}{ r' }  \left\| u_1^\ep\right\|_{\underline{L}^2(B_{2r'} (x))}  \leq  \frac{C}{ r' }  \left\| u_1^\ep - (u_1^\ep)_{B_{ 2r' }(x)}\right\|_{\underline{L}^2(B_{2r'} (x))} \leq C f(x).
\end{equation*}
We again set
\begin{equation*} 
 g(x) := \left( \mathsf{M}_{\ep \X(\frac \cdot\ep)} ( \indc_{\tilde U} \left| \mathbf{F} \right|^2) \right)^\frac12  (x).
\end{equation*}
In view of~\eqref{e.max vs f.global}, we need to establish $L^p$-integrability for $f$. We will apply Lemma~\ref{l.CZ-4realz.global}, and for this, fix $x \in B_{1/2}$ and $r \in \left [\ep \X(\frac x\ep) \wedge 1 ,1 \right]$, and define 
\begin{equation*} 
f_{x,r}(y) := \max_{ r'  \in \left[ \ep \X_s \left(\frac y\ep\right) \wedge  1  ,  1 \right]} \left( \frac{1}{ r' } \left\| v_r^\ep -  (v_r^\ep)_{B_{ r' }(y)}\right\|_{\underline{L}^2(B_{r'} (y))} \right) , 
\end{equation*}
where $v_r^\ep \in u_1^\ep + H_0^1(B_{2r}(x)  \cap \tilde U)$ solves $\nabla\cdot(\a^\ep \nabla v_r^\ep) = 0$ in $B_{2r}(x) \cap \tilde U$, $v_r^\ep$ is extended to be $u_1^\ep$ outside of $B_{2r}(x) \cap \tilde U$. With these definitions, it is analogous to the proof of Lemma~\ref{l.CZ-4realz}, using this time Theorem~\ref{t.Lipschitz.boundary}, to check assumptions~\eqref{e.CZ-v1.global} and~\eqref{e.CZ-v2.global} of Lemma~\ref{l.CZ-4realz.global}. We obtain 
\begin{equation*} 
\left\| f \right\|_{L^p(B_{1/2}) }  \leq C \left(  \left\| g   \right\|_{L^p(B_{1}) } +  \left\| f \right\|_{L^2(B_1) }  \right),
\end{equation*}
and this together with~\eqref{e.max vs f.global} gives the desired estimate for the first maximal function. 

\smallskip

\emph{Step 3.}
The estimate for $u_2^\ep$ follows easily from~\eqref{e.Lipschitz2} and~\eqref{e.Lipschitz.global.grad}. Indeed, since $u_2^\ep$ has zero boundary values on $\partial \tilde U$, we have, for $x\in  \tilde U$, that 
\begin{equation*} 
\max_{ {r'} \in [\ep \X_s(x) \wedge 1, 1]} \left\|\nabla u_2^\ep\right\|_{\underline{L}^2(B_{ {r'}}(x) \cap  \tilde U)} 
\leq 
C \left\|\nabla u_2^\ep\right\|_{\underline{L}^2(\tilde U)} .
\end{equation*}
This readily implies that
\begin{equation*} 
\left\| \mathsf{M}_{\ep \X(\frac \cdot\ep)} \left( \indc_{ \tilde U} |\nabla u_2^\ep|^2 \right) \right\|_{L^{\infty}(\tilde U)} 
\leq 
C \left\|\nabla u_2^\ep\right\|_{\underline{L}^2( \tilde U)} .
\end{equation*}
 Finally, since $u_2^\ep$ agrees with $u^\ep$ on $\partial \tilde  U$, we have that 
 \begin{equation*} 
\left\|\nabla u_2^\ep\right\|_{\underline{L}^2(\tilde U)}  \leq C \left( \left\|\nabla u^\ep\right\|_{\underline{L}^2( \tilde  U)}  + \left\| F \right\|_{\underline{L}^2( \tilde U)} \right),
\end{equation*}
which can in turn be estimated with the right-hand side in~\eqref{e.CZmaxest.global}. The proof is now complete. 
\end{proof}

\begin{proof}[Proof of Theorem~\ref{t.CZ.global}]
Using Lemma~\ref{l.thetheoremreally.global}, the proof is completely analogous to the one of Theorem~\ref{t.CZ}. We omit the details.
\end{proof}
\index{Calder\'on-Zygmund estimate!global|)}

\section{\texorpdfstring{$W^{1,p}$}{W1,p}-type  estimates for the two-scale expansion error}
\label{s.twoscale.p}
\index{two-scale expansion|(}
\index{Dirichlet problem}

As an application the $L^p$-type gradient bounds proved in the previous sections, we are in a position to improve the quantitative estimates on the two-scale expansion obtained in Chapter~\ref{c.twoscale} from $L^2$ to $L^p$, for $p\in (2,\infty)$. The main result of this section is summarized in the following theorem. As usual, we set~$\a^\ep:=\a\left(\frac\cdot\ep\right)$.  

\begin{theorem}[Two-scale expansion, $L^p$ estimates]
\label{t.twoscale.p}
Fix~$s\in (0,2)$, $\alpha,\gamma \in (0,1]$, $p\in [2,\infty)$, $q\in (p,\infty]$, a bounded $C^{1,\gamma}$ domain $U\subseteq\Rd$ and $\ep\in \left(0,\tfrac 12\right]$. Let $\X$ be as in Theorem~\ref{t.CZ.global}.
There exist $\delta (d,\Lambda)>0$, $C(s,\alpha,p,q,U,d,\Lambda)<\infty$ and a random variable $\Y_{\ep,p,q}$ satisfying
\begin{equation} 
\label{e.XepSI.p}
\Y_{\ep,p,q} 
\leq 
\left\{ 
\begin{aligned}
& \O_s\left( C \ep^\alpha \left| \log \ep \right|^{\frac12} \right) & \mbox{if} & \ d=2, \ \alpha \in \left( 0,\tfrac1p \right], \\
& \O_s\left( C \ep^{\frac1p} \left| \log \ep \right|^{\frac1{2p}} \right) & \mbox{if} & \ d=2, \ \alpha \in \left( \tfrac1p,\infty \right), \\
& \O_{2+\delta} \left( C  \ep^{\alpha\wedge \frac1p} \right) & \mbox{if} & \ d>2, \ \alpha \in (0,\infty),
\end{aligned}
\right.
\end{equation}
such that the following holds: 
for every $u \in W^{1+\alpha,q}(\Rd)$, if we let $w^\ep\in H^1(U)$ be defined by~\eqref{e.wepdef} and $u^\ep\in H^1(U)$ be the solution of the Dirichlet problem
\begin{equation}
\label{e.twoscaleDP.p}
\left\{
\begin{aligned}
& -\nabla \cdot \left( \a^\ep\nabla u^\ep \right) = -\nabla \cdot \left( \ahom \nabla u \right) & \mbox{in} & \ U, \\
& u^\ep = u & \mbox{on} & \ \partial U,
\end{aligned}
\right.
\end{equation}
then we have the estimate
\begin{equation} 
\label{e.twoscale.p}
\left\| \mathsf{M}_{\ep \X\left( \frac{\cdot}{\ep}\right) } \left| \nabla \left( u^\ep - w^\ep \right) \right|^2  \right\|_{L^{p/2}(U)}^{\frac12}
\leq 
\Y_{\ep,p,q} \left\| u \right\|_{W^{1+\alpha,q}(\Rd)}.
\end{equation}
\end{theorem}

We first present an improvement of Theorem~\ref{t.twoscaleplugveng} by upgrading the $H^{-1}$ norm on the left side of~\eqref{e.twoscaleplugH2} to a $W^{-1,p}$-type norm. 

\begin{proposition}[{$L^p$-type estimates for the two-scale expansion}]
\label{p.twoscaleplugveng.p}
Fix $s\in (0,2)$, $\alpha \in (0,1]$, $p\in (2,\infty)$, $q\in (p,\infty]$ and $\ep\in \left(0,\tfrac 12\right]$. 
There exist an exponent $\delta (d,\Lambda)>0$, a constant $C(s,\alpha,p,q,d,\Lambda)<\infty$ and a random variable $\mathcal{Z}_{\ep,p,q}$ satisfying
\begin{equation} 
\label{e.twoscaleXepbound2.p}
\mathcal{Z}_{\ep,p,q} \leq 
\left\{ 
\begin{aligned}
& \O_s\left( C  \left| \log \ep \right|^{\frac12} \right) & \mbox{if} & \ d=2,\\
& \O_{2+\delta} \left( C \right) & \mbox{if} & \ d>2,
\end{aligned}
\right.
\end{equation}
such that the following holds: for every~$u\in W^{1+\alpha,q}(\Rd)$, if we define~${w}^\ep\in H_{\mathrm{loc}}^1(\Rd)$ by~\eqref{e.wepdef}, then there exists $\mathbf{F}^\ep \in L^2_{\mathrm{loc}}(\Rd;\Rd)$ satisfying 
\begin{equation} 
\label{e.Fepeq}
\nabla \cdot \mathbf{F}^\ep 
= 
\nabla \cdot \left( 
\a^\ep \nabla {w}^\ep
-\ahom \nabla u 
\right) 
\end{equation}
as well as the estimate
\begin{equation} 
\label{e.twoscaleplugH2.p}
\left\| \left| \mathbf{F}^\ep \right|^2 \ast \zeta_{\ep} \right\|_{L^{\frac p2}(B_2)}^{\frac12}
\leq
\ep^\alpha \mathcal{Z}_{\ep,p,q} \left\| u \right\|_{W^{1+\alpha,q} (\Rd)}.
\end{equation}
\end{proposition}

The statement of Proposition~\ref{p.twoscaleplugveng.p} should indeed be compared to that of Theorem~\ref{t.twoscaleplugveng} since~\eqref{e.Fepeq} implies   
\begin{equation*} \label{}
\left\| \left(\nabla \cdot \left( \a^\ep \nabla {w}^\ep \right)- \nabla \cdot \ahom \nabla u  \right)  \ast \zeta_{\ep}  \right\|_{W^{-1,p}(B_2)}
\leq
\left\| \left| \mathbf{F}^\ep \right|^2 \ast \zeta_{\ep} \right\|_{L^{\frac p2}(B_2)}^{\frac12},
\end{equation*}
and therefore the bound~\eqref{e.twoscaleplugH2.p} gives us a $W^{-1,p}$-type estimate which generalizes~\eqref{e.twoscaleplugH2}, after convolution with $\zeta_\ep$. However, splitting this estimate into~\eqref{e.Fepeq} and~\eqref{e.twoscaleplugH2.p} gives a slightly stronger statement which is required in order for Proposition~\ref{p.twoscaleplugveng.p} to be used in combination with the Calder\'on-Zygmund-type estimate in Theorem~\ref{t.CZ.global}.

\begin{proof}[{Proof of Proposition~\ref{p.twoscaleplugveng.p}}]
%
According to Lemma~\ref{l.letswritetheflux}, we have
\begin{equation}
\label{e.divergencesagree}
\nabla \cdot \left( \a^\ep  \nabla w^\ep 
- \ahom \nabla u \right) 
= \nabla \cdot \mathbf{F}^\ep,
\end{equation}
where $\mathbf{F}^\ep$ is the vector field defined in~\eqref{e.F_i^ep}, $\phi_e^\ep$ and $\mathbf{S}_e^\ep$ are defined in~\eqref{e.phiep.def} and~\eqref{e.Sep.def}, respectively, and $\mathbf{S}_e$ is the flux corrector defined in~\eqref{e.yesfluxcorrector}.

\smallskip

In view of the definition of~$\mathbf{F}^\ep$ in~\eqref{e.F_i^ep}, we have 
\begin{align*} 
\left\| \zeta_{\ep} \ast \left( \left| \mathbf{F}^\ep \right|^2  \right)   \right\|_{L^{\frac p2}(B_2)} 
& \leq C \left\|\zeta_\ep \ast \nabla u - \nabla u  \right\|_{L^{p}(B_2)}^2
\\ & \qquad
+ C \ep^2 
\sum_{k=1}^d \left\| \zeta_{\ep} \ast \left( \left( \left| \mathbf{S}_{e_k}^\ep \right|^2  + \left|  \phi_{e_k}^\ep  \right|^2  \right) \left|\nabla \left(  \zeta_\ep \ast \nabla u \right) \right|^2 \right)   \right\|_{L^{\frac p2}(B_2)} .
\end{align*}
The first term on the right side can be bounded using Lemma~\ref{l.convolution2}: we have
\begin{equation*} 
 \left\|\zeta_\ep \ast \nabla u - \nabla u  \right\|_{L^{p}(B_2)} \leq C \ep^{\alpha} \left\| \nabla u \right\|_{W^{\alpha,p}(\R^d)} .
\end{equation*}
Next, applying Lemma~\ref{l.convolution0} (with $U = B_{2\ep}(z)$), we obtain, for $q>p$, 
\begin{multline*} 
 \ep \left|\zeta_{\ep} \ast \left( \left( \left| \mathbf{S}_{e_k}^\ep \right|^2 + \left| \phi_{e_k}^\ep \right|^2 \right) \left|\nabla \left( \zeta_\ep \ast \nabla u \right) \right|^2 \right)(z) \right|^{\frac 12}  
 \\ \leq C \ep^{\alpha }  \left\| \nabla u \right\|_{\underline{W}^{\alpha,q}\left(z + 3\ep \cu_{0}\right)} 
 \left\| \left| \mathbf{S}_{e_k}^\ep  \right| + \left| \phi_{e_k}^\ep  \right| \right\|_{\underline{L}^2(z + 3\ep \cu_0)}.
\end{multline*}
Hence, by H\"older's inequality,
\begin{multline*} 
\ep^2 
 \sum_{k=1}^d \left\| \zeta_{\ep} \ast \left( \left( \left| \mathbf{S}_{e_k}^\ep    \right|^2   + \left|  \phi_{e_k}^\ep  \right|^2  \right) \left|\nabla \left(  \zeta_\ep \ast \nabla u \right) \right|^2 \right)   \right\|_{L^{\frac p2}(B_2)} 
 \\ \leq
 C \ep^{2\alpha} \left\| \left\| \nabla u \right\|_{\underline{W}^{\alpha,q}\left(\cdot + 3\ep \cu_{0}\right)} \right\|_{L^q(B_2)}^2
 \left\| \left\| \left| \mathbf{S}_{e_k}^\ep \right| + \left| \phi_{e_k}^\ep \right| \right\|_{\underline{L}^2(\cdot + 3\ep \cu_0)}
 \right\|_{L^{\frac{q p}{q-p}}(B_2)}^2. 
\end{multline*}
As explained in Remark~\ref{r.explicitXep}---see~\eqref{e.Xepp.bound1}---we have 
\begin{equation*} 
 \left\| \left\| \left| \mathbf{S}_{e_k}^\ep \right| + \left| \phi_{e_k}^\ep \right| \right\|_{\underline{L}^2(\cdot + 3\ep \cu_0)}
  \right\|_{L^{\frac{q p}{q-p}}(B_3)}
 \leq 
\left\{ 
\begin{aligned}
& \O_s\left( C \left| \log \ep \right|^{\frac12} \right) & \mbox{if} & \ d=2,\\
& \O_{2+\delta} \left( C   \right) & \mbox{if} & \ d>2. 
\end{aligned}
\right.
\end{equation*}
On the other hand, by Jensen's inequality, 
\begin{equation*} 
 \left\| \left\| \nabla u \right\|_{\underline{W}^{\alpha,q}\left(\cdot + 3\ep \cu_{0}\right)} \right\|_{L^q(B_3)} \leq C  \left\| u \right\|_{W^{1+\alpha,q}(\R^d)}.  
\end{equation*}
This completes the proof of~\eqref{e.twoscaleplugH2.p} and of the proposition. 
\end{proof}

By combining Proposition~\ref{p.twoscaleplugveng.p} and the $L^p$-type estimates of Theorem~\ref{t.CZ.global}, we can upgrade the estimates of Theorem~\ref{t.twoscale} for the two-scale expansion of the Dirichlet problem from $H^1$-type bounds to the $W^{1,p}$-type bounds of Theorem~\ref{t.twoscale.p}. The argument begins by formulating an appropriate version of Lemma~\ref{l.bndr v}. 

\begin{lemma} 
\label{l.bndr v.CZ}
Let $s \in (0,2)$, $p \in (2,\infty)$, $q \in (p,\infty)$, $\alpha \in (0,\infty)$,  and 
\begin{equation}  
\label{e.beta bndr v_ep.CZ}
\beta \in \left(0,\tfrac 1p \right] \cap \left( 0 ,  \tfrac1p  + \al - \tfrac1q   \right).
\end{equation}
There exist $\delta(d,\Lambda) > 0$, $C(s,p,q,\alpha,\beta,U, d,\Lambda)<\infty$ and, for every $\eps \in \Ll( 0,\frac 1 2 \Rr]$, a random variable $\mathcal{Z}_\ep$ satisfying
\begin{equation}
\label{e.Xeptrippybounds.CZ}
\mathcal{Z}_\ep \leq \left\{ 
\begin{aligned}
& \O_s\left( C \ep^{\frac 1p} \left| \log \ep \right|^{\frac1{2p}} \right) & \mbox{if} & \ d=2  \mbox{ and } \alpha > \tfrac 1q, \\
& \O_s\left( C \ep^{\beta} \right) & \mbox{if} & \ d=2  \mbox{ and } \alpha \leq \tfrac 1q,  \\  
& \O_{2+\delta} \left( C \ep^\beta \right) & \mbox{if} & \ d>2,
\end{aligned}
\right.
\end{equation}
such that for every $u \in W^{1+\alpha,q}(\Rd)$ and $T_\eps$ as defined in \eqref{e.T_ep},
\begin{equation}
\label{e.Tepboundzz.CZ}
\left\| \zeta_\ep \ast \left| \nabla  T_\ep \right|^2  \right\|_{L^{\frac p2}(U)}^{\frac 12} \leq \mathcal{Z}_\ep \left\|  u  \right\|_{W^{1+\alpha,q}(\R^d)}.
\end{equation}
\end{lemma}

\begin{proof}
As usual we let $\phi_e^\ep$ and $\mathbf{S}_e^\ep$ be defined as in~\eqref{e.phiep.def} and~\eqref{e.Sep.def}, respectively. Recall that $T^\ep$ is defined by the formula 
\begin{equation} 
\label{e.T_ep.CZ}
T^\ep(x) =   \left(\indc_{\R^d \setminus U_{2R(\ep)}} \ast \zeta_{R(\ep)} \right)(x) \sum_{k=1}^d \ep \phi_{e_k}^\ep (x) \partial_k \left( u \ast \zeta_\ep \right)(x).
\end{equation}
Taking the gradient gives us 
\begin{align} 
\label{e.nabla T_ep.CZ}\notag
\left| \nabla T^\ep \right| 
& 
\leq 
\indc_{\R^d \setminus U_{3R(\ep)}} 
\sum_{k=1}^d \left( \frac{C\ep}{R(\ep)} \left| \phi_{e_k}^\ep  \right| + \left| \nabla \phi_{e_k}\left( \tfrac \cdot\ep \right) \right| \right) \left| \partial_k \left( u \ast \zeta_\ep \right) \right| 
\\  & \qquad 
+ \indc_{\R^d \setminus U_{3R(\ep)}} \sum_{k=1}^d \ep \left| \phi_{e_k}^\ep \right| \left|\nabla \partial_k \left( u \ast \zeta_\ep \right) \right|.
\end{align}
To treat the first term on the right, 
Lemma~\ref{l.convolution0} (with $U = B_{\ep}(z)$) yields that 
\begin{multline*} 
 \zeta_\ep \ast \left( \sum_{k=1}^d
 \left( \frac{C\ep}{R(\ep)} \left| \phi_{e_k}^\ep \right| + \left| \nabla \phi_{e_k}\left(\tfrac\cdot\ep\right) \right| \right) \left| \partial_k \left( u \ast \zeta_\ep \right) \right|  \right)^2 (z)
 \\ \leq 
C \sum_{k=1}^d \left(  \frac{ \ep}{R(\ep)} \left\| \phi_{e_k}^\ep  \right\|_{\underline{L}^2(z + 3\ep \cu_0)}    + \left\| \nabla \phi_{e_k}\left( \tfrac \cdot\ep \right)\right\|_{\underline{L}^2(z + 3\ep \cu_0)}   \right)^2 \left\| \nabla u \right\|_{\underline{L}^2(z + 3\ep \cu_0)}^2.
\end{multline*}
Denote
\begin{equation*} 
g(z) := \sum_{k=1}^d \left( \frac{ \ep}{R(\ep)} \left\| \phi_{e_k}^\ep \right\|_{\underline{L}^2(z + 3\ep \cu_0)}    + \left\| \nabla \phi_{e_k}\left(\tfrac\cdot\ep\right)  \right\|_{\underline{L}^2(z + 3\ep \cu_0)} \right).
\end{equation*}
By the H\"older inequality, we get
\begin{align*} 
\left\| g  \left\| \nabla u \right\|_{\underline{L}^2(\cdot + 3\ep \cu_0)}  \right\|_{L^p\left(\R^d \setminus U_{2R(\ep) } \right)}  
\leq C \left\| g \right\|_{L^{\frac{qp}{q-p}}\left(\R^d \setminus U_{3R(\ep)} \right)} \left\| \nabla u\right\|_{L^{q}\left(\R^d \setminus U_{5R(\ep)} \right)} .
\end{align*}
Lemma~\ref{l.Rellichtype} implies that, for $\beta$ as in~\eqref{e.beta bndr v_ep.CZ},
\begin{equation*} 
\left\| \nabla u\right\|_{L^{q}\left(\R^d \setminus U_{5R(\ep)} \right)}   \leq  C R(\ep)^{\beta - \frac1p + \frac1q } \left\| u \right\|_{W^{1+\alpha,q}(\R^d)}.
\end{equation*}
Moreover, defining 
\begin{align*} 
 \hat \X_\ep^1 & := 
 \sum_{k=1}^d   \left\| \left\| \phi_{e_k}^\ep  \right\|_{\underline{L}^2(\cdot+3\ep \cu_0)} \right\|_{\underline{L}^{\frac{q-p}{qp} } \left( \R^d \setminus U_{3R(\ep)}  \right) }
 \\
  \hat \X_\ep^2 & := 
 \sum_{k=1}^d   \left\| \left\|\nabla \phi_{e_k}\left( \tfrac \cdot\ep \right)  \right\|_{\underline{L}^2(\cdot+3\ep \cu_0)}^{\frac{qp}{q-p}}\right\|_{\underline{L}^{\frac{q-p}{qp} } \left( \R^d \setminus U_{3R(\ep)}  \right) },
\end{align*}
we have that 
\begin{equation*} 
\left\| g \right\|_{L^{\frac{qp}{q-p}}\left(\R^d \setminus U_{2R(\ep)} \right)} \leq C R(\ep)^{\frac{1}{p} - \frac{1}{q}} \left(  \frac{ \ep}{R(\ep)}  \hat \X_\ep^1  +  \hat \X_\ep^2    \right).
\end{equation*}
Summarizing, we arrive at the inequality 
\begin{multline} \label{e.nabla T_ep est1.CZ}
\left\| \zeta_\ep \ast \left( \sum_{k=1}^d
\left( \frac{C\ep}{R(\ep)} \left| \phi_{e_k}^\ep  \right| + \left| \nabla \phi_{e_k}\left( \tfrac \cdot\ep \right) \right| \right) \left| \partial_k \left( u \ast \zeta_\ep \right) \right| \right)^2 \right\|_{L^{\frac p2} \left(\R^d \setminus U_{2R(\ep) } \right)}^{\frac 12} 
\\ \leq C R(\ep)^{\beta} \left(  \frac{ \ep}{R(\ep)}  \hat \X_\ep^1  +  \hat \X_\ep^2    \right) \left\| u \right\|_{W^{1+\alpha,q}(\R^d)}.
\end{multline}
Notice that, as in the proof of Lemma~\ref{l.bndr v}, we have
\begin{equation} \label{e.tilde X1 and tilde X2.CZ}
\hat \X_\ep^1 
\leq 
\left\{ 
\begin{aligned}
& \O_s\left( C  \left| \log \ep \right|^{\frac12} \right) & \mbox{if} & \ d=2,\\
& \O_{2+\delta} \left( C \right) & \mbox{if} & \ d>2,
\end{aligned}
\right. 
\qquad \mbox{and} \qquad 
\hat \X_\ep^2 \leq \O_{2+\delta} \left( C \right).
\end{equation}

\smallskip

We now turn our focus on the last term in~\eqref{e.nabla T_ep.CZ}. Details are very similar to Step 2 of the proof of Lemma~\ref{l.bndr v} and the above reasoning, so we only provide a sketch of the proof. We have that 
\begin{equation*} 
 \ep \left|\zeta_{\ep} \ast \left( \left| \phi_{e_k}^\ep \right|^2  \left|\nabla \left( \zeta_\ep \ast \nabla u \right) \right|^2 \right)(z) \right|^{\frac 12}  
 \\ \leq C \ep^{\alpha }  \left\| \nabla u \right\|_{\underline{W}^{\alpha,q}\left(z + 3\ep \cu_{0}\right)}  \left\| \phi_{e_k}^\ep \right\|_{\underline{L}^2(z + 3\ep \cu_0)},
\end{equation*}  
and hence, as above, we get
\begin{equation*} 
 \ep \left\| \zeta_{\ep} \ast \left( \left| \phi_{e_k}^\ep \right|^2  \left| \nabla \left( \zeta_\ep \ast \nabla u \right) \right|^2 \right)  \right\|_{L^{\frac p2} \left(\R^d \setminus U_{2R(\ep) } \right)}  \leq C \ep^{\alpha}  R(\ep)^{\frac{1}{p} - \frac{1}{q}} \hat \X_\ep^1 \left\| \nabla u \right\|_{W^{\alpha,q}(\R^d)}.
\end{equation*}
Connecting the above estimate with~\eqref{e.nabla T_ep.CZ} and~\eqref{e.nabla T_ep est1.CZ} leads to \eqref{e.Tepboundzz.CZ}, similarly to the proof of Lemma~\ref{l.bndr v}, and thus completes the proof. 
\end{proof}

\begin{proof}[Proof of Theorem~\ref{t.twoscale.p}]
To correct for the boundary values, we consider, as in the previous chapter, the solution $v^\ep \in H^1(U)$ of the Dirichlet problem
\begin{equation}
\label{e.twoscaleDP.CZ2}
\left\{
\begin{aligned}
& -\nabla \cdot \left( \a^\ep \nabla v^\ep \right) =  0 & \mbox{in} & \ U, \\
& v^\ep = T^\ep & \mbox{on} & \ \partial U,
\end{aligned}
\right.
\end{equation}
which can be written equivalently as an equation for $v^\ep - T^\ep$:
\begin{equation}
\label{e.twoscaleDP.CZ2.vTep}
\left\{
\begin{aligned}
& -\nabla \cdot \left( \a^\ep \nabla (v^\ep-T^\ep) \right) =   \nabla \cdot \left( \a^\ep \nabla T^\ep \right) & \mbox{in} & \ U, \\
& v^\ep-T^\ep =  0 & \mbox{on} & \ \partial U.
\end{aligned}
\right.
\end{equation} 
Subtracting~\eqref{e.twoscaleDP.CZ2} from~\eqref{e.twoscaleDP.p} and using~\eqref{e.Fepeq}, we also get an equation for the function $u^\ep-w^\ep+v^\ep$:
\begin{equation} \label{e.twoscaleDP.CZ3}
\left\{
\begin{aligned}
&-\nabla \cdot \left( \a^\ep \nabla (u^\ep - w^\ep + v^\ep) \right) = \nabla \cdot \mathbf{F}^\ep & \mbox{in}  & \ U, \\
& u^\ep - w^\ep + v^\ep= 0  & \mbox{on} & \ \partial U,
\end{aligned}
\right.
\end{equation}
where $\mathbf{F}^\ep$ is defined in~\eqref{e.F_i^ep} and estimated in Proposition~\ref{p.twoscaleplugveng.p}. Applying the global Calder\'on-Zygmund estimate, Theorem~\ref{t.CZ.global}, to the previous two displays, we obtain  
\begin{multline}
\label{e.firsttermwhip} 
\left\| \mathsf{M}_{\ep \X\left( \frac{\cdot}{\ep}\right) } \left| \nabla \left( u^\ep - w^\ep + v^\ep \right) \right|^2  \right\|_{L^{p/2}(U)} 
\\
\leq 
C \left\| \nabla \left( u^\ep - w^\ep + v^\ep \right) \right\|_{L^2( U)}^2
+  C \left\| \zeta_{\ep \X\left( \frac{\cdot}{\ep}\right)} 
\ast \left( \indc_{U} \left| \mathbf{F}^\ep \right|^2  
\right) \right\|_{L^{p/2}(U)}
\end{multline}
and
\begin{multline}
\label{e.secondtermwhip}
\left\| \mathsf{M}_{\ep \X\left( \frac{\cdot}{\ep}\right) } \left| \nabla \left(v^\ep - T^\ep\right) \right|^2  \right\|_{L^{p/2}(U)} 
\\
\leq 
C \left\| \nabla \left(v^\ep - T^\ep\right) \right\|_{L^2( U)}^2
 +  C \left\| \zeta_{\ep \X\left( \frac{\cdot}{\ep}\right)} \ast \left( \indc_{U} \left| \nabla T_\ep \right|^2  \right)   \right\|_{L^{p/2}(U)}.
\end{multline}
The first terms on the right sides of~\eqref{e.firsttermwhip} and~\eqref{e.secondtermwhip} are controlled by testing~\eqref{e.twoscaleDP.CZ2.vTep} with $v^\ep -T^\ep$ and~\eqref{e.twoscaleDP.CZ3} with $u^\ep-w^\ep+v^\ep$. We get
\begin{equation*} \label{}
\left\| \nabla \left( u^\ep - w^\ep + v^\ep \right) \right\|_{L^2( U)}
\leq 
C \left\|  \mathbf{F}^\ep \right\|_{L^2( U)}
\end{equation*}
and
\begin{equation*} \label{}
\left\| \nabla \left(v^\ep - T^\ep\right) \right\|_{L^2( U)}
\leq
C\left\| \nabla T^\ep \right\|_{L^2(U)}. 
\end{equation*}
Furthermore,  it is easy to see from the properties of~$\X$ that 
\begin{equation*} 
\left\| \zeta_{\ep \X\left( \frac{\cdot}{\ep}\right)} \ast \left( \indc_{U} \left| \mathbf{F}^\ep \right|^2  \right)   \right\|_{L^{p/2}(U)}  \leq C \left\| \zeta_{\ep} \ast \left( \indc_{U} \left| \mathbf{F}^\ep \right|^2  \right)   \right\|_{L^{p/2}(U)} ,
\end{equation*}
and similarly for $\nabla T^\ep$. Combining these, we therefore obtain  
\begin{equation}
\label{e.firsttermwhip2} 
\left\| \mathsf{M}_{\ep \X\left( \frac{\cdot}{\ep}\right) } \left| \nabla \left( u^\ep - w^\ep + v^\ep \right) \right|^2  \right\|_{L^{p/2}(U)} 
\leq 
C \left\| \zeta_{\ep \X\left( \frac{\cdot}{\ep}\right)} \ast \left( \indc_{U} \left| \mathbf{F}^\ep \right|^2  \right) \right\|_{L^{p/2}(U)}
\end{equation}
and
\begin{equation}
\label{e.secondtermwhip2}
\left\| \mathsf{M}_{\ep \X\left( \frac{\cdot}{\ep}\right) } \left| \nabla v^\ep \right|^2  \right\|_{L^{p/2}(U)} 
\leq 
C \left\| \mathsf{M}_{\ep \X\left( \frac{\cdot}{\ep}\right) } 
 \left| \nabla T^\ep \right|^2    \right\|_{L^{p/2}(U)} .
\end{equation}
Combining these yields
\begin{multline}
\label{e.secondtermwhip3}
\left\| \mathsf{M}_{\ep \X\left( \frac{\cdot}{\ep}\right) } \left| \nabla \left( u^\ep - w^\ep \right) \right|^2  \right\|_{L^{p/2}(U)} 
\\
\leq 
C\left\| \zeta_{\ep \X\left( \frac{\cdot}{\ep}\right)} \ast \left( \indc_{U} \left| \mathbf{F}^\ep \right|^2  \right) \right\|_{L^{p/2}(U)}
+
 \left\| \mathsf{M}_{\ep \X\left( \frac{\cdot}{\ep}\right) } 
 \left| \nabla T^\ep \right|^2    \right\|_{L^{p/2}(U)}.
\end{multline}
We may now apply Proposition~\ref{p.twoscaleplugveng.p} and Lemma~\ref{l.bndr v.CZ} to obtain~\eqref{e.twoscale.p}. This completes the proof.
\end{proof}

%
%
%
%
%
%
%

\index{two-scale expansion|)}
\section*{Notes and references}

Calder\'on-Zygmund-type gradient $L^p$ estimates were obtained in the context of periodic homogenization by Avellaneda and Lin~\cite{AL1}, who derived them from singular integral estimates for the elliptic Green function. Such estimates were first proved in the stochastic case in~\cite{ADan,DGO} using similar arguments as the ones here. Lemma~\ref{l.CZ-4realz} is an efficient formalization of the method introduced by Caffarelli and Peral~\cite{CP} for obtaining~$L^p$ estimates by approximating on small scales in $L^q$ for some $q>p$. Proposition~\ref{p.CZsmallcontrast} is essentially due to~\cite{CP}. The optimal estimates for the two-scale expansion error proved in Section~\ref{s.twoscale.p}, like most of the results of Chapter~\ref{c.twoscale}, are presented here for the first time.



\chapter{Estimates for parabolic problems}
\label{c.parabolic}

We present quantitative estimates for the homogenization of the parabolic equation
\begin{equation} 
\label{e.pde.parab}
\partial_t u - \nabla \cdot \a(x) \nabla u = 0 \quad \mbox{in} \ I\times U \subseteq \R\times \Rd.
\end{equation}
The coefficients~$\a(x)$ are assumed to depend only on the spatial variable~$x$ rather than~$(t,x)$.\footnote{Quantitative homogenization results for parabolic equations with \emph{space-time} random coefficients can also be obtained from the ideas presented in this book: see~\cite{ABM}.} 
The main purpose of this chapter is to illustrate that the parabolic equation~\eqref{e.pde.parab} can be treated satisfactorily using the elliptic estimates we have already obtained in earlier chapters. In particular, we present error estimates for general  Cauchy-Dirichlet problems in bounded domains, two-scale expansion estimates, a parabolic large-scale regularity theory. We conclude, in the last two sections, with $L^\infty$--type estimates for the homogenization error and the two-scale expansion error for both the parabolic and elliptic Green functions. The statements of these estimates are given below in Theorem~\ref{t.PGF.basecase} and Corollary~\ref{c.GtoGbar0}. Like the estimates in Chapter~\ref{c.two}, these estimates are suboptimal in the scaling of the error but optimal in stochastic integrability (i.e., the scaling of the error is given by a small exponent $\alpha>0$ and the stochastic integrability is $\O_{d-}$-type). In the next chapter, we present complementary estimates which are optimal in the scaling of the error and consistent with the bounds on the first-order correctors proved in Chapter~\ref{c.A1}. See Theorem~\ref{t.PtoPbar}.

\smallskip

The arguments in this chapter are deterministic, with the only probabilistic inputs being the regularity theory and the sub-optimal estimates on the first-order correctors presented in Chapter~\ref{c.regularity}. In particular, the optimal bounds on the first-order correctors proved in Chapter~\ref{c.A1} are not needed here.

\smallskip

We begin in the first section by reviewing the basic setup for parabolic equations, defining the parabolic Sobolev spaces and giving the parabolic version of some basic estimates such as the Caccioppoli inequality. In Section~\ref{s.parab.homog}, we prove the quantitative homogenization results. The large-scale regularity theory is presented in Sections~\ref{s.parab.Lip} and~\ref{s.parab.reg}. We conclude with the estimates on the parabolic and elliptic Green functions in Sections~\ref{s.greenie.decay} and~\ref{s.greenie.base}.

\smallskip

In Sections~\ref{s.greenie.decay} and~\ref{s.greenie.base} we will deviate from the convention maintained in the rest of the book by allowing ourselves to use estimates which are available \emph{only for scalar equations}: namely the De Giorgi-Nash $L^\infty$ estimates and the Nash-Aronson pointwise upper bound on the parabolic Green function. Without these estimates to control the small scales, it is not straightforward to even define the Green functions. However, the reader should \emph{not} interpret this choice to mean that the arguments here will not extend to the case of elliptic systems. Indeed, for systems, these deterministic estimates can be replaced by the $C^{0,1}$-type regularity proved in Chapters~\ref{c.regularity} and~\ref{c.parabolic}, which also implies Nash-Aronson-type bounds valid for times larger than a minimal time (i.e., for $\sqrt{t} \geq \X$, where $\X$ is the minimal scale for the $C^{0,1}$-type estimate). We can also replace the $\delta_0$ distribution in the definition of the Green functions by a smooth, positive bump function with characteristic length scale of order one (which makes no essential difference since, as we have explained many times previously, homogenization results ``should not be concerned with scales smaller than the correlation length scale''). This does however complicate the argument on a technical level, and therefore for readability we present only the scalar case here.

\section{Function spaces and some basic estimates}
\label{s.parab.setup}

The purpose of this section is to introduce the basic functional analytic framework for parabolic equations. We also record some basic estimates, such as parabolic versions of the Caccioppoli and Meyers estimates. 

\smallskip

We begin with the definitions of the parabolic Sobolev spaces. For every bounded Lipschitz domain~$U\subseteq\Rd$, Banach space~$X$ and~$p\in[1,\infty)$, we denote by~$L^p(U;X)$ the Banach space of Lebesgue-measurable mappings $u:U \to X$ with (volume-normalized) norm
\begin{equation*} \label{}
\left\| u \right\|_{\underline{L}^p(U;X)} 
:=
\left( \fint_{U} \left\| u(x) \right\|_{X}^p \,dx
\right)^{\frac1p} <\infty. 
\end{equation*}
For every interval $I = (I_-,I_+) \subset \R$ and bounded Lipschitz domain $U \subset \Rd$, 
we define the \emph{parabolic boundary}~$\partial_\sqcup (I\times U)$ of the cylinder $I\times U$ by 
\index{parabolic boundary $\partial_\sqcup$}
\begin{equation*}
\partial_\sqcup (I\times U) 
:= 
\left( \{ I_- \} \times U \right) \cup \left( I \times \partial U \right). 
\end{equation*}
We define the function space
\index{Sobolev space!parabolic $H^1_\pa$ and $W^{1,p}_\pa$} 
\begin{equation}  
\label{e.def.X}
H^1_\pa(I\times U) := \Ll\{ u \in L^2(I;H^1(U))\,:\,\partial_t u \in L^2(I;H^{-1}(U)) \Rr\} ,
\end{equation}
which is the closure of bounded smooth functions on $I\times U$ equipped with the norm
\begin{equation}
\label{e.def.Xnorm}
\|u\|_{\underline{H}^1_\pa(I \times U)} := \|u\|_{\underline{L}^2(I;\underline H^1(U))} + \|\partial_t u\|_{\underline{L}^2(I;\underline H^{-1}(U))}.
\end{equation}
We denote by $H^1_{\pa,\sqcup}(I\times U)$ the closure in $H^1_\pa(I\times U)$ of the set of smooth functions with compact support in $(I\times U) \setminus \partial_\sqcup (I\times U)$. In other words, a function in $H^1_{\pa,\sqcup}(I\times U)$ has zero trace on the lateral boundary $I \times \partial U$ and the initial time $\{ I_-\}\times U$ but does not necessarily vanish at the final time.

\smallskip

We sometimes work with generalizations of $H^1_\pa(I\times U)$ in which the exponent of integrability $p\in (1,\infty)$ is not necessarily equal to~$2$. Therefore we introduce the function space
\begin{equation}  
\label{e.def.W1p.par}
W^{1,p}_\pa(I\times U) := 
\Ll\{ u \in L^p\left(I;W^{1,p}(U)\right)\,:\,\partial_t u \in L^p\left(I;W^{-1,p}(U)\right) \Rr\} ,
\end{equation}
equipped with the (volume-normalized) norm
\begin{equation}
\label{e.def.W1p.par.norm}
\left\|u\right\|_{\underline{W}^{1,p}_\pa\left(I \times U\right)} 
:= \left\|u\right\|_{\underline{L}^p(I;\underline{W}^{1,p}(U))} + \|\partial_t u\|_{\underline{L}^p\left(I;\underline{W}^{-1,p}(U)\right)}.
\end{equation}
Similarly to $H^1_{\pa,\sqcup}(I\times U)$, we denote by $W^{1,p}_{\pa,\sqcup}(I\times U)$ the closure in $W^{1,p}_\pa(I\times U)$ of the set of smooth functions with compact support in $(I_-,I_+] \times U$. Finally, for every parabolic cylinder $V$, we denote by $W^{1,p}_{\pa,\,\mathrm{loc}}(V)$, $H^1_{\pa,\,\mathrm{loc}}(V)$, and so forth, the functions on~$V$ which are, respectively, elements of $W^{1,p}_\pa(W)$ and $H^1_{\pa}(W)$, etc, for every subcylinder $W\subseteq V$ with $\partial_\sqcup W \subseteq V$. 

\smallskip

We give a weak interpretation of the parabolic equation
\begin{equation}
\label{e.weaksoldefparab}
\partial_t u - \nabla \cdot \a\nabla u = u^* \quad \mbox{in} \ I \times U
\end{equation}
for right-hand sides $u^*$ belonging to $L^2(I; H^{-1}(U))$. We say that $u \in H^1_\pa(I\times U)$ is a \emph{weak solution} (or just \emph{solution})  of~\eqref{e.weaksoldefparab} if
\begin{equation*} \label{}
\forall \phi \in L^2(I;H^1_0(U)), \qquad \int_{I\times U} \left( \phi \partial_t u + \nabla \phi \cdot \a\nabla u \right) 
=
\int_{I\times U} \phi u^*.
\end{equation*}
The integral expression $\int_{I\times U} \phi \partial_t u$ is a shorthand for the canonical pairing between $\phi \in L^2(I;H^1_0(U))$ and $\partial_t u\in L^2(I;H^{-1}(U))$ which extends the integral of the product of bounded smooth functions (cf. the explanation just below~\eqref{e.thedefWminus}). 

\smallskip

We next present a parabolic version of the Caccioppoli inequality. For each $r\in (0,\infty]$, we set $I_r:=(-r^2,0]$ and let $Q_r$ be the \emph{parabolic cylinder}
\begin{equation*} \label{}
Q_r := I_r \times B_r. 
\end{equation*}
We also write $Q_r(t,x):= (t,x) + Q_r$ for each $(t,x)\in\R\times\Rd$ and $r>0$.

\index{Caccioppoli inequality}

\begin{lemma}[parabolic Caccioppoli inequality]
\label{l.cacciopp.this}
There exists $C(d,\Lambda)<\infty$ such that, for every 
$u\in H^1_\pa(Q_{2r})$ and $u^* \in L^2(I_{2r};H^{-1}(B_{2r}))$ satisfying
\begin{equation} 
\label{e.cacc.pde.parab}
\partial_t u -\nabla \cdot (\a \nabla u ) 
= u^* \quad \mbox{in} \ Q_{2r}, 
\end{equation}
we have
\begin{equation} 
\label{e.cacciopp.this}
 \left\| \nabla u \right\|_{L^2(Q_r)}  
\leq 
Cr^{-1}\left\|  u \right\|_{L^2(Q_{2r})} 
+C \left\| u^* \right\|_{L^2 \left( I_{2r};H^{-1}(B_{2r}) \right) }
\end{equation}
and
\begin{equation} 
\label{e.Cacciopp.thetime}
\sup_{s\in I_{r}} 
\left\|  u(s,\cdot) \right\|_{L^2(B_r)}  
\leq C\left\|  \nabla u \right\|_{L^2(Q_{2r})} 
+ C \left\| u^* \right\|_{L^2\left(I_{2r};H^{-1}(B_{2r})\right)}.
\end{equation}
\end{lemma}
\begin{proof}
Select $\eta_r \in C^\infty_c(Q_{2r})$ satisfying 
\begin{equation*} \label{}
0\leq \eta \leq 1, \quad
\eta \equiv 1 \ \mbox{on} \ Q_r,  \quad
\left| \partial_t \eta \right| + \left| \nabla \eta \right|^2 \leq Cr^{-2}. 
\end{equation*}
Testing~\eqref{e.cacc.pde.parab} with $\phi := \eta^2_r u \in L^2(I_{2r};H^{1}_0(B_{2r}))$ gives
\begin{equation*} \label{}
\int_{Q_{2r}} \phi \left( u^* - \partial_t u \right) = \int_{Q_{2r}} \nabla \phi\cdot \a\nabla u. 
\end{equation*}
Estimating the right side, we find that 
\begin{align*} \label{}
\int_{Q_{2r}}\nabla \phi\cdot \a \nabla u 
& 
\geq 
\int_{Q_{2r}} \eta_r^2 \left| \nabla u \right|^2 -  C \int_{Q_{2r}} \eta_r \left| \nabla \eta_r \right| |u| \left| \nabla u \right|
\\ &
\geq 
\frac12 \int_{Q_{2r}} \eta_r^2 \left| \nabla u \right|^2 - C\int_{Q_{2r}} \left| \nabla \eta_r \right|^2 |u|^2 
\\ & 
\geq  
\frac12 \int_{Q_{2r}} \eta_r^2 \left| \nabla u \right|^2 - Cr^{-2} \int_{Q_{2r}}  |u|^2. 
\end{align*}
Next we estimate the left side:
\begin{align*}
\int_{Q_{2r}} \eta^2_r u \left( u^*-\partial_tu \right) 
&
\leq - \int_{Q_{2r}} \partial_t \left( \frac12 \eta^2_r u^2 \right) + \int_{Q_{2r}} \eta_r \left| \partial_t\eta_r \right| u^2 
\\ & \qquad
+ \int_{-4r^2}^0 \left\| (\eta^2_r u)(t,\cdot) \right\|_{H^1(B_{2r})} \left\| u^*(t,\cdot) \right\|_{H^{-1}(B_{2r})} \,dt 
\\ & 
\leq - \frac12 \int_{B_{2r}} \eta_r^2(0,x) u^2(0,x)\,dx + Cr^{-2} \int_{Q_{2r}} u^2 
\\ & \qquad
+C  \left\| \eta^2_r u \right\|_{L^2\left(I_{2r};H^1(B_{2r})\right)} \left\| u^* \right\|_{L^2\left(I_{2r};H^{-1}(B_{2r})\right)}.
\end{align*}
Using also that 
\begin{equation*} \label{}
\left\| \eta^2_r u \right\|_{L^2\left(I_{2r};H^1(B_{2r})\right)}
\leq Cr^{-1} \left\| u \right\|_{L^2(I_{2r}\times B_{2r})} + C \left\| \eta_r \nabla u \right\|_{L^2(I_{2r}\times B_{2r})},
\end{equation*}
we obtain 
\begin{multline*} \label{}
 \left\| \eta^2_r u \right\|_{L^2\left(I_{2r};H^1(B_{2r})\right)} \left\| u^* \right\|_{L^2\left(I_{2r};H^{-1}(B_{2r})\right)}
 \\
 \leq r^{-2} \left\| u \right\|_{L^2(I_{2r}\times B_{2r})}^2 + \frac14 \left\| \eta_r \nabla u \right\|_{L^2(I_{2r}\times B_{2r})}^2
 + C  \left\| u^* \right\|_{L^2\left(I_{2r};H^{-1}(B_{2r})\right)}^2.
\end{multline*}
Combining the above, we finally obtain
\begin{equation*} \label{}
\frac12 \int_{B_{2r}} \eta_r^2(0,x) u^2(0,x)\,dx +  \frac14 \int_{Q_{2r}} \eta_r^2 \left| \nabla u \right|^2 
\leq
Cr^{-2} \int_{Q_{2r}}  |u|^2 +  C \left\| u^* \right\|_{L^2\left(I_{2r};H^{-1}(B_{2r})\right)}^2.
\end{equation*}
This yields~\eqref{e.cacciopp.this}.

\smallskip

By repeating the above computation but using instead the test function
$\phi := \eta^2_r u \indc_{\{ t < s\}}$ for fixed $s\in I_{2r}$, and estimating the right side of the weak formulation from below differently, namely
\begin{align*} \label{}
\int_{Q_{2r}}\nabla \phi\cdot \a \nabla u
& 
\geq -C \left\| \eta_r \nabla u \right\|_{L^2(Q_{2r})}^2 -C \left\| \nabla \eta_r \nabla u \right\|_{L^2(Q_{2r})}  \left\| u  \eta_r \right\|_{L^2(Q_{2r})} 
\\ & 
\geq -C  \left\| \nabla u \right\|_{L^2(Q_{2r})}^2 - \frac1{16}r^{-2}  \int_{Q_{2r}} \eta_r^2 u^2
\\ & 
\geq -C  \left\|  \nabla u \right\|_{L^2(Q_{2r})}^2 - \frac14 \sup_{t\in I_{2r}} \int_{B_{2r}} \eta_r^2(t,x) u^2(t,x)\,dx,
\end{align*}
we get 
\begin{multline*} \label{}
\frac12 \int_{B_{2r}} \eta_r^2(s,x) u^2(s,x)\,dx 
\\
\leq C\left\|  \nabla u \right\|_{L^2(Q_{2r})}^2 +  \frac14 \sup_{t\in I_{2r}}  \int_{B_{2r}} \eta_r^2(t,x) u^2(t,x)\,dx + C \left\| u^* \right\|_{L^2\left(I_{2r};H^{-1}(B_{2r})\right)}^2.
\end{multline*}
Taking the supremum over $s\in I_{2r}$ and rearranging yields~\eqref{e.Cacciopp.thetime}. 
\end{proof}

Since the parabolic equation~\eqref{e.pde.parab} does not depend on time, we may upgrade~$L^2$ gradient bounds to~$L^\infty_tL^2_x$ gradient bounds. 

\begin{lemma}[$L^2$ gradient estimate in time slices]
\label{l.slicesvsaverages.PAR}
Fix $r>0$ and suppose that $u\in H^1_\pa(Q_{2r})$ satisfies
\begin{equation*} \label{}
\partial_tu  -\nabla \cdot \a \nabla u  = 0 \quad \mbox{in} \ Q_{2r}. 
\end{equation*}
There exists a constant~$C(d,\Lambda)<\infty$ such that 
\begin{equation} 
\label{e.gradinslices.PAR}
\sup_{t \in (-r^2,0)}   
\left\| \nabla u(t,\cdot) \right\|_{\underline{L}^2(B_r)}  
\leq   
C  \left\|\nabla u\right\|_{\underline{L}^2(Q_{2r})}.
\end{equation}
\end{lemma}
\begin{proof}
Since $\partial_t u$ is also a solution of the equation, the Caccioppoli inequality yields
\begin{equation}
\label{e.cacciop.fort}
\fint_{Q_r} \left| \nabla \partial_t u \right|^2 
\leq 
\frac{C}{r^2} \fint_{Q_{2r}} \left| \partial_t u \right|^2.
\end{equation}
(Strictly speaking, we need to justify that $\partial_t u$ belongs to $H^1_\pa(I\times U)$. One may justify this by considering difference quotients in time, obtaining a version of~\eqref{e.cacciop.fort} by the Caccioppoli inequality for these difference quotients, which then allows to pass to weak limits in the difference quotient parameter to obtain that $\partial_t u$ does indeed belong to $H^1_\pa$ and is a solution of the equation. We leave the details of this standard argument to the reader.) Differentiating gives  
\begin{equation*}
\partial_t \fint_{B_r} |\nabla u|^2 
= 
2 \fint_{B_r} \nabla u \cdot \nabla \partial_t u  ,
\end{equation*}
and thus we have, after an integration in time, that 
\begin{align*}
\sup_{t \in (-r^2,0)} \fint_{B_r} |\nabla u(t,\cdot)|^2 
& 
\leq
\fint_{Q_r} |\nabla u|^2 
+ 2 r^2 \fint_{Q_r}  \left| \nabla u \cdot \nabla \partial_t u \right|
\\ & 
\leq   2 \fint_{Q_r} |\nabla u|^2 + r^4 \fint_{Q_r}  \left| \nabla \partial_t u \right|^2.
\end{align*}
Combining this with~\eqref{e.cacciop.fort} yields 
\begin{equation} 
\label{e.gradinslices2.PAR}
\sup_{t \in (-r^2,0)} \fint_{B_r} |\nabla u(t,\cdot)|^2 
\leq   
2 \fint_{Q_r} |\nabla u|^2 + C r^2 \fint_{Q_{2r}}  \left| \partial_t u \right|^2.
\end{equation}

The next goal is hence to estimate the second term on the right side of~\eqref{e.gradinslices2.PAR}. For this, fix $s<r$ and $\eta \in C_0^\infty(B_{2s})$ such that $\eta = 1$ in $B_s$, $0 \leq \eta \leq 1$, and $|\nabla \eta| \leq \tfrac{2}{s}$. Testing the equation for $u$ with the function $\eta \partial_t u$, we obtain
\begin{equation*} 
\fint_{Q_s} \left| \partial_t u \right|^2 
\leq  2^{d+2} \fint_{Q_{2s}}  \eta \partial_t u \nabla \cdot \left( \a \nabla u\right) 
\leq 2^{d+2} \Lambda \left(\fint_{Q_{2s}} \left| \nabla (\eta \partial_t u) \right|^2\right)^{\frac12}   
\left(\fint_{Q_{2s}} \left| \nabla u \right|^2
\right)^{\frac12}  .
\end{equation*}
Therefore
\begin{equation*} 
\fint_{Q_s} \left| \partial_t u \right|^2 
\leq  
\frac{C}{s}  \left(\fint_{Q_{2s}} \left| \partial_t u \right|^2\right)^{\frac12}   \left(\fint_{Q_{2s}} \left| \nabla u \right|^2\right)^{\frac12}  .
\end{equation*}
Choosing now~$\tfrac 12 \leq \sigma' < \sigma \leq 1$ we obtain by a covering argument (with~$s = (\sigma-\sigma')\tfrac{r}{2}$ and translations of~$Q_s$) that 
\begin{equation*} 
\fint_{Q_{\sigma' 4r}} \left| \partial_t u \right|^2 
\leq  
\frac{C}{(\sigma - \sigma')^{d+3}r}  
\left(\fint_{Q_{\sigma 4r}} \left| \partial_t u \right|^2\right)^{\frac12}   \left(\fint_{Q_{4r}} \left| \nabla u \right|^2\right)^{\frac12}  .
\end{equation*}
Thus, by Young's inequality, 
\begin{equation*} 
\fint_{Q_{\sigma' 4r}} \left| \partial_t u \right|^2 
\leq  \frac12 \fint_{Q_{\sigma 4r}} \left| \partial_t u \right|^2  + \frac{C}{(\sigma - \sigma')^{2(d+3)}}     r^{-2} \fint_{Q_{4r}} \left| \nabla u \right|^2  .
\end{equation*}
Lemma~\ref{l.simpleiter} then implies that 
\begin{equation*} 
\fint_{Q_{2r}} \left| \partial_t u \right|^2 
\leq  C r^{-2} \fint_{Q_{4r}} \left| \nabla u \right|^2  .
\end{equation*}
Combining this with~\eqref{e.gradinslices2.PAR} gives
\begin{equation*} 
\sup_{t \in (-r^2,0)} \fint_{B_r} |\nabla u(t,\cdot)|^2 
\leq C \fint_{Q_{4r}} |\nabla u|^2.
\end{equation*}
We conclude the proof by a covering argument by considering~$r/4$ instead~$r$. 
\end{proof}

\section{Homogenization of the Cauchy-Dirichlet problem}
\label{s.parab.homog}

The purpose of this section is to generalize Theorem~\ref{t.DP} by estimating the homogenization error for the Cauchy-Dirichlet problem. As in that previous theorem, we are interested in an estimate which is suboptimal in the size of the error but optimal in stochastic integrability and allows for relatively rough boundary data. It is also possible (and relatively straightforward) to obtain error estimates for the Cauchy-Dirichlet problem which are sharper in the size of the error and generalize those of Chapter~\ref{c.twoscale} for the Dirichlet problem (e.g., Theorems~\ref{t.twoscale} and~\ref{t.L2EE}), but we do not give such estimates here. 

\index{Cauchy-Dirichlet problem}

\begin{theorem}
\label{t.CDP}
Fix exponents $\delta>0$, $\sigma \in(0,d)$, an interval $I\subseteq \left( -\frac12 ,0\right)$ and a Lipschitz domain $U\subseteq B_{1/2}$. There exist an exponent $\beta(\delta,d,\Lambda)>0$, a constant $C(s,I,U,\delta,d,\Lambda)<\infty$ and a random variable $\X_\sigma$ satisfying 
\begin{equation*} \label{}
\X_\sigma = \O_1(C)
\end{equation*}
such that the following holds. 
For each $\ep \in \left(0,\tfrac12\right]$, $f\in W^{1,2+\delta}_\pa(I\times U)$ if we let $u,u^\ep\in f + H^1_{\pa,\sqcup}(I\times U)$ respectively denote the solutions of the Cauchy-Dirichlet problems
\begin{equation*} \label{}
\left\{
\begin{aligned}
& \left( \partial_t - \nabla \cdot \a^\ep\nabla\right) u^\ep = 0 & \mbox{in} & \ I \times U, \\
& u^\ep = f & \mbox{on} & \ \partial_{\sqcup}(I\times U),
\end{aligned}
\right.
\end{equation*}
and
\begin{equation*} \label{}
\left\{
\begin{aligned}
& \left( \partial_t - \nabla \cdot \ahom \nabla\right) u = 0 & \mbox{in} & \ I \times U, \\
& u = f & \mbox{on} & \ \partial_{\sqcup}(I\times U),
\end{aligned}
\right.
\end{equation*}
then we have the estimate 
\begin{multline}
\label{e.CDPestimates.bb}
\left\| u^\ep - u \right\|_{L^2(I\times U)} 
+ \left\| \nabla (u^\ep - u) \right\|_{L^2(I;H^{-1}(U))}
+ \left\| \a^\ep \nabla u^\ep - \ahom \nabla u \right\|_{L^2(I;H^{-1}(U))}
\\
\leq  
C  \left\| f \right\|_{W^{1,2+\delta}_\pa(I\times U)} \left( \ep^{\beta(d-\sigma)} + \X_\sigma \ep^\sigma \right). 
\end{multline}
\end{theorem}

The proof of Theorem~\ref{t.CDP} closely follows the one of Theorem~\ref{t.DP}, most of which is formalized in Theorem~\ref{t.DP.blackbox}. We henceforth fix $I$, $U$, $\ep$, $\delta$, $f$,~$u$ and~$u^\ep$ as in the statement of the theorem. We also fix a parameter~$r\in [\ep,1)$ to be selected in the last paragraph of the proof, which represents a mesoscopic scale (it will depend on $\ep$ in such a way that $\ep \ll r \ll 1$). We set $I : =(I_-,I_+)$ and denote, for~$s>0$, 
\begin{equation*} \label{}
I_s:= \left( I_- + s^2, I_+\right]
\quad \mbox{and} \quad 
U_s:= \left\{ x\in U \,:\, \dist(x,\partial U) > s \right\}. 
\end{equation*}
We work with the modified two-scale expansion~$w^\ep$ defined by 
\begin{equation} 
\label{e.wepdef.par}
w^\ep(t,x) 
:=u(t,x) + \ep\eta(t,x) \sum_{k=1}^d  \partial_{x_k} u(t,x) \phi_{e_k}^\ep(x),
\end{equation}
where $\phi^\ep_e$ is defined in~\eqref{e.phiep.def} and $\eta \in C^\infty(\overline{ I} \times \overline{U} )$ is a cutoff function satisfying
\begin{equation} 
\label{e.etacutoff.parab}
\left\{ 
\begin{aligned}
& 0\leq \eta \leq 1, \quad \eta \equiv 1 \ \mbox{in} \ I_{2r} \times U_{2r}, \\
& \eta \equiv 0 \ \mbox{in} \ (I\times U) \setminus (I_r\times U_r), \\
& \forall k,l\in\N, \quad \left| \nabla^k \partial_t^l \eta \right| \leq C_{k+2l} r^{-(k+2l)}. 
\end{aligned}
\right.
\end{equation}
In other words, $w^\ep(t,\cdot)$ is the product of the cutoff function~$\eta$ and two-scale expansion of the function~$u(t,\cdot)$, with the time variable playing no role. Note that we suppress the dependence of~$\eta$ on the parameter~$r>0$.

\smallskip

As in the proof of Theorem~\ref{t.DP}, the strategy of the argument is, roughly speaking, to plug $w^\ep$ into the equation for $u^\ep$ and estimate the error which arises, here with respect to~$\| \cdot \|_{L^2(I;H^1(U))}$. 
This estimate is presented in Lemma~\ref{l.CDP.plugineq}, below.  As usual, the error is controlled by a random variable which measures how ``well-behaved'' the first-order correctors are: we define
\begin{equation*} \label{}
\mathcal{E}(\ep)
:= 
\sum_{k=1}^d \left( 
\left\| \nabla \phi_{e_k} \left(\tfrac \cdot\ep\right) \right\|_{H^{-1} (B_1)} 
+
\left\| \a^\ep \left( e_k+ \nabla \phi_{e_k}\left(\tfrac \cdot\ep\right) \right) - \ahom e_k \right\|_{H^{-1} (B_1)} 
\right).
\end{equation*}
We then obtain in Lemma~\ref{l.H1est.parab} an estimate on the difference between $u^\ep$ and $w^\ep$ in the (strong) norm $\| \cdot \|_{H^1_\pa(I\times U)}$, and finally on the homogenization error in Lemma~\ref{l.throwcorrectors}. Note that the argument here is completely deterministic, the only stochastic ingredient being the bounds on the random variable~$\mathcal{E}(\ep)$ we get from estimates on the first-order correctors proved in earlier chapters. 

\smallskip

We recall first some classical pointwise estimates for the homogenized solution~$u$. Up to an affine change of variables, the equation for $u$ is of course the standard heat equation. Therefore at each point $(t,x) \in I_r\times U_r$, we have, for every $k,l\in\N$, 
\begin{equation*} \label{}
\left| \nabla^k\partial_t^l u (t,x)  \right| 
\leq 
C_{k+2l} r^{-(k+2l)+1} 
\left\| \nabla u \right\|_{\underline{L}^2(Q_r(t,x))}
\leq 
C_{k+2l} r^{-(k+2l)+1 - (2+d)/2}  \left\| \nabla u \right\|_{L^2(I_r\times U_r)}.
\end{equation*}
(See~\cite[Theorem 9 in Section 2.3]{Evans} for a proof of this standard pointwise estimate.) We deduce therefore that 
\begin{equation} 
\label{e.upoint.bounds}
\left\| \nabla^k\partial_t^l u \right\|_{L^\infty(I_r\times U_r)} 
\leq C_{k+2l} r^{-(k+2l)+1-(d+2)/2} \left\| \nabla f \right\|_{L^2(I\times U)}.
\end{equation}
\index{Meyers estimate}
We also need the parabolic version of the global Meyers estimate, which states that, if $\delta(d,\Lambda)>0$ is sufficiently small, then  
\begin{equation} 
\label{e.umeyers.parab}
\left\| u \right\|_{W^{1,2+\delta}_\pa(I\times U)} 
\leq 
C \left\| f \right\|_{W^{1,2+\delta}_\pa(I\times U)}.
\end{equation}
Although we do not give the proof of the parabolic version of the Meyers estimate here, it is quite close to that of the elliptic version presented in Appendix~\ref{a.meyers} (cf. Theorem~\ref{t.Meyers appendix global}). A complete proof of the parabolic Meyers estimate can be found in~\cite[Proposition B.2]{ABM}.

\begin{lemma}
\label{l.CDP.plugineq}
There exists $C(\delta,d,\Lambda)<\infty$ such that 
\begin{equation} \label{e.CDP.plugineq}
\left\| \left( \partial_t - \nabla \cdot \a^\ep \nabla \right) w^\ep \right\|_{L^2(I ; H^{-1}( U))}
\leq 
C \left( r^{\frac{\delta}{4+2\delta}} + r^{-3-(2+d)/2} \mathcal{E}(\ep)  \right) \left\| f \right\|_{W^{1,2+\delta}_\pa(I\times U)}. 
\end{equation}
\end{lemma}
\begin{proof}
For the computations below, we note that~$\ep \nabla \phi^\ep  = \nabla \phi \left( \frac\cdot \ep \right)$.
We compute 
\begin{equation*} \label{}
\left\{
\begin{aligned}
\nabla w^\ep & = 
\eta \sum_{k=1}^d \left( e_k+\nabla \phi_{e_k}\left( \tfrac\cdot \ep \right)  \right) \partial_{x_k} u 
+ \sum_{k=1}^d \phi^\ep_{e_k} \nabla \left( \eta \partial_{x_k} u \right) 
+ (1-\eta)\nabla u,
\\
\partial_t w^\ep & = 
\partial_t u 
+ \sum_{k=1}^d \phi^\ep_{e_k} \partial_t \left( \eta \partial_{x_k} u \right).
\end{aligned}
\right.
\end{equation*}
Using the equation for $\phi^\ep_{e_k}$, we therefore obtain
\begin{align*}
\left( \partial_t - \nabla \cdot \a^\ep\nabla \right) w^\ep
& =
\partial_t u + \sum_{k=1}^d \phi^\ep_{e_k} \partial_t \left( \eta \partial_{x_k} u \right)
- \sum_{k=1}^d \nabla \left( \eta \partial_{x_k} u \right) \cdot \a^\ep \left( e_k + \nabla \phi_{e_k}\left(\tfrac \cdot\ep \right) \right) 
\\ & \quad
- \nabla \cdot \left( \a^\ep \left( \sum_{k=1}^d \phi^\ep_{e_k} \nabla \left( \eta \partial_{x_k} u \right) 
+ (1-\eta)\nabla u \right) \right).
\end{align*}
Writing the homogenized equation for $u$ in the form
\begin{equation*} \label{}
\partial_t u 
= \nabla \cdot \ahom \nabla u 
= \sum_{k=1}^d \nabla \left( \eta \partial_{x_k} u\right) \cdot \ahom e_k + \nabla \cdot\left( \left( 1-\eta\right) \ahom \nabla u \right),
\end{equation*}
we obtain the identity
\begin{align*}
\left( \partial_t - \nabla \cdot \a^\ep\nabla \right) w^\ep
& =
 \sum_{k=1}^d \phi^\ep_{e_k} \partial_t \left( \eta \partial_{x_k} u \right)
- \sum_{k=1}^d \nabla \left( \eta \partial_{x_k} u \right) \cdot \left( \a^\ep \left( e_k + \nabla \phi_{e_k}\left(\tfrac \cdot\ep \right) \right) - \ahom e_k \right) 
\\ & \qquad
- \nabla \cdot \left( \a^\ep  \sum_{k=1}^d \phi^\ep_{e_k} \nabla \left( \eta \partial_{x_k} u \right)  \right) 
-\nabla \cdot \left(
\left( \a^\ep - \ahom \right)(1-\eta)\nabla u 
\right).
\end{align*}
We therefore obtain the bound 
\begin{align*}
\lefteqn{
\left\| \left( \partial_t - \nabla \cdot \a^\ep\nabla \right) w^\ep\right\|_{L^2(I ;H^{-1}(U)) }
} \ \ &
\\ & 
\leq
C  \sum_{k=1}^d \ep \left\|  \phi^\ep_{e_k} \partial_t \left( \eta \partial_{x_k} u \right) \right\|_{L^2(I\times U)}
\\ & \quad 
+
\sum_{k=1}^d \left\| \nabla \left( \eta \partial_{x_k} u \right) \cdot \left( \a^\ep \left( e_k + \nabla \phi_{e_k}\left(\tfrac \cdot\ep \right) \right) - \ahom e_k \right) \right\|_{L^2(I;H^{-1}(U))} 
\\ & \quad 
+ C  \sum_{k=1}^d\ep \left\|  \phi^\ep_{e_k} \nabla \left( \eta \partial_{x_k} u \right) \right\|_{L^2(I\times U)} 
+ C \left\| (1-\eta)\nabla u \right\|_{L^2(I\times U)}. 
\end{align*}
We will now estimate each of the four terms appearing on the right side of the previous inequality using~\eqref{e.etacutoff.parab},~\eqref{e.upoint.bounds},~\eqref{e.umeyers.parab} and the definition of~$\mathcal{E}(\ep)$. For the first term, we use~\eqref{e.etacutoff.parab} and~\eqref{e.upoint.bounds} to get, for each $k\in\{1,\ldots,d\}$, 
\begin{align*}
\ep\left\|  \phi^\ep_{e_k} \partial_t \left( \eta \partial_{x_k} u \right) \right\|_{L^2(I\times U)}
&
\leq
C \left\| \partial_t \left( \eta \partial_{x_k} u \right)  \right\|_{L^\infty(I_r\times U_r)} \ep \left\| \phi^\ep_{e_k} \right\|_{L^2(U)}
\\ & 
\leq C r^{-3 -(2+d)/2} \mathcal{E}(\ep) \left\| f \right\|_{H^1_\pa(I \times U)},
\end{align*}
where here we also use that
\begin{equation} 
\label{e.Hminusoneandscaling}
\ep \left\| \phi^\ep \right\|_{L^2(U)} 
\leq C\ep \left\| \nabla \phi^\ep \right\|_{H^{-1}(U)} 
= C \left\| \nabla \phi\left( \tfrac \cdot\ep\right) \right\|_{H^{-1}(U)} 
\leq C \mathcal{E}(\ep).
\end{equation}
We estimate the second term similarly, to find 
\begin{align*} \label{}
\lefteqn{ 
\left\| \nabla \left( \eta \partial_{x_k} u \right) \cdot \left( \a^\ep \left( e_k + \nabla \phi_{e_k}\left(\tfrac \cdot\ep \right) \right) - \ahom e_k \right) \right\|_{L^2(I;H^{-1}( U))} 
} \qquad & 
\\ & 
\leq 
\left\| \nabla \left( \eta \partial_{x_k} u \right) \right\|_{L^{\infty}(I; W^{1,\infty}( U))} 
\left\| \a^\ep \left( e_k + \nabla \phi_{e_k}\left(\tfrac \cdot\ep \right) \right) - \ahom e_k\right\|_{H^{-1}(U)} 
\\ &
\leq 
C r^{-2 -(2+d)/2} \mathcal{E}(\ep) \left\| f \right\|_{H^1_\pa(I \times U)}. 
\end{align*}
For the third term 
\begin{align*}
\ep \left\|  \phi^\ep_{e_k} \nabla \left( \eta \partial_{x_k} u \right) \right\|_{L^2(I\times U)} 
&
\leq 
\left\| \nabla \left( \eta \partial_{x_k} u \right)  \right\|_{L^\infty(I\times U)} 
\ep \left\| \phi^\ep_{e_k} \right\|_{L^2(U)} 
\\ &
\leq 
C r^{-2 -(2+d)/2} \mathcal{E}(\ep) \left\| f \right\|_{H^1_\pa(I \times U)}.
\end{align*}
Finally, we estimate the fourth term using H\"older's inequality and~\eqref{e.umeyers.parab}:
\begin{align*}
\left\| (1-\eta)\nabla u \right\|_{L^2(I\times U)}
&
\leq  
C \left| \left\{ x\in I\times U\,:\, \eta \neq 1 \right\} \right|^{\frac{\delta}{4+2\delta}} \left\| \nabla u \right\|_{L^{2+\delta}(I\times U)} 
\\ &
\leq 
C r^{\frac{\delta}{4+2\delta}}
\left\| f \right\|_{W^{1,2+\delta}_\pa(I\times U)}.
\end{align*}
This completes the proof of the lemma. 
\end{proof}

We next deduce from the previous lemma an estimate on~$\left\| u^\ep-w^\ep\right\|_{H^1_\pa(I\times U)}$. 

\begin{lemma}
\label{l.H1est.parab}
There exists $C(\delta,d,\Lambda)<\infty$ such that
\begin{equation*} \label{}
\left\| u^\ep - w^\ep \right\|_{H^1_\pa(I\times U)} 
\leq 
C \left( r^{\frac{\delta}{4+2\delta}} + r^{-3-(2+d)/2} \mathcal{E}(\ep)  \right) \left\| f \right\|_{W^{1,2+\delta}_\pa(I\times U)}. 
\end{equation*}
\end{lemma}
\begin{proof}
We compute 
\begin{align*}
\left\| \nabla \left(u^\ep - w^\ep \right) \right\|_{L^2(I\times U)}^2
& 
\leq 
C \int_{I\times U} \nabla \left(u^\ep - w^\ep \right) \cdot \a^\ep \nabla \left(u^\ep - w^\ep \right)
\\ & 
\leq
C\int_{I\times U} \left(u^\ep - w^\ep \right) \left( \partial_t - \nabla \cdot \a^\ep \nabla \right) \left(u^\ep - w^\ep \right)
\\ & 
= 
C\int_{I\times U} \left(u^\ep - w^\ep \right) \left( \partial_t - \nabla \cdot \a^\ep \nabla \right) w^\ep 
\\ & 
\leq 
C \left\| u^\ep - w^\ep \right\|_{L^2(I ; H^{1}( U))}
\left\|\left( \partial_t - \nabla \cdot \a^\ep \nabla \right) w^\ep  \right\|_{L^2(I ; H^{-1}( U))}. 
\end{align*}
This yields
\begin{equation*} \label{}
\left\| u^\ep - w^\ep \right\|_{L^2(I ; H^{1}( U))}
\leq C \left\|\left( \partial_t - \nabla \cdot \a^\ep \nabla \right) w^\ep  \right\|_{L^2(I ; H^{-1}( U))}.
\end{equation*}
For the time derivative, we have
\begin{align*}
\lefteqn{
\left\| \partial_t  \left(u^\ep - w^\ep \right) \right\|_{L^2(I ; H^{-1}(U))}
} \quad &
\\ &
\leq 
\left\| \left( \partial_t  - \nabla \cdot \a^\ep\nabla \right) \left(u^\ep - w^\ep \right) \right\|_{L^2(I ; H^{-1}(U))}
+
\left\| \nabla \cdot \a^\ep\nabla \left(u^\ep - w^\ep \right) \right\|_{L^2(I ; H^{-1}(U))}
\\ & 
\leq
\left\| \left( \partial_t  - \nabla \cdot \a^\ep\nabla \right) w^\ep  \right\|_{L^2(I ; H^{-1}(U))}
+
C \left\| u^\ep - w^\ep  \right\|_{L^2(I ; H^{1}(U))}
\\ & 
\leq 
C \left\| \left( \partial_t  - \nabla \cdot \a^\ep\nabla \right) w^\ep  \right\|_{L^2(I ; H^{-1}(U))}.
\end{align*}
Combining these yields
\begin{align*}
\left\| u^\ep - w^\ep \right\|_{H^1_\pa(I\times U)} 
\leq 
C \left\| \left( \partial_t  - \nabla \cdot \a^\ep\nabla \right) w^\ep  \right\|_{L^2(I ; H^{-1}(U))}.
\end{align*}
An appeal to Lemma~\ref{l.CDP.plugineq} completes the proof. 
\end{proof}

We next estimate the difference~$w^\ep - u$ in weak norms. 

\begin{lemma}
\label{l.throwcorrectors}
There exists $C(\delta,d,\Lambda)<\infty$ such that
\begin{multline}
\label{e.throwcorrectors}
\left\| w^\ep - u \right\|_{L^2(I\times U)} 
+ \left\| \nabla (w^\ep - u) \right\|_{L^2(I;H^{-1}(U))}
+ \left\| \a^\ep \nabla w^\ep - \ahom \nabla u \right\|_{L^2(I;H^{-1}(U))}
\\
\leq  
C \left( r^{\frac{\delta}{4+2\delta}} + r^{-2-(2+d)/2} \mathcal{E}(\ep)  \right) \left\| f \right\|_{W^{1,2+\delta}_\pa(I\times U)}.
\end{multline}
\end{lemma}
\begin{proof}
Observe that 
\begin{equation*} \label{}
\nabla w^\ep - \nabla u 
= \ep \sum_{k=1}^d \nabla \left( \eta \partial_{x_k} u \right) \phi_{e_k}^\ep(x)
+ 
\eta 
\sum_{k=1}^d \partial_{x_k} u \nabla \phi_{e_k}\left( \tfrac \cdot \ep \right).
\end{equation*}
Thus we have, by~\eqref{e.etacutoff.parab},~\eqref{e.upoint.bounds}, the definition of~$\mathcal{E}(\ep)$ and~\eqref{e.Hminusoneandscaling},
\begin{align*}
\lefteqn{
\left\| \nabla w^\ep - \nabla u \right\|_{L^2(I;H^{-1}(U))}
} \quad & 
\\ &
\leq  \left\| \nabla (\eta \nabla u) \right\|_{L^\infty(I\times U)} \sum_{k=1}^d \ep \left\| \phi^\ep_{e_k} \right\|_{L^2(U)}
+
C \left\| \eta \nabla u\right\|_{L^\infty(I;W^{1,\infty}(U))} 
\left\|  \nabla \phi_{e_k}\left( \tfrac \cdot \ep \right) \right\|_{H^{-1}(U)}
\\ & 
\leq 
Cr^{-2-(2+d)/2} \left\| f \right\|_{H^1_\pa(I\times U)} \mathcal{E}(\ep). 
\end{align*}
Notice that since~$\left\| w^\ep - u \right\|_{L^2(I\times U)}\leq C \left\| \nabla w^\ep - \nabla u \right\|_{L^2(I; H^{-1} (U))}$, 
this also gives the estimate for $\left\| w^\ep - u \right\|_{L^2(I\times U)}$. 

\smallskip

We turn to the estimate for the fluxes. For clarity, we perform the computations using coordinates. The $i$th component of the flux vector $\a^\ep\nabla w^\ep$ is given by 
\begin{align*}
\left( \a^\ep\nabla w^\ep \right)_i
= 
\sum_{j,k=1}^d \eta \a^\ep_{ij} \partial_{x_k} u \left( \delta_{jk} + \partial_{x_j} \phi_{e_k} \left( \tfrac \cdot \ep \right) \right)
+ 
\ep \sum_{j,k=1}^d \a^\ep_{ij} \partial_{x_j} \left( \eta \partial_{x_k} u \right) \phi_{e_k}^\ep.
\end{align*}
Thus
\begin{align*}
\left( \a^\ep \nabla w^\ep \right)_i - \left( \ahom \nabla u\right)_i
& 
=
\sum_{j,k=1}^d \eta \partial_{x_k} u \left( \a^\ep_{ij} \left( \delta_{jk}+\partial_{x_j} \phi_{e_k}\left( \tfrac \cdot\ep\right)\right) - \ahom_{ik} \right)
\\ & \quad
+ \ep \sum_{j,k=1}^d \a^\ep_{ij} \partial_{x_j} \left( \eta \partial_{x_k} u\right) \phi^\ep_{e_k}
+ \sum_{j,k=1}^d \left( 1- \eta \right) \ahom_{ik}\partial_{x_k} u.
\end{align*}
We next estimate each of the~$\| \cdot\|_{L^2(I;H^{-1}(U))}$ norm of the three terms on the right side of the previous display. For the first term, we have, by~\eqref{e.etacutoff.parab},~\eqref{e.upoint.bounds} and the definition of~$\mathcal{E}(\ep)$,   
\begin{align*}
\lefteqn{
\left\| 
\eta \partial_{x_k} u \left( \a^\ep_{ij} \left( \delta_{jk}+\partial_{x_j} \phi_{e_k}\left( \tfrac \cdot\ep\right)\right) - \ahom_{ik} \right)
\right\|_{L^2(I;H^{-1}(U))}
} \quad & 
\\ &
\leq 
C \left\| \eta \nabla u \right\|_{L^\infty(I;W^{1,\infty}(U))}
\left\|  \a^\ep \left( e_k+\nabla \phi_{e_k}\left( \tfrac \cdot\ep\right)\right) - \ahom e_k \right\|_{L^2(I;H^{-1}(U))}
\\ & 
\leq 
Cr^{-2-(2+d)/2} \mathcal{E}(\ep). 
\end{align*}
For the second term, we use~\eqref{e.etacutoff.parab},~\eqref{e.upoint.bounds},  the definition of~$\mathcal{E}(\ep)$ and~\eqref{e.Hminusoneandscaling} to obtain
\begin{align*}
\ep \left\| \a^\ep_{ij} \partial_{x_j} \left( \eta \partial_{x_k} u\right) \phi^\ep_{e_k} \right\|_{L^2(I; H^{-1}( U))} 
&
\leq
C \ep \left\| \a^\ep_{ij} \partial_{x_j} \left( \eta \partial_{x_k} u\right) \phi^\ep_{e_k} \right\|_{L^2(I\times U)} 
\\ &
\leq 
C 
\left\| \nabla \left( \eta \nabla u\right) \right\|_{L^{\infty}(I\times U)}
\ep \left\| \phi^\ep_{e_k} \right\|_{L^2(U)} 
\\ & 
\leq 
C r^{-2-(2+d)/2} \mathcal{E}(\ep). 
\end{align*}
For the third term, we use H\"older's inequality, ~\eqref{e.etacutoff.parab} and~\eqref{e.umeyers.parab} to find
\begin{align*} \label{}
\left\| (1-\eta) \ahom \nabla u \right\|_{L^2(I; H^{-1}( U))}
& 
\leq 
\left\| (1-\eta) \ahom \nabla u \right\|_{L^2(I\times U)}
\\ & 
\leq 
C \left| \left\{ x\in I\times U\,:\, \eta \neq 1 \right\} \right|^{\frac{\delta}{4+2\delta}} \left\| \nabla u \right\|_{L^{2+\delta}(I\times U)} 
\\ &
\leq 
C r^{\frac{\delta}{4+2\delta}}
\left\| f \right\|_{W^{1,2+\delta}_\pa(I\times U)}.
\end{align*}
Combining the last four displays yields the desired estimate on the fluxes. This completes the proof. 
\end{proof}

We now complete the proof of Theorem~\ref{t.CDP}. 

\begin{proof}[Proof of Theorem~\ref{t.CDP}]
We deduce from Lemmas~\ref{l.H1est.parab} and~\ref{l.throwcorrectors} the estimate
\begin{multline}
\label{e.CDPestimates2.bb}
\left\| u^\ep - u \right\|_{L^2(I\times U)} 
+ \left\| \nabla (u^\ep - u) \right\|_{L^2(I;H^{-1}(U))}
+ \left\| \a^\ep \nabla u^\ep - \ahom \nabla u \right\|_{L^2(I;H^{-1}(U))}
\\
\leq  
C \left( r^{\frac{\delta}{4+2\delta}} + r^{-3-(2+d)/2} \mathcal{E}(\ep)  \right) \left\| f \right\|_{W^{1,2+\delta}_\pa(I\times U)}.
\end{multline}
This is valid for any choice of the mesoscale $r\in [\ep,1)$. In other words, we have proved the parabolic version of Theorem~\ref{t.DP.blackbox}: compare the estimate~\eqref{e.CDPestimates.bb} to the one in~\eqref{e.DPestimates.bb}. The passage from~\eqref{e.CDPestimates2.bb} now follows exactly the same argument as the passage from Theorem~\ref{t.DP.blackbox} to Theorem~\ref{t.DP}. Indeed, there are no differences in the argument, since what is left concerns how to choose the mesoscale~$r$, define the random variable $\X_s$ and estimate it in terms of $\mathcal{E}(\ep)$.
\end{proof}

\section{Parabolic \texorpdfstring{$C^{0,1}$}{C0,1}-type estimate}
\label{s.parab.Lip}

In this section, we prove the following interior $C^{0,1}$-type estimate for the parabolic equation~\eqref{e.pde.parab}, which is a generalization of the elliptic version given in Theorem~\ref{t.Lipschitz}. 

\begin{theorem}[{Parabolic $C^{0,1}$-type estimate}]
\label{t.Lipschitz.parab}
Fix $s\in (0,d)$. There exist a constant  $C(s,d,\Lambda)<\infty$ and a random variable $\X_s :\Omega \to [1,\infty]$ satisfying
\begin{equation}
\label{e.XLipschitz.parab}
\X_s = \O_s\left( C \right),
\end{equation}
such that the following holds: for every $R\geq \X_s$ and weak solution $u\in H^1_\pa(Q_R)$ of
\begin{equation} 
\label{e.wksrt2.parab}
\partial_t u -\nabla \cdot \left( \a\nabla u \right) = 0 \quad \mbox{in} \ Q_R,
\end{equation}
we have, for every  $r\in \left[ \X_s , R\right]$, the estimate 
\begin{equation}
\label{e.Lipschitz.parab}
 \frac1{r} \left\| u - \left( u \right)_{Q_{r}} \right\|_{\underline{L}^2(Q_r)}
\leq
\frac{C}{R} \left\| u - \left( u \right)_{Q_{R}} \right\|_{\underline{L}^2(Q_{R})}.
\end{equation}
\end{theorem}
\begin{remark} 
\label{r.Lipschitz.parab2}
We can use the parabolic versions of the Caccioppoli and Poincar\'e's inequalities to deduce, for $r \in [\X_s,R]$,
\begin{equation} 
\label{e.Lipschitz.parab2}
\left\| \nabla u \right\|_{\underline{L}^2(Q_{r})} 
\leq 
C \left\| \nabla u \right\|_{\underline{L}^2(Q_{R})}.
\end{equation}
Indeed, the Caccioppoli estimate yields that 
\begin{equation*} 
\left\| \nabla u \right\|_{\underline{L}^2(Q_{r})} \leq \frac{C}{r} \left\| u - \left( u \right)_{Q_{2r}} \right\|_{\underline{L}^2(Q_{2r})}.
\end{equation*}
For the other direction we use the equation. Let $\tilde R \in \left(\frac12 R , R \right)$ be such that
\begin{equation*} 
\int_{-R^2}^0 \int_{\partial B_{\tilde R}} |\nabla u(t,x)| \, dH^{d-1}(x) \, dt \leq \frac{C}{R} \int_{-R^2}^0 \int_{B_{R}} |\nabla u(t,x)| \, dx \, dt.
\end{equation*}
Then integration by parts provides us, for $t_1,t_2 \in (-R^2,0)$, 
\begin{equation*} 
\frac1R \left|\fint_{B_{\tilde R}} u(t_1,x)\,dx - \fint_{B_{\tilde R}} u(t_2,x)\,dx \right| \leq C \fint_{Q_{R}} |\nabla u(t,x)| \, dx \, dt \leq C \left\| \nabla u \right\|_{\underline{L}^2(Q_{R})} .
\end{equation*}
It also follows that 
\begin{equation*} 
\sup_{t \in (-R^2,0)} \frac 1R \left|\fint_{B_{\tilde R}} u(t,x)\,dx -   \left( u \right)_{Q_{\tilde R}} \right| \leq C \left\| \nabla u \right\|_{\underline{L}^2(Q_{R})}
\end{equation*}
Therefore Poincar\'e's inequality yields
\begin{equation*} 
\frac{C}{\tilde R} \left\| u - \left( u \right)_{Q_{\tilde R}} \right\|_{\underline{L}^2(Q_{\tilde R})} \leq  C \left\| \nabla u \right\|_{\underline{L}^2(Q_{R})}.
\end{equation*}
In view of Lemma~\ref{l.slicesvsaverages.PAR}, we also obtain from~\eqref{e.Lipschitz.parab} the bound, for every $r\in\left[\X_s, \frac12 R \right]$ and $t\in \left[-\frac12 R^2,0\right]$, 
\begin{equation} 
\label{e.Lipschitz.parab.slice}
\left\| \nabla u(t,\cdot) \right\|_{\underline{L}^2(B_r)}
\leq 
C \left\| \nabla u \right\|_{\underline{L}^2(Q_{R})}.
\end{equation}
\end{remark}

The proof of Theorem~\ref{t.Lipschitz.parab} requires the following parabolic version of Lemma~\ref{l.harmapproxlipschitz}.

\begin{lemma}
\label{l.harmapproxlipschitz.parab}
Fix $\alpha \in (0,1]$,  $K\geq 1$ and $X\geq 1$. 
Let $R\geq 2X$ and $u\in L^2(Q_{R})$ have the property that, for every $r\in \left[ X, \frac12R \right]$, there exists $w\in {H}^1_\pa(Q_r)$ which is a weak solution of 
\begin{equation*} \label{}
\partial_t w -\nabla \cdot \left( \ahom \nabla w \right) = 0 \quad \mbox{in} \ Q_r,
\end{equation*}
satisfying 
\begin{equation} 
\label{e.wharmapproxit.parab}
\left\| u - w \right\|_{L^2(Q_R)} \leq K r^{-\alpha} \left\| u - \left( u \right)_{Q_{r}} \right\|_{\underline{L}^2(Q_{r})}. 
\end{equation}
Then there exists $C(\alpha,K,d,\Lambda)<\infty$ such that, for every $r\in \left[ X, R\right]$, 
\begin{equation}
\label{e.elipsc.parab}
 \frac1{r} \left\| u - \left( u \right)_{Q_{r}} \right\|_{\underline{L}^2(Q_{r})}
\leq
\frac{C}{R} \left\| u - \left( u \right)_{Q_{R}} \right\|_{\underline{L}^2(Q_{R})}.
\end{equation}
\end{lemma}

The proof of Lemma~\ref{l.harmapproxlipschitz.parab} requires only typographical changes to the proof of Lemma~\ref{l.harmapproxlipschitz}, provided we have a suitable parabolic version of Lemma~\ref{l.u:decayestimate}. The latter is given in Lemma~\ref{l.u:decayestimate.parab}, below. Before giving the statement of this lemma, we must first introduce the notion of~\emph{$\ahom$-caloric polynomials} and discuss the pointwise regularity of caloric functions. This requires some new notation and definitions. 

\smallskip

The parabolic analogue of the pointwise estimate~\eqref{e.pointwiseharm} states (see for instance~\cite[Section~2.3.3.c]{Evans}) that, for every~$k,l\in\N$, there exists~$C(k+2l,d,\Lambda)<\infty$ such that, for every~$R>0$ and~$\ahom$-caloric function~$u$ on~$Q_R$, we have the estimate
\begin{equation}
\label{e.pointwisecaloric}
\left\|  \nabla^k \partial_t^l u \right\|_{L^\infty(Q_{R/2})}  \leq \frac{C}{R^{k+2l}} \left\| u \right\|_{\underline{L}^1(Q_{R})}. 
\end{equation}
Just like in Chapter~\ref{c.regularity}, we next write some consequences of~\eqref{e.pointwisecaloric} which tell us that~$\ahom$-caloric functions can be locally well-approximated by polynomials.

\smallskip

We denote polynomials in the variables~$t,x_1,\ldots,x_d$ by $\mathcal{P}(\R\times \Rd)$. The \emph{parabolic order}~$\deg_p(w)$ of an element~$w \in \mathcal{P}(\R\times \Rd)$ is the order of the polynomial~$(t,x) \mapsto w(t^2,x)$. Slightly abusing notation, for each $k\in\N$ we let $\mathcal{P}_k(\R \times\Rd)$ be the subset of~$\mathcal{P}(\R\times\Rd)$ of polynomials with parabolic order at most~$k$. For $\alpha>0$, we say that a function $\phi:\R\times \Rd\to \R$ is \emph{parabolically $\alpha$-homogeneous} if
\index{parabolically homogeneous}
\begin{equation*} \label{}
\forall \lambda \in\R, \quad \phi(\lambda^2 t,\lambda x) = \lambda^\alpha \phi(t,x).
\end{equation*}

It is clear that any element of $\mathcal{P}_k(\R \times\Rd)$ can be written as a sum of parabolically homogeneous polynomials. We denote by $\Ahom_k(Q_\infty)$ the set of~$\ahom$-caloric functions on~$Q_\infty$ with growth which is slower than a polynomial of parabolic degree $k+1$:
\begin{multline*} \label{}
\Ahom_k (Q_\infty) 
:=
\\
\left\{ 
w\in H^1_\pa(Q_\infty) \,:\,
\limsup_{r \to \infty} r^{-(k+1)}  \left\| w \right\|_{\underline{L}^2(Q_r) }  = 0, \ \partial_t w - \nabla \cdot \left( \ahom \nabla w \right) = 0 \ \mbox{in} \ Q_\infty
\right\}.
\end{multline*}
It turns out that $\Ahom_k(Q_\infty)$ coincides with the set of \emph{$\ahom$-caloric polynomials}\footnote{$\ahom$-caloric polynomials are often called \emph{heat polynomials} in the literature, in the case $\ahom=I_d$} of parabolic degree at most $k$. That is, 
\begin{equation} 
\label{e.identifyAhomkQinfty}
\Ahom_k(Q_\infty)
=
\left\{ w\vert_{Q_\infty}  \,:\, w \in \mathcal{P}_k(\R \times \Rd), \ \partial_t w- \nabla \cdot \left( \ahom \nabla w\right) = 0 \ \mbox{in} \ \R\times \Rd \right\}. 
\end{equation}
Let us comment briefly on the case of the homogeneous heat equation. We denote by $v_k(t,s)$ the parabolically $k$-homogeneous heat polynomial \index{heat polynomial} in the one-dimensional case, which satisfies the heat equation,~$\partial_t v_k(t,s) = \partial_s^2 v_k(t,s)$. It can be written as
\begin{equation*} 
v_n(t,s) = \sum_{k=0}^{\lfloor n/2 \rfloor} \frac{s^{n-2k} t^k }{(n-2k)! k!} .
\end{equation*}
It turns out that the vector space of heat polynomials can be characterized as
\begin{multline*} 
\left\{ w\vert_{Q_\infty}  \,:\, w \in \mathcal{P}_k(\R \times \Rd), \ \partial_t w- \Delta w = 0 \ \mbox{in} \ \R\times \Rd \right\}
\\ = \spn \left\{ \prod_{j=1}^d v_{\beta_j}(t,x_j) \, : \, \beta_1 + \ldots + \beta_d \leq k \right\}.
\end{multline*}
To see this, letting $y^\beta = \prod_{j=1}^d y_j^{\beta_j}$, $|\beta|=n \in \N$, be the initial values for the heat equation, we have by the binomial formula and symmetry that the solution has the form, up to a multiplicative constant,  
\begin{align*} 
u(t,x) & = \int_{\R^d} y^\beta t^{-\frac d2} \exp \left( - \frac{|x-y|^2}{4t} \right) \, dy
\\ & =  \int_{\R^d} \left( x+ \sqrt{t}y \right)^\beta  \exp \left( -  \frac{|y|^2}{4} \right) \, dy 
\\ & = \prod_{j=1}^d \int_{\R}  \left( x_j +  \sqrt{t}y_j \right)^{\beta_j} \exp \left( - \frac{|y_j|^2}{4}\right) \, dy_j
\\ & = \prod_{j=1}^d \sum_{k=0}^{\beta_j} \binom{\beta_j}{k}  x_j^{\beta_j - k} t^{\frac{k}{2}}  \int_{\R}  y_j^{k} \exp \left( - \frac{|y_j|^2}{4} \right) \, dy_j 
\\ & = \prod_{j=1}^d \sum_{k=0}^{\lfloor \beta_j /2 \rfloor} \binom{\beta_j}{2k}  x_j^{\beta_j - 2k} t^{k} \int_{\R}  y_j^{2k} \exp \left( - \frac{|y_j|^2}{4} \right) \, dy_j  .
\end{align*}
Since
\begin{equation*} 
\int_{\R}  y_j^{2k} \exp \left( - \frac{|y_j|^2}{4} \right) \, dy_j   =  \frac{(2k)!}{k!},
\end{equation*}
we conclude that $u \in \spn \left\{ \prod_{j=1}^d v_{\beta_j}(t,x_j) \, : \, \beta_1 + \ldots + \beta_d = n \right\}$. 

\smallskip

Conversely, by the uniqueness of the backwards heat equation, see e.g.~\cite{Vessella} and the references therein, in the class
\begin{equation*} 
\mathcal{H}_\delta := \left\{ u \in H^2((0,1) \times  \R^d) \, : \, \left\| u \exp\left(- \delta |\cdot|^2\right) \right\|_{L^\infty} < \infty \right\}
\end{equation*}
for small $\delta \in (0,1)$, we have that if there are two heat polynomials, say $u_1$ and $u_2$, agreeing at the time~$t$, then $u_1 \equiv u_2$. Taking $w$ to be a $n$-homogeneous heat polynomial, we have that $w(0,\cdot) \in   \mathcal{P}_n(\Rd)$, and by the previous computation we find $\tilde w \in \spn \left\{ \prod_{j=1}^d v_{\beta_j}(t,x_j) \, : \, \beta_1 + \ldots + \beta_d = n \right\}$ such that $\tilde w(0,\cdot) = w(0,\cdot)$. Thus the uniqueness implies that $\tilde w \equiv w$ and $w \in \spn \left\{ \prod_{j=1}^d v_{\beta_j}(t,x_j) \, : \, \beta_1 + \ldots + \beta_d = n \right\}$. 
As a consequence, the class of $n$-homogeneous heat polynomials is isomorphic to the $n$-homogeneous polynomials of $\R^d$.

\begin{exercise}
Show that~\eqref{e.identifyAhomkQinfty} is a consequence of~\eqref{e.pointwisecaloric}. Show moreover that, for every $k\in\N$, there exists $C(k,d,\Lambda)<\infty$ such that, for every $r>0$ and $\ahom$-caloric function $w \in H^1_\pa(Q_r)$, there exists an $\ahom$-caloric polynomial~$p$ of parabolic order at most $k$ such that, for every  $s\in \left( 0,\tfrac12r\right]$, 
\begin{equation} 
\label{e.Ahomcaloric}
\left\| w - p \right\|_{L^\infty(Q_s)} 
\leq
C\left( \frac sr \right)^{k+1} \left\| w \right\|_{\underline{L}^2(Q_r)}. 
\end{equation}
Hint: Prove the second statement first. Take~$p$ to be the~$k$th (parabolic) order Taylor polynomial to~$w$ at~$(0,0)$ and use a blow-up argument to show that~$p$ is~$\ahom$-caloric. Apply~\eqref{e.pointwisecaloric} and Taylor's remainder theorem to get the desired estimate. 
\end{exercise}

We next state a parabolic counterpart of Lemma~\ref{l.u:decayestimate}. 

\begin{lemma}
\label{l.u:decayestimate.parab}
Fix $\alpha \in [0,1]$,  $K\geq 1$ and $X\geq 1$. 
Let $R\geq 2X$ and $u\in L^2(Q_{R})$ have the property that, for every $r\in \left[ X, R \right]$, there exists $w_r \in H^1_{\pa}(Q_{r/2})$ which is a solution of 
\begin{equation*} \label{}
\partial_t w_r -\nabla \cdot \left( \ahom \nabla w_r \right) = 0 \quad \mbox{in} \ Q_{r/2}
\end{equation*}
and satisfies 
\begin{equation} 
\label{e.wharmapproxit02.parab}
\left\| u - w_r \right\|_{\underline{L}^2(Q_{r/2})} \leq K r^{-\alpha} \left\| u - \left( u \right)_{Q_{r}} \right\|_{\underline{L}^2(Q_{r})}. 
\end{equation}
Then, for every $k \in \N$, there exists $\theta(\alpha,k,d,\Lambda) \in (0,\frac12)$  and $C(\alpha,k,d,\Lambda)<\infty$ such that, for every $r\in \left[ X, R\right]$, 
\begin{multline}
\label{e.u:decayestimate.parab}
\inf_{p \in \Ahom_k(Q_\infty) } \left\| u - p \right\|_{\underline{L}^2(Q_{\theta r})}
\\
\leq
\frac14 \theta^{k+1 - \alpha/2} \inf_{p \in \Ahom_k(Q_\infty) } \left\| u - p  \right\|_{\underline{L}^2(Q_{r})}   + C K r^{-\alpha} \left\| u - (u)_{Q_r} \right\|_{\underline{L}^2(Q_{r})} .
\end{multline}
\end{lemma}
\begin{proof} 
The proof is similar to that of Lemma~\ref{l.u:decayestimate} and we omit it. 
\end{proof}

\begin{proof}[{Proof of Lemma~\ref{l.harmapproxlipschitz.parab}}]
The proof is similar to that of Lemma~\ref{l.harmapproxlipschitz} and is omitted. 
\end{proof}

\section{Parabolic higher regularity and Liouville theorems}
\label{s.parab.reg}

In this section, we generalize the parabolic Liouville theorem implicit in~\eqref{e.identifyAhomkQinfty} and the pointwise bounds given in~\eqref{e.Ahomcaloric} to $\a(x)$-caloric functions. Denote, for every $r\in (0,\infty]$,
\begin{equation*} \label{}
\A(Q_r):= 
\left\{
w\in H^1_\pa(Q_r)
\,:\,
 \partial_t w - \nabla \cdot \left( \a \nabla w \right) = 0 \ \mbox{in} \ Q_r
\right\},
\end{equation*}
and, for every $k\in\N$,
\begin{equation*} \label{}
\A_k(Q_\infty):= 
\left\{ 
w\in \A(Q_\infty)  \,:\,
\limsup_{r \to \infty} r^{-(k+1)}  \left\| w \right\|_{\underline{L}^2(Q_r) }  = 0
\right\}.
\end{equation*}
The following theorem is a parabolic analogue of Theorem~\ref{t.regularity}. 
\index{Liouville theorem}
\index{$C^{k,1}$ estimate}

\begin{theorem}
[Parabolic higher regularity theory]
\label{t.regularity.parab}
Fix $s \in (0,d)$. There exist an exponent $\delta(s,d,\Lambda)\in \left( 0, \frac12 \right]$ and a random variable $\X_s$ satisfying the estimate
\begin{equation}
\label{e.X.parab}
\X_s \leq \O_s\left(C(s,d,\Lambda)\right)
\end{equation}
such that the following statements hold, for every $k\in\N$:
\begin{enumerate}
\item[{$\mathrm{(i)}_k$}] There exists $C(k,d,\Lambda)<\infty$ such that, for every $u \in \A_k(Q_\infty)$, there exists $p\in \Ahom_k(Q_\infty)$ such that, for every $R\geq \X_s$,
\begin{equation} \label{e.liouvillec.parab}
\left\| u - p \right\|_{\underline{L}^2(Q_R)} \leq C R^{-\delta} \left\| p \right\|_{\underline{L}^2(Q_R)}.
\end{equation}

\item[{$\mathrm{(ii)}_k$}]For every $p\in \Ahom_k(Q_\infty)$, there exists $u\in \A_k(Q_\infty)$ satisfying~\eqref{e.liouvillec.parab} for every $R\geq \X_s$. 

\item[{$\mathrm{(iii)}_k$}]
There exists $C(k,d,\Lambda)<\infty$ such that, for every $R\geq \X_s$ and $u\in \A(Q_R)$, there exists $\phi \in \A_k(Q_\infty)$ such that, for every $r \in \left[ \X_s,  R \right]$, we have the estimate
\begin{equation}
\label{e.intrinsicreg.parab}
\left\| u - \phi \right\|_{\underline{L}^2(Q_r)} \leq C \left( \frac r R \right)^{k+1} 
\left\| u \right\|_{\underline{L}^2(Q_R)}.
\end{equation}
\end{enumerate}
In particular, $\P$-almost surely, we have for every $k\in\N$ that
\begin{equation} 
\label{e.dimensionofAk.parab}
\dim(\A_k(Q_\infty)) = \dim(\Ahom_k(Q_\infty)) =   \binom{d+k}{d}.
\end{equation}
\end{theorem}

The proof of Theorem~\ref{t.regularity.parab} is almost the same as the one of Theorem~\ref{t.regularity}. With Lemma~\ref{l.u:decayestimate.parab} in hand, we just needs to replace balls by parabolic cylinders in the argument. We therefore omit the details.

\begin{exercise}
Formulate and prove a parabolic version of Corollary~\ref{cor.regularity} and the result of Exercise~\ref{ex.regularity}.
\end{exercise}

\section{Decay of the Green functions and their gradients}
\label{s.greenie.decay}

We begin this section by introducing the elliptic and parabolic Green functions and reviewing some of their basic properties. A complete construction of the Green functions is given in Appendix~\ref{a.NA}, where we also prove the Nash-Aronson estimates mentioned below. Here we briefly summarize the deterministic estimates presented there, before presenting gradient decay estimates which are consequences of the (stochastic, large-scale) regularity theory. 

\smallskip

We denote by~$P = P(t,x,y)$ the \emph{parabolic Green function} for the operator $-\nabla \cdot \a\nabla$. We can characterize $P$ by the fact that, for each fixed $y \in \R\times \Rd$, the function $(t,x) \mapsto P(t,x,y)$ is the solution of the initial-value problem
\begin{equation} 
\label{e.Green.parab}
\left\{ 
\begin{aligned}
& \partial_t  P (\cdot,y) - \nabla \cdot \left( \a\nabla  P (\cdot,y) \right) = 0 & \mbox{in} & \ (s,\infty) \times \Rd,\\ 
&  P (0,x,y) = \delta_0(x-y) & \mbox{in} & \ \Rd,
\end{aligned}
\right.
\end{equation}
where $\delta_0$ is the Dirac delta distribution. Likewise, for each fixed $(t,x) \in\R \times\Rd$, the map $(s,y) \mapsto P^*(s,y,x):= P(-s,x,y)$ is the solution of the terminal-value problem for the backwards parabolic equation
\begin{equation} 
\label{e.Green.parab.backasswards}
\left\{ 
\begin{aligned}
& -\partial_t  P^* (\cdot,x) - \nabla \cdot \left( \a\nabla  P^* (\cdot,x) \right) = 0 & \mbox{in} & \ (-\infty,0) \times \Rd,\\ 
&  P^* (0,y,x) = \delta_0(x-y) & \mbox{in} & \ \Rd.
\end{aligned}
\right.
\end{equation}
Since we are working in the case of symmetric coefficients (i.e.,~$\a=\a^t$), we have that $P^*(s,y,x) = P(-s,y,x)$, that is, the parabolic Green function is symmetric in the variables~$x$ and~$y$. If the matrix~$\a$ were not symmetric, then $P^*(\cdot,x)$ would be the solution of~\eqref{e.Green.parab.backasswards} with~$\a^t$ in place of~$\a$. 

\smallskip

The Duhamel formula allows us to write a formula for the solution of the Cauchy problem: given $s\in \R$, we consider 
\begin{equation} 
\label{e.Cauchy}
\left\{ 
\begin{aligned}
& \partial_t  u - \nabla \cdot \left( \a(x)\nabla  u \right) = \nabla \cdot \mathbf{f} & \mbox{in} & \ (s,\infty) \times \Rd,\\ 
&  u (s,\cdot) = g & \mbox{in} & \ \Rd.
\end{aligned}
\right.
\end{equation}
The Duhamel formula for the solution $u$ is
\begin{equation} 
\label{e.duhamel}
u(t,x) = \int_{\Rd} g(y) P(t-s,x,y)\,dy  + \int_s^t \int_{\Rd} \mathbf{f}(s',y)\cdot \nabla_y P(t-s',x,y)\,dy\,ds'.
\end{equation}
It is sometimes convenient to refer to the rescaled parabolic Green function $P^\ep$, defined for $\ep>0$ by
\begin{equation*} \label{}
P^\ep(t,x,y):= \ep^{-d} P\left( \frac {t-s}{\ep^2},\frac x\ep,\frac y\ep  \right), 
\end{equation*}
which is the parabolic Green function for the operator $-\nabla \cdot \a\left( \tfrac \cdot\ep \right)\nabla$. 

\smallskip

Note that the parabolic Green function for $\ahom$ is given by the formula
\begin{equation}
\label{e.heatkernelahom}
\overline P (t,x):=(4\pi t)^{-\frac d2} (\det \ahom)^{\,-\frac12} \exp\left( - \,\frac{x\cdot \ahom^{\,-1} x}{4t}\right).
\end{equation}
It is the solution of
\begin{equation} 
\label{e.heatkernelahomeq}
\left\{
\begin{aligned}
& \partial_t \overline P  - \nabla \cdot \left(\ahom \nabla \overline P \right) = 0 & \mbox{in}  & \ \Rd \times(0,\infty), \\
& \overline P  = \delta_0 & \mbox{on} & \ \Rd \times \{0 \}. 
\end{aligned}
\right.
\end{equation}
The Nash-Aronson estimate 
\index{Nash-Aronson estimate}
states that, for every choice of~$\alpha \in (0,\Lambda^{-1})$, there exist~$0<\beta(d,\Lambda)<\infty$ and constants $0< c(d,\Lambda)\leq C(\alpha,d,\Lambda)<\infty$ such that, for every $t>0$ and $x,y \in\R^d\times\Rd$,
\begin{equation}
\label{e.aronson}
c t^{-\frac d2} \exp\left( -\beta \frac{\left|x-y\right|^2}{4t} \right) \leq  P (t,x,s,y)  \leq C t^{-\frac d2} \exp\left( -\alpha \frac{\left|x-y\right|^2}{4t} \right).
\end{equation}
As the parameters depend only on $d$, $\Lambda$ and the choice of $\alpha$, this estimate is independent of any structural assumption on the coefficients besides uniform ellipticity (such as our probabilistic assumptions). In this chapter (and the rest of the book) we need only the \emph{upper} bound in~\eqref{e.aronson}, a proof of which is presented in Appendix~\ref{a.NA} (and which relies on the Nash H\"older estimate for parabolic equations). We refer the reader to~\cite{FS} for an efficient and direct demonstration of~\eqref{e.aronson}, including the lower bound.

\smallskip

The \emph{elliptic Green function} for~$-\nabla \cdot \a\nabla$ is denoted by $G = G(x,y)$. We can characterize it for every $y\in\Rd$ as the solution of 
\begin{equation} 
\label{e.Green.ellip}
-\nabla \cdot \left( \a\nabla G(\cdot,y) \right) = \delta_0(\cdot-y) \quad \mbox{in} \ \Rd. 
\end{equation}
Recall that the elliptic and parabolic Green functions are related by the formula
\begin{equation} 
\label{e.intGammatoG}
G(x,y) = \int_0^\infty  P (t,x,y)\,dt. 
\end{equation}
Actually, we see from~\eqref{e.aronson} that the integral on the right of~\eqref{e.intGammatoG} is only finite in dimensions~$d>2$. Therefore we need to be slightly more precise in~$d=2$. In this case we can define
\begin{equation} 
\label{e.intGammatoG.d=2}
G(x,y) : = \lim_{T \to \infty}
\int_0^T \left( P (t,x,y) - \left( P(t,\cdot,y) \ast \Phi(1,\cdot) (0) \right)  \right) \,dt. 
\end{equation}
The convergence of this limit is a consequence of Nash's H\"older estimate for parabolic equations, as we show in~\eqref{e.Goscbound} below. The definition~\eqref{e.intGammatoG.d=2} will agree with the usual Green function definition, up to an additive constant. 

\smallskip

Using~\eqref{e.intGammatoG} (and~\eqref{e.intGammatoG.d=2} in the case~$d=2$), decay estimates for the parabolic Green function~$P$ immediately imply corresponding bounds for the elliptic Green function~$G$.

\smallskip

In the rest of this section, we present estimates on the decay of the \emph{gradient} of the parabolic and elliptic Green functions. Such estimates are not valid for general coefficient fields without further assumption. The extra ingredient in our setting is probabilistic, and it is encoded in the \emph{minimal scales} which prescribe the validity of the large-scale regularity theory. Before giving the bounds on the gradient of the Green functions, we pause to record some facts concerning these random scales~$\X_\sigma$, which appear not only in this section but also repeatedly throughout the chapter. We collect them in the following remark for ease of reference.

\smallskip

We emphasize that, throughout the rest of the chapter,~$\X_\sigma(x)$ denotes the random variable defined in the following remark. 

\begin{remark}[{Random variables $\X_\sigma(x)$}]
\label{r.minimalscale.PGF}
Let $\X_\sigma'$ denote, for each fixed~$\sigma \in (0,d)$, the maximum of the minimal scales $\X_\sigma$ appearing in Theorems~\ref{t.regularity} and~\ref{t.regularity.parab}, Lemma~\ref{l.correctorminbounds} and Proposition~\ref{p.fluxcorrectorests}. We let $\X_\sigma'(z)$ denote, for~$z\in\Zd$, its stationary extension. In particular, we have that $\X_\sigma'(z) \leq \O_\sigma(C)$ for a constant~$C(\sigma,d,\Lambda)<\infty$. We then (re-)define, for every~$x \in \Rd$, 
\begin{equation*}
\mathcal{X}_\sigma(x) := \sup_{z \in \Z^d} \left(2^d \X_{\frac{\sigma+d}{2}}'(z) - |x-z| \right).
\end{equation*}
We see then that $\X_\sigma(z)\geq \X_\sigma'(z)$ for every $z\in\Zd$ and that $x\mapsto \X_\sigma(x)$ is~$1$-Lipschitz continuous ($\P$--almost surely), which is very convenient. 
We will argue that, for some~$C(\sigma,d,\Lambda)<\infty$, 
\begin{equation*} \label{}
\mathcal{X}_s(x) \leq \O_{\sigma}(C).
\end{equation*}
Fix $r\in \left[ \tfrac 12,\infty\right)$ and observe that 
\begin{align} \notag
(2^d \X_{\frac{\sigma+d}{2}}'(z) - r)_+ & \leq \sum_{k=0}^\infty 2^k r \indc_{\left\{2^d \X_{\frac{\sigma+d}{2}}'(z) > 2^k r\right\}}
\\ \notag & \leq \sum_{k=0}^\infty 2^k r \left( \frac{2^d \X_{\frac{\sigma+d}{2}}'(z)}{2^k r} \right)^{\frac{\sigma+d}{2\sigma}}
\\ \notag & = \sum_{k=0}^\infty \left( 2^k r \right)^{- \frac{d-\sigma}{2\sigma}} \left( 2^d \X_{\frac{\sigma+d}{2}}'(z) \right)^{\frac{\sigma+d}{2\sigma}}
\\ \notag & \leq C r^{- \frac{d-\sigma}{2\sigma}} \left( \X_{\frac{\sigma+d}{2}}'(z) \right)^{\frac{\sigma+d}{2\sigma}}.
\end{align}
Thus, for $x \in \R^d$ and $z \in \Z^d$, we have
\begin{equation*} 
\left(1 + |x-z|\right)^{\frac{d-\sigma}{2\sigma}} (2^d \X_{\frac{\sigma+d}{2}}(z) - |x-z|)_+ \leq \O_\sigma(C).
\end{equation*}
We can now bound $\mathcal{X}_\sigma(x)$, using Lemma~\ref{l.sum-O} and taking $q = \frac{4d\sigma}{d-\sigma}$, as 
\begin{align} \notag 
\mathcal{X}_\sigma(x)& \leq \left( \sum_{z \in \Z^d } (2^d \X_{\frac{\sigma+d}{2}}'(z) - |x-z|)_+^q \right)^{\frac 1q}
\\ \notag & = \left( \sum_{z \in \Z^d } \left( \left(1 + |x-z|\right)^{\frac{d-\sigma}{2\sigma}} (2^d \X_{\frac{\sigma+d}{2}}(z) - |x-z|)_+ \right)^q \left(1 + |x-z|\right)^{- 2d} \right)^{\frac 1q}
\\ \notag & \leq \O_\sigma(C),
\end{align}
as claimed. Finally, we note that, since $\X_\sigma(z)\geq \X_\sigma'(z)$ for every $z\in\Zd$, each of the results mentioned above (Theorems~\ref{t.regularity} and~\ref{t.regularity.parab}, Lemma~\ref{l.correctorminbounds} and Proposition~\ref{p.fluxcorrectorests}) hold with the $\X_\sigma$ we have defined above replacing the one in the statements.
\end{remark}

\index{Green function!decay estimates|(}

\begin{theorem}
[Parabolic Green function decay bounds]
\label{t.gradientGF}
Fix $\sigma \in (0,d)$ and $\alpha \in (0,\Lambda^{-1})$. There exists a constant $C(\sigma,\alpha,d,\Lambda)<\infty$ such that, for every $x,y\in\Rd$ and $t>0$ with $\sqrt{t} \geq \X_\sigma(x)$, 
\begin{equation}
\label{e.gradPx}
\sup_{r\in \left[ \X_\sigma(x) , \sqrt{ t} \right] } 
\left\| \nabla_x P(t,\cdot,y) \right\|_{\underline{L}^2(B_r(x))}
\leq
Ct^{-\frac d2-\frac12}
\exp\left( -\alpha \frac{ \left|x-y\right|^2}{4t} \right),
\end{equation}
for every $x,y\in\Rd$ and $t>0$ with $\sqrt{t} \geq \X_\sigma(y)$, 
\begin{equation}
\label{e.gradPy}
\sup_{r\in \left[ \X_\sigma(y) , \sqrt{t} \right] } 
\left\| \nabla_y P(t,x,\cdot) \right\|_{\underline{L}^2\left(B_r(y)\right)}
\leq
Ct^{-\frac d2-\frac12}
\exp\left( -\alpha \frac{\left|x-y\right|^2}{4t} \right) ,
\end{equation}
and, for every $x,y\in\Rd$ and $t>0$ with $\sqrt{t} \geq \X_\sigma(x) \vee \X_\sigma(y)$, 
\begin{multline} 
\label{e.gradPxy}
\sup_{r' \in \left[ \X_\sigma(x) ,  \sqrt{t} \right], \, r \in \left[ \X_\sigma(y) ,  \sqrt{ t} \right]} 
\left\| \nabla_x \nabla_y P(t,\cdot,\cdot) \right\|_{\underline{L}^2(B_{r'}(x) \times B_{r}(y))}
\\
\leq
Ct^{-\frac d2-1}
\exp\left( -\alpha\frac{ \left|x-y\right|^2}{4t} \right).
\end{multline}
\end{theorem}

Like many of the statements we have presented previously in this book, the previous theorem gives \emph{deterministic} estimates which are valid on length scales larger than a certain \emph{random} scale. It is easy to see however that this implies estimates which are valid on \emph{every} scale (larger than the unit scale) but have \emph{random} right-hand sides. For instance, the first estimate~\eqref{e.gradPx} implies that, for every~$\sigma\in (0,2)$ and $\alpha \in \left( 0,\Lambda^{-1} \right)$, there exists~$C(\sigma,\alpha,d,\Lambda)<\infty$ such that, for every~$x,y\in\Rd$ and~$t\geq 1$, 
\begin{equation} 
\label{e.random.PGFbound}
\left\| \nabla_x P(t,\cdot,y) \right\|_{\underline{L}^2(B_1(x))} 
\leq 
\O_{\sigma} \left( Ct^{-\frac d2-\frac12} \exp\left( -\alpha \frac{|x-y|^2}{4t}  \right)  \right).
\end{equation}
This will be explained below in Remark~\ref{r.yesrandomGFRHS}.

\smallskip

We also mention that, similarly to the discussion just below the statement of Theorem~\ref{t.regularity}, the gradient bounds above should be thought of as ``almost pointwise'' estimates since we are integrating only over the small scales. In particular, under the additional assumption that the coefficients are uniformly H\"older continuous, then they can be combined with the classical Schauder estimates to obtain true pointwise bounds on $\nabla_xP$, $\nabla_yP$ and $\nabla_x\nabla_yP$.

\smallskip

We next present an analogue of the previous theorem for the elliptic Green function, which is an immediate consequence of Theorem~\ref{t.gradientGF} by formula~\eqref{e.intGammatoG} and an integration in time.

\begin{theorem}
[Elliptic Green function decay bounds]
There exists $C(d,\Lambda)<\infty$ such that, for every $x,y\in\Rd$ with $\frac12 |x-y|>\X_s(x)$,
\label{t.PGF.decay}
\begin{equation}
\label{e.gradGx}
\sup_{r\in \left[ \X_s(x), \frac12|x-y|\right]}
\left\| \nabla_x G(\cdot,y) \right\|_{\underline{L}^2(B_r(x))}
\leq C|x-y|^{1-d},
\end{equation}
for every $x,y\in\Rd$ with $\frac12 |x-y|>\X_s(y)$,
\begin{equation}
\label{e.gradGy}
\sup_{r\in \left[ \X_s(y), \frac12|x-y|\right]}
\left\| \nabla_y G(\cdot,y) \right\|_{\underline{L}^2(B_r(y))}
\leq C|x-y|^{1-d},
\end{equation}
and, for every $x,y\in\Rd$ with $\frac12 |x-y|>\X_s(x) \vee \X_s(y)$,
\begin{equation}
\label{e.gradGxy}
\sup_{r\in \left[ \X_s(x), \frac12|x-y|\right], \, r' \in  \left[ \X_s(y), \frac12|x-y|\right]}
\left\| \nabla_x\nabla_y G \right\|_{\underline{L}^2\left(B_r(y) \times B_{r'}(y)\right)}
\leq C|x-y|^{-d}.
\end{equation}
\end{theorem}

Theorems~\ref{t.gradientGF} and~\ref{t.PGF.decay} follow from the Lipschitz regularity developed in Chapter~\ref{c.regularity}. To apply Theorem~\ref{t.Lipschitz}, we use the Nash-Aronson estimate~\eqref{e.aronson} to bound the oscillation of the Green functions. We begin with the proof of Theorem~\ref{t.gradientGF}.

\begin{proof}
[{Proof of Theorem~\ref{t.gradientGF}}]
We first prove~\eqref{e.gradPx}. For notational convenience, for every $x,y\in\Rd$ and $t > s > 0$, we denote $P(t,x,s,y) := P(t-s,x,y)$.
An application of Theorem~\ref{t.Lipschitz.parab} yields that, for every $x,y\in\Rd$ and $t > s > 0$ with $\tfrac12\sqrt{t-s} \geq \X_s(x)$,
\begin{equation*}
\sup_{r \in \left[ \X_s(x),  \frac12 \sqrt{ t-s} \right] }
\left\| \nabla_x P(t,x,s,y) \right\|_{\underline{L}^2(Q_r(t,x))}
\leq
C(t-s)^{-\frac12}  \left\| P(\cdot,s,y)\right\|_{\underline{L}^2\left(Q_{ \frac12 \sqrt{ t-s}}(t,x)\right)}.
\end{equation*}
The Nash-Aronson estimate~\eqref{e.aronson} gives us the bound for the right side:
\begin{align*} \label{}
 \left\| P(\cdot,s,y)\right\|_{\underline{L}^2\left(Q_{ \frac12 \sqrt{ t-s}}(t,x)\right)}
&
\leq 
C
\sup_{(t',x')\in Q_{ \frac12 \sqrt{ t-s}}(t,x)}
\left( t'-s \right)^{-\frac d2}
\exp\left( -\frac{\alpha\left|x'-y\right|^2}{t'-s} \right)
\\ & 
\leq C\left( t-s \right)^{-\frac d2}
\exp\left( -\frac{\alpha \left|x-y\right|^2}{2(t-s)} \right).
\end{align*}
Combining the previous two displays yields~\eqref{e.gradPx} after redefining~$\alpha$. 

\smallskip

The statement~\eqref{e.gradPy} is an immediate consequence of \eqref{e.gradPx}, by the symmetry $P(t,x,y) = P(t,y,x)$.

\smallskip

To prove~\eqref{e.gradPxy} we also argue in the same way after noticing that, for each index~$j\in\{1,\ldots,d\}$, the function $(s,y) \mapsto \partial_{x_i} P( t,x,s,y)$ also solves the backwards equation. Theorem~\ref{t.Lipschitz.parab} (see Remark~\ref{r.Lipschitz.parab2}) give us that, for every $(t,x)\in (s,\infty)\times\Rd$ and $r' \in \left[ \X_s(y) ,\tfrac12 (t-s)^{\frac12} \right]$,
\begin{align*}
\left\| \nabla_y \partial_{x_i} P(t,x,\cdot)  \right\|_{\underline{L}^2\left( Q_{r'}(s,y) \right)} 
\leq C(t-s)^{-\frac12}  \left\| \partial_{x_i} P(t,x,\cdot)  \right\|_{\underline{L}^2\left( Q_{\frac12 \sqrt{t-s}}(s,y) \right)}.
\end{align*}
Fixing $(t,x)$ and $r>0$ with $\X_s(x)\leq r \leq \frac12 (t-s)^{\frac12}$, integrating over $Q_{r}(t,x)$ and applying~\eqref{e.gradPx} gives
\begin{align*}
\lefteqn{
\left\| \nabla_y \partial_{x_i} P  \right\|_{\underline{L}^2\left( Q_r(t,x) \times Q_{r'}(s,y) \right)}^2
} \qquad & 
\\ & 
=
\fint_{Q_{r}(t,x)} \left\| \nabla_y \partial_{x_i} P(t',x',\cdot)  \right\|_{\underline{L}^2\left( Q_{r'}(s,y) \right)}^2\,dt'\,dx'
\\ & 
\leq 
C(t-s)^{-1} \fint_{Q_r(t,x)}  \left\| \partial_{x_i} P(t',x',\cdot)  \right\|_{\underline{L}^2\left( Q_{\frac12 \sqrt{t-s}}(s,y) \right)}^2 \,dt'\,dx'
\\ & 
= C(t-s)^{-1} \fint_{Q_{r'}(s,y)} \left(  \fint_{Q_r(t,x)} \left| \partial_{x_i} P(t',x',s',y') \right|^2 \, dt'\,dx' \right) \,ds'\,dy'
\\ & 
\leq 
C(t-s)^{-1} \fint_{Q_{r'}(s,y)} \left(  C\left( t-s' \right)^{-d-1}
\exp\left( -\frac{2\alpha \left|x-y'\right|^2}{(t-s')} \right) \right) \, ds'\, dy'
\\ &
\leq C (t-s)^{-d-2} \exp\left( -\frac{2\alpha \left|x-y'\right|^2}{(t-s')} \right).
\end{align*}
Summing over $i\in\{1,\ldots,d\}$ completes the proof of~\eqref{e.gradPxy}. 
\end{proof}

\begin{exercise}
\label{ex.ellipGF}
Deduce Theorem~\ref{t.PGF.decay} from~\eqref{e.intGammatoG} and Theorem~\ref{t.gradientGF}.
\end{exercise}

\index{Nash-Aronson estimate}
\index{De Giorgi-Nash estimate}

In addition to Exercise~\ref{ex.ellipGF}, we offer another, perhaps more direct proof which essentially just interchanges the order of the time integration and application of the~$C^{0,1}$ estimate. Let us first show that the Nash-Aronson estimate~\eqref{e.aronson} for the parabolic Green function can be translated into oscillation bounds on the elliptic Green function via the formula~\eqref{e.intGammatoG}. We obtain that there exist constants $0<c(d,\Lambda) \leq C(d,\Lambda)<\infty$ such that, for every $x\in\Rd\setminus \{0\}$,
\begin{equation} 
\label{e.ellipAronson}
c\overline{G}(x) \leq G(x,0) \leq C \overline{G}(x). 
\end{equation}
Indeed, for every $x \in B_{2r}\setminus B_r$, 
\begin{align*} \label{}
G(x,0)
& 
= \int_0^\infty  P (t,x,0)\,dt
\\ & 
\leq C \int_0^\infty  t^{-\frac d2} \exp\left( -\frac{\alpha\left|x\right|^2}{t} \right)  \,dt
\\ & 
= C \left( \frac{\alpha |x|^2}{r^2} \right)^{1-\frac d2}  \int_0^\infty  t^{-\frac d2} \exp\left( -\frac{r^2}{4t} \right)  \,dt
\\ &
\leq C \overline G(x,0). 
\end{align*}
The lower bound is similar. 

\smallskip

We need a more refined estimate in $d=2$. In order to obtain gradient estimates, the quantity we need to estimate is $\osc_{B_{2r}\setminus B_r} G(\cdot,0)$ and, while we can get an estimate on this quantity from the upper bound in~\eqref{e.ellipAronson}, the decay will be optimal in $r$ only in~$d>2$. In dimension $d=2$, this oscillation is bounded by $C$, as we will discover below, and yet~\eqref{e.ellipAronson} gives us only that it is bounded by $\left| \log r \right|$. 
To get a better estimate, we return to~\eqref{e.aronson} and combine it with Nash's H\"older estimate to obtain that, for every $r>0$, $y\in B_{2r}\setminus B_r$ and $t\geq r^2$,
\begin{equation*} \label{}
\osc_{x\in B_{r/2}(y) }
 P \left( t,x,0 \right) 
\leq C \left( \frac{r^2}{t} \right)^{\alpha}  t^{-\frac d2}  = C r^{2\alpha} t^{-\frac d2-\alpha}. 
\end{equation*}
For small times $t\lesssim r^2$, we will still use the bound~\eqref{e.aronson}. We deduce that, for every $r>0$, $y\in B_{2r}\setminus B_r$ and $t\geq r^2$,
\begin{align*} \label{}
\osc_{x\in B_{r/2}(y) } G(x,0) 
&
\leq \int_0^\infty \osc_{x\in B_{r/2}(y) }
 P \left( t,x,0 \right)  \,dt
\\ &
\leq 
C \int_0^{r^2}  t^{-\frac d2} \exp\left( -\frac{cr^2}{t} \right)  \,dt
+ C r^{2\alpha}  \int_{r^2}^\infty t^{-\frac d2-\alpha} \,dt.
\end{align*}
Both integrals on the right side are bounded by $Cr^{2-d}$ and so we have proved that, with~$r$ and~$y$ as above, 
\begin{equation*} \label{}
\osc_{x\in B_{r/2}(y) } G(x,0) 
\leq 
Cr^{2-d}. 
\end{equation*}
Covering $B_{2r}\setminus B_r$ with $C$ many such balls yields, for every $r>0$,
\begin{equation} 
\label{e.Goscbound}
\osc_{x\in B_{2r}\setminus B_r} G(x,0)  
\leq Cr^{2-d}. 
\end{equation}
Notice that, in $d=2$, the estimate~\eqref{e.Goscbound} is stronger than the upper bound in~\eqref{e.ellipAronson}, because we can sum the former over dyadic scales to obtain the latter.

\begin{proof}[{Proof of Theorem~\ref{t.PGF.decay}}]
We prove only~\eqref{e.gradGx}. Theorem~\ref{t.Lipschitz},~\eqref{e.Goscbound} and the inclusion 
\begin{equation*} \label{}
B_{\frac12|x-y|}(x) 
\subseteq 
B_{|x-y|}(y) \setminus B_{\frac12|x-y|}(y)
\end{equation*}
yield
\begin{align*}
\lefteqn{
\sup_{r\in \left[ \X_s(x), \frac12|x-y|\right]}
\left\| \nabla_x G(\cdot,y) \right\|_{\underline{L}^2(B_{r}(x))}
} \qquad & \\
& \leq C|x-y|^{-1} \Big\| G(\cdot,y)  - \left( G(\cdot,y) \right)_{B_{\frac12|x-y|}(x)} \Big\|_{\underline{L}^2\left( B_{\frac12|x-y|}(x) \right)}
\\ &
\leq C|x-y|^{-1} \osc_{B_{|x-y|} \setminus B_{\frac12|x-y|}} G(\cdot,y)
\\ &
\leq C|x-y|^{1-d}.  \qedhere
\end{align*}
\end{proof}

\begin{remark}
\label{r.yesrandomGFRHS}
In this remark, we demonstrate, as promised, how to obtain~\eqref{e.random.PGFbound} from the statement of Theorem~\ref{t.gradientGF}. Analogous bounds to~\eqref{e.random.PGFbound} can be proved for $\nabla_y P$ and $\nabla_x\nabla_yP$ via a similar computation.

\smallskip

We first use Lemma~\ref{l.slicesvsaverages.PAR}, the parabolic Caccioppoli inequality (Lemma~\ref{l.cacciopp.this}) and the Nash-Aronson estimate~\eqref{e.aronson} and take $\ep(\alpha,\Lambda)>0$ sufficiently small to obtain, for every time~$t$ such that $\sqrt{t} \geq 2\ep^{-1}$, the crude deterministic bound 
\begin{align*} \label{}
\left\| \nabla_x P(t,\cdot,y) \right\|_{\underline{L}^2(B_1(x))}
& 
\leq C t^{\frac d4} \left\| \nabla_x P(t,\cdot,y) \right\|_{\underline{L}^2(B_{\ep \sqrt{t}/2}(x))}
\\ & 
\leq 
C t^{\frac d4}
\left\| \nabla_x P(\cdot,y) \right\|_{\underline{L}^2(Q_{\ep \sqrt{t}}(t,x))}
\\ & 
\leq
Ct^{-\frac d4-\frac12} \exp\left( -\alpha \frac{|x-y|^2}{4t}  \right) .
\end{align*}
Using then that 
\begin{equation*} \label{}
\indc_{\{ \sqrt{t} \leq \X_{\sigma d/2} (x) \} } 
\leq \O_\sigma \left( C t^{-\frac d4} \right),
\end{equation*}
we therefore obtain
\begin{align*} \label{}
\left\| \nabla_x P(t,\cdot,y) \right\|_{\underline{L}^2(B_1(x))}
\indc_{\{ \sqrt{t} \leq \X_{\sigma d/2} (x) \} } 
&
\leq
Ct^{-\frac d4-\frac12} \exp\left( -\alpha \frac{|x-y|^2}{4t}  \right) \indc_{\{ \sqrt{t} \geq \X_{\sigma d/2} (x) \} } 
\\ & 
\leq
\O_{\sigma} \left( Ct^{-\frac d2-\frac12} \exp\left( -\alpha \frac{|x-y|^2}{4t}  \right)  \right).\end{align*}
On the other hand, by the statement of Theorem~\ref{t.gradientGF}, namely~\eqref{e.gradPx}, we have
\begin{align*} \label{}
\left\| \nabla_x P(t,\cdot,y) \right\|_{\underline{L}^2(B_1(x))}
\indc_{\{ \sqrt{t} \geq \X_{\sigma d/2} (x) \} } 
&
\leq 
C \X_{\sigma d/2}(x)^{\frac d2} \left\| \nabla_x P(t,\cdot,y) \right\|_{\underline{L}^2(B_{\X_{\sigma d/2}(x)}(x))}
\\ & 
\leq 
C \X_{\sigma d/2}(x)^{\frac d2} t^{-\frac d2 -\frac12} \exp\left( -\alpha \frac{|x-y|^2}{4t}  \right). 
\end{align*}
Since $ \X_{\sigma d/2}(x)^{\frac d2} \leq \O_\sigma(C)$, this completes the proof of~\eqref{e.random.PGFbound} for times~$t\geq C(\alpha,\Lambda)$. For~$t$ satisfying~$1\leq t\leq C(\alpha,\Lambda)$, no argument is required since Lemmas~\ref{l.cacciopp.this},~\ref{l.slicesvsaverages.PAR} and the Nash-Aronson bound give us~$\left\| \nabla_x P(t,\cdot,y) \right\|_{\underline{L}^2(B_1(x))} \leq C$. 
\end{remark}

We conclude this section with a technical lemma which formalizes a computation that we need to perform several in this chapter and the next one. Since this computation comes in several variations, the statement of this lemma must (unfortunately) be rather general. 

\begin{lemma} 
\label{l.PGF.layercake}
Fix $\sigma \in (0,d)$, $p \in [2,\infty]$, $\alpha \in (0,\Lambda^{-1})$, and $\beta \in (\alpha,\Lambda^{-1})$.  There exists a constant $C(\sigma,p,\alpha,\beta,d,\Lambda)< \infty$ such that, for every $x,y \in \Rd$, $t \in (0,\infty)$, $s \in \left[ \X_\sigma(y)^2,\infty\right)$ and $g \in L_{\textrm{loc}}^p(\R^d)$, we have
\begin{multline} \label{e.PGF.layercake}
\left| \left( \Phi\left( \frac{s}{\beta} , \cdot\right) \ast \left( |g(\cdot)| \left( P(t,x,\cdot) + t^{\frac 12} \left| \nabla_y P(t,x,\cdot) \right| \right) \right) \right)(y) \right| 
\\ 
 \leq C  \Psi_{p}[g] (t+s)^{-\frac d2} \exp \left( - \alpha \frac{|x-y|^2}{4(t+s)} \right),
\end{multline}
where $\Psi_{p}$ is defined, for $p \in [2,\infty)$, as 
\begin{multline*} 
\Psi_{p}[g] :=  
  \left( \left(\frac{t+s}{t} \right)^{\frac d2} \int_{1}^\infty r^{d+1} \exp\left(- p(\beta-\alpha) \frac{r^2}{8} \right) \left\| g \right\|_{\underline{L}^{p} \left( B_{r \sqrt{s}}(y) \right)}^p \, dr \right)^{\frac 1p}  \\
\wedge  \left(  \left(\frac{t+s}{s} \right)^{\frac d2}  \int_{1}^\infty r^{d+1} \exp\left(- p(\beta-\alpha) \frac{r^2}{8} \right) \left\| g \right\|_{\underline{L}^{p} \left( B_{r \sqrt{t}}(x) \right)}^p \, dr \right)^{\frac 1p}. 
\end{multline*}
and, for $p=\infty$, as
\begin{multline*} 
\Psi_{\infty}[g] :=  
 \sup_{r\in [1,\infty)} \left( \exp\left(- (\beta-\alpha) \frac{r^2}{8} \right) \left\| g \right\|_{L^{\infty} \left( B_{r \sqrt{s}}(y)\right) } \right) \\
\wedge  \sup_{r\in [1,\infty)} \left( \exp\left(- (\beta-\alpha) \frac{r^2}{8} \right) \left\| g \right\|_{L^{\infty} \left( B_{r \sqrt{t}}(x)\right) } \right). 
\end{multline*}

\end{lemma}

\begin{proof}
Throughout the proof, we fix $\sigma \in (0,d)$, $p \in [2,\infty]$, $\alpha \in (0,\Lambda^{-1})$, and $\beta \in (\alpha,\Lambda^{-1})$. We decompose the proof into four steps.

\smallskip

\emph{Step 1.}
We show that there exists a constant $C(\beta,d)  \in (0,\infty)$ such that for every $f \in L_{\textrm{loc}}^1(\R^d)$, $t'>0$ and $z \in\Rd$, 
\begin{equation} 
\label{e.PGF.layercake.1}
\left( \Phi\left( \frac{t'}{\beta} , \cdot\right) \ast f(\cdot) \right)(z)
=
C \int_{0}^\infty r^{d+1} \exp\left(- \beta \frac{r^2}{4} \right)  \fint_{ B_{r \sqrt{t'}}(z) }  f(z') \, dz'   \, dr.
\end{equation}
To see this, write first 
\begin{equation*} 
\Phi\left(  \frac{t'}{\beta} , y-z' \right) = C (t')^{-\frac d2} \exp\left(- \beta \frac{|z-z'|^2}{4t'}\right) = C \int_{0}^\infty \indc_{ \left\{ (t')^{-\frac d2} \exp\left(- \beta \frac{|z-z'|^2}{4t'} \right) > \lambda \right\} } \, d\lambda.
\end{equation*}
The change of variables  
\begin{equation*} 
\lambda = (t')^{-\frac d2} \exp\left(- \beta \frac{r^2}{4}\right)\,, \quad d \lambda = - \frac{\beta}{2} r  (t')^{-\frac d2} \exp\left(- \beta \frac{r^2}{4} \right) \, dr,
\end{equation*}
leads to
\begin{equation*} 
(t')^{-\frac d2} \exp\left(- \beta \frac{|z-z'|^2}{4t'}\right) = \frac{\beta}{2}  \int_{0}^\infty  r  \exp\left(- \beta \frac{r^2}{4} \right) (t')^{-\frac d2}  \indc_{ \left\{ |z'-z| < r \sqrt{s}  \right\} } \, dr.
\end{equation*}
Applying Fubini's theorem yields~\eqref{e.PGF.layercake.1}.

\smallskip

\emph{Step 2.} We show that there exists a constant $C(\alpha,\beta,d) < \infty$ such that, for every $s,t > 0$ and $x,y,z \in \R^d$, we have
\begin{equation} \label{e.PGF.layercake.2}
\left\| \frac{\Phi\left( \frac{s}{\beta} , y - \cdot  \right)}{\Phi\left( \frac{s}{\alpha} , y - z \right)} \right\|_{L^\infty(B_{\sqrt{s}}(z))} +  
\left\| \frac{\Phi\left( \frac{t}{\beta}, x-z  \right)}{\Phi\left( \frac{t}{\alpha}, x - \cdot \right)}  \right\|_{L^\infty(B_{\sqrt{t}}(z))} \leq C .
\end{equation}
To see this, observe that, for $\tilde z \in B_{\sqrt{t}}(z)$, we have 
\begin{equation*} 
\left| x- \tilde z   \right|^2 \leq \frac{\beta}{\alpha} \left| x -z \right|^2 + \frac{\beta}{\beta-\alpha} \left|\tilde z - z \right|^2 \leq \frac{\beta}{\alpha}   \left| x - z \right|^2 +  \frac{\beta}{\beta-\alpha}  t
\end{equation*}
and
\begin{equation*} 
\left| y - \tilde z   \right|^2 \geq \frac{\alpha}{\beta}\left| y - z\right|^2 - \frac{\alpha}{\beta-\alpha}  \left|\tilde z - z \right|^2 \geq \frac{\alpha}{\beta} \left| x - z \right|^2 - \frac{\alpha}{\beta-\alpha} t .
\end{equation*}
Therefore,~\eqref{e.PGF.layercake.2} follows easily with $C = 2 \left(\frac{\beta}{\alpha}\right)^{\frac d2} \exp \left(\frac{\beta \alpha}{\beta-\alpha}\right)$.

\smallskip

\emph{Step 3.} 
For the rest of the proof, we also fix $x,y \in \R^d$, $t \in (0,\infty)$, and $s \in \left[ \X_\sigma(y)^2,\infty\right)$.
In this step, we prove the version of~\eqref{e.PGF.layercake} when the ball $B_{r \sqrt{s}}(y)$ appears in the estimate. We consider only the term involving $t^{\frac 12}\nabla_y P(t,x,\cdot)$, since it is easy to check that the term involving $ P(t,x,\cdot)$ has the same estimate using directly the Nash-Aronson bound. By~\eqref{e.PGF.layercake.1}, we obtain, for $q = \tfrac{p}{p-1} \in [1,2]$, 
\begin{align} \label{e.PGF.layercake.3.1}
\notag  
\lefteqn{\left( \Phi\left(  \frac{s}{\beta} , \cdot\right) \ast  \left( |g(\cdot)| \left| \nabla_y P(t,x,\cdot) \right| \right) \right)(y)} \quad  &  
\\ 
\notag & \leq C  \int_{1}^\infty r^{d+1} \exp\left(- \beta \frac{r^2}{4} \right) \left\| g\right\|_{\underline{L}^{p} \left( B_{r \sqrt{s}}(y)\right)} \left\|  \nabla_y P(t,x,\cdot) \right\|_{\underline{L}^q\left( B_{r \sqrt{s}}(y)\right)}  \, dr 
\\ 
\notag & \leq C  \left(  \int_{1}^\infty r^{d+1} \exp\left(- q \alpha \frac{r^2}{4}   \right)  \left\|  \nabla_y P(t,x,\cdot) \right\|_{\underline{L}^{q} \left( B_{r \sqrt{s}}(y)\right)}^q \, dr \right)^{\frac 1q}  
\\ 
 & \qquad \times \left(  \int_{1}^\infty r^{d+1} \exp\left(- p(\beta-\alpha) \frac{r^2}{4}   \right)  \left\| g \right\|_{\underline{L}^{p} \left( B_{r \sqrt{s}}(y)\right)}^p \, dr \right)^{\frac 1p}.
\end{align}
The interpretation for the second term on the right in the case $p=\infty$ is as in~\eqref{e.PGF.layercake}. 
We now have, for $R \geq \sqrt{s} \geq  \X_\sigma(y)$, that
\begin{align} \label{e.PGF.layercake.3.2}
\notag 
 \left\| \nabla_{y} P(t,x,\cdot) \right\|_{\underline{L}^q \left(B_{R}(y)\right)}^q 
 & \leq C \fint_{B_{R}(y)} \left\| \nabla_{y} P(t,x,\cdot) \right\|_{\underline{L}^q\left(B_{R \wedge \sqrt{t}}(z)\right)}^q \, dz
 \\ \notag & \leq C \fint_{B_{R}(y)}   \left\| \nabla_{y} P(t,x,\cdot) \right\|_{\underline{L}^2\left(B_{R \wedge \sqrt{t}}(z)\right)}^{q}  \, dz
 \\  & \leq C t^{-\frac q2-(q-1) \frac d2} \fint_{B_{R}(y)} t^{- \frac d2} \exp \left( - q \alpha \frac{|z-x|^2}{4t} \right) \, dz.
\end{align}
Indeed, the first inequality is a triviality, and the second one we get from H\"older's inequality. The third inequality follows, in the case $\sqrt{t} \le R$, from the Caccioppoli estimate (Lemma~\ref{l.parcacc.NA}) after integration in time, together with \eqref{e.gradientP.NA}, and in the case $R \le \sqrt{t}$,  from~\eqref{e.gradPy}. We also made use of the fact that since $y \mapsto \X_\sigma(y)$ is 1-Lipschitz continuous, we have that
\begin{equation*} \label{}
\sup_{z \in B_R(y)} \X_\sigma(z) \leq  \X_\sigma(y) + R \leq 2R.
\end{equation*}
Using once more~\eqref{e.PGF.layercake.1} and the semigroup property yields 
\begin{align} 
\notag  
\lefteqn{\int_{1}^\infty r^{d+1} \exp\left(- q \alpha \frac{r^2}{4}   \right)  \left\|  \nabla_y P(t,x,\cdot) \right\|_{\underline{L}^{q} \left( B_{r \sqrt{s}}(y)\right)}^q \, dr} \quad  &  
\\ 
\notag & \leq C t^{-\frac q2-(q-1) \frac d2}  \int_{1}^\infty r^{d+1} \exp\left(- q \alpha \frac{r^2}{4}   \right)  \left\|  t^{- \frac d2} \exp \left( - q \alpha \frac{|\cdot-x|^2}{4t} \right)  \right\|_{\underline{L}^{1} \left( B_{r \sqrt{s}}(y)\right)} \, dr
\\ 
\notag & \leq C \left( t^{-\frac 12} \left( \frac{t+s}{t}\right)^{ \frac{d}{2p}} (t+s)^{-\frac d2} \exp \left( -  \alpha \frac{|z-x|^2}{4(t+s)} \right) \right)^{q}.
\end{align}
Combining this with~\eqref{e.PGF.layercake.3.1} gives~\eqref{e.PGF.layercake}. 

\smallskip

\emph{Step 4.}  We finally prove~\eqref{e.PGF.layercake} when now we have the ball $B_{r \sqrt{t}}(x)$ instead of $B_{r \sqrt{s}}(y)$ appearing in the estimate.  Consider again only the term involving $t^{\frac 12}\nabla_y P(t,x,\cdot)$. We set, for $z \in \R^d$, 
\begin{equation*} 
\tilde g(z) := \exp \left(- (\beta-\alpha')\frac{|z-y|^2}{4s}\right)|g(z)| \quad \mbox{and} \quad h(z) := \frac{\Phi\left(  \frac{s}{\alpha'} , y -z \right)}{ \Phi\left(  \frac{t}{\beta} , x -z \right)}\left| \nabla_y P(t,x,z) \right| ,
\end{equation*}
so that 
\begin{equation*} 
 \Phi\left(  \frac{s}{\beta} , y -z \right) |g(z)| \left| \nabla_y P(t,x,z) \right| =  C \Phi\left(  \frac{t}{\beta} , x -z \right)  \tilde g(z) h(z).
\end{equation*}
Then, as in~\eqref{e.PGF.layercake.3.1}, we get by H\"older's inequality that
\begin{align} \label{e.PGF.layercake.4.1}
\notag  
\lefteqn{\left( \Phi\left(  \frac{s}{\beta} , \cdot\right) \ast  \left( |g(\cdot)| \left| \nabla_y P(t,x,\cdot) \right| \right) \right)(y)} \quad  &  
\\ 
\notag & \leq C  \left(  \int_{1}^\infty r^{d+1} \exp\left(- q \alpha \frac{r^2}{4}   \right)  \left\|  h \right\|_{\underline{L}^{q} \left( B_{r \sqrt{t}}(x)\right)}^q \, dr \right)^{\frac 1q}  
\\ 
 & \qquad \times \left(  \int_{1}^\infty r^{d+1} \exp\left(- p(\beta-\alpha) \frac{r^2}{8}   \right)  \left\| \tilde g \right\|_{\underline{L}^{p} \left( B_{r \sqrt{t}}(x)\right)}^p \, dr \right)^{\frac 1p}.
\end{align}
Following the computation in~\eqref{e.PGF.layercake.3.2} and~\eqref{e.PGF.layercake.2}, we have, for $\beta' := \frac12 \left( \Lambda^{-1} + \beta\right)$, $\alpha' = \frac12 \left( \alpha + \beta\right)$ and $R \geq \sqrt{t}$, that
\begin{align} \label{e.PGF.layercake.4.2}
\notag 
 \left\| h \right\|_{\underline{L}^q \left(B_{R}(x)\right)}^q 
 & \leq C \fint_{B_{R}(y)} \left\|\frac{\Phi\left( \frac{s}{\beta} , y - \cdot  \right)}{\Phi\left( \frac{t}{\beta} , x - \cdot \right)}   \right\|_{L^\infty\left(B_{\sqrt{t}/4}(z)\right)}^q \left\| \nabla_{y} P(t,x,\cdot) \right\|_{\underline{L}^q\left(B_{\sqrt{t}/4}(z)\right)}^q \, dz
 \\ \notag &
  \leq C \fint_{B_{R}(x)}  \left( \frac{\Phi\left( \frac{s}{\alpha} , y - z \right)}{\Phi\left( \frac{t}{ \beta} , x - z \right)}\right)^q   \left\| \nabla_{y} P(t,x,\cdot) \right\|_{\underline{L}^2\left(Q_{\sqrt{t}/2}(t,z)\right)}^q   \, dz
 \\  & \leq C t^{-\frac q2} s^{-(q-1) \frac d2}  \fint_{B_{R}(x)} s^{- \frac d2} \exp \left( - q \alpha \frac{|z-y|^2}{4s} \right) \, dz.
\end{align}
The proof is now concluded by the same argument as in the end of Step 3. 
\end{proof}

\index{Green function!decay estimates|)}

\section{Homogenization of the Green functions}
\label{s.greenie.base}

In this section, we prove a quantitative estimate for the convergence of $P$ to $\overline{P}$ at large times which is suboptimal in the scaling (i.e., the rate is roughly $O(t^{-\ep})$) but optimal in stochastic integrability. The precise statement is given in the following theorem, which is the main result of the section. We remark that the results of this section depend only on the bounds for the first-order correctors and flux correctors which are consequences of the suboptimal estimates proved in Chapters~\ref{c.two} and~\ref{c.regularity}, and summarized in Lemma~\ref{l.correctorminbounds}, Proposition~\ref{p.correctorsbasecase} (for the first-order correctors) and~\eqref{e.fluxcorrectorgradbound} and~\eqref{e.correctorminbounds0.flux}  (for the flux correctors). In particular, we are not using here the more involved, precise theory developed in Chapter~\ref{c.A1} which gives the optimal quantitative bounds on the first-order correctors. 

\smallskip

For each~$t>0$ and $x,y \in \R^d$, we define the two-scale expansion of~$\overline{P}$ (in the~$x$-variable) by 
\begin{equation} 
\label{e.def.H}
H(t,x,y) := \overline{P}(t,x-y) + \sum_{k=1}^d \left( \phi_{e_k}(x) - \left( \phi_{e_k} \ast \Phi(t,\cdot) \right)(y) \right)  \partial_{x_k}\overline{P}(t,x-y).
\end{equation}
As in Chapter~\ref{c.twoscale} (for example, see~\eqref{e.phiep.def} and~\eqref{e.wepdef}), we have subtracted the ``constant'' $\left( \phi_{e_k} \ast \Phi(t,\cdot) \right)(y)$ from the first-order corrector~$\phi_{e_k}$ in order to make the definition of~$H$ unambiguous (at least in dimension~$d=2$) and because this is actually needed in order for the bounds we prove below to be optimal in stochastic integrability (even in~$d>2$). This function $y\mapsto \left( \phi_{e_k} \ast \Phi(t,\cdot) \right)(y)$ is obvious not constant, but it has only small, low frequency oscillations which do not harm the analysis. 

\index{Green function!homogenization estimates|(}
\index{Green function!two-scale expansion estimates|(}
\index{Berry-Esseen theorem}

\begin{theorem}[Homogenization of the parabolic Green function]
\label{t.PGF.basecase}
Fix $\sigma \in (0,d)$. There exist~$\delta(d,\Lambda) >0$ and, for each $y\in\Rd$, a random variable~$\Y_\sigma(y)$ satisfying 
\begin{equation} 
\label{e.PGD.basecase.Ysig}
\Y_\sigma(y) \leq \O_\sigma\left( C(\sigma,d,\Lambda) \right)
\end{equation}
and, for each $\alpha \in (0,\Lambda^{-1})$, a constant~$C(\alpha,\sigma,d,\Lambda)<\infty$ such that, 
for every $x,y\in\Rd$ and $\sqrt{t} \geq \Y_\sigma(y)$,
\begin{multline}
\label{e.PGF.basecase}
\left| H(t,x,y)-P(t,x,y) \right|  + \left| P(t,x,y) - \overline{P}(t,x-y) \right| 
\\ \leq 
Ct^{-\delta(d-\sigma)} t^{- \frac d2}  \exp\left(- \alpha \frac{|x-y|^2}{4t}\right)
\end{multline}  
and additionally, if  $r \in \left[  \Y_\sigma(x), \frac12 \sqrt{t} \right]$, then
\begin{equation}
\label{e.PGF.basecase.grad}
 \left\| \nabla_x P(t,\cdot,y) - \nabla_x H(t,\cdot,y) \right\|_{\underline{L}^2(B_r(x))}
\leq 
Ct^{-\delta(d-\sigma)} t^{ - \frac 12 - \frac d2}  \exp\left(- \alpha \frac{|x-y|^2}{4t}\right).
\end{equation}
\end{theorem}

The first estimate~\eqref{e.PGF.basecase} is a true pointwise estimate, which is possible due to our use of (scalar) De Giorgi-Nash~$L^\infty$~estimates in this chapter. The second estimate~\eqref{e.PGF.basecase.grad} for the gradients should also be considered as ``almost pointwise'', similarly to the~$C^{0,1}$-type estimate~\eqref{e.Lipschitz2} and as explained below~\eqref{e.Lipschitz2}, since the integration is only over the correlation length scale. In particular, if we assume the coefficients to be H\"older continuous on the microscopic scale, for instance, then this two-scale expansion estimate can be upgraded to a true pointwise bound in the same manner as~\eqref{e.schaudertopointwise}. A two-scale expansion-type estimate for the mixed second derivatives $\nabla_x\nabla_y P$ is also obtained below in Exercise~\ref{exer.mixedsecond}. 

\smallskip

In view of~\eqref{e.intGammatoG}, we can integrate the result of Theorem~\ref{t.PGF.basecase} in time to obtain the following quantitative homogenization result for the elliptic Green function~$G$. We denote the two-scale expansion of $\overline{G}$ by 
\begin{equation} 
\label{e.def.F}
F(x,y) := \overline{G}(x-y) + \sum_{k=1}^d \left( \phi_{e_k}(x) - \left( \phi_{e_k} \ast \Phi(|x-y|^2,\cdot) \right)(y) \right)  \partial_{x_k}\overline{G}(x-y).
\end{equation}

\begin{corollary}
[Homogenization of the elliptic Green function]
Fix~$\sigma\in (0,d)$ and $\alpha\in \left( 0,\Lambda^{-1} \right)$ and let~$\Y_\sigma(y)$ be as in Theorem~\ref{t.PGF.basecase}. There exist constants~$\delta(d,\Lambda) >0$ and~$C(\sigma,\alpha,d,\Lambda)<\infty$ such that,  
for every~$x,y\in\Rd$ with $|x-y|\geq \Y_\sigma(y)$, 
\label{c.GtoGbar0}
\begin{equation} 
\label{e.GtoGbar0}
\left| G(x,y) - \overline{G}(x-y) \right| 
\leq
C \left| x-y \right|^{-\delta} \left| x-y \right|^{2-d}
\end{equation}
and, if $r \in \left[  \Y_\sigma(x), \frac12 |x-y| \right]$, then we also have
\begin{equation} 
\label{e.GtoGbar.nabla0}
\left\| \nabla_x G(x,y) - \nabla_x F(x,y) \right\|_{\underline{L}^2(B_r(x))} 
\leq
C \left| x-y \right|^{-\delta} \left| x-y \right|^{1-d}.
\end{equation}
\end{corollary}

\index{Green function!two-scale expansion estimates|)}

For future reference, we observe that, by the symmetry of the parabolic Green function, the first estimate~\eqref{e.PGF.basecase} of Theorem~\ref{t.PGF.basecase} also holds if~$\sqrt{t} \geq \Y_\sigma(x) \wedge \Y_\sigma(y)$. That is, for every such $t$, we have
\begin{equation} 
\label{e.PGF.basecase2}
\left| P(t,x,y) - \overline{P}(t,x-y) \right| \leq  Ct^{-\delta(d-\sigma)} t^{- \frac d2}  \exp\left(- \alpha \frac{|x-y|^2}{4t}\right).
\end{equation}

\smallskip

Most of the effort in the proof of Theorem~\ref{t.PGF.basecase} will be focused on obtaining the first estimate~\eqref{e.PGF.basecase}. The second estimate~\eqref{e.PGF.basecase.grad} follows relatively easily from it and the $C^{1,1}$-type regularity estimate (the case $k=1$ of the parabolic version stated Theorem~\ref{t.regularity.parab}). The argument for~\eqref{e.PGF.basecase} is similar in philosophy to the proof of Theorem~\ref{t.DP.blackbox}. We first approximate the parabolic Green function $P(t,x,y)$ by a function $Q(t,x,s,y)$, which is the solution of an initial value problem with a more regular initial condition. For fixed initial time~$s>0$ and~$y\in\Rd$, we take $(t,x) \mapsto Q(t,x,s,y)$ to be the solution of the initial-value problem 
\begin{equation} 
\label{e.Q}
\left\{
\begin{aligned}
& \left( \partial_t - \nabla\cdot \a\nabla \right) Q = 0 & \mbox{in} & \ (s,\infty)\times \Rd, \\ 
& Q(s,\cdot) = \overline{P}(s, \cdot-y) & \mbox{on} & \ \Rd. 
\end{aligned}
\right. 
\end{equation}
The first step in the proof of Theorem~\ref{t.PGF.basecase} is to show that $Q(t,x,s,y)$ is indeed a good approximation to~$P(t,x,y)$ provided that $s \ll t$, that is, whenever $s$ is a mesoscopic time compared to~$t$. (This can be compared to the trimming of the mesoscopic boundary layer of width~$r$ in the proof of Theorem~\ref{t.DP.blackbox}.)

\begin{lemma}
\label{l.PtoQ}
Fix $\sigma \in (0,d)$ and $\alpha \in (0,\Lambda^{-1})$. 
There exists~$C(\sigma,\alpha,d,\Lambda)<\infty$ such that, for every~$x,y\in\Rd$,~$t \in\left[  3 \X_\sigma(y)^2,\infty\right)$ and~$s \in \left[ \X_\sigma(y)^2, \frac 13 t \right]$,
\begin{equation*} 
\left| Q(t,x,s,y) - P(t,x,y) \right|
 \leq C \left( \frac{s}{t}\right)^{\frac12}t^{-\frac d2}  
\exp\left(-  \alpha \frac{|x-y|^2}{4t}\right).
\end{equation*}
\end{lemma}
\begin{proof}
Fix $\sigma \in (0,d)$, $\alpha \in (0,\Lambda^{-1})$,~$x,y\in\Rd$, $t \in\left[  3 \X_\sigma(y)^2,\infty\right)$ and $s \in \left[ \X_\sigma(y)^2, \frac 13 t \right]$. By Duhamel's formula, we have the representation 
\begin{align*}
Q(t,x,s,y) - P(t,x,y)
&
= \int_{\Rd} \left( \overline{P}(s,z-y) - P(s,z,y) \right) P(t-s,x,z) \, dz.
\end{align*}
Since $ \overline{P}(s,\cdot-y) - P(s,\cdot,y)$ has zero mass, we also get, using the Aronson bound~\eqref{e.upperP.NA},
\begin{align} \notag 
\lefteqn{\left| Q(t,x,s,y) - P(t,x,y) \right|} \quad &
\\ \notag & = \left| \int_{\Rd} \left( \overline{P}(s,z-y) - P(s,z,y) \right) \left( P(t-s,x,z) - \left( P(t-s,x,\cdot)  \right)_{B_{\sqrt{s}}(y)} \right) \, dz \right|
\\ \notag & \leq C \int_{\Rd}  \Phi\left( \frac{s}{\alpha} ,z-y \right) \left| P(t-s,x,z) - \left( P(t-s,x,\cdot)  \right)_{B_{\sqrt{s}}(y)} \right| \, dz .
\end{align}
Using the layer-cake formula~\eqref{e.PGF.layercake.1} with $f(z)  : =  \left| P(t-s,x,z) - \left( P(t-s,x,\cdot)  \right)_{B_{\sqrt{s}}(y)} \right|$, we obtain 
\begin{align} \notag 
\lefteqn{\left| Q(t,x,s,y) - P(t,x,y) \right|} \quad & 
\\ \notag & \leq C \int_{0}^\infty r^{d+1} \exp\left(- \alpha \frac{r^2}{4} \right)  \fint_{ B_{r \sqrt{s}}(y) }  \left| P(t-s,x,z) - \left( P(t-s,x,\cdot)  \right)_{B_{\sqrt{s}}(y)} \right| \, dz  \, dr.
\end{align}
Thus, by the triangle inequality and Poincar\'e's inequality, we deduce that
\begin{align} \notag 
\lefteqn{r^d \fint_{ B_{r \sqrt{s}}(y) }  \left| P(t-s,x,z) - \left( P(t-s,x,\cdot)  \right)_{B_{\sqrt{s}}(y)} \right| \, dz} \quad &
\\ \notag & \leq   ( r \vee 1)^{d} \fint_{ B_{(r\vee 1) \sqrt{s}}(y) }  \left| P(t-s,x,z) - \left( P(t-s,x,\cdot)  \right)_{B_{r \sqrt{s}}(y)} \right| \, dz
\\ \notag & \leq   C\sqrt{s} ( r \vee 1)^{d+1} \left\|  \nabla_y  P(t-s,x,z) \right\|_{\underline{L}^1\left(  B_{(r\vee 1) \sqrt{s}}(y) \right)}. 
\end{align}
Applying~\eqref{e.PGF.layercake.3.2}, we then obtain
\begin{equation*} 
\left\|  \nabla_y  P(t-s,x,z) \right\|_{\underline{L}^1\left(  B_{(r\vee 1) \sqrt{s}}(y) \right)} \leq 
 C (t-s)^{-\frac 12- \frac d2} \fint_{B_{(r\vee 1) \sqrt{s}}(y)} \exp \left( -  \alpha \frac{|z-x|^2}{4(t-s)} \right) \, dz. 
\end{equation*}
Putting thus the computations together yields 
\begin{align} \notag 
\lefteqn{\left| Q(t,x,s,y) - P(t,x,y) \right|} \quad & 
\\ \notag & \leq C \left( \frac{s}{t-s} \right)^{\frac 12}\int_{1}^\infty r^{d+1} \exp\left(- \alpha \frac{r^2}{4} \right)  (t-s)^{- \frac d2} \fint_{B_{ r \sqrt{s}}(y)} \exp \left( -  \alpha \frac{|z-x|^2}{4(t-s)} \right) \, dz  \, dr.
\end{align}
Another application of~\eqref{e.PGF.layercake.1} then gives the result, since $t \geq 3s$. 
\end{proof}

The final remaining step in the proof of Theorem~\ref{t.PGF.basecase} is to compare $Q$ to the two-scale expansion for $\overline{P}$, which is denoted by $H$ and defined above in~\eqref{e.def.H}. We first define, for each $x \in\Rd$ and $t \in (s,\infty)$,
\begin{multline} \label{e.PGF.v.def1}
v(t,x,s,y) \\ := \sum_{k=1}^d \int_{\Rd}  \left( \phi_{e_k}(z) - \left( \phi_{e_k} \ast \Phi(s,\cdot) \right)(y) \right)  \partial_{x_k}\overline{P}(s,z-y) P(t-s,x,z)\,dz,
\end{multline}
which  is the solution of
\begin{equation*} \label{}
\left\{
\begin{aligned}
& \left(\partial_{t} - \nabla_x \cdot\a\nabla_x \right)v(\cdot,\cdot,s,y) = 0 & \mbox{in} & \ (s,\infty) \times \Rd, \\
& v(s,\cdot,s,y) = H(s,\cdot,y)-Q(s,\cdot,s,y) & \mbox{on} & \ \Rd. 
\end{aligned} 
\right.
\end{equation*}
The function~$v$ represents the error due to the fact that the initial condition for $H$ is not quite the same as $\overline{P}$, because $H$ has the additional perturbation caused by the correctors. We will show that this function $v$ is small (provided~$s \gg 1$) since the correctors grow sublinearly. 

\begin{lemma} 
Fix $\sigma \in (0,d)$, $\alpha \in [0,\Lambda^{-1})$, and let $v$ be as in \eqref{e.PGF.v.def1}. There exist constants $C(\sigma,\alpha,d,\Lambda)<\infty$ and $\delta(d,\Lambda) > 0$ such that for every $x,y \in \R^d$, $s \in \left[ \X_\sigma(y)^2 ,\infty\right)$, and $t>s$,  we have
\label{l.PGF.v.est}
\begin{equation} \label{e.PGF.v.est}
\left| 
v(t,x,s,y)
\right|
\leq 
Cs^{-\delta(d-\sigma)} t^{-\frac d2} \exp\left( - \alpha \frac{|x-y|^2}{4t} \right).
\end{equation}
\end{lemma}

\begin{proof}
Since $\sqrt{s} \geq \X_{\sigma}(y)$, we have by  Remark~\ref{r.scalar.Linfty} that, for $t' \geq s$ and $z \in \R^d$,  
\begin{equation} \label{e.PGF.corrLinfty}
\sum_{k=1}^d \left| \phi_{e_k}(z) - \left( \phi_{e_k} \ast \Phi(t',\cdot) \right)(y) \right| \leq C (t')^{\frac 12 -\delta(d-\sigma) } \left(1 + \frac{|z-y|}{\sqrt{t'}} \right).
\end{equation} 
Hence, using also the Nash-Aronson bound \eqref{e.upperP.NA}, we have
\begin{multline*} 
\sum_{k=1}^d \left| \left( \phi_{e_k}(z) - \left( \phi_{e_k} \ast \Phi(s,\cdot) \right)(y) \right)  \partial_{x_k}\overline{P}(s,z-y) \right| P(t-s,x,z) \\ \leq 
C s^{ -\delta(d-\sigma) } \Phi\left( \frac{s}{\alpha}, z-y \right) \Phi\left( \frac{t-s}{\alpha}, x-y \right).
\end{multline*}
Now~\eqref{e.PGF.v.est} follows by the semigroup property and the definition of $v$ in~\eqref{e.PGF.v.def1}.
\end{proof}

We now compute the equation for $H$ and estimate the corresponding errors.

\begin{lemma} 
\label{l.PGF.H.eq}
For each $y \in \R^d$, there exist $\mathbf{F}(\cdot,\cdot,y) \in L_{\textrm{loc}}^2((0,\infty) \times \R^d;\R^d)$ and  $f(\cdot,\cdot,y) \in L_{\textrm{loc}}^2( (0,\infty) \times \R^d)$ such that, for every $x \in \R^d$ and $t \in [1,\infty)$,  
\begin{equation}  \label{e.PGF.H.eq.source}
\left(\left( \partial_{t} -\nabla_x \cdot\a\nabla_x \right) H(\cdot,\cdot,y)\right)(t,x) = \left( f(\cdot,\cdot,y) + \nabla_x \cdot \mathbf{F}(\cdot,\cdot,y)  \right) (t,x).
\end{equation}
Moreover, for each $\sigma \in (0,d)$ and $\alpha \in [0,\Lambda^{-1})$, there exist constants $C(\sigma,\alpha,d,\Lambda)<\infty$ and $\delta(d,\Lambda) \in (0,1)$ such that, for every $y \in \R^d$ and $t \geq \X_{\sigma}(y)^2$, 
\begin{equation} \label{e.PGF.H.est.source}
 \left\| \left( \left| \mathbf{F}(t,\cdot,y) \right| + t^{\frac12} \left| f(t,\cdot,y) \right| \right) \exp\left( \alpha \frac{|\cdot-y|^2}{t} \right) \right\|_{L^2(\R^d)}  \leq C t^{-\frac d2- 1 - \delta(d-\sigma)} .
\end{equation}
\end{lemma}
\begin{proof}
By a computation very similar to the one in the proo of~Lemma~\ref{l.letswritetheflux}, we obtain the identity~\eqref{e.PGF.H.eq.source} with
\begin{equation*} 
\mathbf{F}_i(t,x,y) := \sum_{j,k=1}^d \partial_{x_j} \partial_{x_k} \overline{P}(t,x-y) \left(  \mathbf{S}_{e_k,ij}(t,x,y)  - \a_{ij}(x) \phi_{e_k}(t,x,y) \right) 
\end{equation*}
and
\begin{equation*} 
f(t,x,y) := \partial_t \left(  \sum_{k=1}^d  \phi_{e_k}(t,x,y) \partial_{x_k} \overline{P}(t,x-y)  \right) ,
\end{equation*}
where we used the shorthand notation
\begin{equation*}
\left\{ 
\begin{aligned}
&
\phi_{e_k}(t,x,y) :=  \left( \phi_{e_k}(x) - \left( \phi_{e_k} \ast \Phi(t,\cdot) \right)(y) \right),
\\ &
\mathbf{S}_{e_k,ij}(t,x,y) := \mathbf{S}_{e_k,ij}(x) - \left(\mathbf{S}_{e_k,ij} \right)_{B_{\sqrt{t}}(y)}.
\end{aligned}
\right. 
\end{equation*}
We turn to the proof of~\eqref{e.PGF.H.est.source}. By Lemma~\ref{l.correctorminbounds} and Proposition~\ref{p.fluxcorrectorests}, we have, for every $r \geq \sqrt{t} \geq \X_{\sigma}(y)$, 
\begin{equation*} 
\sum_{i,j,k=1}^d  \left\| \mathbf{S}_{e_k,ij}(t,\cdot,y) \right\|_{\underline{L}^2\left(B_{r}(y) \right)} 
+\sum_{k=1}^d  \left\| \phi_{e_k}(t,\cdot,y) \right\|_{\underline{L}^2\left(B_{r}(y) \right)} 
\leq 
C r^{1-\delta}
\end{equation*}
and
 \begin{equation*} 
\sum_{k=1}^d  \left\| \partial_t \phi_{e_k}(t,\cdot,y) \right\|_{\underline{L}^2\left(B_{r}(y) \right)} \leq C r^{-1-\delta}.
\end{equation*}
On the other hand, there exists~$C(\alpha,d,\Lambda) < \infty$  such that 
\begin{equation*} 
\max_{k\in \{1,2,3\}}  t^{k} \left| \nabla^k \overline{P}(t,x-y) \right|^2 \exp\left( 2\alpha \frac{|x-y|^2}{t} \right) \leq Ct^{-\frac d2} \Phi\left(\frac{t}{\Lambda^{-1} - \alpha},x-y\right)
\end{equation*}
Together with~\eqref{e.PGF.layercake.1} from the proof of Lemma~\ref{l.PGF.layercake}, the last three displays yield~\eqref{e.PGF.H.est.source}. The proof is complete.
\end{proof}

We next estimate the error arising from the fact that~$H$ is not an exact solution of the equation (that is, from the error in the two-scale expansion). We define, for each $s,t \in [1,\infty)$ with~$t>s$ and $x,y \in \R^d$, 
\begin{equation}
\label{e.PGF.w.def1}
w(t,x,s,y) := H(t,x,y) -Q(t,x,s,y) - v(t,x,s,y).
\end{equation}
Notice that $w(s,z,s,y) = 0 $ and, according to Lemma~\ref{l.PGF.H.eq}, there exist $f$ and $\mathbf{F}$ satisfying~\eqref{e.PGF.H.est.source} and such that 
\begin{equation}  \label{e.PGF.H.eq.source2}
\left(\left( \partial_{t} -\nabla_x \cdot\a\nabla_x \right) w(\cdot,\cdot,s,y)\right)(t,x) = \left(  f(\cdot,\cdot,y) + \nabla_x \cdot \mathbf{F}(\cdot,\cdot,y)  \right) (t,x) .
\end{equation}
We next present Gaussian-type bounds on~$w$. 

\begin{lemma} 
\label{l.PGF.w.basiclemma}
Fix $\sigma \in (0,d)$, $\alpha \in [0,\Lambda^{-1})$,~$\ep \in \left(0,\tfrac12\right]$ and let $w$ be defined as in~\eqref{e.PGF.w.def1}. There exist~$C(\sigma,\alpha,\ep,d,\Lambda)<\infty$ and~$\delta(d,\Lambda)>0$ such that, 
for every $s \in \left[ \X_\sigma(y)^2 ,\infty\right)$, $t>s$ and $y\in\Rd$, we  have
\begin{equation} \label{e.PGF.w.goal}
\left\| w(t,\cdot,s,y) \exp \left( \alpha \frac{|\cdot - y|^2}{4t}\right)\right\|_{L^2(\R^d)} \leq C s^{-\frac d4 -  \delta(d-\sigma)} \left(\frac{t}{s}\right)^{\ep}.
\end{equation} 
\end{lemma}

\begin{proof}
\emph{Step 1}. We derive a differential inequality for $w$. The claim is that,  for every~$\ep  \in \left(0,\frac12 \right]$ and $\psi \in C^{\infty}(\R^d)$,
\begin{multline}  \label{e.NA.w.basic1}
\frac12 \partial_t \left( w^2 \psi^2 \right)  + \frac{\ep}{2} \psi^2 \a \nabla w \cdot \nabla w  
 \leq - \nabla \cdot \left( w \psi^2 \left( \a \nabla w  +  \mathbf{F}   \right)\right) 
\\  +  w^2 \left( \frac{\Lambda}{1-2 \ep} |\nabla \psi|^2  + \psi \partial_t \psi + \frac{\ep}{t} \psi^2 \right) + \frac{2}{\ep} \psi^2 \left( \left| \mathbf{F}  \right|^2 + t \left| f  \right|^2\right) .
\end{multline}
By a direct computation we first get
\begin{align} \notag 
w \psi^2 \left(\partial_t - \nabla \cdot \a \nabla\right) w   & = \frac12 \partial_t \left( w^2 \psi^2 \right) - \nabla \cdot \left(w \psi^2  \a \nabla w \right) 
\\ \notag & \quad +  \psi^2 \a \nabla w \cdot \nabla w +  2w\psi \a \nabla w \cdot \nabla \psi  - w^2 \psi \partial_t \psi.
\end{align}
Applying Young's inequality, we obtain, for any $\ep \in \left(0,\frac12 \right]$, 
\begin{equation*} 
\psi^2 \a \nabla w \cdot \nabla w +  2w\psi \a \nabla w \cdot \nabla \psi  \geq \ep \psi^2 \a \nabla w \cdot \nabla w - \frac{\Lambda}{1-\ep} w^2 |\nabla \psi|^2 ,
\end{equation*}
and hence we arrive at 
\begin{multline*} 
w \psi^2 \left(\partial_t - \nabla \cdot \a \nabla\right) w    = \frac12 \partial_t \left( w^2 \psi^2 \right) - \nabla \cdot \left(w \psi^2  \a \nabla w \right) 
\\ +  \ep \psi^2 \a \nabla w \cdot \nabla w - w^2 \left( \frac{\Lambda}{1-\ep} |\nabla \psi|^2  + \psi \partial_t \psi \right).
\end{multline*}
On the other hand, again by Young's inequality, we have
\begin{align} \notag 
w \psi^2 \left( \nabla \cdot \mathbf{F} + f \right)  & = \nabla \cdot \left( w \psi^2  \mathbf{F}  \right)  - \psi^2 \mathbf{F} \cdot \nabla w - 2 \psi w \mathbf{F} \cdot \nabla \psi   + w \psi^2f 
\\ \notag & \leq \nabla \cdot \left( w \psi^2  \mathbf{F}  \right) + \frac{\ep}{2} \psi^2 \a \nabla w \cdot \nabla w 
\\ \notag &  \quad  + w^2 \left(\frac{\ep}{1-\ep} |\nabla \psi|^2 +  \frac{\ep}{t}  \psi^2 \right) + \frac{2}{\ep} \psi^2 \left( \left| \mathbf{F}  \right|^2 + t \left| f  \right|^2\right) .
\end{align}
Combining the above two inequalities and using the equation for $w$ yields~\eqref{e.NA.w.basic1}. 

\smallskip

\emph{Step 2}. 
We next claim that, for each~$\ep\in \left(0,\tfrac12\right]$, there exist~$C(\sigma,\ep,d,\Lambda) < \infty $ and~$\delta(d,\Lambda) \in (0,1)$ such that, for every~$t \geq s$ and~$y\in\R^d$, 
\begin{equation} 
\label{e.PGF.w.goal1}
\left\| w(t,\cdot,y) \right\|_{L^2(\R^d)} 
\leq 
C  s^{-\frac d4 -  \delta(d-\sigma)} 
\left(\frac{t}{s}\right)^{\ep}.
\end{equation}
To see this, choose $\psi = 1$  in~\eqref{e.NA.w.basic1} and define 
\begin{equation*} 
\rho(t) :=  \int_{\R^d} w^2(t,x) \, dx 
\quad \mbox{and} \quad
\tau(t) := \frac{4}{\ep}  \int_{\R^d}  \left( \left| \mathbf{F}(t,x,y)  \right|^2 + t \left| f(t,x,y)  \right|^2\right)  \, dx .
\end{equation*}
Then~\eqref{e.NA.w.basic1} reads as
\begin{equation*} 
\partial_t \rho(t) \leq \frac{\ep}{t} \rho(t) + \tau(t) .
\end{equation*}
Since $\sqrt{s} \geq \X_{\sigma}(y)$, we have by~\eqref{e.PGF.H.est.source} that, for every~$t \geq s$,
\begin{equation*} 
\tau(t)  \leq \frac{C}{\ep} t^{-\frac d2 - 1-\delta(d-\sigma)}.
\end{equation*}
Integrating and using the fact that $\rho(s) =0$, we obtain
\begin{equation*} 
\rho(t) \leq \int_{s}^{t} \frac{\ep}{t'} \rho(t') \, dt'  + \frac{C}{\ep} s^{-\frac d2 -\delta(d-\sigma)}.
\end{equation*}
Gr\"onwall's inequality now implies~\eqref{e.PGF.w.goal1}. 

\smallskip

\emph{Step 3}. 
We prove~\eqref{e.PGF.w.goal}. Applying\eqref{e.NA.w.basic1} with 
\begin{equation*} 
\psi(t,x) := \exp\left( \alpha \frac{|x-y|^2}{4t}\right)
\end{equation*}
we obtain, for $\tilde{\ep}(\alpha,d,\Lambda)>0$ sufficiently small (recall also that $0<\alpha<\Lambda^{-1}$), 
\begin{align*}
&
\left( \frac{\Lambda}{1-2 \tilde{\ep}} |\nabla \psi|^2 + \psi \partial_t \psi + \frac{\tilde{\ep}}{t} \psi^2 \right)(x,t)
\\ & \qquad 
= \frac{\psi^2(x,t)}{t}\left( \frac{\Lambda}{1-2 \tilde{\ep}} \alpha^2 \frac{|x-y|^2}{4t} - \alpha \frac{|x-y|^2}{4t} + \tilde{\ep} \right) \leq \frac{C}{t}.
\end{align*}
Thus, after integration with respect to~$x$, we obtain by~\eqref{e.NA.w.basic1},~\eqref{e.PGF.w.goal1} and~\eqref{e.PGF.H.est.source} that 
\begin{equation*} 
\partial_{t} \left\| w(t,\cdot,y) \exp \left( \alpha \frac{|\cdot - y|^2}{4t}\right)\right\|_{L^2(\R^d)}^2 \leq C t^{-1} s^{-\frac d4 -  \delta(d-\sigma)} \left(\frac{t}{s}\right)^{\ep}.
\end{equation*}
Integrating with respect to~$t$ and recalling that $w(s,\cdot,\cdot) = 0$, we obtain~\eqref{e.PGF.w.goal}.
\end{proof}

We now combine the previous three lemmas to obtain Theorem~\ref{t.PGF.basecase}.

\begin{proof}[{Proof of Theorem~\ref{t.PGF.basecase}}]
Throughout the proof, we fix $\alpha \in (0,\Lambda^{-1})$. 
\smallskip

\emph{Step 1.} 
 We complete the proof of~\eqref{e.PGF.basecase}. We set, for $\theta \in \left(0, \frac{d-\sigma}{2\sigma} \wedge 1 \right) $, 
\begin{equation*} 
\Y_\sigma(y) := \X_{\frac{\sigma + d}{2}}(y)^{1+\theta}.
\end{equation*}
 Clearly $\frac{\sigma}{1-\theta} < d$, and hence $\Y_\sigma(y) \leq \O_\sigma(C)$. Taking $s = t^{1-\theta}$ we see that $\sqrt{t} \geq \Y_\sigma(y)$ implies $\sqrt{s} \geq \X_{\frac{\sigma+d}{2}}(y)$. We decompose $P-\overline{P}$ as 
\begin{align} \notag 
P(t,x,y) - \overline{P}(t,x-y)  & = \left( P(t,x,y) - Q(t,x,s,y) \right) 
+ \left(  H(t,x,y) -  \overline{P}(t,x-y) \right) 
\\ \notag & \quad  - \left( w(t,x,s,y) - v(t,x,s,y) \right).
\end{align}
We estimate the three terms in parentheses on the right using Lemmas~\ref{l.PtoQ},~\ref{l.PGF.v.est} and~\ref{l.PGF.w.basiclemma}. 
Notice that these lemmas are applicable since~$\sqrt{s} \geq \X_{\frac{\sigma+d}{2}}(y)$.
Lemma~\ref{l.PGF.v.est} implies that 
\begin{equation*} 
\left\| v(t,\cdot,s,y) \exp \left( \alpha \frac{|\cdot - y|^2}{4t}\right)\right\|_{L^2(\R^d)} \leq C t^{-\frac d4 -  \delta(1-\theta)(d-\sigma)}.
\end{equation*}
Moreover, by Remark~\ref{r.scalar.Linfty}, there exist $C(\sigma,d,\Lambda)<\infty$ and $\gamma(d,\Lambda)\in (0,1)$ such that, for every $r \geq \sqrt{s} \geq \X_\sigma(y)$, we have
\begin{equation*} 
\left\| \phi_{e_{k}} - \left( \phi_{e_{k}} \ast \Phi(s,\cdot) \right)(y)\right\|_{L^{\infty}(B_r(y))}  \leq C r^{1-\gamma(d-\sigma)}.
\end{equation*}
Using this, we also obtain, for $\alpha' := \frac12 \left(\alpha + \Lambda^{-1} \right)$, that
\begin{equation*}
\left| \left( \phi_{e_k}(x) - \left( \phi_{e_k} \ast \Phi(t,\cdot) \right)(y) \right) \partial_{x_k}\overline{P}(t,x-y)\right|  \leq C t^{-\frac{\gamma}{2}  (d-\sigma)} \Phi\left(\frac{t}{\alpha'}, x-y \right),
\end{equation*}
and thus
\begin{equation*} 
\left\| \left( H(t,\cdot,y)-\overline{P}(t,\cdot,y) \right) \exp \left( \alpha \frac{|\cdot - y|^2}{4t}\right)\right\|_{L^2(\R^d)} \leq C t^{-\frac d4 -\frac{\gamma}{2}  (d-\sigma)}.
\end{equation*}
Lemma~\ref{l.PtoQ} yields 
\begin{equation*} 
\left\| \left( P(t,\cdot,y)-Q(t,\cdot,s,y) \right) \exp \left( \alpha \frac{|\cdot - y|^2}{4t}\right)\right\|_{L^2(\R^d)} \leq C t^{-\frac d4 -\frac{\theta}{2} (d-\sigma)}.
\end{equation*}
Lemma~\ref{l.PGF.w.basiclemma} implies
\begin{equation*} 
\left\| w(t,\cdot,s,y) \exp \left( \alpha \frac{|\cdot - y|^2}{4t}\right)\right\|_{L^2(\R^d)} \leq C t^{-\frac d4 -  \frac{\delta}{2}(d-\sigma) + \theta\left( \frac d4 + \delta(d-\sigma) + \ep \right)  } .
\end{equation*} 
Combining the above displays and choosing $\theta$ small enough, we get
\begin{equation} \label{e.P-barP.temp1}
\left\| \left( P(t,\cdot,y)-\overline{P}(t,\cdot-y) \right) \exp \left( \alpha \frac{|\cdot - y|^2}{4t}\right)\right\|_{L^2(\R^d)} \leq C t^{-\frac d4 -\frac{\delta \wedge \gamma \wedge \theta}{4} (d-\sigma)}.
\end{equation}
We next turn this inequality into a pointwise estimate using the semigroup property. We write
\begin{align} \notag \label{e.P-barP.temp2}
 P(t,x,y)-\overline{P}(t,x-y) & = 
\int_{\R^d} \left(P\left(\tfrac t2 ,z,y\right)-\overline{P}\left( \tfrac t2,z-y\right) \right)  \overline{P}\left(\tfrac t2 ,x-z\right) \, dz
\\  & \quad 
+ \int_{\R^d}  \left(P\left(\tfrac t2 ,x,z\right)-\overline{P}\left( \tfrac t2,x-z\right) \right)  P\left(\tfrac t2 ,z,y\right) \, dz .
\end{align}
H\"older's inequality yields
\begin{align} \notag 
\lefteqn{\left| \int_{\R^d} \left(P\left(\tfrac t2 ,z,y\right)-\overline{P}\left( \tfrac t2,z-y\right) \right)  \overline{P}\left(\tfrac t2 ,x-z\right) \, dz \right| } \quad &
\\ \notag & \leq \left\| \left(P\left(\tfrac t2 ,\cdot ,y\right)-\overline{P}\left( \tfrac t2,\cdot-y\right) \right)\exp \left(  \frac{\alpha|\cdot - y|^2}{2t}\right)   \right\|_{L^2(\R^d)} 
\\ \notag & \qquad \times 
\left\| \overline{P}\left(\tfrac t2 ,x-\cdot\right)  \exp \left( -  \frac{\alpha|\cdot - y|^2}{2t}\right)  \right\|_{L^2(\R^d)} .
\end{align}
By the semigroup property and the Nash-Aronson estimate, 
\begin{equation*} 
\left\| \overline{P}\left(\tfrac t2 ,x-\cdot\right)  \exp \left( - \alpha \frac{|\cdot - y|^2}{2t}\right)  \right\|_{L^2(\R^d)} 
\leq C t^{\frac d4} \Phi\left(\frac {t}{\alpha} , x-y\right).
\end{equation*}
The second term on the right in~\eqref{e.P-barP.temp2} can be estimated similarly. Thus, combining the above two displays with~\eqref{e.P-barP.temp1} and~\eqref{e.P-barP.temp2} completes the proof of~\eqref{e.PGF.basecase} by taking $\delta$ smaller in the statement, if necessary.

\smallskip

\emph{Step 2.} We use~\eqref{e.PGF.basecase} and the $C^{1,1}$ regularity estimate to obtain~\eqref{e.PGF.basecase.grad}. Fix $\alpha \in (0,\Lambda^{-1})$ and $\beta:=\frac13\left(\Lambda^{-1}-\alpha\right)>0$ 
so that, by the Nash-Aronson bound, we have, for every $x,y\in\Rd$, and $t>0$,  
\begin{equation} 
\label{e.pickthebetayouass}
\left\| P(\cdot ,y) \right\|_{L^\infty( Q_{\beta \sqrt{t}}(t,x))}
+ t \left\| \nabla^2 \overline{P} (\cdot ,y) \right\|_{L^\infty( Q_{\beta \sqrt{t}}(t,x))}
\leq
Ct^{-\frac d2} \exp\left( -\alpha \frac{|x-y|^2}{4t} \right). 
\end{equation}
Fix $x,y\in\Rd$ and $t \in [1,\infty)$ such that $\beta \sqrt{t} \geq \Y_\sigma(x)$. 
By the parabolic $C^{1,1}$-type regularity estimate (Theorem~\ref{t.regularity.parab} for $k=1$),  we may select~$\xi\in\Rd$ and $a\in\R$ such that, for every $r\in \left[ \X_\sigma(x),\beta \sqrt{t} \right]$, 
\begin{align} \label{}
\left\| P(t,\cdot,y) - \left( \ell_\xi + \phi_\xi +a \right) \right\|_{\underline{L}^2(Q_r(t,x))}
&
\leq
C \left( \frac{r}{\sqrt{t}} \right)^2 \left\|P(t,\cdot,y) \right\|_{\underline{L}^2(Q_{\beta\sqrt{t}}(t,x)}
\notag \\ & 
\leq C r^2 t^{-1 -\frac d2} \exp\left( -\alpha \frac{|x-y|^2}{4t} \right).
\label{e.PGFaffineapprox}
\end{align}
Here $\ell_\xi(x)=\xi\cdot x$ and $\phi_\xi$ is the first-order stationary corrector. 
By the Caccioppoli inequality (Lemma~\ref{l.cacciopp.this}) and Lemma~\ref{l.slicesvsaverages.PAR}, this implies, for every such~$r$,
\begin{equation} 
\label{e.PGFaffineapprox.cacc}
\left\| \nabla_x P(t,\cdot,y) - \left( \xi + \nabla \phi_\xi \right) \right\|_{\underline{L}^2(B_r(x))}
\leq 
C rt^{-1 -\frac d2} \exp\left( -\alpha \frac{|x-y|^2}{4t} \right).
\end{equation}
By~\eqref{e.PGF.basecase},~\eqref{e.PGFaffineapprox} and the triangle inequality, we deduce that, for every such~$r$, 
\begin{equation*} \label{}
\left\| \overline{P}(\cdot,\cdot-y) - \left( \ell_\xi + \phi_\xi +a \right) \right\|_{\underline{L}^2(Q_r(t,x))}
\leq 
C \left( r^2 t^{-1 -\frac d2} +  t^{-\delta(d-\sigma)-\frac d2} \right) \exp\left( -\alpha \frac{|x-y|^2}{4t} \right).
\end{equation*}
Using, in view of~\eqref{e.pickthebetayouass}, that 
\begin{multline*} \label{}
\sup_{(s,z) \in Q_r(t,x)}
\left| \overline{P}(s,z-y) - \left( \overline{P}(t,x-y) + (z-x) \cdot \nabla \overline{P}(t,x-y) \right) \right|
\\
\leq 
Cr^2 \left\| \nabla^2 \overline{P} (\cdot ,y) \right\|_{L^\infty( Q_{\beta \sqrt{t}}(t,x))}
\leq
C r^2t^{-1-\frac d2}\exp\left( -\alpha \frac{|x-y|^2}{4t} \right),
\end{multline*}
we deduce that
\begin{align*} \label{}
\lefteqn{
r \left|  \nabla \overline{P}(t,x-y) - \xi  \right|
} \quad & 
\\ &
\leq \left\| \phi_\xi - \left( \phi_\xi \right)_{B_r(x)} \right\|_{\underline{L}^2(B_r(x))} + C \left( r^2 t^{-1 -\frac d2} +  t^{-\delta(d-\sigma)-\frac d2} \right) \exp\left( -\alpha \frac{|x-y|^2}{4t} \right)
\\ & 
\leq 
Cr^{1-\delta} |\xi| + C \left( r^2 t^{-1 -\frac d2} +  t^{-\delta(d-\sigma)-\frac d2} \right) \exp\left( -\alpha \frac{|x-y|^2}{4t} \right).
\end{align*}
This implies that 
\begin{equation*} \label{}
\left|  \nabla \overline{P}(t,x-y) - \xi  \right|
\leq 
C \left( r t^{-1 -\frac d2} +  r^{-1}t^{-\delta(d-\sigma)-\frac d2} \right) \exp\left( -\alpha \frac{|x-y|^2}{4t} \right).
\end{equation*}
Taking $r = t^{\frac12(1 - \delta(d-\sigma))}$ in the above expression yields, after shrinking $\delta$, 
\begin{equation*} \label{}
\left|  \nabla \overline{P}(t,x-y) - \xi  \right|
\leq 
C t^{-\delta(d-\sigma)-\frac 12- \frac d2} \exp\left( -\alpha \frac{|x-y|^2}{4t} \right).
\end{equation*}
From this,~\eqref{e.PGFaffineapprox.cacc} and the fact that $r\geq \Y_\sigma(x) \geq \X_\sigma(x)$, we obtain 
\begin{align*} \label{}
\lefteqn{
\left\| \nabla_x P(t,\cdot,y) - \left( \nabla \overline{P}(t,x-y) + \nabla \phi_{\nabla \overline{P}(t,x-y)} \right) \right\|_{\underline{L}^2(B_r(x))}
} \ \ & 
\\ & 
\leq
\left\| \nabla_x P(t,\cdot,y) - \left( \xi + \nabla \phi_\xi \right) \right\|_{\underline{L}^2(B_r(x))}
+
 \left|  \nabla \overline{P}(t,x-y) - \xi  \right| \left( 1+ \sup_{e\in B_1} \left\| \nabla \phi_e \right\|_{L^2(B_r(x))} \right)
\\ & 
\leq 
C \left( r t^{-1 -\frac d2} +  t^{-\delta(d-\sigma)-\frac12-\frac d2} \right) \exp\left( -\alpha \frac{|x-y|^2}{4t} \right).
\end{align*}
Since 
\begin{align*} \label{}
\lefteqn{
\left\| \nabla_x H(t,\cdot,y) - \left( \nabla \overline{P}(t,x-y) + \nabla \phi_{\nabla \overline{P}(t,x-y)} \right) \right\|_{\underline{L}^2(B_r(x))}
} \quad \quad  & 
\\ & 
\leq 
\left( \sup_{z\in B_r(x)} \left| \nabla\overline{P}(t,x-y) - \nabla\overline{P}(t,z-y) \right| \right) \left(1 + \sup_{e\in B_1} \left\| \nabla \phi_e \right\|_{L^2(B_r(x))} \right)
\\ & 
\leq
C \left( r t^{-1-\frac d2} \right) \exp\left( -\alpha \frac{|x-y|^2}{4t} \right),
\end{align*}
we obtain
\begin{align*} \label{}
&
\left\| \nabla_x P(t,\cdot,y) - \nabla_x H(t,\cdot,y) \right\|_{\underline{L}^2(B_r(x))}
\\ & \qquad 
\leq 
C \left( r t^{-1-\frac d2} + t^{-\delta(d-\sigma)-\frac12-\frac d2} \right) \exp\left( -\alpha \frac{|x-y|^2}{4t} \right).
\end{align*}
After shrinking $\delta$ again, this completes the proof of~\eqref{e.PGF.basecase.grad} under the additional restriction that 
\begin{equation*} \label{}
r \leq t^{\frac12 - \frac12 \delta(d-\sigma)} \wedge \beta \sqrt{t}.
\end{equation*}
For $r \in \left[ t^{\frac12 - \frac12 \delta(d-\sigma)}, \frac12\sqrt{t} \right]$, however, the desired estimate follows immediately from the Caccioppoli inequality and~\eqref{e.PGF.basecase}. Therefore the proof is complete. 
\end{proof}

We conclude this section with the proof of Corollary~\ref{c.GtoGbar0}.

\begin{proof}[{Proof of Corollary~\ref{c.GtoGbar0}}]
We can immediately pass from~\eqref{e.PGF.basecase} to~\eqref{e.GtoGbar0} by integration in time, using formula~\eqref{e.intGammatoG} and, in $d=2$, the formula~\eqref{e.intGammatoG.d=2}. One can pass from~\eqref{e.PGF.basecase.grad} to~\eqref{e.GtoGbar.nabla0} in the same way, although it is necessary to control the difference between the ``additive constants'' in the definitions of~$H$ and~$F$, i.e., the difference between 
\begin{equation*} \label{}
\int_0^\infty \phi_{e_k} \ast \Phi(t,\cdot)(y) \partial_{x_k} \overline{P}(t,x-y)\,dt
\quad \mbox{and} \quad
\phi_{e_k} \ast \Phi(|x-y|^2,\cdot)(y) \partial_{x_k} \overline{G}(x-y). 
\end{equation*}
Alternatively, we can obtain~\eqref{e.GtoGbar.nabla0} from~\eqref{e.GtoGbar0} by repeating the argument of Step~2 in the proof of Theorem~\ref{t.PGF.basecase} above, using the (elliptic) $C^{1,1}$ regularity estimate. The details are left to the reader. 
\end{proof}

\index{Green function!two-scale expansion of $\nabla_x\nabla_yP$|(}

\begin{exercise}[{Quantitative two-scale expansion of $\nabla_x\nabla_yP$}]
\label{exer.mixedsecond}
In this exercise, we extend the statement of Theorem~\ref{t.PGF.basecase} by proving a bound on the two-scale expansion of the mixed second derivatives of the Green functions. To give the statement, we must first modify the definition~\eqref{e.def.H} of $H$ so that it has the proper oscillations in both variables. We define
\begin{align} 
\label{e.def.H.mixed}
\lefteqn{
\tilde{H}(t,x,y) 
} \ & 
\\  & \notag
:= \overline{P}(t,x-y) 
+ \sum_{k=1}^d \left( \phi_{e_k}(x) - \left( \phi_{e_k} \ast \Phi(t,\cdot) \right)(y) \right)  \partial_{x_k}\overline{P}(t,x-y)
\\ & \  \notag
- \sum_{l=1}^d \left( \phi_{e_l}(y) - \left( \phi_{e_l} \ast \Phi(t,\cdot) \right)(x) \right)  \partial_{x_k}\overline{P}(t,x-y)
\\ & \  \notag
- \sum_{k,l=1}^d 
\left( \phi_{e_k}(x) - \left( \phi_{e_k} \ast \Phi(t,\cdot) \right)(y) \right) 
\left( \phi_{e_l}(y)- \left(\phi_{e_l} \ast \Phi(t,\cdot) \right)(x) \right) 
\partial_{x_k} \partial_{x_l} \overline{P}(t,x-y).
\end{align}
Note that, neglecting terms of lower order, we have
\begin{align*}
\nabla_x \nabla_y \tilde{H}(t,x,y)
&
=
- \sum_{k,l=1}^d \left( e_k + \nabla \phi_{e_k} (x) \right)\left( e_k + \nabla \phi_{e_k} (y) \right) \partial_{x_k}{\partial_{x_l}} \overline{P}(t,x-y)
\\ & \quad 
+ \mbox{lower order terms.} 
\end{align*}
The result is that there exists a random variable $\Y_\sigma$ satisfying~\eqref{e.PGD.basecase.Ysig} such that, for every $\sqrt{t} \geq 2\Y_\sigma(y)\vee \Y_\sigma(x)$ and $r \in \left[ \Y_\sigma(y)\vee \Y_\sigma(x), \frac12 \sqrt{t} \right]$, we have 
\begin{equation*}
\left\| \nabla_x \nabla_y P(t,\cdot,y) - \nabla_x \nabla_y \tilde{H}(t,\cdot) \right\|_{\underline{L}^2(B_r(x)\times B_r(y))}
\leq 
Ct^{-\delta(d-\sigma)} t^{ - 1 - \frac d2}  \exp\left(- \alpha \frac{|x-y|^2}{4t}\right).
\end{equation*}
To prove this, we use the bounds already proved in Theorem~\ref{t.PGF.basecase} and then 
follow the argument in Step~2 of the proof of Theorem~\ref{t.PGF.basecase} above, in the~$y$ variable.
\end{exercise}

\index{Green function!two-scale expansion of $\nabla_x\nabla_yP$|)}

\section*{Notes and references}

In the case of periodic coefficients, estimates for the decay of the gradient of the elliptic Green function (and the Green function itself for systems) was proved by Avellaneda and Lin~\cite{AL1,AL2} and the quantitative estimate for the decay of the error of the two-scale expansion was proved by Kenig, Lin and Shen~\cite{KLS2} for the elliptic Green function and by Geng and Shen~\cite{GS2} for the parabolic Green function. Decay estimates for the gradient of the elliptic Green function in the random setting were first proved by Marahrens and Otto~\cite{MO}, although with suboptimal stochastic integrability (i.e., with finite moments rather than the stronger $\O_{d-}$-type minimal scale bounds of Theorem~\ref{t.PGF.decay}). A two-scale expansion estimate for the elliptic Green function was proved by Bella, Giunti and Otto~\cite{BGO}.

\addtocontents{toc}{\protect\newpage}



\chapter{Decay of the parabolic semigroup}
\label{c.semigroup}

Consider the parabolic initial-value problem
\begin{equation}
\label{e.flow}
\left\{
\begin{aligned}
& \partial_t u -\nabla \cdot \left( \a \nabla u \right) = 0 &  \mbox{in} & \ (0,\infty) \times \Rd, \\
& u(0,\cdot) = \nabla\cdot \g & \mbox{on} & \ \Rd,
\end{aligned}
\right.
\end{equation}
where the vector field~$\g$ is a bounded, stationary random field with a unit range of dependence. The solution of~\eqref{e.flow} can be written by the Green function formula
\begin{equation}  
\label{e.SG.u.def}
u(t,x) := - \int_{\R^d} \g(y) \cdot \nabla_y P(t,x,y) \, dy.
\end{equation}
Since~$\g$ is a stationary random field, the field~$u(t,\cdot)$ remains stationary for all times~$t>0$. As we will see later in \eqref{e.null.average}, since the parabolic equation preserves mass, the function $u$ satisfies, for every $t>0$, 
\begin{equation*}
\E \left[ \int_{ [0,1]^d } u(t,x)\,dx \right] 
= 0. 
\end{equation*}
We can expect therefore that $u(t,x)\to 0$ as~$t\to \infty$. The goal here is to obtain quantitative bounds on the rate of decay of the solution $u(t,x)$ to zero for large~$t$.

\smallskip
 
If the heterogeneous operator~$\partial_t - \nabla \cdot \a\nabla$ is replaced by the heat operator~$\partial_t -\Delta$, then the formula~\eqref{e.SG.u.def} becomes the convolution of the gradient of the standard heat kernel against the random field~$\g$. It is easy to verify that, in this case, 
\begin{equation*}
\left| u(t,x) \right| 
\leq  
\O_2 \left( C t^{-\frac12 - \frac d4} \right) .
\end{equation*}
To see this, we write
\begin{equation*}
u(t,x) = \int_{\Rd} \g(y) \cdot \nabla \Phi(t,x-y) \,dy
=
t^{-\frac12} \int_{\Rd} \g(y) \cdot \left( t^{-\frac d2} \nabla \Phi\left(1,\frac{x-y}{\sqrt{t}}\right)\right) \,dy
\end{equation*}
and then observe that the integral on the right side represents a weighted average of~$\g$ over a length scale of order~$\sqrt{t}$. That is, roughly speaking, it is an average of~$O(t^{\frac d2})$ many independent random variables, each of which is bounded by a constant times the assumed bound on~$\g$. By Lemma~\ref{l.barO.boxes}, this integral should therefore be at most of size~$\O_2\left(Ct^{-\frac d4} \right)$. The extra factor of~$t^{-\frac12}$ in front of the integral therefore gives the claimed bound on~$|u(t,x)|$. 

\smallskip

The first goal of this chapter is to obtain an analogous estimate for the heterogeneous parabolic operator. The statement is given in Theorem~\ref{t.semigroup}, below. Its proof depends on the theory developed in the first three chapters, namely the 
(suboptimal) bounds on the first-order correctors which are summarized in Lemma~\ref{l.correctorminbounds} and the corresponding (still suboptimal) bounds on the flux correctors established in Section~\ref{s.fluxcorrectors}. We also use the two-scale expansion result of Theorem~\ref{t.twoscaleplugveng} as well as the estimates on the parabolic Green function proved earlier in Chapter~\ref{c.parabolic}. We do not, however, use the optimal quantitative estimates on the first-order correctors proved in Chapter~\ref{c.A1}. 

\smallskip

In fact, Theorem~\ref{t.semigroup} can be used to give a second, independent proof of the results of Chapter~\ref{c.A1}. Indeed, by the Duhamel formula, the first-order corrector~$\phi_e$ can be written formally as the integral 
\begin{equation} 
\label{e.duhamel.FOC}
\phi_{e} (x) = \int_0^\infty u_e (t,x) \,dt, 
\end{equation}
where $u_e(t,x)$ is the solution of the parabolic initial-value problem
\begin{equation}
\label{e.flow.FOC}
\left\{
\begin{aligned}
& \partial_t u_e -\nabla \cdot \left( \a \nabla u_e \right) = 0 &  \mbox{in} & \ (0,\infty) \times \Rd, \\
& u(0,\cdot) = \nabla\cdot \left( \a(\cdot) e \right) & \mbox{on} & \ \Rd.
\end{aligned}
\right.
\end{equation}
The integral in~\eqref{e.duhamel.FOC} may not converge, which is related to the fact that the correctors are only defined up to constants and may not exist as stationary objects themselves. Therefore, to be more precise, we should define their gradients by
\begin{equation} 
\label{e.duhamel.FOC.grad}
\nabla\phi_{e} (x) = \int_0^\infty \nabla u_e (t,x) \,dt.
\end{equation}
This integral is convergent, as we will see, and in fact the relation~\eqref{e.duhamel.FOC.grad} allows us to obtain from Theorem~\ref{t.semigroup} all of the main results of Chapter~\ref{c.A1}, in particular the optimal bounds on the first-order correctors given in Theorems~\ref{t.correctors} and~\ref{t.HKtoSob}. Since the proof of Theorem~\ref{t.semigroup} does not rely on the bootstrap arguments of Chapter~\ref{c.A1}, we obtain a second, independent proof of these results.

\smallskip

In Section~\ref{s.PGFopt}, we give quantitative estimates on the homogenization error of the elliptic and parabolic Green functions. For the parabolic Green function, the estimate is roughly that 
\begin{equation} 
\label{e.rough.PGFts}
\left| P(t,x,y) - \overline{P}(t,x-y) \right|
= \O_{2-}\left( Ct^{-\frac12} t^{-\frac d2} 
\exp\left( -\frac{\alpha|x-y|^2}t \right) \right),
\end{equation}
up to a logarithmic correction of~$\log^{\frac12} t$ in dimension two. See Theorem~\ref{t.PtoPbar} below for the precise statement. This improves the scaling of the estimate proved in the previous chapter, which has~$t^{-\delta}$ in place of~$t^{-\frac 12}$ (cf. Theorem~\ref{t.PGF.basecase}). Since the characteristic length scale of~$\overline{P}(t,\cdot)$ is of order~$t^{-\frac12}$, the bound~\eqref{e.rough.PGFts} matches the order of the homogenization error found previously (e.g., in Theorems~\ref{t.correctors} or~\ref{t.L2EE}) and therefore is optimal in scaling of~$t$.

\smallskip

These Green function estimates are closely related to the semigroup decay estimates.  The connection arises because the center of mass of the parabolic Green function, which unlike the mass itself is not preserved by the flow, can be expressed in terms of the semigroup flow. It must be controlled precisely in order to have sharp estimates, and the bounds on the decay of the parabolic semigroup gives exactly what is needed to obtain~\eqref{e.rough.PGFts}. 

\smallskip

Throughout this chapter, as in Sections~\ref{s.greenie.decay} and~\ref{s.greenie.base} (but unlike in the rest of the book), we allow ourselves to use estimates which are only valid for scalar elliptic equations, but not for systems.

\section{Optimal decay of the parabolic semigroup}
\label{s.sgp.opt}

This section is devoted to the proof of the following theorem.

\begin{theorem}[Optimal decay of the parabolic semigroup]
\label{t.semigroup}
For every~$\sigma \in (0,2)$, there exists a constant $C(\sigma,d,\Lambda)<\infty$ such that the following holds. Let~$\g$ be an~$\Rd$-valued, $\Zd$--stationary random field (in the sense of Definition~\ref{def.stat.field}) such that~$\| \g \|_{L^\infty(\Rd)}\leq 1$ and, for every $x\in\Rd$,~$\g(x)$ is $\F(B_1(x))$--measurable. Let~$u$ be the solution of the parabolic initial-value problem~\eqref{e.flow}. Then, for every $t\in [1,\infty)$ and $x\in\Rd$,
\begin{equation} 
\label{e.semidecay}
\left| u(t,x) \right| 
\leq \O_\sigma \left(C  t^{-\frac 12 - \frac d4} \right).
\end{equation}
\end{theorem}

Throughout this section, the function~$u$ is the solution of the initial-value problem~\eqref{e.flow}, where~$\g$ is a random field satisfying the assumptions of Theorem~\ref{t.semigroup}.

\smallskip

We begin by recording a crude, deterministic bound on $u$ which is an easy consequence of the Nash-Aronson estimate. It will be useful below mostly for technical reasons as well as for establishing the base case for the inductive proof of Lemma~\ref{l.SG.pointwise.u} below.

\begin{lemma}
\label{l.crude.ubound}
There exists~$C(d,\Lambda)<\infty$ such that, for every $t\in  (0,\infty)$,  
\begin{equation} 
\label{e.SG.u.decay}
\left\| u(t,\cdot) \right\|_{L^\infty(\Rd)}  
 +t^{\frac12} \left\|\nabla u(t,\cdot)\right\|_{L^\infty(\Rd)}
\leq 
C t^{-\frac 12}. 
\end{equation}
\end{lemma}
\begin{proof}
We combine the bound~\eqref{e.gradientP.NA} with the time-slice gradient estimate for parabolic equations with coefficients which are independent of time contained in Lemma~\ref{l.slicesvsaverages.PAR} to obtain the existence of a constant $C(d,\Lambda)<\infty$ such that for every $t > 0$ and $x,y \in \Rd$,  
\begin{equation} 
\label{e.dumbdecay.greenie}
\left\| \nabla_y P(t,x,\cdot) \right\|_{\underline{L}^2\left(B_{\sqrt{t}} (y) \right) }
\leq 
Ct^{-\frac12-\frac d2} \exp\left( -\frac{|x-y|^2}{Ct} \right). 
\end{equation}
Since we assume that $\left\| \g \right\|_{L^\infty(\Rd)} \leq 1$, we obtain the $L^\infty$ bound on $u$ in~\eqref{e.SG.u.decay} from the previous line and~\eqref{e.SG.u.def}. By the Caccioppoli inequality (Lemma~\ref{l.cacciopp.this}) and Lemma~\ref{l.slicesvsaverages.PAR}, we deduce the gradient bound on $u$ in~\eqref{e.SG.u.decay}.
\end{proof}

To prepare for the use of the independence assumption, we next show that the parabolic Green function is well-approximated by a function with local dependence on the coefficient field. 

\index{Green function!localization of}

\begin{lemma}  
[{Localization of $P$}]
\label{l.semigp.loc}
There exist a constant $C(d,\Lambda) < \infty$ and, for each $r\in [2,\infty)$, a random element of $L^\infty_{\mathrm{loc}}((0,r^2] \times\Rd\times\Rd)$, denoted by 
\begin{equation*}
(t,x,y) \mapsto P'_r(t,x,y)
\end{equation*}
which satisfies, for every $t \in (0,r^2]$ and $x,y \in \Rd$, the following statements: 
\begin{equation}
\label{e.Pprimemeas}
P'_r(t,x,y) \quad \mbox{is $\F(B_r(x))$--measurable}
, 
\end{equation}
\begin{equation}
\label{e.Pprimesupport}
P'_r(t,x,\cdot) \equiv 0 \quad \mbox{in} \ \Rd\setminus B_{r-1}(x),
\end{equation}
\begin{equation}  
\label{e.tildeP.prob}
\int_{\Rd}  P'_r(t,x',y) \, dx' = 1, 
\end{equation}
\begin{equation}
\label{e.Plocalize1}
\Ll| P(t,x,y) - P'_r(t,x,y) \Rr| \le C t^{-\frac d 2} \exp \Ll( -\frac{r^2 + |x-y|^2}{Ct} \Rr) ,
\end{equation}
and
\begin{equation}
\label{e.Plocalize2}
\Ll\| \nabla P(t,x,\cdot) - \nabla P'_r(t,x,\cdot) \Rr\|_{\underline L^2(B_{\sqrt{t}}(y))} \le C t^{-\frac 1 2 - \frac d 2} \exp \Ll( -\frac{r^2 + |x-y|^2}{Ct} \Rr).
\end{equation}
\end{lemma} 
\begin{proof}
In place of~\eqref{e.Pprimemeas} and~\eqref{e.Pprimesupport}, we will construct $P'_r$ so that $P'_r(t,x,y)$ is $\mathcal{F}(x+15r\cu_0)$--measurable and $P'_r(t,x,\cdot)$ is supported in $x+15r\cu_0$. We then obtain the statement of the proposition by reindexing the parameter~$r$ and performing a standard covering argument to get the last estimate~\eqref{e.Plocalize2}. 

\smallskip

For each $x'\in\Rd$ and $r\in[1,\infty)$, define the modified coefficient field 
\begin{equation*}
\tilde{\a}_{x',r} := \a\indc_{x'+2r\cu_0} + \Id \indc_{\Rd \setminus \left( x'+2r\cu_0 \right) }. 
\end{equation*}
Let $\tilde{P}_{x',r} = \tilde{P}_{x',r} (t,x,y)$ denote the parabolic Green function corresponding to~$\tilde{\a}_{x',r}$. It is clear that $\tilde{P}_{x',r}(t,x,y)$ is an~$\F(x'+2r\cu_0)$--measurable random variable for each fixed $x',x,y\in\Rd$, $r\in [2,\infty)$ and~$t>0$.

\smallskip

\emph{Step 1.} 
We show that, for every $x,x',y\in\Rd$,  $r\in [2,\infty)$ and $t\in (0,r^2]$ such that $x \in x' + r\cu_0$, we have
\begin{equation}
\label{e.Ploctildeest1}
\left| P(t,x,y) - \tilde{P}_{x',r} (t,x,y) \right|  
\leq 
C t^{-\frac d2} 
\exp\left( - \frac{r^2+\left|x-y\right|^2 }{Ct} \right).
\end{equation}
To estimate this, we begin from the identity
\begin{multline}
\label{e.duhamelPtileP}
P(t,x,y) - \tilde{P}_{x',r} (t,x,y)
\\
=
\int_0^t \int_{\Rd} 
\nabla_z P(s,y,z) 
\cdot 
\left( \tilde{\a}_{x',r}(z) - \a(z) \right) 
\nabla_z \tilde{P}_{x',r} (s,x,z) \, dz\, ds.
\end{multline}
This can be seen by computing the equation satisfied by~$P(t,x,\cdot)-\tilde{P}_{x',r} (t,x,\cdot)$ and applying the Duhamel formula. Thus 
\begin{align*}
\left| P(t,x,y) - \tilde{P}_{x',r} (t,x,y) \right| 
\leq  
C
\int_0^t \int_{\Rd \setminus \left(x'+2r\cu_0\right)}
\left| \nabla_z P(s,y,z)  \right| 
\left| \nabla_z \tilde{P}_{x',r} (s,x,z) \right| 
\, dz\, ds.
\end{align*}
Using H\"older's inequality, 
Lemma~\ref{l.slicesvsaverages.PAR},
~\eqref{e.gradientP.NA} and $\sqrt{t} \leq r\leq \dist\left(x, \partial (x'+2r\cu_0)\right)$, we find that, for every $s\in (0,t]$ and $z'\in \Rd\setminus \left( x'+2r\cu_0\right)$, 
\begin{align*}
\lefteqn{
\left( \int_{B_{\sqrt{s}/2}(z')}
\left| \nabla_z P(s,y,z)  \right| 
\left| \nabla_z \tilde{P}_{x',r} (s,x,z) \right| 
\, dz \right)^2
} \qquad & 
\\ &
\leq
\int_{B_{\sqrt{s}/2}(z')}
\left| \nabla_z P(s,y,z)  \right|^2\,dz
\int_{B_{\sqrt{s}/2}(z')}
\left| \nabla_z \tilde{P}_{x',r} (s,x,z) \right|^2 \, dz
\\ & 
\leq
C s^{-2} \int_{Q_{\sqrt{s}}(s,z')}
\left| \nabla_z P(s',y,z)  \right|^2\,dz \,ds'
\int_{Q_{\sqrt{s}}(s,z')}
\left| \nabla_z \tilde{P}_{x',r} (s',x,z) \right|^2 \, dz \,ds'
\\ & 
\leq 
C s^{-d-2} \exp\left( - \frac{\left| y-z' \right|^2 + \left| x-z' \right|^2}{Cs} \right)
\\ & 
\leq 
Cs^{-d-2} \exp\left( - \frac{\left| x-z' \right|^2 + \left| x-y \right|^2}{Cs} \right).
\end{align*}
By covering~$\Rd\setminus \left( x'+2r\cu_0\right)$ with such balls $B_{\sqrt{s}/2}(z')$ and summing the previous display over the covering, we obtain, for every $x \in x' + r\cu_0$ and $y\in\Rd$,  
\begin{align*}
\int_{\Rd \setminus \left(x'+2r\cu_0\right)}
\left| \nabla_z P(s,y,z)  \right| 
\left| \nabla_z \tilde{P}_{x',r} (s,x,z) \right| 
\, dz
\leq 
C r^d s^{-d-1} 
\exp\left( - \frac{r^2 + \left|x-y\right|^2 }{Cs} \right).
\end{align*}
Integration over $s\in (0,t]$ yields, in view of~\eqref{e.duhamelPtileP},
\begin{align*}
\left| P(t,x,y) - \tilde{P}_{x',r} (t,x,y) \right| 
&
\leq
C r^d t^{-d} 
\exp\left( - \frac{r^2+\left|x-y\right|^2 }{Ct} \right)
\leq 
C t^{-\frac d2} 
\exp\left( - \frac{r^2+\left|x-y\right|^2 }{Ct} \right).
\end{align*}
This is~\eqref{e.Ploctildeest1}.

\smallskip

\emph{Step 2.} 
We show that, for every $x,x',y\in\Rd$,  $r\in [2,\infty)$ and $t\in (0,r^2]$ such that $x \in x'+r\cu_0$, we have
\begin{equation}
\label{e.Ploctildeest2}
\left\| \nabla P(t,x,\cdot) - \nabla \tilde{P}_{x',r} (t,x,\cdot) \right\|_{\underline{L}^2( B_{\sqrt{t}}(y))}
\leq 
C t^{-\frac12-\frac d2} 
\exp\left( - \frac{r^2+\left|x-y\right|^2 }{Ct} \right).
\end{equation}
By a standard covering argument, we may assume that $\sqrt{t} \leq \frac18r$. 
In the case that $y\in \Rd \setminus B_{r/2}(x)$, we use~\eqref{e.dumbdecay.greenie} to obtain
\begin{align*}
\lefteqn{
\left\| \nabla P(t,x,\cdot) - \nabla \tilde{P}_{x',r} (t,x,\cdot) \right\|_{\underline{L}^2( B_{\sqrt{t}}(y))}
} \quad & 
\\ & 
\leq 
\left\| \nabla P(t,x,\cdot) \right\|_{\underline{L}^2( B_{\sqrt{t}}(y))}
+
\left\| \nabla \tilde{P}_{x',r} (t,x,\cdot) \right\|_{\underline{L}^2( B_{\sqrt{t}}(y))}
\\ & 
\leq
C t^{-\frac12-\frac d2} \exp\left( -\frac{|x-y|^2}{Ct} \right) 
\leq 
C t^{-\frac12-\frac d2} 
\exp\left( - \frac{r^2+\left|x-y\right|^2 }{Ct} \right).
\end{align*}
In the alternative case that $y\in B_{r/2}(x)$, we have that $B_{4\sqrt{t}}(y) \subseteq x'+2r\cu_0$ and therefore $(t,y) \mapsto H(t,y):= P(t,x,y) - \tilde{P}_{x',r} (t,x,y)$ is a solution of 
\begin{equation*}
\partial_t H - \nabla \cdot \a \nabla H = 0 \quad \mbox{in} \ Q_{4\sqrt{t}}. 
\end{equation*}
Thus we may apply the parabolic Caccioppoli inequality (Lemma~\ref{l.cacciopp.this}), the gradient estimate for time slices (Lemma~\ref{l.slicesvsaverages.PAR}) and the estimate~\eqref{e.Ploctildeest1} above to obtain 
\begin{align*}
\left\| \nabla P(t,x,\cdot) - \nabla \tilde{P}_{x',r} (t,x,\cdot) \right\|_{\underline{L}^2(B_{\sqrt{t}}(y))}
&
\leq 
C t^{-\frac12} 
\left\| P(t,x,\cdot) - \tilde{P}_{x',r} (t,x,\cdot) \right\|_{\underline{L}^2(Q_{\sqrt{t}}(y))}
\\ & 
\leq 
C t^{-\frac12-\frac d2} 
\exp\left( - \frac{r^2+\left|x-y\right|^2 }{Ct} \right).
\end{align*}
This completes the proof of~\eqref{e.Ploctildeest2}. 

\smallskip

\emph{Step 3.} We construct the localized approximation~$P'_r$ of $P$. We first select a function $\psi\in C^\infty_c(2r\cu_0)$ such that $0\leq \psi \leq 1$, $\left\| \nabla \psi\right\|_{L^\infty(\Rd)} \leq Cr^{-1}$ and 
\begin{equation*}
\sum_{z\in  r\Zd} 
\psi(x-z) = 1, \quad \forall x\in\Rd. 
\end{equation*}
We can construct~$\psi$ by mollifying~$\indc_{r\cu_0}$, for instance. We also select a smooth cutoff function $\zeta \in C^\infty_c(10r\cu_0)$ such that $\zeta \equiv 1$ on $8r\cu_0$ and $\left\| \nabla \zeta \right\|_{L^\infty(\Rd)} \leq Cr^{-1}$. 

\smallskip

We now define
\begin{equation*}
\tilde{P}'_r(t,x,y)
:=
\sum_{z\in r\Zd} 
\psi(x-z) \zeta(y-z) \tilde{P}_{z,6r} (t,x,y).
\end{equation*}
We finally define $P'_r(t,x,y)$ by normalizing~$\tilde{P}'_r(t,x,y)$ in order to enforce the unit mass condition~\eqref{e.tildeP.prob}:
\begin{equation*}
P'_r(t,x,y) := \left( \int_{\Rd} \tilde{P}'_r(t,x,y) \,dx \right)^{-1} \tilde{P}'_r(t,x,y).
\end{equation*}
It is easy to check from the construction that~$P'_r(t,x,y)$ is $\mathcal{F}(x+15r\cu_0)$--measurable and that $P'_r(t,x,\cdot)$ is supported in $x+15r\cu_0$. 

\smallskip

What remains is to prove~\eqref{e.Plocalize1} and~\eqref{e.Plocalize2}. It is straightforward to obtain from the construction above,~\eqref{e.Ploctildeest1},~\eqref{e.Ploctildeest2} and the triangle inequality that, for every $x,y\in\Rd$, $r\in [2,\infty)$ and $t\in (0,r^2]$,  
\begin{equation}
\label{e.Ploctildeest.prime1}
\left| P(t,x,y) - \tilde{P}'_{r} (t,x,y) \right|  
\leq 
C t^{-\frac d2} 
\exp\left( - \frac{r^2+\left|x-y\right|^2 }{Ct} \right).
\end{equation}
and
\begin{equation}
\label{e.Ploctildeest.prime2}
\left\| \nabla P(t,x,\cdot) - \nabla \tilde{P}'_{r} (t,x,\cdot) \right\|_{\underline{L}^2( B_{\sqrt{t}}(y))}
\leq 
C t^{-\frac12-\frac d2} 
\exp\left( - \frac{r^2+\left|x-y\right|^2 }{Ct} \right).
\end{equation}
Therefore to complete the proof, it suffices by the triangle inequality to obtain the estimates 
\begin{equation}
\label{e.Ploctildeest.toprime1}
\left| P'_r(t,x,y) - \tilde{P}'_{r} (t,x,y) \right|  
\leq 
C t^{-\frac d2} 
\exp\left( - \frac{r^2+\left|x-y\right|^2 }{Ct} \right)
\end{equation}
and
\begin{equation}
\label{e.Ploctildeest.toprime2}
\left\| \nabla P'_r(t,x,\cdot) - \nabla \tilde{P}'_{r} (t,x,\cdot) \right\|_{\underline{L}^2( B_{\sqrt{t}}(y))}
\leq 
C t^{-\frac12-\frac d2} 
\exp\left( - \frac{r^2+\left|x-y\right|^2 }{Ct} \right).
\end{equation}
We have that 
\begin{align*}
\left| P'_r(t,x,y) - \tilde{P}'_{r} (t,x,y) \right| 
&
=
\left| \left( \int_{\Rd} \tilde{P}'_r(t,x,y) \,dx \right)^{-1} - 1 \right|
\tilde{P}'_{r} (t,x,y) 
\\ & 
\leq
\left| \left( \int_{\Rd} \tilde{P}'_r(t,x,y) \,dx \right)^{-1} - 1 \right|
t^{-\frac d2} \exp\left( - \frac{|x-y|^2}{Ct} \right). 
\end{align*}
Similarly, we find that 
\begin{align*}
\lefteqn{
\left\| \nabla P'_r(t,x,y) - \nabla \tilde{P}'_{r} (t,x,\cdot) \right\|_{\underline{L}^2( B_{\sqrt{t}}(y))}
} \quad & 
\\ &
\leq
C \left| \left( \int_{\Rd} \tilde{P}'_r(t,x,y) \,dx \right)^{-1} - 1 \right|
\left\| \nabla P'_r(t,x,y) - \nabla \tilde{P}'_{r} (t,x,\cdot) \right\|_{\underline{L}^2( B_{\sqrt{t}}(y))}
\\ & \quad 
+ 
C \left| \nabla_y \left( \int_{\Rd} \tilde{P}'_r(t,x,y) \,dx \right)^{-1} \right|
\left\| P'_r(t,x,y) - \tilde{P}'_{r} (t,x,\cdot) \right\|_{\underline{L}^2( B_{\sqrt{t}}(y))}
\\ & 
\leq 
C \left| \left( \int_{\Rd} \tilde{P}'_r(t,x,y) \,dx \right)^{-1} - 1 \right|
t^{-\frac12 -\frac d2} \exp\left( - \frac{|x-y|^2}{Ct} \right)
\\ & \quad
+
C \left| \nabla_y \left( \int_{\Rd} \tilde{P}'_r(t,x,y) \,dx \right)^{-1} \right|
t^{-\frac d2} \exp\left( - \frac{|x-y|^2}{Ct} \right).
\end{align*}
Therefore what remains is to show that the normalizing constant is sufficiently close to one, and this follows from~\eqref{e.Ploctildeest.prime1} and~\eqref{e.Ploctildeest.prime2}. Indeed, we have
\begin{align*}
\left| \int_{\Rd} \tilde{P}'_r(t,x,y) \,dx - 1 \right|
&
=
\left| \int_{\Rd} \left( \tilde{P}'_r(t,x,y) - P(t,x,y) \right) \,dx \right|
\\ & 
\leq 
\int_{\Rd} \left| \tilde{P}'_r(t,x,y) - P(t,x,y) \right| \,dx
\\ & 
\leq \int_{\Rd}
C t^{-\frac d2} 
\exp\left( - \frac{r^2+\left|x-y\right|^2 }{Ct} \right)\,dx
\\ & 
= 
C\exp\left( - \frac{r^2}{Ct} \right),
\end{align*}
and, similarly, 
\begin{align*}
\left|  \int_{\Rd} \nabla_y \tilde{P}'_r(t,x,y) \,dx  \right|
&
=
\left| \int_{\Rd} \nabla_y \left( \tilde{P}'_r(t,x,y) - P(t,x,y) \right) \,dx \right|
\\ & 
\leq 
\int_{\Rd} \left| \nabla_y \tilde{P}'_r(t,x,y) - \nabla_y P(t,x,y) \right| \,dx
\\ & 
\leq \int_{\Rd}
C t^{-\frac12-\frac d2} 
\exp\left( - \frac{r^2+\left|x-y\right|^2 }{Ct} \right)\,dx
\\ & 
= 
Ct^{-\frac12} \exp\left( - \frac{r^2}{Ct} \right).
\end{align*}
This completes the proof. 
\end{proof}

Using the localized Green functions, we can localize the function~$u(t,x)$ itself. For each $r \in [2,\infty)$, $t \in (0,r^2]$ and $x \in \Rd$, we define
\begin{equation}  \label{e.SG.u'.def}
u'(r,t,x) := -\int_{\Rd} \g(y) \cdot \nabla_y P'_r(t,x,y) \, dy,
\end{equation}
where $P'_r(t,x,y)$ is the function given by Lemma~\ref{l.semigp.loc}.

\begin{lemma}[Localization of~$u$]
\label{l.SG.loc.u}
There exists a constant $C(d,\Lambda) < \infty$ such that for every $r \in [2,\infty)$, $t \in (0,r^2]$ and $x \in \Rd$, the random variable $u'(r,t,x)$ is $\mcl F(B_r(x))$-measurable and satisfies the estimates
\begin{equation}  
\label{e.SG.loc.u1}
\Ll|u(t,x) - u'(r,t,x) \Rr| \le C t^{-\frac 1 2} \exp \Ll( - \frac{r^2}{Ct}\Rr)
\end{equation}
and 
\begin{equation}  
\label{e.SG.loc.u2}
\Ll|\nabla u(t,x) - \nabla u'(r,t,x) \Rr| \le C t^{-1}\exp \Ll( - \frac{r^2}{Ct}\Rr).
\end{equation}
\end{lemma}
\begin{proof}
Since for every $x \in \Rd$, the random variable $\g(x)$ is $\mcl F(B_1(x))$-measurable, we obtain that $u'(r,t,x)$ is $\mcl F(B_r(x))$-measurable from \eqref{e.Pprimemeas}.
The bounds \eqref{e.SG.loc.u1} and \eqref{e.SG.loc.u2} are immediate consequences of the estimates \eqref{e.Plocalize1} and \eqref{e.Plocalize2} and the fact that $\|\g\|_{L^\infty(\Rd)} \le 1$. 
\end{proof}

The strategy for the proof of Theorem~\ref{t.semigroup} is to propagate information on the size of spatial averages of $u(t,\cdot)$ to larger and larger times. 
The precise statement  to be propagated takes the following form: given $\sigma \in (0,2)$, and $T, \mathsf{K} \in [1,\infty)$, we denote by $\S(T,\mathsf K,  \sigma)$ the statement that, for every $t \in [1,T]$, $s \in \left[  t ,\infty\right)$ and $x \in \Rd$, 
\begin{equation}
\label{e.whatweprop}
\int_{\Rd} u(t,y) \bar P(s,x-y) \, dy = \O_\sigma \Ll( \mathsf K s^{-\frac 1 2 -\frac d 4} \Rr).
\end{equation}

We show first that~$\S(T,\mathsf K, \sigma)$ implies  pointwise estimates on~$u$. 

\begin{lemma}
\label{l.SG.pointwise.u}
Fix $\sigma \in [1,2)$. There exists $C(\sigma,d,\Lambda) < \infty$ such that the following implication holds for every $T, \mathsf K \ge 1$:
\begin{equation*}  
\S(T, \mathsf K,\sigma) 
\implies 
\forall (t,x) \in [1,2T] \times \Rd, \quad 
u(t,x) = \O_\sigma \Ll( C \mathsf K t^{-\frac 1 2 - \frac d 4} \Rr) .
\end{equation*}
\end{lemma}

The following lemma is needed in the proof of Lemma~\ref{l.SG.pointwise.u}, so we present it first. It states roughly that, up to a lower-order error and assuming certain pointwise bounds on~$u$, we can replace the homogenized Green function~$\overline{P}$ in the integral on the left side of~\eqref{e.whatweprop} with the heterogeneous Green function~$P$.

\begin{lemma} 
There exist an exponent $\delta(d,\Lambda) > 0$ and, for each $\sigma \in (0,2)$, a constant $C(\sigma,d,\Lambda) < \infty$ such that the following holds. Assume that for some $\mathsf K \in [1,\infty)$ and $t \in [1,\infty)$,  we have for every $x \in \R^d$ that
\label{l.SG.opt5}
\begin{equation} \label{e.SG.opt3.a1}
\left| u(t,x) \right|  \leq  \mathsf K  t^{-\frac12} \wedge  \O_{\sigma}\left( \mathsf K t^{-\frac12 - \frac d4} \right).
\end{equation}
Then, for every $s \in [t,\infty)$ and $x \in \Rd$, we have the estimate
\begin{equation} \label{e.SG.opt3.r1}
\int_{\R^d} u(t,y) \left( P(s,x,y) - \overline{P}(s,x-y) \right) \, dy = \O_{\sigma}\left( C s^{-\delta(2-\sigma)} \mathsf K t^{-\frac12 - \frac d4}  \right).
\end{equation}
\end{lemma}
\begin{proof}
We fix $\sigma \in (0,2)$, $K, t \in [1,\infty)$ such that \eqref{e.SG.opt3.a1} holds, $\theta :=  \frac{d(\sigma+2)}{4}\in (0,d)$, and denote by $\Y_\theta(y)$ the random variable given by Theorem~\ref{t.PGF.basecase}. This theorem guarantees the existence of an exponent $\delta(d,\Lambda) > 0$ and a constant $C(\sigma,d,\Lambda) < \infty$ such that for every $s \in [1,\infty)$ and $x,y \in \Rd$,
\begin{equation*} 
\left|  P(s,x,y) - \overline{P}(s,x-y)   \right| \indc_{ \left\{ \Y_\theta(y) \le  \sqrt{s}\right\} }  \leq  Cs^{-\delta(d-\sigma)} s^{- \frac d2}  \exp\left(-  \frac{|x-y|^2}{Cs}\right).
\end{equation*}
We thus deduce from \eqref{e.SG.opt3.a1} and Lemma~\ref{l.sum-O} that
\begin{equation*}  
\int_{\R^d} u(t,y) \left( P(s,x,y) - \overline{P}(s,x-y) \right)\indc_{ \left\{ \Y_\theta(y) \le  \sqrt{s}\right\} }  \, dy = \O_{\sigma}\left( C s^{-\delta(d-\sigma)} \mathsf K t^{-\frac12 - \frac d4}  \right).
\end{equation*}
On the other hand, by the Nash-Aronson bound \eqref{e.upperP.NA},  we have
\begin{equation*}  
\left|  P(s,x,y) - \overline{P}(s,x-y)   \right| \indc_{ \left\{ \Y_\theta(y) >  \sqrt{s}\right\} }  \le C \Ll(\frac{\Y_\theta(y)}{\sqrt{s}}\Rr)^{\frac \theta \sigma} s^{- \frac d2}  \exp\left(-  \frac{|x-y|^2}{Cs}\right),
\end{equation*}
and by our choice of $\sigma$, we have $s^{-\frac \theta {2\sigma}} = s^{-\frac d 4 - \frac{2-\sigma}{8\sigma}}$. Using also that $\mcl Y_\theta = \O_\theta(C)$ and the deterministic bound in \eqref{e.SG.opt3.a1}, we thus obtain the result.
\end{proof}

\begin{proof}[Proof of Lemma~\ref{l.SG.pointwise.u}]
For each $T',\msf K' \in [1,\infty)$, we denote by $\S'(T',\mathsf K ',\sigma)$ the statement that, for every $t \in[1,T']$ and $x \in \Rd$, we have
\begin{equation*}  
u(t,x) = \O_\sigma \Ll( \mathsf K' t^{-\frac 1 2 - \frac d 4} \Rr) .
\end{equation*}
We decompose the proof into two steps.

\smallskip

\emph{Step 1.} In this first step, we show that there exists a constant $T_0(\sigma,d,\Lambda) < \infty$ such that under the assumption of $\S(T,\msf K, \sigma)$, we have for every $T' \in [T_0,T]$ and  $\msf K' \in [1,\infty)$ that
\begin{equation}  
\label{e.impl.SG.pointwise.u}
\S'(T',\msf K',\sigma) \implies \S'\Ll(2T',\msf K' \vee \Ll(\msf K + \frac {\msf K'}{2}\Rr),\sigma\Rr) ,
\end{equation}
We thus assume that $\S(T,\msf K, \sigma)$ and $\S'(T',\msf K',\sigma)$ hold, and fix $t \in \left( \tfrac {T'}{2}, T'\right]$ and $x\in\R^d$. The argument is based on the decomposition
\begin{multline} \label{e.decomp.again}
 u(2t,x)   = \int_{\Rd} u(t,y) P(t,x,y) \, dy \\
  = \int_{\Rd} u(t,y) \bar P(t,x-y) \, dy + \int_{\Rd} u(t,y) \Ll( P(t,x,y) - \bar P(t,x-y) \Rr) \, dy.
\end{multline}
The first term on the right side can be easily controlled using the assumption of~$\S(T,\mathsf K,\sigma)$, which gives us that
\begin{equation*} 
\int_{\Rd} u(t,y) \bar P(t,x-y) \, dy = \O_\sigma \Ll( \mathsf K t^{-\frac 1 2 -\frac d 4} \Rr).
\end{equation*}
We now turn to the estimation of the second term on the right side of \eqref{e.decomp.again}. By Lemma~\ref{l.crude.ubound} and the assumption of $\S'(T',\mathsf K',\sigma)$, we have, for  every $z \in \R^d$, that 
\begin{equation*} 
\left| u(t,z) \right| \leq C t^{-\frac 12} \wedge  \O_\sigma \Ll( \mathsf K' t^{-\frac 1 2 - \frac d 4} \Rr) . 
\end{equation*}
By Lemma~\ref{l.SG.opt5}, we deduce that, for some $\delta(d,\Lambda) > 0$ and $C(\sigma,d,\Lambda)<\infty$, 
\begin{equation*} 
\int_{\R^d} u(t,y) \left( P(t,x,y) - \overline{P}(t,x-y) \right) \, dy = \O_{\sigma}\left( C t^{-\delta(2-\sigma)} \mathsf K' t^{-\frac12 - \frac d4}  \right).
\end{equation*}
The implication \eqref{e.impl.SG.pointwise.u} thus follows from these estimates and Lemma~\ref{l.sum-O}, provided that we choose $T_0(\sigma,d,\Lambda) < \infty$ sufficiently large that, for the constant $C(\sigma,d,\Lambda) < \infty$ in the previous display,
\begin{equation*}  
C \Ll( \frac{T_0}{2} \Rr)^{-\delta(2-\sigma)} \le 2^{-\frac 3 2 -\frac d 4}.
\end{equation*}

\smallskip

\emph{Step 2.} We complete the proof. We fix the constant $T_0( \sigma, d,\Lambda) < \infty$ identified in the previous step. By Lemma~\ref{l.crude.ubound}, there exists $\msf K_0(d,\Lambda) <\infty$ such that $\S'(T_0,\msf K_0,\sigma)$ holds. In particular, the statement $\S'(T_0,2 \msf K + \msf K_0,\sigma)$ holds. A recursive application of \eqref{e.impl.SG.pointwise.u} then yields that $\S'(2T,2 \msf K + \msf K_0,\sigma)$, and this implies the announced result.
\end{proof}

In the next lemma, we obtain the optimal decay in $s^{-\frac 1 2 - \frac d 4}$ for spatial averages of $u(t,\cdot)$ on scales of order $s^\frac 1 2$. At this stage we pay no attention to obtaining the right scaling in the $t$ variable. This lemma will serve us to establish the base case for propagating the statement $\S(T,\msf K,\sigma)$ forward in time.

\begin{lemma}
\label{l.SG.opt1}
There exists $C(d,\Lambda) < \infty$ such that for every $t \in [1,\infty)$,  $s  \in [t,\infty)$ and $x \in \Rd$,
\begin{equation}  
\label{e.SG.opt1.r}
\int_{\Rd} u(t,y) \bar P(s,x-y) \, dy = \O_2 \Ll( C t^{ \frac d 4} s^{-\frac 1 2 - \frac d 4} \Rr) .
\end{equation}
\end{lemma}
\begin{proof}
Without loss of generality, we focus on proving \eqref{e.SG.opt1.r} for $x = 0$. For every $t \in [1,\infty)$ and $r \ge t^\frac 1 2$, we write
\begin{equation*} 
w(t,x,y) := \g(y) \cdot \nabla_y P(t,x,y)
\quad \mbox{and} \quad 
w' (r,t,x,y) := \g(y) \cdot \nabla_y P'_r(t,x,y) ,
\end{equation*}
where $P'_r(t,x,y)$ is the function defined in Lemma~\ref{l.semigp.loc}.
By \eqref{e.SG.u.def} and \eqref{e.Plocalize2}, we have
\begin{equation*} 
u(t,x) = - \int_{\R^d} w(t,x,y) \, dy = - \lim_{r \to \infty} \int_\Rd  w'(r,t,x,y)\, dy.
\end{equation*}
Denoting  
\begin{equation*} 
\tilde w'(r,t,x,y) := w'(r,t,x,y) - \E\left[ w'(r,t,x,y) \right],
\end{equation*}
we can thus split the left side of \eqref{e.SG.opt1.r} using the identity 
\begin{align}  
\label{e.SG.opt1.0.1}
\notag
\lefteqn{\int_{\R^d}  u(t,y)  \overline{P}(s,y) \, dy 
}
\quad & 
\\ \notag 
&  
=  \int_{\R^d} \E\left[ u(t,y) \right]  \overline{P}(s,y) \, dy
 -  \int_{\R^d} \int_{\R^d} \tilde w'(t^\frac 1 2,t,y,z)  \overline{P}(s,y) \, dz \, dy
\\ & \quad   - \sum_{k=0}^\infty \int_{\R^d} \int_{\R^d}  \left( \tilde w'(2^{k+1}t^\frac 1 2, t ,y ,z) - \tilde w'(2^{k}t^\frac 1 2,t ,y,z) \right) \overline{P}(s,y) \, dz \, dy.
\end{align}
We decompose the estimation of each of these terms into four steps.

\smallskip

\emph{Step 1.}
We show that there exists $C(d,\Lambda)<\infty$ such that for every $s,t \in [1,\infty)$,
\begin{equation}
\label{e.SG.opt1.1.1}
\Ll|\int_{\Rd} \E[u(t,y)] \, \bar P(s,y) \, dy \Rr| \le C t^{-\frac 12} \exp \Ll( -C^{-1}  s\Rr) .
\end{equation}
Note that the function $y \mapsto \E[u(t,y)]$ is $\Z^d$-periodic, and by Lemma~\ref{l.crude.ubound}, it is bounded by $C t^{-\frac 1 2}$. Arguing as in Exercise~\ref{ex.heatdecay}, we see that that proof of \eqref{e.SG.opt1.1.1} amounts to the verification of the fact that for every $t > 0$,
\begin{equation}  
\label{e.null.average}
\int_{[0,1]^d} \E[u(t,y)] \, dy = 0.
\end{equation}
Using again that the function $y \mapsto u(t,y)$ is $\Z^d$-stationary, the equation for $u$ and~Lemma~\ref{l.crude.ubound}, we see that, for every $t_2>t_1>0$, for $\psi \in C_c^\infty(B_1)$ with mass one and $\psi_R = R^{-d}\psi\left( \tfrac{\cdot}{R}\right)$, we have
\begin{multline*} 
\int_{[0,1]^d} \E\left[ u(t_1,y) \right] \, dy -  \int_{[0,1]^d} \E\left[ u(t_2,y) \right] \, dy
\\  =  \lim_{R \to \infty} \E\left[\int_{t_1}^{t_2} \int_{\R^d}  \nabla u(t',y) \cdot \a(y)\nabla \psi_R(y)\, dy \, dt'\right] = 0. 
\end{multline*}
Hence, the quantity on the left side of \eqref{e.null.average} does not depend on $t > 0$. Since by Lemma~\ref{l.crude.ubound}, this quantity decays to $0$ as $t$ tends to infinity, this proves \eqref{e.null.average}.

\smallskip

\smallskip

\emph{Step 2.}
We next show, using the unit range of dependence assumption, that there exists a constant $C(d,\Lambda) < \infty$ such that for every $t \ge 1$, $r\ge t^\frac 1 2$ and $s \ge t$,
\begin{equation}  \label{e.SG.opt1.2.1}
\left| \int_{\R^d} \int_{\R^d} \tilde w'(r,t,y,z)  \overline{P}(s,y) \, dz \, dy \right| \leq \O_2\left( C \left( \frac{r^2}{t}\right)^{\frac {1}{2}}  \left( \frac{r^2}{s}\right)^{\frac d4} s^{-\frac 12} \right).
\end{equation}
To see this, we rewrite the left side as
\begin{align} \notag \label{e.SG.opt1.2.2}
\lefteqn{\int_{\R^d} \int_{\R^d} \tilde w'(r,t,y,z)  \overline{P}(s,y) \, dz \, dy} \quad &
\\ \notag & = \sum_{z' \in r\Z^d}  \int_{\R^d}  \int_{z' + r \cu_0}   \tilde w'(r,t,y,z) \, dz \, \overline{P}(s,y)    \, dy
\\ \notag & = \sum_{z' \in r\Z^d}  \int_{\R^d}  \int_{z' + r \cu_0}   \tilde w'(r,t,y,z) \, dz \left( \overline{P}(s,y) - \left( \overline{P}(s,\cdot) \right)_{z' + 3 r \cu_0 } \right)    \, dy
\\  & = \sum_{z' \in r\Z^d} \int_{z' + 3 r \cu_0}  \int_{z' + r \cu_0}  \tilde w'(r,t,y,z)  \, dz \left( \overline{P}(s,y) - \left( \overline{P}(s,\cdot) \right)_{z' + 3 r \cu_0 } \right)   \, dy .
\end{align}
Indeed, the first equality is just splitting the integral; the second one follows from~\eqref{e.tildeP.prob}, since
\begin{equation*} 
\int_{\R^d}     w'(r,t,y,z) \, dy =  \int_{\R^d}  \g(z) \cdot \nabla_z P'_r(t,y,z) \, dy  =  \g(z) \cdot \nabla_z \int_{\R^d} P'_r(t,y,z) \, dy  = 0;
\end{equation*}
and the third one follows from the fact that $y \mapsto \nabla P'_r(t,y,z)$ is supported in $B_r(z)$. 
To estimate the resulting sum, we use the independence assumption by way of Lemma~\ref{l.barO.boxes}. To prepare for the application of this lemma, we observe that, for every $z'\in r\Zd$, the random variable 
\begin{multline*}
\int_{z' + 3 r \cu_0}  \int_{z' + r \cu_0}  \tilde w'(r,t,y,z)  \, dz \left( \overline{P}(s,y) - \left( \overline{P}(s,\cdot) \right)_{z' + 3 r \cu_0 } \right)   \, dy
\\ 
\mbox{is $\F(z'+4r\cu_0)$--measurable}
\end{multline*}
and its expectation vanishes since~$\E\left[  \tilde w'(r,t,y,z)  \right] = 0$. 
The measurability follows from~\eqref{e.Pprimemeas} and the definition of $\tilde{w}'$ above. Next we observe that, by Lemma~\ref{l.semigp.loc}, for every $z'\in r\Zd$,
\begin{equation}
\label{e.some.bound.for.w'}
\left| \int_{z' + r \cu_0} \tilde w'(r,t,y,z) \, dz \right| \leq 2 \|\g\|_{L^\infty(\R^d)}  \left\| \nabla_y P'_r(t,y,\cdot) \right\|_{L^1(\R^d)} \leq C t^{-\frac 12} .
\end{equation}
Thus, by the Poincar\'e inequality, we obtain that
\begin{multline*} 
\left| \int_{z' + 3 r \cu_0} \int_{z' + r \cu_0} \tilde w'(r,t,y,z) \, dz \left( \overline{P}(s,y) - \left( \overline{P}(s,\cdot) \right)_{z' + 3 r \cu_0 } \right) \, dy \right|
\\ \leq C \left( \frac{r^2}{t}\right)^{\frac {1}{2}}   \left\| \nabla \overline{P}(s,\cdot) \right\|_{L^1\left(z' + 3 r \cu_0\right) }.
\end{multline*}
We now apply Lemma~\ref{l.barO.boxes} to obtain, for each of the~$3^d$ elements~$z''$ of~$r\Zd\cap (3r\cu_0)$, 
\begin{multline*}
\left| \sum_{z' \in z''+3r\Z^d} \int_{z' + 3 r \cu_0}  \int_{z' + r \cu_0}  \tilde w'(r,t,y,z)  \, dz \left( \overline{P}(s,y) - \left( \overline{P}(s,\cdot) \right)_{z' + 3 r \cu_0 } \right)   \, dy \right|
\\ 
\leq 
\O_2\left( C \left( \frac{r^2}{t}\right)^{\frac {1}{2}}  \left(  \sum_{z' \in z''+3r\Z^d} \left\| \nabla \overline{P}(s,\cdot) \right\|_{L^1\left(z' + 3 r \cu_0\right) }^2 \right)^{\frac 12} \right).
\end{multline*}
Summing over $z''\in r\Zd\cap (3r\cu_0)$ and using the triangle inequality, we obtain
\begin{multline*}
\left| \sum_{z' \in 3r\Z^d} \int_{z' + 3 r \cu_0}  \int_{z' + r \cu_0}  \tilde w'(r,t,y,z)  \, dz \left( \overline{P}(s,y) - \left( \overline{P}(s,\cdot) \right)_{z' + 3 r \cu_0 } \right)   \, dy \right| 
\\
\leq 
\O_2\left( C \left( \frac{r^2}{t}\right)^{\frac {1}{2}}  \left(  \sum_{z' \in r\Z^d} \left\| \nabla \overline{P}(s,\cdot) \right\|_{L^1\left(z' + 3 r \cu_0\right) }^2 \right)^{\frac 12} \right).
\end{multline*}
The expression on the right can be further estimated using H\"older's inequality by
\begin{equation*} 
 \left(  \sum_{z' \in r\Z^d} \left\| \nabla \overline{P}(s,\cdot) \right\|_{L^1\left(z' + 3 r \cu_0\right) }^2 \right)^{\frac 12} \leq C r^{\frac{d}{2}} \left\| \nabla \overline{P}(s,\cdot)  \right\|_{L^2(\R^d)} \leq C \left( \frac{r^2}{s}\right)^{\frac d4} s^{-\frac 12},
\end{equation*}
which gives~\eqref{e.SG.opt1.2.1}. 

\smallskip

\emph{Step 3.}
We show that there exists $C(d,\Lambda)<\infty$ such that for every $t \ge 1$, $r \ge t^\frac 1 2$ and $s \ge t$, 
\begin{multline}  \label{e.SG.opt1.3.1}
\left| \int_{\R^d} \int_{\R^d}  \left( \tilde w'(2r ,t,y ,z) - \tilde w'(r ,t,y,z) \right) \overline{P}(s,y) \, dz \, dy \right| \\\leq \O_2\left( C \left( \frac{r^2}{t}\right)^{\frac {1}{2}}  \left( \frac{r^2}{s}\right)^{\frac d4} s^{-\frac 12} \exp \Ll( -\frac{r^2}{Ct} \Rr) \right).
\end{multline}
The proof is almost identical to that of the previous step, except that we use Lemma~\ref{l.semigp.loc} to improve the bound \eqref{e.some.bound.for.w'} to
\begin{equation*}  
\left| \int_{z' + r \cu_0} \tilde w'(r,t,y,z) \, dz \right|  \le C t^{-\frac 1 2} \exp \Ll( -\frac {r^2}{Ct} \Rr) .
\end{equation*}

\smallskip

\emph{Step 4.} Using the decomposition in \eqref{e.SG.opt1.0.1} with \eqref{e.SG.opt1.1.1}, \eqref{e.SG.opt1.2.1}, \eqref{e.SG.opt1.3.1} and Lemma~\ref{l.sum-O}, we obtain the announced result.
\end{proof}

The propagation of the assumption $\S(T,\msf K, \sigma)$ forward in time is based on the following splitting of the left side of~\eqref{e.whatweprop}.

\begin{lemma} \label{l.SG.basic}
For every $t_1,s \in [1,\infty)$, $t_{2}\in[t_1,\infty)$, and $x \in \R^d$,  we have
\begin{multline} 
 \label{e.SG.basic}
\int_{\R^d} u(t_{2},y) \overline{P}(s,x-y) \, dy 
 = \int_{\R^d} u(t_{1},y) \overline{P}(s+t_{2}-t_{1},x-y) \, dy
\\  
+ \int_{t_{1}}^{t_2} \int_{\R^d} \left(\ahom - \a(y) \right) \nabla u\left(t ,y \right) \cdot \nabla \overline{P}(s + t_{2}-t,x-y) \, dy  \, dt . 
\end{multline}

\end{lemma}
\begin{proof}
We fix $t_1,s \in [1,\infty)$ and $t_{2}\in[t_1,\infty)$. For each $t \in (0,s+t_{2})$ , and $x \in \R^d$,  we set
\begin{equation*} 
f(t) := \int_{\R^d} u( t,y) \overline{P}(s+t_{2}- t ,x-y) \, dy ,
\end{equation*}
and observe that
\begin{align} 
\notag  
\partial_t f(t)  & =   \int_{\R^d} \left( \partial_{t} u( t,y) \overline{P}(s+t_{2}- t ,x-y) - u( t,y) \partial_{t}\overline{P}(s+t_{2} - t ,x-y) \right) \, dy.
\end{align}
Using the equation for $u$ in \eqref{e.flow}, we get
\begin{align} 
\notag  
\lefteqn{ \int_{\R^d}  \partial_{t} u( t,y) \overline{P}(s+t_{2}- t ,x-y) \,dy } \quad  &  
\\ \notag  & = 
\int_{\R^d}  \left( \ahom - \a(y)\right) \nabla u( t,y) \cdot \nabla \overline{P}(s+t_{2}- t ,x-y)  \, dy  
\\ \notag  & \quad - \int_{\R^d}   \nabla u( t,y) \cdot \ahom \nabla \overline{P}(s+t_{2}- t ,x-y)  \, dy
\\ \notag  & = 
\int_{\R^d}  \left( \ahom - \a(y)\right) \nabla u( t,y) \cdot \nabla \overline{P}(s+t_{2}- t ,x-y)  \, dy  
\\ \notag  & \quad + \int_{\R^d}  u( t,y) \partial_{t} \overline{P}(s+t_{2}- t ,x-y)  \, dy.
\end{align}
The result follows by combining the two previous displays and integrating over the interval $(t_{1},t_{2})$. 
\end{proof}

In order to propagate the assumption of $\S(T, \msf K, \sigma)$, we will argue that the last term on the right side of~\eqref{e.SG.basic} is of lower order compared with the first term, which is simply the spatial average of~$u$ at an earlier time. Intuitively, the reason why the last term in \eqref{e.SG.basic} is of lower order is an idea we encountered already many times throughout the book: the spatial averages of $\a \nabla u$ and $\ahom \nabla u$ should be close (see e.g.\ Theorem~\ref{t.DP.blackbox} and Lemmas~\ref{l.commutator} and~\ref{l.brz}). We first show the corresponding result for the parabolic Green function using a two-scale expansion argument.

\begin{lemma} \label{l.SG.Pfluxes}
Fix $\sigma \in (0,d)$ and $\alpha \in (0,\Lambda^{-1})$. There exist an exponent $\delta(d,\Lambda) > 0$, a constant $C(\sigma,\alpha,d,\Lambda) < \infty$, and, for each $x \in \Rd$, a random variable $\Y_\sigma(x)$ satisfying
\begin{equation} \label{e.SG.Y}
\Y_\sigma(x) \leq \O_\sigma\left( C \right)
\end{equation}
such that for every $x,y \in \Rd$, $t \ge 1 \vee \Y^2_\sigma(y)$ and $s \ge 1 \vee \Y^2_\sigma(x)$, we have
\begin{multline} \label{e.SG.Pfluxes}
\left| \int_{\Rd} \left( \a(z) - \ahom\right) \nabla_x P\left( t, z, y \right) \cdot \nabla \overline{P}(s,x-z) \,dz \right| 
\\\leq C t^{-\delta(d-\sigma)- \frac12} s^{-\frac 12}(t+s)^{-\frac d2} \exp\left( - \alpha \frac{|x-y|^2}{4(t+s)} \right).
\end{multline}
\end{lemma}

\begin{proof} We denote by $\X_\sigma$ the maximum of the two random variables appearing in Lemma~\ref{l.correctorminbounds} and Proposition~\ref{p.fluxcorrectorests}, and denote by $(\X_\sigma(x), x \in \Rd)$ its stationary extension, see Remark~\ref{r.Lipminscale}. We define $\Y_\sigma(x)$ to be the maximum between $\X_\sigma(x)$ and the random variable appearing in Theorem~\ref{t.PGF.basecase}. Note that there exists $C(\sigma,d,\Lambda) < \infty$ such that $\Y_\sigma = \O_\sigma(C)$.

\smallskip

Throughout the proof, we fix $\alpha \in (0,\Lambda^{-1})$, $x,y \in \R^d$ and $t,s \in [1,\infty)$ satisfying $t^{\frac 12} \geq \Y_\sigma(y)$ and $s^\frac 1 2 \ge \Y_\sigma(x)$. Let $H$ be the two-scale expansion of~$\overline{P}$ defined in~\eqref{e.def.H}. We write
\begin{align} \notag  \label{e.SG.Pfluxes.1}
\lefteqn{\left| \int_{\Rd}  \left( \a(z) - \ahom\right)  \nabla_x P\left( t, z, y \right) \overline{P}(s,x-z) \,dz \right|} \quad &
\\ \notag & \leq  \left| \int_{\Rd}  \left( \a(z) - \ahom\right) \nabla_x H\left( t, z, y \right) \overline{P}(s,x-z) \,dz \right|
\\ & \quad +  \left| \int_{\Rd}  \left( \a(z) - \ahom\right)  \left(\nabla_x P\left( t, z, y \right) -  \nabla_x H\left( t, z,y \right) \right) \overline{P}(s,x-z) \,dz  \right|.
\end{align}
We estimate the two terms appearing on the right side of \eqref{e.SG.Pfluxes.1} separately in the following two steps.

\smallskip

\emph{Step 1.} We show that there exist $\delta(d,\Lambda)>0$ and $C(\sigma,\alpha,d,\Lambda)<\infty$ such that 
\begin{multline}  \label{e.SG.Pfluxes.2}
\left| \int_{\Rd}  \left( \a(z) - \ahom\right) \nabla_x H\left( t, z, y \right) \cdot \nabla \overline{P}(s,x-z) \,dz \right| \\\leq 
 C t^{ -\delta(d-\sigma)-\frac12} \left( \frac{t+s}{s} \right)^{\frac d4} (t+s)^{-\frac d2} \exp\left( - \alpha \frac{|x-y|^2}{4(t+s)} \right).
\end{multline}
By the definition of $H$, we see that 
\begin{align} 
\notag  
\lefteqn{\left( \a(z) - \ahom\right) \nabla_x H\left( t, z, y \right)} \quad  &  
\\ 
\notag & = \sum_{k=1}^d  \left( \a(z)\left(e_k  + \nabla \phi_{e_k}(z) \right) - \ahom e_k  -  \ahom \nabla \phi_{e_k}(z) \right) \partial_{x_k} \overline{P}(t,z-y)
\\ 
\notag & \quad + \sum_{k=1}^d \left( \phi_{e_k}(z) {-} \left( \phi_{e_k} {\ast} \Phi(t,\cdot) \right)(y) \right) \left( \a(z) - \ahom\right)  \partial_{x_k} \nabla \overline{P}(t,z-y) 
\end{align}
By the definition of the flux correctors, see~\eqref{e.yesfluxcorrector}, we have that
\begin{equation*} 
\a(z)\left(e_k  + \nabla \phi_{e_k}(z) \right) - \ahom e_k = \left(\nabla \cdot \mathbf{S}_{e_{k}}\right)(z),
\end{equation*}
and hence, reorganizing the terms, we get
\begin{align} 
\notag  
\lefteqn{\left( \a(z) - \ahom\right) \nabla_x H\left( t, z, y \right)  = \sum_{k=1}^d    \left[ \nabla \cdot \mathbf{S}_{e_{k}} \right](z)\partial_{x_k} \overline{P}(t,z-y)  } \quad  &  
 \\ \notag &\qquad \qquad \qquad  + \sum_{k=1}^d \left[ \ahom \nabla \left( \left( \phi_{e_k} - \left( \phi_{e_k} {\ast} \Phi(t,\cdot) \right)(y) \right) \partial_{x_k}  \overline{P}(t,\cdot-y) \right) \right](z)
\\ \notag & \qquad \qquad \qquad + \sum_{k=1}^d \left( \phi_{e_k}(z) {-} \left( \phi_{e_k} {\ast} \Phi(t,\cdot) \right)(y) \right) \a(z) \partial_{x_k} \nabla \overline{P}(t,z-y) .
\end{align}
Therefore, after integrating by parts, we obtain, for $\beta := \frac{1}{2}(\Lambda^{-1} + \alpha)$, 
\begin{align} 
\notag  
\lefteqn{\left| \int_{\Rd}  \left( \a(z) - \ahom\right) \nabla_x H\left( t, z, y \right) \cdot \nabla \overline{P}(s,x-z) \,dz \right|} \quad  &  
\\ 
\notag & 
\leq Ct^{-1} s^{-\frac 12}\sum_{k=1}^d \int_{\Rd} \left| \phi_{e_k}(z) - \left( \phi_{e_k} {\ast} \Phi(t,\cdot) \right)(y)\right|  \Phi\left(\frac{t}{\beta},z-y \right) \Phi\left(\frac{s}{\beta},z-x \right) \, dz
\\ 
\notag &
 \quad + Ct^{-1}s^{-\frac 12} \sum_{k=1}^d \int_{\Rd} \left| \mathbf{S}_{e_{k}}(z) - \left( \mathbf{S}_{e_{k}} \right)_{B_{\sqrt{s}}(y)} \right| \Phi\left(\frac{t}{\beta},z-y \right) \Phi\left(\frac{s}{\beta},z-x \right) \, dz .
\end{align}
Writing 
\begin{equation*} 
f(z) := \sum_{k=1}^d  \left( \left| \phi_{e_k}(z) - \left( \phi_{e_k} {\ast} \Phi(t,\cdot) \right)(y)\right|  + \left| \mathbf{S}_{e_{k}}(z) - \left( \mathbf{S}_{e_{k}} \right)_{B_{\sqrt{s}}(y)} \right|  \right),
\end{equation*}
we use Step 1 of the proof of Lemma~\ref{l.PGF.layercake} to obtain
\begin{multline*} 
\left| \int_{\Rd}  \left( \a(z) - \ahom\right) \nabla_x H\left( t, z, y \right) \overline{P}(s,x-z) \,dz \right|
\\ \leq  C t^{-1} s^{-\frac 12}  \int_0^\infty r^{d+1} \exp\left( - \beta \frac{r^2}{4} \right)  \left\| f \Phi\left(\frac{s}{\beta},\cdot-x \right)  \right\|_{\underline{L}^1\left( B_{rt^{1/2}}(y)\right) } \, dr.
\end{multline*}
By Lemma~\ref{l.correctorminbounds}, Proposition~\ref{p.fluxcorrectorests} and the fact that $\Y_\sigma(y) \leq t^{\frac 12}$, we have that 
\begin{equation*} 
\left\| f  \right\|_{\underline{L}^2 \left( B_{rt^{1/2}}(y)\right) } \leq C \left( r t^{\frac 12}\right)^{1-\delta(d-\sigma)}.
\end{equation*}
Thus, by H\"older's inequality,
\begin{equation*} 
\left\| f \Phi\left(\frac{s}{\beta},\cdot-x \right)  \right\|_{\underline{L}^1\left( B_{rt^{1/2}}(y)\right)}  \leq C \left( r t^{\frac 12}\right)^{1-\delta(d-\sigma)} 
\left\| \Phi\left(\frac{s}{\beta},\cdot-x \right)  \right\|_{\underline{L}^2\left( B_{rt^{1/2}}(y)\right)},
\end{equation*}
and, by another application of H\"older's inequality and Step 1 of the proof of Lemma~\ref{l.PGF.layercake},
\begin{align} \notag 
\lefteqn{\int_0^\infty r^{d+1} \exp\left( - \beta \frac{r^2}{4} \right)  \left\| f \Phi\left(\frac{s}{\beta},\cdot-x \right)  \right\|_{\underline{L}^1\left( B_{rt^{1/2}}(y)\right) } \, dr } \quad &
\\ \notag &\leq
C t^{\frac 12(1-\delta(d-\sigma))}  \left( s^{-\frac d2}\int_0^\infty r^{d+1}  \exp\left( - 2\alpha \frac{r^2}{4} \right) \left\| \Phi\left(\frac{s}{2\alpha},\cdot-x \right)  \right\|_{\underline{L}^1\left( B_{rt^{1/2}}(y)\right)} \, dr \right)^{\frac12}
\\ \notag & = 
C t^{\frac 12(1-\delta(d-\sigma))} \left( \frac{t+s}{s} \right)^{\frac d4} \Phi\left(\frac{t+s}{\alpha},x-y \right).
\end{align}
Combining the above displays and redefining $\delta > 0$, we obtain~\eqref{e.SG.Pfluxes.2}. 

\smallskip

\emph{Step 2.} We now estimate the second term in~\eqref{e.SG.Pfluxes.2}, and show that
\begin{multline}  \label{e.SG.Pfluxes.3}
\left| \int_{\Rd}  \left( \a(z) - \ahom\right)  \left(\nabla_x P\left( t, z, y \right) -  \nabla_x H\left( t, z,y \right) \right) \overline{P}(s,x-z) \,dz  \right|  \\\leq 
 C s^{ -\delta(d-\sigma)} t^{-\frac12}  (t+s)^{-\frac d2} \exp\left( - \alpha \frac{|x-y|^2}{4(t+s)} \right).
\end{multline}
This, follows directly from Theorem~\ref{t.PGF.basecase}, because $s^{\frac 12} \geq \Y_\sigma(x) $.
\end{proof}

We are now ready to state and prove the lemma that allows to control the last term in \eqref{e.SG.basic} and therefore propagate the assumption of $\S(T,\msf K,\sigma)$ forward in time. The argument relies on Lemma~\ref{l.SG.Pfluxes} and CLT cancellations.

\begin{lemma} 
The exist an exponent $\delta(d,\Lambda) > 0$ and, for each $\sigma \in (0,2)$, a constant $C(\sigma,d,\Lambda)< \infty$ such that the following holds. Assume that for some $\msf K \in [1,\infty)$ and $t \in [1,\infty)$, we have for every $x \in \Rd$ that
\label{l.SG.opt4}
\begin{equation} \label{e.SG.opt4.a1}
\left| u\Ll(\tfrac t 2,x\Rr) \right| \leq \mathsf K t^{-\frac12} \wedge \O_{\sigma}\left( \mathsf K t^{-\frac12 - \frac d4} \right).
\end{equation}
Then, for every $s \in [t,\infty)$ and $x \in \Rd$, we have the estimate
\begin{equation} \label{e.SG.opt4.r}
 \int_{\R^d} \left(\ahom - \a(y) \right) \nabla u\left(t , y \right)  \cdot \nabla \overline{P}(s,x-y) \, dy  = \O_{\sigma}\left( C \mathsf K t^{-\delta(2-\sigma) -1} s^{-\frac 12 - \frac d4}  \right) . 
\end{equation}

\end{lemma}

\begin{proof}
We fix $\sigma \in (0,2)$ and decompose the proof into two steps.

\smallskip

\emph{Step 1.} We first show that there exist $\delta(d,\Lambda) > 0$ and $C(\sigma,d,\Lambda)<\infty$ such that,  for every $t \in [1,\infty)$, $s \in[t,\infty)$, and $x,x' \in \R^d$,
\begin{multline} 
\label{e.SG.opt4.1.r}
 \int_{x'+r\cu_{0}} \left(\ahom - \a(y) \right) \nabla u\left(t , y \right)  \cdot \nabla \overline{P}(s,x-y) \, dy 
\\=
\O_\sigma\left( C \mathsf K t^{-\delta(2-\sigma)- 1 - \frac d4} s^{-\frac 12} \left\| \Phi\left(\tfrac{t+s}{\alpha},\cdot-x \right) \right\|_{L^1(x'+r\cu_{0})}\right).
\end{multline}
To see this, we write 
\begin{multline*} 
 \int_{x'+r\cu_{0}} \left(\ahom - \a(y) \right) \nabla u\left(t , y \right)  \cdot \nabla \overline{P}(s,x-y) \, dy 
 \\= 
  \int_{x'+r\cu_{0}} u\left(\tfrac t2 , z \right) \left( \int_{\R^d} \left(\ahom - \a(y) \right) \nabla_x P\left(\tfrac t2,y,z\right) \cdot \nabla \overline{P}(s,x-y) \, dy \right) \, dz 
\end{multline*}
Set $\theta :=  \frac{d(\sigma+2)}{4}$ and $\alpha := \frac12 \Lambda^{-1}$. By Lemma~\ref{l.SG.Pfluxes}, there exist $\delta (d,\Lambda) > 0$ and, for every $z \in \R^d$,  a random variable $\Y_\theta(z)$  such that $\Y_\theta(z) \leq \O_\theta(C)$ and
\begin{multline*} 
\left| \int_{\R^d} \left(\ahom - \a(y) \right) \nabla_x P\left(\tfrac t2,y,z\right) \cdot \nabla \overline{P}(s,x-y) \, dy \right| \indc_{\left\{\Y_\theta(y) \vee \Y_\theta(x)  \leq t^{\frac12} \right\}}
\\ \leq
C  t^{-\delta(d-\sigma)- \frac12} s^{-\frac12}\Phi\left(\frac{t+s}{\alpha},z-x \right).
\end{multline*}
Hence the assumption of~\eqref{e.SG.opt4.a1} yields, via Lemma~\ref{l.sum-O}, that
\begin{multline*} 
\left| \int_{x'+r\cu_{0}} u\left(\tfrac t2 , z \right) \int_{\R^d} \left(\ahom {-} \a(y) \right) \nabla_x P\left(\tfrac t2,y,z\right) \cdot \nabla \overline{P}(s,x{-}y) \, dy \indc_{\left\{\Y_\theta(y) \vee \Y_\theta(x)  \leq t^{\frac12} \right\}} \, dz \right|  
\\\leq \O_\sigma\left( C \mathsf K t^{-\delta(d-\sigma)- 1 - \frac d4} s^{-\frac 12} \left\| \Phi\left(\tfrac{t+s}{\alpha},\cdot-x \right) \right\|_{L^1(x'+r\cu_{0})}\right).
\end{multline*}
On the other hand, by Lemma~\ref{l.PGF.layercake} we have the quenched bound
\begin{equation*} 
\left| \int_{\R^d} \left(\ahom - \a(y) \right) \nabla_x P\left(\tfrac t2,y,z\right) \cdot \nabla \overline{P}(s,x-y) \, dy \right|  \leq 
C t^{- \frac12} s^{-\frac12} \Phi\left(\frac{t+s}{\alpha},z-x \right).
\end{equation*}
Since
\begin{equation*} 
\indc_{\left\{\Y_\theta(y) \vee \Y_\theta(x)  > t^{\frac12} \right\}} \leq t^{-\frac d4 - \frac{d}{8\sigma}(2-\sigma)} \left(\Y_\theta(x)^{\frac\theta\sigma} + \Y_\theta(x)^{\frac\theta\sigma} \right),
\end{equation*}
Lemma~\ref{l.sum-O} and the bound $\left| u\Ll(\tfrac t2,z\Rr) \right| \leq \mathsf K t^{-\frac12}$ in~\eqref{e.SG.opt4.a1} imply, after possibly reducing the value of $\delta(d,\Lambda) > 0$,
\begin{multline*} 
\left| \int_{x'+r\cu_{0}} u\left(\tfrac t2 , z \right) \int_{\R^d} \left(\ahom {-} \a(y) \right) \nabla_x P\left(\tfrac t2,y,z\right) \cdot \nabla \overline{P}(s,x{-}y) \, dy \indc_{\left\{\Y_\theta(y) \vee \Y_\theta(x)  > t^{\frac12} \right\}} \, dz \right|  
\\\leq \O_\sigma\left( C \mathsf K t^{-\delta(2-\sigma)- 1 - \frac d4} s^{-\frac 12} \left\| \Phi\left(\tfrac{t+s}{\alpha},\cdot-x \right) \right\|_{L^1(x'+r\cu_{0})}\right).
\end{multline*}
Combining these estimates yields~\eqref{e.SG.opt4.1.r}. 

\smallskip

\emph{Step 2.} We now show~\eqref{e.SG.opt4.r} using the unit range of dependence assumption. 
Set $r := t^{\frac{1+\ep}{2}}$ with $\ep : = \frac{\delta(2-\sigma)}{d}$. Using the notation $u'(r,t,x)$ introduced in \eqref{e.SG.u'.def}, we define
\begin{equation*} 
\f(y) := \left(\ahom - \a(y) \right) \nabla u\left(t,y \right),  \qquad
\f'(r,y) := \left(\ahom - \a(y) \right) \nabla u'\left(r,t,y \right) ,
\end{equation*}
as well as
\begin{equation*} 
\tilde \f(y) :=  \f(t,y) - \E[\f(t,y)]  \quad \text{ and } \quad
\tilde \f'(r,y) := \f'(r,t,y) - \E[\f'(r,t,y)].
\end{equation*}
Recall that, by Lemma~\ref{l.SG.loc.u}, we have 
\begin{equation*} 
\left\| \f - \f'(r,\cdot)\right\|_{L^\infty(\R^d)} \leq C t^{-\frac 12} \exp \left( - C^{-1} \frac{r^2}{t} \right) = C t^{-\frac 12} \exp \left( - C^{-1} t^{\ep} \right) . 
\end{equation*}
Set also, for every $k \in \N$, 
\begin{equation*} 
\tilde \g_{k}'(r,y) := \tilde \f'(2^{k+1} r,t,y)  - \tilde \f'(2^{k}r,t,y) ,
\end{equation*}
and note using Lemma~\ref{l.SG.loc.u} again that
\begin{equation*} 
\left\| \tilde \g_{k}'(r,\cdot) \right\|_{L^\infty(\R^d)} \leq  C t^{-\frac 12} \exp \left( - C^{-1} 4^k t^{\ep} \right) .
\end{equation*}
We use the decomposition
\begin{multline} \label{e.SG.opt4.1}
\int_{\R^d} \tilde \f(y) \cdot  \nabla \overline{P}(s,x-y) \, dy  = \sum_{x' \in r\Z^d}\int_{x ' + r \cu_{0}} \tilde \f'(r,y) \cdot  \nabla \overline{P}(s,x-y) \, dy
\\ + \sum_{k=0}^\infty \sum_{x' \in 2^k r\Z^d}\int_{x ' + 2^k r \cu_{0}} \tilde \g_{k}'(r,y) \cdot \nabla \overline{P}(s,x-y) \, dy . 
\end{multline}
Let us now estimate the sizes of summands. First, we have that 
\begin{equation*} 
 \int_{x'+r\cu_{0}} \tilde \f'(r,y) \cdot \nabla \overline{P}(s,x-y) \, dy = 
 \O_\sigma\left( C \mathsf K t^{-\delta(2-\sigma)- 1 - \frac d4} s^{-\frac 12} \left\| \Phi\left(\tfrac{t+s}{\alpha},\cdot-x \right) \right\|_{L^1(x'+r\cu_{0})}\right) ,
\end{equation*}
which follows by~\eqref{e.SG.opt4.1.r} and
\begin{multline*} 
 \left| \int_{x'+r\cu_{0}} \left( \tilde \f(y) - \tilde \f'(r,y)\right) \cdot \nabla \overline{P}(s,x-y) \, dy \right| 
 \\ \leq C t^{-\frac 12} \exp \left( - C^{-1} t^{\ep} \right) s^{-\frac 12} r^{\frac d2} \left\| \Phi\left(\tfrac{t+s}{\alpha},\cdot-x \right) \right\|_{L^1(x'+r\cu_{0})}  . 
\end{multline*}
We also have that
\begin{multline*} 
 \left| \int_{x'+2^k r\cu_{0}} \tilde \g_{k}'(r,\cdot) \cdot \nabla \overline{P}(s,x-y) \, dy \right| 
 \\ \leq C 2^{-k(d+2)} t^{-\delta(2-\sigma)- 1 - \frac d4} s^{-\frac 12} \left\| \Phi\left(\tfrac{t+s}{\alpha},\cdot-x \right) \right\|_{L^1(x'+2^k r\cu_{0})}  . 
\end{multline*}
For each~$y \in \R^d$, the random variables $\tilde \f'(r,y)$ and $\g_{k}'(r,y) $ are $\F(B_r(y))$--measurable and $\F(B_{2^k r}(y))$--measurable, respectively. We thus obtain from Lemma~\ref{l.barO.boxes}, Lemma~\ref{l.sum-O} and the previous two displays the existence of~$C(\sigma,d,\Lambda)<\infty$ such that
\begin{align} \notag 
\lefteqn{\left| \sum_{z' \in r \Z^d} \int_{z' + r\cu_0} \tilde \f'(r,y) \cdot  \nabla \overline{P}(s,x-y) \, dy \right| } \qquad &
\\ \notag & \leq \O_{\sigma}\left( C \mathsf K t^{-\delta(2-\sigma)- 1 - \frac d4} s^{-\frac 12}  \left(\sum_{z' \in r \Z^d} \left\| \Phi\left( \frac{s}{\alpha} , \cdot - x\right) \right\|_{L^1(z' + r\cu_0)}^2 \right)^{\frac 12} \right)
\\ \notag & = \O_{\sigma}\left( C \mathsf K t^{-\delta(2-\sigma)- 1 - \frac d4} s^{-\frac 12}  r^{\frac d2} \left\|\Phi\left( \frac{s}{\alpha} , \cdot \right) \right\|_{\underline{L}^2(\R^d)} \right)
\\ \notag & = \O_{\sigma}\left( C\mathsf K t^{-\frac{\delta(2-\sigma)}{2} - 1 } s^{-\frac12 - \frac d4} \right).
\end{align} 
Similarly, 
\begin{equation*} 
\sum_{k=0}^\infty \left| \sum_{z' \in 2^k r \Z^d} \int_{z' + 2^k r\cu_0} \tilde v_k'(r,y) \partial_{t} \overline{P}(s,x-y) \, dy \right| \leq  \O_{\sigma}\left( C\mathsf K t^{-\frac{\delta(2-\sigma)}{2} - 1 } s^{-\frac12 - \frac d4} \right).
\end{equation*}
Combining previous two displays with~\eqref{e.SG.opt4.1} and redefining $\delta$ completes the argument. 
\end{proof}

We are now ready to complete the proof of Theorem~\ref{t.semigroup}.

\begin{proof}[Proof of Theorem~\ref{t.semigroup}]
We decompose the proof into two steps.

\smallskip

\emph{Step 1.} In the first step, we show that there exist an exponent $\eps(\sigma,d,\Lambda) > 0$ and a constant $C(\sigma,d,\Lambda) < \infty$ such that for every $T, \mathsf K \ge 1$,
\begin{equation}
\label{e.impl.opt}
\S(T,\mathsf K, \sigma) \implies \S\left(2T,\mathsf K(1+C T^{-\ep}), \sigma\right).
\end{equation}
In order to do so, we will use the decomposition provided by Lemma~\ref{l.SG.basic}. To start, by Lemma~\ref{l.SG.pointwise.u}, there exists a constant $C(\sigma,d,\Lambda) < \infty$ such that, under the assumption of $\S(T, \mathsf K,\sigma)$, we have, for every $t \in [1,2T]$ and $x \in \Rd$, that
\begin{equation}
\label{e.pointwise.utx}
u(t,x) = \O_\sigma \Ll( C \mathsf K t^{-\frac 1 2 - \frac d 4} \Rr) .
\end{equation}
Using also Lemma~\ref{l.crude.ubound}, we obtain, for every $t \in [1,2T]$,
\begin{equation}
\label{e.pointwise.utx2}
\left|u(t,x)\right| \leq C t^{- \frac 12} \wedge \O_\sigma \Ll( C \mathsf K t^{-\frac 1 2 - \frac d 4} \Rr) .
\end{equation}
For each $t \in \left(\tfrac12 T,T\right]$ and $s \ge t$, we use Lemma~\ref{l.SG.basic} in the form of
\begin{multline*} 
\int_{\R^d} u(2t,y) \overline{P}(s,x-y) \, dy 
 = \int_{\R^d} u(t,y) \overline{P}(s+t,x-y) \, dy
\\ 
+ \int_{t}^{2t} \int_{\R^d} \left(\ahom - \a(y) \right) \nabla u\left(t' ,y \right) \cdot \nabla \overline{P}(s+t-t',x-y) \, dy  \, dt'.
\end{multline*}
In view of \eqref{e.pointwise.utx2}, we can apply Lemma~\ref{l.SG.opt4} and get that
\begin{equation*}
 \int_{t}^{2t} \int_{\R^d} \left(\ahom - \a(y) \right) \nabla u\left(t' , y \right)  \cdot \nabla \overline{P}(s+t-t',x-y) \, dy \, dt'  = \O_{\sigma}\left( C \mathsf K t^{-\delta(2-\sigma)} s^{-\frac 12 - \frac d4}  \right) .
\end{equation*}
By our assumption of $\S(T,\mathsf K,  \sigma)$, we also have
\begin{equation*}  
\int_\Rd u(t,y) \bar P(s+t,x-y) \, dy = \O_\sigma \Ll( \msf K s^{-\frac 1 2 - \frac d 4} \Rr) ,
\end{equation*}
and thus
\begin{equation*} 
 \int_{\R^d} u(2t,y) \overline{P}(s,x-y) \, dy  \leq  \O_{\sigma}\left( \left(1 + CT^{-\delta(2-\sigma)}  \right) \mathsf K s^{-\frac 12 - \frac d4}  \right) .
\end{equation*}
This proves \eqref{e.impl.opt} with $\ep :=   \delta(2-\sigma) > 0$. 

\smallskip

\emph{Step 2.} 
We conclude the proof. By Lemma~\ref{l.SG.opt1}, there exists a constant $\msf K_0(d,\Lambda) < \infty$ such that $\S(1,\msf K_0,2)$ holds. Since $\prod_{k = 0}^\infty (1 + C 2^{-k\ep}) < \infty$, 
an induction on \eqref{e.impl.opt} implies the existence of a constant $\msf K(\sigma,d,\Lambda) < \infty$ such that $\S(T,\msf K,\sigma)$ holds for every $T \ge 1$. The conclusion then follows by Lemma~\ref{l.SG.pointwise.u}.
\end{proof}

\section{Homogenization of the Green functions: optimal scaling}
\label{s.PGFopt}

In this section, we prove optimal quantitative estimates on the difference between the heterogeneous parabolic Green function~$P(t,x,y)$ and its homogenized counterpart~$\bar P(t,x-y)$, which are stated in the following theorem.


\smallskip

\index{Berry-Esseen theorem}

\begin{theorem}[Optimal convergence of parabolic Green function]
\label{t.PtoPbar}
\index{Berry-Esseen theorem}
Fix~$\sigma\in (0,2)$ and $\alpha \in (0,\Lambda^{-1})$. There exists~$C(\sigma,\alpha,d,\Lambda)<\infty$ such that, for every~$t\in [1,\infty)$ and~$x,y\in\Rd$, 
\begin{equation} 
\label{e.PtoPbar}
\left| P(t,x,y) - \overline{P}(t,x-y) \right| 
\leq
\left\{
\begin{aligned}
& \O_\sigma\left( C t^{-\frac12 -\frac d2} \exp\left( - \alpha \frac{|x-y|^2}{4t} \right)\right) & \mbox{if} & \ d>2, \\
&  \O_\sigma\left( C \log^{\frac12}(1+t) t^{-\frac12-\frac d2}  \exp\left( - \alpha \frac{|x-y|^2}{4t} \right)\right) & \mbox{if} & \ d=2.
\end{aligned}
\right. 
\end{equation}
\end{theorem}

In view of~\eqref{e.intGammatoG}, Theorem~\ref{t.PtoPbar} implies the following quantitative homogenization result for the elliptic Green function~$G$. 

\begin{corollary}
[Optimal convergence of elliptic Green function]
\label{c.GtoGbar}
For each~$\sigma\in (0,2)$, there exists a constant~$C(\sigma,d,\Lambda)<\infty$ such that, for every~$t\in [1,\infty)$ and~$x,y\in\Rd$ with $|x-y|\geq 1$, 
\begin{equation} 
\label{e.GtoGbar}
\left| G(x,y) - \overline{G}(x-y) \right| 
\leq
\left\{
\begin{aligned}
&\O_\sigma\left( C |x-y|^{1-d} \right) & \mbox{if} & \ d>2, \\
&  \O_\sigma\left( 
C |x-y|^{-1} \log^{\frac12}(2+|x-y|)  \right) & \mbox{if} & \ d=2.
\end{aligned}
\right. 
\end{equation}
\end{corollary}

The proof of Theorem~\ref{t.PtoPbar} is more subtle than the one of Theorem~\ref{t.PGF.basecase}, due to the difficulty in obtaining bounds at the optimal scaling. While it is relatively straightforward to bootstrap the arguments from the previous section to obtain a version of~\eqref{e.rough.PGFts} with $t^{-\frac12+\delta}$ in place of the factor of~$t^{-\frac12}$, getting \emph{on top of the exponent} requires us to ``go to next order'' in our asymptotic expansion of~$P$ (or, at least we must expand to next order the term causing the difficulty in the bootstrap argument). This should not hide the fact that the argument is still a bootstrapping of the suboptimal result of Theorem~\ref{t.PGF.basecase}, since it uses the latter in a crucial way. However, as discussed below, the formal argument will be an induction in the time variable rather than a bootstrap in the exponent.

\smallskip

Estimates for the error in the two-scale expansion for both~$\nabla_xP$ and~$\nabla_x\nabla_yP$ which are ``nearly optimal'' can be obtained from Theorem~\ref{t.PtoPbar} by using an argument similar to the one in Step~2 of the proof of Theorem~\ref{t.PGF.basecase} in the previous section. This argument however leads to a small loss of the exponent in the estimate. The statement one would obtain for~$\nabla_xP$ for instance is that, for every $\alpha\in (0,\Lambda^{-1})$ and $\beta>0$, there exists a constant $C(\beta,\alpha,d,\Lambda)<\infty$ such that  
\begin{equation}
\label{e.PGF.basecase.grad.nearlyopt}
 \left\| \nabla_x P(t,\cdot,y) - \nabla_x H(t,\cdot,y) \right\|_{\underline{L}^2(B_r(x))}
\leq 
\O_2\left( 
Ct^{-\frac12 + \beta} t^{ - \frac 12 - \frac d2}  \exp\left(- \alpha \frac{|x-y|^2}{4t}\right)
\right).
\end{equation}
In fact, this small loss of exponent $\beta>0$ can be removed, but this requires a more subtle argument that is beyond the scope of this book. 

\smallskip

\index{intermediate asymptotics|(}

In order to motivate the analysis underlying the proof of Theorem~\ref{t.PtoPbar}, let us briefly discuss some classical facts concerning \emph{intermediate asymptotics} of solutions to the heat equation. Consider the initial value problem for the heat equation
\begin{equation} 
\label{e.heatflowexample}
\left\{
\begin{aligned}
&\partial_t u -\Delta u = 0  & \mbox{in} & \ (0,\infty) \times \Rd,\\
& u(0,\cdot) = u_0 & \mbox{on} & \ \Rd, 
\end{aligned}
\right.
\end{equation}
for a smooth and compactly supported initial condition $u_0\in C^\infty_c(\Rd)$. It is well-known that, after a long time, the solution~$u$ will be close to a multiple of the heat kernel $\Phi$, namely $m_0\Phi$, where $m = \int_{\Rd} u_0$ is the mass of the initial condition. For instance, in terms of the $L^1(\Rd)$ norm, we have that 
\begin{equation*} \label{}
\lim_{t\to \infty} \left\| u(t,\cdot) - m\Phi(t,\cdot) \right\|_{L^1(\Rd)} = 0. 
\end{equation*}
One may wonder about the rate of convergence in this limit. It turns out to be $O(t^{-\frac12})$, that is, for a constant $C$ depending on~$d$ and the support of $u_0$, 
\begin{equation} 
\label{e.decayfaster}
\left\| u(t,\cdot) - m\Phi(t,\cdot) \right\|_{L^1(\Rd)} \leq Ct^{-\frac12}. 
\end{equation}
In particular, if $u_0$ has zero mass, then $u(t,\cdot)$ decays to zero by a factor of~$t^{-\frac12}$ in $L^1(\Rd)$ and in every other $L^p$ norm by~$t^{-\frac12}$ faster than a solution with nonzero mass. This estimate can be proved directly from the representation formula 
\begin{equation*} \label{}
u(t,x) = \int_{\Rd} u_0(y) \Phi(t,x-y)\,dy
\end{equation*}
and an integration by parts. It is actually sharp, as we can see immediately from considering for instance~$u(t,x) = \partial_{x_i} \Phi(t+1,x)$ for any $i\in\{1,\ldots,d\}$, which has zero mass and has an~$L^1$ norm which decays at exactly this rate. 

\smallskip

Moreover, this suggests a more precise asymptotic expansion for $u_0$ of the form
\begin{equation} 
\label{e.expandtonextorder}
u(t,x) = m \Phi(t,x) + \mathbf{m}\cdot \nabla \Phi(t,x) + o\left( t^{-\frac12} \right) \quad \mbox{as} \ t \to \infty, 
\end{equation}
where the vector $\mathbf{m}$ should be the \emph{center of mass} (or \emph{first moment}) of $u_0$:
\begin{equation*} \label{}
\mathbf{m}:= \int_\Rd x u_0(x) \,dx. 
\end{equation*}
As it turns out, this is indeed correct, and one can show that 
\begin{equation*} \label{}
\left\| u(t,\cdot) - \left(  m \Phi(t,\cdot) + \mathbf{m}\cdot \nabla \Phi(t,\cdot ) \right) \right\|_{L^1(\Rd)}  \leq Ct^{-1}. 
\end{equation*}
In fact, a full asymptotic expansion (to any power of $t^{-\frac12}$) is possible in terms of the moments of $u$ and higher derivatives of the heat kernel, and we refer the reader to~\cite{DZ} for the precise statement and proof of this result.  

\index{intermediate asymptotics|)}
\smallskip

To understand the relevance of this discussion to our context, namely that of proving Theorem~\ref{t.PtoPbar}, notice that the estimate~\eqref{e.PtoPbar} says that the difference between~$P$ and~$\overline{P}$ has relative size of order $t^{-\frac12}$, the same as the error of the leading-order intermediate asymptotic expansion above. The proof strategy for Theorem~\ref{t.PtoPbar} relies on a propagation argument, in other words we assume the validity of the estimate up to time $t$ and try to establish it at time $2t$. To do this, we will study two sources of error: the error that has already been accumulated up to time $t$, which will be pushed forward it time, and the ``new'' error created between time $t$ and $2t$ due to the fact that $P$ and $\overline{P}$ solve different equations (the homogenization error). The problem we face is that even if the latter error were \emph{zero} (no homogenization error), we can see from the discussion above that the error that is propagated forward in time show decay no faster than the heat kernel times~$Ct^{-\frac12}$. This may at first glance appear to be just what we need, but the accumulation of many such errors across all scales will lead to a logarithmic correction (i.e.,~by this argument one would be able to show only~\eqref{e.PtoPbar} with an additional factor of $\log t$ on the right side). 

\smallskip

Therefore the proof strategy must be slightly smarter, and so we keep track of the first moment of~$P$ as we propagate the estimate. This leads us to define
\begin{equation} 
\label{e.def.Qbar}
\overline{Q}(t,x,y)
:=
\overline{P}(t,x-y) - \mathbf{m} (t,y) \cdot \nabla \overline{P} (t,x-y),
\end{equation}
where we denote the \emph{center of mass} (or \emph{first moment}) of~$P(t,\cdot,y)$ relative to~$y$ by
\begin{equation} 
\label{e.def.m}
\mathbf{m}(t,y)
:= 
\int_{\Rd} (x-y) P(t,x,y)\,dx. 
\end{equation}
Note that $(t,x)\mapsto \overline{Q}(t,x,y)$ is not an exact solution of the homogenized equation, since $\mathbf{m}$ depends on~$t$. We also stress that it is also not a deterministic function (despite the use of the ``bar'' notation) since $\mathbf{m}(t,y)$ is random. However, as we will see,~$Q(\cdot,y)$ is \emph{almost} a solution of the $\ahom$-heat equation. We record some basic estimates on~$\mathbf{m}(t,y)$ in Lemma~\ref{l.trash.m} below.

\smallskip

We denote the two-scale expansion of $\overline{Q}$ by 
\begin{equation} 
\label{e.def.K}
K(t,x,y)
:=
\overline{Q}(t,x,y) 
+
\sum_{k=1}^d  \left( \phi_{e_k} (x) -  \left( \phi_{e_k} \ast \Phi(t,\cdot) \right) (y) \right) \partial_{x_k} \overline{Q}(t,x,y). 
\end{equation}
The proof of Theorem~\ref{t.PtoPbar} proceeds by proving bounds on the difference~$P-K$, which we expect to be easier to analyze since, informed by~\eqref{e.expandtonextorder}, we can expect the error propagated forward in time to decay faster than $t^{-\frac12}$. Indeed, the motivation for using~$\overline{Q}$ in place of~$\overline{P}$ in the definition is that it has the same center of mass as~$P$ (as well as the same mass). By the definitions of $\mathbf{m}(t,y)$ and $\overline{Q}(t,x,y)$ and the identity
\begin{equation*} \label{}
\int_{\Rd} -x \partial_{x_i} \overline{P}(t,x)\,dx = \int_{\Rd} e_i \overline{P}(t,x)\,dx = e_i,
\end{equation*}
we have that 
\begin{equation} 
\label{e.firstmoment}
\int_{\Rd} (x-y) \overline{Q}(t,x,y) \,dx = \mathbf{m}(t,y) = \int_{\Rd} (x-y)P(t,x,y)\,dx. 
\end{equation}
It is also easy to check that 
\begin{equation} 
\label{e.checkmass}
\int_{\Rd} \overline{Q}(t,x,y) \,dx = 1 = \int_{\Rd} P(t,x,y)\,dx. 
\end{equation}

\smallskip

In the following lemmas, we check that the difference between~$K$ and~$\overline{P}$ is small enough that appropriate estimates on $P-K$ are sufficient to imply Theorem~\ref{t.PtoPbar}. We do this by showing that $\overline{Q} - \overline{P}$ is sufficiently small, which reduces to showing that stochastic moments of~$\mathbf{m}(t,y)$ are bounded independently of $t$, and then showing that $K - \overline{Q}$ is sufficiently small as well, which is an exercise in applying some familiar bounds on the correctors to estimate the second term on the right of~\eqref{e.def.K}.

\smallskip

Throughout the rest of this section, for each $\sigma\in (0,d)$ and $y\in\Rd$, we let $\Y_\sigma(y)$ denote the random variable in the statement of Theorem~\ref{t.PGF.basecase} which satisfies, for some constant $C(\sigma,\Lambda,d)<\infty$,
\begin{equation} 
\label{e.Ysig.again}
\Y_\sigma(y) \leq \O_\sigma(C).
\end{equation}

\begin{lemma}
\label{l.trash.m}
There exists a constant $C(d,\Lambda)<\infty$ such that, for $t>0$ and $y\in\Rd$, 
\begin{equation}
\label{e.m.basecase}
\left|\mathbf{m}(t,y)\right|  + t \left| \partial_t \mathbf{m}(t,y) \right| \leq C t^{\frac 12} . 
\end{equation}
Moreover, for every $\sigma \in (0,2)$, there exist $\delta(d,\Lambda)>0$ and $C(\sigma,d,\Lambda)<\infty$ such that, for every $t \geq 1$ and $y\in\Rd$, 
\begin{equation} \label{e.m.basecase.O1}
 \left| \partial_t \mathbf{m}(t,y)\right| \leq \O_\sigma\left(C t^{- \frac 12 -\frac d4}\right) 
\end{equation}
and
\begin{equation} \label{e.m.basecase.O}
 \left|  \mathbf{m}(t,y)\right|\leq \left\{
\begin{aligned}
& \O_\sigma\left( C \right) & \mbox{if} & \ d>2, \\
&  \O_\sigma\left( C \log^{\frac 12}(1+t)  \right) & \mbox{if} & \ d=2.
\end{aligned}
\right. 
\end{equation}
\end{lemma}

\begin{lemma}
\label{l.trashK}
For every $\sigma\in (0,2)$ and $\alpha \in (0,\Lambda^{-1})$, there exists $C(\sigma,\alpha,d,\Lambda)<\infty$ such that, for every $x,y \in\R^d$ and $t \in [1,\infty)$, 
\begin{equation*} 
\left| K(t,x,y) - \overline{Q}(t,x,y) \right| 
\leq 
\left\{
\begin{aligned}
& \O_\sigma\left( C t^{-\frac12 -\frac d2} \exp\left( -\alpha \frac{|x-y|^2}{4t} \right)\right) & \mbox{if} & \ d>2, \\
&  \O_\sigma\left( C \log^{\frac12}(1+t)  t^{-\frac12-\frac d2}  \exp\left( -\alpha  \frac{|x-y|^2}{4t} \right)\right) & \mbox{if} & \ d=2.
\end{aligned}
\right. 
\end{equation*}
\end{lemma}

Before we give the proofs of Lemmas~\ref{l.trash.m} and~\ref{l.trashK}, notice that they imply that 
\begin{multline} 
\label{e.QbartoPbar}
\left| K(t,x,y) - \overline{P}(t,x-y) \right| 
\\\leq 
\left\{
\begin{aligned}
& \O_\sigma\left( C t^{-\frac12 -\frac d2} \exp\left( -\alpha \frac{|x-y|^2}{4t} \right)\right) & \mbox{if} & \ d>2, \\
&  \O_\sigma\left( C  \log^{\frac 12}(1+t) t^{-\frac12-\frac d2} \exp\left( -\alpha \frac{|x-y|^2}{4t} \right)\right) & \mbox{if} & \ d=2.
\end{aligned}
\right. 
\end{multline}
Indeed, the bound~\eqref{e.m.basecase.O} implies that, for every $t > 0 $ and $x,y\in\Rd$, 
\begin{align*} \label{}
\left|  \overline{P}(t,x-y) - \overline{Q}(t,x-y) \right| 
&
\leq \left| \mathbf{m}(t,y) \right| \left| \nabla\overline{P}(t,x-y) \right|
\\ &  
\leq \left\{
\begin{aligned}
& \O_\sigma\left( C t^{-\frac12 -\frac d2} \exp\left( -\alpha \frac{|x-y|^2}{4t}\right)\right) & \mbox{if} & \ d>2, \\
&  \O_\sigma\left( C  \log^{\frac 12}(1+t) t^{-\frac12-\frac d2}  \exp\left( -\alpha \frac{|x-y|^2}{4t}\right)\right) & \mbox{if} & \ d=2.
\end{aligned}
\right. 
\end{align*}
Combining this bound with the bound in Lemma~\ref{l.trashK}, we get~\eqref{e.QbartoPbar}. The upshot of~\eqref{e.QbartoPbar} is that we can prove~\eqref{e.PtoPbar} by proving a similar estimate on $P-K$. 

\smallskip

\begin{proof}[{Proof of Lemma~\ref{l.trash.m}}]

We first show~\eqref{e.m.basecase}. 
By the Nash-Aronson bound,
\begin{equation*}
\left|\mathbf{m}(t,y)\right|
\leq 
C \int_{\R^d} t^{  -\frac d2}|x-y| \exp \left( - \alpha \frac{|x-y|^2}{4t}\right) \, dx \leq C t^{\frac 12 }.
\end{equation*}
For the time derivative, observe that the equation for $P$ gives us
\begin{align}
\label{e.m_t} 
\partial_t \mathbf{m}(t,y)
=   \int_{\Rd} (x-y)\partial_t P(t,x,y)\,dx
 = - \int_{\Rd} \a(x) \nabla_x P(t,x,y)\,dx,
\end{align}
and thus \eqref{e.m.basecase} follows since 
\begin{equation*} 
\left\|\nabla_x P(t,\cdot,y) \right\|_{L^1(\R^d)} \leq C t^{-\frac12}. 
\end{equation*}

\smallskip

We turn to the proof of~\eqref{e.m.basecase.O1} and~\eqref{e.m.basecase.O}. For every $e \in \partial B_{1}$ we recognize from~\eqref{e.m_t} that~$u_e:=\partial_t \mathbf{m} \cdot e$ is the solution of
\begin{equation}
\label{e.flow.FOC.again}
\left\{
\begin{aligned}
& \partial_t u_e -\nabla \cdot \left( \a \nabla u_e \right) = 0 &  \mbox{in} & \ (0,\infty) \times \Rd, \\
& u(0,\cdot) = \nabla\cdot \left( \a(\cdot) e \right) & \mbox{on} & \ \Rd.
\end{aligned}
\right.
\end{equation}
Hence, for each $\sigma \in (0,2)$, we have by Theorem~\ref{t.semigroup} that for every $t \geq 1$ and $y \in \R^d$, 
\begin{equation*} 
\left| \partial_t \mathbf{m}(t,y) \right| = \left| \int_{\Rd} \a(x) \nabla_x P(t,x,y) \, dx \right| \leq \O_{\sigma}(C t^{-\frac12 - \frac d4}),
\end{equation*}
giving~\eqref{e.m.basecase.O1}. This also yields by Lemma~\ref{l.sum-O} that if $d>2$, then
\begin{equation*} 
\left| \mathbf{m}(t,y) \right|  \leq \left| \int_1^t \partial_t \mathbf{m}(t',y) \,dt' \right| + \left| \mathbf{m}(1,y) \right| \leq  \O_\sigma\left( C \right).
\end{equation*}
In $d=2$ we instead use an argument from Section~\ref{s.twopoint} using the bound~\eqref{e.phid=2}. Using Lemmas~\ref{l.SG.opt5} and~\ref{l.sum-O} repeatedly, together with  the semigroup property,  we obtain that there exists $\delta(d,\Lambda)>0$ such that, for every $x \in \Rd$, $s \in [0,\infty)$ and $t \in [1,\infty)$, 
\begin{align} 
\notag  
\lefteqn{\int_{\Rd}  u_e(t,y) \overline{P}(s,x-y) \, dy } \quad  &  
\\ 
\notag &= \int_{\R^d} \int_{\R^d} u_e\left(\tfrac t2,z\right) P\left(\tfrac t2,y,z\right) \, dz  \, \overline{P}(s,x-y) \, dy
\\ 
\notag & = \int_{\R^d} \int_{\R^d} u_e\left(\tfrac t2,z\right) \overline{P}\left(\tfrac t2,y-z\right) \, dz \,  \overline{P}(s,x-y) \, dy + \O_\sigma\left( C t^{-\delta(2-\sigma) -1} \right)
\\ 
\notag & = \int_{\R^d}  u_e\left(\tfrac t2,z\right) \overline{P}\left(s+\tfrac t2,x-z\right) \, dz + \O_\sigma\left( C t^{-\delta(2-\sigma) -1} \right)
\\ 
\notag & = \int_{\R^d}  u_e\left(\tfrac t2,z\right) P\left(s+\tfrac t2,x,z\right) \, dz + \O_\sigma\left( C t^{-\delta(2-\sigma) -1} \right) + \O_\sigma\left( C (s+t)^{-\delta(2-\sigma)-1} \right) 
\\ 
\notag & = u_e\left(t+s,x\right) +  \O_\sigma\left( C t^{-\delta(2-\sigma) -1} \right) .
\end{align}
Therefore, by Theorem~\ref{t.semigroup} and Lemma~\ref{l.sum-O}, we have, for every $x \in \Rd$, $s,s' \in (0,\infty)$ and $T>s'\vee s$, 
\begin{align*}
\lefteqn{
\left( \int_1^T u_e(t,\cdot) \ast \left( \overline{P}(s',\cdot) -  \overline{P}(s,\cdot) \right) \right) (x) \,dt 
} \qquad & 
\\ & 
= \int_1^T \left( u_e(t+s',x) - u_e(t+s,x)\right) \, dt + \O_\sigma(C)
\\ & 
= 
\int_{1+s'}^{1+s} u_e(t,x) \, dt +  \O_\sigma\Ll(C(1+ (s'-s) T^{-1})\Rr).
\end{align*}
We also have by~\eqref{e.SG.u.decay} that 
\begin{equation*} 
\left( \int_0^1 u_e(t,\cdot) \ast \left( \overline{P}(s',\cdot) -  \overline{P}(s,\cdot) \right) \right) (x) \,dt \leq C.
\end{equation*}
By sending $T \to \infty$, we obtain that, for every $x \in \Rd$, and $s,s' \in [1,\infty)$,
\begin{equation*} 
\left( \phi_e\ast \overline{P}(s',\cdot) - \phi_e\ast \overline{P}(s,\cdot) \right) (x)  = \int_{s'}^{s} u_e(t,x) \, dt + \O_\sigma(C). 
\end{equation*}
Consequently,~\eqref{e.phid=2} and Remark~\ref{r.Phidoesntmatter} yield that, for every $x \in \Rd$ and $s\ge s' \ge 1$,
\begin{equation*} 
\int_{s'}^{s} u_e(t,x) \, dt  = \O_\sigma\left( C\log^{\frac12} \left(2+\frac{s}{s'}\right) \right).
\end{equation*}
As~$u_{e}(t,x) = \partial_t \mathbf{m}(t,x) \cdot e$, we obtain~\eqref{e.m.basecase.O} in the case~$d=2$.
\end{proof}

Observe that, as a consequence of~\eqref{e.m.basecase}, we obtain, for each~$\alpha \in (0,\Lambda^{-1})$, the existence of a constant~$C(\alpha,d,\Lambda)<\infty$ such that, for every $t>0$ and $x,y \in \Rd$, 
\begin{multline} \label{e.Qbar.bounded}
t \left| \partial_t \overline{Q}(t,x,y) \right| + t \left| \nabla_{x}^2 \overline{Q}(t,x,y) \right| + \sqrt{t} \left| \nabla_{x} \overline{Q}(t,x,y) \right| + \left|\overline{Q}(t,x,y) \right| \\ \leq Ct^{-\frac d2} \exp\left( - \alpha \frac{|x-y|^2}{4t} \right).
\end{multline}

\begin{proof}[{Proof of Lemma~\ref{l.trashK}}]
We have that 
\begin{align} 
\notag  
K(t,x,y) - \overline{Q}(t,x,y) & =   \sum_{k=1}^d  \left( \phi_{e_k} (x) -  \left( \phi_{e_k} \ast \Phi(t,\cdot) \right) (y) \right) \partial_{x_k} \overline{Q}(t,x,y).
\end{align}
In dimensions~$d>2$, Theorem~\ref{t.correctors} and Remark~\ref{r.scalar.Linfty} imply that, for some~$\delta(d,\Lambda)>0$,
\begin{equation*} 
\sum_{k=1}^d  \left|  \phi_{e_k} (x) -  \left( \phi_{e_k} \ast \Phi(t,\cdot) \right) (y) \right| \leq \O_{2+\delta}(C).
\end{equation*}
On the other hand, in dimension~$d=2$, Theorem~\ref{t.correctors} and Remark~\ref{r.scalar.Linfty} yield
\begin{equation*} 
\sum_{k=1}^d  \left|  \phi_{e_k} (x) -  \left( \phi_{e_k} \ast \Phi(t,\cdot) \right) (y) \right| \leq \O_{\sigma}\left(C \log\left(2 + \sqrt{t} + |x-y|\right)^{\frac12}\right).
\end{equation*}
Now the Nash-Aronson bound for $\overline{P}$ and the bound~\eqref{e.Qbar.bounded} for $\partial_{x_k} \overline{Q}(t,x,y)$ yield the result. 
\end{proof}

\smallskip

As discussed above, we prove the desired decay estimate for~$P-K$ by propagating it forward in time. The precise statement we propagate is given in the following definition.
\begin{definition} \label{d.PGF.}
Given constants~$T \in [1,\infty]$, $\mathsf{C}\in [1,\infty)$, $\sigma\in (0,2)$ and $\alpha \in (0,\Lambda^{-1})$, we let~$\mathcal{S}\left( T, \mathsf{C},\sigma,\alpha \right)$ denote the statement that, for every $t\in [1,T] \cap (1,\infty)$ and $x,y\in\Rd$, 
\begin{equation} 
\label{e.veryambitious}
\left| (P-K)(t,x,y) \right|
\leq
\left\{
\begin{aligned}
& \O_\sigma\left( \mathsf{C} t^{-\frac12 -\frac d2} \exp\left( -\alpha \frac{|x-y|^2}{4t} \right)\right) & \mbox{if} & \ d>2, \\
&  \O_\sigma\left( \mathsf{C} \log^{\frac12}(1+t)  t^{-\frac12-\frac d2}  \exp\left( -\alpha \frac{|x-y|^2}{4t} \right)\right) & \mbox{if} & \ d=2.
\end{aligned}
\right. 
\end{equation}
\end{definition}
The main step in the proof of Theorem~\ref{t.PtoPbar} is to show that, for every $\sigma\in (0,2)$ and $\alpha\in \left(0,\Lambda^{-1} \right)$, there exists $C(\sigma, \alpha,d,\Lambda)<\infty$ such that, for every $T,\mathsf{C}\geq C$,
\begin{equation}
\label{e.propagate.T}
\mathcal{S}\left( T, \mathsf{C},\sigma ,\alpha \right)
\implies
\mathcal{S}\left( 2T, \mathsf{C},\sigma , \alpha \right).
\end{equation}
The proof of~\eqref{e.propagate.T} is based on the Duhamel formula. For every~$s,t\in (0,\infty)$ and $x,y \in \R^d$, we write
\begin{align} 
\label{e.Duhamel.P-K}
\lefteqn{ 
(P-K)(t+s,x,y) 
} \quad &
\notag \\ &
=
\int_{\Rd} (P-K) (t,z,y) P(s,x,z) \, dz
\notag \\ & \quad 
-
\int_{t}^{t+s} 
\int_{\Rd}
\left( \partial_{t'} -\nabla_z\cdot \a \nabla_z  \right) \hat K\left( t',z\right)  P(t+s-t',x,z)
\, dz \, dt'. 
\end{align}
Here, to enforce the validity of \eqref{e.Duhamel.P-K}, we could choose $\hat K (\cdot,\cdot) = K(\cdot, \cdot,y)$. However, in order to have \eqref{e.Duhamel.P-K} for fixed $s,t > 0$ and $x,y \in \Rd$, it suffices that
\begin{equation*}  
\hat K(t+s,x) = K(t+s,x,y) \quad \text{and} \quad \hat K(t,\cdot) = K(t,\cdot,y).
\end{equation*}
This simple observation will give us more freedom and allow us to control the interaction between the singularity created by $\nabla P(t+s-t',x,z)$, when $t'$ tends to $t+s$, and correctors and, especially, flux correctors present when computing $\left( \partial_{t'} -\nabla_z\cdot \a \nabla_z  \right) K\left( t',z,y\right)$; see Lemma~\ref{l.PGF.secondterm} below.  

\smallskip

The first integral on the right side of~\eqref{e.Duhamel.P-K} represents the previous error, accumulated by time~$t$, propagated to time~$t+s$, while the second term is due to the error created between times~$t$ and $t+s$ (the homogenization error). We first concentrate on the former, which we rewrite as
\begin{align} 
\label{e.Duhamel.P-K.first}
\notag  
\lefteqn{
\int_{\Rd}  (P-K)(t,z,y)  P(s,x,z)  \,dz 
} \quad & \\
& 
=   \int_{\Rd}  (P-K)(t,z,y)  \left(P(s,x,z)- \overline{P}(s,x-z) \right) \,dz
\notag \\ &  \quad 
+   \int_{\Rd}  (P-\overline{Q})(t,z,y)  \overline{P}(s,x-z)  \,dz
\notag \\ &  \quad 
+   \int_{\Rd}  (\overline{Q}-K)(t,z,y)  \overline{P}(s,x-z)  \,dz.
\end{align}
The last integral on the right side of~\eqref{e.Duhamel.P-K.first} has been already treated in Lemma~\ref{l.trashK}. We will estimate the first two integrals on the right side of~\eqref{e.Duhamel.P-K.first} in Lemmas~\ref{l.PGF.(P-K)vs(P-barP)} and~\ref{l.PGF.PvsbarQ}, respectively. The motivation for using~$\overline{Q}$ in place of $\overline{P}$ is precisely because it allows us to estimate the second integral---this is where we will have to essentially prove a version of~\eqref{e.expandtonextorder}.

\begin{lemma} 
Fix $\sigma\in (0,2)$ and $\alpha \in \left( 0,\Lambda^{-1} \right)$. There exists~$C(\sigma,\alpha,d,\Lambda) < \infty$ such that for every $\msf C,T \in [1,\infty)$, if the statement~$\mathcal{S}\left( T, \mathsf{C},\sigma ,\alpha\right)$ holds, then for every~$x,y \in \Rd$,~$t \in (1,T]$ and $s \in [Ct,\infty)$, we have
\label{l.PGF.PvsbarQ}
\begin{multline*} 
\left| \int_{\Rd} (P-\overline{Q})(t,z,y) \overline{P}(s,x-z) \, dz \right| \\ \leq 
\left\{
\begin{aligned}
& \O_\sigma\left(\frac14 \mathsf{C} (t+s)^{-\frac12 -\frac d2} \exp\left( -\alpha \frac{|x-y|^2}{4(t+s)} \right)\right) & \mbox{if} & \ d>2, \\
&  \O_\sigma\left(\frac14 \mathsf{C} \log^{\frac12}(1+t)  (t+s)^{-\frac12-\frac d2}  \exp\left( -\alpha \frac{|x-y|^2}{4(t+s)} \right)\right) & \mbox{if} & \ d=2.
\end{aligned}
\right. 
\end{multline*}

\end{lemma}

\begin{proof}
We will exploit the fact that, as observed already in~\eqref{e.checkmass} and~\eqref{e.firstmoment},  
\begin{equation}
\label{e.vanishingstuffs}
\int_{\Rd} (P-\overline{Q})(t,z,y) \,dz = 0 \qquad \mbox{and} \qquad \int_{\Rd} (z-y) (P-\overline{Q})(t,z,y) \,dz = 0. 
\end{equation}
We use Taylor's theorem with remainder for $\overline{P}(s,\cdot-x)$ at $y$, that is,
\begin{multline*} 
\overline{P}(s,z-x) = \overline{P}(s,y-x) + \nabla \overline{P}(s,y-x) \cdot (z-y)
\\ + \int_{0}^1 \int_0^1 \theta_1 \nabla^2\overline{P}(s,y-x + \theta_1 \theta_2(z-y)) \, d\theta_1 \, d\theta_2 (z-y)^{\otimes 2}.
\end{multline*}
The last term can be estimated by
\begin{align*} 
\lefteqn{
\left|\int_{0}^1 \int_0^1 \theta_1 \nabla^2\overline{P}(s,y-x + \theta_1 \theta_2(z-y)) \, d\theta_1 \, d\theta_2 (z-y)^{\otimes 2} \right|
} \quad & 
\\ & 
\leq 
C \left(1 + \frac{|y-z|^2}{s} \right)^2 \frac{|y-z|^2}{s} s^{-\frac d2} \exp\left( -\frac{\alpha +\Lambda^{-1}}{2}\frac{|x-y|^2}{4s} + C \frac{|y-z|^2}{s} \right)
\\ & 
\leq C \frac{|y-z|^2}{s} s^{-\frac d2} \exp\left( -\frac{\alpha +\Lambda^{-1}}{2}\frac{|x-y|^2}{4s} + C \frac{|y-z|^2}{s} \right).
\end{align*}
We choose the constant $C$ in the statement of the lemma large enough so that 
\begin{equation*} 
\frac{C}{s} \leq \frac{\alpha}{16t} 
\quad \mbox{and} \quad 
\frac{\alpha +\Lambda^{-1}}{2} \frac 1s \geq  \frac{\alpha}{t+s}.
\end{equation*}
Notice that the first condition ensures that 
\begin{equation*} \label{}
\frac{|y-z|^2}{t} \exp\left( C \frac{|y-z|^2}{s} \right)
\leq
\frac{|y-z|^2}{t}\exp\left( \frac{\alpha}{4} \frac{|y-z|^2}{4t} \right)
\leq 
C\exp\left( \frac{\alpha}{2} \frac{|y-z|^2}{4t} \right).
\end{equation*}
Combining the above bounds and using also that~$t+s\leq Cs$, we obtain
\begin{multline*} 
\left| \overline{P}(s,z-x) - \overline{P}(s,y-x) - \nabla \overline{P}(s,y-x) \cdot (z-y) \right|
\\ \leq C t s^{-\frac12} 
\exp\left( \frac{\alpha}{2} \frac{|y-z|^2}{4t} \right) (t+s )^{-\frac 12 - \frac d2}\exp\left( -\alpha \frac{|x-y|^2}{4(t+s)} \right) .
\end{multline*}
Together with~\eqref{e.vanishingstuffs} this leads to
\begin{align} \notag 
\lefteqn{\left| \int_{\Rd} (P-\overline{Q})(t,z,y) \overline{P}(s,x-z) \, dz \right|}  &
\\ \notag & 
\leq 
C\left(\frac{t}{s}\right)^{\frac 12} 
\left(t+s \right)^{-\frac 12 - \frac d2}
\exp\left( -\alpha \frac{|x-y|^2}{4(t+s)} \right) 
t^{\frac12} \int_{\R^d} \exp\left( \frac{\alpha}{2} \frac{|y-z|^2}{4t} \right) \left|(P-\overline{Q})(t,z,y)\right| \,dz .
\end{align}
By Lemmas~\ref{l.trashK} and~\ref{l.sum-O} and the assumption that~$\mathcal{S}\left( T, \mathsf{C},\sigma ,\alpha\right)$ holds, we have
\begin{align} 
\notag  
\lefteqn{
t^{\frac12}\int_{\R^d} \exp\left( \frac{\alpha}{2} \frac{|y-z|^2}{4t} \right) \left|(P-\overline{Q})(t,z,y)\right| \,dz
} \qquad  &  
\\ 
\notag & \leq t^{\frac12} \int_{\R^d} \exp\left( \frac{\alpha}{2} \frac{|y-z|^2}{4t} \right) \left( \left| \left( P-K\right) (t,z,y) \right|+  \left| (\overline{Q}-K)(t,z,y) \right| \right)  \,dz 
\\ \notag   & \leq \left\{
\begin{aligned}
& \O_\sigma\left(  C \mathsf{C}\right) & \mbox{if} & \ d>2, \\
&  \O_\sigma\left( C \mathsf{C}  \log^{\frac12}(1+t) \right) & \mbox{if} & \ d=2.
\end{aligned}
\right. 
\end{align}
Combining the previous two displays and enlargening our choice of the constant $C$ in the statement of the lemma, if necessary, we obtain the desired estimate. 
\end{proof}

We next give the estimate of the first integral on the right side of~\eqref{e.Duhamel.P-K.first}. 

\begin{lemma} 
Fix $\sigma\in (0,2)$ and $\alpha \in \left( 0,\Lambda^{-1} \right)$. There exist~$\delta(\sigma,d,\Lambda)>0$ and $C(\sigma,\alpha,d,\Lambda) < \infty$ such that for every $\mathsf{C},T\in[1,\infty)$, if~$\mathcal{S}\left( T, \mathsf{C},\sigma ,\alpha\right)$ holds, then
for every~$x,y \in \Rd$, $t \in (1,T]$ and $s\in \left[t,\infty\right)$, we have
\label{l.PGF.(P-K)vs(P-barP)}
\begin{multline} \label{e.PGF.(P-K)vs(P-barP)}
 \int_{\Rd}  \left| (P-K)(t,z,y) \right|   \left| P(s,x,z)- \overline{P}(s,x-z)\right|   \,dz 
\\ 
\leq 
\O_\sigma\left(
C \mathsf{C} \left(\frac{s}{ t} \right)^{\frac12} s^{-\delta} (t+s)^{-\frac12 -\frac d2} \exp\left( -\alpha \frac{|x-y|^2}{4(t+s)} \right) 
\right).
\end{multline}
\end{lemma}
\begin{proof}
Fix $\tau := \frac{d(2+\sigma)}{4}<d$ and let $\X_\tau(z)$ be as in Remark~\ref{r.minimalscale.PGF}.
Using Theorem~\ref{t.PGF.basecase} and the hypothesis that~$\mathcal{S}(T,\mathsf{C},\sigma,\alpha)$ holds, Lemma~\ref{l.sum-O} yields
\begin{align}
\label{e.gosplityourself.term1.in}
\lefteqn{
\int_{\Rd} \left| (P-K)(t,z,y) \right| \left| P(s,x,z) - \overline{P}(s,x-z)\right| 
\indc_{\{ \X_\tau(z)\leq \sqrt{s}  \}}  \,dz
} 
\quad & 
\notag \\ &
\leq 
C 
\int_{\Rd} 
\left| (P-K)(t,z,y) \right|
s^{-\delta-\frac d2} \exp\left( - \alpha \frac{ |x-y|^2}{4s}  \right) \,dz
\notag \\ & 
\leq
\O_\sigma \left( 
C \int_{\Rd}  
\mathsf{C} t^{-\frac12-\frac d2} \log^\frac12(1+t)\exp\left( -\alpha \frac{|z-y|^2}{4t} \right) 
s^{-\delta-\frac d2}  \exp\left( -\alpha \frac{ |z-x|^2}{4s}  \right) 
\,dz
\right) 
\notag \\ & 
\leq 
\O_\sigma \left( C\mathsf{C} \left( \frac{s}{t} \right)^{\frac12} s^{-\frac{\delta}{2}} (t+s)^{-\frac12 -\frac d2}  \exp\left( - \alpha  \frac{ |x-y|^2}{4(t+s)}  \right)  \right).
\end{align}
Recalling the definition of~$\tau$, we have that 
\begin{equation*}
\frac{\tau}{2\sigma} 
= \frac{d}{4} + \frac{d(2-\sigma)} {8\sigma} 
\geq
\frac12 + \frac{2-\sigma}{4\sigma}. 
\end{equation*}
Thus, for $\delta (\sigma):=\frac{2-\sigma}{4\sigma} >0$, 
\begin{equation*} \label{}
\indc_{\{ \X_\tau(z)\geq \sqrt{s}  \}} 
\leq
\left( \frac{X_\tau^2(z)}{s} \right)^{\frac{\tau}{2\sigma}}
\leq 
\O_\sigma\left( C s^{-\frac12 -\delta} \right)
\end{equation*}
and the Nash-Aronson bounds together with 
Lemma~\ref{l.sum-O},
we get
\begin{align*}
\lefteqn{
\int_{\Rd} \left| (P-K)(t,z,y) \right| \left| P(s,x,z) - \overline{P}(s,x-z)\right| \indc_{\{ \X_\tau(z)\geq \sqrt{s}  \}}  \,dz
} 
\quad &
\\ & 
\leq 
C \int_{\Rd}
t^{-\frac d2} \exp\left( -\alpha \frac{|z-y|^2}{4t} \right) 
 s^{-\frac d2} \exp\left( \alpha \frac{|z-y|^2}{4s} \right) \left( \frac{X_\tau^2(z)}{s} \right)^{\frac{\tau}{2\sigma}}
\,dz
\\ & 
\leq 
\O_\sigma \left( C s^{-\delta}(t+s)^{-\frac12 -\frac d2} \exp\left( -\alpha \frac{ |x-y|^2}{4(t+s)}  \right)\right).
\end{align*}
We thus obtain, for some $\delta(\sigma,d,\Lambda)>0$, 
\begin{multline} 
\label{e.gosplityourself.term1.out}
\int_{\Rd} \left| (P-K)(t,z,y) \right| \left| P(s,x,z) - \overline{P}(s,x-z)\right| 
\indc_{\{ \X_\tau(z)\geq \sqrt{s}  \}}  \,dz
\\
\leq 
\O_\sigma \left( C s^{-\delta - \frac 12}(t+s)^{-\frac d2} \exp\left( -\frac{\beta |x-y|^2}{t+s}  \right) \right).
\end{multline}
Combining~\eqref{e.gosplityourself.term1.in} and~\eqref{e.gosplityourself.term1.out} yields~\eqref{e.PGF.(P-K)vs(P-barP)}. The proof is complete. 
\end{proof}

We turn now to the estimate of the second integral on the right side of~\eqref{e.Duhamel.P-K}. In fact, we have already estimated a similar integral in the previous section when we estimated the differences between~$P$ and~$H$ (see Lemmas~\ref{l.PGF.H.eq} and~\ref{l.PGF.w.basiclemma}). The main differences in the argument here compared to that of the previous section are~(i) we replace $\overline{P}$ by $\overline{Q}$ and we conclude differently by using sharper estimates on the correctors (that is, we use Theorem~\ref{t.correctors} instead of Proposition~\ref{p.correctorsbasecase}), and~(ii) the cutoff function in time needs to be chosen more carefully, due to the need for a sharper estimate, and this causes some additional technical difficulties due to the singularity of the Green functions (which is the reason for introducing the modification~$\hat{K}$ of $K$). 

\begin{lemma}
\label{l.PGF.secondterm}
Fix $\sigma\in (0,2)$ and $\alpha \in \left( 0,\Lambda^{-1} \right)$.
There exist~$C(\sigma,\alpha,d,\Lambda)<\infty$ and, for every~$x,y\in\Rd$, $t\in [1,\infty)$ and $s\in [t,\infty)$, 
a function $\hat{K}: [t,t+s] \times \R^d \to \R$ satisfying
\begin{equation*}  
\hat K(t+s,x) = K(t+s,x,y) \quad \text{and} \quad \hat K(t,\cdot) = K(t,\cdot,y)
\end{equation*}
and such that
\begin{multline}
\left| 
\int_{t}^{t+s} \int_{\Rd}
\left( \partial_{t'} -\nabla_z \cdot \a \nabla_z \right)\hat K(t',z) P(t+s-t',x,z) \, dz \, dt'
\right|
\\
\leq 
\left\{
\begin{aligned}
& \O_{\sigma}\left(  C  t^{-\frac 12} (t+s)^{- \frac{d}{2}} \exp\left( - \alpha \frac{|x-y|^2}{4(t+s)}\right) \right) & \mbox{if} & \ d>2, \\
&  \O_\sigma\left( C \log^{\frac12}\left(1+t+s\right)  t^{-\frac 12} (t+s)^{- \frac{d}{2}} \exp\left( - \alpha \frac{|x-y|^2}{4(t+s)}\right) \right) & \mbox{if} & \ d=2.
\end{aligned}
\right. 
\end{multline} 
\end{lemma}

\begin{proof}
Throughout the proof, we fix $\sigma\in (0,2)$, $\alpha \in \left( 0,\Lambda^{-1} \right)$, $\beta := \frac12(\Lambda^{-1} + \alpha)$,~$x,y\in\Rd$, $t,s\in [1,\infty)$ with $s\leq t$.

\smallskip

We modify $K$ close to terminal time $t+s$ by defining, for each~$z \in \R^d$, 
\begin{multline*} 
\tilde K (z) =  \overline{Q}(t+s,x,y) +  \nabla_x \overline{Q}(t+s,x,y) \cdot (z-x) \\+ \sum_{k=1}^d \left( \phi_{e_k}\left(z \right)  -\left( \phi_{e_{k}} \ast  \Phi(t+s,\cdot) \right)(y)  \right) \partial_{x_k} \overline{Q}(t+s,x,y),
\end{multline*}
and, for $t' \in [t,t+s]$ and $z \in \R^d$ we set~$\tau(t') := 1 \wedge \left( \left(t+s - 1 - t' \right) \vee 0\right)$ and
\begin{equation*} 
\hat K(t',z)  = \tau(t') K (t',z,y)  + (1-\tau(t')) \tilde K (z) .
\end{equation*}
Since $\tilde K (\cdot) $ is independent of $t'$ and $ \nabla_z \cdot \a \nabla_z \tilde K (z) = 0$ by the equation for the correctors, we have that 
\begin{align} \notag  \label{e.PGF.Duhamel.hatK}
\lefteqn{\left| 
\int_{t}^{t+s} \int_{\Rd}
\left( \partial_{t'} -\nabla_z \cdot \a \nabla_z \right)\hat K(t',z) P(t+s-t',x,z) \, dz \, dt'
\right|} \quad &
\\ \notag & \leq \left| \int_{t}^{t+s-1}  \int_{\Rd} \left(\partial_{t'} - \nabla_z \cdot \a \nabla_z \right)K (t',z,y)  P(t+s-t',x,z) \, dz \, dt' \right|
\\ &  \quad + \left|  \int_{t+s-2}^{t+s-1}  \int_{\Rd} \left(K (t',z,y) -  \tilde K (z)  \right) P(t+s-t',x,z) \, dz \, dt'  \right|.
\end{align}
We estimate the terms on the right in Steps 1-2 below. 

\smallskip

\emph{Step 1.} We show that there exists~$C(\alpha,\sigma,d,\Lambda)<\infty$ such that
\begin{multline}  \label{e.PGF.K.1}
 \int_{t}^{t+s-1} \left| \int_{\Rd} \left(\partial_{t'} - \nabla_z \cdot \a \nabla_z \right)K (t',z,y)  P(t+s-t',x,z) \, dz \right| \, dt' 
 \\\leq  \left\{
\begin{aligned}
& \O_{2+\delta}\left(  C t^{-\frac12} \Phi\left(\frac{t+s}{\alpha} , x-y \right) \right) & \mbox{if} & \ d>2, \\
&  \O_\sigma\left( C \log^{\frac12}(1+t) t^{-\frac12}  \Phi\left(\frac{t+s}{\alpha} , x-y \right) \right) & \mbox{if} & \ d=2.
\end{aligned}
\right. 
\end{multline}
First, by Lemma~\ref{l.letswritetheflux}, we get
\begin{align} 
\label{e.PGF.K.1.1}
\notag  
\lefteqn{\left| \int_{\Rd} \left(\partial_{t'} - \nabla_z \cdot \a \nabla_z \right)K (t',z,y)  P(t+s-t',x,z) \, dz \right|} \quad  &  
\\  
\notag & \leq   \int_{\Rd} \left( f(t',z,y) +  g(t',z,y) \right) \left| \nabla_{z}^2 \overline{Q}(t',z,y) \right|  \left| \nabla_{y} P(t +s -t',x,z) \right| \, dz
\\  
\notag & \quad +  \left| \int_{\Rd}   \nabla_z \cdot \left( \ahom \nabla^2 \overline{P}(t',z-y) \mathbf{m} (t',y)  \right)  P(t +s-t',x,z)  \, dz \right| 
\\  
 & \quad +   \int_{\Rd}  \left| \partial_{t'} \left( K\left(t',z,y \right) - \overline{P}(t',z-y)\right)  \right|  P(t+s-t',x,z)  \, dz ,
\end{align}
where 
\begin{align} 
\label{e.PGF.K.1.f}
f(t',z,y)  & :=   \sum_{i,j,k=1}^d \left| \mathbf{S}_{e_k,ji}\left(z \right) - \left( \mathbf{S}_{e_k,ji}\right)_{E_{t'}}\right| \\
\label{e.PGF.K.1.g}
g(t',z,y) & :=  \sum_{i,j,k=1}^d \left| \a_{ij}\left( z \right) \left( \phi_{e_k}\left(z \right)  -\left( \phi_{e_{k}} \ast  \Phi(t',\cdot) \right)(y)  \right) \right| .
\end{align}
and 
\begin{equation*} 
E_{t'} : = \left\{ 
\begin{aligned}
& B_{\sqrt{t+s-t'}}(x)  & \mbox{ if } & \  t+s-t' \leq t',\\ 
& B_{\sqrt{t'}}(y)  & \mbox{ if } & \ t+s-t' > t'.
\end{aligned}
\right.
\end{equation*}
Each term appearing in~\eqref{e.PGF.K.1.1} is a convolution with the parabolic Green function $P(t+s-t',x,\cdot)$, or its gradient, and can thus be estimated using Lemma~\ref{l.PGF.layercake}. For the application of that lemma, notice that, on the one hand,  by the Nash-Aronson bound, we have, on the one hand, for all $z\in \R^d$ and $t'>0$, that
\begin{equation*} 
(t')^{\frac 12} \left| \nabla \overline{P}(t',z) \right|  + t' \left| \nabla^2 \overline{P}(t',z) \right| +
(t')^{\frac 32}\left| \nabla^3 \overline{P}(t',z) \right| \leq C \Phi\left(\frac{t'}{\beta} , z \right).
\end{equation*}
On the other hand, by~\eqref{e.Qbar.bounded} and the Nash-Aronson bounds for $\overline{P}$, we also get, for every $\beta \in (0,\Lambda^{-1})$, that there exists a constant $C(\beta,d,\Lambda)<\infty$ such that 
\begin{equation}  \label{e.Qbar.bounded2}
t' \left| \nabla_z^2 \overline{Q}(t',z,y) \right| +  t' \left| \partial_{t'} \overline{Q}(t',z,y) \right|  + (t')^{\frac 32}\left| \nabla_z \partial_{t'}\overline{Q}(t',z,y) \right|  \leq C \Phi\left(\frac{t'}{\beta} , z-y \right). 
\end{equation}

\smallskip

First, by~\eqref{e.m.basecase.O},  we have that
\begin{multline*} 
\left| \nabla_z \cdot \left( \ahom \nabla^2 \overline{P}(t',z-y) \mathbf{m} (t',y)  \right) \right| \\ 
\leq 
\left\{
\begin{aligned}
& \O_\sigma\left(  C (t')^{-\frac 32}  \Phi\left(\frac{t'}{\beta} , z-y \right)   \right) & \mbox{if} & \ d>2, \\
&  \O_\sigma\left( C \log^{\frac12}(1+t') (t')^{-\frac 32}  \Phi\left(\frac{t'}{\beta} , z-y \right)   \right) & \mbox{if} & \ d=2.
\end{aligned}
\right. 
\end{multline*}
Thus, by Lemma~\ref{l.PGF.layercake}, 
\begin{multline} \label{e.PGF.K.1.term2}
 \int_{\Rd} \left| \nabla_z \cdot \left( \ahom \nabla^2 \overline{P}(t',z-y) \mathbf{m} (t',y)  \right)  \right|  P(t+s-t',x,z) \, dz  
 \\
\leq 
\left\{
\begin{aligned}
& \O_\sigma\left(  C (t')^{-\frac 32}  \Phi\left(\frac{t'}{\beta} , z-y \right)   \right) & \mbox{if} & \ d>2, \\
&  \O_\sigma\left( C \log^{\frac12}(1+t')  (t')^{-\frac 32}  \Phi\left(\frac{t'}{\beta} , z-y \right)   \right) & \mbox{if} & \ d=2.
\end{aligned}
\right. 
\end{multline}
Second, we expand the time derivative appearing in~\eqref{e.PGF.K.1.1} as 
\begin{multline} \label{e.PGF.K.1.term2.0}
\partial_{t'} \left( K\left(t',z,y \right) - \overline{P}(t',z-y)\right)  =  
- \partial_{t'} \left( \mathbf{m} (t',y) \cdot \nabla  \overline{P} (t',z-y) \right) 
\\  
+ \partial_{t'}  \left(  \sum_{k=1}^d \left( \phi_{e_k} \left(x \right)  -\left( \phi_{e_{k}} \ast  \Phi(t',\cdot) \right)(y)  \right) \partial_{z_k}\overline{Q}\left(t',z,y \right) \right).
\end{multline}
By~\eqref{e.m.basecase.O1} and~\eqref{e.m.basecase.O},  
\begin{multline} \label{e.PGF.K.1.term2.1}
\left| \partial_{t'} \left( \mathbf{m} (t',y) \cdot \nabla  \overline{P} (t',z-y) \right)  \right|  \\ \leq 
\left\{
\begin{aligned}
& \O_\sigma\left(  C (t')^{-\frac 32}  \Phi\left(\frac{t'}{\beta} , z-y \right)   \right) & \mbox{if} & \ d>2, \\
&  \O_\sigma\left( C \log^{\frac12}(1+t') (t')^{-\frac 32}  \Phi\left(\frac{t'}{\beta} , z-y \right)   \right) & \mbox{if} & \ d=2.
\end{aligned}
\right. 
\end{multline}
For the second term in~\eqref{e.PGF.K.1.term2.0}, we apply Theorem~\ref{t.correctors} and Remark~\ref{r.scalar.Linfty} to obtain, for $k \in \{1,\ldots,d\}$ and $t' \geq 1$, $z \in \R^d$,
\begin{equation} \label{e.PGF.K.1.term2.2}
\left|  \phi_{e_k} (z) -  \left( \phi_{e_k} \ast \Phi(t',\cdot) \right) (y) \right| \leq  \left\{
\begin{aligned}
& \O_{2+\delta}\left(  C \right) & \mbox{if} & \ d>2, \\
&  \O_\sigma\left( C \log^{\frac12}\left(1+\sqrt {t'} + |z-y|\right) \right) & \mbox{if} & \ d=2,
\end{aligned}
\right. 
\end{equation}
and, consequently, by~\eqref{e.Qbar.bounded2}, 
\begin{multline} \label{e.PGF.K.1.term2.3}
\left| \phi_{e_k} (z) - \left( \phi_{e_k} \ast \Phi(t',\cdot) \right) (y) \right| \left| \partial_{t'} \partial_{z_k}\overline{Q}\left(t',z,y \right) \right| 
\\\leq C\left\{
\begin{aligned}
& \O_{2+\delta}\left( C(t')^{-\frac 32} \Phi\left(\frac{t'}{\beta},z-y\right)\right) & \mbox{if} & \ d>2, \\
& \O_\sigma\left( C  \log^{\frac12}\left(1+ t' \right) (t')^{-\frac 32} \Phi\left(\frac{t'}{\beta},z-y\right) \right) & \mbox{if} & \ d=2.
\end{aligned}
\right.
\end{multline} 
We also obtain, by~\eqref{e.PGF.K.1.term2.2}, that
\begin{multline} \label{e.PGF.K.1.term2.4}
\sum_{k=1}^d \left| \partial_{t'}  \left( \phi_{e_{k}} \ast  \Phi(t',\cdot) \right)(y)  
 \right| \left| \partial_{z_k}\overline{Q}\left(t',z,y \right) \right| 
\\
\leq   \left\{
\begin{aligned}
& \O_{2+\delta}\left(  C (t')^{-\frac 32} \Phi\left(\frac{t'}{\beta} , z-y \right) \right) & \mbox{if} & \ d>2, \\
&  \O_\sigma\left( C  \log^{\frac12}(1+t') (t')^{-\frac 32}\Phi\left(\frac{t'}{\beta} , z-y \right) \right) & \mbox{if} & \ d=2.
\end{aligned}
\right. 
\end{multline}
Collecting estimates for the terms appearing in~\eqref{e.PGF.K.1.term2.0}, that is~\eqref{e.PGF.K.1.term2.1},~\eqref{e.PGF.K.1.term2.3} and~\eqref{e.PGF.K.1.term2.4}, we obtain, by the semigroup property and the Nash-Aronson bound for~$P$,  
\begin{multline} \label{e.PGF.K.1.term3} 
 \left| \int_{\Rd}   \partial_{t'} \left( K\left(t',z,y \right) - \overline{P}(t',z-y)\right)    P(t+s-t',x,z)  \, dz \right| \\ \leq
  \left\{
\begin{aligned}
& \O_{2+\delta}\left(  C (t')^{-\frac 32} \Phi\left(\frac{t+s}{\alpha} , x-y \right) \right) & \mbox{if} & \ d>2, \\
&  \O_\sigma\left( C \log^{\frac12}(1+t') (t')^{-\frac 32} \Phi\left(\frac{t+s}{\alpha} , x-y \right) \right) & \mbox{if} & \ d=2.
\end{aligned}
\right. 
\end{multline}

We next concentrate on the terms appearing in the first term on the right side of~\eqref{e.PGF.K.1.1}.
First, by~\eqref{e.PGF.K.1.term2.2},
\begin{equation*} 
g(t',z,y) \left| \nabla_{z}^2 \overline{Q}(t',z,y) \right| \leq  C  \left\{
\begin{aligned}
& \O_{2+\delta}\left(  C (t')^{-1} \Phi\left(\frac{t'}{\beta} , z-y \right) \right) & \mbox{if} & \ d>2, \\
&  \O_\sigma\left( C  \log^{\frac12}(1+t') (t')^{-1} \Phi\left(\frac{t'}{\beta} , z-y \right) \right) & \mbox{if} & \ d=2.
\end{aligned}
\right. 
\end{equation*}
Thus, by Lemmas~\ref{l.PGF.layercake} and~\ref{l.sum-O}, 
\begin{multline} \label{e.PGF.K.1.term1.1} 
\int_{\Rd} g(t',z,y)  \left| \nabla_{z}^2 \overline{Q}(t',z,y) \right|  \left| \nabla_{y} P(t +s -t',x,z) \right| \, dz \\ \leq    \left\{
\begin{aligned}
& \O_{2+\delta}\left(  C (t')^{- 1} (t+s-t')^{-\frac 12} \Phi\left(\frac{t+s}{\alpha} , z-y \right) \right) & \mbox{if} & \ d>2, \\
&  \O_\sigma\left( C \log^{\frac12}(1+t') (t')^{-1}  (t+s-t')^{-\frac 12}\Phi\left(\frac{t+s}{\alpha} , z-y \right) \right) & \mbox{if} & \ d=2.
\end{aligned}
\right. 
\end{multline}
Next, Proposition~\ref{p.fluxcorrectorests} for the flux correctors and Lemma~\ref{l.sum-O} yield, for $t' \leq t+s-1$, 
\begin{multline}  \label{e.PGF.K.1.term1.2} 
\left(  \int_{1}^\infty r^{d+1} \exp\left(- 2(\beta-\alpha) \frac{r^2}{4} \right) \left\| f(t',\cdot,y) \right\|_{\underline{L}^{2} \left( E_{t',r} \right)}^2 \, dr \right)^{\frac 12} 
\\  \leq \left\{
\begin{aligned}
& \O_{2+\delta}\left(  C \right) & \mbox{if} & \ d>2, \\
&  \O_\sigma\left( C \log^{\frac12}(1+t) \right) & \mbox{if} & \ d=2.
\end{aligned}
\right. 
\end{multline}
where
\begin{equation*} 
E_{t',r} : = \left\{ 
\begin{aligned}
& B_{r\sqrt{t+s-t'}}(x)  & \mbox{ if } & \  t+s-t' \leq t',\\ 
& B_{r\sqrt{t'}}(y)  & \mbox{ if } & \ t+s-t' > t'.
\end{aligned}
\right.
\end{equation*}
Let us remark that the above estimate is not valid for $t' \in ( t+s-1,t+s)$, since this would call for the boundedness of flux correctors. This forces us to use the cut-off function $\tau$. 
Now, again by Lemmas~\ref{l.PGF.layercake} and~\ref{l.sum-O}, 
\begin{multline} \label{e.PGF.K.1.term1.3} 
\int_{\Rd} f(t',z,y)  \left| \nabla_{z}^2 \overline{Q}(t',z,y) \right|  \left| \nabla_{y} P(t +s -t',x,z) \right| \, dz \\ \leq    \left\{
\begin{aligned}
& \O_{2+\delta}\left(  C (t')^{- 1} (t+s-t')^{-\frac 12} \Phi\left(\frac{t+s}{\alpha} , z-y \right) \right) & \mbox{if} & \ d>2, \\
&  \O_\sigma\left( C \log^{\frac12}(1+t') (t')^{-1} (t+s-t')^{-\frac 12}  \Phi\left(\frac{t+s}{\alpha} , z-y \right) \right) & \mbox{if} & \ d=2.
\end{aligned}
\right. 
\end{multline}
Therefore, combining~\eqref{e.PGF.K.1.term1.1}  and~\eqref{e.PGF.K.1.term1.3}, we get
\begin{multline} \label{e.PGF.K.1.term1} 
\int_{\Rd} \left( f(t',z,y) +  g(t',z,y) \right) \left| \nabla_{z}^2 \overline{Q}(t',z,y) \right|  \left| \nabla_{y} P(t +s -t',x,z) \right| \, dz
\\ \leq    
\left\{
\begin{aligned}
& \O_{2+\delta}\left(  C (t')^{- 1} (t+s-t')^{-\frac 12} \Phi\left(\frac{t+s}{\alpha} , z-y \right) \right) & \mbox{if} & \ d>2, \\
&  \O_\sigma\left( C \log^{\frac12}(1+t') (t')^{-1} (t+s-t')^{-\frac 12}  \Phi\left(\frac{t+s}{\alpha} , z-y \right) \right) & \mbox{if} & \ d=2.
\end{aligned}
\right. 
\end{multline}
Combining~\eqref{e.PGF.K.1.term1},~\eqref{e.PGF.K.1.term2},~\eqref{e.PGF.K.1.term3}, integrating over time and applying Lemma~\ref{l.sum-O}, we obtain~\eqref{e.PGF.K.1}.

\smallskip

\emph{Step 2.} 
We next show that there is a constant $C(\sigma,\alpha,d,\Lambda)<\infty$ such that 
\begin{multline}  \label{e.PGF.K.2}
\left|  \int_{t+s-2}^{t+s-1}  \int_{\Rd} \left(K (t',z,y) -  \tilde K (z)  \right) P(t+s-t',x,z) \, dz \, dt'  \right| 
\\ 
\leq
\left\{
\begin{aligned}
& \O_{2+\delta}\left(  C  (t+s)^{-\frac 12} \Phi\left(\frac{t+s}{\alpha} , z-y \right) \right) & \mbox{if} & \ d>2, \\
&  \O_\sigma\left( C (t+s)^{-\frac 12} \log^{\frac12}(1+t+s) \Phi\left(\frac{t+s}{\alpha} , z-y \right) \right) & \mbox{if} & \ d=2.
\end{aligned}
\right. 
\end{multline}
To see this, rewrite the difference as
\begin{align} \notag 
K (t',z,y) -  \tilde K (z)  & = \overline{Q}(t',z,y) - \overline{Q}(t+s,x,y) -  \nabla_x \overline{Q}(t+s,x,y) \cdot (z-x) 
\\ \notag & \quad +  \sum_{k=1}^d \left( \phi_{e_k}\left(z \right)  -\left( \phi_{e_{k}} \ast  \Phi(t',\cdot) \right)(y)  \right) \partial_{x_k} \overline{Q}(t',z,y) 
\\ \notag & \quad -  \sum_{k=1}^d \left( \phi_{e_k}\left(z \right)  -\left( \phi_{e_{k}} \ast  \Phi(t+s,\cdot) \right)(y)  \right) \partial_{x_k} \overline{Q}(t+s,x,y).  
\end{align}
It is straightforward to see by~\eqref{e.Qbar.bounded2} that, for $z\in \R^d$ and $t' \in (t+s-2,t+s]$,
\begin{multline*} 
\left| \overline{Q}(t',z,y) - \overline{Q}(t+s,x,y) -  \nabla_x \overline{Q}(t+s,x,y) \cdot (z-x) \right|
\\ \leq C \left(1+|x-z|^2 \right) (t+s)^{-1} \Phi\left(\frac{t+s}{\alpha} , x-y \right).
\end{multline*}
On the other hand, by~\eqref{e.Qbar.bounded2} and~\eqref{e.PGF.K.1.term2.2}, we have, for all $t' \in (t+s-2,t+s]$, that
\begin{multline*} 
\left| \sum_{k=1}^d \left( \phi_{e_k}\left(z \right)  -\left( \phi_{e_{k}} \ast  \Phi(t',\cdot) \right)(y)  \right) \partial_{x_k} \overline{Q}(t',z,y)  \right|
\\ \leq 
\left\{
\begin{aligned}
& \O_{2+\delta}\left(  C (t')^{-\frac12} \Phi\left(\frac{t'}{\beta} , z-y \right)  \right) & \mbox{if} & \ d>2, \\
&  \O_\sigma\left( C \log^{\frac12}\left(1+t' \right) (t')^{-\frac12} \Phi\left(\frac{t'}{\beta} , z-y \right) \right) & \mbox{if} & \ d=2.
\end{aligned}
\right. 
\end{multline*}
Using the semigroup property, together with Lemma~\ref{l.sum-O}, we obtain~\eqref{e.PGF.K.2} after integration in time.

\smallskip

\emph{Step 3.} To conclude the argument, we simply connect~\eqref{e.PGF.K.1} and~\eqref{e.PGF.K.2} with~\eqref{e.PGF.Duhamel.hatK}.  
\end{proof}

We next complete the proof of Theorem~\ref{t.PtoPbar} by combining Lemmas~\ref{l.PGF.PvsbarQ},~\ref{l.PGF.(P-K)vs(P-barP)} and~\ref{l.PGF.secondterm} and an induction argument.

\begin{proof}[{Proof of Theorem~\ref{t.PtoPbar}}]
For the readability of the argument, we will write the proof for~$d>2$. The only modification needed in the case~$d=2$ is that an additional factor of $\log^{\frac12}(t+s)$ should be inserted inside each of the~$\O_\sigma\left( \cdot \right)$ expressions below.

\smallskip

Fix $T,\mathsf{C} \in [1,\infty)$, $\sigma\in (0,2)$ and~$\alpha \in \left(0,\Lambda^{-1} \right)$ and suppose that $\mathcal{S}(T,\mathsf{C},\sigma,\alpha)$ holds. We seek to prove that $\mathcal{S}(2T,\mathsf{C},\sigma,\alpha)$ also holds, provided that $T,\mathsf{C}$ are sufficiently large. 
The strategy is to apply the formula~\eqref{e.Duhamel.P-K} with $t\in [\lambda^{-1}T,T]$ and $s=(\lambda-1) t$, where~$\lambda \in [2,\infty)$ is a fixed parameter to be selected below, and then to estimate both of the terms on the right side. 

\smallskip
Recalling~\eqref{e.Duhamel.P-K}, we have that 
\begin{align} 
\label{e.Duhamel.P-K.again}
\lefteqn{ 
(P-K)(t+s,x,y) 
} \quad &
\notag \\ &
=
\int_{\Rd} (P-K) (t,z,y) P(s,x,z) \, dz
\notag \\ & \quad 
-
\int_{t}^{t+s} 
\int_{\Rd}
\left( \partial_{t'} -\nabla_z\cdot \a \nabla_z  \right) \hat K\left( t',z\right)  P(t+s-t',x,z),
\, dz \, dt'. 
\end{align}
where $\hat K$ is constructed in Lemma~\ref{l.PGF.secondterm}. To estimate the first term, we apply~\eqref{e.Duhamel.P-K.first} with $t\in [\lambda^{-1}T,T]$ and $s:=(\lambda-1) t$ as above, to get
\begin{align} 
\label{e.Duhamel.P-K.second}
\notag  
\lefteqn{
\int_{\Rd}  (P-K)(t,z,y)  P(s,x,z)  \,dz 
} \quad & \\
& 
=   \int_{\Rd}  (P-K)(t,z,y)  \left(P(s,x,z)- \overline{P}(s,x-z) \right) \,dz
\notag \\ &  \quad 
+   \int_{\Rd}  (P-\overline{Q})(t,z,y)  \overline{P}(s,x-z)  \,dz
\notag \\ &  \quad 
+   \int_{\Rd}  (\overline{Q}-K)(t,z,y)  \overline{P}(s,x-z)  \,dz.
\end{align}
We use Lemmas~\ref{l.PGF.(P-K)vs(P-barP)},~\ref{l.PGF.PvsbarQ} and~\ref{l.trashK}, respectively, to bound the three terms on the right side of the previous display. To apply Lemma~\ref{l.PGF.PvsbarQ}, we need that $\lambda$ be sufficiently large, so we take~$\lambda:= C(\sigma,\alpha,d,\Lambda)$ so that $s \geq C t$ and $T$ large enough so that $ C \left(\frac{s}{ t} \right)^{\frac12} s^{-\delta} \leq \frac14$. Observe that, with $\lambda$ and now chosen, we have that~$t\leq s \leq Ct$. As a result, we obtain that 
\begin{align*} \label{}
\lefteqn{
\left| \int_{\Rd}  (P-K)(t,z,y)  P(s,x,z)  \,dz  \right|
} \quad & 
\\ & 
\leq 
\O_\sigma\left(
\left(
C \left( \frac{s}{t} \right)^{\frac12} +
\mathsf{C} \left( \frac14 + C \left(\frac{s}{ t} \right)^{\frac12} s^{-\delta}
\right) 
  \right) (t+s)^{-\frac12 -\frac d2} \exp\left( -\alpha \frac{|x-y|^2}{4(t+s)} \right)
\right) 
\\ & 
\leq 
\O_\sigma\left(
\left(
C +
\frac 12 \mathsf{C}
  \right) (t+s)^{-\frac12 -\frac d2} \exp\left( -\alpha \frac{|x-y|^2}{4(t+s)} \right)
\right). 
\end{align*}
To estimate the second term, we apply Lemma~\ref{l.PGF.secondterm}, which yields, since $s \leq Ct$, 
\begin{multline*} 
\left| 
\int_{t}^{t+s} \int_{t}^{t+s} 
\int_{\Rd}
\left( \partial_{t'} -\nabla_z\cdot \a \nabla_z  \right) \hat K\left( t',z\right)  P(t+s-t',x,z),
\, dz \, dt'.  
\right|
\\\leq 
\O_{\sigma}\left( C (t+s)^{-\frac12-\frac{d}{2}} \exp\left( - \alpha \frac{|x-y|^2}{4(t+s)}\right)\right).
\end{multline*}

\smallskip

Combining the previous two displays with~\eqref{e.Duhamel.P-K} and using the choices of~$\lambda$,~$t$ and $s$, we obtain finally that 
\begin{align*} \label{}
\left| (P-K)( t+s,x,y)  \right| 
\leq 
\O_\sigma\left(
\left(
C +
\frac 12 \mathsf{C} 
  \right) (t+s)^{-\frac12 -\frac d2} \exp\left( -\alpha \frac{|x-y|^2}{4(t+s)} \right)
\right). 
\end{align*}
We now observe that, if $\mathsf{C}$ is sufficiently large, we have that $C\leq \frac 12 \mathsf{C} $, and, in particular, 
\begin{equation*} \label{}
\left| (P-K)( t+s,x,y)  \right| 
\leq 
\O_\sigma\left(
\mathsf{C} (t+s)^{-\frac12 -\frac d2} \exp\left( -\alpha \frac{|x-y|^2}{4(t+s)} \right)
\right). 
\end{equation*}
This holds for every $t\in \left[ \lambda^{-1} T,T \right]$ and therefore we have that the statement $\mathcal{S}(\lambda T,\mathsf{C},\sigma,\alpha)$ is valid. Since~$\lambda\geq 2$, we have proved the implication~\eqref{e.propagate.T}. That is, there exists $C(\sigma,\alpha,d,\Lambda)<\infty$ such that, for every $T,\mathsf{C} \geq C$, 
\begin{equation}
\label{e.propagate.T.yes}
\mathcal{S}\left( T, \mathsf{C},\sigma ,\alpha \right)
\implies
\mathcal{S}\left( 2T, \mathsf{C},\sigma , \alpha \right).
\end{equation}
According to the Nash-Aronson estimate, for each fixed $T\geq 1$, there exists $C(T,\sigma,\alpha,d,\Lambda)<\infty$ such that~$\mathcal{S}\left( T, C,\sigma ,\alpha \right)$. We therefore deduce by induction that, for some constant~$C(\sigma,\alpha,d,\Lambda)<\infty$, the statement~$\mathcal{S}\left( T,C,\sigma , \alpha \right)$ holds for every $T<\infty$. That is, $\mathcal{S}\left( \infty,C,\sigma , \alpha \right)$ is valid. 

\smallskip

We now obtain the desired bound~\eqref{e.PtoPbar} from~$\mathcal{S}\left( \infty ,C,\sigma , \alpha \right)$ and~\eqref{e.QbartoPbar}. This completes the proof of the theorem. 
\end{proof}

\index{Green function!homogenization estimates|)}

\section*{Notes and references}

Optimal bounds on the decay of the parabolic semigroup were first proved by Gloria, Neukamm and Otto~\cite{GNO}, who obtained the result on the discrete lattice with much weaker stochastic integrability (i.e.,~finite moment bounds) using a very different argument than the one here. The full statement of Theorem~\ref{t.semigroup} was first proved by Gloria and Otto~\cite{GO6} using a proof which is more similar to the one presented here, based on a propagation in time. Previously, the first quantitative (but suboptimal) estimates on the semigroup were obtained in~\cite{M}. 

\smallskip

\index{Berry-Esseen theorem}

The optimal quantitative homogenization estimate for the parabolic and elliptic Green functions (Theorem~\ref{t.PtoPbar} and Corollary~\ref{c.GtoGbar}) are obtained here for the first time. Theorem~\ref{t.PtoPbar} can be considered as an ``optimal Berry-Esseen theorem for a diffusion in a random environment''; see also the discussion around \eqref{e.local.clt}. Suboptimal versions were previously obtained by different methods in~\cite{M3,M4}.



\chapter{Linear equations with nonsymmetric coefficients}
\label{c.subadd-fitz}

The analysis of Chapters~\ref{c.two} and~\ref{c.A1} was based on variational methods. The variational formulation of the equation suggested natural subadditive quantities which  turned out to be very convenient to analyze and formed the backbone of these two chapters. Since all subsequent chapters were heavily reliant on the results of Chapter~\ref{c.two}, a quantitative homogenization theory for equations lacking a variational structure must necessarily be very different from the one described in this book. 

\smallskip

On the other hand, it is well known that a general divergence-form elliptic equation 
\begin{equation}
\label{e.pdequasi}
-\nabla \cdot\left( \a(\nabla u, x) \right) = 0
\end{equation}
admits a variational interpretation---at least in the classical sense---if and only if the coefficients $\a(p,x)$ can be written as
\begin{equation*} \label{}
\a(p,x) = \nabla_p L(p,x),
\end{equation*}
for a Lagrangian $L(p,x)$ which is convex in $p$ (cf.~\cite[Chapter 8]{Evans}). In the case that $\a(p,x)$ is linear in $p$, that is, $\a(p,x) = \a(x)p$ for a matrix $\a(x)$, this condition is of course equivalent to~$\a(x)$ being symmetric. 

\smallskip

It may therefore come as a surprise that there is a completely natural and straightforward generalization of the previous chapters to the case of general (not necessarily symmetric) coefficients~$\a(x)$. This is because any divergence-form equation, even a very general quasilinear one like~\eqref{e.pdequasi}, admits a ``double-variable'' variational formulation.\footnote{The need for the extra variable should be thought of as analogous to the fact that a symmetric matrix~$A$ can be determined by evaluating the quadratic form $p\mapsto p\cdot Ap$, while if~$A$ is not symmetric one must consider the full bilinear form~$(p,q) \mapsto q\cdot Ap$. See also Remark~\ref{r.doublevariables}.} 
This variational formulation is no less natural or important than the ``classical'' one and the fact that it is much less well-known seems to be a historical quirk in the development of the calculus of variations. 
It allows the arguments of the previous section to be extended with the only additional difficulty being a notational one: we have to handle twice as many variables and become conversant in the unfamiliar (but beautiful) algebraic structure of the problem. 

\smallskip

The purpose of this chapter is to present this generalization to nonsymmetric coefficient fields, focusing primarily on the results in Chapters~\ref{c.two} and~\ref{c.regularity}. At that point, the reader will agree that the results of the other chapters can also be similarly extended.

\smallskip

We therefore depart from the standing assumptions by enlarging the probability space~$\Omega$. We first modify the uniform ellipticity condition in~\eqref{e.ue} by requiring that  
\begin{equation} 
\label{e.ue.Fitz}
\forall \xi \in \Rd, \quad \left| \xi \right|^2 \leq \xi\cdot \a(x) \xi
 \quad \mbox{and} \quad 
\Ll|\a(x) \xi \Rr| \leq \Lambda \left| \xi \right|.  
\end{equation}
We then define 
\begin{equation*} \label{}
\Omega:= 
\left\{ \a \,:\, \a \ \mbox{is a Lebesgue measurable map from $\Rd$ to $\R^{d\times d}$ satisfying~\eqref{e.ue.Fitz}} \right\}. 
\end{equation*}
The rest of the assumptions read verbatim as in the symmetric case. The transpose of~$\a$ is denoted by~$\a^t$. 

\smallskip

We begin in Section~\ref{s.variationalform} by reviewing this somewhat non-standard variational formulation, as usual in the linear case for simplicity. This variational formulation for non-self adjoint elliptic operators (as well as similar variational principles for parabolic equations, gradient flows and other problems) have been discovered and rediscovered multiple times, in the PDE community (see~\cite{BrEk,Krylov}), by convex analysts~\cite{Fitz,MLT} and by numerical analysts, who use it to efficiently solve non-self adjoint problems~\cite{CLMM,BG}. Recently it has appeared in~\cite{DLPS} in a context quite close to ours, namely the derivation of coarse graining procedures for analyzing singular limits of evolution equations.

\smallskip

In the rest of the chapter, we define appropriate subadditive quantities and then rewrite the arguments of Chapter~\ref{c.two} to obtain the analogous quantitative estimates on their convergence. This chapter is a simplification of~\cite{AM}, which considered general quasilinear equations.

\section{Variational formulation of general linear equations}
\label{s.variationalform}

It is easy to find a variational principle for any PDE: we can just write the equation with all nonzero terms on the left side, square it, and look for a minimizer. For a second-order elliptic equation, this procedure leads to an integral functional which involves second derivatives. One way to avoid this feature is to write the PDE as a system of first-order equations and \emph{then} square it. Or, we can square one of the equations and think of the other as a constraint. For example, the equation
\begin{equation} \label{e.PDEagain}
-\nabla \cdot \left( \a \nabla u \right) = 0 \quad \mbox{in} \ U,
\end{equation}
can be written as 
\begin{equation*} \label{}
\left\{
\begin{aligned}
& \a\nabla u - \g = 0 & \mbox{in} & \ U, \\
& \nabla \cdot \g = 0 & \mbox{in} & \ U.
\end{aligned}
\right. 
\end{equation*}
The latter can also be written as 
\begin{equation*} \label{}
u\in H^1(U) \quad \mbox{is a (null) local minimizer of} \quad \inf_{\g\in \Ls(U)} \frac12 \left\| \a\nabla u - \g \right\|^2_{L^2(U)}
\end{equation*}
where $\Ls(U)$ is the subspace of $L^2(U;\Rd)$ of solenoidal (divergence-free) vector fields defined by
\begin{equation}
\label{e.Ls}
\Ls(U):= \left\{ \g \in L^2(U;\Rd) \,:\, \forall \phi\in H^1_0(U), \ \int_U \g\cdot \nabla \phi = 0 \right\}.
\end{equation}
We can also write this as a minimization problem for the pair $(u,\g) \in H^1(U) \times \Ls(U)$ as 
\begin{equation*} \label{}
(u,\g) \in H^1(U) \times \Ls(U) \quad \mbox{is a minimizer of} \quad \int_U 
\left( \frac12 \left| \a\nabla u \right|^2 +\frac12  \left| \g \right|^2 - \g \cdot \a\nabla u \right). 
\end{equation*}
The drawback we immediately encounter is that this functional is not uniformly convex over $\left( w+H^1_0(U) \right) \times\Ls(U)$ for a given boundary condition $w\in H^1(U)$. In the case that $\a$ is symmetric, we can resolve this problem by replacing $\frac12 \left\| \a\nabla u - \g \right\|^2_{L^2(U)}$ with $\frac12 \left\| \a^{\frac12}\nabla u - \a^{-\frac12} \g \right\|^2_{L^2(U)}$ which leads to the minimization of 
\begin{equation} 
\label{e.convex.crossterm}
\int_U 
\left( \frac12 \nabla u\cdot \a\nabla u +\frac12  \g \cdot \a^{-1} \g - \g \cdot \nabla u \right).
\end{equation}
The integral of $\g \cdot \nabla u$ is actually linear in $(u,\g) \in \left( w+H^1_0(U) \right) \times \Ls(U)$ and does not depend on~$u$ since, for any $w \in H^1(U)$, we have
\begin{equation} 
\label{e.islinear}
\forall  u\in w+H^1_0(U), \quad \int_U  \g \cdot \nabla u  = \int_U \g \cdot \nabla w.  
\end{equation}
The functional in \eqref{e.convex.crossterm} is thus uniformly convex, jointly in $(u,\g) \in \left( w+H^1_0(U) \right) \times \Ls(U)$. 
Of course, the quadratic term in $\g$ does not depend on $u$ either, and so, in the case that $\a$ is symmetric, the situation simplifies: we deduce that 
\begin{align*} \label{}
\lefteqn{
u \ \mbox{is a solution of~\eqref{e.PDEagain}}
}  \qquad & \\
& \iff
\inf_{\g \in \Ls(U)}\int_U  \left( \frac12 \nabla u\cdot \a\nabla u +\frac12  \g \cdot \a^{-1} \g - \g \cdot \nabla u \right) = 0 \\
& \implies \forall w\in u+ H^1_0(U), 
\quad \int_U \frac12 \nabla u\cdot \a\nabla u  \leq  \int_U  \frac12 \nabla w\cdot \a\nabla w
\end{align*}
and we recover that $u$ is a minimizer of the classical Dirichlet energy. 

\smallskip

In the nonsymmetric case, we can generalize the above calculation by replacing $\frac12 \left\| \a^{\frac12}\nabla u - \a^{-\frac12} \g \right\|^2_{L^2(U)}$ with
\begin{equation*} \label{}
\fint_U \frac12 \left( \a\nabla u - \g \right) \cdot \s^{-1}\left( \a\nabla u - \g \right) 
\end{equation*}
where 
\begin{equation*}
\a = \s +\m  \quad \mbox{such that} \quad \s  = \s^t,\  \m  = -\m^t,
\end{equation*}
is the unique decomposition of $\a(x)$ into its symmetric part $\s (x):= \frac12(\a+\a^t)(x) $ and its skew-symmetric part $\m(x):=\frac12(\a-\a^t)(x)$. Expanding the square, we can write the integrand as
\begin{equation*} \label{}
\frac12 \left( \a\nabla u - \g \right) \cdot \s^{-1}\left( \a\nabla u - \g \right) 
= A(\nabla u,\g,x) - \g\cdot \nabla u,
\end{equation*}
where 
\begin{equation}
\label{e.def.A}
A(p,q,x):= \frac12 p \cdot \s (x)p + \frac12 (q-\m(x)p) \cdot \s ^{-1}(x) (q-\m(x)p).
\end{equation}
We obtain therefore that 
\begin{equation} 
\label{e.PDEAequiv}
u \in H^1(U) \ \text{is a solution of~\eqref{e.PDEagain}}
 \iff
\inf_{\g \in \Ls(U)}\int_U \left( A(\nabla u, \g,x) -\nabla u\cdot \g\right) = 0.
\end{equation}
Notice that $A$ satisfies:
\begin{itemize}  
\item the mapping $(p,q) \mapsto A(p,q,x)$ is jointly uniformly convex and $C^{1,1}$ in the sense that there exists $C(d,\Lambda)<\infty$ such that, for every $p_1,p_2,q_1,q_2\in\Rd$,
\begin{multline}
 \label{e.AconvexC11}
\frac1C\left( \left| p_1-p_2 \right|^2+ \left| q_1-q_2 \right|^2 \right) 
\\
\leq
\frac12A(p_1,q_1,x) + \frac12A(p_2,q_2,x) - A\left( \frac12p_1+\frac12p_2,\frac12q_1+\frac12q_2,x \right) 
\\
\leq
C\left( \left| p_1-p_2 \right|^2+ \left| q_1-q_2 \right|^2 \right),
\end{multline}
Indeed, the mapping~$(p,q) \mapsto A(p,q,x)$ is a uniformly convex quadratic form on~$\R^{2d}$. 

\item for every $p,q \in \Rd$, 
\begin{equation}
\label{e.fitz.positivity}
A(p,q,x) \ge p\cdot q \;.
\end{equation}
\item for every $p, q \in \Rd$, we have the equivalence
\begin{equation}
\label{e.fitz.equality}
A(p,q,x) = p\cdot q \quad \iff \quad q = \a(x) p.
\end{equation}
\end{itemize}

\begin{remark}
\label{r.doublevariables}
Note that~\eqref{e.fitz.equality} implies that the quadratic form~$(p,q) \mapsto A(p,q,x)$ completely characterizes the matrix $\a(x)$. We are thus in a similar situation as in the case of symmetric~$\a$, in which the simpler quadratic form $p\mapsto p\cdot \a(x)p$ characterized~$\a(x)$, although there are twice as many variables. 
\end{remark}

The above properties of~$A$ can be summarized in convex analytic language as the statement that the function $(p,q)\mapsto A(p,q,x)$ is a \emph{variational representation} of the linear mapping $p\mapsto \a(x)p$. 

\begin{definition}
\label{def.variationalrep}
\index{variational representative}
Let $\phi:\Rd \to \Rd$. We say that $F:\Rd\times\Rd \to \R$ \emph{variationally represents} $\phi$ if $F$ is convex and $F(p,q) \geq p\cdot q$ for every $p,q\in\Rd$ with equality  if and only if $q =\phi(p)$. 
\end{definition}

The canonical example of Definition~\ref{def.variationalrep} is, for a given convex function $L$, that the mapping $(p,q) \mapsto L(p)+L^*(q)$ represents $p\mapsto \nabla L(p)$. It turns out that a map $\phi$ has a representative if and only if it is \emph{monotone}, that is, 
\index{monotone map}
\begin{equation*} \label{}
\forall p_1,p_2\in\Rd,\quad (p_1-p_2) \cdot \left( \phi(p_1)-\phi(p_2)\right) \geq 0. 
\end{equation*}
This discovery is usually attributed to Fitzpatrick~\cite{Fitz}, but actually first appeared a few years earlier in a paper of Krylov~\cite{Krylov}. We remark that representatives are \emph{not} unique. They can be split into a sum of a function of $p$ and a function of $q$ if and only if the monotone map is the gradient of a convex function. A uniformly convex representative of $\phi$ can be found if and only if $\phi$ is \emph{uniformly monotone}, that is, there exists $\lambda>0$ such that 
\begin{equation*} \label{}
\forall p_1,p_2\in\Rd,\quad (p_1-p_2) \cdot \left( \phi(p_1)-\phi(p_2)\right) \geq  \lambda \left| p_1 - p_2 \right|^2. 
\end{equation*}
We mention these facts (proofs of which can be found in~\cite{AM}) so that the reader is aware that the variational setting described here is part of a much larger and richer convex analytic framework. 

\smallskip

In view of~\eqref{e.islinear},~\eqref{e.PDEAequiv} and~\eqref{e.AconvexC11}, we have therefore succeeded in recasting the PDE as a uniformly convex optimization problem. The minimizer of this problem give us the solution $u$ as well as its flux, $\g=\a\nabla u$. 

\smallskip

For our reference, we summarize what we have shown above in the following proposition. 

\begin{proposition}
\label{p.yesvariational}
For every $u\in H^1(U)$ and $u^*\in H^{-1}(U)$, we have that $u$ satisfies the equation
\begin{equation*}
-\nabla \cdot \left( \a(x) \nabla u \right) = u^*\quad \mbox{in} \ U
\end{equation*}
if and only if 
\begin{equation}
\label{e.yesvariational}
0 = \inf \left\{  \int_U \left( A(\nabla u,\g,\cdot) - \g\cdot \nabla u \right) \,:\, \g \in L^2(U;\Rd), \ -\nabla \cdot \g = u^* \right\}. 
\end{equation}
In this case, the infimum on the right side above is attained for $\g = \a \nabla u$.
\end{proposition}
\begin{proof}
By the direct method, there exists $\g$ realizing the infimum in \eqref{e.yesvariational}. The equivalence is then immediate from \eqref{e.fitz.positivity} and \eqref{e.fitz.equality}.
\end{proof}

The variational formulation can be used to prove the existence of minimizers of the integral function over $w+H^1_0(U)$ via the direct method, and these minimizers can be shown to be null using convex analytic arguments. This alternative approach to well-posedness of boundary value problems for linear elliptic equations is explored in the following exercise. 

\begin{exercise}  
Let $X$ be a reflexive Banach space, $X^*$ denote its dual and $B$ be a continuous and coercive bilinear form on $X$. The assumption of coercivity means that there exists a constant $\lambda > 0$ such that
\begin{equation*}  
\forall u \in X, \quad B(u,u) \ge \lambda \|u\|_{X}^2. 
\end{equation*}
The \emph{Lax-Milgram lemma} states that, for every $u^* \in X^*$, there exists $u \in X$ satisfying
\begin{equation*}  
\forall v \in X, \quad B(u,v) = \langle v,u^* \rangle,
\end{equation*}
where $\langle \cdot,\cdot\rangle$ denotes the duality pairing between $X$ and $X^*$. In this exercise, we outline a proof of this result using a variational, convex analytic argument. In analogy with \eqref{e.def.A}, we define a quadratic form $S$ on $X$ by
\begin{equation*}  
S(u) : = \frac 1 2 B(u,u),
\end{equation*}
denote by $M$ the linear map $M : X \to X^*$ characterized by
\begin{equation*}  
\forall u, v \in X, \quad \langle v, M u \rangle = \frac 1 2 \Ll( B(u,v) - B(v,u) \Rr) ,
\end{equation*}
and let $A$ be the bilinear form $A:X \times X^*\to \R$ such that
\begin{equation*}  
\forall u \in X, u^* \in X^*, \quad A(u,u^*) = S(u) + S^*(u^* - M u),
\end{equation*}
where $S^*:\X^*\to\R$ denotes the convex dual of $S$, that is,
\begin{equation}  
\label{e.def.S*}
S^*(u^*) := \sup_{u \in X} \Ll( \langle u, u^*\rangle - S(u) \Rr) .
\end{equation}
We break the proof into five sub-exercises. 

\begin{enumerate}  

\item Let $S'$ be the symmetric bilinear form on $X \times X$ defined by
\begin{equation*}  
S'(u,v) := \frac 1 2 \Ll( B(u,v) + B(v,u) \Rr) .
\end{equation*}
Using the coercivity assumption, show that there exists a unique $u \in X$ achieving the supremum in \eqref{e.def.S*} which is characterized by
\begin{equation*}  
\forall v \in X, \quad \langle v, u^* \rangle = S'(u,v).
\end{equation*}

\item Deduce that the conclusion of the Lax-Milgram lemma is equivalent to the statement that, for every $u^* \in X^*$,
\begin{equation}  
\label{e.null.min}
\inf_{u \in X} \Ll( A(u,u^*) - \langle u,u^* \rangle \Rr) = 0.
\end{equation}

\item Fix $u^* \in X^*$ and define $G:X^* \to \R$ by 
\begin{equation*}  
G(v^*) := \inf_{u \in X} \Ll(A(u,u^* + v^*) -  \langle u,u^* \rangle\Rr) .
\end{equation*}
Deduce from the joint convexity of $A$ that $G$ is itself convex. Using the continuity and coercivity of $A$, show that $G$ is locally bounded above, and thus that $G^{**} = G$ (see \cite[Chapter~1]{ET} if needed).

\item Argue that, for every $v \in X$,
\begin{align*}  
G^*(v) & = \sup \Ll\{ \langle v,v^* \rangle - A(u,u^* + v^*) + \langle u,u^* \rangle, \ u \in X, \ v^* \in X^* \Rr\}  \\
& = S(v) + S^*(u^* - Mv) - \langle v,u^* \rangle \\
& \ge 0.
\end{align*}

\item Combine the previous steps to obtain that $G(0) \le 0$ and deduce~\eqref{e.null.min}.

\end{enumerate}
\end{exercise}

\smallskip

In the rest of this chapter, we adopt the variational point of view and apply it to the homogenization of our equation~\eqref{e.PDEagain}. 

\section{The double-variable subadditive quantities}

Rather than attempting to homogenize~\eqref{e.PDEagain} directly, we seek to homogenize the integral functional $(\nabla u,\g) \mapsto  \int_U A(\nabla u, \g,\cdot)$ by finding an effective, constant-coefficient $\overline{A}$ so that this functional behaves, on large scales, like $(\nabla u,\g) \mapsto  \int_U \overline{A}(\nabla u, \g)$. We will discover later that the $\overline{A}$ that we produce will be a representative for an homogenized matrix $\ahom$, but it is the integral functional that is the primary object of study here. 

\smallskip

We begin by defining the subadditive quantities which are analogous to~$\nu$ and~$\nu^*$ and are coarsened representations of the mapping $(\nabla u,\g) \mapsto \int_U A(\nabla u, \g)$, in the same way that~$\nu$ and $\nu^*$ are coarsened versions of the mapping $\nabla u \mapsto \int_U \frac 1 2 \nabla u \cdot \a \nabla u$. 

\smallskip

For each $p,q \in \Rd$, we set
\begin{equation}
\label{e.def.mu}
\mu(U,p,q) := \inf \left\{  \fint_U A(\nabla v, \h, \cdot) \,:\, v \in \ell_p + H^1_0(U), \ \h \in q + \Lso(U)\right\}.
\end{equation}
Recall that $\Lso(U)$ is the set of solenoidal vector fields with null flux at the boundary defined in \eqref{e.def.Ls0}.
The approximately dual quantity is defined, for each $q^*, p^* \in \Rd$, by
\begin{equation}
\label{e.def.mu*}
\mu^*(U,q^*,p^*) := \sup_{u \in H^1(U), \, \g \in \Ls(U)} \fint_U \Ll(q^* \cdot \nabla u + p^* \cdot \g - A(\nabla u, \g, \cdot)  \Rr).
\end{equation} 
We next show that, in the symmetric case, the quantities $\mu$ and $\mu^*$ defined above reduce to the quantities $\nu$ and $\nu^*$ introduced in \eqref{e.def.nu} and \eqref{e.def.nu*}. In fact, they are both equal to $J(U,p,q) + p\cdot q$. 

\begin{lemma}
\label{l.mu.when.sym} 
Suppose $\a(x)=\a^t(x)$ for a.e.\ $x \in U$. Then, for every $p, q \in \Rd$,
\begin{equation}  \label{e.mu.sym}
\mu(U,p,q) =  \nu(U,p) + \nu^*(U,q) = \mu^*(U,q,p).
\end{equation}
\end{lemma}
\begin{proof}
Since we assume $\a$ to be symmetric in this proposition, the function $F$ simplifies into 
\begin{equation*}  
A(p,q,x) = \frac 1 2 p \cdot \a(x) p + \frac 1 2 q \cdot \a^{-1}(x) q ,
\end{equation*}
and therefore,
\begin{equation} 
 \label{e.mu.split}
\mu(U,p,q) \\
= \inf_{v \in \ell_p+H^1_0(U)} \fint_U \frac 1 2 \nabla v \cdot \a \nabla v + \inf_{\h \in q+\Lso(U)} \fint_U \frac 1 2 \h \cdot \a^{-1} \h.
\end{equation}
The first term above is $\nu(U,p)$. We denote by $\h(\cdot,U,p,q)$ the solenoidal field achieving the infimum in the last term above. Computing the first variation, we find that
\begin{equation*} \label{}
\mbox{for every} \ \g \in \Lso(U), \quad \int_U \g \cdot \a^{-1} \h = 0.
\end{equation*}
In other words, $\a^{-1} \h \in \Lso(U)^\perp$. 
In view of the definition of $\Lso(U)$ and since $\nabla (H^1(U))$ is a closed subspace of $L^2(U;\Rd)$, we have $\Lso(U)^\perp = \nabla (H^1(U))$. Therefore, there exists $u \in H^1(U)$ such that 
\begin{equation*}  
\h = \a \nabla u. 
\end{equation*}
Since $\h \in q+ \Lso(U)$, we have
\begin{equation*}  
\forall \phi \in H^1(U), \quad \int_U \nabla \phi \cdot (\a \nabla u - q) = 0.
\end{equation*}
That is, the function $u$ solves the Neumann problem \eqref{e.Neumann}. In particular,~$u$ satisfies~\eqref{e.Jfirstvar} for $p=0$, and is therefore the maximizer in the definition of $\nu^*(U,q)$. By the first variation for $\nu^*(U,q)$, see \eqref{e.firstvarJ2}, we have
\begin{equation*}  
\nu^*(U,q) = \fint_U \frac 1 2 \nabla u \cdot \a \nabla u.
\end{equation*}
This coincides indeed with the last term in \eqref{e.mu.split}, and thus completes the proof of the first equality of~\eqref{e.mu.sym}. The proof of the second equality is similar.
\end{proof}

\begin{exercise} 
Give a detailed proof of the second equality of~\eqref{e.mu.sym}.
\end{exercise}

As in the symmetric case, we are going to combine the two subadditive quantities into one master quantity which we denote by~$J$. In order to write this in the most efficient way, we first introduce some further notation. First, we identify the quadratic form $(p,q) \mapsto A(p,q,x)$ with the matrix $\Abf(x)$ such that 
\begin{equation*} \label{}
A(p,q,x) = \frac12 \begin{pmatrix}  p \\ q \end{pmatrix} \cdot \Abf(x) \begin{pmatrix} p \\ q \end{pmatrix}.
\end{equation*}
The matrix is explicitly given by
\begin{equation*} \label{}
\Abf(x) : = \begin{pmatrix} 
\s (x)-\m(x)\s^{-1}(x)\m(x) & \m(x) \s ^{-1}(x) \\
- \s ^{-1}(x) \m(x) & \s ^{-1}(x)
\end{pmatrix}.
\end{equation*}
For each domain $U\subseteq\Rd$, we define a vector space~$\S(U) \subseteq \Lpot(U)\times\Ls(U)$ by
\begin{multline} 
\label{e.def.S}
\S(U) := 
\Big\{ S \in \Lpot(U)\times\Ls(U) \ : \\
 \forall S' \in \Lpoot(U)\times\Lso(U), \ \int_U S \cdot \Abf S' = 0 
\Big\}.
\end{multline}
This definition can be compared with that for~$\A(U)$, see~\eqref{e.def.A(U)}. Recall that the spaces~$\Lpot(U)$,~$\Ls(U)$,~$\Lpoot(U)$ and~$\Lso(U)$ are defined  in~\eqref{e.Lpot.not},~\eqref{e.Ls.not},~\eqref{e.def.Lp0} and~\eqref{e.def.Ls0}. 

\smallskip

The master quantity, which is a variant of the quantity introduced in~\eqref{e.Jsplitting} (see also~\eqref{e.formulafornuJpq}) and which we also denote by~$J$, is defined for every bounded open Lipschitz domain $U\subseteq\Rd$ and $X,X^* \in \R^{2d}$ by
\begin{equation}
\label{e.def.bigJ.compressed}
J\left( U, X, X^* \right) := 
\\
\sup_{S \in \S(U)} 
\fint_U\left( 
-\frac12 S \cdot \Abf S
- X \cdot \Abf S 
+ X^* \cdot S
\right).
\end{equation}
We denote the maximizer in the definition of $J$ by 
\begin{equation*} \label{}
S(\cdot,U,X,X^*):= \mbox{unique maximizer in the definition of $J(U, X, X^*)$.}
\end{equation*}
In order to emphasize the ``double-variable'' nature of the function $J$, we note that the definition \eqref{e.def.bigJ.compressed} can be rewritten, for every $p,q,p^*,q^* \in \Rd$, as
\begin{multline}
\label{e.def.bigJ}
J\left( U,  \begin{pmatrix}  p \\ q \end{pmatrix} ,  \begin{pmatrix}  q^* \\ p^* \end{pmatrix}  \right) = 
\\
\sup_{(\nabla v, \g) \in \S(U)} 
\fint_U\left( 
-\frac12  \begin{pmatrix}  \nabla v \\ \g \end{pmatrix}  \cdot \Abf  \begin{pmatrix} \nabla v \\ \g \end{pmatrix} 
-  \begin{pmatrix}  p \\ q \end{pmatrix} \cdot \Abf  \begin{pmatrix}  \nabla v \\ \g \end{pmatrix}  
+  \begin{pmatrix}  q^* \\ p^* \end{pmatrix} \cdot  \begin{pmatrix}  \nabla v \\ \g \end{pmatrix}
\right).
\end{multline}
\smallskip
We first prove the analogous statement to Lemma~\ref{l.Jsplitting}.

\begin{lemma}
\label{l.Jsplitting.bigJ}
For every $X,X^*\in\R^{2d}$,
\begin{equation} 
\label{e.Jsplitting.bigJ}
J(U,X,X^*)= \mu(U,X)+\mu^*(U,X^*) - X\cdot X^*.  
\end{equation}
Moreover, the maximizer $S(\cdot,U,X,X^*)$ is the difference between the maximizer of~$\mu^*(U,X^*)$ in~\eqref{e.def.mu*} and the minimizer of $\mu(U,X)$ in~\eqref{e.def.mu}.
\end{lemma}
\begin{proof}

We first show that, for every $X^*\in \R^{2d}$, 
\begin{equation}
\label{e.mu*.S}
\mu^*(U,X^*) = \sup_{S \in \S(U)} \fint_U \Ll( X^* \cdot S - \frac 1 2 S \cdot \Abf S \Rr) .
\end{equation}
We may rewrite the definition of $\mu^*(X^*)$ in \eqref{e.def.mu*} as
\begin{equation}  
\label{e.redef.mu*}
\mu^*(U,X^*) = \sup_{S \in \Lpot(U) \times \Ls(U)} \fint_U \Ll( X^* \cdot S - \frac 1 2 S \cdot \Abf S \Rr) .
\end{equation}
Note that for every $X^* \in \R^{2d}$ and $S \in \Lpoot(U) \times \Lso(U)$, we have
\begin{equation}  
\label{e.null.Lag}
\fint_U X^* \cdot S = 0.
\end{equation}
Let $S^*$ denote the maximizer in~\eqref{e.redef.mu*}. Combining the first variation for~$\mu^*$ and~\eqref{e.null.Lag}, we deduce that 
\begin{equation*}  
\forall\, S' \in \Lpoot(U) \times \Lso(U), \qquad \int_U S' \cdot \Abf S^* = 0.
\end{equation*}
That is, $S^*\in \S(U)$, and thus~\eqref{e.mu*.S} holds.

\smallskip

We can now repeat the proof of Lemma~\ref{l.Jsplitting}, with only minor changes to the notation. Let $X = (p,q) \in \R^{2d}$, and let 
\begin{equation}
\label{e.def.S0}
S_0 = (\nabla v, \h) \in X + \Lpoot(U) \times \Lso(U)
\end{equation}
denote the minimizing pair in the definition of $\mu(U,p,q) = \mu(U,X)$, see \eqref{e.def.mu}. For every $S \in \S(U)$, we have
\begin{multline}
\label{e.mu.nu.bigJ}
\mu(U,X) + \fint_U  \Ll( X^* \cdot S - \frac 1 2 S \cdot \Abf S \Rr)  - X \cdot X^* \\
 = \fint_U \Ll( \frac 1 2 S_0 \cdot \Abf S_0 - \frac 1 2 S \cdot \Abf S + X^* \cdot S \right) - X \cdot X^*.
\end{multline}
In view of \eqref{e.def.S0}, we have 
\begin{equation*}  
X  = \fint_U S_0.
\end{equation*}
Since $S \in \S(U)$, we also have
\begin{equation*}  
\fint_U S \cdot \Abf S_0 = \fint_U S \cdot \Abf X,
\end{equation*}
and this last identity holds true in particular for $S = S_0$. 
We deduce that the left side of \eqref{e.mu.nu.bigJ} equals
\begin{equation*}  
\fint_U \Ll( -\frac 1 2 \Ll( S - S_0 \Rr) \cdot \Abf \Ll( S - S_0 \Rr) 
- X \cdot \Abf\Ll( S - S_0 \Rr) + X^* \cdot (S - S_0)\Rr).
\end{equation*}
Comparing with \eqref{e.mu*.S} and \eqref{e.def.bigJ.compressed}, we obtain the announced result.
\end{proof}
We now show an estimate for elements of $\S(U)$ in the spirit of the interior Caccioppoli inequality for elements of $\A(U)$ proved in~Lemma~\ref{l.Caccioppoli appendix}.  

\begin{lemma}[Caccioppoli inequality in $\S(U)$]
\label{l.cacciop.S}
Let $V \subset \Rd$ be a domain such that $\bar V \subset U$. There exists a constant $C(U,V,d,\Lambda) < \infty$ such that, for every $S \in \S(U)$,
\begin{equation*}  
\|S\|_{L^2(V)} \le C \|S\|_{\hat H^{-1}(U)}.
\end{equation*}
\end{lemma}
\begin{proof}
Let $V', V''\subseteq\Rd$ be two Lipschitz domains such that $\bar V \subset V'$, $\bar {V'} \subset V''$ and $\bar {V''} \subset U$. Let $\cu \subset \Rd$ be a cube such that $\bar U \subset \cu$. Denote by $H^1_{\mathrm{per}}(\cu)$ the space of periodic functions over $\cu$ with square-integrable gradient, and by $L^2_\mathrm{sol,per}(\cu)$ the orthogonal space to $\nabla H^1_{\mathrm{per}}(\cu)$ (compare with \eqref{e.Ls.not}). We also let $\eta \in C^\infty_c(\Rd)$ be such that $\eta \equiv 1$ on $V$ and $\eta \equiv 0$ outside of $V'$, and $\eta' \in C^\infty_c(\Rd)$ be such that $\eta' \equiv 1$ on $V''$ and $\eta' \equiv 0$ outside of $U$. Throughout the proof, the constant $C$ is allowed to depend on the choice of $V'$, $V''$, $\eta$ and $\eta'$, and may change from one occurence to another. We fix 
\begin{equation*}  
S = \begin{pmatrix}
\nabla u \\ \g
\end{pmatrix}
\in \S(U),
\end{equation*}
with $u \in H^1(U)$ of mean zero. We decompose the rest of the proof into three steps.

\smallskip

\emph{Step 1.} 
In this step, we define a Helmholtz-Hodge decomposition for $\eta'\g$. Let $h \in H^1_{\per}(\cu)$ be the unique mean-zero periodic solution of
\begin{equation}  
\label{e.caccio.S.eq.h}
\Delta h = \nabla \cdot (\eta' \g),
\end{equation}
and for each $i,j \in \{1,\ldots,d\}$, let $\mathbf T_{ij} \in H^1_{\per}(\cu)$ be the unique mean-zero periodic solution of
\begin{equation}
\label{e.caccio.S.eq.T}
\Delta \mathbf T_{ij} = \partial_j (\eta'\g_i) - \partial_i (\eta'\g_j).
\end{equation}
These quantities provide us with the decomposition 
\begin{equation}
\label{e.caccio.S.decomp.g}
\eta' \g = \nabla h + \nabla \cdot \mathbf T + (\eta'\g)_\cu,
\end{equation}
where $\nabla \cdot \mathbf T$ is the vector field whose $i$-th component is given by
\begin{equation*}  
\sum_{j = 1}^d \partial_j \mathbf T_{ij},
\end{equation*}
and where we recall that $(\eta' \g)_\cu = \fint_\cu \eta'\g$. 
Indeed, by an explicit computation, we verify that each coordinate of the vector field
\begin{equation*}  
\eta'\g - \nabla h - \nabla \cdot \mathbf T
\end{equation*}
is harmonic. Since this vector field is also periodic in $\cu$, it must be constant, and the constant is identified as $(\eta' \g)_\cu$ since $\nabla h$ and $\nabla \cdot \mathbf T$ are of mean zero.

\smallskip 

\emph{Step 2.} In this step, we show the following three estimates:
\begin{equation}
\label{e.caccio.S.est.h}
\|\nabla h\|_{L^2(V')} \le C \|\g\|_{\hat H^{-1}(U)},
\end{equation}
\begin{equation}
\label{e.caccio.S.est.T}
\|\mathbf{T}\|_{L^2(V')} \le C \|\g\|_{\hat H^{-1}(U)},
\end{equation}
and 
\begin{equation}
\label{e.caccio.S.est.u}
\|u\|_{L^2(V')} \le C \|\nabla u\|_{\hat H^{-1}(U)}.
\end{equation}
By the weak formulation of the equation \eqref{e.caccio.S.eq.h}, we see that
\begin{equation*}  
\|\nabla h\|_{H^{-1}_\per(\cu)} \le C \|\eta'\g\|_{H^{-1}_\per(\cu)} \le C \|\g\|_{\hat H^{-1}(U)},
\end{equation*}
where we write $H^{-1}_\per(\cu)$ for the dual space to $H^{1}_\per(\cu)$.
Denoting by $\phi \in H^1_\per(\cu)$ the unique mean-zero periodic solution of
\begin{equation*}  
-\Delta \phi = h,
\end{equation*}
we get via an integration by parts that
\begin{equation*}  
\int_\cu |\nabla^2 \phi|^2 = \int_\cu |\Delta \phi|^2 = \int_\cu h^2.
\end{equation*}
Moreover,
\begin{equation*}  
\int_{\cu} h^2 = \int_\cu \nabla \phi \cdot \nabla h \le \|\nabla h\|_{H^{-1}_\per(\cu)} \, \|\nabla \phi\|_{H^1_\per(\cu)}.
\end{equation*}
Combining the last three displays yields
\begin{equation}  
\label{e.caccio.S.intermed}
\|h\|_{L^2(\cu)} \le \|\nabla h\|_{H^{-1}_\per(\cu)} \le C \|\g\|_{\hat H^{-1}(U)}.
\end{equation}
Moreover, the function $h$ is harmonic on $V''$. By the Caccioppoli inequality for harmonic functions, we deduce that 
\begin{equation}  
\label{e.use.some.standard.caccio}
\|\nabla h\|_{L^2(V')} \le C \|h\|_{L^2(V'')} \le C \|\g\|_{\hat H^{-1}(U)},
\end{equation}
which is \eqref{e.caccio.S.est.h}. The proofs of \eqref{e.caccio.S.est.T} and \eqref{e.caccio.S.est.u} are very similar to the argument for \eqref{e.caccio.S.intermed}, so we omit the details.

\smallskip

\emph{Step 3.} We now conclude the argument.
Note first that, since $\eta' \equiv 1$ on the support of $\eta$, we have
\begin{align*}  
C^{-1} \int_U \eta^2 \Ll( |\nabla u|^2 + |\g|^2 \Rr) 
& \le 
\int_U 
\eta^2\begin{pmatrix}
\nabla u \\ \g
\end{pmatrix}
\cdot \Abf 
\begin{pmatrix}
\nabla u \\ \g
\end{pmatrix}
\\
& = 
\int_U 
\eta^2\begin{pmatrix}
\nabla u \\ \nabla h + \nabla \cdot \mathbf T + (\eta'\g)_\cu
\end{pmatrix}
\cdot \Abf 
\begin{pmatrix}
\nabla u \\ \g
\end{pmatrix}
.
\end{align*}
Moreover, since $\mathbf T$ is skew-symmetric, we have
\begin{equation*}  
\begin{pmatrix}
\nabla (\eta^2 u) \\ \nabla \cdot (\eta^2 \mathbf T)
\end{pmatrix}
\in \Lpoot(U)\times\Lso(U),
\end{equation*}
and the definition of the space $\S(U)$ implies
\begin{align*}  
0 & = 
\int_U 
\begin{pmatrix}
\nabla (\eta^2 u) \\ \nabla \cdot (\eta^2 \mathbf T)
\end{pmatrix}
\cdot \Abf 
\begin{pmatrix}
\nabla u \\ \g
\end{pmatrix}
\\
& = 
\int_U \eta^2
\begin{pmatrix}
\nabla u \\ \nabla \cdot \mathbf T
\end{pmatrix}
\cdot \Abf 
\begin{pmatrix}
\nabla u \\ \g
\end{pmatrix}
+ 2 \int_U \eta 
\begin{pmatrix}
u \nabla \eta \\ \mathbf T \nabla \eta
\end{pmatrix}
\cdot \Abf 
\begin{pmatrix}
\nabla u \\ \g
\end{pmatrix}
.
\end{align*}
Combining the last displays yields
\begin{align*}  
\int_U \eta^2 \Ll( |\nabla u|^2 + |\g|^2 \Rr) & \le C \int_U \eta 
\begin{pmatrix}
-2u \nabla \eta \\  \eta \nabla h - 2 \mathbf T \nabla \eta + \eta (\eta'\g)_\cu
\end{pmatrix}
\cdot \Abf 
\begin{pmatrix}
\nabla u \\ \g
\end{pmatrix}
.
\end{align*}
By H\"older's inequality, \eqref{e.caccio.S.est.h}, \eqref{e.caccio.S.est.T}, \eqref{e.caccio.S.est.u} and $|(\eta'\g)_\cu| \le C \|\g\|_{\hat H^{-1}(U)}$, we obtain
\begin{equation*}  
\Ll(\int_U \eta^2 \Ll( |\nabla u|^2 + |\g|^2 \Rr) \Rr)^\frac 1 2\le C \Ll( \|\nabla u\|_{\hat H^{-1}(U)} + \|\g\|_{\hat H^{-1}(U)} \Rr) .
\end{equation*}
Since $\eta \equiv 1$ on $V$, this completes the proof.
\end{proof}

We show in the next lemma that the space $\S(U)$ can alternatively be described in terms of solutions of the original equation as well as its adjoint (although we remark that this is not used until Section~\ref{s.FitzDP}). We denote by $\A^*(U)\subseteq H^1(U)$ the set of weak solutions of the adjoint equation
\begin{equation}
\label{e.adjoint.equation}
-\nabla \cdot \left( \a^t \nabla u^* \right) = 0 \quad \mbox{in} \ U, 
\end{equation}
that is, 
\begin{equation*} \label{}
\A^*(U):= \left\{ u^* \in H^1(U) \,:\, \forall w\in H^1_0(U), \  \int_U \nabla u^* \cdot \a\nabla w = 0  \right\}.
\end{equation*}
\begin{lemma}
\label{l.whatisS}
We have
\begin{equation}  
\label{e.thisisS}
\S(U) = \left\{ \left( \nabla u + \nabla u^* , \a\nabla u - \a^t\nabla u^* \right)  \,:\, 
u \in \A(U), \, u^* \in \A^*(U)
\right\}.
\end{equation}
\end{lemma}
\begin{proof}
Denote by $\S'(U)$ the set on the right side of \eqref{e.thisisS}. By definition, an element $S$ of $\Lpot(U) \times \Ls(U)$ belongs to $\S(U)$ if and only if
\begin{equation*}  
\text{for every } S' \in \Lpoot(U) \times \Lso(U), \qquad \int_U S \cdot \Abf S' = 0.
\end{equation*}
In more explicit notation, a pair $(\nabla v, \g) \in \Lpot(U) \times \Ls(U)$ belongs to $\S(U)$ if and only if
\begin{multline}  
\label{e.cond.explicit.S}
\text{for every } \phi \in H^1_0(U) \text{ and } \f \in \Lso(U), \\
\int_U \Ll[ \nabla \phi \cdot \s  \nabla v + (\f - \m \nabla \phi) \cdot \s ^{-1} (\g - \m \nabla v) \Rr] = 0.
\end{multline}
This implies in particular that 
\begin{equation*}  
\text{for every } \f \in \Lso(U), \qquad \int_U \f \cdot \s ^{-1} (\g - \m  \nabla v) = 0.
\end{equation*}
The latter condition is equivalent to the existence of $w \in H^1(U)$ such that 
\begin{equation*}  
\s ^{-1} (\g - \m \nabla v) = \nabla w. 
\end{equation*}
We deduce that \eqref{e.cond.explicit.S} is equivalent to the existence of~$w \in H^1(U)$ such that
\begin{equation}  
\label{e.getting.to.S}
\Ll\{
\begin{aligned}  
& \g = \s  \nabla w + \m \nabla v, 
\quad \mbox{and} \ \ 
\\
& \forall \phi \in H^1_0(U), \quad \int_U \nabla \phi \cdot \Ll(\s  \nabla v + \m \nabla w\Rr) = 0.
\end{aligned}
\Rr.
\end{equation}
The last line above is the weak formulation of the equation
\begin{equation*}  
-\nabla \cdot \Ll( \s  \nabla v + \m \nabla w \Rr)  = 0 \qquad \text{ in } U,
\end{equation*}
and we recall that, since $\g \in \Ls(U)$, the first line of \eqref{e.getting.to.S} ensures that
\begin{equation*}  
-\nabla \cdot \Ll( \s  \nabla w + \m \nabla v \Rr)  = 0 \qquad \text{ in } U.
\end{equation*}
Whenever the two displays above are satisfied, we may define
\begin{equation*}  
u := \frac 1 2 (v + w), \qquad u^* := \frac 1 2 (v - w),
\end{equation*}
and verify that $u \in \A(U)$, $u^* \in \A^*(U)$, and 
\begin{equation*}  
\g = \a \nabla u - \a^t \nabla u^*,
\end{equation*}
thereby proving that $\S(U) \subset \S'(U)$.
Conversely, given $u \in \A(U)$ and $u^* \in \A^*(U)$, if we set 
\begin{equation*}  
v := u + u^*, \qquad w := u - u^*, \qquad \g := \a \nabla u - \a^t \nabla u^* = \s  \nabla w + \m \nabla v,
\end{equation*}
then we can verify that the conditions in \eqref{e.getting.to.S} are satisfied, since 
\begin{equation*}  
\s  \nabla v + \m \nabla w = \a \nabla u + \a^t \nabla u^* \in \Ls(U).
\end{equation*}
This proves the converse inclusion $\S'(U) \subset \S(U)$, and hence completes the proof that $\S(U) = \S'(U)$.
\end{proof}

We now give the analogue of Lemma~\ref{l.basicJ} for our double-variable quantity.

\begin{lemma}[Basic properties of $J$]
\label{l.basicJ.bigJ}
Fix a bounded Lipschitz domain $U\subseteq \Rd$. 
The quantity $J(U,X,X^*)$ and its maximizer $S(\cdot,U,X,X^*)$ satisfy the following:

\begin{itemize}

\item \emph{Representation as quadratic form} The mapping $(X,X^*) \mapsto J(U,X,X^*)$ is a quadratic form and there exist $2d$-by-$2d$ matrices $\Abf(U)$ and $\Abf_*(U)$ and a constant $C(\Lambda) < \infty$ such that 
\begin{equation}  
\label{e.JunifconvC11p.bigJ}
C^{-1} \, \Id \le \Abf_*(U) \le \Abf(U) \le C \, \Id
\end{equation}
and 
\begin{equation}  
\label{e.bigJ.matrix}
J(U,X,X^*) = \frac 1 2 X \cdot \Abf(U) X + \frac 1 2 X^* \cdot \Abf_*^{-1}(U) X^* - X \cdot X^*.
\end{equation}
These matrices are characterized by the fact that, for every $X, X^* \in \R^{2d}$,
\begin{equation}
\label{e.Jderp.bigJ}
\Abf(U) X = - \fint_U \Abf S(\cdot,U,X,X^*),
\end{equation}
and
\begin{equation} 
\label{e.Jderq.bigJ}
\Abf_*^{-1}(U) X^*  =  \fint_U  S(\cdot,U,X,X^*).
\end{equation}

\item \emph{Subadditivity.} Let $U_1, \ldots, U_N \subset U$ be bounded Lipschitz domains that form a partition of~$U$, in the sense that $U_i \cap U_j = \emptyset$ if $i\neq j$ and 
\begin{equation*} \label{}
\left| U \setminus \bigcup_{i=1}^N U_i \right| = 0. 
\end{equation*}
Then, for every $X, X^*\in\R^{2d}$,
\begin{equation}
\label{e.Jsubadd.bigJ}
J(U,X,X^*) \leq \sum_{i=1}^N \frac{\left|U_i\right|}{|U|} J(U_i,X,X^*). 
\end{equation}

\item \emph{First variation for $J$.} For $X,X^*\in\R^{2d}$, the function $S(\cdot,U,X,X^*)$ is characterized as the unique element of $\S(U)$ which satisfies 
\begin{equation}
\label{e.Jfirstvar.bigJ} 
\fint_U  T \cdot \Abf S(\cdot,U,X,X^*) = \fint_U \left( - X\cdot \Abf T + X^*\cdot T \right),
\quad 
\forall\, T \in \S(U). 
\end{equation}

\item \emph{Quadratic response.} For every $X,X^*\in\R^{2d}$ and $T\in \S(U)$, 
\begin{multline} 
\label{e.Jquadresponse.bigJ}
 \fint_U \frac12 \left(   T -   S(\cdot,U,X,X^*) \right) \cdot \Abf\left(   T -   S(\cdot,U,X,X^*) \right) 
 \\
=
J(U,X,X^*) - \fint_U \left( -\frac12  T \cdot \Abf  T  -X\cdot\Abf   T + X^*\cdot   T  \right).
\end{multline}
\end{itemize}
\end{lemma}
\begin{proof}
The proof is essentially identical to that of Lemma~\ref{l.basicJ} and so is omitted. 
\end{proof}

\section{Convergence of the subadditive quantities}

In this section, we state and prove the main result of this chapter on the quantitative convergence of the double-variable master quantity~$J$. In order to state it, we must first define the effective matrix which describes the limit of $\mu(U,X)$.  We take it to be the matrix $\Abfh \in \R^{2d\times 2d}$ which satisfies, for every $X \in \R^{2d}$,
\begin{equation}
\label{e.def.bigAbar}
\frac12 X \cdot \Abfh  X = \lim_{n\to \infty} \E \left[ \mu(\cu_n,X) \right]. 
\end{equation}
The right side of~\eqref{e.def.bigAbar} exists and defines a positive quadratic form by the monotonicity of the sequence~$n\mapsto \E \left[ \mu(\cu_n,X) \right]$, which follows from subadditivity and stationarity (in analogy to Chapter~\ref{c.one}) and Lemma~\ref{l.basicJ.bigJ}. Using the matrix notation introduced in this lemma gives that
\begin{equation*}  
\Abfh = \lim_{n \to \infty} \E \Ll[ \Abf(\cu_n) \Rr] ,
\end{equation*}
and, by~\eqref{e.JunifconvC11p.bigJ},
\begin{equation} 
\label{e.Abfhbounds}
C^{-1} \, \Id \leq \Abfh \leq C\, \Id. 
\end{equation}
We define, for every $X, X^* \in \R^{2d}$,
\begin{equation*} \label{}
\overline{J}(X,X^*):= 
\frac12 X \cdot \Abfh  X +  \frac12 X^* \cdot \Abfh ^{-1} X^* - X \cdot X^*
\end{equation*}

The following result is the analogue of Theorem~\ref{t.subadd} for nonsymmetric coefficients and. By Lemma~\ref{l.mu.when.sym}, it is equivalent to Theorem~\ref{t.subadd} in the symmetric case. 

\begin{theorem}
\label{t.subadd.bigJ}
Fix $s\in (0,d)$. 
There exist $\alpha(d,\Lambda)\in \left(0,\frac12\right]$ and $C(s,d,\Lambda)<\infty$ such that, for every $X,X^*\in B_1$ and $n\in\N$,  
\begin{equation}
\label{e.subadderror.bigJ}
 \left| 
 J(\cu_n,X,X^*)
 -
\overline{J}(X,X^*)
  \right| 
\leq C 3^{-n \alpha(d-s)} + \O_1\left( C3^{-ns} \right). 
\end{equation}
\end{theorem}

As in Chapter~\ref{c.two}, the strategy here consists of showing that $\E\left[ J\left(\cu_n,X,\Abfh  X\right) \right]$ decays like a negative power of the length scale~$3^n$ by an iteration of the scales. This then allows us to use the following analogue of Lemma~\ref{l.minimalset} to control~$J(\cu_n,X,X^*)$ for general $(X,X^*)$.

\begin{lemma}
\label{l.minimalset.bigJ}
Fix $\Gamma\geq 1$. There exists a constant~$C(\Gamma,d,\Lambda)<\infty$ such that, for every symmetric matrix $\tilde{\Abf} \in \R^{2d\times 2d}$ satisfying
\begin{equation*} \label{}
\Gamma^{-1} I_{2d} \leq \tilde{\Abf} \leq \Gamma I_{2d}
\end{equation*}
and every bounded Lipschitz domain $U\subseteq\Rd$, we have 
\begin{equation}
\label{e.minimalset.bigJ}
\Ll| \Abf(U) - \tilde \Abf \Rr| + \Ll| \Abf_*(U) - \tilde \Abf \Rr| 
 \leq C \sup_{X  \in B_1} 
 \Ll(J \left(U,X, \tilde{\Abf} X \right) \Rr)^\frac 12 . 
\end{equation}
\end{lemma}
\begin{proof}
The proof of this lemma is almost identical to that of Lemma~\ref{l.minimalset}, with only minor notational differences. 
\end{proof}

We prove Theorem~\ref{t.subadd.bigJ} by following the same steps as in the proof of Theorem~\ref{t.subadd}. Several of the lemmas follow from arguments almost identical to the ones in Chapter~\ref{c.two}, and for these we will omit the details and refer the reader to the appropriate argument. 

\smallskip

We define the subadditivity defect by 
\begin{equation}
\label{e.def.taun.bigJ}
\tau_n := \sup_{X, X^* \in B_1} \Ll( \E[J(\cu_n,X, X^*)]  - \E[J(\cu_{n+1},X, X^*)] \Rr) .
\end{equation}

\begin{lemma}
\label{l.quadresponse.bigJ}
Fix a bounded Lipschitz domain $U\subseteq \Rd$ and let $\{U_1,\ldots,U_k\}$ be a partition of $U$ into Lipschitz subdomains. Then, for every $X,X^*\in \R^{2d}$,
\begin{multline*} \label{}
\sum_{j=1}^k  \frac{|U_j|}{|U|} \frac12 \left\| \Abf^{\frac12} \left(  S (\cdot,U,X,X^*) - S(\cdot,U_j,X,X^*) \right) \right\|_{\underline{L}^2(U_j)}^2 \\
= \sum_{j=1}^k \frac{|U_j|}{|U|}\left( J(U_j,X,X^*) -  J(U,X,X^*) \right).
\end{multline*}
\end{lemma}
\begin{proof}
See the proof of Lemma~\ref{l.quadresponse}. 
\end{proof}

The following lemma gives us control of the spatial averages of~$S(\cdot,\cu_n,X,X^*)$ and $\Abf(\cdot)S(\cdot,\cu_n,X,X^*)$, analogously to Lemma~\ref{l.spatavg}. 

\begin{lemma}
\label{l.spatavg.bigJ}
There exist $\kappa(d)>0$ and~$C(d,\Lambda)<\infty$ such that, for every $m\in\N$ and $X,X^*\in B_1$, 
\begin{equation*}
\var \left[  \fint_{\cu_m} S(\cdot,\cu_{m},X,X^*) \right]  
\leq 
C3^{-m\kappa} + C \sum_{n=0}^m 3^{-\kappa(m-n)} \tau_n.
\end{equation*}
\end{lemma}
\begin{proof}
See the proof of Lemma~\ref{l.spatavg}.
\end{proof}

We next proceed to the definition of the coarsened matrix $\Abfh _U\in \R^{2d\times 2d}$, which is the analogue of the matrix $\ahom_U$ introduced in Definition~\ref{def.ahomU}. 
\begin{definition}
\label{def.ahomU.bigJ}
We define the symmetric matrix $\Abfh_{U} \in \R^{2d\times 2d}$ by
\begin{equation*} \label{}
\Abfh_U := \E \Ll[ \Abf^{-1}(U) \Rr] ^{-1}.
\end{equation*}
We also denote $\Abfh_n:= \Abfh_{\cu_n}$ for short.
\end{definition}
By Lemma~\ref{l.basicJ.bigJ}, there exists a constant $C(\Lambda) < \infty$ such that 
\begin{equation} 
\label{e.ahomUbounds.bigJ}
C^{-1} \, \Id \leq \Abfh_U \leq C \, \Id,
\end{equation}
and moreover,
\begin{equation*}  
\Abfh_U^{-1} X^* = \E \Ll[ \fint_U S(\cdot,U,0,X^*) \Rr] .
\end{equation*}
We now proceed with the double-variable version of Lemma~\ref{l.iterstep}.
\begin{lemma}
\label{l.iterstep.bigJ}
There exist~$\kappa(d)>0$ and~$C(d,\Lambda)<\infty$ such that, for every $n\in\N$ and $X \in B_1$,
\begin{equation}
\label{e.itersteprealz.bigJ}
\E \left[  J(\cu_n, X,\Abfh _n X)\right]  
\leq 
C3^{-n\kappa} + C \sum_{m=0}^n 3^{-\kappa(n-m)} \tau_m.
\end{equation}
\end{lemma}
\begin{proof}
Fix $n\in\N$ and $X\in B_1$. 
We first observe that, by the analogue of~\eqref{e.firstvarJ2}, 
\begin{equation*} \label{}
J(\cu_n, X,\Abfh _n X)  
\leq 
C\fint_{\cu_n} \left| S(\cdot,\cu_n, X,\Abfh _n X) \right|^2.
\end{equation*}
Therefore it suffices to prove that 
\begin{equation} 
\label{e.iterstep.bigJ.wts}
\E \left[ \fint_{\cu_n} \left| S(\cdot,\cu_n, X,\Abfh _n X) \right|^2 \right] 
\leq 
C3^{-n\kappa} + C \sum_{m=0}^n 3^{-\kappa(n-m)} \tau_m.
\end{equation}
Applying the multiscale Poincar\'e inequality (Proposition~\ref{p.mspoin}), we get
\begin{align*}  
& \Ll\| S(\cdot,\cu_{n+1},X,X^*)-  \Abfh_n^{-1} X^* + X \Rr\|_{\Hminus(\cu_{n+1})}^2 \\
& \quad \le  C  \fint_{\cu_{n+1}} \left| S(\cdot,\cu_{n+1},X,X^*) -  \Abfh_n^{-1} X^* + X \right|^2 \notag \\
& \qquad \quad
+C\left(  \sum_{m=0}^n 3^m \left( |Z_{m}|^{-1}\sum_{y\in Z_{m}} \left| \fint_{y+\cu_m} S(\cdot,\cu_{n+1},X,X^*)  -  \Abfh_n^{-1}X^* + X \right|^2 \right)^{\frac12}  \right)^2. 
\end{align*}
The right side above can be compared with that of \eqref{e.multiscaleapp}. Following then Steps~2 and~3 of the proof of Lemma~\ref{l.flatness}, we obtain that
\begin{multline*}  
\E \Ll[\Ll\| S(\cdot,\cu_{n+1},X,X^*) -  \Abfh_n^{-1}X^* + X\Rr\|_{\Hminus(\cu_{n+1})}^2\Rr] 
\\
\leq 
C 3^{2n} \left( C3^{-n\kappa} + C \sum_{m=0}^n 3^{-\kappa(n-m)} \tau_m \right).
\end{multline*}
According to Lemma~\ref{l.cacciop.S}, 
\begin{multline*} \label{}
\Ll\| S(\cdot,\cu_{n+1},X,X^*) -  \Abfh_n^{-1}X^* + X\Rr\|_{\underline{L}^2(\cu_n)}
\\
\leq C 3^{-n}  \Ll\| S(\cdot,\cu_{n+1},X,X^*) -  \Abfh_n^{-1}X^* + X\Rr\|_{\Hminus(\cu_{n+1})}.
\end{multline*}
We therefore deduce that
\begin{equation*}  
\E \Ll[\Ll\| S(\cdot,\cu_{n+1},X,X^*) -  \Abfh_n^{-1}X^* + X\Rr\|_{\underline{L}^2(\cu_n)}^2\Rr]
\leq
C \left( C3^{-n\kappa} + C \sum_{m=0}^n 3^{-\kappa(n-m)} \tau_m \right),
\end{equation*}
and by quadratic response and stationarity,
\begin{equation*}  
\E \Ll[\Ll\| S(\cdot,\cu_{n+1},X,X^*) - S(\cdot,\cu_{n},X,X^*) \Rr\|_{\underline{L}^2(\cu_n)}^2\Rr] \le C \tau_n.
\end{equation*}
The last two displays imply \eqref{e.iterstep.bigJ.wts}, so the proof is complete.
\end{proof}

We now complete the proof of Theorem~\ref{t.subadd.bigJ}. 

\begin{proof}[{Proof of Theorem~\ref{t.subadd.bigJ}}]
Applying Lemma~\ref{l.iterstep.bigJ} iteratively as in the proof of Proposition~\ref{p.subaddE}, we obtain an exponent $\alpha(d,\Lambda) > 0$ and a constant $C(d,\Lambda) < \infty$ such that, for every $n \in \N$ and $X \in B_1$,
\begin{equation*}  
\E \left[ J(\cu_n,X,\Abfh  X)\right] \leq C3^{-n\alpha}. 
\end{equation*}
We use this estimate and Lemma~\ref{l.barO.boxes} to control the stochastic fluctuations of $J(\cu_n, X, \Abfh X)$, and then conclude using Lemma~\ref{l.minimalset.bigJ}, just as in the proof of Theorem~\ref{t.subadd}.
\end{proof}

\index{homogenized coefficients}

We conclude this section by identifying the effective matrix~$\ahom$. We need to show that~$\Abfh$ is a representative in the sense of Definition~\ref{def.variationalrep} of some linear, uniformly monotone mapping~$p\mapsto \ahom p$. This is accomplished in the following lemma, which gives us our definition of~$\ahom$. 

\begin{lemma}
\label{l.werepresent} 
There exist $C(d,\Lambda)<\infty$ and a matrix $\ahom \in \R^{d\times d}$ satisfying
\begin{equation} 
\label{e.ahombounds.bigJ}
\forall \xi \in \Rd, \quad \left| \xi \right|^2 \leq \xi\cdot \ahom \xi
 \quad \mbox{and} \quad 
\Ll|\ahom \xi \Rr| \leq \Lambda \left| \xi \right|
\end{equation}
such that the mapping
\begin{equation}
\label{e.mappingAbar}
\bar A : (p,q) \mapsto 
\frac 12 \begin{pmatrix} p \\ q  \end{pmatrix} \cdot \Abfh \begin{pmatrix} p \\ q  \end{pmatrix} 
\end{equation}
is a variational representative of $p\mapsto \ahom p$. 
\end{lemma}
\begin{proof}
The fact that $(p,q) \mapsto \overline{A}(p,q)$ is uniformly convex is clear from~\eqref{e.JunifconvC11p.bigJ} and~\eqref{e.def.bigAbar}. 
It is also evident from~\eqref{e.fitz.positivity} that 
\begin{equation} 
\label{e.lwbd.mu}
\mu\left(\cu_n, p,q \right) \geq p\cdot q,
\end{equation}
and thus~\eqref{e.def.bigAbar} yields that
\begin{equation} 
\label{e.primal.lower.bound}
\overline{A}(p,q) \geq p\cdot q. 
\end{equation}
It remains to check that, for each $p\in\Rd$,  
\begin{equation}
\label{e.finalcheck.bigJ}
\inf_{q\in\Rd} \left( \overline{A}(p,q) - p\cdot q \right) = 0
\end{equation}
and that the mapping which sends $p$ to the~$q$ achieving the infimum is linear. The linearity is obvious from the fact that $\overline{A}$ is quadratic, which is itself immediate from~\eqref{e.def.bigAbar} and the fact that $(p,q) \mapsto \E \left[\mu(\cu_n,p,q)\right]$ is quadratic. The bounds~\eqref{e.ahombounds.bigJ} follow from~\eqref{e.Abfhbounds}. 

\smallskip

We decompose the proof of~\eqref{e.finalcheck.bigJ} into two steps.

\smallskip

\emph{Step 1.} We first show that 
\begin{equation}
\label{e.dual.lower.bound}
\frac 1 2  \begin{pmatrix} q^* \\ p^*  \end{pmatrix} \Abfh^{-1}  \begin{pmatrix} q^* \\ p^*  \end{pmatrix}  \geq p^*\cdot q^*. 
\end{equation}
Recall that, in view of Proposition~\ref{p.yesvariational} and the solvability of the Dirichlet problem, we have, for every bounded Lipschitz domain $U\subseteq\Rd$, 
\begin{equation*}  
\inf \left\{  \int_U \left( A(\nabla u,\g,\cdot) - \g\cdot \nabla u \right) \,:\, u \in \ell_{p^*} + H^1_0(U), \ \g \in \Ls(U) \right\} = 0.
\end{equation*}
We deduce that, for every $n \in \N$ and $q^*, p^* \in \Rd$,
\begin{equation*}  
\mu^*(\cu_n,q^*, p^*) \ge q^* \cdot p^*. 
\end{equation*}
By Lemma~\ref{l.Jsplitting.bigJ} and Theorem~\ref{t.subadd.bigJ}, we have
\begin{equation*}  
\frac 1 2  \begin{pmatrix} q^* \\ p^*  \end{pmatrix} \Abfh^{-1}  \begin{pmatrix} q^* \\ p^*  \end{pmatrix}  = \lim_{n \to \infty} \E[\mu^*(\cu_n,q^*,p^*)],
\end{equation*}
and therefore \eqref{e.dual.lower.bound} follows.

\smallskip

\emph{Step 2.} 
Fix $p \in \Rd$, select $q \in \Rd$ achieving the infimum in \eqref{e.finalcheck.bigJ}, so that
\begin{equation*}  
 \Abfh \begin{pmatrix} p \\ q  \end{pmatrix} = \begin{pmatrix} \ast \\ p  \end{pmatrix},
\end{equation*}
and then let $q^* \in \Rd$ be such that 
\begin{equation}  
\label{e.define.q*}
\Abfh \begin{pmatrix} p \\ q  \end{pmatrix} = \begin{pmatrix} q^* \\ p  \end{pmatrix}.
\end{equation}
We then have that 
\begin{align*}  
\frac 1 2  \begin{pmatrix} q^* \\ p  \end{pmatrix} \Abfh^{-1}  \begin{pmatrix} q^* \\ p  \end{pmatrix} & = \sup_{p',q' \in \Rd} \Ll(  \begin{pmatrix} q^* \\ p  \end{pmatrix} \cdot  \begin{pmatrix} p' \\ q'  \end{pmatrix} - \frac 12 \begin{pmatrix} p' \\ q'  \end{pmatrix} \cdot \Abfh \begin{pmatrix} p' \\ q'  \end{pmatrix}  \Rr)  \\
& = \begin{pmatrix} q^* \\ p  \end{pmatrix} \cdot  \begin{pmatrix} p \\ q  \end{pmatrix} - \frac 12 \begin{pmatrix} p \\ q  \end{pmatrix} \cdot \Abfh \begin{pmatrix} p \\ q  \end{pmatrix}  
\end{align*}
where we used \eqref{e.define.q*} for the second equality. Together with~\eqref{e.primal.lower.bound} and \eqref{e.dual.lower.bound}, this implies that 
\begin{equation*}  
\bar A (p,q) = p\cdot q,
\end{equation*}
and therefore completes the proof of \eqref{e.finalcheck.bigJ}.
\end{proof}

We conclude this section with the reassuring remark that the identification of~$\ahom$ by Lemma~\ref{l.werepresent} does not conflict, in the symmetric case, with the identification of~$\ahom$ in Chapter~\ref{c.two}. This is immediate from Lemma~\ref{l.mu.when.sym}. 

\begin{exercise}
Show that homogenization commutes with the map $\a \mapsto \a^t$.
\end{exercise}

\section{Quantitative homogenization of the Dirichlet problem}
\label{s.FitzDP}

In this section we give a generalization of Theorem~\ref{t.DP}, demonstrating that the subadditive quantities~$\mu$ and~$\mu^*$ give a quantitative control of the homogenization error for the Dirichlet problem. 

\smallskip

We begin by quantifying the weak convergence of $S$ and $\Abf S$ to constants.  The following proposition is an analogue of~Theorem~\ref{t.Jmaximizers}. We also use the notation 
\begin{equation*} \label{}
\overline{S}(X,X^*) := X - \Abfh^{-1}X^*. 
\end{equation*}

\begin{proposition}[{Weak convergence of $(S,\Abf S)$}]
\label{p.controlofS}
There exist $\alpha(d,\Lambda)\in \left(0,\frac1d\right]$ and $C(s,d,\Lambda)<\infty$ such that, for every $X,X^*\in B_1$ and $n\in\N$,
\begin{multline} 
\label{e.HminusoneS}
\left\| S(\cdot,\cu_m,X,X^*) - \overline{S}(X,X^*) \right\|_{\Hminus(\cu_m)}^2
+
\left\| \Abf S(\cdot,\cu_m,X,X^*) - \Abfh \overline{S}(X,X^*) \right\|_{\Hminus(\cu_m)}^2
\\
\leq
C 3^{m -m \alpha(d-s)} + \O_1\left( C3^{m-ms} \right).
\end{multline}
\end{proposition}
\begin{proof}
See Step~1 of the proof of Proposition~\ref{p.qualweakconv} and the proof of Theorem~\ref{t.Jmaximizers}. 
\end{proof}

We next give the construction of (finite-volume) correctors $\phi_{m,e}$ which play the same role as the functions defined in~\eqref{e.FVC} in the proof of Theorem~\ref{t.DP.blackbox}. These functions will be obtained directly from~$S(\cdot,X,X^*)$ and~$\Abf S(\cdot,X,X^*)$ and Lemma~\ref{l.whatisS}. The estimates we need for $\phi_{m,e}$ to run the proof of Theorem~\ref{t.DP.blackbox} will be immediate from~\eqref{e.HminusoneS}.

\smallskip

We define $\phi_{e,m}$ from the maximizers of $J(\cu_m,X,0)$ for an appropriate choice of $X$, depending on $e$. The selection of~$X$ is a linear algebra exercise using Lemma~\ref{l.werepresent}. We set
\begin{equation*} \label{}
X_e:= 
- \begin{pmatrix} e \\ \ahom e \end{pmatrix}
\end{equation*}
and observe that we have 
\begin{equation} 
\label{e.freedom}
\Abfh X_e = 
- \begin{pmatrix} \ahom e \\ e \end{pmatrix}.
\end{equation}
To confirm~\eqref{e.freedom}, we note (see Lemma~\ref{l.werepresent}) that 
\begin{equation*} \label{}
q \mapsto \frac12 \begin{pmatrix} e \\ q \end{pmatrix} \cdot \Abfh \begin{pmatrix} e \\ q \end{pmatrix} - e\cdot q
\quad \mbox{attains its minimum at} \ q=\ahom e
\end{equation*}
and
\begin{equation*} \label{}
p \mapsto \frac12 \begin{pmatrix} p \\ \ahom e \end{pmatrix} \cdot \Abfh \begin{pmatrix} p \\ \ahom e \end{pmatrix} - p\cdot \ahom e
\quad \mbox{attains its minimum at} \ p= e.
\end{equation*}
Differentiating in $q$ and $p$, respectively, gives~\eqref{e.freedom}. 
%
%
%

\smallskip


We next take~$u_{e,n}\in H^1_\pa(\cu_{n+1})$ to be the element $u\in \mathcal{A}(\cu_{n+1})$ in the representation of~$S(\cdot,\cu_{n+1},X_e,0)$ given in Lemma~\ref{l.whatisS}, with additive constant chosen so that~$\left(u_{e,n}\right)_{\cu_n} = 0$. Equivalently, we can define~$u_{e,n}$ to be the function on $\cu_{n+1}$ with mean zero on~$\cu_n$ with gradient given by
\begin{equation} 
\label{e.deconformugrad}
\nabla u_{e,n} = \frac12 \left( \pi_1 S\left( \cdot,\cu_{n+1},X_e,0\right) + \pi_2 \Abf S\left( \cdot,\cu_{n+1},X_e,0 \right) \right),
\end{equation}
where $\pi_1$ and $\pi_2$ denote the projections $\R^{2d}\to \Rd$ given by $\pi_1(x,y) = x$ and $\pi_2(x,y)=y$ for $x,y\in\Rd$. Note that, by Lemma~\ref{l.whatisS}, we also have
\begin{equation}
\label{e.deconformuflux}
\a \nabla u_{e,n} =  \frac12 \left( \pi_2 S\left( \cdot,\cu_{n+1},X_e,0\right) + \pi_1 \Abf S\left( \cdot,\cu_{n+1},X_e,0 \right) \right).
\end{equation}
By Proposition~\ref{p.controlofS},~\eqref{e.freedom},~\eqref{e.deconformugrad} and~\eqref{e.deconformuflux}, we have that 
\begin{multline} 
\label{e.weakucontrol}
3^{-n} \left\| \nabla u_{e,n} - e \right\|_{\Hminus(\cu_{n+1})} 
+ 
3^{-n} \left\| \a \nabla u_{e,n} - \ahom e \right\|_{\Hminus(\cu_{n+1})} 
\\
\leq
C 3^{-n \alpha(d-s)} + \O_1\left( C3^{-ns} \right).
\end{multline}
As~$u_{e,n}\in \mathcal{A}(\cu_{n+1})$, we have that~$u_{e,n}$ is a solution of 
\begin{equation} 
\label{e.eqnforuen}
 - \nabla \cdot \left( \a \nabla u_{e,n}\right) = 0 \quad \mbox{in} \ \cu_{n+1}.
\end{equation}
The approximate first-order corrector $\phi_{e,n}$ is defined by subtracting the affine function $x\mapsto e\cdot x$ from $u_{e,n}$:
\begin{equation*} \label{}
\phi_{e,n}(x) := u_{e,n}(x) - e\cdot x. 
\end{equation*}
Summarizing, we therefore have that $\phi_{e,n}$ is a solution of 
\begin{equation} 
\label{e.eq.phien}
 - \nabla \cdot \left( \a \left(e+\nabla \phi_{e,n}\right) \right) = 0 \quad \mbox{in} \ \cu_{n+1},
\end{equation}
and satisfies the estimates
\begin{multline} 
\label{e.weakphicontrol}
3^{-n} \left( \left\| \nabla \phi_{e,n} \right\|_{\Hminus(\cu_{n+1})} 
+ 
 \left\| \a \left(e+ \nabla \phi_{e,n} \right)- \ahom e \right\|_{\Hminus(\cu_{n+1})} 
 \right)
\\
\leq
C 3^{-n \alpha(d-s)} + \O_1\left( C3^{-ns} \right).
\end{multline}
By the previous two displays, Lemma~\ref{l.integrate.H-1} and $\left( \phi_{e,n} \right)_{\cu_n} = 0$, we also have
\begin{equation} 
\label{e.L2phicontrol}
3^{-n} \left\|  \phi_{e,n} \right\|_{\underline{L}^2(\cu_n)} 
\leq
C 3^{-n \alpha(d-s)} + \O_1\left( C3^{-ns} \right).
\end{equation}

\smallskip

We next recall the definition of the random variable~$\mathcal{E}'(\ep)$ from~\eqref{e.E'm}:
\begin{equation*} 
\mathcal{E}'(\ep):= \sum_{k=1}^d \left( \ep \left\|  \phi_{m,e_k}\left( \tfrac \cdot \ep \right) \right\|_{L^2(\ep\cu_m)}
+
\left\| \a\left( \tfrac \cdot\ep\right)\left( e_k + \nabla \phi_{m,e_k}\left(\tfrac\cdot\ep\right) \right) - \ahom e_k \right\|_{H^{-1}(\ep\cu_m)}  \right)^2,
\end{equation*}
where $m=m_\ep:= \left\lfloor \left| \log \ep \right| / \log 3 \right\rfloor$. The estimates~\eqref{e.weakphicontrol} and~\eqref{e.L2phicontrol} above immediately yield that 
\begin{equation*} \label{}
\mathcal{E}'(\ep) \leq C\ep^{\alpha(d-s)} + \O_1\left( C \ep^s \right). 
\end{equation*}
Theorem~\ref{t.DP.blackbox} now applies---the reader will notice that the proof makes no use of the symmetry of $\a(\cdot)$---and combining this argument with the more routine bookkeeping in the proof of Theorem~\ref{t.DP}, we obtain the following analogue of the latter.

\begin{theorem}
\label{t.DP.nosymm}
Fix $s\in (0,d)$, an exponent $\delta>0$ and a bounded Lipschitz domain $U\subseteq B_1$. There exist~$\beta(\delta,d,\Lambda)>0$,~$C(s,U,\delta,d,\Lambda)<\infty$ and a random variable~$\X_s$ satisfying
\begin{equation}
\label{e.sizeX.nosymm}
\X_s = \O_1\left(C \right)
\end{equation}
such that the following statement holds. For each $\ep\in (0,1]$ and $f\in W^{1,2+\delta}(U)$, let $u^\ep, u \in f+H^1_0(U)$ respectively denote the solutions of the Dirichlet problems 
\begin{equation}
\label{e.uep.and.ubar.nosymm}
\left\{
\begin{aligned}
& -\nabla \cdot \left( \a\left(\tfrac x\ep\right) \nabla u^\ep \right) = 0 &  \mbox{in} & \ U, \\
& u^\ep = f & \mbox{on} & \ \partial U,
\end{aligned}
\right.
\quad \mbox{and} \quad 
\left\{
\begin{aligned}
& -\nabla \cdot \left( \ahom \nabla  u  \right) = 0 &  \mbox{in} & \ U, \\
&  u  = f & \mbox{on} & \ \partial U.
\end{aligned}
\right.
\end{equation}
Then we have the estimate
\begin{multline}
\label{e.DPestimates.nosymm}
\left\| u^\ep -   u  \right\|_{L^{2}(U)}^2
+ \left\| \nabla u^\ep - \nabla  u  \right\|_{\Hminus(U)}^2 
+ \left\| \a\left(\tfrac\cdot\ep\right) \nabla u^\ep - \ahom \nabla  u  \right\|_{\Hminus(U)}^2 
\\
\leq
C \left\| \nabla f \right\|_{L^{2+\delta}(U)}^2 \left( 
\ep^{\beta(d-s)}
+ \X_s \ep^s
\right).  
\end{multline}
\end{theorem}

\begin{remark}[Regularity theory for nonsymmetric coefficient fields]
\label{r.regularity.nonsymm}
As Theorem~\ref{t.DP.nosymm} gives a full generalization of Theorem~\ref{t.DP} to the case of nonsymmetric coefficient fields, the complete regularity theory of Chapter~\ref{c.regularity} is now available, in particular the statements of Theorems~\ref{t.Lipschitz} and~\ref{t.regularity}, Corollary~\ref{cor.regularity} and Exercise~\ref{ex.regularity}. Indeed, the arguments of Chapter~\ref{c.regularity} make no use of symmetry and the only input from Chapters~\ref{c.one} and~\ref{c.two} is Theorem~\ref{t.DP}, which can therefore by replaced by Theorem~\ref{t.DP.nosymm}.

\smallskip

Similarly, one also obtains in a straightforward way the boundary regularity of Section~\ref{s.boundaryreg} and well as the Calder\'on-Zygmund estimates of Chapter~\ref{c.CZ}: in particular, the statements of Theorems~\ref{t.Lipschitz.boundary},~\ref{t.Liouville.halfspace},~\ref{t.CZ} and~\ref{t.CZ.global}.
\end{remark}

Most of the rest of the results in the book, in particular the optimal bounds for the first-order correctors proved in Chapter~\ref{c.A1} and their scaling limits in Chapter~\ref{c.gff}, also have straightforward extensions to the case of nonsymmetric coefficient fields. Most of the work for obtaining such a generalization consists of extending Theorem~\ref{t.additivity} to a ``double-variable'' version of~$J_1$, using the ideas of the present chapter and following the arguments of Chapter~\ref{c.A1}. We will not present the details of such argument here.

\section*{Notes and references}

An adaptation to the case of nonsymmetric coefficients of the subadditive arguments of~\cite{AS} was first obtained in~\cite{AM} in the more general context of general \emph{nonlinear} equations in divergence form (i.e., in which the coefficients are uniformly monotone maps). The arguments in this chapter follows these ideas, although the presentation here is new and much simpler than that of~\cite{AM}. More on variational principles for non--self adjoint linear operators and evolutionary equations can be found in~\cite{GhBook,Vis} and in the context of periodic homogenization for uniformly monotone maps in~\cite{GMS,Vis2}. A related ``double variable'' variational principle was previously used in the context of homogenization of linear, nonsymmetric elliptic operators in~\cite{FP1,FP2}.



\chapter{Nonlinear equations}
\label{c.nonlinear}

In this chapter, we adapt the approach developed in earlier chapters to the case of \emph{nonlinear} elliptic equations. The main focus is on extending the arguments of Chapters~\ref{c.two} and~\ref{c.regularity}. 

\smallskip

In the first section, we present the assumptions and some basic facts concerning nonlinear elliptic equations. In Section~\ref{s.Lsubadd}, we introduce the subadditive quantities and adapt the arguments of Chapter~\ref{c.two} to prove their convergence. We give a quantitative result on the homogenization of the Dirichlet problem in Section~\ref{s.Ldirichlet}. In Section~\ref{s.LLip.NL} we prove a large-scale~$C^{0,1}$-type estimate. 

\section{Assumptions and preliminaries}

The most general elliptic equation we wish to consider takes the form
\begin{equation} 
\label{e.generalPDE}
-\nabla \cdot \a(\nabla u,x ) = 0.
\end{equation}
Here the coefficient field~$\a(p,x)$ is a map $\a:\Rd\times\Rd\to \Rd$ which may depend nonlinearly in its first variable. The uniform ellipticity assumption is the condition that, for every $p_1,p_2,x\in\Rd$,
\begin{equation} 
\label{e.uniformlymonotone}
\left( \a(p_1,x) - \a(p_2,x) \right) \cdot (p_1-p_2) \geq 
\left| p_1-p_2 \right|^2 
\end{equation}
and
\begin{equation}
\label{e.aLipschitz}
\left| \a(p_1,x) - \a(p_2,x) \right|
\leq \Lambda \left| p_1-p_2 \right|.
\end{equation}
We also require that, for some constant $\msf K \in [1,\infty)$ and every $p\in\Rd$, 
\begin{equation} 
\label{e.abounded}
\left\| \a(p,\cdot) \right\|_{L^\infty(\Rd)} \leq \mathsf K + \Lambda |p|. 
\end{equation}
The condition~\eqref{e.uniformlymonotone} says that the map $p\mapsto \a(p,x)$ is \emph{uniformly monotone}, while~\eqref{e.aLipschitz} and~\eqref{e.abounded} require~$\a(\cdot,\cdot)$ to be Lipschitz in its first variable and bounded uniformly in its second variable. We make no assumption on the regularity of $\a(p,\cdot)$ beyond measurability and boundedness. In the case that~\eqref{e.generalPDE} is linear, that is, $\a(p,x) = \tilde\a(x)p$ for a $d\times d$ matrix $\tilde \a(x)$, the assumptions~\eqref{e.uniformlymonotone} and~\eqref{e.aLipschitz} are equivalent to~\eqref{e.ue.Fitz} and~\eqref{e.abounded} is superfluous.

\smallskip

To simplify the presentation, we assume throughout the chapter that 
\begin{equation} 
\label{e.classicalform}
\a(p,x) =  D_p L(p,x),
\end{equation}
where $L=L(p,x)$ is a Lagrangian which is convex in $p$. (We denote the derivatives of~$L$ in the variable~$p$ by~$D_p$, and reserve the symbol~$\nabla$ for spatial variables.) Note that~\eqref{e.classicalform} is a generalization of the assumption from the linear case that the coefficient matrix is symmetric, and it gives~\eqref{e.generalPDE} a (classical) variational formulation:
\begin{multline*} \label{}
u\in H^1(U) \ \mbox{is a solution of~\eqref{e.generalPDE}}
\\
\iff
u \ \mbox{is the minimizer in $u+H^1_0(U)$ of} \ w\mapsto \int_U L(\nabla w,\cdot). 
\end{multline*}
The assumptions~\eqref{e.uniformlymonotone}, \eqref{e.aLipschitz} and \eqref{e.abounded} can be rewritten in terms of $L$ as the statements that, for every $p_1,p_2, x \in \Rd$,
\begin{equation} 
\label{e.Lunifconv3}
\left(  D_p L(p_1,x) -  D_p L(p_2,x) \right) \cdot \left(p_1 - p_2 \right) 
\geq 
 \left|p_1-p_2\right|^2,
\end{equation}
\begin{equation} 
\label{e.Lunifconv5}
\left|  D_p L(p_{1},x) -  D_p L(p_{2},x) \right|  \leq \Lambda \left|p_{1}-p_{2}\right|,
\end{equation}
and
\begin{equation} 
\label{e.Lnormalize}
\left\|  D_pL(p,\cdot) \right\|_{L^\infty(\Rd)} \leq \msf K + \Lambda |p|. 
\end{equation}
Note that \eqref{e.Lunifconv3} and \eqref{e.Lunifconv5} are equivalent to the statement that for every $x \in \Rd$, the mapping $p \mapsto L(p,x)$ is $C^{1,1}$ and uniformly convex:
\begin{equation}
\label{e.Lunifconv}
\forall x\in\Rd \ \mbox{and a.e.} \, p\in\Rd, 
\quad 
\Id \leq D^2_pL(p,x) \leq \Lambda \, \Id.
\end{equation}
This assumption can be written equivalently as follows: for every $p_1,p_2,x\in\Rd$,
\begin{equation} 
\label{e.Lunifconv2}
\frac14 \left| p_1-p_2 \right|^2 
\leq
L(p_1,x) + L(p_2,x) - 2 L\left(\frac{p_1+p_2}2,x\right) 
\leq 
\frac{\Lambda}4\left|p_1-p_2\right|^2. 
\end{equation}
Note that adding a function of $x$ only to $L$ does not change the definition of $\a$ in \eqref{e.classicalform} or any of the assumptions in \eqref{e.Lunifconv3}-\eqref{e.Lnormalize}. We can therefore assume without loss of generality that, for every $p,x \in \Rd$,
\begin{equation} 
\label{e.Lorigincontrol}
\frac 12 |p|^2 -\msf K |p|
\leq
 L(p,x)
\leq 
\frac \Lambda2 |p|^2 + \msf K |p|.
\end{equation}
In this chapter, we therefore modify the standing assumptions by re-defining  
\begin{multline*} \label{}
\Omega:= 
\big\{ 
L\,:\, L \ \mbox{is a measurable map from} \  \Rd \times \Rd \to \R 
\\
 \mbox{satisfying~\eqref{e.Lnormalize}, \eqref{e.Lunifconv2} and \eqref{e.Lorigincontrol}} 
\big\}.
\end{multline*}
The family~$\left\{ \F_U \right\}$ indexed by the Borel subsets $U \subseteq\Rd$ is defined by
\begin{align*} \label{}
\F_U & :=  \mbox{the $\sigma$-algebra generated by} \ 
\left\{ L \mapsto \int_{\Rd} L(p,x) \varphi(x)\,dx  :  \varphi \in C^\infty_c(U), \ p\in\Rd \right\}.
\end{align*}
We take $\F:=\F_{\Rd}$. The translation operator $T_y$ acts on $\Omega$ by
\begin{equation} 
\label{e.Ty.L}
\left( T_y L \right)(p,x):= L(p,x+y),
\end{equation}
and, as usual, we extend this to elements of~$\F$ by defining $T_y E:= \left\{ T_yL\,:\, L \in E \right\}$. 
With these modifications to the definitions, the rest of the assumptions are in force. In particular, we take $\P$ to be a probability measure on $(\Omega,\F)$ satisfying stationarity with respect to $\Zd$ translations and the unit range of dependence assumption. We also take $\overline\Omega$ to be the subset of $\Omega$ consisting of those~$L$'s which do not depend on~$x$:
\begin{multline} 
\label{e.Omegabar.L}
\overline\Omega:= 
\big\{ 
L\,:\, L \ \mbox{is a measurable map from} \  \Rd \to \R 
\\
\mbox{such that $(p,x) \mapsto L(p)$ belongs to $\Omega$} 
\big\}.
\end{multline}

\smallskip

We denote by $\mathcal{L}(U)$ the subset of $H^1_{\mathrm{loc}}(U)$ consisting of local minimizers of the integral functional, or equivalently, weak solutions $w$ of the equation
\begin{equation*} \label{}
-\nabla \cdot \left(  D_pL\left( \nabla w,x \right) \right) = 0 \quad \mbox{in} \ U. 
\end{equation*}

The nonlinear version of Caccioppoli's inequality is the following (compare the statement and its proof with the linear case, Lemma~\ref{l.Caccioppoli appendix}). 

\begin{lemma}[{Caccioppoli inequality}]
\label{l.nonlinear.caccioppoli}
There exists~$C(d,\Lambda)<\infty$ such that, for every $r>0$ and $u,v\in \mathcal{L}(B_r)$, 
\begin{equation} 
\label{e.Caccioppoli.NL}
\left\| \nabla u - \nabla v \right\|_{\underline{L}^2(B_{r/2})} 
\leq
\frac{C}{r} \left\| u - v - (u-v)_{B_{r}} \right\|_{\underline{L}^2(B_r)}.
\end{equation}
\end{lemma}
\begin{proof}
The proof can be compared to that of Lemma~\ref{l.Caccioppoli appendix}. 
By subtracting a constant from $u$, we may suppose that~$\left( u - v \right)_{B_r} = 0$. Select $\phi \in C^\infty_c(B_r)$ to satisfy
\begin{equation} 
\label{e.phicutoffCacc.NL}
0\leq \phi \leq 1, \quad 
\phi = 1 \quad \mbox{in} \ B_{r/2}, \quad
\left| \nabla \phi \right| \leq 4r^{-1},  
\end{equation}
and test the equations for $u$ and $v$ each with $\phi^2(u-v)$ to get
\begin{equation*} \label{}
\fint_{B_{r}} \nabla \left( \phi^2 (u-v) \right) \cdot  D_pL(\nabla u, \cdot) = 0 = \fint_{B_{r}} \nabla \left( \phi^2 (u-v) \right) \cdot  D_pL(\nabla v, \cdot). 
\end{equation*}
After rearranging, this gives
\begin{multline*} \label{}
\fint_{B_r} \phi^2 \left( \nabla u - \nabla v \right) \cdot \left(  D_pL(\nabla u, \cdot)-  D_pL(\nabla v, \cdot) \right) 
\\
=
- \fint_{B_r} 2\phi(u-v) \nabla \phi \cdot \left(  D_pL(\nabla u, \cdot)-  D_pL(\nabla v, \cdot) \right).
\end{multline*}
Using uniform ellipticity in~\eqref{e.Lunifconv3}, as well as the Lipschitz bound~\eqref{e.Lunifconv5}, we obtain
\begin{equation*} \label{}
\fint_{B_r} \phi^2 \left| \nabla u - \nabla v \right|^2 
\leq 
C \fint_{B_r} \phi \left| \nabla \phi \right| \left| u-v \right| \left| \nabla u - \nabla v \right|.
\end{equation*}
By Young's inequality,
\begin{equation*} \label{}
C \fint_{B_r} \phi \left| \nabla \phi \right| \left| u-v \right| \left| \nabla u - \nabla v \right|
\leq \fint_{B_r}  \frac 12 \left( \phi^2  \left| \nabla u - \nabla v \right|^2  + C\left| \nabla \phi \right|^2 \left| u-v \right|^2\right),
\end{equation*}
and therefore we obtain, after reabsorption,
\begin{equation*} \label{}
\fint_{B_r} \phi^2 \left| \nabla u - \nabla v \right|^2 
\leq
C\fint_{B_r} \left| \nabla \phi \right|^2 \left| u-v \right|^2.
\end{equation*}
Using~\eqref{e.phicutoffCacc.NL}, we deduce
\begin{equation*} \label{}
\fint_{B_{r/2} } \left| \nabla u - \nabla v \right|^2 
\leq 
C r^{-2} \fint_{B_r}  \left| u-v \right|^2.
\end{equation*}
This completes the proof. 
\end{proof}

\section{Subadditive quantities and basic properties}
\label{s.Lsubadd}

In this section, we introduce the subadditive quantities~$\nu$ and~$\nu^*$ and show that, compared to the linear case, many of the basic properties are preserved.

\smallskip

Recall that $\ell_p$ is the affine function $\ell_p(x)=p\cdot x$. 
Given $L\in\Omega$, $p,q\in\Rd$ and a bounded Lipschitz domain $U\subseteq \Rd$, we define
\begin{equation} 
\label{e.nu.NL}
\nu(U,p):= \inf_{u\in \ell_p+H^1_0(U)}  \fint_{U} L\left(\nabla u,\cdot \right)
\end{equation}
and
\begin{equation}
\label{e.def.nu*.NL}
\nu^*(U,q):= \sup_{u\in H^1(U)} \fint_{U} \left( -L\left(\nabla u,\cdot \right) + q\cdot \nabla u \right).
\end{equation}
We also denote
\begin{equation*} \label{}
v(\cdot,U,p) := \mbox{minimizing element of $\ell_p+H^1_0(U)$ in~\eqref{e.nu.NL}}
\end{equation*}
and
\begin{equation*} \label{}
u(\cdot,U,q) := \mbox{maximizing element of $H^1(U)$ in~\eqref{e.def.nu*.NL}.}
\end{equation*}
As in the linear case, see~\eqref{e.obvioustesting}, we may compare $\nu^*(U,q)$ to~$\nu(U,p)$ by testing the minimizer of $\nu(U,p)$ in the definition of $\nu^*(U,q)$. This gives
\begin{equation} 
\label{e.obvioustesting.NL}
\nu(U,p) + \nu^*(U,q) \geq p\cdot q. 
\end{equation}
Thus if we define $J$ by
\begin{equation*} \label{}
J(U,p,q):= \nu(U,p) + \nu^*(U,q) - p\cdot q,
\end{equation*}
we see that $J$ is nonnegative. We cannot however hope to write $J$ as an optimization problem as in~\eqref{e.formulafornuJpq}. This turns out to be only slightly inconvenient.


\smallskip

We begin by recording some basic properties of~$\nu$ and $\nu^*$, extending Lemma~\ref{l.basicJ} to the nonlinear case.

\begin{lemma}[Basic properties of $\nu$ and $\nu^*$]
\label{l.basicJ.NL}
Fix a bounded Lipschitz domain $U\subseteq \Rd$. 
The quantities $\nu(U,p)$ and $\nu^*(U,q)$ and the functions $v(\cdot,U,p)$ and $u(\cdot,U,q)$ satisfy the following properties:

\begin{itemize}

\item \emph{Uniformly convex and $C^{1,1}$.} There exists $C(\Lambda)<\infty$ such that, for every  $p_1,p_2\in\Rd$, 
\begin{multline} 
\label{e.JunifconvC11p.NL}
 \frac14 \left| p_1-p_2\right|^2
\leq
\nu\left( U,p_1 \right) +  \nu\left( U,p_2 \right)  - 2 \nu\left( U,\frac{p_1+p_2}2 \right)
\leq
C \left| p_1-p_2\right|^2
\end{multline}
and, for every $q_1,q_2,p\in\Rd$, 
\begin{multline} 
\label{e.JunifconvC11q.NL}
\frac 1C \left| q_1 - q_2\right|^2
\leq
 \nu^*\left( U,q_1 \right) + \nu^*\left( U,q_2 \right)  - 2\nu^*\left( U,\frac{q_1+q_2}2 \right)
\leq
 \frac14 \left| q_1 - q_2\right|^2.
\end{multline}

\item \emph{Subadditivity.} Let $U_1, \ldots, U_N \subset U$ be bounded Lipschitz domains that form a partition of~$U$, in the sense that $U_i \cap U_j = \emptyset$ if $i\neq j$ and 
\begin{equation*} \label{}
\left| U \setminus \bigcup_{i=1}^N U_i \right| = 0. 
\end{equation*}
Then, for every $p,q\in\Rd$,
\begin{equation}
\label{e.Jsubadd.NL}
\nu(U,p) \leq \sum_{i=1}^N \frac{\left|U_i\right|}{|U|} \nu(U_i,p)
\quad \mbox{and} \quad 
\nu^*(U,q) \leq \sum_{i=1}^N \frac{\left|U_i\right|}{|U|} \nu^*(U_i,q). 
\end{equation}

\item \emph{First variations.} For each~$q\in\Rd$, the maximizer~$u(\cdot,U,q)$ in the definition~\eqref{e.def.nu*.NL} of~$\nu^*(U,q)$ is characterized as the unique element of $H^1(U) \cap \mathcal{L}(U)$ which satisfies
\begin{equation}
\label{e.nu*firstvar.NL} 
\forall\, w \in H^1(U), 
\quad
\fint_U \nabla w \cdot  D_pL \left( \nabla u(\cdot,U,q),\cdot \right) = q\cdot \fint_U\nabla w.
\end{equation}
For each $p\in\Rd$, the minimizer $v(\cdot,U,p)$ is characterized as the unique element of $\mathcal{L}(U)$, up to a constant, such that
\begin{equation}
\label{e.nufirstvar.NL} 
\forall\, w \in \mathcal{L} (U), 
\quad
\fint_U \nabla v \cdot  D_pL \left( \nabla w, \cdot \right) = p\cdot \fint_U  D_pL \left( \nabla w, \cdot \right).
\end{equation}

\item \emph{Quadratic response.} There exists $C(\Lambda)<\infty$ such that, for every $q\in\Rd$ and $w\in H^1(U)$, 
\begin{multline} 
\label{e.nu*quadresponse.NL}
\frac12 \fint_U \left| \nabla w - \nabla u(\cdot,U,q)\right|^2
\leq
\nu^*(U,q) - \fint_U \left( - L(\nabla w,\cdot) + q\cdot \nabla w  \right)
\\
\leq 
\frac{\Lambda}{2} \fint_U \left| \nabla w - \nabla u(\cdot,U,q)\right|^2,
\end{multline}
and, similarly, for every $p\in\Rd$ and $w\in \ell_p + H^1_0(U)$, 
\begin{multline} 
\label{e.nuquadresponse.NL}
\frac12 \fint_U \left| \nabla w - \nabla v(\cdot,U,p)\right|^2
\leq
\nu(U,p) - \fint_U L(\nabla w, \cdot) 
\\
\leq 
\frac{\Lambda}{2}  \fint_U \left| \nabla w - \nabla v(\cdot,U,p)\right|^2.
\end{multline}

\item \emph{Formulas for derivatives of $\nu$ and $\nu^*$.} For every $p,q\in\Rd$,
\begin{equation} 
\label{e.nuderp.NL}
 D_p\nu(U,p)  = \fint_U  D_pL\left( \nabla v(\cdot,U,p),\cdot\right) 
\end{equation}
and
\begin{equation} 
\label{e.nu*derq.NL}
{D}_q\nu^*(U,q)  =  \fint_U  \nabla u(\cdot,U,q).
\end{equation}
\end{itemize}
\end{lemma}
\begin{proof}
The proof of subadditivity~\eqref{e.Jsubadd.NL} is nearly identical to the proof of~\eqref{e.Jsubadd} and thus we omit it.

\smallskip

For the rest of the proof, we fix $p,q \in \R^d$ and $U \subset \R^d$, and use the shorthand notation $u := u(\cdot,U,q)$ and $v := v(\cdot,U,p)$.  

\smallskip

\emph{Step 1.} We first prove the first variations. Fix $w \in H^1(U)$ and $t \in \R \setminus \{0\}$, and set $u_t := u + t w$. We have
\begin{align*} 
\nu^*(U,q) & \geq \fint_{U} \left(- L(\nabla u_t ,\cdot) + q \cdot \nabla u_t  \right) 
\\ & = \nu^*(U,q) + t \fint_{U} \left( q \cdot \nabla w - \nabla w \cdot \int_0^1  D_p L\left(\nabla u + st \nabla w,\cdot \right)\,ds   \right).
\end{align*}
Dividing by $t$ and sending $t \to 0$ proves~\eqref{e.nu*firstvar.NL}. Clearly $u \in \mathcal{L}(U)$. 

\smallskip

To show~\eqref{e.nufirstvar.NL}, fix $\phi \in H_0^1(U)$ and $t \in \R \setminus \{0\}$, and set $v_t = v + t \phi$. Then 
\begin{align*} 
\nu(U,p) & \leq \fint_{U} L(\nabla v_t ,x) \,dx 
\\ & = \nu(U,p) + t \fint_{U} \int_0^1  D_p L(\nabla v(x,U,p) + t s \nabla \phi(x) ,x ) \, ds \cdot \nabla \phi(x) \, dx.
\end{align*}
Dividing by $t$ and sending $t \to 0$ gives that $v \in \mathcal L(U)$. Since $v - p \cdot x \in H_0^1(U)$, we deduce that~\eqref{e.nufirstvar.NL} holds. To prove that $v$ is the unique element satisfying~\eqref{e.nufirstvar.NL}, take $v_1$ and $v_2$ satisfying~\eqref{e.nufirstvar.NL} and observe that we obtain
\begin{equation*} 
\fint_{U} \left( D_p L(\nabla v_1 ,\cdot ) -  D_p L(\nabla v_2 ,\cdot ) \right) \cdot \left(\nabla v_1 - \nabla v_2 \right) = 0.
\end{equation*}
This implies, by~\eqref{e.Lunifconv3}, that $v_1 - v_2$ is a constant, proving the claim.  

\smallskip

\emph{Step 2.} Quadratic response. We only prove~\eqref{e.nu*quadresponse.NL} since the proof of~\eqref{e.nuquadresponse.NL} is very similar. 
By the first variation~\eqref{e.nu*firstvar.NL}, we obtain that 
\begin{align*} 
\lefteqn{\nu^\ast(U,q) - \fint_{U} \left( - L\left(\nabla w,\cdot \right) + q \cdot \nabla w \right) } \quad & \\
& = \fint_U \left( \int_0^1  D_p L \left( \nabla u + t \left( \nabla w - \nabla u\right),\cdot \right) \,dt - q \right) \cdot \left( \nabla w - \nabla u\right) 
\\ & = \fint_U \int_0^1 \left(  D_p L \left( \nabla u + t \left( \nabla w - \nabla u\right),\cdot \right) -  D_p L \left( \nabla u , \cdot \right) \right) \cdot \left( \nabla w - \nabla u \right) \,dt .
\end{align*}
Now~\eqref{e.nu*quadresponse.NL} follows from~\eqref{e.Lunifconv3}.

\smallskip

\emph{Step 3.} The proof of~\eqref{e.nuderp.NL}. We will show more, namely the existence of a constant $C(\Lambda)<\infty$ such that, for every $p,p'\in\Rd$,  
\begin{equation} 
\label{e.quant.diffnu.NL}
\frac12 \left| p-p' \right|^2 \leq  \nu(U,p) - \nu(U,p') - \left( p -p' \right) \cdot \fint_U  D_pL(\nabla v,\cdot)  \leq C \left| p-p' \right|^2. 
\end{equation}
Fix $p,p'\in\Rd$ and set $v:= v(\cdot,U,p)$ and $v':=v(\cdot,U,p')$. Since $v - v' - (p-p')\cdot x \in H_0^1(U)$, we have that 
\begin{equation*} 
0 = \fint_U \left( D_p L(\nabla v,\cdot) -  D_p L(\nabla v',\cdot)   \right) \cdot \left( \nabla v - \nabla v' - (p-p') \right).
\end{equation*}
Using~\eqref{e.Lunifconv3}, we thus obtain
\begin{equation*} 
\fint_{U} \left| \nabla v - \nabla v' \right|^2  \leq C \left| p - p'\right|^2 .
\end{equation*}
On the other hand, by Jensen's inequality,
\begin{equation*} 
 \left| p - p'\right|^2  =  \left| \fint_{U} \left(\nabla v - \nabla v' \right) \right|^2 \leq \fint_{U} \left| \nabla v - \nabla v' \right|^2 .
\end{equation*}
Applying~\eqref{e.nufirstvar.NL}, we get 
\begin{align*}
\lefteqn{
\nu(U,p) - \nu(U,p') - \left( p -p' \right) \cdot \fint_U  D_pL(\nabla v,\cdot) 
} \quad & 
\\ & 
=  \fint_U \left( L(\nabla v,\cdot) - L\left(\nabla v',\cdot\right)  -  D_pL(\nabla v,\cdot) \cdot \left( \nabla v - \nabla v' \right) \right) 
\\ & 
=  \fint_U  \left(\int_0^1  D_p L(\nabla v' + t (\nabla v - \nabla v') ,\cdot) \, dt -  D_pL(\nabla v,\cdot) \right) \cdot \left( \nabla v - \nabla v' \right)  .
\end{align*}
Combining the three displays above and applying~\eqref{e.Lunifconv3} yields~\eqref{e.quant.diffnu.NL}.

\smallskip

\emph{Step 4.} The proof of~\eqref{e.nu*derq.NL}. We will show that there exists~$C(\Lambda)<\infty$ such that, for every $q,q'\in\Rd$, 
\begin{equation} 
\label{e.quant.diffnu*.NL}
\frac 1C \left| q- q' \right|^2 \leq \nu^*(U,q) - \nu^*(U,q') - \left( q-q' \right) \cdot \fint_U \nabla u  \leq  \frac12 \left| q- q' \right|^2. 
\end{equation}
Fix $q,q'\in\Rd$ and set $u:=u(\cdot,U,q)$ and $u':=u(\cdot,U,q')$. First,~\eqref{e.nu*firstvar.NL} gives
\begin{equation*} 
0 = \fint_U \left( \left(  D_p L(\nabla u,\cdot) -  D_p L(\nabla u',\cdot) \right) \cdot \left(\nabla u - \nabla u' \right)  - (q-q') \cdot \left(\nabla u - \nabla u' \right) \right),
\end{equation*}
and hence, by~\eqref{e.Lunifconv3},
\begin{equation*} 
\fint_{U} \left| \nabla u - \nabla u' \right|^2  \leq \left| q - q'\right|^2.
\end{equation*}
On the other hand, testing~\eqref{e.nu*firstvar.NL} with $w(x) := (q-q')\cdot x$ yields by~\eqref{e.Lunifconv5} that
\begin{equation*} 
\left| q - q'\right|^2 \leq  \left| q - q'\right| \fint_U \left|   D_p L(\nabla u,\cdot) -  D_p L(\nabla u',\cdot) \right| \leq C\left| q - q'\right| \left( \fint_{U} \left| \nabla u - \nabla u' \right|^2 \right)^{\frac 12}.
\end{equation*}
From this we obtain
\begin{equation*} 
\left| q - q'\right|^2 \leq C \fint_{U} \left| \nabla u - \nabla u' \right|^2.
\end{equation*}
Using the first variation~\eqref{e.nu*firstvar.NL} again, we find that 
\begin{align*} \label{}
\lefteqn{ 
\nu^*(U,q) - \nu^*(U,q') - \left( q-q' \right) \cdot \fint_U \nabla u 
} \quad & 
\\ &
=
 \fint_U \left( L(\nabla u',\cdot) - L(\nabla u,\cdot) + q'\cdot \left( \nabla u - \nabla u' \right)  \right) 
\\ & 
=  \fint_U  \left( L(\nabla u',\cdot) - L(\nabla u,\cdot)  -  D_pL\left(\nabla u',\cdot \right)  \cdot \left( \nabla u' - \nabla u\right)   \right) 
\\ & =
 \fint_U  \left(\int_0^1  D_p L(\nabla u + t (\nabla u' - \nabla u) ,\cdot) \, dt -  D_pL(\nabla u',\cdot) \right) \cdot \left( \nabla u' - \nabla u \right) .
\end{align*}
From the last three displays and~\eqref{e.Lunifconv3} we get~\eqref{e.quant.diffnu*.NL}.  

\smallskip

\emph{Step 5.} Uniform convexity and $C^{1,1}$:~\eqref{e.JunifconvC11p.NL} and~\eqref{e.JunifconvC11q.NL} follow from~\eqref{e.quant.diffnu.NL} and~\eqref{e.quant.diffnu*.NL}, respectively.
\end{proof}

By~\eqref{e.JunifconvC11p.NL} and~\eqref{e.JunifconvC11q.NL} it follows, as we show in the following exercise, that there exists a constant~$C(\Lambda)<\infty$ such that, for $U \subset \R^d$ and every $p_1,p_2 \in \R^d$, 
\begin{equation} 
\label{e.Lip nu and nu^ast.NL}
\left| D_p \nu(U,p_1) -  D_p \nu(U,p_2)  \right| + \left| D_q \nu^\ast(U,p_1) -  D_q \nu^\ast(U,p_2)  \right|  \leq C |p_1 - p_2|.
\end{equation}

\begin{exercise}
\label{ex.C11}
Suppose that $K>0$, $U\subseteq \Rd$ is convex  and $f:U \to \R$ is continuous and satisfies, for every $x_1,x_2\in U$, 
\begin{equation*} \label{}
\left| \tfrac12 f(x_1) + \tfrac12 f(x_2) - f\left( \tfrac12x_1+\tfrac12x_2 \right) \right| 
\leq
K \left| x_1 - x_2 \right|^2. 
\end{equation*}
Prove that $f \in C^{1,1}(U)$ and that, for every $x_1,x_2\in U$,
\begin{equation*} \label{}
\left| \nabla f(x_1) - \nabla f(x_2) \right| \leq CK \left| x_1-x_2\right|. 
\end{equation*}
\end{exercise}

\index{coarsened coefficient field}
We next introduce ``coarsenings'' $L_U$ of the Lagrangian~$L$ which are defined in terms of the Legendre-Fenchel transform of $\nu^*(U,\cdot)$. For every $p \in \Rd$ and bounded Lipschitz domain~$U\subseteq \Rd$, define
\begin{equation} 
\label{e.LU.def}
L_U(p):= \sup_{q\in\Rd} \left(p\cdot q -  \nu^*(U,q) \right).
\end{equation}
Since $q \mapsto \nu^*(U,q)$ is uniformly convex and bounded above by a quadratic function, we have (see~\cite[Section 3.3.2]{Evans}) 
\begin{equation*} 
\sup_{ p \in\Rd} \left(p\cdot q -  L_U(p) \right) = \nu^*(U,q). 
\end{equation*}
Since the mapping $p \mapsto \nu^*(u,p)$ is $C^{1,1}$ and uniformly convex, we deduce that, for every~$q\in\R^d$,
\begin{equation*} 
DL_U\left(D\nu^*(U,q)\right) = q \quad \mbox{and} \quad D\nu^*\left(U,DL_U(q)\right) = q. 
\end{equation*}
In view of~\eqref{e.nu*derq.NL}, we obtain the identity, for every~$p \in \R^d$, 
\begin{equation*} \label{}
p = \fint_U \nabla u \left( \cdot, U, D L_U(p) \right). 
\end{equation*}
We also introduce deterministic versions of $L_U$ by setting 
\begin{equation} 
\label{e.LbarU.def}
\overline{L}_U(p):= \sup_{q\in\Rd} \left(p\cdot q - \E \left[ \nu^*(U,q) \right] \right). 
\end{equation}
That is, $\overline{L}_U$ is the Legendre-Fenchel transform of $q \mapsto \E \left[ \nu^*(U,q) \right]$. This quantity is analogous to the coarsened matrix $\ahom_U$ introduced in Definition~\ref{def.ahomU}. For notational convenience, we also set
\begin{equation*}  
\bar \nu^*(U,q) := \E \Ll[ \nu^*(U,q) \Rr] ,
\end{equation*}
and observe that 
\begin{equation*}  
\sup_{p \in \Rd} \Ll( p\cdot q - \bar \nu^*(U,q) \Rr)  = \bar L_U(p).
\end{equation*}
It follows from \eqref{e.nu*derq.NL} that
\begin{equation} 
\label{e.barnu*derq.NL}
D \bar \nu^*(U,q) = \E \left[ \fint_{U} \nabla u(\cdot,U,q) \right],
\end{equation}
and from~\eqref{e.JunifconvC11q.NL} that, for every $q_1,q_2 \in \R^d$, 
\begin{equation*} 
\left|q_1 - q_2\right|  \leq C \left|D \bar \nu^*(U,q_1) - D \bar \nu^*(U,q_2) \right|.
\end{equation*}
Since $D \overline{L}_U$ and $D \bar \nu^*(U,\cdot)$ are inverse functions, we infer from \eqref{e.barnu*derq.NL} that, for every $p \in \Rd$,
\begin{equation} 
\label{e.ptoqforL}
p = \E \left[ \fint_{U} \nabla u\left(\cdot,U, D \overline{L}_U(p)  \right) \right]. 
\end{equation}
We also use the shorthand notation $\overline{L}_n:=\overline{L}_{\cu_n}$ and $\bar \nu^*_n := \bar \nu^*(\cu_n,\cdot)$.
We next observe that, by \eqref{e.barnu*derq.NL}, \eqref{e.ptoqforL}, stationarity and quadratic response, we have
\begin{align*}
\lefteqn{
\left| p - D \bar \nu^*_{n+1} \left( D \overline{L}_n(p) \right) \right|^2
} \ \ &
\\ & 
= 
\left| \E \left[    
\fint_{\cu_{n}} \nabla u\left(\cdot,\cu_{n}, D \overline{L}_{n}(p) \right)  \right]
-
\E\left[ \fint_{\cu_{n+1}} \nabla u\left(\cdot,\cu_{n+1}, D \overline{L}_{n}(p) \right)  
\right]    \right|^2 
\\ & 
\leq  \sum_{z\in 3^n\Zd\cap \cu_{n+1}} 
\E\left[ 
\fint_{z+\cu_n} 
\left| \nabla u\left(\cdot,z+\cu_n,D \overline{L}_{n}(p) \right)
- \nabla u\left(\cdot,\cu_{n+1}, D \overline{L}_{n}(p) \right) \right|^2
\right]
\\ & 
\leq C \left( \E \left[ \nu^*(\cu_{n},D \overline{L}_{n}(p)) \right] - \E \left[ \nu^*(\cu_{n+1},D \overline{L}_{n}(p)) \right]\right).
\end{align*}
It follows that 
\begin{multline} 
\label{e.nablaLhom.comparescales}
\left| D \overline{L}_{n+1} (p) - D \overline{L}_{n} (p) \right|^2 \\ \leq C\left( \E \left[ \nu^*(\cu_{n},D \overline{L}_{n}(p)) \right] - \E \left[ \nu^*(\cu_{n+1},D \overline{L}_{n}(p)) \right] \right).
\end{multline}
We also note that, similarly to~\eqref{e.mumap}, the uniformly convex function
\begin{equation}
\label{e.mumap.NL}
q\mapsto \, \E \left[ J\left(U, p,q \right) \right] = \E \left[ \nu\left(U, p \right) \right] + \E \left[ \nu^*\left(U,q \right) \right] -  p\cdot q
\end{equation}
achieves its minimum at the point $q$ such that $p =  D \bar \nu^*(U,q)$, which is equivalent to $q= D \overline{L}_U(p)$. It therefore follows from~\eqref{e.JunifconvC11q.NL} that, for some constant $C(\Lambda)<\infty$ and every $p,q \in \Rd$, 
\begin{equation}
\label{e.quadresponseLhom.NL}
\frac1C \left| q -  D\overline{L}_U(p) \right|^2
 \leq \E \left[ J\left(U, p,q \right) \right] 
- \E \left[ J\left(U, p, D\overline{L}_U(p) \right) \right] 
\leq
C \left| q -  D\overline{L}_U(p) \right|^2. 
\end{equation}
Finally, we define the effective Lagrangian $\overline{L}$ by setting, for every $p \in \Rd$,
\index{homogenized coefficients}
\begin{equation} 
\label{e.def.barL}
\overline{L}(p):=  \lim_{n \to \infty} \E \Ll[ \nu( \cu_n, p) \Rr] .
\end{equation}
This is an extension to the nonlinear setting of the quantity~$\ahom$ introduced in~Definition~\ref{d.ahom}. 
Recalling the definition of~$\overline{\Omega}$ in~\eqref{e.Omegabar.L}, we observe that $\bar L \in \bar \Omega$.  
We denote the convex dual of $\bar L$ by $\bar L\!\,^*$: for every $q \in \Rd$,
\begin{equation}  
\label{e.def.barL*}
\bar L\!\,^*(q) := \sup_{p \in \Rd} \Ll( p \cdot q - \bar L(p) \Rr) .
\end{equation}

\section{Convergence of the subadditive quantities}

The main result of this section is the following generalization of Theorem~\ref{t.subadd}. 

\begin{theorem}
\label{t.subadd.NL}
Fix $s\in (0,d)$ and $\mathsf{M} \in [1,\infty)$. There exist $\alpha(d,\Lambda)\in \left(0,\frac12\right]$ and $C(s,d,\Lambda)<\infty$ such that, for every~$p,q\in\Rd$ and $n\in\N$,  
\begin{multline} 
\label{e.subadderror.NL}
\sup_{p \in B_{\mathsf{M}}}  \left| \nu(\cu_n,p) - \overline{L}(p) \right| + \sup_{q \in B_{\mathsf{M}}}  \left| \nu^*(\cu_n,q) - \overline{L}\!\,^*(q) \right|   \\ 
\leq 
C   3^{-n \alpha(d-s)}  \left( \mathsf{K}  +  \mathsf{M} \right)^2  + \O_1\left( C  \left( \mathsf{K}  +  \mathsf{M} \right)^2 3^{-ns}  \right). 
\end{multline}
\end{theorem}

The strategy is similar to the one of the proof of Theorem~\ref{t.subadd}. One difference is that we have to perform the iteration ``$p$-by-$p$'' rather than for all $p$'s at once. This leads to a different definition of the increment~$\tau_n$ compared to the one from Chapter~\ref{c.two}. 
We define, for each $p,q\in\Rd$ and $n\in\N$,
\begin{equation} 
\label{e.your.taus}
\left\{ 
\begin{aligned}
& \tau_n^\ast(q):= \E \left[ \nu^*(\cu_{n},q) \right] - \E \left[ \nu^*(\cu_{n+1},q) \right],
\\ & 
\tau_n(p):= \E \left[ \nu(\cu_{n},p) \right] - \E \left[ \nu(\cu_{n+1},p) \right],
\end{aligned}
\right.
\end{equation}
which are both nonnegative by subadditivity and stationarity. 
We begin with the generalization of Lemma~\ref{l.spatavg}. 

\begin{lemma}
\label{l.spatavg.NL}
There exist $\kappa(d)>0$ and $C(d,\Lambda)<\infty$ such that, for every $p,q\in \Rd$ and $m\in\N$,
\begin{equation} 
\label{e.spatavg.NL}
\var \left[ \fint_{\cu_m} \nabla u(\cdot,\cu_m,q)  \right] 
\leq  C (\mathsf{K}+|q|)^2 3^{-m\kappa} + C \sum_{n=0}^m 3^{-\kappa(m-n)} \tau_n^\ast(q).
\end{equation}
\end{lemma}
\begin{proof}
The proof is almost the same as the one of Lemma~\ref{l.spatavg}. We fix $q\in\Rd$. 

\smallskip
 
\emph{Step 1.} We show that
\begin{multline} 
\label{e.vardownscales.NL}
\var  \left[  \fint_{\cu_m} \nabla u(\cdot,\cu_m,q)\right]
\\
\leq \var  \left[ 3^{-d(m-n)}\sum_{z\in 3^n\Zd\cap \cu_m}  \fint_{z+\cu_n} \nabla u(\cdot,z+\cu_n,q)\right] 
+ C\sum_{k=n}^{m-1} \tau_k^\ast(q). 
\end{multline}
Simplify the notation by setting
\begin{equation*} \label{}
u:= u(\cu_m,q)\quad \mbox{and} \quad u_z:=u(z+\cu_n,q).
\end{equation*}
Using the identity 
\begin{equation*} \label{}
\fint_{\cu_m} \nabla u = 3^{-d(m-n)} \sum_{z\in 3^n\Zd \cap \cu_m} \fint_{z+\cu_n} \left( \nabla u - \nabla u_z \right)
+ 3^{-d(m-n)} \sum_{z\in 3^n\Zd\cap \cu_m} \fint_{z+\cu_n} \nabla u_z,
\end{equation*}
we find that 
\begin{multline*} \label{}
\var \left[ \fint_{\cu_m} \nabla u \right]
\leq 
2 \E \left[ 3^{-d(m-n)} \sum_{z\in 3^n\Zd \cap \cu_m} \left| \fint_{z+\cu_n} \left( \nabla u - \nabla u_z \right) \right|^2 \right] 
\\
+  2 \var \left[ 3^{-d(m-n)}\sum_{z\in 3^n\Zd\cap \cu_m} \fint_{z+\cu_n} \nabla u_z \right]. 
\end{multline*}
By~\eqref{e.Jsubadd.NL}, we have
\begin{align*} \label{}
\sum_{z\in 3^n\Zd\cap \cu_m} \left| \fint_{z+\cu_n} \left(  \nabla u - \nabla u_z \right) \right|^2
& \leq 
\sum_{z\in 3^n\Zd\cap \cu_m}
\fint_{z+\cu_n} \left|  \nabla u - \nabla u_z \right|^2
\\ & 
\leq C \sum_{z\in 3^n\Zd\cap \cu_m} \left(\nu^*(z+\cu_n,q) - \nu^*(\cu_m,q) \right)
\end{align*}
and taking the expectation of this yields, by stationarity, 
\begin{align*} \label{}
\E \left[ 3^{-d(m-n)} \sum_{z\in 3^n\Zd \cap \cu_m} \left| \fint_{z+\cu_n} \left( \nabla v - \nabla v_z \right) \right|^2 \right]
&
\leq C \E \left[ \nu^*(z+\cu_n,q) - \nu^*(\cu_m,q)  \right] 
\\ & 
= C \sum_{k=n}^{m-1} \tau_k^\ast(q). 
\end{align*}
Combining the two previous displays yields~\eqref{e.vardownscales.NL}.

\smallskip

\emph{Step 2.} We show that there exists~$\theta(d) \in (0,1)$ such that, for every $n\in\N$,
\begin{equation} 
\label{e.varsumcontract.NL}
\var \left[ 3^{-d}\sum_{z\in 3^n\Zd\cap \cu_{n+1}} \fint_{z+\cu_n} \nabla u(\cdot,z+\cu_n,q)\right] 
\\
\leq 
\theta \var  \left[\fint_{\cu_n} \nabla u(\cdot,\cu_n,q) \right]. 
\end{equation}
The proof of~\eqref{e.varsumcontract.NL} is exactly the same as that of Step~2 of the proof of Lemma~\ref{l.spatavg}. 

\smallskip

\emph{Step 3.} Iteration and conclusion. Denote
\begin{equation*} \label{}
\sigma_n^2(q) := \var  \left[\fint_{\cu_n} \nabla u(\cdot,\cu_n,q) \right].
\end{equation*}
By~\eqref{e.vardownscales} and~\eqref{e.varsumcontract}, there exists~$\theta(d) \in (0,1)$ such that, for every $n\in\N$,
\begin{equation} 
\sigma_{n+1}^2(q) 
\leq 
\theta \sigma_n^2(q) 
+ C \tau_n^\ast(q). 
\end{equation}
An iteration yields
\begin{equation*} \label{}
\sigma_m^2(q) \leq \theta^m \sigma_0^2(q) + C \sum_{n=0}^m \theta^{m-n} \tau_n^\ast(q).
\end{equation*}
Since $\sigma_0^2(q) \leq C(\mathsf{K}+|q|)^2$, this inequality implies
\begin{equation*} \label{}
\sigma_m^2 \leq C (\mathsf{K}+|q|)^2 3^{-m\kappa} + C \sum_{n=0}^m 3^{-\kappa(m-n)} \tau_n^\ast(q)
\end{equation*}
for $\kappa:= \log \theta / \log 3$, which completes the proof. 
\end{proof}

Recall that the analogue of the lemma above for~$v(\cdot,\cu_m,p)$ is trivial: indeed, since $v(\cdot,\cu_m,p) \in \ell_p + H^1_0(\cu_m)$, we have
\begin{equation}
\label{e.Stokesie}
p = \fint_{\cu_m} \nabla v(\cdot,\cu_m,p)
\end{equation}
and therefore
\begin{equation} 
\label{e.spatavg.NL.p}
\var \left[ \fint_{\cu_m} \nabla v(\cdot,\cu_m,p)  \right] 
=0.
\end{equation}

\smallskip

We next generalize Lemma~\ref{l.flatness}. Recall that~$\overline{L}_U$ is defined in~\eqref{e.LbarU.def} and satisfies~\eqref{e.ptoqforL}. That is,~$q=D \overline{L}_U(p)$ is the correct value of the parameter~$q$ so that the expectation of $\fint_U \nabla u(\cdot,U,q)$ is equal to~$p$.

\begin{lemma}
\label{l.flatness.NL}
There exist $\kappa(d)>0$ and $C(d,\Lambda)<\infty$ such that, for every $n\in\N$ and $p \in \Rd$,
\begin{multline}
\label{e.flatness.NL}
\E \left[\fint_{\cu_{n+1}} \left| u\left(x,\cu_{n+1}, D \overline{L}_n(p) \right) - p\cdot x \right|^2\,dx \right]  \\
 \leq 
 C3^{2n} \left( \left(\mathsf{K}+|p|\right)^2 3^{-\kappa n}+ \sum_{m=0}^n 3^{-\kappa (n-m)} \tau_m^\ast \left(  D\overline{L}_m(p)  \right)  \right)
\end{multline}
and
\begin{multline}
\label{e.flatness.NL.p}
\E \left[\fint_{\cu_{n+1}} \left| v(x,\cu_{n+1},p) - p\cdot x
\right|^2\,dx \right]  \\
 \leq 
 C3^{2n} \left( \left(\mathsf{K}+|p|\right)^2 3^{-\kappa n}+ \sum_{m=0}^n 3^{-\kappa (n-m)} \tau_m(p)  \right).
\end{multline}
\end{lemma}
\begin{proof}
We omit the proof since the argument has no important differences from the proof of Lemma~\ref{l.flatness}.
\end{proof}

We can now generalize Lemma~\ref{l.iterstep} in a nonlinear fashion by comparing the maximizer of $\nu^*(\cu_n,q)$ to the minimizer of $\nu(\cu_n,D \overline{L}\!\,^*_n(q))$, and then using the Caccioppoli inequality (Lemma~\ref{l.nonlinear.caccioppoli}) to estimate the difference between $\nu^*(\cu_n,q)$ and $p \cdot q -\nu(\cu_n, p)$ for $q=D \overline{L}_n(p)$.

\begin{lemma}
\label{l.iterstep.NL}
There exist $\kappa(d) > 0$ and $C(d,\Lambda)<\infty$ such that, for every~$p \in \Rd$ and~$n\in\N$,
\begin{multline}
\label{e.itersteprealz.NL}
\E \left[  J\left(\cu_n,p,D \overline{L}_n(p) \right)\right]  
\\
\leq 
C \left( \left(\mathsf{K}+|p|\right)^2 3^{-\kappa n}+ \sum_{m=0}^n 3^{-\kappa (n-m)} \left( \tau_m^\ast \left( D \overline{L}_m(p)\right)  + \tau_m(p) \right) \right).
\end{multline}
\end{lemma}
\begin{proof}
Fix $p\in \Rd$. Lemma~\ref{l.flatness.NL} and the triangle inequality give  
\begin{multline*}
\E \left[ \fint_{\cu_{n+1}} \left| u\left(\cdot,\cu_{n+1},D \overline{L}_n(p) \right) - v\left(\cdot,\cu_{n+1}, p \right) \right|^2 \right] 
\\
\leq 
C 3^{2n}\left( \left(\mathsf{K}+|p|\right)^2 3^{-\kappa n}+ \sum_{m=0}^n 3^{-\kappa (n-m)} \left( \tau_m^\ast \left( D \overline{L}_m(p)\right)  + \tau_m(p) \right) \right).
\end{multline*}
Applying the Caccioppoli inequality, Lemma~\ref{l.nonlinear.caccioppoli}, we find that
\begin{multline}
\label{e.NLdiffest}
\E \left[ \fint_{\cu_{n}}  \left| \nabla u\left(\cdot,\cu_{n+1},D \overline{L}_n(p) \right) - \nabla v\left(\cdot,\cu_{n+1}, p \right) \right|^2  \right] 
\\
\leq 
C \left( \left(\mathsf{K}+|p|\right)^2 3^{-\kappa n}+ \sum_{m=0}^n 3^{-\kappa (n-m)} \left( \tau_m^\ast \left( D \overline{L}_m(p)\right)  + \tau_m(p) \right) \right).
\end{multline}
To conclude, it suffices to show that, for every $p,q\in\Rd$,
\begin{multline}
\label{e.toconclude.NL}
\E \left[  J\left(\cu_n,p,q \right)\right] 
\\
\leq
C\E \left[ \fint_{\cu_{n}}  \left| \nabla u\left(\cdot,\cu_{n+1},q \right) - \nabla v\left(\cdot,\cu_{n+1}, p \right) \right|^2  \right] 
+ C \left( \tau_n^\ast(q)  + \tau_n(p) \right) .
\end{multline}
By~\eqref{e.Stokesie}, we have  
\begin{equation*} 
\fint_{\cu_n} \left( - L(\nabla v(\cdot,\cu_n,p),\cdot) + q\cdot \nabla v(\cdot,\cu_n,p)  \right) = - \nu(\cu_n,p ) + p\cdot q.
\end{equation*}
By~\eqref{e.nu*quadresponse.NL}, we obtain  
\begin{align} \notag 
J\left(\cu_n, p ,q\right) & = \nu(\cu_n,p) + \nu^*(\cu_n,p) - p\cdot q 
\\ \notag & =
\nu^*(\cu_n,p) - \fint_{\cu_n} \left( - L(\nabla v(\cdot,\cu_n,p),x) + q\cdot \nabla v(\cdot,\cu_n,p)  \right) 
\\ \notag & \leq 
\frac{\Lambda}{2} \fint_{\cu_n} \left| \nabla v(\cdot,\cu_n,p) - \nabla u(\cdot,\cu_n,q)\right|^2 . 
\end{align}
On the other hand, by quadratic response, we have 
\begin{equation*} \label{}
\fint_{\cu_n} \left| 
\nabla u\left(\cdot,\cu_{n+1},q\right)-
\nabla u\left(\cdot,\cu_{n},q\right)
\right|^2 
\leq C \sum_{z\in 3^n\Zd\cap \cu_{n+1}} \left( \nu^*(z+\cu_n,q) - \nu^*(\cu_{n+1},q) \right) 
\end{equation*}
and
\begin{equation*} \label{}
\fint_{\cu_n} \left| 
\nabla v\left(\cdot,\cu_{n+1},p\right)-
\nabla v\left(\cdot,\cu_{n},p\right)
\right|^2 
\leq C \sum_{z\in 3^n\Zd\cap \cu_{n+1}} \left( \nu(z+\cu_n,p) - \nu(\cu_{n+1},p) \right).
\end{equation*}
Combining these, applying the triangle inequality and taking expectations yields the claim~\eqref{e.toconclude.NL}, finishing the proof. 
\end{proof}

We now perform the iteration to complete the generalization of Proposition~\ref{p.subaddE}. 

\begin{proposition}
\label{p.subaddE.NL}
There exist~$\alpha(d,\Lambda)\in\left(0,\tfrac 12\right]$ and~$C(d,\Lambda)<\infty$ such that, for every~$p\in\Rd$ and~$n\in\N$,
\begin{equation} 
\label{e.subaddE.NL}
\E \left[  
J\left(\cu_n, p, D \overline{L}(p)\right)
\right] 
\leq 
C \left(\mathsf{K}+|p|\right)^2 3^{-n \alpha}.
\end{equation}
\end{proposition}
\begin{proof}
For each $p\in\Rd$, we define the quantity 
\begin{equation*} \label{}
F_n(p) := \E \left[ J\left(\cu_n,p,D \overline{L}_n(p)\right) \right]. 
\end{equation*}
Notice the difference between~$F_n$ defined here compared to the one defined in~\eqref{e.defDn}, which is that we are not tracking the convergence of all the~$p$'s at once, just a single~$p$. We also define
\begin{equation*} \label{}
\tilde{F}_n(p):= 3^{-\frac \kappa 2n} \sum_{m=0}^n 3^{\frac \kappa2m}  F_m(p),
\end{equation*}
where $\kappa(d)>0$ is the exponent in the statement of Lemma~\ref{l.iterstep.NL}. It is also convenient to define 
\begin{equation*} \label{}
T_n(p) := \tau_n^\ast \left( D \overline{L}_n (p) \right) + \tau_n(p).
\end{equation*}
Recall that~$\tau_n^\ast$ and  $\tau_n$ are defined in~\eqref{e.your.taus}.

\smallskip

\emph{Step 1.} We prove that
\begin{equation} 
\label{e.octopus}
F_n(p) - F_{n+1}(p) \geq T_n(p). 
\end{equation}
Since the map in~\eqref{e.mumap.NL} achieves its minimum at $q= D \overline{L}_U(p)$, we have 
\begin{align*}
F_n(p) - F_{n+1}(p)
& 
=  \E \left[ J\left(\cu_n,p,D \overline{L}_n(p)\right) \right] - \E \left[ J\left(\cu_{n+1},p,D \overline{L}_{n+1}(p) \right) \right]
\\ & 
\geq \E \left[ J\left(\cu_n,p,D \overline{L}_n(p)\right) \right] - \E \left[ J\left(\cu_{n+1},p,D \overline{L}_{n}(p)\right) \right]
\\ & 
= T_n(p). 
\end{align*}

\smallskip

\emph{Step 2.} We show that there exists $\theta(d,\Lambda) \in \left[\tfrac12,1\right)$ and $C(d,\Lambda)<\infty$ such that 
\begin{equation} 
\label{e.tilderocking}
\tilde{F}_{n+1}(p) 
\leq 
\theta \tilde{F}_n(p) + C\left(\mathsf{K}+|p|\right)^2 3^{-\frac\kappa 2n}. 
\end{equation}
This is essentially the same as the proof of~\eqref{e.tildeiterstep}. We have by~\eqref{e.octopus} and
\begin{equation} 
\label{e.basegiveups}
F_0(p)\leq C\left(\mathsf{K}+|p|\right)^2
\end{equation}
that
\begin{align*} \label{}
\tilde{F}_n(p) - \tilde{F}_{n+1} (p)
&
= 
3^{-\frac \kappa 2n} \sum_{m=0}^n 3^{\frac \kappa 2m} \left( F_m(p) - F_{m+1}(p) \right) - 
 3^{-\frac \kappa 2(n+1)}  F_{0} (p) 
\\ & 
\geq  
 3^{-\frac \kappa 2n} \sum_{m=0}^n 3^{\frac \kappa 2m} T_m(p) - C\left(\mathsf{K}+|p|\right)^2 3^{-\frac \kappa 2n} . 
\end{align*}
By a similar computation as in the display~\eqref{e.octopusalign}, we apply Lemma~\ref{l.iterstep.NL} to obtain 
\begin{align*} \notag
\tilde{F}_{n+1} (p)
 \leq C\left(\mathsf{K}+|p|\right)^2 3^{-\frac \kappa 2n} + C3^{-\frac \kappa 2n}  \sum_{m=0}^n 3^{\frac \kappa 2 m} T_m(p).
\end{align*}
Combining the previous two displays yields~\eqref{e.tilderocking}.

\smallskip

\emph{Step 3.} The conclusion. By an iteration of the inequality~\eqref{e.tilderocking} and using~\eqref{e.basegiveups}, we obtain, for some $\alpha(d,\Lambda)>0$,
\begin{equation*} \label{}
\tilde{F}_n (p)
\leq \theta^n \tilde{F}_0 (p) + C\left(\mathsf{K}+|p|\right)^23^{-\frac \kappa 2n}
\leq C\left(\mathsf{K}+|p|\right)^2 3^{-n\alpha} .
\end{equation*}
Since $F_n(p) \leq \tilde{F}_n(p)$, we deduce that
\begin{equation*} \label{}
F_n(p) \leq C\left(\mathsf{K}+|p|\right)^2 3^{-n\alpha} .
\end{equation*}
By~\eqref{e.octopus}, we also have
\begin{equation*} \label{}
T_n(p) \leq C\left(\mathsf{K}+|p|\right)^2 3^{-n\alpha}.
\end{equation*}
Hence by~\eqref{e.nablaLhom.comparescales}, we obtain that, for every $m,n\in\N$ with $n\leq m$, 
\begin{align*} 
\left| D \overline{L}_n(p) - D \overline{L}_m(p) \right|^2 
\leq 
\left( \sum_{k=n}^m \left| D \overline{L}_k(p) - D \overline{L}_{k+1} (p) \right|\right)^2
&
\leq
C \left( \sum_{k=n}^m T_k(p)^{\frac12} \right)^2
\\ & 
\leq C \left(\mathsf{K}+|p|\right)^2 \left( \sum_{k=n}^m 3^{-k\alpha/2} \right)^2
\\ & 
\leq C\left(\mathsf{K}+|p|\right)^2 3^{-n\alpha}. 
\end{align*}
For each $p \in \Rd$, we define $\bar L(p)$ by
\begin{equation}  
\label{e.alt.def.barL}
\bar L (p) := \lim_{n \to \infty} \bar L_n(p).
\end{equation}
We postpone the verification of the fact that this definition coincides with that given in \eqref{e.def.barL} to the next step. 
By the two previous displays and~\eqref{e.quadresponseLhom.NL}, we get
\begin{equation} 
\label{e.looks.like.concl.NL}
\E \left[ J\left( \cu_n,p,D \overline{L}(p)\right) \right] 
\leq F_n(p) + C\left|  D\overline{L}_n(p) -  D\overline{L}(p) \right|^2
\leq C\left(\mathsf{K}+|p|\right)^2 3^{-n\alpha}. 
\end{equation}

\smallskip

\emph{Step 4.} 
There remains to verify that the definitions in \eqref{e.alt.def.barL} and \eqref{e.def.barL} coincide. 
By \eqref{e.obvioustesting.NL}, we have, for every $p \in \Rd$,
\begin{equation}  
\label{e.obviousineq.NL}
\E \Ll[ \nu(\cu_n,p) \Rr]  \ge \bar L_n(p),
\end{equation}
and we recall that, for every $q \in \Rd$,
 \begin{equation*}  
 \bar L_n(p) = \sup_{q \in \Rd} \Ll(p\cdot q - \E \Ll[ \nu^*(\cu_n,q) \Rr] \Rr).
 \end{equation*}
Selecting $q = D \bar L(p)$ and using \eqref{e.looks.like.concl.NL}, we infer that
\begin{equation*}  
\bar L_n(p) \ge \E \Ll[ \nu(\cu_n,p) \Rr] - C \left(\mathsf{K}+|p|\right)^2  3^{-n\alpha}.
\end{equation*}
This and \eqref{e.obviousineq.NL} imply that \eqref{e.def.barL} holds.
\end{proof}
\begin{exercise}  
Show that for $\bar L\!\,^*$ defined in \eqref{e.def.barL*} and every $q \in \Rd$,
\begin{equation}  
\label{e.lim.barL*}
\bar L\!\,^*(q) = \lim_{n \to \infty} \E \Ll[ \nu^*(\cu_n,q) \Rr] .
\end{equation}
\end{exercise}

We are now ready to complete the proof of~Theorem~\ref{t.subadd.NL}. 

\begin{proof}[{Proof of Theorem~\ref{t.subadd.NL}}]
We decompose the proof into four steps. 


\smallskip

\emph{Step 1.} We show that, for every $p\in\Rd$ and $n \in \N$,
\begin{equation} 
\label{e.crazyeasy.p}
0 
\leq \E \left[ \nu(\cu_n,p) \right] - \overline{L}(p) 
\leq C \left(\mathsf{K}+|p|\right)^2 3^{-n\alpha}
\end{equation}
and, for every $q\in\Rd$ and $n \in \N$,
\begin{equation} 
\label{e.crazyeasy.q}
0 
\leq \E \left[ \nu^\ast(\cu_n,q) \right] - \overline{L}\!\,^*(q) 
\leq C \left(\mathsf{K}+|q|\right)^2 3^{-n\alpha}. 
\end{equation}
To check~\eqref{e.crazyeasy.p}, we observe that the leftmost inequality is an immediate consequence of subadditivity and stationarity.  
For the rightmost inequality, we observe that by subadditivity and stationarity, for every $p\in\Rd$,
\begin{align*} \label{}
\E\left[ \nu(\cu_n,p) \right]
&
= \E \left[ J\left(\cu_n,p,D \overline{L}(p)\right) \right]
+ p\cdot D \overline{L}(p) 
- \E \left[ \nu^*\left(\cu_n,  D\overline{L}(p) \right) \right]
\\ & 
\leq \E \left[ J\left(\cu_n,p,D \overline{L}(p)\right) \right]
+ p\cdot D \overline{L}(p) 
- \overline{L}\!\,^* \left(  D\overline{L}(p) \right)
\\ & 
= \E \left[ J\left(\cu_n,p,D \overline{L}(p)\right) \right]
+ \overline{L}(p). 
\end{align*}
Together with Proposition~\ref{p.subaddE.NL}, this completes the proof of~\eqref{e.crazyeasy.p}. The proof of~\eqref{e.crazyeasy.q} is similar using \eqref{e.lim.barL*}, so we omit it.

\smallskip

\emph{Step 2.} We upgrade the stochastic integrability of~\eqref{e.crazyeasy.p} and~\eqref{e.crazyeasy.q}. The claim is that there exists $C(d,\Lambda)<\infty$ such that, for every $p,q \in \Rd$ and $m,n\in\N$ with $m\leq n$,
\begin{equation}
\label{e.silyyasyp}
\left| \nu(\cu_n,p) - \overline{L}(p) \right| 
\leq C\left(\mathsf{K}+|p|\right)^2  3^{-m\alpha} + \O_1\left( C\left(\mathsf{K}+|p|\right)^2 3^{-d(n-m)} \right)
\end{equation}
and
\begin{equation}
\label{e.silyyasyq}
\left|\nu^*(\cu_n,q) - \overline{L}\!\,^*(q) \right| 
\leq
C\left(\mathsf{K}+|q|\right)^2 3^{-m\alpha} + \O_1\left( C\left(\mathsf{K}+|q|\right)^2 3^{-d(n-m)} \right).
\end{equation}
Fix $p\in\Rd$ and compute for each $m,n\in\N$ with $m\leq n$ and $\mathcal{Z}_m:=3^m\Zd\cap \cu_n$, using subadditivity and Lemmas~\ref{l.change-s} and~\ref{l.barO.boxes}, to obtain
\begin{align*}
\nu(\cu_n,p) - \E \left[ \nu(\cu_m,p) \right] 
&
\leq
\frac1{\left|\mathcal{Z}_m\right|} \sum_{z\in \mathcal{Z}_m} 
\left( 
\nu(z+\cu_m,p) - \E \left[ \nu(z+\cu_m,p) \right] \right)
\\ &
\leq 
\O_2 \left( C\left(\mathsf{K}+|p|\right)^2 3^{-\frac d2(n-m)} \right)\wedge C\left(\mathsf{K}+|p|\right)^2  
\\ &
\leq 
\O_1\left( C\left(\mathsf{K}+|p|\right)^2  3^{-d(n-m)} \right).
\end{align*}
By~\eqref{e.crazyeasy.p}, we obtain
\begin{equation}
\label{e.silyyasypup}
\nu(\cu_n,p) - \overline{L}(p) 
\leq C\left(\mathsf{K}+|p|\right)^2 3^{-m\alpha} + \O_1\left( C\left(\mathsf{K}+|p|\right)^2 3^{-d(n-m)} \right).
\end{equation}
Similarly, for every $q\in\Rd$, 
\begin{equation*}
\nu^*(\cu_n,q) - \overline{L}\!\,^*(q) 
\leq
C\left(\mathsf{K}+|q|\right)^23^{-m\alpha} + \O_1\left( C\left(\mathsf{K}+|q|\right)^23^{-d(n-m)} \right).
\end{equation*}
By~\eqref{e.obvioustesting.NL} with~$q = D_p \overline{L}(p)$ and the previous display, we have, for every $p \in \Rd$, 
\begin{align*}
\overline{L}(p) - \nu(\cu_n,p) 
&
\leq 
\nu^*(\cu_n,D_p \overline{L}(p))
- \left( p\cdot  D_p \overline{L}(p) - \overline{L}(p) \right)
\\ & 
=
\nu^*(\cu_n,D_p \overline{L}(p))
-
\overline{L}\!\,^*(D_p \overline{L}(p))
\\ & 
\leq 
C\left(\mathsf{K}+|p|\right)^2 3^{-m\alpha} 
+ \O_1\left( C\left(\mathsf{K}+|p|\right)^2 3^{-d(n-m)} \right).
\end{align*}
Combining this inequality with~\eqref{e.silyyasypup} yields~\eqref{e.silyyasyp}. Arguing similarly and using that~$\overline{L}\!\,^{**} = \overline{L}$,  we also obtain~\eqref{e.silyyasyq}. 

\smallskip

\emph{Step 3.} We complete the proof of the theorem. Fix $\mathsf{M}\in [1,\infty)$. By~\eqref{e.Lnormalize} and~\eqref{e.nuderp.NL}, we have
\begin{equation*}
\sup_{p\in B_{\mathsf{M}}} 
\left( \left| D_p \nu(U,p) \right| + \left| D_p \overline{L}(p) \right| \right) \leq C\left( \mathsf{K} + \mathsf{M} \right).
\end{equation*}
Thus, for every $L>0$ and $m,n\in\N$ with $m\leq n$, 
\begin{align*}
\lefteqn{
\sup_{p\in B_{\mathsf{M}}}  
\left| \nu(\cu_n,p) - \overline{L}(p) \right| 
} \quad & 
\\ &
\leq C \left( \mathsf{K} + \mathsf{M} \right)L 
+ \sum_{p\in L\Zd\cap B_{\mathsf{M}}} 
\left| \nu(\cu_n,p) - \overline{L}(p) \right|
\\ & 
\leq C \left( \mathsf{K} + \mathsf{M} \right)L 
+ \left| L\Zd\cap B_{\mathsf{M}} \right| \cdot\left( C\left( \mathsf{K} + \mathsf{M} \right)^23^{-m\alpha} + \O_1\left( C\left( \mathsf{K} + \mathsf{M} \right)^2 3^{-d(n-m)} \right) \right) 
\\ & 
= C \left( \mathsf{K} + \mathsf{M} \right)L + C\left( \mathsf{K} + \mathsf{M} \right)^2 \left( \frac{ \mathsf{M}}{L} \right)^d 3^{-m\alpha}
+ \O_1\left(  C\left( \mathsf{K} + \mathsf{M} \right)^2 \left( \frac{ \mathsf{M}}{L} \right)^d 3^{-d(n-m)} \right).
\end{align*}
Taking $L:= \left( \mathsf{K} + \mathsf{M} \right) 3^{-m\alpha/2d}$ yields
\begin{equation*}
\sup_{p\in B_{\mathsf{M}}}  
\left| \nu(\cu_n,p) - \overline{L}(p) \right| 
\leq 
C\left( \mathsf{K} + \mathsf{M} \right)^2 3^{-m\alpha/2d}
+ \O_1\left(  C\left( \mathsf{K} + \mathsf{M} \right)^2  3^{m\alpha/2-d(n-m)} \right).
\end{equation*}
Given~$s\in (0,d)$, the choice of $m$ to be the closest integer to $(d-s)n/(d+\alpha/2)$ and then shrinking~$\alpha$ yields the desired estimate for the first term in~\eqref{e.subadderror.NL}. The estimate for the second term follows similarly. 
\end{proof}

\section{Quantitative homogenization of the Dirichlet problem}
\label{s.Ldirichlet}

In this section, we prove an estimate on the homogenization error for general Dirichlet boundary value problems. This is a generalization of Theorem~\ref{t.DP} to the nonlinear setting. The proof is based on an adaptation of the two-scale expansion argument used in the proof of Theorem~\ref{t.DP.blackbox}.

\smallskip

\begin{theorem}
\label{t.DP.nonlinear}
Fix $s \in (0,d)$, $\delta>0$, $\mathsf{M} \in [1,\infty)$, and a bounded Lipschitz domain~$U\subseteq \cu_0$. There exist constants $\alpha(\delta,d,\Lambda),\gamma(\delta,U,d,\Lambda)>0$ and $C(\delta,U,s,d,\Lambda)<\infty$, and a random variable $\X_{s} : \Omega \to [0,\infty)$ satisfying
\begin{equation*} 
\X_{s} \leq \O_{1}\left( C (\mathsf{K} + \mathsf{M}) \right)
\end{equation*}
such that the following statement holds.
For every~$\ep \in \left( 0,\tfrac12\right]$ and~$f\in W^{1,2+\delta}(U)$ satisfying
\begin{equation}  \label{e.fnormalize.NL}
\left\| \nabla f \right\|_{\underline{L}^{2+\delta \wedge \gamma}(U)} \leq \mathsf{M},
\end{equation}
and solutions~$u^\ep, u \in f+H^1_0(U)$ of the Dirichlet problems 
\begin{equation}
\label{e.uep.and.ubar.NL}
\left\{
\begin{aligned}
& -\nabla \cdot \left( D_{p}L \left(\nabla u^\ep , \tfrac \cdot \ep\right)  \right) = 0 &  \mbox{in} & \ U, \\
& u^\ep = f & \mbox{on} & \ \partial U,
\end{aligned}
\right.
 \quad 
\left\{
\begin{aligned}
& -\nabla \cdot D\overline{L}\left( \nabla  u\right)   = 0 &  \mbox{in} & \ U, \\
&  u  = f & \mbox{on} & \ \partial U,
\end{aligned}
\right.
\end{equation}
we have the estimate
\begin{multline} \label{e.DP.nonlinear}
\left\| \nabla u^\ep - \nabla u \right\|_{H^{-1}(U)} 
+ 
\left\| D_{p} L\left( \nabla u^\ep , \tfrac \cdot\ep\right) - D\overline{L}\left( \nabla u \right) \right\|_{H^{-1}(U)} 
\\ \leq 
C (\mathsf{K} + \mathsf{M}) \ep^{\alpha(d-s)} 
+  \X_{s} \ep^s .
\end{multline}
\end{theorem}

We begin the proof of Theorem~\ref{t.DP.nonlinear} by using Theorem~\ref{t.subadd.NL} to obtain bounds on finite-volume correctors. These are built on the functions~$v(\cdot,\cu_m,p)$ which are the minimizers of~$\nu(\cu_m,p)$, in other words, the solutions of the Dirichlet problem in~$\cu_m$ with affine boundary data~$\ell_p$.

\begin{lemma}
\label{l.approx.correctors.NL}
Fix $s\in (0,d)$ and $\mathsf{M} \in [1,\infty)$. There exist constants $\alpha(d,\Lambda)>0$ and~$C(s,d,\Lambda)<\infty$ such that, for every~$p\in\Rd$ and~$m\in\N$
\begin{multline} 
\label{e.approx.correctors.NL}
3^{-m}\sup_{p \in B_\mathsf{M}} \left( \left\| \nabla v(\cdot,\cu_m,p) - p \right\|_{\Hminusul(\cu_m)}
+  \left\| D_pL(\nabla v(\cdot,\cu_m,p),\cdot) - D\overline{L}(p) \right\|_{\Hminusul(\cu_m)} \right)
\\
\leq C  (\mathsf{K} + \mathsf{M}) 3^{-m\alpha(d-s)}  + \O_1\left( C    (\mathsf{K} + \mathsf{M})  3^{-ms} \right) .
\end{multline}
\end{lemma}
\begin{proof}
We follow the outline of the proof of Proposition~\ref{p.qualweakconv} using Theorem~\ref{t.subadd.NL}. Fix $s\in (0,d)$. We recall from \eqref{e.subadderror.NL} that, for every $k \in \N_0$ such that $k \leq m$, $z \in 3^k \Z^d$ and $p,q \in \R^d$, 
\begin{equation} 
\label{e.subadderror.NL.applied}
\sup_{p \in B_\mathsf{M}}\left| \nu(z + \cu_k,p) - \overline{L}(p) \right|   \leq C \left( \mathsf{K}  +  \mathsf{M} \right)^2  3^{-n \alpha(d-s)} + \O_1\left( C \left( \mathsf{K}  +  \mathsf{M} \right)^23^{-ns} \right) . 
\end{equation}

\smallskip

\emph{Step 1.} We show that there exist~$\beta (d,\Lambda)>0$ and~$C(s,d,\Lambda)<\infty$ such that, for every~$m \in \N$, we have
\begin{equation}  \label{e.approx.correctors.NL.goal1}
3^{-m} \sup_{p \in B_\mathsf{M}} \left\| \nabla v(\cdot,\cu_m,p) - p \right\|_{\Hminusul(\cu_m)} \leq C \left( \mathsf{K}  +  \mathsf{M} \right)  3^{-m \beta(d-s)} + \O_1\left( C \left( \mathsf{K}  +  \mathsf{M} \right) 3^{-ms}  \right). 
\end{equation}
Denote, in short,
\begin{equation*} 
v(x) := v(x,\cu_m,p), \quad v_n(x) :=  v(x, z + \cu_n,p), \quad x \in z + \cu_n,
\end{equation*}
and
\begin{equation*} 
\mathcal{Z}_n := 3^n \Z^d \cap \cu_m \quad \mbox{and} \quad \left| \mathcal{Z}_n\right| := 3^{d(m-n)}.
\end{equation*}
Then, since $ \left( \nabla v - p \right)_{z +\cu_{n} } =  \left( \nabla v -\nabla  v_n  \right)_{z +\cu_{n} }$,  we have by the multiscale Poincar\'e inequality (Proposition~\ref{p.mspoin}) that
\begin{equation*} 
 \| \nabla v - p  \|_{\Hminusul(\cu_m)}
\leq C \left\| \nabla v- p    \right\|_{\underline{L}^2(\cu_m)}  
+  C \sum_{n=0}^{m-1} 3^n    \left(  \frac{1}{|\mathcal{Z}_n|}\sum_{z \in \mathcal{Z}_n }\left| \left( \nabla v -\nabla  v_n \right)_{z +\cu_{n} }  \right|^2 \right)^{\frac12}.
\end{equation*}
The first term is bounded by $C(1+|p|)$. For the second term, we get from Jensen's inequality that 
\begin{equation*} \label{}
 \frac{1}{|\mathcal{Z}_n|} \sum_{z\in \mathcal{Z}_n} 
 \left| \left( \nabla v -\nabla  v_n  \right)_{z +\cu_{n} } \right|^2 
  \leq 
   \left\| \nabla v - \nabla v_n \right\|_{\underline{L}^2(\cu_m)}^2.
\end{equation*}
As in the case of~\eqref{e.vvtildesnap}, we have by~\eqref{e.nuquadresponse.NL} that  
\begin{align*} 
\left\| \nabla v - \nabla v_n \right\|_{\underline{L}^2(\cu_m)}^2
\leq 
 \frac{2}{|\mathcal{Z}_n|} \sum_{z\in \mathcal{Z}_n} \left( \nu(z+\cu_n,p) - \nu(\cu_m,p) \right),
\end{align*}
and hence
\begin{equation*} \label{}
 \frac{1}{|\mathcal{Z}_n|} \sum_{z\in \mathcal{Z}_n} 
 \left| \left( \nabla v \right)_{z+\cu_n}  - p \right|^2 
  \leq 
\frac{2}{|\mathcal{Z}_n|} \sum_{z\in \mathcal{Z}_n} \left( \nu(z+\cu_n,p) - \nu(\cu_m,p) \right).  
\end{equation*}
Thus~\eqref{e.subadderror.NL.applied} and Lemma~\ref{l.sum-O} yields that 
\begin{multline*} 
\sup_{p \in B_{\mathsf{M}}} \frac{1}{|\mathcal{Z}_n|} \sum_{z\in \mathcal{Z}_n} 
 \left| \left( \nabla v(\cdot, \cu_m,p) \right)_{z+\cu_n}  - p \right|^2 
 \\
\leq 
C \left( \mathsf{K}  +  \mathsf{M} \right)^2   3^{-n \alpha(d-s')} + \O_1\left(C \left( \mathsf{K}  +  \mathsf{M} \right)^2  3^{- s' n} \right) . 
\end{multline*}
Using~\eqref{e.sqrtelem} as in the proof of Theorem~\ref{t.subadd} gives that 
\begin{equation*} 
\sup_{p \in B_{\mathsf{M}}} \left( \frac{1}{|\mathcal{Z}_n|} \sum_{z\in \mathcal{Z}_n} 
 \left| \left( \nabla v \right)_{z+\cu_n}  - p \right|^2 \right)^{\frac12}
\leq C  \left( \mathsf{K}  +  \mathsf{M} \right)   3^{-n \beta(d-s)} + \O_1\left( C  \left( \mathsf{K}  +  \mathsf{M} \right)  3^{- s n} \right) . 
\end{equation*}
Applying Lemma~\ref{l.sum-O} once more yields~\eqref{e.approx.correctors.NL.goal1}. 

\smallskip

\emph{Step 2.}
We show that there exists $\beta (d,\Lambda)>0$ and $C(d,\Lambda)<\infty$ such that 
\begin{multline}  \label{e.approx.correctors.NL.goal2}
3^{-m}  \sup_{p \in B_\mathsf{M}}\| D_pL(\nabla v,\cdot)  - D\overline{L}(p)   \|_{\Hminusul(\cu_m)} 
\\ \leq C  (\mathsf{K} + \mathsf{M}) 3^{-m\beta(d-s)}  + \O_1\left( C (\mathsf{K} + \mathsf{M})3^{-ms} \right) .
\end{multline}
Applying the multiscale Poincar\'e inequality yields
\begin{multline*} 
 \| D_pL(\nabla v,\cdot)  - D\overline{L}(p)   \|_{\Hminusul(\cu_m)}
 \leq C \left\| D_pL(\nabla v,\cdot)  - D\overline{L}(p)      \right\|_{\underline{L}^2(\cu_m)}  
\\+  C \sum_{n=0}^{m-1} 3^n    \left(  \frac{1}{|\mathcal{Z}_n|}\sum_{z \in \mathcal{Z}_n }\left| \left( D_pL(\nabla v,\cdot) \right)_{z +\cu_{n} } - D\overline{L}(p)  \right|^2 \right)^{\frac12}.
\end{multline*}
On the one hand, using~\eqref{e.Lnormalize} we get 
\begin{equation*} 
\left\| D_pL(\nabla v,\cdot)  - D\overline{L}(p)      \right\|_{\underline{L}^2(\cu_m)}   \leq C(\mathsf{K}+|p|) . 
\end{equation*}
On the other hand, by~\eqref{e.Lunifconv5}, we obtain
\begin{multline*} 
\left| \left( D_pL(\nabla v,\cdot) \right)_{z +\cu_{n} } - D\overline{L}(p)  \right|^2
\\
\leq C\left\|\nabla v - \nabla v_n \right\|_{\underline{L}^2(z + \cu_n)}^2 +  2 \left| \left( D_pL(\nabla v_n,\cdot) \right)_{z +\cu_{n} } - D\overline{L}(p)\right|^2 .
\end{multline*}
The first term on the right can be estimated exactly as in the first step. For the second term we instead use~\eqref{e.nuderp.NL}, and the inequality~\eqref{e.quant.diffnu.NL} from the proof of Lemma~\ref{l.basicJ.NL}.  We have that 
\begin{equation*} 
\left(D_pL(\nabla v_n,\cdot) \right)_{z +\cu_{n} }  = D_p \nu(z+\cu_n,p). 
\end{equation*}
Applying~\eqref{e.quant.diffnu.NL}, we get, for $z \in 3^n \Z^d$, 
\begin{equation*} 
\frac12 \left| p-p' \right|^2 \leq  \nu(z+\cu_n,p) - \nu(z+\cu_n,p') - \left( p -p' \right) \cdot D_p \nu(z+\cu_n,p) \leq C \left| p-p' \right|^2. 
\end{equation*}
and, by taking expectation and sending $n \to \infty$ for $z=0$,
\begin{equation*} 
\frac12 \left| p-p' \right|^2 \leq  \overline{L}(p) - \overline{L}(p') - \left( p -p' \right) \cdot D\overline{L}(p) \leq C \left| p-p' \right|^2. 
\end{equation*}
Applying this with $p' = p \pm \delta e_k$, for $\delta>0$, we get 
\begin{equation*} 
\left|D_p \nu(z+\cu_n,p) -  D\overline{L}(p) \right| \leq C \left( \delta +  \delta^{-1}  \sum_{e \in \{-1,0,1\}^d }  \left| \nu(\cu_n,p+\delta e) -  \overline{L}(p + \delta e)\right| \right).
\end{equation*}
We choose $\delta = 3^{-n \frac\alpha4 (d-s)}(\mathsf{K}+\mathsf{M})$ and obtain by~\eqref{e.subadderror.NL.applied} (applied with $s = \frac{s'+d}{2}$) that
\begin{multline*}
\sup_{p \in B_\mathsf{M}}\left( \delta +  \delta^{-1}  \sum_{e \in \{-1,0,1\}^d }  \left| \nu(\cu_n,p+\delta e) -  \overline{L}(p + \delta e)\right| \right)^2
 \\ \leq C \left( \mathsf{K}  +  \mathsf{M} \right)^2   3^{-n\frac\alpha2 (d-s')} + \O_{1}\left( C \left( \mathsf{K}  +  \mathsf{M} \right)^2   3^{- s' n }  \right). 
\end{multline*}
The rest of the proof is analogous to Step 1. The proof is complete. 
\end{proof}

We next record necessary basic regularity estimates for $\overline{L}$-minimizers. 

\begin{lemma} \label{l.refreg.NL}
Fix $\delta \in (0,\infty]$ and a Lipschitz domain $U$. Let $f \in W^{1,2+\delta}(U)$.
Suppose that $u \in f + H_0^1(B_R)$ is $\overline{L}$-minimizer in $B_R$.  Then there exist constants $\gamma(U,d,\Lambda)>0$ and $C(U,d,\Lambda)< \infty$ such that 
\begin{equation} \label{e.meyers.NL}
 \left\| \nabla u \right\|_{\underline{L}^{2+\delta \wedge \gamma}(U)} 
 \leq 
 C  \left\| \nabla f \right\|_{\underline{L}^{2+\delta \wedge \gamma}(U)}  .
\end{equation}
Moreover, there exist $\beta(d,\Lambda)>0$ and $C(d,\Lambda)< \infty$ such that if $B_R \subset U$, then 
 \begin{equation} \label{e.refreg.NL}
   \left\| \nabla u - (\nabla u)_{B_{R/2}}  \right\|_{L^\infty(B_{R/2})} + 
R^{\beta} \left[ \nabla u \right]_{C^\beta(B_{R/2})}   \leq  \frac{C}{R} \inf_{p \in \mathcal{P}_1 }\left\| u - p \right\|_{\underline{L}^{2}(B_R)} .
\end{equation}
\end{lemma}

\begin{proof}
The proof of~\eqref{e.meyers.NL} is analogous to the one of Theorem~\ref{t.Meyers appendix}. Indeed, by defining 
\begin{equation*} 
\mathbf{b}(p) := D \overline{L}(p) - D \overline{L}(0),
\end{equation*}
we see, by the first variation, that $\nabla \cdot \mathbf{b}(\nabla u) = 0$. We thus get a Caccioppoli inequality, for $B_r(y) \subset U$,
\begin{equation*} 
\left\| \nabla u \right\|_{\underline{L}^{2}(B_{r/2}(y))} \leq C  \frac1r \left\| u - (u)_{B_{r}(y)} \right\|_{\underline{L}^{2}(B_{r})},
\end{equation*}
and at the boundary, that is, for $B_r(y)$ such that  $B_{r/2}(y) \cap U \neq \emptyset$,
\begin{equation*} 
\left\| \nabla (u-f) \right\|_{\underline{L}^{2}(B_{r/2}(y) \cap U)} \leq C \left( \frac1r \left\| u -f \right\|_{\underline{L}^{2}(B_{r}(y) \cap U)}  + \left\| \nabla f \right\|_{\underline{L}^{2}(B_{r}(y) \cap U)}  \right).
\end{equation*}
We consequently deduce a reverse H\"older inequality using Sobolev-Poincar\'e, as in Corollary~\ref{c.Meyers reverse}, and then~\eqref{e.meyers.NL} follows by Gehring's lemma, Lemma~\ref{l.Meyers really}. 

\smallskip 

The regularity result~\eqref{e.refreg.NL} is standard, if $\overline{L}$ would be smooth. Indeed, then each component of $\nabla u$ would satisfy a uniformly elliptic equation with coefficients $D^2\overline{L}(\nabla u)$, which by the uniform ellipticity would be bounded from above and below. Thus one may apply De Giorgi-Nash-Moser theory to obtain~\eqref{e.refreg.NL}. In our case, as $\overline{L}$ is only $C^{1,1}$, the above reasoning is not, strictly speaking, possible. There are several ways to overcome this difficulty. One, which we choose to do here, is to mollify $\overline{L}$ and then obtain~\eqref{e.refreg.NL} for the regularized solution with boundary values $u$. Indeed, we have, for $\delta>0$ and $\zeta_\delta$ defined in~\eqref{e.standardmollifier}, that, for $z \in \R^d$, 
\begin{equation*} 
C^{-1} \Id \leq D^2 \left(\zeta_\delta \ast \overline{L}\right)(z) \leq C \Id .
\end{equation*}
We then pass to the limit by Arzel\'a-Ascoli theorem in the equation using $C^1$ regularity of $D\overline{L}$. By the uniqueness, the original solution $u$ inherits the regularity estimate~\eqref{e.refreg.NL} from the approximating sequence. The proof is complete. 
\end{proof}

We turn to the proof of the quantitative homogenization result, Theorem~\ref{t.DP.nonlinear}.

\begin{proof}[Proof of Theorem~\ref{t.DP.nonlinear}] Fix $\ep\in \left(0,\tfrac12\right]$, $s \in (0,d)$, $\delta>0$, $\mathsf{M} \in [1,\infty)$, and a bounded Lipschitz domain~$U\subseteq \cu_0$. Fix also~$f\in W^{1,2+\delta}(U)$ satisfying~\eqref{e.fnormalize.NL}. We divide the proof into eight steps. As usual, we may restrict to the case that $\ep \leq c$ for any~$c(\delta,U,s,d,\Lambda) >0$.

\smallskip

\emph{Step 1.} We setup the argument. 
Fix $m \in \N$ such that $\ep \in \left[ 3^{-m}, 3^{-m+1} \right)$. We also fix a mesoscale $k , \ell \in \N$ with $k < \ell < m$ and 
\begin{equation*} 
m - \ell \leq \ell - k. 
\end{equation*}
These mesoscales will be selected below in terms of $m$. We also set $r := 3^{\ell-m}$ and
\begin{equation*} 
U_r := \left\{ x \in U \, : \,  \dist(x , \partial U) > r \right\},
\end{equation*}
and let $\zeta_r $ be a smooth cut-off function such that $\zeta_r = 1$ in $U_{3r}$ and it vanishes outside of $U_{2r}$, and 
\begin{equation}  \label{e.gradzetabounded.NL}
r \| \nabla \zeta_r \|_{L^\infty(\R^d)} + r^2  \| \nabla^2 \zeta_r \|_{L^\infty(\R^d)}  \leq C \indc_{U_{3r} \setminus U_{2r}}.
\end{equation}
We also define, using standard mollifier defined in~\eqref{e.standardmollifier}, 
\begin{equation*} 
\psi(x):= \left( \zeta_{3^{k}} \ast \indc_{\cu_{k}}\right)(x),
\end{equation*}
so that 
\begin{equation*}
\sum_{z\in 3^{k} \Zd} \psi(\cdot-z) = 1. 
\end{equation*}
Write $\psi_z^\ep := \psi\left(\tfrac{\cdot-z}{\ep}\right)$. We clearly have that 
\begin{equation} \label{e.gradpsibounded.NL}
3^k \ep  \|\nabla \psi_z^\ep \|_{L^\infty(\R^d)} + (3^{k} \ep)^2 \|\nabla^2 \psi_z^\ep \|_{L^\infty(\R^d)} \leq C.
\end{equation}
Furthermore, for $z \in 3^{k} \Z^d$ and $p \in\R^d$, we denote in short
\begin{equation*} 
v_z(x,p) := v\left(x , z + \cu_{k+2} , p \right) - \left( v\left(x , z + \cu_{k+2} , p \right) \right)_{ z + \cu_{k+2} } 
\end{equation*}
and
\begin{equation*} 
 \phi_z(x,p) := v_z(x,p) - p \cdot (x-z) .
\end{equation*}
Then both $v_z$ and $\phi_{z}$ have zero mean in $z + \cu_{k+2}$.  We also set, for $z \in 3^{k} \Z^d$, 
\begin{equation*} 
v_z^\ep(x,p) := v_z\left( \tfrac{x}{\ep},p \right) 
\quad  \mbox{and} \quad 
\phi_z^\ep(x,p) := \phi_z\left( \tfrac{x}{\ep},p \right)  .
\end{equation*}
For simplicity, we use the short-hand notation, for each~$y \in 3^{k}\Z^d$, 
\begin{equation*} 
\sum_{z} =  \sum_{z \in 3^{k} \Z^d } 
\quad \mbox{and} \quad \sum_{z \sim y} =  \sum_{z \in 3^k \Z^d \cap \left( y + \cu_{k+2}\right)}. 
\end{equation*}
We define the two-scale expansion of $u$ as 
\begin{equation*}
w^\ep(x)
:=
u(x)
+
\ep \zeta_r(x) \sum_{z} \psi_z^\ep(x) \phi_z^\ep \left(x,\nabla u( \ep z)\right) .
\end{equation*}

\smallskip 

\emph{Step 2.} We collect necessary regularity estimates for $u^\ep$ and $u$. 
By Meyers estimate~\eqref{e.meyers.NL} we have that there exist~$\gamma(d,\Lambda)>0$ and~$C(d,\Lambda)<\infty$ such that 
\begin{equation}  \label{e.regumeyers.NL}
  \left\| \nabla u^\ep \right\|_{\underline{L}^{2+\delta \wedge \gamma}(U)} 
  + \left\| \nabla u \right\|_{\underline{L}^{2+\delta \wedge \gamma}(U)} 
 \leq 
 C\left( \mathsf{K} + \mathsf{M} \right) . 
\end{equation}
The proof for $u^\ep$ is similar to the one for~\eqref{e.meyers.NL}, but now using also~\eqref{e.Lnormalize}.  On the other hand, for any $y \in U_{2r}$, we have by~\eqref{e.refreg.NL} that there exist constants $\beta(d,\Lambda)>0$ and $C(d,\Lambda)<\infty$ such that 
\begin{equation}  \label{e.regutwoscale.NL}
\left\| \nabla u \right\|_{L^\infty(U_{r})} + r^{\beta} \left[ \nabla u \right]_{C^\beta(U_r)} 
 \leq   C  r^{-\frac d2}\mathsf{M}  =:  \mathsf{M}_r. 
\end{equation}

\smallskip 

\emph{Step 3.} 
We define a quantity, which measures all the errors appearing in the estimates below. Set, for $s' = \frac{d+s}{2}$, 
\begin{equation*} 
\mathcal{E} (m) :=  (\mathsf{K} + \mathsf{M})
\left( 3^{-\beta(\ell-k) + \frac{d}{2}(m-\ell)}   +  (3^{\ell - m})^{\frac{\delta \wedge \gamma}{2+\delta \wedge \gamma}} \right)
+  \mathcal{E}_{1}(m)  + \mathcal{E}_{2}(m) + \mathcal{E}_{3}(m) ,
\end{equation*}
where the random variables $ \mathcal{E}_{1}(m)$, $ \mathcal{E}_{2}(m)$, and $ \mathcal{E}_{3}(m)$ are defined as
\begin{equation*} 
\mathcal{E}_{1}(m) :=  \sum_{j=k}^{k+2} 
\left( 3^{(j-n)d} \sum_{z \in  3^{j} \Z^d \cap \cu_{n} }\sup_{p \in B_{  \mathsf{M}_r} } \left| \nu(z + \cu_j,p) - \overline{L}(p) \right| \right)^{\frac 12},
\end{equation*}
\begin{equation*} 
\mathcal{E}_{2}(m) := 
 \sum_{j=k}^{k+2} 
\left( 3^{(j-n)d -j}   \sum_{z \in  3^{j} \Z^d \cap \cu_{n} }  \sup_{p \in B_{  \mathsf{M}_r} }    \left\| \nabla \phi_z \left(\cdot, p\right)   \right\|_{\underline{H}^{-1}(z + \cu_{j})}^2  \right)^{\frac 12},
\end{equation*}
and
\begin{equation*} 
\mathcal{E}_{3}(m) := 
 \sum_{j=k}^{k+2} 
\left( 3^{(j-n)d - j}   \sum_{z \in  3^{j} \Z^d \cap \cu_{n} }   \sup_{p \in B_{  \mathsf{M}_r} } 
  \left\|  D_pL\left( \nabla v_z \left(\cdot, p \right) , \cdot \right)- D \overline{L} \left( p  \right) \right\|_{\underline{H}^{-1}(z + \cu_{j})}^2  \right)^{\frac 12} .
\end{equation*}
In Step~9 we show that the parameters can be selected so that the error above becomes small.

\smallskip 

\emph{Step 4.}  We show that the gradient of $w^\ep$ can be written as
\begin{equation} \label{e.gradw.nonlinear}
\nabla w^\ep(x)
 =   \zeta_{r}(x) \sum_{y \in 3^k \Z^d \cap \ep^{-1} U_{2r}} \nabla v_y^\ep \left(x,\nabla u( \ep y) \right) \indc_{ \ep (y + \cu_{k})}(x)  + E(x),
\end{equation}
where~$E(x)$ is an error defined below in~\eqref{e.Exterm}.  To see this, observe first that by a direct computation we get
\begin{equation*} 
\nabla w^\ep(x)  =
 \zeta_r(x) \sum_{z} \psi_z(x) \nabla v_z^\ep \left(x,\nabla u( \ep z)\right)  + E_1(x) ,
\end{equation*}
where
\begin{multline} \label{e.E1.NL}
E_1(x)  :=  (1-\zeta_r(x))  \nabla u(x) 
\\ 
+ \zeta_r(x) \sum_{z} \psi_z^\ep(x) \left( \nabla u(x) - \nabla u( \ep z)\right) 
+ \ep  \sum_{z} \nabla ( \zeta_r \psi_z^\ep)(x)  \phi_z^\ep(x,\nabla u( \ep z)) .
\end{multline}
Furthermore, for $x \in \ep (y + \cu_k)$, taking $y \in 3^{k} \Z^d $ such that $\ep y \in U_{2r} $, we see, for $x \in \ep (y + \cu_k)$, that 
\begin{multline*} 
 \zeta_r(x) \sum_{z} \psi_z^\ep (x)  \nabla v_z^\ep\left(x , \nabla u( \ep z) \right)
\\ =    \zeta_r(x)\sum_{y \in 3^k \Z^d \cap \ep^{-1} U_{2r}} \left( \nabla v_y^\ep \left(x,\nabla u( \ep y) \right)   + E_2(x,y) \right) \indc_{ \ep (y + \cu_{k})}(x),
\end{multline*}
where 
\begin{multline} \label{e.E2.NL}
 E_2(x,y) =  \zeta_r(x) \sum_{z \sim y} \psi_{z}^\ep(x)  \left( \nabla v_z^\ep \left(x,\nabla u( \ep z) \right)  - \nabla v_z^\ep \left(x,\nabla u( \ep y) \right) 
 \right) 
\\ + \zeta_r(x) \sum_{z \sim y} \psi_{z}^\ep(x)  \left( \nabla v_z^\ep \left(x,\nabla u( \ep y) \right)  - \nabla v_y^\ep \left(x,\nabla u( \ep y) \right) \right).
\end{multline}
Indeed, notice that $\sum_{z \sim y} \psi_{z}^\ep(x) = 1$ for $x \in \ep (y + \cu_k)$.
 This proves~\eqref{e.gradw.nonlinear} with 
 \begin{equation}
 \label{e.Exterm}
E(x) := E_{1}(x) + \sum_{y \in 3^k \Z^d \cap \ep^{-1} U_{2r}} \indc_{ \ep (y + \cu_{k})}(x) E_{2}(x,y). 
\end{equation}

\smallskip

\emph{Step 5.}
In this step we show that $E$ in~\eqref{e.gradw.nonlinear} is small in $L^2(U)$ norm, namely
\begin{equation}  \label{e.Eestimate.NL}
\left\| E \right\|_{L^2(U)}  \leq C\mathcal{E} (m).
\end{equation}
We first estimate $E_{1}$ defined in~\eqref{e.E1.NL}, and claim that
\begin{equation} \label{e.E1estimate.NL}
\left\| E_{1} \right\|_{L^2(U)}  \leq C \mathcal{E}(m). 
\end{equation}
For the first term in $E_{1}$ we use the Meyers estimate in~\eqref{e.regumeyers.NL}, and get
\begin{equation*} 
\left\| (1-\zeta_{r}) \nabla u \right\|_{L^2(U)} 
\leq 
C
\left| U \setminus U_{r} \right|^{\frac{\delta \wedge \gamma}{2+\delta \wedge \gamma}} \left( \mathsf{K} + \mathsf{M} \right)  
\leq  
C (3^{\ell - m})^{\frac{\delta \wedge \gamma}{2+\delta \wedge \gamma}} \left( \mathsf{K} + \mathsf{M} \right)  
\leq 
C\mathcal{E}(m).
\end{equation*}
For the second term we apply~\eqref{e.regutwoscale.NL} and obtain
\begin{equation*} 
\left\| \zeta_r \sum_{z} \psi_z^\ep \left( \nabla u(\cdot) - \nabla u( \ep z)\right) \right\|_{L^2(U)} \leq C 3^{-\beta(\ell-k)} \mathsf{M}_{r}  \leq \mathcal{E}(m).
\end{equation*} 
For the third term we first have
\begin{equation*} 
\left\|  \ep  \sum_{z} \nabla ( \zeta_r \psi_z^\ep)  \phi_z^\ep(\cdot,\nabla u( \ep z)) \right\|_{L^2(U)}^2 
\leq C 3^{-2k}   \sum_{z}  \left\|  \phi_z^\ep(\cdot,\nabla u( \ep z)) \right\|_{L^2(\ep(z + \cu_{k+2}))}^2 ,
\end{equation*}
and then, by Lemma~\ref{l.integrate.H-1} (using the zero mean of~$\phi_{z}$), and~\eqref{e.regutwoscale.NL}, we deduce that
\begin{align} 
\notag  
3^{-2k}\left\|  \phi_z^\ep(\cdot,\nabla u( \ep z)) \right\|_{L^2(\ep(z + \cu_{k+2}))}^2  &  
\leq
C 3^{(j-n)d}  3^{-2k} \left\|   \phi_z(\cdot,\nabla u( \ep z)) \right\|_{\underline{L}^2(z + \cu_{k+2})}^2
\\ 
\notag &
\leq 
C 3^{(j-n)d}  3^{-2k} \sup_{p \in B_{\mathsf{M}_{r} }} \left\|  \nabla  \phi_z^\ep(\cdot,p) \right\|_{\Hminusul(z + \cu_{k+2})}^2 .
\end{align}
Collecting last three displays yields~\eqref{e.E1estimate.NL}, after summation, recalling the definitions of $\mathcal{E}_{2}(m)$ and $\mathcal{E}(m)$. 

\smallskip

We next claim that we have
\begin{equation}  \label{e.sumE_2(x,y).NL}
\left\|  \sum_{y \in 3^k \Z^d \cap \ep^{-1} U_{2r}}  \indc_{ \ep (y + \cu_{k})} E_2(\cdot,y)  \right\|_{L^2 \left(\ep\left( y + \cu_k \right) \right) } 
\leq 
C\mathcal{E} (m).
\end{equation}
Together with~\eqref{e.E1estimate.NL} this yields~\eqref{e.Eestimate.NL}. 
To see~\eqref{e.sumE_2(x,y).NL}, take $y \in 3^k \Z^d \cap \ep^{-1} U_{2r}$ and estimate as
\begin{align} \notag 
\fint_{\ep (y + \cu_k)} \left| E_2(x,y)\right|^2 \, dx  
&  
\leq  
\sum_{z \sim y} \left\|\nabla v_z \left(\cdot,\nabla u( \ep z) \right) - \nabla v_z \left(\cdot,\nabla u( \ep y) \right) \right\|_{\underline{L}^2(y + \cu_{k})}^2
\\ \notag &  \quad 
+ \sum_{z \sim y} \left\|\nabla v_z \left(\cdot,\nabla u( \ep y) \right) - \nabla v_y \left(\cdot,\nabla u( \ep y) \right) \right\|_{\underline{L}^2(y + \cu_{k})}^2.
\end{align}
For the first term on the right we have, as in the proof of~\eqref{e.nuderp.NL}, by~\eqref{e.regutwoscale.NL} the deterministic bound
\begin{equation*} 
\left\|\nabla v_z \left(\cdot,\nabla u( \ep z) \right) - \nabla v_z \left(\cdot,\nabla u( \ep y) \right) \right\|_{\underline{L}^2(y + \cu_{k})}^2
 \leq 
 C \left| \nabla u( \ep z) - \nabla u( \ep y)\right|^2 \leq C 3^{2\beta(k-\ell)} \mathsf{M}_r^2 .
\end{equation*}
For the second term, on the other hand, we compare to patched $\nu$-minimizers in $(z' + \cu_k)$'s with the slope $\nabla u(\ep y)$ and use quadratic response as in~\eqref{e.nuquadresponse.NL} to obtain
\begin{align} 
\notag  
\lefteqn{\sum_{z \sim y}  \left\|\nabla v_z \left(\cdot,\nabla u( \ep y) \right) - \nabla v_y \left(\cdot,\nabla u( \ep y) \right) \right\|_{\underline{L}^2(y + \cu_{k})}^2 } \quad  &  
\\ 
\notag & \leq
C 
\sum_{z' \in 3^k \Z^d \cap (z + \cu_{k+3}) } \left( \nu(z' + \cu_{k},\nabla u( \ep y) )  - \nu(z + \cu_{k+2}, \nabla u( \ep y))     \right).
\\ 
\notag & \leq C
\sum_{j=k}^{k+2} \sum_{z' \in 3^k \Z^d \cap (z + \cu_{k+3}) } \sup_{p \in B_{  \mathsf{M}_r} } \left| \nu(z + \cu_k,p) - \overline{L}(p) \right|
\end{align}
Summation yields~\eqref{e.sumE_2(x,y).NL} by the definitions of $\mathcal{E}_{1}(m)$ and $\mathcal{E}(m)$.

\smallskip 

\emph{Step 6.} 
We next show that 
\begin{equation}  \label{e.twoscalegoal3.NL}
\left\| D_pL\left(\nabla w^\ep(\cdot), \tfrac \cdot\ep \right) - \zeta_{r} \sum_{z} \psi_z^\ep D_pL\left( \nabla v_z^\ep \left(\cdot, \nabla u( \ep z) \right) , \tfrac{\cdot}{\ep} \right) \right\|_{L^2(U)} 
\leq 
C\mathcal{E}(m). 
\end{equation}
We again split the analysis to an interior estimate and a boundary layer estimate. Our first claim is that 
\begin{equation}  \label{e.twoscalegoal3pre1.NL}
\left\| D_pL\left(\nabla w^\ep(\cdot), \tfrac \cdot\ep \right) - \sum_{z} \psi_z^\ep D_pL\left( \nabla v_z^\ep \left(\cdot, \nabla u( \ep z) \right) , \tfrac{\cdot}{\ep} \right) \right\|_{L^2(U_{4r})}
\leq 
C\mathcal{E}(m).
\end{equation}
For this, notice that, by~\eqref{e.gradw.nonlinear} and~\eqref{e.Lunifconv5}, we get, for $x \in U_{4r}$, that
\begin{equation*} 
\left| D_pL\left(\nabla w^\ep(x), \tfrac x\ep \right) - \sum_{y \in 3^k \Z^d \cap \ep^{-1} U_{2r}} \indc_{ \ep (y + \cu_{k})}(x)  D_pL\left(\nabla v_y^\ep \left(x,\nabla u( \ep y) \right) , \tfrac x\ep \right)  \right|
\leq 
C \left| E(x) \right|.
\end{equation*}
We then further decompose, for $y \in 3^k \Z^d \cap \ep^{-1} U_{2r}$ and $x\in\ep (y + \cu_{k})$,  as
\begin{align} \notag 
\lefteqn{D_pL\left( \nabla v_y^\ep \left(x,\nabla u( \ep y) \right),\tfrac{x}{\ep} \right) 
-
 \sum_{z  \sim y} \psi_z^\ep(x) D_pL\left( \nabla v_z^\ep \left(x,\nabla u( \ep z) \right),\tfrac{x}{\ep} \right) } \quad &
\\ \notag &  
= 
\sum_{z\sim y} \psi_z^\ep(x) \left( D_pL\left( \nabla v_y^\ep \left(x,\nabla u( \ep y) \right),\tfrac{x}{\ep} \right) 
- D_pL\left( \nabla v_z^\ep \left(x,\nabla u( \ep y) \right),\tfrac{x}{\ep} \right) 
\right) 
\\ \notag &  \quad
+ \sum_{z\sim y} \psi_z^\ep(x) \left( D_pL\left( \nabla v_z^\ep \left(x,\nabla u( \ep y) \right),\tfrac{x}{\ep} \right) 
- D_pL\left( \nabla v_z^\ep \left(x,\nabla u( \ep z) \right),\tfrac{x}{\ep} \right) 
\right) ,
\end{align}
and the second term on the left is the one appearing on the left in~\eqref{e.twoscalegoal3.NL}. Using~\eqref{e.Lunifconv5} we notice that both terms on the right can be estimated exactly as in the case of $E_{2}$.  This proves~\eqref{e.twoscalegoal3pre1.NL}.  

\smallskip 

We then claim that 
\begin{equation}  \label{e.twoscalegoal3pre2.NL}
\left\| D_pL\left(\nabla w^\ep(\cdot), \tfrac \cdot\ep \right) - \zeta_r \sum_{z} \psi_z^\ep D_pL\left( \nabla v_z^\ep \left(\cdot, \nabla u( \ep z) \right) , \tfrac{\cdot}{\ep} \right) \right\|_{L^2(U \setminus U_{4r})}
\leq 
C\mathcal{E}(m). 
\end{equation}
This is straightforward. Indeed, by~\eqref{e.gradw.nonlinear} and~\eqref{e.Lnormalize} we have
\begin{multline*} 
\left\| D_pL\left(\nabla w^\ep(\cdot), \tfrac \cdot\ep \right) \right\|_{L^2(U \setminus U_{4r})}  
 \leq 
C \left\| E \right\|_{L^2(U)} + C \left| U \setminus U_{4r} \right|^{\frac 12} \mathsf{K} 
\\ +C \left( \sum_{y \in 3^k \Z^d \cap \ep^{-1} (U_{2r} \setminus U_{5r})} \left\| \nabla v_y^\ep \left(x,\nabla u( \ep y) \right) \right\|_{L^2(\ep (y + \cu_{k}))}^2 \right)^{\frac12},
\end{multline*}
and since $v_y$ is $\nu$-minimizer with the slope $\nabla u( \ep y)$, the last term can be bounded using~\eqref{e.Lnormalize} and~\eqref{e.regutwoscale.NL} as
\begin{align} 
\notag  
\lefteqn{\sum_{y \in 3^k \Z^d \cap \ep^{-1} (U_{2r} \setminus U_{5r})} \left\| \nabla v_y^\ep \left(x,\nabla u( \ep y) \right) \right\|_{L^2(\ep (y + \cu_{k}))}^2 } \qquad  &  
\\ 
\notag &  \leq 
C \sum_{y \in 3^k \Z^d \cap \ep^{-1} (U_{2r} \setminus U_{5r})} \left| \ep\cu_{k}  \right|  \left(\mathsf{K} +  \left| \nabla u( \ep y) \right| \right)^2 .  
\\ 
\notag &  \leq 
C \sum_{y \in 3^k \Z^d \cap \ep^{-1} (U_{2r} \setminus U_{5r})} \left| \ep\cu_{k}  \right| \left(\mathsf{K}^2 +  \left\| \nabla u \right\|_{\underline{L}^2(y + \ep\cu_{k})}^2 
+ \left( 2^k \ep \right)^{2\beta} \left[ \nabla u \right]_{C^{\beta}(U_{2r})}^2 \right)
\\ 
\notag &  \leq 
C  \left| U \setminus U_{6r} \right|  \mathsf{K}^2  + C \left| U \setminus U_{r} \right|^{\frac{2 (\delta \wedge \gamma)}{2+\delta \wedge \gamma}}  \left(\mathsf{K} +  \mathsf{M} \right)^2  + C \left(\frac{2^k \ep}{r}\right)^{2\beta} \mathsf{M}_{r}^2
\\ 
\notag &  
\leq 
C\mathcal{E}^2(m) .
\end{align}
Similarly,
\begin{equation*} 
\left\| \zeta_r \sum_{z} \psi_z^\ep D_pL\left( \nabla v_z^\ep \left(\cdot, \nabla u( \ep z) \right) , \tfrac{\cdot}{\ep} \right) \right\|_{L^2(U \setminus U_{4r})} 
\leq 
C\mathcal{E}(m),
\end{equation*}
and~\eqref{e.twoscalegoal3pre2.NL} follows. Combining~\eqref{e.twoscalegoal3pre1.NL} and~\eqref{e.twoscalegoal3pre2.NL} proves~\eqref{e.twoscalegoal3.NL}.

\smallskip 

\emph{Step 7.} 
In this step we show that 
\begin{equation}  \label{e.twoscalegoal4.NL}
\left\| \nabla \cdot  D_pL\left(\nabla w^\ep(\cdot), \tfrac \cdot\ep \right) \right\|_{H^{-1}(U)} 
\leq 
C\mathcal{E}(m). 
\end{equation}
To see this, we first estimate as
\begin{align} \notag 
\lefteqn{\left\| \nabla \cdot  D_pL\left(\nabla w^\ep(\cdot), \tfrac \cdot\ep \right) \right\|_{H^{-1}(U)} } \quad &
\\ \notag & \leq 
\left\| D_pL\left(\nabla w^\ep(\cdot), \tfrac \cdot\ep \right) - \zeta_{r} \sum_{z} \psi_z^\ep D_pL\left( \nabla v_z^\ep \left(\cdot, \nabla u( \ep z) \right) , \tfrac{\cdot}{\ep} \right) \right\|_{L^2(U)} 
 \\ \notag & \quad + 
\left\| \nabla \cdot \left( \zeta_{r} \sum_{z} \psi_z^\ep \left( D_pL\left( \nabla v_z^\ep \left(\cdot, \nabla u( \ep z) \right) , \tfrac{\cdot}{\ep} \right) \right) 
\right) 
\right\|_{H^{-1}(U)} . 
\end{align}
Thus, recalling~\eqref{e.twoscalegoal3.NL}, we only need to estimate the last term on the right.  Using the fact that
\begin{equation*} 
 \nabla \cdot D_pL\left( \nabla v_z^\ep \left(\cdot, \nabla u( \ep z) \right) , \tfrac{\cdot}{\ep}\right) =  \nabla \cdot D \overline{L} \left(\nabla u\right) = \nabla \cdot D \overline{L} \left(\nabla u( \ep z)   \right) =  0,
\end{equation*}
giving also
\begin{equation*} 
 \nabla \cdot \left( \sum_{z}  \zeta_r   \psi_z^\ep  D \overline{L} \left(\nabla u \right)  \right) = \nabla \zeta_r \cdot D \overline{L} \left(\nabla u \right)  , 
\end{equation*}
we expand the divergence as follows:
\begin{align} \notag  \label{e.twoscalegoal4pre.NL}
\lefteqn{\nabla \cdot \left( \zeta_{r} \sum_{z} \psi_z^\ep \left( D_pL\left( \nabla v_z^\ep \left(\cdot, \nabla u( \ep z) \right) , \tfrac{\cdot}{\ep} \right) \right) 
 \right) } \quad &
\\ \notag & = \sum_{z} \nabla  \left(\zeta_r   \psi_z^\ep \right) \cdot 
\left(  D_pL\left( \nabla v_z^\ep \left(\cdot, \nabla u( \ep z) \right) , \tfrac{\cdot}{\ep} \right) 
- D \overline{L} \left(\nabla u( \ep z)   \right) \right)
\\ & \quad + \nabla \cdot \left( \sum_{z}  \zeta_r   \psi_z^\ep  
\left(  D \overline{L} \left(\nabla u( \ep z)   \right) - D \overline{L} \left(\nabla u \right) \right) - (1-\zeta_r)  D \overline{L} \left(\nabla u \right)  \right) .
\end{align}
Therefore, by the triangle inequality, 
\begin{align} \notag 
\lefteqn{\left\| \nabla \cdot \left( \zeta_{r} \sum_{z} \psi_z^\ep \left( D_pL\left( \nabla v_z^\ep \left(\cdot, \nabla u( \ep z) \right) , \tfrac{\cdot}{\ep} \right) \right) 
 \right) 
 \right\|_{H^{-1}(U)}
 } \quad &
\\ \notag & \leq  \left\| \sum_{z}  \nabla  \left(\zeta_r   \psi_z^\ep \right) \cdot 
\left(  D_pL\left( \nabla v_z^\ep \left(\cdot, \nabla u( \ep z) \right) , \tfrac{\cdot}{\ep} \right) 
- D \overline{L} \left(\nabla u( \ep z)   \right) \right)  \right\|_{H^{-1}(U)}
\\ \notag & \quad + \left\| \sum_{z}  \zeta_r   \psi_z^\ep  
\left(  D \overline{L} \left(\nabla u( \ep z)   \right) - D \overline{L} \left(\nabla u \right) \right) - (1-\zeta_r)  D \overline{L} \left(\nabla u \right)  \right\|_{L^{2}(U)}. 
\end{align}
On the one hand, we have, for any $\theta \in H_0^1(U)$ with $\|\nabla \theta\|_{L^2(U)} \leq 1$, that
\begin{align} \notag 
\lefteqn{ \sum_z \int_{U}\nabla  \left(\zeta_r   \psi_z^\ep \right) \cdot 
\left(  D_pL\left( \nabla v_z^\ep \left(\cdot, \nabla u( \ep z) \right) , \tfrac{\cdot}{\ep} \right) 
- D \overline{L} \left(\nabla u( \ep z)   \right) \right)  \theta } \quad &
\\ \notag & \leq 
C \left(  \sum_z 
\frac{\left|\ep \cu_k \right| }{(2^k\ep)^{2}} \left\| 
D_pL\left( \nabla v_z^\ep \left(\cdot, \nabla u( \ep z) \right) , \tfrac{\cdot}{\ep} \right) - D \overline{L} \left(\nabla u( \ep z)   \right) \right\|_{\underline{H}^{-1}(\ep(z + \cu_{k+2}))}^2 
\right)^{\frac12} 
\\ \notag & \qquad  \times \left(  \sum_z (2^k\ep)^{2} \left\| \nabla \left( \nabla  \left(\zeta_r   \psi_z^\ep \right) \theta \right) \right\|_{L^2(\ep(z + \cu_{k+2}))}^2 \right)^{\frac 12} . 
\end{align}
By~\eqref{e.notation.Sobolev.scaling-},~\eqref{e.regutwoscale.NL} and the definition of $\mathcal{E}_{3}(m)$, we have
\begin{equation*} 
\sum_z 
\frac{\left|\ep \cu_k \right| }{(2^k\ep)^{2}} \left\| 
D_pL\left( \nabla v_z^\ep \left(\cdot, \nabla u( \ep z) \right) , \tfrac{\cdot}{\ep} \right) - D \overline{L} \left(\nabla u( \ep z)   \right) \right\|_{\underline{H}^{-1}(\ep(z + \cu_{k+2}))}^2 
\\ \leq  
C\mathcal{E}^2(m),
\end{equation*}
and, by~\eqref{e.gradzetabounded.NL} and~\eqref{e.gradpsibounded.NL},
\begin{equation*} 
\sum_z (2^k\ep)^{2} \left\| \nabla \left( \nabla  \left(\zeta_r   \psi_z^\ep \right) \theta \right) \right\|_{L^2(\ep(z + \cu_{k+2}))}^2 \leq C \left(1 + \frac{2^k \ep}{r} \right)^2 \leq C  .
\end{equation*}
Thus we get
\begin{equation*} 
\left\| \sum_{z}  \nabla  \left(\zeta_r   \psi_z^\ep \right) \cdot 
\left(  D_pL\left( \nabla v_z^\ep \left(\cdot, \nabla u( \ep z) \right) , \tfrac{\cdot}{\ep} \right) 
- D \overline{L} \left(\nabla u( \ep z)   \right) \right)  \right\|_{H^{-1}(U)} 
\leq 
C\mathcal{E}(m). 
\end{equation*}
On the other hand, for the other term, we use a similar reasoning as in the case of $E_1$ in Step 5, using now properties of $\overline{L}$, and also obtain 
\begin{equation*} 
\left\| \sum_{z}  \zeta_r   \psi_z^\ep  
\left(  D \overline{L} \left(\nabla u( \ep z)   \right) - D \overline{L} \left(\nabla u \right) \right) - (1-\zeta_r)  D \overline{L} \left(\nabla u \right)  \right\|_{L^{2}(U)} \leq 
C\mathcal{E}(m).
\end{equation*}
Now~\eqref{e.twoscalegoal4.NL} follows by~\eqref{e.twoscalegoal3.NL}. 

\smallskip 

\emph{Step 8.} 
We prove that
\begin{equation} \label{e.twoscalegoal5.NL}
\left\| \nabla u^\ep - \nabla u \right\|_{H^{-1}(U)} + \left\| D_{p} L\left( \nabla u^\ep , \tfrac \cdot\ep\right) - D\overline{L}\left( \nabla u \right) \right\|_{H^{-1}(U)} \leq 
C\mathcal{E} (m). 
\end{equation}
Since $w^\ep - u^\ep \in H_{0}^1(U)$, we have by~\eqref{e.twoscalegoal4.NL} and~\eqref{e.Lunifconv3} that 
\begin{equation}  \label{e.twoscalegoal6.NL}
\left\| \nabla u^\ep - \nabla w^\ep \right\|_{L^2(U)} \leq 
C\mathcal{E}(m).  
\end{equation}
In particular, by~\eqref{e.Lunifconv5} it follows that 
\begin{equation*} 
\left\| \nabla u - \nabla w^\ep \right\|_{H^{-1}(U)} + \left\| D_{p} L\left( \nabla u^\ep , \tfrac \cdot\ep\right) - D_{p} L\left( \nabla w^\ep , \tfrac \cdot\ep\right) \right\|_{H^{-1}(U)} 
\leq 
C\mathcal{E}(m).
\end{equation*}
Therefore it is enough to estimate the average and the flux of $\nabla w^\ep$ against $\nabla u$. The average is easy to compute directly from the definition of $w^\ep$, as in Step 5, namely
\begin{align} 
\notag  
\left\| \nabla u - \nabla w^\ep \right\|_{H^{-1}(U)}  & =  \ep  \left\|   \nabla \left( \zeta_r \sum_{z} \psi_z^\ep \phi_z^\ep \left(\cdot,\nabla u( \ep z)\right) \right) \right\|_{H^{-1}(U)} 
\\ 
\notag & 
\leq  C \ep \left\|   \zeta_r \sum_{z} \psi_z^\ep \phi_z^\ep \left(\cdot,\nabla u( \ep z)\right)  \right\|_{L^2(U)}
\leq 
C\mathcal{E}(m).
\end{align}
To estimate the fluxes, we have by~\eqref{e.twoscalegoal3.NL} that 
\begin{equation*} 
\left\| D_pL\left(\nabla w^\ep(\cdot), \tfrac \cdot\ep \right) - \zeta_{r} \sum_{z} \psi_z^\ep D_pL\left( \nabla v_z^\ep \left(\cdot, \nabla u( \ep z) \right) , \tfrac{\cdot}{\ep} \right) \right\|_{H^{-1}(U)} 
\leq 
C\mathcal{E}(m). 
\end{equation*}
Similarly to~\eqref{e.twoscalegoal4pre.NL}, we write
\begin{align} \notag  
\lefteqn{
\zeta_{r} \sum_{z} \psi_z^\ep D_pL\left( \nabla v_z^\ep \left(\cdot, \nabla u( \ep z) \right) , \tfrac{\cdot}{\ep} \right)    - D \overline{L} \left(\nabla u \right)  } \quad &
\\ \notag & = \sum_{z} \zeta_r   \psi_z^\ep \ \cdot 
\left(  D_pL\left( \nabla v_z^\ep \left(\cdot, \nabla u( \ep z) \right) , \tfrac{\cdot}{\ep} \right) 
- D \overline{L} \left(\nabla u( \ep z)   \right) \right)
\\ & \notag \quad +  \sum_{z}  \zeta_r   \psi_z^\ep  \left(  D \overline{L} \left(\nabla u( \ep z)   \right) - D \overline{L} \left(\nabla u \right) \right) - (1-\zeta_r)  D \overline{L} \left(\nabla u \right) .
\end{align}
The estimates for the two terms on the right are analogous to Step 7, and we omit the details. Now~\eqref{e.twoscalegoal4.NL} follows. 

\smallskip 

\emph{Step 9.} We estimate the error term~$\mathcal{E}(m)$. Similar to~\eqref{e.sumthisupO} in the proof of Theorem~\ref{t.DP}, what we need to show is that, with appropriate choices of the parameters, we have, for some $\beta(\delta,d,\Lambda)>0$,
\begin{equation} 
\label{e.sumthisupO.NLcase}
\sup_{m\in\N} 3^{ms}\left(  \mathcal{E}(m)   - C(\mathsf{K} + \mathsf{M})3^{-m\beta(d-s)} \right)_+ 
\leq 
\O_1\left(C(\mathsf{K} + \mathsf{M})\right). 
\end{equation}
Applying Theorem~\ref{t.subadd.NL}, together with~\eqref{e.sqrtelem}, and Lemma~\ref{l.approx.correctors.NL}, we obtain that, for $s' := \frac{d+s}{2}$, there exists constants $\tilde \alpha(d,\Lambda)>0$ and $C(s,d,\Lambda)<\infty$ such that, for every~$m\in\N$, 
\begin{equation*} 
\mathcal{E}_{1}(m) + \mathcal{E}_{2}(m) + \mathcal{E}_{3}(m) 
\leq 
C (\mathsf{K} + \mathsf{M}) 
r^{-\frac d2}3^{-\tilde \alpha(d-s')k}  
+ \O_{1}\left(C  (\mathsf{K} + \mathsf{M}) r^{-\frac d2} 3^{-s' k}  \right).
\end{equation*}
Now the desired estimate~\eqref{e.sumthisupO.NLcase} is obtained by choosing $r$ appropriately and arguing as in the proof of~\eqref{e.sumthisupO}. The proof is complete.
\end{proof}

\section{\texorpdfstring{$C^{0,1}$}{C0,1}-type estimate}
\label{s.LLip.NL}

The goal of this section is to provide the counterpart for Theorem~\ref{t.Lipschitz} in Chapter~\ref{c.regularity}.
\begin{theorem}
[{Quenched $C^{0,1}$-type estimate}]
\label{t.Lipschitz.NL}
Fix $s\in (0,d)$ and $\mathsf{M} \in [1,\infty)$. There exist a constant  $C(\mathsf{M},s,\mathsf{K},d,\Lambda)<\infty$ and a random variable $\X_s :\Omega \to [1,\infty]$ satisfying
\begin{equation}
\label{e.XLipschitz.NL}
\X_s = \O_s\left( C \right),
\end{equation}
such that the following holds: for every $R \geq \X_s$ and a weak solution $u\in H^1(B_{R})$ of
\begin{equation} 
\label{e.wksrt2.NL}
-\nabla \cdot \left( D_pL\left(\nabla u,x\right) \right) = 0 \quad \mbox{in} \ B_{R},
\end{equation}
satisfying 
\begin{equation*} 
\frac1R \left\| u - (u)_{B_R}\right\|_{\underline{L}^2(B_R)} \leq \mathsf{M} 
\end{equation*}
we have, for every  $r\in \left[ \X_s , R\right]$, the estimate
\begin{equation}
\label{e.Lipschitz.NL}
 \left\| \nabla u \right\|_{\underline{L}^2(B_{r})}
\leq
C \left( \frac{1}{R} \left\| u - \left( u \right)_{B_{R}} \right\|_{\underline{L}^2(B_{R})}  + \mathsf{K} \right).
\end{equation}
\end{theorem}

We begin the proof of Theorem~\ref{t.Lipschitz.NL} with the observation that Theorem~\ref{t.DP.nonlinear} yields a corresponding version of  Proposition~\ref{p.harmonicapproximation} for nonlinear equations.

\begin{proposition}[{Harmonic approximation}]
\label{p.harmonicapproximation.NL}
Fix $s\in (0,d)$ and $\mathsf{M} \in [1,\infty)$. There exist constants $\alpha(d,\Lambda)>0$, $C(s,d,\Lambda)<\infty$, and a random variable $\X_s:\Omega \to [1,\infty]$ satisfying 
\begin{equation}
\label{e.Xharmonicapproximation.NL}
\X_s = \O_s\left( C\left(\mathsf{K}+ \mathsf{M} \right) \right),
\end{equation}
such that the following holds: for every $R\geq \X_s$ and weak solution $u\in H^1(B_{R})$ of
\begin{equation} 
\label{e.wksrt.NL}
-\nabla \cdot \left( D_pL\left(\nabla u,x\right) \right) = 0 \quad \mbox{in} \ B_{R},
\end{equation}
satisfying 
\begin{equation}  \label{e.unormalization.NL}
\frac1R \left\| u - (u)_{B_R}\right\|_{\underline{L}^2(B_R)} \leq \mathsf{M} ,
\end{equation}
there exists a solution $\overline{u} \in H^1(B_{R/2})$ of
\begin{equation*} \label{}
-\nabla \cdot \left( D\overline{L}\left(\nabla \overline{u}\right) \right) = 0 \quad \mbox{in} \ B_{R/2}
\end{equation*}
satisfying 
\begin{equation}
\label{e.harmonicapproximation.NL}
\left\| u - \overline{u}  \right\|_{\underline{L}^2(B_{R/2})} \leq CR^{-\alpha(d-s)}\left(\mathsf{K}+ \mathsf{M} \right).
\end{equation}
\end{proposition}

\begin{proof}
By the Meyers estimate~\eqref{e.regumeyers.NL} and the Caccioppoli inequality,
\begin{equation*}  
  \left\| \nabla u \right\|_{\underline{L}^{2 + \gamma}(B_{R/2})} 
 \leq 
 C\left( \frac1R \left\| u - (u)_{B_R}\right\|_{\underline{L}^{2}(B_R)} + \mathsf{K} \right) \leq C \left(\mathsf{K}+ \mathsf{M} \right).
\end{equation*}
After this the proof is analogous to the one of  Proposition~\ref{p.harmonicapproximation} using Theorem~\ref{t.DP.nonlinear}. 
Indeed, for each $s\in (0,d)$, we let $\tilde{\X}_s$ denote the random variable in the statement of Theorem~\ref{t.DP.nonlinear} with $\delta$ the exponent in the Meyers estimate, as above, such that $\tilde{\X}_s \leq \O_1 \left(C\left(\mathsf{K}+ \mathsf{M} \right)\right)$. 
Setting 
\begin{equation*} \label{}
\X_s := \tilde{\X}_{(s+d)/2}^{1/s} \vee 1,
\end{equation*}
we have that~$\X_s = \O_s\left(C\left(\mathsf{K}+ \mathsf{M} \right)\right)$. Applying Theorem~\ref{t.DP.nonlinear} with~$\ep = R^{-1}$ and~$\tilde s = \frac{s+d}{2}$ gives us the desired conclusion for every~$R\geq \X_s$ if we take~$\overline{u}$ to be the solution of the Dirichlet problem for the homogenized equation with Dirichlet boundary condition~$u$ on~$\partial B_{R/2}$ as in the statement. 
\end{proof}

\begin{proof}[Proof of Theorem~\ref{t.Lipschitz.NL}]
We may replace the application of Lemma~\ref{l.u:decayestimate} by using instead~\eqref{e.refreg.NL} and Proposition~\ref{p.harmonicapproximation.NL}.  Let $H,\tilde H$ be large constants to be fixed in the course of the proof. Let $\X_s = \O_s\left( C\left( \mathsf{K} + \tilde H \mathsf{M} \right) \right)$ be as in Proposition~\ref{p.harmonicapproximation.NL}. Without loss of generality we may assume that $R \geq H$, since otherwise the statement follows easily by giving up a volume factor.  Also assume that $\X_s \leq R$. 

\smallskip

Arguing by induction, we assume next that, for some~$r \in [\X_s \vee H,R]$, we have
\begin{equation}  \label{e.indassMr.NL}
\mathsf{M}_r := \sup_{t \in[r,R]} \frac1t \left\| u - (u)_{B_t}\right\|_{\underline{L}^{2}(B_t)} \leq \frac KR \left\| u - (u)_{B_R}\right\|_{\underline{L}^2(B_R)} \leq \tilde H \left( \mathsf{M} + \mathsf{K} \right). 
\end{equation}
Notice that this assumption is valid for ~$r = R$ by~\eqref{e.unormalization.NL}. The goal is to show that the assumption is also valid for~$r/2$ in place of~$r$. 

\smallskip

Using Proposition~\ref{p.harmonicapproximation.NL} gives that if $t \geq r \geq \X_s$, then we find $\overline{u}_t \in H^1(B_t)$ such that  it solves 
\begin{equation*} \label{}
-\nabla \cdot \left( D\overline{L}\left(\nabla \overline{u}_t\right) \right) = 0 \quad \mbox{in} \ B_{t}
\end{equation*}
and satisfies, by the induction assumption~\eqref{e.indassMr.NL},  
\begin{equation*}
\left\| u - \overline{u}_t  \right\|_{\underline{L}^2(B_{t})} \leq C\tilde H t^{1-\alpha(d-s)} \left(\mathsf{K}+ \mathsf{M} \right) .
\end{equation*}
Moreover, by~\eqref{e.refreg.NL} and the Caccioppoli inequality, $\overline{u}_t$ satisfies the decay estimate 
 \begin{equation} \label{e.refregpre.NL}
  \inf_{p \in \mathcal{P}_1} \left\| \overline{u}_t - p \right\|_{\underline{L}^{2}(B_{\theta t})}
 \leq C \theta^{1+\beta} \left\| \overline{u}_t - p \right\|_{\underline{L}^{2}(B_{t})} . 
\end{equation}
By combining the previous two displays with the aid of the triangle inequality and choosing $\theta$ so that $C \theta^\beta = \frac12$, we have that, for $t \in \left[r,R\right]$, 
\begin{equation*} 
 \frac{1}{\theta t} \inf_{p \in \mathcal{P}_1} \left\| u - p \right\|_{\underline{L}^{2}(B_{\theta t})} \leq \frac12 
\frac{1}{t} \inf_{p \in \mathcal{P}_1}  \left\| u - p \right\|_{\underline{L}^{2}(B_{t})} +  Ct^{-\alpha(d-s)}  \tilde H  \left(\mathsf{K}+ \mathsf{M}\right) .
\end{equation*}
Varying $t$, summing and reabsorbing yields that, for $k \in \N$ such that $2^{-k} R \leq r  \leq 2^{-k+1}R$, we have that 
\begin{equation*} 
\sum_{j=0}^k  \frac{1}{2^{-j} R} \inf_{p \in \mathcal{P}_1} \left\| u - p \right\|_{\underline{L}^{2}(B_{2^{-j} R})} \leq 
\frac{C}{R} \left\| u - (u)_{B_R} \right\|_{\underline{L}^{2}(B_{R})} +  CH^{-\alpha(d-s)}  \tilde H  \left(\mathsf{K}+ \mathsf{M} \right) .
\end{equation*}
Denoting $p_j \in \mathcal{P}_1$ to be the affine function realizing the infimum above for $j \in \N$, and setting $p_0 = (u)_{B_R}$, we see by the triangle inequality that 
\begin{equation*} 
\max_{j \in \{0,\ldots,k\}} \left| \nabla p_{j} \right| \leq  \sum_{j=0}^{k-1}  \left| \nabla p_{j+1} - \nabla p_{j} \right| \leq \frac{C}{R} \left\| u - (u)_{B_R} \right\|_{\underline{L}^{2}(B_{R})} +  CH^{-\alpha(d-s)}  \tilde H \left(\mathsf{K}+ \mathsf{M}\right)  ,
\end{equation*}
which also yields, for $j \in \{0,\ldots,k\}$, 
\begin{multline*} 
\frac{1}{2^{-j}R} \left\| u - (u)_{B_{2^{-j} R}} \right\|_{\underline{L}^{2}(B_{2^{-j} R})}  \leq  
\frac{1}{2^{-j}R}\inf_{p \in \mathcal{P}_1} \left\| u - p \right\|_{\underline{L}^{2}(B_{2^{-j} R})}  + \left| \nabla p_{j} \right| 
\\ \leq \frac{C}{R} \left\| u - (u)_{B_R} \right\|_{\underline{L}^{2}(B_{R})} +  CH^{-\alpha(d-s)}  \tilde H \left(\mathsf{K}+ \mathsf{M} \right) . 
\end{multline*}
Therefore, since
\begin{equation*} 
\mathsf{M}_{r/2} \leq C \max_{j \in \{0,\ldots,k\}} \frac{1}{2^{-j}R} \left\| u - (u)_{B_{2^{-j} R}} \right\|_{\underline{L}^{2}(B_{2^{-j} R})} ,
\end{equation*}
we deduce that 
\begin{equation*} 
\mathsf{M}_{r/2} \leq \frac{C}{R} \left\| u - (u)_{B_R} \right\|_{\underline{L}^{2}(B_{R})} +  CH^{-\alpha(d-s)} \tilde H \left(\mathsf{K}+ \mathsf{M} \right) . 
\end{equation*}
Choosing thus $H$ large enough so that $CH^{-\alpha(d-s)} = \frac12$, we obtain 
\begin{equation*} 
\mathsf{M}_{r/2} \leq \frac{C}{R} \left\| u - (u)_{B_R} \right\|_{\underline{L}^{2}(B_{R})} + \frac12 \tilde H \left(\mathsf{K}+ \mathsf{M} \right) , 
\end{equation*}
which proves the induction step with $\tilde H=2C$, and completes the proof. 
\end{proof}

\section*{Notes and references}
The results in this chapter first appeared in~\cite{AS} and were extended to the case of general monotone maps (that is, for~$\a$ not necessarily of the form~\eqref{e.classicalform}) in~\cite{AM}. The higher regularity theory in the nonlinear case is still incomplete (although see~\cite{AFK} for the first steps in this direction). The problem of adapting the arguments here to nonlinear functionals with more general, possibly non-quadratic growth conditions is also open.

\appendix



\chapter{The~\texorpdfstring{$\O_s$}{Os} notation}
\label{a.bigO}

Throughout the book, we use the symbol  $\O_s(\cdot)$ as a way to express inequalities for random variables with respect to what is essentially the Orlicz norm~$L^\varphi(\Omega,\P)$ with $\varphi(t) = \exp \left( t^s \right)$. In this appendix, we collect some basic properties implied by this notation for the convenience of the reader. 

\smallskip

Let us recall the definition given in~\eqref{e.Os.not}. For a random variable $X$ and an exponent~$s, \theta \in (0,\infty)$, we write
\begin{equation}
\label{e.def.Os.notation}
X \leq \O_s(\theta)
\end{equation}
to mean that 
\begin{equation}
\label{e.def.Os}
\E \Ll[ \exp\left(\left(\theta^{-1} X_+\right)^s\right) \Rr] \le 2,
\end{equation}
where $X_+ := X \vee 0 = \max\{ X,0\}$. The notation is obviously homogeneous: 
\begin{equation*}  
X \le \O_s(\theta) \quad \iff \quad \theta^{-1} X \le \O_s(1).
\end{equation*}
We also write 
\begin{equation*}  
X = \O_s(\theta) \quad \iff \quad X \le \O_s(\theta) \ \ \text{and} \ \ -X \le \O_s(\theta), 
\end{equation*}
and by extension,
\begin{equation*}
X\leq Y +  \O_s(\theta)  \quad \iff \quad X-Y \leq \O_s(\theta),
\end{equation*}
\begin{equation*}
X = Y + \O_s(\theta) \quad\iff \quad X-Y \leq \O_s(\theta) \ \ \mbox{and} \ \ Y-X \leq \O_s(\theta). 
\end{equation*}
We start by relating the property \eqref{e.def.Os.notation} with the behavior of the tail probability of the random variable.
\begin{lemma}[Tail probability]
\label{l.bigO.vs.tail}
For every random variable $X$ and $s,\theta \in (0,\infty)$, we have
\begin{equation*}  
X \le \O_s(\theta) \quad \implies \quad \forall x \ge 0, \quad \P\Ll[ X \ge \theta x \Rr] \le  2 \exp \Ll( - x^s \Rr) ,
\end{equation*}
and 
\begin{equation*}  
\forall x \ge 0, \quad \P\Ll[ X \ge \theta x \Rr] \le \exp \Ll( - x^s \Rr) \quad \implies \quad X \le \O_s \Ll(2^\frac 1 s \, \theta\Rr).
\end{equation*}
\end{lemma}
\begin{proof}
The first implication is immediate from Chebyshev's inequality. For the second one, we apply Fubini's theorem to get
\begin{align*}  
\E \Ll[ \exp\left(2^{-1}\left( \theta^{-1} X_+\right)^s\right) \Rr] &
= 1 + \E \Ll[ \int_0^{\theta^{-1} X_+} 2^{-1} s x^{s-1} \exp \Ll(2^{-1} x^s \Rr) \, dx  \Rr] \\
& = 1 + \int_0^{\infty} 2^{-1} s x^{s-1} \exp \Ll( 2^{-1} x^s \Rr) \P \Ll[ X \ge \theta x \Rr]  \, dx \\
& \le 1 + \int_0^{\infty} 2^{-1} s x^{s-1} \exp \Ll( - 2^{-1} x^s \Rr) \,  dx  \\
& = 2.
\qedhere
\end{align*}
\end{proof}
\begin{remark}  
\label{r.shift.x0}
If we only assume that, for some $x_0 \in [0,\infty)$, we have
\begin{equation*}  
\forall x \ge x_0, \quad \P\Ll[ X \ge \theta x \Rr] \le \exp \Ll( - x^s \Rr)
\end{equation*}
then we can apply Lemma~\ref{l.bigO.vs.tail} to the random variable $\theta^{-1} X \1_{\theta^{-1} X \ge x_0}$ and get
\begin{equation}  
\label{e.first.ambiguity}
X \le (\theta x_0) \vee \O_s \Ll( 2^{\frac 1 s} \theta \Rr) ,
\end{equation}
where we understand this inequality to mean
\begin{equation*}  
\text{there exists a r.v. } Y \text{ s.t. } Y \le \O_s(2^{\frac 1 s} \, \theta) \text{ and } X \le (\theta x_0) \vee Y.
\end{equation*}
\end{remark}

In the next lemma, we record simple interpolation and multiplication results for $\O_s$-bounded random variables.
\begin{lemma}
\label{l.change-s}
\emph{(i)} 
For every random variable $X$ and $s < s' \in (0,\infty)$,
$$
\Ll\{
\begin{array}{l}
X \mbox{ takes values in }  [0,1] \\
X \le \O_s(\theta)
\end{array}
\Rr.
\quad  \implies \quad X \le \O_{s'}(\theta^{\frac s {s'}}).
$$

\noindent \emph{(ii)} 
If $0 \le X_i \le \O_{s_i}(\theta_i)$ for $i \in\{1,2\}$, then
\begin{equation*} 
X_1 X_2 \le \O_{\frac{s_1 s_2}{s_1 + s_2}} \left( \theta_1\theta_2 \right)\,.
\end{equation*}

\end{lemma}
\begin{proof}
The first statement follows from
$$
\E\Ll[\exp( (\theta^{-\frac s {s'}} X)^{s'}) \Rr] \le \E\Ll[\exp( (\theta^{-1} X)^{s})\Rr] \le 2.
$$
For the second statement, we apply Young's and H\"older's inequalities to get
\begin{align*}
\lefteqn{
\E\Ll[\exp\Ll(\Ll[(\theta_1 \theta_2)^{-1} X_1 X_2\Rr]^\frac{s_1 s_2}{s_1 + s_2}\Rr) \Rr] 
} \qquad & \\
& \le \E\Ll[\exp\Ll(\frac{s_2}{s_1 + s_2}(\theta_1^{-1} X_1)^{s_1} +  \frac{s_1}{s_1 + s_2} (\theta_2^{-1} X_2)^{s_2} \Rr) \Rr] \\
& \le \E\Ll[\exp\Ll(\theta_1^{-1} X_1)^{s_1}\Rr)\Rr]^{\frac{s_2}{s_1 + s_2}}\, \E\Ll[\exp\Ll(\theta_2^{-1} X_2)^{s_2}\Rr)\Rr]^{\frac{s_1}{s_1 + s_2}} \\
& \le 2. \qedhere
\end{align*}
\end{proof}

We next show that an average of random variables bounded by $\O_s(\theta)$ is bounded by $\O_s(C \theta)$, with $C = 1$ when $s \ge 1$.
\begin{lemma}
\label{l.sum-O}
For each $s \in (0,\infty)$, there exists a constant $C_s < \infty$ such that the following holds. Let $\mu$ be a measure over an arbitrary measurable space $E$, let $\theta : E \to (0,\infty)$ be a measurable function and $(X(x))_{x \in E}$ be a jointly measurable family of nonnegative random variables such that, for every $x \in E$, $X(x) \le \O_s(\theta(x))$. We have
\begin{equation*}  
\int X \, d\mu \le  \O_s \Ll(C_s \int \theta \, d \mu \Rr).
\end{equation*}
Moreover, for each $s \ge 1$, the statement holds for $C_s = 1$.
\end{lemma}
\begin{proof}
Without loss of generality, we can assume that $\int \theta \, d \mu < \infty$, and by homogeneity, we can further assume that $\int \theta \, d \mu = 1$. We first consider the case $s \ge 1$. 
Since the function $x \mapsto \exp(x^s)$ is convex on $\R_+$, Jensen's inequality gives
\begin{equation*} 
\E\left[\exp\left(  \left(\int X \, d\mu \right)^{s} \right) \right]  \leq \E\left[\int \exp \left( (\theta^{-1} X)^s  \right) \, \theta \, d \mu \right] 
 \leq 2,
\end{equation*}
as announced. For $s \in (0,1)$, we define $t_s := \Ll( \frac{1-s} s \Rr)^{\frac 1 s}$ so that the function $x \mapsto \exp((x + t_s)^s)$ is convex on $\R_+$. By Jensen's inequality,
\begin{align*} 
\E\left[\exp\left(  \left(\int X \, d\mu \right)^{s} \right) \right]  & \leq \E\left[\exp \left( \left( \int X  \, d\mu  + t_s \right)^{s}  \right) \right]  
\\ & \leq \E\left[\int \exp \left( (\theta^{-1} X + t_s)^s \right) \, \theta d \mu \right] 
\\ & \leq \E\left[\int \exp \left( (\theta^{-1} X)^s + t_s^s \right) \, \theta d \mu \right] 
\\ & \leq 2 \exp\left( t_s^s\right)\,.
\end{align*}
For $\sigma \in (0,1)$ sufficiently small in terms of $s$, we thus have
\begin{equation*}
\E\left[\exp\left( \sigma \left(\int X \, d\mu \right)^{s}  \right) \right] \leq \E\left[\exp\left(\left(\int X \, d\mu \right)^{s} \right) \right]^\sigma
\leq  \Ll[2 \exp \Ll( t_s^s\Rr)\Rr]^\sigma \le 2\,. \qedhere
\end{equation*}
\end{proof}

We next explore a correspondence between the assumption of $X \le \O_s(1)$, for $s \in (1,\infty)$, and the behavior of the Laplace transform of $X$. 

\begin{lemma}[Laplace transform]
\label{l.lapl.bigO}
For every $s \in (1,\infty)$, there exists a constant $C(s) < \infty$ such that
\begin{equation}
\label{e.stoch-equiv1}
X \le \O_s(1)  \quad \implies \quad \forall \lambda \ge 1, \ \log \E[\exp(\lambda X)] \le C \lambda^{\frac s {s-1}},
\end{equation}
and
\begin{equation}
\label{e.stoch-equiv2}
\forall \lambda \ge 1, \ \log \E[\exp(\lambda X)] \le  \lambda^{\frac s {s-1}} \quad \implies \quad X \le \O_s(C).
\end{equation}
\end{lemma}
\begin{proof}
For each $\lambda \ge 1$, we apply Fubini's theorem and then Lemma~\ref{l.bigO.vs.tail} to get
\begin{align*}  
\E[\exp(\lambda X)] & = 1 + \lambda \int_0^\infty \exp(\lambda x) \, \P[X \ge x] \, dx \\
& \le 1 + 2\lambda \int_0^\infty \exp(\lambda x - x^s) \, dx.
\end{align*}
Decomposing the integral into two parts:
\begin{equation*}  
\lambda \int_0^{ (2\lambda)^{\frac 1 {s-1}} } \exp \Ll( \lambda x - x^s \Rr) \, dx \le \lambda \int_0^{ (2\lambda)^{\frac 1 {s-1}} } \exp \Ll( \lambda x  \Rr) \, dx \le  \exp \Ll( C \lambda^{\frac {s}{s-1}}\Rr) ,
\end{equation*}
and 
\begin{equation*}  
\int_{ (2\lambda)^{\frac 1 {s-1}} }^\infty \exp \Ll( \lambda x - x^s \Rr) \, dx \le  \int_{ (2\lambda)^{\frac 1 {s-1}} }^\infty \exp \Ll( - 2^{-1} x^s \Rr) \, dx \le C ,
\end{equation*}
we obtain \eqref{e.stoch-equiv1}. As for \eqref{e.stoch-equiv2}, by Chebyshev's inequality, for every $x \ge 2$ and $\lambda \ge 1$, we have
\begin{equation*}  
\P \Ll[ X \ge x \Rr]  \le \exp \Ll(\lambda^{\frac s {s-1}} -\lambda x \Rr) .
\end{equation*}
Selecting $\lambda = (x/2)^{s-1}$, we get
\begin{equation*}  
\P \Ll[ X \ge x \Rr]  \le \exp \Ll( - 2^{-s} x^s \Rr) .
\end{equation*}
The conclusion thus follows from Lemma~\ref{l.bigO.vs.tail} and Remark~\ref{r.shift.x0}.
\end{proof}
\begin{exercise}  
Show that there exists a constant $C \in (0, \infty)$ such that
\begin{equation*}  
X \le \O_1(1)  \quad \implies \quad  \log \E[\exp(C^{-1} X)] \le 1,
\end{equation*}
and
\begin{equation*}  
\log \E[\exp(X)] \le 1 \quad \implies \quad  X \le \O_1(C).
\end{equation*}
\end{exercise}
One classical way to prove the central limit theorem for sums of i.i.d.\ and sufficiently integrable random variables is based on the computation of Laplace transforms. The scaling of the central limit theorem follows from the fact that the Laplace transform of a centered random variable is quadratic near the origin, as recalled in the next lemma. 
\begin{lemma}
\label{l.quad.lapl}
There exists a constant $C < \infty$ such that, if a random variable $X$ satisfies
\begin{equation} 
\label{e.integ.2X.app}
\E \Ll[ \exp(2|X|) \Rr]  \le 2 \quad \text{and} \quad \E[X] = 0,
\end{equation}
then for every $\lambda \in [-1,1]$,
\begin{equation*}  
\Ll|\log \E \Ll[ \exp \Ll( \lambda X \Rr)  \Rr]  - \frac {\E[X^2]}{2} \lambda^2 \Rr| \le C \lambda^3.
\end{equation*}
\end{lemma}
\begin{proof}
For $\lambda \in (-2,2)$, let $\psi(\lambda) := \log \E\left[\exp(\lambda X)\right]$, and let
\begin{equation*}  
\E_{\lambda}[\cdot] := \frac{\E\left[\, \cdot \, \exp(\lambda X)\right]}{\E\left[\exp(\lambda X) \right]}
\end{equation*}
be the ``Gibbs measure'' associated with $\lambda X$. Using the formal derivation rule 
$$
\partial_\lambda \E_\lambda[F] = \E_\lambda[FX] - \E_\lambda[F] \E_\lambda[X]
$$
and the assumption of \eqref{e.integ.2X.app}, we get that, for every $|\lambda| \le 1$,
\begin{align*}  
\psi'(\lambda) & = \E_\lambda[X], \\
\psi''(\lambda) & = \E_\lambda\left[X^2\right] - \left(\E_\lambda[X]\right)^2, \\
\psi'''(\lambda) & = \E_\lambda \left[X^3 \right] -3\E_\lambda\left[X^2\right] \E_\lambda[X] + 2 \left(\E_\lambda[X]\right)^3.
\end{align*}
In particular, by \eqref{e.integ.2X.app}, there exists a constant $C < \infty$ such that, for every $|\lambda| \le 1$, we have $|\psi'''(\lambda)| \le C$. Noting also that $\psi(0) = \psi'(0) = 0$ and that $\psi''(0) = \E[X^2]$, we obtain the result by a Taylor expansion of $\psi$ around $0$.
\end{proof}
Inspired by this, we now introduce new notation to measure the size of centered random variables. For every random variable $X$ and $s \in (1,2]$, $\theta \in (0,\infty)$, we write
\begin{equation}
\label{e.def.Ost.appendix}
X = \bar \O_s(\theta)
\end{equation}
to mean that
$$
\forall \lambda \in \R, \quad \log \E \left[ \exp\left( \lambda \theta^{-1}  X \right) \right] \leq  \lambda^2 \vee   | \lambda|^{\frac s {s-1}}  .
$$
\begin{lemma}
\label{l.bigO.barO}
For every $s \in (1,2]$, there exists a constant $C(s) < \infty$ such that, for every random variable $X$, 
\begin{equation*}  
X = \bar \O_s(1) \quad \implies \quad X = \O_s(C),
\end{equation*}
and
\begin{equation*}  
X = \O_s(1) \quad \text{and} \quad \E[X] = 0 \quad  \implies \quad X = \bar \O_s(C).
\end{equation*}
\end{lemma}
\begin{proof}
The first statement is immediate from  Lemma~\ref{l.lapl.bigO}, while the second one follows from Lemmas~\ref{l.lapl.bigO} and \ref{l.quad.lapl}.
\end{proof}
We now state a result analogous to Lemma~\ref{l.sum-O}, but for $\bar \O$ in place of $\O$.
\begin{lemma}	
\label{l.sum.barO}
For each $s \in (1,2]$, there exists $C(s) < \infty$ such that the following holds.
Let $\mu$ be a measure over an arbitrary measurable space $E$, let $\theta : E \to \R_+$ be a measurable function and $(X(x))_{x \in E}$ be a jointly measurable family of random variables such that, for every $x\in E$, $X(x) = \bar \O_s(\theta(x))$. We have
\begin{equation*} 
\int X \, d\mu = \bar \O_s \Ll(C   \int \theta \, d\mu  \Rr) .
\end{equation*}
\end{lemma}
\begin{proof}
By Lemma~\ref{l.bigO.barO}, there exists a constant $C(s)$ such that, for every $x \in E$, 
\begin{equation*}  
X(x) = \O_s(C \theta(x)).
\end{equation*}
Applying Lemmas~\ref{l.sum-O} and \ref{l.bigO.barO} thus yields the result.
\end{proof}
The key motivation for the introduction of the notation $\bar \O$ is that a sum of $k$ independent $\bar \O_s(\theta)$ random variables is $\bar \O_s(\sqrt{k} \theta)$, in agreement with the scaling of the central limit theorem.
\begin{lemma}
\label{l.barO}
For every $s \in (1,2]$, there exists $C(s) < \infty$ such that the following holds.
Let $\theta_1, \ldots, \theta_k \ge 0$ and $X_1, \ldots, X_k$ be \emph{independent} random variables such that, for  every $i \in \{1,\ldots,k\}$, $X_i = \bar \O_s (\theta_i)$. We have
\begin{equation}
\label{e.barO}
\sum_{i = 1}^k X_i = \bar \O_s \Ll( C \Big(\sum_{i = 1}^k \theta_i^2\Big)^{\frac 1 2} \Rr) .
\end{equation}
Moreover, if $\theta_i  = \theta_j$ for every $i, j \in \{1,\ldots,k\}$, then the constant $C$ in \eqref{e.barO} can be chosen equal to $1$. 
\end{lemma}
\begin{proof}[Proof of Lemma~\ref{l.barO}]
We define 
$$
\bar \theta := \Ll(\sum_{i = 1}^k \theta_i^2\Rr)^{\frac 1 2}.
$$ 
By Lemmas~\ref{l.lapl.bigO} and~\ref{l.bigO.barO}, in order to prove \eqref{e.barO}, it suffices to show that there exists $C(s) < \infty$ such that, for  every $\lambda \in \R$,
\begin{equation}
\label{e.barO2}
\log \E \Ll[ \exp \Ll(  \bar \theta^{\, -1} \lambda \sum_{i = 1}^k X_i \Rr)  \Rr] \le C\Ll(1+|\lambda|^{\frac s {s-1}}\Rr).
\end{equation}
We use independence and then the assumption $X_i = \bar \O_s(\theta_i)$ to bound the term on the left side by
\begin{align*}  
\sum_{i = 1}^k \log \E \Ll[ \exp \Ll( \bar \theta^{\, -1}  \lambda   X_i \Rr)  \Rr] & \le \sum_{i = 1}^k \Ll( \bar \theta^{\, -1} \theta_i \lambda \Rr)^2 \vee \Ll|\bar \theta^{\, -1} \theta_i\lambda \Rr|^{\frac s {s-1}} \\
& \le \lambda^2  + |\lambda|^{\frac s {s-1}} \bar \theta^{\, -\frac s {s-1}}  \sum_{i = 1}^k \theta_i^{\frac s {s-1}}.
\end{align*}
Since $s \le 2$, we have $\frac s{s-1} \ge 2$, and thus \eqref{e.barO2} follows from the observation that
\begin{equation*}  
\bar \theta \ge \Ll(\sum_{i = 1}^k \theta_i^{\frac s {s-1}}\Rr)^{\frac {s-1} s}.
\end{equation*}
When $\theta_i  = \theta_j$ for every $i, j \in \{1,\ldots,k\}$, without loss of generality we may set $\theta_i = 1$, and observe that
\begin{align*}
\sum_{i = 1}^k \log \E \Ll[ \exp \Ll( k^{-\frac 1 2} \lambda   X_i \Rr) \Rr]
& \le \sum_{i = 1}^k \Ll(k^{-\frac 1 2 } \lambda\Rr)^2 \vee  \left| k^{-\frac 1 2 } \lambda\right|^{\frac s {s-1}}  \\
& \le \lambda^2 \vee |\lambda|^{\frac s {s-1}},
\end{align*}
where in the last step we used the fact that $\frac s {s-1} \ge 2$.
\end{proof}
Since we also often encounter families of random variables with a finite range of dependence, we also provide a version of Lemma~\ref{l.barO} tailored to this situation.

\begin{lemma}
\label{l.barO.boxes}
For every $s \in (1,2]$, there exists a constant $C(s) < \infty$ such that the following holds. 
Let $\theta > 0$, $R \ge 1$, $\mcl Z \subset R\Z^d$, and for each $x \in \mcl Z$, let $X(x)$ be an $\mcl F(\cu_{2R}(x))$-measurable random variable such that $X(x) = \bar \O_s(\theta(x))$. We have
$$
\sum_{x \in \mcl Z} X(x) = \bar \O_s\Ll(C\, \Big( \sum_{x \in \mcl Z} \theta(x)^2 \Big)^{\frac 1 2} \Rr).
$$
\end{lemma}
\begin{proof}
We partition $\mcl Z$ into $\mcl Z^{(1)}, \ldots, \mcl Z^{(3^d)}$ in such a way that for every $j \in \{1,\ldots, 3^d\}$, if $x \neq x' \in \mcl Z^{(j)}$, then $|x-x'| \ge 3R \ge 2R + 1$. (That is, we define $\mcl Z^{(1)} = (3R\Z^d) \cap \mcl Z$, and so on with translates of $(3R\Z)^d$.) For each $j$, the random variables $(X(x))_{x \in \mcl Z^{(j)}}$ are independent. By Lemma~\ref{l.barO},
$$
\sum_{x \in \mcl Z^{(j)}} X(x) = \bar \O_s \Ll(C\, \Big( \sum_{x \in \mcl Z^{(j)}} \theta(x)^2 \Big)^{\frac 1 2}  \Rr) .
$$
The conclusion follows by summing over $j$ and applying Lemma~\ref{l.sum.barO}.
\end{proof}



\chapter{Function spaces and elliptic equations on Lipschitz domains}
\label{a.sobolev}

In this appendix, we collect some standard facts about Sobolev and Besov spaces on Lipschitz domains $U\subseteq \Rd$, including negative and fractional spaces. 

\begin{definition}[$C^{k,\alpha}$ domain]
\label{def.Lipdomain}
\index{Lipschitz domain}
\index{$C^{k,\alpha}$ domain}
Let $k\in\N$ and $\alpha\in (0,1]$. We say that a domain $U\subseteq\Rd$ is a \emph{$C^{k,\alpha}$ domain} if every point of~$\partial U$ has a neighborhood~$\mathcal{N}$ such that~$\partial U \cap \mathcal{N}$ can be represented, up to a change of variables, as the graph of a~$C^{k,\alpha}$ function of~$d-1$ of the variables. A~$C^{0,1}$ domain is also called a~\emph{Lipschitz domain}.
\end{definition}

We proceed by defining the Sobolev spaces $W^{\alpha,p}(U)$ for $\alpha\in\R$ and $p\in [1,\infty]$.

\index{Sobolev space|(}
\begin{definition}[Sobolev space $W^{k,p}(U)$]
\index{Sobolev space~$W^{\alpha,p}$}
\label{def.Sobolev}
For every $\be =(\be_1,\ldots,\be_d) \in \N^d$, we write $\partial^\beta := \partial^{\be_1}_{x_1} \cdots \partial^{\be_d}_{x_d}$, and $|\be|:= \sum_{i = 1}^d \be_i$. For every $k \in \N$ and $p \in [1,\infty]$, the Sobolev space $W^{k,p}(U)$ is defined by
\begin{equation*}  
W^{k,p}(U) := \{v \in L^p(U) \ : \ \forall |\beta| \le k, \ \partial^\beta u \in L^p(U)\},
\end{equation*}
endowed with the norm
\begin{equation}  \label{e.Wnorm}
\|v\|_{W^{k,p}(U)}  := \sum_{0 \le |\be|\le k} \|\partial^\beta v\|_{L^p(U)}.
\end{equation}
For $|U| < \infty$, we also define the rescaled norm
\begin{equation*}  
\|v\|_{\underline W^{k,p}(U)} := \sum_{0 \le |\be|\le k} |U|^{\frac{|\beta|-k}{d}} \|\partial^\beta v\|_{\underline L^p(U)}.
\end{equation*}
We denote by $W_0^{k,p}(U)$ the closure of $C_0^\infty(U)$ with respect to the norm in~\eqref{e.Wnorm}. We say that $v \in W_{\textrm{loc}}^{k,p}(U)$ if $v \in W^{k,p}(V)$ whenever~$\bar V$ is a compact subset of~$U$. 
\end{definition}

\begin{definition}
[Sobolev space $W^{\alpha,p}(U)$]
\index{Sobolev space~$W^{\alpha,p}$}
\label{def.sobolev.fractional}
\index{Sobolev space!fractional}
Let $\alpha\in (0,1)$ and $p\in [1,\infty]$. The fractional Sobolev $W^{\alpha,p}(U)$ seminorm is defined, for $p < \infty$, as 
\begin{equation} \label{e.Wseminorm.fractional}
\left[ v \right]_{W^{\alpha,p}(U)}^p :=(1-\alpha) \int_U \int_U \frac{|v(x)-v(y)|^p}{|x-y|^{d+\alpha p}} \, dx \,dy 
\end{equation}
and, for $p = \infty$,
\begin{equation*}  
\Ll[ v \Rr]_{W^{\al,\infty}(U)} := \esssup_{x,y \in U} \frac{|v(x) - v(y)|}{|x-y|^\al}.
\end{equation*}
For each $k \in \N$, we define the Sobolev space $W^{k+\al,p}(U)$ by
\begin{equation*} 
W^{k+\alpha,p}(U):= \left\{ v \in W^{k,p}(U)\,:\,   \left\| \nabla^k v \right\|_{W^{\alpha,p}(U)} < \infty\right\},
\end{equation*}
with the norm
\begin{equation} \label{e.Wnorm.fractional}
\left\| v \right\|_{W^{k+\alpha,p}(U)} := \left\| v \right\|_{W^{k,p}(U)} + \sum_{|\beta| = k} \left[ \partial^\beta v \right]_{W^{\alpha,p}(U)}.
\end{equation}
For $|U| < \infty$, we also define the rescaled norm
\begin{equation*}  
\|v\|_{\underline W^{k+\alpha,p}(U)} := \sum_{0 \le |\be|\le k} |U|^{\frac{|\beta|-k-\alpha}{d}} \|\partial^\beta v\|_{\underline L^p(U)} + \sum_{|\beta| = k} |U|^{-\frac 1 p}\left[ \partial^\beta v \right]_{W^{\alpha,p}(U)} .
\end{equation*}
\end{definition}

The space $W^{\alpha,p}(U)$ is a Banach space for every $\alpha\in (0,\infty)$ and $p\in [1,\infty]$. If $p \in [1,\infty)$, then the space $W^{\alpha,p}(U)$ is the closure of the set of smooth functions with respect to the norm in~\eqref{e.Wnorm}, a result first proved for integer~$\alpha$ by Meyers and Serrin~\cite{H=W}.

\smallskip

When $p=2$ and $\alpha>0$, we denote~$H^\alpha(U) : = W^{\alpha,2}(U)$. 

\smallskip

Note that the normalizing constant~$(1-\alpha)$ in~\eqref{e.Wseminorm.fractional}  is chosen so that, for every $v \in W^{1,p}(U)$, 
\begin{equation*} 
\lim_{\alpha \to 1} \left[ v \right]_{W^{\alpha,p}(U)}^p = c(U) \left\|\nabla v \right\|_{L^p(U)}^p.
\end{equation*}
Thus the normalizing constant provides continuity in seminorms. In the case $U = \R^d$, this can be complemented by  the fact that, for $v \in  \bigcup_{\alpha \in (0,1)} W^{\alpha,p}(\R^d)$, $$\lim_{\alpha \to 0}  \alpha \left[ v \right]_{W^{\alpha,p}(\R^d)}^p = \tfrac{1}{p} c(d) \left\| v \right\|_{L^p(\R^d)}^p.$$  Both of these properties are established in~\cite[Chapter 10]{Mazya}.

\smallskip

We denote by $W^{\alpha,p}_0(U)$ the closure of $C^\infty_c(U)$ with respect to~$\left\| \cdot\right\|_{W^{\alpha,p}(U)}$ .

\smallskip

We denote by $W^{\alpha,p}_{\mathrm{loc}}(U)$ the vector space of measurable functions on $U$ which belong to $W^{\alpha,p}(V)$ for every bounded open set $V$ with $\overline{V} \subseteq U$.

\begin{definition}[Sobolev space $W^{-\alpha,p}(U)$]
\label{def.negsob}
\index{negative Sobolev space~$W^{-\alpha,p}$}
Let $\al > 0$, $p \in [1,\infty]$, and let $p'$ be the conjugate exponent of $p$, that is, $p' = \tfrac{p}{p-1} \in [1,\infty]$. For every distribution $u$, we define 
\begin{equation}  \label{e.thedefWminus}
\|u\|_{W^{-\al,p}(U)} := \sup  \Ll\{ \int_U u v \ : \ v \in C^\infty_c(U), \ \|v\|_{W^{\al,p'}(U)} \le 1 \Rr\},
\end{equation}
where we abuse notation slightly and denote by $(u,v) \mapsto \int_U uv$ the canonical duality pairing. The Sobolev space $W^{-\al,p}(U)$ is the space of distributions for which this norm is finite. Equipped with this norm, it is a Banach space.
\end{definition}
\begin{remark}
\label{r.non.separable.sobolev}
For every $p\in(1,\infty]$, the space $W^{-\al,p}(U)$ may be identified with the dual of $W^{\al,p'}_0(U)$. 
When $p = 1$, a difficulty arises since the dual of $W^{\al,\infty}(U)$ is not canonically embedded into the space of distributions. To wit, it is elementary to associate a distribution with each element of the dual of $W^{\al,\infty}(U)$, by restriction; but this mapping is not injective. This problem is caused by the fact that the space $C^\infty_c(U)$ is not dense in $W^{\al,\infty}(U)$.
\end{remark}
\index{Sobolev space|)}

\begin{remark}
\label{r.RRT}
For every~$p\in (1,\infty]$, we can identify~$W^{-1,p}(U)$ with~$\nabla\cdot L^p(U;\Rd)$. That is, 
each element of $f\in W^{-1,p}(U)$ can be associated with a vector field $\mathbf{F} \in L^p(U;\Rd)$ such that
\begin{equation} 
\label{e.RRT}
\forall  w  \in W^{1,p'}_0(U), 
\quad
\int_U w f
= 
\int_U \mathbf{F} \cdot\nabla w
\end{equation}
and satisfying, for some $C(U,d)<\infty$, 
\begin{equation*} \label{}
\frac1C \left\| f \right\|_{W^{-1,p}(U)}
\leq
\left\| \mathbf{F} \right\|_{L^p(U)} 
\leq
C \left\| f \right\|_{W^{-1,p}(U)}. 
\end{equation*}
Note that in \eqref{e.RRT}, the left side is interpreted as the duality pairing between $w \in W^{1,p'}_0(U)$ and $f \in W^{-1,p}(U)$, while the right side is a standard Lebesgue integral.
The existence of such $\mathbf F$ is a consequence of the Riesz representation theorem and the Poincar\'e inequality (the latter tells us that $\| w \|_{W^{1,p}(U)}$ is equivalent to $\| \nabla w\|_{L^p(U)}$). Note that such an identification does not hold for~$p=1$, since, as explained above,~$W^{-1,1}(U)$ cannot be identified with the dual of~$W^{1,\infty}(U)$.

\smallskip

This identification can be generalized to~$W^{-k,p}(U)$ with~$k\in\N$ and~$p\in (1,\infty]$. Indeed, again by the Riesz representation theorem, we can identify each~$f \in W^{-k,p}(U)$ with $\left(f_\beta\right)_{\beta}\subseteq L^p(U)$, indexed by the multiindices $\beta$ with $|\beta| \leq k$, such that 
\begin{equation*} 
\forall  w  \in W^{k,p'}_0(U),  \quad
\int_U w f
= 
\sum_{|\beta| \leq k}\int_U f_\beta  \partial^\beta w.
\end{equation*}
\end{remark}

\smallskip

We also need to define the Besov space $B^{p,p}_\alpha(\partial U)$, where~$\partial U$ is the boundary of a Lipschitz domain~$U$. This is because we need a relatively sharp version of the Sobolev trace theorem for some results in Section~\ref{s.blayers}. We denote the $(d-1)$-dimensional Hausdorff measure on~$\partial U$ by~$\sigma$. 

\begin{definition}[Besov space $B_\alpha^{p,p}(\partial U)$]
\index{Besov space~$B_\alpha^{p,p}(\partial U)$}
Let $U\subseteq \Rd$ be a Lipschitz domain, $\alpha\in (0,1)$ and $p\in [1,\infty]$. The Besov $B_\alpha^{p,p}(\partial U)$ seminorm is defined as 
\begin{equation*} 
\left[ v \right]_{B_\alpha^{p,p}(\partial U)}^p := (1-\alpha) \int_{\partial U} \int_{\partial U} \frac{|v(x)-v(y)|^p}{|x-y|^{d-1+\alpha p}} \, d\sigma(x) \,d\sigma(y) ,
\end{equation*}
and the Besov space $B_\alpha^{p,p}(\partial U)$ by
\begin{equation*} 
B_\alpha^{p,p}(\partial U):= \left\{ v \in L^p(\partial U)\,:\,  \left[ v \right]_{B_\alpha^{p,p}(\partial U)}  < \infty\right\}.
\end{equation*}
This is a Banach space under the $B_\alpha^{p,p}(\partial U)$ norm, which is defined by
\begin{equation*} \label{}
\left\| v \right\|_{B_\alpha^{p,p}(\partial U)} := \left\| v \right\|_{L^p(\partial U)} + \left[ v \right]_{B_\alpha^{p,p}(\partial U)}.
\end{equation*}
\end{definition}

The next four statements give us various versions of the Sobolev and Sobolev-Poincar\'e inequalities which are used in the book. They can be found in~\cite{Mazya}. See also~\cite{Adams,Gilbarg-Trudinger}.

\begin{theorem}
[Global Sobolev inequality]
Fix $k \in \N$ and~$p \in \left(1, \frac{d}{k}\right)$. There exists~$C(k,p,d)<\infty$ such that, for every~$v \in W^{k,p}(\R^d)$,
\begin{equation} 
\label{e.Sobolev.global}
\left\| v \right\|_{L^{ \frac{dp}{d-kp}}(\R^d)} 
\leq 
C \left\|v\right\|_{W^{k,p}(\R^d)}.
\end{equation}
\end{theorem}

Recall that~$\mathcal{P}_k$ denotes the set of real-valued polynomials of order at most~$k$.

\begin{theorem}[Local Sobolev-Poincar\'e inequality]
Let $k \in \N$, $p \in \left[1, \frac{d}{k}\right)$ and $r>0$. There exists~$C(k,p,d)<\infty$ such that, for every~$v \in W^{k,p}(B_r)$,
\begin{equation} \label{e.Sobolev-Poincare.local}
\inf_{w \in \mathcal{P}_{k-1}}\left\| v -w \right\|_{\underline{L}^{\frac{dp}{d-kp}}(B_r)} \leq C r^k  \left\| \nabla^k v\right\|_{\underline{L}^{p}(B_r)}. 
\end{equation}
Moreover, for~$v \in W_0^{k,p}(B_r)$, then we may take $w=0$ on the left side of~\eqref{e.Sobolev-Poincare.local}.
\end{theorem}

For the borderline exponent $p = \frac dk$, there is a version of Sobolev's inequality, attributed to Trudinger, with $L^q$-norm on the left replaced by a certain exponential integral. However, since we do not use this result in the book, it will not be stated. For exponents~$p>\frac{d}{k}$, we have the following. 

\begin{theorem}
[Morrey's inequality]
Let $k \in \N$, $p  \in \left(\frac{d}{k},\infty\right]$ and $r>0$. There exists~a constant $C(k,p,d)$ such that for every $v \in W^{k,p}(B_r)$
\begin{equation} \label{e.Morrey.local}
\inf_{w \in \mathcal{P}_{k-1}} \left\| v - w \right\|_{L^\infty(B_r)} \leq C r^{k} \left\|\nabla^k v\right\|_{\underline{L}^{p}(B_r)} .
\end{equation}
Moreover, if $v \in W_0^{k,p}(B_r)$, then $w$ in the infimum above can be taken to be zero.  Finally, if $v \in W^{k,p}(\R^d)$, then
\begin{equation} \label{e.Morrey.global}
\left\| v \right\|_{L^\infty(\R^d)} \leq C \left\|v\right\|_{W^{k,p}(\R^d)} .
\end{equation}
\end{theorem}

We next record a version of the Sobolev-Poincar\'e-Wirtinger inequality. In most of the book, we just need the statement in the case $\alpha =1$. However, we also need a fractional version for cubes in Lemmas~\ref{l.convolution2} and~\ref{l.convolution0}. The former case is classical and can be found for example in~\cite[Section 6.3.4]{Mazya}. The latter case can be found in~\cite[Section 10.2.2]{Mazya}.

\begin{proposition}[{Sobolev-Poincar\'e-Wirtinger inequality}]
\label{p.sobolevpoincare}
\index{{Sobolev-Poincar\'e inequality}}
Fix a bounded Lipschitz domain $U\subseteq\Rd$, $\alpha \in (0,1]$, $p\in \left(1, \tfrac{d}{\alpha} \right)$, and denote 
\begin{equation*} \label{}
p^*=p^*(\alpha,p):= \frac{dp} {d-\alpha p}.
\end{equation*}
Suppose also that~$\alpha = 1$ or else that~$U$ is a cube. Then there exists $C(\alpha,p,U,d)<\infty$ such that, for every $u\in W^{\alpha,p}(U)$, 
\begin{equation*} \label{}
\left\| u - \left( u \right)_U \right\|_{L^{p^*}(U)}
\leq 
C  \left[ u \right]_{W^{\alpha,p}(U)}. 
\end{equation*}
\end{proposition}

\begin{remark}
In the case $U$ is either a ball or a cube, for which the above inequality is mainly applied in this book, the constant $C$ above has the form
\begin{equation*} 
C =  \frac{C(d,p)}{(d-\alpha p)^{1-1/p}}.
\end{equation*}
See~\cite[Section 10.2.3]{Mazya}. 
\end{remark}

There is also an embedding between fractional Sobolev spaces, recorded in the following proposition, which follows from~\cite[Theorem 7.58]{Adams} and the extension theorem, Proposition~\ref{p.extension} below. 

\begin{proposition}[{Embedding between fractional Sobolev spaces}]
\label{p.sobolevembedding}
\index{{Sobolev embedding!fractional}}
Fix a Lipschitz domain $U\subseteq\Rd$, $\alpha \in (0,1)$, $\beta \in (0,\alpha)$,~$p\in (1,\infty)$ and~$q\in [1,\infty]$ such that
\begin{equation*} 
\beta \leq \alpha - d \left( \frac 1p - \frac 1q \right).
\end{equation*}
Then there exists $C(\alpha,\beta,p,q,k,U,d)<\infty$ such that, for every $u\in W^{\alpha,p}(U)$, 
\begin{equation*} \label{}
\left[ u \right]_{W^{\beta,q}(U)} \leq C  \left[ u \right]_{W^{\alpha,p}(U)}. 
\end{equation*}
\end{proposition}

A proof of the following Sobolev extension theorem for $W^{\alpha,p}$ spaces in Lipschitz domains can be found in~\cite[Chapter VI, Theorem 5]{Stein} for $\alpha \in\N$ and, for $\alpha\in (0,1)$, in~\cite[Theorem 5.4]{Hitch}. The general statement follows from these. 

\begin{proposition}[Extension theorem]
\label{p.extension}
\index{Sobolev extension theorem}
Let $U\subseteq \Rd$ be a bounded Lipschitz domain, $\alpha \in (0,\infty)$ and $p \in (1,\infty)$. The restriction operator $W^{\alpha,p}(\Rd) \to W^{\alpha,p}(U)$ has a bounded linear right inverse. That is, there exists $C(U,\alpha,p,d)<\infty$ and a linear operator 
\begin{equation*} \label{}
\Ext: W^{\alpha,p}(U) \to W^{\alpha,p}(\Rd)
\end{equation*}
such that, for every $u \in W^{\alpha,p}(U)$, 
\begin{equation*} \label{}
\Ext(u) = u \quad \mbox{a.e. in} \ U
\end{equation*}
and
\begin{equation*} \label{}
\left\| \Ext(u) \right\|_{W^{\alpha,p}(\Rd)} \leq C \left\| u \right\|_{W^{\alpha,p}(U)}. 
\end{equation*}
\end{proposition}

The following version of the trace theorem for Lipschitz domains is a special case of~\cite[Chapter VII, Theorem 1]{JW} (see also~\cite{Marschall}). \index{Sobolev trace theorem}

\begin{proposition}[Sobolev trace theorem]
\label{p.trace}
\index{trace operator}
Let $U\subseteq \Rd$ be a bounded Lipschitz domain, $\alpha\in (0,\infty)$ and $p \in (1,\infty)$ be such that $\alpha> \frac1p$. The linear operator $C^\infty(\overline{U}) \to \mathrm{Lip}(\partial U)$ that restricts a smooth function on $\overline{U}$ to $\partial U$ has an extension to a bounded linear mapping $W^{\alpha,p}(U) \to B_{\alpha - \frac1p}^{p,p}(\partial U)$.
That is, there exists 
$C(\alpha,p,U,d)<\infty$ and a linear operator 
\begin{equation*} \label{}
\Tr: W^{\alpha,p}(U) \to B_{\alpha - \frac1p}^{p,p}(\partial U)
\end{equation*}
such that, for every $u\in  W^{\alpha,p}(U)$, 
\begin{equation} \label{e.trace embedding}
\left\| \Tr(u) \right\|_{B_{\alpha - \frac1p}^{p,p}(\partial U)} \leq C \left\| u \right\|_{W^{\alpha,p}( U)}
\end{equation}
and, for every $u\in C^\infty(\overline{U})$,
\begin{equation*} \label{}
\Tr(u) = u \quad \mbox{on} \ \partial U. 
\end{equation*} 
\end{proposition}

We next show that the trace operator~$\Tr$ is surjective and has a bounded linear right inverse. This is a consequence of the solvability of the Dirichlet problem for the Poisson equation in Lipschitz domains, which was proved in~\cite{JeKe}. For another proof, see~\cite[Theorem 10.1]{FMM}). Note that in~$C^{1,1}$ domains, the result is true for $p\in (1,\infty)$ by an argument which is simpler (cf.~\cite{ADN}). 

\begin{proposition}
\label{p.DPpoisson}
Let $U\subseteq \Rd$ be a bounded Lipschitz domain and $p \in \left[ \tfrac 32,3\right]$. For every $g \in B_{1-\frac1p}^{p,p}(\partial U)$ and $f\in W^{-1,p}(U)$, there exists a unique solution of the Dirichlet problem
\begin{equation}
\label{e.DirichletPoisson}
\left\{ 
\begin{aligned}
& -\Delta u = f & \mbox{in} & \ U, \\
& u = g & \mbox{on} & \ \partial U,
\end{aligned}
\right. 
\end{equation}
where the boundary condition is interpreted to mean that $\Tr(u) = g$. Moreover, there exists $C(U,p,d)<\infty$ such that 
\begin{equation*} \label{}
\left\| u \right\|_{W^{1,p}(U)} \leq C \left(  \left\| g \right\|_{B_{1-\frac1p}^{p,p}(\partial U)} + \left\| f \right\|_{W^{-1,p}(U)}\right). 
\end{equation*}
\end{proposition}

We next state some~$H^2$ estimates for solutions of the Poisson equation. The interior estimate is essentially just the Caccioppoli inequality (see Lemma~\ref{l.Caccioppoli appendix} below), while a global version necessarily requires some regularity of boundary of the domain. 

\begin{lemma}[Interior $H^2$ estimate]
\label{l.H2estimate.interior}
There exists $C(d) < \infty$ such that if $f\in L^2(B_1)$ and $u \in H^1(B_1)$ satisfy the Poisson equation
\begin{equation*} \label{}
-\Delta u = f \quad \mbox{in} \ B_1, 
\end{equation*}
then $u\in H^2_{\mathrm{loc}}(B_1)$ and
\begin{equation} 
\label{e.interior.H2}
\left\| u \right\|_{H^2(B_{1/2})}
\leq 
C \left( \left\| u \right\|_{L^2(B_1)} + \left\| f \right\|_{L^2(B_1)}  \right). 
\end{equation}
\end{lemma}
\begin{proof}
Since $\partial_j u$ satisfies the equation $-\Delta \partial_j u= \partial_j f$, the Caccioppoli estimate (Lemma~\ref{l.Caccioppoli appendix}) gives the result. 
\end{proof}

Proposition~\ref{p.DPpoisson} can be improved if~$U$ is more regular than just Lipschitz. We next record a global~$H^2$ estimate for $C^{1,1}$ and for convex domains. The proofs can be found in~\cite[Theorems 2.4.2.5 and 3.1.2.1]{Grisvard}.

\begin{proposition}[Global $H^2$ estimate]
\label{p.H2estimate}
Let $U\subseteq \Rd$ be a bounded domain which is convex or~$C^{1,1}$. There exists $C(U,d) < \infty$ such that if~$f\in L^2(U)$ and~$u\in H^1(U)$ solves the Dirichlet problem 
\begin{equation}
\label{e.Dirichlet.forH2}
\left\{ 
\begin{aligned}
& -\Delta u = f & \mbox{in} & \ U, \\
& u = 0 & \mbox{on} & \ \partial U,
\end{aligned}
\right. 
\end{equation}
then~$u\in H^2(U)$ and
\begin{equation}
\label{e.global.H2}
\left\| u \right\|_{H^2(U)} 
\leq
C \left\| f \right\|_{L^2(U)}. 
\end{equation}
\end{proposition}

\smallskip

The proof of Corollary~\ref{c.mP.gradient} requires the following global $H^2$ estimate for the Neumann problem in a cube. This result is true in more generality (similar to the Dirichlet case given in Proposition~\ref{p.H2estimate}), but the case of a cube has a short proof by the reflection principle which we include here for the convenience of the reader. 

\begin{lemma}
\label{l.H2.est}
There exists $C(d) < \infty$ such that if $\cu\subseteq \Rd$ is a cube, $f \in L^2(\cu)$ is such that $\int_\cu f = 0$ and~$w \in H^1(\cu)$ solves the Neumann problem
\begin{equation}
\label{e.H2.Neumann}
\left\{ 
\begin{aligned}
& -\Delta w = f & \quad \mbox{in} & \ \cu, \\
&  \partial_\nu w = 0 &\quad \mbox{on} & \ \partial \cu,
\end{aligned} 
\right.
 \end{equation}
then $w \in H^2(\cu)$ and
\begin{equation}
\label{e.H2.H2}
 \left\| \nabla^2 w \right\|_{L^2(\cu)}  \leq C   \left\| f \right\|_{L^2(\cu)}. 
\end{equation} 
\end{lemma}
\begin{proof}
By translation and scaling, we may assume that $\cu := (0,1)^d$. Denote $\td \cu := (-1,1)^d$ and define $\td f \in L^2(\td \cu)$ by
\begin{equation*}  
\td f(x_1,\ldots,x_d) := f(|x_1|,\ldots, |x_d|).
\end{equation*}
Let $\td w \in H^1_{\per}(\td \cu)$ be the unique $\td \cu$-periodic function with mean zero solving
\begin{equation*}  
-\Delta \td w = \td f.
\end{equation*}
We have that $\td w \in H^2_\per(\td \cu)$ as well as the estimate
\begin{equation*}  
\int_{\td \cu} |\nabla^2 \td w|^2 \le \int_{\td \cu} \td f^2. 
\end{equation*}
This can be obtained by testing, for each $j\in\{1,\ldots,d\}$, the equation $-\Delta \partial_j \td w = \partial_j f$ with $\partial_j \td w$, summing over $j$ and then applying Young's inequality (or, alternatively, applying Lemma~\ref{l.H2estimate.interior} and using periodicity). 
To conclude, we check that the restriction of $\td w$ to $\cu$ solves the Neumann problem~\eqref{e.H2.Neumann}. It only remains to check that the boundary condition holds, and for this we argue by symmetry. The function $\td w$ is invariant under the change of coordinates $x_i \mapsto -x_i$, for each $i \in \{1,\ldots,d\}$.  We also recall that the trace mapping $v \mapsto   \partial_\nu v$  is continuous operator from $H^2_\per(\td \cu) \to  L^2(\partial \cu)$. Moreover, the function $\td w$ can be approximated by smooth functions in $H^2_\per(\td \cu)$ that are also invariant under these changes of coordinates. These approximations all have a null image under $T$, so the proof is complete.
\end{proof}

\chapter{The Meyers \texorpdfstring{$L^{2+\delta}$}{L(2+delta)} estimate}
\label{a.meyers}
\index{Meyers estimate|(}

In this appendix, we give a mostly self-contained proof of the Meyers improvement of integrability estimate for gradients of solutions of uniformly elliptic equations. The estimate states roughly that a solution of a uniformly elliptic equation must have slightly better regularity than a typical~$H^1$ function: it belongs to the space $W^{1,2+\delta}$ for some positive exponent~$\delta$. We include both a local version (Theorem~\ref{t.Meyers appendix}) as well as a global version in Lipschitz domains with Dirichlet boundary conditions (Theorem~\ref{t.Meyers appendix global}).

\smallskip

The Meyers estimate (in any form) is a consequence of the Caccioppoli and Sobolev inequalities, which immediately yield a reverse H\"older inequality for solutions, and Gehring's lemma, a measure-theoretic fact stating that a reverse H\"older inequality implies an improvement of integrability. 

\smallskip

We first give the statement of the interior Meyers estimate. Throughout we assume that $\a\in\Omega$ is fixed, that is, $\a$ is a measurable map from $\Rd$ to the set of symmetric matrices with eigenvalues belonging to~$[1,\Lambda]$. 

\begin{theorem}[Interior Meyers estimate]
\label{t.Meyers appendix}
Fix $r>0$ and $p\in(2,\infty)$. Suppose that $h \in W^{-1,p}(B_r)$ and that $u \in H^1(B_r)$ satisfy
\begin{equation} \label{e.Meyers eq}
-\nabla \cdot \left( \a(x) \nabla u\right) = h 
\quad \mbox{in} \ B_r.
\end{equation}
Then there exist~$\delta(d,\Lambda) >0$ and $C(d,\Lambda)<\infty$ such that
\begin{equation} \label{e.Meyers appendix}
\left\| \nabla u \right\|_{\underline{L}^{(2+\delta) \wedge p}(B_{r/2})} \leq C \left( \left\| \nabla u \right\|_{\underline{L}^{2}(B_{r})}  + \left\| h \right\|_{\underline{W}^{-1,(2+\delta)\wedge p}(B_{r})}  \right).
\end{equation}
 \end{theorem}

We begin the proof of Theorem~\ref{t.Meyers appendix} with the Caccioppoli inequality, perhaps the most basic elliptic regularity estimate which is also used many times for other purposes throughout the text. 

\begin{lemma}
[Interior Caccioppoli inequality] 
\label{l.Caccioppoli appendix}
\index{Caccioppoli inequality}
Fix $r>0$. Suppose that $h \in H^{-1}(B_r)$ and $u \in H^1(B_r)$ satisfy 
\begin{equation} \label{e.Caccioppoli appendix pde}
-\nabla \cdot \left( \a(x) \nabla u\right) = h \quad \mbox{in } B_r.
\end{equation}
Then there exists $C(d,\Lambda)<\infty$ such that 
\begin{equation} \label{e.Caccioppoli appendix}
\left\| \nabla u \right\|_{\underline{L}^2(B_{r/2})}
\leq 
C
\left( 
\frac{1}{r} \left\| u - \left( u \right)_{B_r} \right\|_{\underline{L}^2(B_r)} 
+
\left\| h \right\|_{\underline{H}^{-1}(B_r)} 
\right) .
\end{equation}
\end{lemma}
\begin{proof}
By subtracting a constant from $u$, we may suppose~$\left( u \right)_{B_r} = 0$. Fix a cutoff function $\phi \in C^\infty_c(B_r)$ satisfying
\begin{equation} 
\label{e.phicutoffCacc}
0\leq \phi \leq 1, \quad 
\phi = 1 \quad \mbox{in} \ B_{r/2}, \quad
\left| \nabla \phi \right| \leq 4r^{-1},  
\end{equation}
and test the equation with $\phi^2u$ to get
\begin{equation*} \label{}
\fint_{B_r} \phi^2 \nabla u\cdot \a\nabla u
= 
\fint h \phi^2 u  
- 
\fint_{B_r} 2 \phi u \nabla \phi \cdot \a\nabla u.
\end{equation*}
Recall that we use the expression $\fint h \phi^2 u$ to denote the (normalized) pairing between $h\in H^{-1}(B_r)$ and $\phi^2u \in H^1_0(B_r)$, as explained below~\eqref{e.thedefWminus}. We have that 
\begin{equation*} \label{}
\left|\fint h \phi^2 u    \right|
\leq 
\left\| h \right\|_{\underline{H}^{-1}(B_r)} 
\left\| \phi^2 u \right\|_{\underline{H}^1(B_r)}
\end{equation*}
and, by the Poincar\'e inequality and~\eqref{e.phicutoffCacc}, 
\begin{align*} \label{}
\left\| \phi^2 u \right\|_{\underline{H}^1(B_r)}
\leq 
C \left\| \nabla \left( \phi^2 u \right) \right\|_{\underline{L}^2(B_r)}
\leq
\frac{C}{r} \left\| u \right\|_{\underline{L}^2(B_r)} 
+ C \left\| \phi \nabla u \right\|_{\underline{L}^2(B_r)}.
\end{align*}
By Young's inequality and the upper bound $\left| \a \right| \leq\Lambda$, we have 
\begin{align*} \label{}
\left| \fint_{B_r} 2 \phi u \nabla \phi \cdot \a\nabla u \right|
&
\leq 
\frac12 \fint_{B_r} \phi^2 \left| \nabla u \right|^2 
+ C \fint_{B_r} \left| \nabla \phi \right|^2 u^2
\\ &
\leq 
\frac12 \left\| \phi \nabla u \right\|_{\underline{L}^2(B_r)}^2
+ \frac{C}{r^2} \left\| u \right\|_{\underline{L}^2(B_r)}^2. 
\end{align*}
By the uniform ellipticity assumption, $\a \geq \Id$, and thus we have 
\begin{equation*} \label{}
\fint_{B_r} \phi^2 \nabla u\cdot \a\nabla u
\geq
\fint_{B_r} \phi^2 \left| \nabla u \right|^2 = 
\left\| \phi \nabla u \right\|_{\underline{L}^2(B_r)}^2. 
\end{equation*}
Combining the five previous displays, we obtain
\begin{equation*} \label{}
\left\| \phi \nabla u \right\|_{\underline{L}^2(B_r)}^2
\leq 
\frac{C}{r^2} \left\| u \right\|_{\underline{L}^2(B_r)}^2
+ C\left\| h \right\|_{\underline{H}^{-1}(B_r)} 
\left( 
\frac{C}{r} \left\| u \right\|_{\underline{L}^2(B_r)} 
+ C \left\| \phi \nabla u \right\|_{\underline{L}^2(B_r)}\right).
\end{equation*}
Young's inequality then yields
\begin{equation*} \label{}
\left\| \phi \nabla u \right\|_{\underline{L}^2(B_r)}
\leq 
\frac{C}{r} \left\| u \right\|_{\underline{L}^2(B_r)} 
+
C\left\| h \right\|_{\underline{H}^{-1}(B_r)}. 
\end{equation*}
As $\left\| \phi \nabla u \right\|_{\underline{L}^2(B_r)} \geq \left\| \nabla u \right\|_{\underline{L}^2(B_{r/2})}$ by~\eqref{e.phicutoffCacc}, the proof is complete. 
\end{proof}

The Caccioppoli inequality can be combined with the Sobolev-Poincar\'e inequality to obtain a reverse H\"older inequality for the gradients of solutions. Recall that we denote~$2_\ast:=  \frac{2d}{d+2}$. 

\begin{corollary} 
\label{c.Meyers reverse}
Under the assumptions of Lemma~\ref{l.Caccioppoli appendix}, there exists~$C(d,\Lambda) < \infty$ such that 
\begin{equation} 
\label{e.rev holder appendix}
\left\| \nabla u \right\|_{\underline{L}^2(B_{r/2})}
\leq
C \left( 
\left\| \nabla u \right\|_{\underline{L}^{2_\ast}(B_r)}
+
\left\| h \right\|_{\underline{H}^{-1}(B_r)}
\right). 
\end{equation}
\end{corollary}
\begin{proof}
By the Sobolev-Poincar\'e inequality,
\begin{equation*} 
\left\| u - \left( u \right)_{B_r} \right\|_{\underline{L}^2(B_r)}
\leq 
Cr \left\| \nabla u \right\|_{\underline{L}^{2_\ast}(B_r)}.
\end{equation*}
Plugging this into~\eqref{e.Caccioppoli appendix} gives the result. 
\end{proof}

It is a basic real analytic fact that a reverse H\"older inequality like~\eqref{e.rev holder appendix} implies a small improvement of integrability. 
This is formalized in the following lemma, which is usually attributed to Gehring. \index{Gehring lemma}

\begin{lemma}[Gehring's lemma]
\label{l.Meyers really}
Fix~$p\in(2,\infty)$, $q \in (1,2)$, $K\geq 1$, and $R>0$. Suppose that $f \in L^2(B_R)$, $g \in L^p(B_R)$ and that $f$ satisfies the following reverse H\"older inequality, whenever $B_r(z) \subset B_R$, 
\begin{equation} 
\label{e.revHol}
\left\| f \right\|_{\underline{L}^2(B_{r/2}(z))} \leq K \left( \left\| f \right\|_{\underline{L}^q(B_{r}(z))} + \left\| g \right\|_{\underline{L}^2(B_{r}(z))}  \right).
\end{equation}
Then there are constants $\delta(q,K,d) \in (0,1)$ and $C(q,K,d)<\infty$ such that
\begin{equation*} 
\left\| f \right\|_{\underline{L}^{(2+\delta) \wedge p}(B_{R/2})} \leq C \left( \left\| f \right\|_{\underline{L}^2(B_{R})} + \left\| g \right\|_{\underline{L}^{(2+\delta) \wedge p}(B_{R})}  \right).
\end{equation*}
\end{lemma}

The proof of Gehring's lemma requires a measure theoretic, geometric covering tool. There are several choices that could be made and here we use the Vitali covering lemma, which we state next (for a proof, see~\cite[Theorem 1.5.1]{Evans-Gariepy}). 

\begin{lemma}[{Vitali covering lemma}] 
\label{l.vitali}
Let $\{B^\alpha\}_{\alpha\in\Gamma}$ be a 
family of balls in~$\R^d$ with uniformly bounded radii. Then there exists a countable, pairwise disjoint subfamily $\left\{B^{\alpha_j}\right\}_{j \in \N}$ such that 
\begin{equation*} 
\bigcup_{\alpha\in\Gamma} B^\alpha \subset \bigcup_{j \in \N} 5 B^{\alpha_j} .
\end{equation*}
\end{lemma}

We next proceed with the proof of Gehring's lemma. 

\begin{proof}[Proof of Lemma~\ref{l.Meyers really}] 
Fix $ \frac12 \leq s < t \leq 1$. For any $z \in B_{s R}$ and $r \in \left(0,(t-s)R\right]$, set
\begin{equation*} 
E(z,r) :=  \left\| f \right\|_{\underline{L}^2(B_{r/2}(z))} +  \left\| g \right\|_{\underline{L}^2(B_{r}(z))} .
\end{equation*}
Take 
\begin{equation*} 
\lambda_0 := \left(\frac{20}{t-s}\right)^{\frac d2} \left( \left\| f \right\|_{\underline{L}^2(B_{R})} +  \left\| g \right\|_{\underline{L}^2(B_{R})} \right), 
\end{equation*}
and let $\mathcal{L}$ be the set of Lebesgue points of both $f$ and $g$, that is 
\begin{equation*} 
\mathcal{L} := \left\{ x \in B_{sR} \, : \,  \lim_{r \to 0} \left\| f -f(x)\right\|_{\underline{L}^2(B_{r}(x))} = 0 
\ \ \mbox{and} \ \
\lim_{r \to 0} \left\| g -g(x) \right\|_{\underline{L}^2(B_{r}(x))} = 0\right\}.
\end{equation*} 
Then $\left| B_{sR} \setminus \mathcal{L} \right| = 0$. Now, for each $\lambda > \lambda_0$ and $z \in \mathcal{L} \cap \{|f|>\lambda\}$ there exists~$r_z \in \left(0,\tfrac{t-s}{10}R\right]$ such that 
\begin{equation} \label{e.MAstop}
E(z,r_z) = \lambda , \quad \mbox{while } \sup_{r \in \left(r_z,\frac{t-s}{10}R\right]} E(z,r) \leq \lambda. 
\end{equation}
Indeed, the existence of such $r_z$ is clear since $z \in \mathcal{L} \cap \{|f|>\lambda\}$, giving
\begin{equation*} 
\lim_{r \to 0} E(z,r) = |f(z)| +  |g(z)| > \lambda,
\end{equation*}
the mapping $r \mapsto E(z,r)$ is continuous, and 
\begin{equation*} 
E\left(z,\tfrac{t-s}{10}R\right) \leq \left(\frac{20}{t-s}\right)^{\frac d2} \left( \left\| f \right\|_{\underline{L}^2(B_{R})} +  \left\| g \right\|_{\underline{L}^2(B_{R})} \right) \leq \lambda_0<\lambda.
\end{equation*}
By~\eqref{e.MAstop} and the reverse H\"older inequality~\eqref{e.revHol} we obtain
\begin{equation} \label{e.MAlambda}
\lambda = E(z,r_z) \leq K \left\| f \right\|_{\underline{L}^q(B_{r}(z))} + 2K \left\| g \right\|_{\underline{L}^2(B_{r}(z))}.
\end{equation}
Set, for any Borel set $A \subset B_{tR}$,
\begin{equation*} 
\mu_f(A) := \int_{B_{tR} \cap A} |f(x)|^q \, dx \quad \mbox{and} \quad \nu_g(B_{tR} \cap A) := \int_{A} |g(x)|^2 \, dx,
\end{equation*}
and similarly for $\nu_f(A)$. We have from~\eqref{e.MAlambda} that
\begin{equation*} 
|B_{r_z}| \leq \left(\frac{8K}{\lambda}\right)^q \mu_f\left( B_{r_z}(z) \cap \left\{|f| > \frac{\lambda}{8K} \right\}\right) 
+ \left(\frac{8K}{\lambda}\right)^2 \nu_g \left( B_{r_z}(z) \cap \left\{|g| > \frac{\lambda}{8K} \right\}\right) .
\end{equation*}
On the other hand, by~\eqref{e.MAstop} we get
\begin{equation*} 
\nu_f\left( B_{5r_z}(z) \cap \left\{|f| \geq \lambda \right\}\right) \leq \left\| f \right\|_{L^2(B_{5r_z}(x))}^2 \leq 5^d \lambda^2 |B_{r_z}|.
\end{equation*}
It follows by the two previous displays that 
\begin{multline*} 
\nu_f\left( B_{5r_z}(z) \cap \left\{|f| \geq \lambda \right\}\right) 
\\ \leq
5^d \left(8K\right)^q \lambda^{2-q}\mu_f\left( B_{r_z}(z) \cap \left\{|f| \geq \frac{\lambda}{8K} \right\}\right) 
+ 5^d \left(8K\right)^2 \nu_g \left( B_{r_z}(z) \cap \left\{|g| \geq \frac{\lambda}{8K} \right\}\right) .
\end{multline*}
Applying then Vitali's covering lemma with $\{B_{r_z}(z)\}$ yields, taking into account that $B_{r_z}(z)\subset B_{tR}$, that
\begin{multline} \label{e.MAlevels}
\nu_f\left( B_{sR} \cap \left\{|f| \geq \lambda \right\}\right) 
\\ \leq 
5^d \left(8K\right)^q \lambda^{2-q}\mu_f\left(  \left\{|f| \geq \frac{\lambda}{8K} \right\}\right) 
+ 5^d \left(8K\right)^2 \nu_g \left(  \left\{|g| \geq \frac{\lambda}{8K} \right\}\right) .
\end{multline}
Let then $\ep := (2+\delta) \wedge p$ and set, for $m> \lambda_0$, $f_m := |f| \wedge m$. Integration and the definition of $\lambda_0$ give that 
\begin{multline*} 
\left\| f_m \right\|_{L^{2+\ep}(B_{sR})}^{2+\ep} = \ep \int_{0}^m \lambda^{\ep-1}  \nu_f\left( B_{sR} \cap \left\{|f| \geq \lambda \right\}\right) \, d\lambda
\\ \leq \lambda_0^{2+\ep}|B_{tR}| + \ep \int_{\lambda_0}^m \lambda^{\ep-1}  \nu_f\left( B_{sR} \cap \left\{|f| \geq \lambda \right\}\right) \, d\lambda.
\end{multline*}
For the latter term we may use~\eqref{e.MAlevels} and obtain, after change of variables, that
\begin{multline*} 
\ep \int_{\lambda_0}^m \lambda^{\ep-1}  \nu_f\left( B_{sR} \cap \left\{|f| \geq \lambda \right\}\right) \, d\lambda
\\ \leq 
5^d (8K)^{2+\ep} \left( \ep \int_{\lambda_0/(8K)}^{m/(8K)} \lambda^{1+\ep-q} \mu_f\left(  \left\{|f| \geq \lambda \right\}\right)  \, d\lambda + \left\| g \right\|_{L^{2+\ep}(B_{tR})}^{2+\ep} \right).
\end{multline*}
The first term on the right can be further estimated as
\begin{equation*} 
\ep \int_{\lambda_0/(8K)}^{m/(8K)} \lambda^{1+\ep-q} \mu_f\left(  \left\{|f| \geq \lambda \right\}\right)  \, d\lambda 
\leq 
\frac{\delta}{2-q+\delta} \left\| f_m \right\|_{L^{2+\ep}(B_{tR})}^{2+\ep}.
\end{equation*}
Now we can exploit the smallness of $\delta$. Indeed, we choose it so that 
\begin{equation*} 
\frac{\delta}{2-q+\delta}5^d (8K)^{2+\delta} = 2^{-d-1} ,
\end{equation*}
which is possible since $q<2$. Collecting all the estimates gives, after taking averages,
\begin{equation*} 
\left\| f_m \right\|_{\underline{L}^{2+\ep}(B_{sR})} \leq \frac12 \left\| f_m \right\|_{\underline{L}^{2+\ep}(B_{tR})} + \frac{C}{(t-s)^{\frac d2}} \left(\left\| f \right\|_{\underline{L}^{2}(B_{R})} +  \left\| g \right\|_{\underline{L}^{p}(B_{R})}  \right).
\end{equation*}
An application of Lemma~\ref{l.simpleiter} below then yields that 
\begin{equation*} 
\left\| f_m \right\|_{\underline{L}^{2+\ep}(B_{R/2})} \leq  C \left(\left\| f \right\|_{\underline{L}^{2}(B_{R})} +  \left\| g \right\|_{\underline{L}^{p}(B_{R})}  \right). 
\end{equation*}
We conclude by sending~$m\to \infty$, using the monotone convergence theorem.
\end{proof}

\begin{lemma} 
\label{l.simpleiter}
Suppose that $A,\xi \geq 0$ and $\rho : \left[\tfrac12, 1\right) \to [0,\infty)$ satisfies 
\begin{equation*} \label{}
\sup_{t \in \left[\frac12,1\right)} (1-t)^\xi \rho(t) < \infty
\end{equation*}
and, for all $\frac12 \leq s < t < 1$,
\begin{equation} \label{e.iterineq}
\rho(s) \leq \frac12 \rho(t) + (t-s)^{-\xi} A.
\end{equation}
Then there exists~a constant $C(\xi)<\infty$ such that $\rho\left(\tfrac 12 \right) \leq C A.$ 
\end{lemma}
\begin{proof}
Denote $M:= \sup_{t \in \left[\frac12,1\right)} (1-t)^\xi \rho(t)$. Fix $s \in \left[\frac12,1\right)$ and take~$\delta := 1 - \left(\tfrac 23\right)^{1/\xi}$ and~$t := 1 - (1-\delta)(1-s)$ in~\eqref{e.iterineq} and multiply by~$(1-s)^{\xi}$ to obtain
\begin{equation*} 
(1-s)^{\xi} \rho(s) \leq \frac12 \left(\frac{1-s}{1-t} \right)^\xi (1-t)^{\xi} \rho(t) + (1-s)^\xi (t-s)^{-\xi} A = \frac 34 (1-t)^{\xi} \rho(t)  + \delta^{-\xi} A. 
\end{equation*}
This yields
\begin{equation*} \label{}
(1-s)^{\xi} \rho(s) \leq \frac34 M + \delta^{-\xi} A.
\end{equation*}
Taking the supremum of the left side over $s\in \left[ \frac12,1\right)$ gives 
\begin{equation*} \label{}
M \leq \frac 34 M + \delta^{-\xi} A.
\end{equation*}
Thus $M \leq 4 \delta^{-\xi} A \leq CA$. Since $\rho(\frac12) \leq 2^\xi M \leq CM$, we obtain the lemma. 
\end{proof}

\begin{proof}[{Proof of Theorem~\ref{t.Meyers appendix}}]
We first consider the case that $h$ can be written in the form 
\begin{equation*} \label{}
h = \nabla \cdot \mathbf{H}
\end{equation*}
for a vector field $\mathbf{H} \in L^p(B_r)$ satisfying 
\begin{equation*} \label{}
\left\| \mathbf{H} \right\|_{\underline{L}^{(2+\delta)\wedge p}(B_r)} 
\leq 
C\left\| h \right\|_{\underline{W}^{-1,(2+\delta)\wedge p}(B_r)}.
\end{equation*}
If this is the case, then the result follows by  Corollary~\ref{c.Meyers reverse} and Lemma~\ref{l.Meyers really}, applied with $f = \left|\nabla u\right|$, $g = \left|\mathbf{\mathbf{H}}\right|$, and $q = \tfrac{2d}{d+2}$. 

\smallskip

The general case follows from the observation that the special case described above is actually general. To see this, we first notice that it suffices to consider the case that $p\in \left[ \frac32,3\right]$, since we may assume that~$\delta \leq 1$. Then, given $h\in W^{-1,p}(U)$, we solve the Dirichlet problem
\begin{equation} 
\label{e.DPLippy}
\left\{
\begin{aligned}
& - \Delta w = h & \mbox{in} & \ U, \\
& w = 0 & \mbox{on} & \ \partial U. 
\end{aligned}
\right.
\end{equation}
According to Proposition~\ref{p.DPpoisson}, there exists a solution $w$ belonging to $W^{1,p}(U)$ with the estimate
\begin{equation*} \label{}
\left\| w \right\|_{\underline{W}^{1,p}(B_r)} 
\leq 
C \left\| h \right\|_{\underline{W}^{-1,p}(B_r)}. 
\end{equation*}
(We leave it to the reader to check that we have scaled this estimate correctly.) We then set $\mathbf{H}:= -\nabla w$, which completes the proof of the claim and of the theorem. 
\end{proof}

We state and prove a global version of the Meyers estimate in bounded Lipschitz domains with Dirichlet boundary conditions. It is useful to define, for a given bounded Lipschitz domain~$U\subseteq \Rd$, the geometric constant
\begin{equation} \label{e.kappa_U appendix}
\kappa_U := 
\sup_{z \in U, \, r>0} \frac{\left| B_r(z) \cap U \right|}{\left| B_{r/2}(z) \cap U \right|} + \sup_{z \in \partial U,\, r>0} \frac{\left| B_r \right|}{\left| B_{r}(z) \setminus U \right|}.
\end{equation}
The assumption that $U$ is Lipschitz ensures that $\kappa_U<\infty$. It is the finiteness of $\kappa_U$ and the solvability of~\eqref{e.DPLippy} which are the main reasons for assuming~$U$ is Lipschitz.

\begin{theorem}[Global Meyers estimate]
\label{t.Meyers appendix global}
\index{Meyers estimate}
Fix $p\in (2,\infty)$ and let~$U\subseteq \Rd$ be a bounded Lipschitz domain. Suppose that $h \in W^{-1,p}(U)$, $f \in W^{1,p}(U)$, and that $u \in f + H_0^1(U)$ is the solution of 
\begin{equation} 
\label{e.Meyers eq global}
\left\{ 
\begin{aligned} 
& -\nabla \cdot \left( \a(x) \nabla u\right) = h & \mbox{in} & \ U, \\
& u = f & \mbox{on} & \ \partial U. 
\end{aligned}
\right.
\end{equation}
There exist~$\delta(U,d,\Lambda)  >0$ and $C(U,d,\Lambda)<\infty$ such that $u \in W^{1,(2+\delta)\wedge p}(U)$ and 
\begin{equation} \label{e.Meyers appendix global}
\left\| \nabla u \right\|_{\underline{L}^{(2+\delta) \wedge p}(U)} \leq C \left( \left\| \nabla f \right\|_{\underline{L}^{(2+\delta)\wedge p}(U)}  +  \left\| h\right\|_{\underline{W}^{-1,(2+\delta)\wedge p}(U)}  \right).
\end{equation}
 \end{theorem}
 
Before embarking on the proof of Theorem~\ref{t.Meyers appendix global}, we mention one simple reduction, which is that we may assume that $f=0$. Indeed, by defining 
\begin{equation*} \label{}
\tilde{u}:= u-f
\quad \mbox{and} \quad 
\tilde{h}:=h + \nabla\cdot\left( \a(x) \nabla f(x)\right),
\end{equation*}
we see that $\tilde{u} \in H_0^1(U)$ and $\tilde{u}$  solves the equation
\begin{equation*} 
-\nabla \cdot \left( \a(x) \nabla \tilde{u} \right) = \tilde{h} \quad \mbox{in} \ U. 
\end{equation*}
We then observe that the conclusion of the theorem for $u$ follows from the one for~$\tilde{u}$. Therefore, we may, without loss of generality, consider only the case $f = 0$.

\begin{lemma}[Global Caccioppoli estimate] 
\label{l.Caccioppoli appendix glob}
\index{Caccioppoli inequality}
Let $U$ be a bounded Lipschitz domain. Suppose that $h \in H^{-1}(U)$ and that $u \in H_0^1(U)$ solves 
\begin{equation} 
\label{e.Caccioppoli appendix glob.pde}
-\nabla \cdot \left( \a(x) \nabla u\right) = h
\quad \mbox{in} \ U.
\end{equation}
There exists~$C(U,d,\Lambda) < \infty$ such that
\begin{equation} 
\label{e.Caccioppoli appendix glob 0}
\left\|\nabla u \right\|_{\underline{L}^2(U)}
\leq
C \left\| h \right\|_{\underline{H}^{-1}(U)},
\end{equation}
and, for every $z \in \partial U$ and $r>0$,
\begin{equation} 
\label{e.Caccioppoli appendix glob}
\left\| \nabla u \right\|_{\underline{L}^2(B_{r/2}(z)\cap U)}
\leq 
C
\left( \frac1r
\left\| u \right\|_{\underline{L}^2(B_r(z) \cap U)}
+
\left\| h \right\|_{\underline{H}^{-1}(B_r(z) \cap U)}
\right).
\end{equation}
\end{lemma}
\begin{proof}
The estimate~\eqref{e.Caccioppoli appendix glob 0} follows simply by testing the equation~\eqref{e.Caccioppoli appendix glob.pde} with~$u$.  The proof of~\eqref{e.Caccioppoli appendix glob} is similar to the one of Lemma~\ref{l.Caccioppoli appendix}. Note that we have control on the volume factor~$|B_{r}(z) \cap U|/|B_{r/2}(z) \cap U|$ by~$\kappa_U < \infty$ since~$U$ is Lipschitz.  
\end{proof}

The global version of the Caccioppoli inequality naturally leads to a global version of the reverse H\"older inequality. 

\begin{corollary} 
\label{c.Meyers reverse glob}
Assume the hypothesis of Lemma~\ref{l.Caccioppoli appendix glob} and consider $u$ to belong to the space~$H^1(\Rd)$ by extending it to be zero in~$\Rd\setminus U$. 
Then there exists~$C(U,d,\Lambda) < \infty$ such that, for every $z \in U$ and $r>0$,
\begin{equation} 
\label{e.rev holder appendix glob}
\left\| \nabla u\right\|_{\underline{L}^2 (B_{r/2}(z)) }
\leq 
C
 \left(
\left\| \nabla u\right\|_{\underline{L}^{2_\ast} (B_{r}(z) \cap U)}
+
\left\| h \right\|_{\underline{H}^{-1}(B_r(z)\cap U)}
\right).
\end{equation}
\end{corollary}
\begin{proof}
If $B_r(z) \subseteq U$, then the result follows from the interior estimate given in Corollary~\ref{c.Meyers reverse}. If $B_r(z) \not\subseteq U$, then we can find a point $z' \in \partial U \cap B_r(z)$, which implies that $B_r(z) \subseteq B_{2r}(z')$. We may then apply the result of Lemma~\ref{l.Caccioppoli appendix glob} in the ball $B_{4r}(z')$ and follow the rest of the proof of Corollary~\ref{c.Meyers reverse glob}. The only difference is that the version of the Sobolev-Poincar\'e inequality we use is different: rather than the version for mean-zero functions, we use the one which applies to functions vanishing on a set of positive measure (cf.~\cite[Lemma 4.8]{HL}). It is valid here because~$u$ vanishes on $B_{4r}(z') \cap U$, which has measure at least $\kappa_U^{-1} \left|B_{4r} \right| \geq c\left|B_{4r} \right|$. 
\end{proof}

\begin{proof}[{Proof of Theorem~\ref{t.Meyers appendix global}}]
The result is a consequence of~\eqref{e.Caccioppoli appendix glob 0}, Corollary~\ref{c.Meyers reverse glob} and Lemma~\ref{l.Meyers really}. The details are almost identical to those of the proof of Theorem~\ref{t.Meyers appendix}, so we leave them to the reader. 
\end{proof}

\index{Meyers estimate|)}

\chapter{Sobolev norms and heat flow}
\label{a.MSP}

In this appendix, we review the characterization of Sobolev norms in terms of weighted integral expressions involving spatial averages of the function with respect to the standard heat kernel. These inequalities are used in Chapters~\ref{c.A1} and~\ref{c.gff} in order to obtain estimates in~$W^{-\alpha,p}$ for the gradients and fluxes of the correctors as well as Gaussian free fields, and can be compared to Propositions~\ref{p.mspoin} and~\ref{l.mspoincare.masks}. 

\smallskip

It turns out that we need to distinguish two cases, depending on whether~$\alpha$ is an integer or not. We first consider the case that $\alpha \in\N$ in the following proposition. The fractional case in which~$\alpha\not\in\N$ is presented in Proposition~\ref{p.MSP.Walpha}, below.

\begin{proposition} \label{p.MSP_int}
Let $\alpha \in \N$ and $p \in (1,\infty)$. There exists~$C (\alpha,p,d) < \infty$ such that for every $u \in W^{-\alpha,p}(\R^d)$,
\begin{equation}
\label{e.triebel}
\frac1C \left\| u \right\|_{W^{-\alpha,p}(\R^d)} 
\leq 
\left( \int_{\Rd}   
\left( \int_0^1 \left( t^{\frac \alpha2} \left| \left( \Phi\left(t,\cdot\right)  \ast u\right)(x) \right| \right)^2 \, \frac{dt}{t} \right)^{\frac p2}   dx
\right)^{\frac1p}  
\leq C \left\| u \right\|_{W^{-\alpha,p}(\R^d)} .
\end{equation}
\end{proposition}

Inequalities similar to~\eqref{e.triebel} are often presented as consequences of singular integrals, multiplier theory, Littlewood-Paley decompositions and functional calculus. We give here a more direct proof based on the following property of the heat equation.

\begin{lemma} 
\label{l.MSP_pre}
Let $q \in (1,\infty)$ and $k\in \N$. There exists~$C(k,q,d)<\infty$ such that for every $f \in L^q(\R^d)$,
\begin{equation} \label{e.MSP_pre}
 \left( \int_{\Rd} \left(   \int_0^\infty t^k \left| \nabla^k 
 \left( \Phi\left(t,\cdot\right)  \ast f\right)(x) \right|^2 \, \frac{dt}{t} \right)^{\frac q2} \,dx \right)^{\frac1q} 
 \leq C \left\| f \right\|_{L^q(\R^d)}.
\end{equation}
\end{lemma}

\begin{proof} 
By density, it is enough to show the estimate for compactly supported smooth function~$f$. In Steps 2 and 3, we consider the case $q \in (1,2]$, in Step 4 the case $q \in [4,\infty)$, and finally in Step 5 we use interpolation to obtain the result for $q \in (2,4)$. We denote 
\begin{equation*} \label{}
v(t,x):=  \left( \Phi\left(t,\cdot\right)  \ast f\right)(x).
\end{equation*}

\smallskip

\emph{Step 1.} We first prove a pointwise bound for certain maximal functions appearing below. We will prove that, for $k \in \N$ and a constant~$C(k,d)<\infty$,
\begin{equation} \label{e.maxfct heat vs HL}
\sup_{t > 0} \left| \left( g_k\left(t,\cdot\right)  \ast | f | \right)(x) \right| \leq C_k \mathsf{M}f(x).
\end{equation}
Here $\mathsf{M}f$ is the centered Hardy-Littlewood maximal function defined by
\begin{equation*} \label{}
\mathsf{M}f(x):= \sup_{r>0} \fint_{B_r(x)} \left| f(y) \right|\,dy. 
\end{equation*}
and
\begin{equation} \label{e.MSP.g_k}
g_k(t,x) := \left(1 + \frac{|x|^2}{t} \right)^{k-1} \Phi(t,x).
\end{equation}
To prove~\eqref{e.maxfct heat vs HL}, choose~$t_x$ realizing the supremum and estimate as
\begin{align*} 
\sup_{t > 0} \left| \left( g\left(t,\cdot\right)  \ast | f | \right)(x) \right|  & \leq 2 \left|  \left( g(t_x,\cdot) \ast f \right)(x)  \right|  \\ 
 & \leq 2 \int_{\R^d}  \left(1 + \frac{|x-z|^2}{t_x} \right)^{k-1}  \Phi(t_x,x-z) \left| f(z)\right| \, dz 
\\ & \leq C \sum_{j=0}^\infty 4^{ (k-1) j d} \exp(-4^{j-2}) \fint_{B_{2^j \sqrt{t_x}}(x)} \left| f(z)\right|  \,dz 
\\ & \leq C_k \mathsf{M}f(x).
\end{align*}
Therefore~\eqref{e.maxfct heat vs HL} follows.

\smallskip

\emph{Step 2.} 
We now assume that $q \in (1,2]$. 
Compute formally
\begin{equation*} 
\partial_t |v|^q = q |v|^{q-2} v \partial_t v \quad \mbox{and} \quad  \Delta |v|^q = q(q-1) |v|^{q-2} \left|\nabla v\right|^2 +   q |v|^{q-2} v \Delta v. 
\end{equation*}
To give a rigorous meaning for negative powers of $|v|$ above, and in what follows, one needs to take $(|v|+\ep)^{q-2} v$ instead of $|v|^{q-2} v$, and then pass $\ep \to 0$ as soon as it is appropriate. Rearranging the terms and multiplying with $|v|^{2-q}$ leads to
\begin{equation*} 
\left|\nabla v\right|^2 = \frac{|v|^{2-q}}{q-1} \left( \frac 1q \left( \Delta |v|^q - \partial_t |v|^q \right) - |v|^{q-2} v \left( \Delta v - \partial_t v \right) \right).
\end{equation*}
In particular, since $v$ solves the heat equation, we have
\begin{equation} \label{e.MSP-grad v squared}
\left|\nabla v\right|^2 = \frac{|v|^{2-q}}{q(q-1)} \left( \Delta |v|^q - \partial_t |v|^q \right),
\end{equation}
and thus the term on the right is nonnegative.

\smallskip

\emph{Step 3.}  We prove~\eqref{e.MSP_pre} in the case $q\in (1,2]$. Fix $k \in \N \cup \{0\}$. 
Notice that, by integration by parts,
\begin{equation*} 
t^{\frac{k}{2}} |\nabla^{k+1}  v(t,x)| \leq C\left(  t^{\frac{k}{2}}  \left|\nabla^{k} \Phi \left(\tfrac t2,\cdot\right)\right| \ast \left| \nabla v \left(\tfrac t2,\cdot\right) \right| \right)(x) . 
\end{equation*}
Therefore, by~\eqref{e.MSP-grad v squared}, we have that
\begin{multline*} 
t^{\frac{k}{2}} \left|\nabla^{k+1}  v(t,x) \right| \\ \leq C\left(  t^{\frac k2}  \left|\nabla^k \Phi \left(\tfrac t2,\cdot\right)\right| \ast \left(  \left| v\left(\tfrac t2,\cdot\right) \right|^{2-q}\left( \Delta \left|v\left(\tfrac t2,\cdot\right)\right|^q - \partial_t \left|v\left(\tfrac t2,\cdot\right)\right|^q \right) \right)^{\frac12} \right)(x).
\end{multline*}
Applying H\"older's inequality yields 
\begin{multline*} 
\left(  t^{\frac k2}  \left|\nabla^k \Phi \left(\cdot,t\right)\right| \ast \left(  \left| v\left(t,\cdot\right) \right|^{2-q}\left( \Delta \left|v\left(t,\cdot\right)\right|^q - \partial_t \left|v\left(t,\cdot\right)\right|^q \right) \right)^{\frac12} \right)(x)
\\ \leq C \left( \left( g_k \left(t,\cdot\right) \ast  \left| v\left(t,\cdot\right) \right|^{2-q} \right)^{\frac12} \left( \Phi\left(t,\cdot\right) \ast \left(  \Delta \left|v\left(t,\cdot\right)\right|^q - \partial_t \left|v\left(t,\cdot\right)\right|^q \right)  \right) ^{\frac 12} \right) (x),
\end{multline*}
where $g_k$ is defined in~\eqref{e.MSP.g_k}. By Step 2, 
\begin{equation*} 
\sup_{t>0} \left( g_k \left(t,\cdot\right)  \ast  \left| v\left(t,\cdot\right) \right|^{2-q} \right)(x) \leq C\left( \mathsf{M}f(x) \right)^{2-q},
\end{equation*}
and combining the displays above gives the pointwise bound
\begin{equation*} 
t^{k} \left|\nabla^{k+1}  v(t,x) \right|^2 \leq C \left(  \mathsf{M}f(x) \right)^{2-q} \left( \Phi\left(t,\cdot\right) \ast \left( \Delta \left|v\left(t,\cdot\right)\right|^q - \partial_t \left|v\left(t,\cdot\right)\right|^q \right)   \right) (x).
\end{equation*}
 It follows, again by H\"older's inequality, that
\begin{multline*} 
\int_{\R^d} \left( \int_0^\infty  t^{k} \left|\nabla^{k+1}  v(t,x) \right|^2  \, dt \right)^{\frac q2} \,dx  \\ 
\leq C \left\| \mathsf{M}f \right\|_{L^q(\R^d)}^{\frac{q}{2}(2-q)} \left( \int_{\R^d} \int_{0}^\infty \Phi\left(t,\cdot\right) \ast \left( \Delta \left|v\left(t,\cdot\right)\right|^q - \partial_t \left|v\left(t,\cdot\right)\right|^q \right)(x) \, dt \, dx    \right)^{\frac q2}.
\end{multline*}
But since $\Delta \left|v\left(t,\cdot\right)\right|^q - \partial_t \left|v\left(t,\cdot\right)\right|^q \geq 0$ and the convolution contracts, we have that  
\begin{multline*} 
 \int_{\R^d} \int_{0}^\infty \Phi\left(t,\cdot\right) \ast \left( \Delta \left|v\left(t,\cdot\right)\right|^q - \partial_t \left|v\left(t,\cdot\right)\right|^q \right)(x) \, dt \, dx
 \\ \leq  \int_{\R^d} \int_{0}^\infty \left( \Delta \left|v\left(t,x\right)\right|^q - \partial_t \left|v\left(t,x\right)\right|^q \right) \, dx \, dt \leq 2 
 \left\| f \right\|_{L^q(\R^d)}^{q} .
\end{multline*}
Consequently, 
\begin{equation*} 
\int_{\R^d} \left( \int_0^\infty  t^{k} \left|\nabla^{k+1}  v(t,x) \right|^2  \, dt \right)^{\frac q2} \,dx  \leq C \left\| \mathsf{M}f \right\|_{L^q(\R^d)}^{\frac{q}{2}(2-q)}  \left\| f \right\|_{L^q(\R^d)}^{\frac{q^2}{2}} ,
\end{equation*}
and the result follows by the strong type $(q,q)$-estimate for maximal functions, see e.g. Lemma~\ref{l.maximalfct}. This proves~\eqref{e.MSP_pre} when $q\in (1,2]$.

\smallskip

\emph{Step 4.} We now proceed in proving~\eqref{e.MSP_pre} when $q\in [4,\infty)$. By duality, it is enough to show that, for $k \in \N \cup \{0\}$,
\begin{equation*} 
\sup_{g \in L^p(\R^d), \left\| g \right\|_{L^{p}(\R^d)}\leq 1 }\int_{\R^d} \int_0^\infty t^k \left| \nabla^{k+1} v(t,x) \right|^2 g(x) \,dt \, dx \leq C \left\| f \right\|_{L^q(\R^d)}^2 ,
\end{equation*}
where $p$ is the conjugate of $\tfrac{q}{2}$, that is $p = \tfrac{q}{q-2}$. Notice that $p \in (1,2]$ if and only if $q \in [4,\infty)$. Pick nonnegative $g \in L^p(\R^d)$ such that $\left\| g \right\|_{L^{p}(\R^d)}\leq 1$. We first observe that
\begin{equation*} 
\partial_t  \left| \nabla^{k+1} v \right|^2  -\Delta  \left| \nabla^{k+1} v \right|^2 = - \left| \nabla  \left| \nabla^{k+1} v \right| \right|^2,
\end{equation*}
since each component of $\nabla^{k+1} v$ satisfies the heat equation. In particular, $\left| \nabla^{k+1} v \right|^2$ is a smooth subsolution. By the comparison principle we thus have that 
\begin{equation*} 
\left| \nabla^{k+1} v(t,x) \right|^2 \leq \left( \Phi \left(\tfrac t2,\cdot \right) \ast \left| \nabla^{k+1} v\left(\tfrac t2,\cdot \right) \right|^2 \right) (x),
\end{equation*}
since the latter term is a solution to the heat equation in $(t/2,\infty) \times \R^d$ with initial values $\left| \nabla^{k+1} v\left(\tfrac t2,\cdot \right) \right|^2$. 
Denote $u(t,x) = \left(\Phi \left(t,\cdot\right) \ast g\right) (x)$. By Fubini it is thus enough to show that 
\begin{equation*} 
\int_{\R^d} \int_0^\infty t^k \left| \nabla^{k+1} v(t,x) \right|^2 u(t,x) \,dt \, dx \leq C \left\| f \right\|_{L^q(\R^d)}^2 \left\| g \right\|_{L^{p}(\R^d)} .
\end{equation*}
Compute now, for any solution $w$ to the heat equation, 
\begin{align*} 
\left| \nabla w \right|^2 u & = \nabla \cdot \left( w \nabla w u \right) - \left( w \Delta w \right) u - w \nabla w \cdot \nabla u 
\\ \notag & = \nabla \cdot \left( w \nabla w u \right) - \frac12 \left(\partial_t w^2\right) u - w \nabla w \cdot \nabla u 
\\ \notag & = \nabla \cdot \left( w \nabla w u \right) - \frac12 \partial_t \left( w^2 u \right) + \frac12 \left(\partial_t v\right) w^2 - w \nabla w \cdot \nabla u 
\\ \notag & = \nabla \cdot \left( w \nabla w u  +\frac12 w^2 \nabla u \right) - \frac12 \partial_t \left( w^2 u \right)  - 2 w \nabla w \cdot \nabla u . 
\end{align*}
Therefore,
\begin{align*} 
t^k \left| \nabla w \right|^2 u & = t^k \nabla \cdot \left( w \nabla w u  +\frac12 w^2 \nabla u \right)  - \frac12 \partial_t \left( t^k w^2 u \right) - 2 t^k  w \nabla w \cdot \nabla u + \frac{k}{2} t^{k-1} w^2 u  .
\end{align*}
Then, by integration by parts, 
\begin{align*} 
\lefteqn{ \int_{\R^d} \int_0^\infty  t^k \left| \nabla w(t,x) \right|^2 u(t,x) \,dx \, dt } \qquad & \\ 
\notag & \leq \frac{k}{2} \int_{\R^d} \int_0^\infty t^{k-1} \left|w(t,x)\right|^2 u(t,x) \, dx \, dt
\\ \notag & \quad + \frac12 \lim_{t \to 0} \int_{\R^d} t^k |w(t,x)|^2 |u(t,x)| \, dx 
\\ \notag & \quad + 2 \int_{\R^d} \int_0^\infty t^k |w(t,x)| | \nabla w(t,x)| |\nabla u(t,x)| \, dx \, dt.
\end{align*}
We will apply the above inequality with components of $w :=  \nabla^{k} v $.
By H\"older's inequality we get, 
\begin{equation*} 
 \lim_{t \to 0} \int_{\R^d} t^k |w(t,x)|^2 |u(t,x)| \,dx   \leq \int_{\R^d} \left| f \right|^2 |g| \leq \left\|  f \right\|_{L^q(\R^d)}^2 \left\| g \right\|_{L^p(\R^d)}.
\end{equation*}
Notice that by the assumed smoothness of $f$, the term on the left is zero for $k>0$. By H\"older's inequality we also obtain
\begin{multline*} 
\int_{\R^d} \int_0^\infty t^k |w(t,x)| | \nabla w(t,x)| |\nabla u(t,x)| \, dx \, dt  \leq C \left\| \sup_{t>0} \left( t^{\frac k2}  |\nabla^k v\left(t,\cdot\right) | \right) \right\|_{L^q(\R^d)} 
\\ \times \left\| \left( \int_0^\infty t^k | \nabla^{k+1} v \left(t,\cdot\right) |^2 \, dt  \right)^{\frac12} \right\|_{L^q(\R^d)} \left\| 
\left( \int_0^\infty | \nabla u\left(t,\cdot\right) |^2 \, dt  \right)^{\frac12}  \right\|_{L^p(\R^d)}.
\end{multline*}
Furthermore, by Steps 1 and 3, and the strong $(q,q)$-type estimate, we obtain
\begin{equation*} 
\left\|  \sup_{t>0} \left( t^{\frac k2}  |\nabla^k v \left(t,\cdot\right) |  \right)  \right\|_{L^q(\R^d)}  \leq C \left\| f \right\|_{L^q(\R^d)} 
\end{equation*}
and
\begin{equation*} 
\left\|  \left( \int_0^\infty | \nabla u\left(t,\cdot\right) |^2 \, dt  \right)^{\frac12}  \right\|_{L^p(\R^d)} \leq C \left\| g \right\|_{L^p(\R^d)}.
\end{equation*}
We deduce that 
\begin{align} \label{e.MSP.q>4}
\lefteqn{\int_{\R^d} \int_0^\infty  t^k \left| \nabla^{k+1} v (t,x) \right|^2 u(t,x) \,dx \, dt} \qquad &  \\ 
\notag & \leq \frac{k}{2} \int_{\R^d} \int_0^\infty  t^{k-1} \left| \nabla^{k} v(t,x) \right|^2 u(t,x) \,dx \, dt \\
\notag & \quad + C  \left\|  f \right\|_{L^q(\R^d)}^2 \left\| g \right\|_{L^p(\R^d)}\\
\notag & \quad +C \left\| f \right\|_{L^q(\R^d)}  \left\| g \right\|_{L^p(\R^d)} \left\| \left( \int_0^\infty t^k | \nabla^{k+1} v\left(t,\cdot\right) |^2 \, dt  \right)^{\frac12} \right\|_{L^q(\R^d)}.
\end{align}
Assume now inductively that, for nonnegative $g \in L^p(\R^d)$ with $\| g \|_{L^p(\R^d)} \leq 1$ and $m \in \{0,\ldots,k\}$,
\begin{equation*} 
\frac{m}{2} \int_{\R^d} \int_0^\infty  t^{m-1} \left| \nabla^m v(t,x) \right|^2 u(t,x) \,dx \, dt \leq C_m \left\| f \right\|_{L^q(\R^d)}^2 ,
\end{equation*}
where we set $C_0 = 0$ for $m=0$. Taking supremum over nonnegative $g \in L^p(\R^d)$ with $\| g \|_{L^p(\R^d)} \leq 1$ and recalling that 
\begin{equation*} 
 \int_{\R^d}  \int_0^\infty  t^k \left| \nabla^{k+1} v(t,x) \right|^2  g(x) \,dt \, dx  \leq 2^k \int_{\R^d} \int_0^\infty  t^k \left| \nabla^{k+1} v(t,x) \right|^2 u(t,x) \,dx \, dt ,
\end{equation*}
we obtain by the induction assumption and~\eqref{e.MSP.q>4}  that
\begin{multline*} 
\left\| \left( \int_0^\infty t^k \left| \nabla^{k+1} v\left(t,\cdot\right)  \right|^2 \, dt  \right)^{\frac12} \right\|_{L^q(\R^d)}^2 
\\ \leq C \left\|  f \right\|_{L^q(\R^d)} \left( \left\|  f \right\|_{L^q(\R^d)} + \left\| \left( \int_0^\infty t^k \left| \nabla^{k+1} v\left(t,\cdot\right)  \right|^2 \, dt  \right)^{\frac12} \right\|_{L^q(\R^d)}\right).
\end{multline*}
Hence, after applying Young's inequality and reabsorbing terms,  
\begin{equation*} 
\left\| \left( \int_0^\infty t^k \left| \nabla^{k+1} v\left(t,\cdot\right)  \right|^2  \right)^{\frac12} \right\|_{L^q(\R^d)} \leq C  \left\|  f \right\|_{L^q(\R^d)} .
\end{equation*}
This proves the estimate~\eqref{e.MSP_pre} for $q \in [4,\infty)$.

\smallskip

\emph{Step 5.} Define, for $f \in L^q(\R^d)$, the operator
\begin{equation*} 
Tf(x) := \left( \int_0^\infty  t^k \left| \nabla^{k} \left( \Phi\left(t,\cdot\right)  \ast f(\cdot) \right)(x) \right|^2  \, \frac{dt}{t} \right)^{\frac 12},
\end{equation*}
which is clearly quasilinear. By the previous steps we have that 
\begin{equation*} 
\left\| Tf\right\|_{L^{3/2}(\R^d)} \leq C \left\| f\right\|_{L^{3/2}(\R^d)} \quad \mbox{and} \quad \left\| Tf\right\|_{L^{5}(\R^d)} \leq C \left\| f\right\|_{L^{5}(\R^d)}.
\end{equation*}
By the Marcinkiewicz interpolation theorem (see \cite[Appendix~D]{Taylor}), we deduce that, for every $q \in [2,4]$, 
\begin{equation*} 
\left\| Tf \right\|_{L^{q}(\R^d)} \leq C \left\| f\right\|_{L^{q}(\R^d)} .
\end{equation*}
The proof is complete.
\end{proof}

Above we made use of the following continuity property of maximal functions. 
\begin{lemma} \label{l.maximalfct}
Let $p \in (1,\infty]$ and $f \in L^p(\R^d)$. Let
\begin{equation*} \label{e.maximalfct}
\mathsf{M}f(x):= \sup_{r>0} \fint_{B_r(x)} \left| f(y) \right|\,dy. 
\end{equation*}
Then we have so-called strong $(p,p)$-type estimate
\begin{equation} \label{e.strong-(p,p)}
\left\| \mathsf{M}f \right\|_{L^p(\R^d)} \leq  2 \left(\frac{ 5^d p}{p-1}\right)^{\frac 1p}  \left\| f \right\|_{L^p(\R^d)}.
\end{equation}
\end{lemma}

\begin{proof} The statement is clear if $p=\infty$. Assume thus that $p \in (1,\infty)$. It suffices to show that, for $\lambda>0$, 
\begin{equation} \label{e.pre-weak-(1,1)}
\left| \left\{ \mathsf{M} f > 2 \lambda \right\} \right|  \leq  \frac{5^d}{\lambda} \int_{\left\{ |f| >  \lambda \right\} } |f(y)|\,dy.
\end{equation}
Indeed, this can be seen immediately from the formulas
\begin{equation} \label{e.Mf in L^p}
\left\| \mathsf{M} f \right\|_{L^p(\R^d)}^p = 2^{p} p \int_0^\infty \lambda^{p-1} \left| \left\{ \mathsf{M} f  > 2 \lambda \right\} \right| \, d\lambda
\end{equation}
and
\begin{equation} \label{e.f in L^p}
\left\| f \right\|_{L^p(\R^d)}^p = (p-1) \int_0^\infty \lambda^{p-2} \int_{\left\{ |f| > \lambda \right\} } |f(y)|\,dy \, d\lambda.
\end{equation}
To prove~\eqref{e.pre-weak-(1,1)}, we fix $\lambda>0$ and observe that, for every~$x \in \left\{ \mathsf{M} f  > 2 \lambda \right\}$, there exists $r_x>0$ such that 
\begin{equation*} 
2\lambda <   \fint_{B_{r_x}(x)} |f(y)| \, dy \leq \mathsf{M} f(x). 
\end{equation*}
From this it follows that 
\begin{equation} \label{e.pre-weak-(1,1)-pre}
\left| B_{r_x}(x) \right| \leq \frac{1}{\lambda} \int_{\left\{ |f| >  \lambda \right\} \cap B_{r_x}(x) } |f(y)|\,dy.
\end{equation}
Observe that $\sup_{x \in \left\{ \mathsf{M} f  > 2 \lambda \right\}} r_x < \infty$, since, by H\"older's inequality,
\begin{equation*} 
2\lambda < \fint_{B_{r_x}(x)} |f(y)| \, dy \leq \left|B_{r_x}\right|^{-\frac 1p} \| f\|_{L^p(\R^d)}  \quad \implies \quad r_x \leq C \lambda^{-\frac pd}. 
\end{equation*}
Thus, by the Vitali covering lemma~\ref{l.vitali}, there exists a countable set $\left\{x_j\right\}_j$ such that the family $\{B_{r_{x_j}}(x_j)\}_j$ is pairwise disjoint and
\begin{equation*} 
\bigcup_{x \in \left\{ \mathsf{M} f > 2 \lambda \right\}} B_{r_x}(x) \subset \bigcup_{j} B_{5r_{x_j}}(x_j).
\end{equation*}
Using the fact the balls are disjoint, we obtain by~\eqref{e.pre-weak-(1,1)-pre} that
\begin{equation*} 
\left| \left\{ \mathsf{M} f > 2\lambda \right\}\right| \leq \left|\bigcup_{x \in \left\{ \mathsf{M} f > 2 \lambda \right\}} B_{r_x}(x) \right| \leq 5^d \sum_j \left| B_{r_{x_j}}(x_j)\right| \leq
 \frac{5^d}{\lambda} \int_{\left\{ |f| >  \lambda \right\} } |f(y)|\,dy.
\end{equation*}
This completes the proof of~\eqref{e.pre-weak-(1,1)} and thus the lemma. 
\end{proof}

\begin{proof}[Proof of Proposition~\ref{p.MSP_int}]
Throughout we fix $\alpha \in \N$ and $q \in (1,\infty)$ to be $q := \tfrac{p}{p-1}$.

\smallskip

\emph{Step 1.} One direction is easy using Lemma~\ref{l.MSP_pre}. By Remark~\ref{r.RRT}, every distribution $u$ in $W^{-\alpha,p}(\R^d)$ can be written as 
\begin{equation*} 
u = \sum_{0\leq |\beta| \leq \alpha} \partial^\beta f_\beta , \quad f_\beta \in L^p(\R^d).
\end{equation*}
Here we denote $\partial^\beta = \partial_{x_1}^{\beta_1} \ldots \partial_{x_d}^{\beta_d}$ with $\sum_{j=1}^d \beta_j = |\beta|$. 
Then 
\begin{equation*} 
 \left| \Phi\left(t,\cdot\right)  \ast u \right|  \leq \sum_{0\leq |\beta| \leq \alpha}  \left| \partial^\beta \Phi\left(t,\cdot\right)  \ast f_\beta \right| , 
\end{equation*}
and hence Lemma~\ref{l.MSP_pre} implies that 
\begin{equation*} 
\left\|  \left( \int_0^1 \left( t^{\frac \alpha2} \left| \Phi\left(t,\cdot\right)  \ast u \right| \right)^2 \, \frac{dt}{t}  \right)^{\frac 12}   \right\|_{L^p(\R^d)}  \leq C 
 \sum_{0\leq |\beta| \leq \alpha}  \left\| f_\beta \right\|_{L^p(\R^d)} \leq C \left\| u \right\|_{W^{-\alpha,p}(\R^d)}.
\end{equation*}

\smallskip

\emph{Step 2.} We then prove the other direction. For this, select $g \in W^{\alpha,q}(\Rd)$. We study the solution $w$ of the Cauchy problem for the heat equation:
\begin{equation*} 
\left\{ 
\begin{aligned}
& \partial_t w - \Delta w = 0 & \mbox{in}  & \ \R^d \times (0,\infty), \\ 
& w\left(0,\cdot\right)  = g & \mbox{in} & \ \Rd. 
\end{aligned}
\right.
\end{equation*}
Note that the solution is given by convolution against the heat kernel:
\begin{equation*} 
w(t,x) = (g \ast \Phi\left(t,\cdot\right))(x). 
\end{equation*}
The goal is to estimate $\left| \rho(0) \right|$, where  we define
\begin{equation*} \label{}
\rho(t):= \int_{\Rd} u(x) w(t,x)\,dx. 
\end{equation*}
By Taylor's theorem we first get that 
\begin{equation} \label{e.rho(0) stupid}
\rho(0) = \sum_{j=0}^{k-1} \frac{(-1)^j}{j!} \rho^{(j)}(1) +  \frac{(-1)^{k}}{(k-1)!} \int_0^1 t^{k-1} \rho^{(k)}(t) \, dt ,
\end{equation}
where $\rho^{(j)}$ denotes $j^{th}$ derivative of $\rho$.  Computing then the $j^{th}$ derivative gives 
\begin{align} 
\label{e.rho(j)}
\rho^{(j)}(t)
& 
= \int_{\Rd} u(x) \, \partial_t^j w(t,x)\,dx
\\ \notag &
 = \int_{\Rd} u(x)  \!\left( g \ast \partial_t \Delta^{j-1}\Phi\left(t,\cdot\right)  \right) (x) \,dx
\\ \notag  & 
= \int_{\R^d} u(x) \partial_t \left( \Delta^{j-1} g \ast \Phi\left(t,\cdot\right) \right)(x) \,dx
\\ \notag  & 
= \int_{\R^d} \left( u \ast \Phi\left(\tfrac t2,\cdot\right) \right)(x) \partial_t \left( \Delta^{j-1} g \ast \Phi\left(\tfrac t2,\cdot\right) \right)(x) \,dx.
\end{align}
Choosing $ k \in \N$ and $\beta \in \{ 0,1\}$ such that $\alpha = 2k - \beta$, we have 
\begin{align} \label{e.taylorresidualrho_int}
\lefteqn{ \left| \int_0^1 t^{k-\beta} \rho^{(k+1-\beta)}(t) \, dt  \right| } \qquad & \\ 
& \notag = \left|  \int_{\R^d} \int_0^1 t^{k-\beta}  \left( u \ast \Phi\left(\tfrac t2,\cdot\right) \right)(x) \partial_t \left( \Delta^{k-\beta} g \ast \Phi\left(\tfrac t2,\cdot\right) \right)(x)  \, dt \,dx \right| 
\\ & \notag \leq \int_{\R^d} \left( \int_0^1 \left( t^{\frac{\alpha}{2}} \left| \left(u\ast \Phi\left(\tfrac t2,\cdot\right)\right)(x) \right| \right)^2 \, \frac{dt}{t} \right)^{\frac12} 
\\ &\notag  \qquad \times  \left( \int_0^1 \left( t^{1-\frac{\beta}{2}} \left|\left( \Delta^{k-\beta} g \ast \partial_t  \Phi\left(\tfrac t2,\cdot\right) \right) (x) \right| \right)^2 \, \frac{dt}{t} \right)^{\frac12} \, dx
\\ \notag & \leq 
 \left\| \left( \int_0^1 \left( t^{\frac{\alpha}{2}} \left| \left( u \ast   \Phi\left(t,\cdot\right) \right) (\cdot) \right| \right)^2 \, \frac{dt}{t} \right)^{\frac12}  \right\|_{L^p(\Rd)}
\\ \notag & \qquad \times \left\| \left( \int_0^1 \left( t^{1-\frac{\beta}{2}} \left| \left( \Delta^{k-\beta} g \ast \partial_t  \Phi\left(t,\cdot\right)  \right) (\cdot)   \right| \right)^2 \, \frac{dt}{t} \right)^{\frac12}  \right\|_{L^q(\Rd)}.
\end{align}
A similar computation for $j \leq k$ gives
\begin{equation*} 
\left| \rho^{(j)}(1) \right| \leq  C \left\| g \right\|_{W^{\alpha,q}(\R^d)} \left\| u \ast \Phi\left(\tfrac 12,\cdot\right) \right\|_{L^p(\Rd)}.
\end{equation*}
To estimate the last term, observe that we can find $\tau \in \left( \tfrac14 ,\tfrac12 \right)$ such that 
\begin{equation*} 
\left\| u \ast \Phi\left(\tau,\cdot\right)   \right\|_{L^p(\Rd)} \leq  C \left\| \left( \int_0^1 \left( t^{\frac{\alpha}{2}} \left| \left( u \ast   \Phi\left(t,\cdot\right) \right) (\cdot) \right| \right)^2 \, \frac{dt}{t} \right)^{\frac12}  \right\|_{L^p(\Rd)}.
\end{equation*}
Since the heat kernel is contracting $L^p$-norms, we hence have that 
\begin{align*} 
\left\| u \ast \Phi(1,\cdot) \right\|_{L^p(\Rd)} & = \left\| u \ast \Phi\left(\tau,\cdot\right)   \ast \Phi(1-\tau,\cdot) \right\|_{L^p(\Rd)} \\
\notag & \leq \left\| u \ast \Phi\left(\tau,\cdot\right)   \right\|_{L^p(\Rd)} \\
\notag & \leq C \left\| \left( \int_0^1 \left( t^{\frac{\alpha}{2}} \left| \left( u \ast   \Phi\left(t,\cdot\right) \right) (\cdot) \right| \right)^2 \, \frac{dt}{t} \right)^{\frac12}  \right\|_{L^p(\Rd)}.
\end{align*}
Therefore, 
\begin{equation*} 
\max_{j\in \{0,\ldots,\alpha\} }\left| \rho^{(j)}(1) \right| \leq  C \left\| g \right\|_{W^{\alpha,p}(\R^d)}  \left\| \left( \int_0^1 \left( t^{\frac{\alpha}{2}} \left| \left( u \ast   \Phi\left(t,\cdot\right) \right) (\cdot) \right| \right)^2 \, \frac{dt}{t} \right)^{\frac12}  \right\|_{L^p(\Rd)}
\end{equation*}
follows. Thus we are left to estimate the last term on the right in~\eqref{e.taylorresidualrho_int}. On the one hand, in the case $\beta = 0$, we apply Lemma~\ref{l.MSP_pre} 
for $f =  \Delta^{k} g \in L^q(\R^d)$ to get
\begin{equation*} 
\left\| \left( \int_0^1 \left( t  \left|  \Delta \left( \Phi\left(t,\cdot\right) \ast  \Delta^{k} g  \right)  \right| \right)^2 \, \frac{dt}{t} \right)^{\frac12}  \right\|_{L^q(\Rd)} \leq C \left\|\Delta^{k} g \right\|_{L^q(\R^d)} . 
\end{equation*}
On the other hand,  in the case $\beta = 1$ we apply the same lemma for $f = \partial_{j} \Delta^{k-1} g$ together with
\begin{equation*} 
\left| \left( \Delta^{k-1} g \ast \partial_t  \Phi\left(t,\cdot\right)  \right) (x)   \right|  \leq \sum_{j=1}^d \left| \left( \partial_{j} \Delta^{k-1} g \ast  \nabla \Phi\left(t,\cdot\right)  \right) (x)   \right|,
\end{equation*}
which follows directly from the triangle inequality, to obtain
\begin{equation*} 
\sum_{j=1}^d \left\| \left( \int_0^1 \left| \nabla \left(\Phi\left(t,\cdot\right) \ast  \partial_{j} \Delta^{k-1} g \right) \right|^2 \, dt \right)^{\frac12}  \right\|_{L^q(\Rd)} 
\leq C \left\|\nabla \Delta^{k-1} g \right\|_{L^q(\R^d)}.
\end{equation*}
This finishes the proof. 
\end{proof}

We next give the main result of this appendix for the case $\alpha \not\in\N$.

\begin{proposition}
\label{p.MSP.Walpha}
\index{Poincar\'e inequality!multiscale}
Fix $\alpha\in (0,\infty) \setminus \N$ and $p\in[1,\infty]$.
There exists~$C(\alpha,p,d)<\infty$
such that, for every $u\in W^{-\alpha,p}(\Rd)$,
\begin{equation} 
\label{e.MSP.Walpha}
\frac{1}{C} \left\| u \right\|_{W^{-\alpha,p}(\Rd)} 
\leq 
 \left( \int_0^1 \left( t^{\frac\alpha2} \left\| u \ast \Phi\left(t,\cdot\right)  \right\|_{L^p(\Rd)} \right)^p \,\frac{dt}t \right)^{\frac1p} \leq C \left\| u \right\|_{W^{-\alpha,p}(\Rd)} .
\end{equation}
For $p = \infty$, the inequality above is interpreted as
\begin{equation} 
\label{e.MSP.Walphainfty}
\frac{1}{C} \left\| u \right\|_{W^{-\alpha,\infty}(\Rd)} 
\leq \sup_{t \in (0,1)} t^{\frac\alpha2} \left\| u \ast \Phi\left(t,\cdot\right)  \right\|_{L^\infty(\Rd)}  \leq C \left\| u \right\|_{W^{-\alpha,\infty}(\Rd)} .
\end{equation}
\end{proposition}

We will only prove the inequality 
\begin{equation}
\label{e.onlywhatweshow}
 \left\| u \right\|_{W^{-\alpha,p}(\Rd)} 
\leq  C
 \left( \int_0^1 \left( t^{\frac\alpha2} \left\| u \ast \Phi\left(t,\cdot\right)  \right\|_{L^p(\Rd)} \right)^p \,\frac{dt}t \right)^{\frac1p}.
\end{equation}
For the converse inequality (which we only use in this book in Exercises~\ref{ex.unbounded.critical.gff} and~\ref{ex.unbounded.critical}), we refer to \cite[Section 2.6.4]{Triebel2}.

\begin{proof}[{Proof of~\eqref{e.onlywhatweshow}}]
Throughout, we fix $\alpha \in (0,\infty)$, $p\in (1,\infty]$ and let $q:= \frac{p}{p-1} \in [1,\infty)$ denote the H\"older conjugate of~$p$. The result for $p = 1$ can be obtained by uniformity of constants by sending $p \to 1$. 

\smallskip

\emph{Step 1.} We proceed as in the proof of Proposition~\ref{p.MSP_int}. Choosing this time $ k \in \N$ such that $\alpha = 2(k-1) + \beta$, $\beta  \in (0,1) \cup(1,2)$, H\"older's inequality  implies
\begin{align} \label{e.taylorresidualrho}
\lefteqn{ \left| \int_0^1 t^{k-1} \rho^{(k)}(t) \, dt  \right| } \qquad & \\ 
& \notag \leq \int_0^1 t^{k-1} \left| \int_{\R^d} \left( u \ast \Phi\left(\tfrac t2,\cdot\right) \right)(x) \partial_t \left( \Delta^{k-1} g \ast \Phi\left(\tfrac t2,\cdot\right) \right)(x) \,dx \right| \, dt
\\ & \notag \leq \int_0^1 t^{\frac \alpha2 - \frac{\beta}{2}} \left\| u\ast \Phi\left(\tfrac t2,\cdot\right)\right\|_{L^p(\Rd)} \left\| \Delta^{k-1} g \ast \partial_t \Phi\left(\tfrac t2,\cdot\right) \right\|_{L^q(\Rd)} \, dt
\\ \notag & \leq \left( \int_0^1 \left( t^{\frac \alpha2 } \left\| u\ast \Phi\left(\tfrac t2,\cdot\right)\right\|_{L^p(\Rd)}\right)^{p} \, \frac{dt}{t} \right)^{\frac 1p} 
\\ \notag & \qquad \times \left( \int_{0}^1 \left( t^{1-\frac{\beta}{2} } \left\| \Delta^{k-1} g\ast \partial_t \Phi\left(\tfrac t2,\cdot\right) \right\|_{L^q(\Rd)} \right)^{q} \, \frac{dt}{t} \right)^{\frac 1q}. 
\end{align}
Notice that for $p = \infty$ we have that 
\begin{align}  \label{e.taylorresidualrhoinfty}
\left| \int_0^1 t^{k-1} \rho^{(k)}(t) \, dt  \right| & \leq 
\sup_{t \in (0,1)} \left( t^{\frac \alpha2 } \left\| u\ast \Phi\left(\tfrac t2,\cdot\right)\right\|_{L^\infty(\Rd)} \right)  
\\ \notag  & \qquad \times \int_{0}^1 \left( t^{1-\frac{\beta}{2} } \left\| \Delta^{k-1} g\ast \partial_t \Phi\left(\tfrac t2,\cdot\right) \right\|_{L^1(\Rd)} \right) \, \frac{dt}{t} .
\end{align}
As in the proof of Proposition~\ref{p.MSP_int}, we also have that 
\begin{equation*} 
\max_{j\in{0,\ldots,k-1} }\left| \rho^{(j)}(1) \right| \leq  C \left\| g \right\|_{W^{2k-1,q}(\R^d)}  \left( \int_0^1 \left( t^{\frac\alpha2} \left\| u \ast \Phi\left(t,\cdot\right)  \right\|_{L^p(\Rd)} \right)^p \,\frac{dt}t \right)^{\frac1p},
\end{equation*}
or
\begin{equation*} 
\max_{j\in{0,\ldots,k-1} }\left| \rho^{(j)}(1) \right| \leq  C \left\| g \right\|_{W^{2k-1,1}(\R^d)} \sup_{t \in (0,1)} \left( t^{\frac \alpha2 } \left\| u\ast \Phi\left(\tfrac t2,\cdot\right)\right\|_{L^\infty(\Rd)} \right) .
\end{equation*}
In view of~\eqref{e.rho(0) stupid},~\eqref{e.taylorresidualrho} and~\eqref{e.taylorresidualrhoinfty}, it thus suffices to show that
\begin{equation} 
\label{e.MSPgoal1}
\int_{0}^1 \left( t^{1-\frac{\beta}{2} }  \left\| \Delta^{k-1} g\ast \partial_t \Phi\left(\tfrac t2,\cdot\right) \right\|_{L^q(\Rd)} \right)^{q} \, \frac{dt}{t} \leq C^q \left\| g \right\|_{W^{\alpha,q}(\R^d)}^q.
\end{equation}
Indeed, then we would have by the definition of $\rho(0)$ that
\begin{equation*} 
\left| \int_{\R^d} u(x) g(x) \, dx \right|  \leq C \left\| g \right\|_{W^{\alpha,q}(\R^d)}  \left( \int_0^1 \left( t^{\frac\alpha2} \left\| u \ast \Phi\left(t,\cdot\right)  \right\|_{L^p(\Rd)} \right)^p \,\frac{dt}t \right)^{\frac1p},
\end{equation*}
which yields~\eqref{e.onlywhatweshow}. We devote the rest of the argument to the proof of~\eqref{e.MSPgoal1}. 

\smallskip

Denote, for simplicity, $f := \Delta^{k-1} g$. Since $\partial_t \Phi\left(t,\cdot\right) = \Delta \Phi\left(t,\cdot\right)$, an integration by parts gives that, for any affine function~$\ell_x(\cdot)$,
\begin{equation*} 
\partial_t \int_{\R^d} \Phi\left(t,x-y\right) f (y) \, dy = \int_{\R^d} \Delta_y \Phi\left(t,x-y\right)  \left(f(y) - \ell_x(y) \right) \, dy = 0.
\end{equation*}
The subscript $x$ emphasizes the fact that it may depend on $x$.  We will apply this with $\ell_x(y) = f(x)$ when $\beta \in (0,1)$, and with $\ell_x(y) = \nabla f(x) \cdot (y-x)$ when $\beta \in (1,2)$. Let us consider these two cases separately. 

\smallskip

\emph{Step 2.} The case $\beta \in (0,1)$. Decomposing the integrand into annular domains and then using Jensen's inequality gives
\begin{align*} 
\lefteqn{t^{q} \left| ( f - f(x))  \ast \partial_t \Phi\left(t,\cdot\right)(x) \right|^q } \qquad &
\\ & \leq  \left( C \int_{\R^d} t^{-\frac d2} \left(1+\frac{|x-y|^2}{t} \right) \exp\left(-\frac{|x-y|^2}{4t} \right)  \left| f(y) - f(x) \right| \, dy \right)^q
\\ & \leq C^q \left( \sum_{k=0}^\infty 2^{k(d+2)} \exp(-4^{k-2}) \fint_{B_{2^{k} \sqrt{t} }(x)} \left| f(y) - f(x) \right| \,dy \right)^q
\\ & \leq C^q  \sum_{k=0}^\infty 2^{k(d+2)} \exp(-4^{k-2}) \fint_{B_{2^{k}\sqrt{t}}(x)} \left| f(y) - f(x) \right|^q \,dy
\\ & \leq C^q  \int_{0}^\infty s^{d+1} \exp\left(-\tfrac{s^{2}}{16}\right) \fint_{B_{s \sqrt{t}}(x)} \left| f(y) - f(x) \right|^q \,dy \, ds,
\end{align*}
with a constant $C(d)<\infty$. 
It thus follows that
\begin{multline*} 
\int_{0}^1 \left( t^{1-\frac{\beta}{2} }  \left\| \Delta^{k-1} g\ast \partial_t \Phi\left(\tfrac t2,\cdot\right) \right\|_{L^q(\Rd)} \right)^{q} \, \frac{dt}{t}  
 \\ 
\leq C^q \int_{0}^1 \int_{0}^\infty s^{d+1} \exp\left(-\tfrac{s^{2}}{16}\right) t^{-q \frac{\beta}{2}}\int_{\R^d }\fint_{B_{s \sqrt{t}}(x)} \left| f(y) - f(x) \right|^q \,dy \,dx \, ds \, \frac{dt}{t} .
\end{multline*}
The integral on the right can be equivalently written as 
\begin{equation*} 
\int_{0}^\infty s^{d+1 + q \beta} \exp\left(-\tfrac{s^{2}}{16}\right) \int_{0}^1 \int_{\R^d }  \fint_{B_{s \sqrt{t}}(x)} \frac{\left| f(y) - f(x) \right|^q}{\left(s \sqrt t\right)^{q \beta}} \,dy \,dx  \, \frac{dt}{t} \, ds . 
\end{equation*}
Changing variable $t =  s^{-2} r^2 $ then leads to 
\begin{multline*} 
\int_{0}^1 \left( t^{1-\frac{\beta}{2} }  \left\| \Delta^{k-1} g\ast \partial_t \Phi\left(\tfrac t2,\cdot\right) \right\|_{L^q(\Rd)} \right)^{q} \, \frac{dt}{t}  
 \\ 
\leq C^q \int_{0}^\infty s^{d+1 + q \beta} \exp\left(-\tfrac{s^{2}}{16}\right) \int_{0}^{s} \int_{\R^d }  \fint_{B_{r}(x)} \frac{\left| f(y) - f(x) \right|^q}{r^{q \beta}} \,dy \,dx  \, \frac{dr}{r} \, ds.
\end{multline*}
An elementary computation shows that
\begin{equation*} 
\int_{0}^{\infty} \int_{\R^d }  \fint_{B_{r}(x)} \frac{\left| f(y) - f(x) \right|^q}{r^{q \beta}} \,dy \,dx  \, \frac{dr}{r}  \leq \frac{C^q}{| \beta - [\beta]|} \left\| f \right\|_{W^{\beta,q}(\R^d) }^q. 
\end{equation*}
This validates~\eqref{e.MSPgoal1} in the first case. Note also that here $C$ depends only on $d$.

\smallskip

\emph{Step 3.} The case $\beta \in (1,2)$. The computation is similar to the case $\beta \in (0,1)$. The difference is that integration by parts gives 
\begin{equation*} 
\left| ( f - \ell_x)  \ast \partial_t \Phi\left(t,\cdot\right)(x)\right| = \left(\left| \nabla f - \nabla f(x)\right|  \ast \left| \nabla \Phi\left(t,\cdot\right)\right| \right)(x)
\end{equation*}
and then, consequently, 
\begin{align*} 
\lefteqn{t^{\frac q2} \left| ( f - \ell_x )  \ast \partial_t \Phi\left(t,\cdot\right)(x) \right|^q } \qquad &
\\ & \leq   \left( C \int_{\R^d} t^{-\frac d2} \frac{|x-y|}{\sqrt{t}} \exp\left(-\frac{|x-y|^2}{4t} \right)  \left| \nabla f(y) - \nabla f(x) \right| \, dy \right)^q.
\end{align*}
Starting from the above estimate, the rest of the argument is completely analogous to the previous step. The proof is complete. 
\end{proof}

The following proposition is immediate from Proposition~\ref{p.MSP_int} and Fubini's theorem. However, since it has an easier proof than that of Proposition~\ref{p.MSP_int}, we include it here. 

\begin{proposition}
\label{p.MSP.Walpha2}
\index{Poincar\'e inequality!multiscale}
Fix $\alpha \in \N$. Then there exists $C(k,d)<\infty$ such that, for every $u\in H^{-\alpha}(\Rd)$, 
\begin{equation}
\label{e.MSP.Walpha_2k}
\left\| u \right\|_{H^{-\alpha}(\R^d)} \leq C  \left( \int_0^1 t^\alpha  \left\| u \ast \Phi\left(t,\cdot\right)\right\|_{L^2(\Rd)}^2  \, \frac{dt}{t} \right)^{\frac 12}. 
\end{equation}
\end{proposition}

\begin{proof}
The proof is similar to the one of Proposition~\ref{p.MSP.Walpha}. However, now we will show that 
\begin{equation} \label{e.MSPgoal2}
\int_{0}^1 t^{1- \beta}  \left\| \Delta^{k-1} g\ast \partial_t \Phi\left(\tfrac t2,\cdot\right) \right\|_{L^2(\Rd)}^{2} \, dt \leq C \left\| g \right\|_{H^{\alpha}(\R^d)}^2,
\end{equation}
with $\beta \in \{0,1\}$. Testing the equation of $w(t,x) := \left(\Delta^{k-1} g \ast \Phi\left(\tfrac t2,\cdot\right)\right)(x)$ with $\partial_t w$ we get
\begin{align*} 
\int_0^1 \int_{\R^d} \left| \partial_t w(t,x) \right|^2 \, dx \, dt & = \int_0^1 \int_{\R^d} \partial_t w(t,x) \Delta w(t,x) \, dx \, dt 
\\ \notag & = - \frac12 \int_0^1 \int_{\R^d} \partial_t \left| \nabla w(t,x) \right|^2 \, dx \, dt 
\\ \notag & = \frac12 \int_{\R^d}  \left| \nabla w(x,0) \right|^2 \, dx  - \frac12 \int_{\R^d}  \left| \nabla w(x,1) \right|^2 \, dx .
\end{align*}
But now
\begin{equation*} 
 \int_{\R^d}  \left| \nabla w(x,0) \right|^2 \, dx = \left\| \Delta^{k-1} \nabla g \right\|_{L^2(\R^d)}^2 \leq \left\| g \right\|_{H^{\alpha}(\R^d)}^2.
\end{equation*}
This finishes the proof when $\beta = 1$. If, on the other hand, $\beta=0$, testing with $t \partial_t w$ leads to 
\begin{align*} 
\int_0^1 \int_{\R^d} t \left| \partial_t w(t,x) \right|^2 \, dx \, dt & = \int_0^1 \int_{\R^d} t \partial_t w(t,x) \Delta w(t,x) \, dx \, dt 
\\ \notag & = - \frac12 \int_0^1 \int_{\R^d} t \partial_t \left| \nabla w(t,x) \right|^2 \, dx \, dt 
\\ \notag & = \frac12  \int_0^1 \int_{\R^d} \left| \nabla w(t,x) \right|^2 \, dx \, dt    - \frac12 \int_{\R^d}  \left| \nabla w(x,1) \right|^2 \, dx ,
\end{align*}
and since testing with $w$ gives 
\begin{equation*} 
\frac12 \int_{\R^d} \left| w(x,1) \right|^2 \, dx +  \int_0^1 \int_{\R^d} \left| \nabla w(t,x) \right|^2 \, dx \, dt =  \frac12 \int_{\R^d} \left| w(x,0) \right|^2 \, dx ,
\end{equation*}
we obtain 
\begin{equation*} 
\int_0^1 \int_{\R^d} t \left| \partial_t w(t,x) \right|^2 \, dx \, dt \leq \frac14 \int_{\R^d} \left| \Delta^{k-1} g(x) \right|^2 \, dx. 
\end{equation*}
The proof is complete. 
\end{proof}

Throughout the book, we often need to measure the size of stationary random fields (see Definition~\ref{def.stat.field}). Since these fields do not decay in space, it is natural to use local versions of Sobolev spaces. For every $\al \in \R$ and $p \in [1,\infty]$, the local Sobolev space $W^{\al,p}_\mathrm{loc}(\Rd)$ can be defined as the space of distributions $u$ such that for every $r > 0$, the semi-norm $\|u\|_{W^{\al,p}(B_r)}$ is finite. Alternatively, since the spaces $W^{\al,p}(\Rd)$ are stable under multiplication, one may define $W^{\al,p}_\mathrm{loc}(\Rd)$ as the space of distributions $u$ such that
\begin{equation*}  
\forall \phi \in C^\infty_c(\Rd) \ : \ \phi u \in W^{\al,p}(\Rd).
\end{equation*}
We often use this characterization to measure the size of a distribution $u \in W^{\al,p}_\mathrm{loc}(\Rd)$, through the evaluation of $\|\phi u\|_{W^{\al,p}(\Rd)}$ for $\phi \in C^\infty_c(\Rd)$ satisfying suitable constraints on its size and the size of its support; upper bounds on $\|\phi u\|_{W^{\al,p}(\Rd)}$ are typically obtained using Propositions~\ref{p.MSP_int} or~\ref{p.MSP.Walpha}. An exception is found in the proof of Theorem~\ref{t.HKtoSob}, where we use instead the following remark which has the advantage of relying only on spatial averages of $u$ against the heat kernel.

\begin{remark}[Local versions] \label{r.MSPlocal}
In Propositions~\ref{p.MSP_int},~\ref{p.MSP.Walpha} and~\ref{p.MSP.Walpha2} one can replace the norm on the left by~$\left\| u \right\|_{W^{-\alpha,p}(B_1)}$ and the $L^p(\R^d)$ norms on the right sides of the inequalities by norms $L^p(\Psi)$ with $\Psi(z) = \exp(-|z|)$. In particular, for $\alpha \in \N$ and $p \in (2,\infty)$ we have
\begin{equation}
\label{e.triebel_loc}
\left\| u \right\|_{W^{-\alpha,p}(B_1)} \leq C \left( \int_{\R^d} \exp(-|x|)  \left( \int_0^1 \left( t^{\frac \alpha2} \left| u \ast \left(\Phi\left(t,\cdot\right)   \right)(x)\right| \right)^2 \, \frac{dt}{t}  \right)^{\frac p2}  \,dx \right)^{\frac1p} .
\end{equation}
For $\alpha \in (0,\infty) \setminus \N$, on the other hand, we have
\begin{equation} 
\label{e.MSP.Walpha_loc}
\left\| u \right\|_{W^{-\alpha,p}(B_1)} 
\leq 
\frac{C}{| \alpha - [\alpha]|} \left( \int_0^1 \left( t^{\frac\alpha2} \left\| u \ast \Phi\left(t,\cdot\right)  \right\|_{L^p(\Psi)} \right)^p \,\frac{dt}t \right)^{\frac1p},
\end{equation} 
and, finally, for $p=2$ and $\alpha \in (0,\infty)$, 
\begin{equation}
\label{e.MSP.Walpha_2k_loc}
\left\| u \right\|_{H^{-\alpha}(B_1)} \leq C  \left( \int_0^1 t^\alpha  \left\| u \ast \Phi\left(t,\cdot\right)\right\|_{L^2(\Psi)}^2  \, \frac{dt}{t} \right)^{\frac 12}. 
\end{equation}

\end{remark}

\begin{proof} The main difference compared to~\eqref{p.MSP_int} is that now $g \in W_0^{\alpha,q}(B_1)$. This allows us to localize in~\eqref{e.rho(j)}. In particular, we have that
\begin{equation*} 
\sup_{m \in \{1,\ldots,k\}} \sup_{x \in (\R^{d}\setminus B_2 )} \sup_{(y,t) \in B_1\times (0,1)} \left( t^{-2\alpha} \exp(|x|) \left(\frac{|x-y|^2}{t^2}\right)^{-m} \Phi\left(t,x-y\right) \right) \leq C(\alpha,k,d),
\end{equation*}
and hence 
\begin{equation*} 
\sup_{m \in \{1,\ldots,k\}}  \sup_{(t,x) \in (\R^d \setminus B_2) \times (0,1)}  \exp(|x|) \left| t^{-2\alpha} \partial_t \left( \Delta^{j-1} g \ast \Phi\left(\tfrac t2,\cdot\right) \right)(x) \right| \leq C \left\| g \right\|_{L^p(B_1)}.
\end{equation*}
Therefore,~\eqref{e.rho(j)} yields
\begin{align*}
\rho^{(j)}(t)
&  \leq \left| \int_{B_2} \left( u \ast \Phi\left(\tfrac t2,\cdot\right) \right)(x) \partial_t \left( \Delta^{j-1} g \ast \Phi\left(\tfrac t2,\cdot\right) \right)(x) \,dx \right|
\\ & \quad + C \left\| g \right\|_{L^p(B_1)} \int_{\R^d \setminus B_2} t^{2\alpha}\exp(-|x|)\left| \left( u \ast \Phi\left(\tfrac t2,\cdot\right) \right)(x)\right|  \,dx
\end{align*}
Now the analysis of the first term is completely analogous to the proofs of global estimates, and the second term can be estimated by means of H\"older's inequality appropriately depending on the case. We leave the details to the reader. 
\end{proof}



\chapter{Parabolic Green functions}
\label{a.NA}

\index{Green function!construction of|(}

In this appendix, we give a mostly self-contained construction of the parabolic Green functions for a uniformly elliptic operators and a proof of the Nash-Aronson upper bounds. We do not make any structural assumption on the coefficient field~$\a(x)$ other than the uniform ellipticity condition~\eqref{e.ue}, and, for convenience, the symmetry condition $\a(x) = \a^t(x)$.
The parabolic Green function $P : (0,\infty) \times \R^d \times \R^d \mapsto \R$ is a nonnegative function that, for each $t>0$ and $x, y \in\R^d$, solves the equations
\begin{equation}
\label{e.Peq.NA}
\left\{
\begin{aligned}
& \left( \partial_t - \nabla \cdot \a \nabla \right)P(\cdot,\cdot,y) = 0 & \mbox{in} & \ (0,\infty) \times \R^d, \\ 
& P(0,\cdot,y) = \delta_y, & & 
\end{aligned}
\right.
\end{equation}
and
\begin{equation} 
\label{e.Peqdual.NA}
\left\{
\begin{aligned}
& \left( \partial_t - \nabla \cdot \a \nabla \right)P(\cdot,x,\cdot) = 0 & \mbox{in} & \ (0,\infty) \times \R^d, \\ 
& P(0,x,\cdot) = \delta_x. & & 
\end{aligned}
\right.
\end{equation}
Since $\a(\cdot)$ is symmetric, also $P$ is symmetric in $x,y$:
\begin{equation} \label{e.Psymm.NA}
P(t,x,y) = P(t,y,x).
\end{equation}
Moreover, since the equations preserve mass, we have
\begin{equation}  \label{e.Pmass.NA}
 \left\| P(t,\cdot,y) \right\|_{L^1(\R^d)} =  \left\| P(t,x,\cdot) \right\|_{L^1(\R^d)}= 1.
\end{equation}
General solutions of parabolic Cauchy problems can be represented using the Green function. Indeed, if $f \in L_{\textrm{loc}}^1(\R^d)$ is such that, for some $k\geq 0$,
\begin{equation}  \label{e.initintegrability2.NA}
\int_{\R^d} \left|f(x)\right| \exp\left(-  k |x|^2 \right) \, dx < \infty,
\end{equation}
then the function
\begin{equation*} 
w := \int_{\R^d} P(\cdot,\cdot,y) f(y) \, dy
\end{equation*}
solves $(\partial_{t} - \nabla \a \nabla )w = 0$ in $\left(0, \left( 4\Lambda k \right)^{-1}\right) \times \R^d$, and $w(t,\cdot) \to f$ almost everywhere as $t \to 0$.  
The fundamental solution satisfies a semigroup property, that is, for $t,s>0$ and $x,y \in \R^d$, 
\begin{equation*} 
P(s+t,x,y) = \int_{\R^d} P(t,x,z) P(s,z,y) \, dz.
\end{equation*}
Finally, for $g \in L_{\textrm{loc}}^1(\R\times \R^d )$ satisfying, for $k \geq 0$, 
\begin{equation}  \label{e.initintegrability3.NA}
\sup_{t>0}\int_{\R^d} \left|g(t,x)\right| \exp\left(-  k |x|^2 \right) \, dx < \infty,
\end{equation}
we have that, for $t\in\left(0,\left( 4\Lambda k \right)^{-1}\right)$ and $x \in \R^d$, the function
\begin{equation*} 
w(t,x) := \int_{0}^t \int_{\R^d} g(s,z) P(t-s,x,z) \, dz \, ds
\end{equation*}
solves $(\partial_{t} - \nabla \a \nabla )w = g$ in $\left(0, \left( 4\Lambda k \right)^{-1}\right) \times \R^d$ and $w(0,\cdot) = 0$.

\smallskip

The next proposition summarizes the results proved in this appendix. It states the basic properties of~$P$ described above, as well as Gaussian-type upper bounds in the Nash-Aronson estimate.
\index{Nash-Aronson estimate}
Recall that we denote, for every $r \in (0,\infty)$, $t > 0$ and $x \in \Rd$, the parabolic cylinder
\begin{equation*}  
Q_r(t,x) := \Ll( t-r^2,t \Rr) \times B_r(x).
\end{equation*}

\begin{proposition} 
\label{p.GFexistence.NA} 
There exists a continuous nonnegative function $P:(0,\infty) \times \R^d \times \R^d \to \R$ such that 
$P(\cdot,\cdot,y)$ solves~\eqref{e.Peq.NA}, $P(\cdot,x,\cdot)$ solves~\eqref{e.Peqdual.NA}, and~\eqref{e.Psymm.NA} and~\eqref{e.Pmass.NA} are valid. 
Moreover, for each $\alpha \in (0,\Lambda^{-1})$, there exists  a constant $C(\alpha,d,\Lambda)<\infty$ such that we have, for every $t >0$ and $x,y \in \Rd$, 
\begin{equation} 
\label{e.upperP.NA}
\left|P(t,x,y)\right| \leq C t^{-\frac d2} \exp\left(- \alpha \frac{|x-y|^2}{4 t} \right),
\end{equation}
\begin{multline} 
\label{e.gradientP.NA}
\sup_{r \in \left(0,\frac12 \sqrt{t}\right]} r \left( \left\| \nabla_x P(\cdot,\cdot,y)\right\|_{\underline{L}^2(Q_r(t,x))} + \left\| \nabla_y P(\cdot,x,\cdot)\right\|_{\underline{L}^2(Q_r(t,y))}\right) \\ \leq C t^{-\frac d2} \exp\left(- \alpha \frac{|x-y|^2}{4 t} \right),
\end{multline}
and
\begin{equation}  
\label{e.mixedP.NA}
\sup_{r \in \left(0,\frac12 \sqrt{t}\right]} r^2 \left\| \nabla_x \nabla_y P(\cdot,\cdot,\cdot)\right\|_{\underline{L}^2\left((t-r^2,t) \times B_r(x) \times B_r(y)\right)} \leq C t^{-\frac d2} \exp\left(- \alpha \frac{|x-y|^2}{4 t} \right).
\end{equation}
Finally, if $f \in L_{\textrm{loc}}^1(\R^d)$ is such that, for some $k\geq 0$,
\begin{equation}  \label{e.initcondP.NA}
\int_{\R^d} \left|f(x)\right| \exp\left(- k |x-y|^2 \right) \, dx < \infty,
\end{equation}
and $y \in \R^d$ is a Lebesgue point of $f$ in the sense $ \lim_{h \to 0} \left\| f(\cdot) - f(y) \right\|_{\underline{L}^1(B_{h}(y))} = 0$ and $|f(y)| < \infty$, then
\begin{equation} 
\label{e.initP.NA}
\lim_{t \to 0} \int_{\R^d} P(t,x,y) f(x) \, dx = f(y).
\end{equation}
\end{proposition}

\begin{remark} 
\label{r.slicesGF.NA}
By Lemma~\ref{l.slicesvsaverages.PAR}, one can replace $(t,t-r^2) \times B_r(x) \times B_r(y)$ with  $\{t\} \times B_r(x) \times B_r(y)$ in~\eqref{e.mixedP.NA}, changing integral averages appropriately.
\end{remark}

The first result we recall is the basic regularity result for parabolic equations: solutions are bounded and H\"older continuous. We refer to \eqref{e.def.X} for the definition of the function space $H^1_\pa(I\times U)$. For each $\beta \in (0,1]$, we define the parabolic $C_{\pa}^{0,\beta}$ H\"older seminorm by
\begin{equation*} 
\left[ u  \right]_{C_{\pa}^{0,\beta} \left(I \times {U} \right)} := \sup_{(t,x),(s,y) \in I \times {U} } \frac{|u(t,x) - u(s,y)|}{(|x-y| + |t-s|^{1/2})^\beta} ,
\end{equation*}
and set
\begin{equation*} 
C_{\pa}^{0,\beta} \left(I \times {U} \right) := \left\{u \in C\left(I \times {U} \right)  \, : \, 
\| u\|_{L^\infty\left(I \times {U} \right)}
+ 
 \left[ u \right]_{C_{\pa}^{0,\beta} \left(I \times {U} \right)} < \infty \right\}  .
\end{equation*}
The parabolic boundary $\partial_\sqcup \left((t_1,t_2) \times U \right)$ is defined as 
\begin{equation*} 
\partial_\sqcup \left((t_1,t_2) \times U \right) := \left( [t_1 , t_2]  \times \partial U\right) \cup \left( \{ t_1 \} \times U \right) . 
\end{equation*}

\begin{proposition} \label{p.par.reg.NA}
There exist $\gamma(d,\Lambda) > 0$ and, for each bounded open connected interval $I\subset \R$ and Lipschitz domain $U\subset \Rd$, a constant $C(U,d,\Lambda) < \infty$ such that the following holds. Let $\beta \in (0,1]$,  $\psi \in H_{\pa}^1(I \times U) \cap C_{\pa}^{0,\beta} \left( {I}  \times {U} \right)$ and $u \in H_{\pa}^1(I \times U)$ be such that
\begin{equation} 
\label{e.weaksol.bndr}
\left\{
\begin{aligned}
& \left( \partial_t  - \nabla \cdot \a \nabla  \right) u = 0 & \mbox{in}  & \ I \times U, \\
& u  = \psi & \mbox{on}  & \ \partial_{\sqcup} \left( I \times U \right).
\end{aligned}
\right.
\end{equation}
Then we have, for every $(t,x) \in  I \times {U}$ and $r>0$, 
\begin{multline*} 
\left\| u \right\|_{L^\infty(Q_r(t,x) \cap \left( I \times U \right))} +  r^{\gamma \wedge \beta} \left[ u \right]_{C_{\pa }^{0,\gamma \wedge \beta} \left(Q_{r}(t,x) \cap ( I \times {U}) \right) } 
  \\ \leq C \left( \left\| u  \right\|_{\underline{L}^2({Q_{2r}(t,x) \cap \left( I \times U \right)} )} + \left\| \psi \right\|_{C^{0,\beta} \left(Q_{2r}(t,x) \cap ( I \times {U} )\right) } \right).
\end{multline*}
\end{proposition}

\begin{proof}
The proof is based on the parabolic De Giorgi-Nash-Moser theory, see Nash~\cite{Nash58} and Moser~\cite{Moser64,Moser67}.
\end{proof}

The next lemma gives us the basic energy estimate. It is a variant of the Caccioppoli inequality. 
\index{Caccioppoli inequality}

\begin{lemma} \label{l.parcacc.NA}
Let $u \in H^1_{\pa}((t_1,t_2)\times B_r)$ be such that
\begin{equation*}  
(\partial_{t} -\nabla \cdot \a \nabla )u = 0 \qquad \text{in } (t_1,t_2) \times B_r,
\end{equation*}
and let $\psi \in C^\infty((t_1,t_2) \times B_{r})$ be such that 
\begin{equation*}
u\psi = 0 \qquad  \text{on } (t_1,t_2) \times \partial B_{r}. 
\end{equation*}
For every $\ep \in \left(0,\frac12 \right]$ and almost every $t \in (t_1,t_2)$,  we have	
\begin{multline}  \label{e.parcacc.NA}
\frac12 \partial_t  \int_{B_{r}} |u(t,x) \psi(t,x) |^2 \, dx +  \ep \int_{B_{r}} |\nabla u(t,x)  |^2 \psi^2(t,x)
  \, dx 
 \\ \leq   \int_{B_{r}}  |u(t,x) |^2 \left(\frac{\Lambda}{1-\ep} \left|\nabla \psi(t,x)\right|^2 + \psi(t,x) \partial_t \psi(t,x) \right) \, dx .
\end{multline}
\end{lemma}

\begin{proof}
By testing the equation for $u$ with $u\psi^2$, we will show that
\begin{equation}  \label{e.NA.u.basic1}
\frac12 \partial_t \left( u^2 \psi^2 \right)  + \ep \psi^2 \nabla u \cdot \a \nabla u
 \leq -  \nabla \cdot \left( u \psi^2 \a \nabla u \right) 
+  u^2 \left( \frac{\Lambda}{1- \ep} |\nabla \psi|^2  + \psi \partial_t \psi \right). 
\end{equation}
The desired inequality~\eqref{e.parcacc.NA} then follows by integrating over $B_r$, because we assume that $u \psi = 0 $ on $\partial B_r$. 
Compute, formally, 
\begin{align} \notag 
0 = u \psi^2 \left(\partial_t - \nabla \cdot \a \nabla\right) u   & = \frac12 \partial_t \left( u^2 \psi^2 \right) - \nabla \cdot \left(u \psi^2  \a \nabla u \right) 
\\ \notag & \quad +  \psi^2 \nabla u \cdot \a \nabla u +  2u\psi \nabla \psi \cdot \a \nabla u  - u^2 \psi \partial_t \psi.
\end{align}
Applying Young's inequality, we obtain, for any $\ep \in \left(0,\frac12 \right]$, 
\begin{equation*} 
\psi^2 \nabla u \cdot \a \nabla u +  2u \psi \nabla \psi \cdot \a \nabla u  \geq \ep \psi^2 \nabla u \cdot \a \nabla u - \frac{\Lambda}{1-\ep} u^2 |\nabla \psi|^2 ,
\end{equation*}
and, consequently, 
\begin{equation*} 
0 \geq \frac12 \partial_t \left( u^2 \psi^2 \right) - \nabla \cdot \left( u \psi^2  \a \nabla u \right) 
+  \ep \psi^2 \nabla u \cdot \a \nabla u - u^2 \left( \frac{\Lambda}{1-\ep} |\nabla \psi|^2  + \psi \partial_t \psi \right),
\end{equation*}
which is~\eqref{e.NA.u.basic1}. 
\end{proof}

The next lemma provides a reverse H\"older inequality for solutions. 
\begin{lemma} \label{l.RevHolder1.NA}
For each $q \in (0,2]$ and $\ep>0$, there exists a constant $C(q,\ep,d,\Lambda) < \infty$ such that the following holds. Let $r > 0$, $V := (-  \ep r^2 ,0) \times B_{r}$. and let~$u \in H^1_\pa(V)$ satisfy 
\begin{equation*}
(\partial_{t} -\nabla \cdot \a \nabla )u = 0 \qquad \text{in } V.
\end{equation*} 
Then we have the estimate
\begin{equation} \label{e.RevHolder1.NA}
 \left\| u(0,\cdot) \right\|_{\underline{L}^{2}(B_{r/2})} \leq C \left\| u \right\|_{\underline{L}^{q}(V)}.
\end{equation}
\end{lemma}

\begin{proof}
Fix $q \in (0,2]$ and $\ep>0$. Let $\sigma',\sigma \in \R$ be such that $\frac12 \leq \sigma' < \sigma \leq 1$, and $\sigma V := (-  \ep \sigma r^2,0) \times B_{\sigma r}$. We will show that there exists an exponent $p(d)>2$ and a constant $C(\ep,d,\Lambda)<\infty$ such that 
\begin{equation} \label{e.RevHolder2.NA}
 \left\| u \right\|_{\underline{L}^{p}(\sigma' V)} \leq \frac{C}{\sigma-\sigma'}  \left\| u \right\|_{\underline{L}^{2}(\sigma V)} . 
\end{equation}
Having proved this, we obtain by Lemma~\ref{l.revholderimprove} below (with $d\mu := |V|^{-1}  dx\,dt$) that
\begin{equation*} 
 \left\| u \right\|_{\underline{L}^{2}(\frac23 V)} \leq C \left\| u \right\|_{\underline{L}^{q}(V)}.
\end{equation*}
This then yields~\eqref{e.RevHolder1.NA} by Lemma~\ref{l.parcacc.NA}, cf.~\eqref{e.RevHolder1.NA.temp1} below. 

\smallskip

To prove~\eqref{e.RevHolder2.NA}, fix a smooth test function $\eta$ vanishing on the parabolic boundary of $ \sigma V$ and such that $\eta = 1$ in~$\sigma' V$, $0\leq \eta \leq 1$, and $|\nabla \eta|^2 + |\partial_{t} \eta| \leq \tfrac{16}{\ep r^2 (\sigma-\sigma')^2}$. From Lemma~\ref{l.parcacc.NA}, by integrating in time and using that $u \eta = 0$ on $\partial_{\sqcup}  \left( \sigma V \right)$,  we obtain that
\begin{multline} \label{e.RevHolder1.NA.temp1}
\sup_{s \in (- \ep \sigma r^2,0)} \fint_{B_{\sigma r}} |u(s,x) \eta(s,x) |^2 \, dx + \ep r^2 \fint_{\sigma V} |\nabla (u(t,x) \eta(t,x)) |^2 \, dx \, dt 
\\ \leq \frac{C}{(\sigma-\sigma')^2}
\fint_{\sigma V}  |u(t,x) |^2 \, dx \,dt .
\end{multline}
Letting $2^* = \frac{2d}{d-2}$ when $d>2$ and $ 2^* = 4$ otherwise, we deduce by H\"older's and Sobolev's inequalities that
\begin{align} 
\notag 
\lefteqn{\fint_{\sigma' V} \left| u(t,x) \right|^{4  - \frac{4}{2^*}} \, dx \, dt} \quad & 
\\ 
\notag & \leq C \fint_{-\ep \sigma r^2}^0 \left( \fint_{B_{\sigma r}} |u(t,x) \eta(t,x)|^{2^*} \, dx \, dt  \right)^{\frac 2{2^*}} \left( \fint_{B_{\sigma r}} |u(t,x) \eta(t,x)|^{2} \, dx \, dt  \right)^{1-\frac 2{2^*}} 
\\ 
\notag & \leq C r^2 \fint_{\sigma V}  |\nabla (u(t,x) \eta(t,x))|^{2} \, dx \, dt  \left(  \sup_{s \in (-\ep \sigma r^2,0)} \fint_{B_{\sigma r}} |u(s,x) \eta(s,x)|^{2} \, dx   \right)^{1-\frac 2{2^*}} 
\\ 
\notag & \leq  \frac C\ep \left(\frac{1}{(\sigma-\sigma')^2} 
\fint_{\sigma V} |u(t,x) |^2 \, dx \,dt \right)^{2 - \frac 2{2^*}}.
\end{align}
This yields~\eqref{e.RevHolder2.NA} with $p = 2\left(1 + \tfrac{2^*-2}{2^*} \right)>2$. The proof is complete. 
\end{proof}

Above we employed a very general statement about reverse H\"older inequalities. The result shows that, under general assumptions, they improve themselves. 
\begin{lemma} \label{l.revholderimprove}
Fix $p,q,r \in \R$ such that $0<p<q<r \leq \infty$, and fix $A,a > 0$.  Let $\mu$ be a non-negative Borel measure in $\R^n$. Let $U$ be a Borel set in $\R^n$ and let $\{\sigma U\}_{\sigma \in (0,1]}$ be a collection of Borel subsets of $U$. Assume that $f \in L^p(U;\mu) \cap L^r(U;\mu)$ satisfies a reverse H\"older inequality, that is, for every $\sigma' <\sigma \in \Ll[ \frac12,1\Rr]$, we have
\begin{equation} \label{e.revholderimproveass}
\left\| f  \right\|_{L^r(\sigma' U;\mu)} \leq A (\sigma-\sigma')^{-a} \left\| f  \right\|_{L^q(\sigma U;\mu)} .
\end{equation}
Then there exists a constant $C(a,A,p,q,r)<\infty$ such that 
\begin{equation*} 
\left\| f  \right\|_{L^r\left(\frac23 U;\mu\right)} \leq C  \left\| f  \right\|_{L^p(U;\mu)} .
\end{equation*}
\end{lemma}
\begin{proof}
Fix $p \in (0,q]$. Setting $h := \tfrac{p(r-q)}{q(r-p)} \in (0,1)$, we may rewrite $q$ as
\begin{equation*} 
q =  \left( \frac{hq}{p} \right) p +  \left( \frac{(1-h)q}{r} \right) r \qquad \mbox{and} \qquad \frac{hq}{p} + \frac{(1-h)q}{r} = 1.
\end{equation*}
Thus, by H\"older's inequality we get 
\begin{equation*} 
\left\| f \right\|_{L^{q}(\sigma U;\mu)} 
\leq 
\left\| f \right\|_{L^{p}(\sigma U;\mu)}^{h }  \left\| f \right\|_{L^{r}(\sigma U;\mu)}^{1-h } .
\end{equation*}
By Young's inequality we then obtain
\begin{equation*} 
A (\sigma-\sigma')^a \left\| f  \right\|_{L^q(\sigma U;\mu)}  \leq \frac12 \left\| f \right\|_{L^{r}(\sigma U;\mu)} + 2^{\frac{h}{1-h}} \left(A (\sigma-\sigma')^a \right)^{\frac{1}{1-h} }  \left\| f \right\|_{L^{p}(U;\mu)}.
\end{equation*}
Hence~\eqref{e.revholderimproveass} implies that
\begin{equation*} 
\left\| f  \right\|_{L^r(\sigma' U;\mu)} \leq \frac12 \left\| f \right\|_{L^{r}(\sigma U;\mu)} + 2^{\frac{p(r-q)}{r(q-p)}} \left(A (\sigma-\sigma')^a \right)^{\frac {q(r-p)}{r(q-p)}  }  \left\| f \right\|_{L^{p}(U;\mu)}.
\end{equation*}
An application of Lemma~\ref{l.simpleiter} then proves the result. 
\end{proof}

%
%
%
%

The next lemma is the key to the exponential behavior of fundamental solutions.

\begin{lemma} \label{l.scaled.NA}
There exists a constant $\delta_0(d,\Lambda) > 0$, and, for each $\al \in (0,\Lambda^{-1})$, a constant $C(\alpha,d,\Lambda) < \infty$ such that the following holds for every $\delta \in (0,\delta_0]$. Let $u \in H^1_\pa((0,\infty)\times B_{1/\delta})$ be such that
\begin{equation*}  
\begin{aligned}
 (\partial_t - \nabla \cdot \a \nabla)u & = 0 &  \text{in } & (0,\infty) \times B_{1/\delta}, \\
 u & = 0 &  \text{on } & (0,\infty) \times \partial B_{1/\delta}, \\
 \supp u(0,\cdot) & \subset B_\delta,
\end{aligned}
\end{equation*}
and let $M \in [1,\infty)$ be such that
\begin{equation}  
\label{e.ass.MBde}
0 \le u(0,\cdot) \le M |B_{\delta}|^{-1}.
\end{equation}
For every $t > 0$, we have
\begin{equation}  \label{e.preNash.NA}
 \left\|   \psi_{\alpha,\delta}\nabla u \right\|_{L^2((0,t) \times B_{1/\delta})}  + \left\|  u(t,\cdot) \psi_{\alpha,\delta}(t,\cdot)  \right\|_{L^2(B_{1/\delta})} 
  \leq C M (t+ \delta^2)^{\frac d4} , 
\end{equation}
where $\psi_{\alpha,\delta}$ is defined, for $(t,x) \in [0,\infty) \times \R^d$, by
\begin{equation}  \label{e.psi.NA}
\psi_{\alpha,\delta}(t,x) : = (t+\delta^2)^{\frac d2} \exp\left(\alpha \frac{|x|^2}{4 (t+\delta^2)} \right). 
\end{equation}

\end{lemma}

\begin{proof}
By scaling, we can assume that $M=1$.  We suppress $\alpha$ and $\delta$ in the notation for $\psi$.  

\smallskip

\emph{Step 1.} We first show that there exists a constant $C(\alpha,d,\Lambda)<\infty$ such that, for every $t>0$, 
\begin{multline}  \label{e.preNash.goal1.NA}
\left\| \left(u \psi \right)(t,\cdot) \right\|_{L^2(B_{1/\delta})} + \left\|  \psi  \nabla u \right\|_{L^2((0,t) \times B_{1/\delta})}     
\\ \leq C (t + \delta^2)^{\frac d4} \left(1+  \sup_{s \in (0,t)}  \left\| u(s,\cdot) \right\|_{L^1\left(B_{Cs^{1/2} \wedge \delta^{-1}} \right)}  \right)  .
\end{multline}
To show~\eqref{e.preNash.goal1.NA}, we employ Lemma~\ref{l.parcacc.NA} with $\psi$ defined in~\eqref{e.psi.NA}.  Notice that $u \psi = 0$ on $(0,\infty)\times \partial B_{1/\delta}$. By a direct computation we get
\begin{equation*} 
\frac{\Lambda}{1- \ep} |\nabla \psi|^2  + \psi \partial_t \psi 
=
\frac1{t+\delta^{2}} \psi^2 \left( \frac{\Lambda}{1- \ep} \alpha^2 \frac{|x|^2}{4(t+\delta^{2})}  - \alpha \frac{|x|^2}{4(t+\delta^{2})} + \frac{d}{2} \right) 
.
\end{equation*}
In particular, since $\alpha < \Lambda^{-1}$, there exist constants $\ep(\alpha,\Lambda)>0$ and $C(\alpha,d,\Lambda)<\infty$  such that 
\begin{equation*} 
\frac{\Lambda}{1- \ep} |\nabla \psi|^2  + \psi \partial_t \psi  \leq C (t+\delta^{2})^{d-1} \indc_{\{|x|^2 \leq C(t+\delta^{2})\}}. 
\end{equation*}
Plugging $\psi$ into~\eqref{e.parcacc.NA} and using this upper bound, we obtain that, for every $t > 0$, 
\begin{multline} \label{e.NA.u.basic111}
\partial_t \int_{B_{1/\delta}} \left( u^2 \psi^2 \right)(t,x) \, dx + 
\ep \int_{B_{1/\delta}} \left( |\nabla u|^2 \psi^2 \right)(t,x) \, dx  
\\ \leq C (t+\delta^{2})^{d-1}  \int_{B_{C (t^{1/2}+\delta) \wedge \delta^{-1}} } u^2(t,x) \, dx .
\end{multline}
We next show that there exists $c(\alpha,d,\Lambda) > 0$ such that, for every $t \in (0,c\delta^{-2})$, 
\begin{equation}  \label{e.NA.u.basic112}
\int_{B_{C (t^{1/2}+\delta) \wedge \delta^{-1}} } u^2(t,x) \, dx  \leq C (t + \delta^{2})^{-\frac d2}   (m_t+1)^2 ,
\end{equation}
where we denote
\begin{equation*} 
m_t := \sup_{s \in (0,t)}  \left\| u(s,\cdot) \right\|_{L^1\left(B_{8Cs^{1/2} } \right)}.
\end{equation*}
To see this, we first use Lemma~\ref{l.RevHolder1.NA} to obtain, for every $t \in  (\delta^{2}, (4C)^{-2} \delta^{-2})$, 
\begin{equation*} 
\int_{B_{2 C t^{1/2} } } u^2(t,x) \, dx \leq Ct^{\frac d2} \left\| u \right\|_{\underline{L}^1\left((t/2,t) \times B_{4Ct^{1/2} } \right)}^2 \leq C t^{-\frac d2} m_t^2 . 
\end{equation*}
On the other hand, by Lemma~\ref{l.parcacc.NA} (with $\psi \equiv 1$ there) the fact that $u$ vanishes on $(0,\infty) \times \partial B_{1/\delta}$, that $u(0,\cdot)$ is supported in $B_\delta$ and the assumption of \eqref{e.ass.MBde} (with $M = 1$), we have 
\begin{equation}  
\label{e.NA.u.basic113}
\sup_{t>0} \int_{B_{1/\delta} } u^2(t,x) \, dx \leq \int_{B_{1/\delta} } u^2(0,x) \, dx  \leq |B_\delta|^{-1} \leq 
C \delta^{-d} .
\end{equation}
Combining the last two displays, we obtain~\eqref{e.NA.u.basic112}. Observe that, using again that $u(0,\cdot)$ is supported in $B_\delta$ and \eqref{e.ass.MBde}, we also get
\begin{equation*} 
\int_{B_{1/\delta}} \left( u^2 \psi^2 \right) (0,x) \, dx \leq  \delta^{2d} \exp(\alpha)\int_{B_{1/\delta} } u^2(0,x) \, dx \leq C \delta^{d} .    
\end{equation*}
We then integrate~\eqref{e.NA.u.basic111} in time, using the previous display together with~\eqref{e.NA.u.basic112}, and conclude that up to a redefinition of $C(\alpha,d,\Lambda) < \infty$, we have for every $t \in (0,c\delta^{-2})$ that
\begin{equation*} 
\left\| \left(u \psi \right)(t,\cdot) \right\|_{L^2(B_{1/\delta})}^2 + \left\|  \psi  \nabla u \right\|_{L^2((0,t) \times B_{1/\delta})}^2     \leq C (t + \delta^{2})^{\frac d2} (m_t+1)^2  . 
\end{equation*}
This yields~\eqref{e.preNash.goal1.NA} for every $t \in (0,c\delta^{-2})$. To prove the estimate for $t \geq c\delta^{-2}$, we use Poincar\'e's inequality and Lemma~\ref{l.parcacc.NA} (with $\psi \equiv 1$ there) to get that there exists a constant $c'(d)>0$ such that
\begin{equation*} 
\partial_t  \int_{B_{1/\delta}} |u(t,x)|^2  \, dx \leq - c' \delta^{2} \int_{B_{1/\delta}} |u(t,x)  |^2 \, dx.
\end{equation*}
Using also \eqref{e.NA.u.basic113}, we thus obtain that for every $t \ge c \delta^{-2}$,
\begin{equation*} 
\int_{B_{1/\delta}} |u(t,x)|^2  \, dx \leq \exp\left( - c'(\delta^2 t - c) \right) \int_{B_{1/\delta}} |u(c\delta^{-2},x)|^2  \, dx 
 \leq 
C \exp\left( - c'\delta^{2} t \right) \delta^{-d}.
\end{equation*}
From this, \eqref{e.preNash.goal1.NA} follows easily also for $t \geq c\delta^{-2}$. 

\smallskip

\emph{Step 2.} 
We show that 
\begin{equation}  \label{e.preNash.goal2.NA}
\sup_{s \in (0,\infty)}  \left\| u(s,\cdot) \right\|_{L^1\left(B_{1/\delta} \right)} \leq 1,
\end{equation}
which, together with~\eqref{e.preNash.goal1.NA}, yields~\eqref{e.preNash.NA}. 
By the minimum principle, $u \geq 0$. Testing thus the equation of $u$ with $u/(\ep + u)$, $\ep>0$, and integrating over $(0,t) \times \R^d$ yields
\begin{multline*} 
\int_{B_{1/\delta}} \int_{0}^{u(t,z)} \frac{\theta}{\ep+\theta} \, d\theta \,dz + \ep \int_0^t \int_{B_{1/\delta}} \frac{\nabla u(t',z) \cdot \a(z) \nabla u(t',z)}{(\ep + |u(t',z)|)^2} \, dz \, dt'
\\ = \int_{B_{1/\delta}} \int_{0}^{u(0,z)} \frac{\theta}{\ep+\theta} \, d\theta \,dz \leq \int_{B_{1/\delta}} u(0,z) \,dz \leq 1.
\end{multline*}
Hence~\eqref{e.preNash.goal2.NA} follows after sending $\ep$ to zero, completing the proof. 
\end{proof}

The next lemma will be used to construct the parabolic Green function.  
\begin{lemma} \label{l.parfundsol.GF}
Fix $y \in \R^d$ and $t>0$. There exist a nonnegative solution $u$ of 
\begin{equation} 
\label{e.parfundsol.GF}
\left\{
\begin{aligned}
& \left( \partial_t  - \nabla \cdot \a \nabla  \right) u = 0 & \mbox{in}  & \ (0,\infty) \times \Rd, \\
& u  = \delta_y, & & 
\end{aligned}
\right.
\end{equation}
and there exists $\beta(d,\Lambda) \in (0,1)$ such that, for any $\ep>0$, we have
\begin{multline} \label{e.uinthreespaces.NA}
u \in H^1_{\pa}\left( (\ep^2,\infty) \times \R^d   \right) \cap H^1_{\pa}\left( (0,\infty) \times (\R^d \setminus B_\ep(y)   \right) 
\\\cap C_{\pa}^{0,\beta}\left( \left( (0,\infty) \times \R^d  \right) \setminus \left( (0,\ep^2) \times B_\ep(y) \right) \right) .
\end{multline}
Moreover, $u$ satisfies the following estimates. For all $t>0$, 
\begin{equation}  \label{e.masspreserved.NA}
 \left\| u(t,\cdot) \right\|_{L^1(\R^d)} = 1,
\end{equation}
and, for every $\alpha \in (0,\Lambda^{-1})$, there exists  a constant $C(\alpha,d,\Lambda)<\infty$ such that, for every $x \in \R^d$, we have
\begin{equation} \label{e.upperNA.GF}
|u(t,x)| \leq C t^{-\frac d2} \exp\left(- \alpha \frac{|x-y|^2}{4 t} \right)
\end{equation}
and
\begin{equation} \label{e.gradientNA.GF}
\sup_{r \in \left(0,\frac12 \sqrt{t}\right]} r \left\| \nabla u(t,\cdot)\right\|_{\underline{L}^2(B_r(x))} \leq C t^{-\frac d2} \exp\left(- \alpha \frac{|x-y|^2}{4 t} \right).
\end{equation}
Finally, if $f \in L^1(\R^d)$ is such that, for some $k>0$,
\begin{equation}  \label{e.initintegrability.NA}
\int_{\R^d} \left|f(x)\right| \exp\left(- k |x-y|^2 \right) \, dx < \infty
\end{equation}
and $y \in \R^d$ is a Lebesgue point of $f$ in the sense that $\lim_{h \to 0} \left\| f(\cdot) - f(y) \right\|_{\underline{L}^1(B_{h}(y))} = 0$ and $|f(y)| < \infty$, then
\begin{equation} \label{e.uinit.NA}
\lim_{t \to 0} \int_{\R^d} u(t,x) f(x) \, dx = f(y).
\end{equation}
\end{lemma}

\begin{proof}
\emph{Step 1.} We first show the existence of a nonnegative solution to \eqref{e.parfundsol.GF}. For each $\delta \in (0,\frac 1 2]$, let $u_\delta$ solve the equation 
\begin{equation} 
\label{e.parfundsoldelta.GF}
\left\{
\begin{aligned}
& \left( \partial_t - \nabla \cdot \a \nabla  \right)u_\delta = 0 & \mbox{in} & \ (0,\infty) \times B_{1/\delta}(y), \\
& u_\delta  = 0& \mbox{on} & \ (0,\infty) \times \partial B_{1/\delta}(y), \\ 
& u_\delta(0,\cdot) = \gamma_\delta(\cdot-y), & & 
\end{aligned}
\right.
\end{equation}
where $\gamma_\delta$ is a smooth function with unit mass supported in $B_{\delta}$ such that $0 \leq \gamma_\delta \leq 2 |B_{\delta}|^{-1}$. Extend $u_\delta$ to be zero on $(0,\infty) \times (\R^d \setminus B_{1/\delta}(y))$.  By the maximum principle, $0 \leq u_\delta \leq 2 |B_{\delta}|^{-1}$, and, by Proposition~\ref{p.par.reg.NA}, $u_\delta$ is continuous on $[0,\infty) \times \R^d$.
Lemma~\ref{l.scaled.NA} then yields that, for every $t > 0$, 
\begin{equation} 
\label{e.preNash.NA.applied}
 \left\| \psi_{\alpha,\delta} \nabla u_{\delta}  \right\|_{L^2((0,t) \times B_{1/\delta})}  + \left\|  u_{\delta}(t,\cdot) \psi_{\alpha,\delta}(t,\cdot)  \right\|_{L^2(B_{1/\delta})} 
  \leq C  (t+ \delta^2)^{\frac d4} . 
\end{equation}
This, in turn, gives by Lemma~\ref{l.parcacc.NA} and Proposition~\ref{p.par.reg.NA} that the norm of $u_\delta$ in the spaces indicated in~\eqref{e.uinthreespaces.NA} is uniformly bounded in $\delta$. This implies that there is $u$ satisfying~\eqref{e.uinthreespaces.NA} such that $u_\delta$ converges uniformly to $u$ away from $(0,y)$ and, for any $\ep>0$, $u_\delta \to u$ weakly in $H_{\pa}^1((\ep^2,\infty) \times \R^d)$, and hence $u$, extended as zero in $(-\infty,0)$ is a continuous weak solution outside of $(0,y)$.

\smallskip

\emph{Step 2.} We next show that $u$ satisfies~\eqref{e.upperNA.GF}. By~\eqref{e.preNash.NA.applied} and the pointwise convergence, we obtain that, for every $t>0$,  
\begin{equation} \label{e.preNash.NA.applied2}
 \left\|  u(t,\cdot) \exp \left(  \alpha \frac{|\cdot-y|^2}{4t} \right) \right\|_{L^2(\R^d)} \leq C t^{-\frac d4} .
 \end{equation}
 Replacing $\alpha$ with $\alpha' := \frac12 \left(\alpha + \Lambda^{-1}\right)$ above, we can find a small enough $\ep(\alpha,\Lambda)>0$ such that
 \begin{equation*} 
 \left\| u \right\|_{\underline{L}^2\left((t-\ep t,t) \times B_{(\ep t)^{1/2}}\right) } \leq C t^{-\frac d2} \exp \left(  - \alpha \frac{|x-y|^2}{4t} \right) .
\end{equation*}
Now, the $L^\infty$--$L^2$ estimate provided by Proposition~\ref{p.par.reg.NA} implies~\eqref{e.upperNA.GF}.

\smallskip

\emph{Step 3.} We now prove~\eqref{e.gradientNA.GF}. By the Caccioppoli estimate and Lemma~\ref{l.slicesvsaverages.PAR} we have, for every $t > 0$ and $r \in \left( 0,  \tfrac18(\Lambda^{-1} - \alpha)\sqrt{t} \right)$, that 
\begin{align*} \notag 
r^2 \left\| \nabla u(t,\cdot)\right\|_{\underline{L}^2(B_r(z))}^2  & \leq C \fint_{t- (2r)^2}^t \fint_{B_{3r}(z)} |u(s,z')|^2 \,dz' \, ds
\\ & \leq C t^{-d } \exp \left( - (\Lambda^{-1} + \alpha) \frac{|z-y|^2}{t} \right).
\end{align*}
On the other hand, for $z \in B_{\sqrt{t}}(x)$, we have by the triangle inequality and convexity that, for any $\ep \in (0,1)$, 
\begin{equation*} 
\left| z-y \right|^2 \geq (1-\ep) \left| x-y  \right|^2  - \frac{1-\ep}{\ep}  \left| x - z  \right|^2 \geq (1-\ep) \left| x-y  \right|^2 - \frac{1-\ep}{\ep} t.
\end{equation*}
Therefore, taking $\ep(\alpha,\Lambda) \in (0,1)$ small enough, we obtain, for $z \in B_{\sqrt{t}}(x)$ and $r \in \left( 0,  \tfrac14(\Lambda^{-1} - \alpha)\sqrt{t} \right)$,
\begin{equation*} 
r^2 \left\| \nabla u(t,\cdot)\right\|_{\underline{L}^2(B_r(z))}^2 \leq C  t^{-d} \exp \left( - 2\alpha \frac{|x-y|^2}{4t} \right).
\end{equation*}
A covering argument then shows the gradient estimate~\eqref{e.gradientNA.GF}. 

\smallskip

\emph{Step 4.} We next show~\eqref{e.masspreserved.NA}. Recall the definition of the standard mollifier $\zeta_R$ in \eqref{e.standardmollifer.delta}. We test the equation for $u$ with $\zeta_R$ and use~\eqref{e.uinthreespaces.NA} and \eqref{e.gradientNA.GF} to get that, for every $t>0$, 
\begin{align} \notag 
\int_{\R^d} u(t,x)\,dx - 1 & = \lim_{\delta \to 0} \int_{\R^d} (u_\delta(t,x) - u_\delta(0,x)) \zeta_{1/\delta}(x-y) \,dx   
\\ \notag & \leq C \limsup_{\delta \to 0} \left| \int_0^t \int_{B_{1/\delta}(y)} \left|\nabla u_\delta(x,z) \right| \left|\nabla \zeta_{1/\delta}(x-y)\right| \, dx \, dt \right|
\\ \notag & = 0, 
\end{align}
as desired.
\smallskip

\emph{Step 5.} We finally prove~\eqref{e.uinit.NA}. Let $f\in L_{\textrm{loc}}^1(\R^d)$ be such that  it satisfies~\eqref{e.initintegrability.NA} for some $k>0$ and that $\lim_{h \to 0} \left\| f(\cdot) - f(y) \right\|_{\underline{L}^1(B_{h}(y))} = 0$.
By~\eqref{e.masspreserved.NA}, \eqref{e.upperNA.GF} and a layer-cake formula, we obtain, for every $t \in \left(0,\frac{\alpha}{4k}\right)$, that
\begin{multline*} 
\left| \int_{\R^d} u(t,x) f(x) \, dx -  f(y)  \right| = \left| \int_{\R^d} u(t,x) (f(x)- f(y)) \, dx \right| 
\\ \leq C \int_{0}^\infty r^{d+1} \exp\left( - \frac{\alpha}{4} r^2 \right) \left\| f(\cdot) - f(y) \right\|_{\underline{L}^1(B_{r \sqrt{t}}(y))} \, dr,
\end{multline*}
and the right hand side tends to zero as $t \to 0$. The proof is complete.
\end{proof}

\begin{proof}[Proof of Proposition~\ref{p.GFexistence.NA}]
\emph{Step 1.} We first construct $P(t,x,y)$. Fix $x,y \in \R^d$.  By Lemma~\ref{l.parfundsol.GF} we find a solution $u(\cdot,\cdot,y)$ (we add $y$ in arguments to emphasize the location of the pole at $y$) of
\begin{equation} 
\label{e.parfundsol2.GF}
\left\{
\begin{aligned}
& \left( \partial_t  - \nabla \cdot \a \nabla  \right) u(\cdot,\cdot,y) = 0 & \mbox{in}  & \ (0,\infty) \times \Rd, \\
& u(0,\cdot,y) = \delta_y. & & 
\end{aligned}
\right.
\end{equation}
On the other hand, for each fixed $t > 0$, denoting $v(s,x,z) := u(t-s,z,x)$, we see that $v$ solves the dual equation 
\begin{equation} 
\label{e.parfundsol2dual.GF}
\left\{
\begin{aligned}
& \left[ \left( \partial_t  + \nabla \cdot \a \nabla  \right)  \right] v(\cdot,\cdot,x) = 0 & \mbox{in}  & \ (-\infty,t) \times \Rd, \\
& v(t,\cdot,x)  = \delta_x. & & 
\end{aligned}
\right.
\end{equation}
Fixing $t_1,t_2 \in \R$ such that $0<t_1<t_2<t$, we have that since $u,v \in H_\pa^1((t_1,t_2) \times \R^d)$,
then $u$ is an admissible test function for the equation of $v$, and vice versa,  on $(t_1,t_2 )\times \R^d$. We obtain
\begin{align} \notag \label{e.Greenformula.GF}
\lefteqn{\int_{\R^d} \left( (uv)(t_2,z) -(uv)(t_1,z) \right)\, dz } \quad &
\\ \notag & = \int_{t_1}^{t_2} \int_{\R^d} \partial_t  \left( u v \right) (s,z) \, dz \, ds 
\\ \notag& = \int_{t_1}^{t_2} \int_{\R^d} \left( \left( u \partial_t v + \a \nabla u \cdot \nabla v \right) + \left( v \partial_t u - \a \nabla v \cdot \nabla u \right) \right)(s,z) \, dz \,ds 
\\ \notag& = 0.
\end{align}
Letting $t_2 \to t$ and $t_1\to 0$ yields that 
\begin{equation*} 
u(t,x,y) = v(0,y,x) = u(t,y,x).
\end{equation*}
Now, for $u$ and $v$ defined as above, we set, for $t>0$ and $x,y \in \R^d$, 
\begin{equation*} 
P(t,x,y) = u(t,x,y),
\end{equation*}
and by the symmetry, for all $t>0$ and $x,y \in \R^d$, 
\begin{equation*} 
P(t,x,y) = P(t,y,x).
\end{equation*}
Now~\eqref{e.Peqdual.NA} is valid, and~\eqref{e.Psymm.NA} and~\eqref{e.Pmass.NA}  follow by Lemma~\ref{l.parfundsol.GF}.

\smallskip

\emph{Step 2.} 
We show~\eqref{e.mixedP.NA}. Notice that~\eqref{e.upperP.NA},~\eqref{e.gradientP.NA} and~\eqref{e.initP.NA} are all consequences of Lemma~\ref{l.parfundsol.GF}. 
For the mixed derivatives we first observe that $w_i := \partial_{y_i} P(\cdot,\cdot,y)$ still solves $(\partial_t - \nabla \cdot \a \nabla)w_i = 0$. Therefore, Lemma~\ref{l.parcacc.NA}  yields, for $r \in \left(0, \tfrac12 \sqrt{t} \right]$, 
\begin{equation*} 
\left\| \nabla_x \nabla_y P(\cdot,\cdot,y) \right\|_{\underline{L}^2(Q_r(t,x))} \leq C r^{-1} \left\|  \nabla_y P(\cdot,\cdot,y)  \right\|_{\underline{L}^2\left(Q_{\frac32 r}(t,x)\right)} 
\end{equation*}
On the other hand, by~\eqref{e.Psymm.NA}, $\nabla_y P(t,x,y) = \nabla_x P(t,y,x)$. Thus, applying~\ref{l.parcacc.NA} once more, after integrating in $y$, gives~\eqref{e.mixedP.NA}. 
\end{proof}

\index{Green function!construction of|)}

%
%
%

\backmatter


\bibliographystyle{plain}
{\small
\bibliography{lecturenotes}
}

\cleardoublepage
\phantomsection
\addcontentsline{toc}{chapter}{\indexname}
\printindex




\end{document}